\documentclass[sn-mathphys,Numbered]{sn-jnl}% Math and Physical Sciences Reference Style

\usepackage{graphicx}
\usepackage{multirow}
\usepackage{amsmath,amssymb,amsfonts}
\usepackage{amsthm}
\usepackage{mathrsfs}
\usepackage[title]{appendix}
\usepackage{xcolor}
\usepackage{textcomp}
\usepackage{manyfoot}
\usepackage{booktabs}
\usepackage{algorithm}
\usepackage{algorithmicx}
\usepackage{algpseudocode}
\usepackage{listings}

\usepackage{subcaption}
\usepackage{url}
\usepackage{adjustbox}
\usepackage{morefloats}

\graphicspath{{figs/}}

\theoremstyle{thmstyleone}

\theoremstyle{thmstyletwo}

\theoremstyle{thmstylethree}

\begin{document}

\title[ADER-DG method with local DG predictor for DAE systems]{%
	High order ADER-DG method with local DG predictor for 
	solutions of differential-algebraic systems of equations
}

\author*[1]{\fnm{Ivan S.} \sur{Popov}}\email{diphosgen@mail.ru, popovis@omsu.ru}
	\affil*[1]{%
	\orgdiv{Department of Theoretical Physics}, 
	\orgname{Dostoevsky Omsk State University}, 
	\orgaddress{\street{Mira prospekt}, \city{Omsk}, \postcode{644077}, \country{Russia}}
}

\abstract{%
A numerical method ADER-DG with a local DG predictor for solving a DAE system has been developed, which was based on the formulation of numerical methods ADER-DG using a local DG predictor for solving ODE and PDE systems. The basis functions were chosen in the form of Lagrange interpolation polynomials with nodal points at the roots of the right Radau polynomials, which differs from the classical formulations of the ADER-DG method, where it is customary to use the roots of Legendre polynomials. It was shown that the use of this basis leads to $A$-stability and $L_{1}$-stability in the case of using the DAE solver as ODE solver. The numerical method ADER-DG allows one to obtain a highly accurate numerical solution even on very coarse grids, with a step greater than the main characteristic scale of solution variation. The local discrete time solution can be used as a numerical solution of the DAE system between grid nodes, thereby providing subgrid resolution even in the case of very coarse grids. The classical test examples, including simple DAE systems of indices 1, 2, 3, were solved by developed numerical method ADER-DG. It was found that in the case of DAE systems of index 1, the empirical convergence orders $p$ approximately correspond to the expected convergence orders that occur in the case of solving ODE systems. With increasing index of the DAE system, a decrease in the empirical convergence orders $p$ is observed. Using the method of decreasing the index of the DAE system led to an increase in the empirical convergence orders --- the situation when the empirical convergence orders decrease with increasing index of the DAE system is expected. An unexpected result was obtained in the numerical solution of the stiff DAE system --- the empirical convergence orders of the numerical solution obtained using the developed method turned out to be significantly higher than the values expected for this method in the case of stiff problems. It turns out that the use of Lagrange interpolation polynomials with nodal points at the roots of the right Radau polynomials is much better suited for solving stiff problems. This result is of fundamental importance for the choice of basis functions when using the ADER-DG method with a local DG predictor to solve stiff problems. Using the method of decreasing the index of the DAE system led to an increase in the empirical convergence orders. Estimates of the computational costs showed that the computational costs of the ADER-DG numerical method with a local DG predictor are approximately comparable to the computational costs of implicit Runge-Kutta methods used to solve DAE systems. Methods were proposed to reduce the computational costs of the ADER-DG numerical method.
}



\maketitle

\section*{Introduction}
\addcontentsline{toc}{section}{Introduction}
\label{intro}

In this paper, a study of the arbitrary high order ADER discontinuous Galerkin (DG) method with local DG predictor, which is frequently used to solve systems of partial differential equations, is carried out to solve a system of differential-algebraic equations (DAE). The DAE system was chosen in the following form~\cite{Hairer_book_2}:
\begin{equation}\label{eq:dae_chosen_form}
\begin{split}
\frac{d\mathbf{u}}{dt} = \mathbf{F}(\mathbf{u}, \mathbf{v}, t),\qquad &t\in\Omega = \left\{t\, |\, t \in [t_{0},\ t_{f}]\right\},\\
\mathbf{0} = \mathbf{G}(\mathbf{u}, \mathbf{v}, t),\qquad &\mathbf{u}(t_{0}) = \mathbf{u}_{0},\quad \mathbf{v}(t_{0}) = \mathbf{v}_{0},\\
\end{split}
\end{equation}
where the functions $\mathbf{u}: \Omega \rightarrow \mathbb{R}^{D_{\rm u}}$, $\mathbf{v}: \Omega \rightarrow \mathbb{R}^{D_{\rm v}}$ are desired functions, the functions $\mathbf{F}: \mathbb{R}^{D_{\rm u}}\times\mathbb{R}^{D_{\rm v}}\times\Omega\rightarrow\mathbb{R}^{D_{\rm u}}$, $\mathbf{G}: \mathbb{R}^{D_{\rm u}}\times\mathbb{R}^{D_{\rm v}}\times\Omega\rightarrow\mathbb{R}^{D_{\rm v}}$ are right side functions, which are given. In terms of DAE systems the desired functions $\mathbf{u}$ are a set of differential variables, the desired functions $\mathbf{v}$ are a set of algebraic variables, $t\in\Omega$ is an independent variable. $D_{\rm u}$ is the number of differential variables, $D_{\rm v}$ is the number of algebraic variables, $D = D_{\rm u} + D_{\rm v}$ is the total number of variables. The initial conditions $\mathbf{u}(t_{0}) = \mathbf{u}_{0}$, $\mathbf{v}(t_{0}) = \mathbf{v}_{0}$ are given at point $t_{0}$ and must be consistent with the DAE system (\ref{eq:dae_chosen_form}), however, technically the initial value for the algebraic variables $\mathbf{v}_{0}$ can be formally discarded if its value is not required in solving the problem. An important numerical characteristic of a DAE system (\ref{eq:dae_chosen_form}) is the index, the value of which usually determines the difficulty of obtaining a numerical solution and high order convergence for a numerical method. There are several definitions of the DAE system index~\cite{Ascher_book_1998, Kunkel_book_2006, Riaza_book_2008, Gerdts_book_2012, handbook_num_anal_daes_2002, Hairer_book_2}, among which it is necessary to note Kronecker index for linear constant coefficient DAE systems, differentiation index, perturbation index, tractability index, geometric index, and strangeness index. In this paper, the index of a DAE system will be understood as the differentiation index~\cite{Hairer_book_2}, defined as the minimum number of analytical differentiations required to obtain an initial value problem for an ODE system from a DAE system. The main methods for solving DAE systems, especially those with high indices, are index decreasing methods~\cite{Kunkel_book_2006}, which allow one to obtain DAE systems with indices no higher than 3, or even index 0 (i.e. ODE systems), based on the use of analytical transformations and automatic differentiation (an example of a modeling language developed for these purposes is APMonitor~\cite{APMonitor_ref}).

The problem of numerical solution of DAE systems arises in a wide range of scientific and technical problems. Among these problems, one can note the problems of describing mechanical systems, electrical circuits, systems with chemical reactions, within the framework of the description of the behavior of which DAE systems arise naturally (see classical books~\cite{Hairer_book_2} and~\cite{Ascher_book_1998, Kunkel_book_2006}). At the moment, there are many numerical methods for solving DAE systems~\cite{Ascher_book_1998, Kunkel_book_2006, Riaza_book_2008, Gerdts_book_2012, handbook_num_anal_daes_2002}, many of which are essentially related to numerical methods for solving ODE systems~\cite{Butcher_book_2016, Hairer_book_1, Hairer_book_2}. In recent years, the study of DG methods and their superconvergence has been an active research field in numerical analysis, see books~\cite{Babuska_book_2001, Wahlbin_lectures_1995} and recent articles~\cite{dg_ivp_ode_4, dg_ivp_ode_5, dg_ivp_ode_6}. DG methods occupy a special place among the numerical methods for solving differential equations, since they allow obtaining a numerical solution of very high accuracy, while being sufficiently simple and scalable to solve complex and large problems. DG methods were created in the work~\cite{lasl_rep_dg_1973} for numerical solution of the neutron transport equations. In 1981 Delfour \textit{et al} \cite{Delfour_1981} constructed a Runge-Kutta-like DG method for solving the initial value problem for the ODE system that demonstrates superconvergence up to order $2N+2$ for polynomials with a degree $N$. Later, in 1986 Delfour \textit{et al} \cite{Delfour_1986} constructed a general theory of one-step, hybrid, and multi-step methods using discontinuous polynomial approximations for solving the initial value problem for the ODE system. Cockburn, Shu \textit{et al} in a series of works~\cite{Cockburn_base_1, Cockburn_base_2, Cockburn_base_3, Cockburn_base_4, Cockburn_base_5} created an accurate and thoroughly developed mathematical basis of DG methods, which stimulated their further development and use for solving a wide class of problems. DG methods are used to solve the initial value problem for ODE systems~\cite{ader_dg_ivp_ode, dg_ivp_ode_4, dg_ivp_ode_5, dg_ivp_ode_6, dg_ivp_ode_1, dg_ivp_ode_2, dg_ivp_ode_3}, the boundary value problem for ODE systems~\cite{dg_bvp_ode_1, dg_bvp_ode_2}, to solve stochastic differential equations~\cite{dg_stoch_ode_1, dg_stoch_ode_2}, and especially widely to solve PDE systems~\cite{ader_dg_ideal_flows, ader_dg_diss_flows, ader_dg_ale, ader_dg_grmhd, ader_dg_gr_prd, ader_dg_gr_z4_2024, ader_dg_simple_mod, ader_dg_PNPM, fron_phys, ader_dg_axioms, exahype, PNPM_DG_2009, ader_dg_semiexpl, ader_dg_hyperelastic, ader_dg_seiemic, ader_dg_seiemic_underwater, ader_eno_fv_blood_2022, ader_dg_gr_z4_2024, ader_dg_wb_shwater_2022, ader_eff, dg_entropy, dg_entropy_add, ader_weno_sph, ader_dg_mod_1, ader_dg_mod_2, ader_rev_2024}.

In this work, a numerical method ADER-DG with a local DG predictor for solving a DAE system (\ref{eq:dae_chosen_form}) has been developed, which was based on the formulation of numerical methods ADER-DG using a local DG predictor for solving ODE systems~\cite{ader_dg_ivp_ode} and PDE systems~\cite{ader_dg_ideal_flows, ader_dg_diss_flows, ader_dg_PNPM, fron_phys, ader_dg_axioms}. In 2008 Dumbser \textit{et al}~\cite{ader_stiff_1} proposed a paradigm of local solution and DG predictor and showed the possibility of using the numerical method to obtain a sufficiently accurate local solution of IVP for a linear scalar ODE in comparison with using the Cauchy-Kovalewski procedure (using Taylor series)~\cite{ader_init_1, ader_init_2}, which was further used in~\cite{ader_stiff_1, ader_stiff_2} and in subsequent works related to the use of a local DG predictor. In 2010 Dumbser~\cite{PNPM_DG_2010} proposed an effective way to use the ADER-DG method to solve initial value and boundary value problems for ODE systems, which was demonstrated by solving the boundary layer equation. In the work~\cite{ader_dg_ivp_ode}, the numerical method ADER-DG with a local DG predictor was adapted to solve the initial value problem for an ODE system, where the superconvergence inherent to other DG methods for solving ODE systems was demonstrated. The implementation of the method~\cite{ader_dg_ivp_ode} was based on the use of the nodal basis of Lagrange interpolation polynomials with nodal points at the nodes of the Gauss-Legendre quadrature formula. An important feature of the variant of the DG method~\cite{ader_dg_ivp_ode} was the significant simplicity of its algorithmic and software implementations, which in no way reduces the accuracy and efficiency of the method; and it was shown that the proposed version of the ADER-DG numerical method is $A$-stable and $L$-stable, and implements superconvergence with convergence order $p = 2N+1$ for the solution at grid nodes, while the local solution obtained using the local DG predictor has convergence order $p = N+1$, where $N$ is the basis polynomials degree, which is usually expected from classical implementations of DG methods for solving ODE. Numerical methods similar to the ADER-DG methods for DAE systems (\ref{eq:dae_chosen_form}) have not yet been developed, so \textit{this work is the first paper} where the ADER-DG method with a local DG predictor was developed for solving the DAE system (\ref{eq:dae_chosen_form}).

The numerical method ADER-DG developed in this paper for solving the DAE system is based on and is a development of the method proposed in the work~\cite{ader_dg_ivp_ode}. The Lagrange interpolation polynomials passing through the nodes of the Gauss-Radau quadrature formula based on the right Radau polynomials were chosen as the basis functions. It was shown that the main properties of the numerical method ADER-DG presented in the work~\cite{ader_dg_ivp_ode} are also preserved for the method presented in this paper in the case of solving the initial value problem for the ODE system. Baccouch in~\cite{dg_ivp_ode_1, dg_ivp_ode_3} obtained rigorous proofs for the convergence orders $2N+1$ and $N+1$ for DG methods presented in a modal functional representation based on Legendre polynomials, which is expected to be correct for the ADER-DG numerical method studied in this work.

It should be noted that in 2021 Han Veiga \textit{et al}~\cite{dec_vs_ader_2021} showed that the numerical methods of the ADER family are significantly interconnected with numerical methods based on the deferred correction (DeC) paradigm (see also~\cite{dec_vs_ader_2023}). DeC methods has a long history and its application to initial value problems for ODE systems goes back to~\cite{dec_src_1968}, and are effectively used to solve both ODE systems~\cite{dec_dutt_2000, dec_dutt_2000, dec_minion_2003, dec_shu_2008} and partial differential equations~\cite{dec_abgrall_2017, dec_abgrall_2019}. DeC methods, like the methods of the ADER family, allow one to obtain an arbitrarily high order and are characterized by a high accuracy of the numerical solution. The current state of research on DeC methods~\cite{dec_vs_ader_2021} shows that they compete with methods of the ADER family. The work notes~\cite{dec_vs_ader_2021} that the explicit ADER method can be considered as a special interpretation of the DeC paradigm presented in the work; it is also noted that the implementation of numerical methods of DeC is easier to implement than ADER methods.

The numerical methods of ADER-DG, based on the use of the ADER paradigm, allow the creation of numerical methods of an arbitrarily high order. Among the existing numerical methods, we can distinguish high-order numerical methods based on the Taylor expansion for solving ODE systems~\cite{ivp_ode_taylor_series_2017, ivp_ode_taylor_series_soft_2005}, which allow obtaining an arbitrary high order. High-order numerical Taylor methods, in this context, in their general structure are close to the use of the Cauchy-Kovalevskaya procedure, which was used in the original versions of the development of the ADER paradigm~\cite{ader_init_1, ader_init_2}, where this procedure was used to obtain a local solution. Modern implementations of high-precision methods based on the ADER paradigm use a local solution obtained by a local DG predictor~\cite{ader_stiff_1, ader_stiff_2, ader_weno_lstdg_ideal, ader_weno_lstdg_diss, ader_dg_dev_1, ader_dg_dev_2, ader_dg_ideal_flows, ader_dg_diss_flows, ader_eff, dg_entropy, dg_entropy_add}. Approximate Taylor methods allow arbitrarily high order, and there are efficient high-performance software implementations for using them~\cite{ivp_ode_taylor_series_soft_2005}. However, in comparison with numerical methods of ADER-DG, approximate Taylor methods have their own disadvantages, which led to the transition from using the Cauchy-Kovalevskaya procedure to using a local DG predictor to obtain a local solution in numerical methods based on the use of the ADER paradigm. Also among the existing numerical methods, high-order classical extrapolation methods for solving ODE~\cite{Hairer_book_1, Hairer_book_2, comp_rk_extr_dec} can be distinguished, which also make it possible to obtain an arbitrary high order. In contrast to the numerical methods of ADER-DG and the approximate Taylor methods, the use of extrapolation methods for solving ODE usually allow obtaining solutions only at grid nodes, while the numerical methods of ADER-DG and approximate Taylor methods allow obtaining a high-order numerical solution and in domains between grid nodes.

An important feature of the DG numerical methods, in particular, the version of the ADER-DG method with a local DG predictor presented in this work, is an arbitrarily high order of convergence and superconvergence, reaching convergence of orders $2N+1$ when using polynomials of degree $N$ in the representation of the numerical solution. Guaranteed achievement of such a high superconvergence is a rather difficult problem, especially for stiff problems~\cite{ader_dg_ivp_ode} and significantly non-uniform grids~\cite{siac_rev, siac_ref_1, siac_ref_3}. To solve this problem and related problems, there are currently post-processing methods for the solution that allow achieving high superconvergence~\cite{siac_rev, siac_ref_1}. In particular, a relevant approach to increasing the convergence of the order are the Smoothness-Increasing Accuracy-Conserving (SIAC) filtering methods~\cite{siac_ref_5}. The work~\cite{siac_rev} presents a detailed review of the development of this approach and describes the main emerging problems, the details of the description and solution of which are presented in the works~\cite{siac_rev, siac_ref_1, siac_ref_3, siac_ref_5, siac_ref_11, siac_ref_20, siac_ref_12, siac_ref_15, siac_ref_16}. The SIAC filter proposed in the work~\cite{siac_rev} allows obtaining superconvergence of the numerical solution on non-uniform grids. In this work, the methods of post-processing of the numerical solution were not used explicitly --- a ``pure'' ADER-DG method with a local DG predictor was developed for solving DAE systems. However, it is clear that the use of post-processing to obtain guaranteed high orders of convergence and superconvergence is an interesting and relevant problem for ADER-DG numerical methods. This will become one of the further directions of development of this numerical method.

An important applied feature of this implementation of the ADER-DG numerical method with a local DG predictor for solving a DAE system is the possibility of using the local solution as a solution with a subgrid resolution, which makes it possible to obtain a detailed solution even on very coarse coordinate grids. The scale of the error of the local solution, when calculating using standard representations of single or double precision floating point numbers, using large values of the degree $N$, practically does not differ from the error of the solution at the grid nodes. It should be noted that the local solution is not some kind of reconstruction of the solution from the set of values at the grid nodes, it is precisely the representation of the solution in the form of a finite element. Further in the text of this paper, a detailed description of the ADER-DG numerical method, analysis of its stability, a demonstration example of the possibility of using a local solution as output numerical solution, calculation of the convergence orders $p$ and estimations of the computational costs with comparison to the computational costs of implicit Runge-Kutta methods used to solve DAE systems are presented. It should be noted that this numerical method is fully one-step ODE solver, with uniform steps in terms of formula implementation.

\section{General description of the numerical method}
\label{sec:1}

\subsection{Formulation of the numerical method}

This Subsection presents the mathematical apparatus of the numerical method ADER-DG for solving the DAE system (\ref{eq:dae_chosen_form}). The developed numerical method is a development of the numerical method for solving the ODE system presented in the work, which was an adaptation of the classical form of the ADER-DG method for solving systems of partial differential equations presented in the works~\cite{ader_dg_ideal_flows, ader_dg_diss_flows, ader_dg_ale, ader_dg_grmhd, ader_dg_gr_prd, ader_dg_gr_z4_2024, ader_dg_simple_mod, ader_dg_PNPM, fron_phys, ader_dg_axioms, exahype, PNPM_DG_2009}.

The numerical solution of the DAE system (\ref{eq:dae_chosen_form}) is presented on the discretization of the domain of definition $\Omega$ of the desired functions by a finite number of non-overlapping discretization domains $\Omega_{n} = \{t\, |\, t\in[t_{n},\ t_{n+1}]\}$ that completely cover the domain of definition $\Omega = \cup_{n} \Omega_{n}$, where $t_{n}$ and $t_{n+1} = t_{n} + {\Delta t_{n}}$ are discretization nodes, and $\Delta t_{n}$ is the discretization step, which is not assumed to be constant for different discretization domains $\Omega_{n}$. Therefore, the proposed ADER-DG numerical method is a method that allows implementation for a variable step $\Delta t_{n}$. The set of discretization domains $\Omega_{n}$ covering the domain of definition $\Omega$ of the desired functions is a one-dimensional grid. Further points $t_{n}$ will be denoted by the grid nodes, and the space between nodes $t\in[t_{n},\ t_{n+1}]$ will be denoted by the domain between nodes.

The numerical method ADER-DG proposed in this paper for solving the DAE system is based on the use of an integral expression following from the ODE subsystem of the original DAE system (\ref{eq:dae_chosen_form}) on the discretization domains $\Omega_{n}$:
\begin{equation}\label{eq:basic_int_rel}
\mathbf{u}_{n+1} = \mathbf{u}_{n} + \int\limits_{t_{n}}^{t_{n+1}} \mathbf{F}(\mathbf{u}(t), \mathbf{v}(t), t) dt,
\end{equation}
where $\mathbf{u}_{n}$ and $\mathbf{u}_{n+1}$ are the solution to the problem at the grid nodes $t_{n}$ and $t_{n+1}$, which will hereinafter be called simply the ``solution at the nodes''; the functions $\mathbf{u}(t)$ and $\mathbf{v}(t)$ included in the integrand represent the solution to the problem in the space between the nodes $t\in[t_{n},\ t_{n+1}]$. Along with the values of the differential variables at the nodes $\mathbf{u}_{n}$, the complete numerical solution at the nodes also contains the values of the algebraic variables at the nodes $\mathbf{v}_{n}$. Further, the solution at the nodes will be understood as the complete set of values of the differential variables $\mathbf{u}_{n}$ and algebraic variables $\mathbf{v}_{n}$ at the nodes --- $(\mathbf{u}_{n}, \mathbf{v}_{n})$.

The main difference between the proposed numerical method ADER-DG from the existing implementations of this numerical method for solving ODE and PDE systems is the inclusion of an additional method for calculating algebraic variables $\mathbf{v}$ that are not present in classical ODE and PDE systems problems~\cite{ader_dg_ideal_flows, ader_dg_diss_flows, ader_dg_ale, ader_dg_grmhd, ader_dg_gr_prd, ader_dg_gr_z4_2024, ader_dg_simple_mod, ader_dg_PNPM, fron_phys, ader_dg_axioms, exahype, PNPM_DG_2009, ader_dg_ivp_ode}.

The formula apparatus of the numerical method ADER-DG was formulated in terms of the local variable $\tau$:
\begin{equation}\label{eq:tau_mapping}
t(\tau) = t_{n} + \tau\cdot{\Delta t_{n}},\qquad t\in\Omega_{n},\quad \tau\in\omega,
\end{equation}
using which the discretization domain $\Omega_{n}$ was mapped onto the reference domain $\omega = \{\tau\, |\, \tau\in[0,\ 1]\}$, which was associated with the convenience of determining the basis functions. Using the local variable $\tau$, the integral expression (\ref{eq:basic_int_rel}) was written in the following form:
\begin{equation}\label{eq:basic_int_rel_by_tau}
\mathbf{u}_{n+1} = \mathbf{u}_{n} + \int\limits_{0}^{1} \mathbf{f}(\mathbf{u}(t(\tau)), \mathbf{v}(t(\tau)), t(\tau)) d\tau,
\end{equation}
where $\mathbf{f} = {\Delta t_{n}} \cdot \mathbf{F}$ is the rescaled to the reference domain $\omega$ function of the right side of the ODE subsystem (\ref{eq:dae_chosen_form}), expression $t = t(\tau)$ is defined by the relation (\ref{eq:tau_mapping}). The solution between the nodes $(\mathbf{u}(t(\tau)), \mathbf{v}(t(\tau)))$ included in the integrand was chosen in the form of a local discrete time solution:
\begin{equation}\label{eq:qr_def_exp}
\begin{split}
\mathbf{u}(t(\tau)) \mapsto \mathbf{q}(\tau) = \sum\limits_{p} \hat{\mathbf{q}}_{p}\cdot\varphi_{p}(\tau),\qquad
\mathbf{v}(t(\tau)) \mapsto \mathbf{r}(\tau) = \sum\limits_{p} \hat{\mathbf{r}}_{p}\cdot\varphi_{p}(\tau),
\end{split}
\end{equation}
where $\varphi_{p}$ is the basis functions, in the form of an expansion of which the local discrete time solution was presented; $\hat{\mathbf{q}}_{p}$ and $\hat{\mathbf{r}}_{p}$ are the coefficients ($D_{\rm u}$-component and $D_{\rm v}$-component sets, respectively, in the implementation represented as $D_{\rm u}$-dimensional and $D_{\rm v}$-dimensional vectors, respectively) of the expansion of the local discrete time solution in the selected basis functions $\varphi_{p}$.

The basis functions $\varphi_{p}$, based on the expansion of which the local discrete time solution (\ref{eq:qr_def_exp}) representation was constructed in this work, were chosen in polynomial form $\varphi_{p}(\tau) = \sum_{l = 0}^{N} \varphi_{p, l}\tau^{l}$, where $N$ is the degree of the polynomials. The basis functions were Lagrange interpolation polynomials with nodal points at the roots $\tau_{k}$ of the right Radau polynomials $R_{N}(\tau)$. The right Radau polynomials $R_{N}(\tau)$ were chosen in the form (see~\cite{ref_to_radau_polynomials}, Section 2.4):
\begin{equation}
R_{N}(\tau) = \frac{(-1)^{N}}{2}\Big[L_{N}(\tau) - L_{N-1}(\tau)\Big],\quad N \geqslant 1,
\end{equation}
where $L_{N}(\tau)$ are the shifted Legendre polynomials. The coefficients $\varphi_{p, l}$ of the Lagrange interpolation polynomials $\varphi_{p}$ were determined from the solution of a system of linear algebraic equations $\sum_{l = 0}^{N} \varphi_{p, l}\tau_{k}^{l} = \delta_{p, k}$, where $\delta_{p, k}$ is the Kronecker delta symbol, which was technically implemented in the form of taking the inverse matrix to the Vandermonde matrix $||\tau_{k}^{l}||$. The right Radau polynomials $R_{N}(\tau)$ do not have multiple roots $\tau_{k}$, which proves the existence and uniqueness of the set of basis polynomials $\varphi_{p}$, $0 \leqslant p \leqslant N$.

It is necessary to make an important remark about the choice of basis functions $\varphi_{p}$ for the ADER-DG method with a local DG predictor for solving DAE systems. The use of Lagrange polynomials with nodal points at the roots of Legendre polynomials $L_{N}(\tau)$ as basis functions $\varphi_{p}$ showed that such a basis is not very suitable for use in the ADER-DG method with a local DG predictor for solving the DAE system (\ref{eq:dae_chosen_form}), despite the fact that this basis was chosen in works~\cite{ader_stiff_1, PNPM_DG_2010, ader_dg_ivp_ode} for the numerical solution of the initial value problem for the ODE system and in works~\cite{ader_dg_ideal_flows, ader_dg_diss_flows, ader_dg_ale, ader_dg_grmhd, ader_dg_gr_prd, ader_dg_gr_z4_2024, ader_dg_simple_mod, ader_dg_PNPM, fron_phys, ader_dg_axioms, exahype, PNPM_DG_2009} for solving PDE problems. In this case, the local solution demonstrated very strong discontinuities at the grid nodes for DAE systems with a high index (2 and more). In the case of $N=1$, the absolute value of the discontinuities increased with decreasing discretization step ${\Delta t_{n}}$, which led to negative empirical convergence orders $p$ for local solution between grid nodes. Such behavior of the local discrete time solution led to the fact that it could not be used as a correct solution between nodes. The choice of basis functions in the form of Lagrange polynomials with nodal points at the roots of right Radau polynomials $R_{N}(\tau)$ made it possible to eliminate this drawback of the ADER-DG method with a local DG predictor for solving DAE systems (\ref{eq:dae_chosen_form}).

Substituting a local discrete time solution (\ref{eq:qr_def_exp}) into the original DAE system (\ref{eq:dae_chosen_form}) generally does not allow one to obtain an exact identity, so a local solution error $\boldsymbol{\sigma}$ occurs:
\begin{equation}\label{eq:lst_dg_error_cond}
\begin{split}
&\boldsymbol{\sigma}_{\rm u}(\tau) = \frac{d\mathbf{q}(\tau)}{d\tau} - \mathbf{f}(\mathbf{q}(\tau), \mathbf{r}(\tau), t(\tau)),\quad\quad
	\mathbf{q}(0) = \mathbf{u}_{n},\\
&\boldsymbol{\sigma}_{\rm v}(\tau) = \mathbf{g}(\mathbf{q}(\tau), \mathbf{r}(\tau), t(\tau)),\quad\quad\quad\quad\quad\,\,\,\,\, \mathbf{r}(0) = \mathbf{v}_{n},
\end{split}
\end{equation}
where $\mathbf{g}\equiv\mathbf{G}$ is introduced only to match the notations with the rescaled function $\mathbf{f}$. The condition for selecting the functions of the local discrete time solution $(\mathbf{q}, \mathbf{r})$ was based on the $L_{2}$-orthogonality of the local error $\boldsymbol{\sigma}$ to all selected basis functions $\varphi_{p}$:
\begin{equation}\label{eq:dae_weak_form}
\begin{split}
\int\limits_{0}^{1} \varphi_{p}(\tau)\cdot\boldsymbol{\sigma}_{\rm u}(\tau) d\tau = \mathbf{0};&\quad\Rightarrow\quad
	\int\limits_{0}^{1} \varphi_{p}(\tau)\left[
		\frac{d\mathbf{q}(\tau)}{d\tau} - 
		\mathbf{f}(\mathbf{q}(\tau), \mathbf{r}(\tau), t(\tau))
	\right] d\tau = \mathbf{0},\\
\int\limits_{0}^{1} \varphi_{p}(\tau)\cdot\boldsymbol{\sigma}_{\rm v}(\tau) d\tau = \mathbf{0};&\quad\Rightarrow\quad
	\int\limits_{0}^{1} \varphi_{p}(\tau)\cdot\mathbf{g}(\mathbf{q}(\tau), \mathbf{r}(\tau), t(\tau)) d\tau = \mathbf{0},
\end{split}
\end{equation}
as a result of which a DAE system in weak form was obtained. The functions $\mathbf{f}$ and $\mathbf{g}$ of the weak form of the ODE system (\ref{eq:dae_weak_form}) were also presented in the form of an expansion in terms of a set of polynomials $\varphi_{p}$:
\begin{equation}\label{eq:fg_pw_repr}
\begin{split}
&\mathbf{f}(\mathbf{q}, \mathbf{r}, t(\tau)) = \sum\limits_{p = 0}^{N} \hat{\mathbf{f}}_{p}\cdot\varphi_{p}(\tau) \,\,\,\,\mapsto\,\,
	\sum\limits_{p = 0}^{N} \mathbf{f}(\hat{\mathbf{q}}_{p}, \hat{\mathbf{r}}_{p}, t(\tau_{p}))\cdot\varphi_{p}(\tau),\\
&\mathbf{g}(\mathbf{q}, \mathbf{r}, t(\tau)) = \sum\limits_{p = 0}^{N} \hat{\mathbf{g}}_{p}\cdot\varphi_{p}(\tau) \,\,\mapsto\,\,
	\sum\limits_{p = 0}^{N} \mathbf{g}(\hat{\mathbf{q}}_{p}, \hat{\mathbf{r}}_{p}, t(\tau_{p}))\cdot\varphi_{p}(\tau),
\end{split}
\end{equation}
where the point-wise evaluation was used~\cite{Jackson_2017} to determine the $L_{2}$-projections $\hat{\mathbf{f}}_{p}$ and $\hat{\mathbf{g}}_{p}$ of the functions $\mathbf{f}$ and $\mathbf{g}$, respectively, for the subsequent substitution into the integral relations (\ref{eq:dae_weak_form}), which is an effective property of nodal bases, which is formally based on calculating the integrals in $L_{2}$-projections using the Gauss-Radau quadrature formula. The Gauss-Radau quadrature formula in this work was chosen in the following form:
\begin{equation}\label{eq:gr_rule}
\int\limits_{0}^{1}d\xi \cdot f(\xi) \approx \sum\limits_{k = 0}^{N} w_{k} f(\xi_{k}),\qquad
w_{k} = \int\limits_{0}^{1}d\xi \cdot \varphi_{k}^{2}(\xi),\qquad
\sum_{k = 0}^{N} w_{k} = 1,
\end{equation}
where $w_{k}$ is the weights of the quadrature formula, normalized to one; the weights $w_{k}$ numerically coincide with the diagonal elements of the mass matrix of the basis $\varphi_{p}$; this formula gives the exact value of the integral of a polynomial function with a degree up to $2N+1$.

The weak form (\ref{eq:dae_weak_form}) of the DAE system, after integration by parts in ODE subsystem and substitution of representations (\ref{eq:qr_def_exp}) and (\ref{eq:fg_pw_repr}), was rewritten as follows:
\begin{equation}
\sum\limits_{q = 0}^{N} \Big[
	\mathrm{K}_{pq} \hat{\mathbf{q}}_{q} - 
	\mathrm{M}_{pq} \mathbf{f}\left(\hat{\mathbf{q}}_{q},\, \hat{\mathbf{r}}_{p},\, t\left(\tau_{q}\right)\right)
\Big] = \varphi_{p}(0)\, \mathbf{u}_{n},\qquad
\mathbf{g}(\hat{\mathbf{q}}_{p}, \hat{\mathbf{r}}_{p}, t(\tau_{p})) = \mathbf{0},
\end{equation}
where the matrices $\mathrm{K} = ||\mathrm{K}_{pq}||$ and $\mathrm{M} = ||\mathrm{M}_{pq}||$ were expressed and calculated in terms of the coefficients of the basis polynomials:
\begin{equation}
\mathrm{K}_{pq} = \varphi_{p}(1)\varphi_{q}(1) - \int\limits_{0}^{1} \frac{d\varphi_{p}(\tau)}{d\tau}\, \varphi_{q}(\tau)\, d\tau,\quad
\mathrm{M}_{pq} = \int\limits_{0}^{1} \varphi_{p}(\tau)\, \varphi_{q}(\tau)\, d\tau.
\end{equation}
It is interesting to note that the fulfillment of algebraic constraints $\mathbf{g} = \mathbf{0}$ in the original DAE system (\ref{eq:dae_chosen_form}) is reduced to their fulfillment only at the nodes $\tau_{k}$ of the quadrature formula, which is characteristic of DG methods with a nodal functional basis. The matrix coefficients were analytically expressed in terms of the coefficients $\varphi_{p, k}$ of the basis polynomials $\varphi_{p}$ and were pre-computed in the code. It should be noted that the functional basis $\varphi_{p}(\tau)$ is orthogonal in $L_{2}$, so the mass matrix $\mathrm{M}$ is diagonal: $\mathrm{M}_{pq} = \mathrm{M}_{p}\delta_{pq}$. The final expression, convenient for the implementation of the local DG predictor, was obtained in the following form:
\begin{equation}
\begin{split}\label{eq:lst_dg_snae}
&\hat{\mathbf{q}}_{p} - 
\sum\limits_{q = 0}^{N} \mathrm{A}_{pq} \cdot \mathbf{f}\left(\hat{\mathbf{q}}_{q},\, \hat{\mathbf{r}}_{p},\, t\left(\tau_{q}\right)\right) = \mathbf{u}_{n},\\
&\mathbf{g}(\hat{\mathbf{q}}_{p}, \hat{\mathbf{r}}_{p}, t(\tau_{p})) = \mathbf{0},
\end{split}
\end{equation}
where the notation for the matrix $\mathrm{A} = \mathrm{K}^{-1} \cdot \mathrm{M}$ was introduced, and the property
\begin{equation}\label{eq:some_int_prop}
\sum\limits_{q = 0}^{N} \left[\mathrm{K}^{-1}\right]_{pq} \varphi_{q}(0) = 1,\quad 0 \leqslant p \leqslant N,
\end{equation}
was also used. 

The first subsystem of nonlinear algebraic equations of the system (\ref{eq:lst_dg_snae}) of the local DG predictor practically coincides with the system of nonlinear algebraic equations of the local DG predictor of the ADER-DG method for solving the problem of the initial value of the ODE system, which was presented in the work~\cite{ader_dg_ivp_ode}. The difference is the additional dependence of the function $\mathbf{f}$ on the algebraic variables $\hat{\mathbf{r}}_{p}$, which was not present in the problems of the ODE systems. The second subsystem of nonlinear algebraic equations of the system (\ref{eq:lst_dg_snae}) of the local DG predictor arises only in the case of using the ADER-DG method for DAE systems.

After obtaining the local discrete time solution $(\mathbf{q}(\tau), \mathbf{r}(\tau))$ at the local DG predictor stage by solving a system of nonlinear algebraic equations (\ref{eq:lst_dg_snae}), the value of the differential variables $\mathbf{u}_{n+1}$ at the next node $t_{n+1} = t_{n} + {\Delta t}_{n}$ was calculated using the formula (\ref{eq:basic_int_rel_by_tau}), where the integral by $\tau$ was calculated using the Gauss-Radau quadrature formula (\ref{eq:gr_rule}):
\begin{equation}\label{eq:u_sol_in_node}
\mathbf{u}_{n+1} = \mathbf{u}_{n} + \sum\limits_{p = 0}^{N} w_{p}\, \mathbf{f}\left(\hat{\mathbf{q}}_{p},\, \hat{\mathbf{r}}_{p},\, t\left(\tau_{p}\right)\right),
\end{equation}
where the property $\varphi_{p}(\tau_{k}) = \delta_{p, k}$ of the basis functions $\varphi_{p}$ was taken into account.

There is no expression for algebraic variables $\mathbf{v}$ similar to the relation (\ref{eq:basic_int_rel_by_tau}) for differential variables $\mathbf{u}$. Any attempt to differentiate the original DAE system (\ref{eq:dae_chosen_form}) will lead to the formulation of a different problem, and the algebraic equations $\mathbf{G} = \mathbf{0}$ of the original DAE may not explicitly contain algebraic variables $\mathbf{v}$. Therefore, in the numerical method ADER-DG with a local DG predictor developed in this paper, the local solution $\mathbf{r}(\tau)$ at point $\tau = 1$ (this corresponds to the point $t(\tau = 1) = t_{n+1}$ in (\ref{eq:tau_mapping})) was chosen as the solution $\mathbf{v}_{n+1}$ at the node $t_{n+1}$:
\begin{equation}\label{eq:v_sol_in_node}
\mathbf{v}_{n+1} = \sum\limits_{p = 0}^{N} \hat{\mathbf{r}}_{p}\cdot\varphi_{p}(1),
\end{equation}
where the values $\varphi_{p}(1)$ of the polynomials were pre-computed in the code. 

It would be possible to add an additional matching stage instead of using a local solution (\ref{eq:v_sol_in_node}), at which the value of $\mathbf{v}_{n+1}$ is calculated and the new value of the differential variables $\mathbf{u}_{n+1}$ is recalculated (then the value of $\mathbf{u}_{n+1}$ (\ref{eq:u_sol_in_node}) would have to be considered as candidate solution), however, for this it would be necessary to obtain with sufficient accuracy the values of not only $\mathbf{u}$, but also the derivatives of $d\mathbf{u}/dt$ (to calculate the $\mathbf{f}$), which would require additional differentiation of the ODE subsystem of the original DAE system (\ref{eq:dae_chosen_form}) and the solution of a ODE system of twice the higher order --- this would be a different problem than the original problem (\ref{eq:dae_chosen_form}).

The software implementation of the numerical method ADER-DG with a local DG predictor was developed using the \texttt{python} programming language. The error of the numerical solution obtained using this numerical method becomes very small even on fairly coarse grids. Therefore, in the software implementation, floating-point numbers of arbitrarily high precision were used within module \texttt{mpmath} (with an additional import of module \texttt{gmpy2}) of the \texttt{python} programming language, the value \texttt{mpmath.mp.dps = 500} was chosen. All calculations presented in this paper were performed using this developed software implementation.

\subsection{Solution of the local DG predictor system of equations}

The resulting system of equations (\ref{eq:lst_dg_snae}) of the local DG predictor can be solved using the classical Newton numerical method for solving a system of nonlinear algebraic equations. In this paper, a version of the classical Newton numerical method was used to perform test calculations, the detailed structure of which was as follows. The solution at each subsequent iteration $(\hat{\mathbf{q}}_{p}^{i+1}, \hat{\mathbf{r}}_{p}^{i+1})$ was calculated through the increment:
\begin{equation}
\begin{split}
\hat{\mathbf{q}}_{p}^{i+1} = \hat{\mathbf{q}}_{p}^{i} + \Delta\hat{\mathbf{q}}_{p}^{i},\qquad
\hat{\mathbf{q}}_{r}^{i+1} = \hat{\mathbf{r}}_{p}^{i} + \Delta\hat{\mathbf{r}}_{p}^{i},
\end{split}
\end{equation}
where $i$ and $i+1$ are the numbers of previous and next iterations, $\Delta\hat{\mathbf{q}}_{p}^{i}$ and $\Delta\hat{\mathbf{r}}_{p}^{i}$ are the increments of the desired quantities $\hat{\mathbf{q}}_{p}^{i+1}$ and $\hat{\mathbf{r}}_{p}^{i+1}$. The criterion for stopping the iterative process of the numerical Newton method was chosen to be the condition of smallness of the increment: $|\Delta\hat{\mathbf{q}}_{p}^{i}| < \varepsilon_{0}$ and $|\Delta\hat{\mathbf{r}}_{p}^{i}| < \varepsilon_{0}$ for all $p = 0, \ldots, N$, where the operation $|\ldots|$ is understood as taking the maximum absolute value of the vector element. The increments $(\Delta\hat{\mathbf{q}}_{p}^{i}, \Delta\hat{\mathbf{r}}_{p}^{i})$ at each iteration were calculated based on the solution of the system of linear algebraic equations obtained as a result of the linearization of the system of nonlinear algebraic equations (\ref{eq:lst_dg_snae}) of the local DG predictor by the desired quantities $(\hat{\mathbf{q}}_{p}, \hat{\mathbf{r}}_{p})$ in the neighborhood of the value of the desired quantity at the previous iteration $(\hat{\mathbf{q}}_{p}^{i}, \hat{\mathbf{r}}_{p}^{i})$. The linearization of the system of nonlinear algebraic equations (\ref{eq:lst_dg_snae}) was obtained in the following form:
\begin{equation}
\begin{split}
&\sum\limits_{q = 0}^{N}\left[
	\mathbb{I}_{D_{\rm u} \times D_{\rm u}}\cdot\delta_{pq} - 
	\mathrm{A}_{pq}\cdot\frac{\partial\mathbf{f}}{\partial\mathbf{u}}
	\left(\hat{\mathbf{q}}_{q}^{i}, \hat{\mathbf{r}}_{q}^{i}, t(\tau_{q})\right)
\right]\cdot\Delta\hat{\mathbf{q}}_{q}^{i}\\
&-\sum\limits_{q = 0}^{N}\mathrm{A}_{pq}\cdot\frac{\partial\mathbf{f}}{\partial\mathbf{v}}
	\left(\hat{\mathbf{q}}_{q}^{i}, \hat{\mathbf{r}}_{q}^{i}, t(\tau_{q})\right)\cdot\Delta\hat{\mathbf{r}}_{q}^{i} = 
	\mathbf{u}_{n} + \sum\limits_{q = 0}^{N} \mathrm{A}_{pq}\cdot\mathbf{f}\left(
		\hat{\mathbf{q}}_{q}^{i}, \hat{\mathbf{r}}_{q}^{i}, t(\tau_{q})
	\right) - \hat{\mathbf{q}}_{p}^{i},\\
&\frac{\partial\mathbf{g}}{\partial\mathbf{u}}
	\left(\hat{\mathbf{q}}_{p}^{i}, \hat{\mathbf{r}}_{p}^{i}, t(\tau_{p})\right)\cdot\Delta\hat{\mathbf{q}}_{p}^{i} +
\frac{\partial\mathbf{g}}{\partial\mathbf{v}}
	\left(\hat{\mathbf{q}}_{p}^{i}, \hat{\mathbf{r}}_{p}^{i}, t(\tau_{p})\right)\cdot\Delta\hat{\mathbf{r}}_{p}^{i}
= -\mathbf{g}\left(\hat{\mathbf{q}}_{p}^{i}, \hat{\mathbf{r}}_{p}^{i}, t(\tau_{p})\right),
\end{split}
\end{equation}
where $\partial\mathbf{f}/\partial\mathbf{u}\equiv{\Delta t_{n}}\partial\mathbf{F}/\partial\mathbf{u}$, $\partial\mathbf{f}/\partial\mathbf{v}\equiv{\Delta t_{n}}\partial\mathbf{F}/\partial\mathbf{v}$, $\partial\mathbf{g}/\partial\mathbf{u}\equiv\partial\mathbf{G}/\partial\mathbf{u}$, $\partial\mathbf{g}/\partial\mathbf{v}\equiv\partial\mathbf{G}/\partial\mathbf{v}$ are the Jacobian matrices of the vector-functions $\mathbf{F}$ and $\mathbf{G}$ included in the original DAE system (\ref{eq:dae_chosen_form}): $\partial\mathbf{f}/\partial\mathbf{u}$ is the matrix of size $D_{\rm u} \times D_{\rm u}$, $\partial\mathbf{f}/\partial\mathbf{v}$ is the matrix of size $D_{\rm u} \times D_{\rm v}$, $\partial\mathbf{g}/\partial\mathbf{u}$ is the matrix of size $D_{\rm u} \times D_{\rm v}$, $\partial\mathbf{g}/\partial\mathbf{v}$ is the matrix of size $D_{\rm v} \times D_{\rm v}$; and $\mathbb{I}_{D_{\rm u} \times D_{\rm u}}$ is the identity matrix of size $D_{\rm u} \times D_{\rm u}$. To use standard solvers, the resulting system of equations was rewritten in matrix form of the following form:
\begin{equation}\label{eq:lst_dg_newton_slae}
\left[
\begin{array}{cccc|cccc}
\mathrm{P}_{0,0}^{i} & \mathrm{P}_{0,1}^{i} & \cdots & \mathrm{P}_{0,N}^{i} & 
\mathrm{Q}_{0,0}^{i} & \mathrm{Q}_{0,1}^{i} & \cdots & \mathrm{Q}_{0,N}^{i} \\
\mathrm{P}_{1,0}^{i} & \mathrm{P}_{1,1}^{i} & \cdots & \mathrm{P}_{1,N}^{i} & 
\mathrm{Q}_{1,0}^{i} & \mathrm{Q}_{1,1}^{i} & \cdots & \mathrm{Q}_{1,N}^{i} \\
\vdots & \vdots & \ddots & \vdots & \vdots & \vdots & \ddots & \vdots\\
\mathrm{P}_{N,0}^{i} & \mathrm{P}_{N,1}^{i} & \cdots & \mathrm{P}_{N,N}^{i} & 
\mathrm{Q}_{N,0}^{i} & \mathrm{Q}_{N,1}^{i} & \cdots & \mathrm{Q}_{N,N}^{i} \\[1mm]
\hline \vphantom{1^{1^{1}}}
\mathrm{R}_{0}^{i} & {0}_{R} & \cdots & {0}_{R} & 
\mathrm{S}_{0}^{i} & {0}_{S} & \cdots & {0}_{S} \\
{0}_{R} & \mathrm{R}_{1}^{i} & \cdots & {0}_{R} & 
{0}_{S} & \mathrm{S}_{1}^{i} & \cdots & {0}_{S} \\
\vdots & \vdots & \ddots & \vdots & \vdots & \vdots & \ddots & \vdots\\
{0}_{R} & {0}_{R} & \cdots & \mathrm{R}_{N}^{i} & 
{0}_{S} & {0}_{S} & \cdots & \mathrm{S}_{N}^{i}
\end{array}
\right]\left[
\begin{array}{c}
\Delta\mathbf{q}_{0}^{i}\\
\Delta\mathbf{q}_{1}^{i}\\
\vdots\\
\Delta\mathbf{q}_{N}^{i}\\[1mm]
\hline \vphantom{1^{1^{1}}}
\Delta\mathbf{r}_{0}^{i}\\
\Delta\mathbf{r}_{1}^{i}\\
\vdots\\
\Delta\mathbf{r}_{N}^{i}
\end{array}
\right] = \left[
\begin{array}{c}
\mathbf{b}_{0}^{i}\\
\mathbf{b}_{1}^{i}\\
\vdots\\
\mathbf{b}_{N}^{i}\\[1mm]
\hline \vphantom{1^{1^{1}}}
\mathbf{c}_{0}^{i}\\
\mathbf{c}_{1}^{i}\\
\vdots\\
\mathbf{c}_{N}^{i}
\end{array}
\right],
\end{equation}
where ${0}_{R}$ is the zero matrix of size $D_{\rm v} \times D_{\rm u}$, ${0}_{S}$ is the zero matrix of size $D_{\rm v} \times D_{\rm v}$, and the following notations were introduced:
\begin{equation}
\begin{split}
&\mathrm{P}_{p, q}^{i} = \mathbb{I}_{D_{\rm u} \times D_{\rm u}}\cdot\delta_{pq} - 
	\mathrm{A}_{pq}\cdot\frac{\partial\mathbf{f}}{\partial\mathbf{u}}
	\left(\hat{\mathbf{q}}_{q}^{i}, \hat{\mathbf{r}}_{q}^{i}, t(\tau_{q})\right),\quad
\mathrm{Q}_{p, q}^{i} = \mathrm{A}_{pq}\cdot\frac{\partial\mathbf{f}}{\partial\mathbf{v}}
	\left(\hat{\mathbf{q}}_{q}^{i}, \hat{\mathbf{r}}_{q}^{i}, t(\tau_{q})\right),\\
&\mathrm{R}_{p}^{i} = \frac{\partial\mathbf{g}}{\partial\mathbf{u}}\left(\hat{\mathbf{q}}_{q}^{i}, \hat{\mathbf{r}}_{q}^{i}, t(\tau_{q})\right),\quad
\hspace{33mm}\mathrm{S}_{p}^{i} = \frac{\partial\mathbf{g}}{\partial\mathbf{v}}\left(\hat{\mathbf{q}}_{q}^{i}, \hat{\mathbf{r}}_{q}^{i}, t(\tau_{q})\right),\\
&\mathbf{b}_{p}^{i} = \mathbf{u}_{n} + \sum\limits_{q = 0}^{N} \mathrm{A}_{pq}\cdot\mathbf{f}\left(
		\hat{\mathbf{q}}_{q}^{i}, \hat{\mathbf{r}}_{q}^{i}, t(\tau_{q})
	\right) - \hat{\mathbf{q}}_{p}^{i},\quad
\hspace{7.1mm}\mathbf{c}_{p}^{i} = -\mathbf{g}\left(\hat{\mathbf{q}}_{q}^{i}, \hat{\mathbf{r}}_{q}^{i}, t(\tau_{q})\right),
\end{split}
\end{equation}
which made it possible to formulate the matrix of the system of equations (\ref{eq:lst_dg_newton_slae}) in block form. Such a structure of the matrix of the system of equations makes it convenient to use continuous slices (with a unit step along one of the directions of a two-dimensional array, which will limit the step to only one leading dimension \texttt{lda}) in programming languages such as \texttt{fortran} and \texttt{python}, as well as when using \texttt{std::valarray} in \texttt{C++}, in software implementation. 

The resulting system of linear algebraic equations (\ref{eq:lst_dg_newton_slae}) has matrix of size $[(N+1)D] \times [(N+1)D]N$. To solve the resulting system of linear algebraic equations, a standard method based on LU decomposition was used. This approach is characterized by a computational complexity of $O(D^{3}(N+1)^{3})$.

\begin{figure}[h!]
\captionsetup[subfigure]{%
	position=bottom,
	font+=smaller,
	textfont=normalfont,
	singlelinecheck=off,
	justification=raggedright
}
\centering
\begin{subfigure}{0.24\textwidth}
\includegraphics[width=\textwidth]{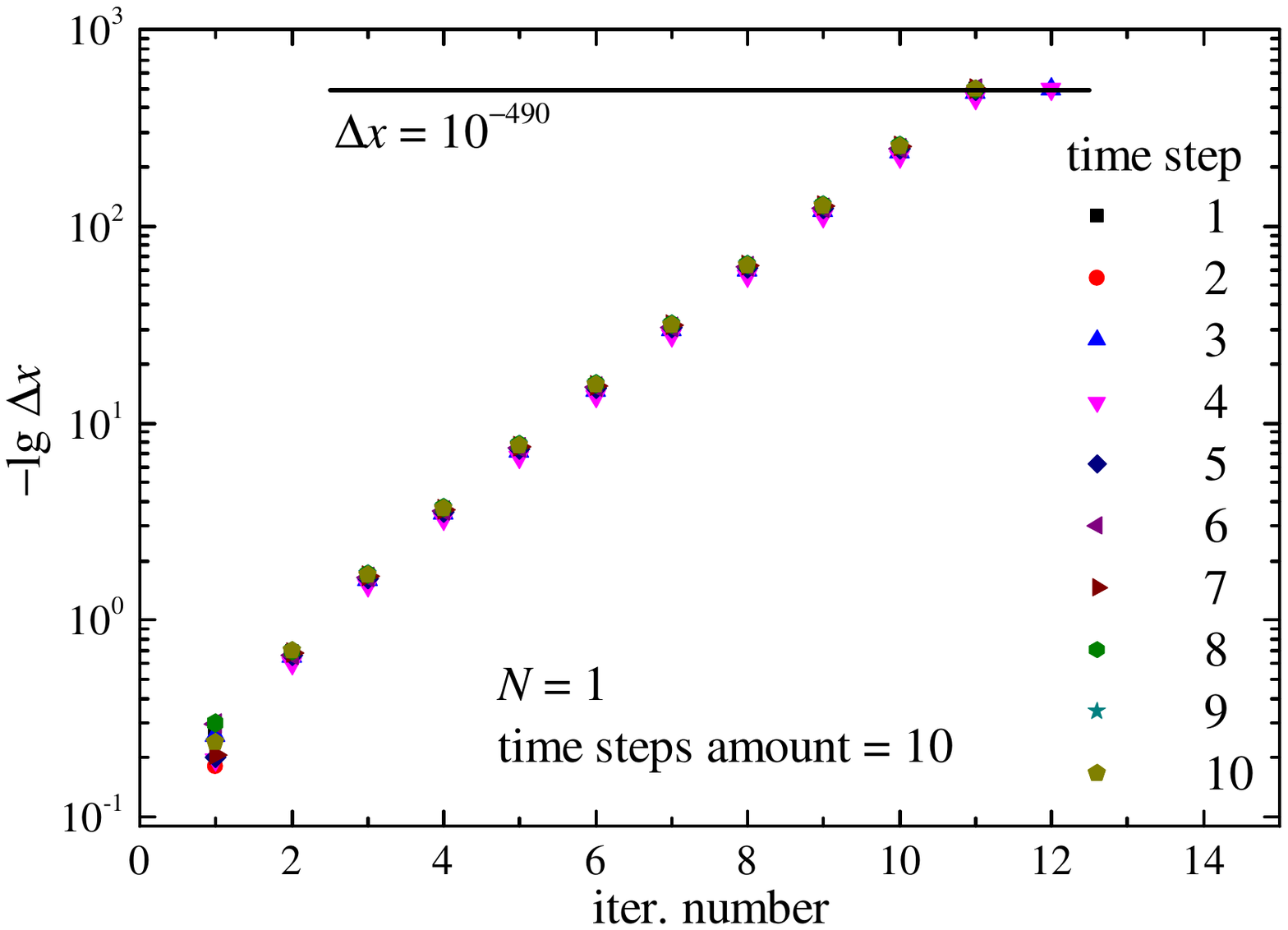}
\vspace{-8mm}\caption{\label{fig:newton_errs:a}\vspace{2mm}}
\end{subfigure}
\begin{subfigure}{0.24\textwidth}
\includegraphics[width=\textwidth]{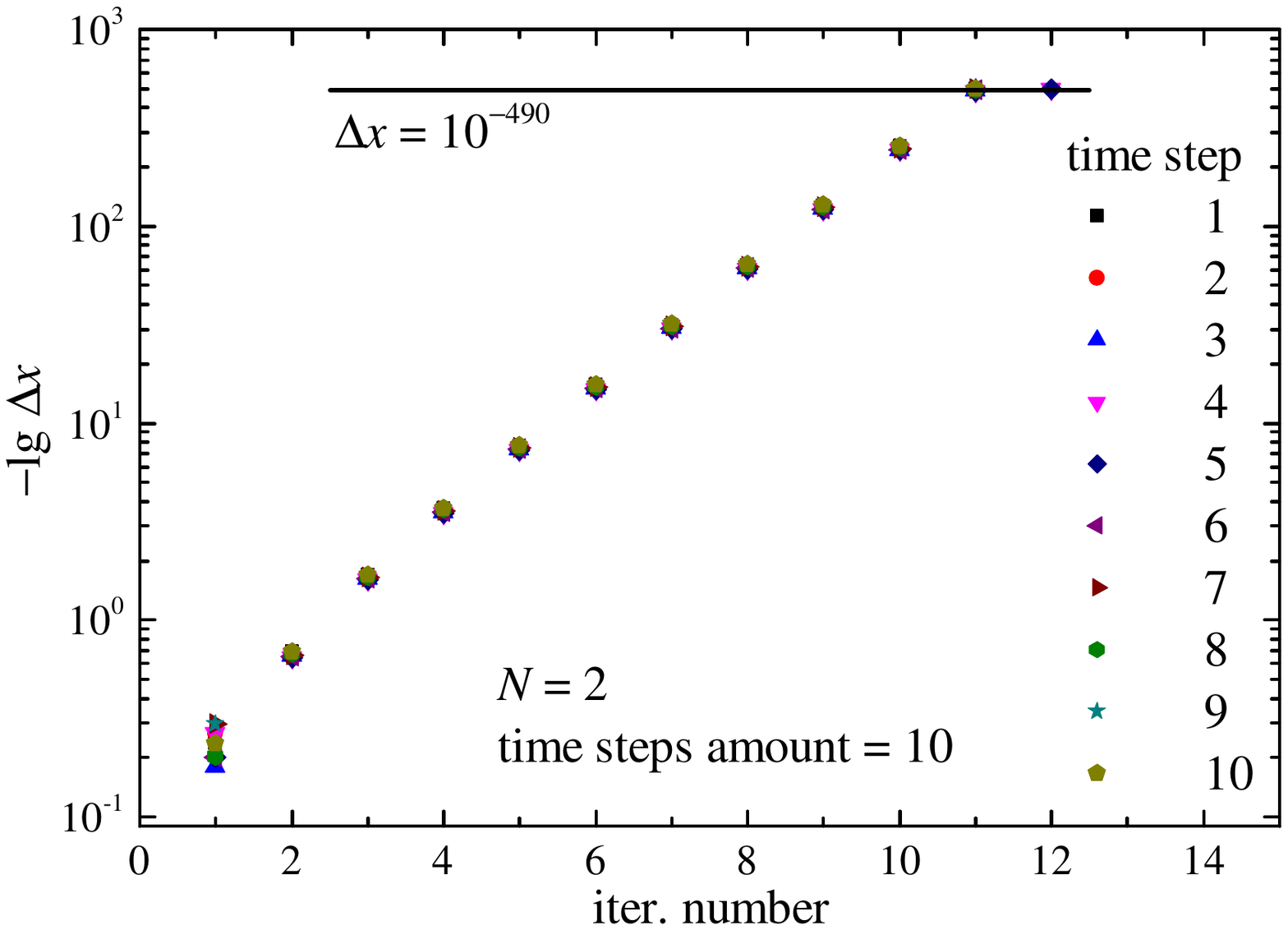}
\vspace{-8mm}\caption{\label{fig:newton_errs:b}\vspace{2mm}}
\end{subfigure}
\begin{subfigure}{0.24\textwidth}
\includegraphics[width=\textwidth]{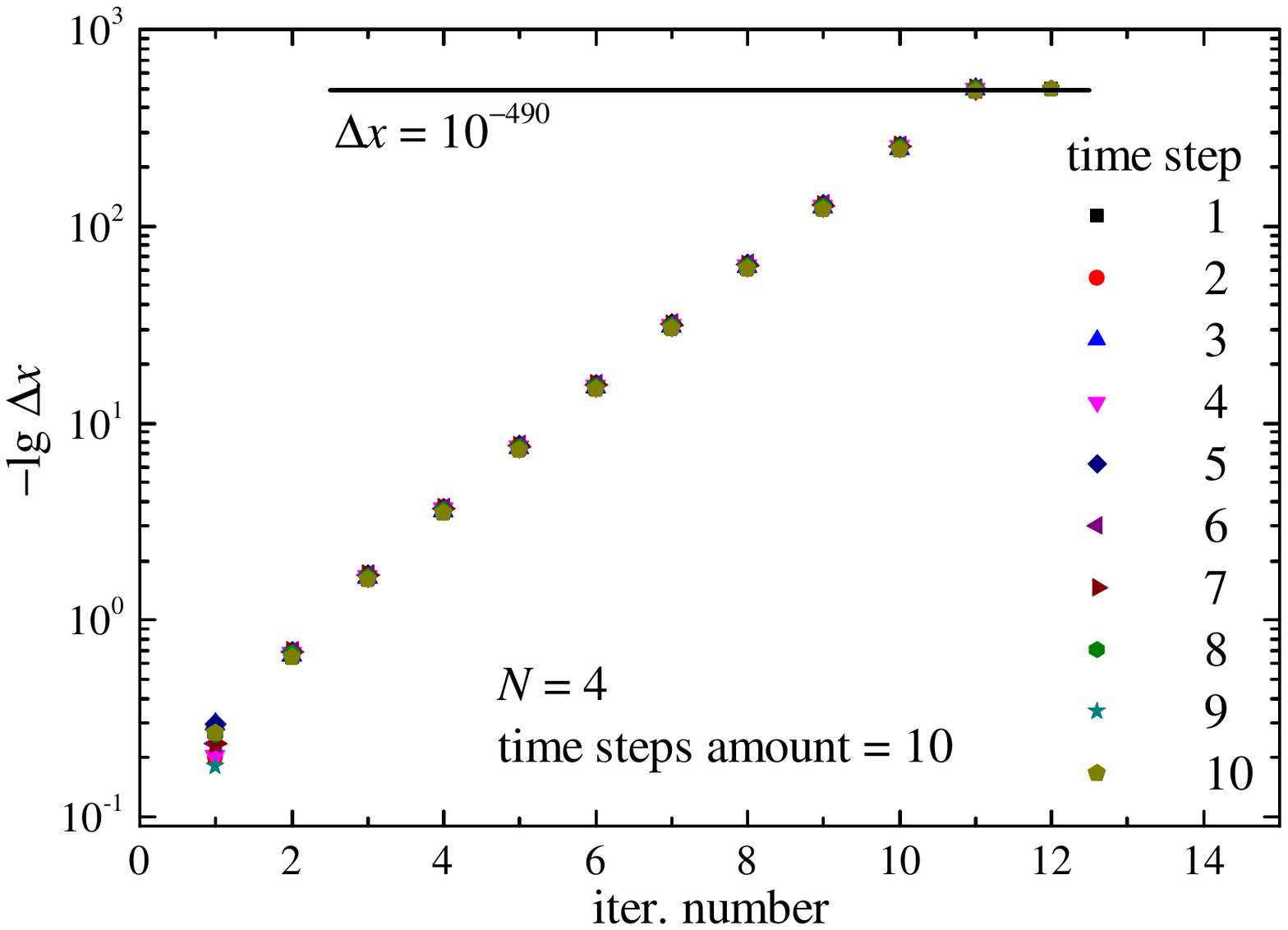}
\vspace{-8mm}\caption{\label{fig:newton_errs:c}\vspace{2mm}}
\end{subfigure}
\begin{subfigure}{0.24\textwidth}
\includegraphics[width=\textwidth]{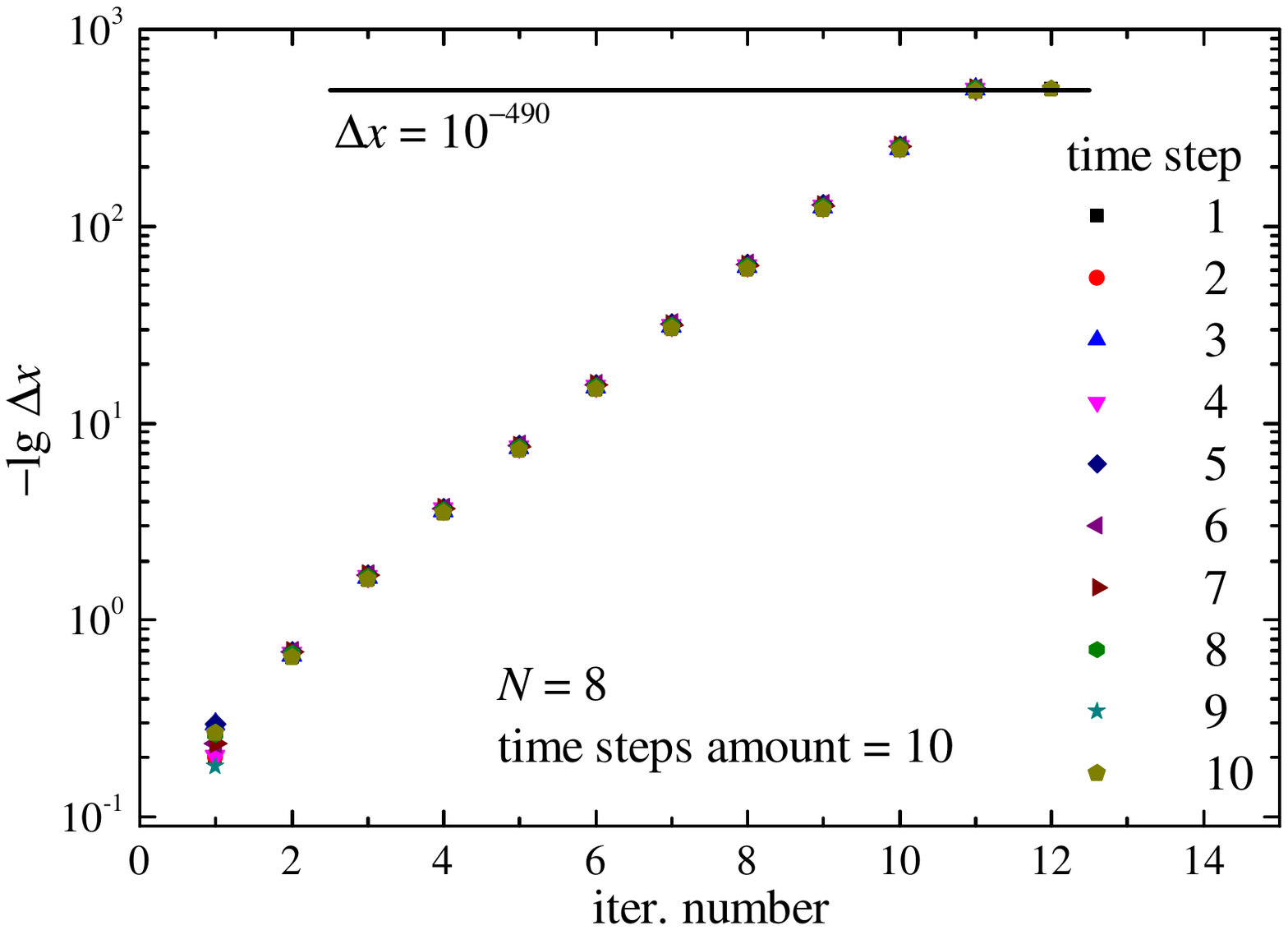}
\vspace{-8mm}\caption{\label{fig:newton_errs:d}\vspace{2mm}}
\end{subfigure}
\begin{subfigure}{0.24\textwidth}
\includegraphics[width=\textwidth]{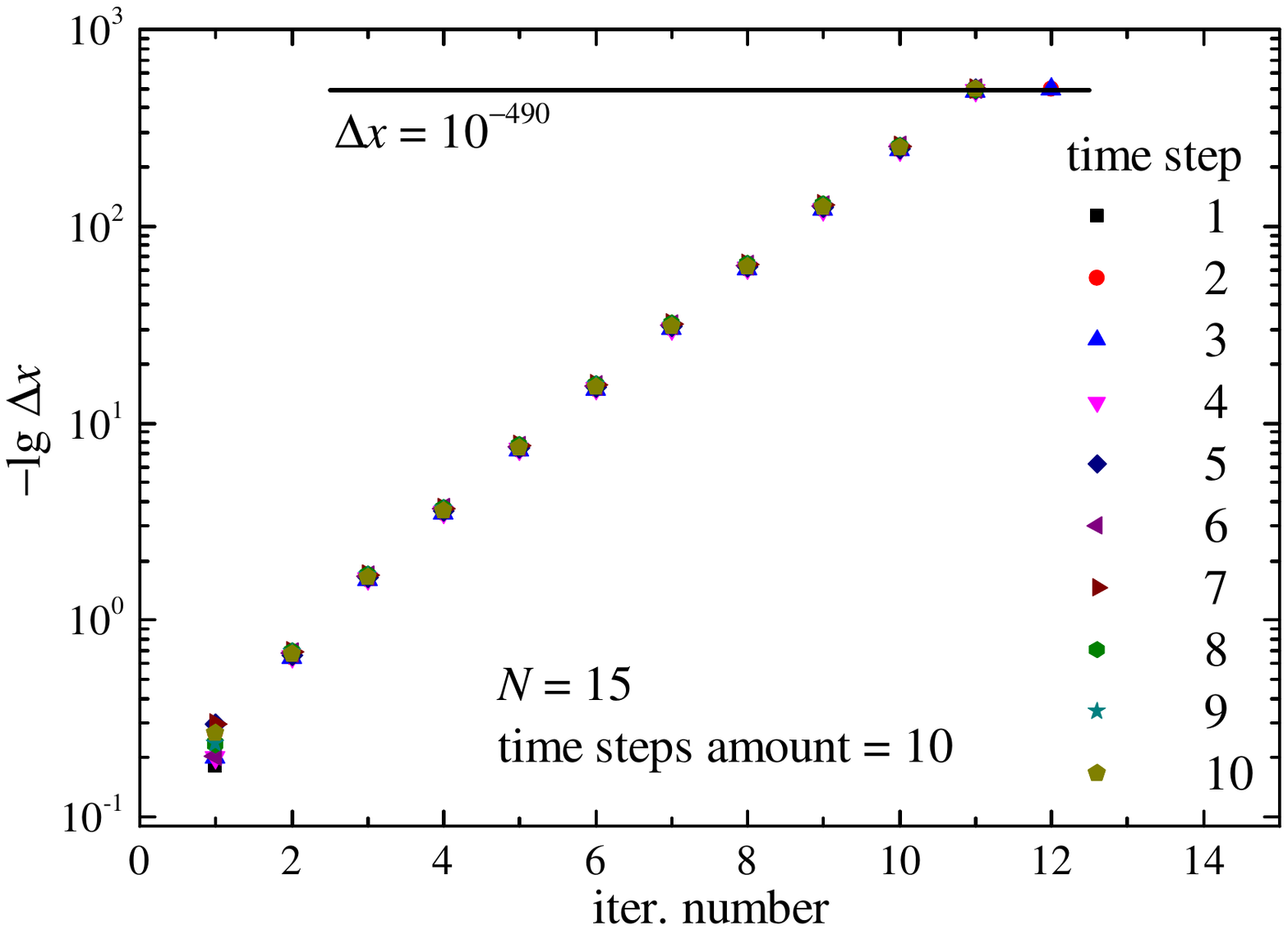}
\vspace{-8mm}\caption{\label{fig:newton_errs:e}\vspace{2mm}}
\end{subfigure}
\begin{subfigure}{0.24\textwidth}
\includegraphics[width=\textwidth]{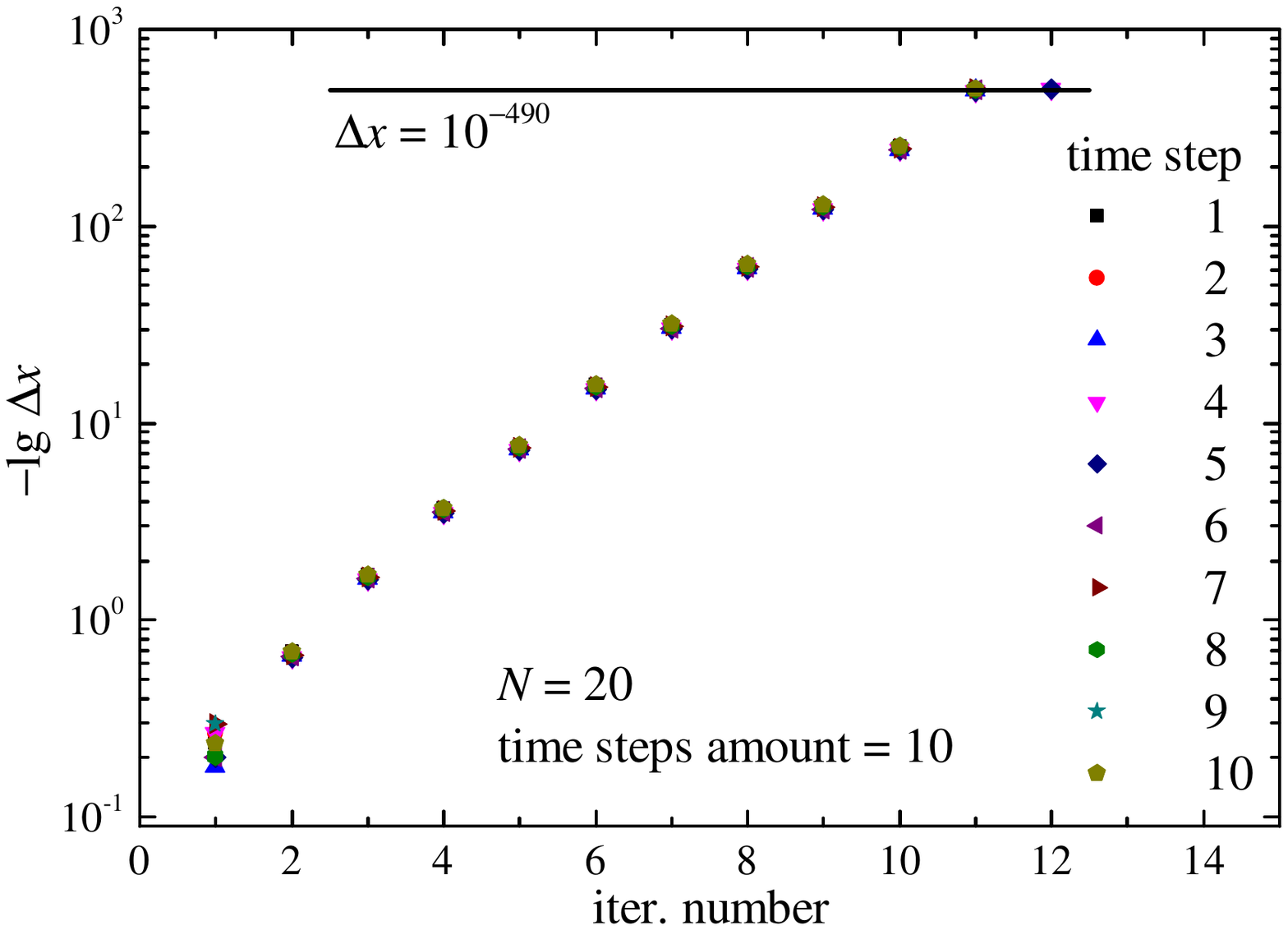}
\vspace{-8mm}\caption{\label{fig:newton_errs:f}\vspace{2mm}}
\end{subfigure}
\begin{subfigure}{0.24\textwidth}
\includegraphics[width=\textwidth]{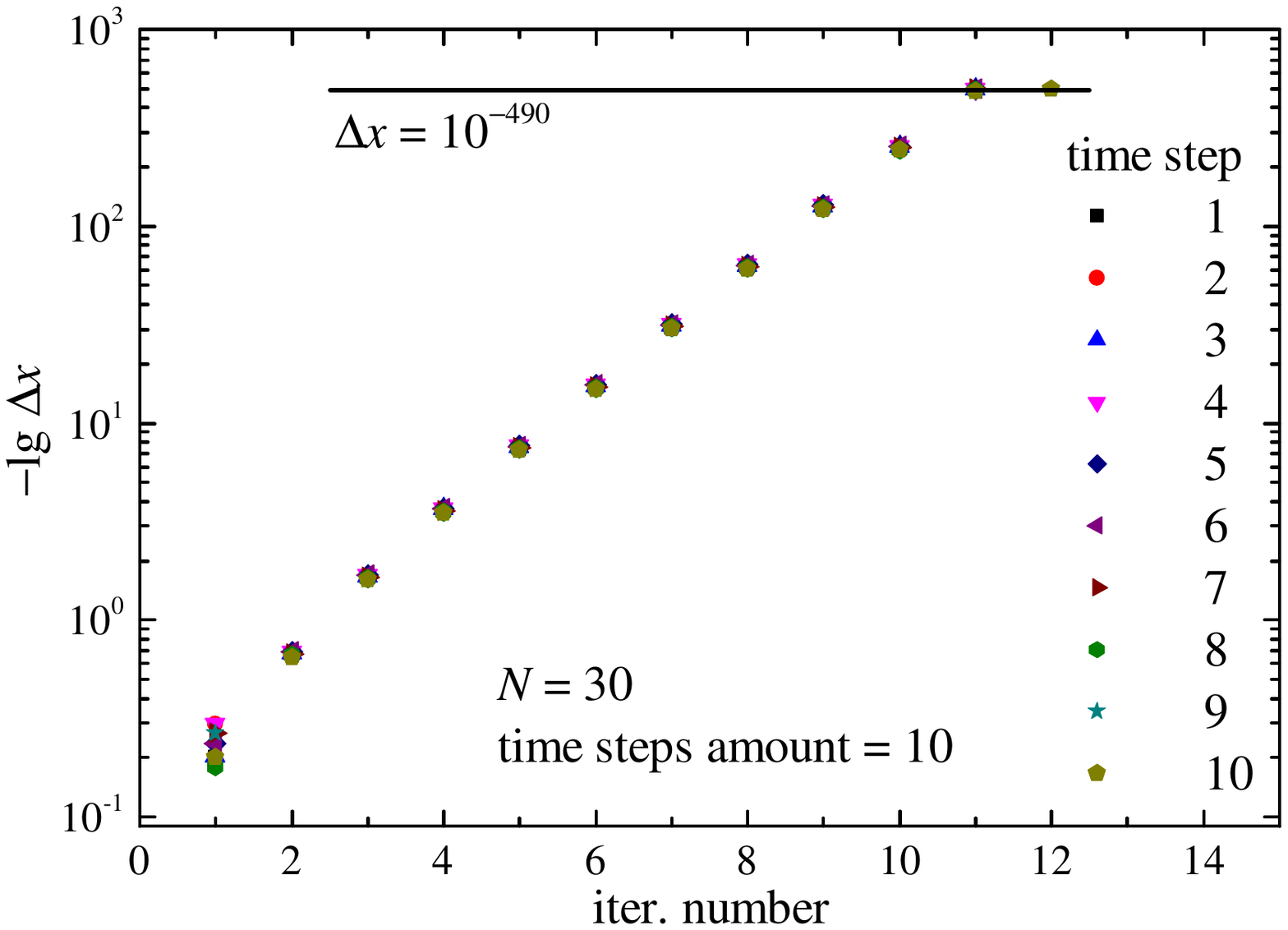}
\vspace{-8mm}\caption{\label{fig:newton_errs:g}\vspace{2mm}}
\end{subfigure}
\begin{subfigure}{0.24\textwidth}
\includegraphics[width=\textwidth]{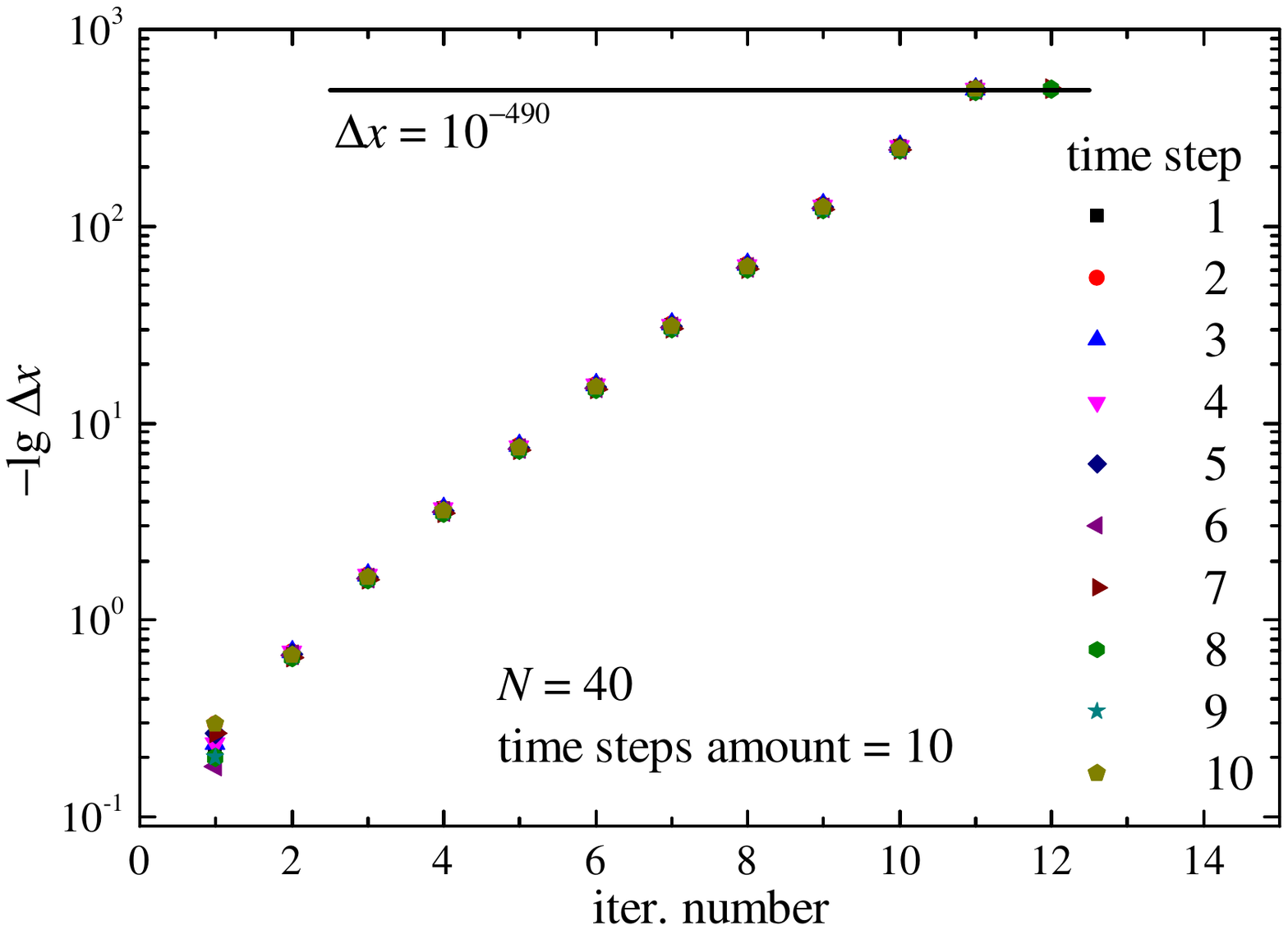}
\vspace{-8mm}\caption{\label{fig:newton_errs:h}\vspace{2mm}}
\end{subfigure}\\
\caption{%
Convergence of iterations of the classical numerical Newton method for solving a system of nonlinear algebraic equations (\ref{eq:lst_dg_newton_slae}) of a local DG predictor: log plot of the negative logarithm $-\lg(\Delta x)$ of the maximum increment $\Delta x = \max(|\Delta\hat{\mathbf{q}}_{p}^{i}|, |\Delta\hat{\mathbf{r}}_{p}^{i}|)$ on the iteration number, in the case of a numerical solution of the problem for a DAE system (\ref{eq:dae_test_newton_conv}) by the ADER-DG method with polynomial degrees $N = 1$ (\subref{fig:newton_errs:a}), $2$ (\subref{fig:newton_errs:b}), $4$ (\subref{fig:newton_errs:c}), $8$ (\subref{fig:newton_errs:d}), $15$ (\subref{fig:newton_errs:e}), $20$ (\subref{fig:newton_errs:f}), $30$ (\subref{fig:newton_errs:g}) and $40$ (\subref{fig:newton_errs:h}), using $10$ time steps (${\Delta t} = 2\pi/10$).
}
\label{fig:newton_errs}
\end{figure}

As a demonstration of the convergence properties of the iterative process of the numerical Newton method for solving the system of nonlinear algebraic equations of the local DG predictor, the following simple DAE system of index 1 was solved:
\begin{equation}\label{eq:dae_test_newton_conv}
\begin{split}
&\ddot{x} + x = z - 1,\quad \ddot{y} + y = 1 - z,\quad x^{2} + y^{2} = z^{2},\\
&x(0) = 1,\quad \dot{x}(0) = 0,\quad y(0) = 0,\quad \dot{y}(0) = 1,\quad z(0) = 1,\quad t\in[0, 2\pi],
\end{split}
\end{equation}
for which the exact analytical solution has the following form: $x = \cos(t)$, $\dot{x} = -\sin(t)$, $y = \sin(t)$, $\dot{y} = \cos(t)$, $z \equiv 1$. For the numerical solution of the presented problem by the ADER-DG method with a local DG predictor, a uniform grid containing $10$ discretization domains was selected --- ${\Delta t} = 2\pi/10$. The value $\varepsilon_{0} = 10^{-490}$ was chosen for the criterion for stopping the iterative process of the numerical Newton method. Fig.~\ref{fig:newton_errs} shows the log plot of the negative logarithm $-\lg(\Delta x)$ of the maximum increment $\Delta x = \max(|\Delta\hat{\mathbf{q}}_{p}^{i}|, |\Delta\hat{\mathbf{r}}_{p}^{i}|)$ on the iteration number, in the case of a numerical solution of the problem for a DAE system (\ref{eq:dae_test_newton_conv}) by the ADER-DG method with polynomial degrees $N = 1$, $2$, $4$, $8$, $15$, $20$, $30$ and $40$. From the data dependencies presented in Fig.~\ref{fig:newton_errs} it is clear that in the case of this problem the Newton method demonstrates the classical second convergence order --- the convergence is clearly ``superlinear'', while before reaching the neighborhood of $\varepsilon_{0}$, a dependence of type $\log(\Delta x) \sim -10^{a \cdot i}$ is observed, where $i$ is the iteration number and the dependence coefficient $a \approx 0.30$--$0.33$. In all the presented cases --- for all the presented degrees $N$ of polynomials and all time steps, convergence up to $\varepsilon_{0}$ was achieved in $11$--$12$ iterations. The accuracy typical for double-precision floating-point numbers $\sim 10^{-15}$ is achieved in $5$--$6$ iterations. It is clear that in the general case the use of the Newton method is characterized by all its well-known disadvantages associated with a decrease in the convergence order in the case of multiple roots, a possible failure of convergence, and the complexity of choosing the initial approximation. In this work, the Newton method was used in all tests, and the solution at the previous node $t_{n}$ was chosen as the initial iteration.

It is necessary to note an important feature of the obtained system of linear algebraic equations (\ref{eq:lst_dg_newton_slae}). The second set of equations of the system (\ref{eq:lst_dg_newton_slae}), which has a block-diagonal structure, can be used to reduce the dimensionality of the problem if one of the sets of matrices --- $\{\mathrm{R}_{p}^{i}\, |\, p = 0,\ldots,N\}$ or $\{\mathrm{S}_{p}^{i}\, |\, p = 0,\ldots,N\}$, contains invertible matrices. Matrices $\mathrm{S}_{p}^{i}$ have a size of $D_{\rm v} \times D_{\rm v}$ and are always square, however, the functions $\mathbf{G}$ may not contain an explicit dependence of algebraic variables $\mathbf{v}$ (in sufficient quantity), therefore the matrices $\mathrm{S}_{p}^{i}$ can be singular. Matrices $\mathrm{R}_{p}^{i}$ have a size of $D_{\rm v} \times D_{\rm u}$ and are usually non-singular in case $D_{\rm v} = D_{\rm u}$, however, in case $D_{\rm v} \neq D_{\rm u}$ they are not square. Therefore, the use of preliminary partial elimination of variables is not a general universal method, but rather refers to problem-specific methods of algorithmic optimization. 

In the first case, when all matrices of the matrix set $\{\mathrm{R}_{p}^{i}\, |\, p = 0,\ldots,N\}$ are invertible, variables can be preliminarily eliminated, which leads to the following relation:
\begin{equation}
\Delta\mathbf{q}_{p}^{i} = \left[\mathrm{R}_{p}^{i}\right]^{-1}\cdot\left[\mathbf{c}_{p}^{i} - \mathrm{S}_{p}^{i}\cdot\Delta\mathbf{r}_{p}^{i}\right].
\end{equation}
Substituting this relation between variables $\Delta\mathbf{q}_{p}^{i}$ and $\Delta\mathbf{r}_{p}^{i}$ into the system of equations (\ref{eq:lst_dg_newton_slae}) leads to the following system of linear algebraic equations:
\begin{equation}\label{eq:lst_dg_newton_slae_v_opt}
\left[
\begin{array}{cccc}
\mathrm{U}_{0,0}^{i} & \mathrm{U}_{0,1}^{i} & \cdots & \mathrm{U}_{0,N}^{i} \\
\mathrm{U}_{1,0}^{i} & \mathrm{U}_{1,1}^{i} & \cdots & \mathrm{U}_{1,N}^{i} \\
\vdots & \vdots & \ddots & \vdots \\
\mathrm{U}_{N,0}^{i} & \mathrm{U}_{N,1}^{i} & \cdots & \mathrm{U}_{N,N}^{i}
\end{array}
\right]\left[
\begin{array}{c}
\Delta\mathbf{r}_{0}^{i}\\
\Delta\mathbf{r}_{1}^{i}\\
\vdots\\
\Delta\mathbf{r}_{N}^{i}
\end{array}
\right] = \left[
\begin{array}{c}
\mathbf{d}_{0}^{i}\\
\mathbf{d}_{1}^{i}\\
\vdots\\
\mathbf{d}_{N}^{i}
\end{array}
\right],
\end{equation}
where the following notations were introduced:
\begin{equation}
\begin{split}
\mathrm{U}_{p, q}^{i} = \mathrm{Q}_{p, q}^{i} - \mathrm{P}_{p, q}^{i}\cdot\left[\mathrm{R}_{p}^{i}\right]^{-1}\cdot\mathrm{S}_{p}^{i},\qquad
\mathbf{d}_{p}^{i} = \mathbf{b}_{p}^{i} - \mathrm{P}_{p, q}^{i}\cdot\left[\mathrm{R}_{p}^{i}\right]^{-1}\cdot\mathbf{c}_{p}^{i}.
\end{split}
\end{equation}
The matrix size of the resulting reduced system of linear algebraic equations is $[(N+1)D_{\rm v}] \times [(N+1)D_{\rm v}]$. In the case of using direct numerical solution methods based on LU decomposition, the computational complexity will be $O(D_{\rm v}^{3}(N+1)^{3})$, which can be significantly less than $O(D^{3}N^{3})$ (it is necessary to take into account that matrices $\mathrm{R}_{p}^{i}$ are square only in the case when the number of differential $\mathbf{u}$ and algebraic variables $\mathbf{v}$ is the same --- $D_{\rm v} = D_{\rm u} = D/2$, then a decrease in the matrix size by $2$ times can asymptotically decrease the computational complexity by $8$ times). However, in this case it is also necessary to take into account the computational costs of calculating the inverse matrices by $\mathrm{R}_{p}^{i}$ with a size of $D_{\rm v} \times D_{\rm v}$ each, the computational costs of calculating $\mathrm{U}_{p, q}^{i}$ the matrices included in the block structure of the system of linear equations, as well as the computational costs of calculating the variables $\Delta\mathbf{q}_{p}^{i}$ from the obtained variables $\Delta\mathbf{r}_{p}^{i}$.

In the second case, when all matrices of the matrix set $\{\mathrm{S}_{p}^{i}\, |\, p = 0,\ldots,N\}$ are invertible, variables can be preliminarily eliminated, which leads to the following relation:
\begin{equation}
\Delta\mathbf{r}_{p}^{i} = \left[\mathrm{S}_{p}^{i}\right]^{-1}\cdot\left[\mathbf{c}_{p}^{i} - \mathrm{R}_{p}^{i}\cdot\Delta\mathbf{q}_{p}^{i}\right].
\end{equation}
Substituting this relation between variables $\Delta\mathbf{r}_{p}^{i}$ and $\Delta\mathbf{q}_{p}^{i}$ into the system of equations (\ref{eq:lst_dg_newton_slae}) leads to the following system of linear algebraic equations:
\begin{equation}\label{eq:lst_dg_newton_slae_u_opt}
\left[
\begin{array}{cccc}
\mathrm{V}_{0,0}^{i} & \mathrm{V}_{0,1}^{i} & \cdots & \mathrm{V}_{0,N}^{i} \\
\mathrm{V}_{1,0}^{i} & \mathrm{V}_{1,1}^{i} & \cdots & \mathrm{V}_{1,N}^{i} \\
\vdots & \vdots & \ddots & \vdots \\
\mathrm{V}_{N,0}^{i} & \mathrm{V}_{N,1}^{i} & \cdots & \mathrm{V}_{N,N}^{i}
\end{array}
\right]\left[
\begin{array}{c}
\Delta\mathbf{q}_{0}^{i}\\
\Delta\mathbf{q}_{1}^{i}\\
\vdots\\
\Delta\mathbf{q}_{N}^{i}
\end{array}
\right] = \left[
\begin{array}{c}
\mathbf{h}_{0}^{i}\\
\mathbf{h}_{1}^{i}\\
\vdots\\
\mathbf{h}_{N}^{i}
\end{array}
\right],
\end{equation}
where the following notations were introduced:
\begin{equation}
\begin{split}
\mathrm{V}_{p, q}^{i} = \mathrm{P}_{p, q}^{i} - \mathrm{Q}_{p, q}^{i}\cdot\left[\mathrm{S}_{p}^{i}\right]^{-1}\cdot\mathrm{R}_{p}^{i},\qquad
\mathbf{h}_{p}^{i} = \mathbf{b}_{p}^{i} - \mathrm{Q}_{p, q}^{i}\cdot\left[\mathrm{S}_{p}^{i}\right]^{-1}\cdot\mathbf{c}_{p}^{i},
\end{split}
\end{equation}
The matrix size of the resulting reduced system of linear algebraic equations is $[(N+1)D_{\rm u}] \times [(N+1)D_{\rm u}]$. In the case of using direct numerical methods of solution based on LU-decomposition, the computational complexity will be $O(D_{\rm u}^{3}(N+1)^{3})$, which can be significantly less than $O(D^{3}N^{3})$ in case $D_{\rm u} \leqslant D_{\rm v}$. However, in this case it is also necessary to take into account the computational costs of calculating the inverse matrices by $\mathrm{S}_{p}^{i}$ with the size $D_{\rm u} \times D_{\rm u}$ each, the computational costs of calculating $\mathrm{V}_{p, q}^{i}$ the matrices included in the block structure of the system of linear equations, as well as the computational costs of calculating the variables $\Delta\mathbf{r}_{p}^{i}$ from the obtained variables $\Delta\mathbf{q}_{p}^{i}$.

It should also be noted that these approaches can be used partially and jointly, in the case of the presence of singular or non-square matrices in the sets $\{\mathrm{R}_{p}^{i}\, |\, p = 0,\ldots,N\}$ or $\{\mathrm{S}_{p}^{i}\, |\, p = 0,\ldots,N\}$. These approaches relate to problem-specific methods of algorithmic optimization, therefore, like any other optimization methods, they require quantitative estimates of the gain in computational costs. Often the number of variables $D$ and polynomial degrees $N$ are small in applied calculations, and computer architectures can support vector extensions instruction sets, and in such cases the direct gain in computational costs by such algorithmic optimizations may be small or nonexistent.

\subsection{Stability of the numerical method}

In the present work, the numerical method ADER-DG with a local DG predictor is investigated, the functional representations of which are based on the expansion of a local discrete time solution $\mathbf{q}$ over a set of Lagrange interpolation polynomials $\varphi_{p}$ with nodal points at the roots of the right Radau polynomials. In the work~\cite{ader_dg_ivp_ode}, where the numerical method ADER-DG with a local DG predictor was investigated for solving initial value problems for ODE systems, Lagrange interpolation polynomials with nodal points at the roots of shifted Legendre polynomials were used, and it was shown that this version of the numerical method has $A$-stability and $L$-stability (more precisely, $L_{1}$-stability). Therefore, it is of certain interest to determine the linear stability properties of the ADER-DG method with a local DG predictor, in which the node of the functional basis is constructed at the nodes of the Gauss-Radau quadrature formula (\ref{eq:gr_rule}).

The study of the linear stability of the numerical method ADER-DG with a local DG predictor, the functional representations of which are based on the expansion of a local discrete time solution by a set of Lagrange interpolation polynomials with nodal points at the roots of the right Radau polynomials, was carried out similarly to the study in the work~\cite{ader_dg_ivp_ode}. 

The linear stability was studied based on the analysis of the Dahlquist's test equation $\dot{u} = \lambda u$. The numerical solution of the Dahlquist's test equation has the form $u_{n+1} = R(\lambda\cdot\Delta t_{n})\, u_{n}$, where $R = R(z)$ is called the stability function of the numerical method. The stability function argument $z = \lambda\cdot\Delta t_{n}$ is a complex number: $z = \Re(z) + i\cdot\Im(z) \in \mathbb{C}$. Using the equation (\ref{eq:lst_dg_snae}) and formula (\ref{eq:u_sol_in_node}), the following result was obtained:
\begin{equation}
\begin{split}
&\sum\limits_{q = 0}^{N}\left[\delta_{pq} - \left(\lambda\cdot\Delta t_{n}\right)\mathrm{A}_{pq}\right]\cdot\hat{q}_{q} = u_{n},\\
&u_{n+1} = u_{n} + \left(\lambda\cdot\Delta t_{n}\right) \sum\limits_{p = 0}^{N} w_{p}\, \hat{q}_{p},
\end{split}
\end{equation}
which allows using the well-known expression used to analyze the stability of implicit Runge-Kutta methods~\cite{Butcher_book_2016, Hairer_book_1, Hairer_book_2} using matrix coefficients $\mathrm{A} = ||\mathrm{A}_{pq}||$ and vector weights $\mathbf{w} = [w_{0} \ldots w_{p}]^{T}$ of the Gauss-Radau quadrature formula (\ref{eq:gr_rule}):
\begin{equation}\label{eq:stab_func_expr}
R(z) = 1 + \mathbf{w}^{T} \cdot \left[\mathbb{I}_{d \times d} - z\, \mathrm{A}\right]^{-1} \cdot \mathbf{1} 
=\frac{\left|\mathrm{E} - z\, \mathrm{A} + z\, \mathbf{1}\otimes\mathbf{w}\right|}{\left|\mathrm{E} - z\mathrm{A}\right|},
\end{equation}
where $\mathbb{I}_{d \times d}$ is the identity matrix with sizes $(N+1)\times(N+1)$, $\mathbf{1} = [1 \ldots 1]^{T}$ is a vector of ones with size $N+1$, $\otimes$ is the tensor product.

\begin{figure}[h!]
\captionsetup[subfigure]{%
	position=bottom,
	font+=smaller,
	textfont=normalfont,
	singlelinecheck=off,
	justification=raggedright
}
\centering
\begin{subfigure}{0.24\textwidth}
\includegraphics[width=\textwidth]{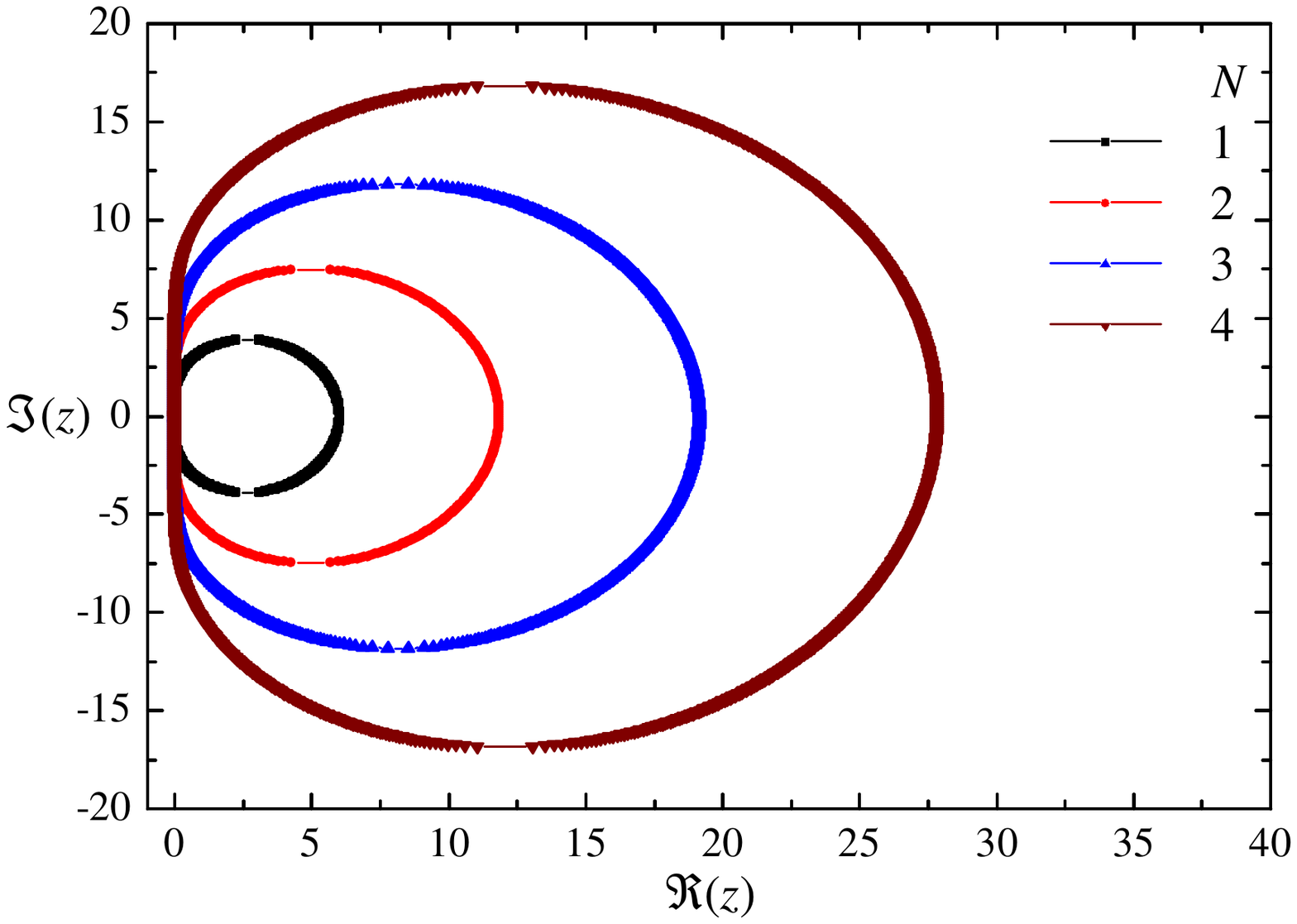}
\vspace{-8mm}\caption{\label{fig:stab_domain:a}\vspace{2mm}}
\end{subfigure}
\begin{subfigure}{0.24\textwidth}
\includegraphics[width=\textwidth]{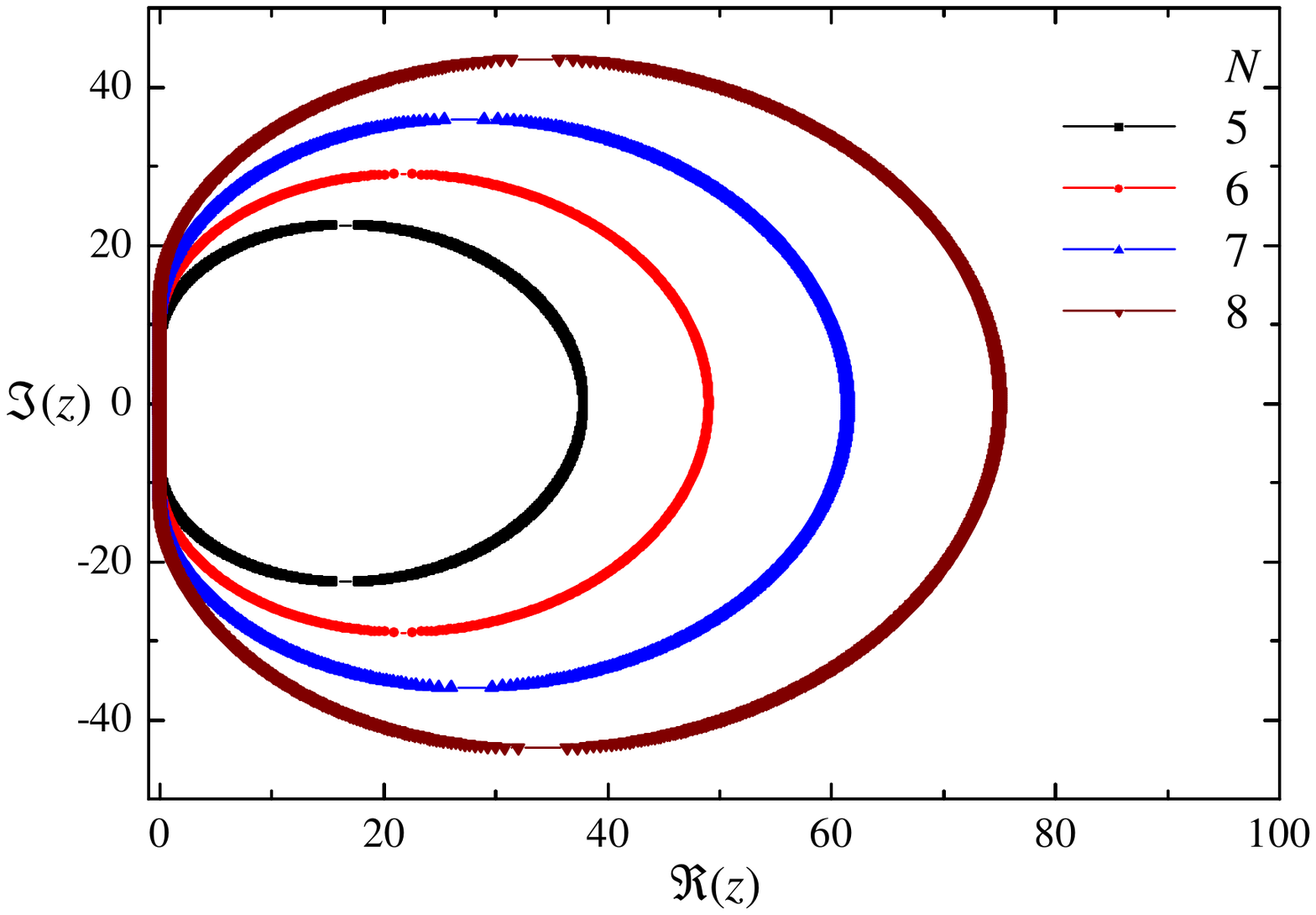}
\vspace{-8mm}\caption{\label{fig:stab_domain:b}\vspace{2mm}}
\end{subfigure}
\begin{subfigure}{0.24\textwidth}
\includegraphics[width=\textwidth]{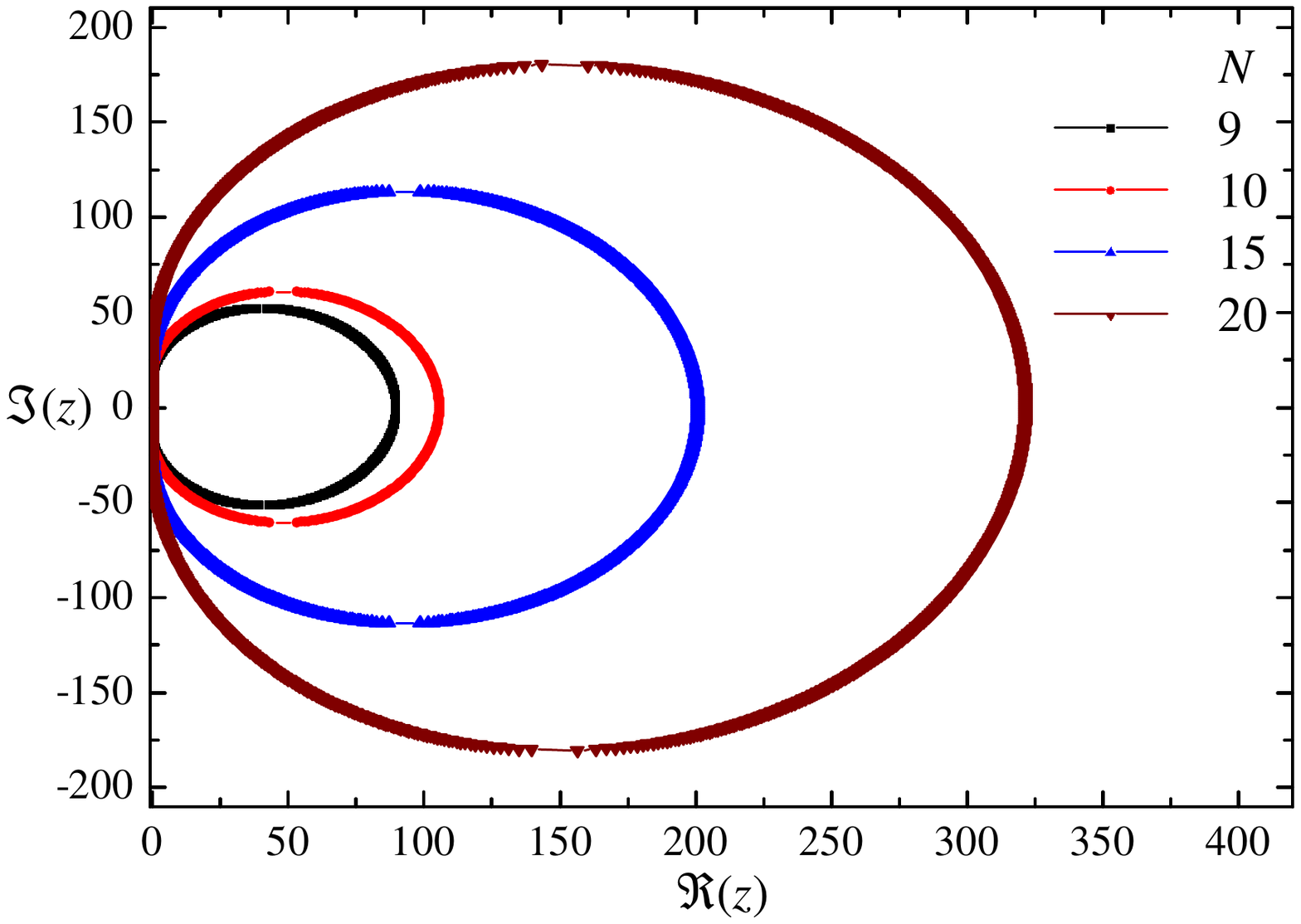}
\vspace{-8mm}\caption{\label{fig:stab_domain:c}\vspace{2mm}}
\end{subfigure}
\begin{subfigure}{0.24\textwidth}
\includegraphics[width=\textwidth]{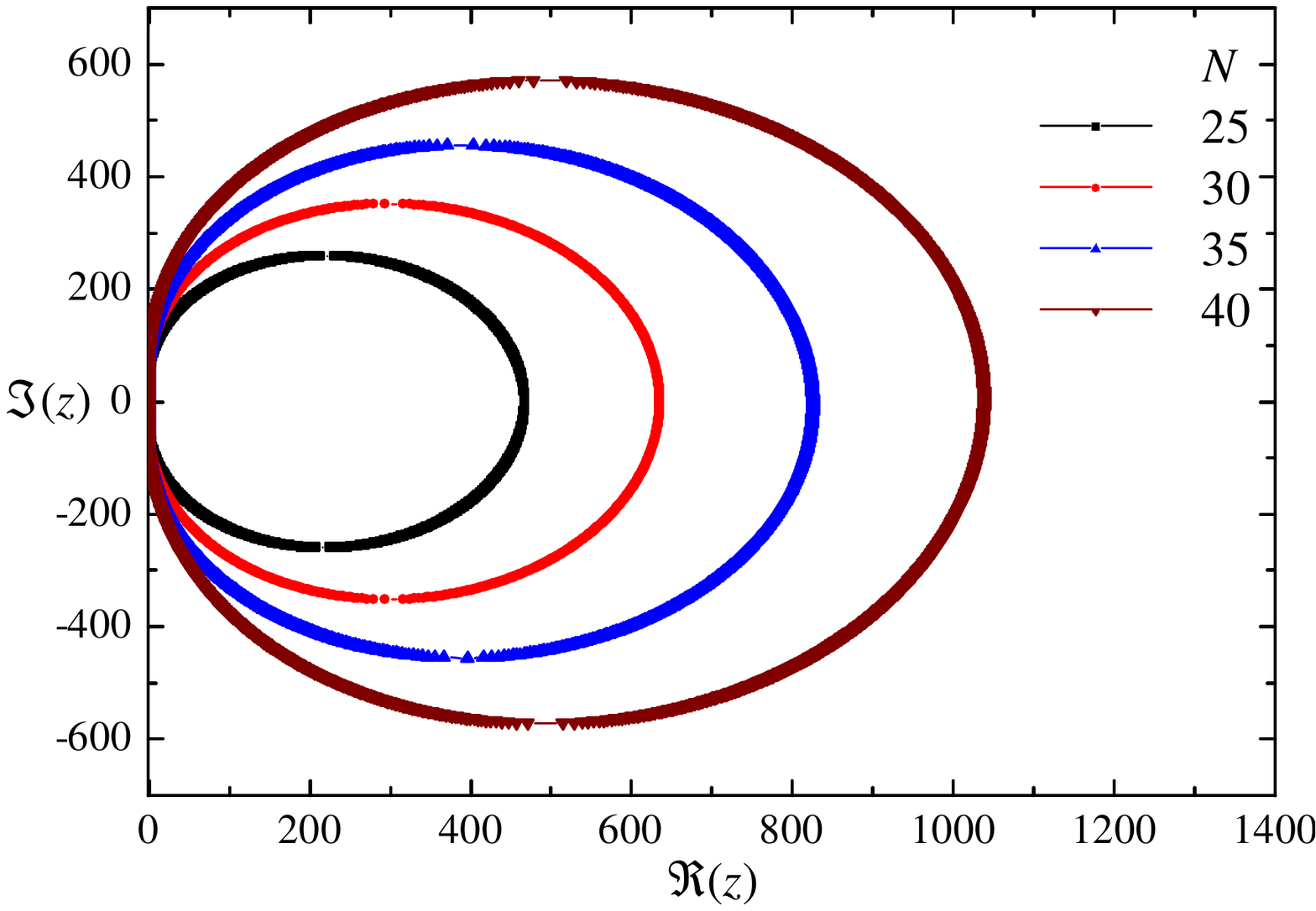}
\vspace{-8mm}\caption{\label{fig:stab_domain:d}\vspace{2mm}}
\end{subfigure}
\caption{%
Regions of absolute stability $|R(z)| < 1$ (outer part of closed curves) in the complex plane $z$ of the ADER-DG numerical method with a local DG predictor for the degrees $1 \leqslant N \leqslant 40$ of polynomials in the DG representation: $N = 1,\, 2,\, 3,\, 4$ (a), $N = 5,\, 6,\, 7,\, 8$ (b), $N = 9,\, 10,\, 15,\, 20$ (c), $N = 25,\, 30,\, 35,\, 40$ (d). Horizontal and vertical axes represent real $\Re(z)$ and imaginary $\Im(z)$ parts of complex number $z = \lambda\cdot\mathrm{\Delta}t_{n}$. In the region $\Re(z) < 0$ the inequality $|R(z)| < 1$ is satisfied, therefore the ADER-DG numerical method with a local DG predictor is $A$-stable.
}
\label{fig:stab_domain}
\end{figure}

\begin{figure}[h!]
\captionsetup[subfigure]{%
	position=bottom,
	font+=smaller,
	textfont=normalfont,
	singlelinecheck=off,
	justification=raggedright
}
\centering
\begin{subfigure}{0.24\textwidth}
\includegraphics[width=\textwidth]{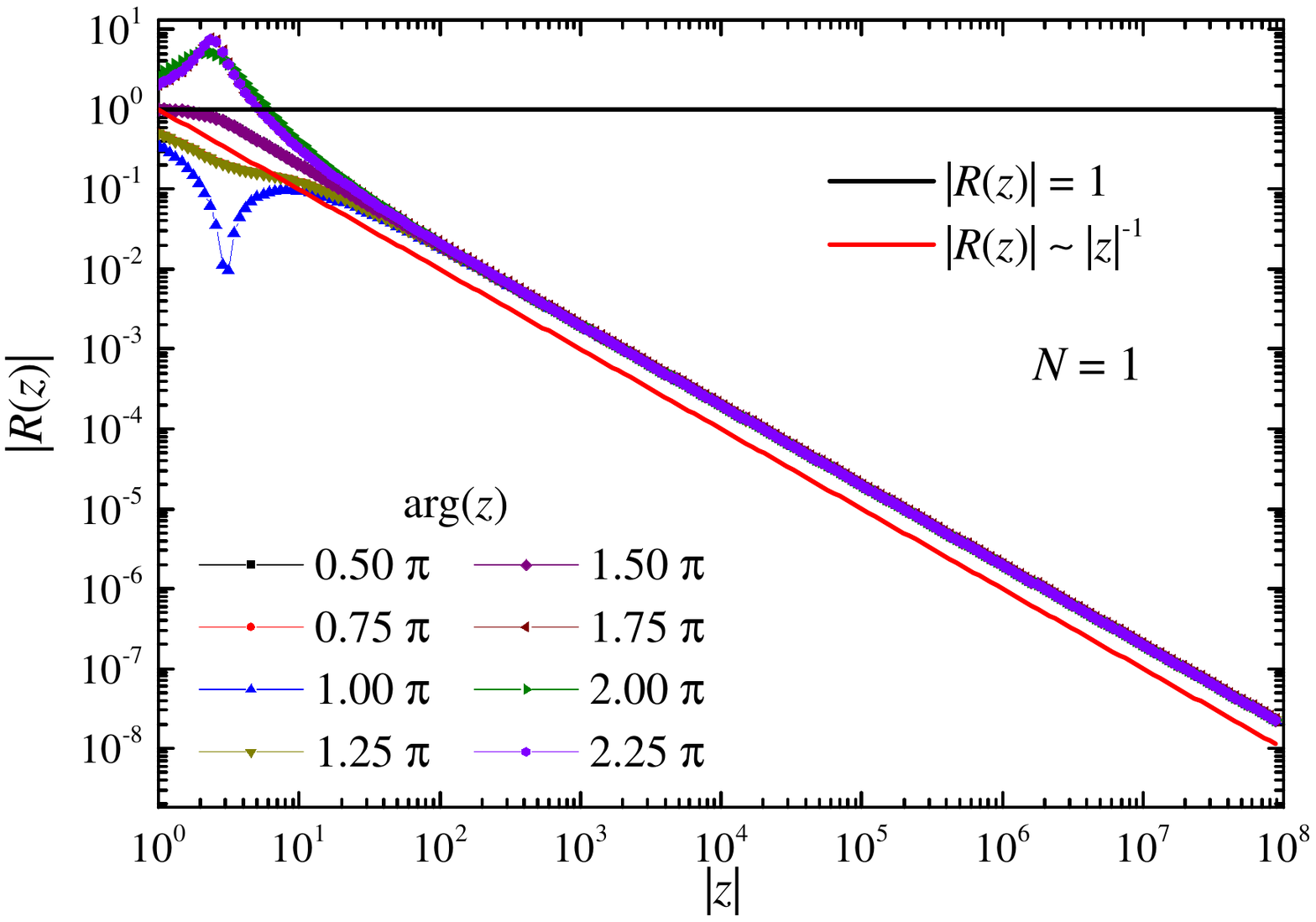}
\vspace{-8mm}\caption{\label{fig:stab_function_rays:a}\vspace{2mm}}
\end{subfigure}
\begin{subfigure}{0.24\textwidth}
\includegraphics[width=\textwidth]{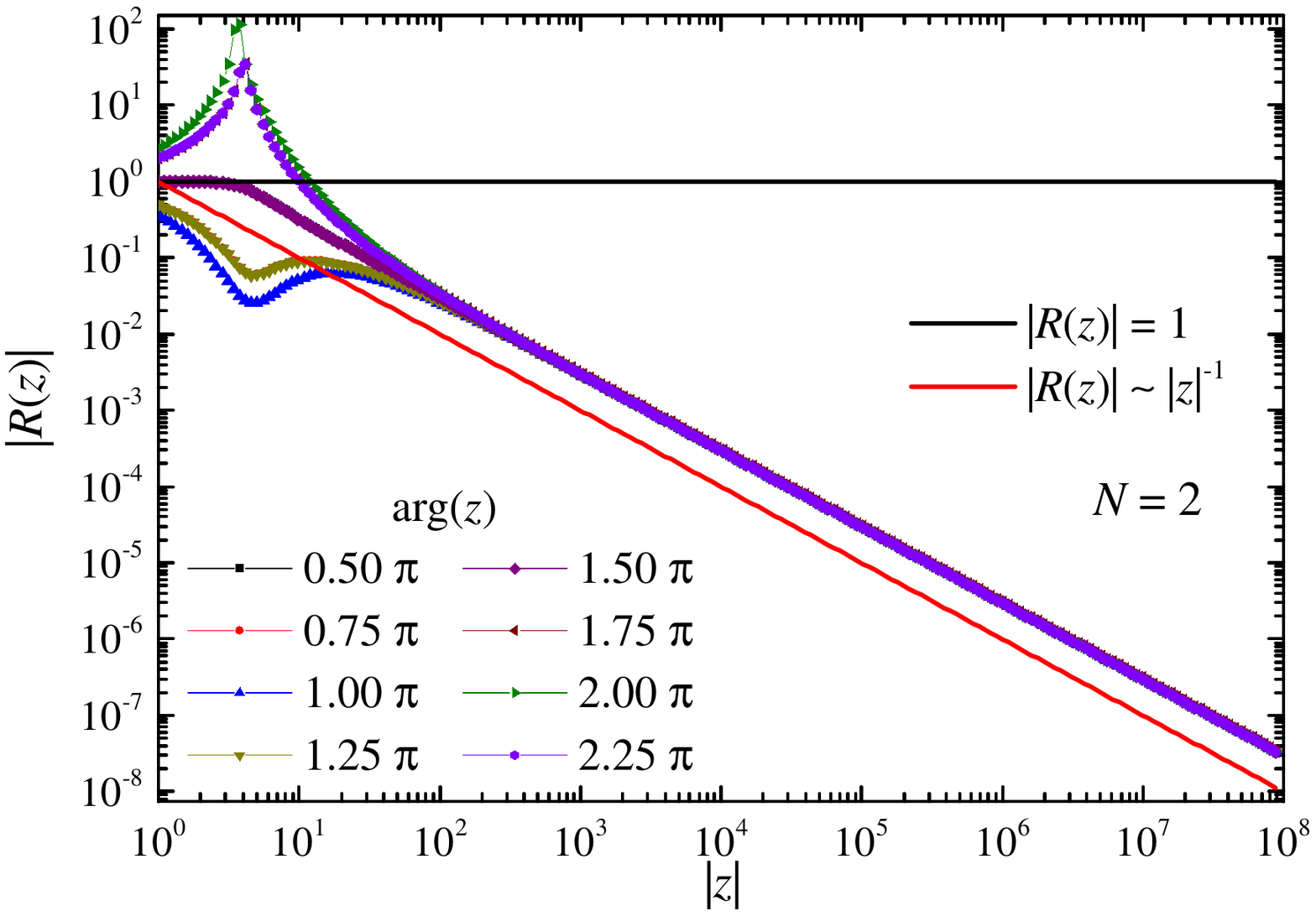}
\vspace{-8mm}\caption{\label{fig:stab_function_rays:c}\vspace{2mm}}
\end{subfigure}
\begin{subfigure}{0.24\textwidth}
\includegraphics[width=\textwidth]{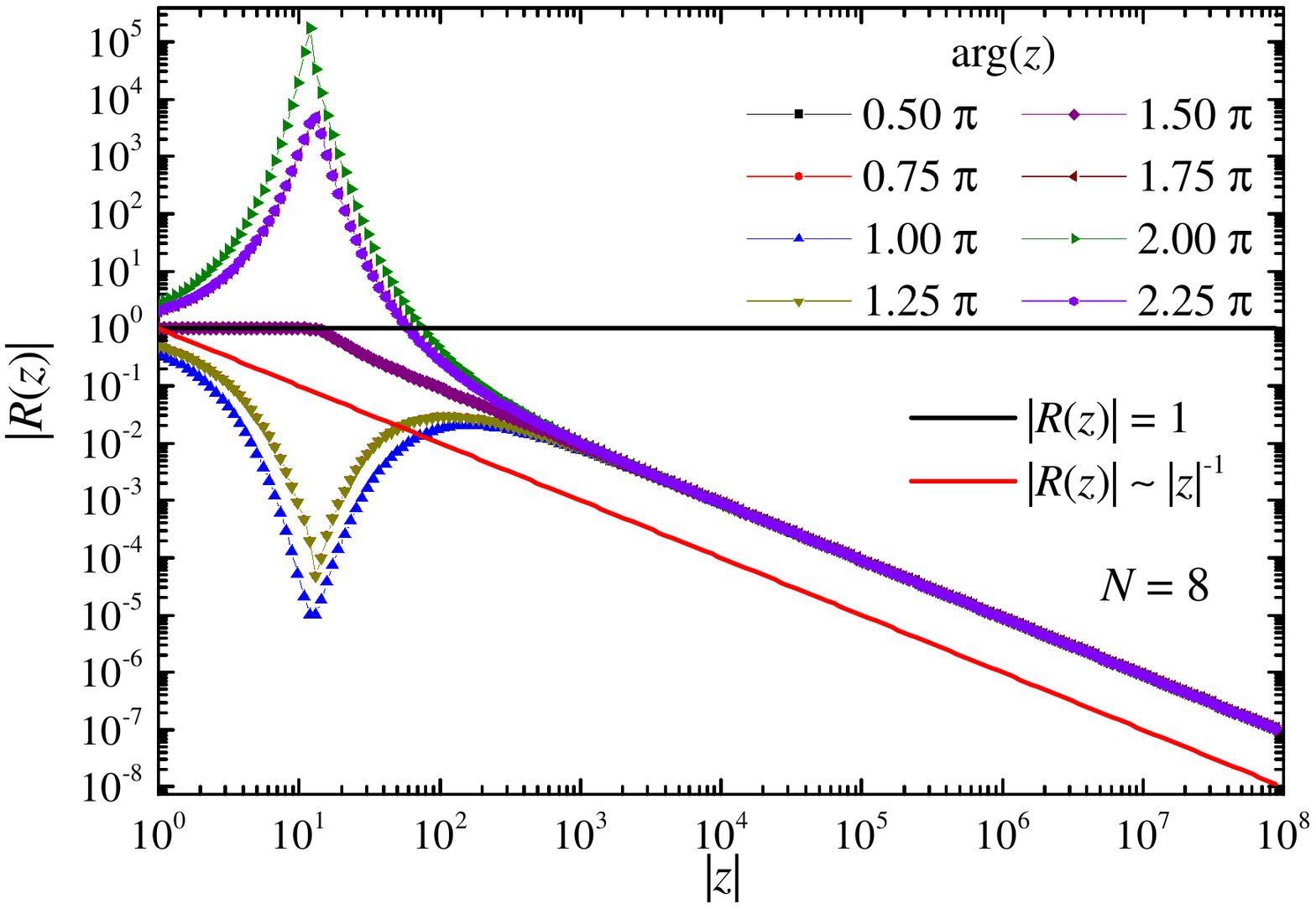}
\vspace{-8mm}\caption{\label{fig:stab_function_rays:d}\vspace{2mm}}
\end{subfigure}
\begin{subfigure}{0.24\textwidth}
\includegraphics[width=\textwidth]{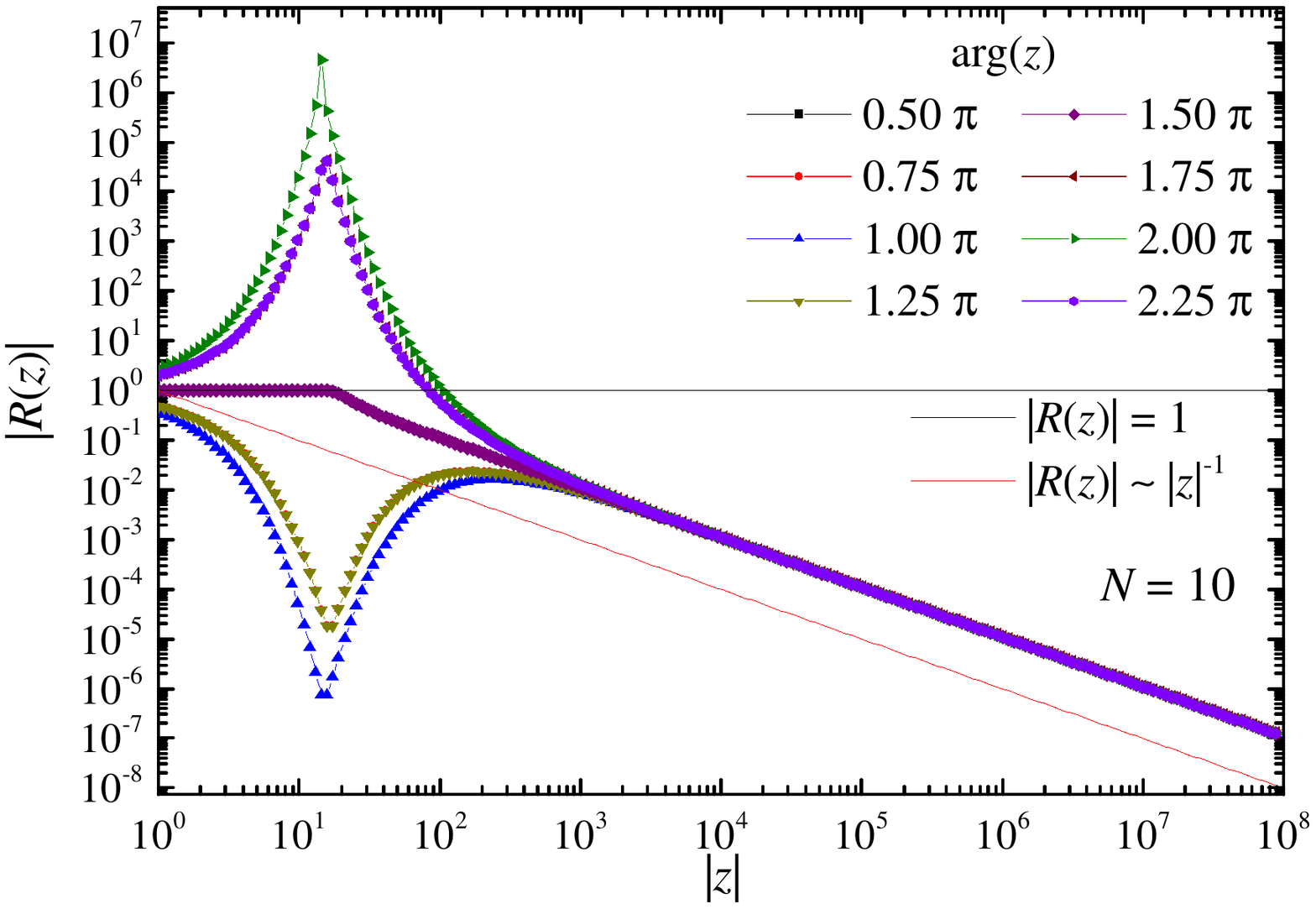}
\vspace{-8mm}\caption{\label{fig:stab_function_rays:e}\vspace{2mm}}
\end{subfigure}\\
\begin{subfigure}{0.24\textwidth}
\includegraphics[width=\textwidth]{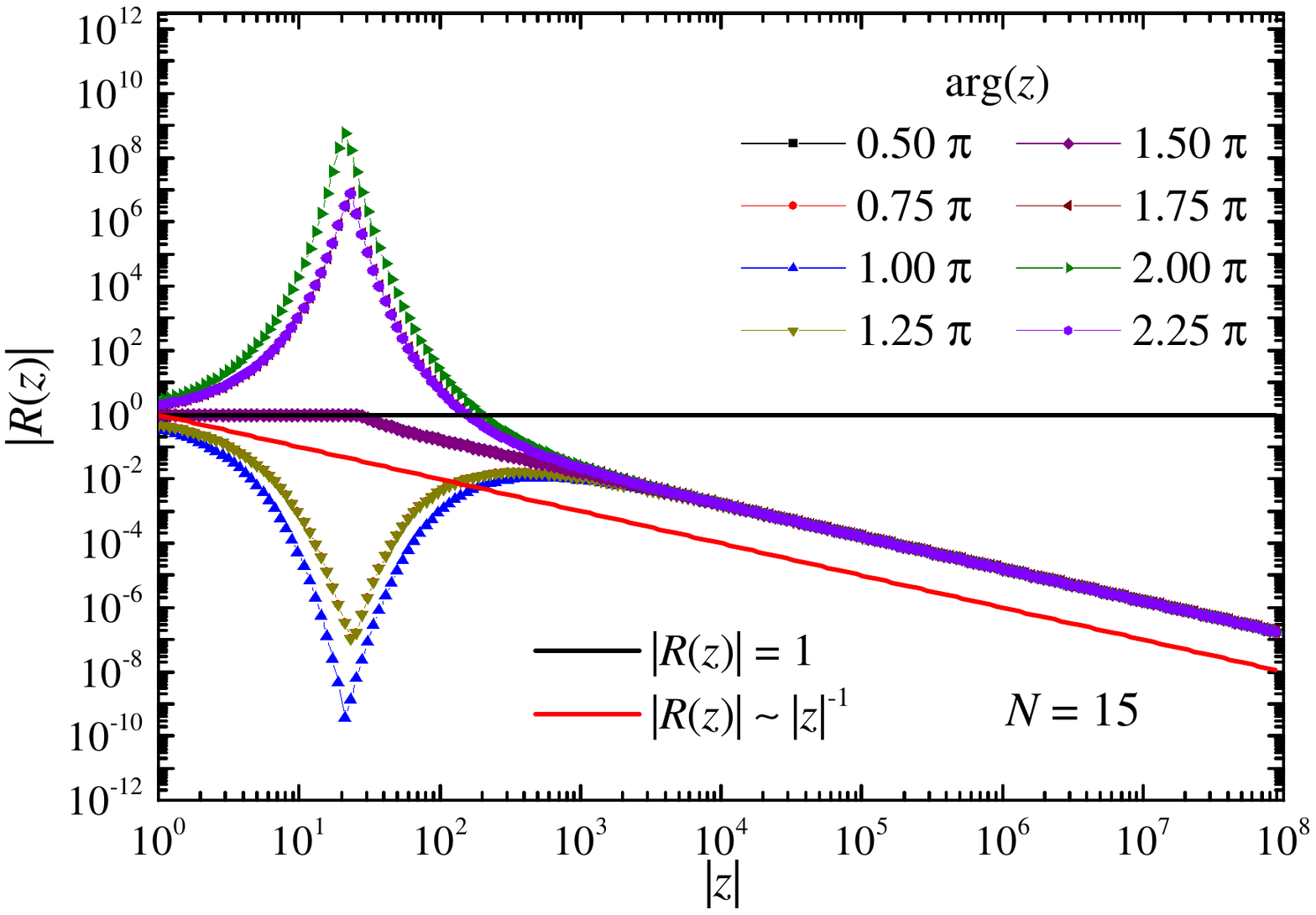}
\vspace{-8mm}\caption{\label{fig:stab_function_rays:f}\vspace{2mm}}
\end{subfigure}
\begin{subfigure}{0.24\textwidth}
\includegraphics[width=\textwidth]{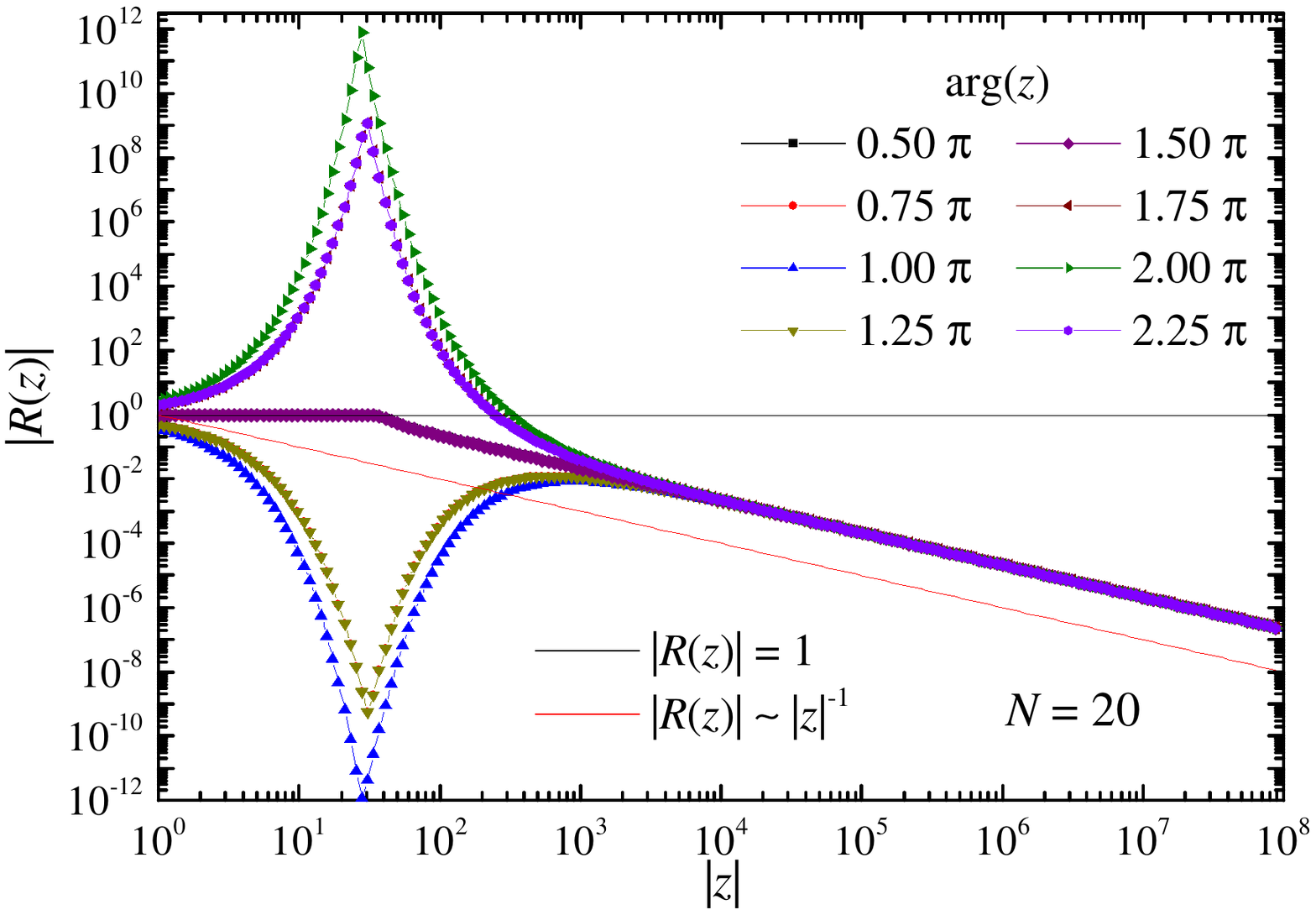}
\vspace{-8mm}\caption{\label{fig:stab_function_rays:g}\vspace{2mm}}
\end{subfigure}
\begin{subfigure}{0.24\textwidth}
\includegraphics[width=\textwidth]{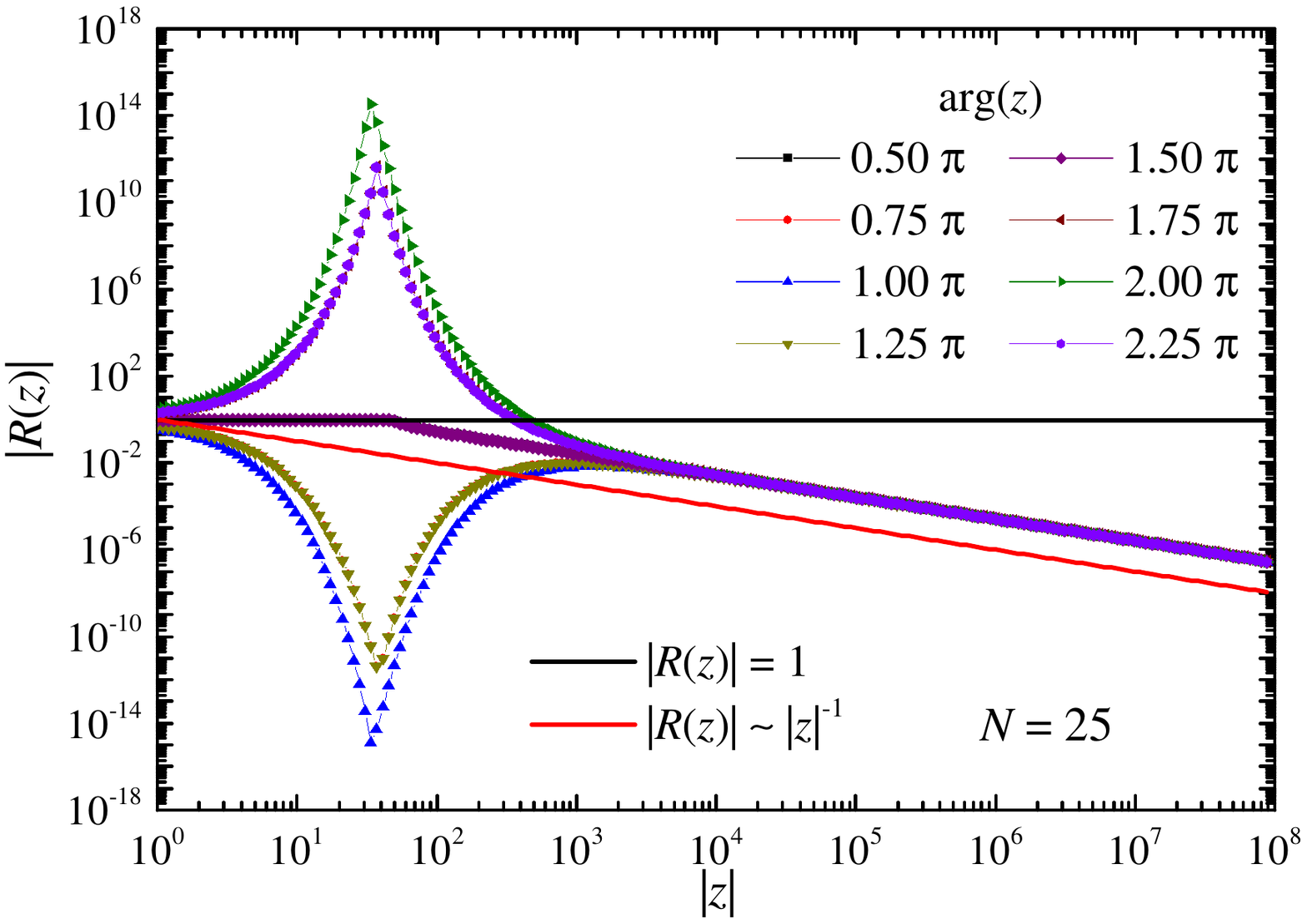}
\vspace{-8mm}\caption{\label{fig:stab_function_rays:h}\vspace{2mm}}
\end{subfigure}
\begin{subfigure}{0.24\textwidth}
\includegraphics[width=\textwidth]{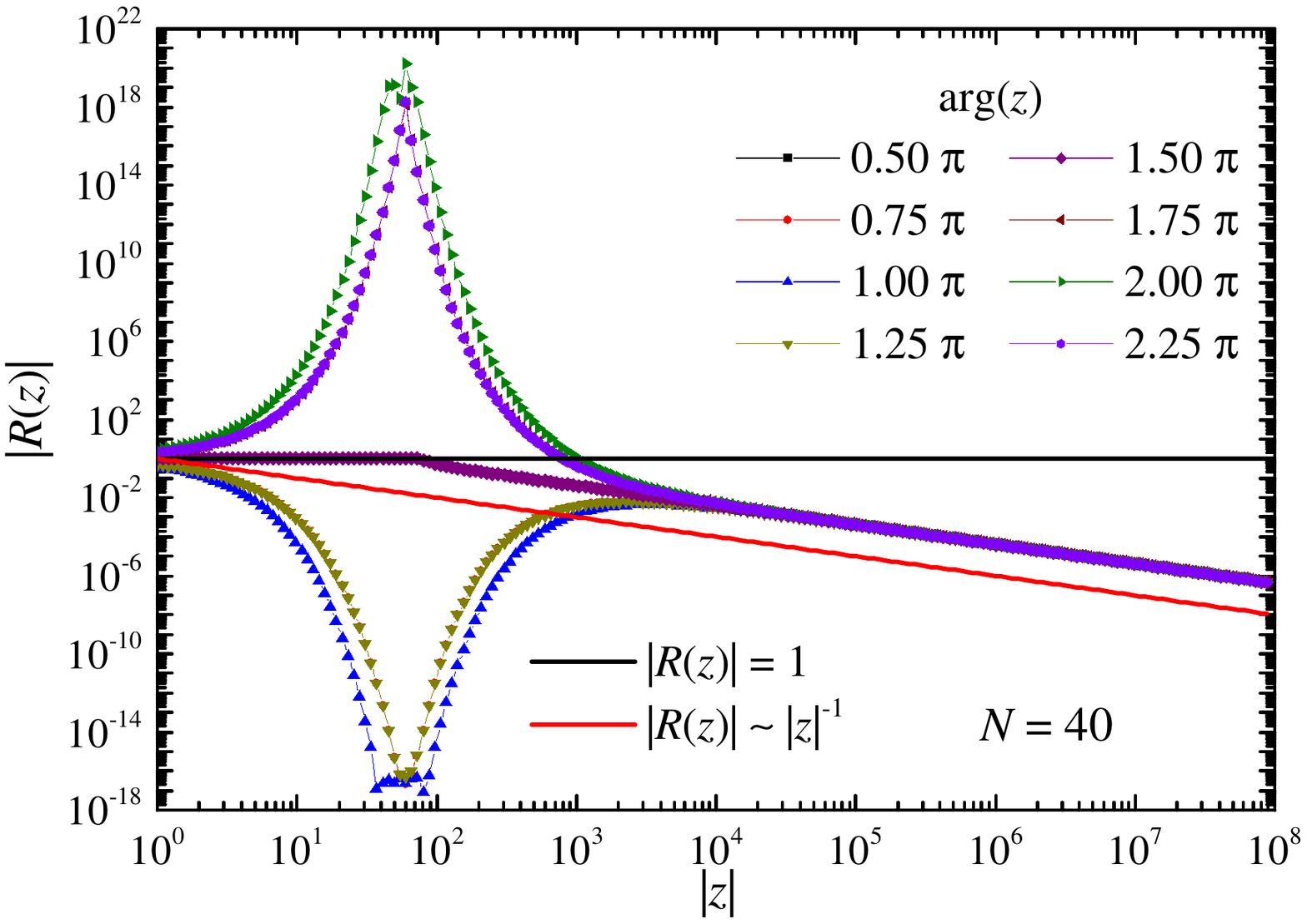}
\vspace{-8mm}\caption{\label{fig:stab_function_rays:i}\vspace{2mm}}
\end{subfigure}\\
\caption{%
The absolute values of the stability function $|R(z)|$ of the ADER-DG numerical method with a local DG predictor for the degrees $N$ on a set of radial rays $z = |z|\exp(i\,\arg(z))$ in the complex plane $z$, with $0.5\pi \leqslant \arg(z) < 2.5\pi$: $N = 1$~(a), $2$~(b), $8$~(c), $10$~(d), $15$~(e), $20$~(g), $25$~(h), $40$~(i). The line $|R(z)| \sim |z|^{-1}$ demonstrates the asymptotic behavior appropriate in the $L$-stability. The line $|R(z)| = 1$ shows the region of absolute stability; the rays $\arg(z) = 1.75\pi,\, 2.0\pi,\, 2.25\pi$ have unstable regions with $|R(z)| > 1$ corresponding to the regions inside the closed curves in Fig.~\ref{fig:stab_domain}. The stability function $|R(z)|$ on the rays $\arg(z) = 0.5\pi,\, 0.75\pi, 1.75\pi$ are the same as $|R(z)|$ on the rays $\arg(z) = 1.5\pi,\, 1.25\pi,\, 2.25\pi$, respectively, which is due to complex conjugation.
}
\label{fig:stab_function_rays}
\end{figure}

The calculated regions of absolute stability $|R(z)| < 1$ are presented in Fig.~\ref{fig:stab_domain} for the degrees $1 \leqslant N \leqslant 40$ of polynomials in the DG representation. In the range $\Re(z) < 0$ the numerical method is absolutely stable. Therefore, it was concluded that the ADER-DG numerical method with a local DG predictor, the functional representations of which are based on the expansion of a local discrete time solution by a set of Lagrange interpolation polynomials with nodal points at the roots of the right Radau polynomials, is $A$-stable. The obtained result completely corresponds to the result of the work~\cite{ader_dg_ivp_ode}, where the ADER-DG with a local DG predictor was investigated for solving initial value problems for ODE systems, Lagrange interpolation polynomials with nodal points at the roots of shifted Legendre polynomials were used.

It should be noted that the determinants in the numerator and in the denominator of the final expression (\ref{eq:stab_func_expr}) for the stability function $R$ represent characteristic polynomials for matrices $[\mathrm{A} - \mathbf{1}\otimes\mathbf{w}]$ and $\mathrm{A}$, in which numbers $z$ that are reciprocals $z = 1/\nu$ of the classical definition of a characteristic polynomial $|\mathrm{A} - \nu\, \mathrm{E}|$ are presented as arguments (taking out the value $z^{N+1}$ in the numerator and denominator does not change the final expression).

The direct calculations showed that the rank of the matrix $[\mathrm{A} - \mathbf{1}\otimes\mathbf{w}]$ equals $N$, and the rank of the matrix $\mathrm{A}$ equals $N+1$ for the values of the degrees $1 \leqslant N \leqslant 40$, so the rank of the matrix $[\mathrm{A} - \mathbf{1}\otimes\mathbf{w}]$ is one less than the rank of the matrix $\mathrm{A}$. The expansion of the determinant in the numerator of the expression (\ref{eq:stab_func_expr}) always has a factor of the form $z = 1/\nu$. Therefore, after reducing the final fraction by a factor $z^{N+1}$, the numerator represents a polynomial of degree $N$ from the argument $z$, and the denominator represents a polynomial of degree $N+1$ from the argument $z$. Therefore, it is expected that the absolute value of the stability function $R(z)$ will have an asymptotic behavior of the form $|R(z)| \sim |z|^{-1}$ in the domain $|z|\rightarrow\infty$ for any directions in the complex plane (this was expected --- the stability function $R$ is a rational function). Therefore, it was concluded that the ADER-DG numerical method with a local DG predictor, the functional representations of which are based on the expansion of a local discrete time solution by a set of Lagrange interpolation polynomials with nodal points at the roots of the right Radau polynomials, is $L$-stable; more precisely --- $L_{1}$-stable. The obtained result completely corresponds to the result of the work~\cite{ader_dg_ivp_ode}, where the ADER-DG with a local DG predictor was investigated for solving initial value problems for ODE systems, and Lagrange interpolation polynomials with nodal points at the roots of shifted Legendre polynomials were used.

In Fig.~\ref{fig:stab_function_rays} additionally presented are the absolute values of the stability function $|R(z)|$ for several basic values of degrees $N$ of polynomials in the DG representation on a set of radial rays $z = |z|\exp(i\,\arg(z))$ in the complex plane $z$, with $0.5\pi \leqslant \arg(z) < 2.5\pi$; the reference point $0.5\pi$ is chosen for the boundary ray completely passing in the region of absolute stability (see Fig.~\ref{fig:stab_domain}). The line $|R(z)| \sim |z|^{-1}$ demonstrates the asymptotic behavior appropriate in the $L$-stability. The values on the rays with $\arg(z) = 0.5\pi,\, 0.75\pi, 1.75\pi$ are the same as the values on the rays with $\arg(z) = 1.5\pi,\, 1.25\pi,\, 2.25\pi$, respectively, which is due to complex conjugation. The rays $\arg(z) = 1.75\pi,\, 2.0\pi,\, 2.25\pi$ have unstable regions with $|R(z)| > 1$ corresponding to the regions inside the closed curves in Fig.~\ref{fig:stab_domain}.

Therefore, it can be concluded that using the functional representations which are based on the expansion of a local discrete time solution by a set of Lagrange interpolation polynomials with nodal points at the roots of the right Radau polynomials instead of Lagrange interpolation polynomials with nodal points at the roots of shifted Legendre polynomials does not change the linear stability properties of the ADER-DG numerical method with a local DG predictor~\cite{ader_dg_ivp_ode}.

\subsection{Numerical solution at nodes and local solution}

As a result of using the numerical ADER-DG method with a local DG predictor for solving the DAE system (\ref{eq:dae_chosen_form}), two types of numerical solutions are obtained: a solution at the grid nodes $(\mathbf{u}_{n}, \mathbf{v}_{n})$, referred to simply as a ``solution at the nodes'' in the text, and a local solution $(\mathbf{u}_{L}, \mathbf{v}_{L})$, which allows one to determine the numerical solution between the grid nodes. This allows one to use the numerical solution of the ADER-DG method with a local DG predictor in the context of problems that require a solution not only at a discrete set of grid nodes $t_{n}$, but also in the entire domain of definition $\Omega$, for example, for an accurate calculation of functionals of the numerical solution without using any methods for reconstructing the numerical solution between nodes from the solution at the nodes. In the present paper, the accuracy and convergence of the numerical solution obtained by the ADER-DG method with a local DG predictor for the DAE system were investigated separately for the solution at the nodes $(\mathbf{u}_{n}, \mathbf{v}_{n})$ and the local solution $(\mathbf{u}_{L}, \mathbf{v}_{L})$. The solution at the nodes for differential variables was further designated as $\mathbf{u}_{n}$, for algebraic variables --- as $\mathbf{v}_{n}$; the local solution for differential variables was further designated as $\mathbf{u}_{L}$, for algebraic variables --- as $\mathbf{v}_{L}$.

The solution at the nodes $(\mathbf{u}_{n}, \mathbf{v}_{n})$ was defined only for a discrete set $\{t_{n}\}$ of coordinate values representing the grid nodes $t_{n}$. The local solution $(\mathbf{u}_{L}, \mathbf{v}_{L})$ was defined by a piecewise-continuous function in each individual discretization domain $\Omega_{n}$, in the sense that the local function was represented by sets of continuous functions $\{\mathbf{u}_{L, n}\}$ and $\{\mathbf{v}_{L, n}\}$, where $\mathbf{u}_{L, n}: \Omega_{n} \rightarrow \mathbb{R}^{D_{\rm u}}$ and $\mathbf{v}_{L, n}: \Omega_{n} \rightarrow \mathbb{R}^{D_{\rm v}}$. The local solution between nodes $(\mathbf{u}_{L}, \mathbf{v}_{L})$ was defined as the local discrete time solution that is computed in the local DG predictor stage and was calculated from the following expression:
\begin{equation}\label{eq:def_local_sol_t}
\begin{split}
&\mathbf{u}_{L}(t) = \mathbf{q}_{n}\left(\frac{t - t_{n}}{\Delta t_{n}}\right) = 
\sum\limits_{p = 0}^{N} \mathbf{q}_{n, p}\cdot\varphi_{p}\left(\frac{t - t_{n}}{\Delta t_{n}}\right),\\
&\mathbf{v}_{L}(t) = \mathbf{r}_{n}\left(\frac{t - t_{n}}{\Delta t_{n}}\right) = 
\sum\limits_{p = 0}^{N} \mathbf{r}_{n, p}\cdot\varphi_{p}\left(\frac{t - t_{n}}{\Delta t_{n}}\right),
\end{split}
\end{equation}
where $\mathbf{q}_{n} = \mathbf{q}_{n}(\tau)$ and $\mathbf{r}_{n} = \mathbf{r}_{n}(\tau)$ denotes the local solution (\ref{eq:qr_def_exp}) obtained at the discretization domain $\Omega_{n}$, $\mathbf{q}_{n, p}$ and $\mathbf{r}_{n, p}$ denote the expansion coefficients of the local solution $\mathbf{q}_{n} = \mathbf{q}_{n}(\tau)$ and $\mathbf{r}_{n} = \mathbf{r}_{n}(\tau)$ in the DG representation (\ref{eq:qr_def_exp}). Discretization $\xi_{m}\in[0, 1]$, $0 \leqslant m < M$, in units of the local coordinate $\tau$, was introduced by grid sub-nodes located in the space between grid nodes. The values of the basis polynomials $\varphi_{p}$ were calculated in sub-nodes $\xi_{m}$ of discretization, which made it possible to calculate the tabulated local solution $\mathbf{u}_{L}(t_{n, m})$ in the form of multiplication the polynomials values matrix, with sizes $M \times (N+1)$, on the matrix of coefficients $\mathbf{q}_{p}$, with sizes $(N+1) \times K$.

\begin{figure}[h!]
\captionsetup[subfigure]{%
	position=bottom,
	font+=smaller,
	textfont=normalfont,
	singlelinecheck=off,
	justification=raggedright
}
\centering
\begin{subfigure}{0.320\textwidth}
\includegraphics[width=\textwidth]{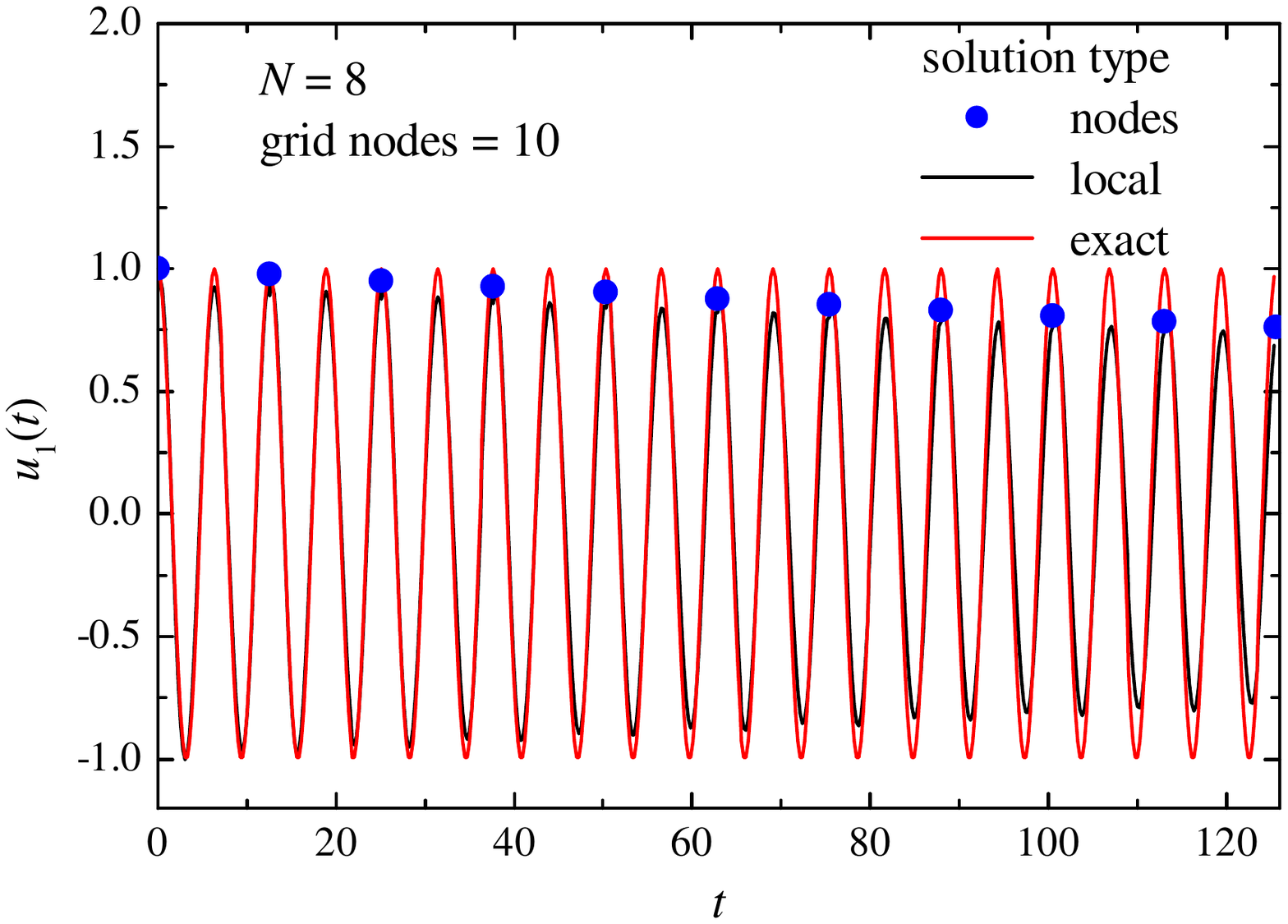}
\vspace{-8mm}\caption{\label{fig:demo_sols_uv_eps:a1}}
\end{subfigure}
\begin{subfigure}{0.320\textwidth}
\includegraphics[width=\textwidth]{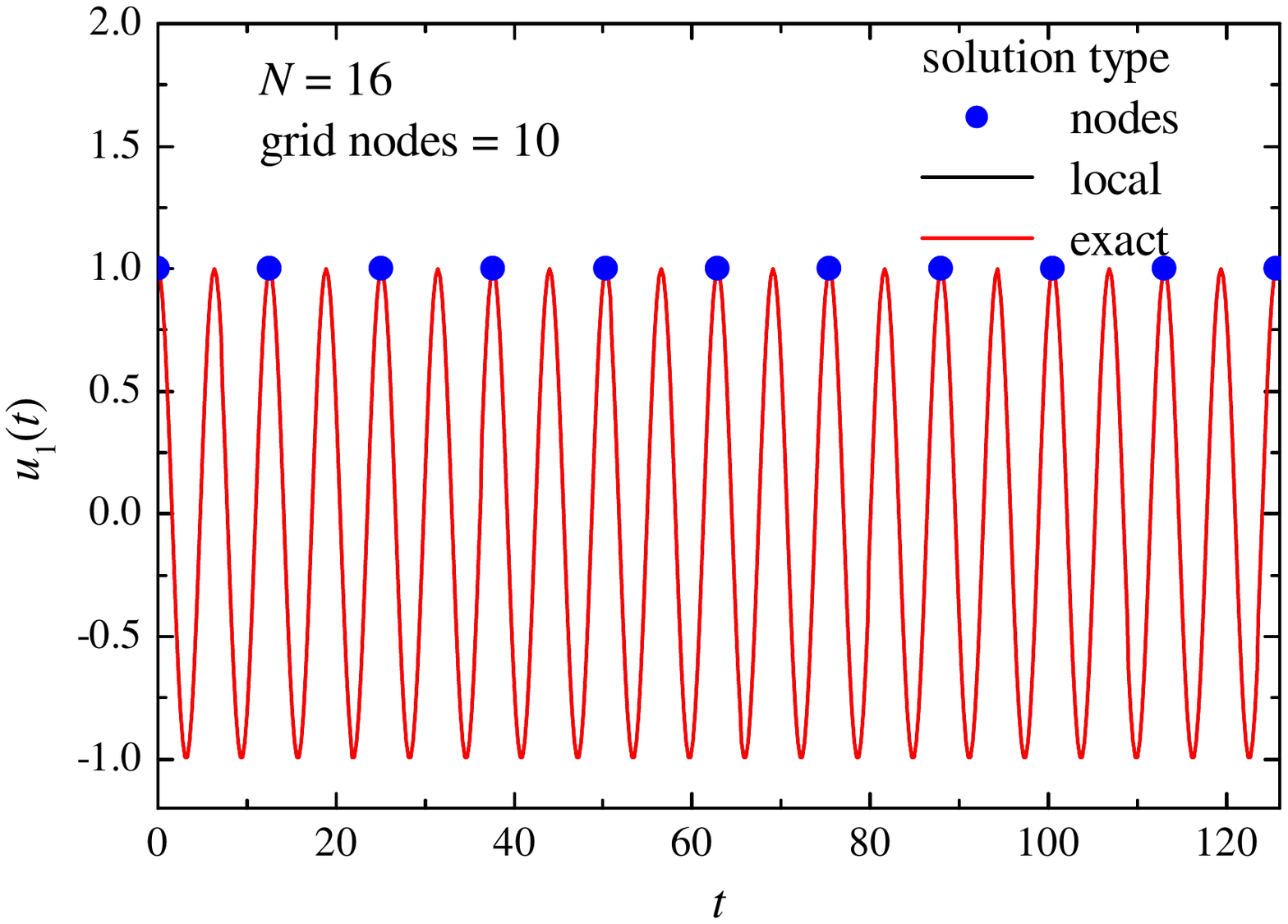}
\vspace{-8mm}\caption{\label{fig:demo_sols_uv_eps:a2}}
\end{subfigure}
\begin{subfigure}{0.320\textwidth}
\includegraphics[width=\textwidth]{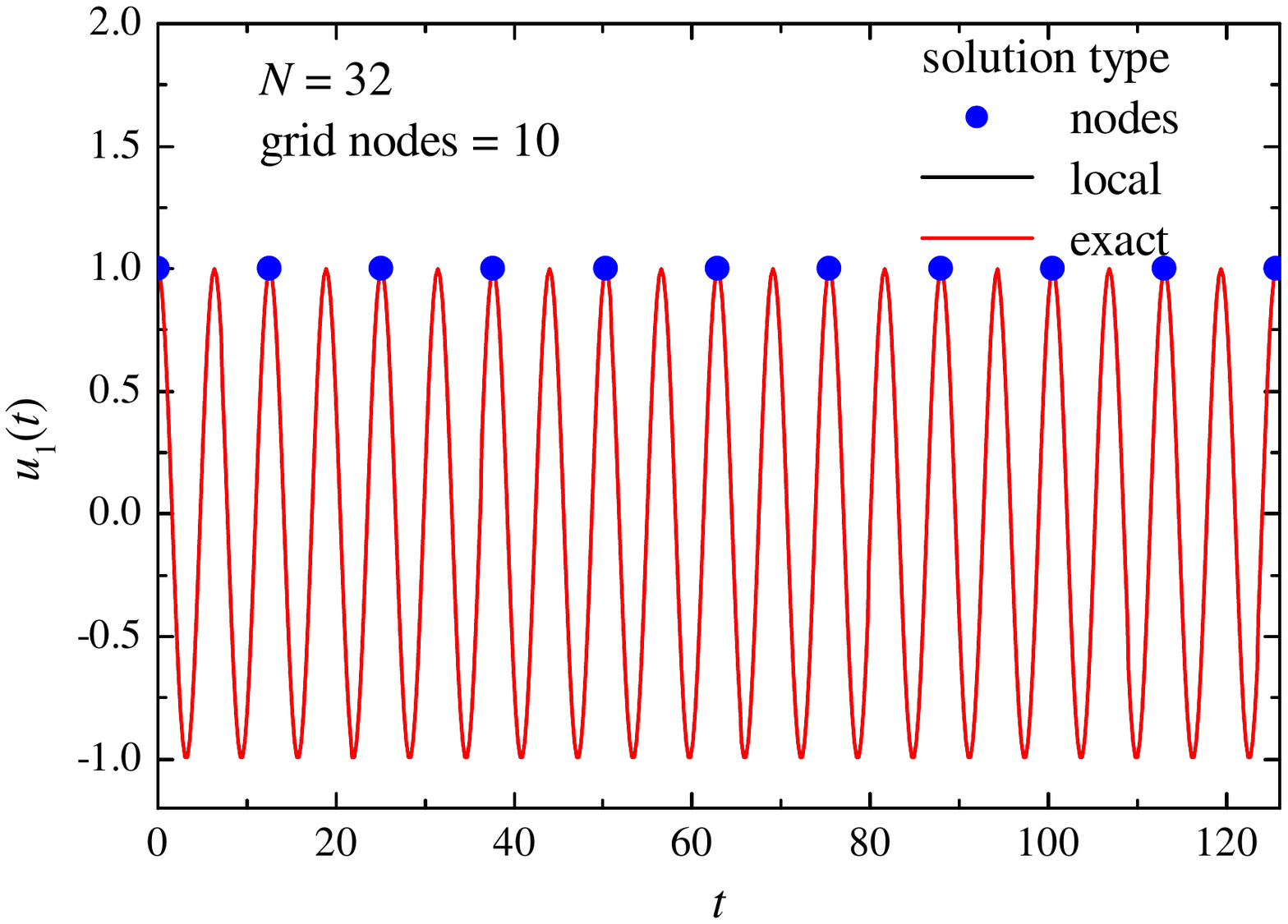}
\vspace{-8mm}\caption{\label{fig:demo_sols_uv_eps:a3}}
\end{subfigure}\\[-2mm]
\begin{subfigure}{0.320\textwidth}
\includegraphics[width=\textwidth]{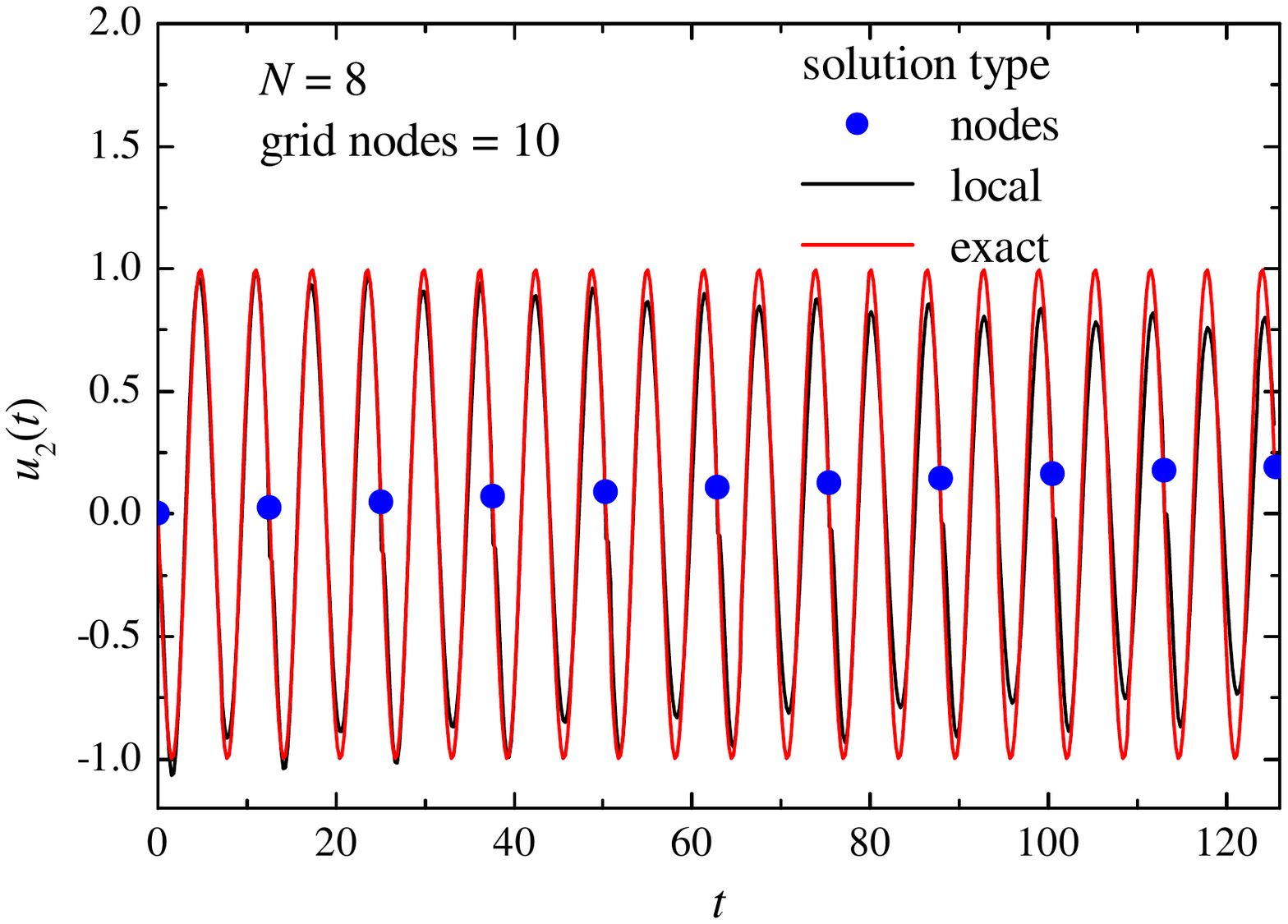}
\vspace{-8mm}\caption{\label{fig:demo_sols_uv_eps:b1}}
\end{subfigure}
\begin{subfigure}{0.320\textwidth}
\includegraphics[width=\textwidth]{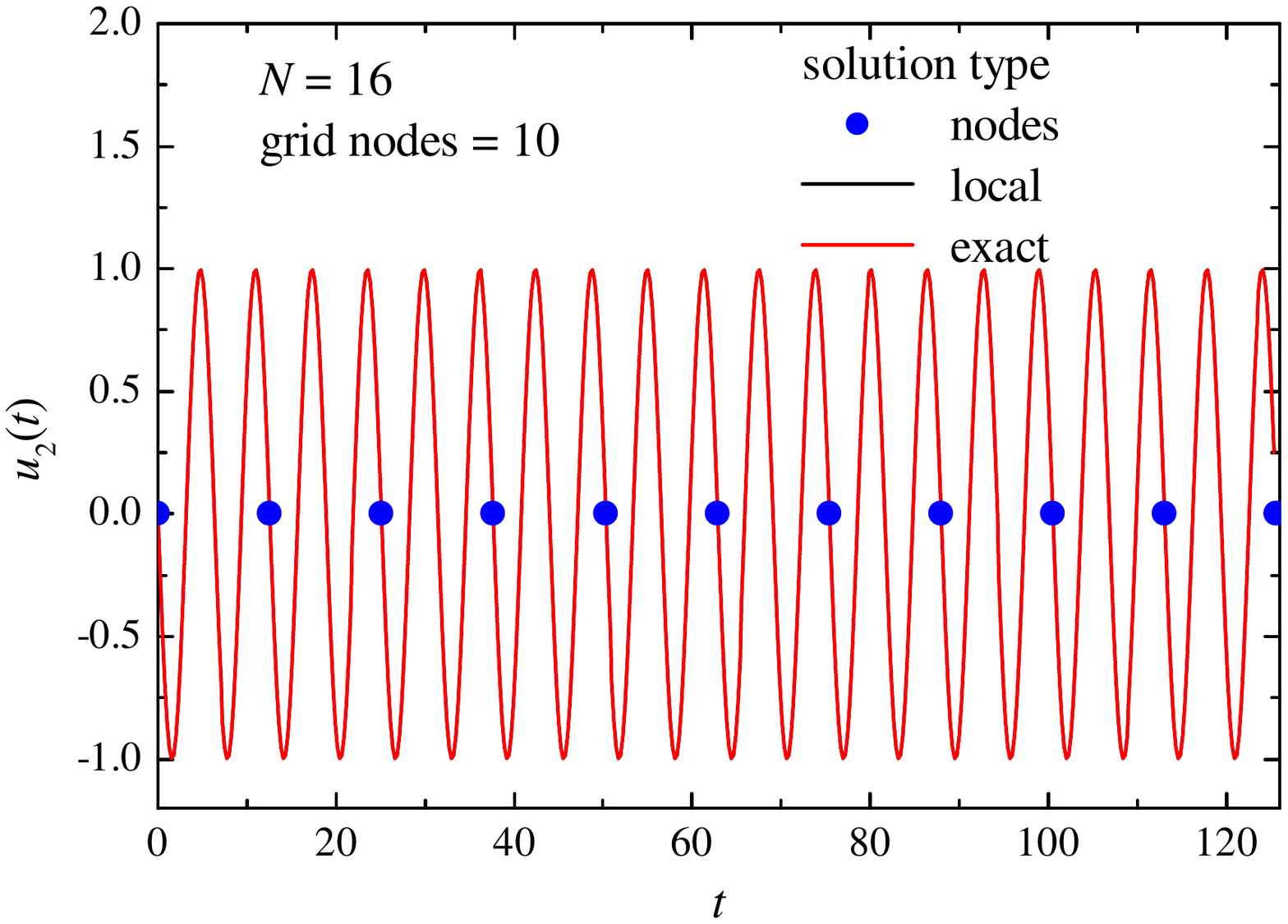}
\vspace{-8mm}\caption{\label{fig:demo_sols_uv_eps:b2}}
\end{subfigure}
\begin{subfigure}{0.320\textwidth}
\includegraphics[width=\textwidth]{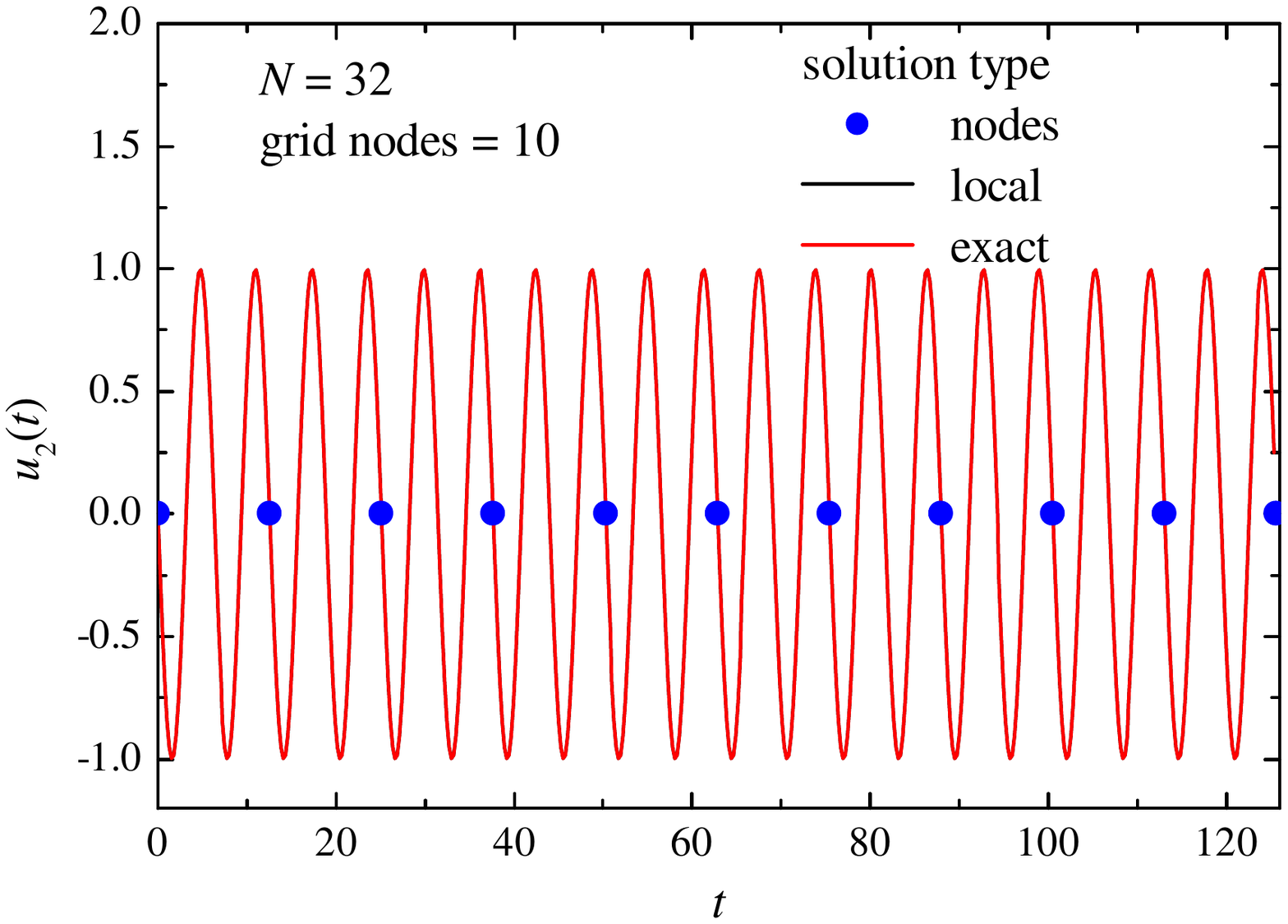}
\vspace{-8mm}\caption{\label{fig:demo_sols_uv_eps:b3}}
\end{subfigure}\\[-2mm]
\begin{subfigure}{0.320\textwidth}
\includegraphics[width=\textwidth]{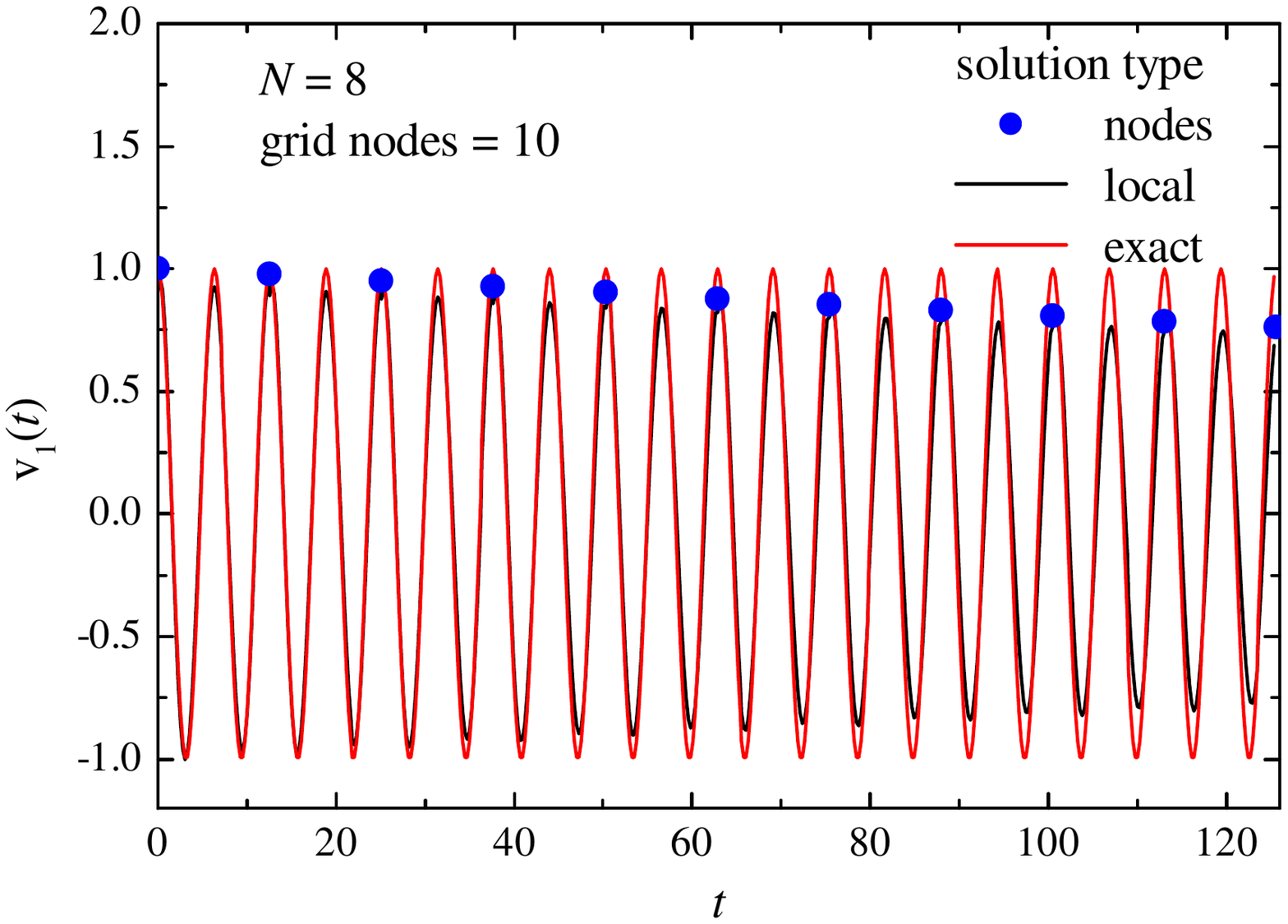}
\vspace{-8mm}\caption{\label{fig:demo_sols_uv_eps:c1}}
\end{subfigure}
\begin{subfigure}{0.320\textwidth}
\includegraphics[width=\textwidth]{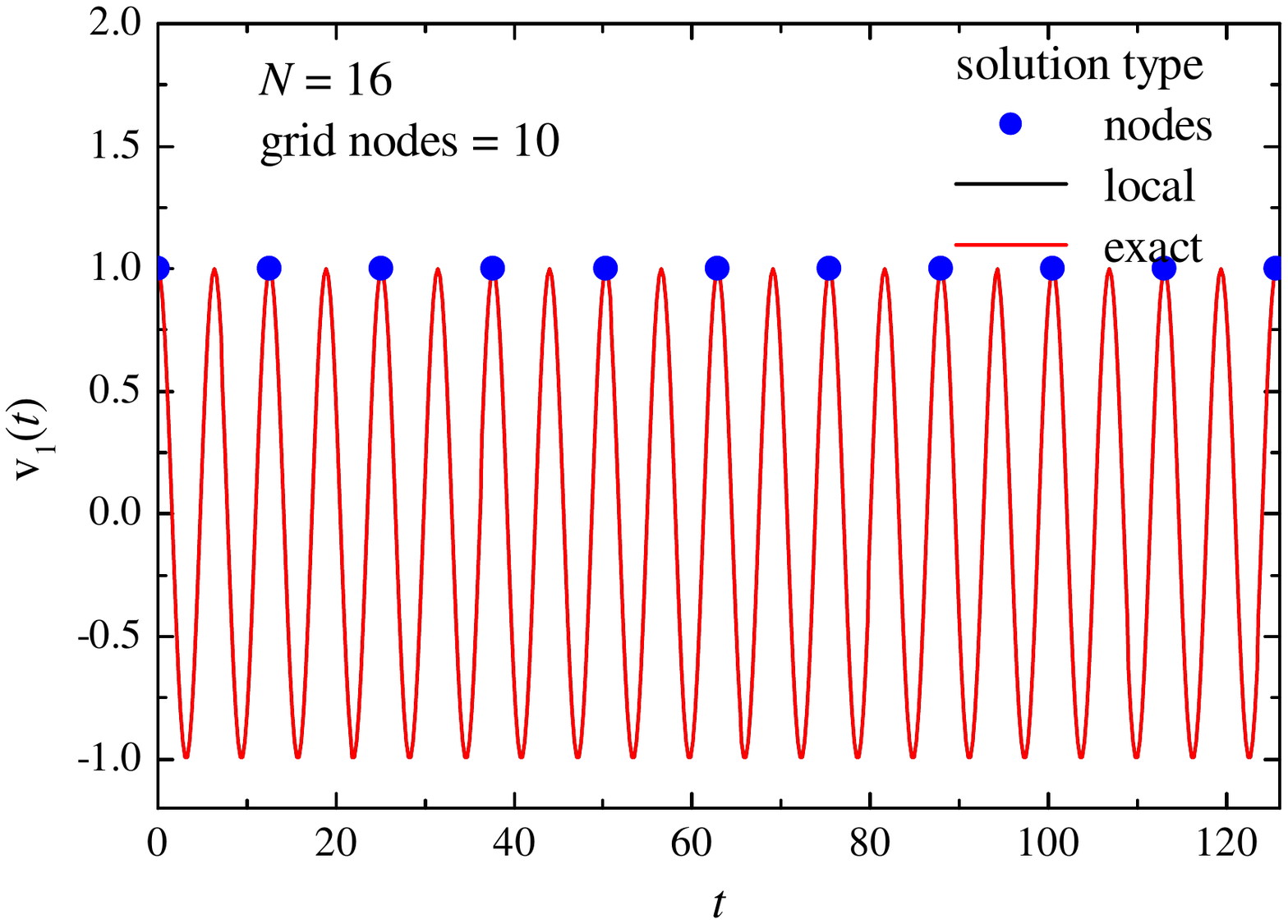}
\vspace{-8mm}\caption{\label{fig:demo_sols_uv_eps:c2}}
\end{subfigure}
\begin{subfigure}{0.320\textwidth}
\includegraphics[width=\textwidth]{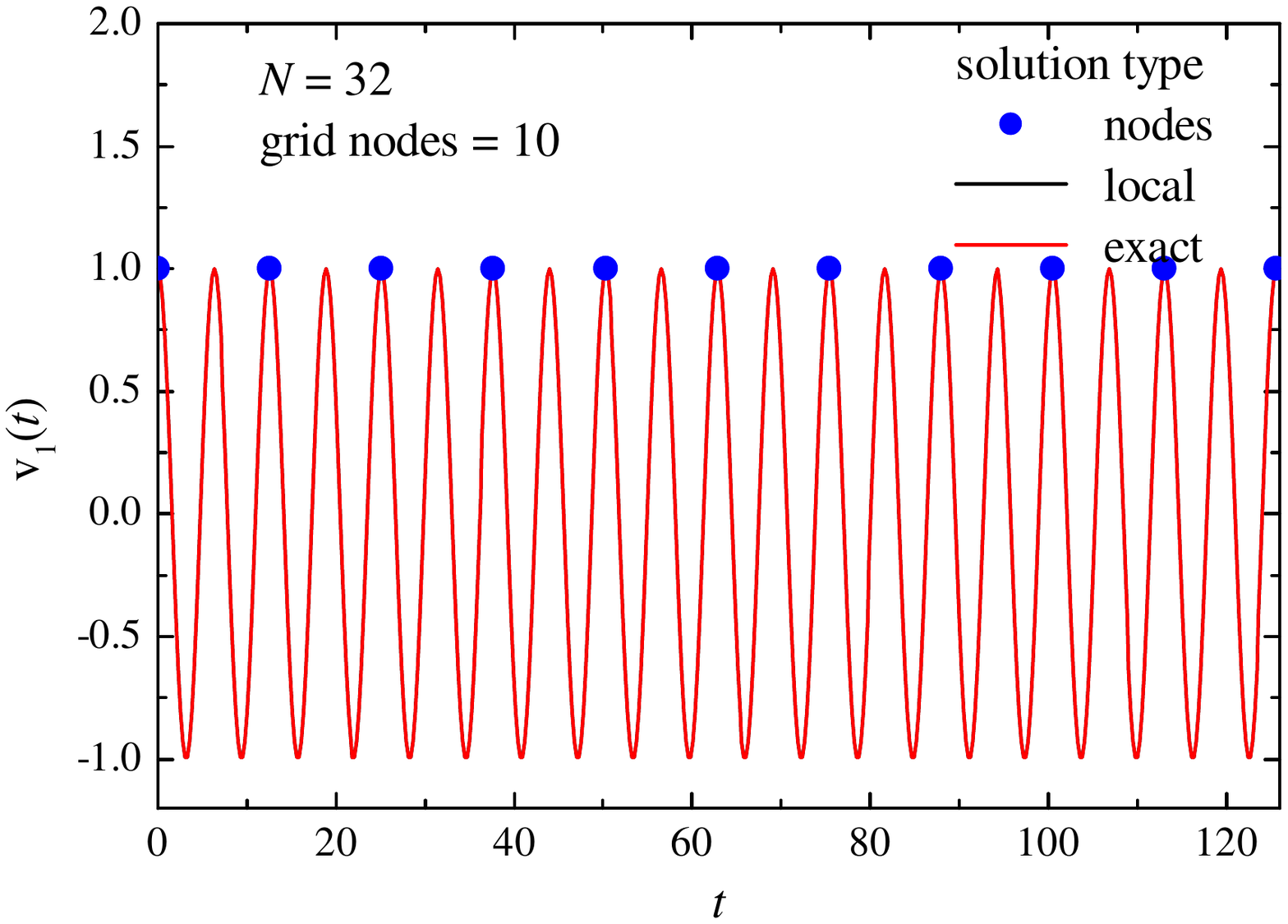}
\vspace{-8mm}\caption{\label{fig:demo_sols_uv_eps:c3}}
\end{subfigure}\\[-2mm]
\begin{subfigure}{0.320\textwidth}
\includegraphics[width=\textwidth]{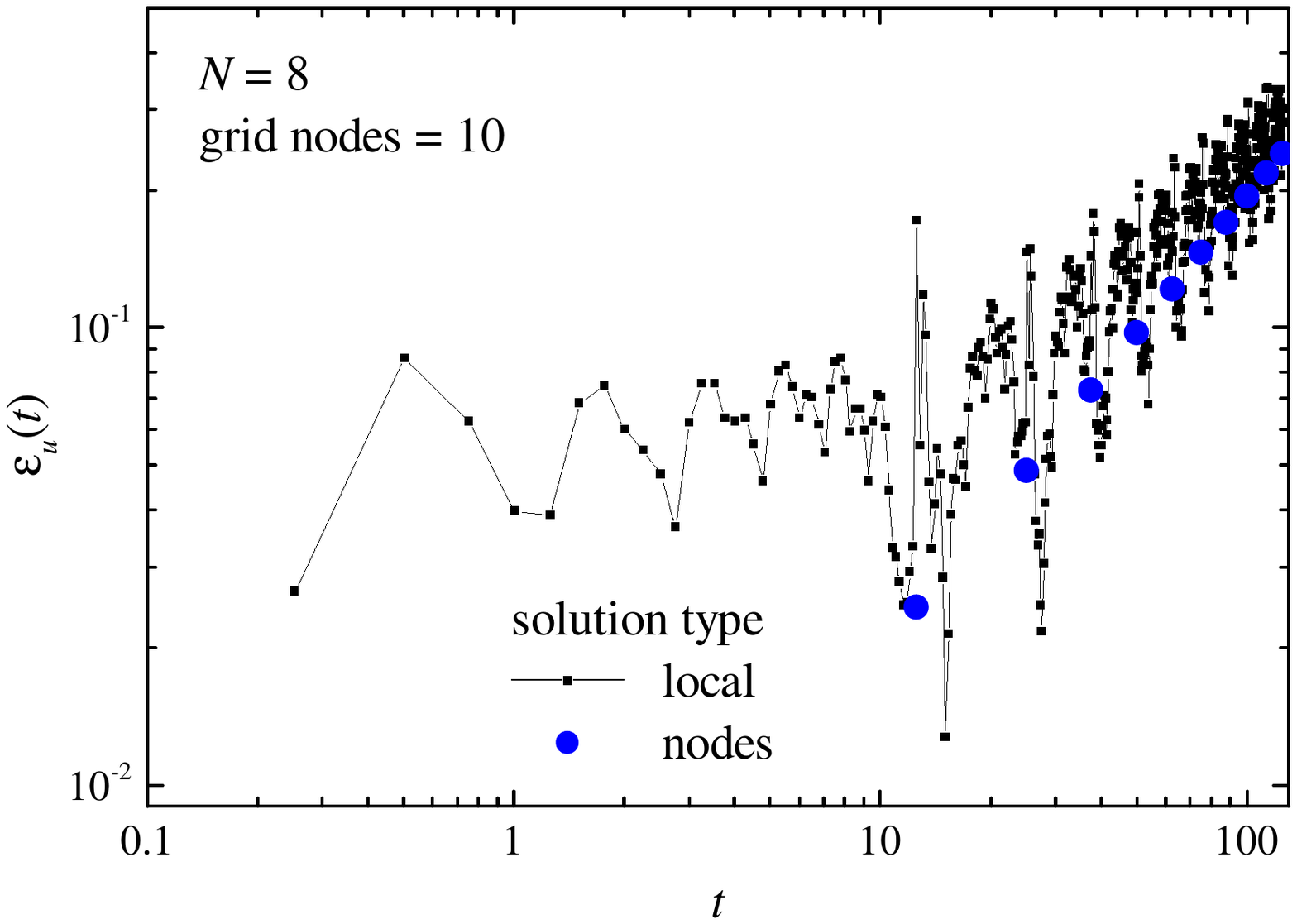}
\vspace{-8mm}\caption{\label{fig:demo_sols_uv_eps:d1}}
\end{subfigure}
\begin{subfigure}{0.320\textwidth}
\includegraphics[width=\textwidth]{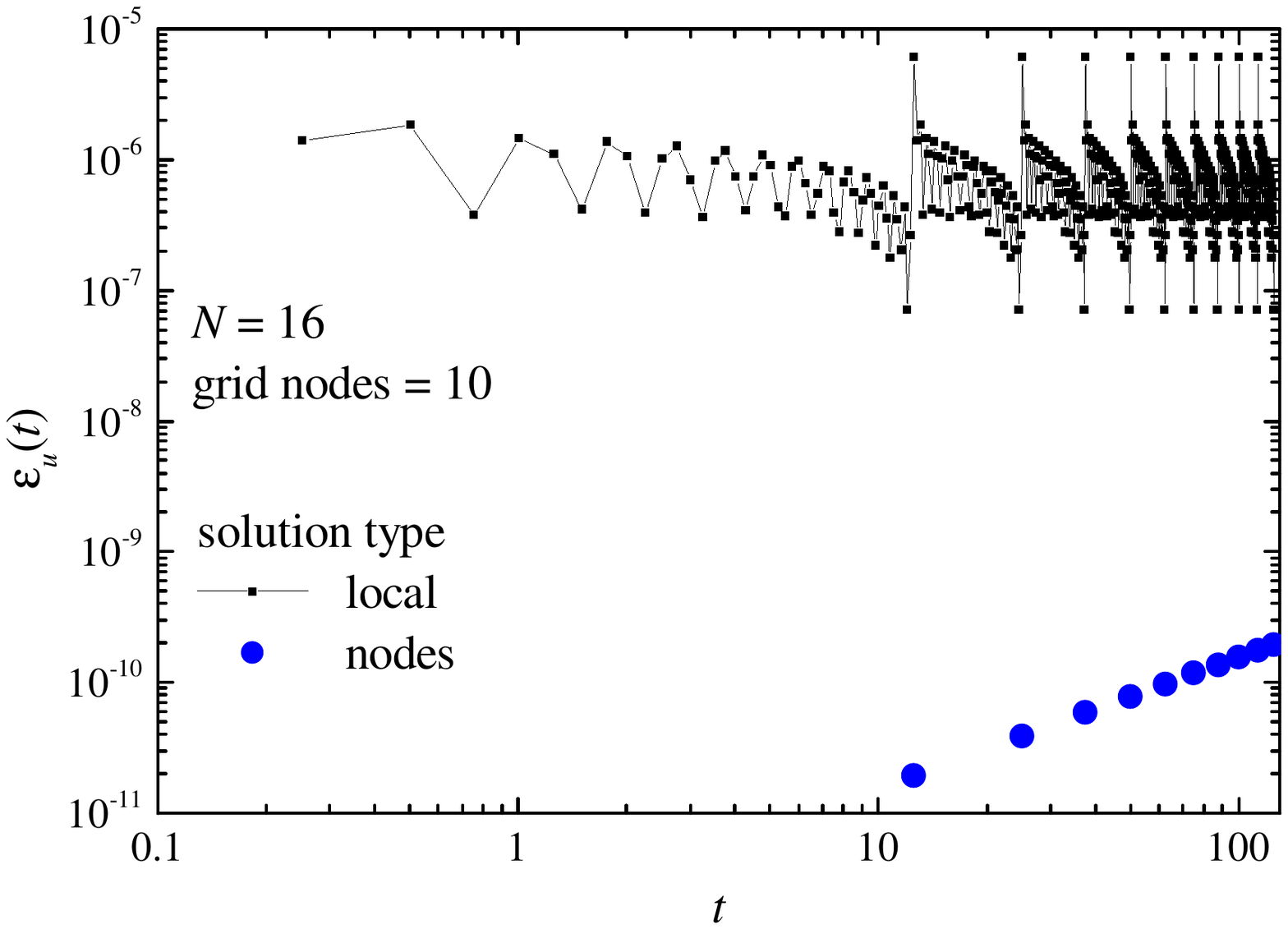}
\vspace{-8mm}\caption{\label{fig:demo_sols_uv_eps:d2}}
\end{subfigure}
\begin{subfigure}{0.320\textwidth}
\includegraphics[width=\textwidth]{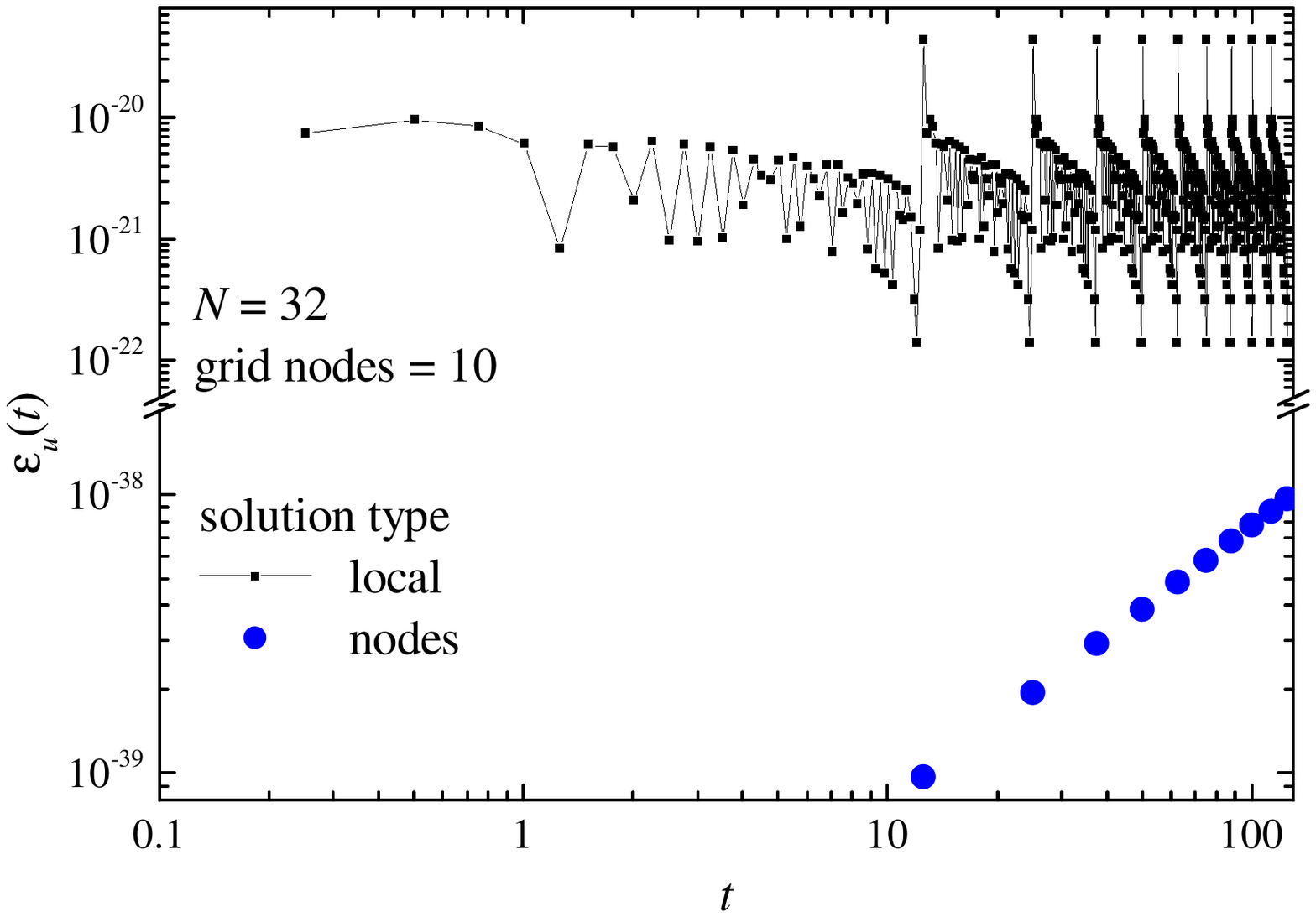}
\vspace{-8mm}\caption{\label{fig:demo_sols_uv_eps:d3}}
\end{subfigure}\\[-2mm]
\begin{subfigure}{0.320\textwidth}
\includegraphics[width=\textwidth]{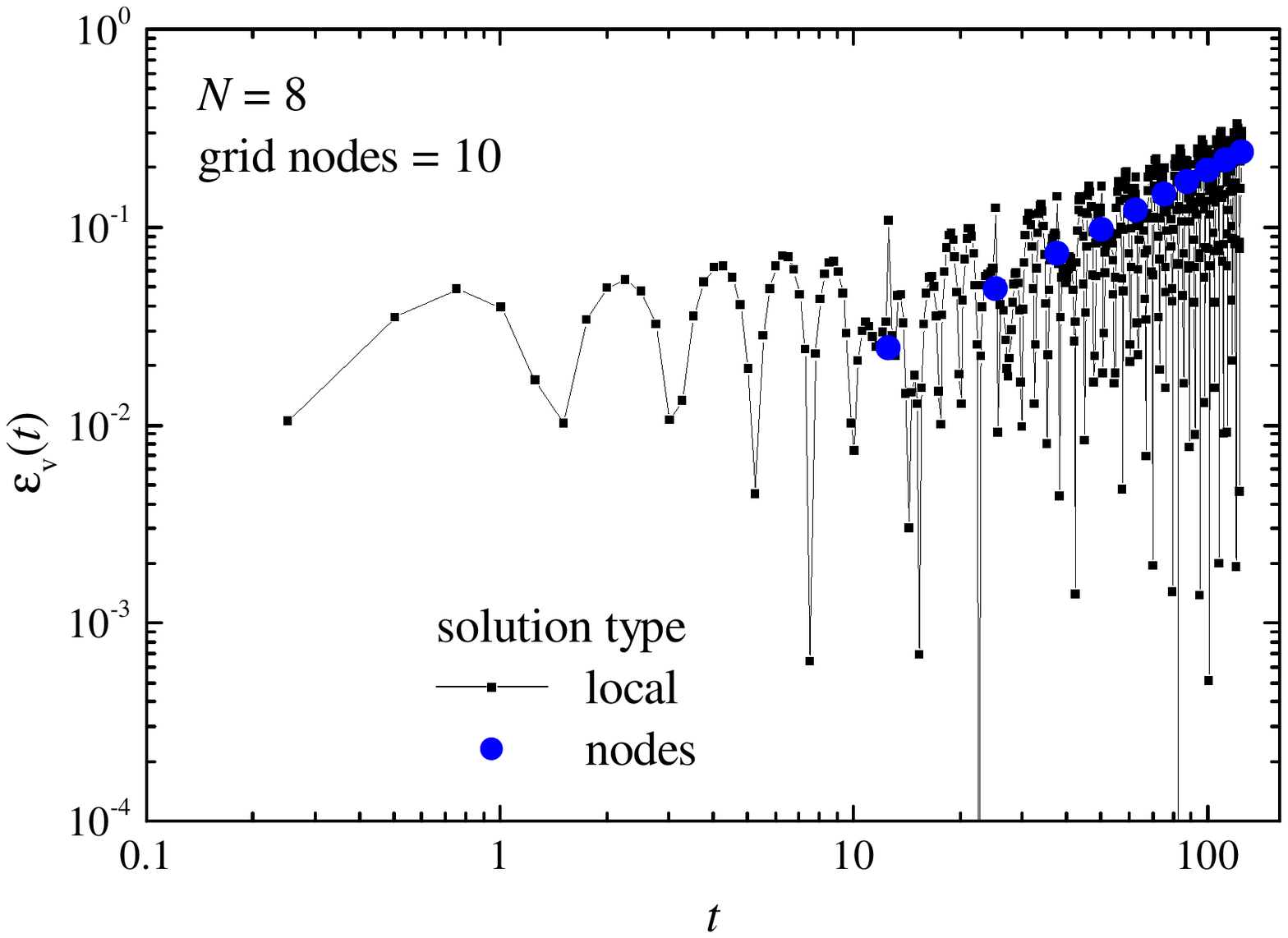}
\vspace{-8mm}\caption{\label{fig:demo_sols_uv_eps:e1}}
\end{subfigure}
\begin{subfigure}{0.320\textwidth}
\includegraphics[width=\textwidth]{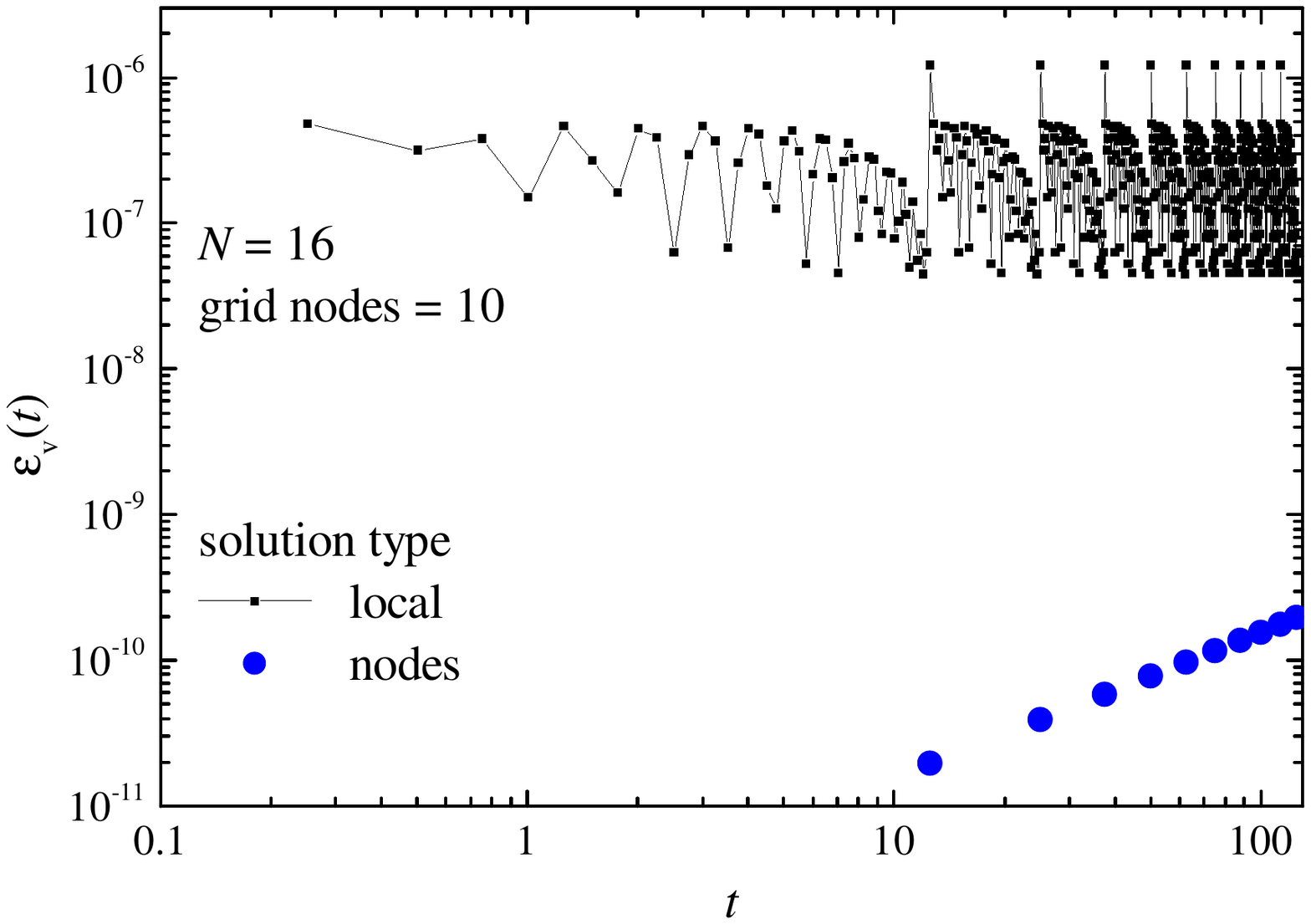}
\vspace{-8mm}\caption{\label{fig:demo_sols_uv_eps:e2}}
\end{subfigure}
\begin{subfigure}{0.320\textwidth}
\includegraphics[width=\textwidth]{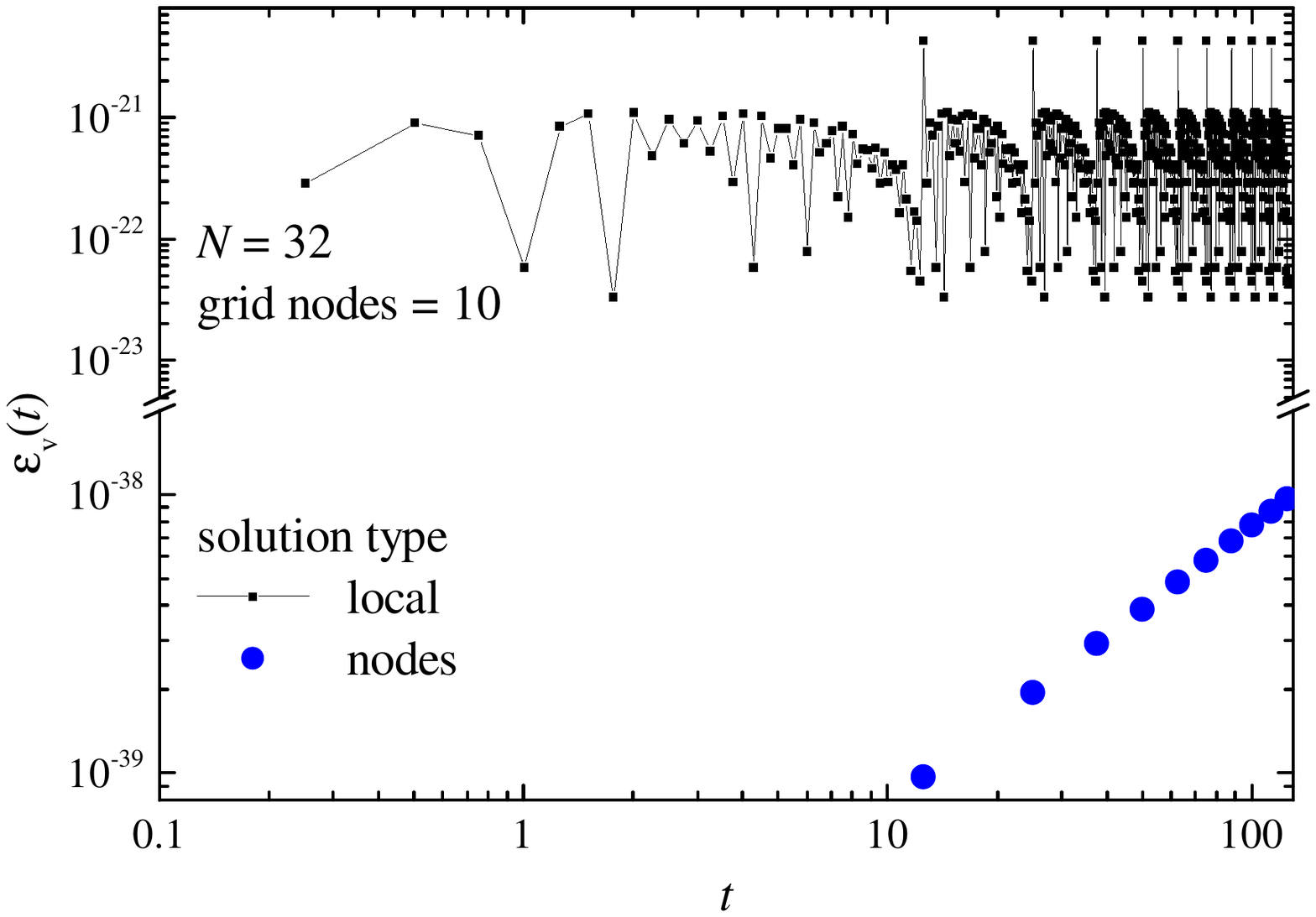}
\vspace{-8mm}\caption{\label{fig:demo_sols_uv_eps:e3}}
\end{subfigure}\\[-2mm]z
\caption{%
Numerical solution of the problem (\ref{eq:demo_test}), in the domain $0 \leqslant t \leqslant 40\pi$ with step $\mathrm{\Delta}t = 4\pi$. Comparison of the solution at nodes $\mathbf{u}_{n}$, the local solution $\mathbf{u}_{L}(t)$ and the exact solution $\mathbf{u}^{\rm ex}(t)$ for components $u_{1}$ (\subref{fig:demo_sols_uv_eps:a1}, \subref{fig:demo_sols_uv_eps:a2}, \subref{fig:demo_sols_uv_eps:a3}), $u_{2}$ (\subref{fig:demo_sols_uv_eps:b1}, \subref{fig:demo_sols_uv_eps:b2}, \subref{fig:demo_sols_uv_eps:b3}), $v_{1}$ (\subref{fig:demo_sols_uv_eps:c1}, \subref{fig:demo_sols_uv_eps:c2}, \subref{fig:demo_sols_uv_eps:c3}), the dependencies of the point-wise errors $\varepsilon_{u}(t)$ (\subref{fig:demo_sols_uv_eps:d1}, \subref{fig:demo_sols_uv_eps:d2}, \subref{fig:demo_sols_uv_eps:d3}), $\varepsilon_{v}(t)$ (\subref{fig:demo_sols_uv_eps:e1}, \subref{fig:demo_sols_uv_eps:e2}, \subref{fig:demo_sols_uv_eps:e3}), obtained using polynomials with degrees $N = 8$ (\subref{fig:demo_sols_uv_eps:a1}, \subref{fig:demo_sols_uv_eps:b1}, \subref{fig:demo_sols_uv_eps:c1}, \subref{fig:demo_sols_uv_eps:d1}, \subref{fig:demo_sols_uv_eps:e1}), $N = 16$ (\subref{fig:demo_sols_uv_eps:a2}, \subref{fig:demo_sols_uv_eps:b2}, \subref{fig:demo_sols_uv_eps:c2}, \subref{fig:demo_sols_uv_eps:d2}, \subref{fig:demo_sols_uv_eps:e2}) and $N = 32$ (\subref{fig:demo_sols_uv_eps:a3}, \subref{fig:demo_sols_uv_eps:b3}, \subref{fig:demo_sols_uv_eps:c3}, \subref{fig:demo_sols_uv_eps:d3}, \subref{fig:demo_sols_uv_eps:e3}).
}
\label{fig:demo_sols_uv_eps}
\end{figure}

An important feature of the ADER-DG numerical method with a local DG predictor is the possibility of using the local solution $\mathbf{u}_{L}(t)$ as the final solution of the problem even in the case of very coarse grids with a large discretization step $\Delta t$. As a demonstration example a numerical solution of the simple DAE system of index 1 for a linear one-dimensional harmonic oscillator $\ddot{x} + x = 0$, $x(0) = 1$, $\dot{x}(0) = 0$:
\begin{equation}\label{eq:demo_test}
\begin{split}
&\dfrac{du_{1}}{dt} = u_{2},\quad \dfrac{du_{2}}{dt} = -v_{1},\quad u_{1} - v_{1} = 0,\\ 
&u_{1}(0) = 1,\quad u_{2}(0) = 0,\quad v_{1}(0) = 1,
\end{split}
\end{equation}
with an analytical exact solution $\mathbf{u}^{\rm ex} = [\cos(t),\, -\sin(t)]^{T}$ and $\mathbf{v}^{\rm ex} = [\cos(t)]$ is presented in Fig.~\ref{fig:demo_sols_uv_eps}. The solution was obtained in the domain $0 \leqslant t \leqslant 40\pi$ with a discretization step $\Delta t = 4\pi$ equal to two harmonic oscillation period --- on the entire domain, the ADER-DG numerical method makes only $10$ steps. To obtain a tabular local solution in all spaces between grid nodes uniform in the local coordinate $\tau$, grids with $M = 50$ sub-nodes were additionally included.

The presented results allow us to conclude that in the case of degree $N = 8$, for such a coarse grid, non-physical dissipation of the numerical solution begins to manifest itself significantly, which is associated with a very high stability of the numerical method. In the case of degree $N = 16$, for the selected values of the discretization step and the domain of definition of the solution, dissipative errors are practically not observed. The point-wise errors $\varepsilon_{u}(t) = |\mathbf{u}(t) - \mathbf{u}^{\rm ex}(t)|$ and $\varepsilon_{v}(t) = |\mathbf{v}(t) - \mathbf{v}^{\rm ex}(t)|$ of the numerical solution in the cases $N = 16$ and $32$ are so small that they cannot be visually observed on the solution plots. The errors $\varepsilon_{u}(t)$ and $\varepsilon_{v}(t)$ of the local solution are $10^{-7}$--$10^{-5}$ for degree $N = 16$ and $10^{-22}$--$10^{-19}$ for degree $N = 32$, which in the case $N = 32$ reach the round-off error of representing real numbers in the double-precision floating-point format. The error of the solution at the nodes is $10^{-11}$--$10^{-10}$ for degree $N = 16$ and $10^{-40}$--$10^{-38}$ for degree $N = 32$. For the point-wise errors $\varepsilon_{u}(t)$ and $\varepsilon_{v}(t)$ of the solution at the nodes $\mathbf{u}_{n}$ and $\mathbf{v}_{n}$, an approximately linear increase is observed by the argument $t$.

This is a demonstration example, and it shows that the real accuracy of the numerical solution obtained using the local DG predictor within the ADER-DG numerical method can already reach the highest applied calculation accuracy.

\subsection{Accuracy, convergence and norms}

A quantitative study of the accuracy and convergence of the numerical method ADER-DG with a local DG predictor for solving DAE systems was carried out on a set of test examples, among which DAE systems of indices 1, 2 and 3 were selected, for which exact analytical or reference solutions exist. In all cases of DAE systems with an index higher than 1, systems were studied whose index was lowered down to 1. Among these examples, such classical tests as the mathematical pendulum and double pendulum were selected. A test example with a stiff problem was studied, including in the area of problem parameters where extremely high stiffness occurs. The results of the study are presented in the following Section~\ref{sec:2}.

This Subsection presents the mathematical apparatus of the norms used in this work and methods for their calculation, on the basis of which the errors of the numerical solution and the convergence orders of the numerical method were calculated. Accuracy and convergence were separately investigated for the solution at the nodes $(\mathbf{u}_{n}, \mathbf{v}_{n})$ and the local solution $(\mathbf{u}_{L}, \mathbf{v}_{L})$ in the space between nodes $[t_{n},\, t_{n+1}]$. The solution in nodes represents the values $(\mathbf{u}_{L}, \mathbf{v}_{L})$ of the desired functions $(\mathbf{u}, \mathbf{v})$ at individual grid points $t_{n}$. The local solution $\mathbf{u}_{L}(t) = \mathbf{q}_{n}(\tau(t))$ and $\mathbf{v}_{L}(t) = \mathbf{r}_{n}(\tau(t))$ are continuous functions (see formula (\ref{eq:def_local_sol_t})) of local coordinate $\tau$ determined by the expansion in basis polynomials $\varphi_{p}(\tau)$, the coefficients $(\hat{\mathbf{q}}_{n, p}, \hat{\mathbf{r}}_{n, p})$ of which are determined for each separate space $[t_{n},\, t_{n+1}]$ between adjacent grid nodes. In this Subsection, the local solution $(\mathbf{q}_{n}, \mathbf{r}_{n})$ will also be matched with the index $n$ of the grid node, which is the left grid node $t_{n}$ for the discretization domain $\Omega_{n} = \{t\, |\, t \in [t_{n},\ t_{n+1}]\}$ in which the local solution is defined.

The accuracy of the numerical solution was determined point-wise, based on the values of the local error of the numerical solution:
\begin{equation}\label{eq:local_uv_eps_def}
\begin{split}
&\varepsilon_{u}(t) = \max\limits_{1 \leqslant k \leqslant D_{\rm u}}|u_{k}(t) - u^{\rm ex}_{k}(t)|,\\
&\varepsilon_{v}(t) = \max\limits_{1 \leqslant k \leqslant D_{\rm v}}|v_{k}(t) - v^{\rm ex}_{k}(t)|,
\end{split}
\end{equation}
where $u_{k}(t)$, $v_{k}(t)$ and $u^{\rm ex}_{k}(t)$, $v^{\rm ex}_{k}(t)$ are the $k$-th components of the vector functions of the numerical solution $\mathbf{u}(t)$, $\mathbf{v}(t)$ and the exact (or reference) solutions $\mathbf{u}^{\rm ex}(t)$, $\mathbf{v}^{\rm ex}(t)$ at a point $t$, respectively. Among the results presented below are graphs of the dependence of the error on $t$. The errors $(\varepsilon_{u}(t_{n}), \varepsilon_{v}(t_{n}))$ of the numerical solution in nodes $(\mathbf{u}, \mathbf{v})$ will be presented only in the nodes $t_{n}$ of the grid. The errors $(\varepsilon_{u}(t_{n, m}), \varepsilon_{v}(t_{n, m}))$ of the local solution $(\mathbf{u}_{L}, \mathbf{v}_{L})$ will be represented in sub-nodes $t_{n, m}$ located between the nodes $t_{n}$ of the grid. In this paper, $M = 50$ points were uniformly added in the spaces $\Omega_{n}$ between the nodes, at which estimates of the local numerical solution were made. Similarly, an error in the execution of algebraic constraints included in the original DAE system (\ref{eq:dae_chosen_form}) was introduced:
\begin{equation}\label{eq:local_g_eps_def}
\begin{split}
&\varepsilon_{g}(t) = \max\limits_{1 \leqslant k \leqslant D_{\rm v}}|g_{k}(t)|,
\end{split}
\end{equation}
where $g_{k}(t)$  are the $k$-th components of the vector functions $\mathbf{g} = \mathbf{g}(\mathbf{u},\, \mathbf{v},\, t)$. Unlike numerical solutions $(\mathbf{u}, \mathbf{v})$ and $(\mathbf{u}_{L}, \mathbf{v}_{L})$, in the case of algebraic constraints the exact solution is $\mathbf{0}$.

The calculation of the convergence orders $p$ of the numerical solution was carried out on the basis of the analysis of the global error $e$ of the numerical solution, integrally defined for the entire domain $\Omega$ of the definition of the desired functions $\mathbf{u}$ and $\mathbf{v}$. The classical functional norms that were used to calculate global errors were defined by the following functionals:
\begin{equation}\label{eq:norms_def}
\begin{split}
&||\boldsymbol{\psi}||_{L_{1}} = \int\limits_{\Omega} \left|\boldsymbol{\psi}(t)\right| dt,\hspace{10mm}
 ||\boldsymbol{\psi}||_{L_{2}}^{2} = \int\limits_{\Omega} \left|\boldsymbol{\psi}(t)\right|^{2} dt,\\
&||\boldsymbol{\psi}||_{L_{\infty}} = \operatorname{ess}\sup\limits_{\hspace{-5mm}t\in\Omega}  \left|\boldsymbol{\psi}(t)\right|,\hspace{5.2mm}
 ||\boldsymbol{\psi}||_{f} = \left|\boldsymbol{\psi}(t_{f})\right|,
\end{split}
\end{equation}
where $\boldsymbol{\psi}: \Omega \rightarrow \mathbb{R}^{D_{\psi}}$ is a function, and under the symbols $|\ldots|$ is meant the maximum absolute value taken over all components of the vector function: $|\boldsymbol{\psi}| = \max_{k} |\psi_{k}|$, as in (\ref{eq:local_uv_eps_def}). 

Using these functional norms (\ref{eq:norms_def}), the following types of global errors of the numerical solution were determined:
\begin{equation}\label{eq:global_eps_def}
\begin{split}
&e_{L}^{u} = ||\mathbf{u}_{L} - \mathbf{u}^{\rm ex}||,\quad 
 e_{L}^{v} = ||\mathbf{v}_{L} - \mathbf{v}^{\rm ex}||,\quad 
 e_{L}^{g} = ||\mathbf{g}(\mathbf{u}_{L}(.),\, \mathbf{v}_{L}(.),\, .)||,\\
&e_{G}^{u} = ||\mathbf{u}_{n} - \mathbf{u}^{\rm ex}||,\quad 
 e_{G}^{v} = ||\mathbf{v}_{n} - \mathbf{v}^{\rm ex}||,\quad 
 e_{G}^{g} = ||\mathbf{g}(\mathbf{u}_{n}(.),\, \mathbf{v}_{n}(.),\, .)||,
\end{split}
\end{equation}
where index $L$ means that this error $e$ is calculated for a local solution, index $G$ means that this error $e$ is calculated for the solution at the nodes, the dot $.$ means the substitution of the argument $t$. The calculation of global errors (\ref{eq:global_eps_def}) for the solution at the nodes $(\mathbf{u}_{n}, \mathbf{v}_{n})$ was carried out on the basis of the replacement of integrals by finite sums multiplied by discretization steps; the operation $\operatorname{ess}\sup$ by $t\in\Omega$ has been replaced by the operation $\max$ by nodes $t_{n}$. The calculation of global errors (\ref{eq:global_eps_def}) for the local solution $(\mathbf{u}_{L}(t), \mathbf{v}_{L}(t))$ in the space between nodes $\Omega_{n}$ was carried out on the basis of the replacement of integrals by finite sums over sub-nodes $t_{n, m}$, multiplied by the discretization steps between sub-nodes $\Delta t_{n}/M$; the operation $\sup$ by $t\in\Omega$ has been replaced by the operation $\max$ by sub-nodes $t_{n, m}$. The ``final'' norm $||\ldots||_{f}$ was calculated only for the solution in the nodes $\mathbf{u}_{n}$ and $\mathbf{v}_{n}$ in the last node of the grid, which corresponds to the coordinate $t_{f}$.

The convergence orders $p$ were calculated based on the approximation of the dependence of the global error $e$ of the numerical solution on the discretization step $\Delta t$ by the power dependence $e(\Delta t) \sim (\Delta t)^{p}$. The calculation of the convergence order was carried out on the basis of the analysis of the numerical solution, performed with a constant discretization step $\Delta t$. The convergence orders $p$ in this work were calculated separately for the solution at the nodes and the local solution in the space between the nodes, separately for the functional norms (\ref{eq:norms_def}): $p_{L_{1}}$, $p_{L_{2}}$, $p_{L_{\infty}}$, $p_{f}$.

The ``final'' norm $||\ldots||_{f}$ (\ref{eq:norms_def}) was used in this work to calculate the global error $e$ and empirical convergence orders $p_{f}$ in the case of test examples of DAE systems (\ref{eq:dae_chosen_form}) for which there is no exact analytical solution in the entire domain of definition $\Omega$ of the desired functions $\mathbf{u}$ and $\mathbf{v}$, however, there is a high accuracy reference solution in the final coordinate $t_{f}$ obtained by a numerical method of incomparably higher accuracy and order than the method under study. In cases where for test examples for the DAE systems (\ref{eq:dae_chosen_form}) there was an exact analytical solution, even expressed in special functions, classical norms $||\ldots||_{L_{1}}$, $||\ldots||_{L_{2}}$, $||\ldots||_{L_{\infty}}$ (\ref{eq:norms_def}) were used to calculate the global error $e$ for calculating the empirical convergence orders $p_{L_{1}}$, $p_{L_{2}}$, $p_{L_{\infty}}$.

According to the results of the works~\cite{dg_ivp_ode_1, dg_ivp_ode_3, ader_dg_ivp_ode}, it is expected that the ADER-DG numerical method with a local DG predictor with degrees of polynomials $N$ considered in this work can demonstrate superconvergence with the convergence order $p_{G}$ for the solution at grid nodes $\mathbf{u}_{n}$ and the convergence order $p_{L}$ for the local solution $\mathbf{u}_{L}(t)$:
\begin{equation}\label{eq:expect_orders}
p_{\rm nodes} = 2N+1;\quad
p_{\rm local} = N+1;
\end{equation}
which is typical for DG methods for solving the initial value problem for ODE systems. It should be noted that these values of the convergence orders are valid only for initial value problem for the ODE systems. It is expected that in the case of solving DAE systems, especially high index (2 and more), the convergence orders of the numerical method will be lower than the convergence orders achieved for the initial value problem for ODE systems.

One important note should be made regarding the expected convergence orders. Convergence orders $p_{\rm local} = N+1$ for local solution $(\mathbf{u}_{L}(t), \mathbf{v}_{L}(t))$ are quite expected in the case of using polynomials $\varphi_{p}$ of degree $N$ in the representation of the local solution. Convergence orders $p_{\rm nodes} = 2N+1$ for the solution at grid nodes for differential variables $\mathbf{u}_{n}$, which was calculated using formula (\ref{eq:u_sol_in_node}), were expected due to the use of the quadrature formula (\ref{eq:gr_rule}), as was already shown in the work~\cite{ader_dg_ivp_ode}. However, due to the use of formula (\ref{eq:v_sol_in_node}) for the solution at grid nodes for algebraic variables $\mathbf{v}_{n}$, one would expect convergence orders $N+1$. However, as subsequent calculations for the DAE system of index 1 showed (see Subsection~\ref{sec:2:ct:ex1} ``Example 1: simple DAE system of index 1'' in the Section~\ref{sec:2:ct} ``Classical tests'' below in the text), convergence orders $p_{\rm nodes} = 2N+1$ were obtained. This is due to the superconvergence of the DG methods at downwind points, which was revealed, in particular, in the work~\cite{dg_ivp_ode_3} in the case of DG methods for solving the initial value problem for second-order ODE.

It should be noted that currently there are approaches that allow obtaining a guaranteed superconvergence of the solution with convergence orders $2N+1$ for DG methods with polynomials of degree $N$ not only at nodal and downwind points, but also in continuous norms in ranges such as the $L_{2}$ norm. These approaches are based on post-processing~\cite{siac_ref_1} using filters to improve convergence orders, such as Smoothness-Increasing Accuracy-Conserving (SIAC) filtering. The term SIAC filtering was introduced in~\cite{siac_ref_5}. An important feature of early works on the use of SIAC filters was the difficulty of generalizing to non-uniform grids, which was noted in~\cite{siac_ref_3}, where the post-processing method was extended to improve the convergence orders $2N+1$ for linear hyperbolic equations. For uniform grids, it was shown~\cite{siac_ref_3} that when applying the SIAC filtering to the DG approximation at the final time, the order of accuracy improves from $N+1$ to $2N+1$ for linear hyperbolic equations with periodic boundary conditions. The superconvergence of the convergence orders $2N+1$ is promising, but it was limited to uniform grids~\cite{siac_ref_3}. Only for a certain family of non-uniform grids, smoothly varying grids, it was proven that the filtered solutions have superconvergence of order $2N+1$~\cite{siac_ref_20}. A series of studies on various aspects of SIAC filters are presented in references~\cite{siac_ref_5, siac_ref_11, siac_ref_20, siac_ref_12, siac_ref_15, siac_ref_16}. In the work~\cite{siac_rev}, an SIAC filter for non-uniform grids was proposed and developed. In the present work, the numerical ADER-DG method with a local DG predictor is developed for initially non-uniform grids, and possibly using post-processing based on SIAC filtering will allow to achieve a more formally guaranteed superconvergence for the node solution for the algebraic variables $\mathbf{v}_{n}$ and for the local solution $(\mathbf{u}_{L}(t), \mathbf{v}_{L}(t))$, however, this remains for future research. Probably, this extension of the ADER-DG method with a local DG predictor will be based on the SIAC filter developed in the work~\cite{siac_rev}.

\subsection{DAE solver as ODE solver}

In the present work, the numerical method ADER-DG with a local DG predictor is investigated, the functional representations of which are based on the expansion of a local discrete time solution $\mathbf{q}$ over a set of Lagrange interpolation polynomials $\varphi_{p}$ with nodal points at the roots of the right Radau polynomials. In the work~\cite{ader_dg_ivp_ode}, where the numerical method ADER-DG with a local DG predictor was investigated for solving initial value problems for ODE systems, Lagrange interpolation polynomials with nodal points at the roots of shifted Legendre polynomials were used, and it was shown that this version of the numerical method has very high convergence orders in solving ODE systems. The values of the convergence orders were comparable and sometimes exceeded the expected values (\ref{eq:expect_orders}). Therefore, it is of certain interest to calculate convergence orders of the ADER-DG method with a local DG predictor, in which the node of the functional basis is constructed at the nodes of the Gauss-Radau quadrature formula (\ref{eq:gr_rule}) in solving the initial value problem for ODE systems. Otherwise, it would be unclear what the decrease in the convergence order in solving the DAE system compared to the ODE system is associated with --- suddenly the numerical method ADER-DG with a local DG predictor with such a functional basis initially did not have the expected convergence orders. 

The solution to the initial value problem for the ODE system was obtained using the developed ADER-DG method with a local DG predictor, the original problem for which was reduced to an ODE system, which is completely equivalent to the numerical method studied in the work~\cite{ader_dg_ivp_ode}. The basis functions were chosen to be Lagrange interpolation polynomials with nodal points at the roots of the right Radau polynomials, which is the main difference from the formulation of the ADER-DG method with a local DG predictor in the work~\cite{ader_dg_ivp_ode}.

The initial value problem for the ODE system, describing a linear one-dimensional harmonic oscillator, was chosen as a test problem (Example 1 in~\cite{ader_dg_ivp_ode}): $\ddot{x} + x = 0$, $x(0) = 1$, $\dot{x}(0) = 0$, which was reformulated as initial value problem for the ODE of the first-order system:
\begin{equation}\label{eq:harm_osc_ode}
\frac{du_{1}}{dt} = u_{2};\quad \frac{du_{2}}{dt} = -u_{1};\quad 
u_{1}(0) = 1;\quad u_{2}(0) = 0;
\end{equation}
where $\mathbf{u} = [u_{1},\, u_{2}]^{T} = [x,\, \dot{x}]^{T}$ is a desired two component vector function. The exact analytical solution of this problem has the form $\mathbf{u}^{\rm ex} = [\cos(t),\, -\sin(t)]^{T}$. The solution definition domain $\Omega$ was chosen to be equal to two oscillation periods as $0 \leqslant t \leqslant 4\pi$.

The study of the convergence of the numerical solution was carried out on the basis of the analysis of the solution for the set of $6$ different uniform discretizations of the domain of definition $\Omega$ of the desired function, with the number of grid nodes $L = 10$, $12$, $14$, $16$, $18$, $20$; note that the zero node, in which the initial condition was defined, also refers to grid nodes, so the discretization step was determined by the expression $\Delta t = 4\pi/(L-1)$.

\begin{figure}[h!]
\captionsetup[subfigure]{%
	position=bottom,
	font+=smaller,
	textfont=normalfont,
	singlelinecheck=off,
	justification=raggedright
}
\centering
\begin{subfigure}{0.320\textwidth}
\includegraphics[width=\textwidth]{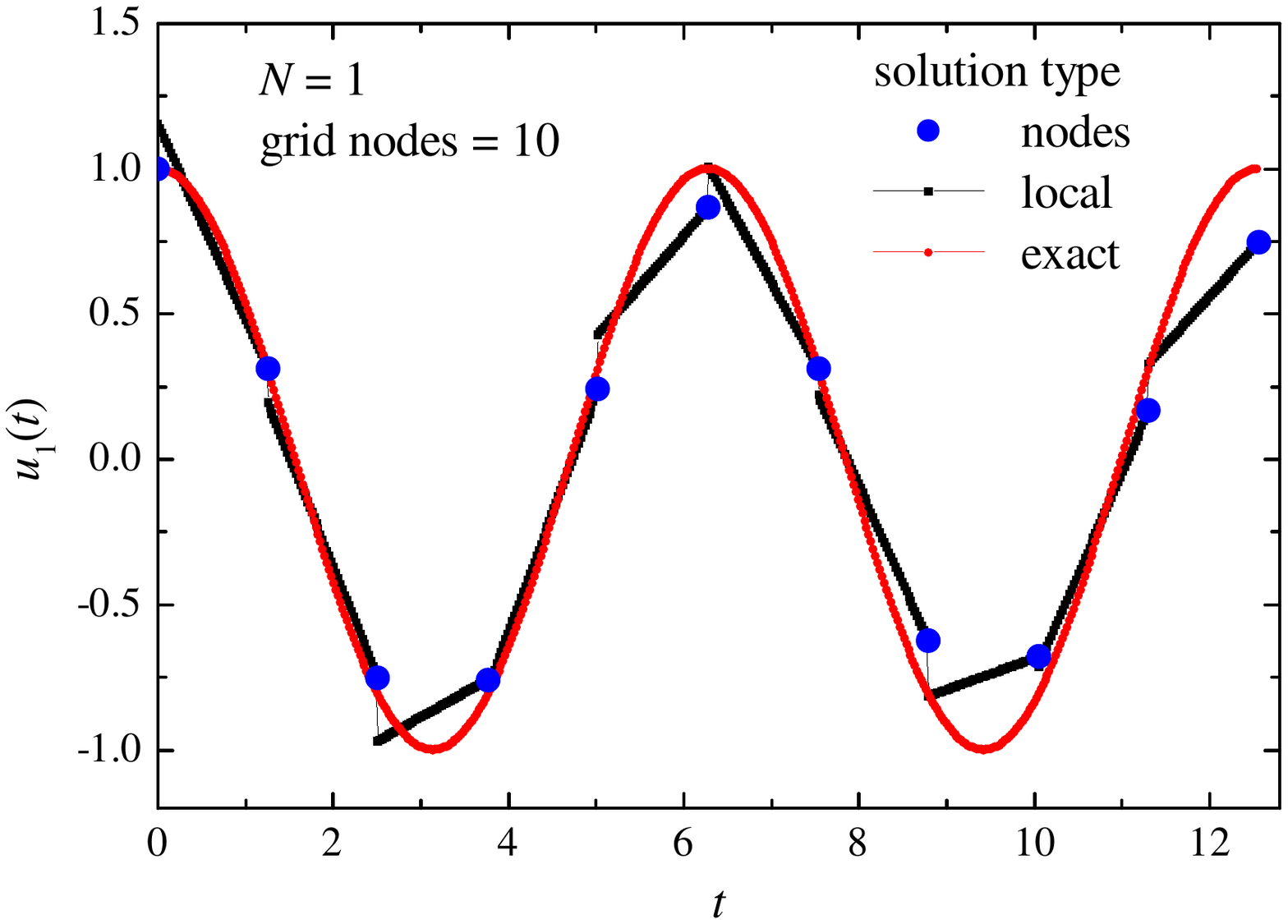}
\vspace{-10mm}\caption{\label{fig:harm_osc:a1}}
\end{subfigure}
\begin{subfigure}{0.320\textwidth}
\includegraphics[width=\textwidth]{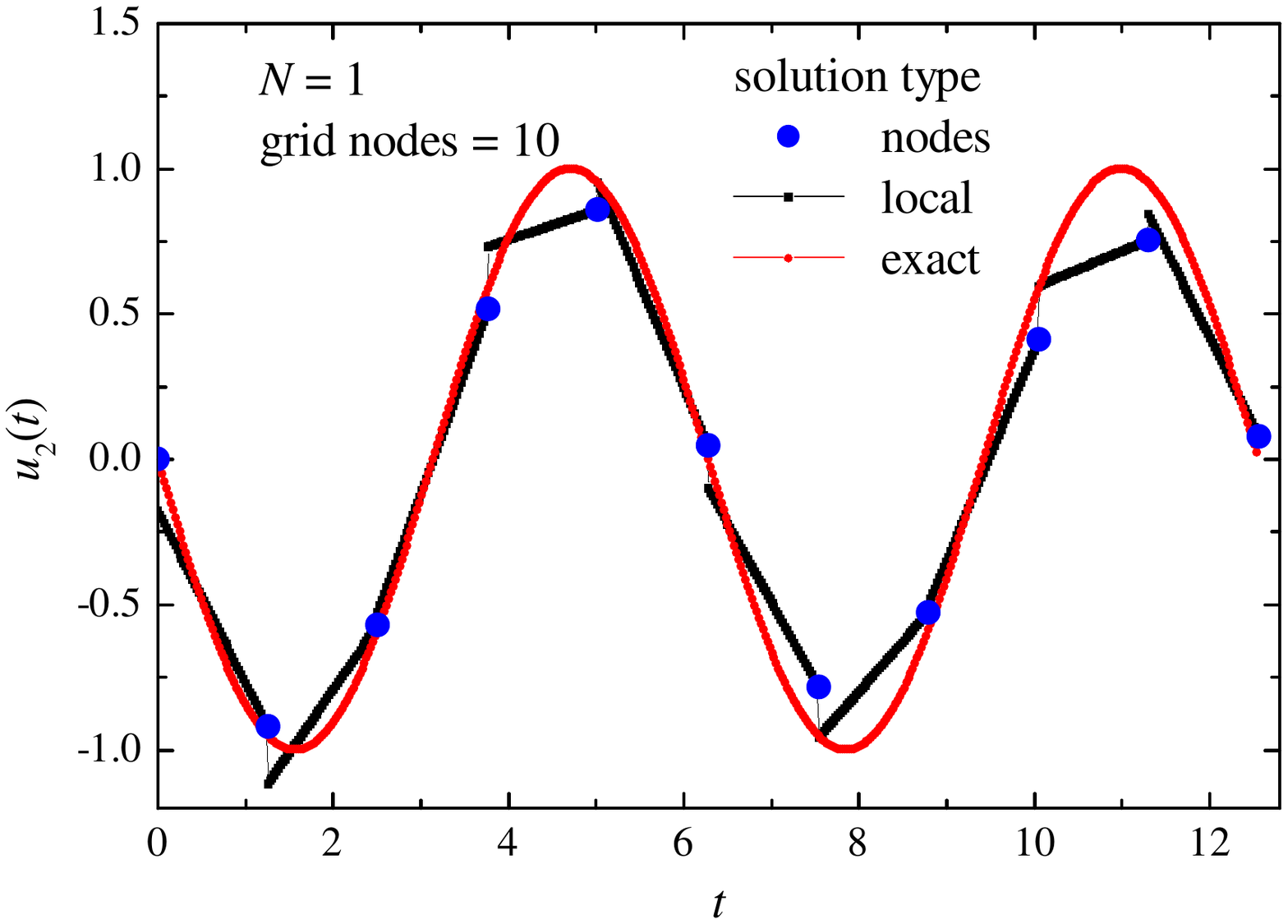}
\vspace{-10mm}\caption{\label{fig:harm_osc:a2}}
\end{subfigure}
\begin{subfigure}{0.320\textwidth}
\includegraphics[width=\textwidth]{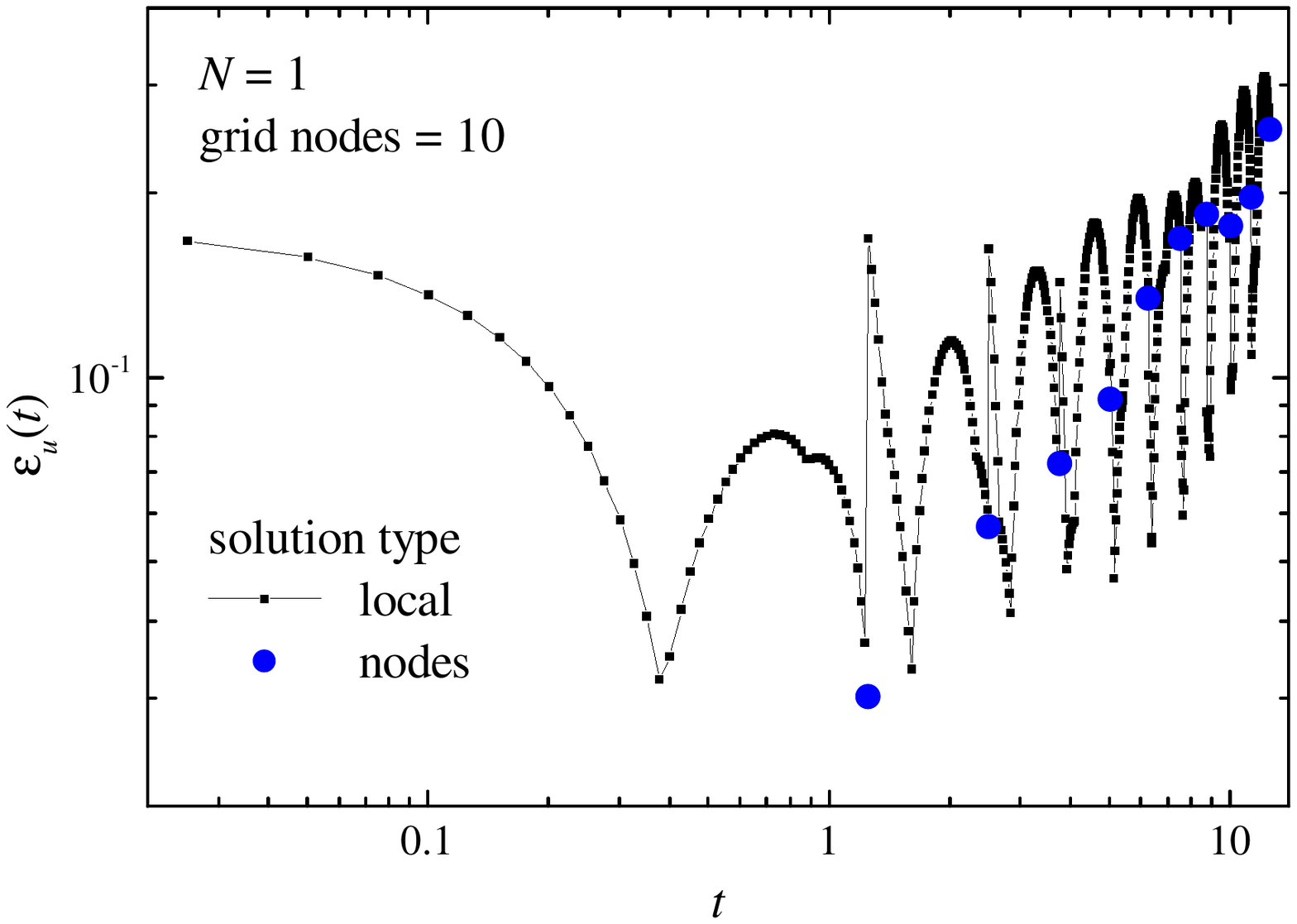}
\vspace{-10mm}\caption{\label{fig:harm_osc:a3}}
\end{subfigure}\\
\begin{subfigure}{0.320\textwidth}
\includegraphics[width=\textwidth]{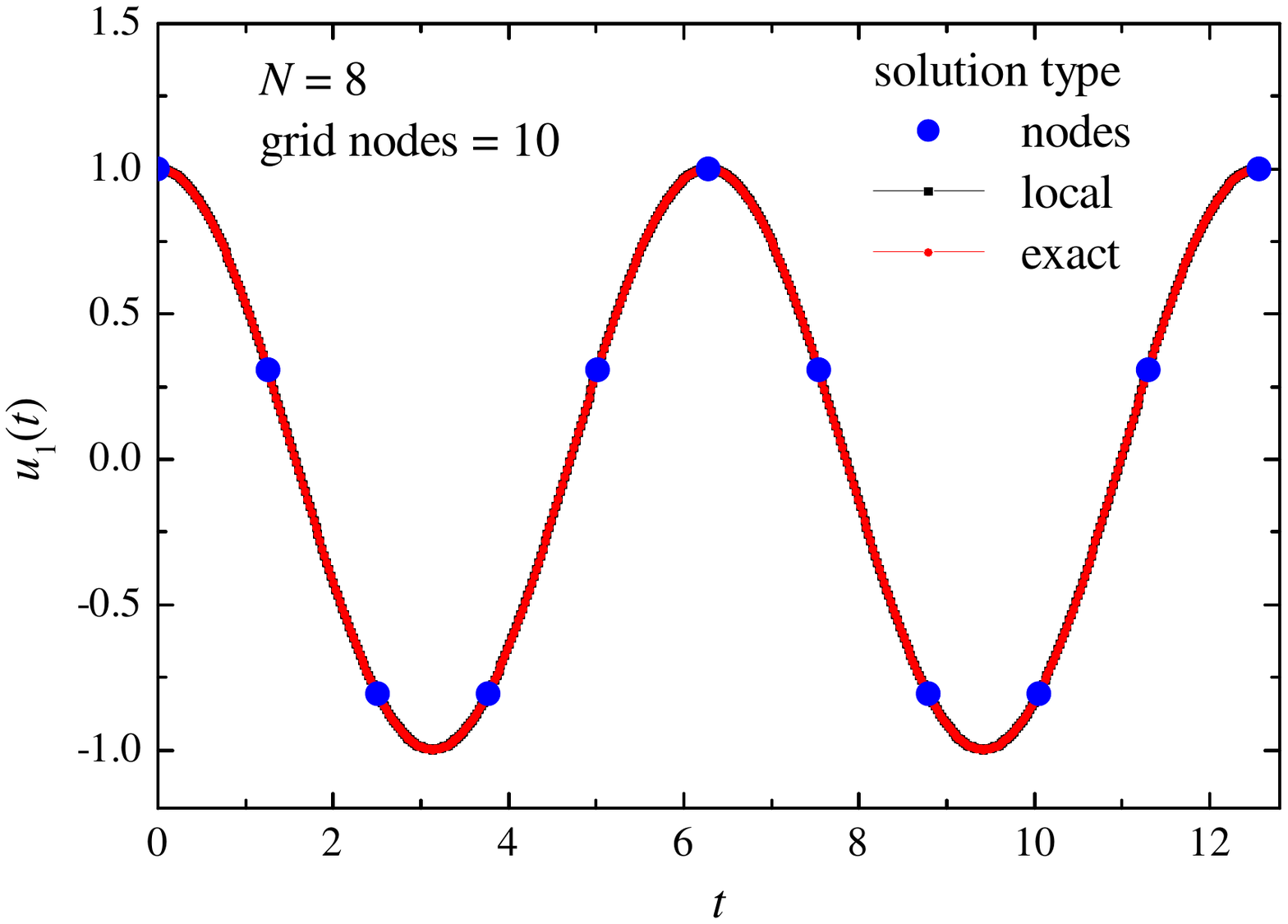}
\vspace{-10mm}\caption{\label{fig:harm_osc:b1}}
\end{subfigure}
\begin{subfigure}{0.320\textwidth}
\includegraphics[width=\textwidth]{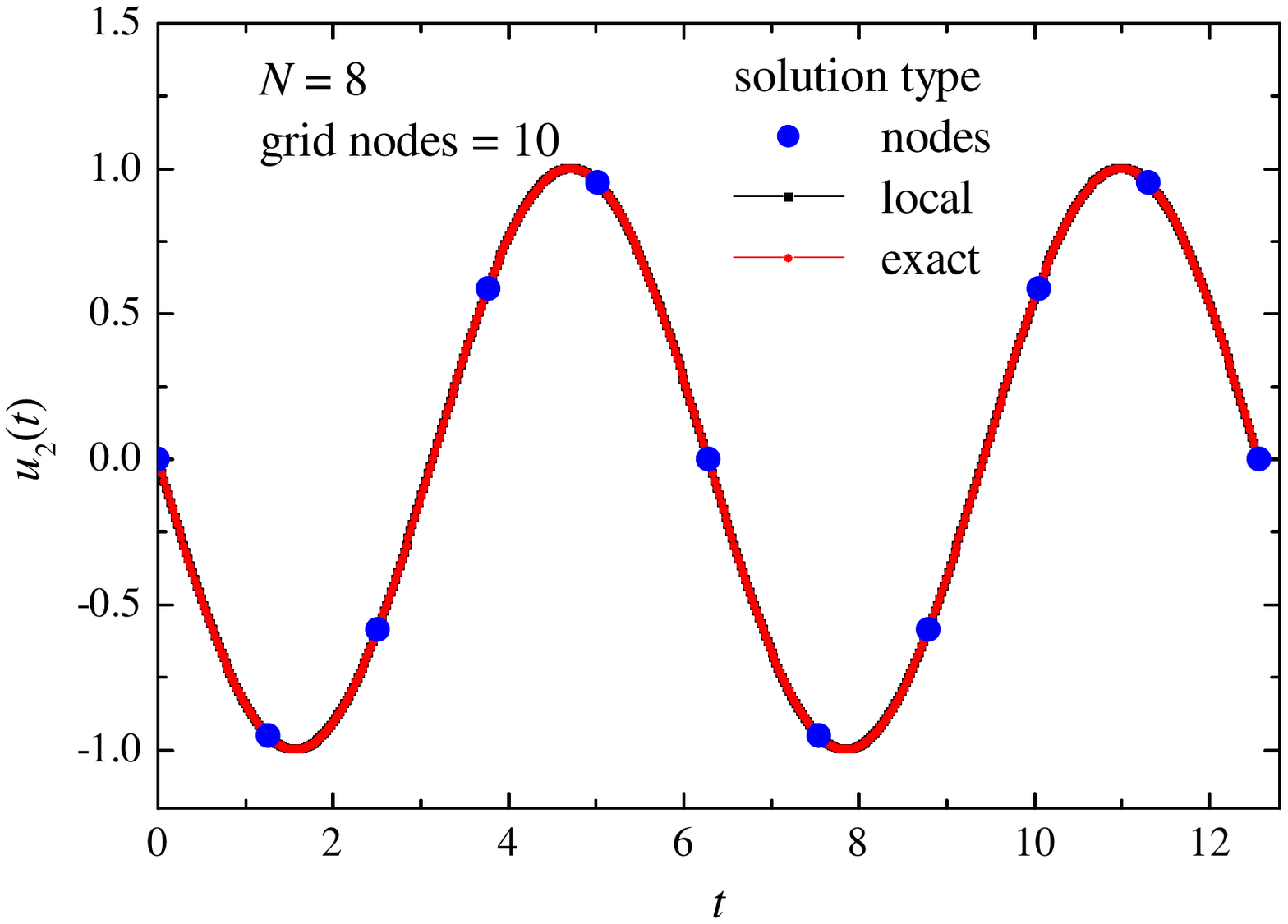}
\vspace{-10mm}\caption{\label{fig:harm_osc:b2}}
\end{subfigure}
\begin{subfigure}{0.320\textwidth}
\includegraphics[width=\textwidth]{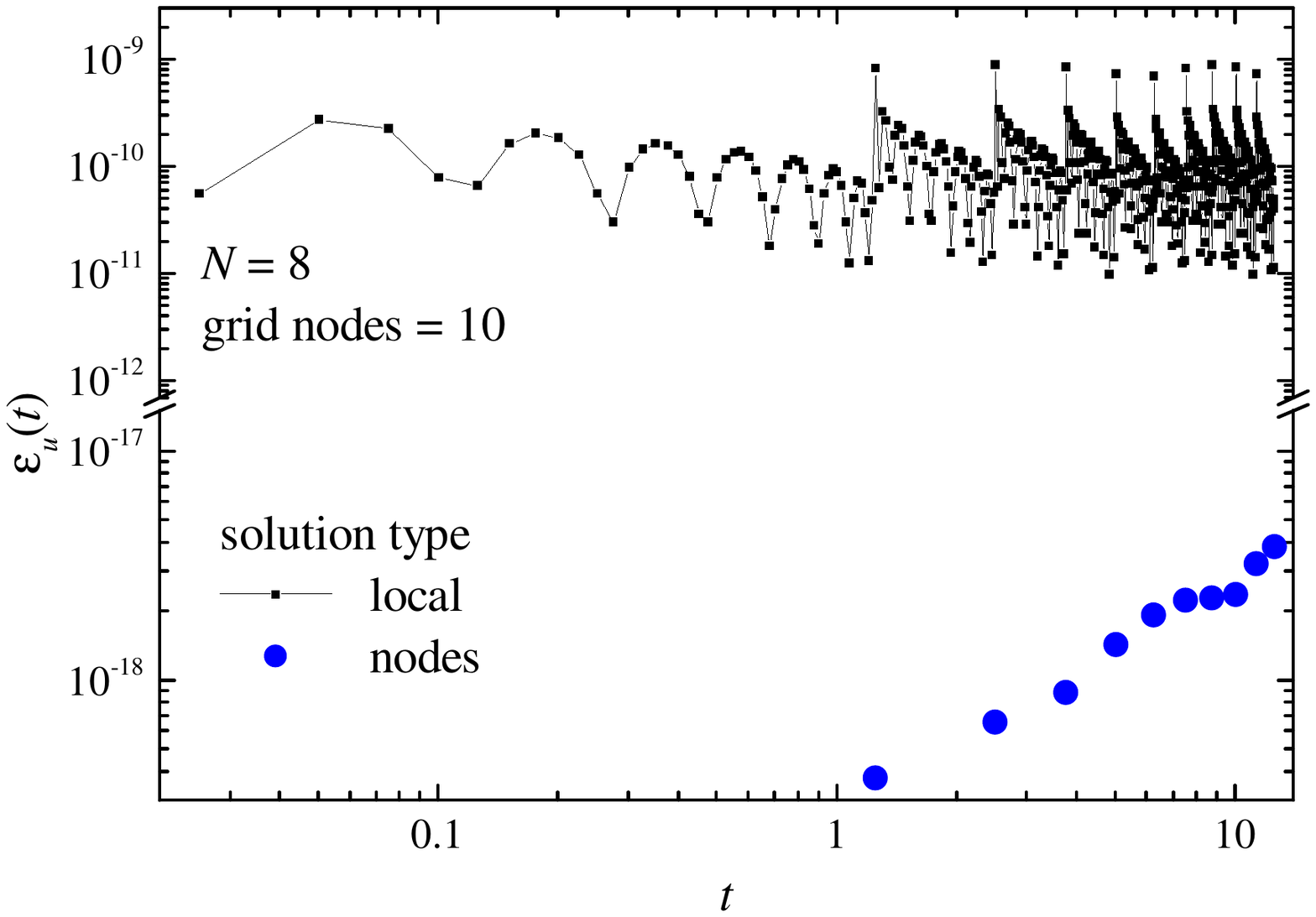}
\vspace{-10mm}\caption{\label{fig:harm_osc:b3}}
\end{subfigure}\\
\begin{subfigure}{0.320\textwidth}
\includegraphics[width=\textwidth]{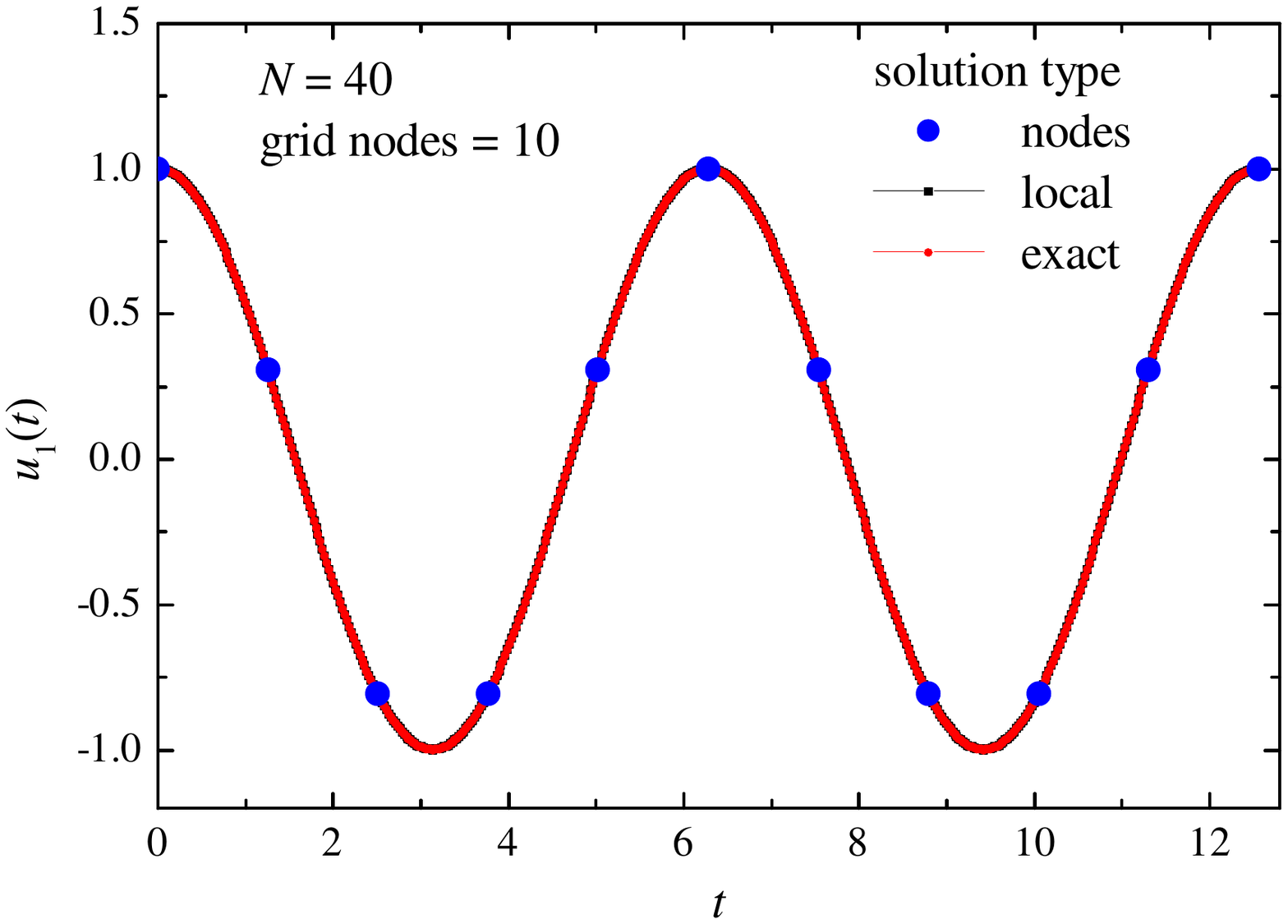}
\vspace{-10mm}\caption{\label{fig:harm_osc:c1}}
\end{subfigure}
\begin{subfigure}{0.320\textwidth}
\includegraphics[width=\textwidth]{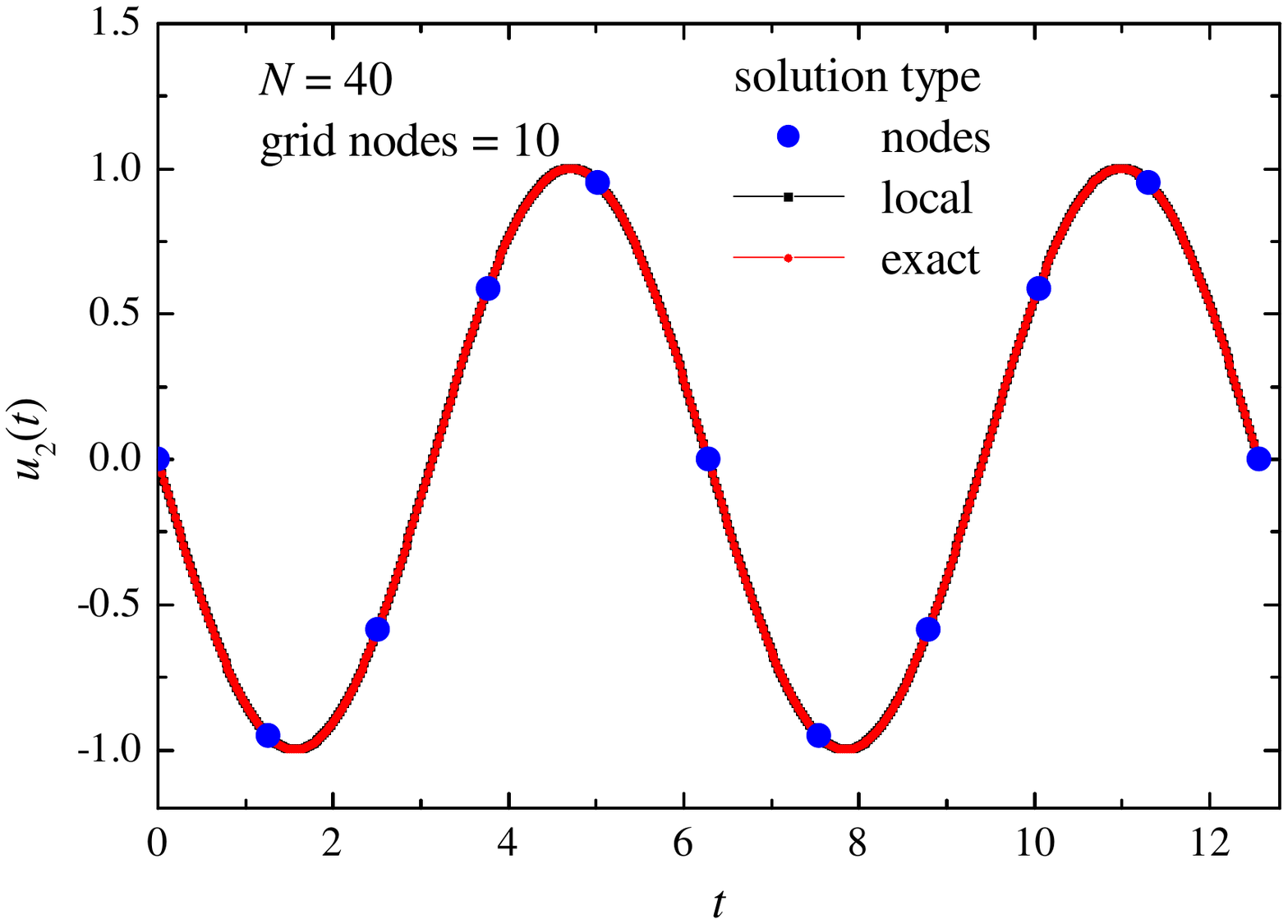}
\vspace{-10mm}\caption{\label{fig:harm_osc:c2}}
\end{subfigure}
\begin{subfigure}{0.320\textwidth}
\includegraphics[width=\textwidth]{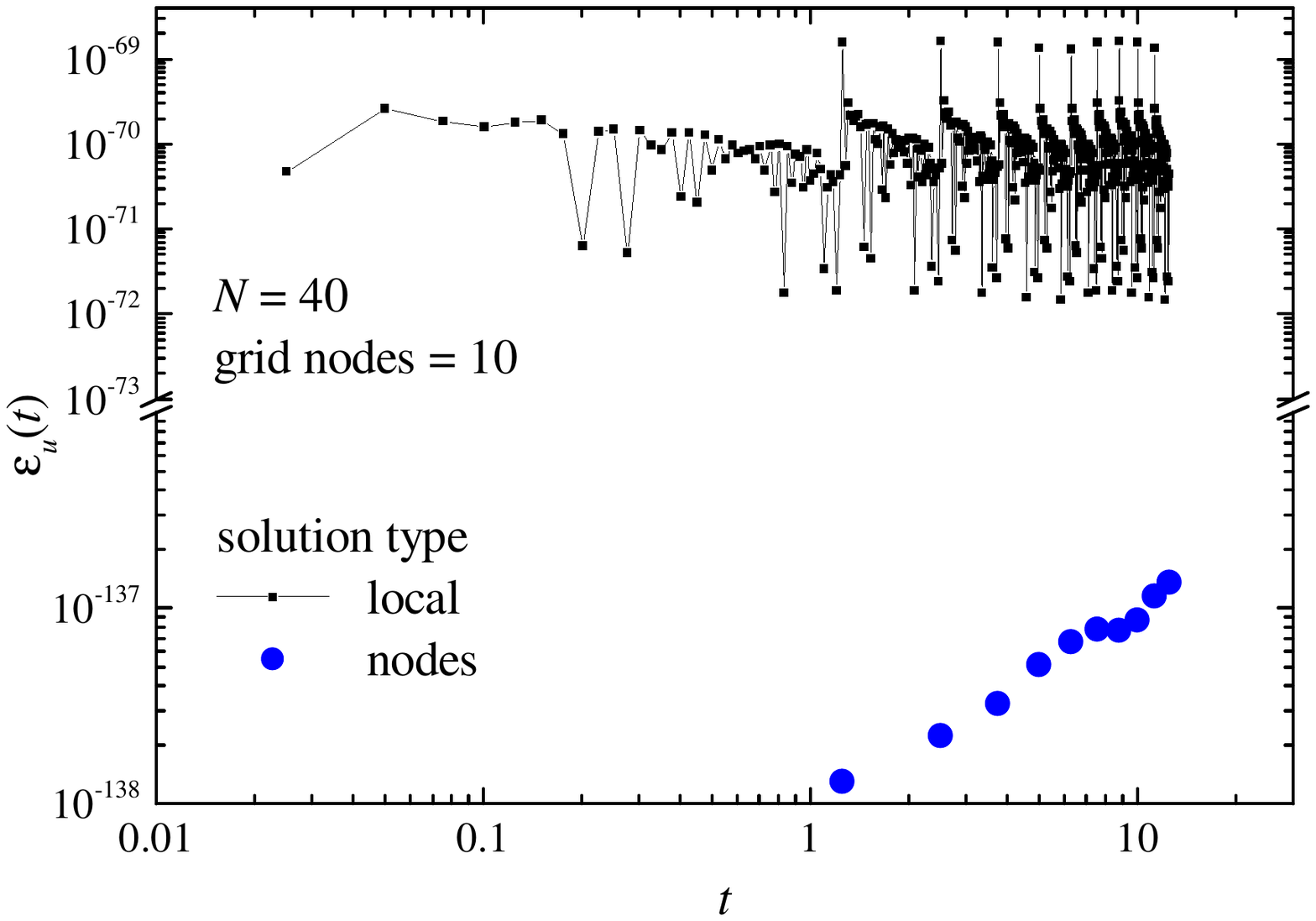}
\vspace{-10mm}\caption{\label{fig:harm_osc:c3}}
\end{subfigure}\\
\begin{subfigure}{0.320\textwidth}
\includegraphics[width=\textwidth]{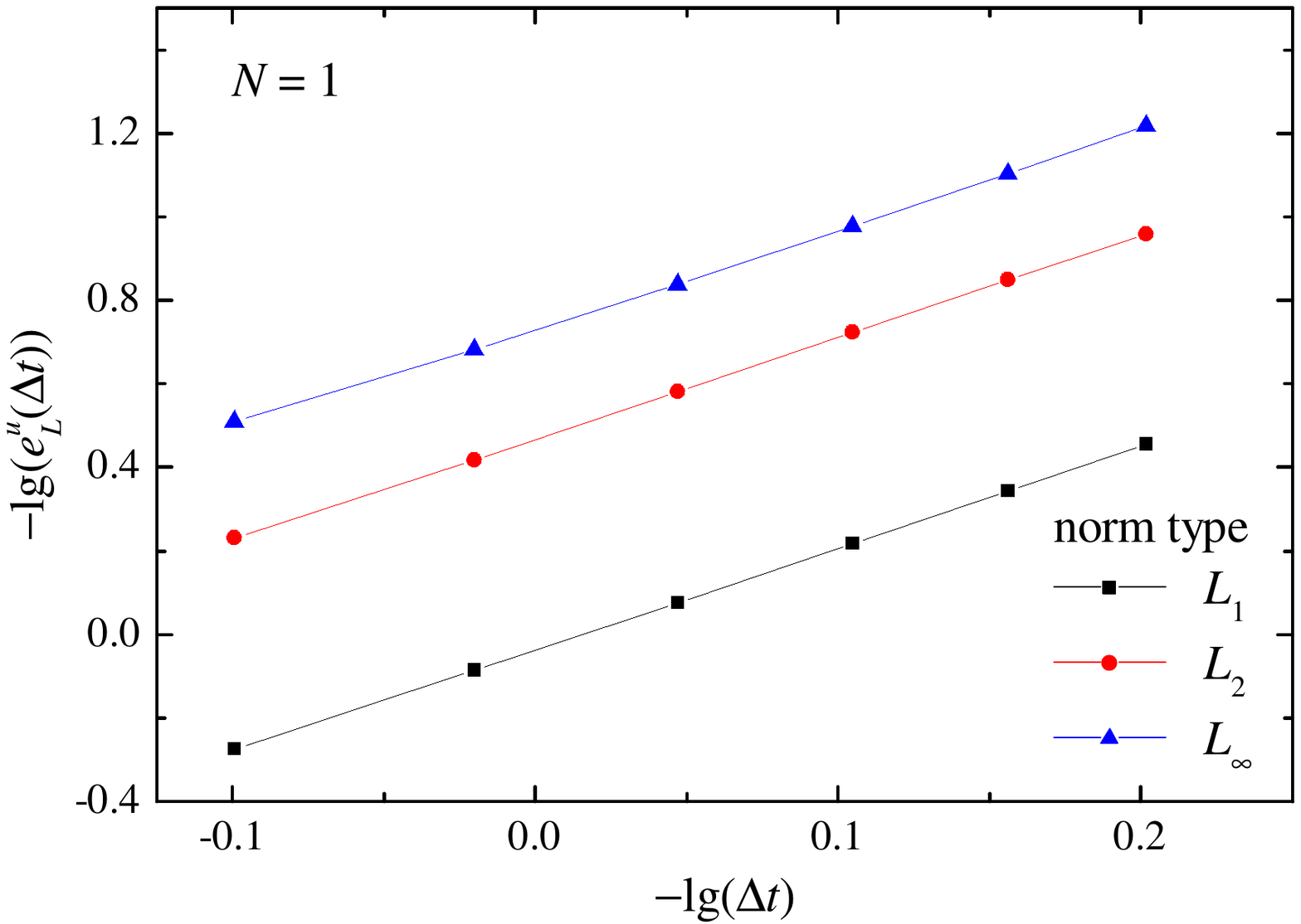}
\vspace{-10mm}\caption{\label{fig:harm_osc:d1}}
\end{subfigure}
\begin{subfigure}{0.320\textwidth}
\includegraphics[width=\textwidth]{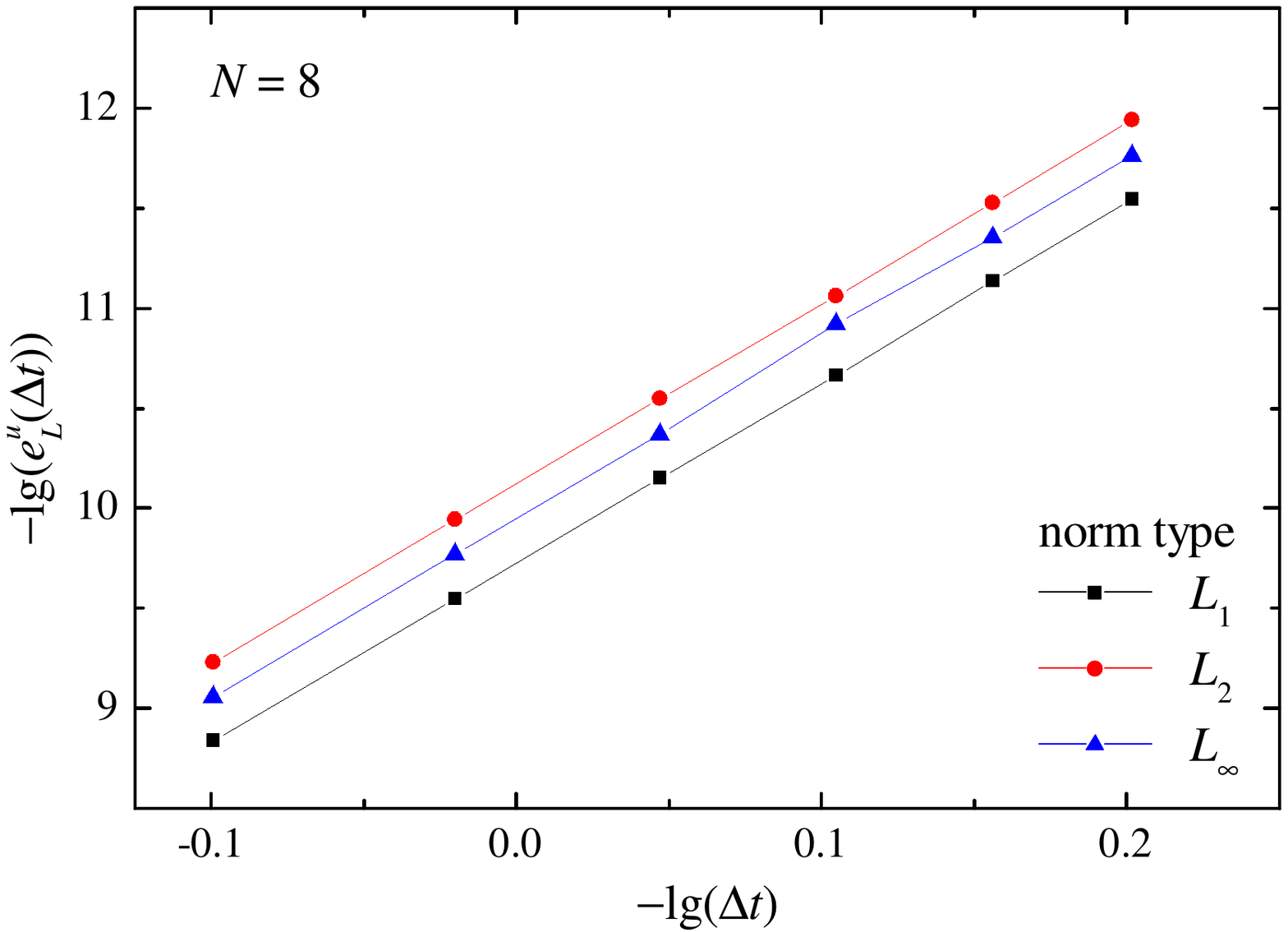}
\vspace{-10mm}\caption{\label{fig:harm_osc:d2}}
\end{subfigure}
\begin{subfigure}{0.320\textwidth}
\includegraphics[width=\textwidth]{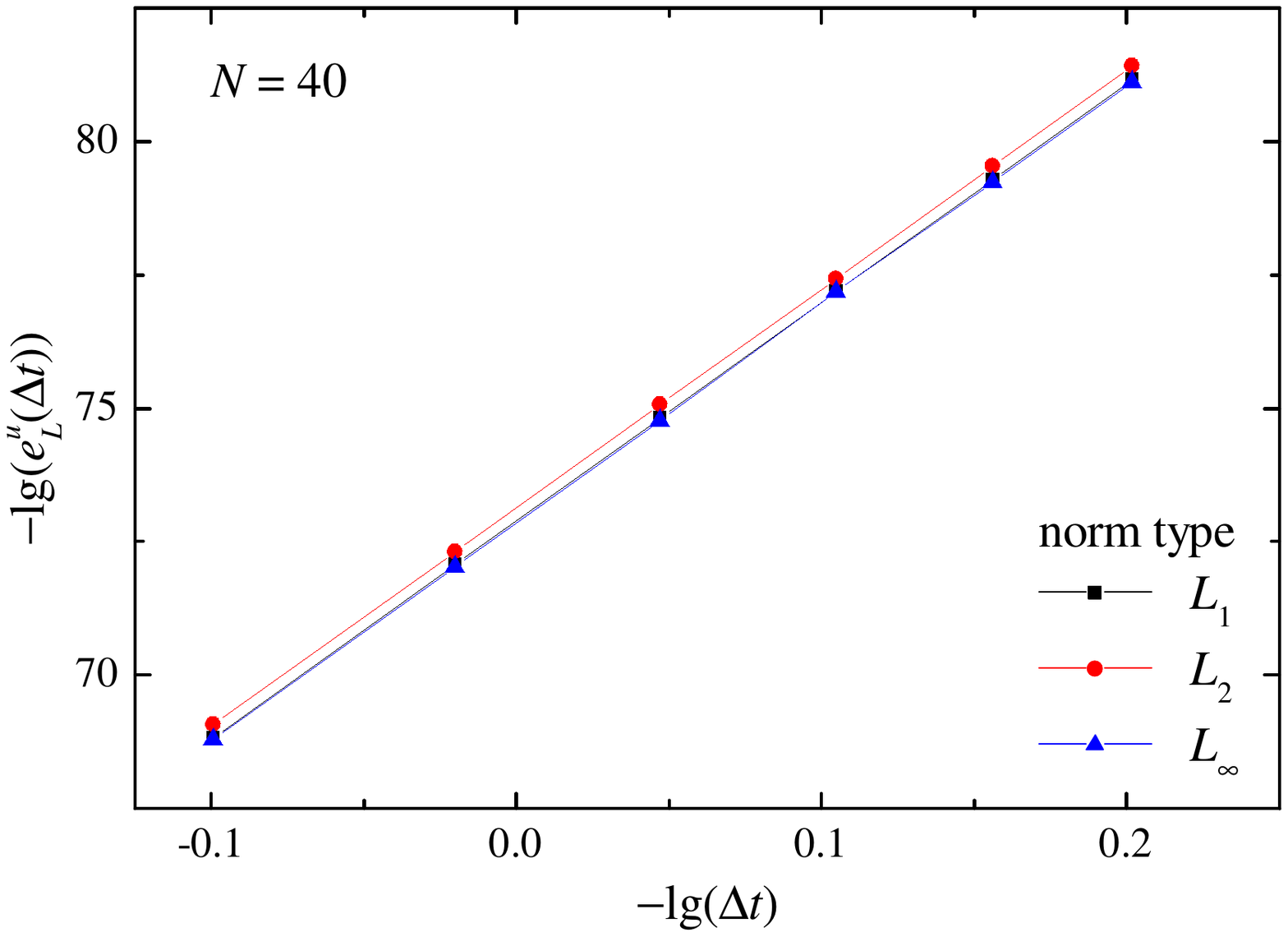}
\vspace{-10mm}\caption{\label{fig:harm_osc:d3}}
\end{subfigure}\\
\begin{subfigure}{0.320\textwidth}
\includegraphics[width=\textwidth]{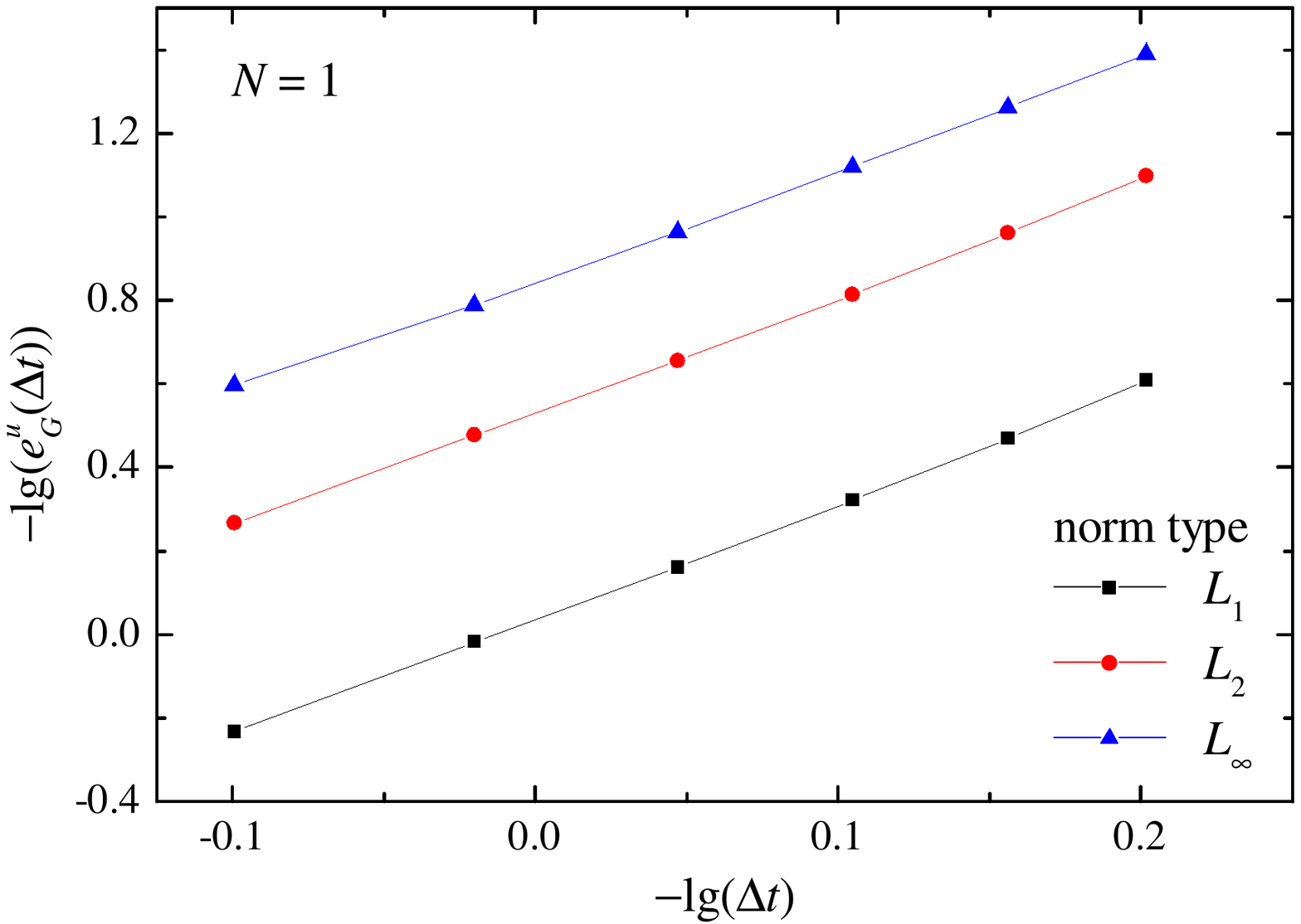}
\vspace{-10mm}\caption{\label{fig:harm_osc:e1}}
\end{subfigure}
\begin{subfigure}{0.320\textwidth}
\includegraphics[width=\textwidth]{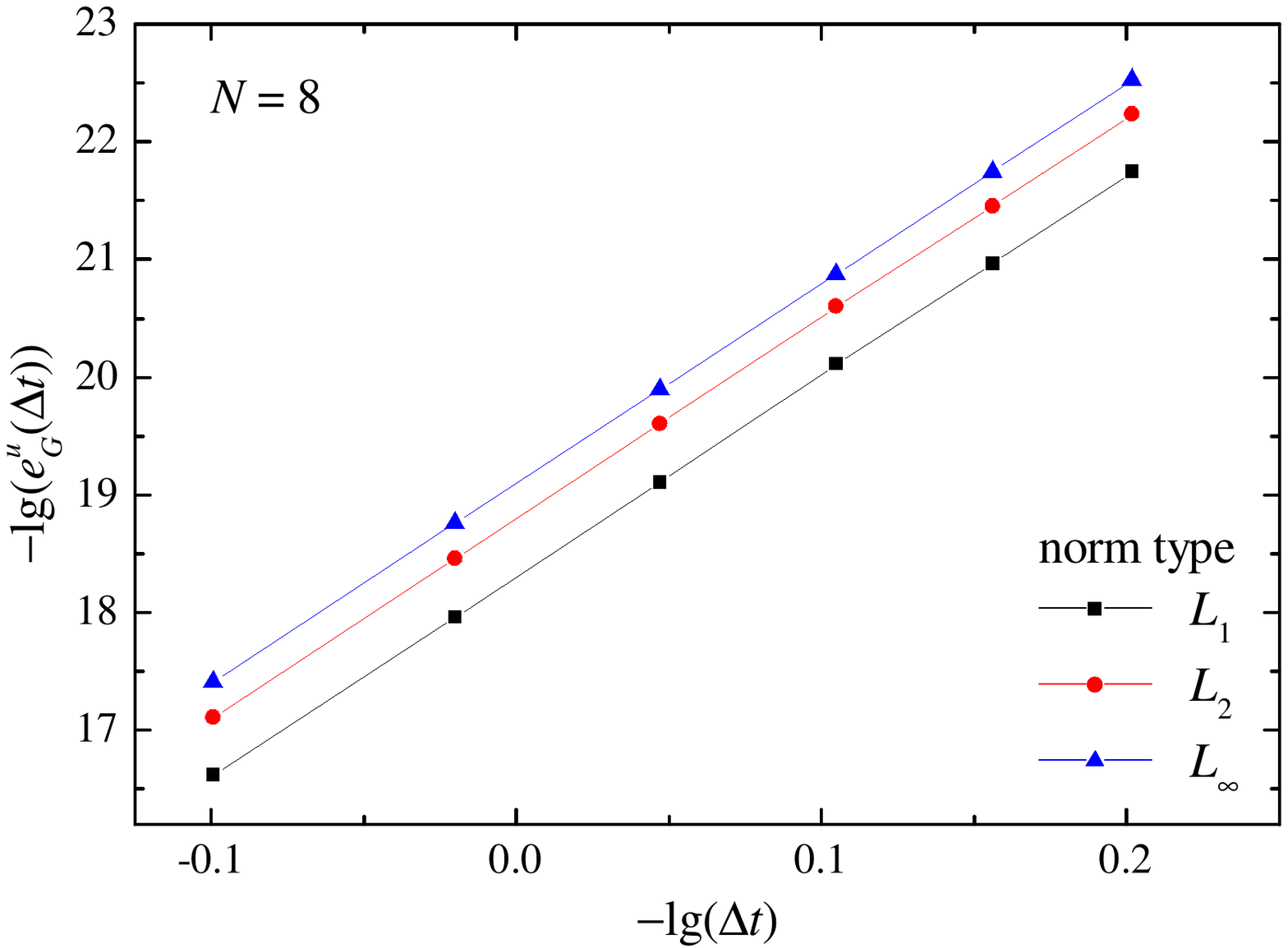}
\vspace{-10mm}\caption{\label{fig:harm_osc:e2}}
\end{subfigure}
\begin{subfigure}{0.320\textwidth}
\includegraphics[width=\textwidth]{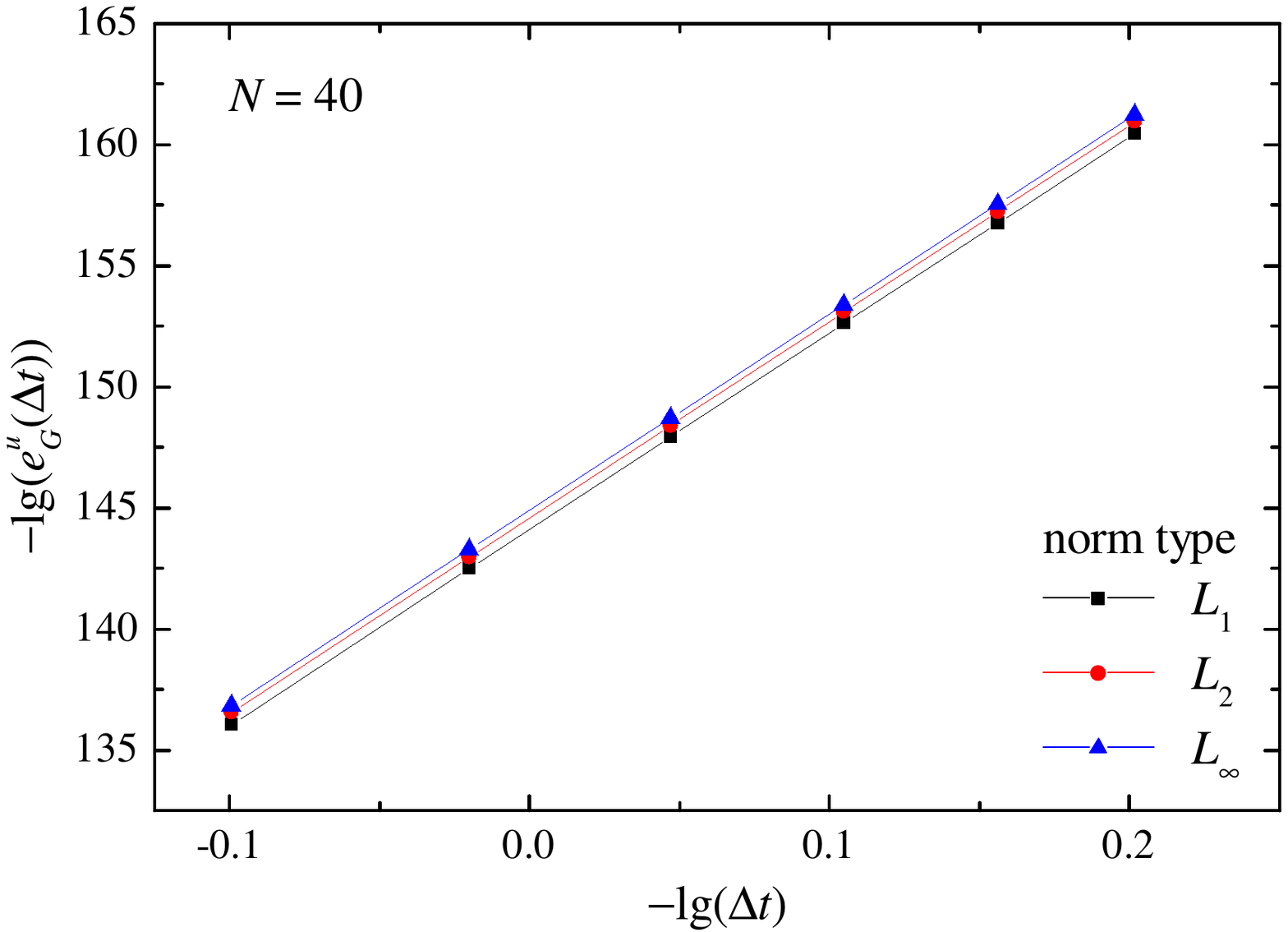}
\vspace{-10mm}\caption{\label{fig:harm_osc:e3}}
\end{subfigure}\\
\caption{%
Numerical solution of the ODE system (\ref{eq:harm_osc_ode}). Comparison of the solution at nodes $\mathbf{u}_{n}$, the local solution $\mathbf{u}_{L}(t)$ and the exact solution $\mathbf{u}^{\rm ex}(t)$ for components $u_{1}$ (\subref{fig:harm_osc:a1}, \subref{fig:harm_osc:b1}, \subref{fig:harm_osc:c1}) and $u_{2}$ (\subref{fig:harm_osc:a2}, \subref{fig:harm_osc:b2}, \subref{fig:harm_osc:c2}), the error $\varepsilon_{u}(t)$ (\subref{fig:harm_osc:a3}, \subref{fig:harm_osc:b3}, \subref{fig:harm_osc:c3}), obtained using polynomials with degrees $N = 1$ (\subref{fig:harm_osc:a1}, \subref{fig:harm_osc:a2}, \subref{fig:harm_osc:a3}), $N = 8$ (\subref{fig:harm_osc:b1}, \subref{fig:harm_osc:b2}, \subref{fig:harm_osc:b3}) and $N = 40$ (\subref{fig:harm_osc:c1}, \subref{fig:harm_osc:c2}, \subref{fig:harm_osc:c3}). Log-log plot of the dependence of the global error for the local solution $e_{L}^{u}$ (\subref{fig:harm_osc:d1}, \subref{fig:harm_osc:d2}, \subref{fig:harm_osc:d3}) and the solution at nodes $e_{G}^{u}$ (\subref{fig:harm_osc:e1}, \subref{fig:harm_osc:e2}, \subref{fig:harm_osc:e3}) on the discretization step $\mathrm{\Delta}t$, obtained in the norms $L_{1}$, $L_{2}$ and $L_{\infty}$, obtained using polynomials with degrees $N = 1$ (\subref{fig:harm_osc:d1}, \subref{fig:harm_osc:e1}), $N = 8$ (\subref{fig:harm_osc:d2}, \subref{fig:harm_osc:e2}) and $N = 40$ (\subref{fig:harm_osc:d3}, \subref{fig:harm_osc:e3}).
}
\label{fig:harm_osc}
\end{figure}

The results of the study of the application of the numerical method for this problem are presented in Fig.~\ref{fig:harm_osc} and in Table~\ref{tab:conv_orders_ode_harm_osc}. The numerical solution obtained using the ADER-DG numerical method with a local DG predictor accurately reproduces all the main features of the solution. In the case of degree $N = 1$, the discontinuity of the local solution $\mathbf{u}_{L}$ at the nodes is clearly observed. A similar behavior is also observed in cases of other values of the degrees $N > 1$, but this is not clearly visible in these figures. An analysis of the error $\varepsilon(t)$ for $N > 1$ shows that the error scales for the local solution $\mathbf{u}_{L}$ and the solution at nodes $\mathbf{u}_{n}$ differ by several orders of magnitude, reaching a value of $60$--$65$ orders of magnitude for a degree of $N = 40$. It can be noted that for the error of the numerical solution at the nodes $\mathbf{u}_{n}$, there is an approximately linear increase in the error $\varepsilon(t)$ with an increase in the argument $t$. The calculated convergence orders $p_{L_{1}}$, $p_{L_{2}}$, $p_{L_{\infty}}$ for various functional norms are presented in Table~\ref{tab:conv_orders_ode_harm_osc}.

All obtained convergence orders $p$ sufficiently correspond to the expected theoretical values $p_{\rm (G)}$ and $p_{\rm (L)}$. There is a slight downward mismatch for low degrees of polynomials $N = 1,\, 2$. As a result, it becomes clear that the solution at nodes $\mathbf{u}_{n}$ obtained by the ADER-DG numerical method with a local DG predictor demonstrates the classical superconvergence $2N+1$ expected from such numerical methods, at least in problem (\ref{eq:harm_osc_ode}). The local solution $\mathbf{u}_{L}$ obtained by the DG predictor demonstrates the classical convergence characteristic of DG methods.

\begin{table*}[h!]
\centering
\caption{%
Convergence orders $p_{L_{1}}$, $p_{L_{2}}$, $p_{L_{\infty}}$, calculated in norms $L_{1}$, $L_{2}$, $L_{\infty}$, of the numerical solution of the ADER-DG method for the ODE problem; $N$ is the degree of the basis polynomials $\varphi_{p}$. Orders $p^{n, u}$ are calculated for \textit{the numerical solution at the nodes} $\mathbf{u}_{n}$; orders $p^{l, u}$ --- for \textit{the local solution} $\mathbf{u}_{L}$. The theoretical values of convergence orders $p_{\rm th.}^{n} = 2N+1$ and $p_{\rm th.}^{l} = N+1$ are presented for comparison.
}
\label{tab:conv_orders_ode_harm_osc}
\begin{tabular}{@{}|l|lll|c|lll|c|@{}}
\toprule
$N$ & $p_{L_{1}}^{n, u}$ & $p_{L_{2}}^{n, u}$ & $p_{L_{\infty}}^{n, u}$  & $p_{\rm th.}^{n}$ & $p_{L_{1}}^{l, u}$ & $p_{L_{2}}^{l, u}$ & $p_{L_{\infty}}^{l, u}$ & $p_{\rm th.}^{l}$ \\
\midrule
$1$	&	$2.78$	&	$2.74$	&	$2.65$	&	$3$	&	$2.43$	&	$2.42$	&	$2.36$	&	$2$\\
$2$	&	$4.97$	&	$4.95$	&	$4.89$	&	$5$	&	$3.14$	&	$3.03$	&	$2.94$	&	$3$\\
$3$	&	$7.00$	&	$6.98$	&	$6.93$	&	$7$	&	$4.00$	&	$3.98$	&	$3.99$	&	$4$\\
$4$	&	$9.02$	&	$9.00$	&	$8.95$	&	$9$	&	$5.00$	&	$4.98$	&	$4.99$	&	$5$\\
$5$	&	$11.03$	&	$11.02$	&	$10.96$	&	$11$	&	$6.00$	&	$5.98$	&	$5.99$	&	$6$\\
$6$	&	$13.04$	&	$13.03$	&	$12.97$	&	$13$	&	$6.99$	&	$6.99$	&	$7.00$	&	$7$\\
$7$	&	$15.05$	&	$15.03$	&	$14.97$	&	$15$	&	$7.99$	&	$7.99$	&	$8.00$	&	$8$\\
$8$	&	$17.05$	&	$17.04$	&	$16.97$	&	$17$	&	$8.99$	&	$8.99$	&	$9.01$	&	$9$\\
$9$	&	$19.06$	&	$19.04$	&	$18.98$	&	$19$	&	$9.99$	&	$9.99$	&	$10.01$	&	$10$\\
$10$	&	$21.06$	&	$21.05$	&	$20.98$	&	$21$	&	$10.99$	&	$10.99$	&	$11.01$	&	$11$\\
$11$	&	$23.06$	&	$23.05$	&	$22.98$	&	$23$	&	$11.99$	&	$11.99$	&	$12.01$	&	$12$\\
$12$	&	$25.07$	&	$25.05$	&	$24.98$	&	$25$	&	$12.99$	&	$12.99$	&	$13.01$	&	$13$\\
$13$	&	$27.07$	&	$27.05$	&	$26.98$	&	$27$	&	$13.99$	&	$13.99$	&	$14.01$	&	$14$\\
$14$	&	$29.07$	&	$29.06$	&	$28.98$	&	$29$	&	$14.99$	&	$14.99$	&	$15.01$	&	$15$\\
$15$	&	$31.07$	&	$31.06$	&	$30.99$	&	$31$	&	$15.99$	&	$15.99$	&	$16.01$	&	$16$\\
$16$	&	$33.07$	&	$33.06$	&	$32.99$	&	$33$	&	$17.00$	&	$16.99$	&	$17.01$	&	$17$\\
$17$	&	$35.07$	&	$35.06$	&	$34.99$	&	$35$	&	$17.99$	&	$17.99$	&	$18.01$	&	$18$\\
$18$	&	$37.08$	&	$37.06$	&	$36.99$	&	$37$	&	$19.00$	&	$18.99$	&	$19.01$	&	$19$\\
$19$	&	$39.08$	&	$39.06$	&	$38.99$	&	$39$	&	$20.00$	&	$19.99$	&	$20.01$	&	$20$\\
$20$	&	$41.08$	&	$41.06$	&	$40.99$	&	$41$	&	$21.00$	&	$20.99$	&	$21.02$	&	$21$\\
\midrule
$25$	&	$51.08$	&	$51.07$	&	$50.99$	&	$51$	&	$26.00$	&	$26.00$	&	$26.02$	&	$26$\\
$30$	&	$61.08$	&	$61.07$	&	$60.99$	&	$61$	&	$31.00$	&	$31.00$	&	$31.02$	&	$31$\\
$35$	&	$71.08$	&	$71.07$	&	$70.99$	&	$71$	&	$36.00$	&	$36.00$	&	$36.02$	&	$36$\\
$40$	&	$81.09$	&	$81.07$	&	$80.99$	&	$81$	&	$41.00$	&	$41.00$	&	$41.02$	&	$41$\\
\bottomrule
\end{tabular}
\end{table*}

Therefore, it can be concluded that using the functional representations which are based on the expansion of a local discrete time solution by a set of Lagrange interpolation polynomials with nodal points at the roots of the right Radau polynomials instead of Lagrange interpolation polynomials with nodal points at the roots of shifted Legendre polynomials does not change very high convergence orders in solving ODE systems of the ADER-DG numerical method with a local DG predictor~\cite{ader_dg_ivp_ode}.

\subsection{On computational cost estimates}

The numerical method ADER-DG with a local DG predictor makes it possible to obtain an arbitrarily high order, which is well achieved for DAE systems (\ref{eq:dae_chosen_form}). Of course, it is of some interest to determine the computational costs that accompany obtaining a numerical solution by this method.

Estimates of computational costs were made based on determining the number of calculations of functions $\mathbf{F}$ and $\mathbf{G}$ included in the right-hand side of the DAE system (\ref{eq:dae_chosen_form}) (the so-called evolutions) and the number of calculations of the Jacobi matrices $\partial\mathbf{F}/\partial\mathbf{u}$, $\partial\mathbf{F}/\partial\mathbf{v}$, $\partial\mathbf{G}/\partial\mathbf{u}$, $\partial\mathbf{G}/\partial\mathbf{v}$. The numerical method ADER-DG with a local DG predictor proposed in this paper requires solving a system of nonlinear algebraic equations (\ref{eq:lst_dg_snae}) at each step, the solution of which leads to obtaining a local discrete time solution $(\mathbf{q},\, \mathbf{r})$, using which, according to formulas (\ref{eq:u_sol_in_node}) and (\ref{eq:v_sol_in_node}), the solution was calculated at the next grid node $t_{n+1}$. The solution of the system of nonlinear algebraic equations (\ref{eq:lst_dg_snae}) was carried out using the classical numerical Newton method for solving systems of nonlinear algebraic equations, based on the linearization of the system of nonlinear equations (\ref{eq:lst_dg_snae}) near the solution approximation $(\mathbf{q}^{i},\, \mathbf{r}^{i})$. The resulting linearized form was presented as a system of linear algebraic equations (\ref{eq:lst_dg_newton_slae}) with a system matrix of size $[(N+1)D] \times [(N+1)D]$, the solution to which was obtained using the LU decomposition (it is clear that in the case of large $D$ and a special form of the system matrix, it is possible to use special solvers that effectively use the properties of special matrices or sparse matrices). To solve the resulting system of linear algebraic equations, a standard method based on LU decomposition was used. This approach is characterized by a computational complexity of $O(D^{3}(N+1)^{3})$. Calculation of the right-hand side of the system of linear algebraic equations (\ref{eq:lst_dg_newton_slae}) at each iteration requires $(N+1)^{2}$ evolutions of functions $\mathbf{F}$ and $\mathbf{G}$. Thus, it can be summarized that the direct algorithmic implementation of the numerical ADER-DG method with a local DG predictor requires at each iteration of the Newton method for solving the system of nonlinear algebraic equations to perform $1$ LU decomposition of a $[(N+1)D] \times [(N+1)D]$ matrix, for the calculation of which it is necessary to perform $N+1$ calculations of the Jacobian matrix $\partial\mathbf{F}/\partial\mathbf{u}$, $\partial\mathbf{F}/\partial\mathbf{v}$, $\partial\mathbf{G}/\partial\mathbf{u}$, $\partial\mathbf{G}/\partial\mathbf{v}$, as well as $(N+1)^{2}$ evolutions of functions $\mathbf{F}$ and $\mathbf{G}$. 

The main one-step numerical methods for solving DAE systems of general form are the $s$-stage implicit Runge-Kutta methods. These computational costs are practically not different from the use of the $s$-stage classical implicit Runge-Kutta methods (stiffly accurate) for solving DAE systems (with $s = N+1$). However, the main difference is that in the case of implicit Runge-Kutta methods, only $s$ evolutions of functions $\mathbf{F}$ and $\mathbf{G}$ are required for one iteration of the Newton method, versus $(N+1)^{2}$ evolutions of functions $\mathbf{F}$ and $\mathbf{G}$ in the case of the ADER-DG method with a local DG predictor. The main computational costs are usually associated not with the calculation of functions $\mathbf{F}$ and $\mathbf{G}$, but with the calculation of the Jacobian matrices, their multiplication by the constants of the method, and with the LU decomposition of the system matrix, which has a size of $sD \times sD$. This approach is characterized by a computational complexity of $O(s^{3}D^{3})$. Calculating the solution $(\mathbf{u}_{n+1},\, \mathbf{v}_{n+1})$ ((\ref{eq:u_sol_in_node}), (\ref{eq:v_sol_in_node})) at the next grid node $t_{n+1}$ based on the local discrete time solution in the ADER-DG method with a local DG predictor requires $N+1$ calculations of functions $\mathbf{F}$ and $\mathbf{G}$, however, these computational costs can be ignored, since these calculations are performed at the stage of solving the system of nonlinear algebraic equations (\ref{eq:lst_dg_snae}) of the DG predictor, and their results can be saved for calculating the solution at the next grid node $t_{n+1}$.

It should be noted that the presented algorithmic optimizations (\ref{eq:lst_dg_newton_slae_v_opt}) and (\ref{eq:lst_dg_newton_slae_u_opt}) can reduce the computational costs of the ADER-DG method in special cases. However, they are applicable not only to the ADER-DG numerical method, but also to the $s$-stage implicit Runge-Kutta methods. Therefore, this aspect of reducing computational costs cannot be attributed only to the ADER-DG numerical method with a local DG predictor.

Therefore, it can be concluded that the presented ADER-DG numerical method with a local DG predictor does not differ significantly from the $s$-stage implicit Runge-Kutta numerical methods with a number of stages $s = N+1$ in terms of computational costs. An important difference is the need to perform $(N+1)^{2}$ evolutions of functions $\mathbf{F}$ and $\mathbf{G}$ in the case of the ADER-DG method with a local DG predictor versus $s$ evolutions of functions $\mathbf{F}$ and $\mathbf{G}$ in the case of the $s$-stage  implicit numerical Runge-Kutta methods, which can make a significant contribution to computational costs in the case where the right-hand side of the DAE system contains ``heavyweight'' functions $\mathbf{F}$ and $\mathbf{G}$ that have ``lightweight'' Jacobian matrices $\partial\mathbf{F}/\partial\mathbf{u}$, $\partial\mathbf{F}/\partial\mathbf{v}$, $\partial\mathbf{G}/\partial\mathbf{u}$, $\partial\mathbf{G}/\partial\mathbf{v}$.

\section{Applications of the numerical method}
\label{sec:2}

This Section presents applications of the numerical method ADER-DG with a local DG predictor, using Lagrange interpolation polynomials with nodes at the roots of the right Radau polynomials $R_{k}$ as basis functions $\varphi_{p}$, for solving DAE systems (\ref{eq:dae_chosen_form}) with indices 1, 2, 3. In cases of problems with indices 2 and 3, numerical solutions of problems (\ref{eq:dae_chosen_form}) with a reduced index down to 1 are also presented. Examples of applications to classical problems and stiffness problems are presented, in particular, to a problem with very high stiffness.

Numerical solutions of the problems are obtained and compared with exact analytical solutions or high-precision reference solutions in the case when there is no exact analytical solution. Global errors $e$ of the numerical solution in various functional norms (\ref{eq:norms_def}) for grids with different numbers of nodes $L$ are calculated, based on which empirical convergence orders $p$ were obtained, which were compared with the expected values.

\subsection{Classical tests}
\label{sec:2:ct}

This Subsection presents Applications of the numerical method ADER-DG with local DG predictor to classical DAE systems (\ref{eq:dae_chosen_form}) with indices 1, 2, 3. Solutions of a simple DAE system of index 1, Hessenberg DAE systems of index 1 and index 2, a solution of the mathematical pendulum problem representing a DAE system of index 3, a solution of the double pendulum problem also representing a DAE system of index 3 are presented. In the case of a double pendulum there is no exact analytical solution to the problem, so a high-precision numerical solution was obtained using the formulation of the problem in the form of an initial value problem for an ODE system, which was solved numerically with very high accuracy.

\subsubsection{Example 1: simple DAE system of index 1}
\label{sec:2:ct:ex1}

The first example of application of the numerical method ADER-DG with local DG predictor consisted in solving a simple DAE of index 1 system, which was previously used in this work for a demonstration example (\ref{eq:demo_test}):
\begin{equation}\label{eq:simple_dae_ind_1}
\begin{split}
&\ddot{x} + x = z - 1,\hspace{17.5mm} x(0) = 1,\quad \dot{x}(0) = 0,\\
&\ddot{y} + y = 1 - z,\hspace{18mm} y(0) = 0,\quad \dot{y}(0) = 1,\\
&g_{1} = x^{2} + y^{2} - z^{2} = 0,\hspace{5mm} z(0) = 1,\quad t\in[0,\, 2\pi].
\end{split}
\end{equation}
The exact analytical solution of this problem was obtained in the following form: $x = \cos(t)$, $\dot{x} = -\sin(t)$, $y = \sin(t)$, $\dot{y} = \cos(t)$, $z \equiv 1$. This problem was rewritten in a form consistent with formulation of the original DAE system (\ref{eq:dae_chosen_form}):
\begin{equation}
\begin{split}
&\frac{du_{1}}{dt} = u_{3},\hspace{26.4mm} u_{1}(0) = 1,\\
&\frac{du_{2}}{dt} = u_{4},\hspace{26.4mm} u_{2}(0) = 0,\\
&\frac{du_{3}}{dt} = -u_{1} + v_{1} - 1,\hspace{10mm} u_{3}(0) = 0,\\
&\frac{du_{4}}{dt} = -u_{2} + 1 - v_{1},\hspace{10mm} u_{4}(0) = 1,\\[1.9mm]
&g_{1} = u_{1}^{2} + u_{2}^{2} - v_{1}^{2} = 0,\hspace{7.5mm} v_{1}(0) = 1,
\end{split}
\end{equation}
where sets of differential variables $\mathbf{u} = [x,\, y,\, \dot{x},\, \dot{y}]^{T}$ and algebraic variables $\mathbf{v} = [z]$ were defined. The vector of algebraic equations (constraints) in the case of this DAE system consisted of only one equation $g_{1} = 0$. The results also present dependencies for the error in fulfilling of ``algebraic equation'' $g_{2} = v_{1} - 1 = 0$, which represents the accuracy of the obtained algebraic variable $v_{1}$. The total number of variables in the chosen problem is $D = 5$, of which $D_{\rm u} = 4$ are differential variables and $D_{\rm v} = 1$ is an algebraic variable. The full-component exact analytical solution of the problem was written in the following form: $\mathbf{u}^{\rm ex} = [\cos(t),\, -\sin(t),\, \sin(t),\, \cos(t)]^{T}$, $\mathbf{v}^{\rm ex} = [1]$.

The domain of definition $[0,\, 2\pi]$ of the desired functions $\mathbf{u}$ and $\mathbf{v}$ was discretized into $L = 10$, $12$, $14$, $16$, $18$, $20$ discretization domains $\Omega_{n}$ with the same discretization steps ${\Delta t}_{n} \equiv {\Delta t} = 2\pi/(L-1)$, which was done to be able to calculate the empirical convergence orders $p$. Therefore, the empirical convergence orders $p$ were calculated by least squares approximation of the dependence of the global error $e$ on the grid discretization step ${\Delta t}$ at 6 data points.

\begin{figure}[h!]
\captionsetup[subfigure]{%
	position=bottom,
	font+=smaller,
	textfont=normalfont,
	singlelinecheck=off,
	justification=raggedright
}
\centering
\begin{subfigure}{0.320\textwidth}
\includegraphics[width=\textwidth]{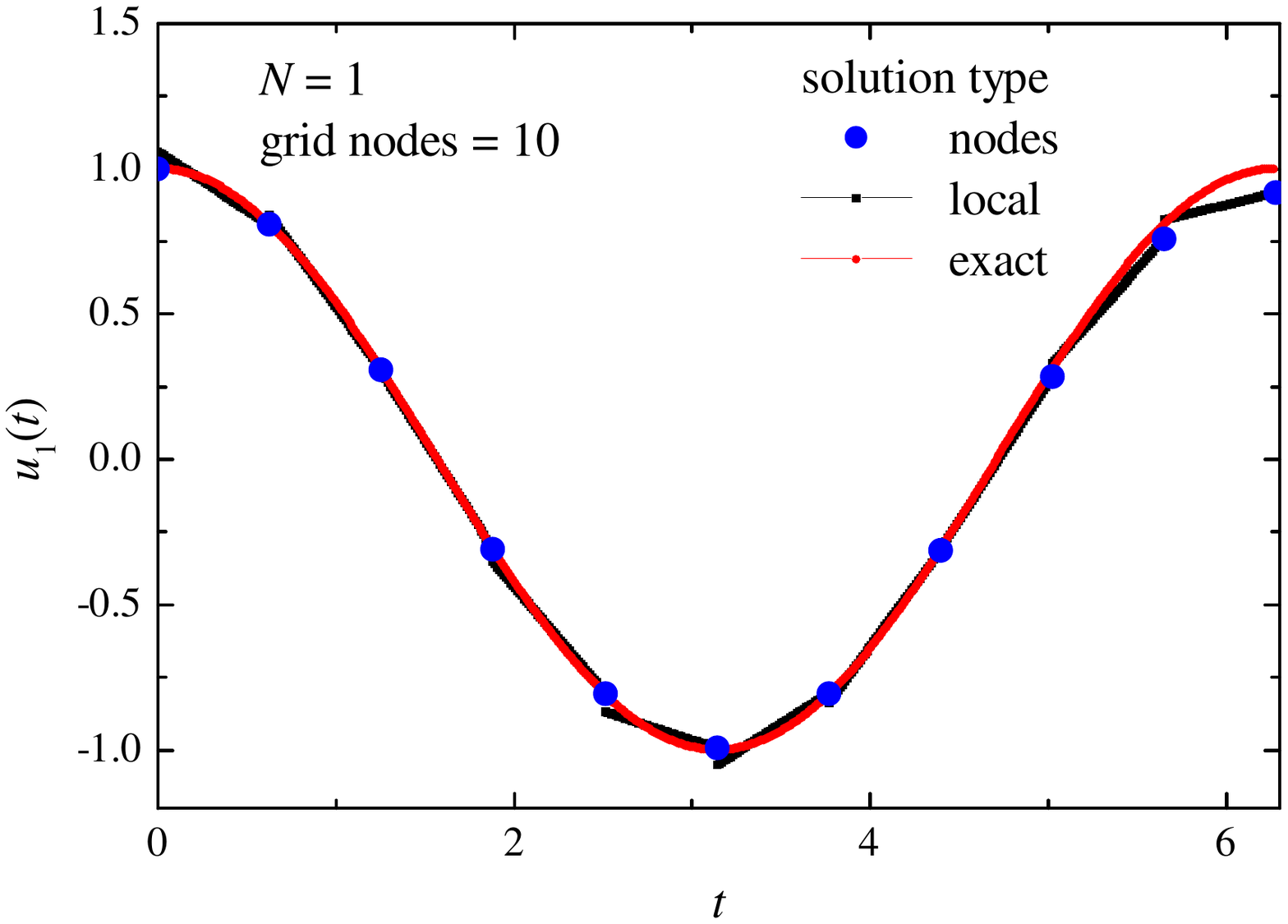}
\vspace{-8mm}\caption{\label{fig:simple_test_sol_uv:a1}}
\end{subfigure}
\begin{subfigure}{0.320\textwidth}
\includegraphics[width=\textwidth]{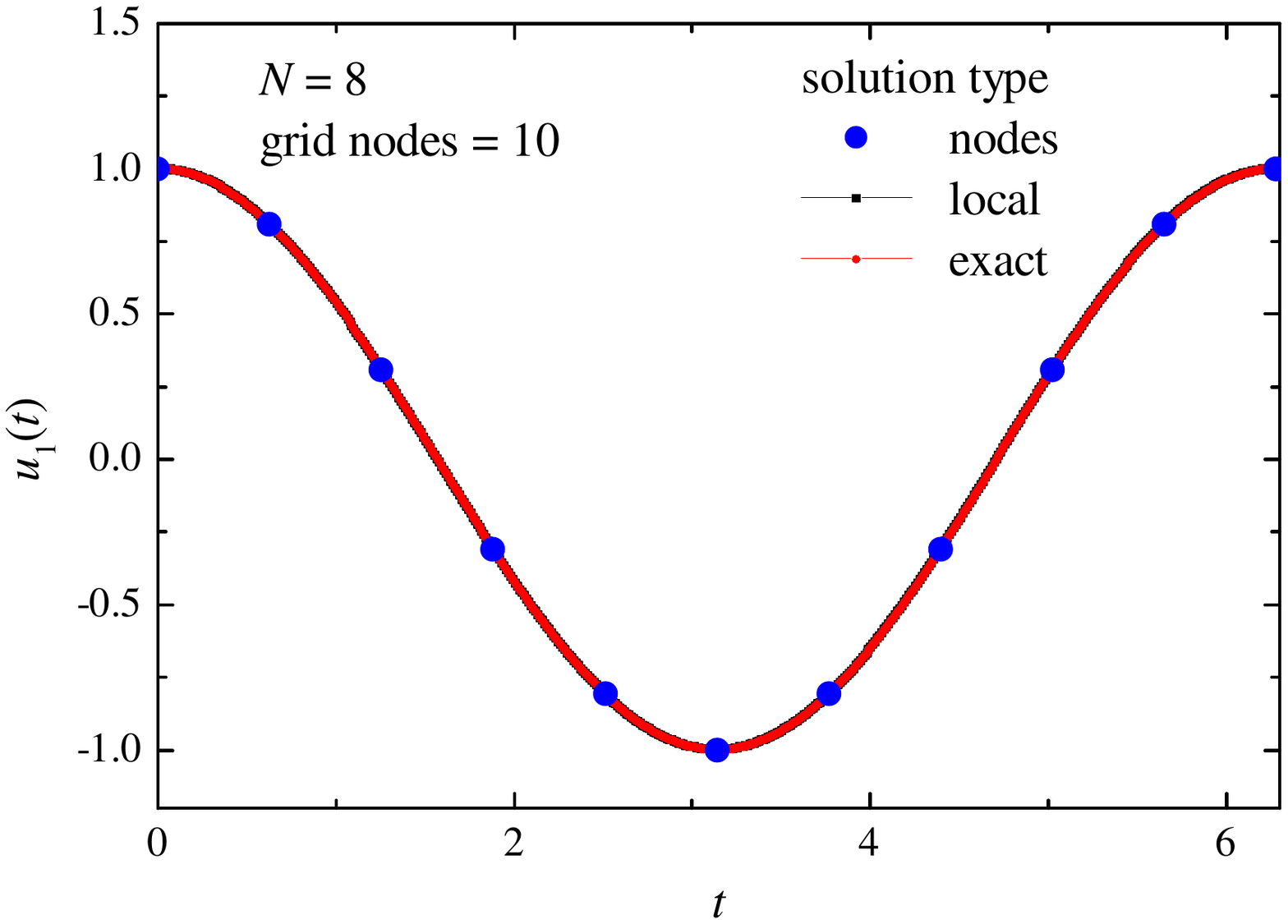}
\vspace{-8mm}\caption{\label{fig:simple_test_sol_uv:a2}}
\end{subfigure}
\begin{subfigure}{0.320\textwidth}
\includegraphics[width=\textwidth]{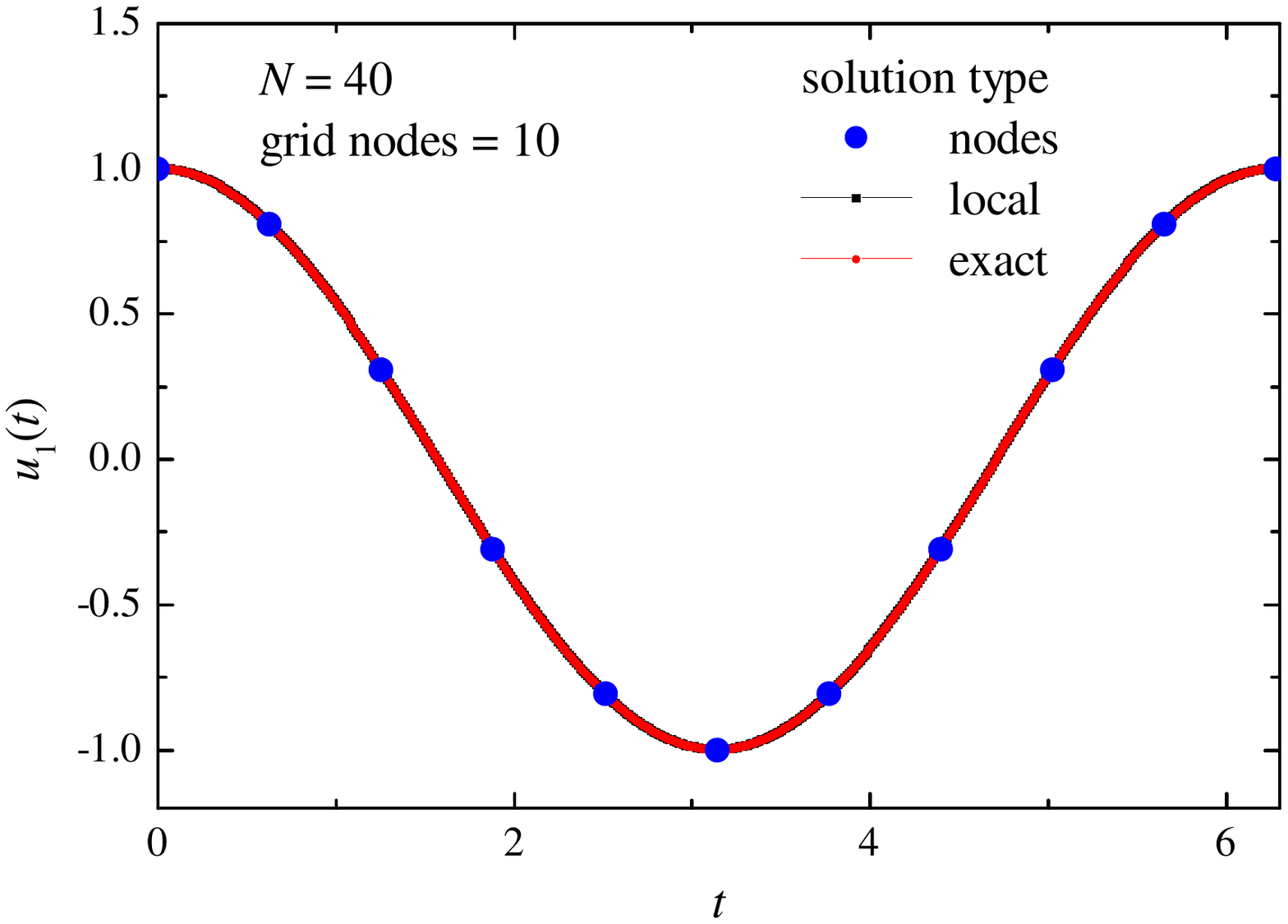}
\vspace{-8mm}\caption{\label{fig:simple_test_sol_uv:a3}}
\end{subfigure}\\
\begin{subfigure}{0.320\textwidth}
\includegraphics[width=\textwidth]{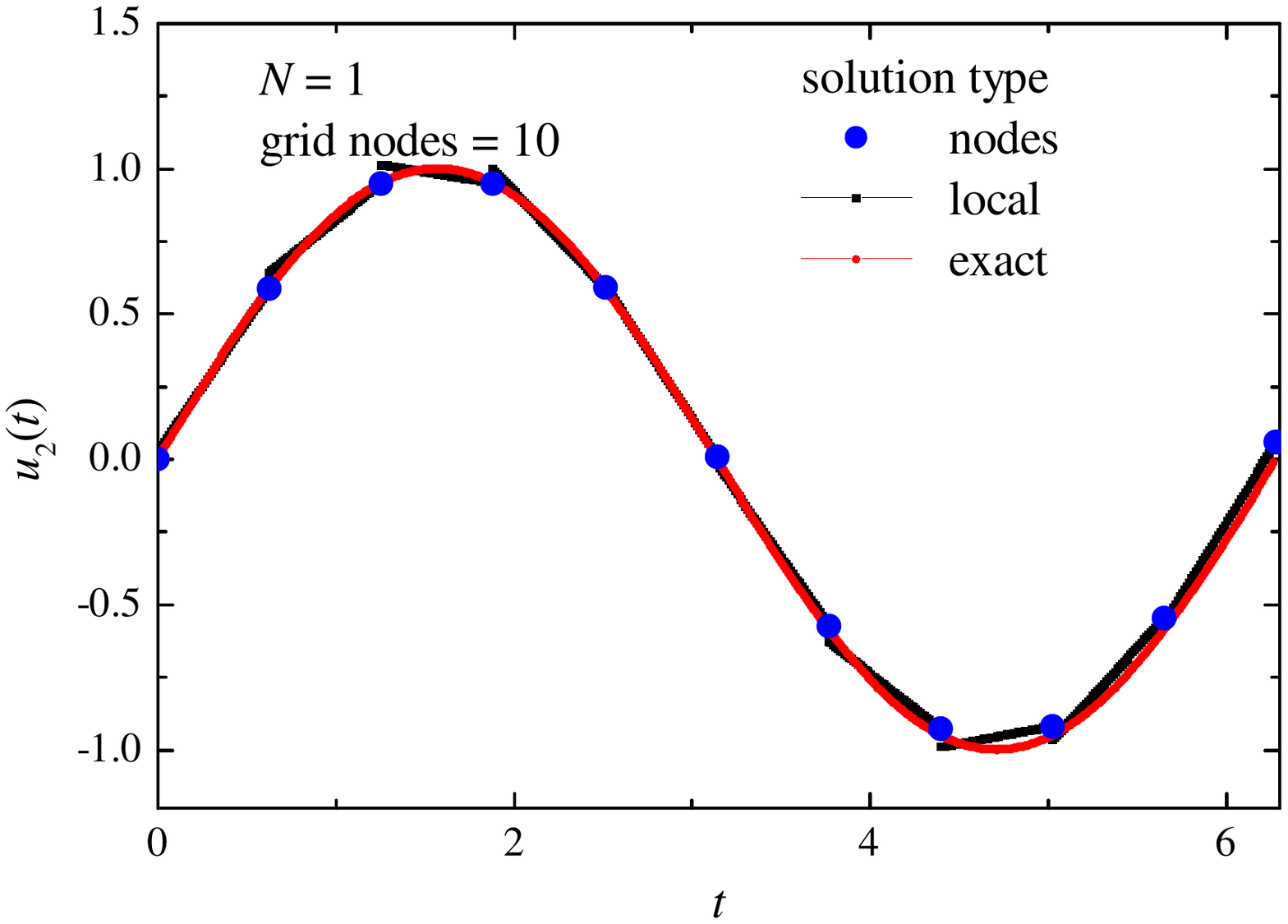}
\vspace{-8mm}\caption{\label{fig:simple_test_sol_uv:b1}}
\end{subfigure}
\begin{subfigure}{0.320\textwidth}
\includegraphics[width=\textwidth]{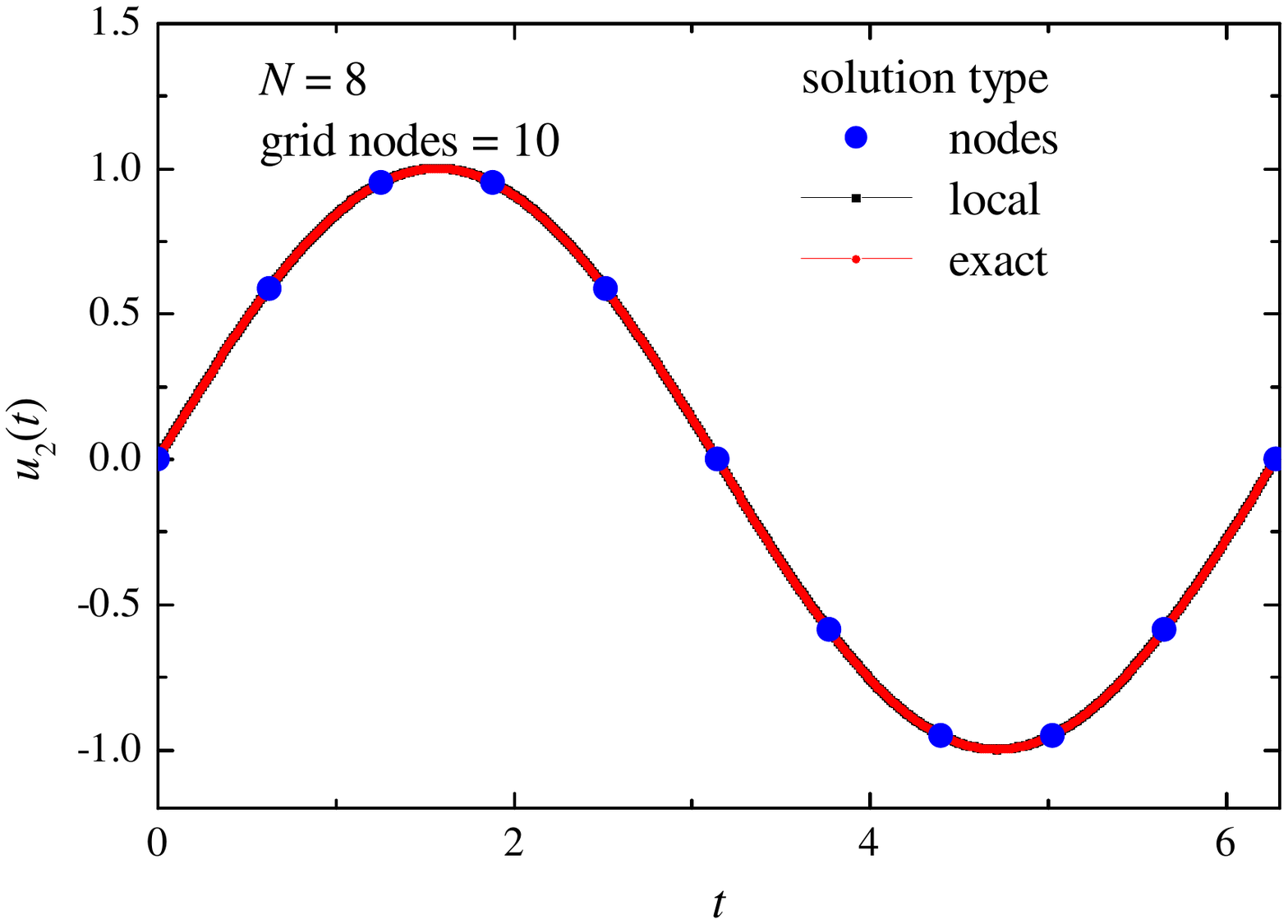}
\vspace{-8mm}\caption{\label{fig:simple_test_sol_uv:b2}}
\end{subfigure}
\begin{subfigure}{0.320\textwidth}
\includegraphics[width=\textwidth]{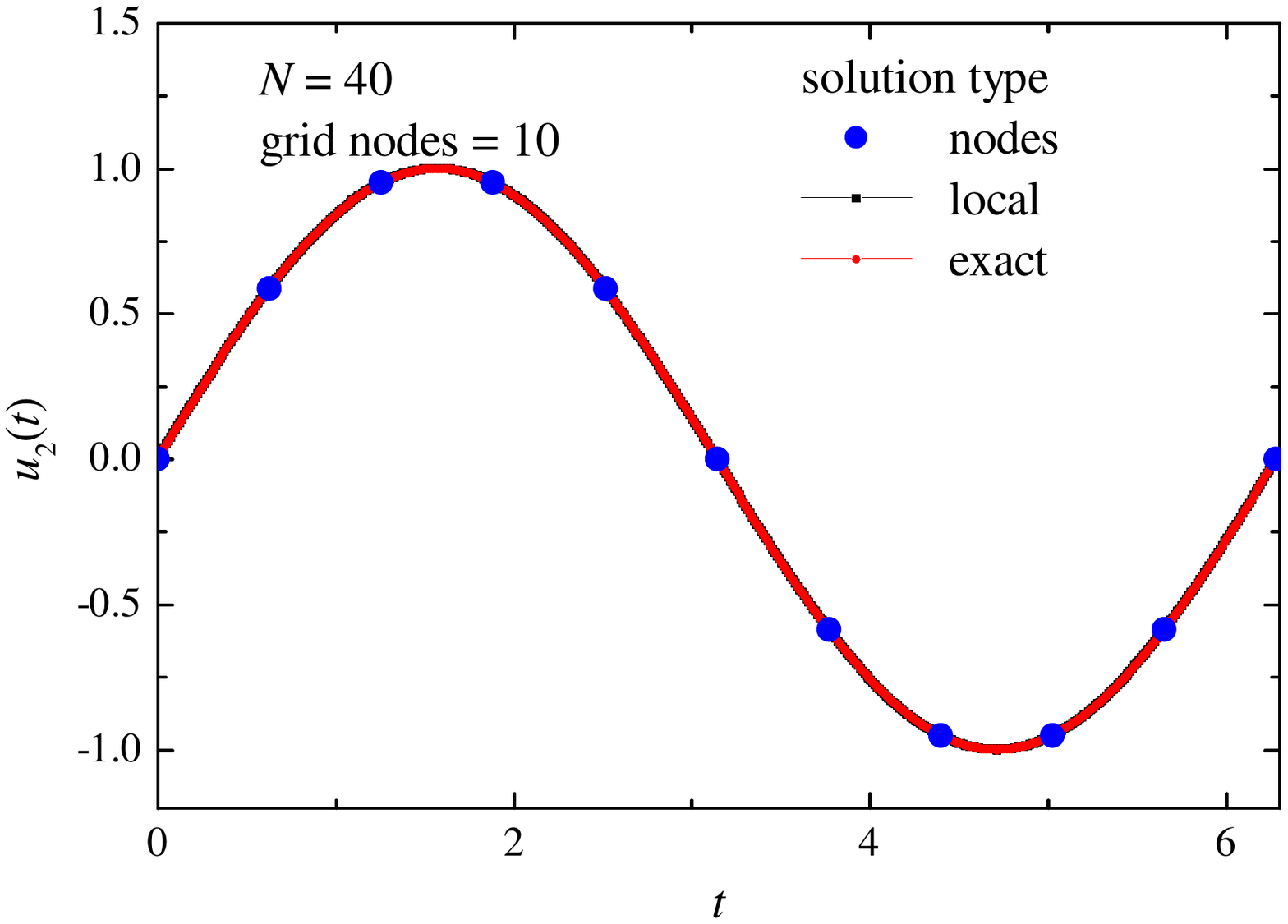}
\vspace{-8mm}\caption{\label{fig:simple_test_sol_uv:b3}}
\end{subfigure}\\
\begin{subfigure}{0.320\textwidth}
\includegraphics[width=\textwidth]{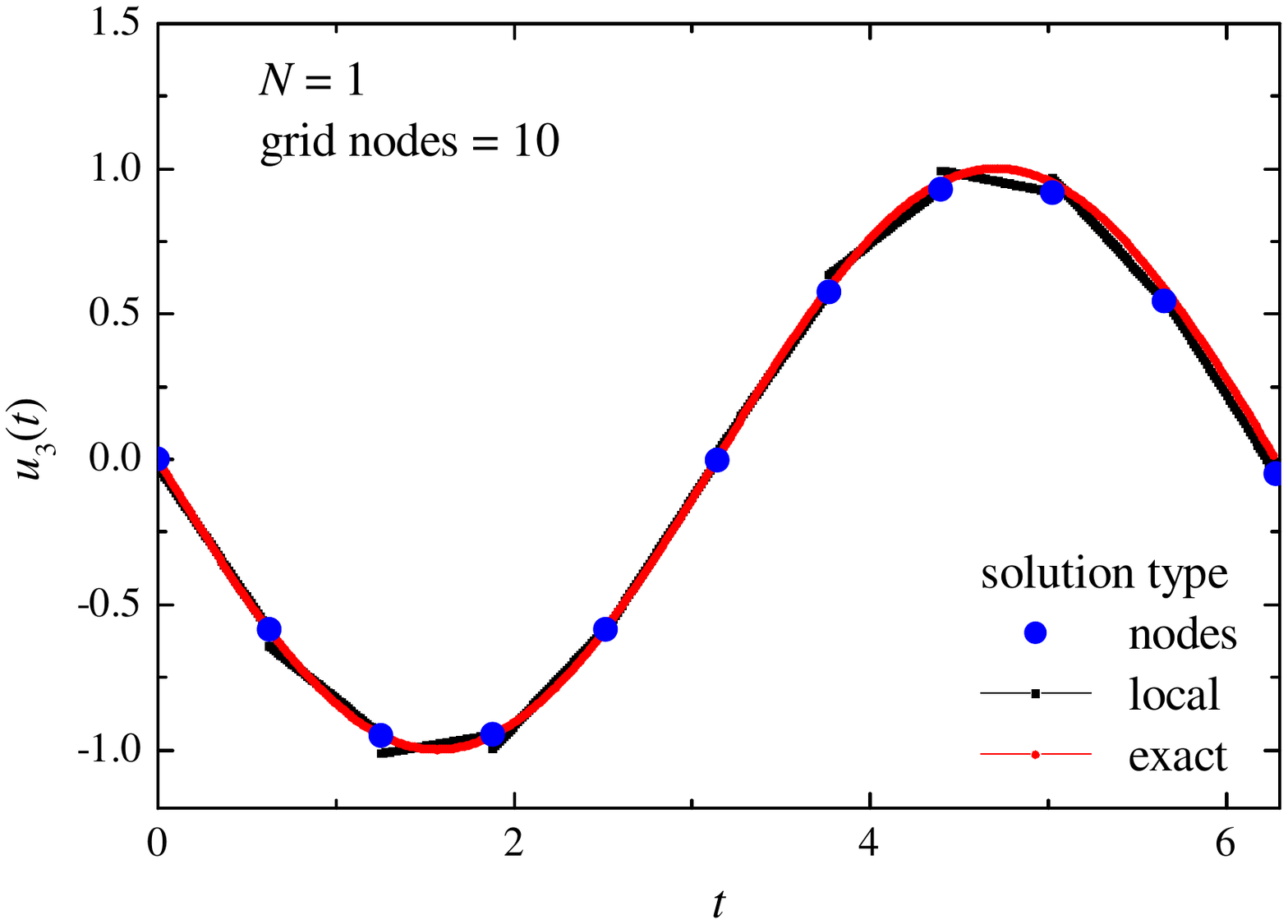}
\vspace{-8mm}\caption{\label{fig:simple_test_sol_uv:c1}}
\end{subfigure}
\begin{subfigure}{0.320\textwidth}
\includegraphics[width=\textwidth]{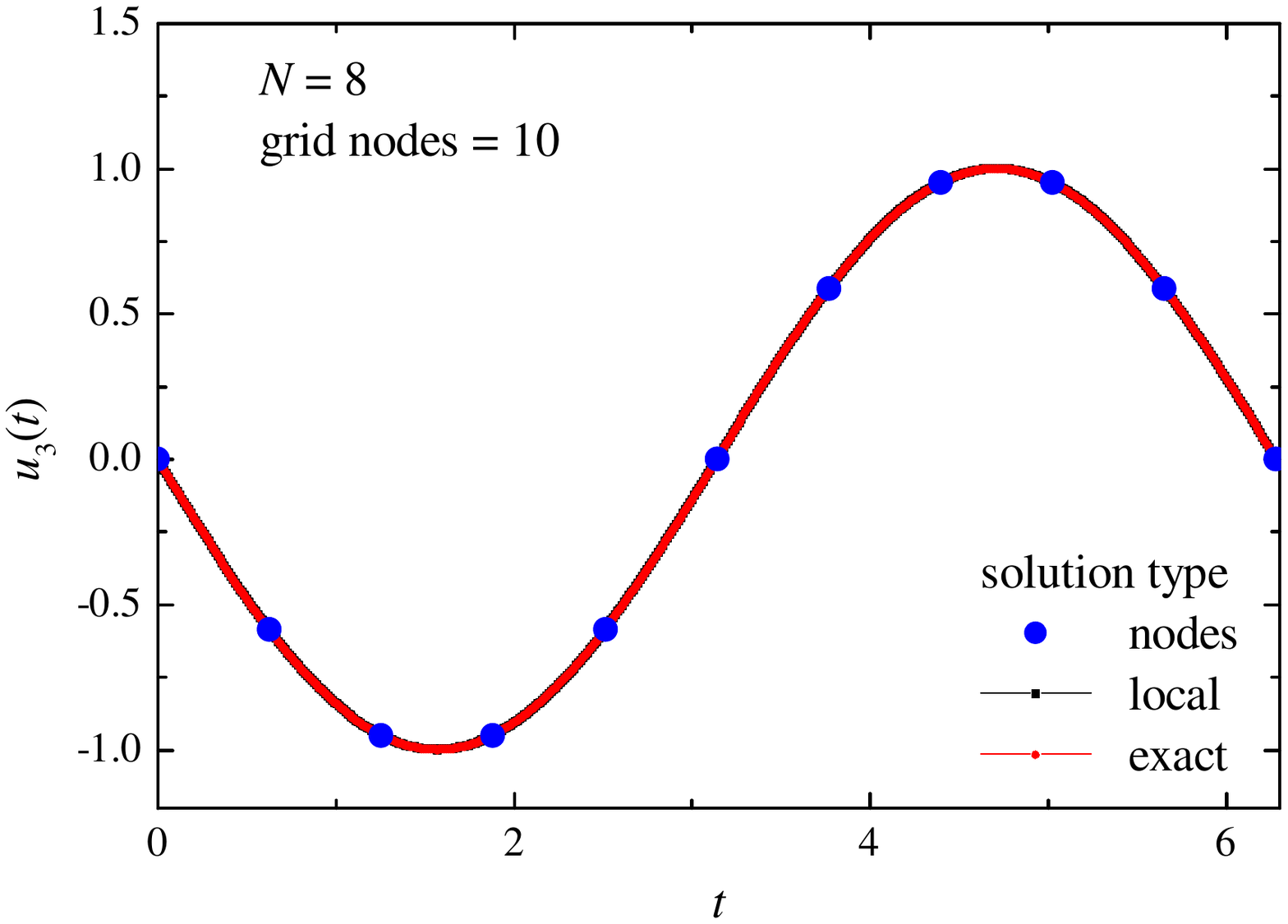}
\vspace{-8mm}\caption{\label{fig:simple_test_sol_uv:c2}}
\end{subfigure}
\begin{subfigure}{0.320\textwidth}
\includegraphics[width=\textwidth]{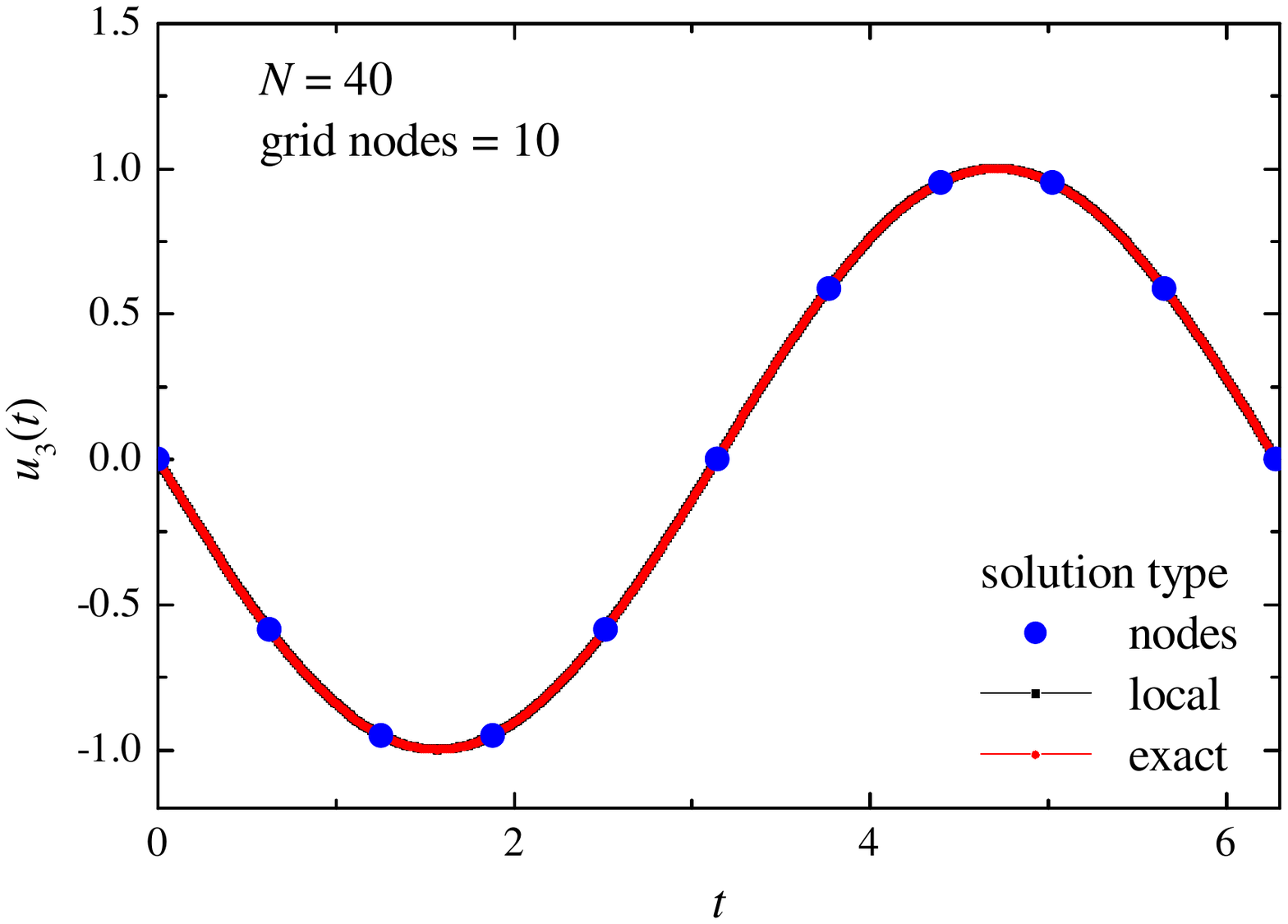}
\vspace{-8mm}\caption{\label{fig:simple_test_sol_uv:c3}}
\end{subfigure}\\
\begin{subfigure}{0.320\textwidth}
\includegraphics[width=\textwidth]{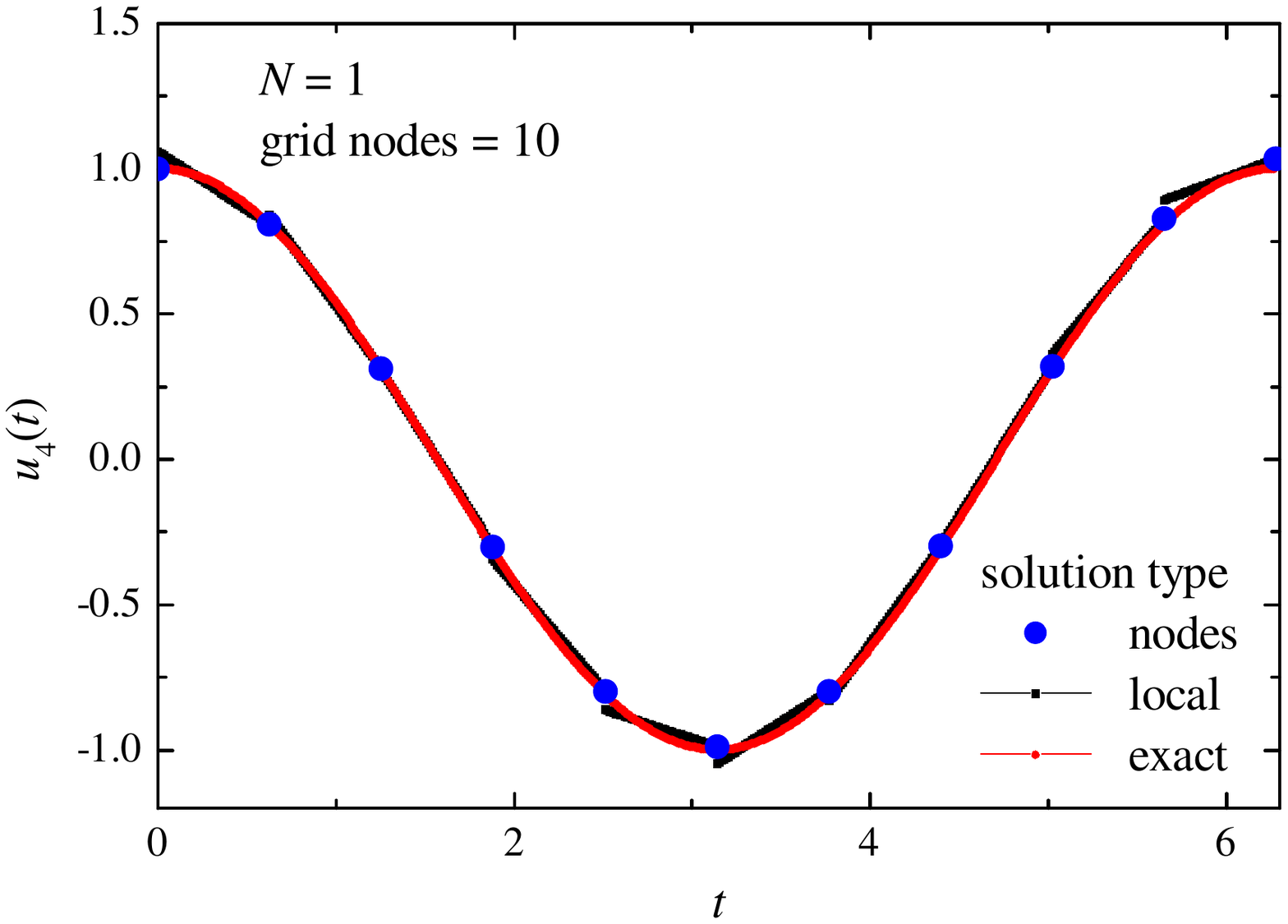}
\vspace{-8mm}\caption{\label{fig:simple_test_sol_uv:d1}}
\end{subfigure}
\begin{subfigure}{0.320\textwidth}
\includegraphics[width=\textwidth]{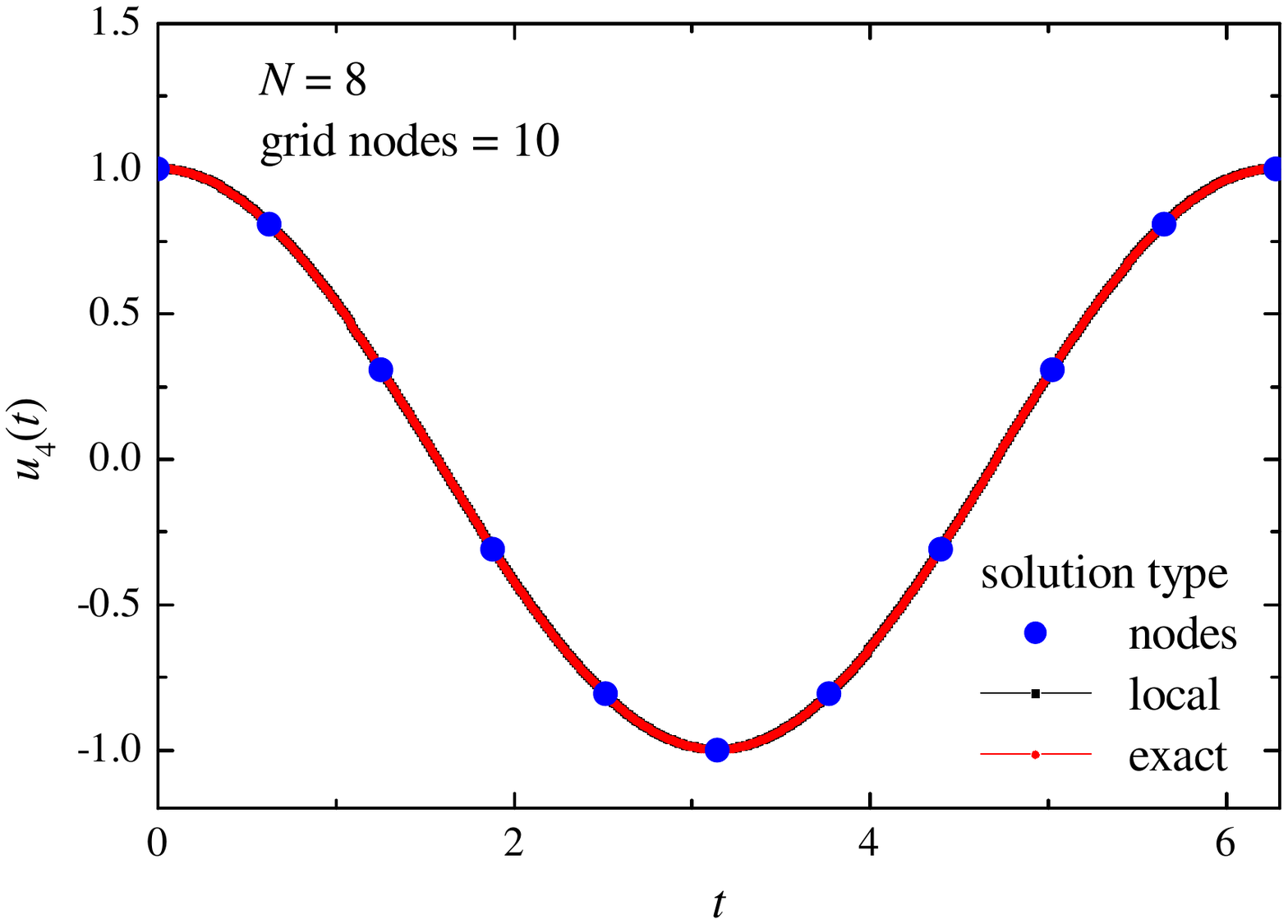}
\vspace{-8mm}\caption{\label{fig:simple_test_sol_uv:d2}}
\end{subfigure}
\begin{subfigure}{0.320\textwidth}
\includegraphics[width=\textwidth]{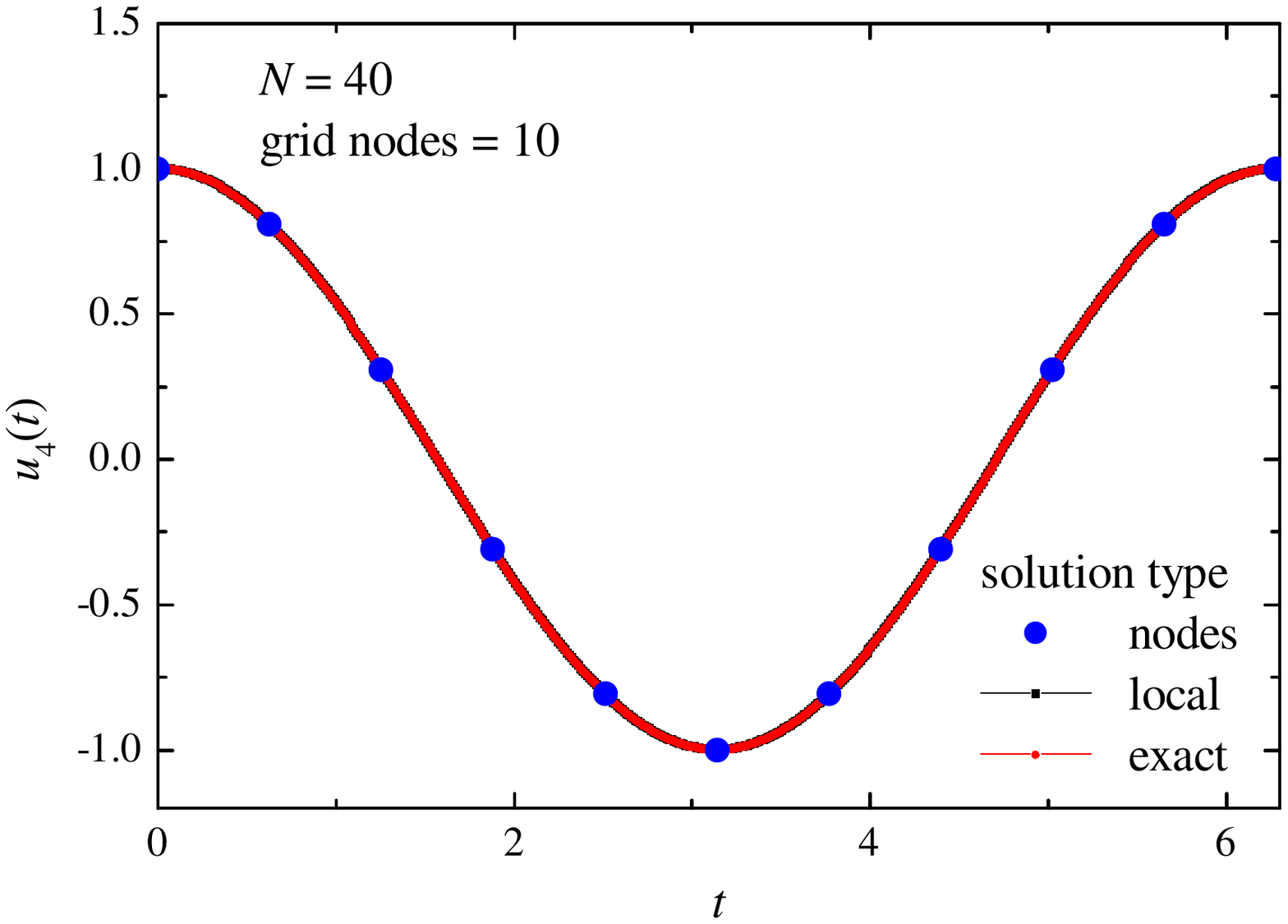}
\vspace{-8mm}\caption{\label{fig:simple_test_sol_uv:d3}}
\end{subfigure}\\
\begin{subfigure}{0.320\textwidth}
\includegraphics[width=\textwidth]{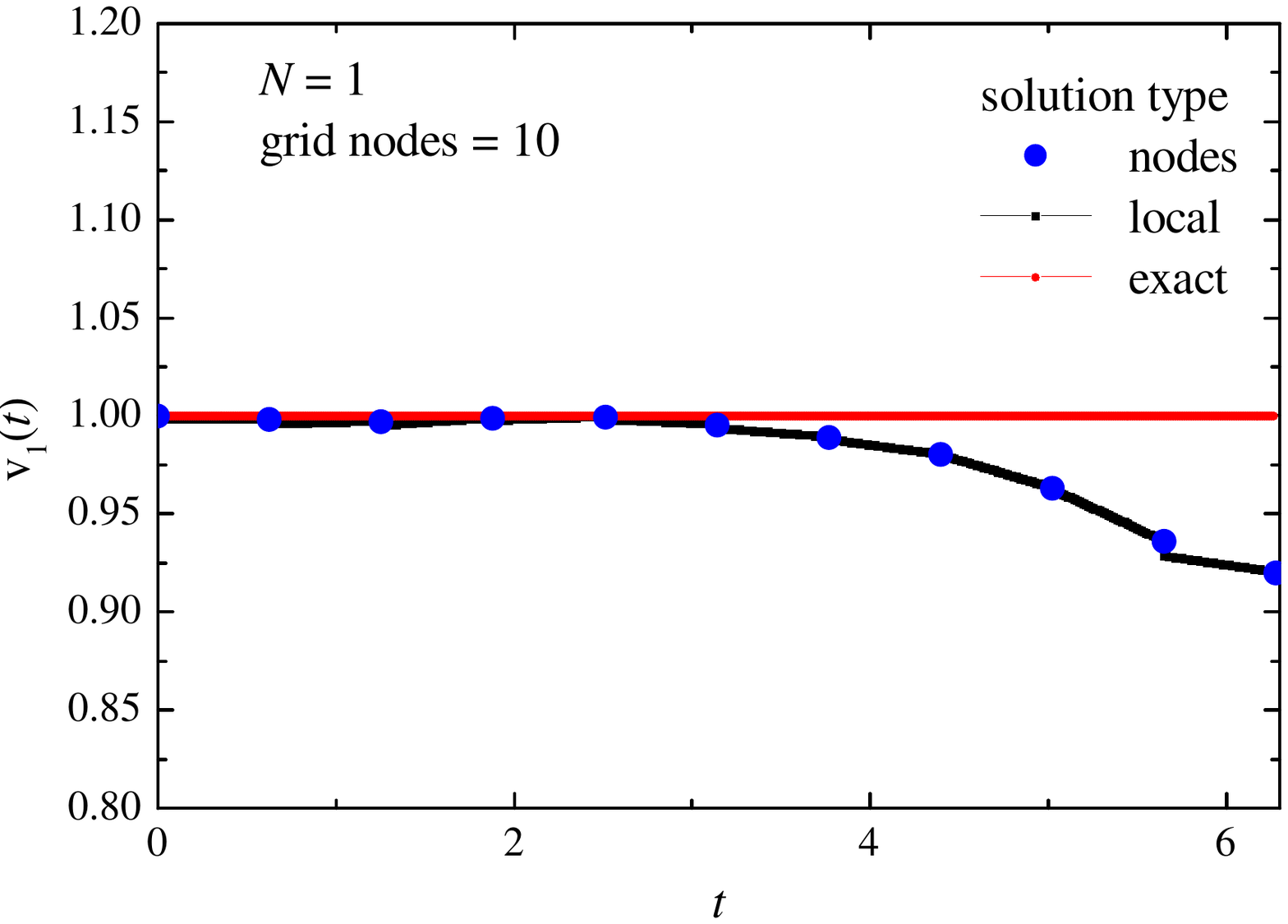}
\vspace{-8mm}\caption{\label{fig:simple_test_sol_uv:e1}}
\end{subfigure}
\begin{subfigure}{0.320\textwidth}
\includegraphics[width=\textwidth]{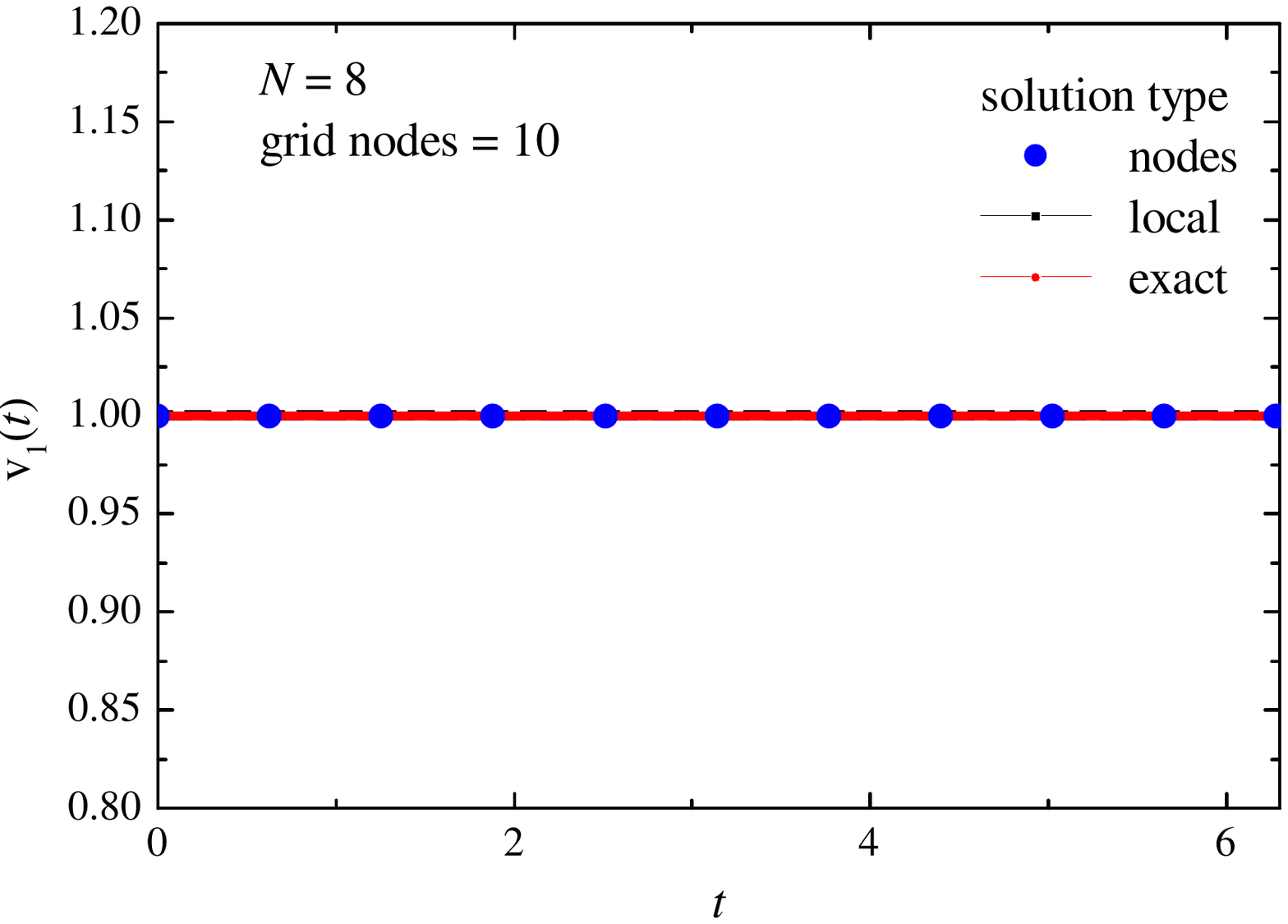}
\vspace{-8mm}\caption{\label{fig:simple_test_sol_uv:e2}}
\end{subfigure}
\begin{subfigure}{0.320\textwidth}
\includegraphics[width=\textwidth]{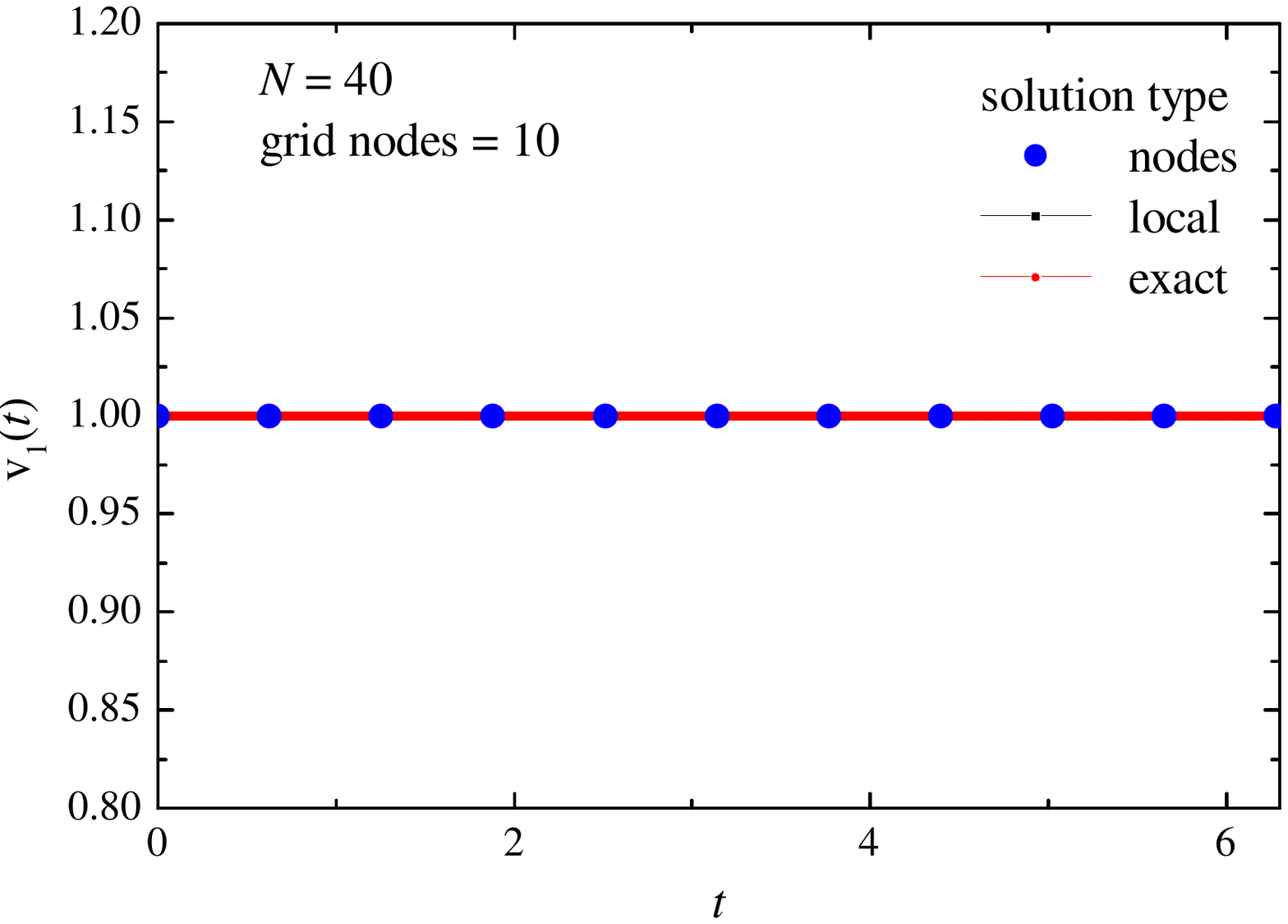}
\vspace{-8mm}\caption{\label{fig:simple_test_sol_uv:e3}}
\end{subfigure}\\
\caption{%
Numerical solution of the problem (\ref{eq:simple_dae_ind_1}). Comparison of the solution at nodes $\mathbf{u}_{n}$, the local solution $\mathbf{u}_{L}(t)$ and the exact solution $\mathbf{u}^{\rm ex}(t)$ for components $u_{1}$ (\subref{fig:simple_test_sol_uv:a1}, \subref{fig:simple_test_sol_uv:a2}, \subref{fig:simple_test_sol_uv:a3}), $u_{2}$ (\subref{fig:simple_test_sol_uv:b1}, \subref{fig:simple_test_sol_uv:b2}, \subref{fig:simple_test_sol_uv:b3}), $u_{3}$ (\subref{fig:simple_test_sol_uv:c1}, \subref{fig:simple_test_sol_uv:c2}, \subref{fig:simple_test_sol_uv:c3}), $u_{4}$ (\subref{fig:simple_test_sol_uv:d1}, \subref{fig:simple_test_sol_uv:d2}, \subref{fig:simple_test_sol_uv:d3}) and $v_{1}$ (\subref{fig:simple_test_sol_uv:e1}, \subref{fig:simple_test_sol_uv:e2}, \subref{fig:simple_test_sol_uv:e3}), obtained using polynomials with degrees $N = 1$ (\subref{fig:simple_test_sol_uv:a1}, \subref{fig:simple_test_sol_uv:b1}, \subref{fig:simple_test_sol_uv:c1}, \subref{fig:simple_test_sol_uv:d1}, \subref{fig:simple_test_sol_uv:e1}), $N = 8$ (\subref{fig:simple_test_sol_uv:a2}, \subref{fig:simple_test_sol_uv:b2}, \subref{fig:simple_test_sol_uv:c2}, \subref{fig:simple_test_sol_uv:d2}, \subref{fig:simple_test_sol_uv:e2}) and $N = 40$ (\subref{fig:simple_test_sol_uv:a3}, \subref{fig:simple_test_sol_uv:b3}, \subref{fig:simple_test_sol_uv:c3}, \subref{fig:simple_test_sol_uv:d3}, \subref{fig:simple_test_sol_uv:e3}).
}
\label{fig:simple_test_sols_uv}
\end{figure}

\begin{figure}[h!]
\captionsetup[subfigure]{%
	position=bottom,
	font+=smaller,
	textfont=normalfont,
	singlelinecheck=off,
	justification=raggedright
}
\centering
\begin{subfigure}{0.320\textwidth}
\includegraphics[width=\textwidth]{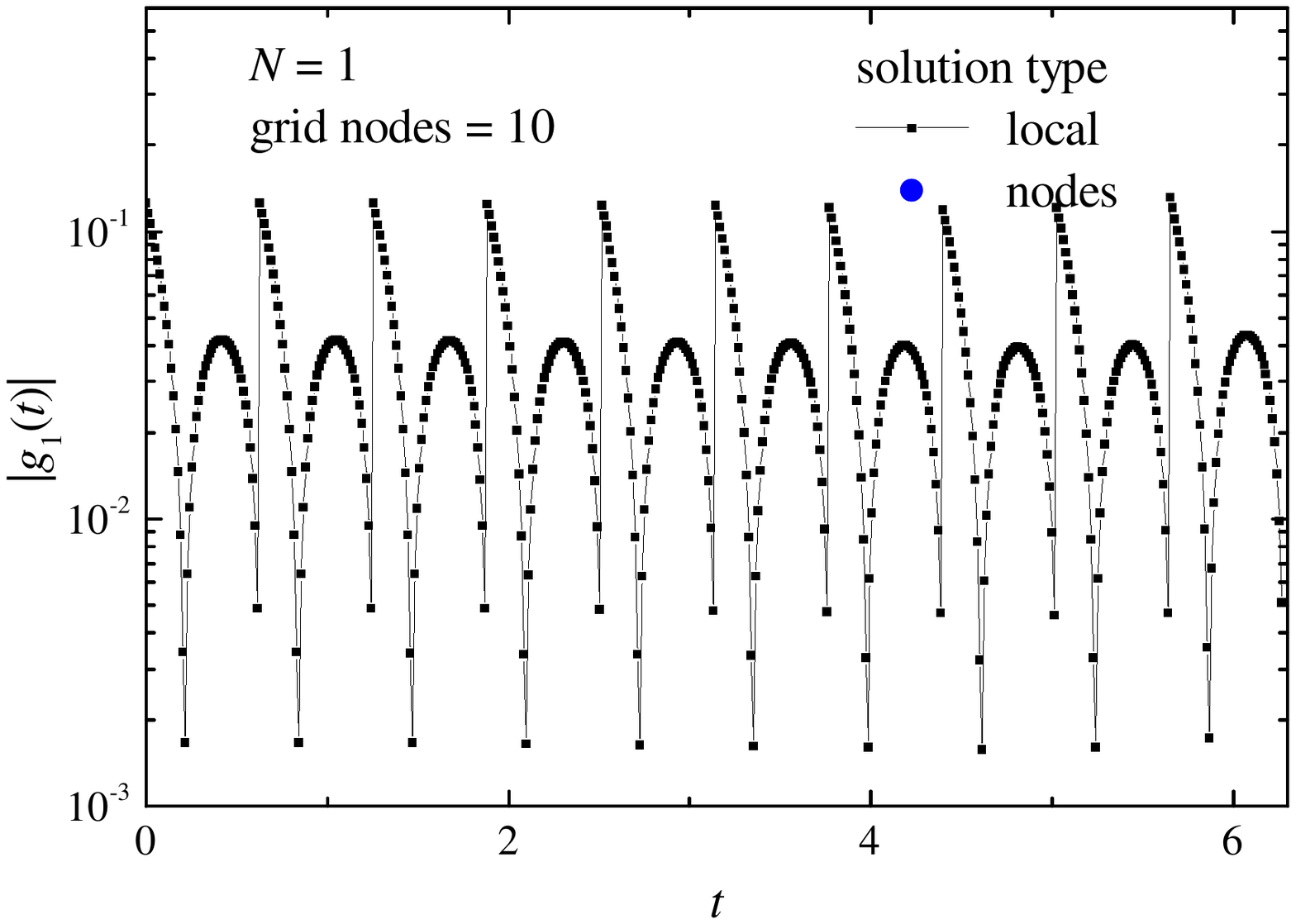}
\vspace{-8mm}\caption{\label{fig:simple_test_sol_g_eps:a1}}
\end{subfigure}
\begin{subfigure}{0.320\textwidth}
\includegraphics[width=\textwidth]{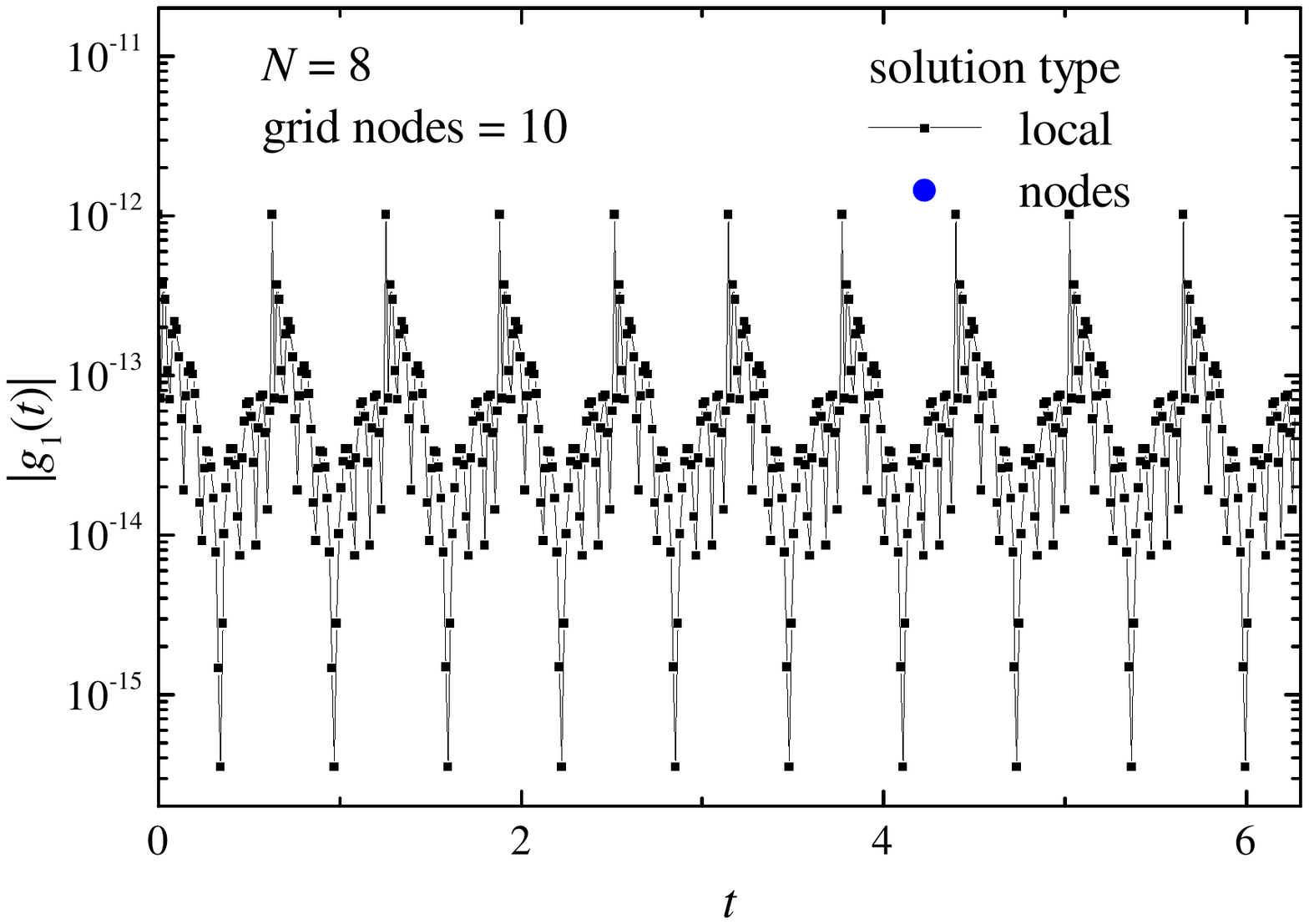}
\vspace{-8mm}\caption{\label{fig:simple_test_sol_g_eps:a2}}
\end{subfigure}
\begin{subfigure}{0.320\textwidth}
\includegraphics[width=\textwidth]{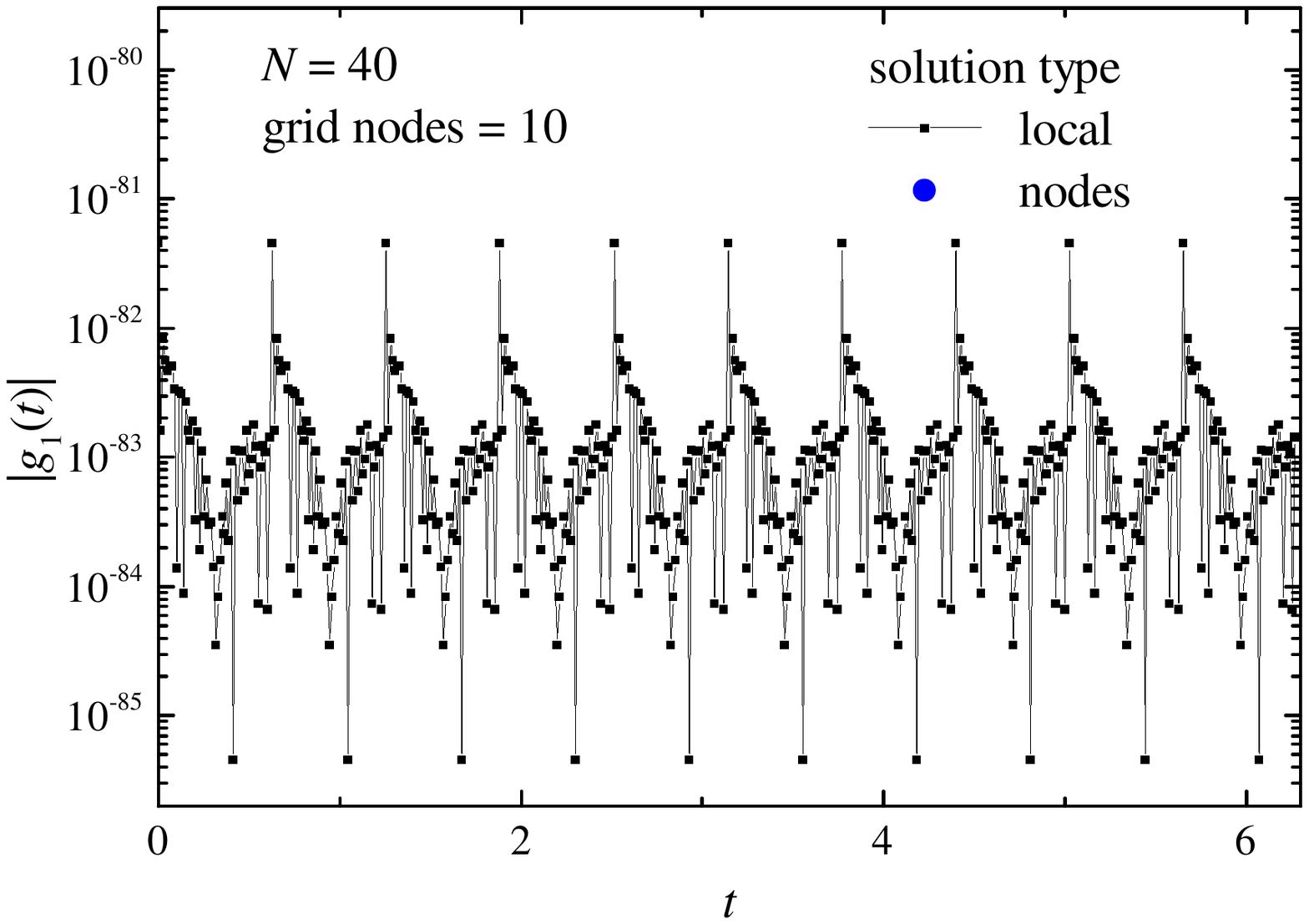}
\vspace{-8mm}\caption{\label{fig:simple_test_sol_g_eps:a3}}
\end{subfigure}\\
\begin{subfigure}{0.320\textwidth}
\includegraphics[width=\textwidth]{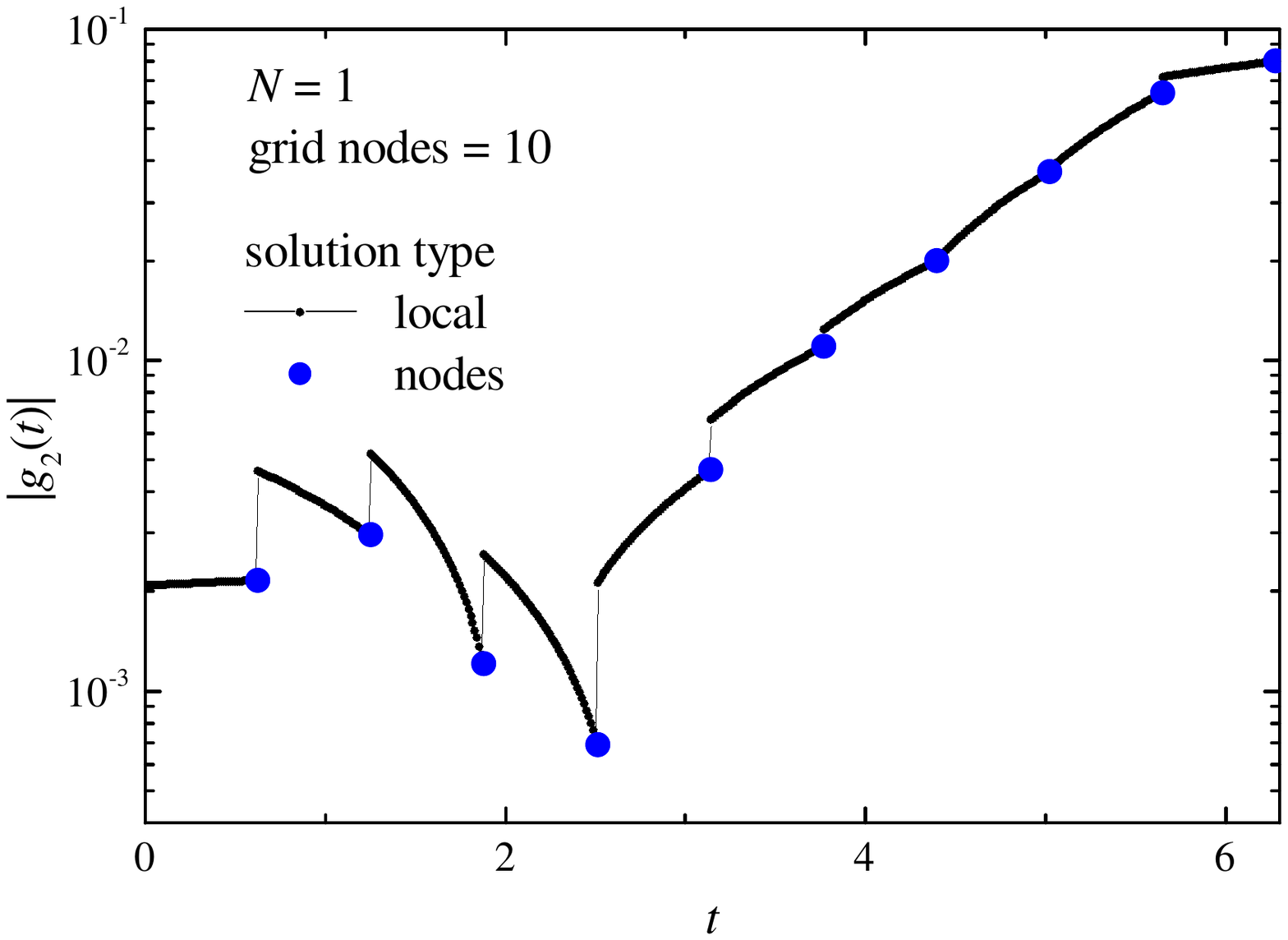}
\vspace{-8mm}\caption{\label{fig:simple_test_sol_g_eps:b1}}
\end{subfigure}
\begin{subfigure}{0.320\textwidth}
\includegraphics[width=\textwidth]{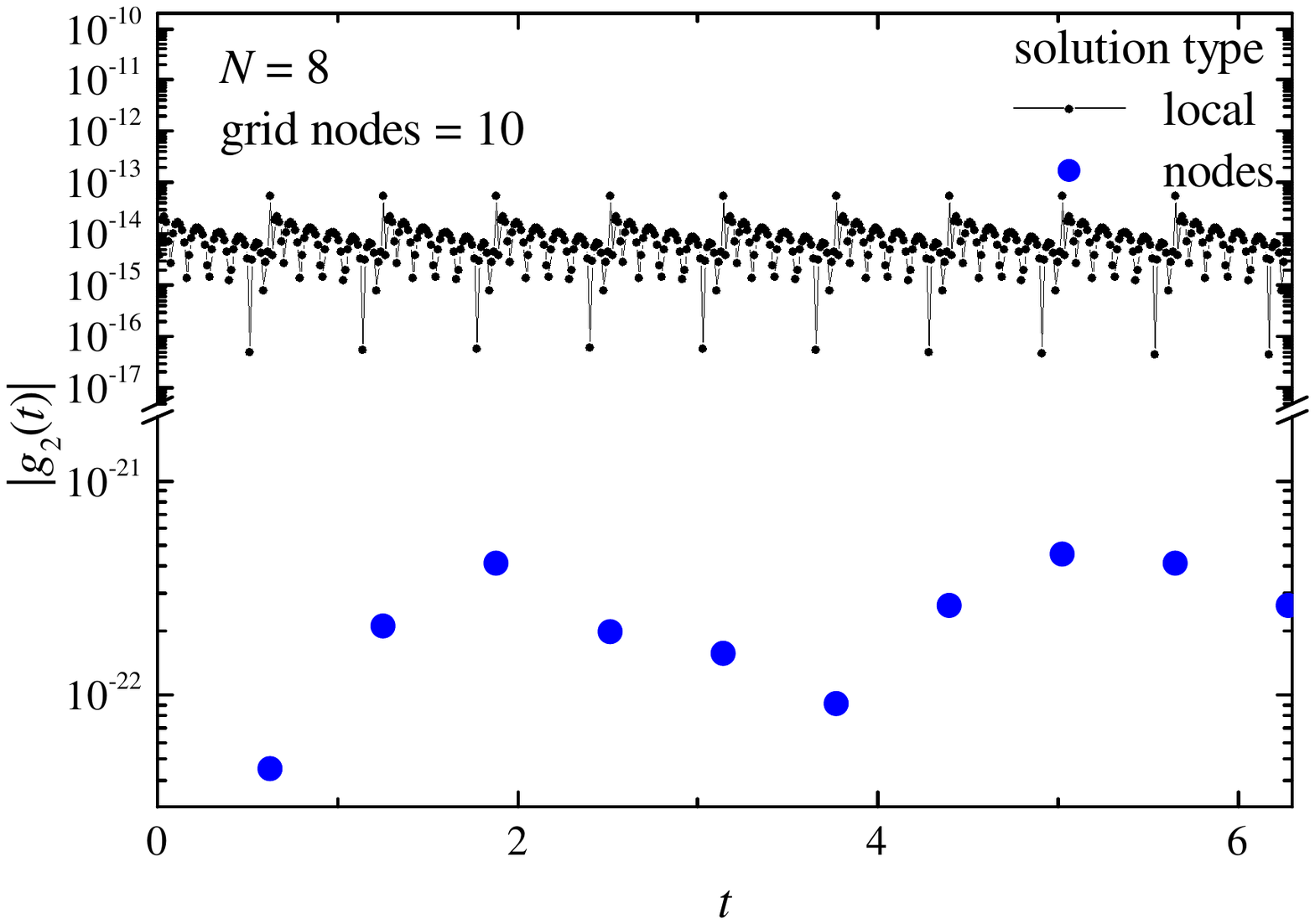}
\vspace{-8mm}\caption{\label{fig:simple_test_sol_g_eps:b2}}
\end{subfigure}
\begin{subfigure}{0.320\textwidth}
\includegraphics[width=\textwidth]{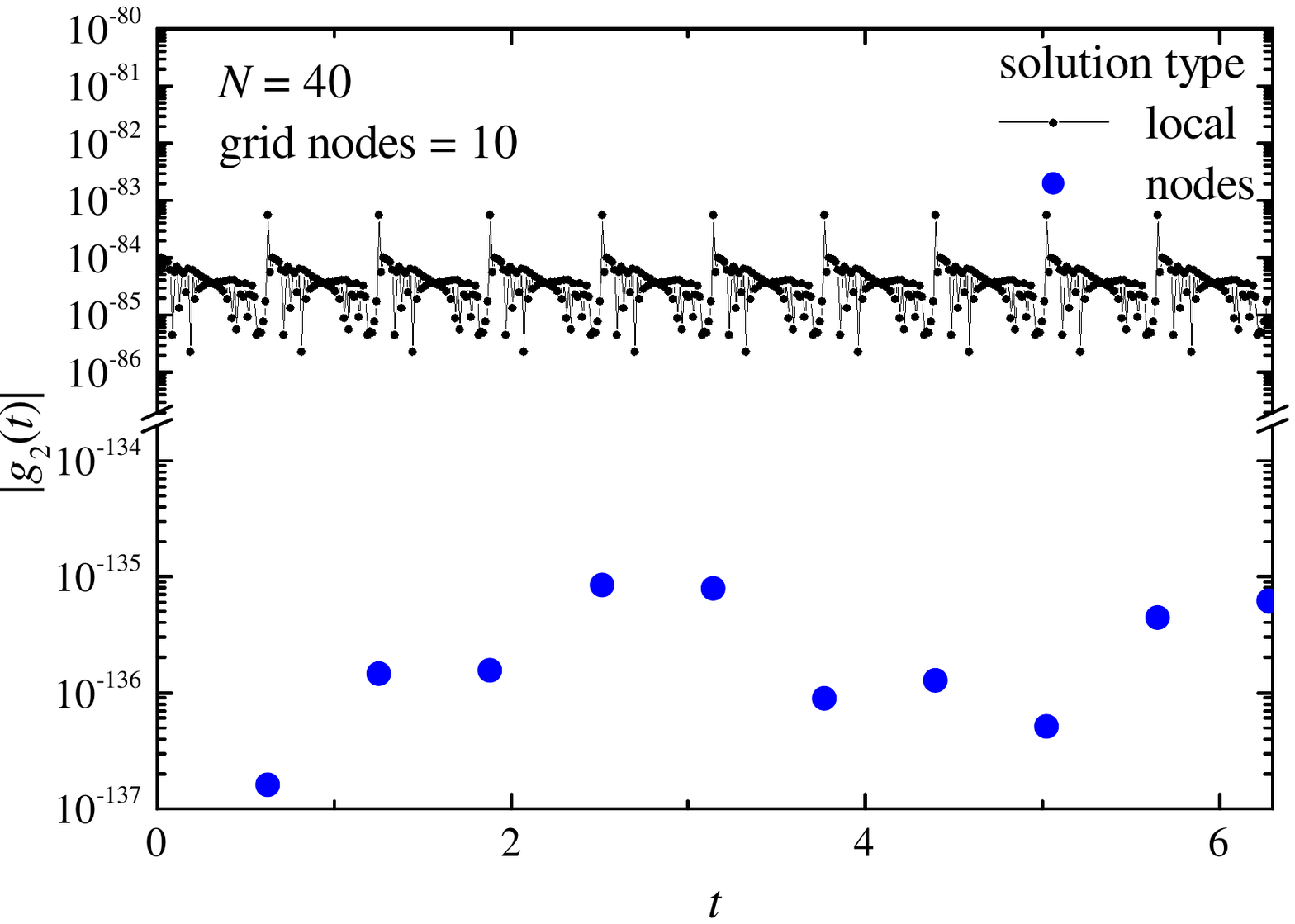}
\vspace{-8mm}\caption{\label{fig:simple_test_sol_g_eps:b3}}
\end{subfigure}\\
\begin{subfigure}{0.320\textwidth}
\includegraphics[width=\textwidth]{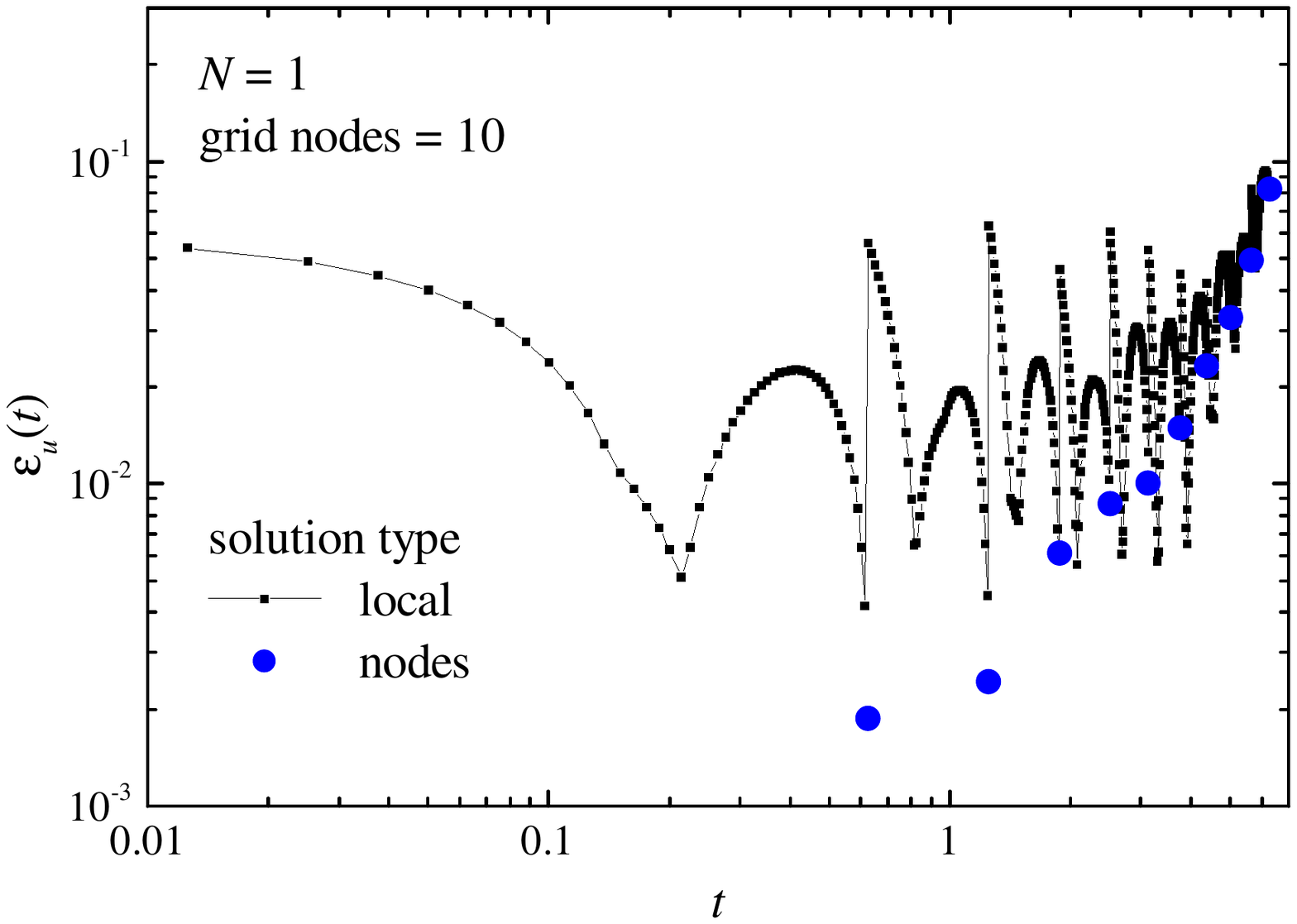}
\vspace{-8mm}\caption{\label{fig:simple_test_sol_g_eps:c1}}
\end{subfigure}
\begin{subfigure}{0.320\textwidth}
\includegraphics[width=\textwidth]{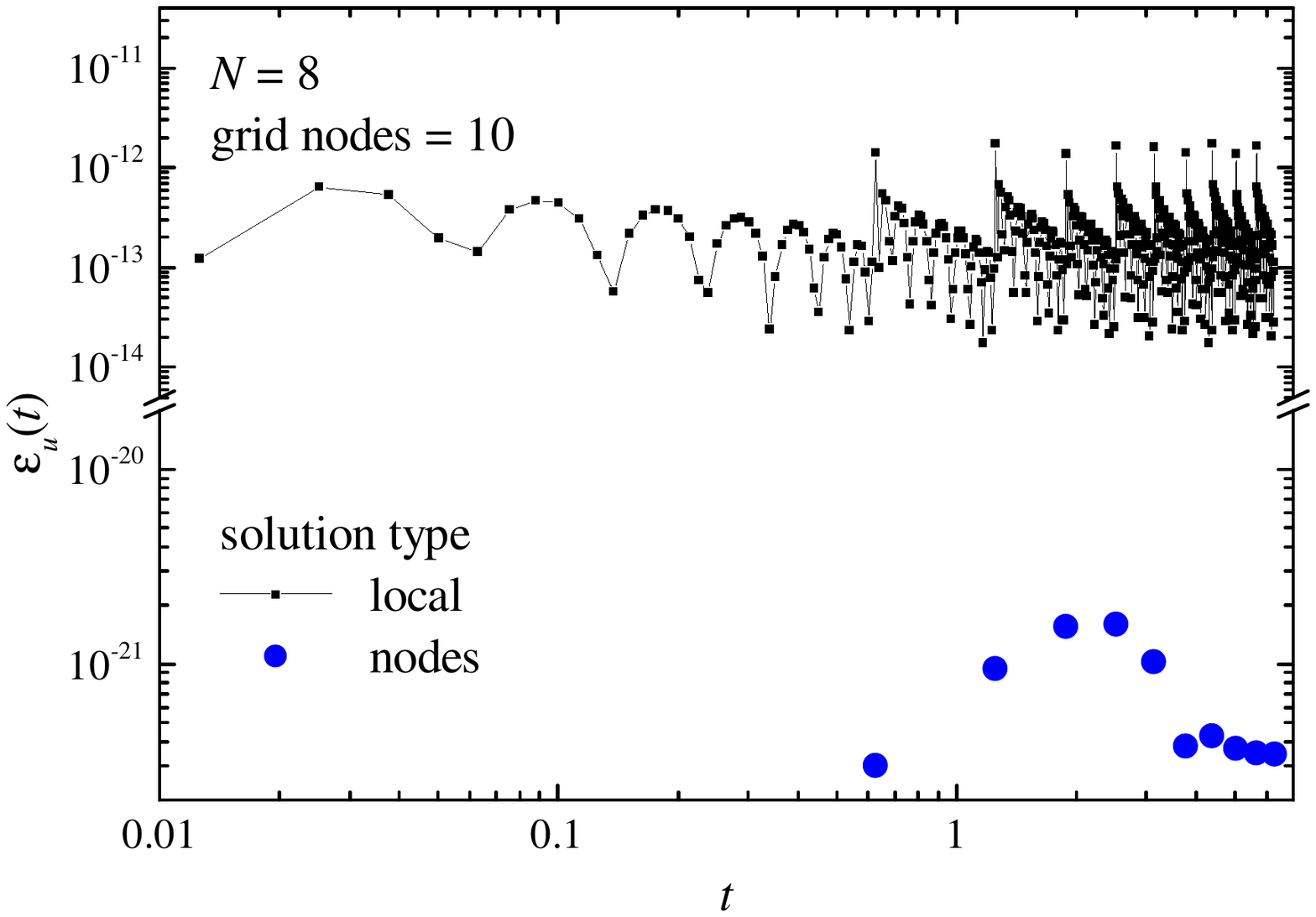}
\vspace{-8mm}\caption{\label{fig:simple_test_sol_g_eps:c2}}
\end{subfigure}
\begin{subfigure}{0.320\textwidth}
\includegraphics[width=\textwidth]{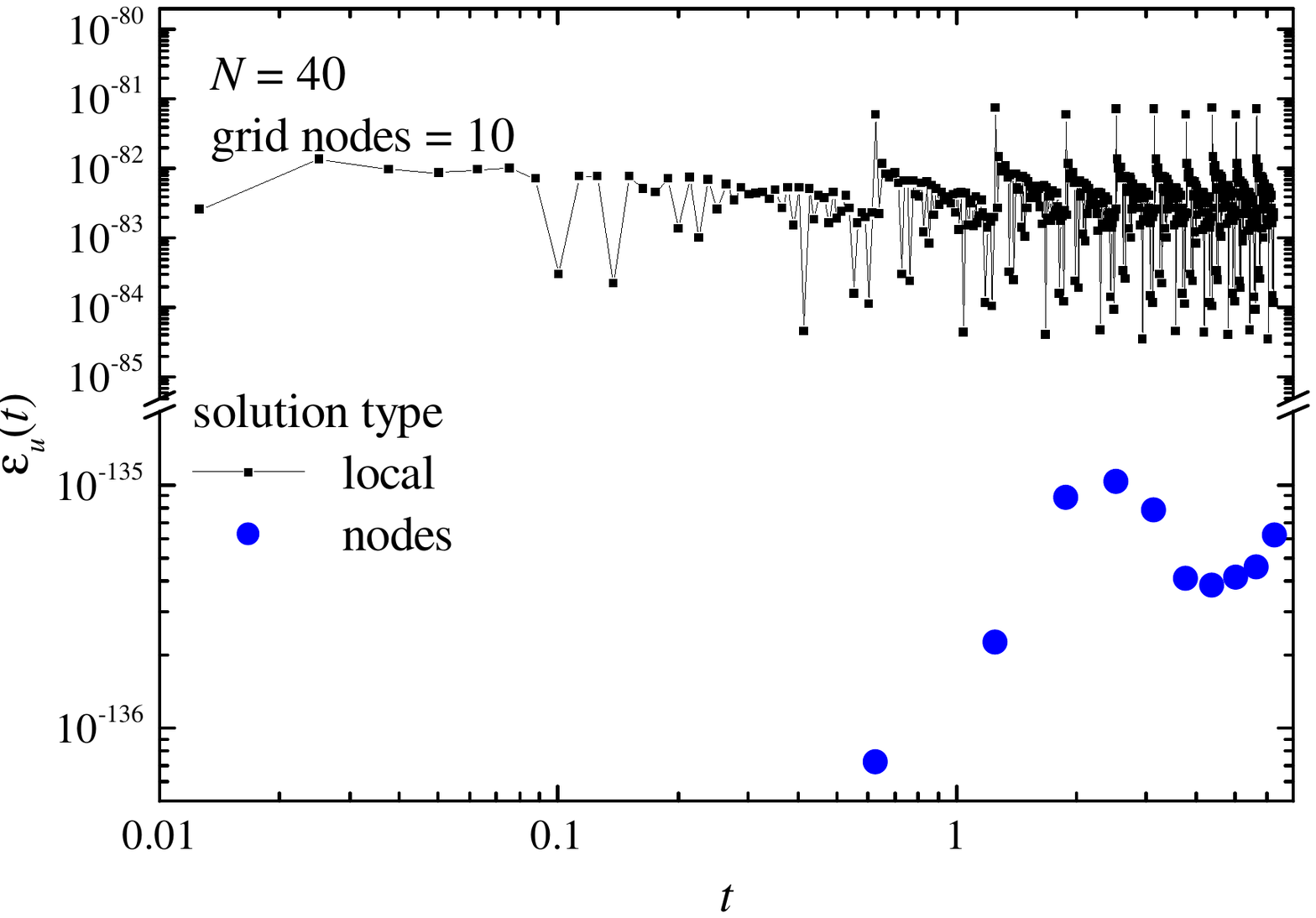}
\vspace{-8mm}\caption{\label{fig:simple_test_sol_g_eps:c3}}
\end{subfigure}\\
\begin{subfigure}{0.320\textwidth}
\includegraphics[width=\textwidth]{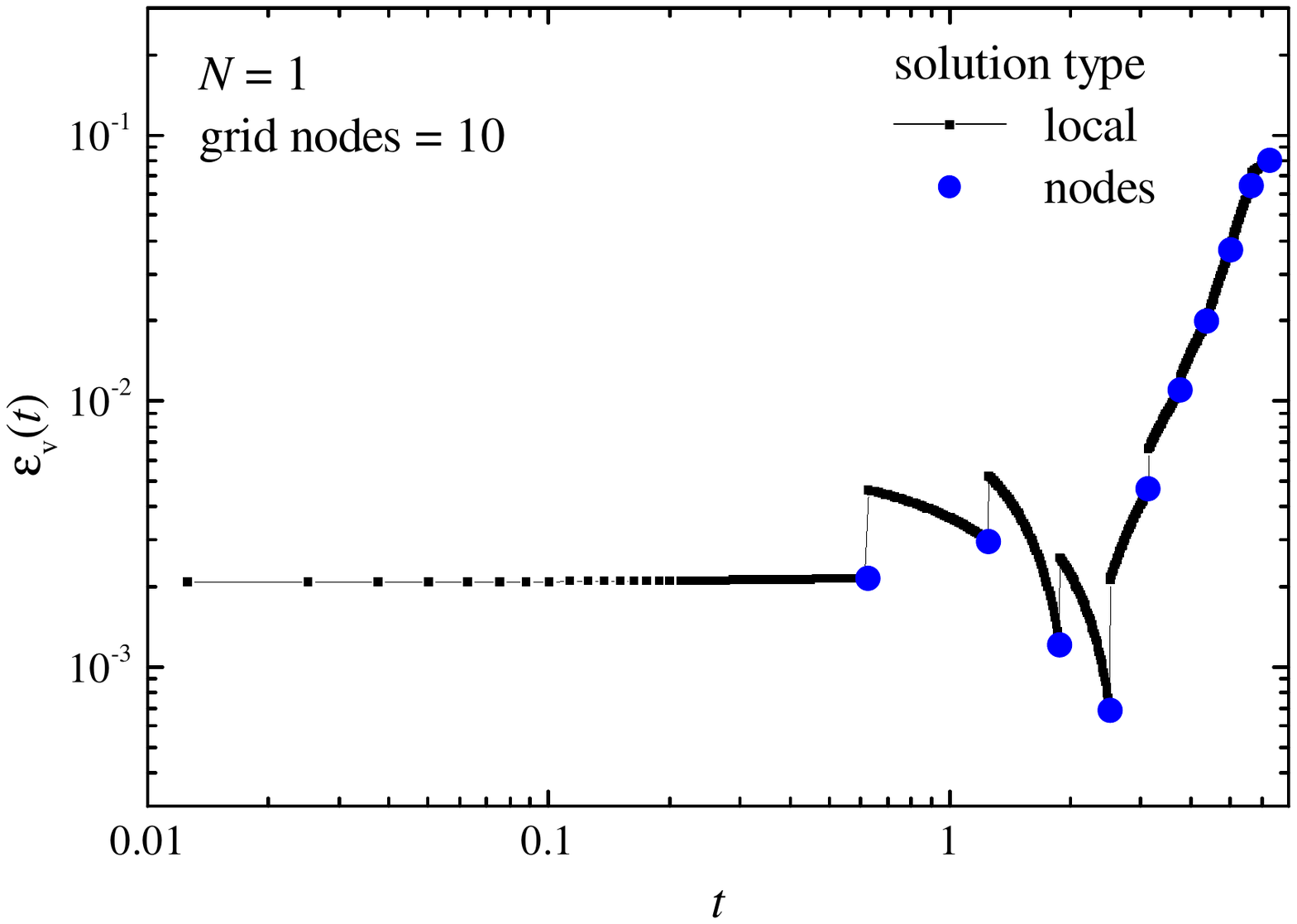}
\vspace{-8mm}\caption{\label{fig:simple_test_sol_g_eps:d1}}
\end{subfigure}
\begin{subfigure}{0.320\textwidth}
\includegraphics[width=\textwidth]{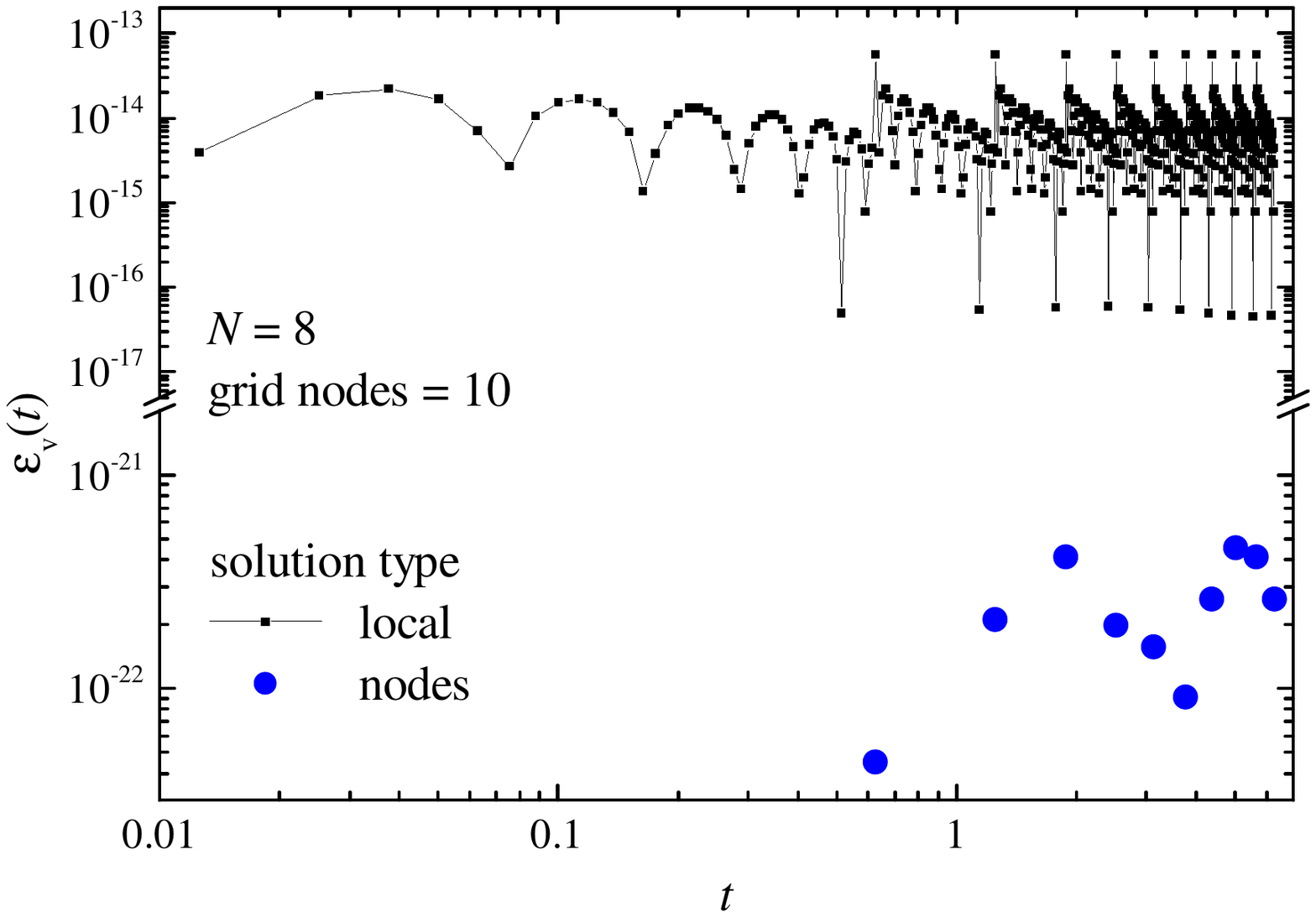}
\vspace{-8mm}\caption{\label{fig:simple_test_sol_g_eps:d2}}
\end{subfigure}
\begin{subfigure}{0.320\textwidth}
\includegraphics[width=\textwidth]{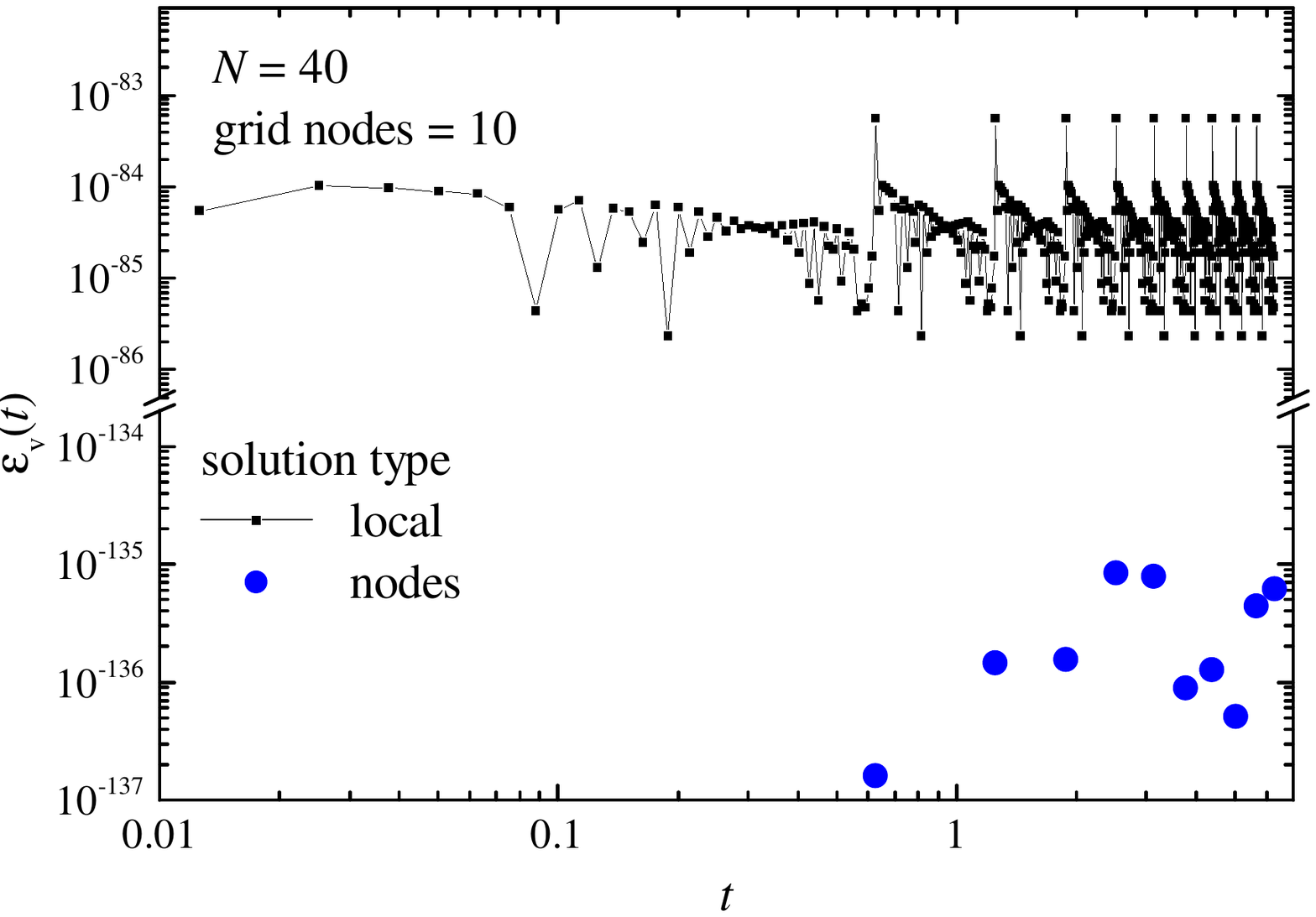}
\vspace{-8mm}\caption{\label{fig:simple_test_sol_g_eps:d3}}
\end{subfigure}\\
\begin{subfigure}{0.320\textwidth}
\includegraphics[width=\textwidth]{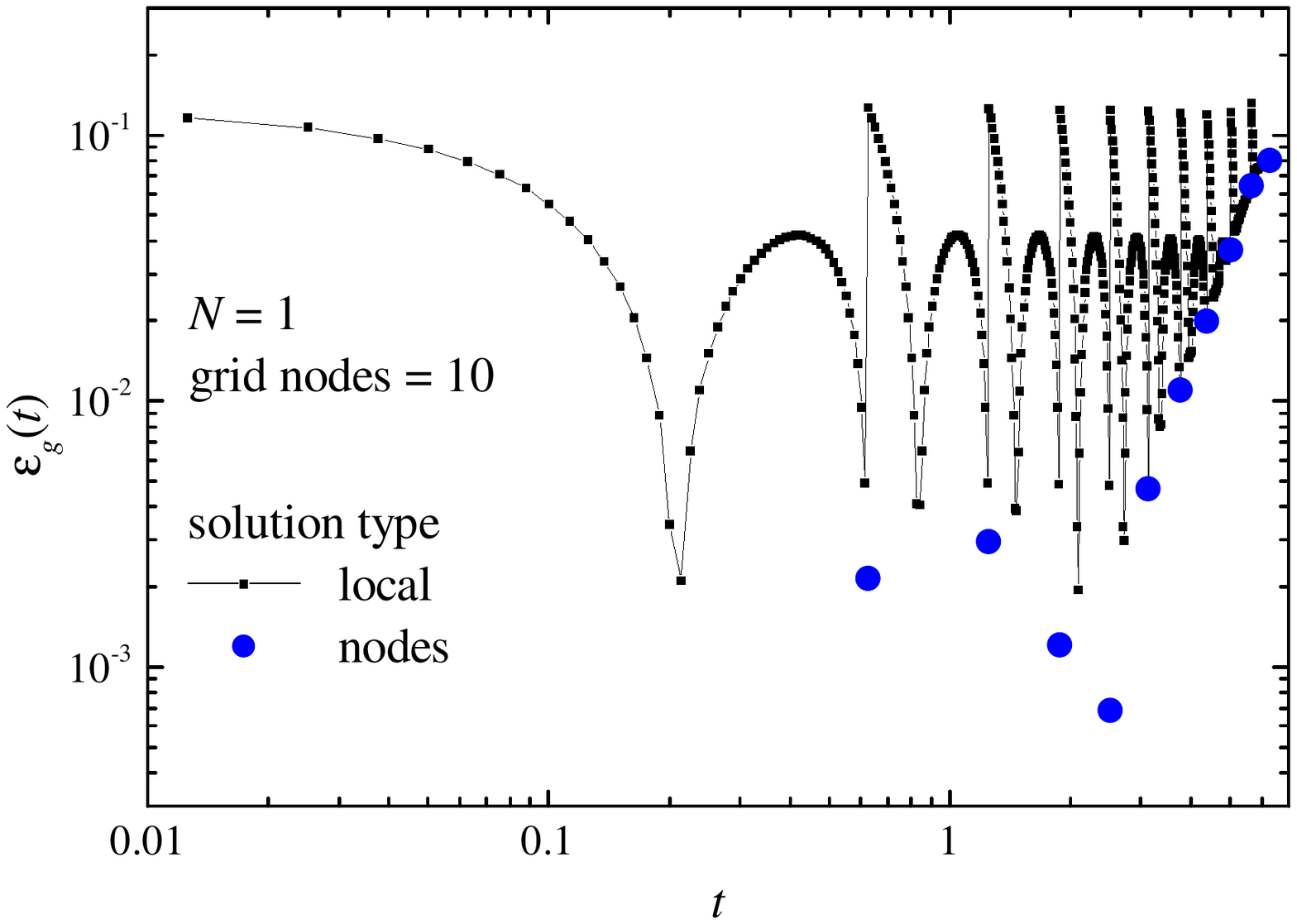}
\vspace{-8mm}\caption{\label{fig:simple_test_sol_g_eps:e1}}
\end{subfigure}
\begin{subfigure}{0.320\textwidth}
\includegraphics[width=\textwidth]{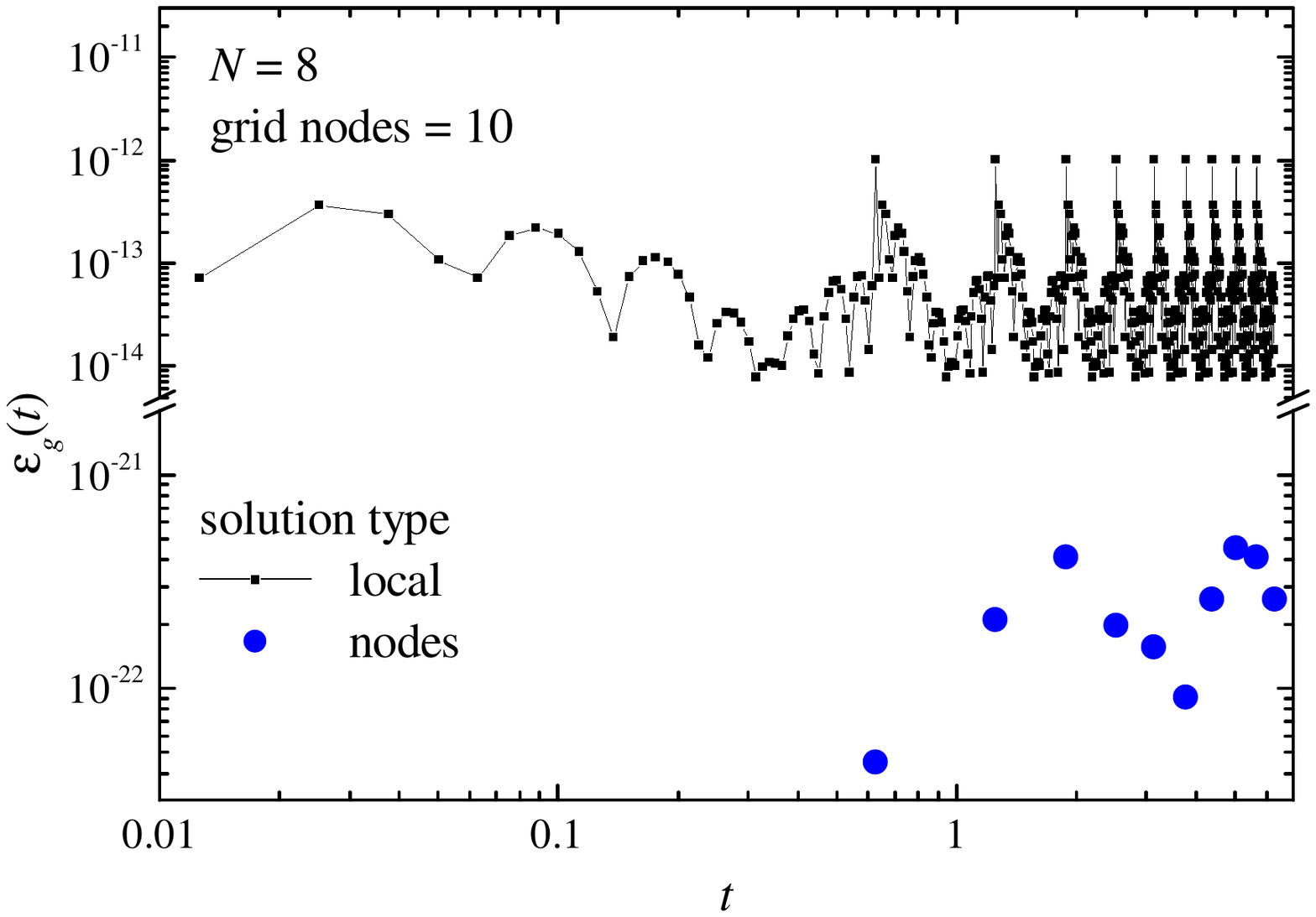}
\vspace{-8mm}\caption{\label{fig:simple_test_sol_g_eps:e2}}
\end{subfigure}
\begin{subfigure}{0.320\textwidth}
\includegraphics[width=\textwidth]{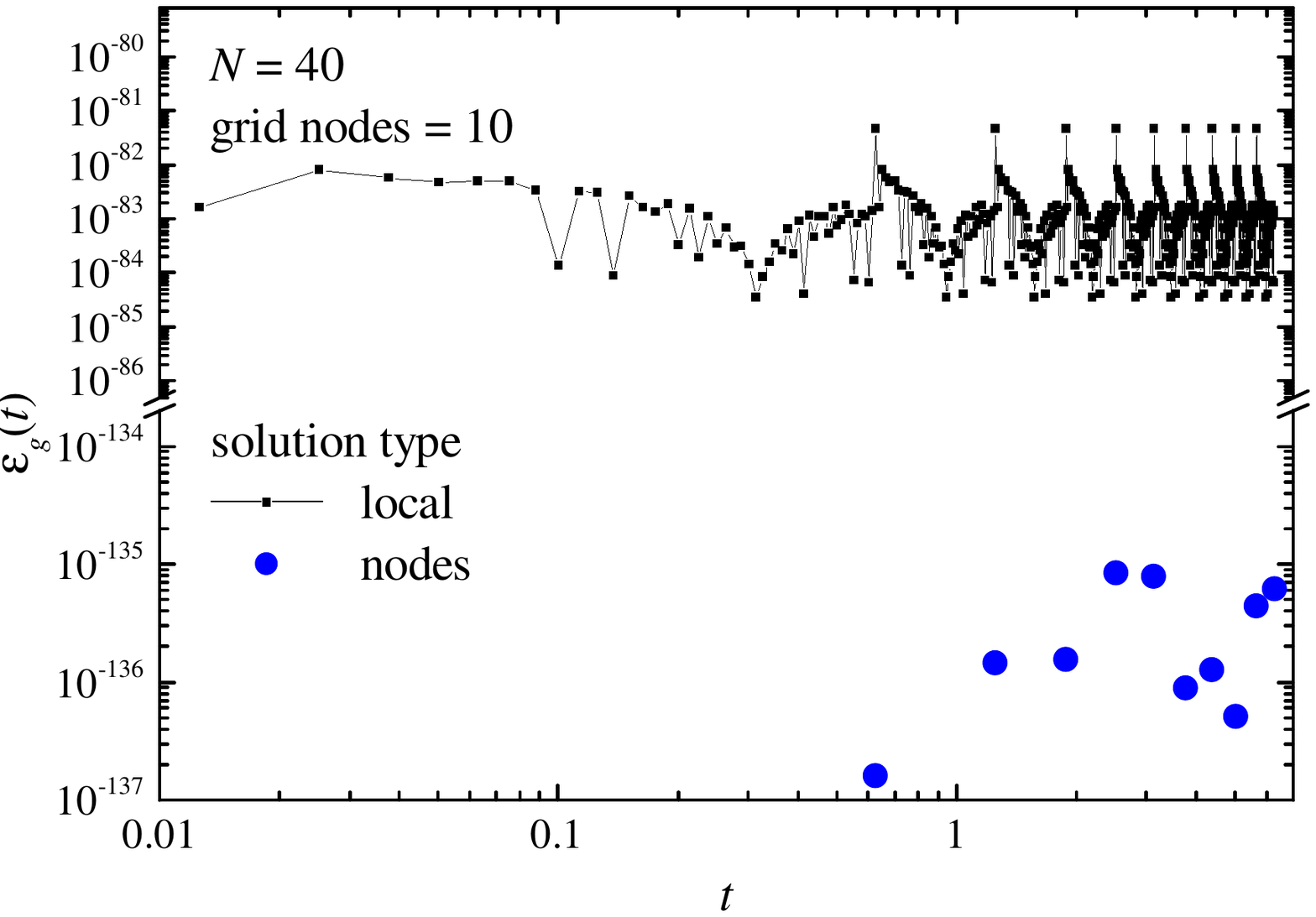}
\vspace{-8mm}\caption{\label{fig:simple_test_sol_g_eps:e3}}
\end{subfigure}\\
\caption{%
Numerical solution of the problem (\ref{eq:simple_dae_ind_1}). Comparison of quantitative satisfiability of the conditions $g_{1} = 0$ (\subref{fig:simple_test_sol_g_eps:a1}, \subref{fig:simple_test_sol_g_eps:a2}, \subref{fig:simple_test_sol_g_eps:a3}) and $g_{2} = 0$ (\subref{fig:simple_test_sol_g_eps:b1}, \subref{fig:simple_test_sol_g_eps:b2}, \subref{fig:simple_test_sol_g_eps:b3}), the errors $\varepsilon_{u}(t)$ (\subref{fig:simple_test_sol_g_eps:c1}, \subref{fig:simple_test_sol_g_eps:c2}, \subref{fig:simple_test_sol_g_eps:c3}), $\varepsilon_{v}(t)$ (\subref{fig:simple_test_sol_g_eps:d1}, \subref{fig:simple_test_sol_g_eps:d2}, \subref{fig:simple_test_sol_g_eps:d3}), $\varepsilon_{g}(t)$ (\subref{fig:simple_test_sol_g_eps:e1}, \subref{fig:simple_test_sol_g_eps:e2}, \subref{fig:simple_test_sol_g_eps:e3}), for numerical solution obtained using polynomials with degrees $N = 1$ (\subref{fig:simple_test_sol_g_eps:a1}, \subref{fig:simple_test_sol_g_eps:b1}, \subref{fig:simple_test_sol_g_eps:c1}, \subref{fig:simple_test_sol_g_eps:d1}, \subref{fig:simple_test_sol_g_eps:e1}), $N = 8$ (\subref{fig:simple_test_sol_g_eps:a2}, \subref{fig:simple_test_sol_g_eps:b2}, \subref{fig:simple_test_sol_g_eps:c2}, \subref{fig:simple_test_sol_g_eps:d2}, \subref{fig:simple_test_sol_g_eps:e2}) and $N = 40$ (\subref{fig:simple_test_sol_g_eps:a3}, \subref{fig:simple_test_sol_g_eps:b3}, \subref{fig:simple_test_sol_g_eps:c3}, \subref{fig:simple_test_sol_g_eps:d3}, \subref{fig:simple_test_sol_g_eps:e3}).
}
\label{fig:simple_test_sol_g_eps}
\end{figure}

\begin{figure}[h!]
\captionsetup[subfigure]{%
	position=bottom,
	font+=smaller,
	textfont=normalfont,
	singlelinecheck=off,
	justification=raggedright
}
\centering
\begin{subfigure}{0.275\textwidth}
\includegraphics[width=\textwidth]{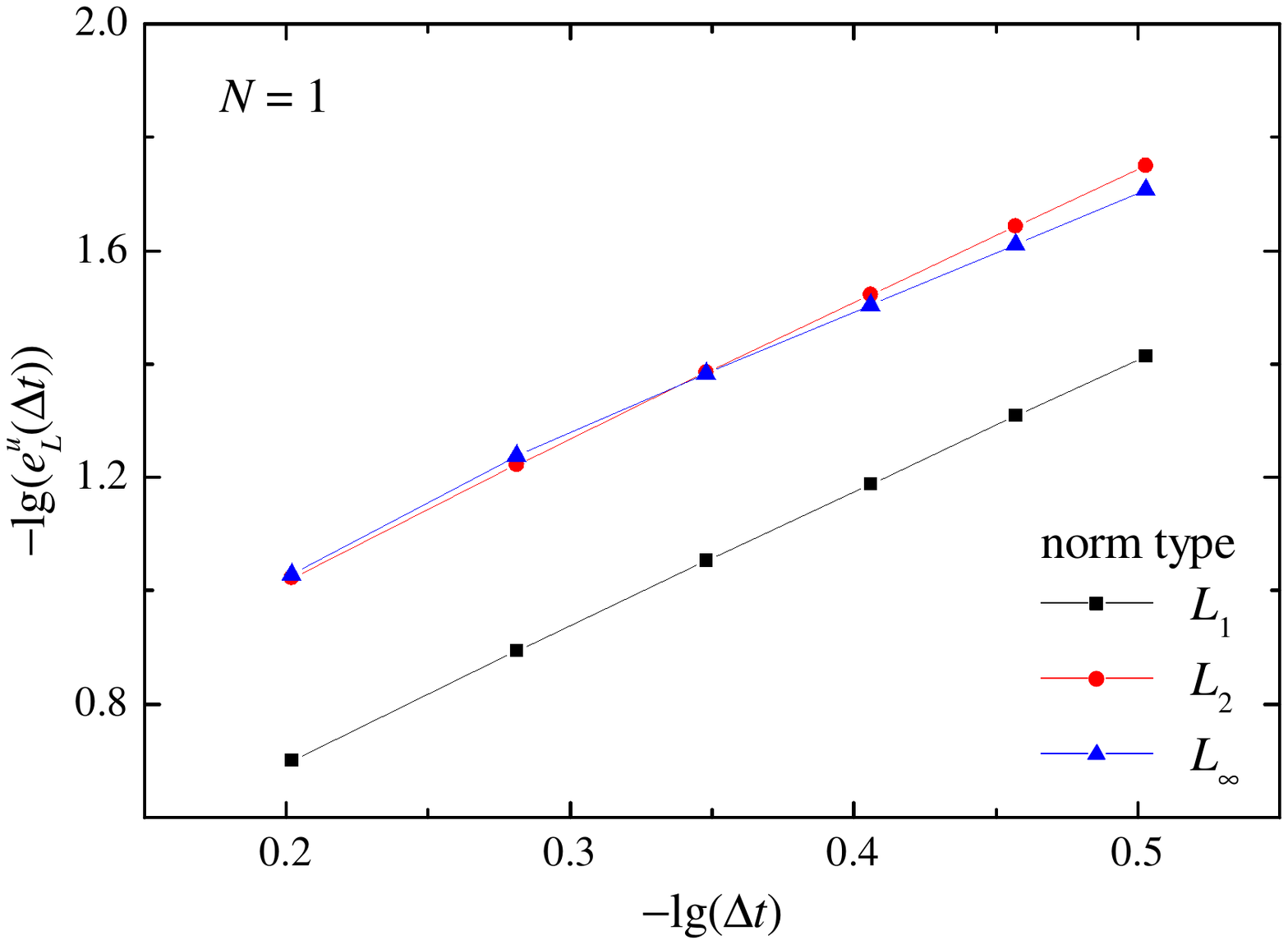}
\vspace{-8mm}\caption{\label{fig:simple_test_errors:a1}}
\end{subfigure}\hspace{6mm}
\begin{subfigure}{0.275\textwidth}
\includegraphics[width=\textwidth]{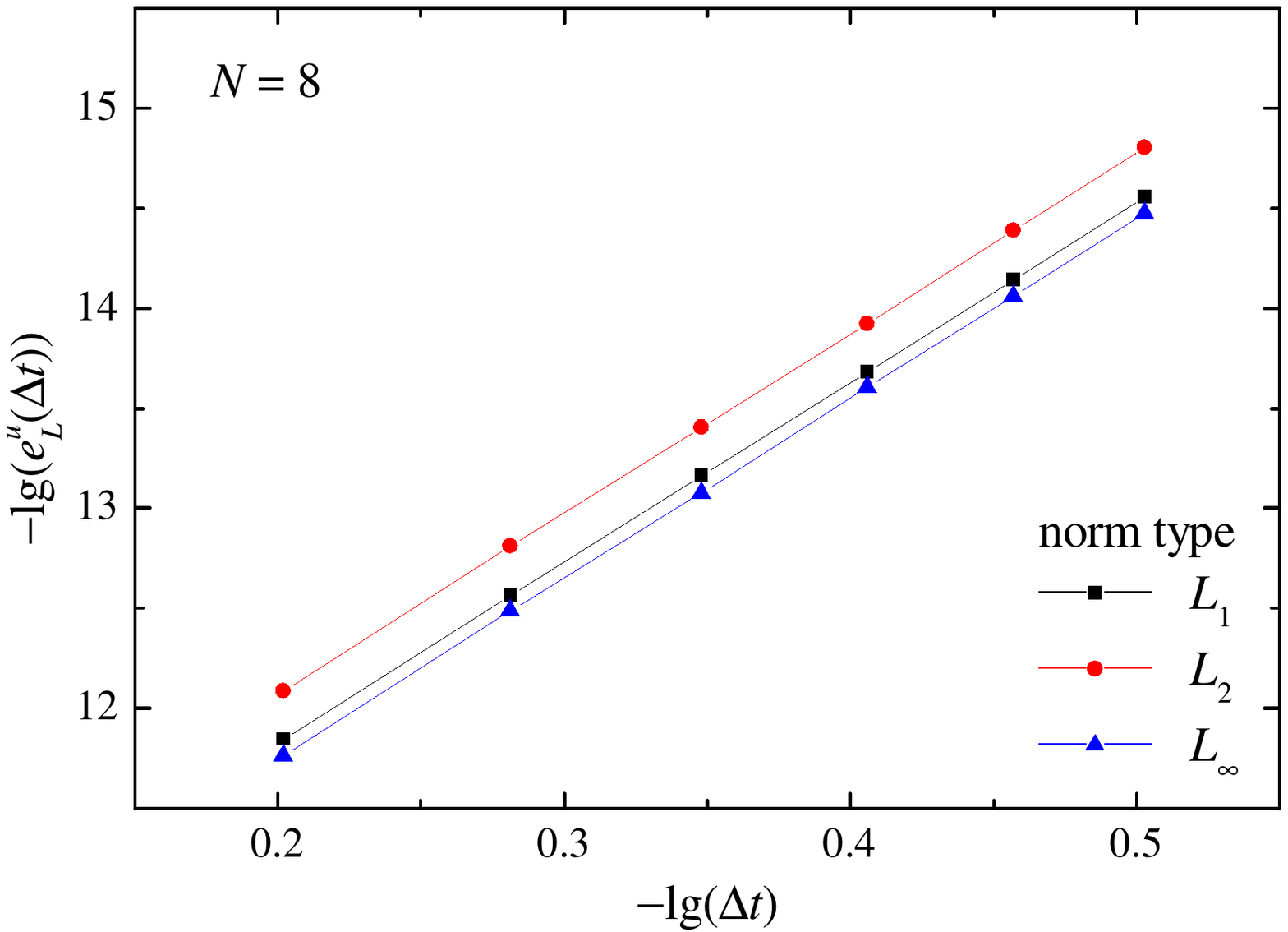}
\vspace{-8mm}\caption{\label{fig:simple_test_errors:a2}}
\end{subfigure}\hspace{6mm}
\begin{subfigure}{0.275\textwidth}
\includegraphics[width=\textwidth]{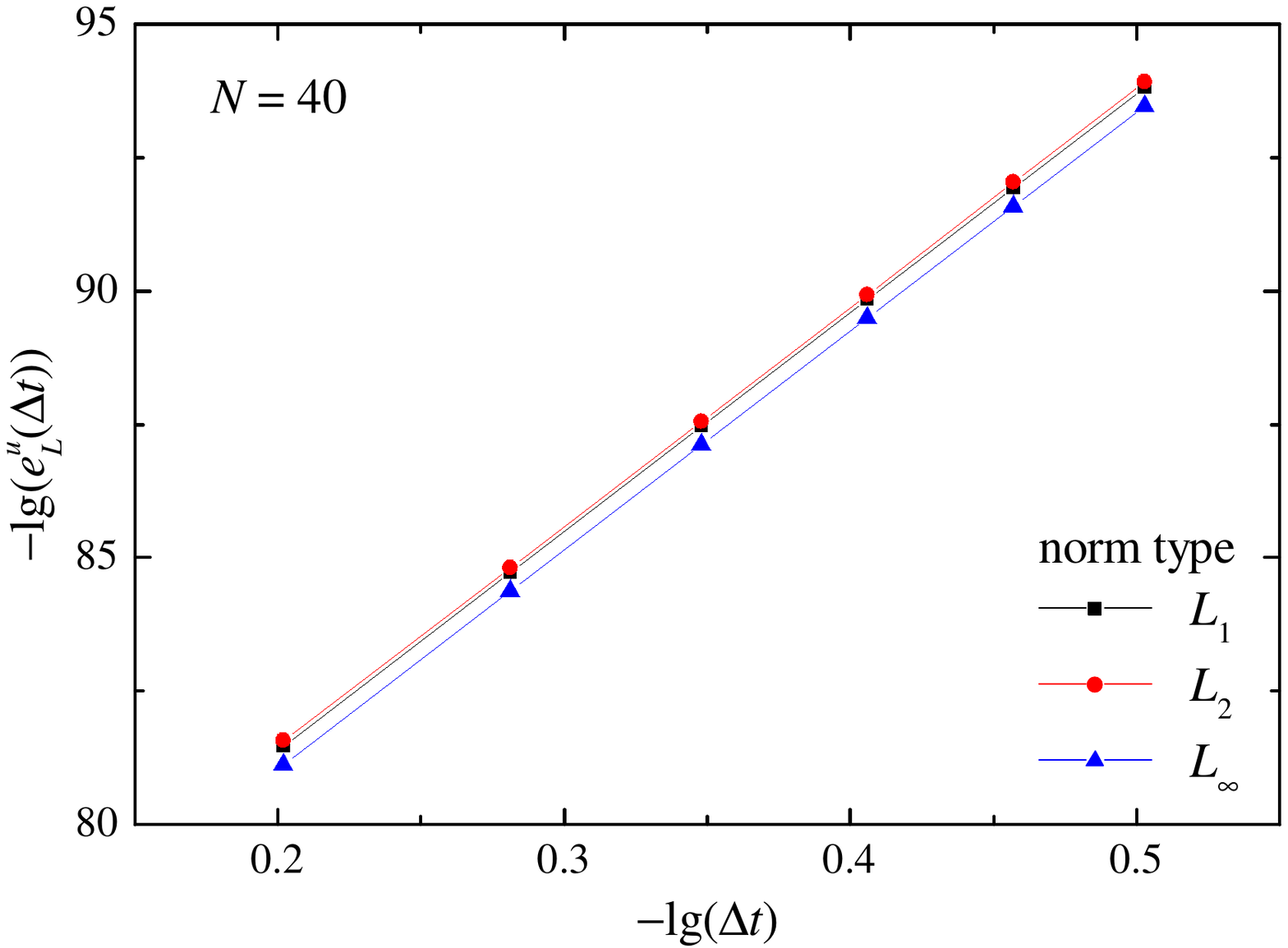}
\vspace{-8mm}\caption{\label{fig:simple_test_errors:a3}}
\end{subfigure}\\[-2mm]
\begin{subfigure}{0.275\textwidth}
\includegraphics[width=\textwidth]{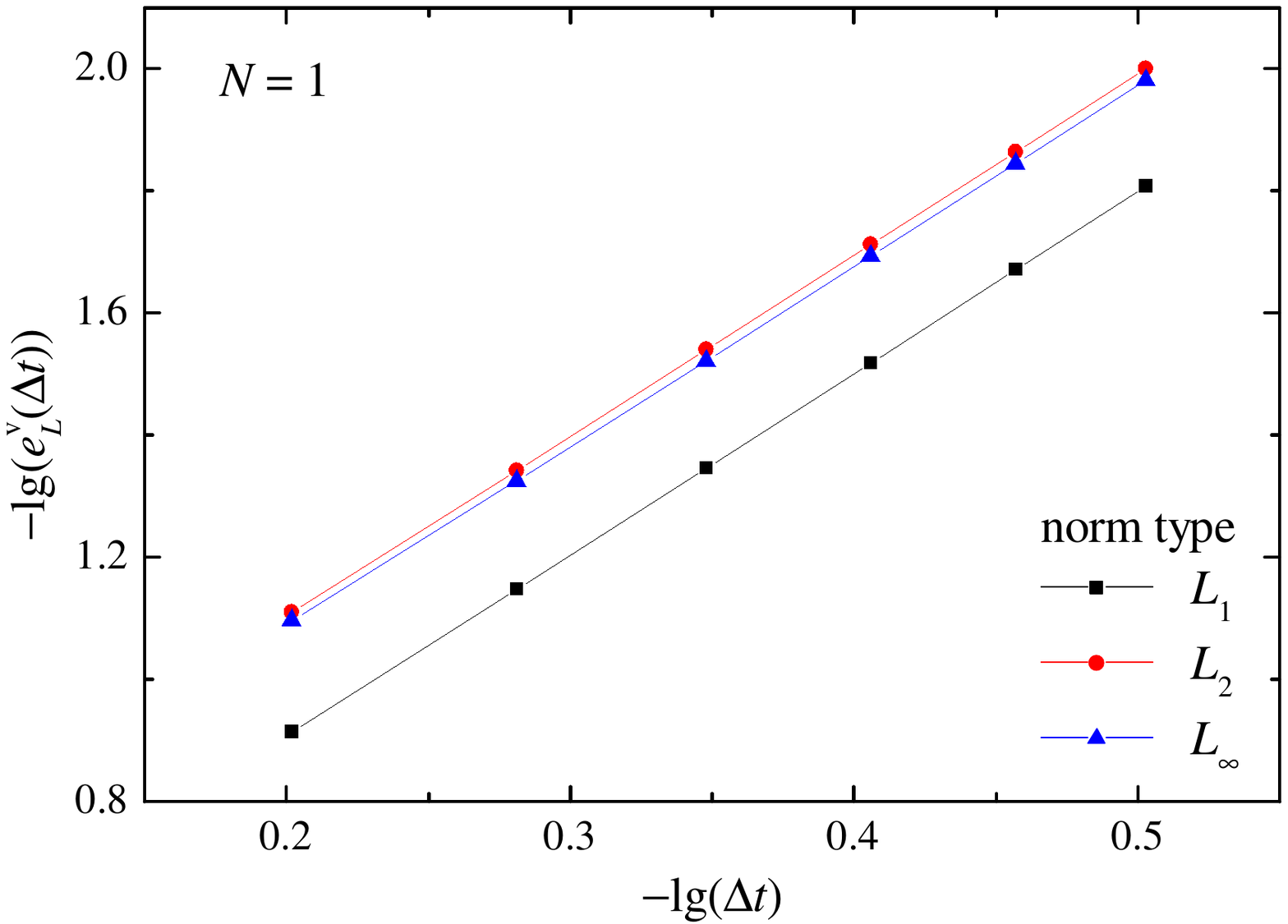}
\vspace{-8mm}\caption{\label{fig:simple_test_errors:b1}}
\end{subfigure}\hspace{6mm}
\begin{subfigure}{0.275\textwidth}
\includegraphics[width=\textwidth]{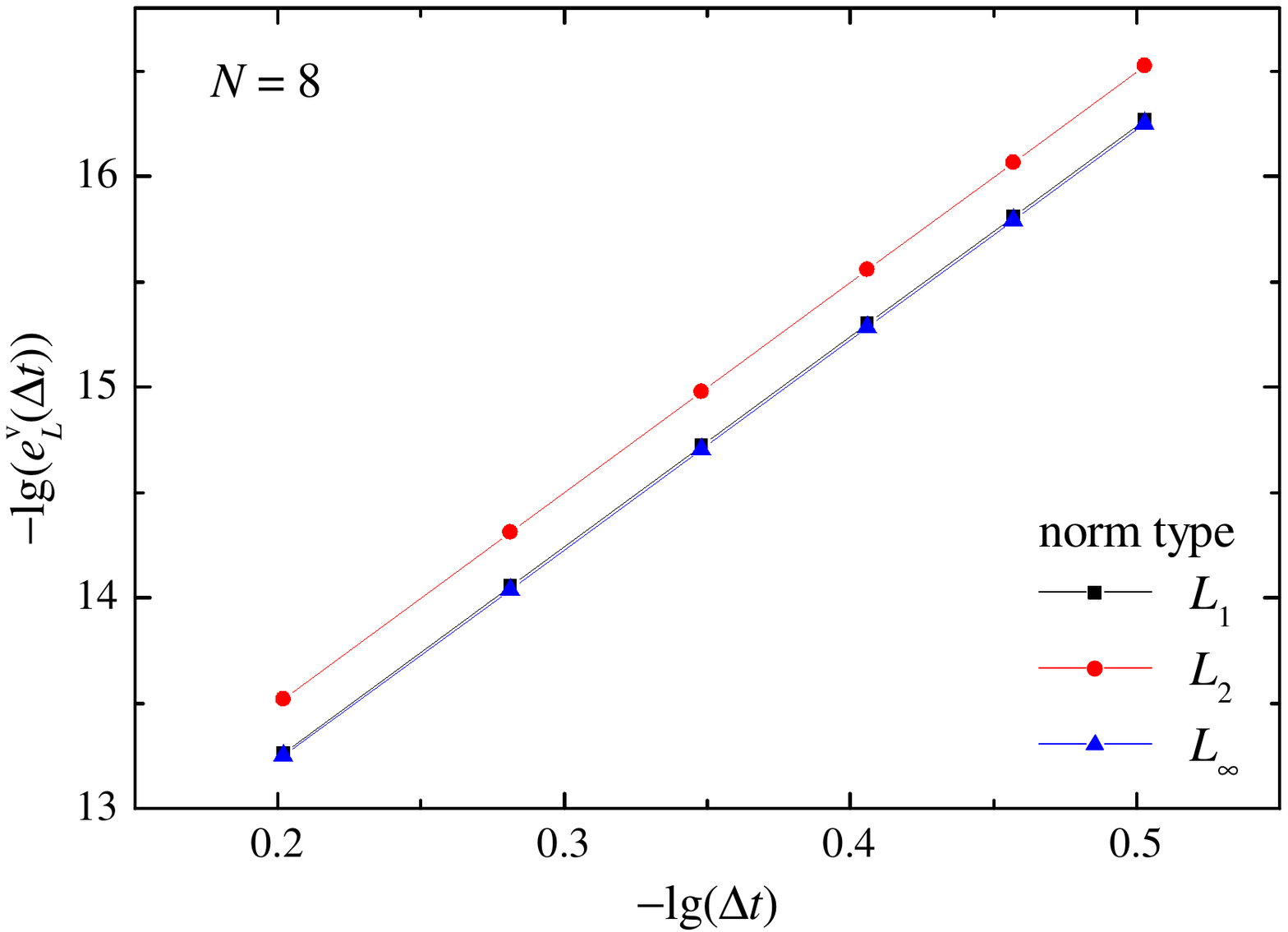}
\vspace{-8mm}\caption{\label{fig:simple_test_errors:b2}}
\end{subfigure}\hspace{6mm}
\begin{subfigure}{0.275\textwidth}
\includegraphics[width=\textwidth]{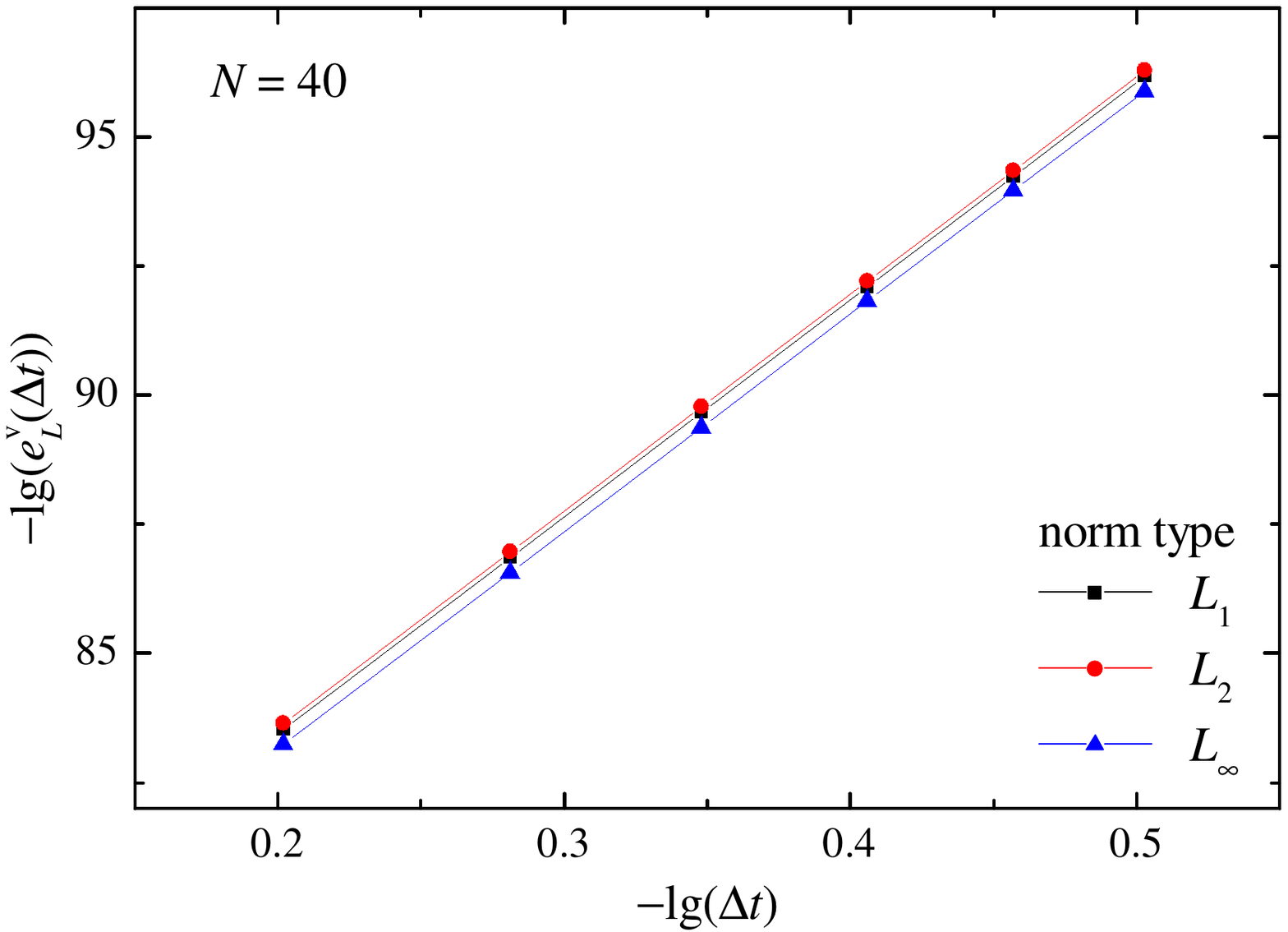}
\vspace{-8mm}\caption{\label{fig:simple_test_errors:b3}}
\end{subfigure}\\[-2mm]
\begin{subfigure}{0.275\textwidth}
\includegraphics[width=\textwidth]{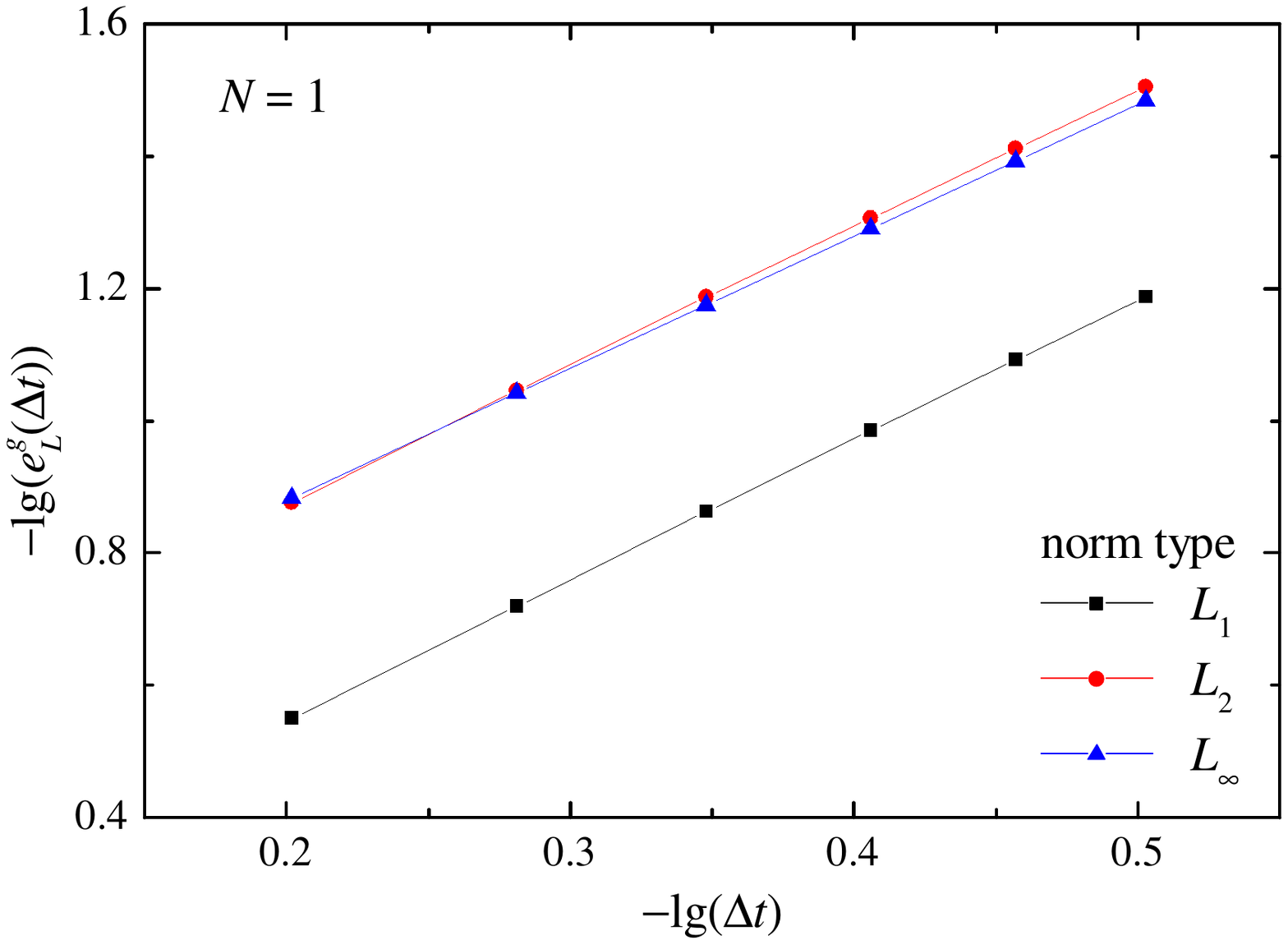}
\vspace{-8mm}\caption{\label{fig:simple_test_errors:c1}}
\end{subfigure}\hspace{6mm}
\begin{subfigure}{0.275\textwidth}
\includegraphics[width=\textwidth]{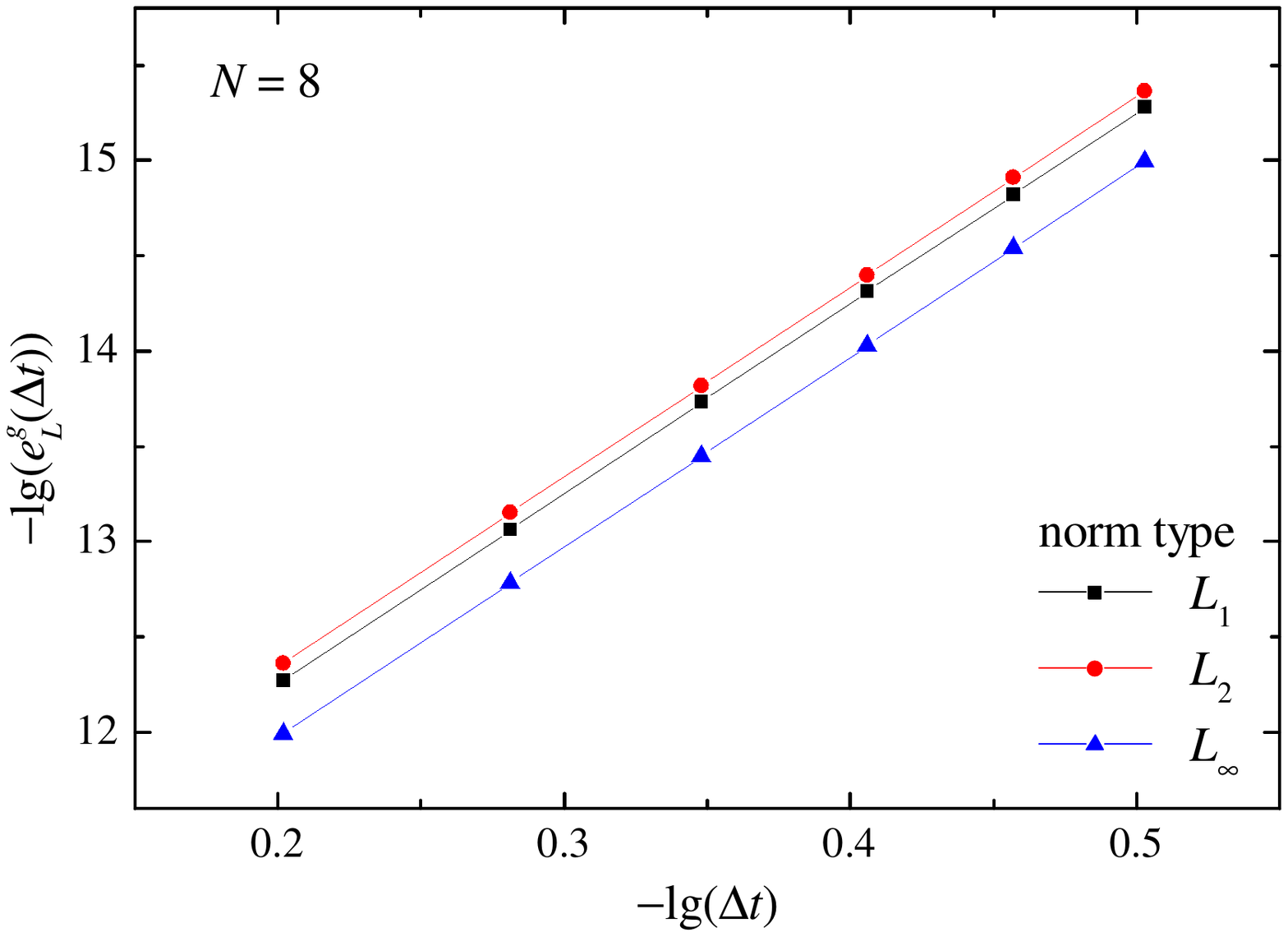}
\vspace{-8mm}\caption{\label{fig:simple_test_errors:c2}}
\end{subfigure}\hspace{6mm}
\begin{subfigure}{0.275\textwidth}
\includegraphics[width=\textwidth]{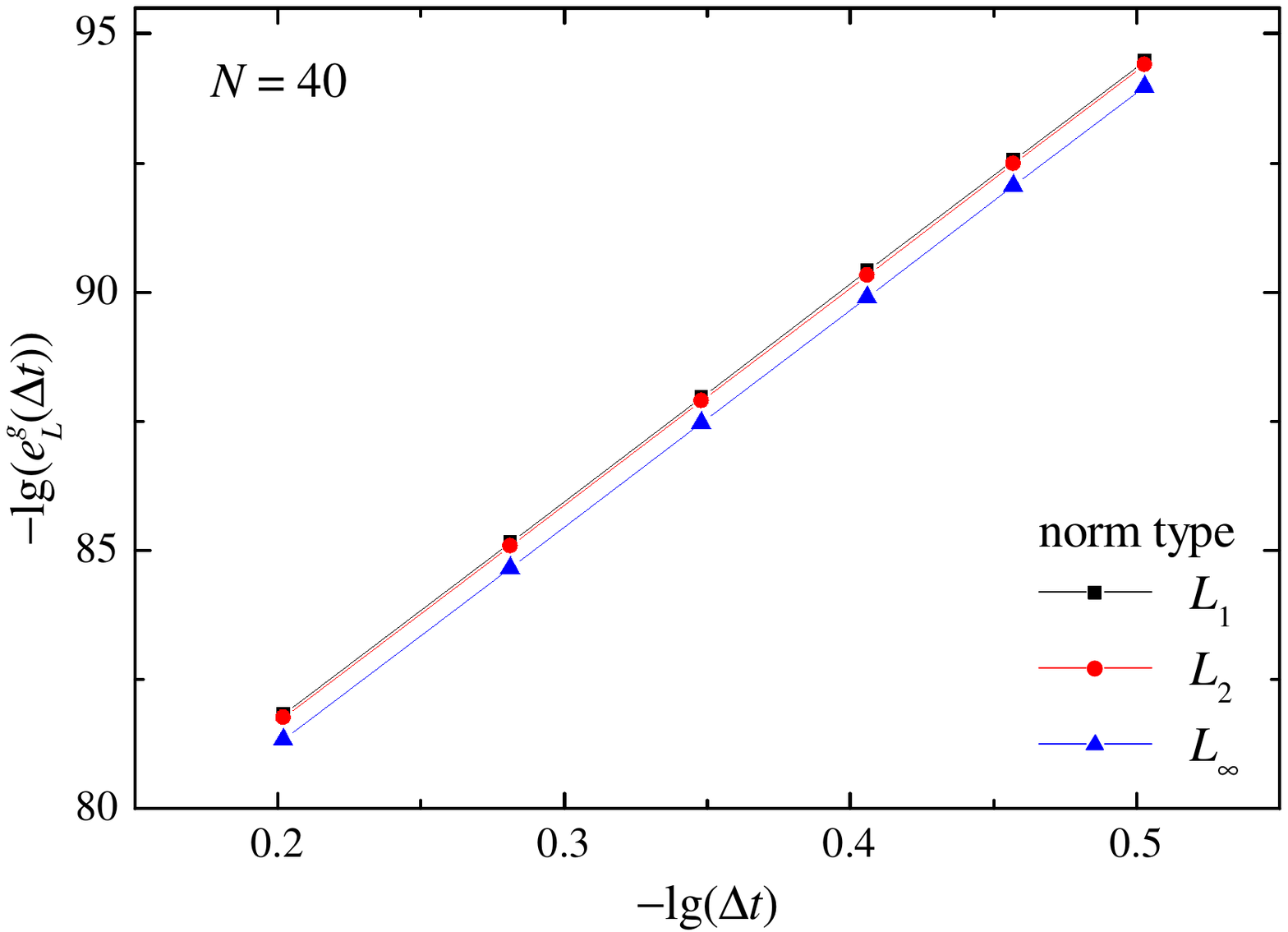}
\vspace{-8mm}\caption{\label{fig:simple_test_errors:c3}}
\end{subfigure}\\[-2mm]
\begin{subfigure}{0.275\textwidth}
\includegraphics[width=\textwidth]{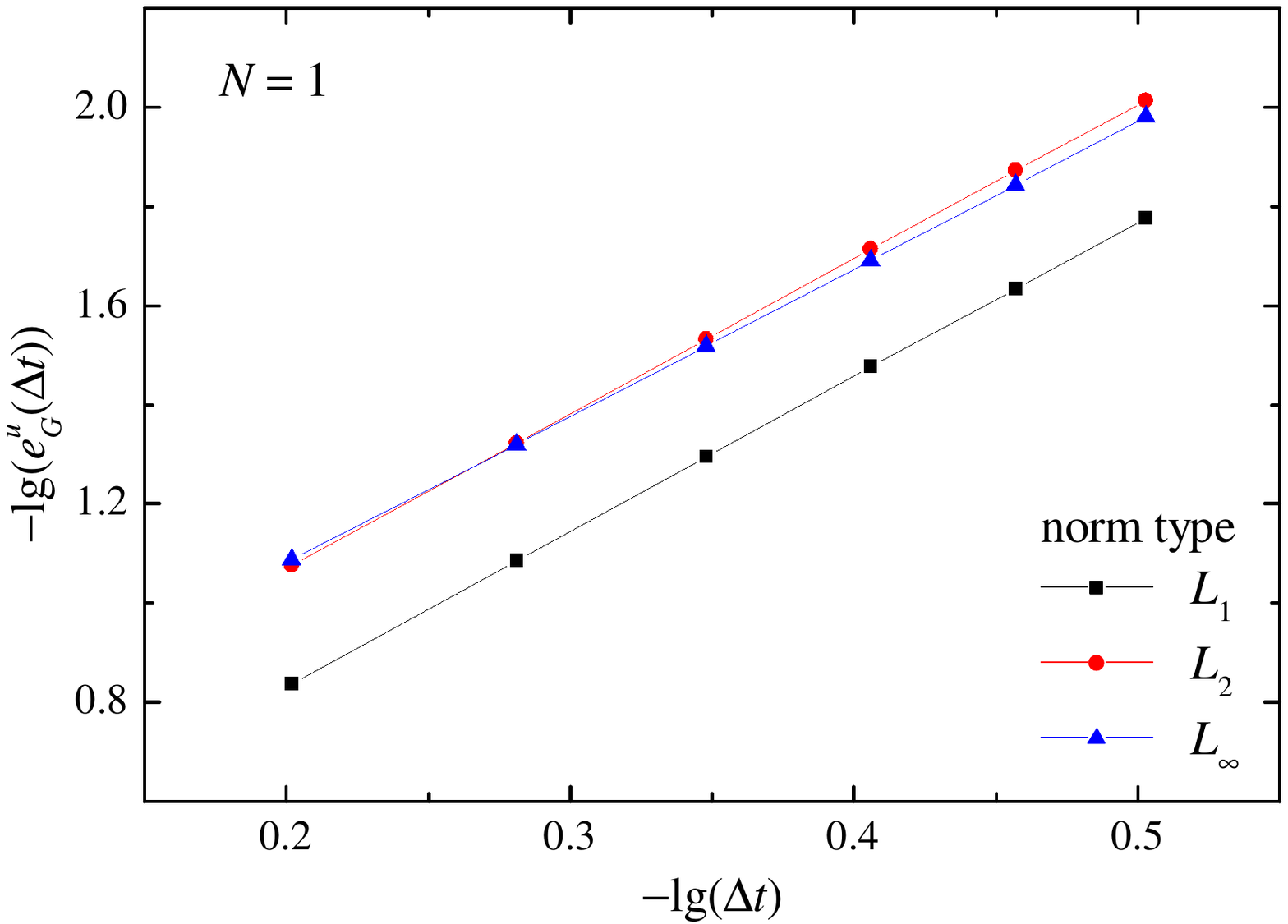}
\vspace{-8mm}\caption{\label{fig:simple_test_errors:d1}}
\end{subfigure}\hspace{6mm}
\begin{subfigure}{0.275\textwidth}
\includegraphics[width=\textwidth]{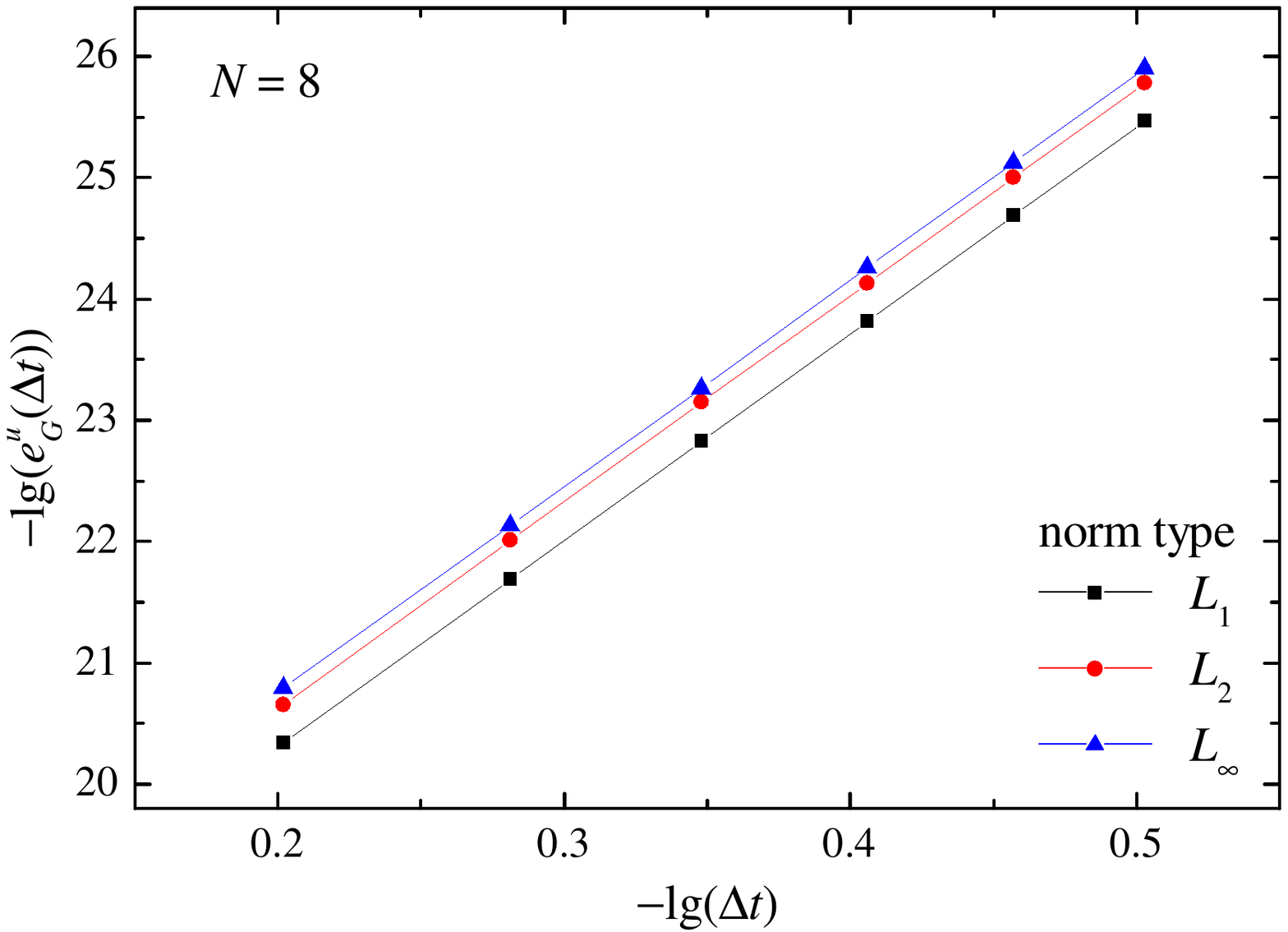}
\vspace{-8mm}\caption{\label{fig:simple_test_errors:d2}}
\end{subfigure}\hspace{6mm}
\begin{subfigure}{0.275\textwidth}
\includegraphics[width=\textwidth]{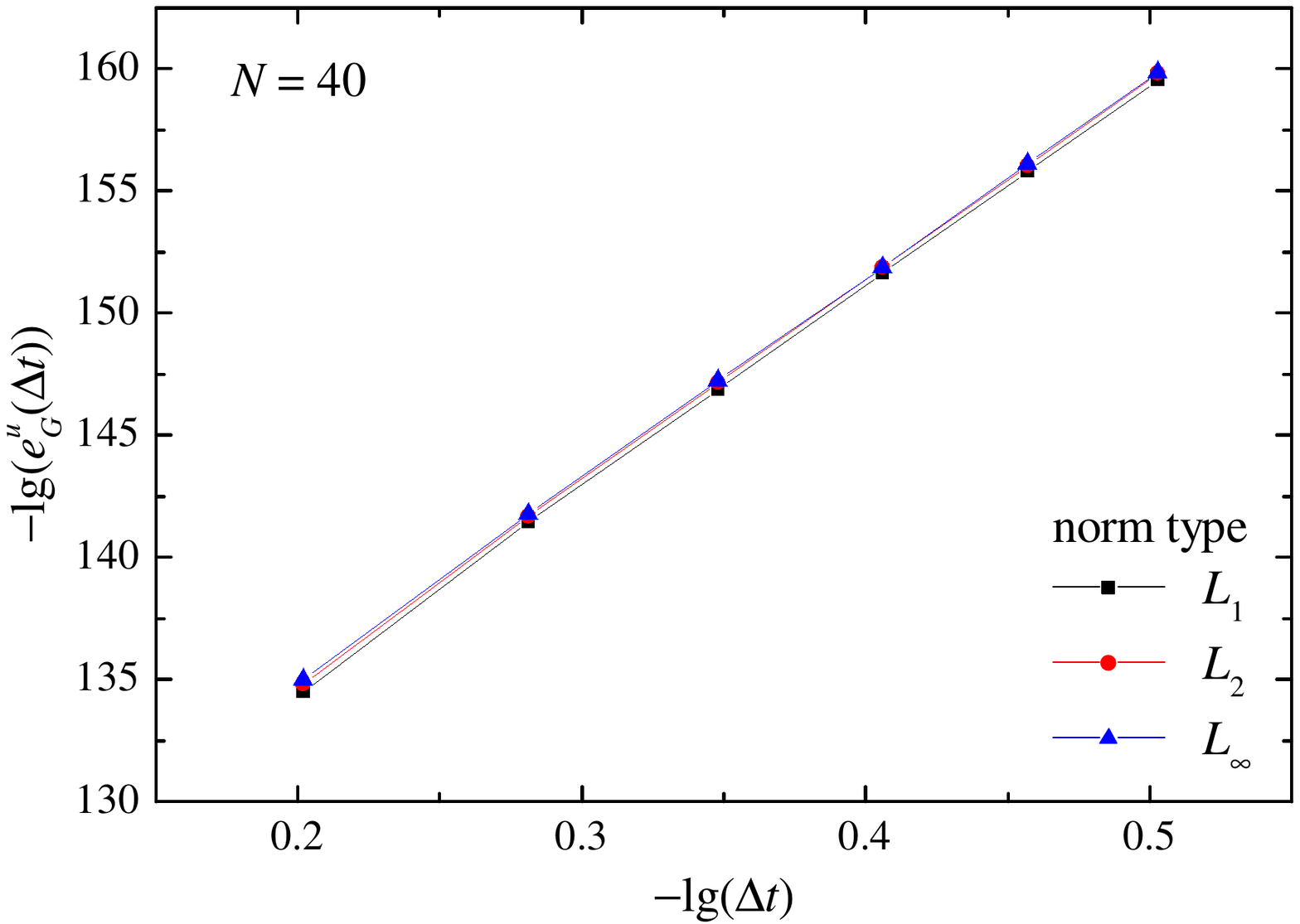}
\vspace{-8mm}\caption{\label{fig:simple_test_errors:d3}}
\end{subfigure}\\[-2mm]
\begin{subfigure}{0.275\textwidth}
\includegraphics[width=\textwidth]{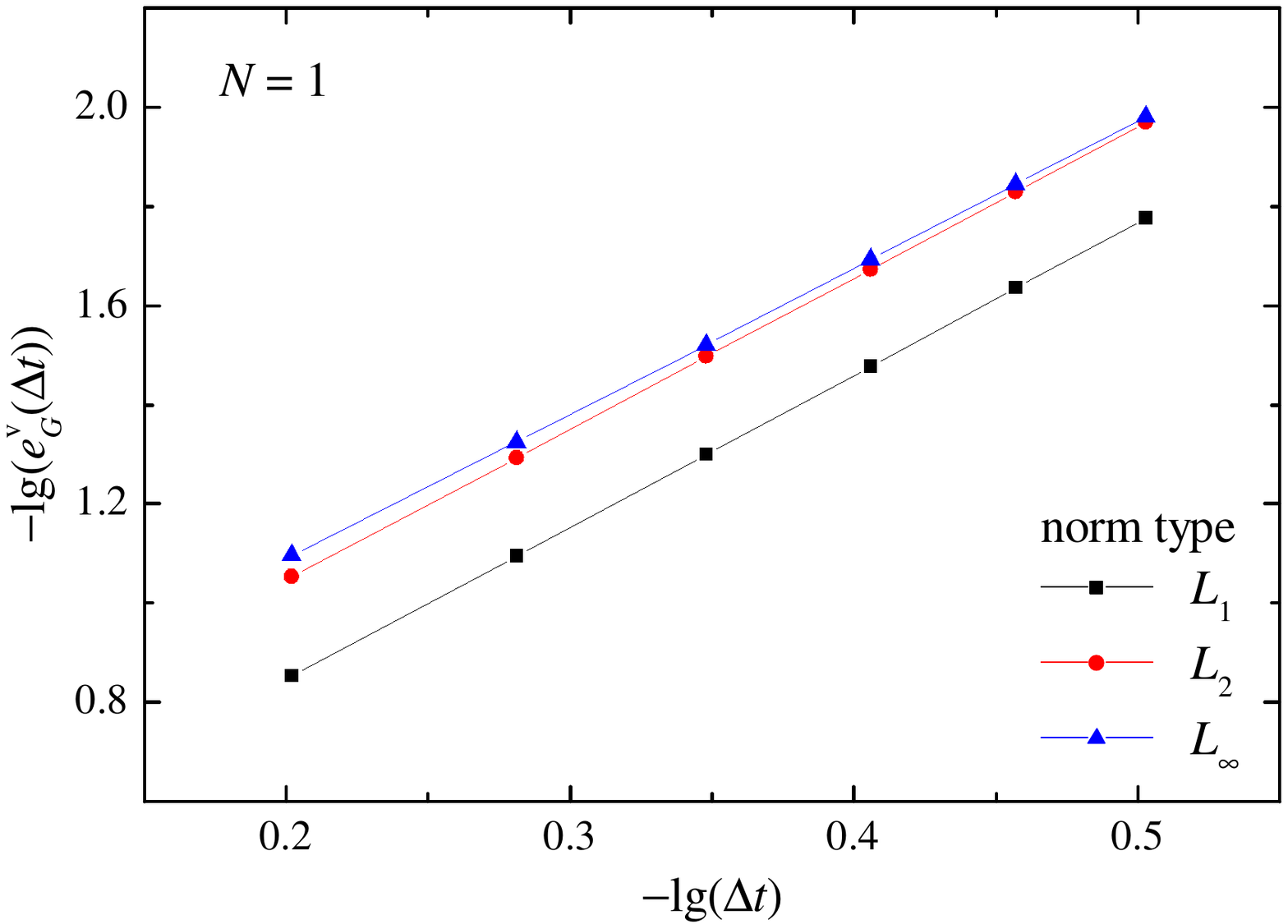}
\vspace{-8mm}\caption{\label{fig:simple_test_errors:e1}}
\end{subfigure}\hspace{6mm}
\begin{subfigure}{0.275\textwidth}
\includegraphics[width=\textwidth]{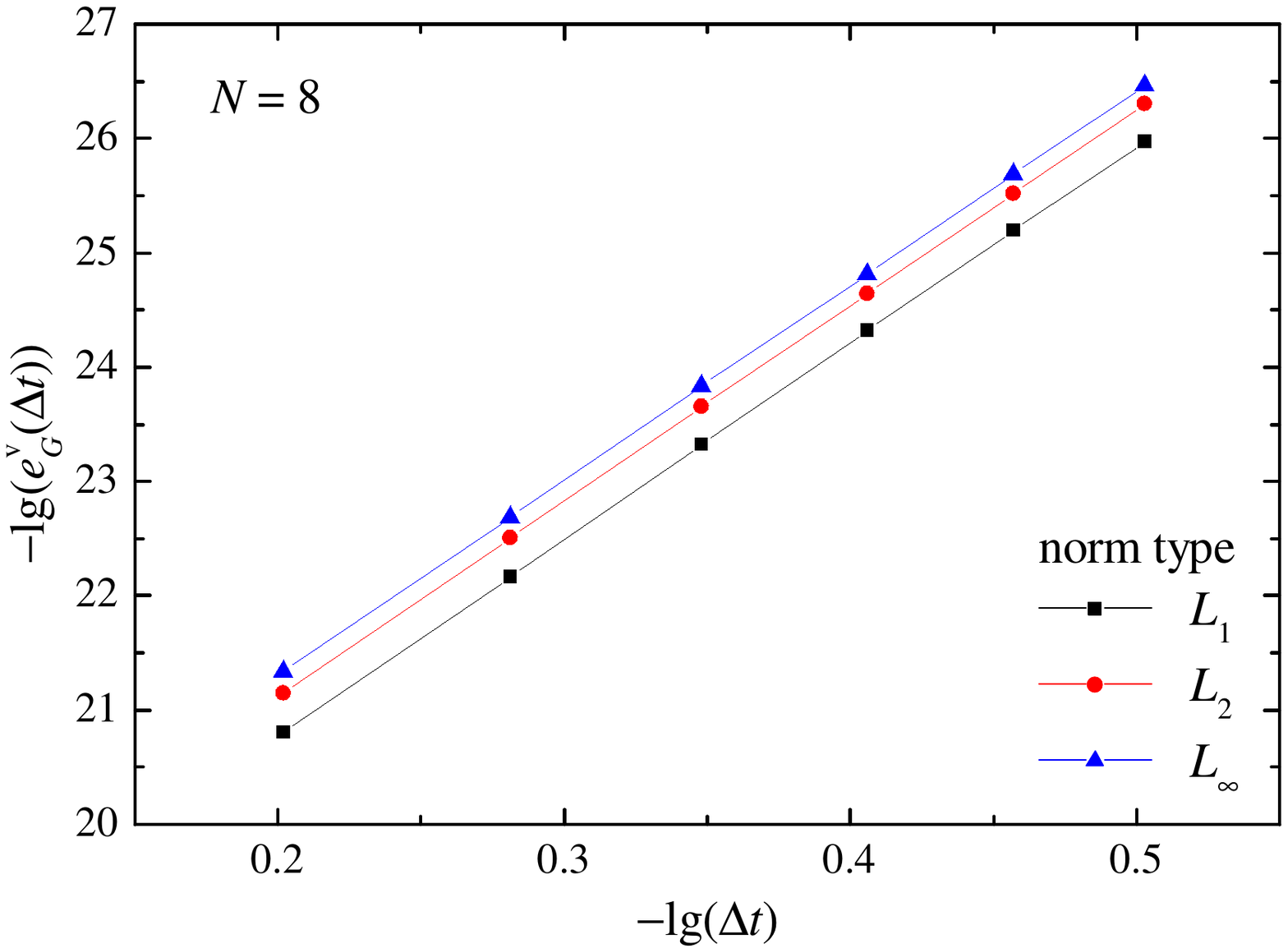}
\vspace{-8mm}\caption{\label{fig:simple_test_errors:e2}}
\end{subfigure}\hspace{6mm}
\begin{subfigure}{0.275\textwidth}
\includegraphics[width=\textwidth]{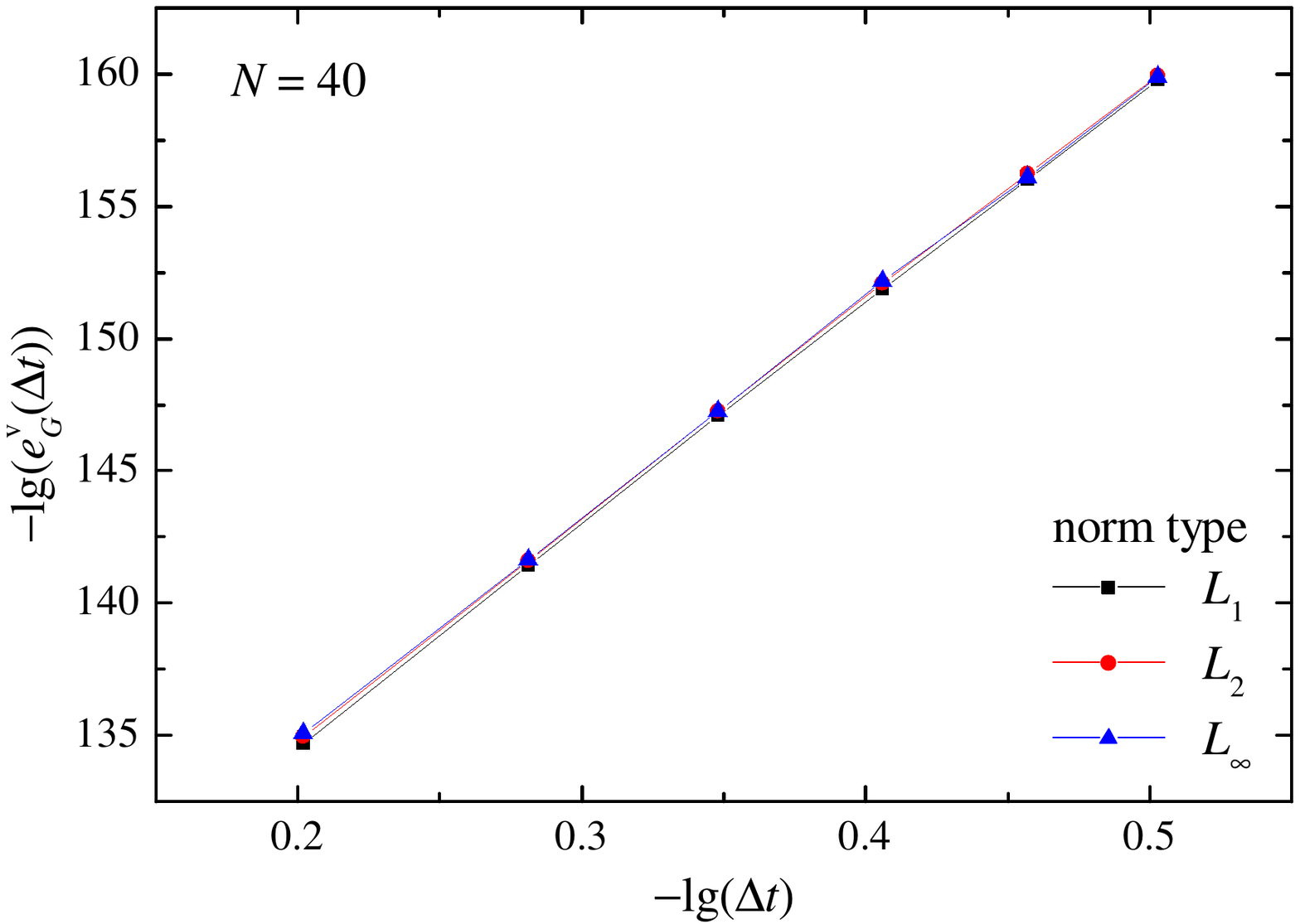}
\vspace{-8mm}\caption{\label{fig:simple_test_errors:e3}}
\end{subfigure}\\[-2mm]
\begin{subfigure}{0.275\textwidth}
\includegraphics[width=\textwidth]{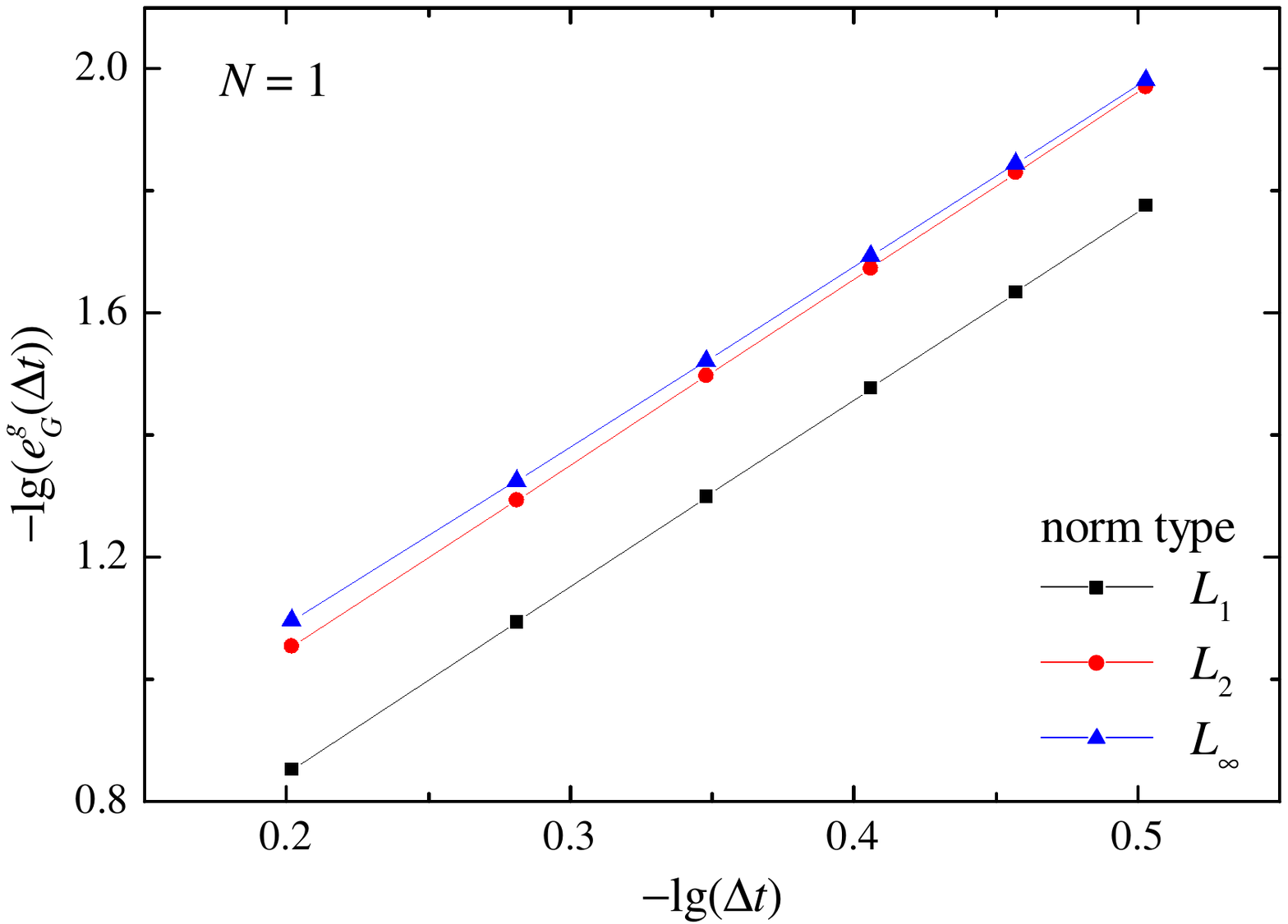}
\vspace{-8mm}\caption{\label{fig:simple_test_errors:f1}}
\end{subfigure}\hspace{6mm}
\begin{subfigure}{0.275\textwidth}
\includegraphics[width=\textwidth]{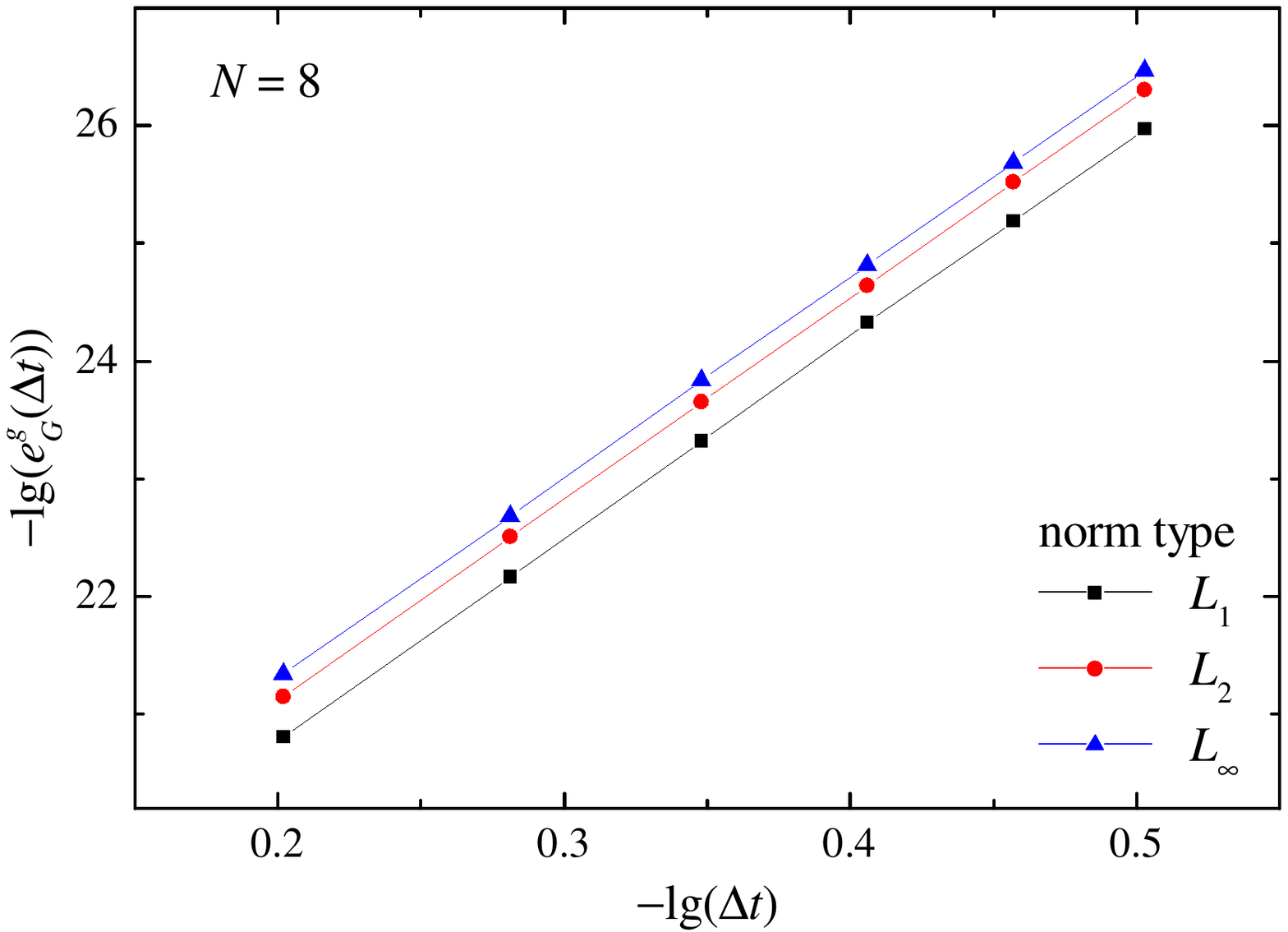}
\vspace{-8mm}\caption{\label{fig:simple_test_errors:f2}}
\end{subfigure}\hspace{6mm}
\begin{subfigure}{0.275\textwidth}
\includegraphics[width=\textwidth]{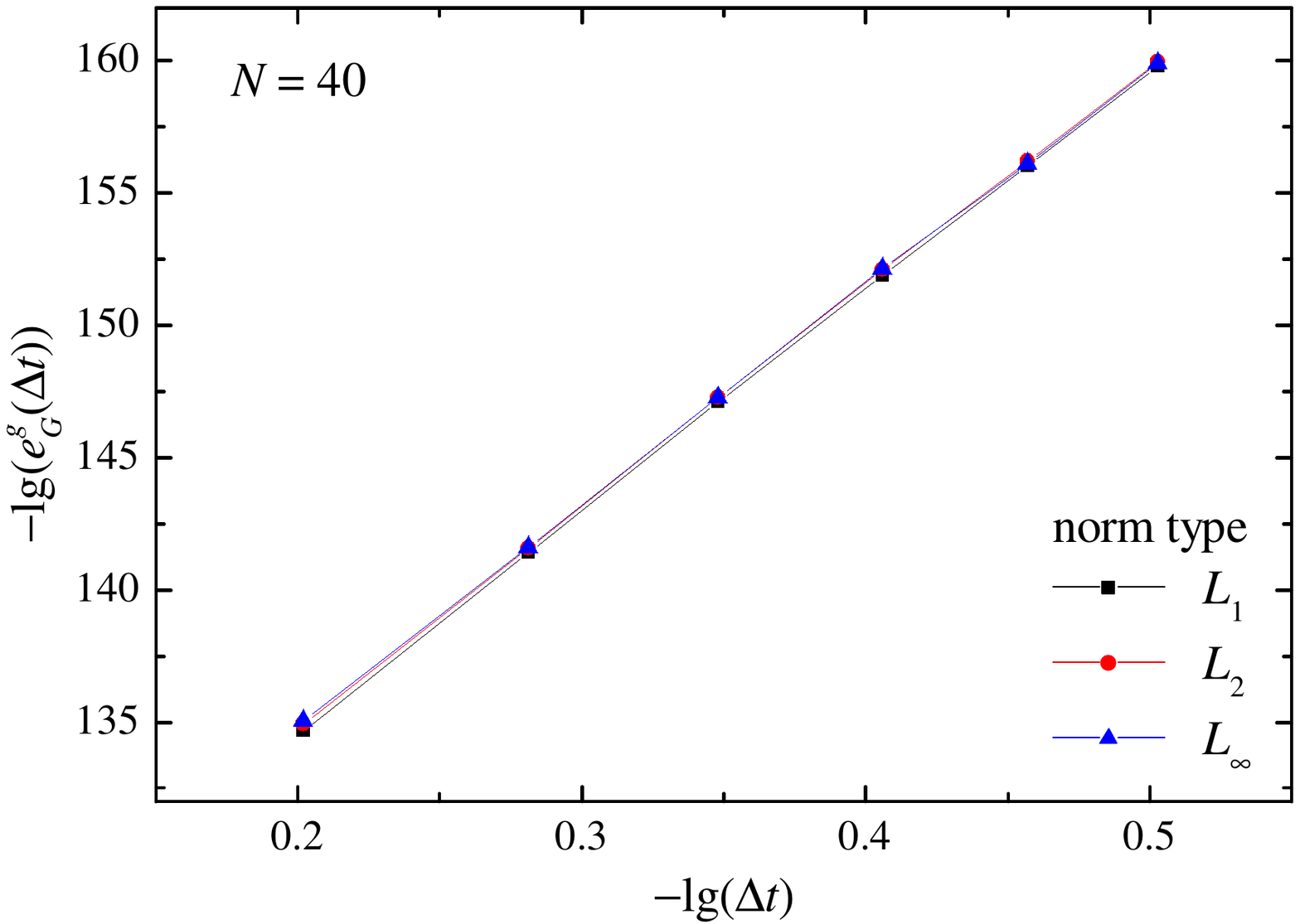}
\vspace{-8mm}\caption{\label{fig:simple_test_errors:f3}}
\end{subfigure}\\[-2mm]
\caption{%
Log-log plot of the dependence of the global errors for the local solution $e_{L}^{u}$ (\subref{fig:simple_test_errors:a1}, \subref{fig:simple_test_errors:a2}, \subref{fig:simple_test_errors:a3}), $e_{L}^{v}$ (\subref{fig:simple_test_errors:b1}, \subref{fig:simple_test_errors:b2}, \subref{fig:simple_test_errors:b3}), $e_{L}^{g}$ (\subref{fig:simple_test_errors:c1}, \subref{fig:simple_test_errors:c2}, \subref{fig:simple_test_errors:c3}) and the solution at nodes $e_{G}^{u}$ (\subref{fig:simple_test_errors:d1}, \subref{fig:simple_test_errors:d2}, \subref{fig:simple_test_errors:d3}), $e_{G}^{v}$ (\subref{fig:simple_test_errors:e1}, \subref{fig:simple_test_errors:e2}, \subref{fig:simple_test_errors:e3}), $e_{G}^{g}$ (\subref{fig:simple_test_errors:f1}, \subref{fig:simple_test_errors:f2}, \subref{fig:simple_test_errors:f3}) on the discretization step $\mathrm{\Delta}t$, obtained in the norms $L_{1}$, $L_{2}$ and $L_{\infty}$, by numerical solution of the problem (\ref{eq:simple_dae_ind_1}) obtained using polynomials with degrees $N = 1$ (\subref{fig:simple_test_errors:a1}, \subref{fig:simple_test_errors:b1}, \subref{fig:simple_test_errors:c1}, \subref{fig:simple_test_errors:d1}, \subref{fig:simple_test_errors:e1}, \subref{fig:simple_test_errors:f1}), $N = 8$ (\subref{fig:simple_test_errors:a2}, \subref{fig:simple_test_errors:b2}, \subref{fig:simple_test_errors:c2}, \subref{fig:simple_test_errors:d2}, \subref{fig:simple_test_errors:e2}, \subref{fig:simple_test_errors:f2}) and $N = 40$ (\subref{fig:simple_test_errors:a3}, \subref{fig:simple_test_errors:b3}, \subref{fig:simple_test_errors:c3}, \subref{fig:simple_test_errors:d3}, \subref{fig:simple_test_errors:e3}, \subref{fig:simple_test_errors:f3}).
}
\label{fig:simple_test_errors}
\end{figure}

The obtained results of the numerical solution of this problem are presented in Figs.~\ref{fig:simple_test_sols_uv}, \ref{fig:simple_test_sol_g_eps}, \ref{fig:simple_test_errors} and in Tables~\ref{tab:conv_orders_nodes_simple_test}, \ref{tab:conv_orders_local_simple_test}. Fig.~\ref{fig:simple_test_sols_uv} shows a comparison of the numerical solution at the nodes $(\mathbf{u}_{n}, \mathbf{v}_{n})$, the numerical local solution $(\mathbf{u}_{L}, \mathbf{v}_{L})$ and the exact analytical solution $(\mathbf{u}^{\rm ex}, \mathbf{v}^{\rm ex})$ separately for each differential $\mathbf{u}$ and algebraic $\mathbf{v}$ variable. Fig.~\ref{fig:simple_test_sol_g_eps} shows the dependencies of the feasibility of conditions $g_{1} = 0$ and $g_{2} = 0$ on the coordinate $t$, as well as the dependencies of the local errors $\varepsilon_{u}$, $\varepsilon_{v}$, $\varepsilon_{g}$ on the coordinate $t$, which allows us to quantitatively estimate the accuracy of the numerical solution, especially taking into account the fact that the numerical solution obtained by the ADER-DG method with a local DG predictor with a high degree of polynomials $N$ very accurately corresponds to the exact analytical solution, and it is impossible to visually identify the error from the solution plots separately (in Fig.~\ref{fig:simple_test_sols_uv}). Fig.~\ref{fig:simple_test_errors} shows the dependencies of the global errors $e^{u}$, $e^{v}$, $e^{g}$ of the numerical solution at the nodes $(\mathbf{u}_{n}, \mathbf{v}_{n})$ and the local solution $(\mathbf{u}_{L}, \mathbf{v}_{L})$ on the discretization step ${\Delta t}$, separately for each differential $\mathbf{u}$ and algebraic $\mathbf{v}$ variable and the algebraic equations $\mathbf{g} = \mathbf{0}$, on the basis of which the empirical convergence orders $p$ were calculated.

\begin{table*}[h!]
\centering
\caption{%
Convergence orders $p_{L_{1}}^{n}$, $p_{L_{2}}^{n}$, $p_{L_{\infty}}^{n}$, calculated in norms $L_{1}$, $L_{2}$, $L_{\infty}$, of \textit{the numerical solution at the nodes} $(\mathbf{u}_{n}, \mathbf{v}_{n})$ of the ADER-DG method for the DAE problem (\ref{eq:simple_dae_ind_1}); $N$ is the degree of the basis polynomials $\varphi_{p}$. Orders $p^{n, u}$ are calculated for solution $\mathbf{u}_{n}$; orders $p^{n, v}$ --- for solution $\mathbf{v}_{n}$; orders $p^{n, g}$ --- for the conditions $\mathbf{g} = 0$ on the numerical solution at the nodes $(\mathbf{u}_{n}, \mathbf{v}_{n})$. The theoretical value of convergence order $p_{\rm th.}^{n} = 2N+1$ is applicable for the ADER-DG method for ODE problems and is presented for comparison.
}
\label{tab:conv_orders_nodes_simple_test}
\begin{tabular}{@{}|l|lll|lll|lll|c|@{}}
\toprule
$N$ & $p_{L_{1}}^{n, u}$ & $p_{L_{2}}^{n, u}$ & $p_{L_{\infty}}^{n, u}$ & $p_{L_{1}}^{n, v}$ & $p_{L_{2}}^{n, v}$ & $p_{L_{\infty}}^{n, v}$ & $p_{L_{1}}^{n, g}$ & $p_{L_{2}}^{n, g}$ & $p_{L_{\infty}}^{n, g}$ & $p_{\rm th.}^{n}$ \\
\midrule
$1$	&	$3.12$	&	$3.11$	&	$2.97$	&	$3.07$	&	$3.05$	&	$2.94$	&	$3.07$	&	$3.05$	&	$2.94$	&	$3$\\
$2$	&	$5.10$	&	$5.11$	&	$4.99$	&	$5.10$	&	$5.08$	&	$4.99$	&	$5.10$	&	$5.08$	&	$4.99$	&	$5$\\
$3$	&	$7.07$	&	$7.06$	&	$6.99$	&	$7.10$	&	$7.07$	&	$6.99$	&	$7.10$	&	$7.07$	&	$6.99$	&	$7$\\
$4$	&	$9.04$	&	$9.02$	&	$8.99$	&	$9.05$	&	$9.04$	&	$8.97$	&	$9.05$	&	$9.04$	&	$8.97$	&	$9$\\
$5$	&	$11.02$	&	$11.01$	&	$10.99$	&	$11.01$	&	$10.99$	&	$10.93$	&	$11.01$	&	$10.99$	&	$10.93$	&	$11$\\
$6$	&	$13.01$	&	$13.01$	&	$12.99$	&	$12.92$	&	$12.96$	&	$13.14$	&	$12.92$	&	$12.96$	&	$13.14$	&	$13$\\
$7$	&	$15.02$	&	$15.02$	&	$14.99$	&	$15.37$	&	$15.28$	&	$15.08$	&	$15.37$	&	$15.28$	&	$15.08$	&	$15$\\
$8$	&	$17.03$	&	$17.02$	&	$16.98$	&	$17.19$	&	$17.12$	&	$17.03$	&	$17.19$	&	$17.12$	&	$17.03$	&	$17$\\
$9$	&	$19.06$	&	$19.04$	&	$18.99$	&	$19.09$	&	$19.08$	&	$19.15$	&	$19.09$	&	$19.08$	&	$19.15$	&	$19$\\
$10$	&	$21.05$	&	$21.06$	&	$21.02$	&	$21.04$	&	$21.08$	&	$21.15$	&	$21.04$	&	$21.08$	&	$21.15$	&	$21$\\
$11$	&	$23.05$	&	$23.08$	&	$23.06$	&	$23.04$	&	$23.11$	&	$23.14$	&	$23.04$	&	$23.11$	&	$23.14$	&	$23$\\
$12$	&	$25.11$	&	$25.11$	&	$25.11$	&	$25.21$	&	$25.15$	&	$25.12$	&	$25.21$	&	$25.15$	&	$25.12$	&	$25$\\
$13$	&	$27.13$	&	$27.14$	&	$27.14$	&	$27.21$	&	$27.16$	&	$27.13$	&	$27.21$	&	$27.16$	&	$27.13$	&	$27$\\
$14$	&	$29.14$	&	$29.16$	&	$29.16$	&	$29.20$	&	$29.16$	&	$29.05$	&	$29.20$	&	$29.16$	&	$29.05$	&	$29$\\
$15$	&	$31.20$	&	$31.22$	&	$31.17$	&	$31.18$	&	$31.17$	&	$30.97$	&	$31.18$	&	$31.17$	&	$30.97$	&	$31$\\
$16$	&	$33.25$	&	$33.26$	&	$33.17$	&	$33.18$	&	$33.19$	&	$32.97$	&	$33.18$	&	$33.19$	&	$32.97$	&	$33$\\
$17$	&	$35.24$	&	$35.28$	&	$35.16$	&	$35.18$	&	$35.23$	&	$35.01$	&	$35.18$	&	$35.23$	&	$35.01$	&	$35$\\
$18$	&	$37.26$	&	$37.29$	&	$37.15$	&	$37.25$	&	$37.29$	&	$37.13$	&	$37.25$	&	$37.29$	&	$37.13$	&	$37$\\
$19$	&	$39.29$	&	$39.29$	&	$39.15$	&	$39.45$	&	$39.35$	&	$39.29$	&	$39.45$	&	$39.35$	&	$39.29$	&	$39$\\
$20$	&	$41.31$	&	$41.31$	&	$41.14$	&	$41.46$	&	$41.38$	&	$41.35$	&	$41.46$	&	$41.38$	&	$41.35$	&	$41$\\
\midrule
$25$	&	$51.45$	&	$51.50$	&	$51.40$	&	$51.42$	&	$51.53$	&	$51.49$	&	$51.42$	&	$51.53$	&	$51.49$	&	$51$\\
$30$	&	$61.66$	&	$61.72$	&	$61.54$	&	$61.80$	&	$61.86$	&	$61.77$	&	$61.80$	&	$61.86$	&	$61.77$	&	$61$\\
$35$	&	$72.21$	&	$72.14$	&	$71.91$	&	$72.50$	&	$72.40$	&	$72.24$	&	$72.50$	&	$72.40$	&	$72.24$	&	$71$\\
$40$	&	$82.95$	&	$82.72$	&	$82.31$	&	$83.32$	&	$83.13$	&	$82.57$	&	$83.32$	&	$83.13$	&	$82.57$	&	$81$\\
\bottomrule
\end{tabular}
\end{table*} 

\begin{table*}[h!]
\centering
\caption{%
Convergence orders $p_{L_{1}}^{l}$, $p_{L_{2}}^{l}$, $p_{L_{\infty}}^{l}$, calculated in norms $L_{1}$, $L_{2}$, $L_{\infty}$, of \textit{the local solution} $(\mathbf{u}_{L}, \mathbf{v}_{L})$ (represented between the nodes) of the ADER-DG method for the DAE problem (\ref{eq:simple_dae_ind_1}); $N$ is the degree of the basis polynomials $\varphi_{p}$. Orders $p^{l, u}$ are calculated for solution $\mathbf{u}_{L}$; orders $p^{l, v}$ --- for solution $\mathbf{v}_{L}$; orders $p^{l, g}$ --- for the conditions $\mathbf{g} = 0$ on the local solution $(\mathbf{u}_{L}, \mathbf{v}_{L})$. The theoretical value of convergence order $p_{\rm th.}^{l} = N+1$ is applicable for the ADER-DG method for ODE problems and is presented for comparison.
}
\label{tab:conv_orders_local_simple_test}
\begin{tabular}{@{}|l|lll|lll|lll|c|@{}}
\toprule
$N$ & $p_{L_{1}}^{l, u}$ & $p_{L_{2}}^{l, u}$ & $p_{L_{\infty}}^{l, u}$ & $p_{L_{1}}^{l, v}$ & $p_{L_{2}}^{l, v}$ & $p_{L_{\infty}}^{l, v}$ & $p_{L_{1}}^{l, g}$ & $p_{L_{2}}^{l, g}$ & $p_{L_{\infty}}^{l, g}$ & $p_{\rm th.}^{l}$ \\
\midrule
$1$	&	$2.37$	&	$2.41$	&	$2.22$	&	$2.97$	&	$2.96$	&	$2.94$	&	$2.12$	&	$2.08$	&	$2.00$	&	$2$\\
$2$	&	$3.05$	&	$3.02$	&	$2.99$	&	$4.23$	&	$4.27$	&	$4.15$	&	$4.06$	&	$4.00$	&	$3.98$	&	$3$\\
$3$	&	$4.00$	&	$4.00$	&	$3.99$	&	$6.07$	&	$6.03$	&	$6.00$	&	$3.96$	&	$3.96$	&	$3.94$	&	$4$\\
$4$	&	$5.00$	&	$5.00$	&	$4.99$	&	$5.98$	&	$5.98$	&	$5.96$	&	$5.98$	&	$5.98$	&	$5.97$	&	$5$\\
$5$	&	$6.00$	&	$5.99$	&	$5.99$	&	$7.99$	&	$7.98$	&	$7.98$	&	$5.97$	&	$5.96$	&	$5.94$	&	$6$\\
$6$	&	$7.00$	&	$7.00$	&	$6.99$	&	$7.99$	&	$7.98$	&	$7.97$	&	$7.98$	&	$7.98$	&	$7.98$	&	$7$\\
$7$	&	$8.00$	&	$8.00$	&	$7.99$	&	$9.99$	&	$9.99$	&	$9.98$	&	$7.97$	&	$7.96$	&	$7.94$	&	$8$\\
$8$	&	$9.00$	&	$9.00$	&	$8.99$	&	$9.99$	&	$9.98$	&	$9.97$	&	$9.98$	&	$9.98$	&	$9.98$	&	$9$\\
$9$	&	$10.00$	&	$10.00$	&	$9.99$	&	$11.99$	&	$11.99$	&	$11.98$	&	$9.97$	&	$9.96$	&	$9.94$	&	$10$\\
$10$	&	$11.00$	&	$11.00$	&	$10.99$	&	$11.99$	&	$11.98$	&	$11.97$	&	$11.99$	&	$11.98$	&	$11.98$	&	$11$\\
$11$	&	$12.00$	&	$12.00$	&	$11.99$	&	$13.99$	&	$13.99$	&	$13.98$	&	$11.97$	&	$11.96$	&	$11.94$	&	$12$\\
$12$	&	$13.00$	&	$13.00$	&	$12.99$	&	$13.99$	&	$13.98$	&	$13.97$	&	$13.99$	&	$13.98$	&	$13.98$	&	$13$\\
$13$	&	$14.00$	&	$14.00$	&	$13.99$	&	$15.99$	&	$15.99$	&	$15.98$	&	$13.97$	&	$13.96$	&	$13.94$	&	$14$\\
$14$	&	$15.00$	&	$15.00$	&	$15.00$	&	$15.99$	&	$15.98$	&	$15.97$	&	$15.99$	&	$15.98$	&	$15.98$	&	$15$\\
$15$	&	$16.00$	&	$16.00$	&	$16.00$	&	$17.99$	&	$17.99$	&	$17.98$	&	$15.97$	&	$15.96$	&	$15.94$	&	$16$\\
$16$	&	$17.00$	&	$17.00$	&	$17.00$	&	$17.99$	&	$17.98$	&	$17.97$	&	$17.99$	&	$17.98$	&	$17.98$	&	$17$\\
$17$	&	$18.00$	&	$18.00$	&	$18.00$	&	$19.99$	&	$19.99$	&	$19.98$	&	$17.97$	&	$17.96$	&	$17.94$	&	$18$\\
$18$	&	$19.00$	&	$19.00$	&	$19.00$	&	$19.99$	&	$19.98$	&	$19.97$	&	$19.99$	&	$19.98$	&	$19.98$	&	$19$\\
$19$	&	$20.00$	&	$20.00$	&	$20.00$	&	$21.99$	&	$21.99$	&	$21.99$	&	$19.97$	&	$19.96$	&	$19.94$	&	$20$\\
$20$	&	$21.00$	&	$21.00$	&	$21.00$	&	$21.99$	&	$21.98$	&	$21.97$	&	$21.99$	&	$21.98$	&	$21.98$	&	$21$\\
\midrule
$25$	&	$26.00$	&	$26.00$	&	$26.00$	&	$27.99$	&	$27.99$	&	$27.99$	&	$25.97$	&	$25.95$	&	$25.95$	&	$26$\\
$30$	&	$31.00$	&	$31.00$	&	$31.00$	&	$31.99$	&	$31.98$	&	$31.97$	&	$31.99$	&	$31.98$	&	$31.98$	&	$31$\\
$35$	&	$36.00$	&	$36.00$	&	$36.00$	&	$37.99$	&	$37.99$	&	$37.99$	&	$35.97$	&	$35.95$	&	$35.95$	&	$36$\\
$40$	&	$41.00$	&	$41.00$	&	$41.00$	&	$41.99$	&	$41.98$	&	$41.97$	&	$41.99$	&	$41.98$	&	$41.98$	&	$41$\\
\bottomrule
\end{tabular}
\end{table*}

The numerical solution presented in Fig.~\ref{fig:simple_test_sols_uv} demonstrates very high accuracy when compared with the exact analytical solution. In the case of polynomials of degree $N = 1$, discontinuities in the grid nodes are visually observed, which is due to the linearity of the local solution representation in this case. At the same time, a one-sided deviation of the local solution curve and solution points at the nodes for differential variables is also observed, which indicates dissipative errors in the numerical solution in the case of polynomials of degree $N = 1$, which is most likely due to the superstability of the ADER-DG numerical method with a local DG predictor and insufficient accuracy of the method in case $N = 1$ to compensate for dissipative errors on a grid with $10$ discretization domains. A systematic decrease in the value of the algebraic variable $\mathbf{v}$ from the exact constant value $\mathbf{v}^{\rm ex} = 1$ is also observed, which is also due to dissipative errors. In the case of polynomials of degrees $N = 8$ and $40$, these errors in the numerical solution are not observed --- the accuracy of the numerical solution is very high even on a coarse grid, and it is impossible to visually distinguish between the numerical solution and the exact analytical solution on the presented plots. The dependencies of the numerical solution errors on the coordinate shown in Fig.~\ref{fig:simple_test_sol_g_eps} also demonstrate the properties of the numerical solution described above. In the case of polynomials of degree $N = 1$, the error $|g_{1}|$ of algebraic equation $g_{1} = 0$ systematically practically does not exceed $10^{-1}$, and does not demonstrate growth with the coordinate $t$, which indicates high efficiency of the ADER-DG numerical method for solving the DAE system even in the case of low polynomials of degree $N$ and a coarse grid. The error in performing algebraic equation $g_{2} = 0$ demonstrates approximately exponential growth, which is associated precisely with the dissipative errors of the numerical method on a coarse grid in the case of polynomials of degree $N = 1$. The local errors $\varepsilon_{u}$, $\varepsilon_{v}$, $\varepsilon_{g}$ of the numerical solution in the case of polynomials of degree $N = 1$ remain sufficiently small in the selected domain of definition of the desired function $[0,\, 2\pi]$. In the case of polynomial degrees $N = 8$ and $40$, the numerical solution errors $\varepsilon_{u}$, $\varepsilon_{v}$, $\varepsilon_{g}$ become very small. The algebraic equation $g_{1} = 0$ is satisfied exactly for the solution at the nodes, so the errors for the solution at the nodes are not presented on the plots $|g_{1}|$; therefore, the algebraic constraints $\mathbf{g} = \mathbf{0}$ in the solution of the DAE system are strictly observed. The errors of the solution of algebraic equation $g_{1} = 0$ for the local solution are not zero, but they are very small and do not demonstrate a systematic increase with the coordinate $t$. The error of the algebraic equation $g_{2} = 0$ corresponds to a very high accuracy of the calculation of the algebraic variable $v_{1}$. It is necessary to note a significant difference in the errors $\varepsilon_{u}$, $\varepsilon_{v}$, $\varepsilon_{g}$ of the local solution $(\mathbf{u}_{L}, \mathbf{v}_{L})$ and the solution at the nodes $(\mathbf{u}_{n}, \mathbf{v}_{n})$ --- in the case of the execution of algebraic equation $|g_{2}|$ and the local errors $\varepsilon_{u}$, $\varepsilon_{v}$, $\varepsilon_{g}$, this difference reaches $10^{8}$--$10^{9}$ times for polynomials of degree $N = 8$ and $10^{50}$--$10^{55}$ times for polynomials of degree $N = 40$. The characteristic values of the errors of the local solution are $10^{-12}$--$10^{-15}$ in the case of $N = 8$ and $10^{-81}$--$10^{-84}$ in the case of $N = 40$. The characteristic values of the errors of the solution at the nodes are $10^{-21}$--$10^{-23}$ in the case of $N = 8$ and $10^{-134}$--$10^{-137}$ in the case of $N = 40$. These error values are in good agreement in order of magnitude with the first point (with the largest discretization step ${\Delta t}$) of the global errors $e$ dependencies on the discretization step ${\Delta t}$, presented in Fig.~\ref{fig:simple_test_errors}. These results also demonstrate a clearly expressed power dependence of the errors on the discretization step $e(\Delta t) \sim (\Delta t)^{p}$. In this case, in the case of polynomials of degree $N = 1$, the error values at the selected step values reach $10^{-2}$, for polynomials of degree $N = 8$ they reach $10^{-15}$ and $10^{-27}$ for the local solution $(\mathbf{u}_{L}, \mathbf{v}_{L})$ and the solution at the nodes $(\mathbf{u}_{n}, \mathbf{v}_{n})$, respectively, and in the case of polynomials of degree $N = 40$ --- $10^{-95}$ and $10^{-161}$, respectively.

The empirical convergence orders $p^{n, u}$, $p^{n, v}$, $p^{n, g}$ for the numerical solution at the nodes $(\mathbf{u}_{n}, \mathbf{v}_{n})$ are presented in Table~\ref{tab:conv_orders_nodes_simple_test} and the empirical convergence orders $p^{l, u}$, $p^{l, v}$, $p^{l, g}$ for the local solution $(\mathbf{u}_{L}, \mathbf{v}_{L})$ --- in Table~\ref{tab:conv_orders_local_simple_test}. The convergence orders were calculated separately for the norms $L_{1}$, $L_{2}$, $L_{\infty}$ (\ref{eq:norms_def}). Also for comparison, the expected values of the convergence orders, which the ADER-DG numerical method with the local DG predictor provides for solving the initial value problem for the ODE system, are presented in Table~\ref{tab:conv_orders_nodes_simple_test}: $p_{\rm nodes} = 2N+1$, $p_{\rm local} = N+1$. Comparison of the obtained convergence orders $p^{n, u}$, $p^{n, v}$, $p^{n, g}$ for the solution at the nodes clearly demonstrates not only the expected superconvergence $2N+1$, but also higher values of the convergence orders --- the obtained empirical values are often greater than the expected values by $1$ and by $2$. Comparison of the obtained convergence orders $p^{l, u}$, $p^{l, v}$, $p^{l, g}$ for the local solution demonstrates similar behavior --- the obtained empirical values of the convergence orders are in good agreement or exceed the expected values $N+1$. An interesting pattern is observed --- the empirical convergence orders for algebraic variables $p^{l, v}$ and algebraic constraints $p^{l, g}$ exceed the expected values by $1$ for odd degrees of polynomials $N$. 

Therefore, it can be concluded that in the case of this example of a DAE system, the ADER-DG method with a local DG predictor allows one to obtain convergence of the numerical solution and convergence orders $p$ similar to or exceeding those expected (\ref{eq:expect_orders}).

\subsubsection{Example 2: Hessenberg DAE system of index 1}
\label{sec:2:ct:ex2}

The second example of application of the numerical method ADER-DG with local DG predictor consisted in solving a Hessenberg DAE system of index 1, which was considered in the work~\cite{hess_daes_examples_2016}:
\begin{equation}\label{eq:hess_dae_ind_1}
\begin{split}
&\ddot{x} + x(4z + 1) + y(3t + 1) = 0,\hspace{15mm} x(0) = 1,\quad \dot{x}(0) = 1,\\
&\ddot{y} + y(4z + 1) - 4\cos(z) = 0,\hspace{17.5mm} y(0) = 0,\quad \dot{y}(0) = 2,\\
&g_{1} = 4x\cos(z) + ty^{2} - 4(z - t^{2}) = 0,\hspace{5.5mm} z(0) = 0,\quad t\in[0, 1].
\end{split}
\end{equation}
The exact analytical solution of this problem was obtained in the following form: $x = t\cos(t^{2} + t)$, $\dot{x} = \cos(t^{2} + t) - t(2t + 1)\sin(t^{2} + t)$, $y = 2\sin(t^{2} + t)$, $\dot{y} = 2(2t + 1)\cos(t^{2} + t)$, $z = t^{2} + t$. This problem was rewritten in a form consistent with the formulation of the original DAE system (\ref{eq:dae_chosen_form}):
\begin{equation}
\begin{split}
&\frac{du_{1}}{dt} = u_{3},\hspace{50.6mm} u_{1}(0) = 0,\\
&\frac{du_{2}}{dt} = u_{4},\hspace{50.6mm} u_{2}(0) = 0,\\
&\frac{du_{3}}{dt} = -u_{1}(4v_{1} + 1) - u_{2}(3t + 1),\hspace{14.4mm} u_{3}(0) = 0,\\
&\frac{du_{4}}{dt} = -u_{2}(4v_{1} + 1) + 4\cos(v_{1}),\hspace{16.7mm}  u_{4}(0) = 1,\\[1.9mm]
&g_{1} = 4x\cos(v_{1}) + t u_{2}^{2} - 4(v_{1} - t^{2}) = 0,\hspace{7.3mm} v_{1}(0) = 1,
\end{split}
\end{equation}
where sets of differential variables $\mathbf{u} = [x,\, y,\, \dot{x},\, \dot{y}]^{T}$ and algebraic variables $\mathbf{v} = [z]$ were defined. The vector of algebraic equations (constraints) in the case of this DAE system consisted of only one equation $g_{1} = 0$. The results also present dependencies for the error in fulfilling of ``algebraic equation'' $g_{2} = v_{1} - t^{2} - t = 0$, which represents the accuracy of the obtained algebraic variable $v_{1}$. The total number of variables in the chosen problem is $D = 5$, of which $D_{\rm u} = 4$ are differential variables and $D_{\rm v} = 1$ is an algebraic variable. The full-component exact analytical solution of the problem was written in the following form:
\begin{equation}
\begin{split}
&\mathbf{u}^{\rm ex} = \left[
\begin{array}{l}
t\cos(t^{2} + t)\\
2\sin(t^{2} + t)\\
\cos(t^{2} + t) - t(2t + 1)\sin(t^{2} + t)\\
2(2t + 1)\cos(t^{2} + t)
\end{array}
\right],\\
&\mathbf{v}^{\rm ex} = \Big[t^{2} + t\Big].
\end{split}
\end{equation}
The domain of definition $[0,\, 1]$ of the desired functions $\mathbf{u}$ and $\mathbf{v}$ was discretized into $L = 8$, $10$, $12$, $14$, $16$, $18$ discretization domains $\Omega_{n}$ with the same discretization steps ${\Delta t}_{n} \equiv {\Delta t} = 1/(L-1)$, which was done to be able to calculate the empirical convergence orders $p$. Therefore, the empirical convergence orders $p$ were calculated by least squares approximation of the dependence of the global error $e$ on the grid discretization step ${\Delta t}$ at 6 data points.

\begin{figure}[h!]
\captionsetup[subfigure]{%
	position=bottom,
	font+=smaller,
	textfont=normalfont,
	singlelinecheck=off,
	justification=raggedright
}
\centering
\begin{subfigure}{0.320\textwidth}
\includegraphics[width=\textwidth]{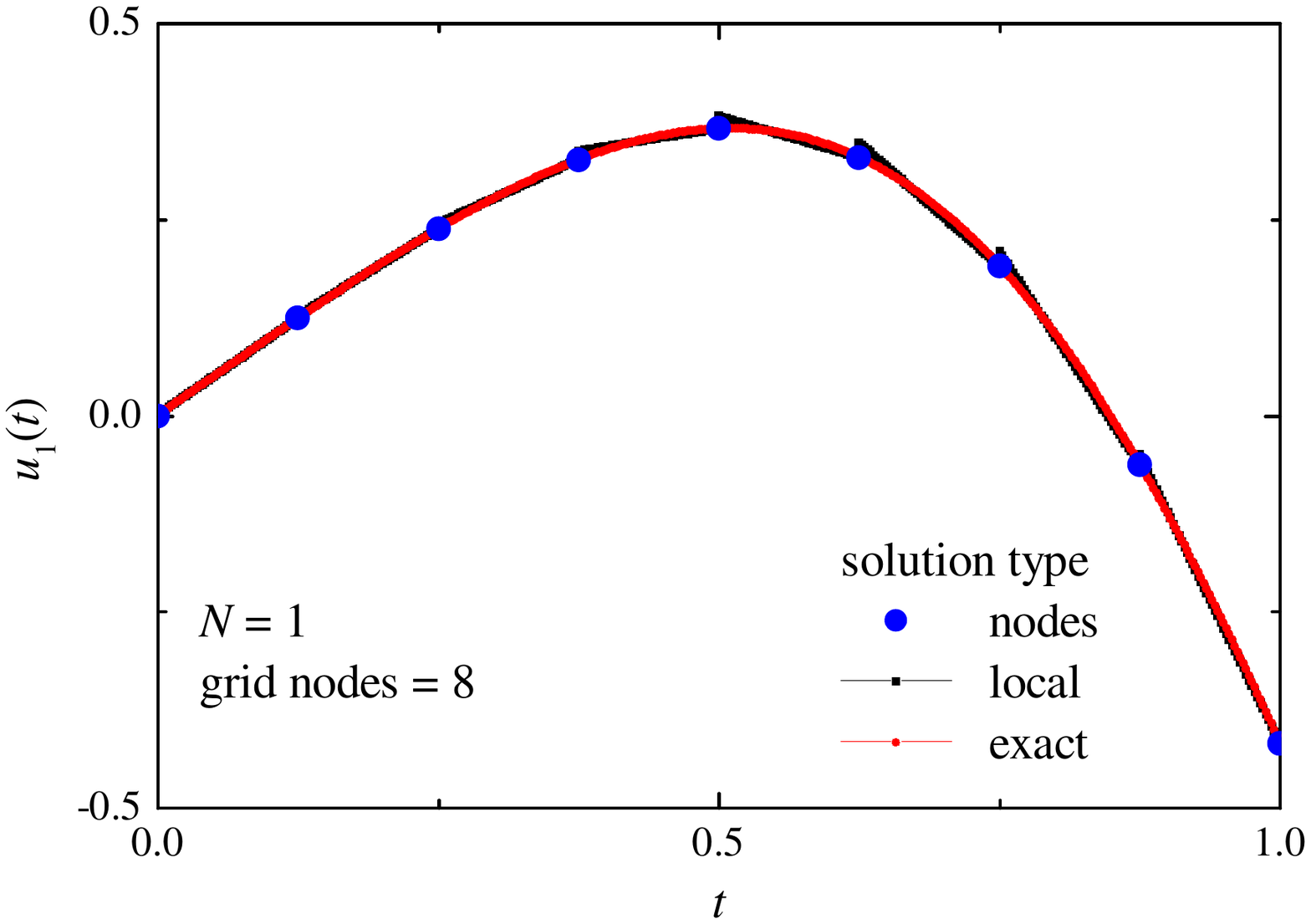}
\vspace{-8mm}\caption{\label{fig:hess_1_sol_uv:a1}}
\end{subfigure}
\begin{subfigure}{0.320\textwidth}
\includegraphics[width=\textwidth]{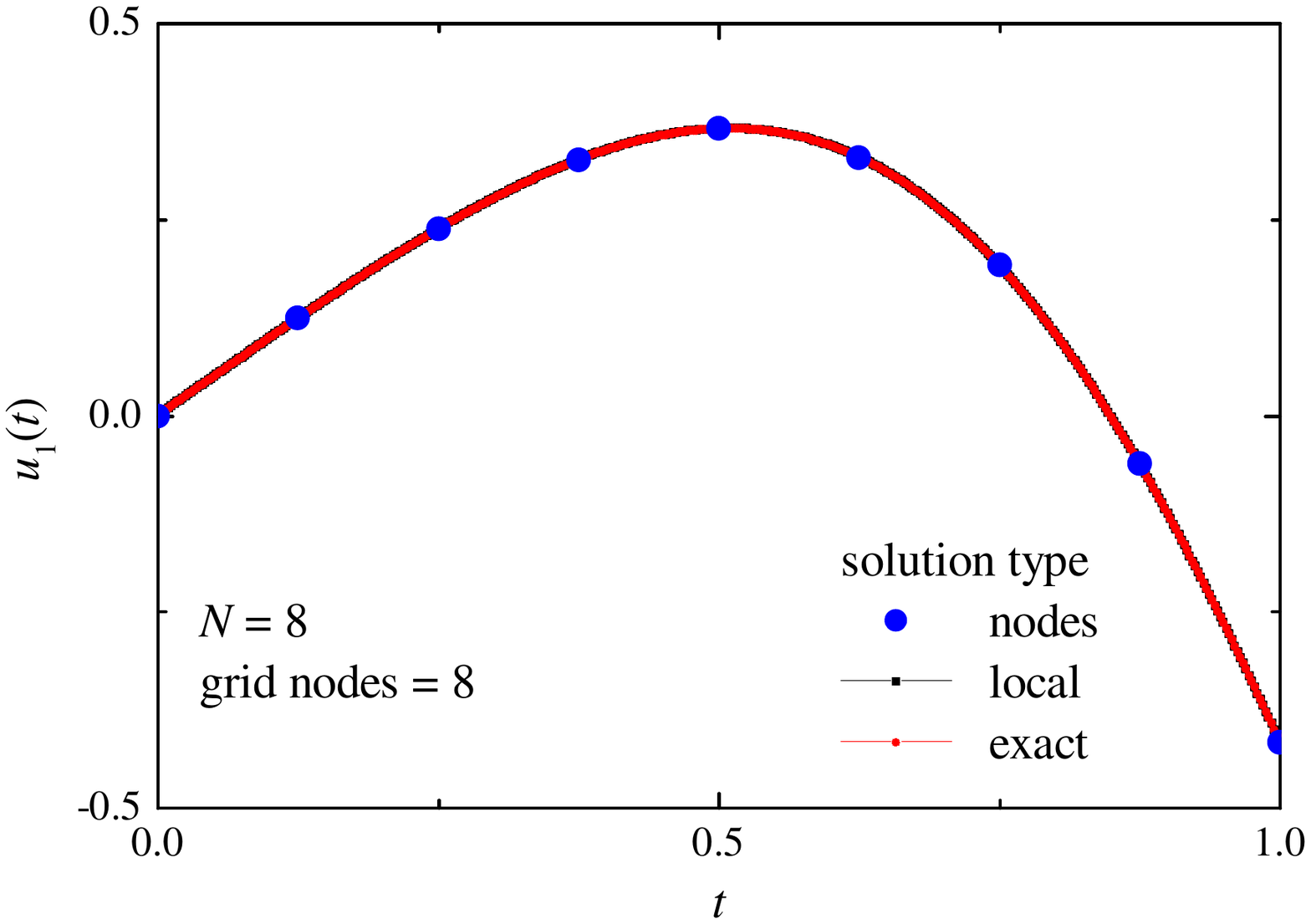}
\vspace{-8mm}\caption{\label{fig:hess_1_sol_uv:a2}}
\end{subfigure}
\begin{subfigure}{0.320\textwidth}
\includegraphics[width=\textwidth]{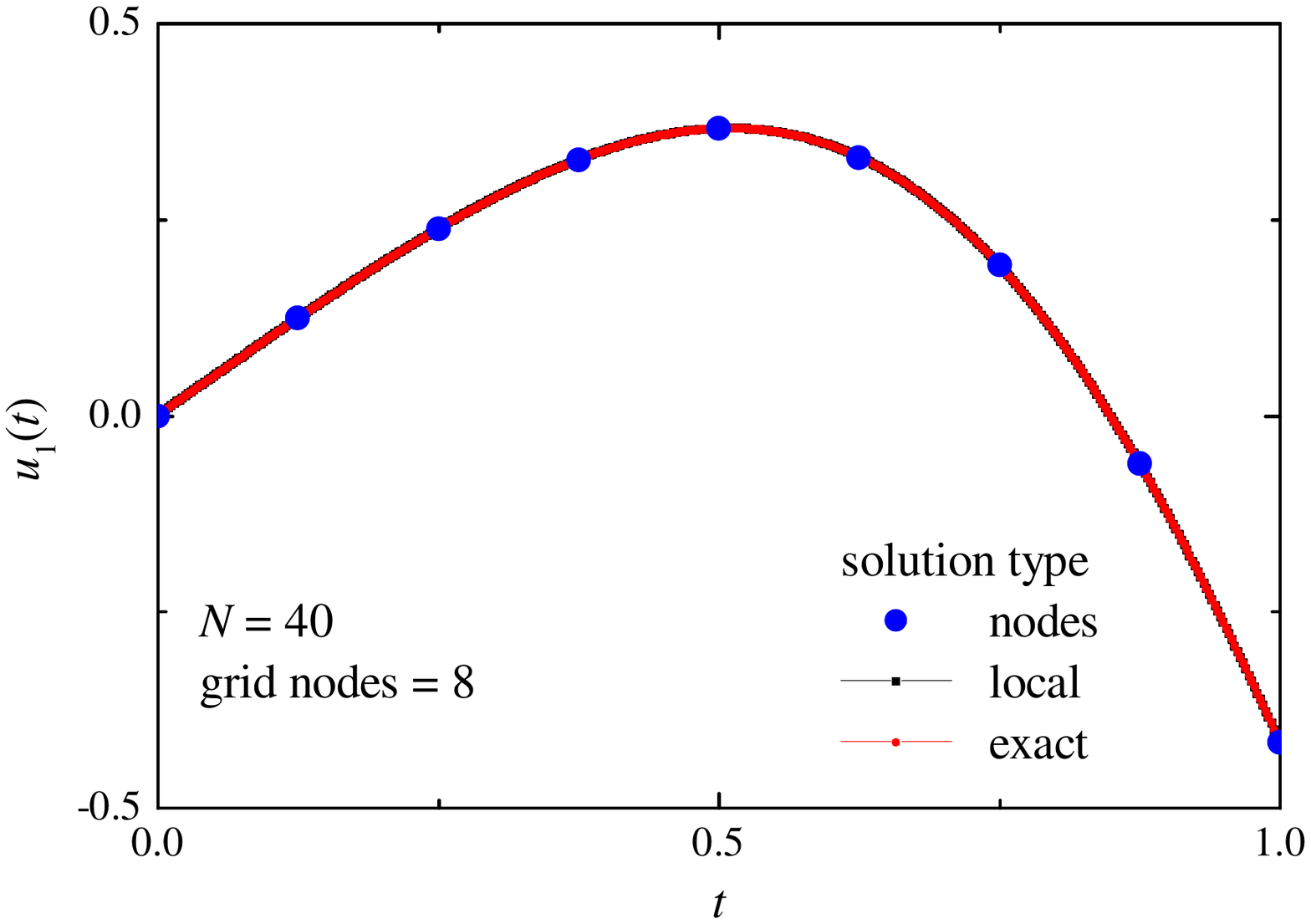}
\vspace{-8mm}\caption{\label{fig:hess_1_sol_uv:a3}}
\end{subfigure}\\
\begin{subfigure}{0.320\textwidth}
\includegraphics[width=\textwidth]{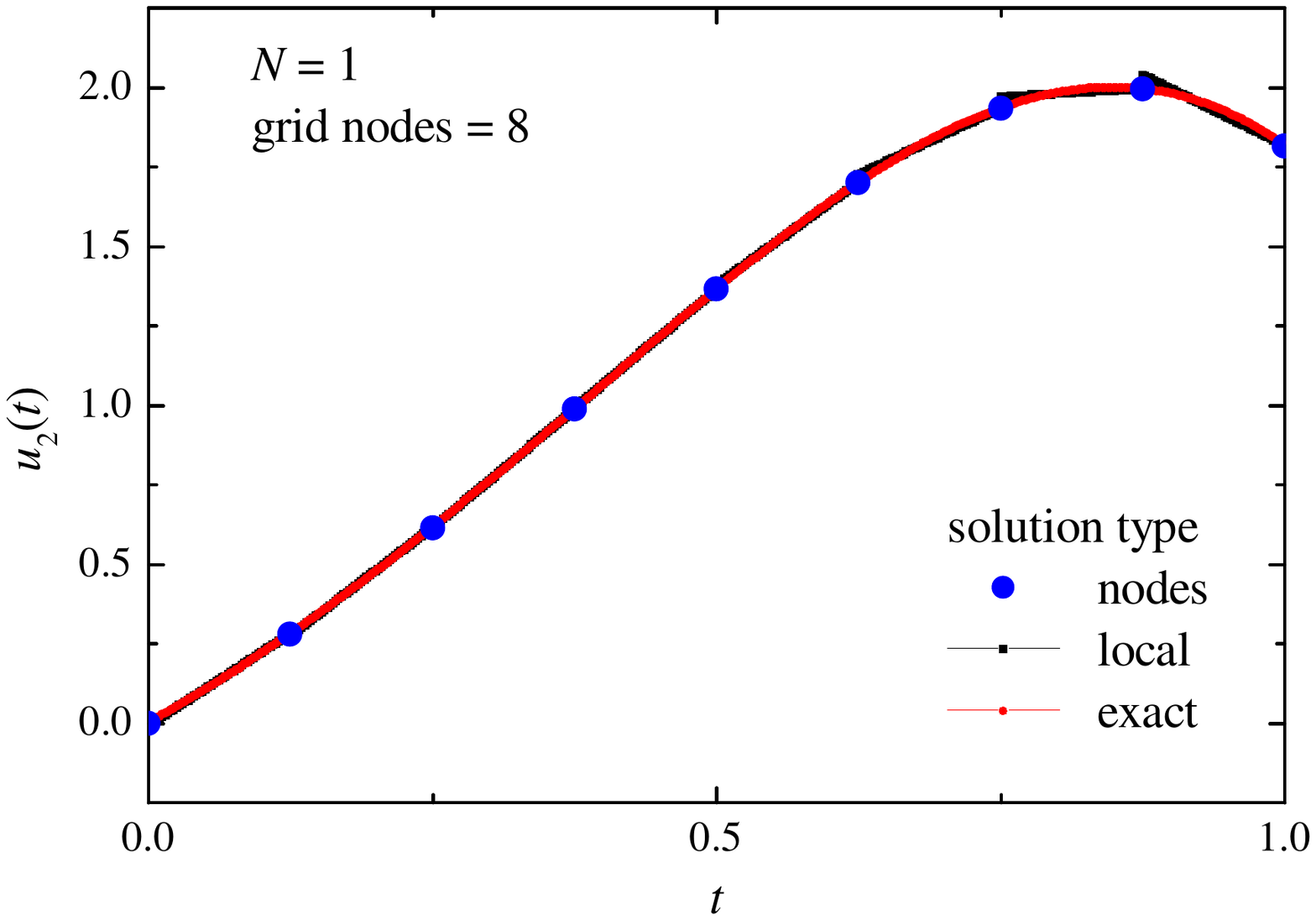}
\vspace{-8mm}\caption{\label{fig:hess_1_sol_uv:b1}}
\end{subfigure}
\begin{subfigure}{0.320\textwidth}
\includegraphics[width=\textwidth]{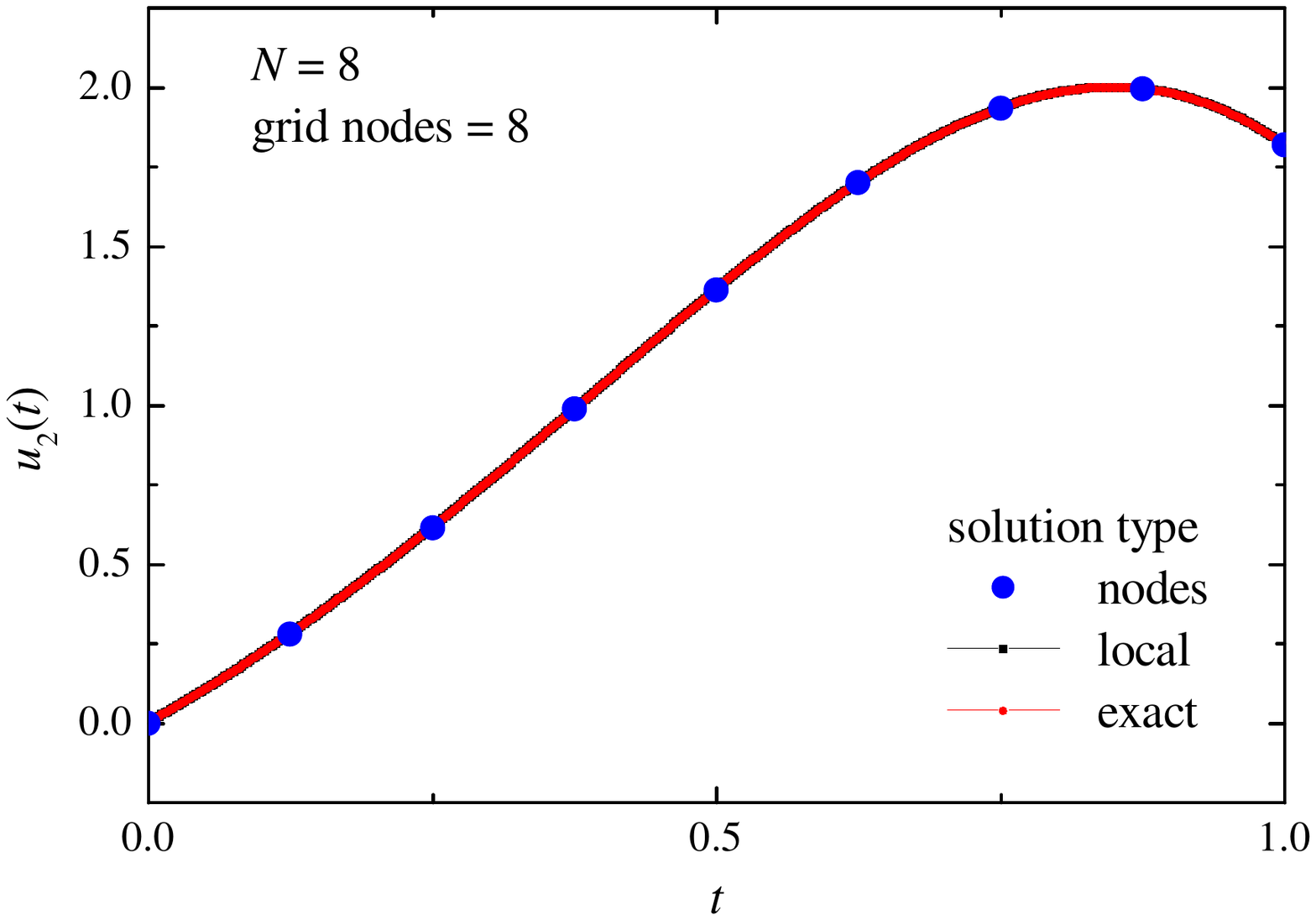}
\vspace{-8mm}\caption{\label{fig:hess_1_sol_uv:b2}}
\end{subfigure}
\begin{subfigure}{0.320\textwidth}
\includegraphics[width=\textwidth]{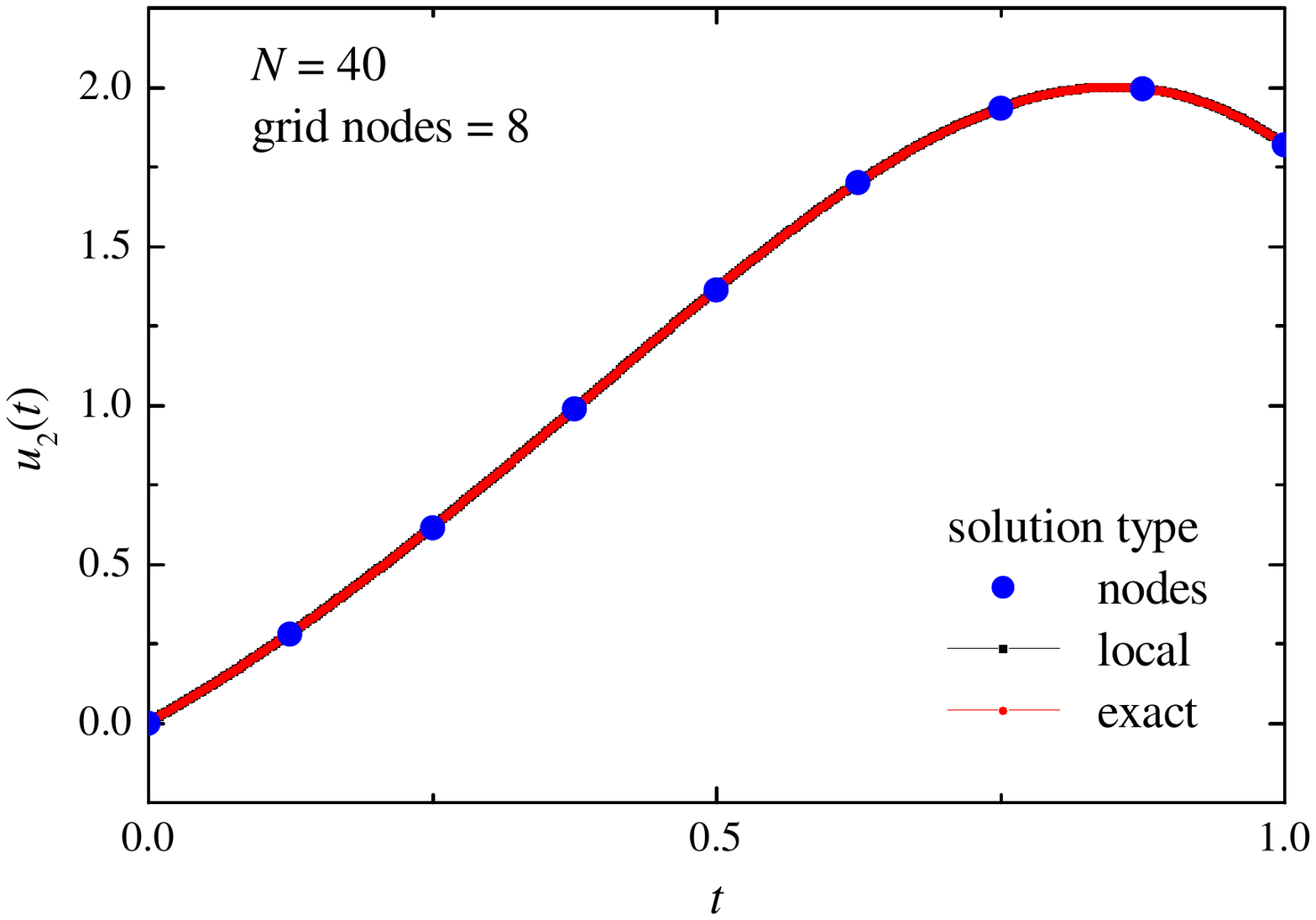}
\vspace{-8mm}\caption{\label{fig:hess_1_sol_uv:b3}}
\end{subfigure}\\
\begin{subfigure}{0.320\textwidth}
\includegraphics[width=\textwidth]{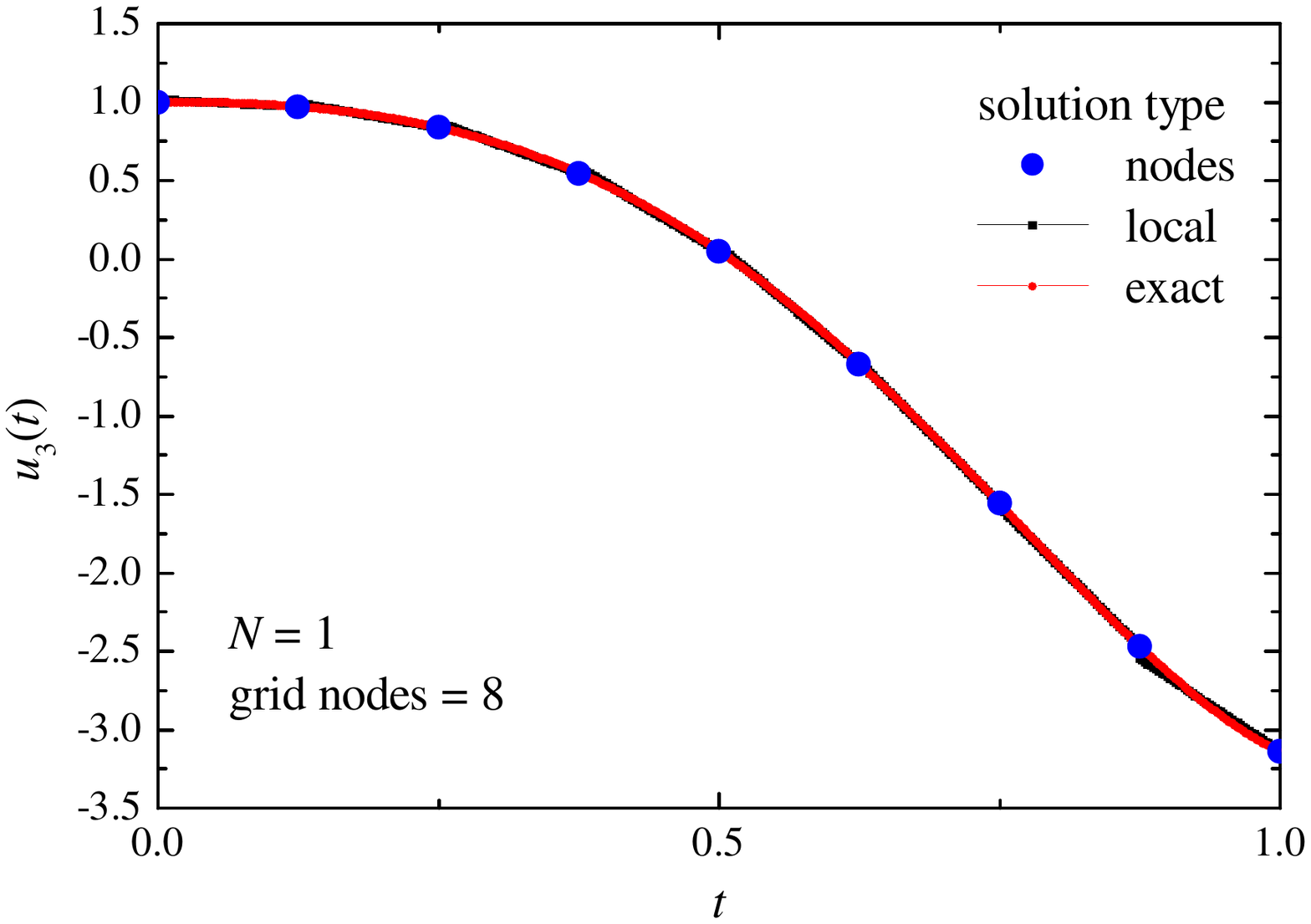}
\vspace{-8mm}\caption{\label{fig:hess_1_sol_uv:c1}}
\end{subfigure}
\begin{subfigure}{0.320\textwidth}
\includegraphics[width=\textwidth]{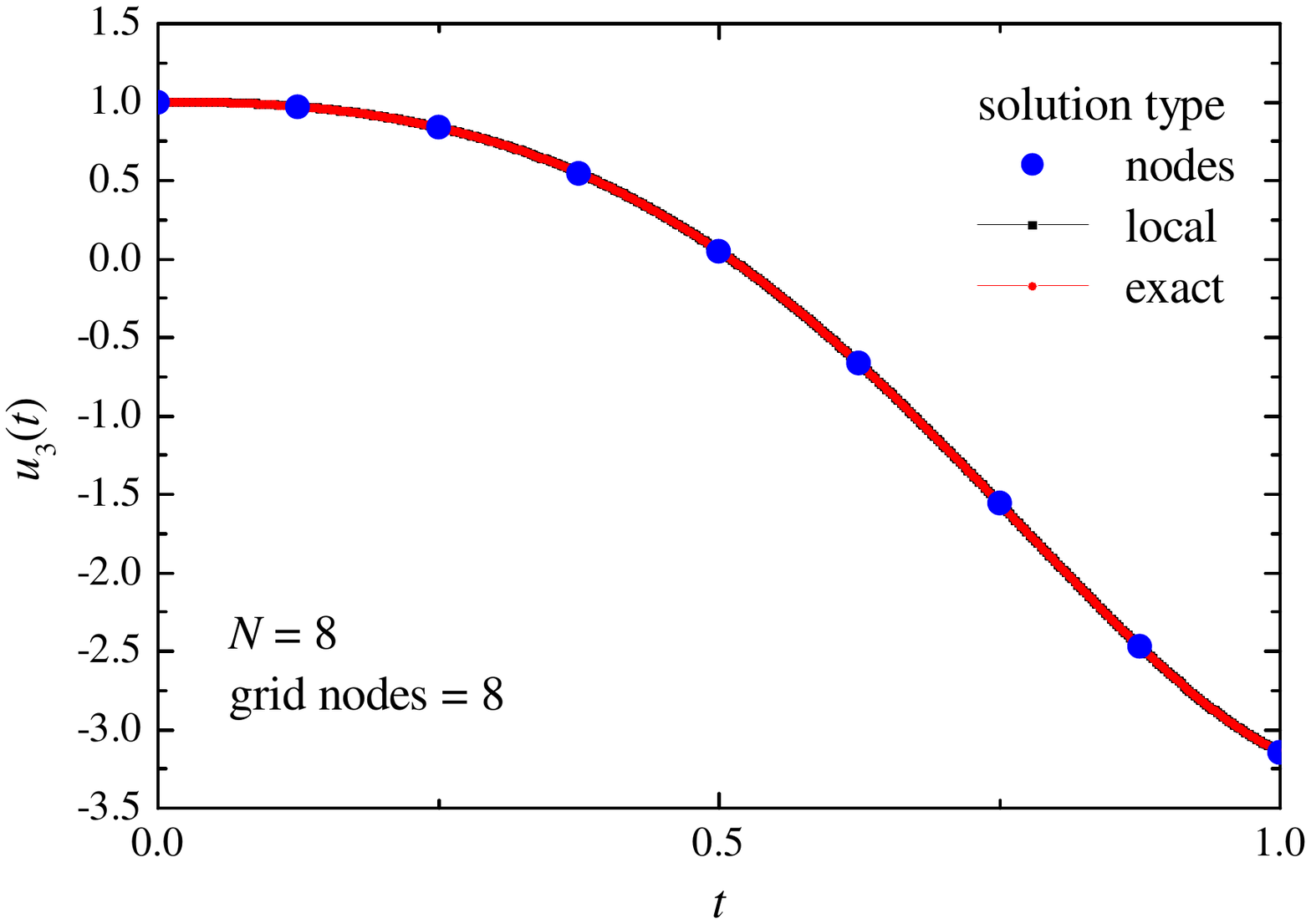}
\vspace{-8mm}\caption{\label{fig:hess_1_sol_uv:c2}}
\end{subfigure}
\begin{subfigure}{0.320\textwidth}
\includegraphics[width=\textwidth]{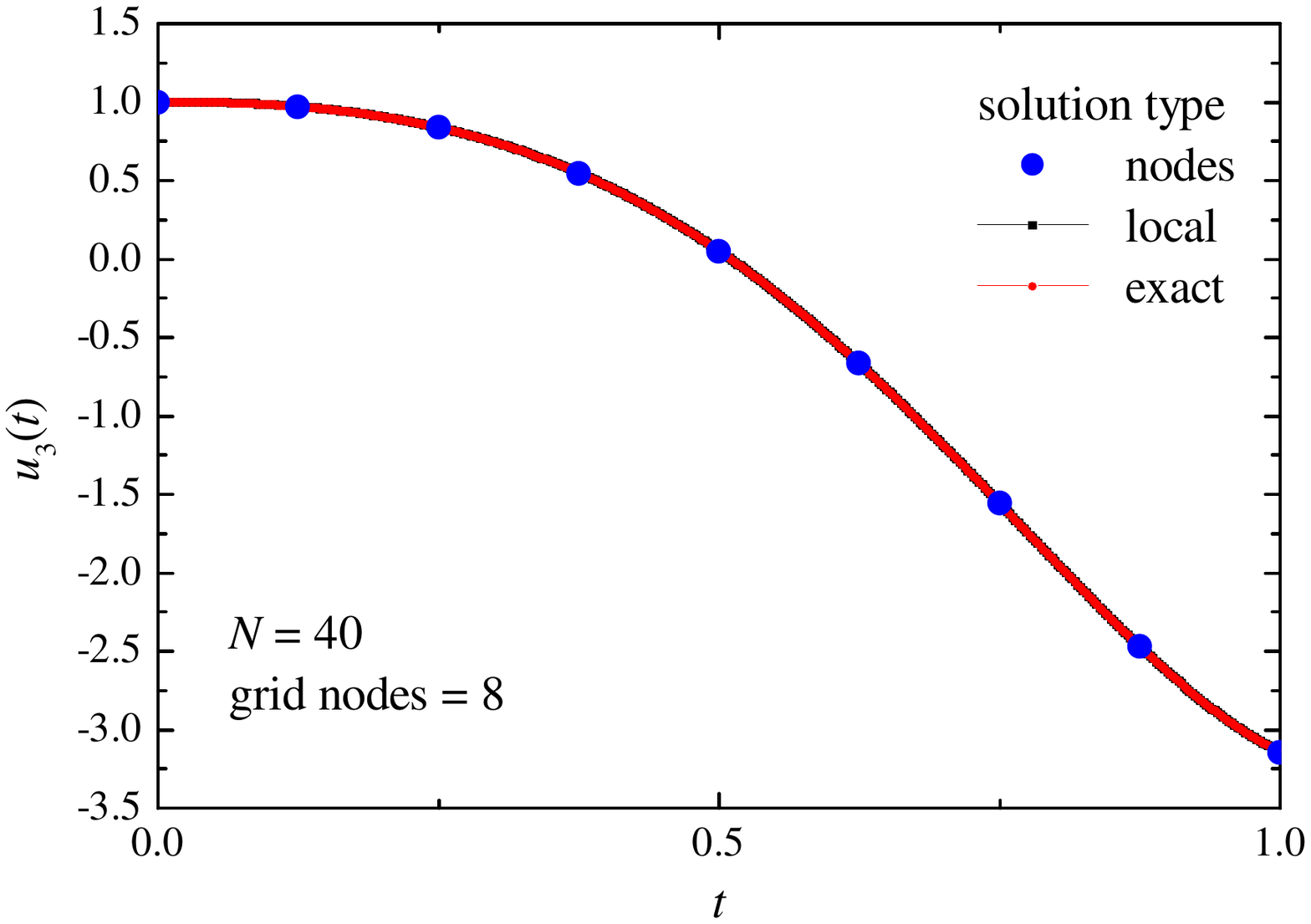}
\vspace{-8mm}\caption{\label{fig:hess_1_sol_uv:c3}}
\end{subfigure}\\
\begin{subfigure}{0.320\textwidth}
\includegraphics[width=\textwidth]{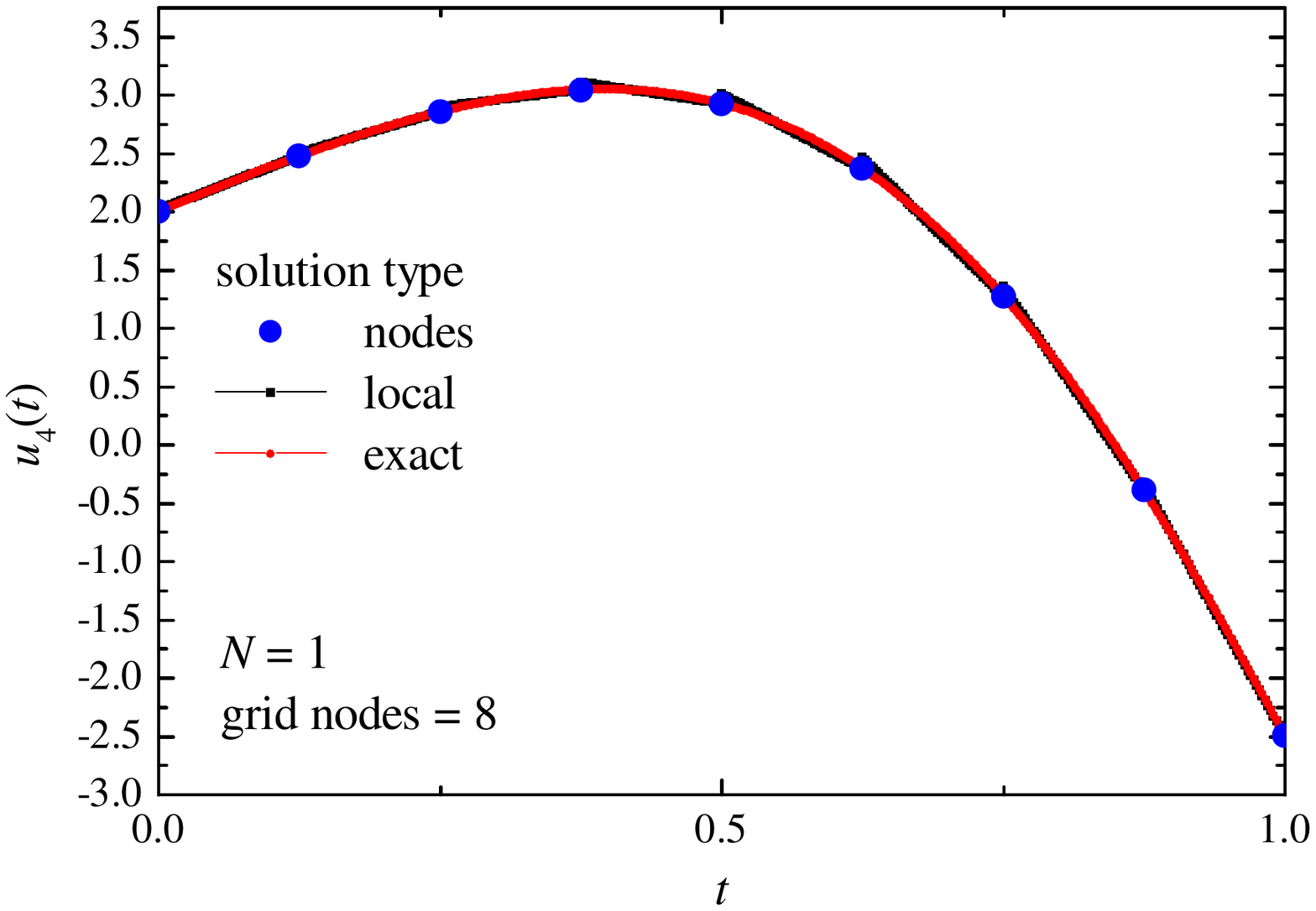}
\vspace{-8mm}\caption{\label{fig:hess_1_sol_uv:d1}}
\end{subfigure}
\begin{subfigure}{0.320\textwidth}
\includegraphics[width=\textwidth]{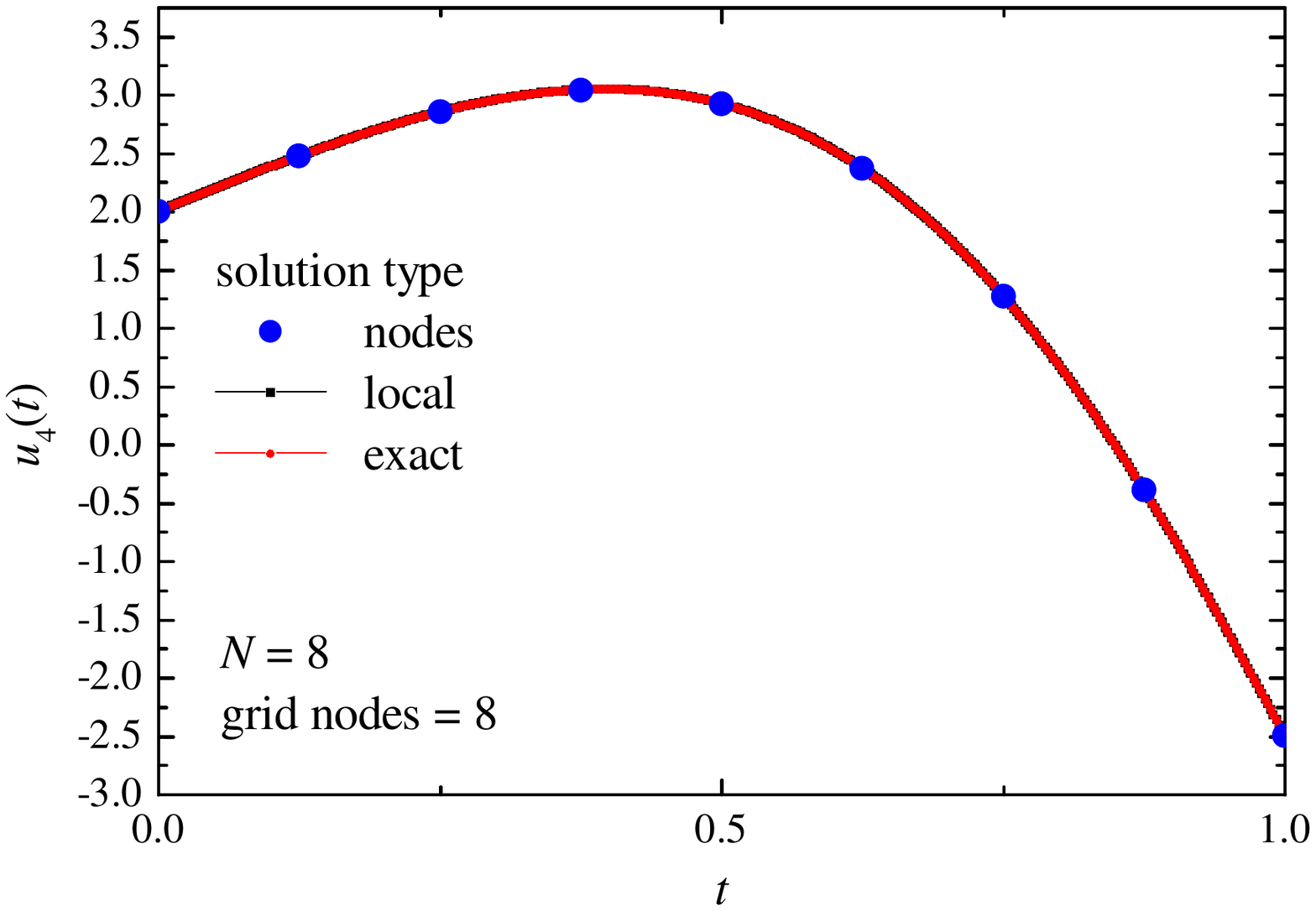}
\vspace{-8mm}\caption{\label{fig:hess_1_sol_uv:d2}}
\end{subfigure}
\begin{subfigure}{0.320\textwidth}
\includegraphics[width=\textwidth]{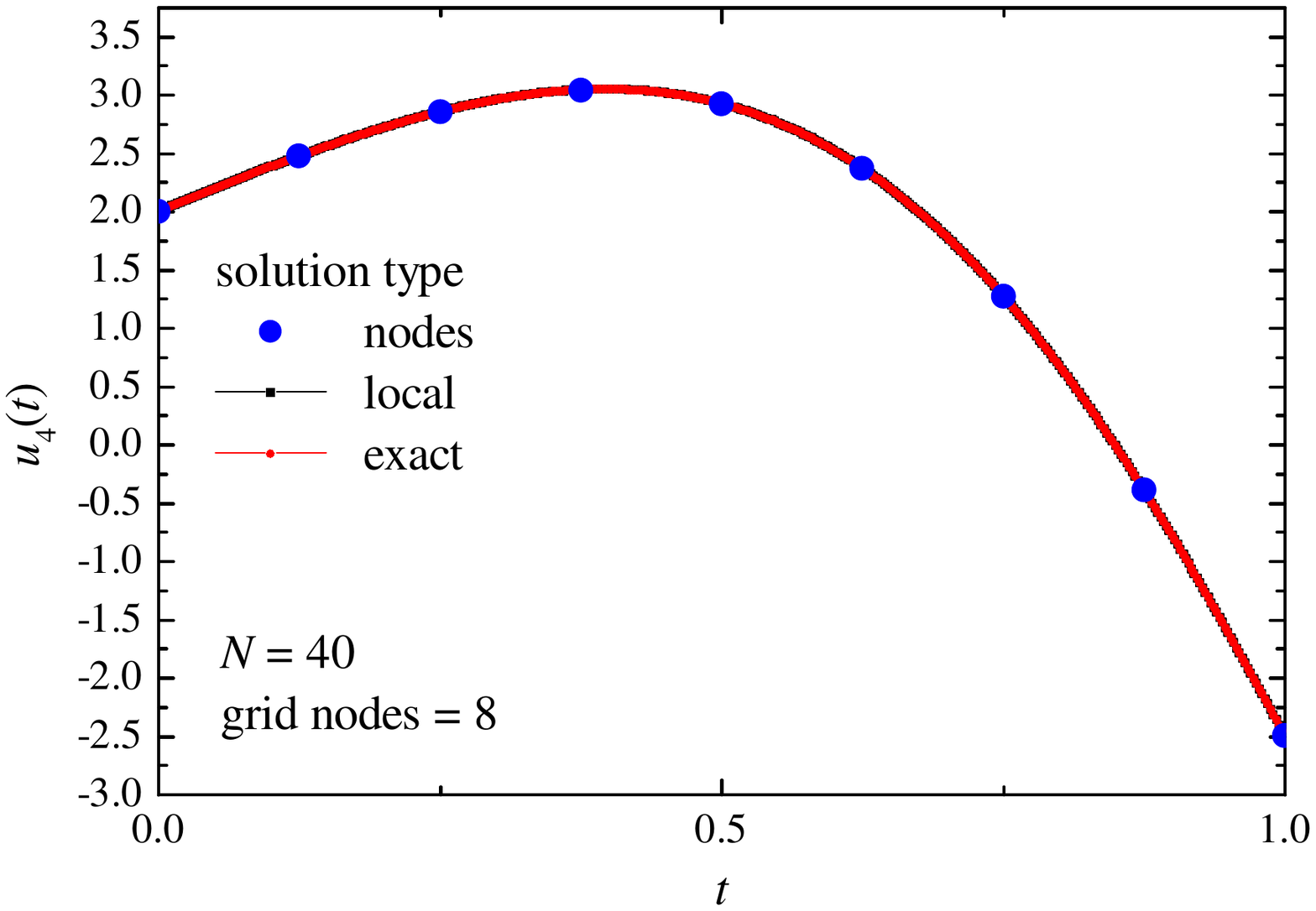}
\vspace{-8mm}\caption{\label{fig:hess_1_sol_uv:d3}}
\end{subfigure}\\
\begin{subfigure}{0.320\textwidth}
\includegraphics[width=\textwidth]{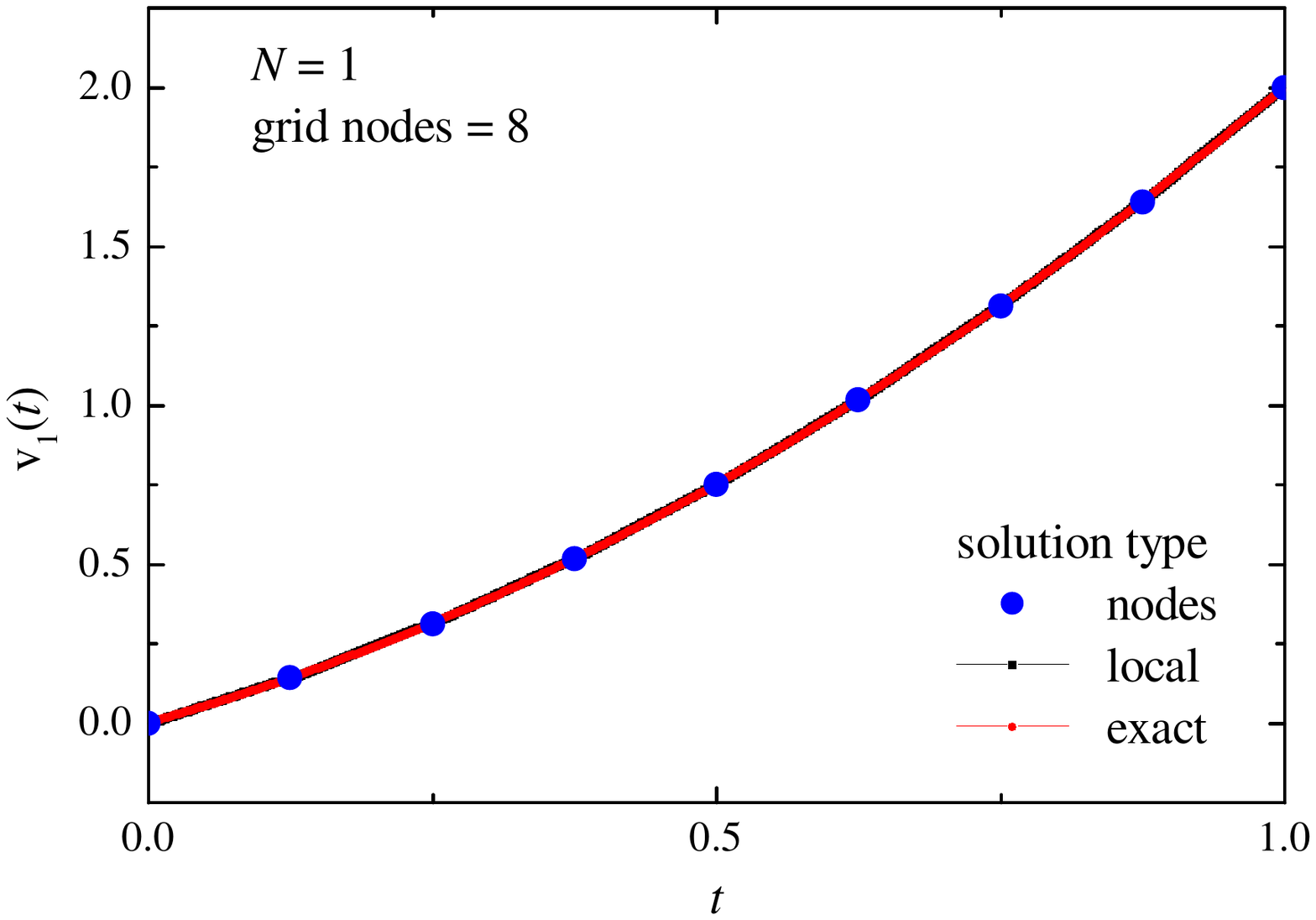}
\vspace{-8mm}\caption{\label{fig:hess_1_sol_uv:e1}}
\end{subfigure}
\begin{subfigure}{0.320\textwidth}
\includegraphics[width=\textwidth]{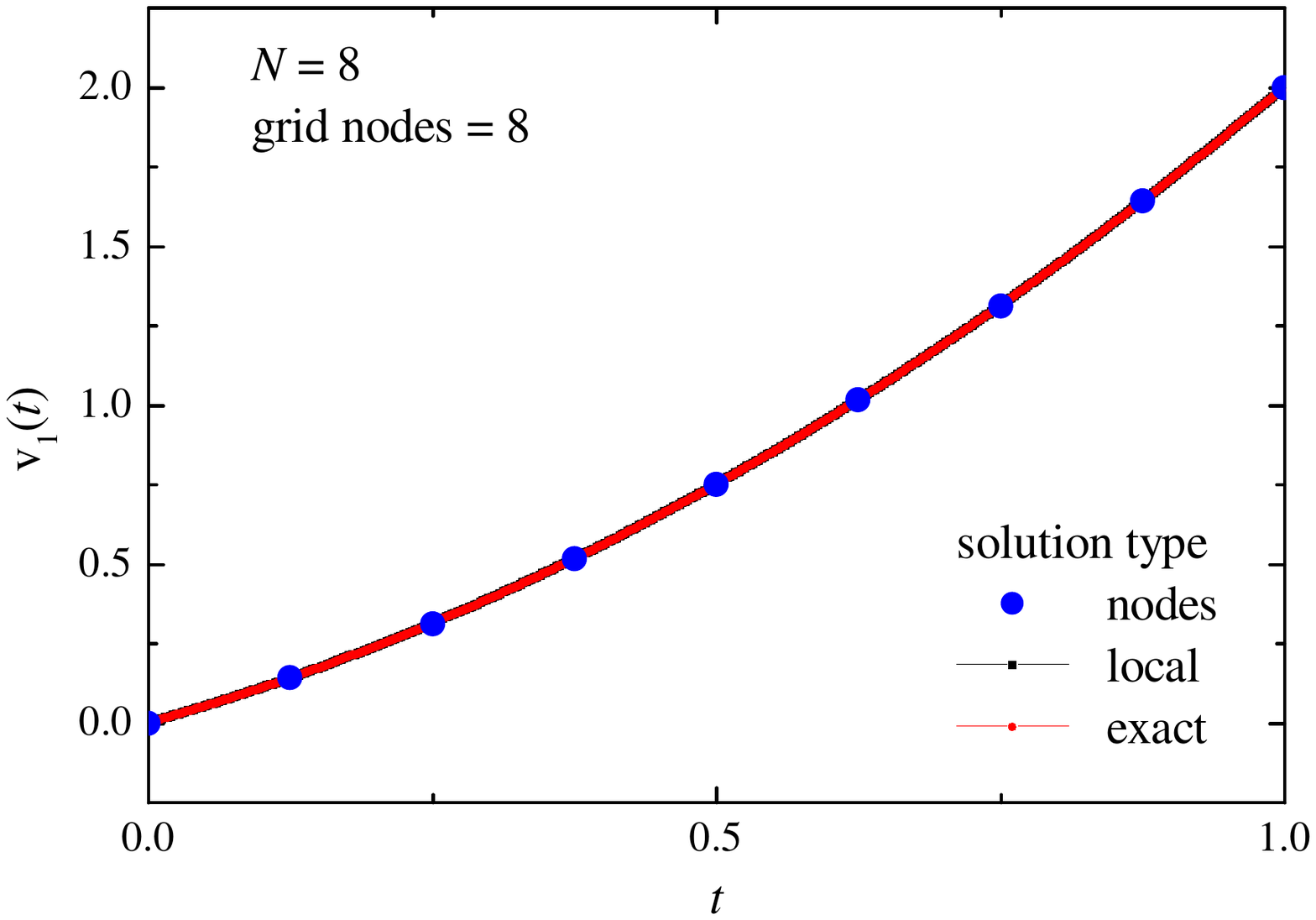}
\vspace{-8mm}\caption{\label{fig:hess_1_sol_uv:e2}}
\end{subfigure}
\begin{subfigure}{0.320\textwidth}
\includegraphics[width=\textwidth]{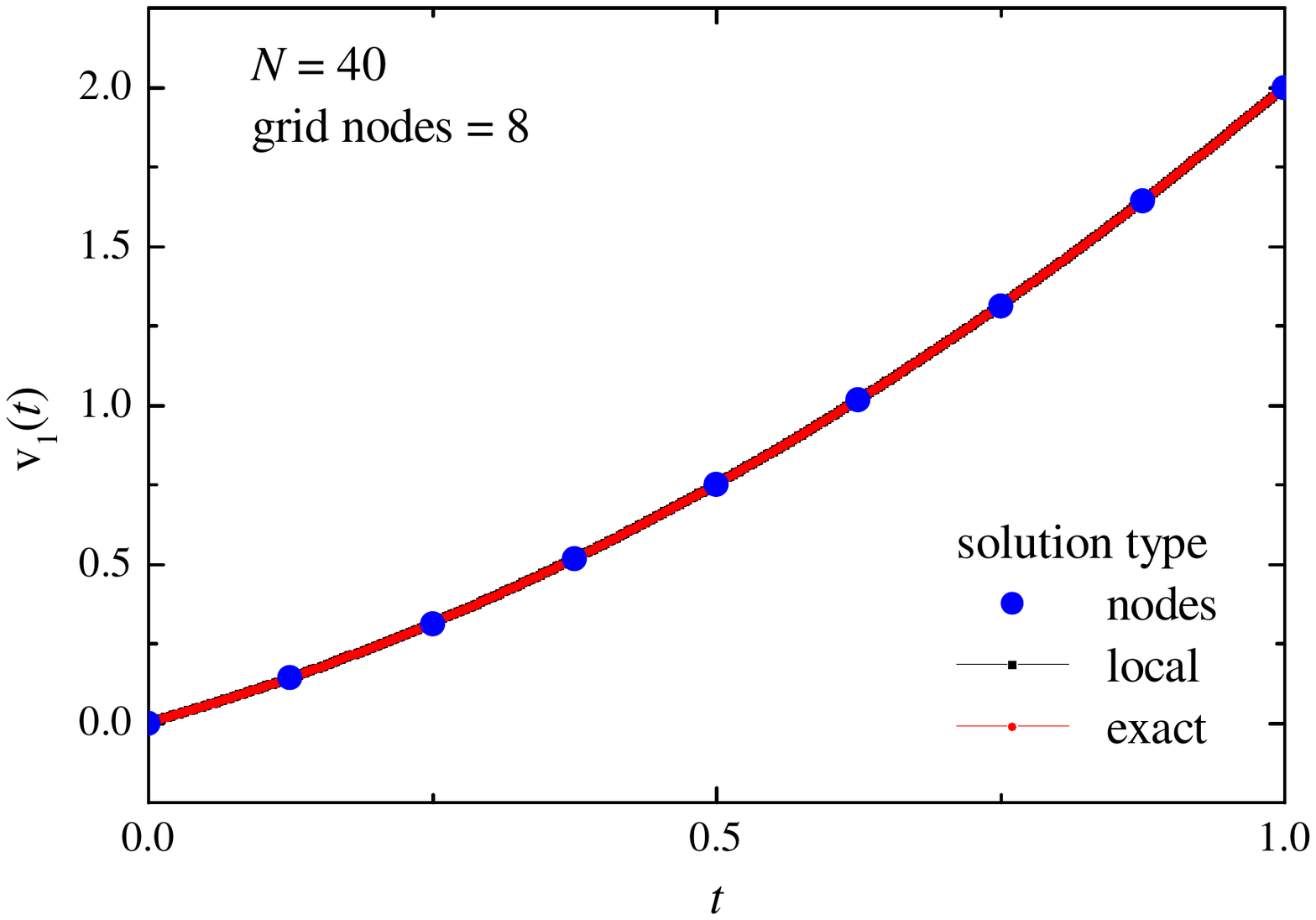}
\vspace{-8mm}\caption{\label{fig:hess_1_sol_uv:e3}}
\end{subfigure}\\
\caption{%
Numerical solution of the problem (\ref{eq:hess_dae_ind_1}). Comparison of the solution at nodes $\mathbf{u}_{n}$, the local solution $\mathbf{u}_{L}(t)$ and the exact solution $\mathbf{u}^{\rm ex}(t)$ for components $u_{1}$ (\subref{fig:hess_1_sol_uv:a1}, \subref{fig:hess_1_sol_uv:a2}, \subref{fig:hess_1_sol_uv:a3}), $u_{2}$ (\subref{fig:hess_1_sol_uv:b1}, \subref{fig:hess_1_sol_uv:b2}, \subref{fig:hess_1_sol_uv:b3}), $u_{3}$ (\subref{fig:hess_1_sol_uv:c1}, \subref{fig:hess_1_sol_uv:c2}, \subref{fig:hess_1_sol_uv:c3}), $u_{4}$ (\subref{fig:hess_1_sol_uv:d1}, \subref{fig:hess_1_sol_uv:d2}, \subref{fig:hess_1_sol_uv:d3}) and $v_{1}$ (\subref{fig:hess_1_sol_uv:e1}, \subref{fig:hess_1_sol_uv:e2}, \subref{fig:hess_1_sol_uv:e3}), obtained using polynomials with degrees $N = 1$ (\subref{fig:hess_1_sol_uv:a1}, \subref{fig:hess_1_sol_uv:b1}, \subref{fig:hess_1_sol_uv:c1}, \subref{fig:hess_1_sol_uv:d1}, \subref{fig:hess_1_sol_uv:e1}), $N = 8$ (\subref{fig:hess_1_sol_uv:a2}, \subref{fig:hess_1_sol_uv:b2}, \subref{fig:hess_1_sol_uv:c2}, \subref{fig:hess_1_sol_uv:d2}, \subref{fig:hess_1_sol_uv:e2}) and $N = 40$ (\subref{fig:hess_1_sol_uv:a3}, \subref{fig:hess_1_sol_uv:b3}, \subref{fig:hess_1_sol_uv:c3}, \subref{fig:hess_1_sol_uv:d3}, \subref{fig:hess_1_sol_uv:e3}).
}
\label{fig:hess_1_sols_uv}
\end{figure} 

\begin{figure}[h!]
\captionsetup[subfigure]{%
	position=bottom,
	font+=smaller,
	textfont=normalfont,
	singlelinecheck=off,
	justification=raggedright
}
\centering
\begin{subfigure}{0.320\textwidth}
\includegraphics[width=\textwidth]{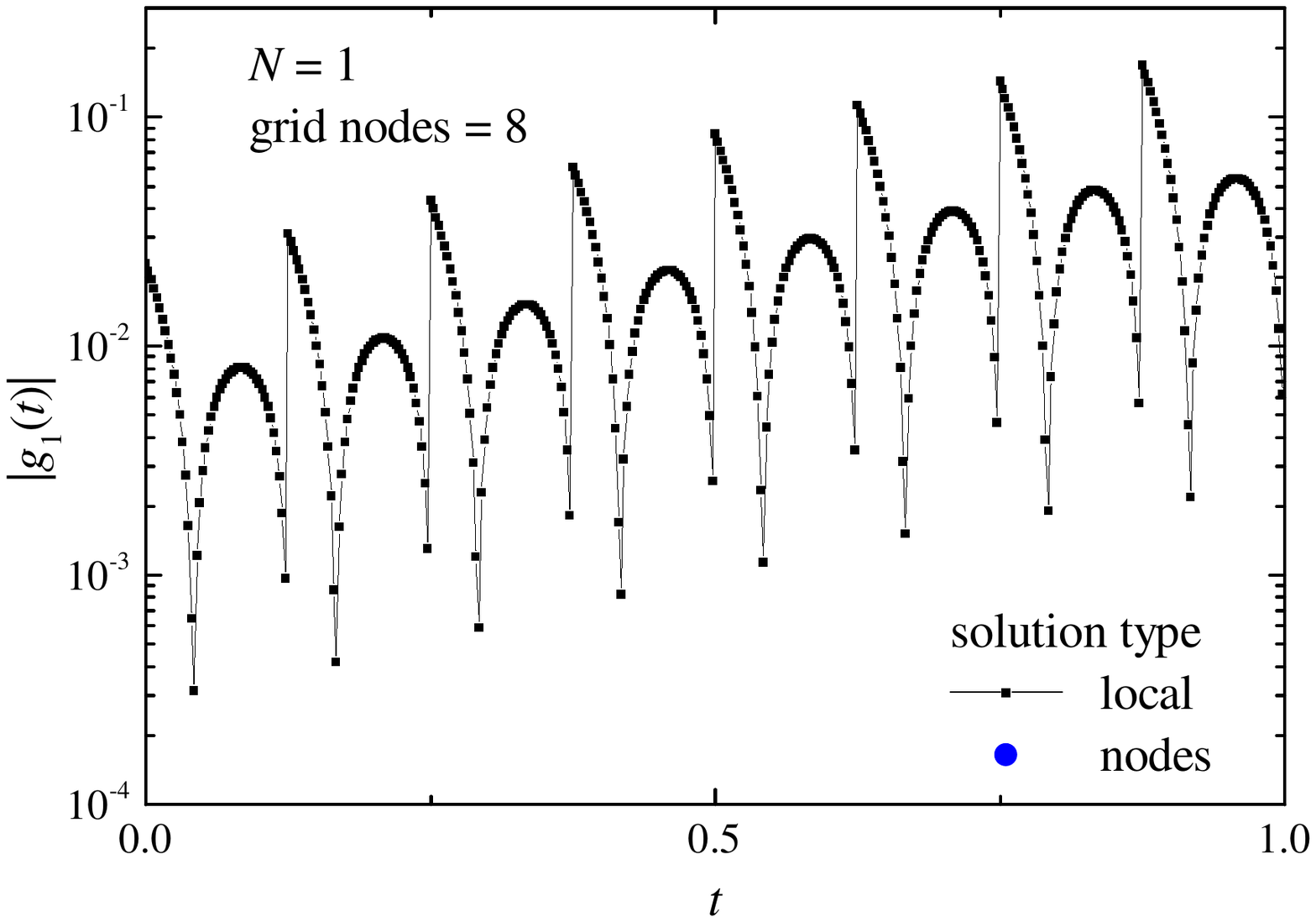}
\vspace{-8mm}\caption{\label{fig:hess_1_sol_g_eps:a1}}
\end{subfigure}
\begin{subfigure}{0.320\textwidth}
\includegraphics[width=\textwidth]{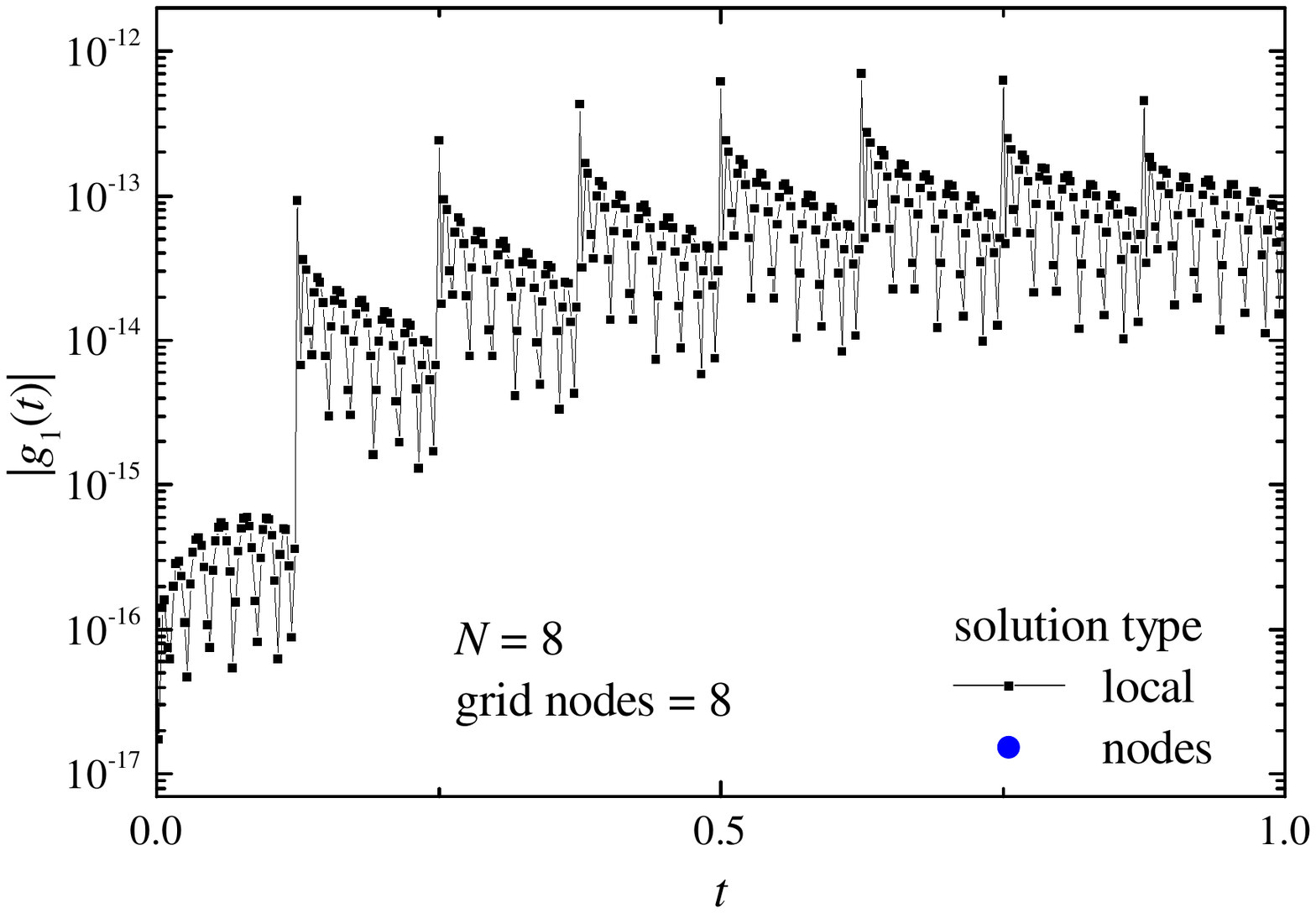}
\vspace{-8mm}\caption{\label{fig:hess_1_sol_g_eps:a2}}
\end{subfigure}
\begin{subfigure}{0.320\textwidth}
\includegraphics[width=\textwidth]{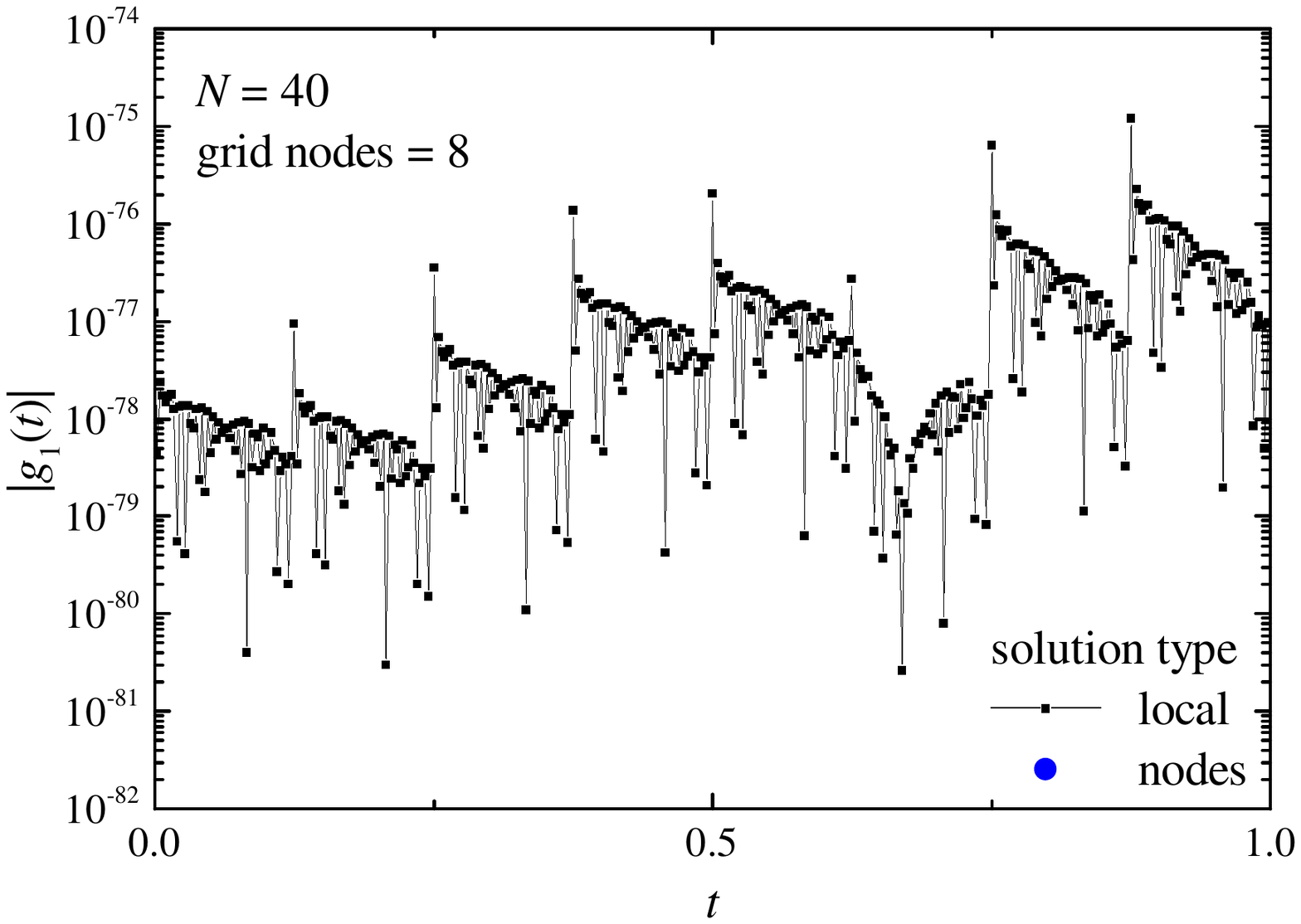}
\vspace{-8mm}\caption{\label{fig:hess_1_sol_g_eps:a3}}
\end{subfigure}\\
\begin{subfigure}{0.320\textwidth}
\includegraphics[width=\textwidth]{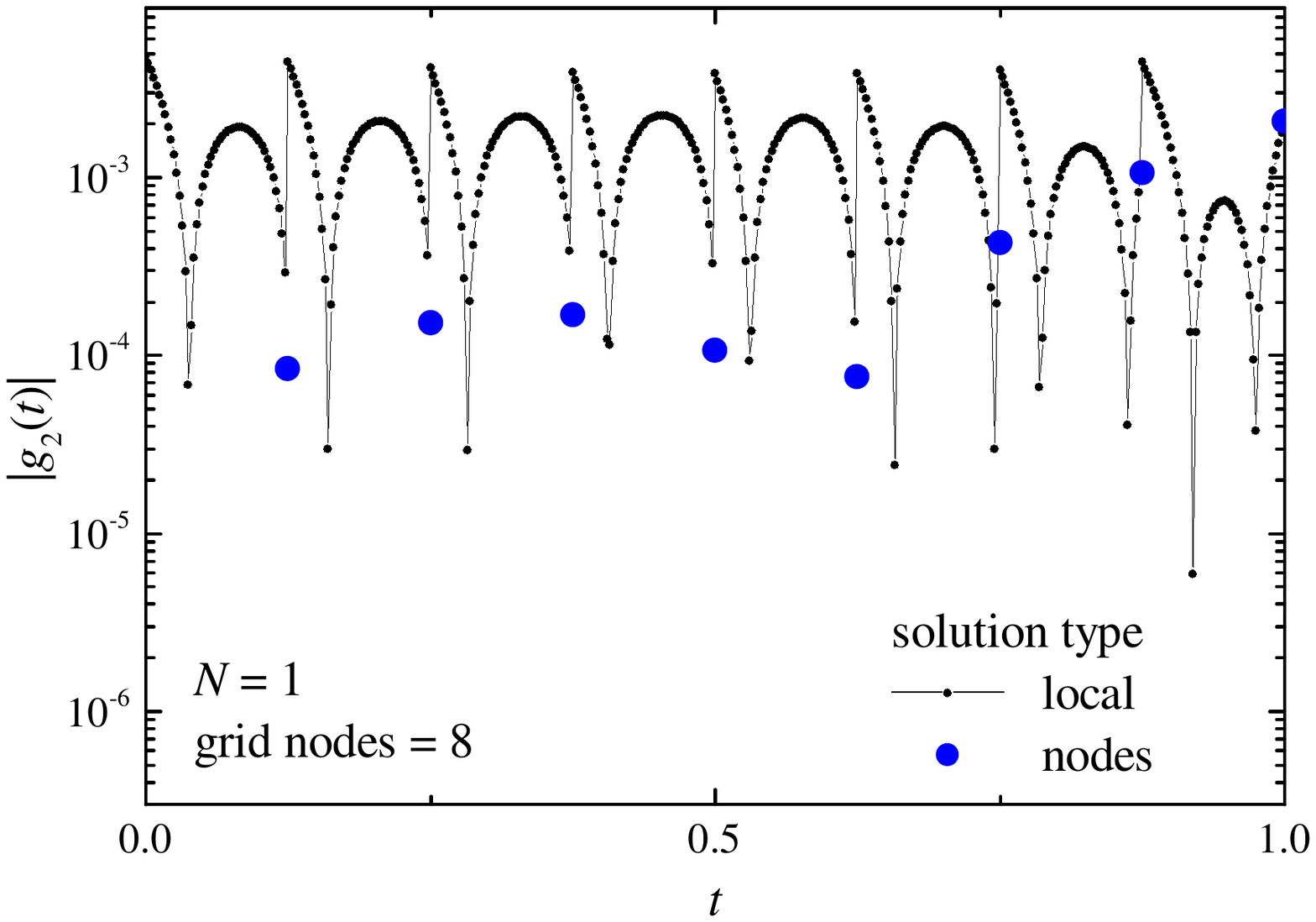}
\vspace{-8mm}\caption{\label{fig:hess_1_sol_g_eps:b1}}
\end{subfigure}
\begin{subfigure}{0.320\textwidth}
\includegraphics[width=\textwidth]{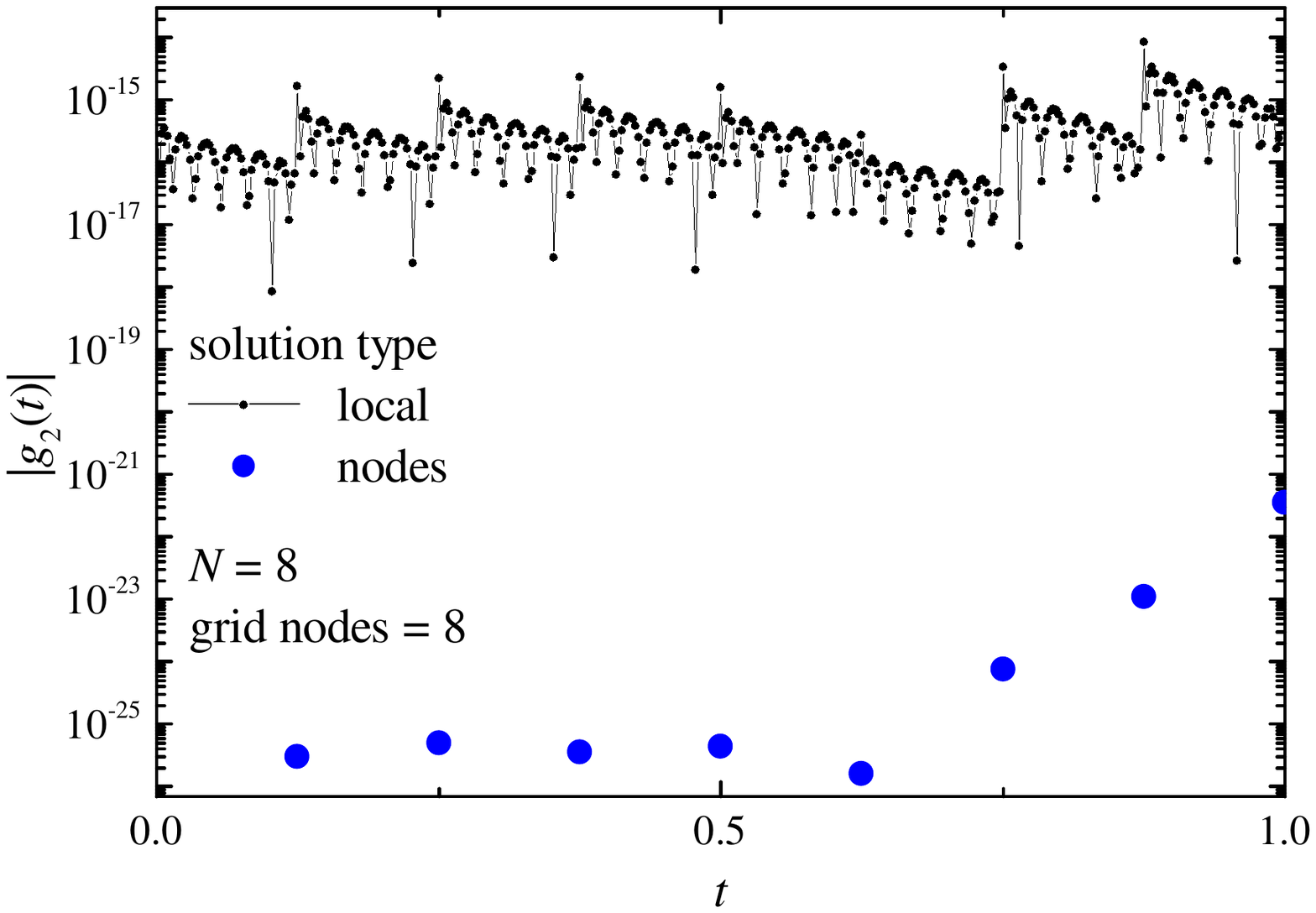}
\vspace{-8mm}\caption{\label{fig:hess_1_sol_g_eps:b2}}
\end{subfigure}
\begin{subfigure}{0.320\textwidth}
\includegraphics[width=\textwidth]{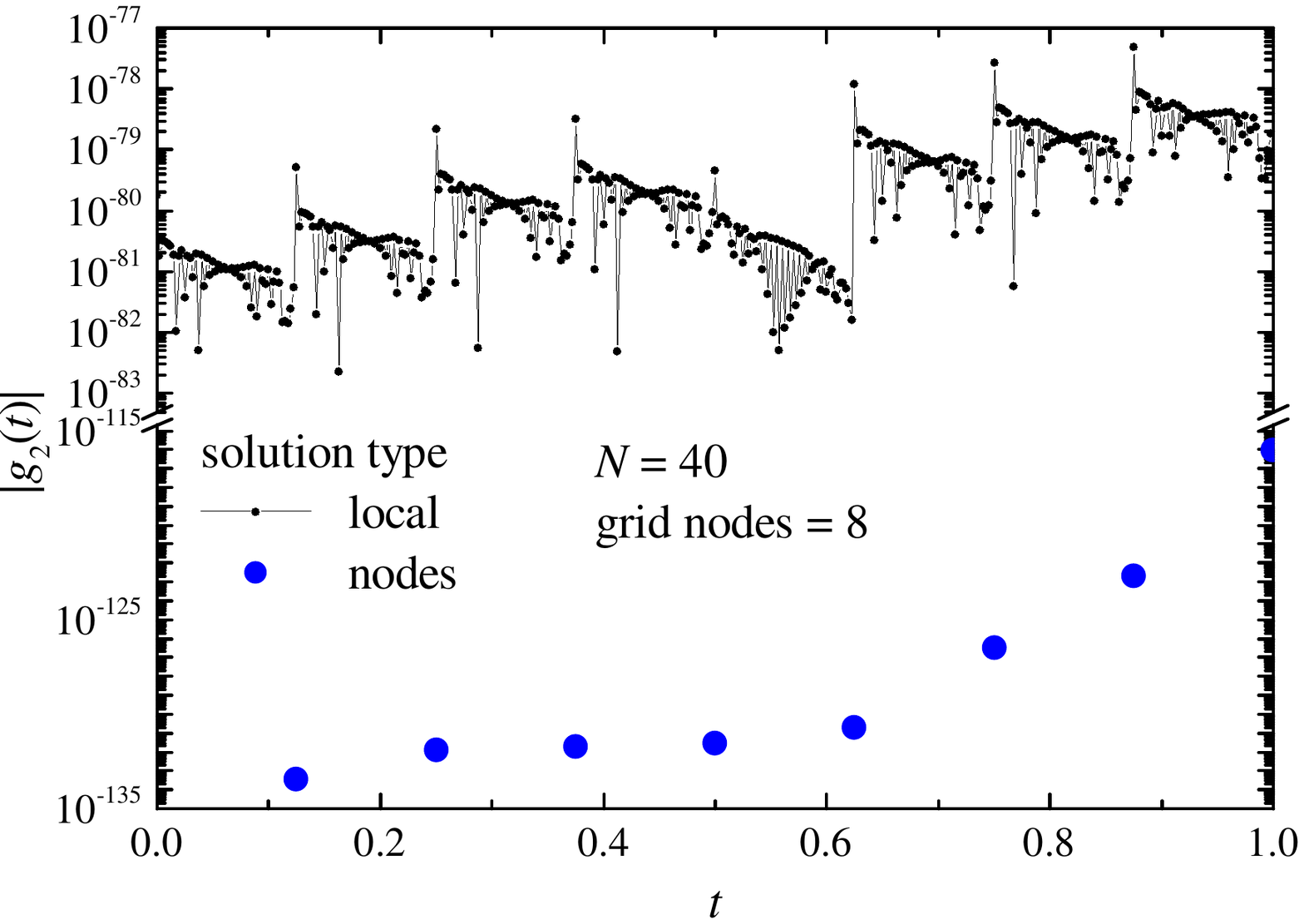}
\vspace{-8mm}\caption{\label{fig:hess_1_sol_g_eps:b3}}
\end{subfigure}\\
\begin{subfigure}{0.320\textwidth}
\includegraphics[width=\textwidth]{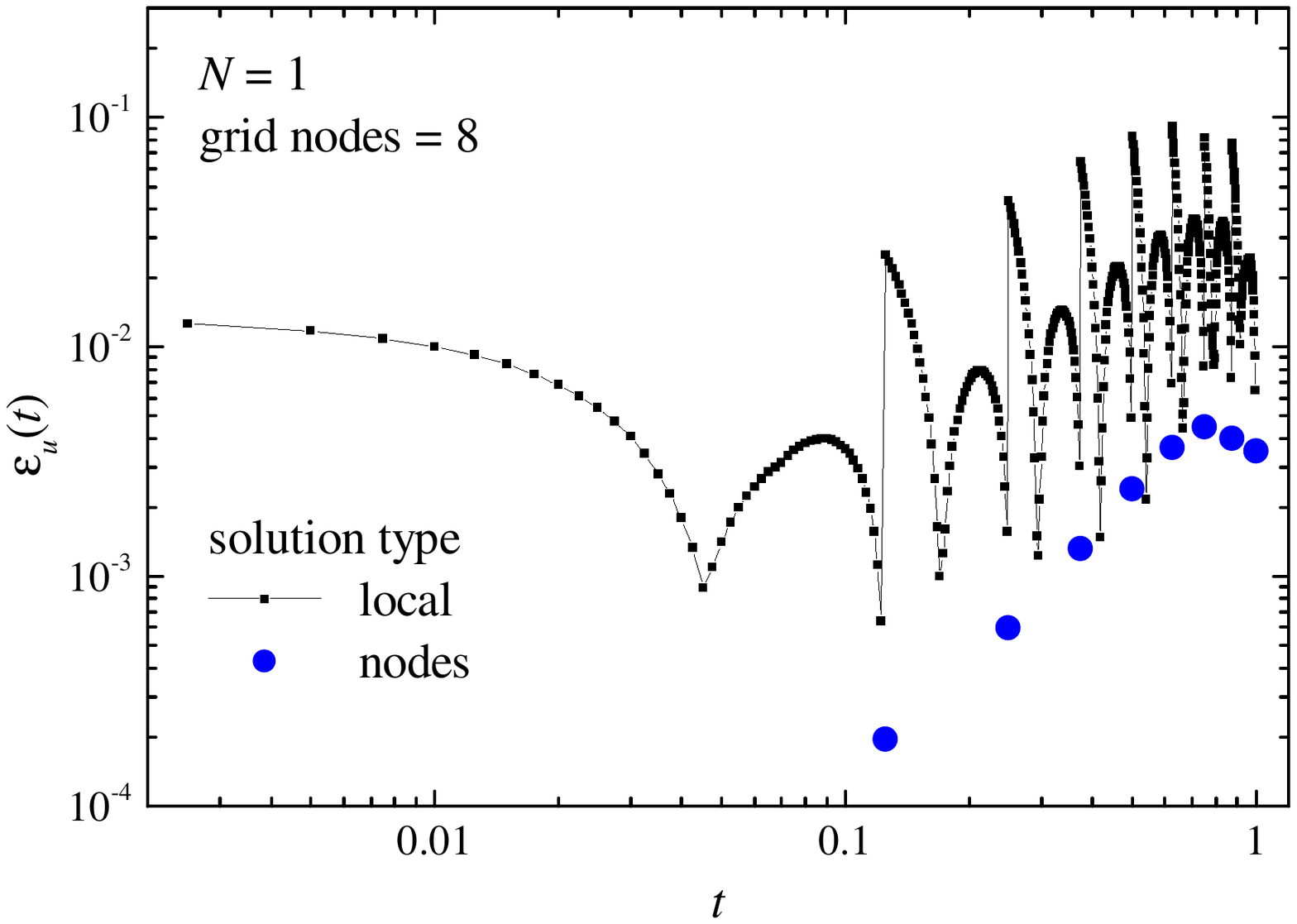}
\vspace{-8mm}\caption{\label{fig:hess_1_sol_g_eps:c1}}
\end{subfigure}
\begin{subfigure}{0.320\textwidth}
\includegraphics[width=\textwidth]{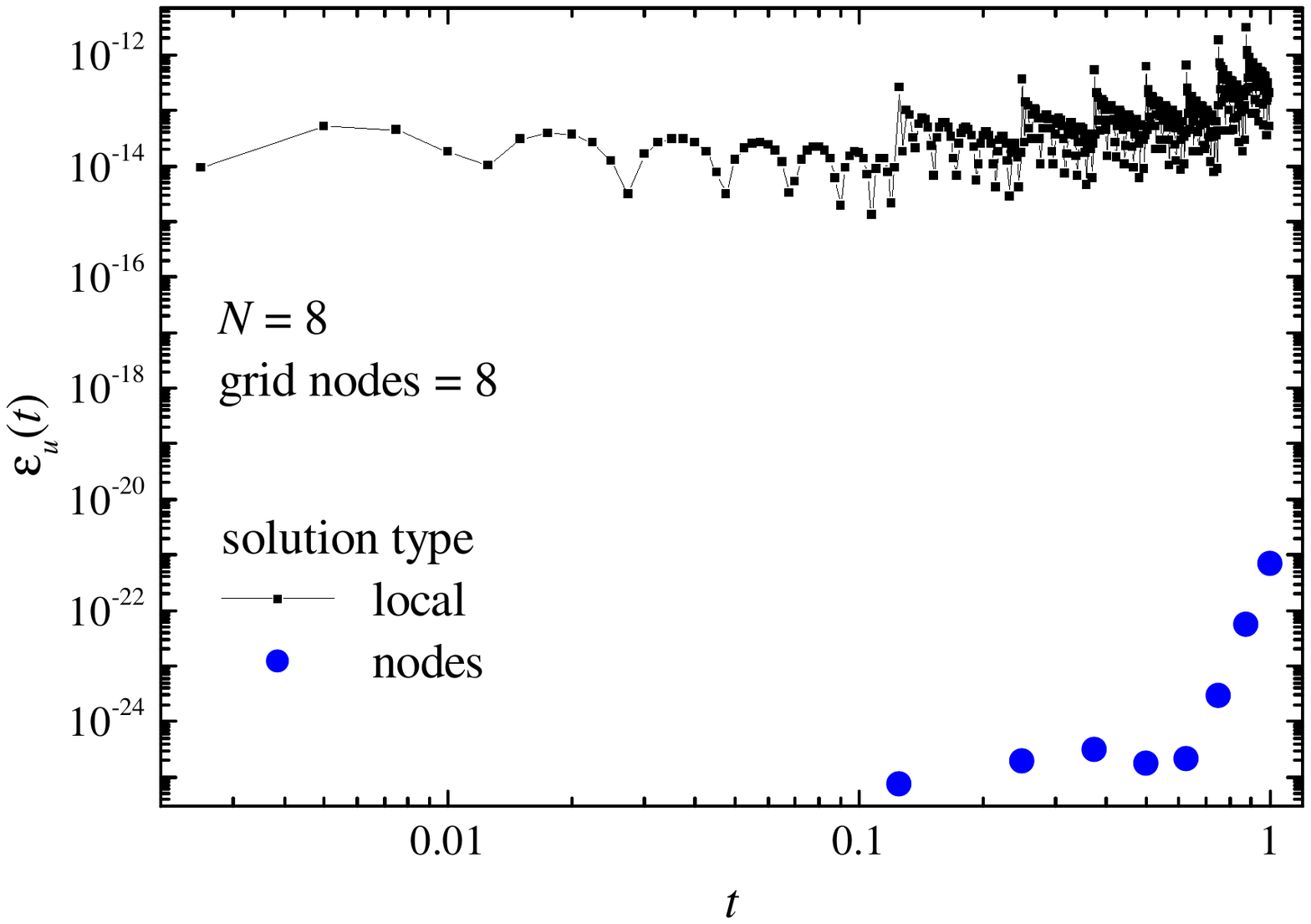}
\vspace{-8mm}\caption{\label{fig:hess_1_sol_g_eps:c2}}
\end{subfigure}
\begin{subfigure}{0.320\textwidth}
\includegraphics[width=\textwidth]{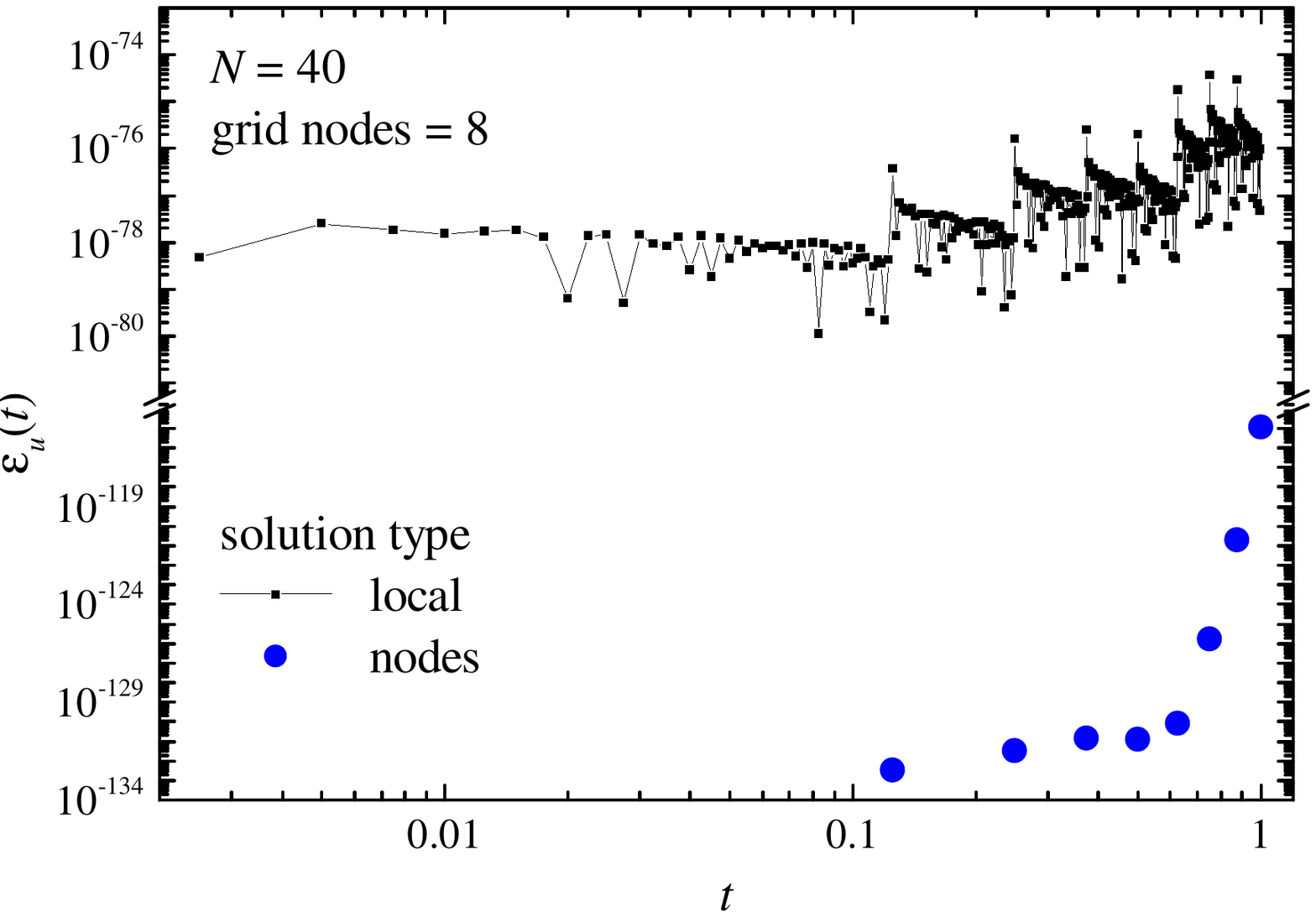}
\vspace{-8mm}\caption{\label{fig:hess_1_sol_g_eps:c3}}
\end{subfigure}\\
\begin{subfigure}{0.320\textwidth}
\includegraphics[width=\textwidth]{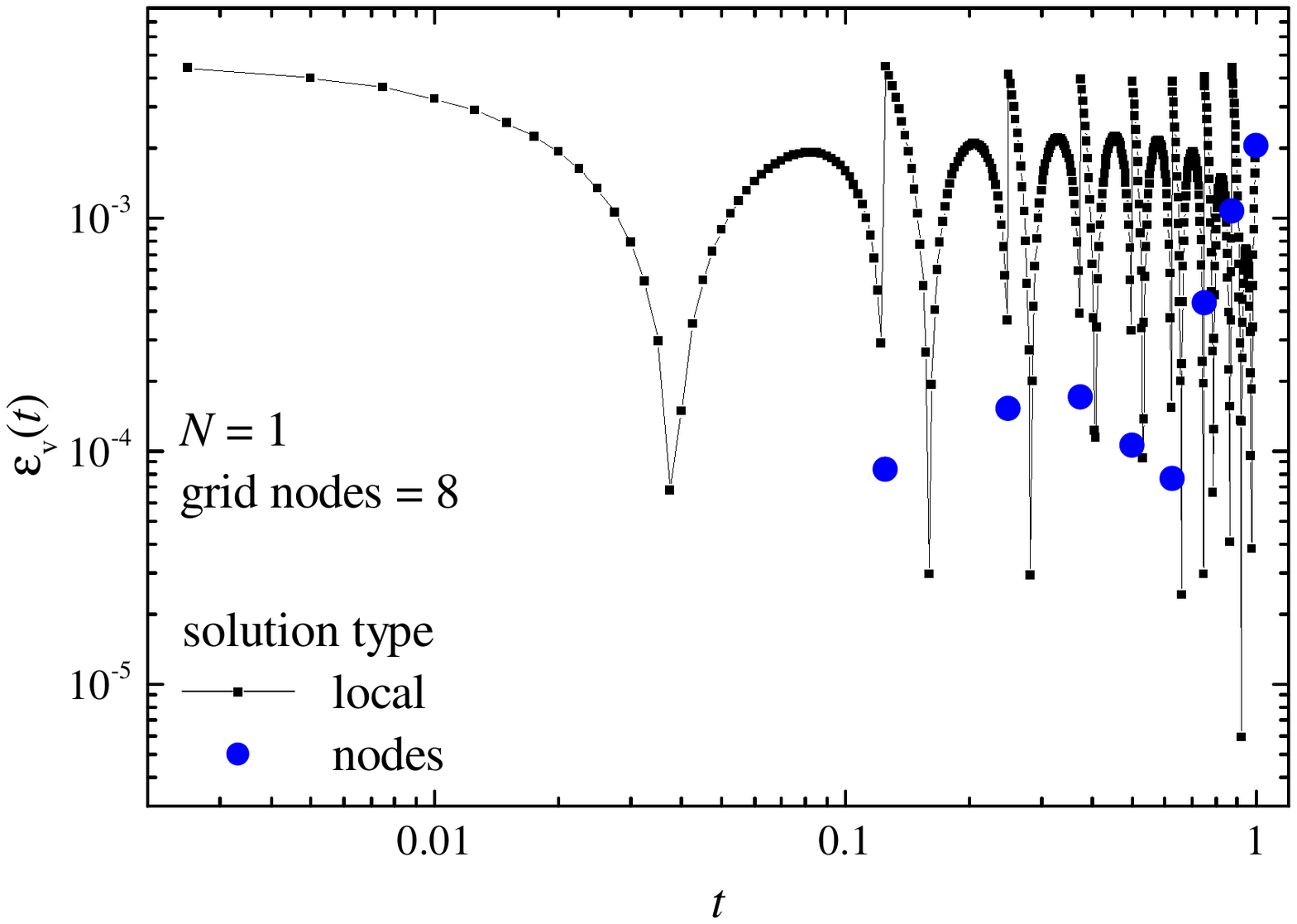}
\vspace{-8mm}\caption{\label{fig:hess_1_sol_g_eps:d1}}
\end{subfigure}
\begin{subfigure}{0.320\textwidth}
\includegraphics[width=\textwidth]{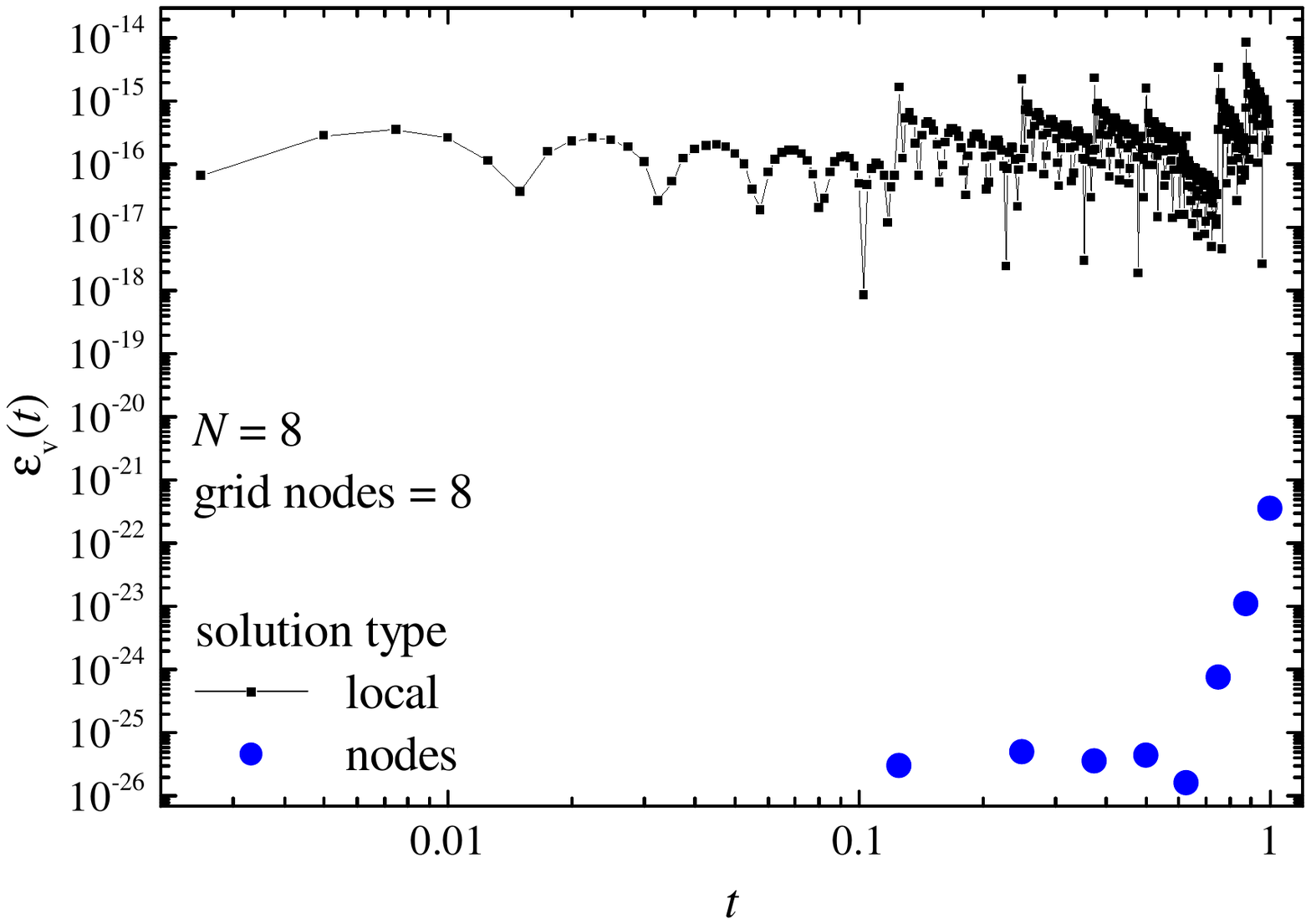}
\vspace{-8mm}\caption{\label{fig:hess_1_sol_g_eps:d2}}
\end{subfigure}
\begin{subfigure}{0.320\textwidth}
\includegraphics[width=\textwidth]{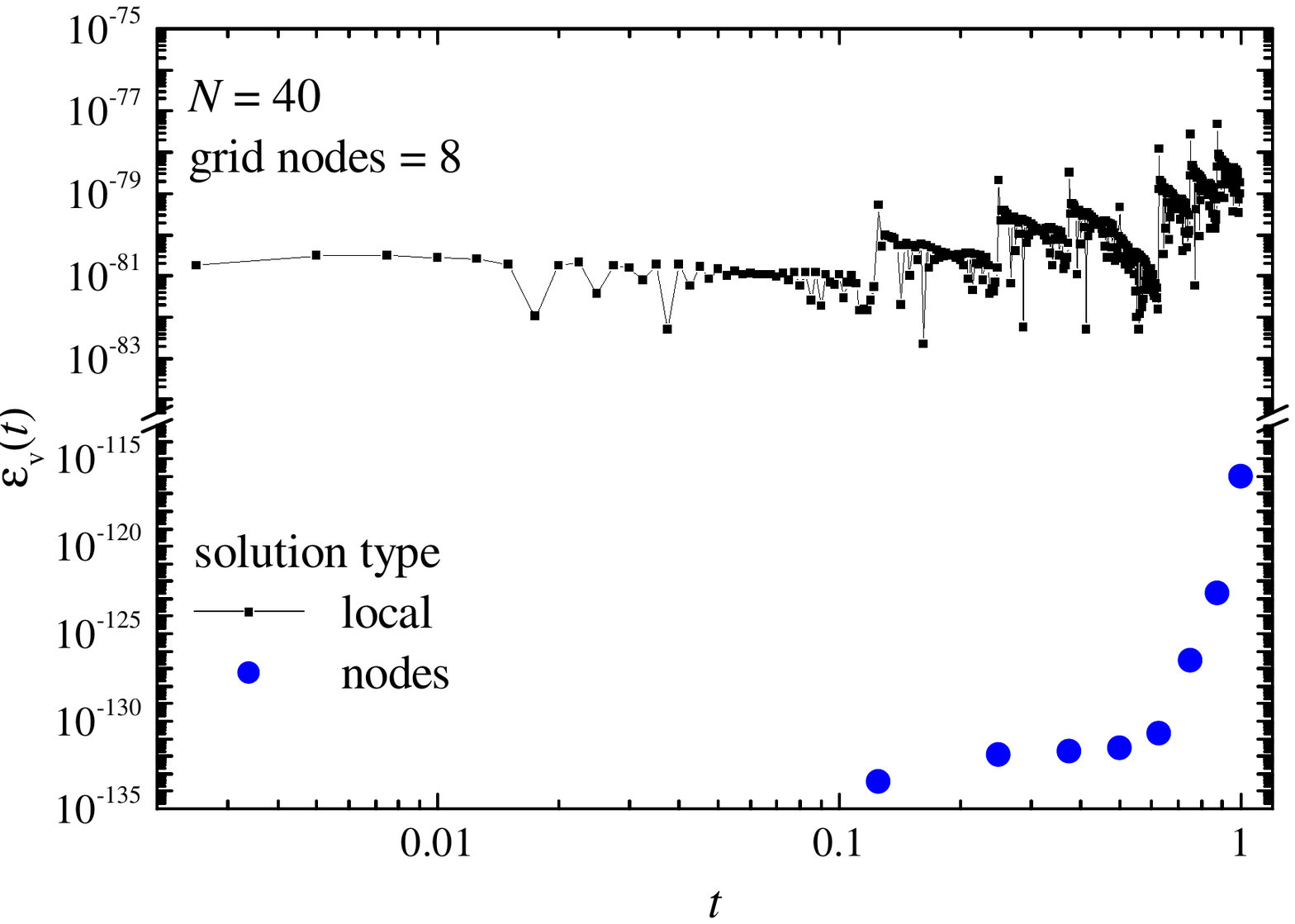}
\vspace{-8mm}\caption{\label{fig:hess_1_sol_g_eps:d3}}
\end{subfigure}\\
\begin{subfigure}{0.320\textwidth}
\includegraphics[width=\textwidth]{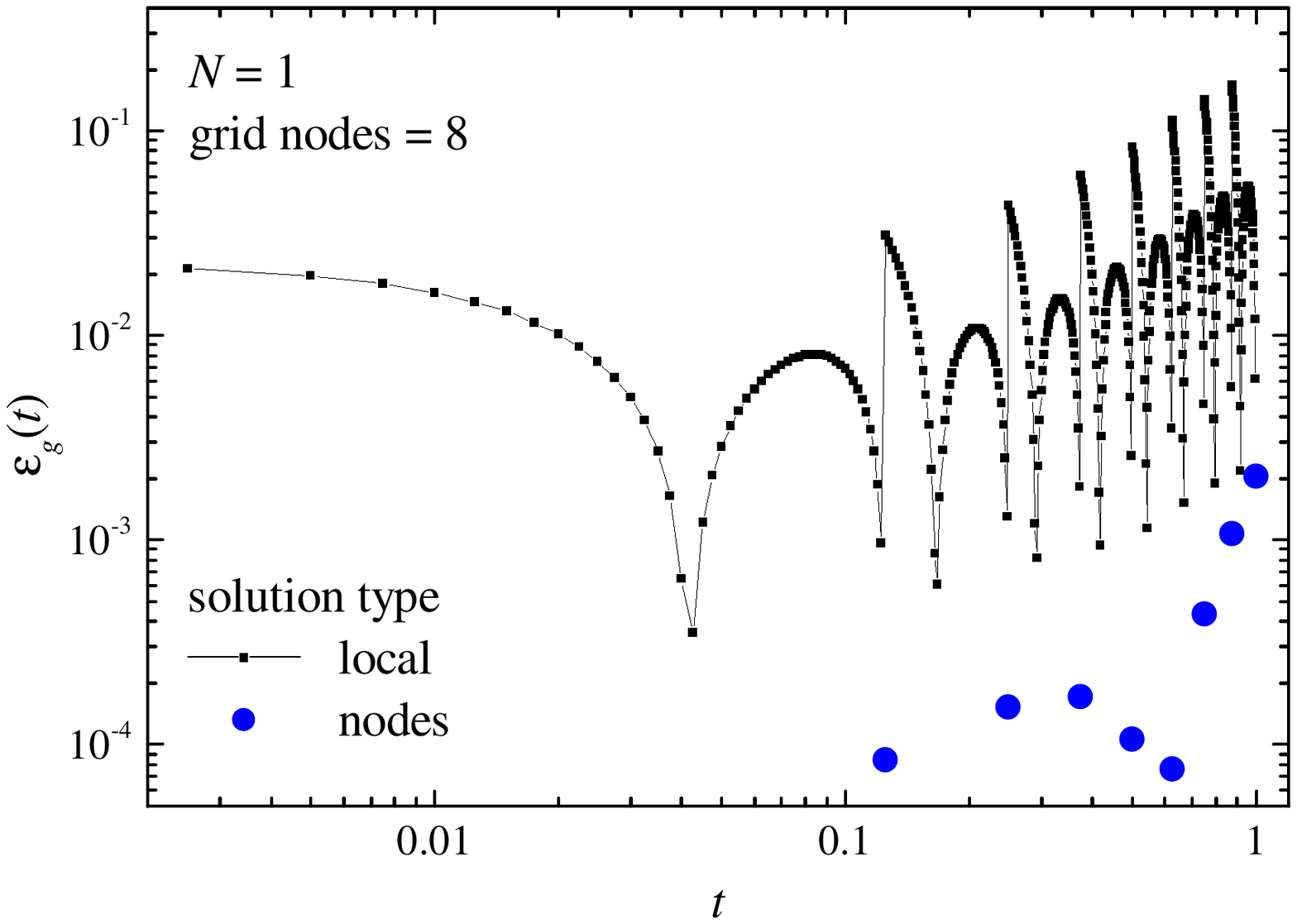}
\vspace{-8mm}\caption{\label{fig:hess_1_sol_g_eps:e1}}
\end{subfigure}
\begin{subfigure}{0.320\textwidth}
\includegraphics[width=\textwidth]{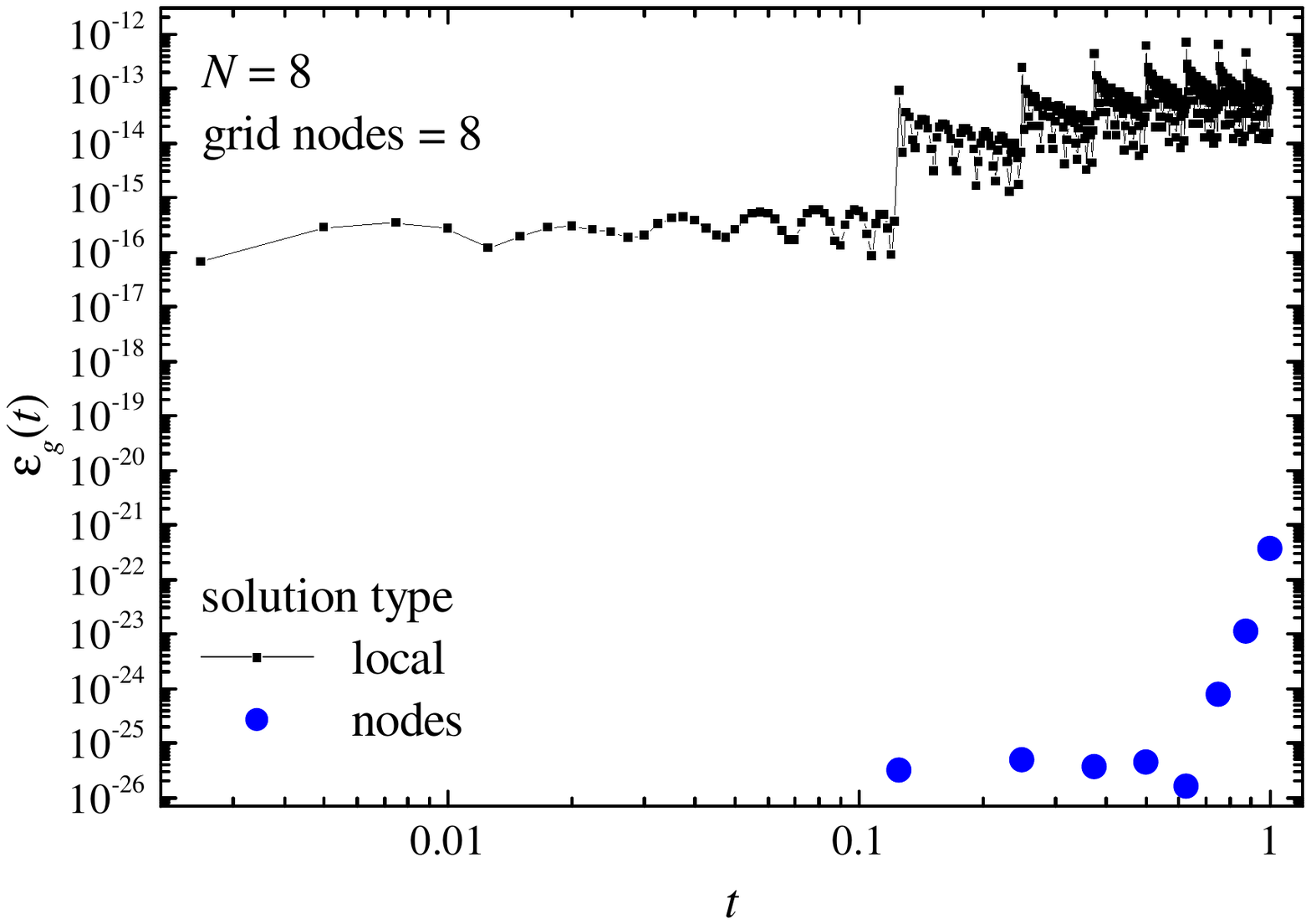}
\vspace{-8mm}\caption{\label{fig:hess_1_sol_g_eps:e2}}
\end{subfigure}
\begin{subfigure}{0.320\textwidth}
\includegraphics[width=\textwidth]{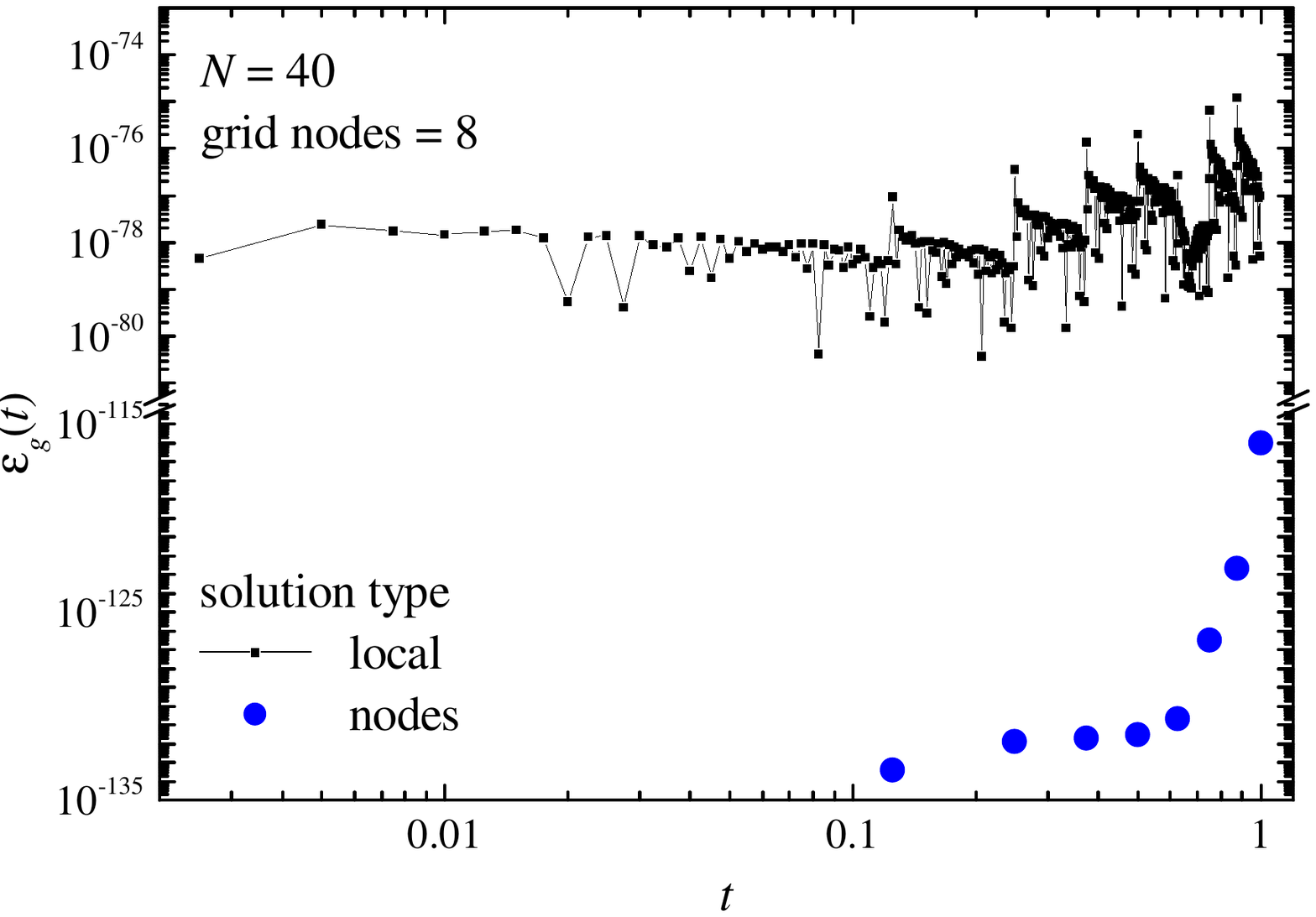}
\vspace{-8mm}\caption{\label{fig:hess_1_sol_g_eps:e3}}
\end{subfigure}\\
\caption{%
Numerical solution of the problem (\ref{eq:hess_dae_ind_1}). Comparison of quantitative satisfiability of the conditions $g_{1} = 0$ (\subref{fig:hess_1_sol_g_eps:a1}, \subref{fig:hess_1_sol_g_eps:a2}, \subref{fig:hess_1_sol_g_eps:a3}) and $g_{2} = 0$ (\subref{fig:hess_1_sol_g_eps:b1}, \subref{fig:hess_1_sol_g_eps:b2}, \subref{fig:hess_1_sol_g_eps:b3}), the errors $\varepsilon_{u}(t)$ (\subref{fig:hess_1_sol_g_eps:c1}, \subref{fig:hess_1_sol_g_eps:c2}, \subref{fig:hess_1_sol_g_eps:c3}), $\varepsilon_{v}(t)$ (\subref{fig:hess_1_sol_g_eps:d1}, \subref{fig:hess_1_sol_g_eps:d2}, \subref{fig:hess_1_sol_g_eps:d3}), $\varepsilon_{g}(t)$ (\subref{fig:hess_1_sol_g_eps:e1}, \subref{fig:hess_1_sol_g_eps:e2}, \subref{fig:hess_1_sol_g_eps:e3}), for numerical solution obtained using polynomials with degrees $N = 1$ (\subref{fig:hess_1_sol_g_eps:a1}, \subref{fig:hess_1_sol_g_eps:b1}, \subref{fig:hess_1_sol_g_eps:c1}, \subref{fig:hess_1_sol_g_eps:d1}, \subref{fig:hess_1_sol_g_eps:e1}), $N = 8$ (\subref{fig:hess_1_sol_g_eps:a2}, \subref{fig:hess_1_sol_g_eps:b2}, \subref{fig:hess_1_sol_g_eps:c2}, \subref{fig:hess_1_sol_g_eps:d2}, \subref{fig:hess_1_sol_g_eps:e2}) and $N = 40$ (\subref{fig:hess_1_sol_g_eps:a3}, \subref{fig:hess_1_sol_g_eps:b3}, \subref{fig:hess_1_sol_g_eps:c3}, \subref{fig:hess_1_sol_g_eps:d3}, \subref{fig:hess_1_sol_g_eps:e3}).
}
\label{fig:hess_1_sol_g_eps}
\end{figure} 

\begin{figure}[h!]
\captionsetup[subfigure]{%
	position=bottom,
	font+=smaller,
	textfont=normalfont,
	singlelinecheck=off,
	justification=raggedright
}
\centering
\begin{subfigure}{0.275\textwidth}
\includegraphics[width=\textwidth]{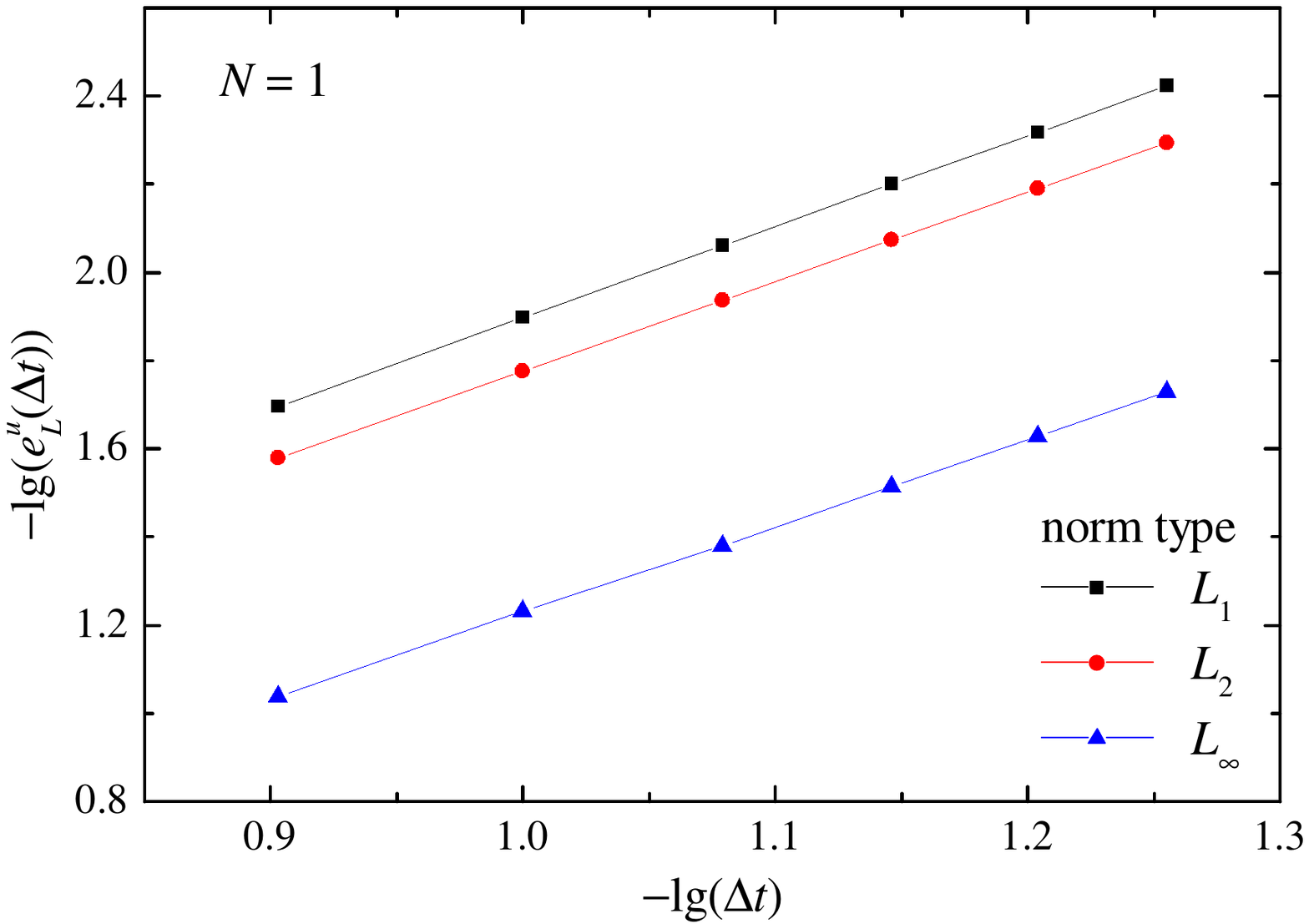}
\vspace{-8mm}\caption{\label{fig:hess_1_errors:a1}}
\end{subfigure}\hspace{6mm}
\begin{subfigure}{0.275\textwidth}
\includegraphics[width=\textwidth]{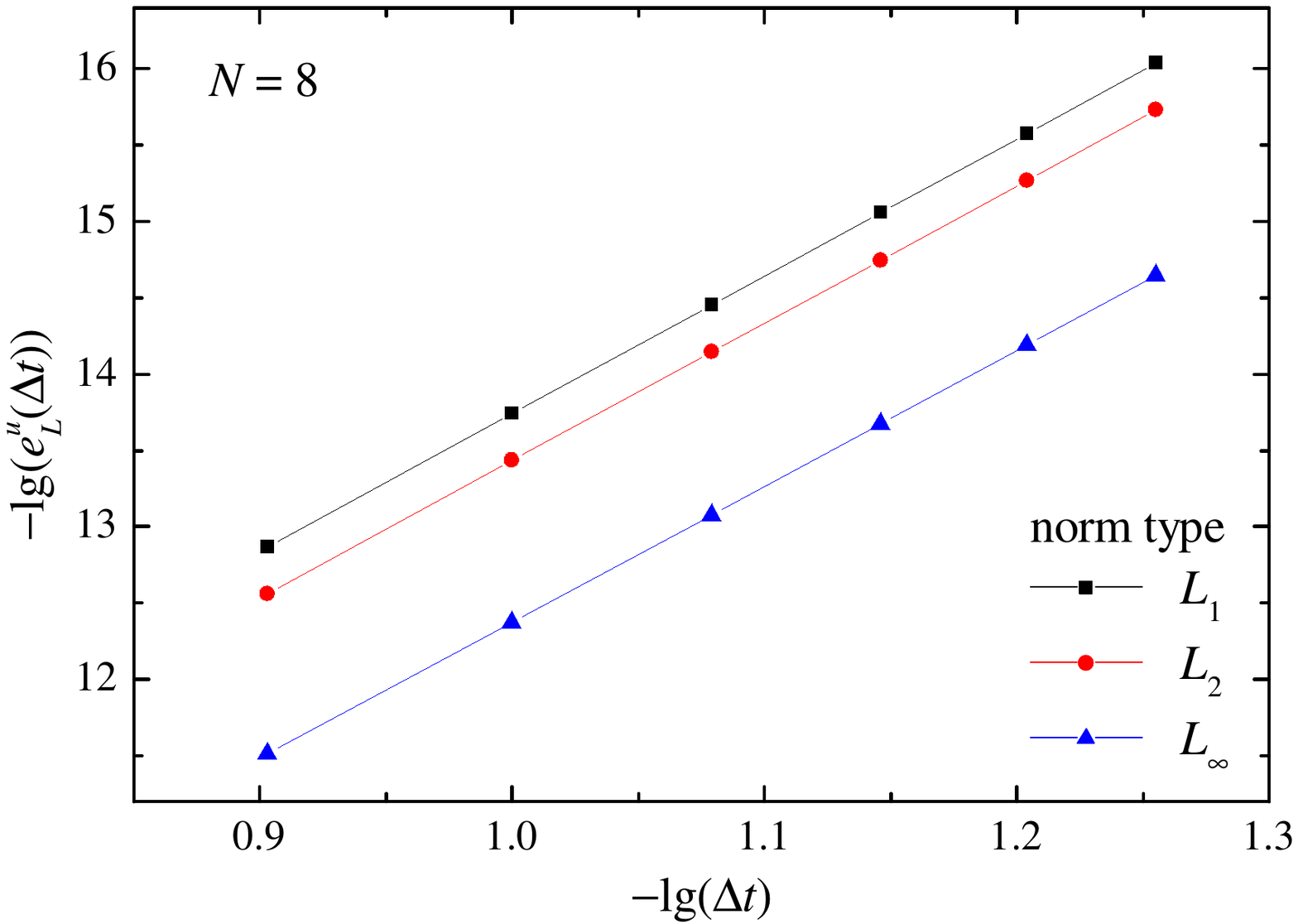}
\vspace{-8mm}\caption{\label{fig:hess_1_errors:a2}}
\end{subfigure}\hspace{6mm}
\begin{subfigure}{0.275\textwidth}
\includegraphics[width=\textwidth]{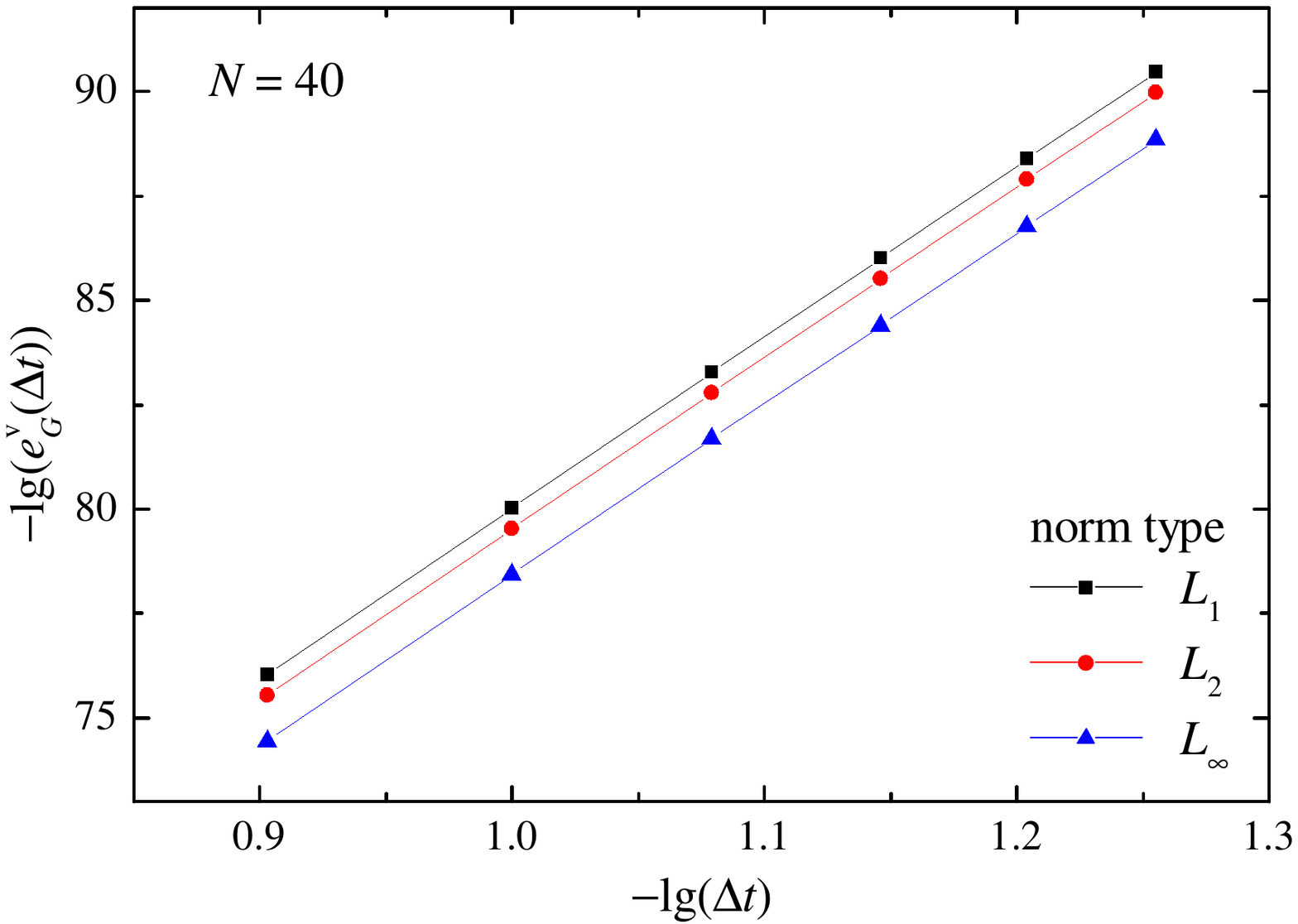}
\vspace{-8mm}\caption{\label{fig:hess_1_errors:a3}}
\end{subfigure}\\[-2mm]
\begin{subfigure}{0.275\textwidth}
\includegraphics[width=\textwidth]{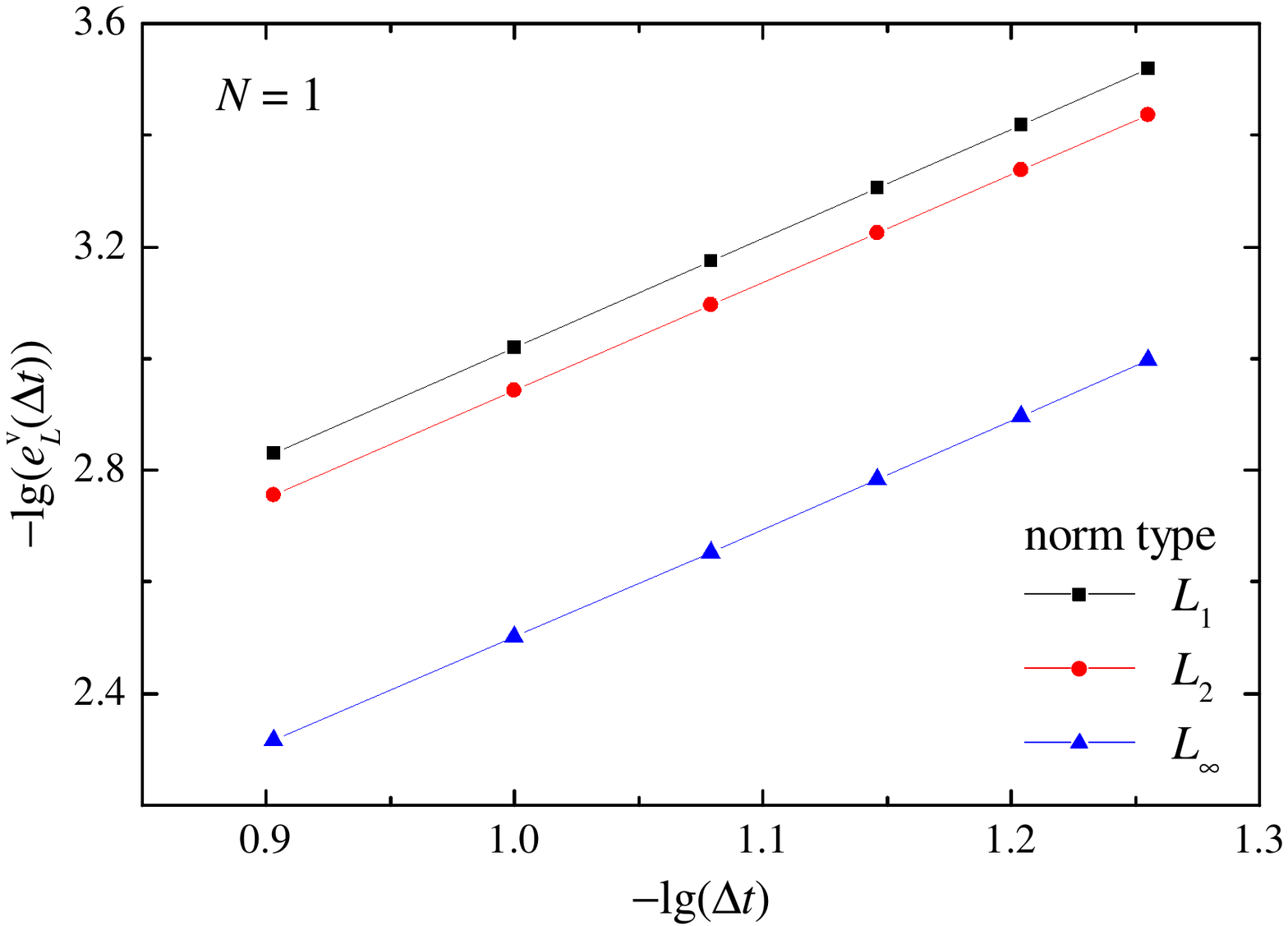}
\vspace{-8mm}\caption{\label{fig:hess_1_errors:b1}}
\end{subfigure}\hspace{6mm}
\begin{subfigure}{0.275\textwidth}
\includegraphics[width=\textwidth]{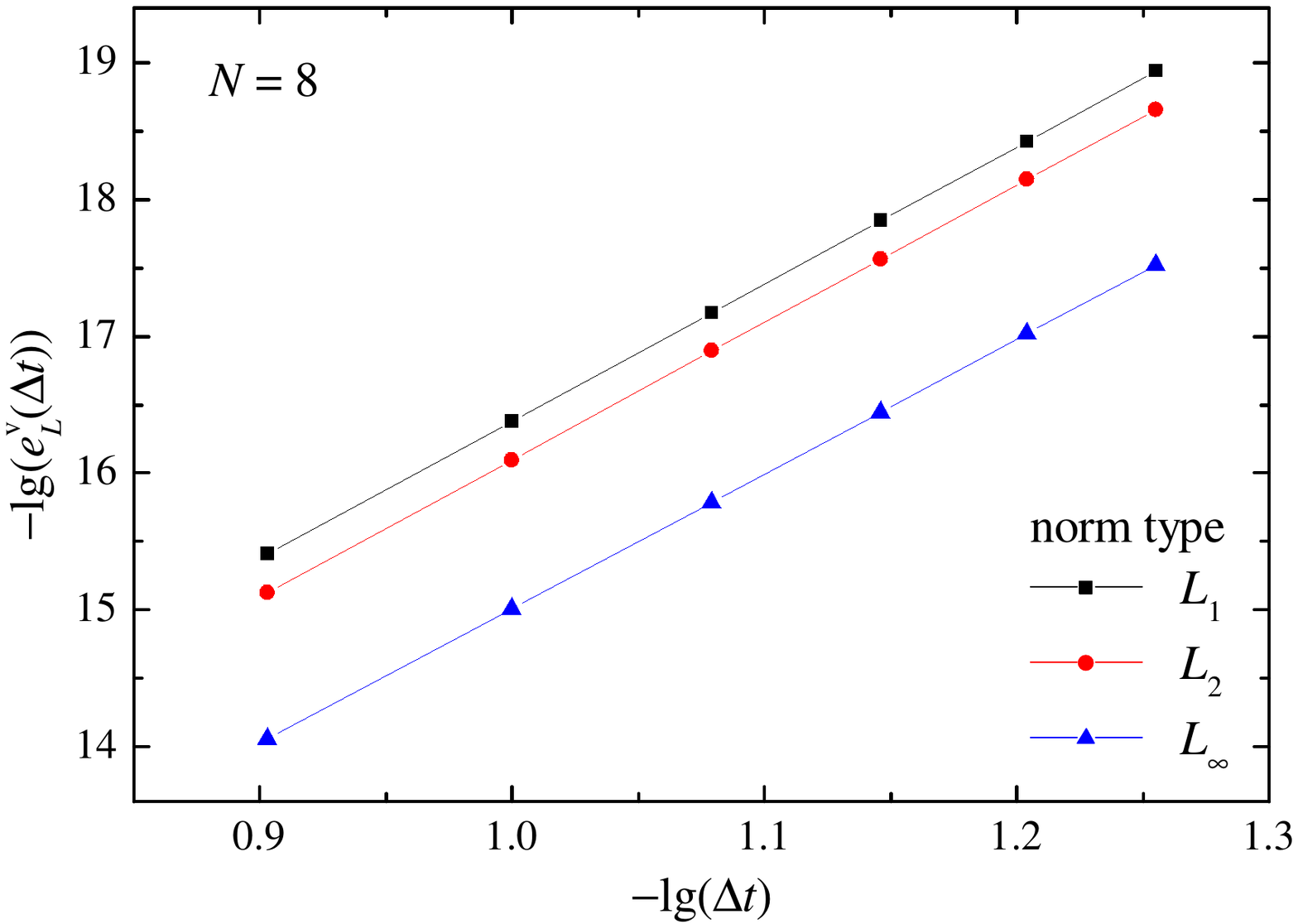}
\vspace{-8mm}\caption{\label{fig:hess_1_errors:b2}}
\end{subfigure}\hspace{6mm}
\begin{subfigure}{0.275\textwidth}
\includegraphics[width=\textwidth]{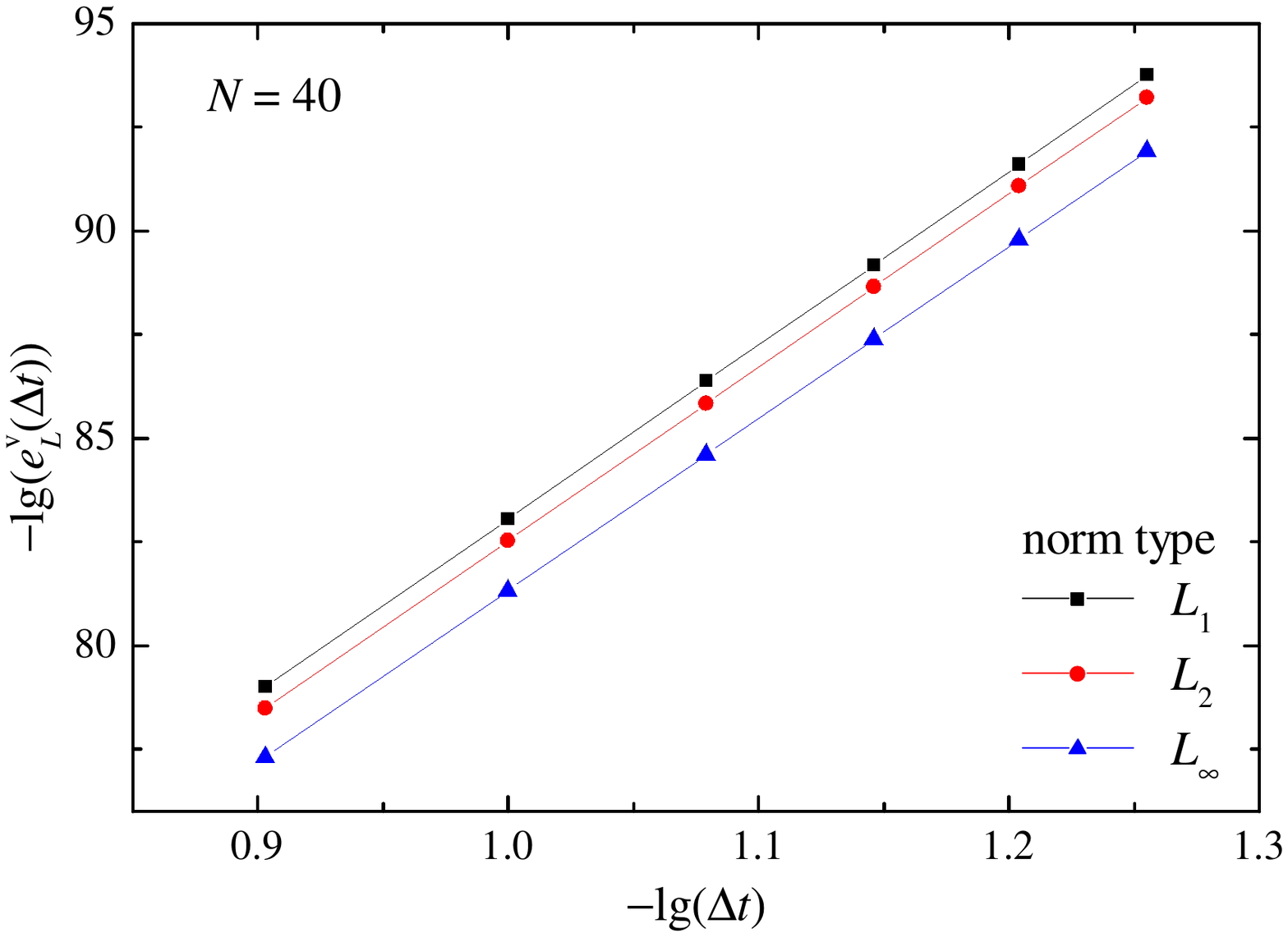}
\vspace{-8mm}\caption{\label{fig:hess_1_errors:b3}}
\end{subfigure}\\[-2mm]
\begin{subfigure}{0.275\textwidth}
\includegraphics[width=\textwidth]{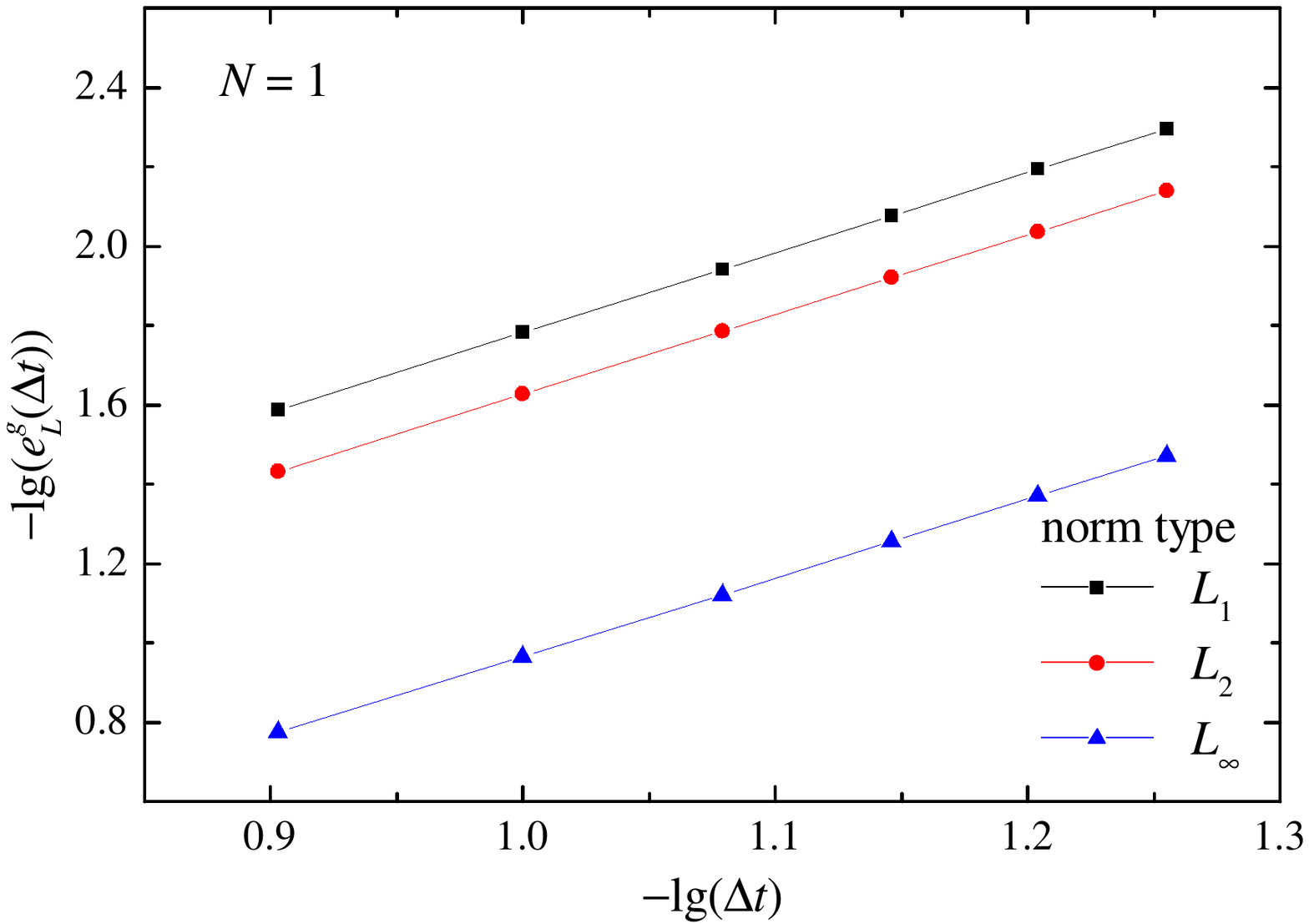}
\vspace{-8mm}\caption{\label{fig:hess_1_errors:c1}}
\end{subfigure}\hspace{6mm}
\begin{subfigure}{0.275\textwidth}
\includegraphics[width=\textwidth]{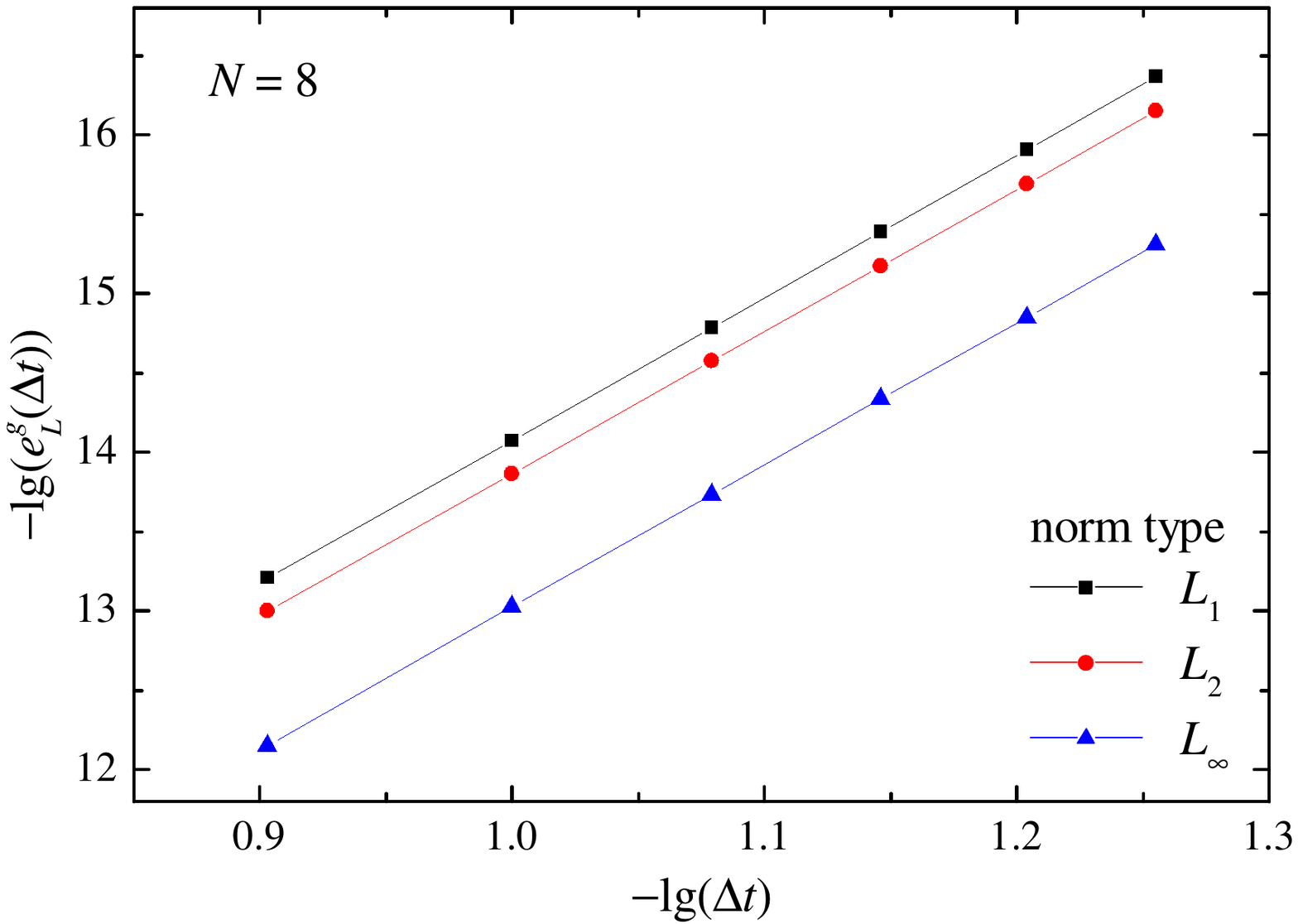}
\vspace{-8mm}\caption{\label{fig:hess_1_errors:c2}}
\end{subfigure}\hspace{6mm}
\begin{subfigure}{0.275\textwidth}
\includegraphics[width=\textwidth]{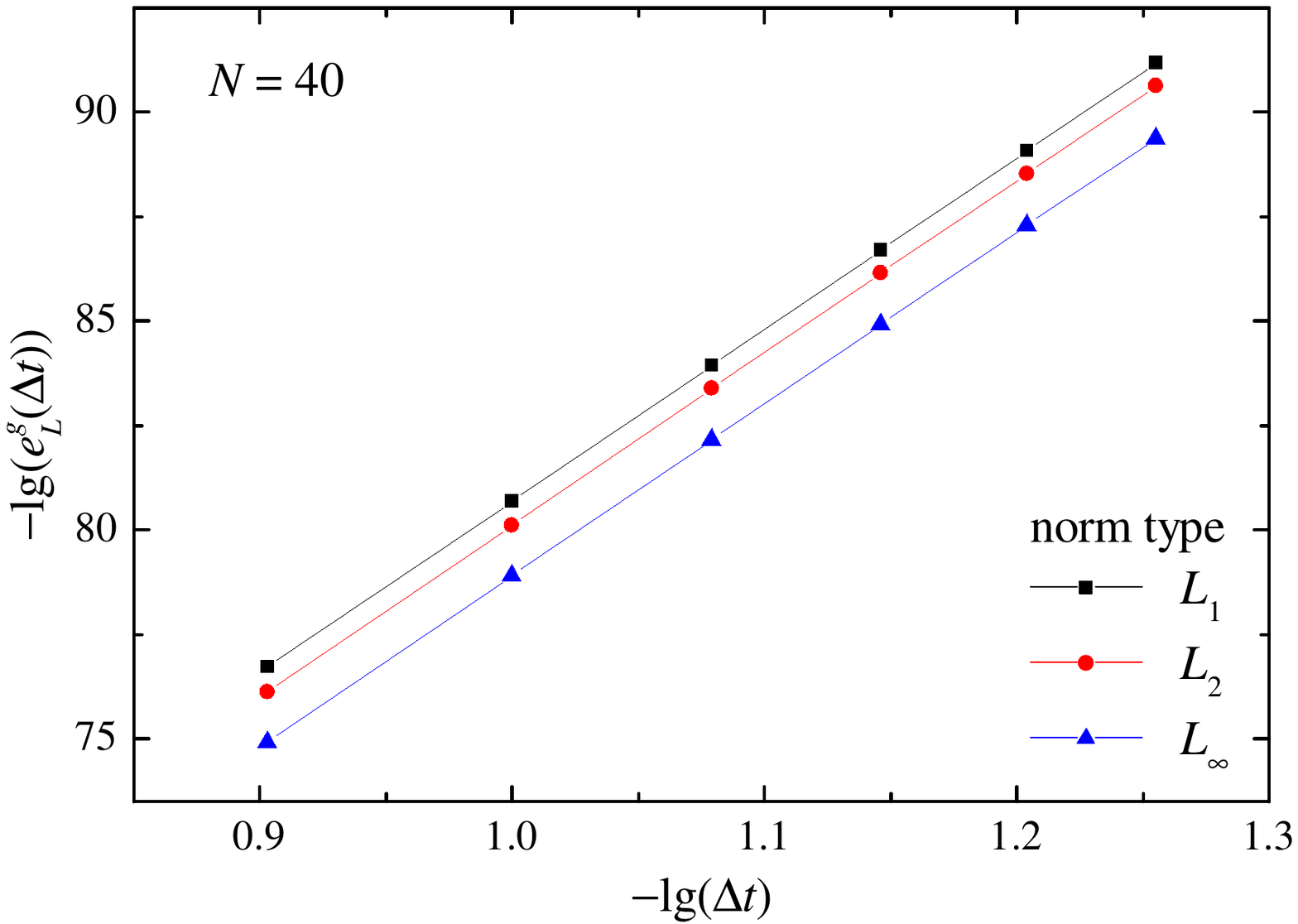}
\vspace{-8mm}\caption{\label{fig:hess_1_errors:c3}}
\end{subfigure}\\[-2mm]
\begin{subfigure}{0.275\textwidth}
\includegraphics[width=\textwidth]{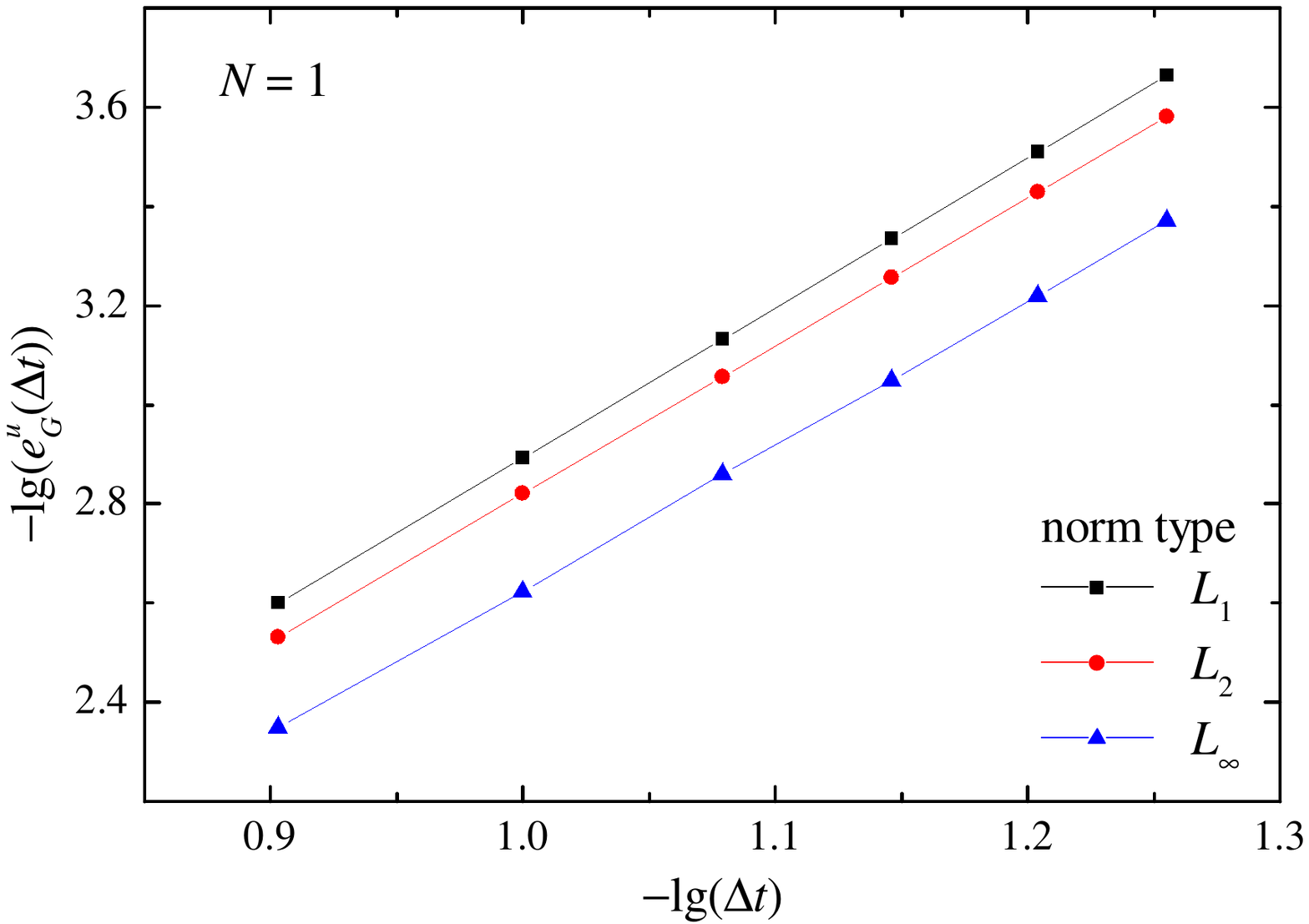}
\vspace{-8mm}\caption{\label{fig:hess_1_errors:d1}}
\end{subfigure}\hspace{6mm}
\begin{subfigure}{0.275\textwidth}
\includegraphics[width=\textwidth]{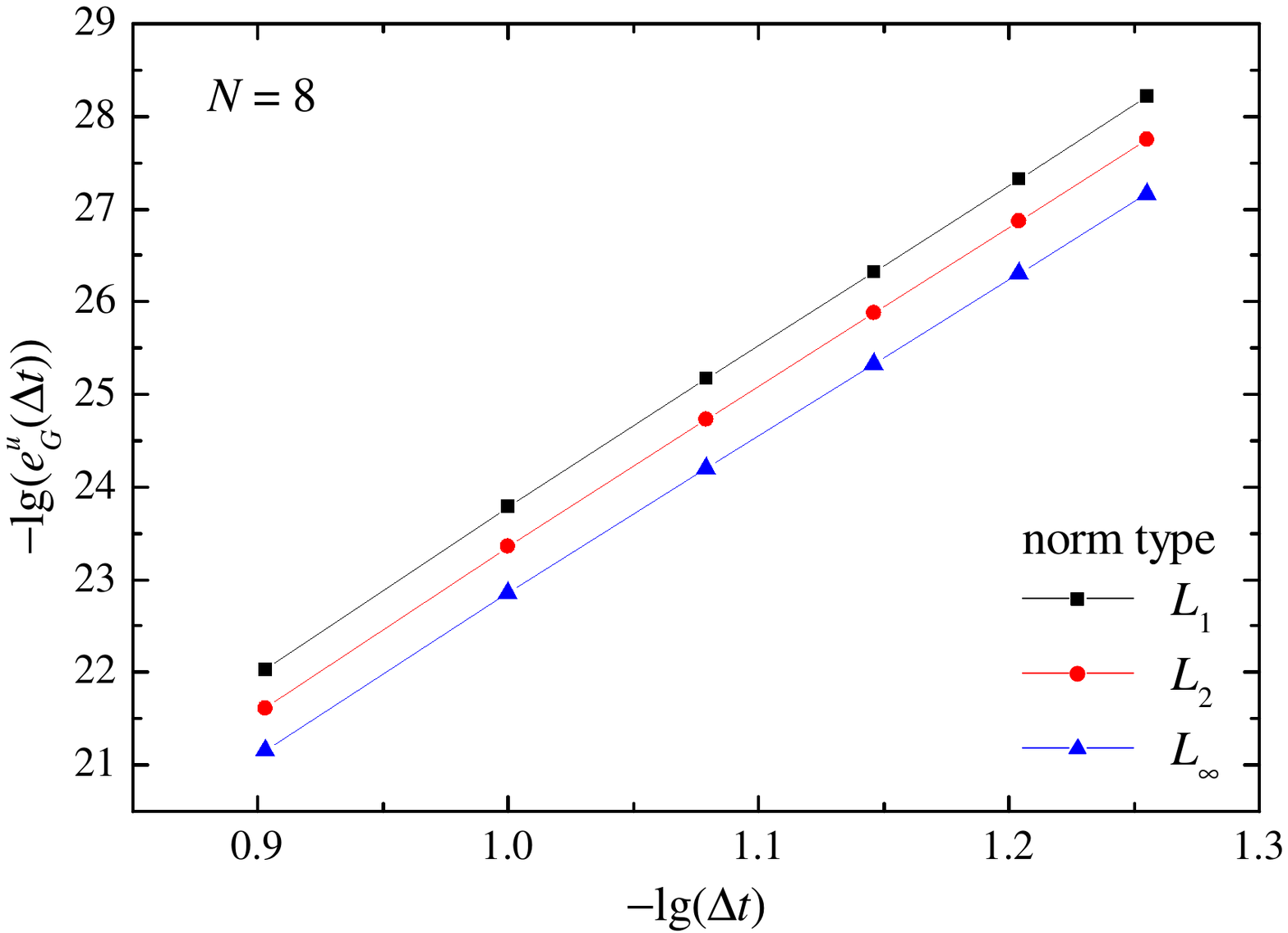}
\vspace{-8mm}\caption{\label{fig:hess_1_errors:d2}}
\end{subfigure}\hspace{6mm}
\begin{subfigure}{0.275\textwidth}
\includegraphics[width=\textwidth]{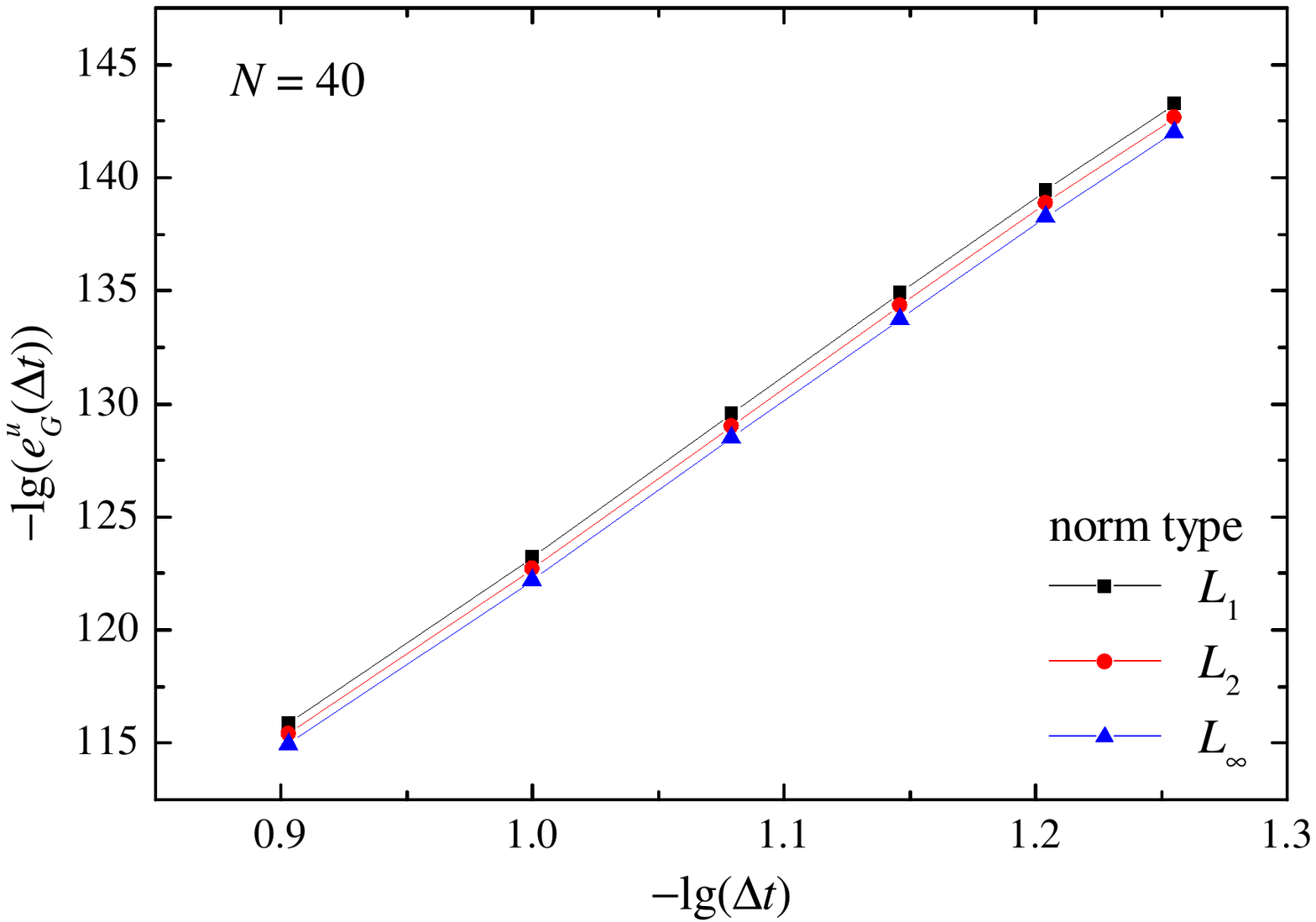}
\vspace{-8mm}\caption{\label{fig:hess_1_errors:d3}}
\end{subfigure}\\[-2mm]
\begin{subfigure}{0.275\textwidth}
\includegraphics[width=\textwidth]{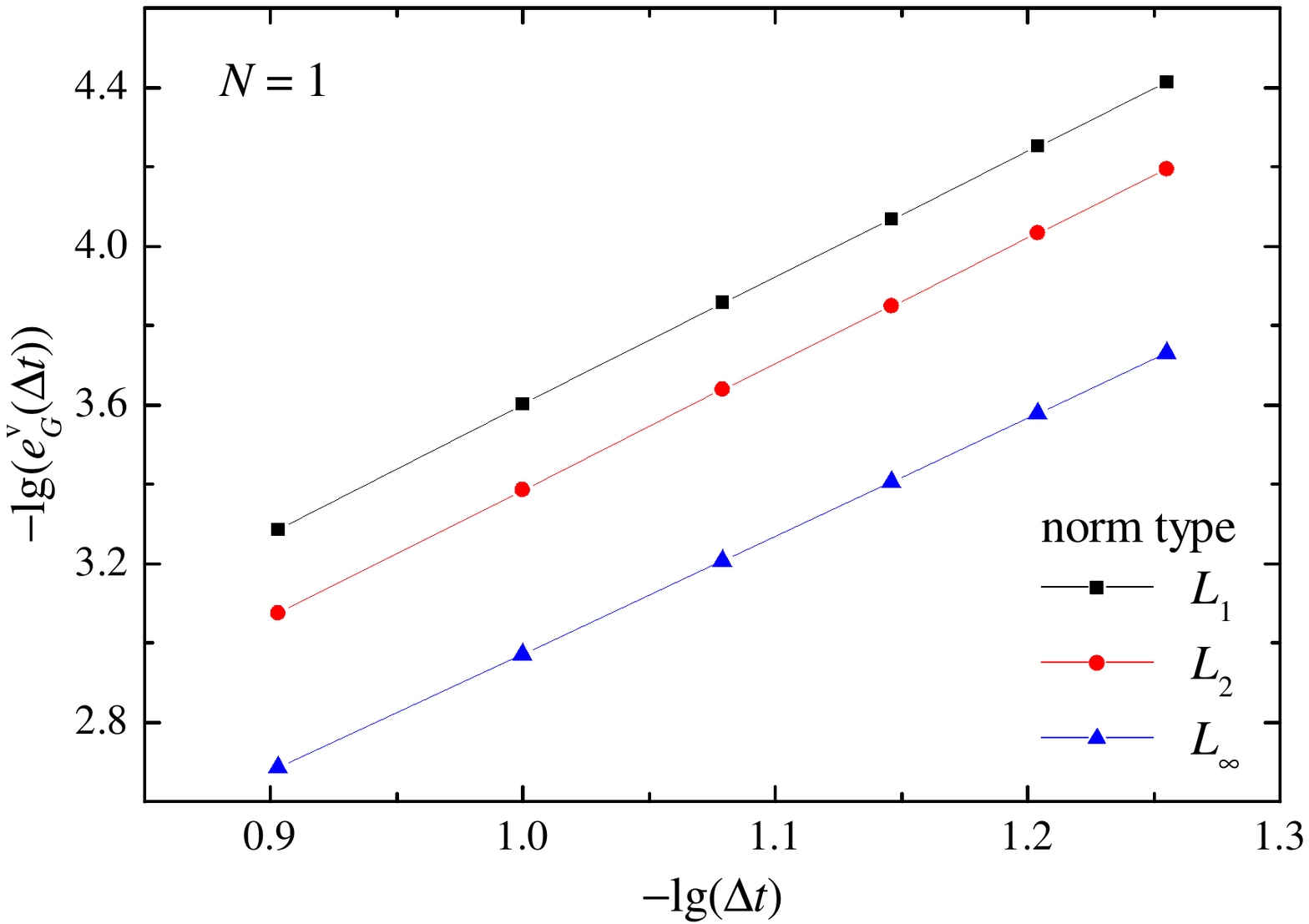}
\vspace{-8mm}\caption{\label{fig:hess_1_errors:e1}}
\end{subfigure}\hspace{6mm}
\begin{subfigure}{0.275\textwidth}
\includegraphics[width=\textwidth]{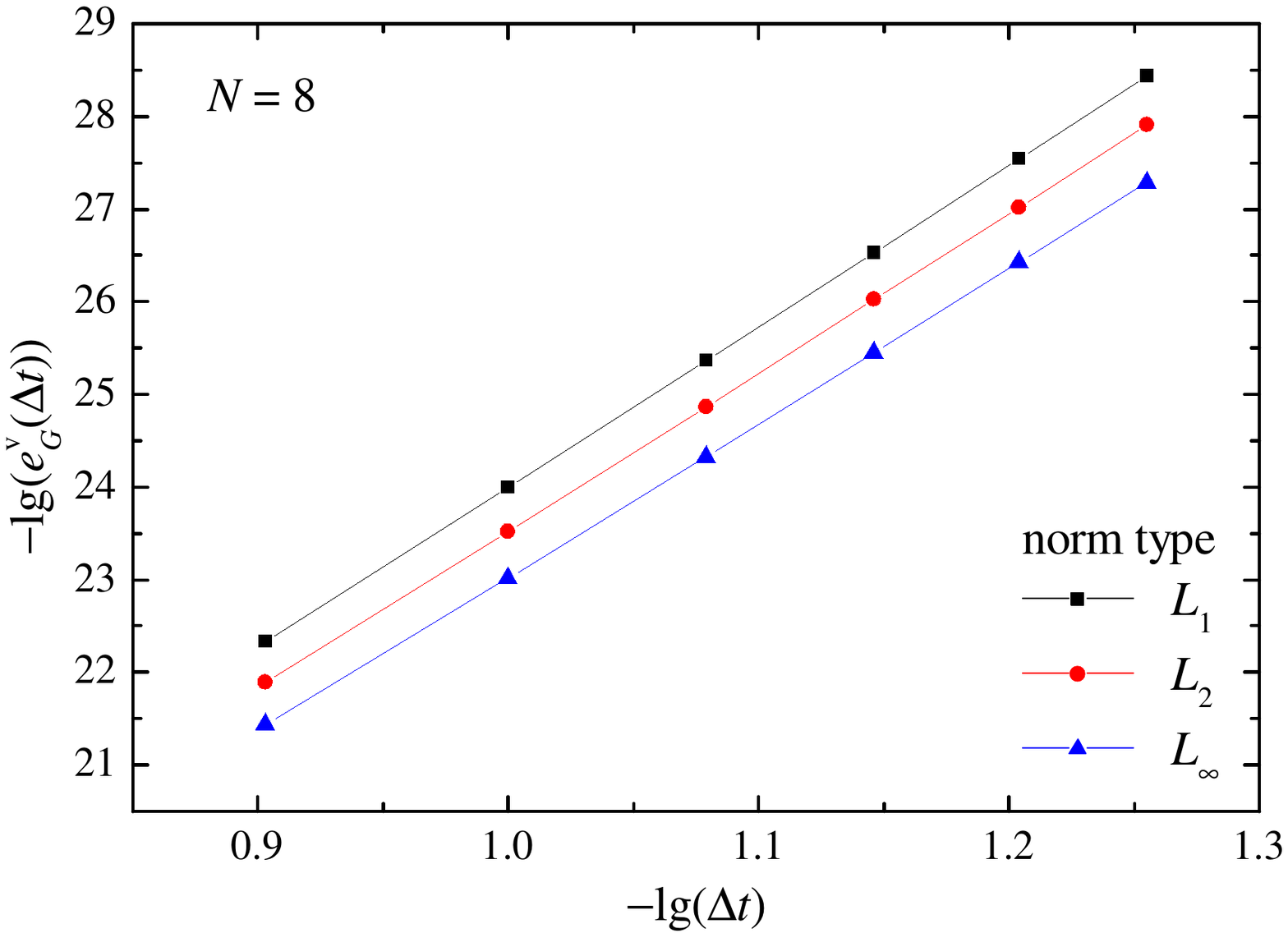}
\vspace{-8mm}\caption{\label{fig:hess_1_errors:e2}}
\end{subfigure}\hspace{6mm}
\begin{subfigure}{0.275\textwidth}
\includegraphics[width=\textwidth]{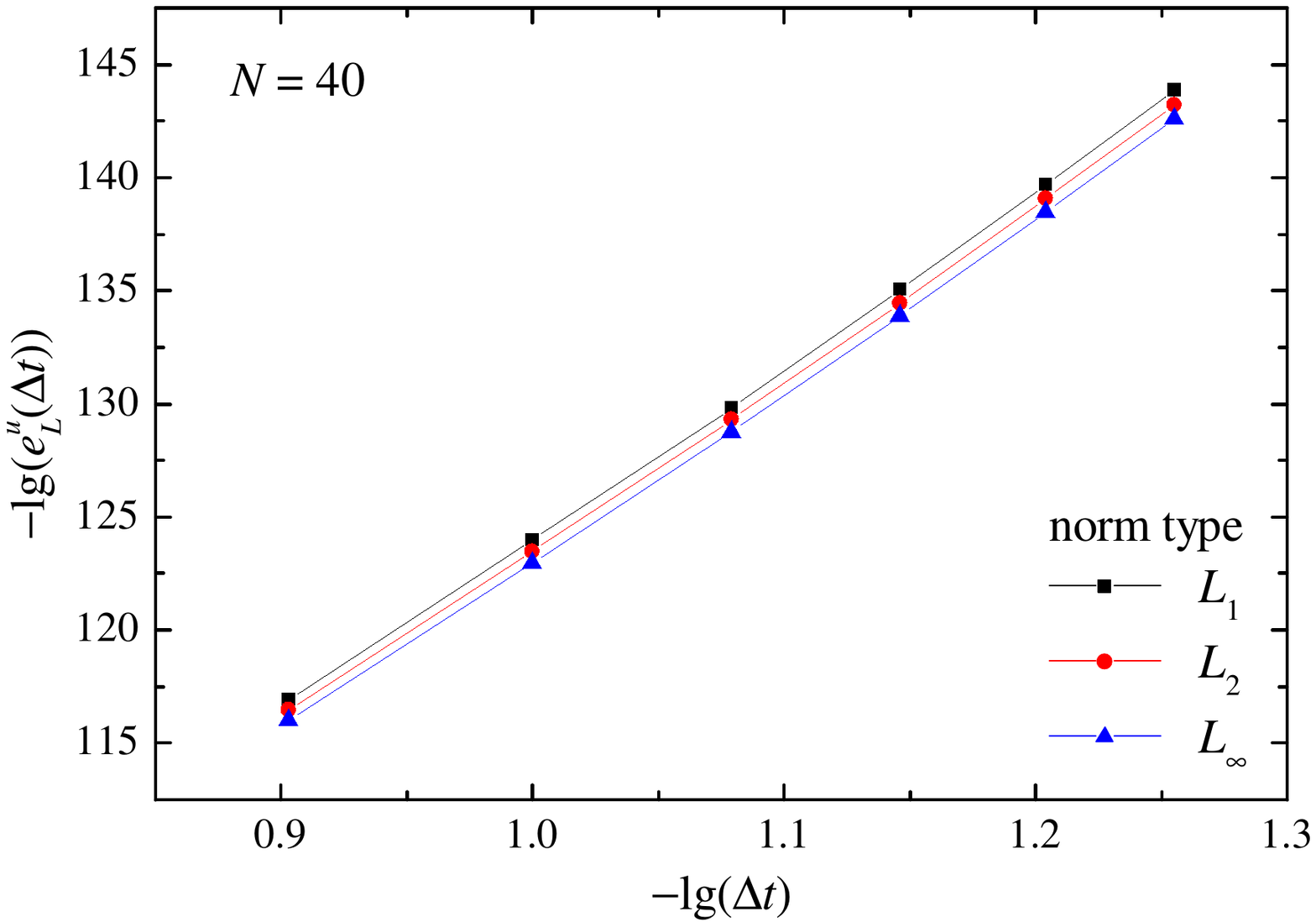}
\vspace{-8mm}\caption{\label{fig:hess_1_errors:e3}}
\end{subfigure}\\[-2mm]
\begin{subfigure}{0.275\textwidth}
\includegraphics[width=\textwidth]{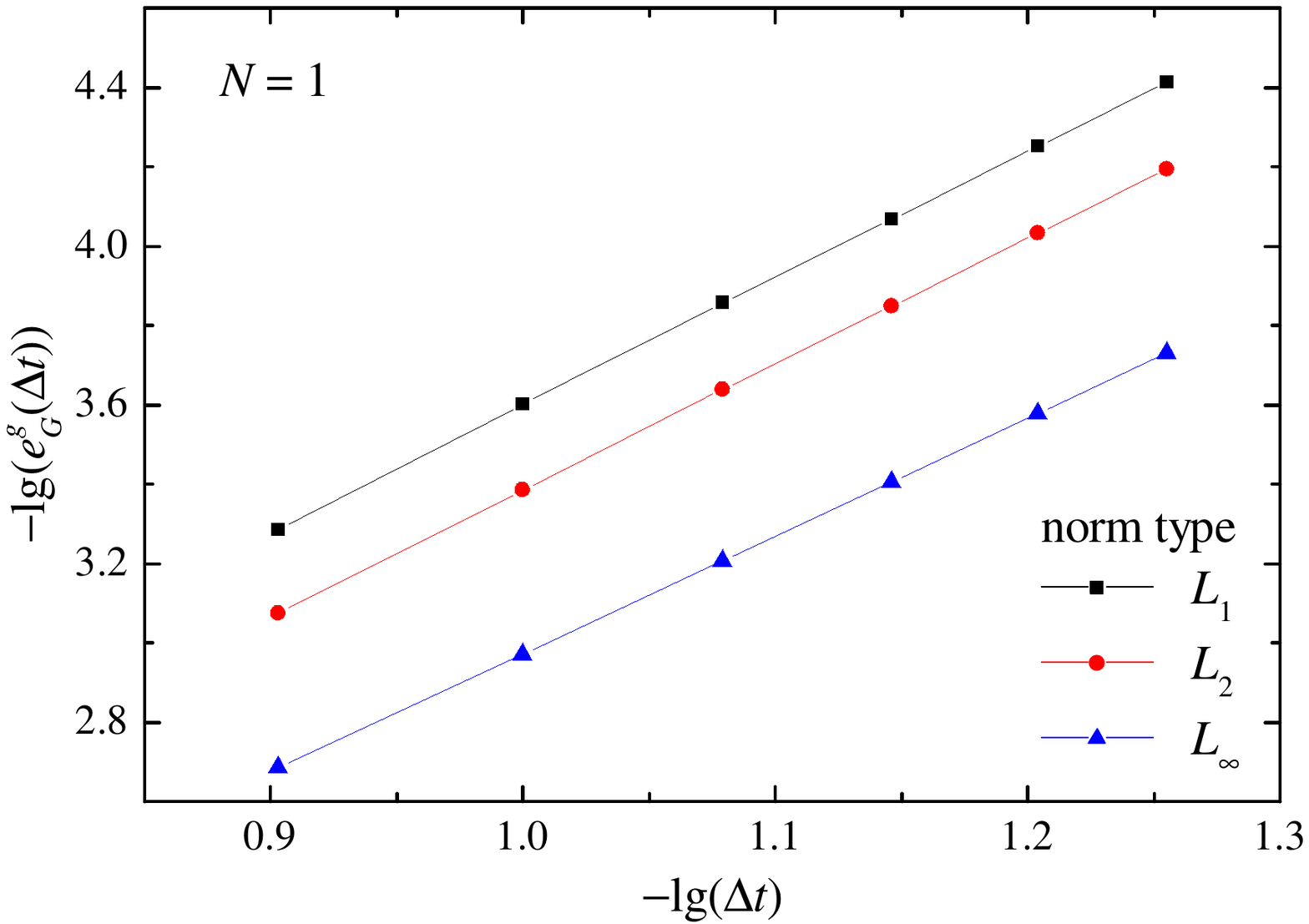}
\vspace{-8mm}\caption{\label{fig:hess_1_errors:f1}}
\end{subfigure}\hspace{6mm}
\begin{subfigure}{0.275\textwidth}
\includegraphics[width=\textwidth]{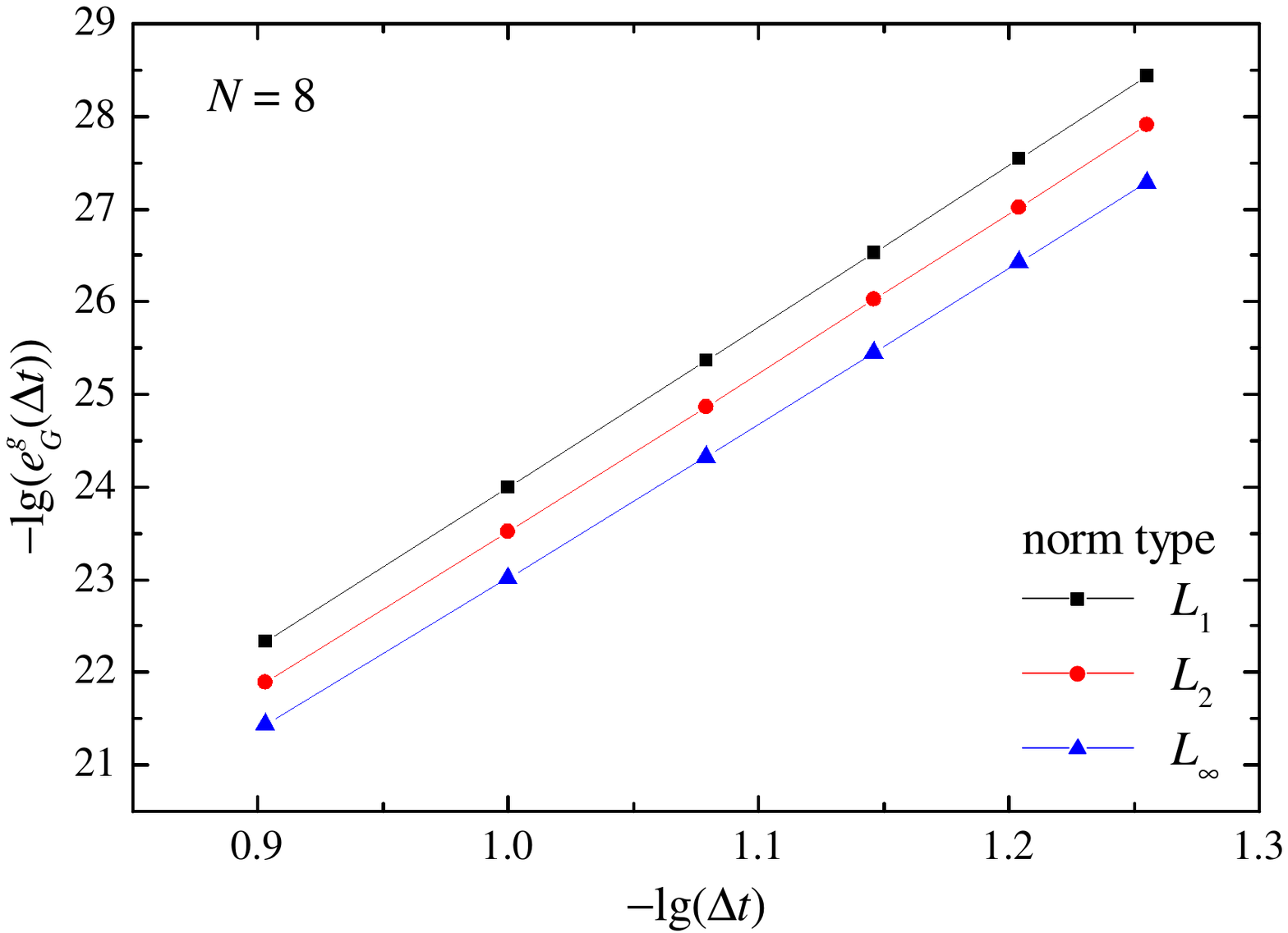}
\vspace{-8mm}\caption{\label{fig:hess_1_errors:f2}}
\end{subfigure}\hspace{6mm}
\begin{subfigure}{0.275\textwidth}
\includegraphics[width=\textwidth]{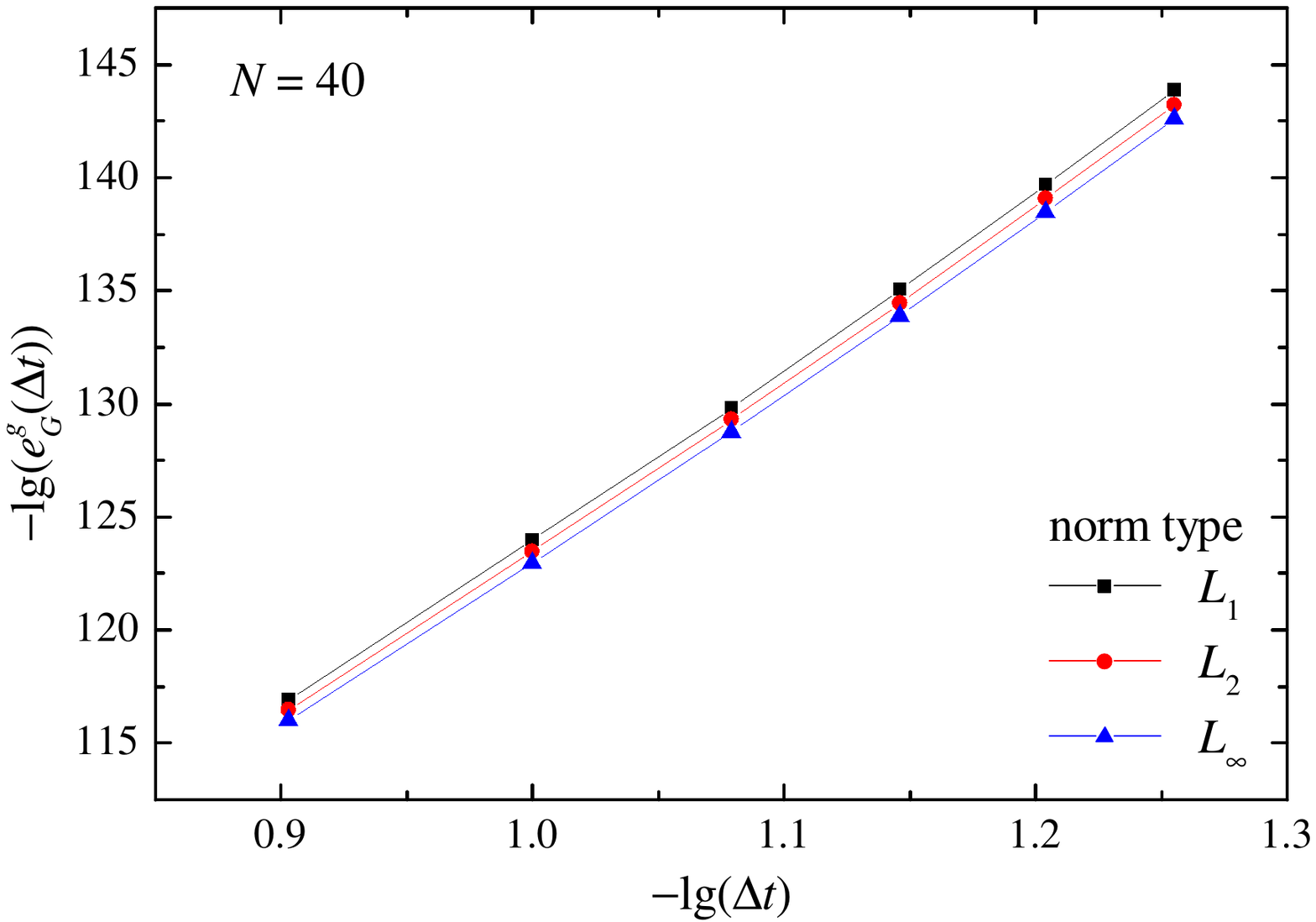}
\vspace{-8mm}\caption{\label{fig:hess_1_errors:f3}}
\end{subfigure}\\[-2mm]
\caption{%
Log-log plot of the dependence of the global errors for the local solution $e_{L}^{u}$ (\subref{fig:hess_1_errors:a1}, \subref{fig:hess_1_errors:a2}, \subref{fig:hess_1_errors:a3}), $e_{L}^{v}$ (\subref{fig:hess_1_errors:b1}, \subref{fig:hess_1_errors:b2}, \subref{fig:hess_1_errors:b3}), $e_{L}^{g}$ (\subref{fig:hess_1_errors:c1}, \subref{fig:hess_1_errors:c2}, \subref{fig:hess_1_errors:c3}) and the solution at nodes $e_{G}^{u}$ (\subref{fig:hess_1_errors:d1}, \subref{fig:hess_1_errors:d2}, \subref{fig:hess_1_errors:d3}), $e_{G}^{v}$ (\subref{fig:hess_1_errors:e1}, \subref{fig:hess_1_errors:e2}, \subref{fig:hess_1_errors:e3}), $e_{G}^{g}$ (\subref{fig:hess_1_errors:f1}, \subref{fig:hess_1_errors:f2}, \subref{fig:hess_1_errors:f3}) on the discretization step $\mathrm{\Delta}t$, obtained in the norms $L_{1}$, $L_{2}$ and $L_{\infty}$, by numerical solution of the problem (\ref{eq:hess_dae_ind_1}) obtained using polynomials with degrees $N = 1$ (\subref{fig:hess_1_errors:a1}, \subref{fig:hess_1_errors:b1}, \subref{fig:hess_1_errors:c1}, \subref{fig:hess_1_errors:d1}, \subref{fig:hess_1_errors:e1}, \subref{fig:hess_1_errors:f1}), $N = 8$ (\subref{fig:hess_1_errors:a2}, \subref{fig:hess_1_errors:b2}, \subref{fig:hess_1_errors:c2}, \subref{fig:hess_1_errors:d2}, \subref{fig:hess_1_errors:e2}, \subref{fig:hess_1_errors:f2}) and $N = 40$ (\subref{fig:hess_1_errors:a3}, \subref{fig:hess_1_errors:b3}, \subref{fig:hess_1_errors:c3}, \subref{fig:hess_1_errors:d3}, \subref{fig:hess_1_errors:e3}, \subref{fig:hess_1_errors:f3}).
}
\label{fig:hess_1_errors}
\end{figure} 

\begin{table*}[h!]
\centering
\caption{%
Convergence orders $p_{L_{1}}^{n}$, $p_{L_{2}}^{n}$, $p_{L_{\infty}}^{n}$, calculated in norms $L_{1}$, $L_{2}$, $L_{\infty}$, of \textit{the numerical solution at the nodes} $(\mathbf{u}_{n}, \mathbf{v}_{n})$ of the ADER-DG method for the DAE problem (\ref{eq:hess_dae_ind_1}); $N$ is the degree of the basis polynomials $\varphi_{p}$. Orders $p^{n, u}$ are calculated for solution $\mathbf{u}_{n}$; orders $p^{n, v}$ --- for solution $\mathbf{v}_{n}$; orders $p^{n, g}$ --- for the conditions $\mathbf{g} = 0$ on the numerical solution at the nodes $(\mathbf{u}_{n}, \mathbf{v}_{n})$. The theoretical value of convergence order $p_{\rm th.}^{n} = 2N+1$ is applicable for the ADER-DG method for ODE problems and is presented for comparison.
}
\label{tab:conv_orders_nodes_hess_1}
\begin{tabular}{@{}|l|lll|lll|lll|c|@{}}
\toprule
$N$ & $p_{L_{1}}^{n, u}$ & $p_{L_{2}}^{n, u}$ & $p_{L_{\infty}}^{n, u}$ & $p_{L_{1}}^{n, v}$ & $p_{L_{2}}^{n, v}$ & $p_{L_{\infty}}^{n, v}$ & $p_{L_{1}}^{n, g}$ & $p_{L_{2}}^{n, g}$ & $p_{L_{\infty}}^{n, g}$ & $p_{\rm th.}^{n}$ \\
\midrule
$1$	&	$3.02$	&	$2.98$	&	$2.90$	&	$3.20$	&	$3.17$	&	$2.97$	&	$3.20$	&	$3.17$	&	$2.97$	&	$3$\\
$2$	&	$5.01$	&	$4.99$	&	$5.00$	&	$5.27$	&	$5.22$	&	$5.00$	&	$5.27$	&	$5.22$	&	$5.00$	&	$5$\\
$3$	&	$7.01$	&	$6.95$	&	$6.73$	&	$7.35$	&	$7.29$	&	$7.01$	&	$7.35$	&	$7.29$	&	$7.01$	&	$7$\\
$4$	&	$9.06$	&	$8.99$	&	$8.69$	&	$9.42$	&	$9.35$	&	$9.01$	&	$9.42$	&	$9.35$	&	$9.01$	&	$9$\\
$5$	&	$11.05$	&	$10.93$	&	$10.55$	&	$11.46$	&	$11.34$	&	$10.95$	&	$11.46$	&	$11.34$	&	$10.95$	&	$11$\\
$6$	&	$13.33$	&	$13.22$	&	$12.82$	&	$13.43$	&	$13.26$	&	$12.82$	&	$13.43$	&	$13.26$	&	$12.82$	&	$13$\\
$7$	&	$15.35$	&	$15.15$	&	$14.71$	&	$15.40$	&	$15.17$	&	$14.71$	&	$15.40$	&	$15.17$	&	$14.71$	&	$15$\\
$8$	&	$17.53$	&	$17.41$	&	$17.02$	&	$17.36$	&	$17.09$	&	$16.62$	&	$17.36$	&	$17.09$	&	$16.62$	&	$17$\\
$9$	&	$19.70$	&	$19.55$	&	$19.12$	&	$19.12$	&	$18.78$	&	$18.30$	&	$19.12$	&	$18.78$	&	$18.30$	&	$19$\\
$10$	&	$21.60$	&	$21.34$	&	$20.87$	&	$20.64$	&	$20.19$	&	$19.69$	&	$20.64$	&	$20.19$	&	$19.69$	&	$21$\\
$11$	&	$23.68$	&	$23.42$	&	$22.95$	&	$22.55$	&	$22.05$	&	$21.55$	&	$22.55$	&	$22.05$	&	$21.55$	&	$23$\\
$12$	&	$25.41$	&	$25.08$	&	$24.59$	&	$23.81$	&	$23.29$	&	$22.79$	&	$23.81$	&	$23.29$	&	$22.79$	&	$25$\\
$13$	&	$26.84$	&	$26.40$	&	$25.90$	&	$26.68$	&	$26.34$	&	$25.84$	&	$26.68$	&	$26.34$	&	$25.84$	&	$27$\\
$14$	&	$28.89$	&	$28.44$	&	$27.94$	&	$28.33$	&	$27.97$	&	$27.48$	&	$28.33$	&	$27.97$	&	$27.48$	&	$29$\\
$15$	&	$30.87$	&	$30.41$	&	$29.91$	&	$30.49$	&	$30.14$	&	$29.65$	&	$30.49$	&	$30.14$	&	$29.65$	&	$31$\\
$16$	&	$31.46$	&	$30.98$	&	$30.48$	&	$34.43$	&	$34.17$	&	$33.70$	&	$34.43$	&	$34.17$	&	$33.70$	&	$33$\\
$17$	&	$33.40$	&	$32.94$	&	$32.44$	&	$36.95$	&	$36.72$	&	$36.28$	&	$36.95$	&	$36.72$	&	$36.28$	&	$35$\\
$18$	&	$35.47$	&	$35.00$	&	$34.50$	&	$39.78$	&	$39.72$	&	$39.45$	&	$39.78$	&	$39.72$	&	$39.45$	&	$37$\\
$19$	&	$37.50$	&	$37.04$	&	$36.54$	&	$39.13$	&	$38.69$	&	$38.19$	&	$39.13$	&	$38.69$	&	$38.19$	&	$39$\\
$20$	&	$40.42$	&	$39.97$	&	$39.47$	&	$40.54$	&	$40.05$	&	$39.55$	&	$40.54$	&	$40.05$	&	$39.55$	&	$41$\\
\midrule
$25$	&	$51.06$	&	$50.59$	&	$50.09$	&	$49.42$	&	$48.93$	&	$48.43$	&	$49.42$	&	$48.93$	&	$48.43$	&	$51$\\
$30$	&	$57.88$	&	$57.39$	&	$56.89$	&	$59.49$	&	$58.99$	&	$58.49$	&	$59.49$	&	$58.99$	&	$58.49$	&	$61$\\
$35$	&	$67.63$	&	$67.13$	&	$66.63$	&	$69.24$	&	$68.74$	&	$68.24$	&	$69.24$	&	$68.74$	&	$68.24$	&	$71$\\
$40$	&	$78.34$	&	$77.84$	&	$77.34$	&	$76.37$	&	$75.87$	&	$75.37$	&	$76.37$	&	$75.87$	&	$75.37$	&	$81$\\
\bottomrule
\end{tabular}
\end{table*} 

\begin{table*}[h!]
\centering
\caption{%
Convergence orders $p_{L_{1}}^{l}$, $p_{L_{2}}^{l}$, $p_{L_{\infty}}^{l}$, calculated in norms $L_{1}$, $L_{2}$, $L_{\infty}$, of \textit{the local solution} $(\mathbf{u}_{L}, \mathbf{v}_{L})$ (represented between the nodes) of the ADER-DG method for the DAE problem (\ref{eq:hess_dae_ind_1}); $N$ is the degree of the basis polynomials $\varphi_{p}$. Orders $p^{l, u}$ are calculated for solution $\mathbf{u}_{L}$; orders $p^{l, v}$ --- for solution $\mathbf{v}_{L}$; orders $p^{l, g}$ --- for the conditions $\mathbf{g} = 0$ on the local solution $(\mathbf{u}_{L}, \mathbf{v}_{L})$. The theoretical value of convergence order $p_{\rm th.}^{l} = N+1$ is applicable for the ADER-DG method for ODE problems and is presented for comparison.
}
\label{tab:conv_orders_local_hess_1}
\begin{tabular}{@{}|l|lll|lll|lll|c|@{}}
\toprule
$N$ & $p_{L_{1}}^{l, u}$ & $p_{L_{2}}^{l, u}$ & $p_{L_{\infty}}^{l, u}$ & $p_{L_{1}}^{l, v}$ & $p_{L_{2}}^{l, v}$ & $p_{L_{\infty}}^{l, v}$ & $p_{L_{1}}^{l, g}$ & $p_{L_{2}}^{l, g}$ & $p_{L_{\infty}}^{l, g}$ & $p_{\rm th.}^{l}$ \\
\midrule
$1$	&	$2.06$	&	$2.03$	&	$1.96$	&	$1.96$	&	$1.93$	&	$1.94$	&	$2.01$	&	$2.00$	&	$1.98$	&	$2$\\
$2$	&	$3.02$	&	$3.02$	&	$2.96$	&	$4.07$	&	$4.04$	&	$3.79$	&	$2.98$	&	$2.95$	&	$2.90$	&	$3$\\
$3$	&	$4.01$	&	$4.01$	&	$3.98$	&	$5.00$	&	$4.97$	&	$4.76$	&	$4.03$	&	$4.06$	&	$4.00$	&	$4$\\
$4$	&	$5.01$	&	$5.00$	&	$4.75$	&	$5.98$	&	$5.95$	&	$5.67$	&	$5.00$	&	$5.00$	&	$4.95$	&	$5$\\
$5$	&	$6.02$	&	$6.01$	&	$5.95$	&	$7.02$	&	$7.01$	&	$6.80$	&	$5.97$	&	$5.95$	&	$5.94$	&	$6$\\
$6$	&	$7.01$	&	$7.01$	&	$6.97$	&	$7.99$	&	$7.96$	&	$7.75$	&	$7.04$	&	$7.09$	&	$7.04$	&	$7$\\
$7$	&	$8.00$	&	$7.98$	&	$7.68$	&	$8.97$	&	$8.94$	&	$8.64$	&	$8.00$	&	$8.01$	&	$7.94$	&	$8$\\
$8$	&	$9.00$	&	$9.00$	&	$8.91$	&	$10.02$	&	$10.03$	&	$9.83$	&	$8.97$	&	$8.95$	&	$8.96$	&	$9$\\
$9$	&	$10.01$	&	$10.01$	&	$9.93$	&	$10.99$	&	$10.96$	&	$10.74$	&	$10.04$	&	$10.14$	&	$10.11$	&	$10$\\
$10$	&	$11.00$	&	$10.96$	&	$10.56$	&	$11.97$	&	$11.93$	&	$11.63$	&	$11.02$	&	$11.03$	&	$10.96$	&	$11$\\
$11$	&	$12.01$	&	$12.00$	&	$11.86$	&	$13.06$	&	$13.08$	&	$12.92$	&	$11.97$	&	$11.95$	&	$11.95$	&	$12$\\
$12$	&	$13.00$	&	$12.99$	&	$12.88$	&	$13.99$	&	$13.96$	&	$13.74$	&	$13.07$	&	$13.17$	&	$13.25$	&	$13$\\
$13$	&	$13.92$	&	$13.89$	&	$13.61$	&	$14.97$	&	$14.92$	&	$14.64$	&	$14.03$	&	$14.06$	&	$13.99$	&	$14$\\
$14$	&	$15.00$	&	$14.98$	&	$14.78$	&	$16.12$	&	$16.06$	&	$15.99$	&	$14.98$	&	$14.97$	&	$14.81$	&	$15$\\
$15$	&	$16.02$	&	$16.01$	&	$16.01$	&	$17.00$	&	$16.98$	&	$16.77$	&	$16.08$	&	$16.08$	&	$15.97$	&	$16$\\
$16$	&	$16.99$	&	$17.00$	&	$16.97$	&	$17.97$	&	$17.92$	&	$17.65$	&	$17.05$	&	$17.11$	&	$17.04$	&	$17$\\
$17$	&	$17.97$	&	$17.95$	&	$17.68$	&	$18.97$	&	$18.91$	&	$18.86$	&	$18.00$	&	$17.99$	&	$17.84$	&	$18$\\
$18$	&	$19.00$	&	$19.01$	&	$18.94$	&	$20.01$	&	$20.03$	&	$19.83$	&	$18.93$	&	$18.94$	&	$19.00$	&	$19$\\
$19$	&	$20.02$	&	$20.02$	&	$19.97$	&	$20.98$	&	$20.93$	&	$20.67$	&	$20.05$	&	$20.19$	&	$20.14$	&	$20$\\
$20$	&	$21.00$	&	$20.92$	&	$20.52$	&	$21.96$	&	$21.87$	&	$21.46$	&	$21.02$	&	$21.03$	&	$20.91$	&	$21$\\
\midrule
$25$	&	$26.01$	&	$26.02$	&	$26.02$	&	$27.00$	&	$27.00$	&	$26.78$	&	$26.08$	&	$26.09$	&	$25.98$	&	$26$\\
$30$	&	$30.93$	&	$30.79$	&	$30.27$	&	$31.94$	&	$31.84$	&	$31.45$	&	$31.03$	&	$31.08$	&	$30.93$	&	$31$\\
$35$	&	$36.01$	&	$36.00$	&	$35.88$	&	$37.11$	&	$37.15$	&	$36.97$	&	$35.90$	&	$35.88$	&	$36.01$	&	$36$\\
$40$	&	$40.97$	&	$40.97$	&	$40.96$	&	$41.91$	&	$41.83$	&	$41.48$	&	$41.06$	&	$41.20$	&	$41.08$	&	$41$\\
\bottomrule
\end{tabular}
\end{table*}

The obtained results of the numerical solution of this problem are presented in Figs.~\ref{fig:hess_1_sols_uv}, \ref{fig:hess_1_sol_g_eps}, \ref{fig:hess_1_errors} and in Tables~\ref{tab:conv_orders_nodes_hess_1}, \ref{tab:conv_orders_local_hess_1}. Fig.~\ref{fig:hess_1_sols_uv} shows a comparison of the numerical solution at the nodes $(\mathbf{u}_{n}, \mathbf{v}_{n})$, the numerical local solution $(\mathbf{u}_{L}, \mathbf{v}_{L})$ and the exact analytical solution $(\mathbf{u}^{\rm ex}, \mathbf{v}^{\rm ex})$ separately for each differential $\mathbf{u}$ and algebraic $\mathbf{v}$ variable. Fig.~\ref{fig:hess_1_sol_g_eps} shows the dependencies of the feasibility of conditions $g_{1} = 0$ and $g_{2} = 0$ on the coordinate $t$, as well as the dependencies of the local errors $\varepsilon_{u}$, $\varepsilon_{v}$, $\varepsilon_{g}$ on the coordinate $t$, which allows one to quantitatively estimate the accuracy of the numerical solution, especially taking into account that the numerical solution obtained by the ADER-DG method with a local DG predictor with a high degree of polynomials $N$ very accurately corresponds to the exact analytical solution, and it is impossible to visually identify the error from the solution plots separately (in Fig.~\ref{fig:hess_1_sols_uv}). Fig.~\ref{fig:hess_1_errors} shows the dependencies of the global errors $e^{u}$, $e^{v}$, $e^{g}$ of the numerical solution at the nodes $(\mathbf{u}_{n}, \mathbf{v}_{n})$ and the local solution $(\mathbf{u}_{L}, \mathbf{v}_{L})$ on the discretization step ${\Delta t}$, separately for each differential $\mathbf{u}$ and algebraic $\mathbf{v}$ variable and the algebraic equations $\mathbf{g} = \mathbf{0}$, on the basis of which the empirical convergence orders $p$ were calculated.

The numerical solution presented in Fig.~\ref{fig:hess_1_sols_uv} demonstrates very high accuracy when compared with the exact analytical solution. In the case of polynomials of degree $N = 1$, small discontinuities of the local solution in the grid nodes $(\mathbf{u}_{n}, \mathbf{v}_{n})$ are visually observed, which is due to the linearity of the local solution representation in this case, while the amplitude of the discontinuities is significantly smaller than in the previous test example. Dissipative errors in the numerical solution, clearly observed in the previous test example, are practically not observed in this case. In the case of polynomials of degrees $N = 8$ and $40$, these errors in the numerical solution are not observed --- the accuracy of the numerical solution is very high even on a coarse grid, and it is not visually distinguished between the numerical solution and the exact analytical solution on the presented plots. The dependencies of the numerical solution errors $\varepsilon_{u}$, $\varepsilon_{v}$, $\varepsilon_{g}$ on the coordinate, presented in Fig.~\ref{fig:hess_1_sol_g_eps}, also demonstrate the properties of the numerical solution described above. In the case of polynomials of degree $N = 1$, the execution error of algebraic equation $g_{2} = 0$ systematically does not exceed $10^{-1}$ and does not demonstrate growth with the coordinate, which indicates high efficiency of the ADER-DG numerical method for solving the DAE system even in the case of low degree polynomials $N$ and a sufficiently coarse grid. The satisfaction error $|g_{1}|$ of algebraic equation $g_{1} = 0$ demonstrates approximately exponential growth, which is associated with small dissipative errors of the numerical method on a coarse grid in the case of polynomials of degree $N = 1$. The numerical solution errors in the case of polynomials of degree $N = 1$ remain sufficiently small on the selected domain of definition $[0,\, 1]$ of the desired functions. In the case of polynomials of degree $N = 8$ and $40$, the numerical solution errors $\varepsilon_{u}$, $\varepsilon_{v}$, $\varepsilon_{g}$ become very small. The satisfaction error $|g_{1}|$ of algebraic equation $g_{1} = 0$ occurs exactly for the solution at the nodes $(\mathbf{u}_{n}, \mathbf{v}_{n})$, therefore, these errors for the solution at the nodes are not presented on the plots $|g_{1}|$; therefore, the algebraic constraints in solving the DAE system are strictly observed. The errors of satisfying the algebraic equation $g_{1} = 0$ for the local solution are not zero, but they are very small and do not show a systematic increase with the coordinate. Satisfaction of the algebraic equation $g_{2} = 0$ corresponds to a very high accuracy of calculating the algebraic variable $v_{1}$. This behavior is essentially the opposite of that observed in the previous example. It is necessary to note a significant difference in the error of the local solution and the solution at the nodes --- in the case of the execution of algebraic equation $g_{2} = 0$ and errors $\varepsilon_{u}$, $\varepsilon_{v}$, $\varepsilon_{g}$, this difference reaches $10^{8}$--$10^{9}$ times for polynomials of degree $N = 8$ and $10^{35}$--$10^{55}$ times for polynomials of degree $N = 40$. The characteristic values of the errors of the local solution are $10^{-17}$--$10^{-15}$ in the case of $N = 8$ and $10^{-78}$--$10^{-83}$ in the case of $N = 40$. The characteristic values of the errors of the solution at the nodes are $10^{-25}$--$10^{-21}$ in the case of $N = 8$ and $10^{-115}$--$10^{-135}$ in the case of $N = 40$. These error values are in good agreement in order of magnitude with the first point (with the largest discretization step ${\Delta t}$) of the error dependencies $e$ on the discretization step ${\Delta t}$, presented in Fig.~\ref{fig:hess_1_errors}. These results also demonstrate a clearly expressed power-law dependence of the errors on the discretization step $e(\Delta t) \sim (\Delta t)^{p}$. In this case, in the case of polynomials of degree $N = 1$, the error values at the selected step values reach $10^{-2}$--$10^{-3}$, for polynomials of degree $N = 8$ error values reach $10^{-16}$ and $10^{-29}$ for the local solution and the solution at the nodes, respectively, and in the case of polynomials of degree $N = 40$ --- $10^{-91}$ and $10^{-145}$, respectively.

The empirical convergence orders $p^{n, u}$, $p^{n, v}$, $p^{n, g}$ for the numerical solution at the nodes $(\mathbf{u}_{n}, \mathbf{v}_{n})$ are presented in Table~\ref{tab:conv_orders_nodes_hess_1} and empirical convergence orders $p^{l, u}$, $p^{l, v}$, $p^{l, g}$ for the local solution $(\mathbf{u}_{L}, \mathbf{v}_{L})$ --- in Table~\ref{tab:conv_orders_local_hess_1}. The convergence orders were calculated separately for the norms $L_{1}$, $L_{2}$, $L_{\infty}$ (\ref{eq:norms_def}). Also for comparison, the expected values of the convergence orders, which the ADER-DG numerical method with the local DG predictor provides for solving the initial value problem for the ODE system, are presented in Tables~\ref{tab:conv_orders_nodes_hess_1} and~\ref{tab:conv_orders_local_hess_1}: $p_{\rm nodes} = 2N+1$, $p_{\rm local} = N+1$. Comparison of the obtained convergence orders $p^{n, u}$, $p^{n, v}$, $p^{n, g}$ for the solution at the nodes clearly demonstrates the expected superconvergence $2N+1$. However, for high values of polynomial degrees $N \geqslant 30$, the empirical convergence orders $p^{n, u}$, $p^{n, v}$, $p^{n, g}$ are lower than the expected value $2N+1$, while in the case of $N = 40$, the empirical values are $3$--$6$ lower than the expected value $p_{\rm nodes} = 81$. Comparison of the obtained convergence orders $p^{l, u}$, $p^{l, v}$, $p^{l, g}$ for the local solution demonstrates good agreement with expected values $N+1$. At the same time, no obvious decrease in the convergence orders for the local solution is observed. 

Therefore, it can be concluded that in the case of this example of a DAE system, the ADER-DG method with a local DG predictor allows one to obtain convergence of the numerical solution and convergence orders $p$ similar to expected. 

\subsubsection{Example 3: Hessenberg DAE system of index 2}
\label{sec:2:ct:ex3}

The third example (which contains two sub-examples) of application of the numerical method ADER-DG with local DG predictor consisted in solving a Hessenberg DAE system of index 2, which was considered in the work~\cite{hess_daes_examples_2016}:
\begin{equation}\label{eq:hess_dae_ind_2}
\begin{split}
&\ddot{x} - x(4z - 1) - 2(1 - 3t)y = 0,\hspace{19.1mm} x(0) = 0,\quad \dot{x}(0) = 0,\\
&\ddot{y} - y(4z - 1) - 2\sin(z) = 0,\hspace{23.6mm} y(0) = 1,\quad \dot{y}(0) = 0,\\
&g_{1} = x^{2} + t^{2}(y^{2} - 1) = 0\hspace{11mm} (\text{index}\ 2),\quad z(0) = 0,\quad t\in[0, 1],\\
&g_{2} = x\dot{x} + t^{2}y\dot{y} + t(y^{2} - 1) = 0\hspace{1.3mm} (\text{index}\ 1).
\end{split}
\end{equation}
The presented DAE system of index 2 contains only one algebraic equation (constraint) $g_{1} = 0$. This test example is the first example in this paper for which the numerical solution of the DAE of index 2 was investigated. 

The exact analytical solution of this problem was obtained in the following form: $x = t\sin(t - t^{2})$, $\dot{x} = \sin(t - t^{2}) + t(1 - 2t)\cos(t - t^{2})$, $y = \cos(t - t^{2})$, $\dot{y} = -(1 - 2t)\sin(t - t^{2})$, $z = t - t^{2}$. This problem was rewritten in a form consistent with the formulation of the original DAE system (\ref{eq:dae_chosen_form}):
\begin{equation}
\begin{split}
&\frac{du_{1}}{dt} = u_{3},\hspace{60.6mm} u_{1}(0) = 0,\\
&\frac{du_{2}}{dt} = u_{4},\hspace{60.6mm} u_{2}(0) = 1,\\
&\frac{du_{3}}{dt} = u_{1}(4v_{1} - 1) + 2(1 - 3t)u_{2},\hspace{25.5mm} u_{3}(0) = 0,\\
&\frac{du_{4}}{dt} = u_{2}(4v_{1} - 1) + 2\sin(v_{1}),\hspace{30mm}  u_{4}(0) = 0,\\[1.9mm]
&g_{1} = u_{1}^{2} + t^{2}(u_{2}^{2} - 1) = 0\hspace{20.2mm} (\text{index}\ 2),\quad v_{1}(0) = 0,\\[2.4mm]
&g_{2} = u_{1}u_{3} + t^{2}u_{2}u_{4} + t(u_{2}^{2} - 1) = 0\hspace{4mm} (\text{index}\ 1),
\end{split}
\end{equation}
where sets of differential variables $\mathbf{u} = [x,\, y,\, \dot{x},\, \dot{y}]^{T}$ and algebraic variables $\mathbf{v} = [z]$ were defined. The full-component exact analytical solution of the problem was written in the following form:
\begin{equation}
\begin{split}
&\mathbf{u}^{\rm ex} = \left[
\begin{array}{l}
t\sin(t - t^{2})\\
\cos(t - t^{2})\\
\sin(t - t^{2}) + t(1 - 2t)\cos(t - t^{2})\\
-(1 - 2t)\sin(t - t^{2})
\end{array}
\right],\\
&\mathbf{v}^{\rm ex} = \Big[t - t^{2}\Big].
\end{split}
\end{equation}
The domain of definition $[0,\, 1]$ of the desired functions $\mathbf{u}$ and $\mathbf{v}$ was discretized into $L = 8$, $10$, $12$, $14$, $16$, $18$ discretization domains $\Omega_{n}$ with the same discretization steps ${\Delta t}_{n} \equiv {\Delta t} = 1/(L-1)$, which was done to be able to calculate the empirical convergence orders $p$. Therefore, the empirical convergence orders $p$ were calculated by least squares approximation of the dependence of the global error $e$ on the grid discretization step ${\Delta t}$ at 6 data points.

As a result of differentiation of algebraic equation $g_{1} = 0$ by independent variable $t$, the DAE index of the system decreases by 1 and becomes equal to 1, while the algebraic equation $g_{1} = 0$ in the DAE system is replaced by the algebraic equation $g_{2} = 0$ --- this procedure is a decrease in the index of the DAE system. In this example, in the present work, two separate DAE systems were studied: a DAE system with an algebraic equation $g_{1} = 0$ (index 2), which was called the DAE system (\ref{eq:hess_dae_ind_2}) of index 2, and a DAE system with an algebraic equation $g_{2} = 0$ (index 1), which was called the DAE system (\ref{eq:hess_dae_ind_2}) of index 1. With an increase in the DAE index of the system, deterioration of the accuracy and convergence properties of the ADER-DG numerical method with a local DG predictor was expected, this property is well known for other numerical methods for solving DAE systems of equations --- usually a decrease in the empirical convergence order is observed compared to the expected convergence order of the numerical method, which occurs when solving the initial value problem for the ODE system~\cite{Hairer_book_2}. Therefore, interest arose in a quantitative study of the convergence of the numerical solution for the original DAE system of index 2 and the DAE system with a decreased index.

It is necessary to note an interesting property of the numerical solution, which was already observed in the examples presented above --- an algebraic equation, which is solved numerically explicitly, is satisfied exactly (within the limits of the accuracy of representing real numbers by floating-point numbers) for the solution at the nodes $(\mathbf{u}_{n}, \mathbf{v}_{n})$. Therefore, when separately solving the DAE system (\ref{eq:hess_dae_ind_2}) of index 2 and the DAE system (\ref{eq:hess_dae_ind_2}) of index 1, the errors were calculated that determine the satisfaction of both algebraic equation $g_{1} = 0$ and algebraic equation $g_{2} = 0$ were calculated --- it was expected that the numerical solution at the nodes $(\mathbf{u}_{n}, \mathbf{v}_{n})$ would exactly satisfy only that algebraic equation, which is explicitly included in the solved DAE system.

\begin{figure}[h!]
\captionsetup[subfigure]{%
	position=bottom,
	font+=smaller,
	textfont=normalfont,
	singlelinecheck=off,
	justification=raggedright
}
\centering
\begin{subfigure}{0.320\textwidth}
\includegraphics[width=\textwidth]{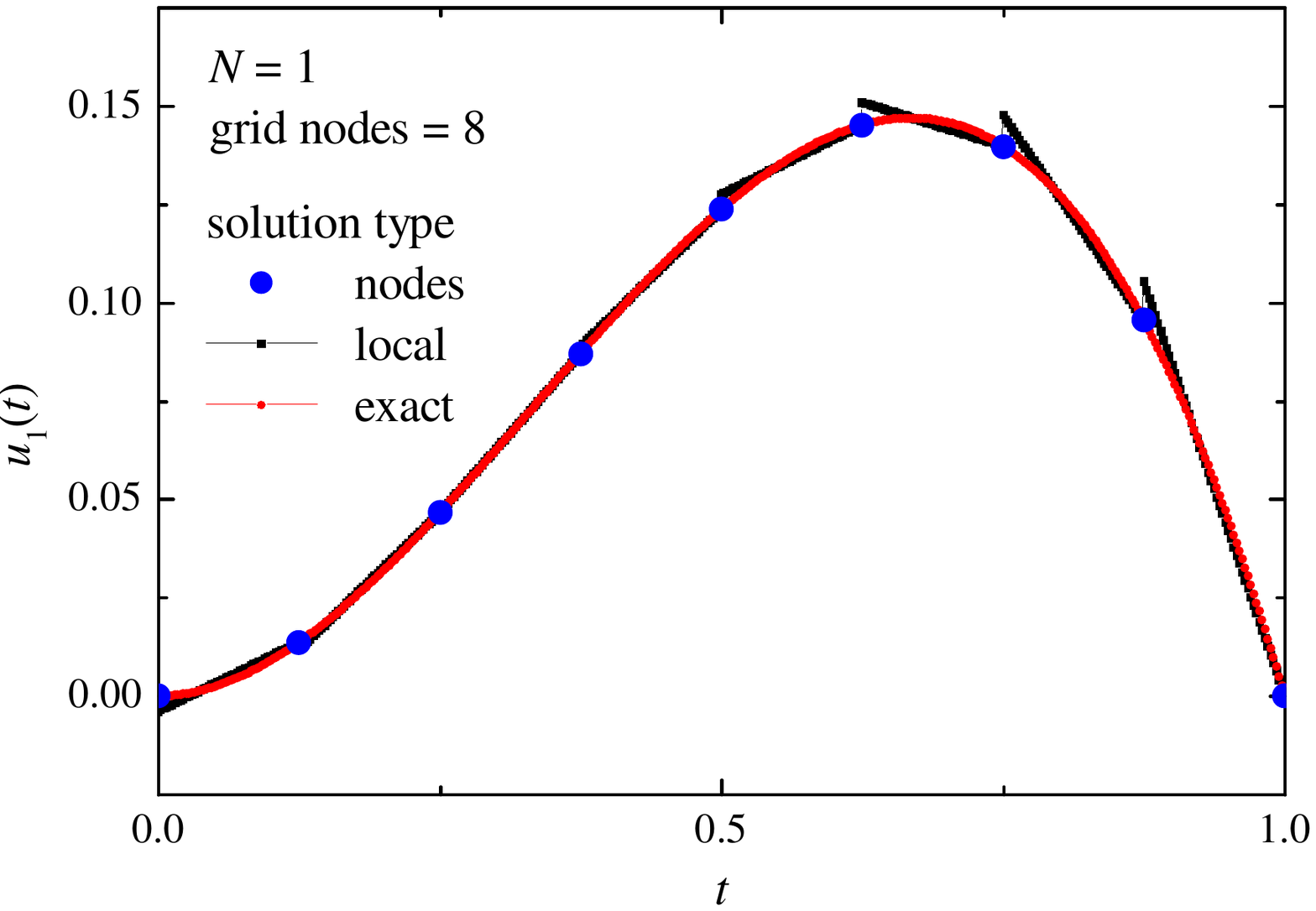}
\vspace{-8mm}\caption{\label{fig:hess_2_ind2_sol_uv:a1}}
\end{subfigure}
\begin{subfigure}{0.320\textwidth}
\includegraphics[width=\textwidth]{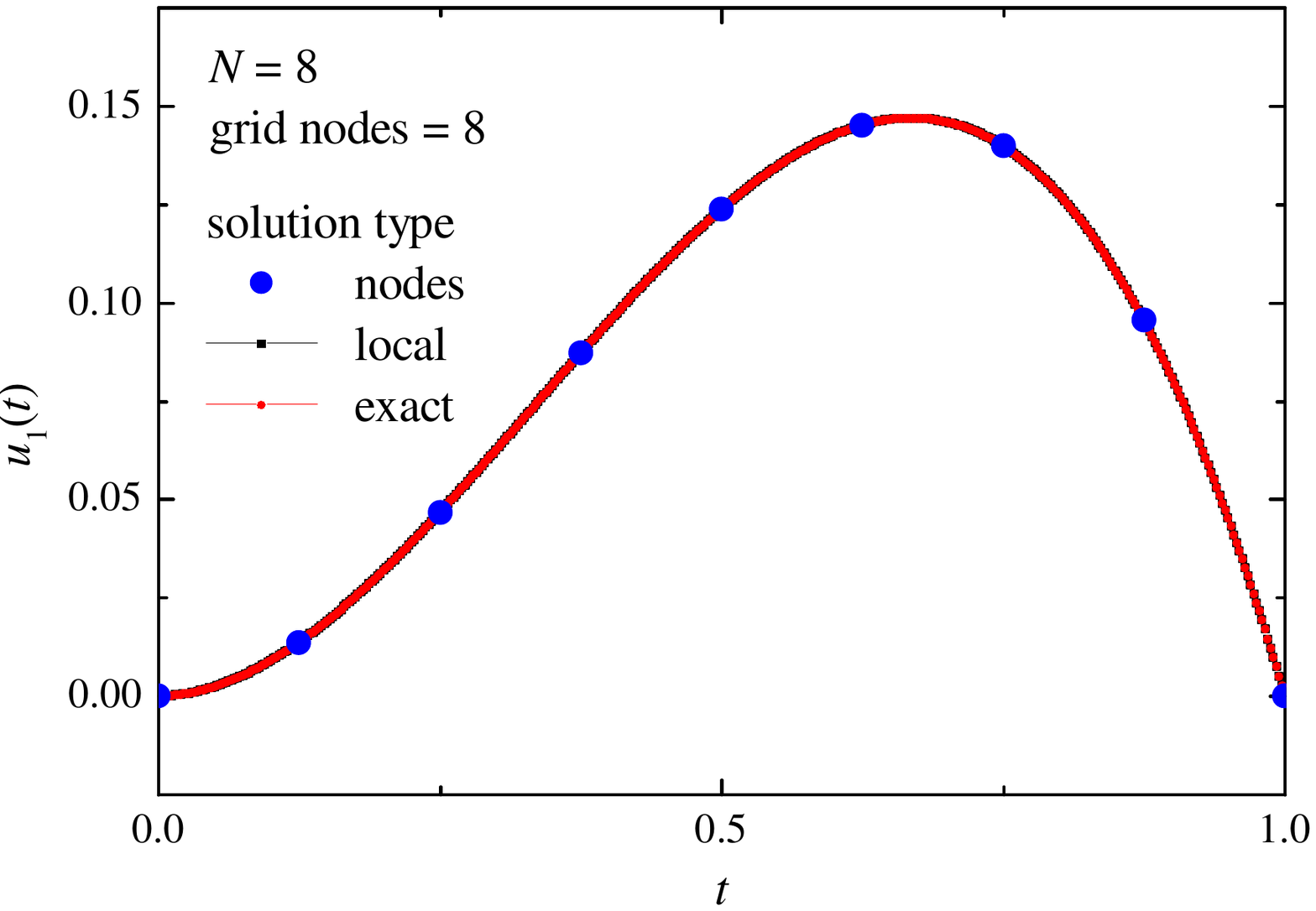}
\vspace{-8mm}\caption{\label{fig:hess_2_ind2_sol_uv:a2}}
\end{subfigure}
\begin{subfigure}{0.320\textwidth}
\includegraphics[width=\textwidth]{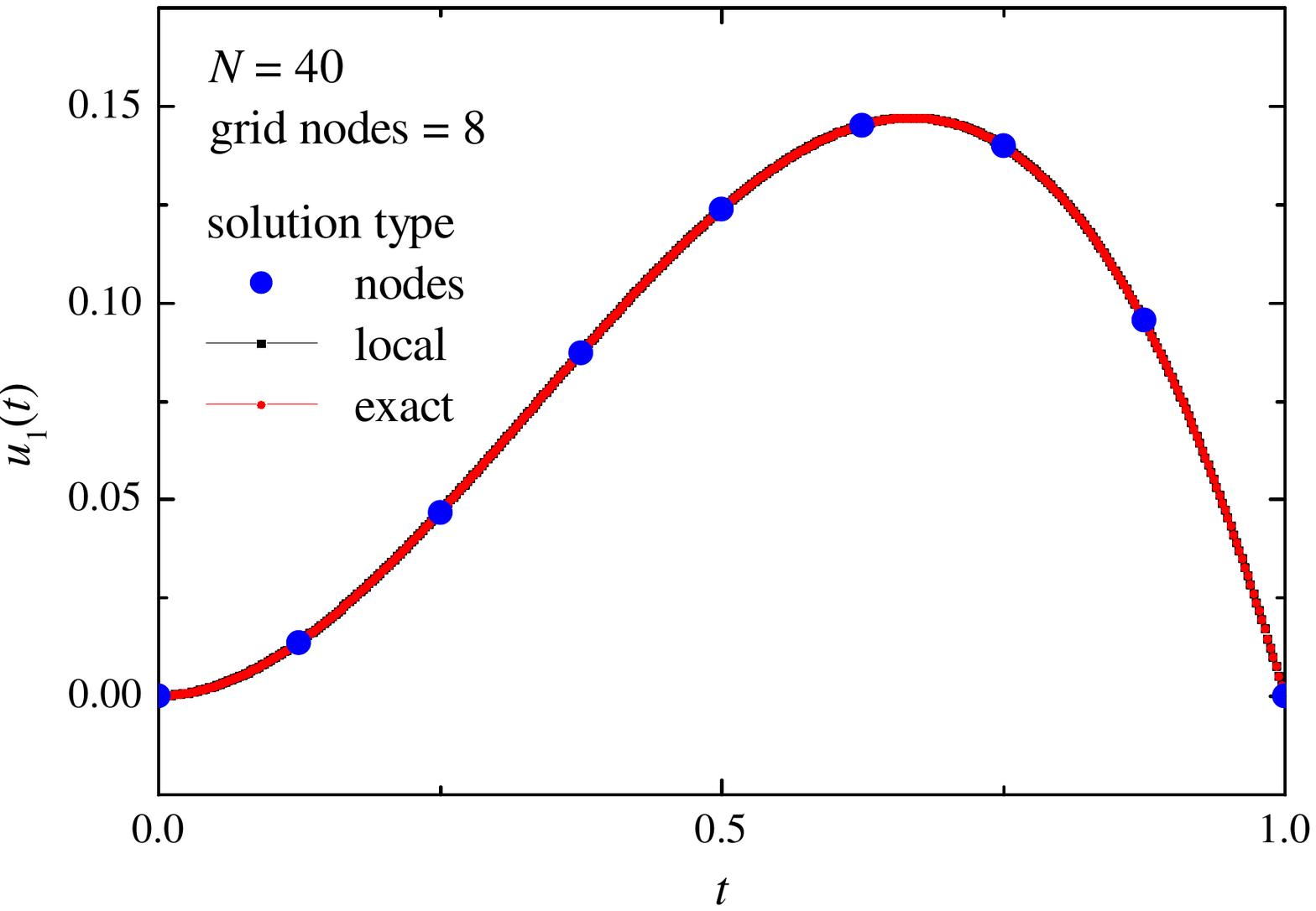}
\vspace{-8mm}\caption{\label{fig:hess_2_ind2_sol_uv:a3}}
\end{subfigure}\\
\begin{subfigure}{0.320\textwidth}
\includegraphics[width=\textwidth]{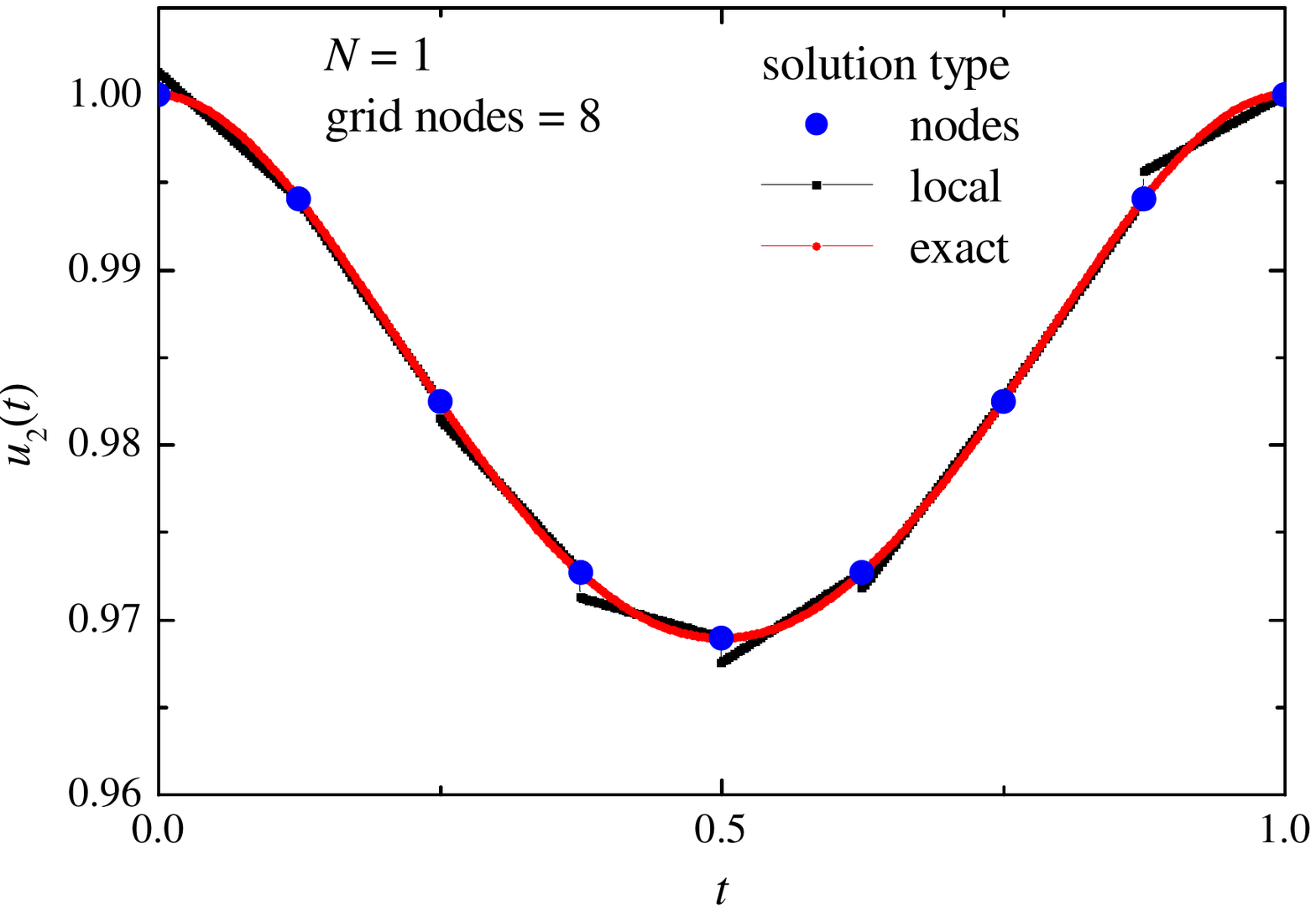}
\vspace{-8mm}\caption{\label{fig:hess_2_ind2_sol_uv:b1}}
\end{subfigure}
\begin{subfigure}{0.320\textwidth}
\includegraphics[width=\textwidth]{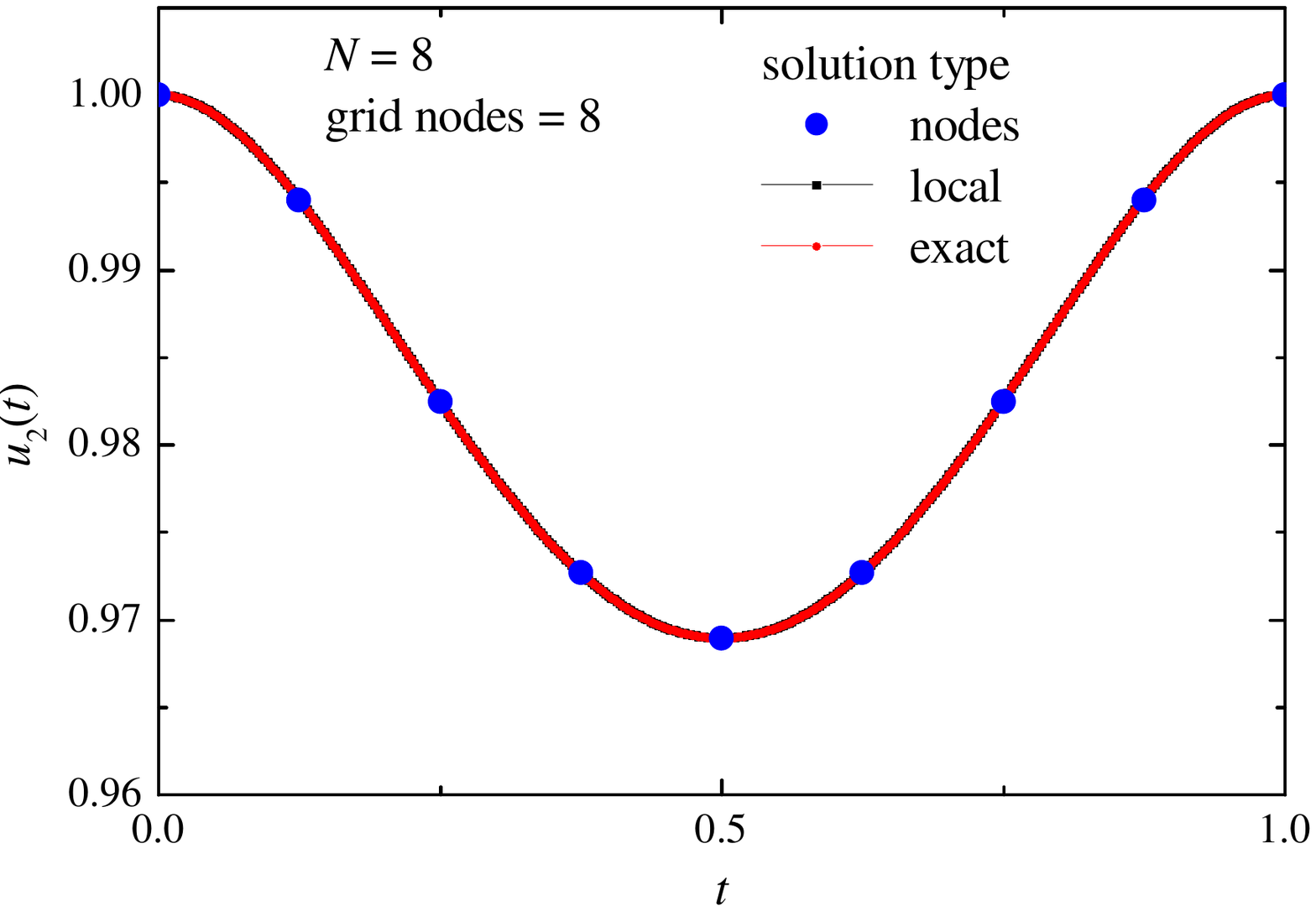}
\vspace{-8mm}\caption{\label{fig:hess_2_ind2_sol_uv:b2}}
\end{subfigure}
\begin{subfigure}{0.320\textwidth}
\includegraphics[width=\textwidth]{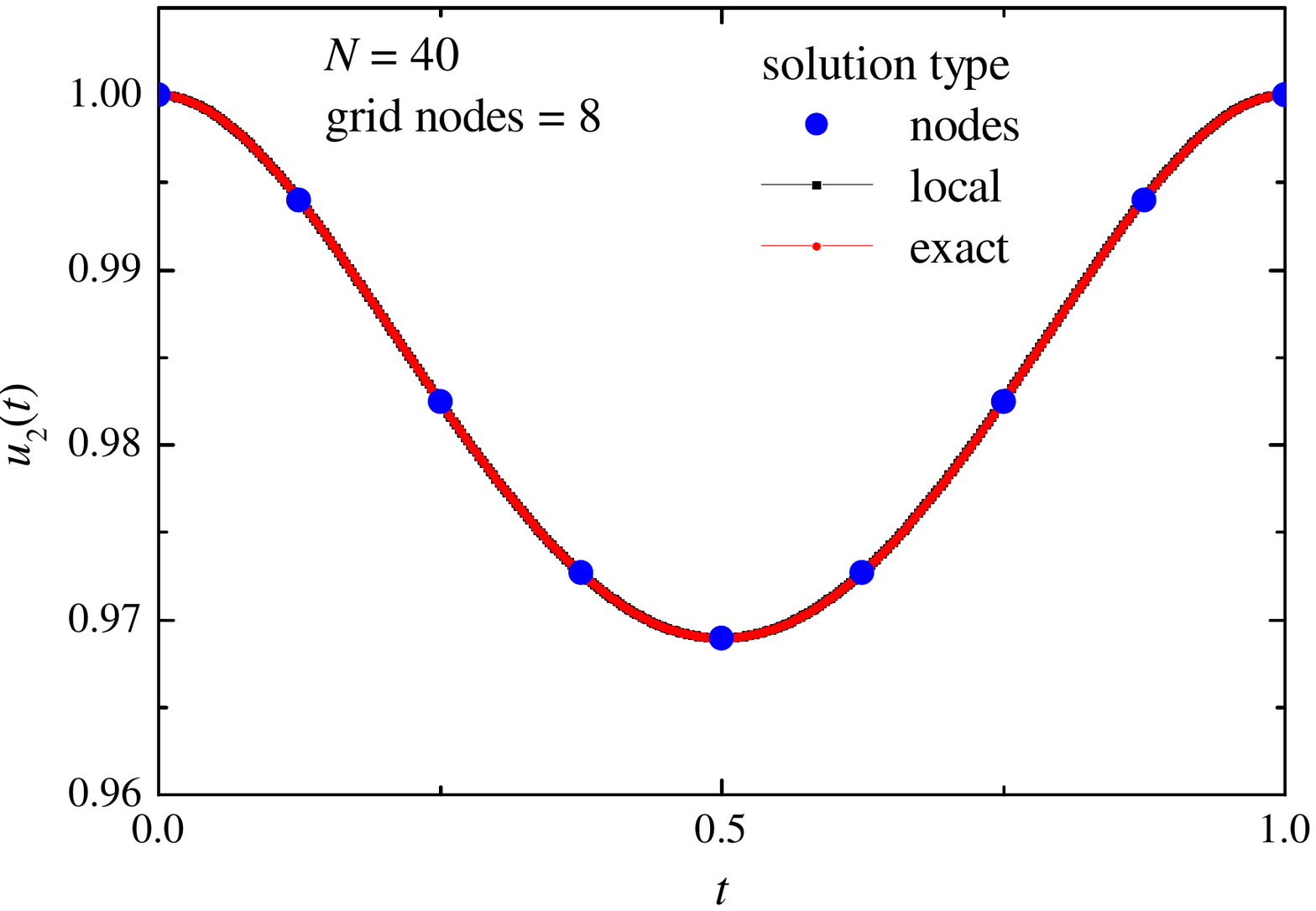}
\vspace{-8mm}\caption{\label{fig:hess_2_ind2_sol_uv:b3}}
\end{subfigure}\\
\begin{subfigure}{0.320\textwidth}
\includegraphics[width=\textwidth]{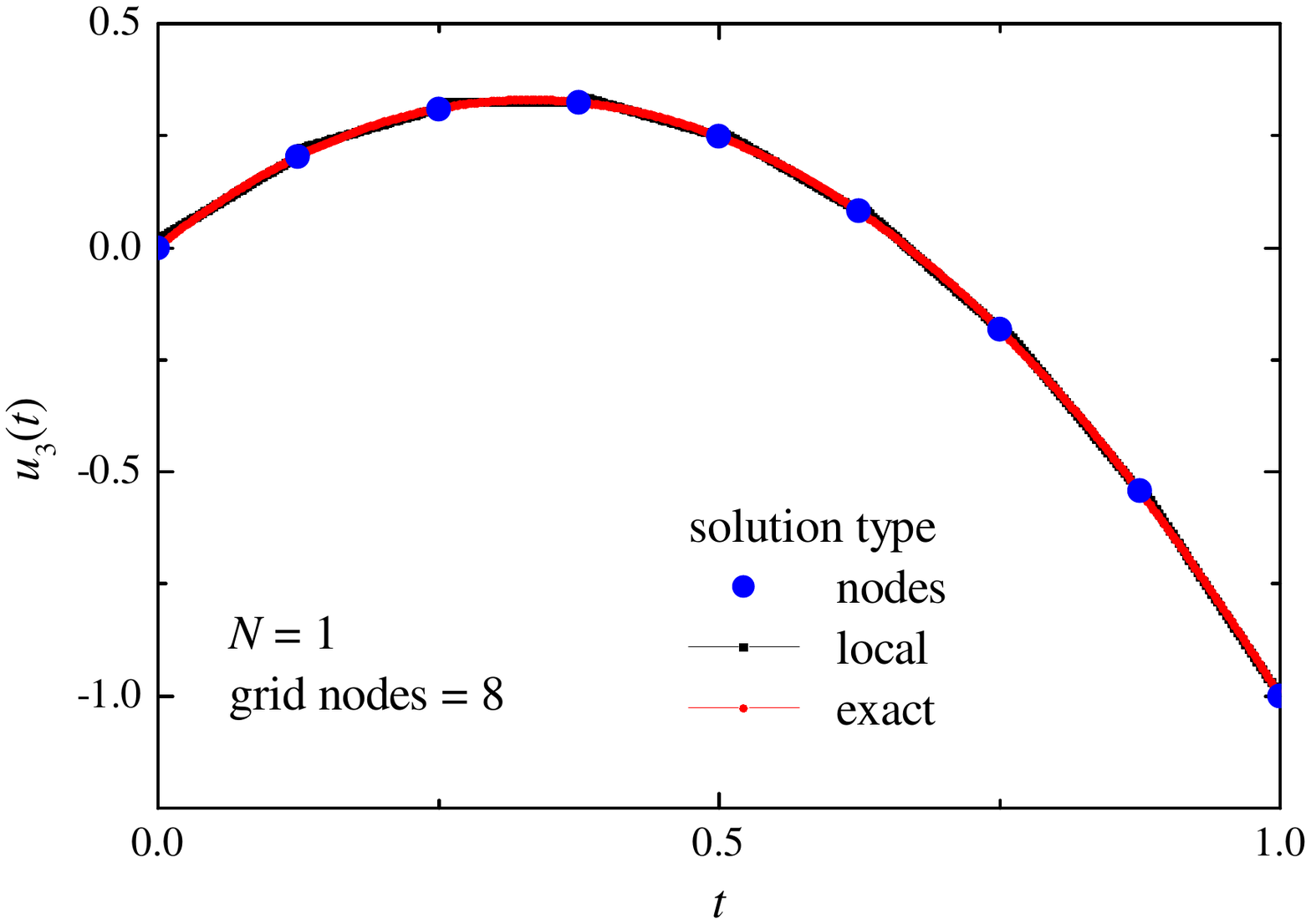}
\vspace{-8mm}\caption{\label{fig:hess_2_ind2_sol_uv:c1}}
\end{subfigure}
\begin{subfigure}{0.320\textwidth}
\includegraphics[width=\textwidth]{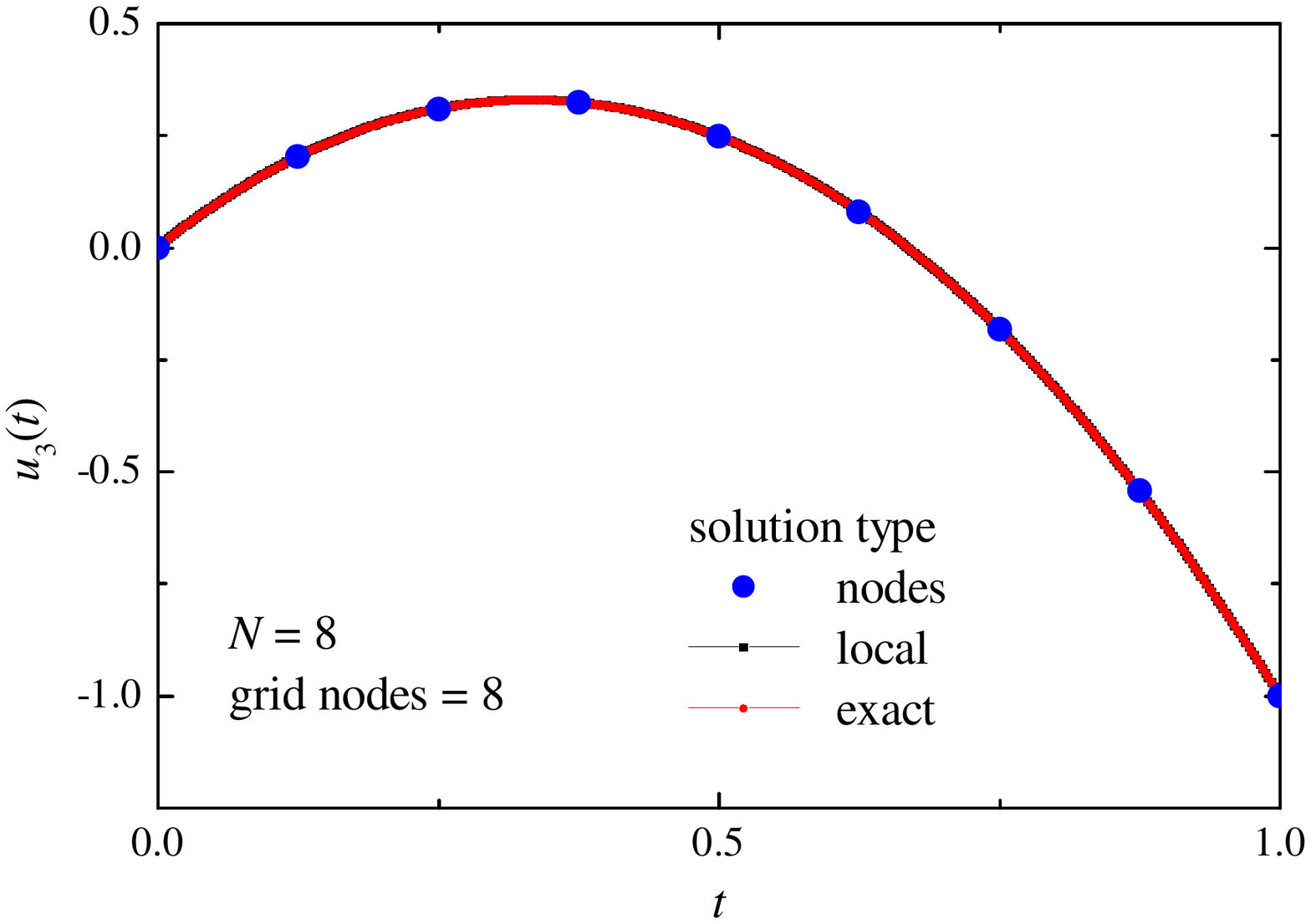}
\vspace{-8mm}\caption{\label{fig:hess_2_ind2_sol_uv:c2}}
\end{subfigure}
\begin{subfigure}{0.320\textwidth}
\includegraphics[width=\textwidth]{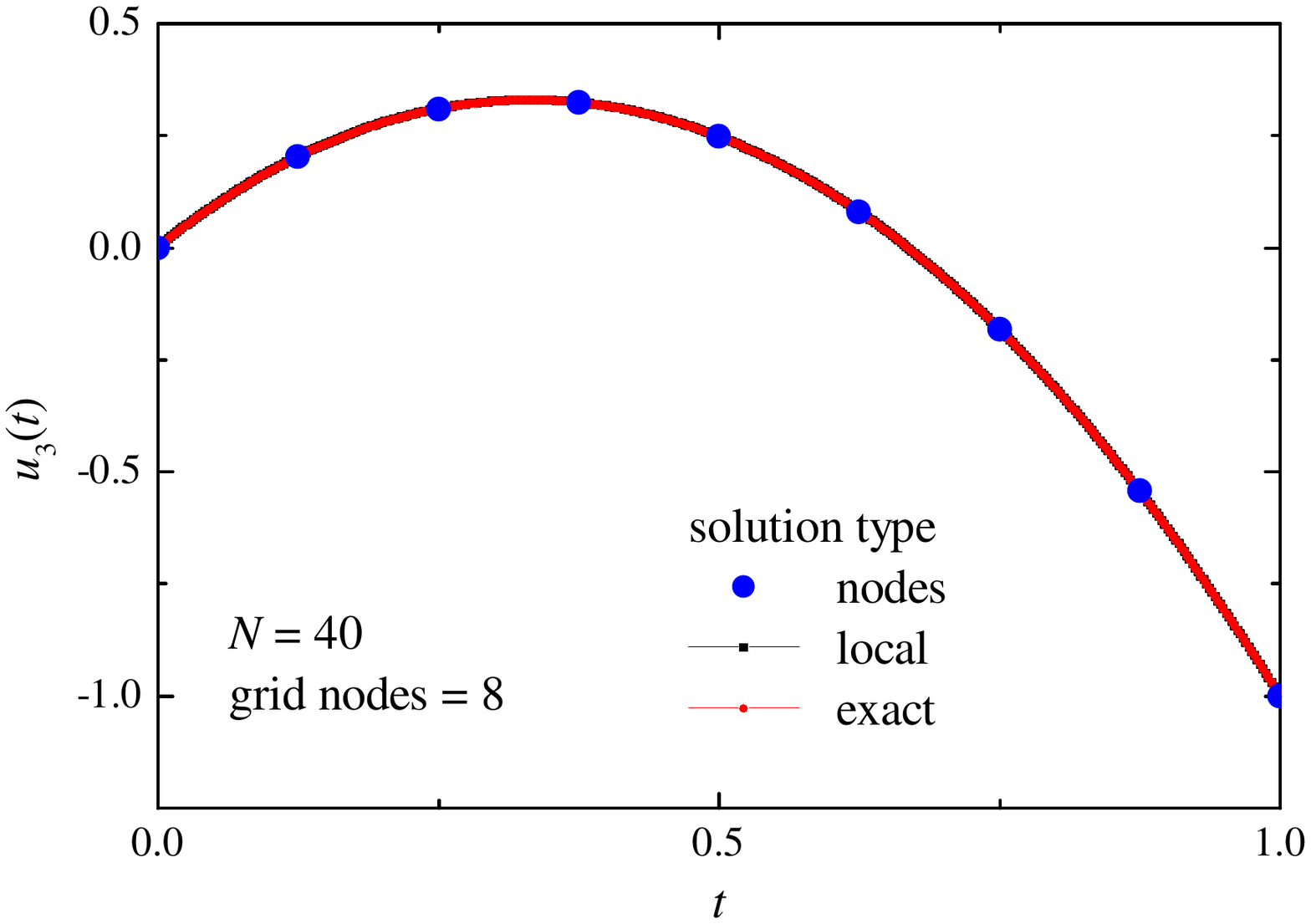}
\vspace{-8mm}\caption{\label{fig:hess_2_ind2_sol_uv:c3}}
\end{subfigure}\\
\begin{subfigure}{0.320\textwidth}
\includegraphics[width=\textwidth]{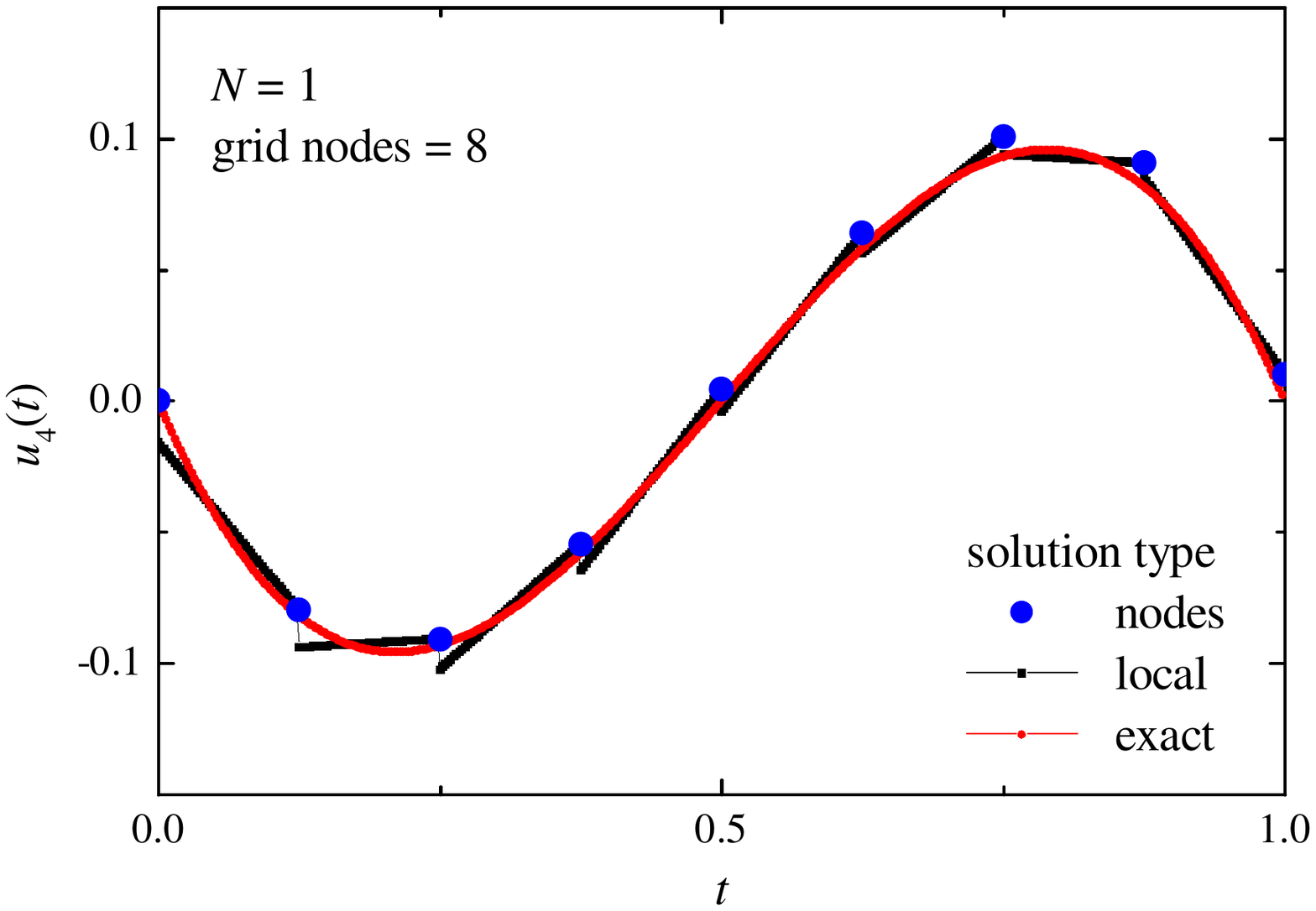}
\vspace{-8mm}\caption{\label{fig:hess_2_ind2_sol_uv:d1}}
\end{subfigure}
\begin{subfigure}{0.320\textwidth}
\includegraphics[width=\textwidth]{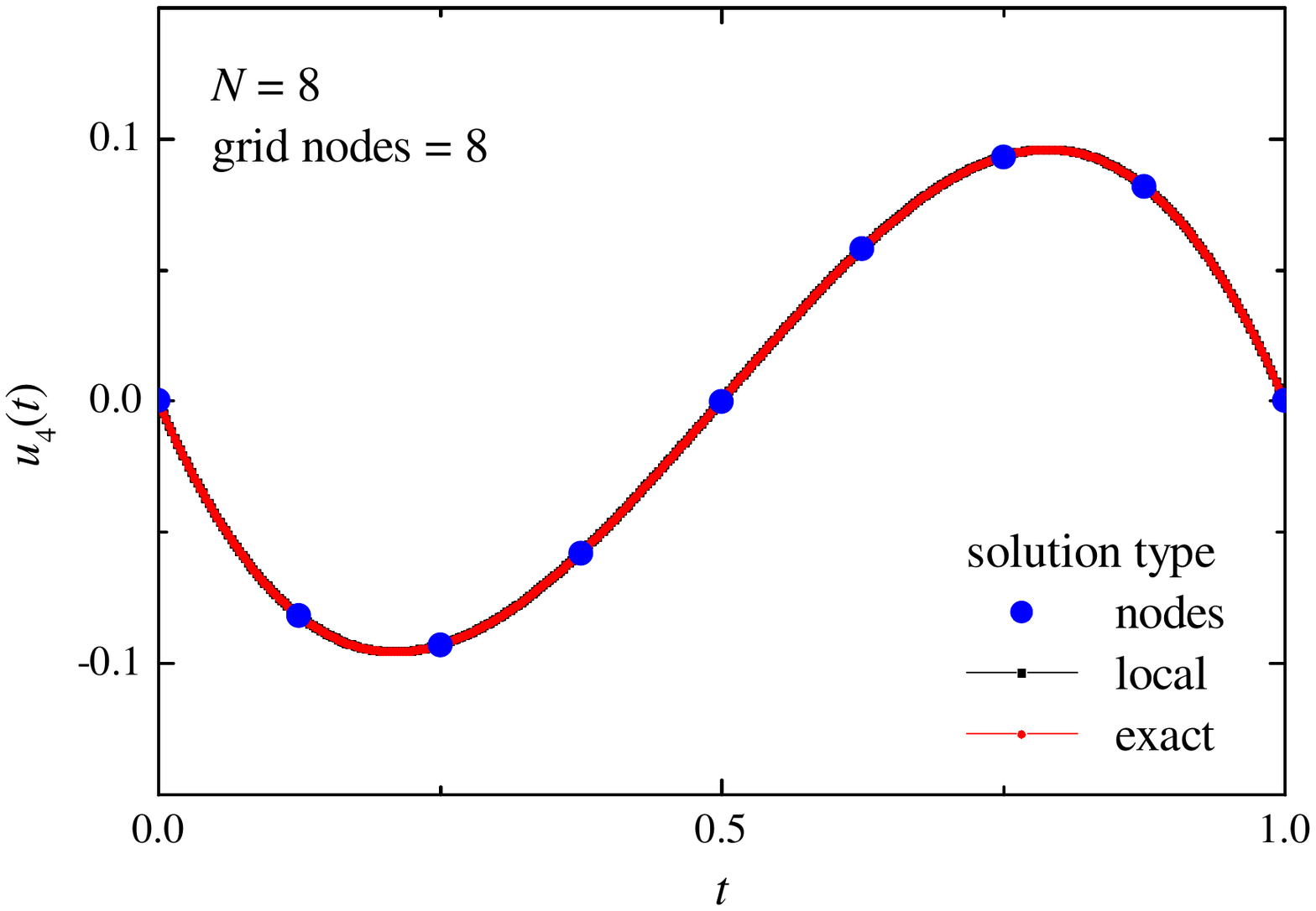}
\vspace{-8mm}\caption{\label{fig:hess_2_ind2_sol_uv:d2}}
\end{subfigure}
\begin{subfigure}{0.320\textwidth}
\includegraphics[width=\textwidth]{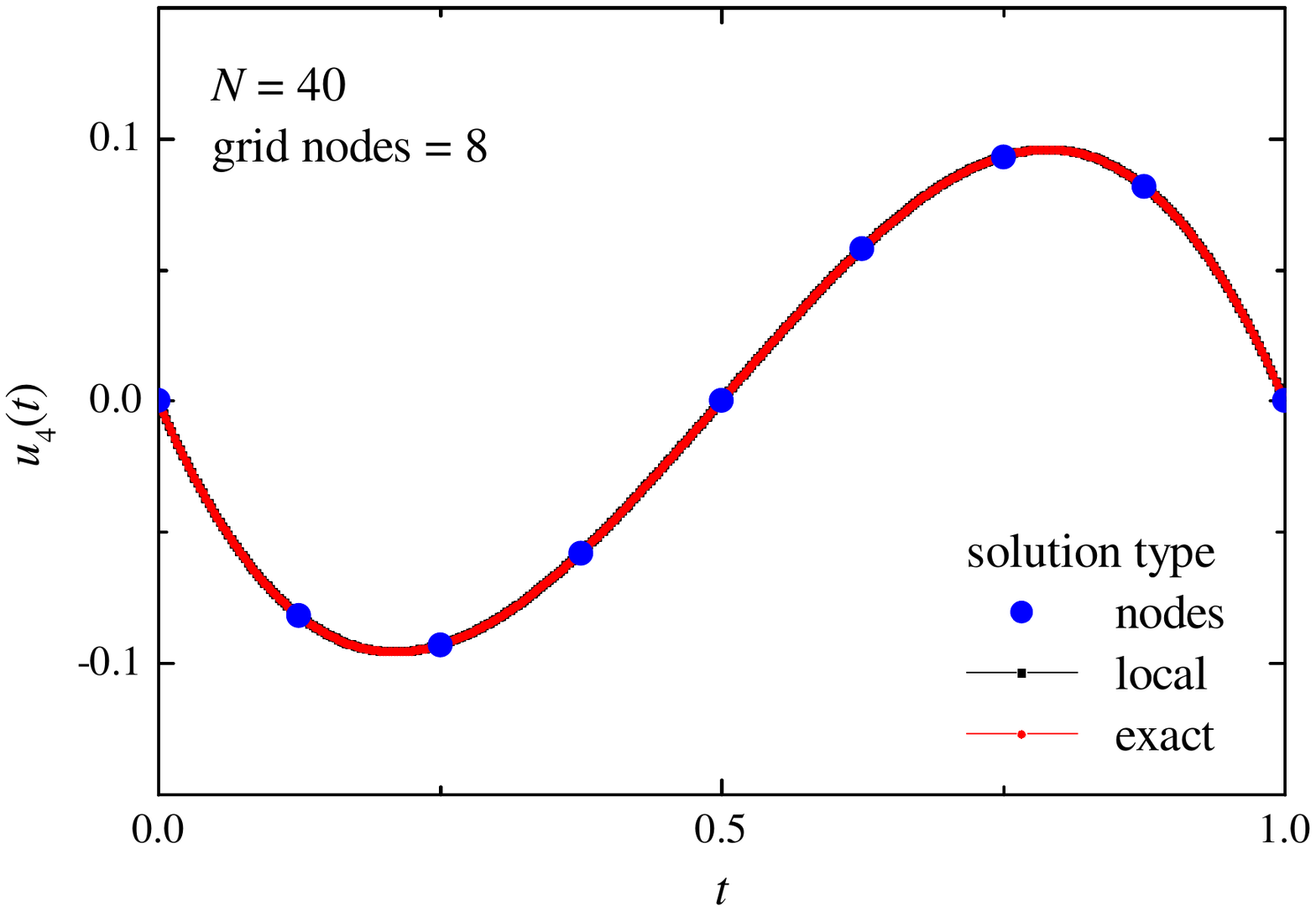}
\vspace{-8mm}\caption{\label{fig:hess_2_ind2_sol_uv:d3}}
\end{subfigure}\\
\begin{subfigure}{0.320\textwidth}
\includegraphics[width=\textwidth]{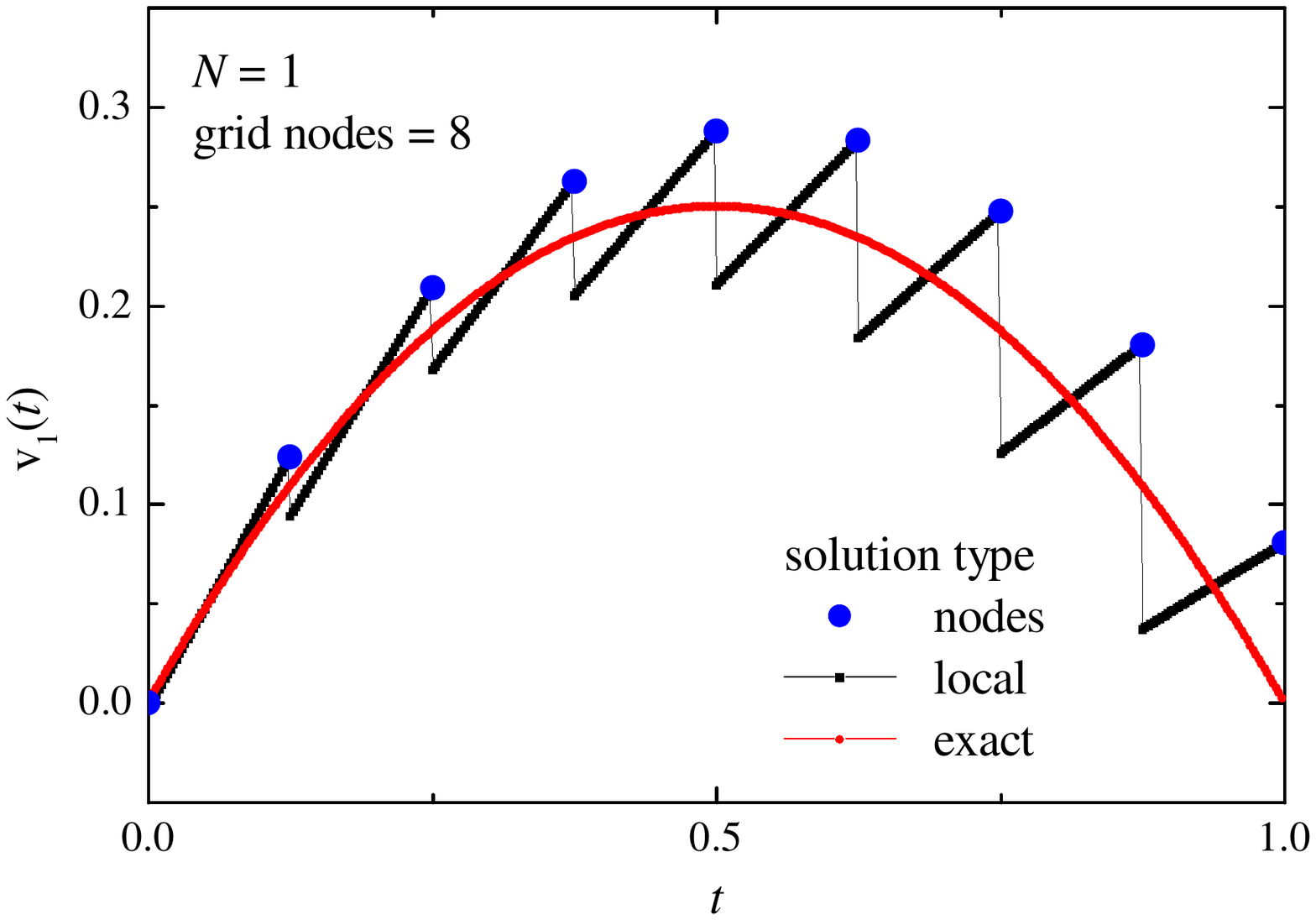}
\vspace{-8mm}\caption{\label{fig:hess_2_ind2_sol_uv:e1}}
\end{subfigure}
\begin{subfigure}{0.320\textwidth}
\includegraphics[width=\textwidth]{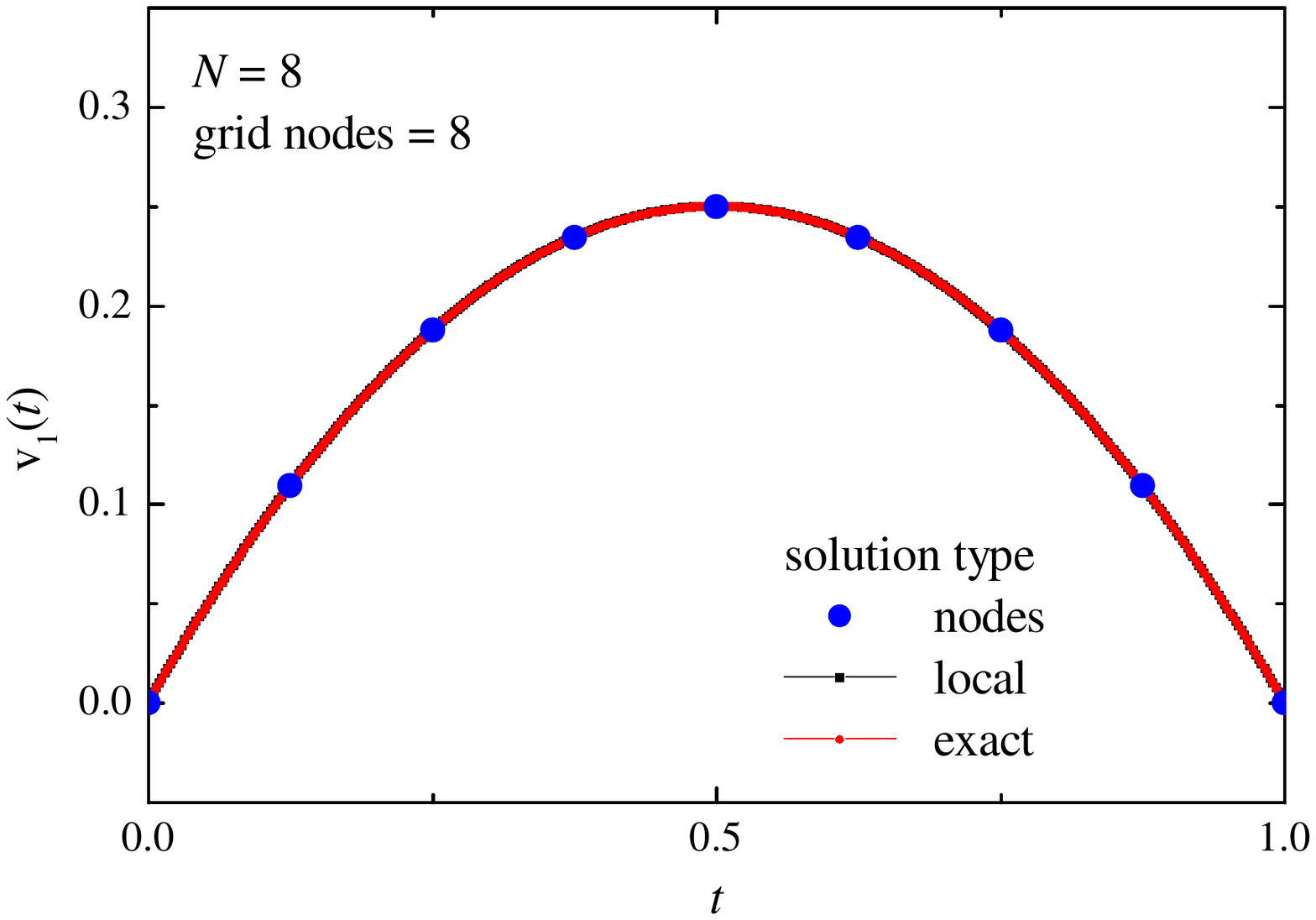}
\vspace{-8mm}\caption{\label{fig:hess_2_ind2_sol_uv:e2}}
\end{subfigure}
\begin{subfigure}{0.320\textwidth}
\includegraphics[width=\textwidth]{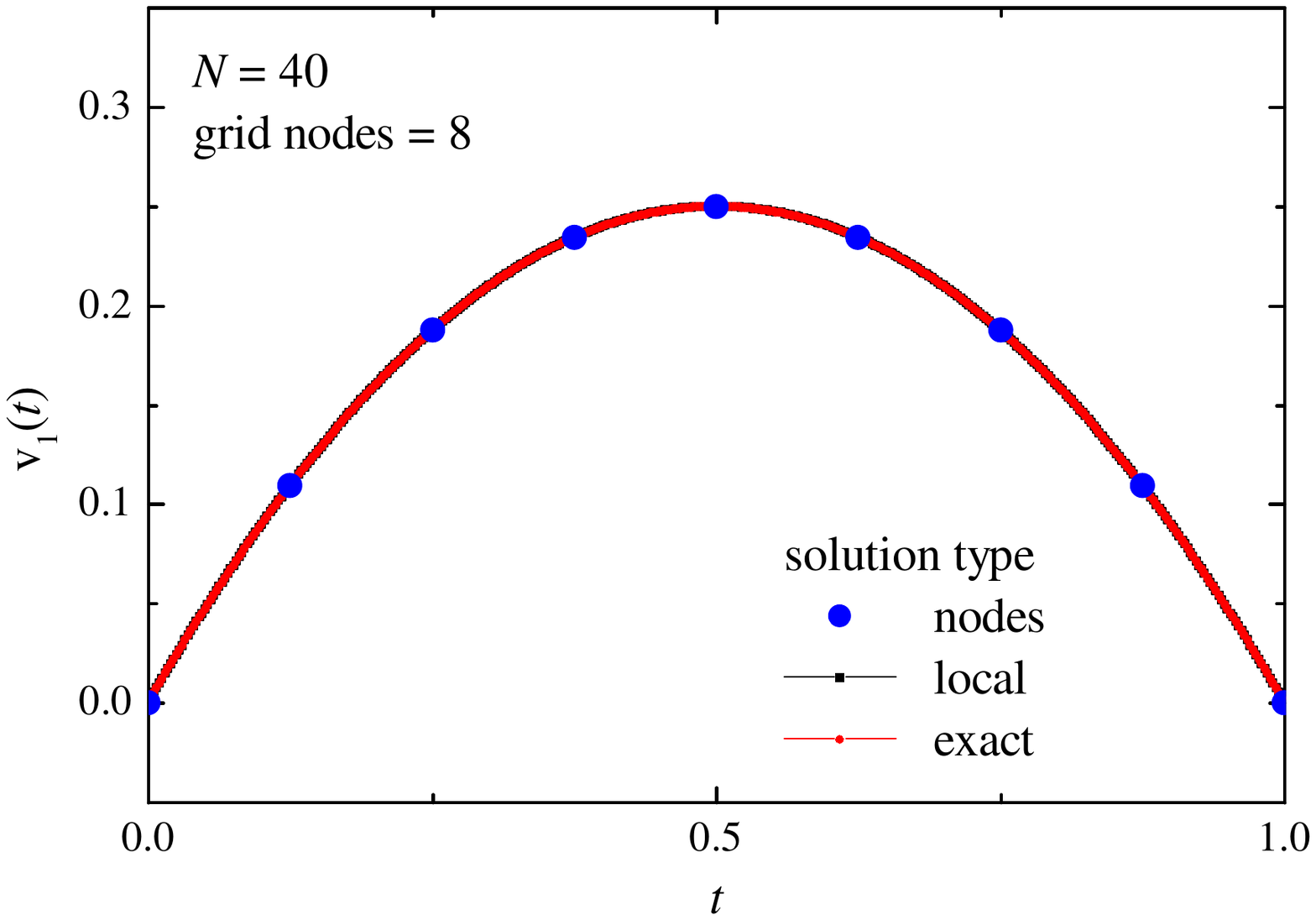}
\vspace{-8mm}\caption{\label{fig:hess_2_ind2_sol_uv:e3}}
\end{subfigure}\\
\caption{%
Numerical solution of the DAE system (\ref{eq:hess_dae_ind_2}) of index 2. Comparison of the solution at nodes $\mathbf{u}_{n}$, the local solution $\mathbf{u}_{L}(t)$ and the exact solution $\mathbf{u}^{\rm ex}(t)$ for components $u_{1}$ (\subref{fig:hess_2_ind2_sol_uv:a1}, \subref{fig:hess_2_ind2_sol_uv:a2}, \subref{fig:hess_2_ind2_sol_uv:a3}), $u_{2}$ (\subref{fig:hess_2_ind2_sol_uv:b1}, \subref{fig:hess_2_ind2_sol_uv:b2}, \subref{fig:hess_2_ind2_sol_uv:b3}), $u_{3}$ (\subref{fig:hess_2_ind2_sol_uv:c1}, \subref{fig:hess_2_ind2_sol_uv:c2}, \subref{fig:hess_2_ind2_sol_uv:c3}), $u_{4}$ (\subref{fig:hess_2_ind2_sol_uv:d1}, \subref{fig:hess_2_ind2_sol_uv:d2}, \subref{fig:hess_2_ind2_sol_uv:d3}) and $v_{1}$ (\subref{fig:hess_2_ind2_sol_uv:e1}, \subref{fig:hess_2_ind2_sol_uv:e2}, \subref{fig:hess_2_ind2_sol_uv:e3}), obtained using polynomials with degrees $N = 1$ (\subref{fig:hess_2_ind2_sol_uv:a1}, \subref{fig:hess_2_ind2_sol_uv:b1}, \subref{fig:hess_2_ind2_sol_uv:c1}, \subref{fig:hess_2_ind2_sol_uv:d1}, \subref{fig:hess_2_ind2_sol_uv:e1}), $N = 8$ (\subref{fig:hess_2_ind2_sol_uv:a2}, \subref{fig:hess_2_ind2_sol_uv:b2}, \subref{fig:hess_2_ind2_sol_uv:c2}, \subref{fig:hess_2_ind2_sol_uv:d2}, \subref{fig:hess_2_ind2_sol_uv:e2}), $N = 40$ (\subref{fig:hess_2_ind2_sol_uv:a3}, \subref{fig:hess_2_ind2_sol_uv:b3}, \subref{fig:hess_2_ind2_sol_uv:c3}, \subref{fig:hess_2_ind2_sol_uv:d3}, \subref{fig:hess_2_ind2_sol_uv:e3}).
}
\label{fig:hess_2_ind2_sol_uv}
\end{figure} 

\begin{figure}[h!]
\captionsetup[subfigure]{%
	position=bottom,
	font+=smaller,
	textfont=normalfont,
	singlelinecheck=off,
	justification=raggedright
}
\centering
\begin{subfigure}{0.320\textwidth}
\includegraphics[width=\textwidth]{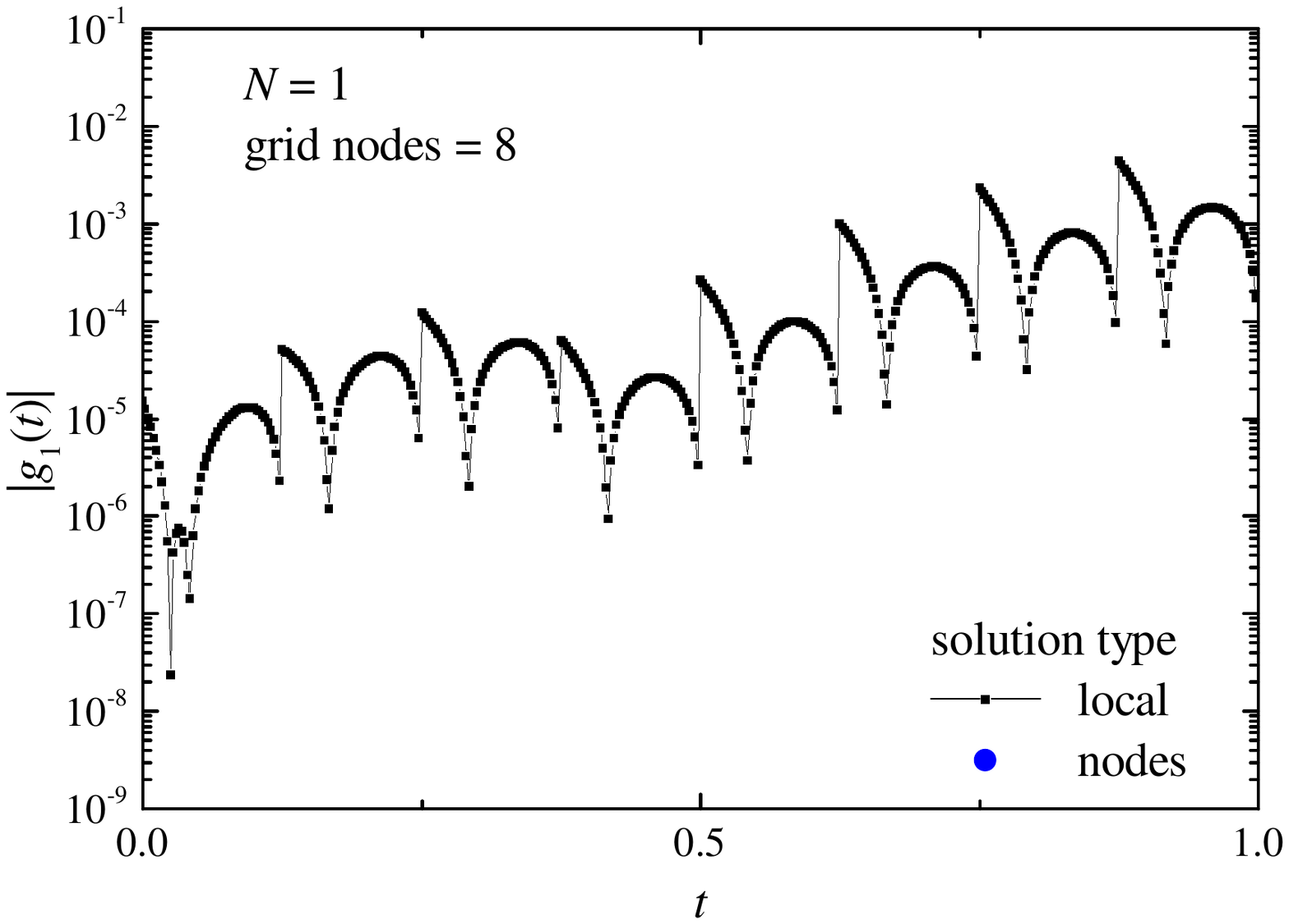}
\vspace{-8mm}\caption{\label{fig:hess_2_ind2_sol_g_eps:a1}}
\end{subfigure}
\begin{subfigure}{0.320\textwidth}
\includegraphics[width=\textwidth]{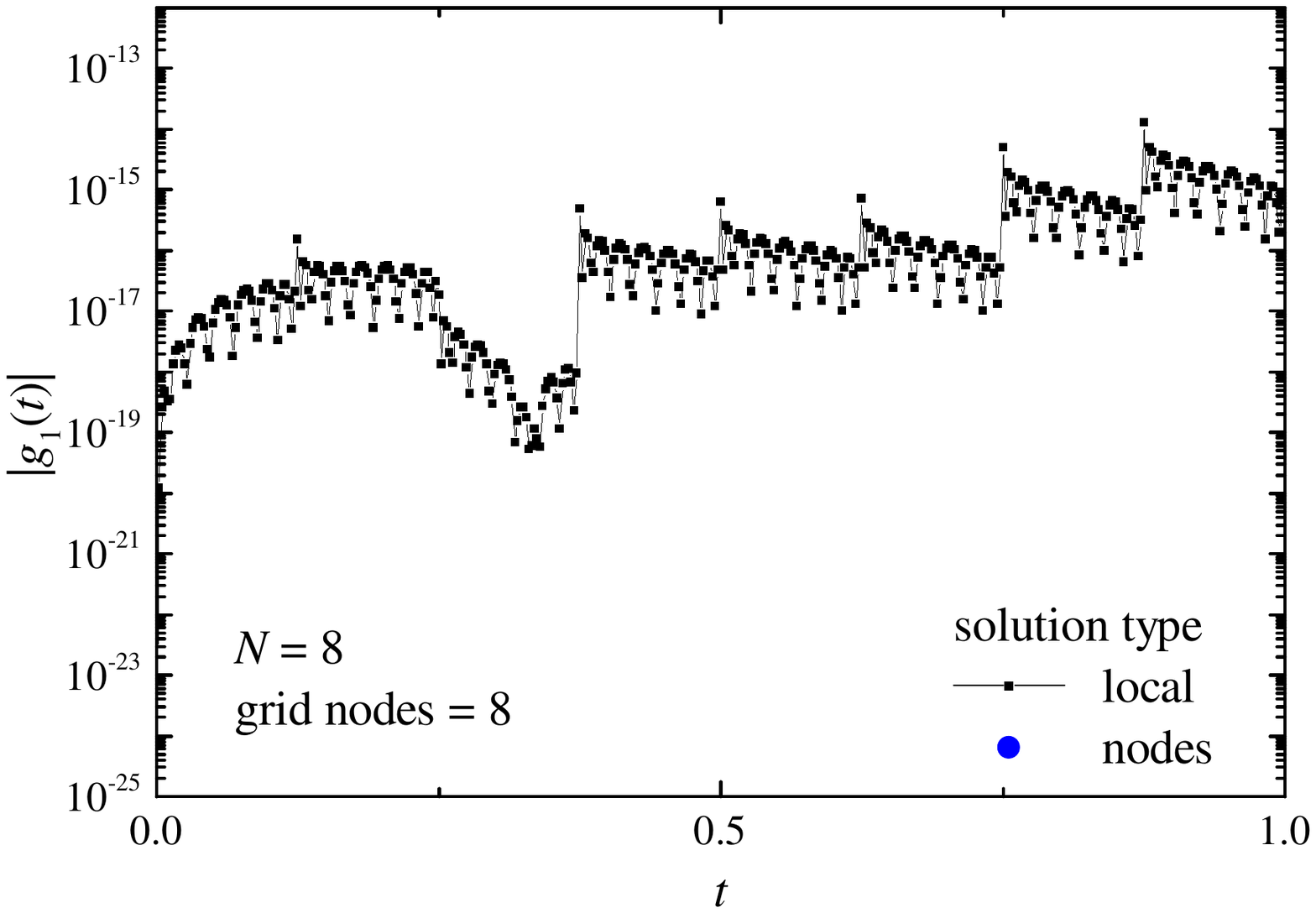}
\vspace{-8mm}\caption{\label{fig:hess_2_ind2_sol_g_eps:a2}}
\end{subfigure}
\begin{subfigure}{0.320\textwidth}
\includegraphics[width=\textwidth]{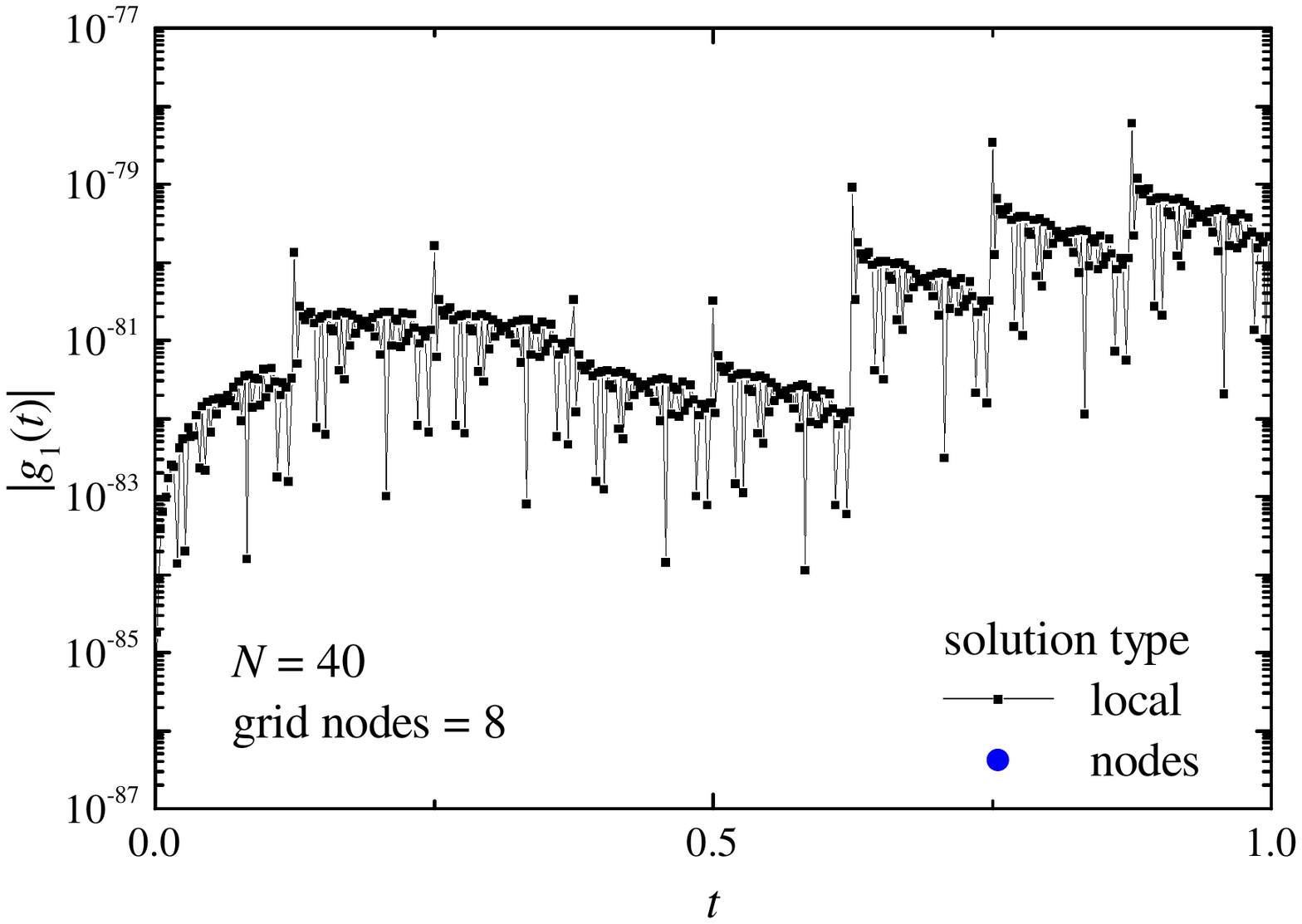}
\vspace{-8mm}\caption{\label{fig:hess_2_ind2_sol_g_eps:a3}}
\end{subfigure}\\
\begin{subfigure}{0.320\textwidth}
\includegraphics[width=\textwidth]{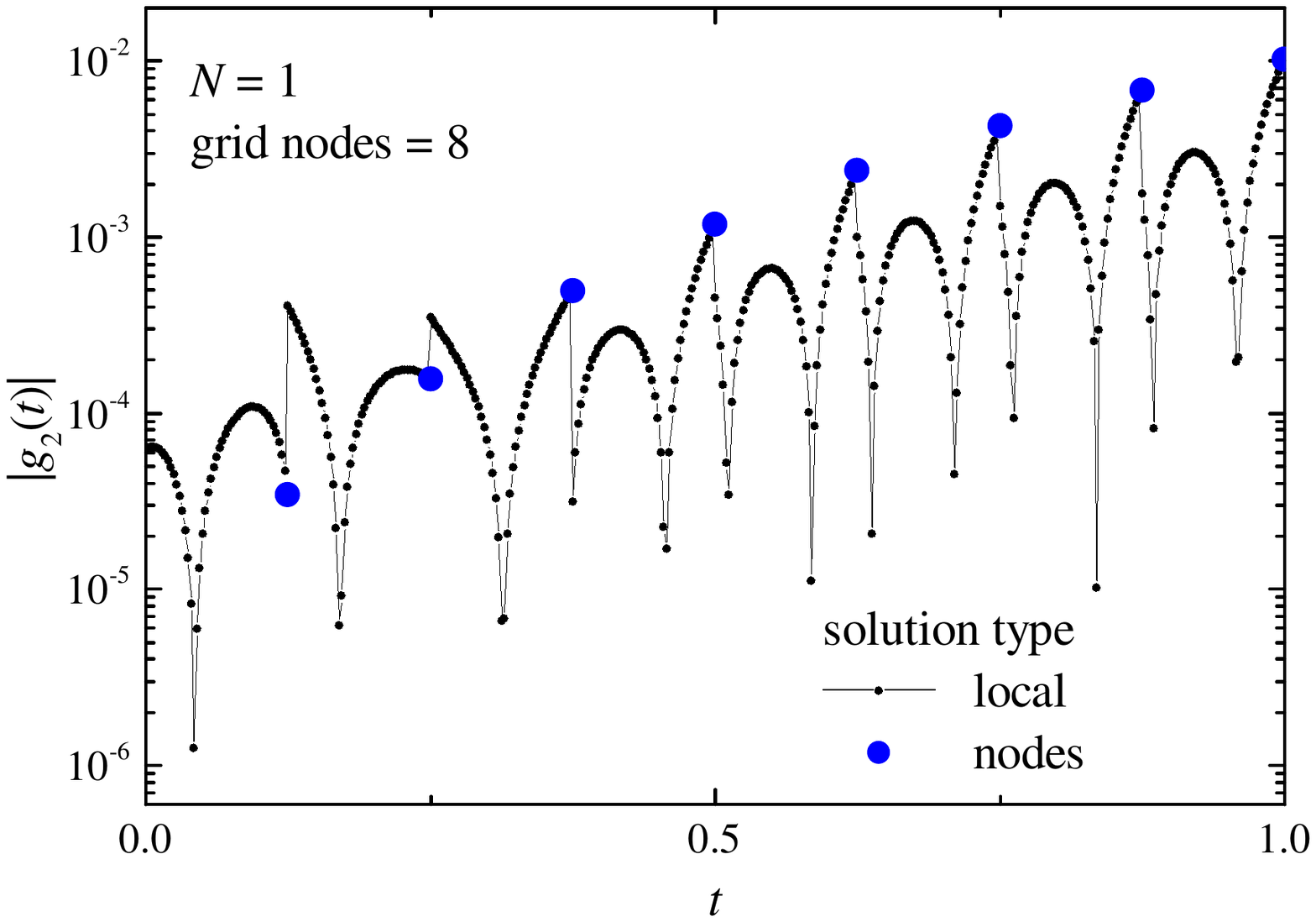}
\vspace{-8mm}\caption{\label{fig:hess_2_ind2_sol_g_eps:b1}}
\end{subfigure}
\begin{subfigure}{0.320\textwidth}
\includegraphics[width=\textwidth]{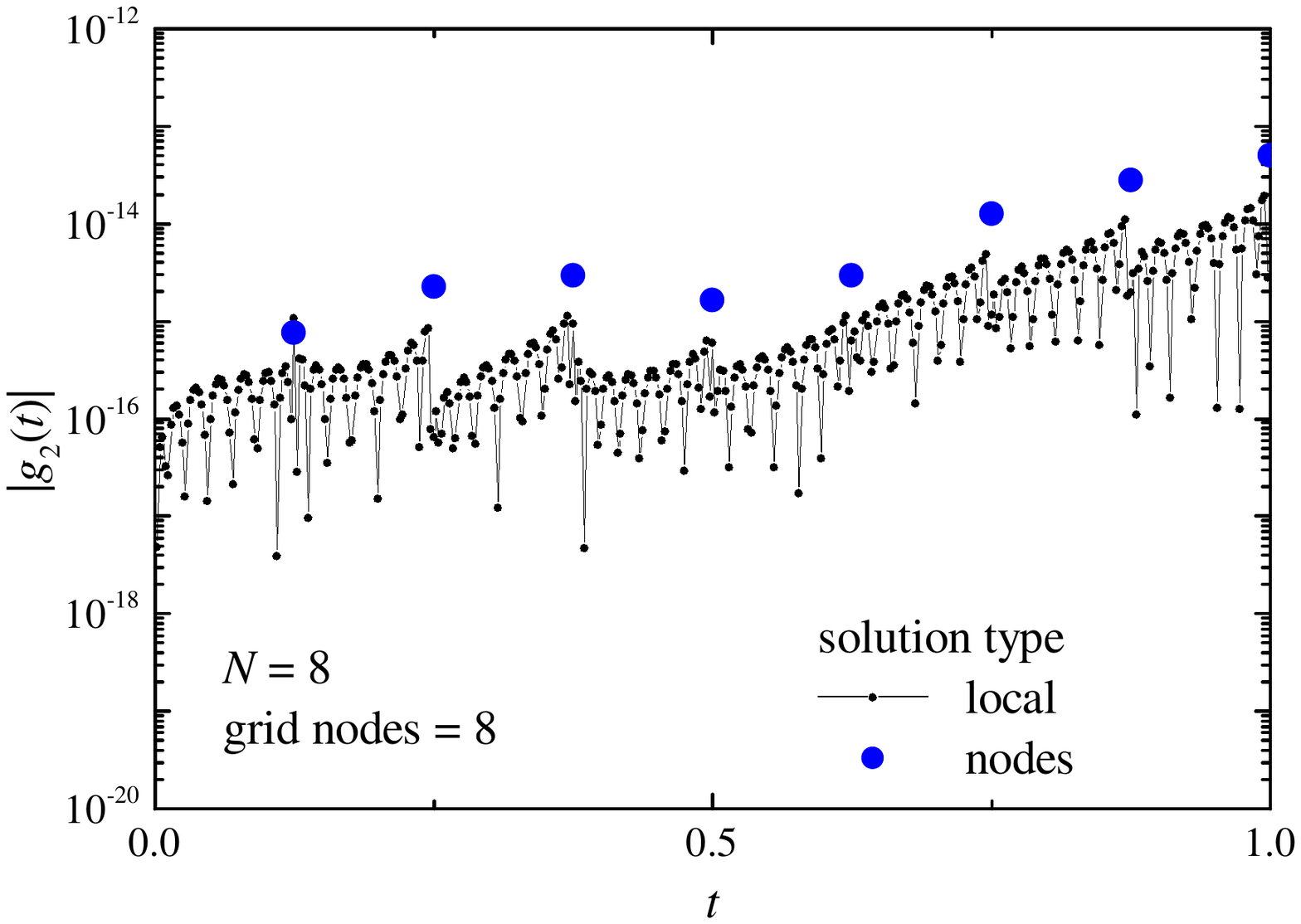}
\vspace{-8mm}\caption{\label{fig:hess_2_ind2_sol_g_eps:b2}}
\end{subfigure}
\begin{subfigure}{0.320\textwidth}
\includegraphics[width=\textwidth]{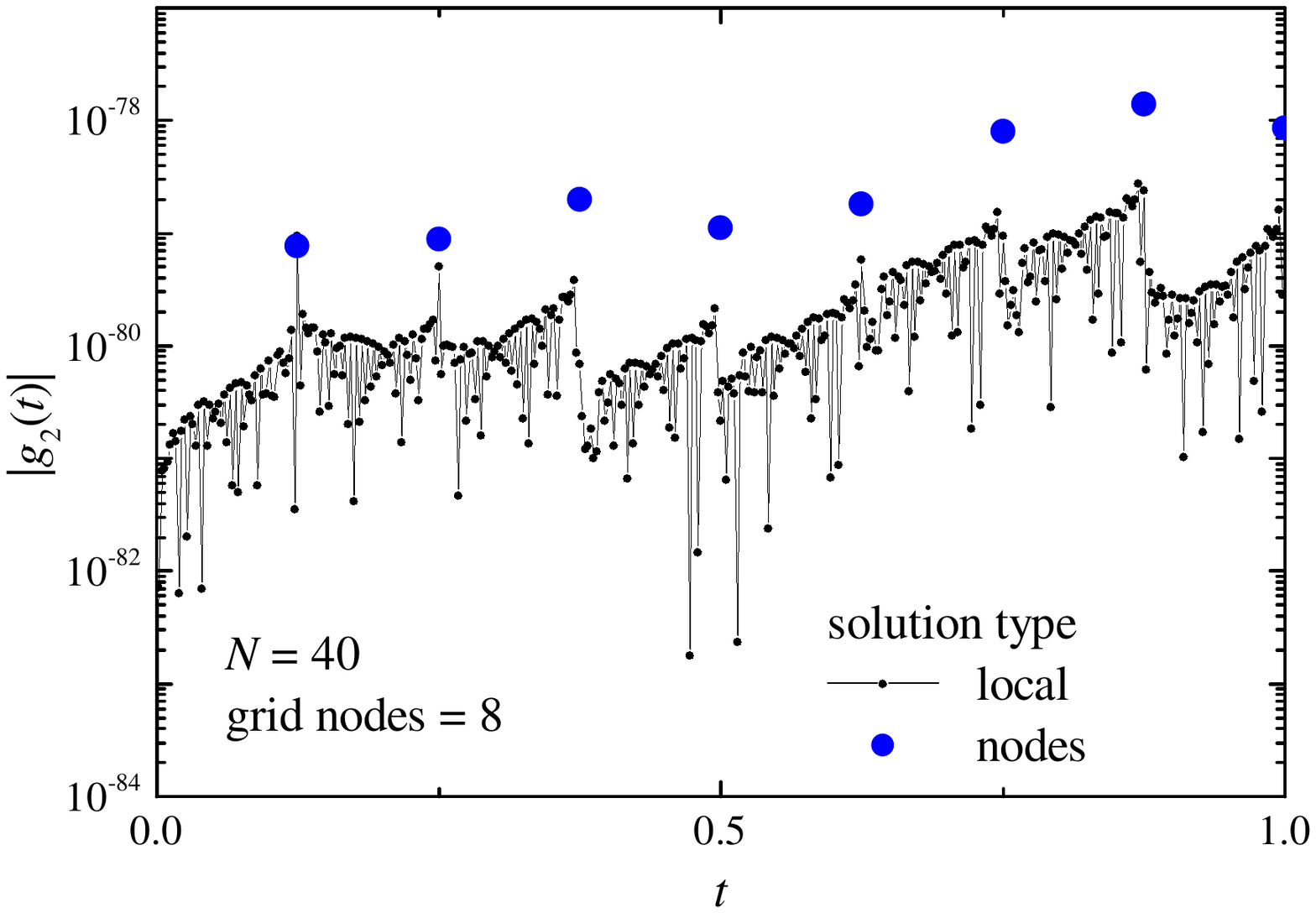}
\vspace{-8mm}\caption{\label{fig:hess_2_ind2_sol_g_eps:b3}}
\end{subfigure}\\
\begin{subfigure}{0.320\textwidth}
\includegraphics[width=\textwidth]{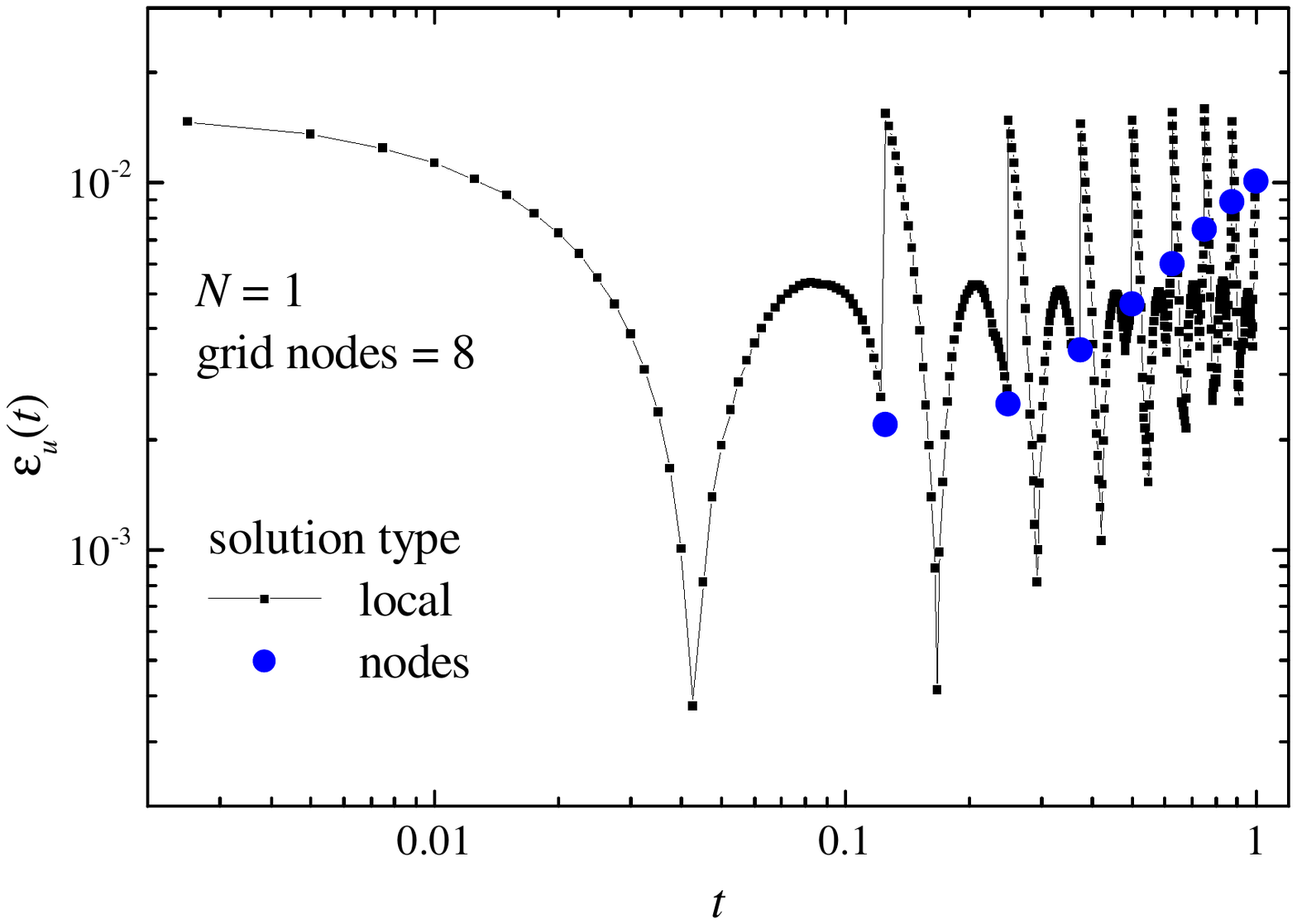}
\vspace{-8mm}\caption{\label{fig:hess_2_ind2_sol_g_eps:c1}}
\end{subfigure}
\begin{subfigure}{0.320\textwidth}
\includegraphics[width=\textwidth]{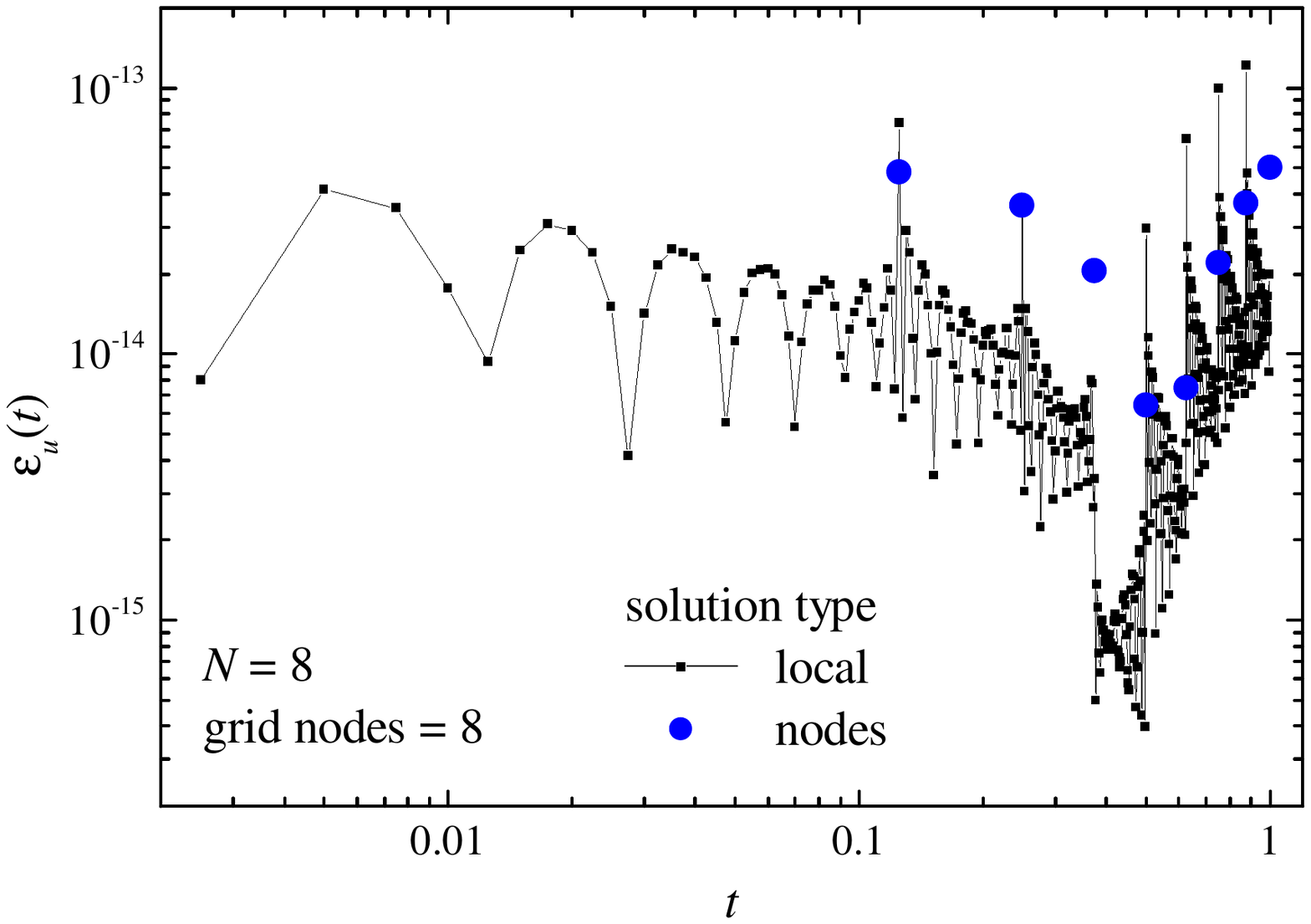}
\vspace{-8mm}\caption{\label{fig:hess_2_ind2_sol_g_eps:c2}}
\end{subfigure}
\begin{subfigure}{0.320\textwidth}
\includegraphics[width=\textwidth]{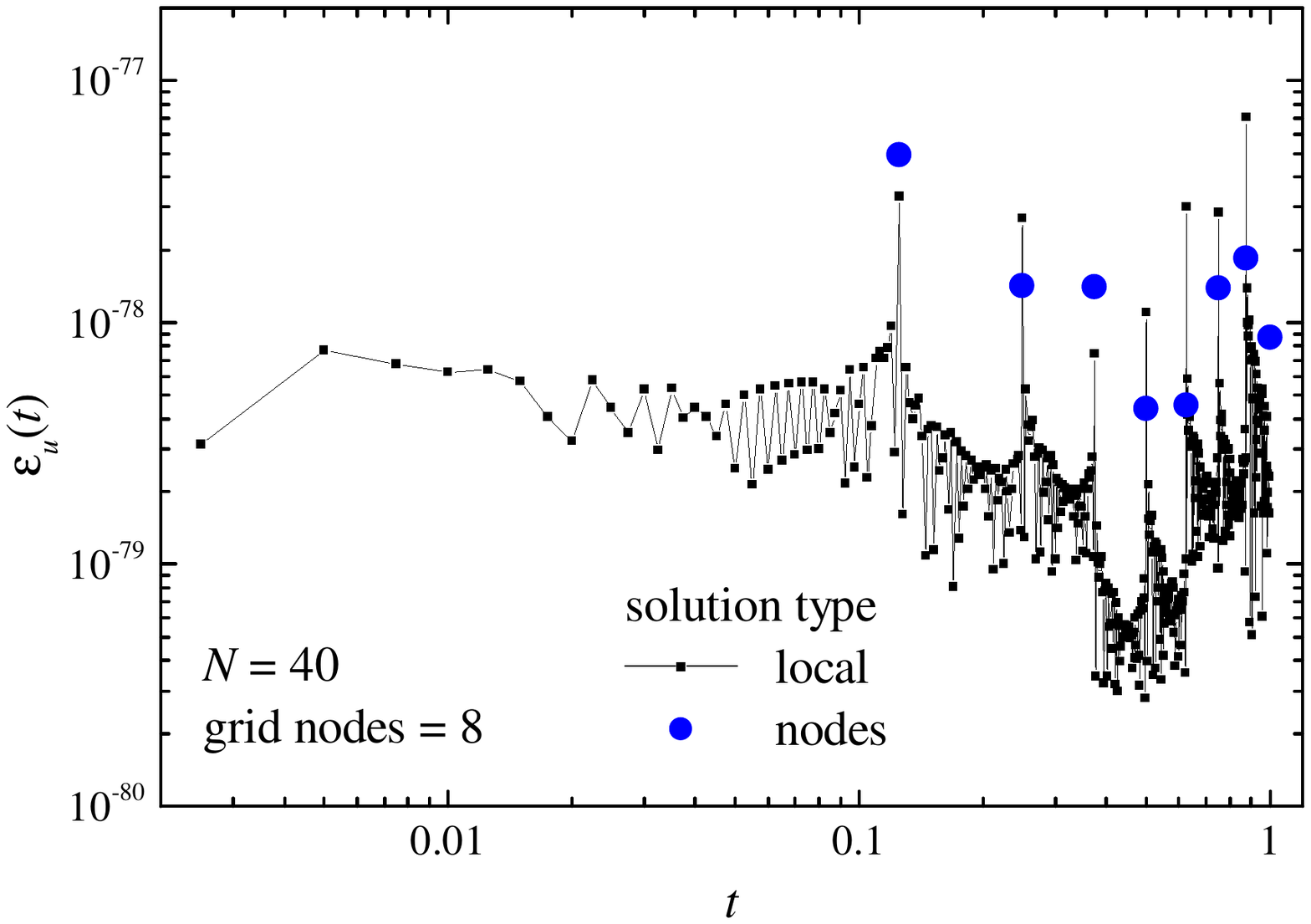}
\vspace{-8mm}\caption{\label{fig:hess_2_ind2_sol_g_eps:c3}}
\end{subfigure}\\
\begin{subfigure}{0.320\textwidth}
\includegraphics[width=\textwidth]{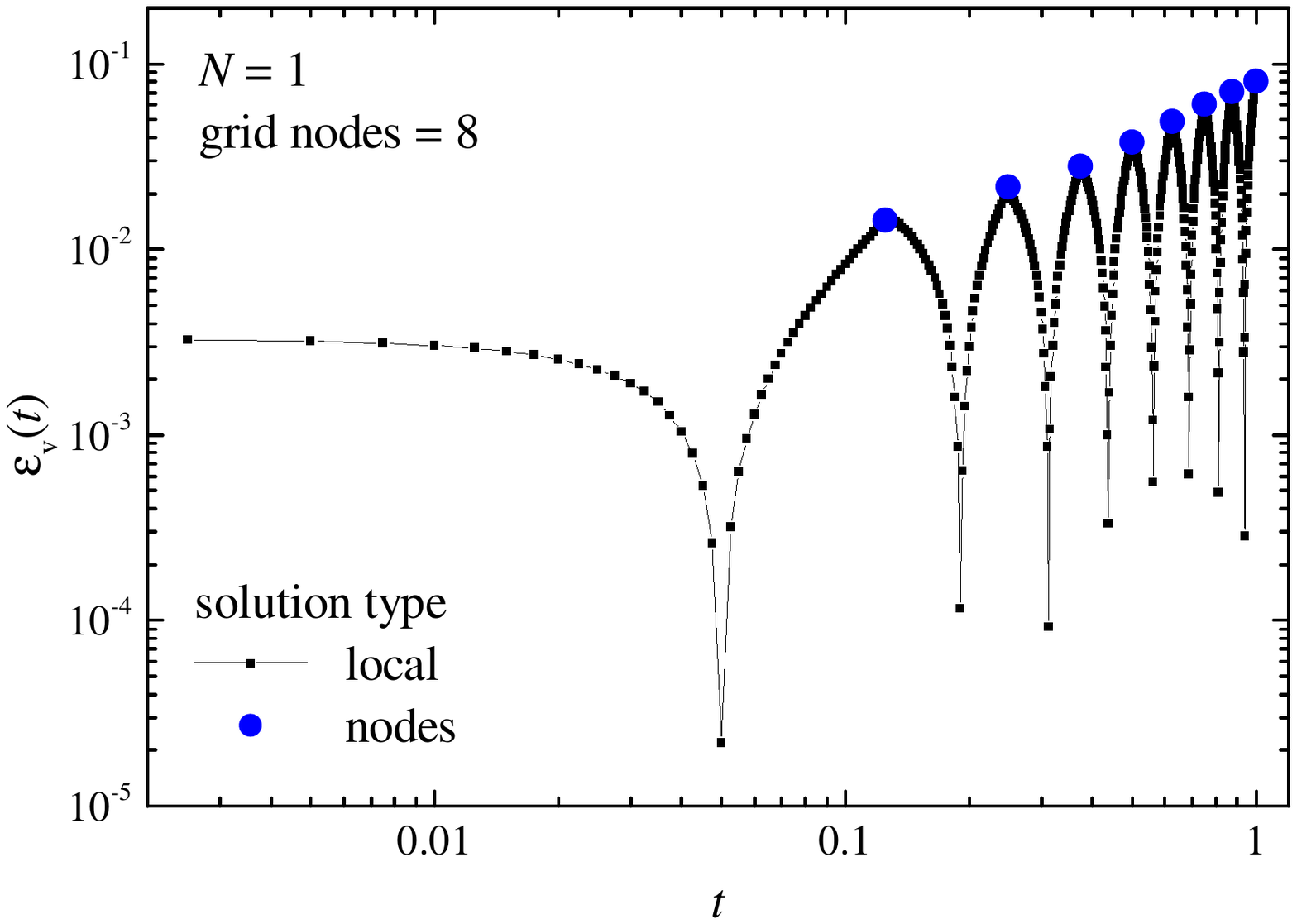}
\vspace{-8mm}\caption{\label{fig:hess_2_ind2_sol_g_eps:d1}}
\end{subfigure}
\begin{subfigure}{0.320\textwidth}
\includegraphics[width=\textwidth]{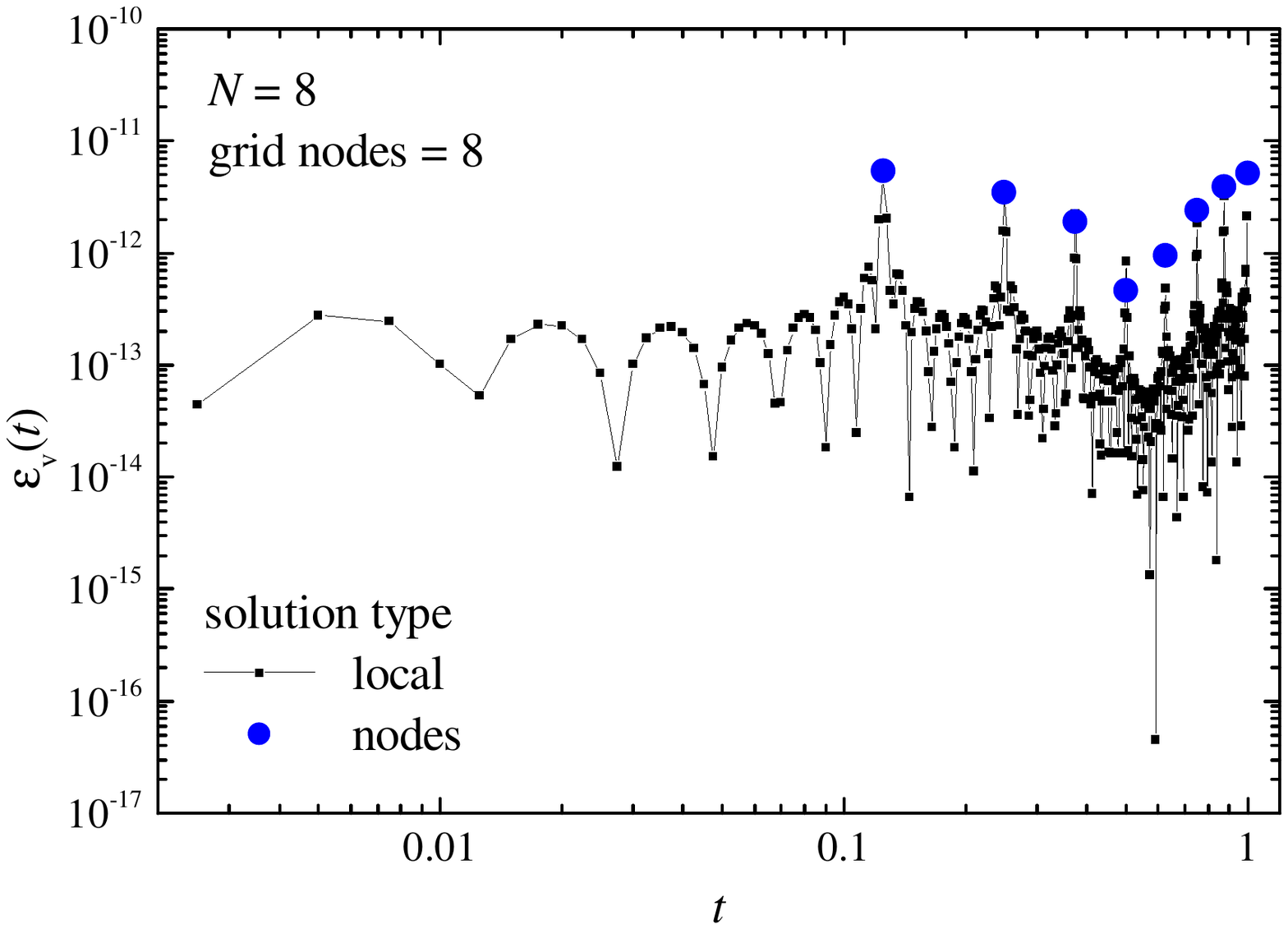}
\vspace{-8mm}\caption{\label{fig:hess_2_ind2_sol_g_eps:d2}}
\end{subfigure}
\begin{subfigure}{0.320\textwidth}
\includegraphics[width=\textwidth]{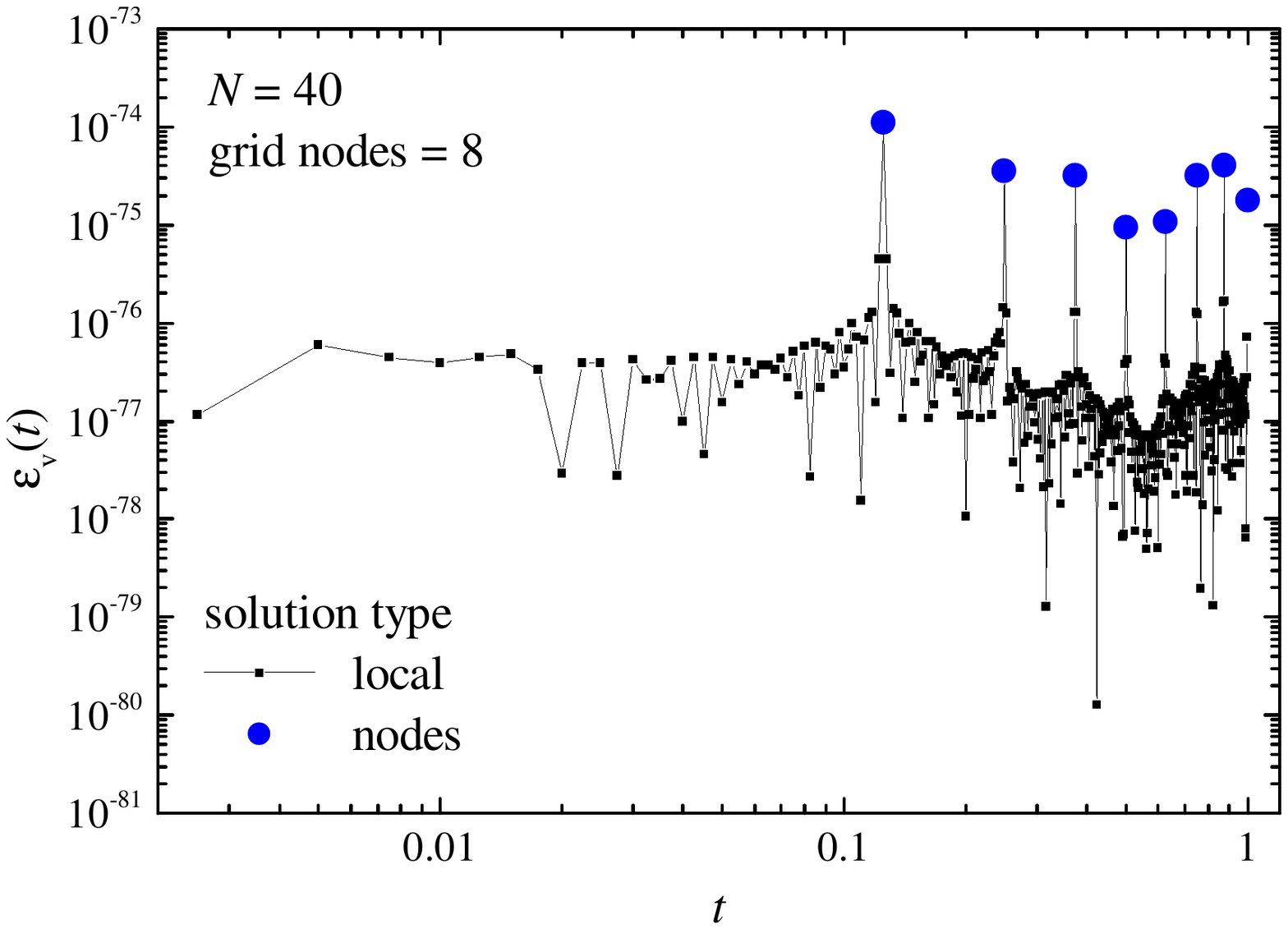}
\vspace{-8mm}\caption{\label{fig:hess_2_ind2_sol_g_eps:d3}}
\end{subfigure}\\
\begin{subfigure}{0.320\textwidth}
\includegraphics[width=\textwidth]{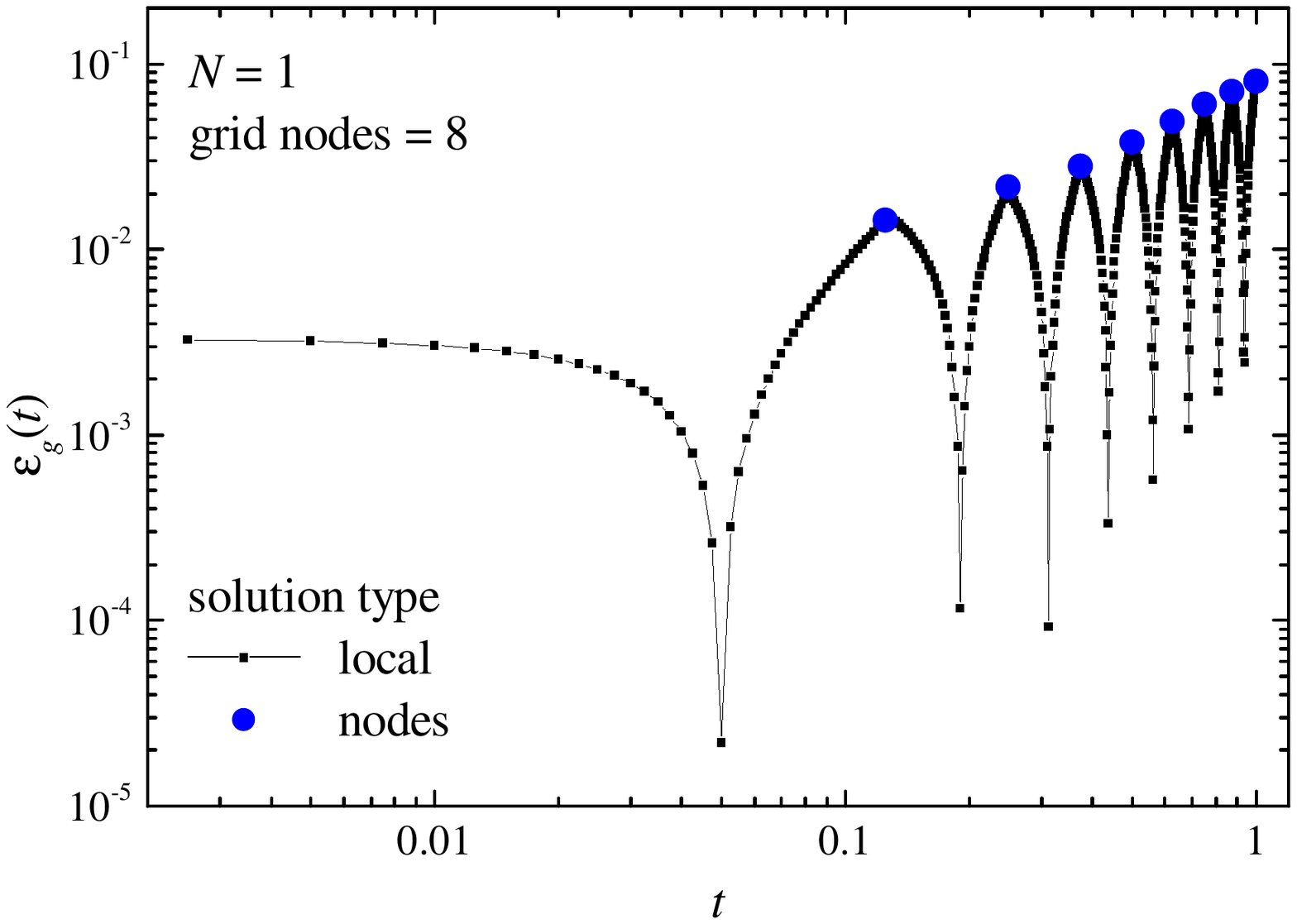}
\vspace{-8mm}\caption{\label{fig:hess_2_ind2_sol_g_eps:e1}}
\end{subfigure}
\begin{subfigure}{0.320\textwidth}
\includegraphics[width=\textwidth]{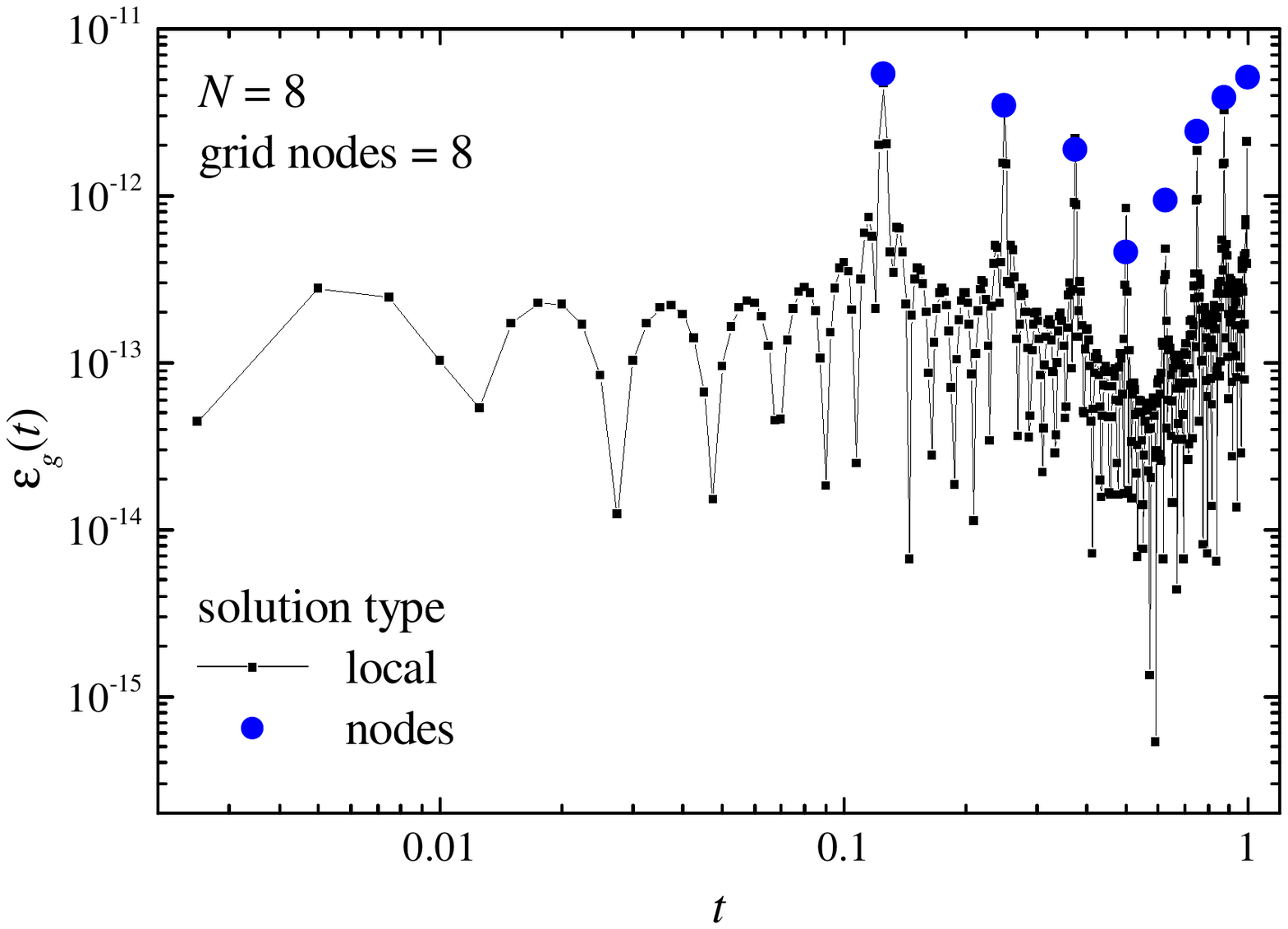}
\vspace{-8mm}\caption{\label{fig:hess_2_ind2_sol_g_eps:e2}}
\end{subfigure}
\begin{subfigure}{0.320\textwidth}
\includegraphics[width=\textwidth]{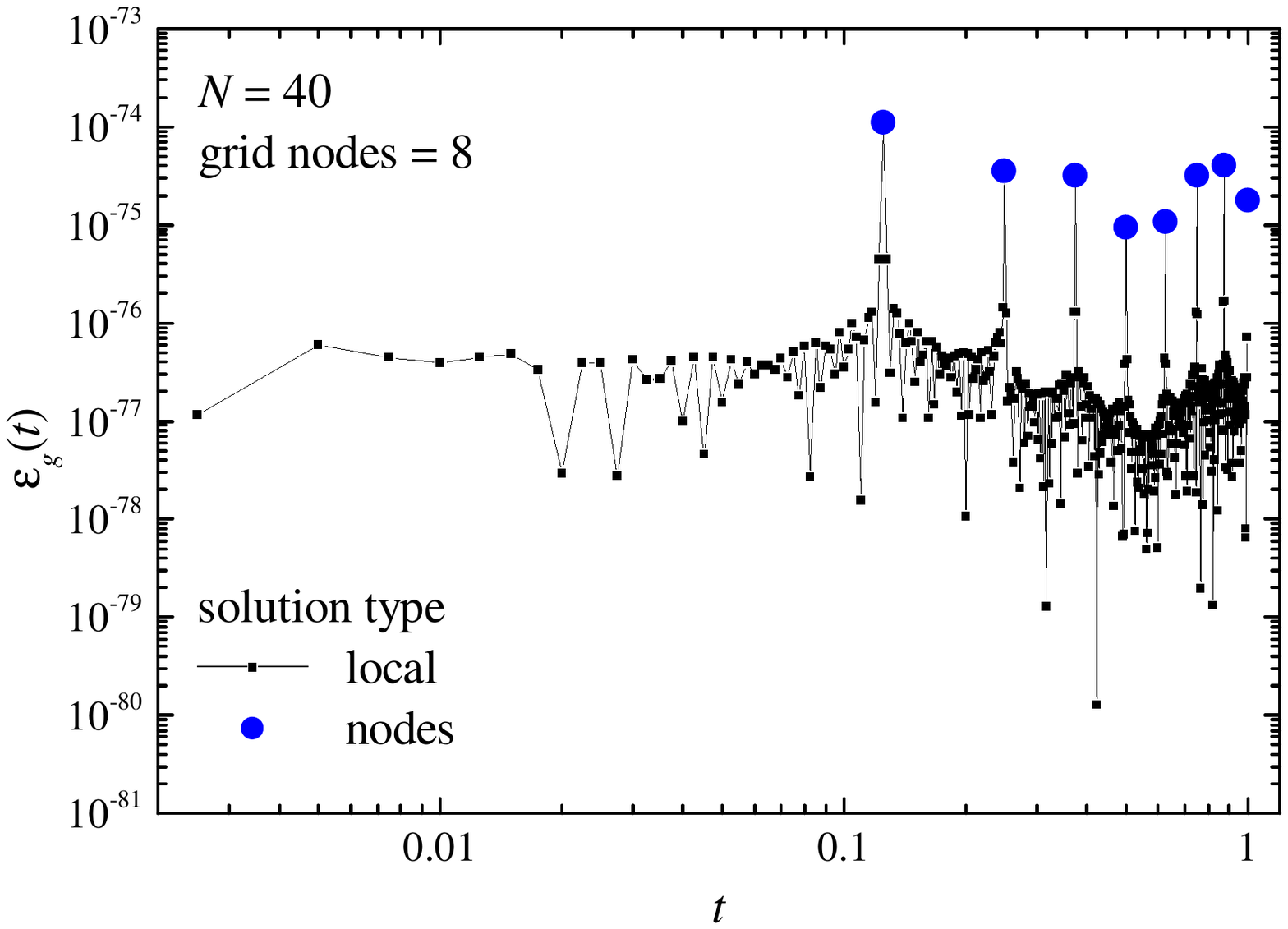}
\vspace{-8mm}\caption{\label{fig:hess_2_ind2_sol_g_eps:e3}}
\end{subfigure}\\
\caption{%
Numerical solution of the DAE system (\ref{eq:hess_dae_ind_2}) of index 2. Comparison of quantitative satisfiability of the conditions $g_{1} = 0$ (\subref{fig:hess_2_ind2_sol_g_eps:a1}, \subref{fig:hess_2_ind2_sol_g_eps:a2}, \subref{fig:hess_2_ind2_sol_g_eps:a3}) and $g_{2} = 0$ (\subref{fig:hess_2_ind2_sol_g_eps:b1}, \subref{fig:hess_2_ind2_sol_g_eps:b2}, \subref{fig:hess_2_ind2_sol_g_eps:b3}), the errors $\varepsilon_{u}(t)$ (\subref{fig:hess_2_ind2_sol_g_eps:c1}, \subref{fig:hess_2_ind2_sol_g_eps:c2}, \subref{fig:hess_2_ind2_sol_g_eps:c3}), $\varepsilon_{v}(t)$ (\subref{fig:hess_2_ind2_sol_g_eps:d1}, \subref{fig:hess_2_ind2_sol_g_eps:d2}, \subref{fig:hess_2_ind2_sol_g_eps:d3}), $\varepsilon_{g}(t)$ (\subref{fig:hess_2_ind2_sol_g_eps:e1}, \subref{fig:hess_2_ind2_sol_g_eps:e2}, \subref{fig:hess_2_ind2_sol_g_eps:e3}), for numerical solution obtained using polynomials with degrees $N = 1$ (\subref{fig:hess_2_ind2_sol_g_eps:a1}, \subref{fig:hess_2_ind2_sol_g_eps:b1}, \subref{fig:hess_2_ind2_sol_g_eps:c1}, \subref{fig:hess_2_ind2_sol_g_eps:d1}, \subref{fig:hess_2_ind2_sol_g_eps:e1}), $N = 8$ (\subref{fig:hess_2_ind2_sol_g_eps:a2}, \subref{fig:hess_2_ind2_sol_g_eps:b2}, \subref{fig:hess_2_ind2_sol_g_eps:c2}, \subref{fig:hess_2_ind2_sol_g_eps:d2}, \subref{fig:hess_2_ind2_sol_g_eps:e2}), $N = 40$ (\subref{fig:hess_2_ind2_sol_g_eps:a3}, \subref{fig:hess_2_ind2_sol_g_eps:b3}, \subref{fig:hess_2_ind2_sol_g_eps:c3}, \subref{fig:hess_2_ind2_sol_g_eps:d3}, \subref{fig:hess_2_ind2_sol_g_eps:e3}).
}
\label{fig:hess_2_ind2_sol_g_eps}
\end{figure} 

\begin{figure}[h!]
\captionsetup[subfigure]{%
	position=bottom,
	font+=smaller,
	textfont=normalfont,
	singlelinecheck=off,
	justification=raggedright
}
\centering
\begin{subfigure}{0.275\textwidth}
\includegraphics[width=\textwidth]{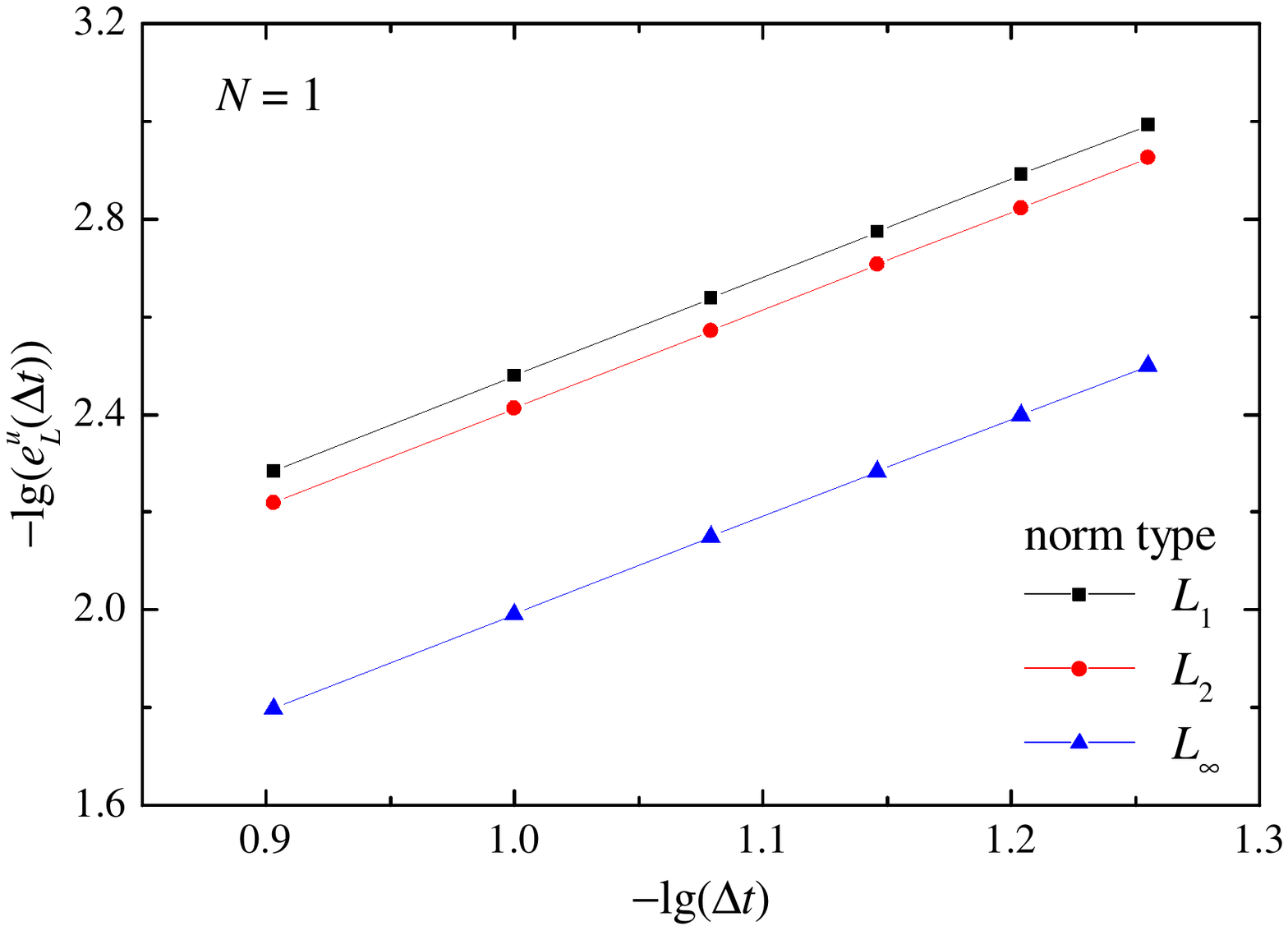}
\vspace{-8mm}\caption{\label{fig:hess_2_ind2_errors:a1}}
\end{subfigure}\hspace{6mm}
\begin{subfigure}{0.275\textwidth}
\includegraphics[width=\textwidth]{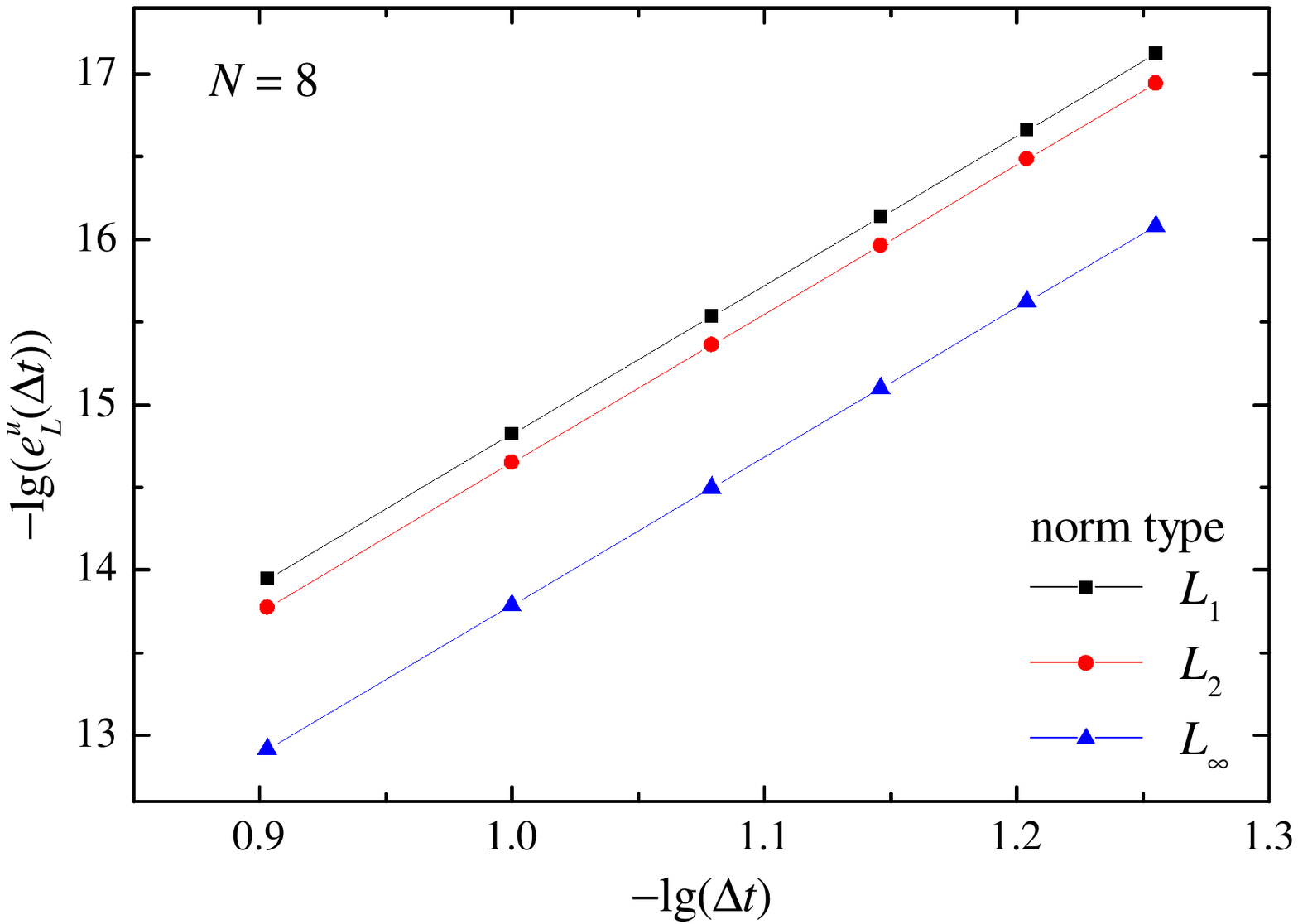}
\vspace{-8mm}\caption{\label{fig:hess_2_ind2_errors:a2}}
\end{subfigure}\hspace{6mm}
\begin{subfigure}{0.275\textwidth}
\includegraphics[width=\textwidth]{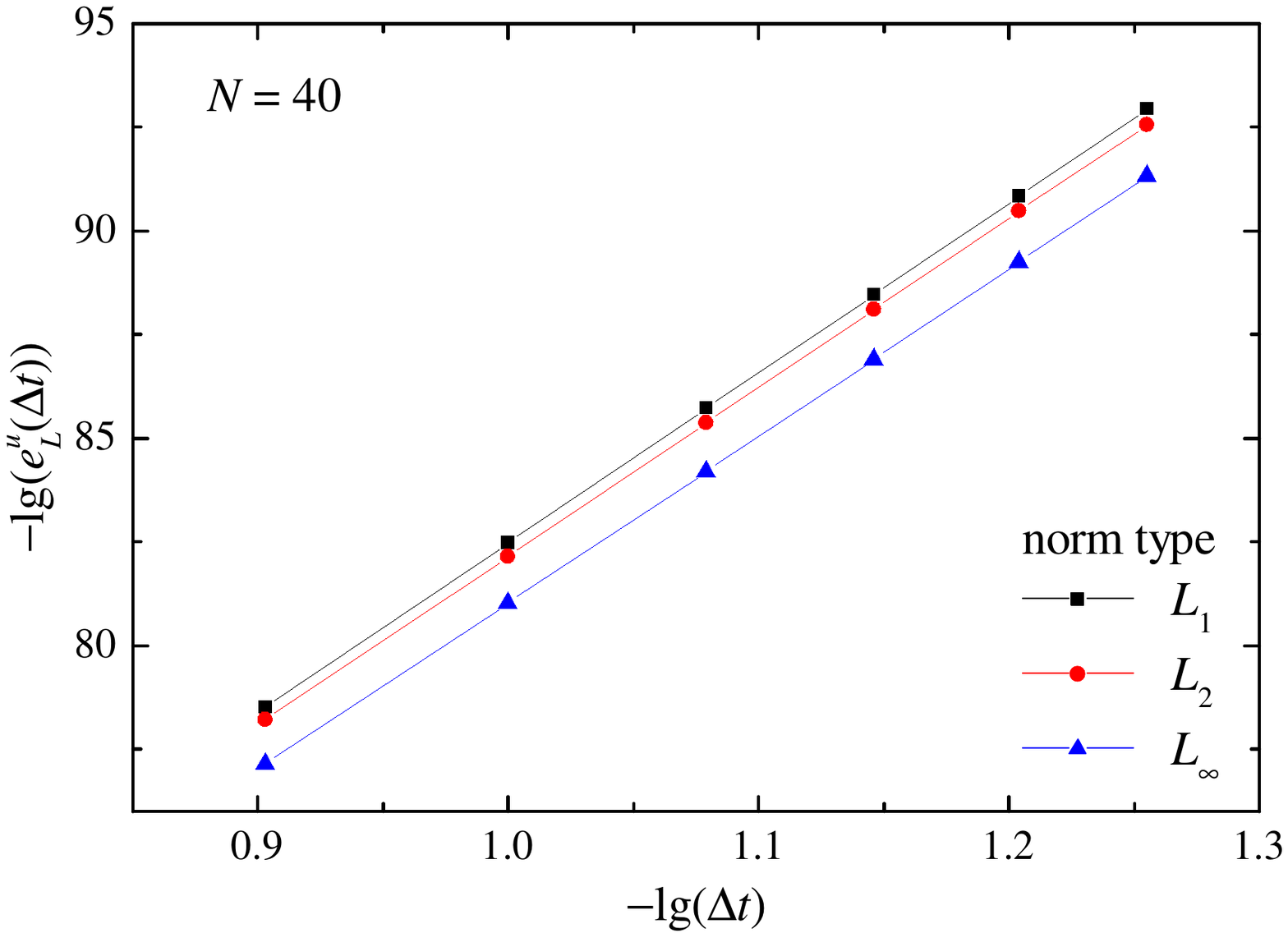}
\vspace{-8mm}\caption{\label{fig:hess_2_ind2_errors:a3}}
\end{subfigure}\\[-2mm]
\begin{subfigure}{0.275\textwidth}
\includegraphics[width=\textwidth]{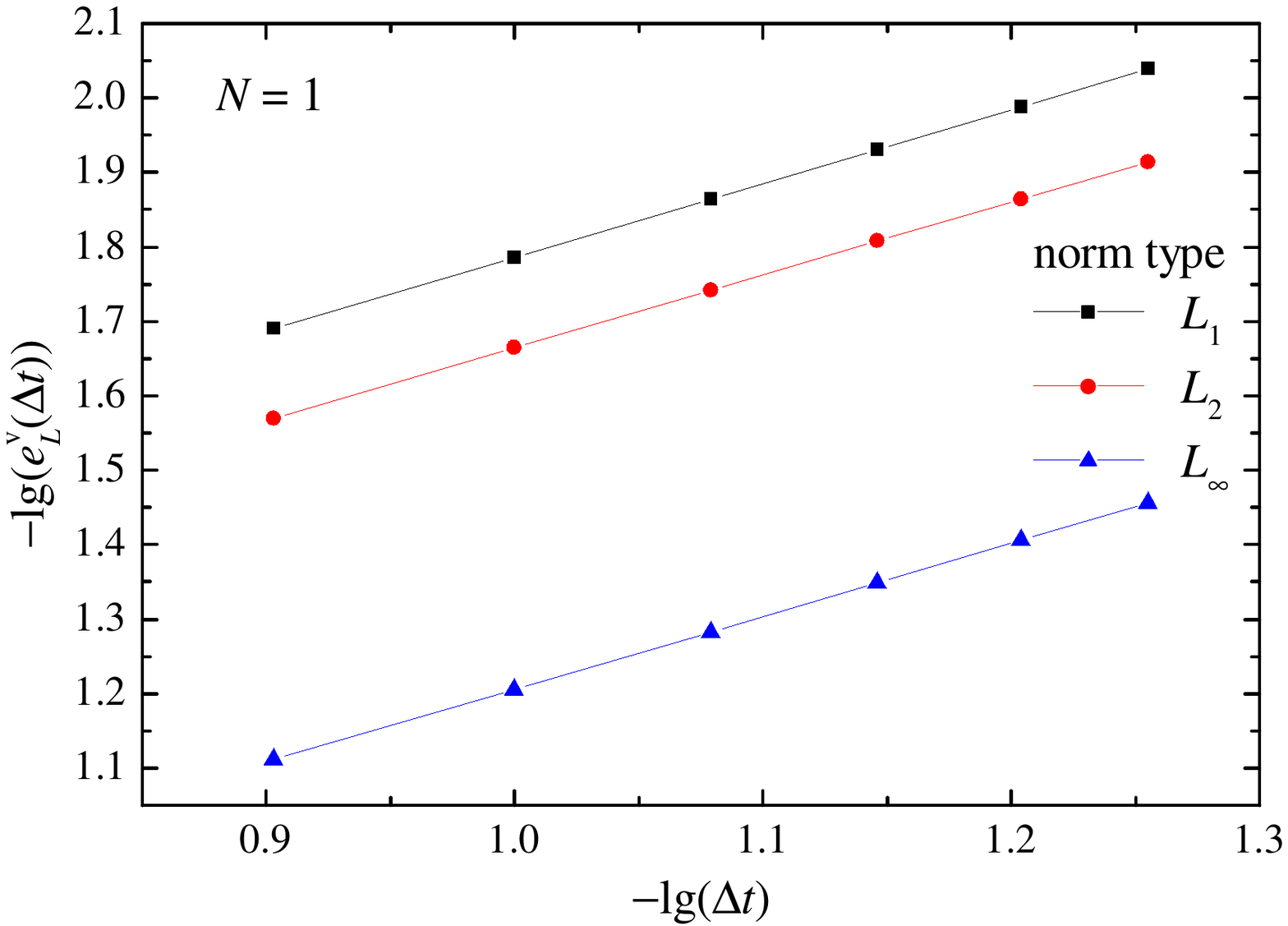}
\vspace{-8mm}\caption{\label{fig:hess_2_ind2_errors:b1}}
\end{subfigure}\hspace{6mm}
\begin{subfigure}{0.275\textwidth}
\includegraphics[width=\textwidth]{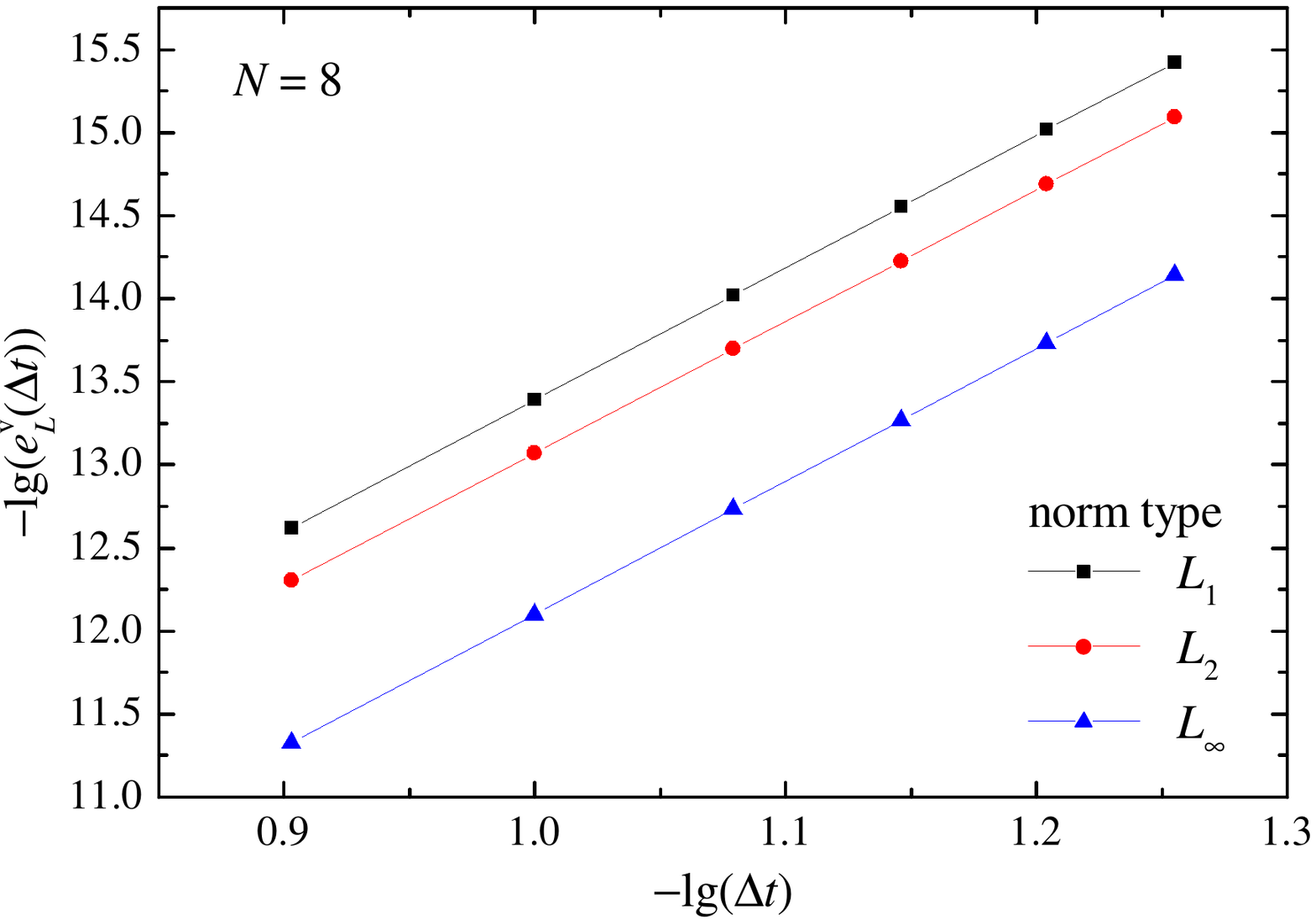}
\vspace{-8mm}\caption{\label{fig:hess_2_ind2_errors:b2}}
\end{subfigure}\hspace{6mm}
\begin{subfigure}{0.275\textwidth}
\includegraphics[width=\textwidth]{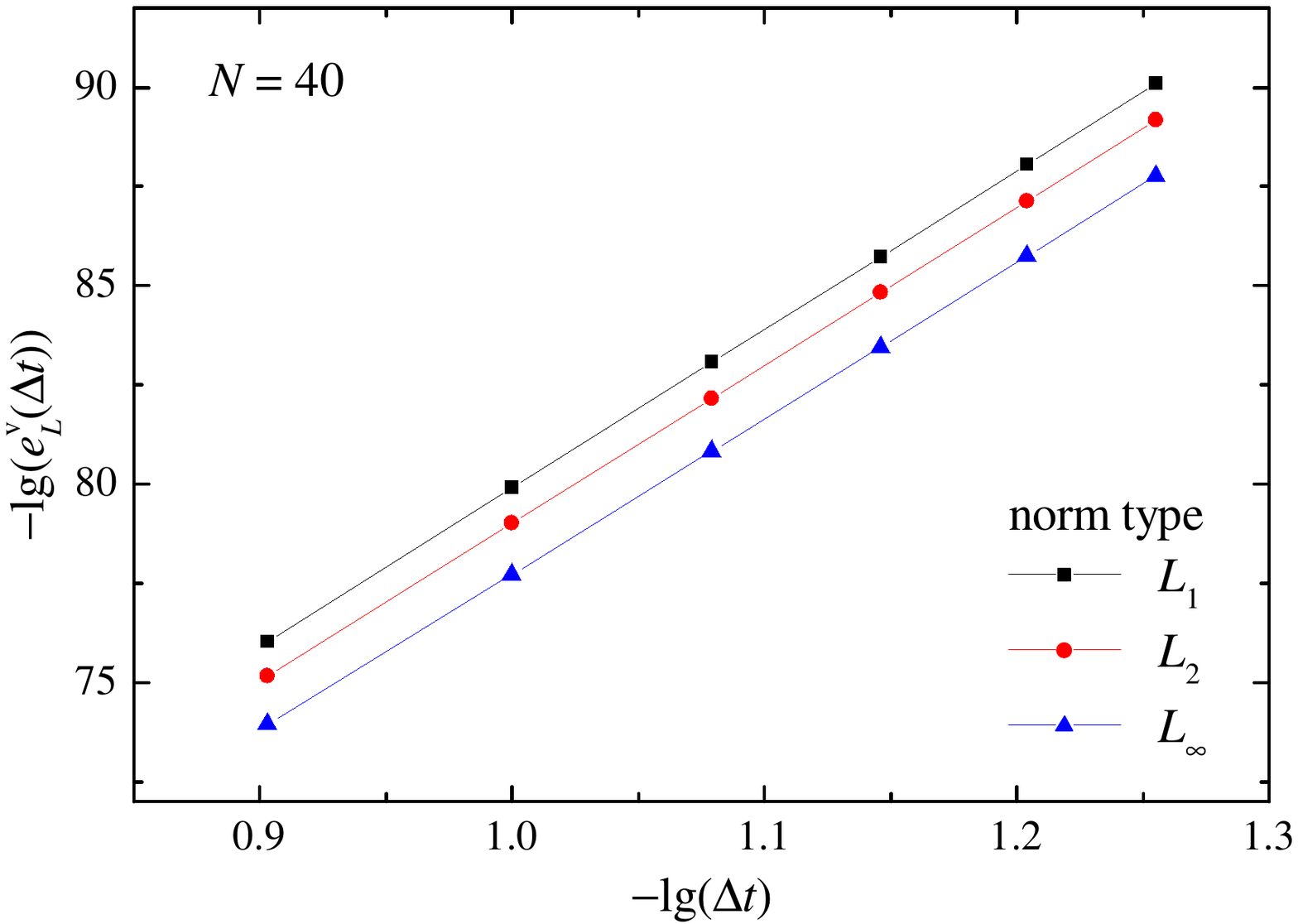}
\vspace{-8mm}\caption{\label{fig:hess_2_ind2_errors:b3}}
\end{subfigure}\\[-2mm]
\begin{subfigure}{0.275\textwidth}
\includegraphics[width=\textwidth]{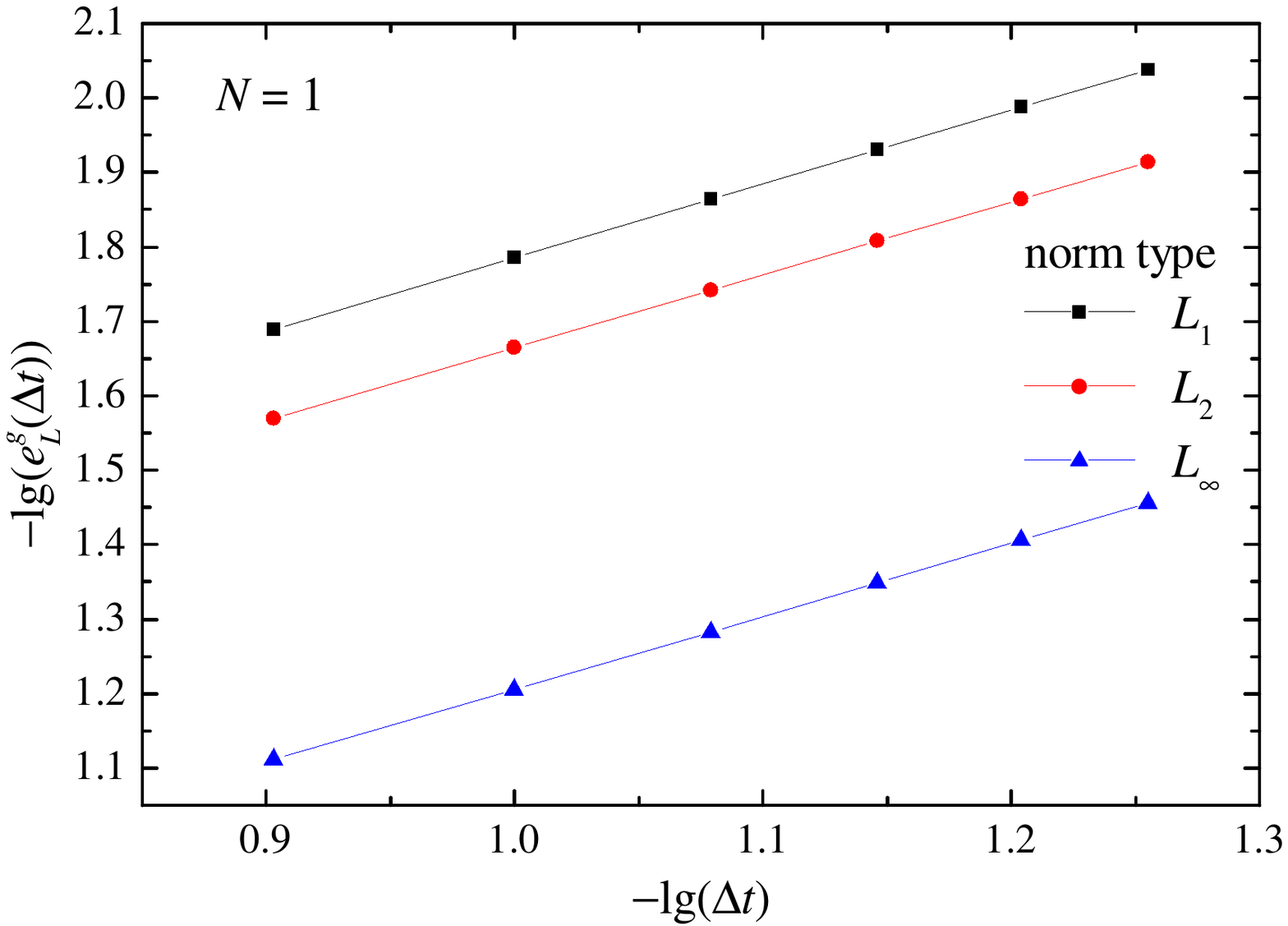}
\vspace{-8mm}\caption{\label{fig:hess_2_ind2_errors:c1}}
\end{subfigure}\hspace{6mm}
\begin{subfigure}{0.275\textwidth}
\includegraphics[width=\textwidth]{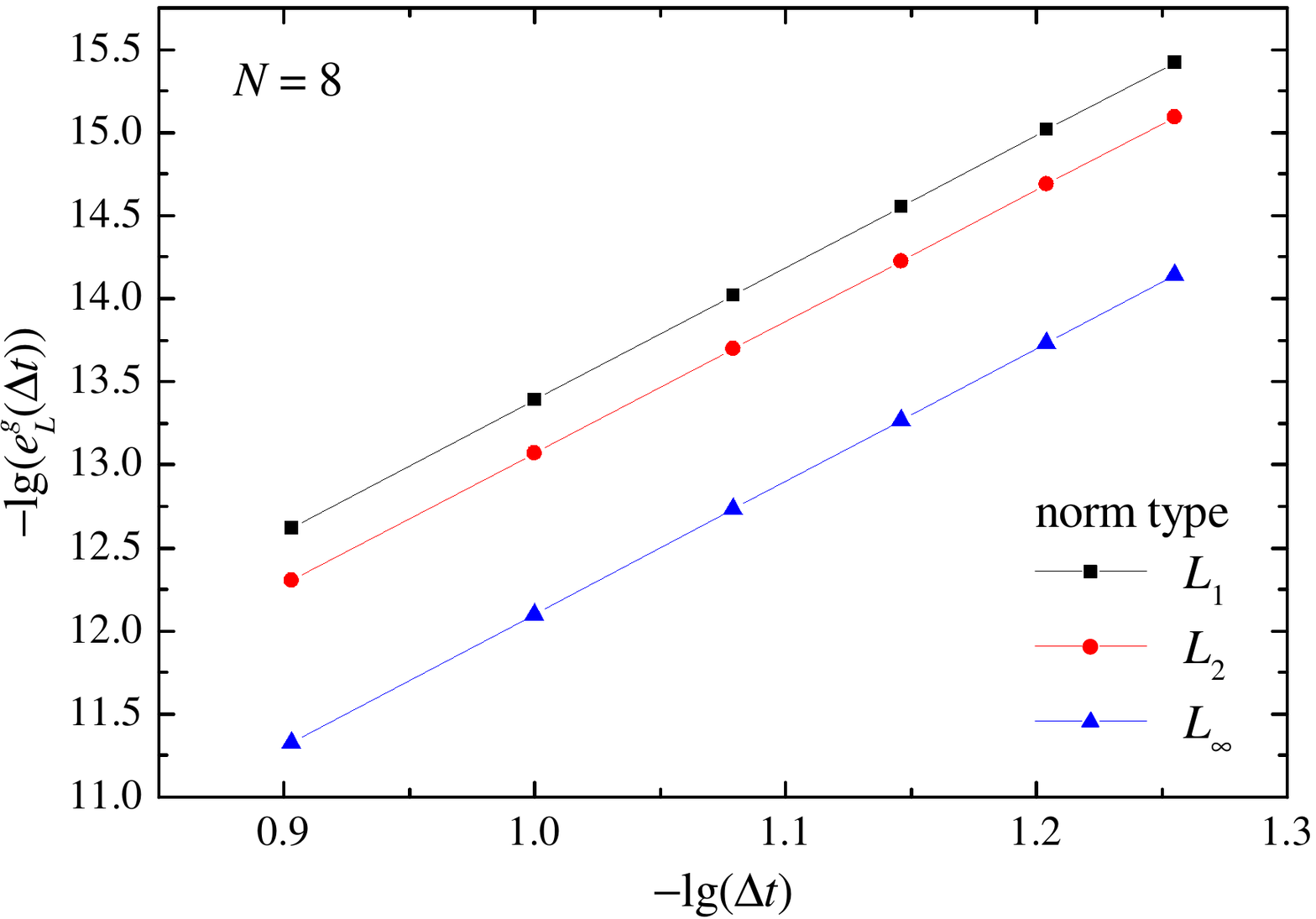}
\vspace{-8mm}\caption{\label{fig:hess_2_ind2_errors:c2}}
\end{subfigure}\hspace{6mm}
\begin{subfigure}{0.275\textwidth}
\includegraphics[width=\textwidth]{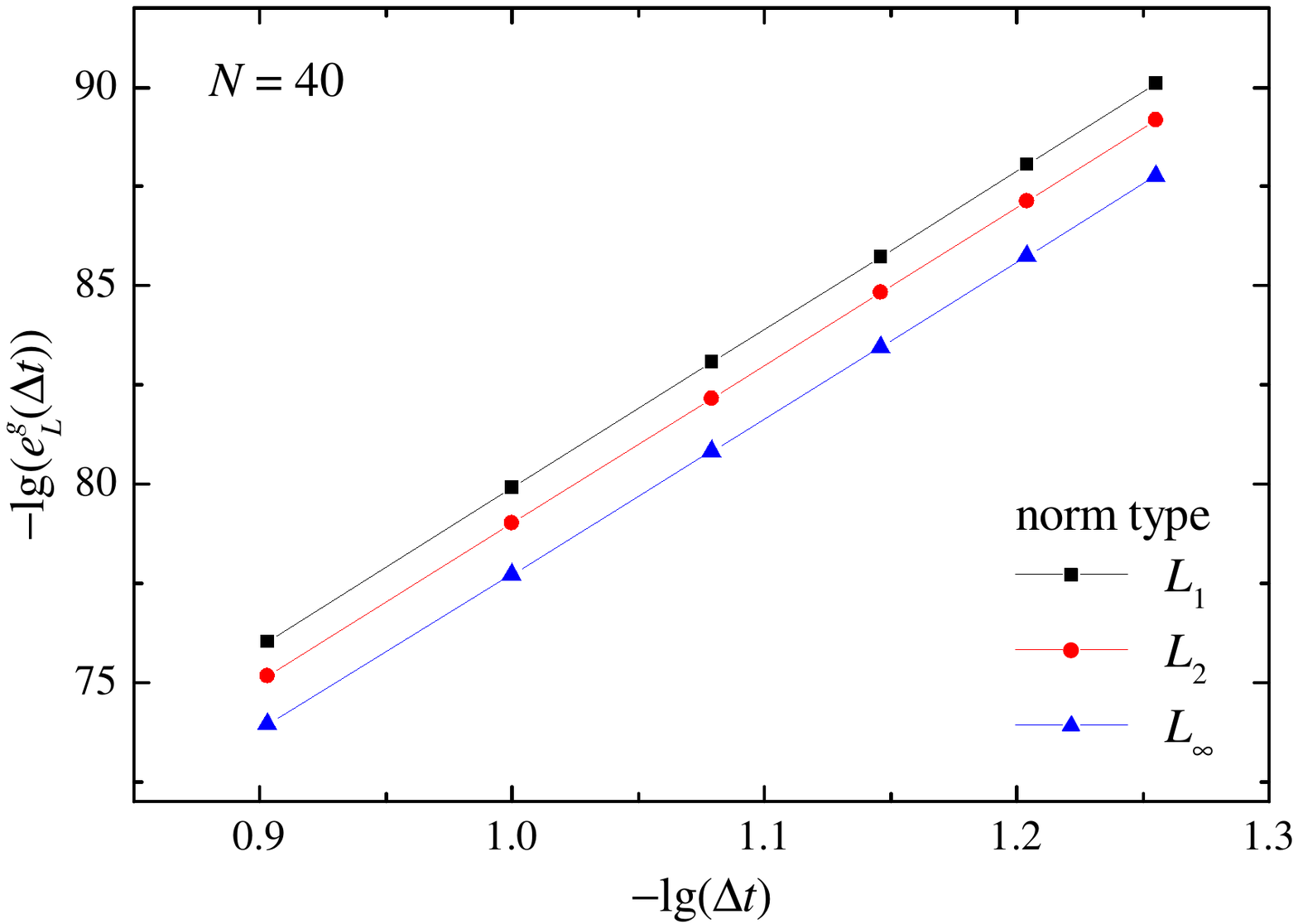}
\vspace{-8mm}\caption{\label{fig:hess_2_ind2_errors:c3}}
\end{subfigure}\\[-2mm]
\begin{subfigure}{0.275\textwidth}
\includegraphics[width=\textwidth]{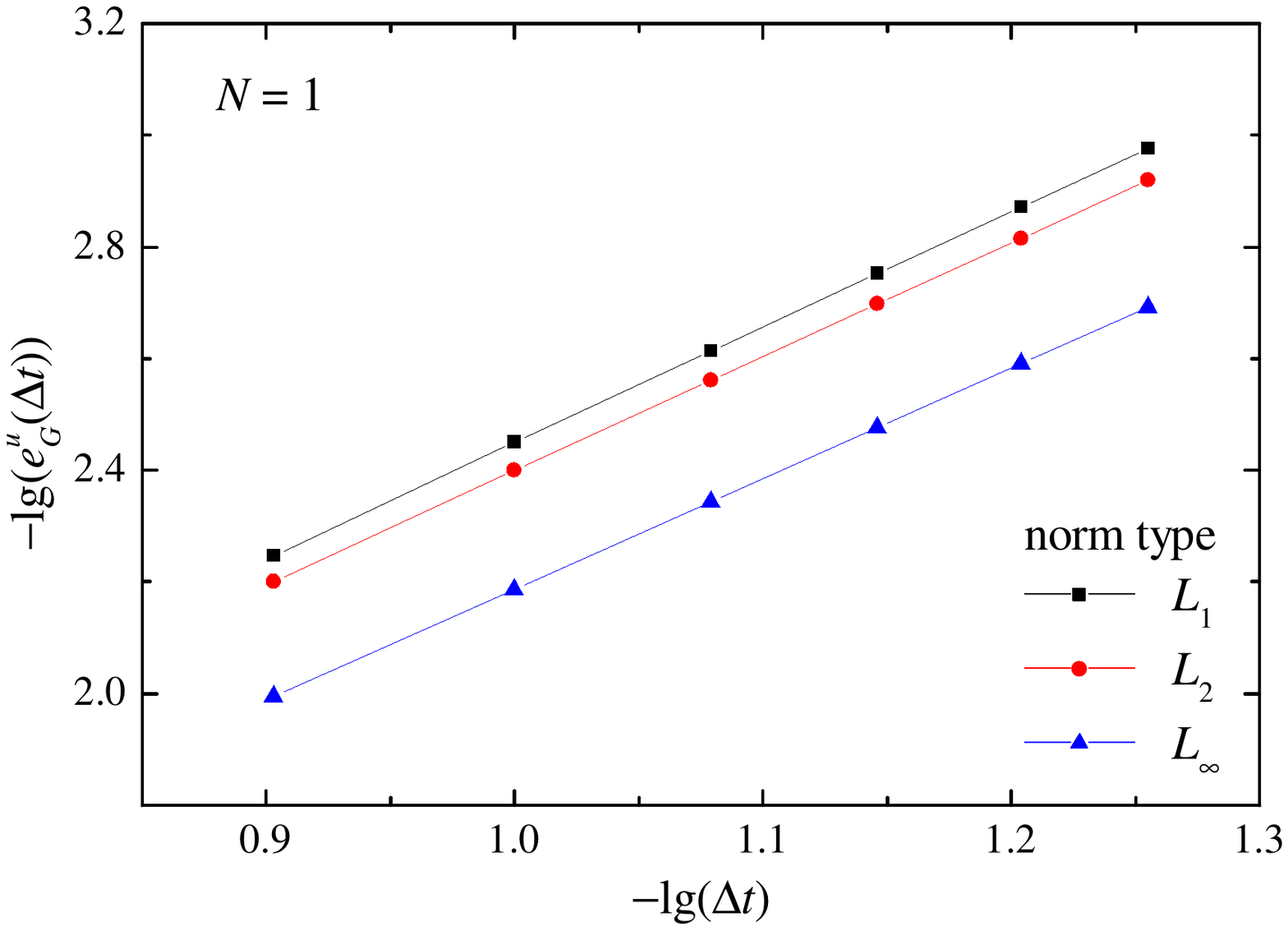}
\vspace{-8mm}\caption{\label{fig:hess_2_ind2_errors:d1}}
\end{subfigure}\hspace{6mm}
\begin{subfigure}{0.275\textwidth}
\includegraphics[width=\textwidth]{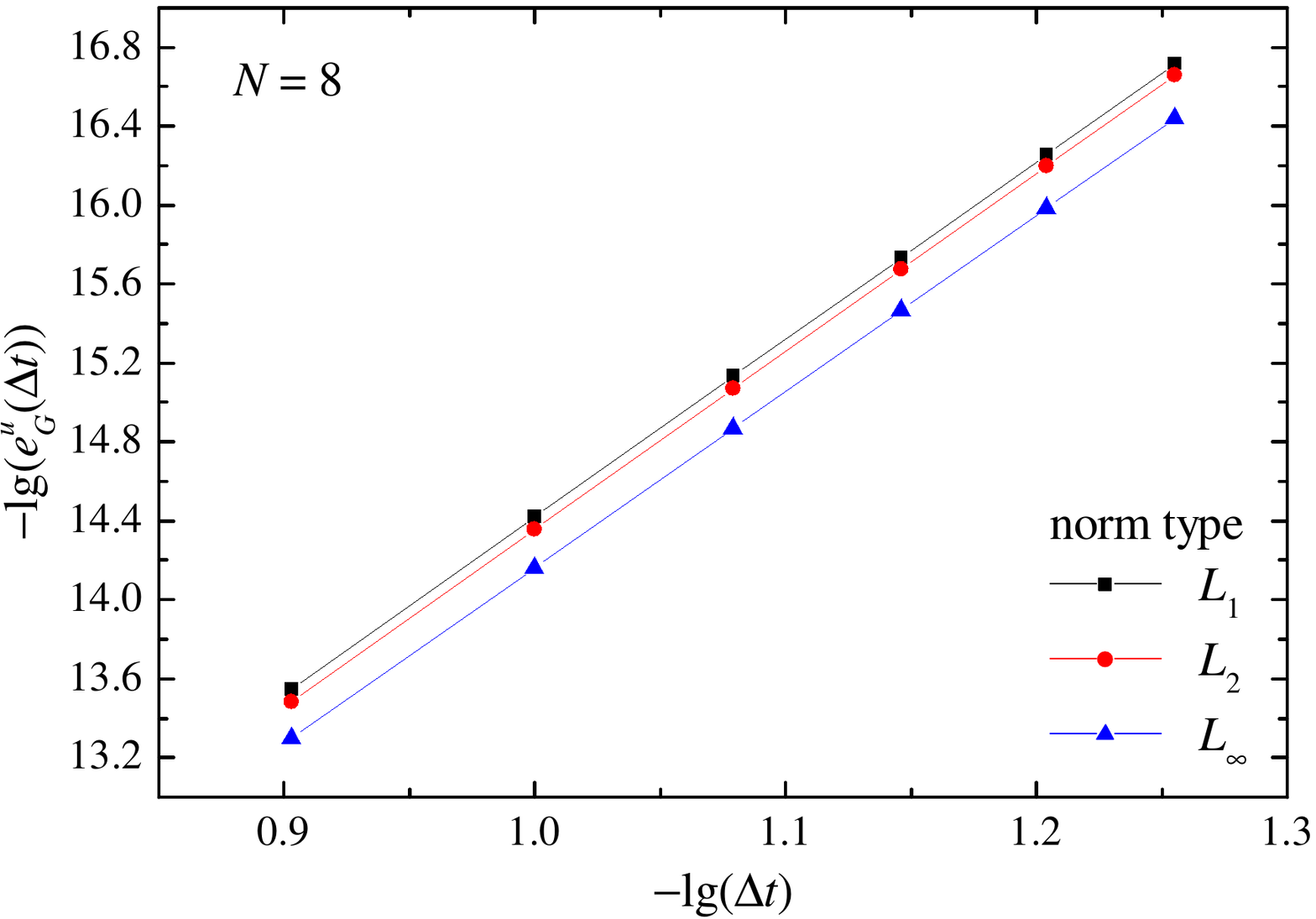}
\vspace{-8mm}\caption{\label{fig:hess_2_ind2_errors:d2}}
\end{subfigure}\hspace{6mm}
\begin{subfigure}{0.275\textwidth}
\includegraphics[width=\textwidth]{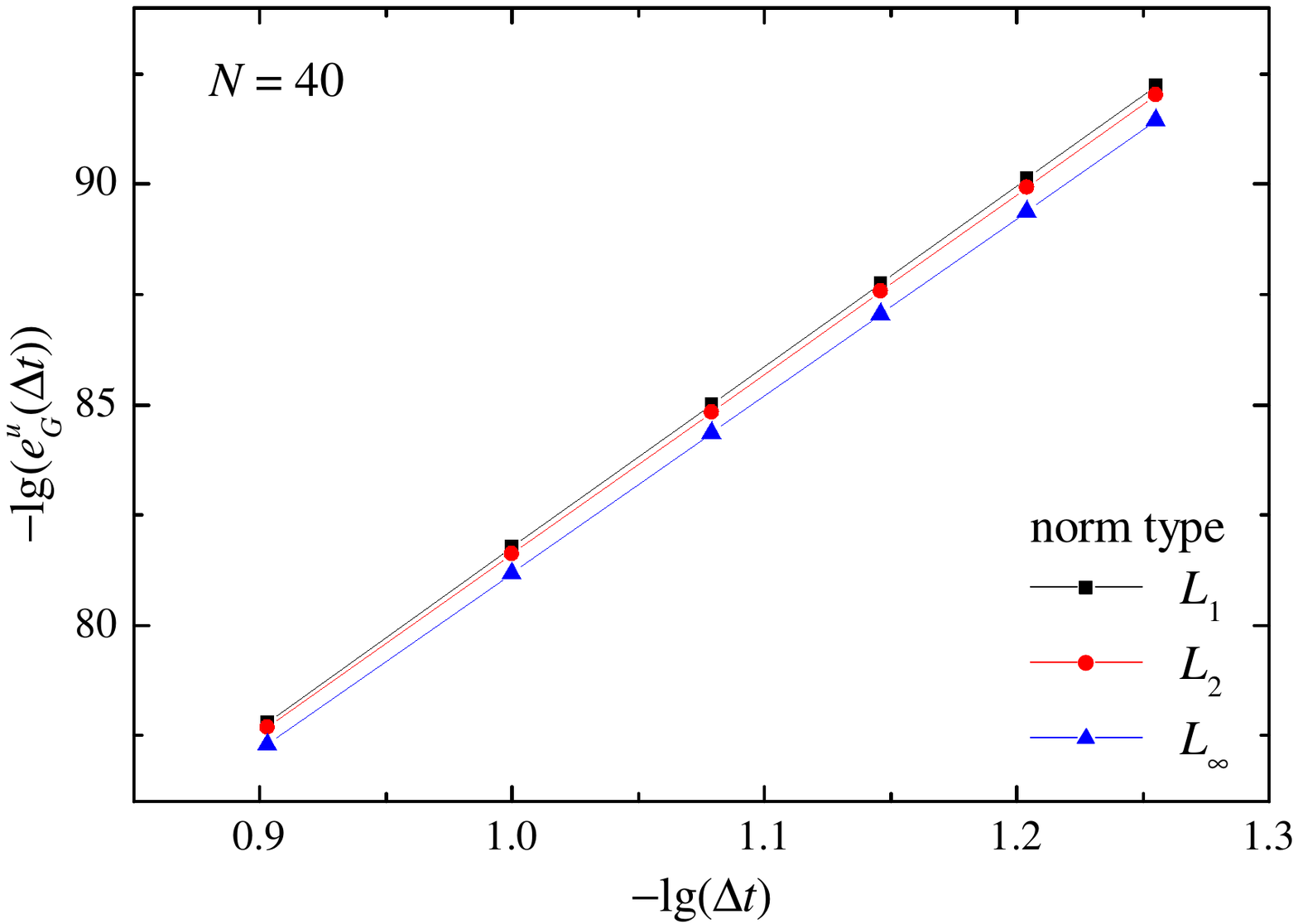}
\vspace{-8mm}\caption{\label{fig:hess_2_ind2_errors:d3}}
\end{subfigure}\\[-2mm]
\begin{subfigure}{0.275\textwidth}
\includegraphics[width=\textwidth]{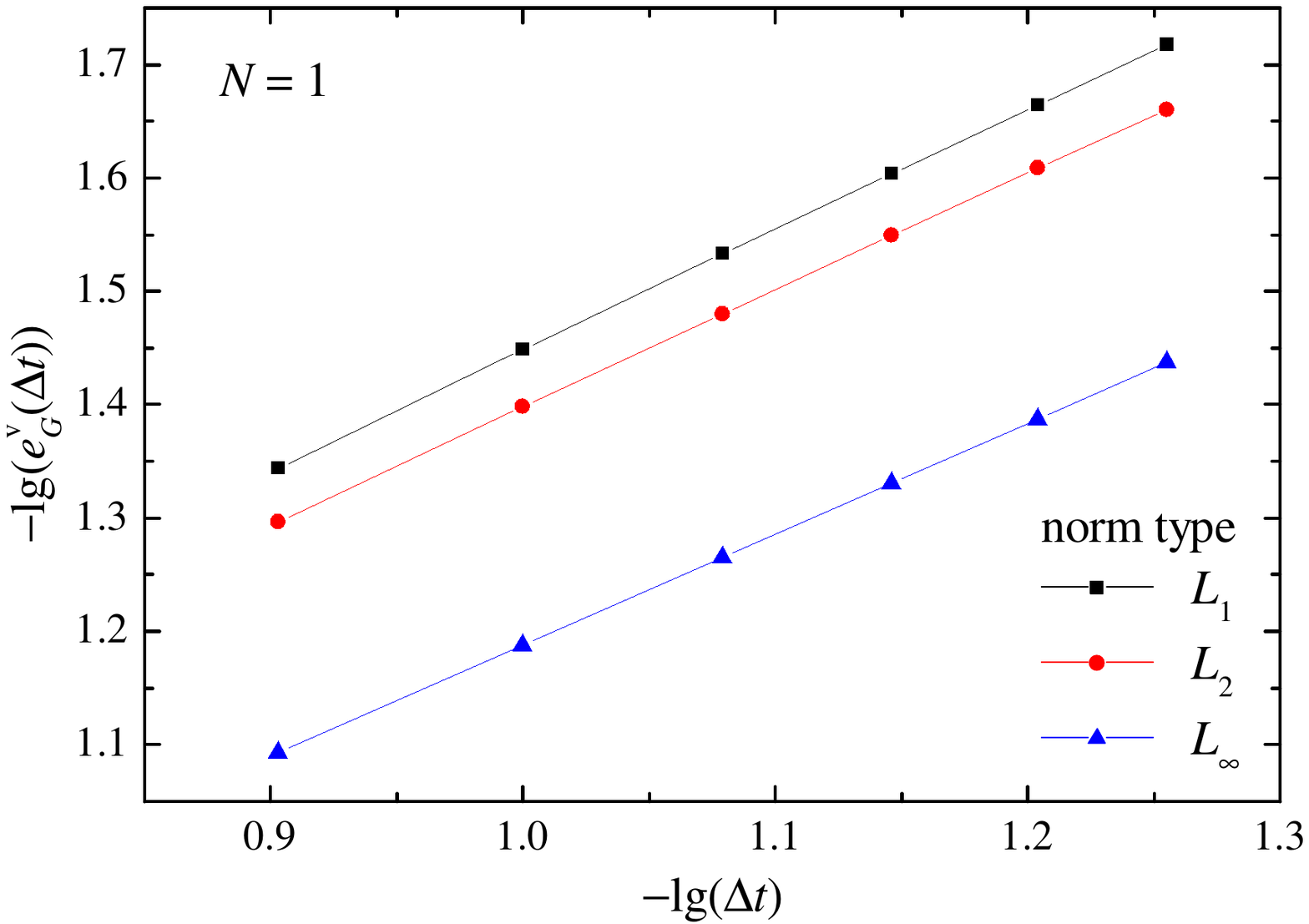}
\vspace{-8mm}\caption{\label{fig:hess_2_ind2_errors:e1}}
\end{subfigure}\hspace{6mm}
\begin{subfigure}{0.275\textwidth}
\includegraphics[width=\textwidth]{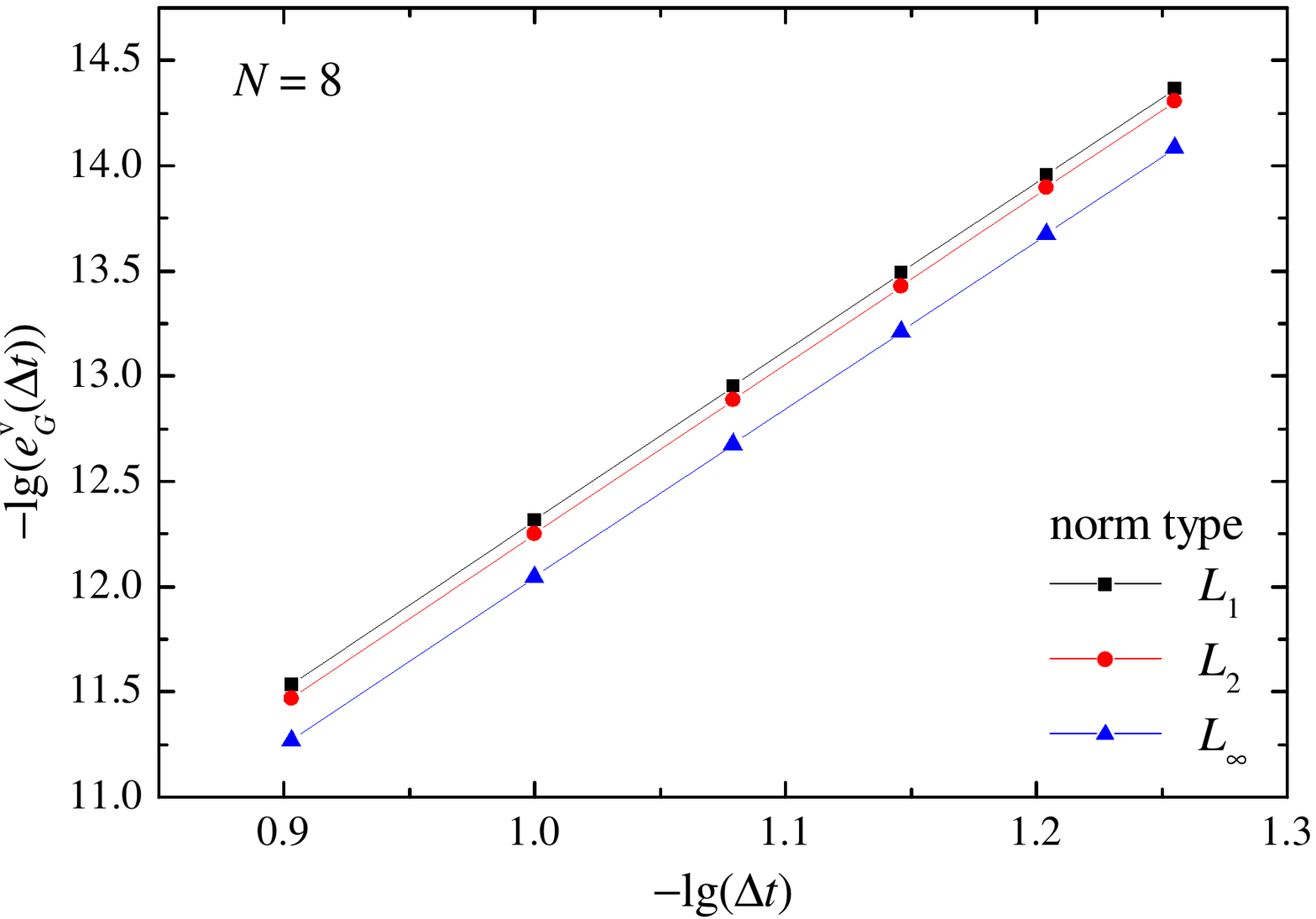}
\vspace{-8mm}\caption{\label{fig:hess_2_ind2_errors:e2}}
\end{subfigure}\hspace{6mm}
\begin{subfigure}{0.275\textwidth}
\includegraphics[width=\textwidth]{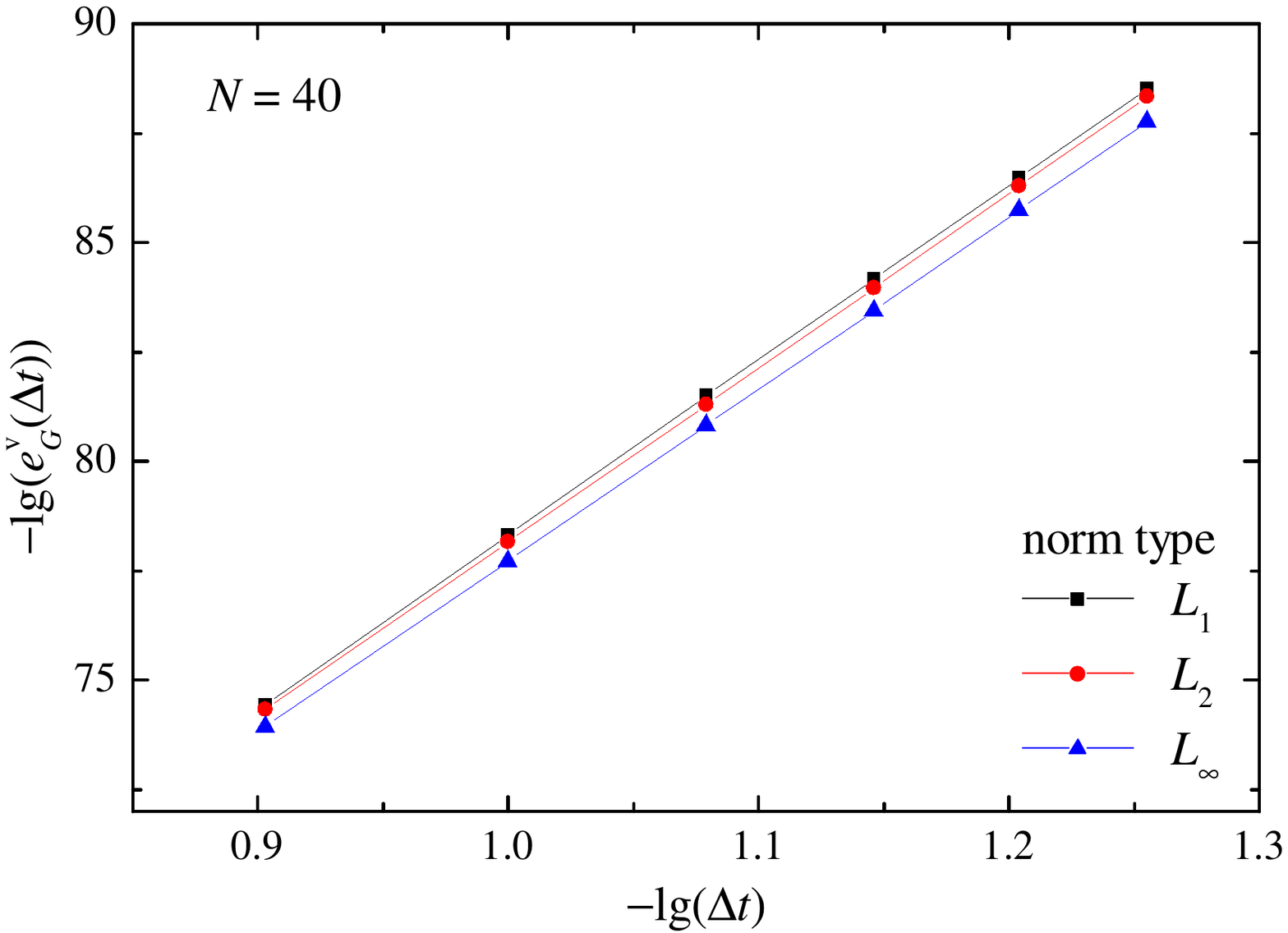}
\vspace{-8mm}\caption{\label{fig:hess_2_ind2_errors:e3}}
\end{subfigure}\\[-2mm]
\begin{subfigure}{0.275\textwidth}
\includegraphics[width=\textwidth]{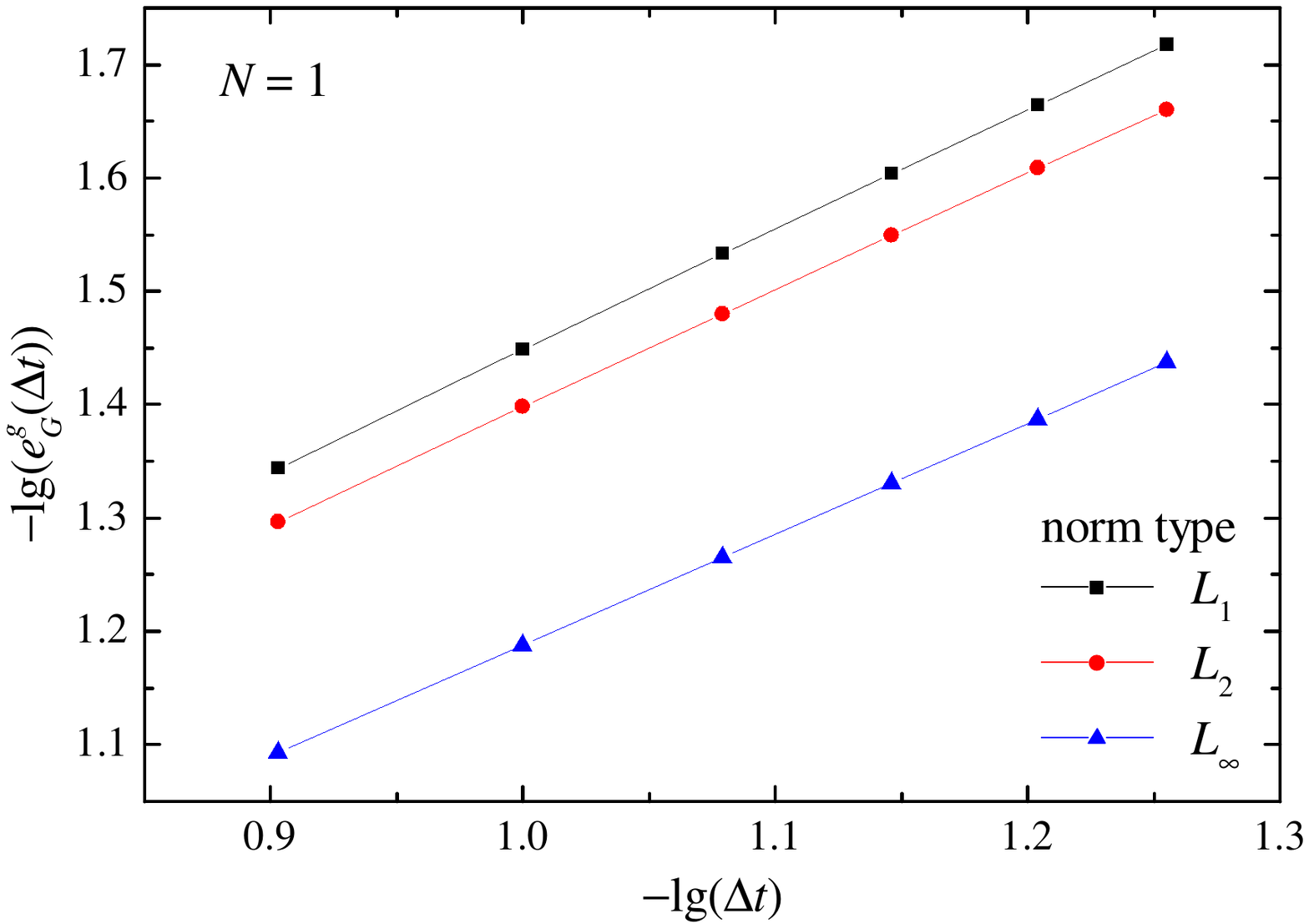}
\vspace{-8mm}\caption{\label{fig:hess_2_ind2_errors:f1}}
\end{subfigure}\hspace{6mm}
\begin{subfigure}{0.275\textwidth}
\includegraphics[width=\textwidth]{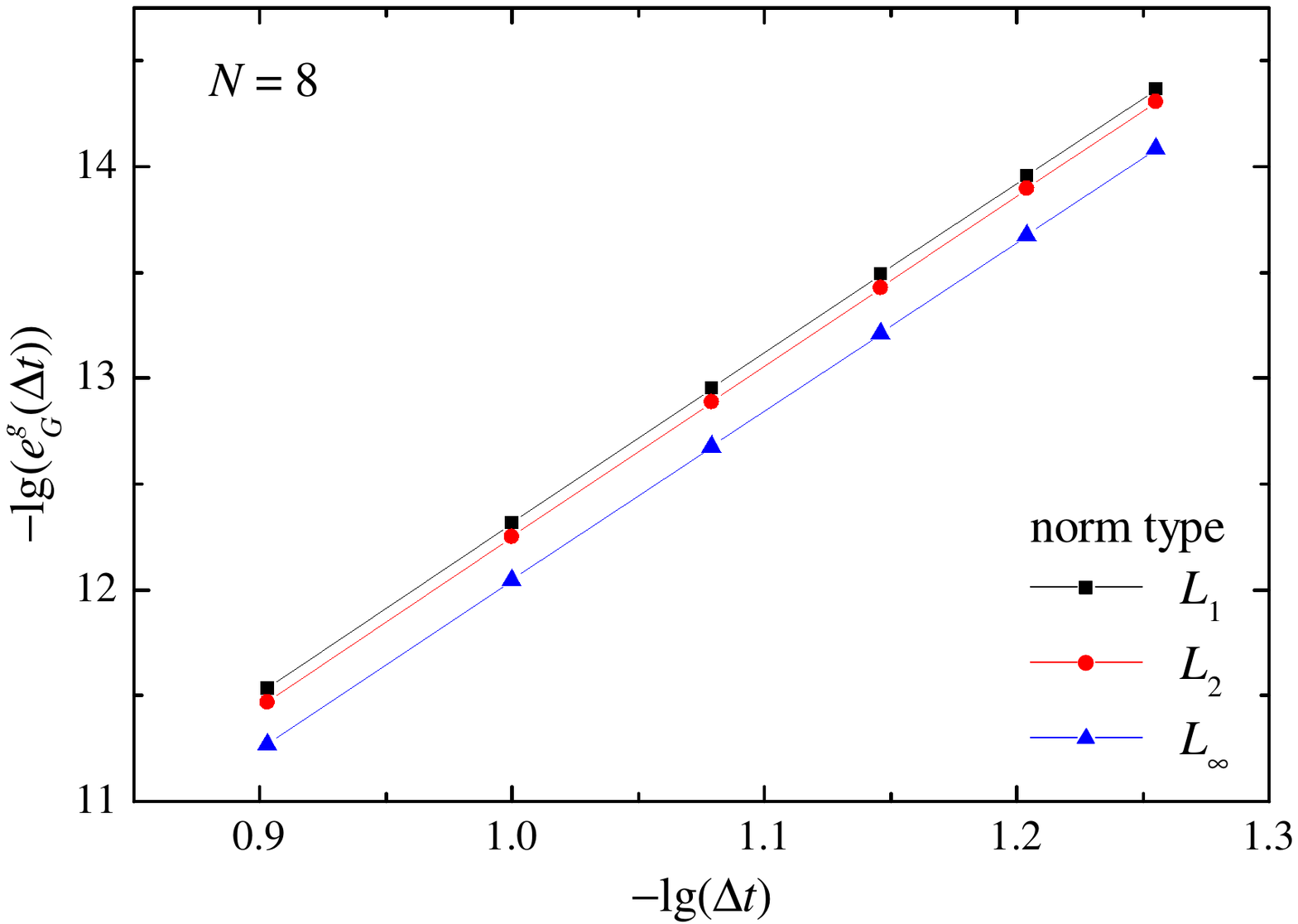}
\vspace{-8mm}\caption{\label{fig:hess_2_ind2_errors:f2}}
\end{subfigure}\hspace{6mm}
\begin{subfigure}{0.275\textwidth}
\includegraphics[width=\textwidth]{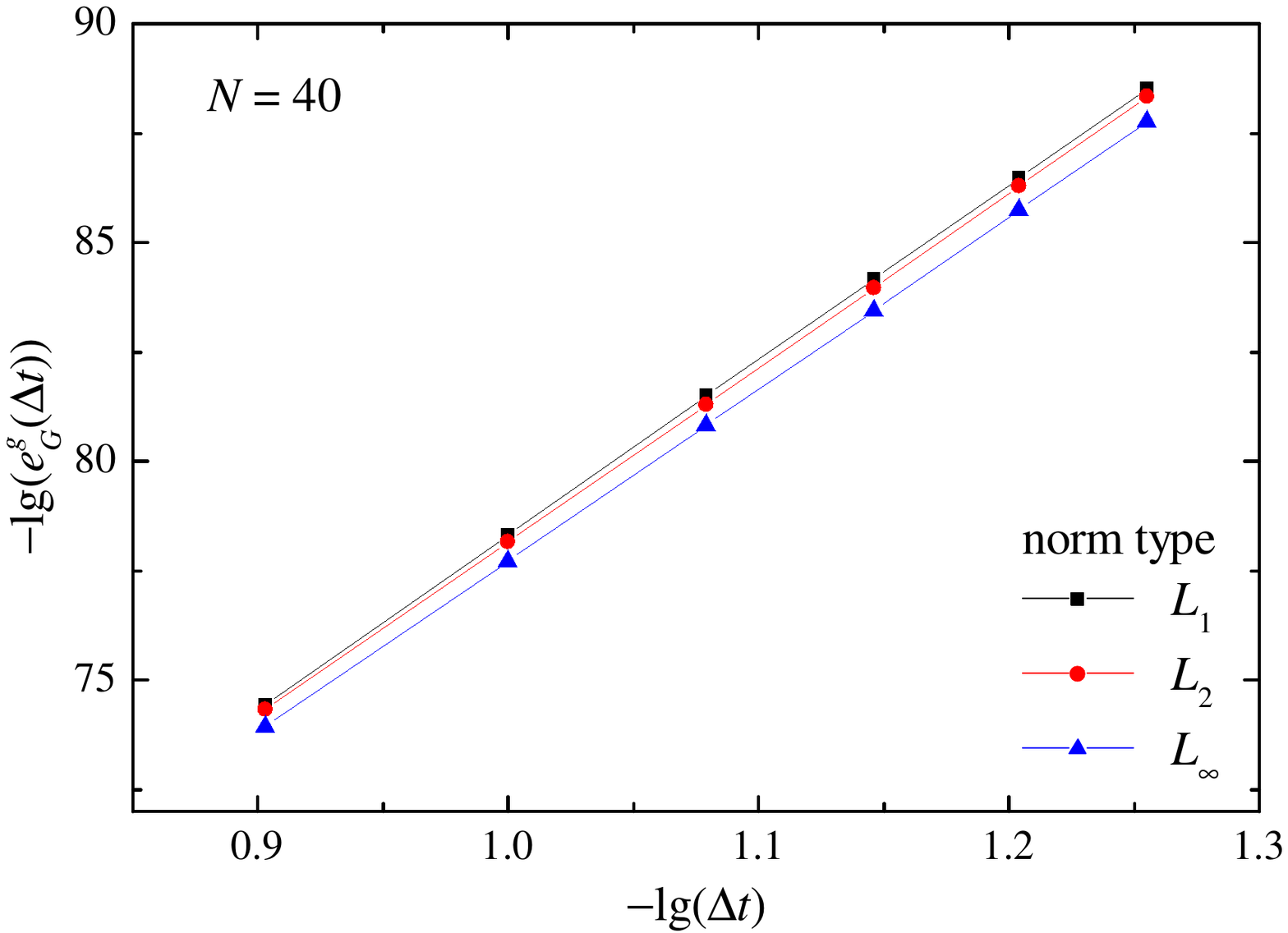}
\vspace{-8mm}\caption{\label{fig:hess_2_ind2_errors:f3}}
\end{subfigure}\\[-2mm]
\caption{%
Log-log plot of the dependence of the global errors for the local solution $e_{L}^{u}$ (\subref{fig:hess_2_ind2_errors:a1}, \subref{fig:hess_2_ind2_errors:a2}, \subref{fig:hess_2_ind2_errors:a3}), $e_{L}^{v}$ (\subref{fig:hess_2_ind2_errors:b1}, \subref{fig:hess_2_ind2_errors:b2}, \subref{fig:hess_2_ind2_errors:b3}), $e_{L}^{g}$ (\subref{fig:hess_2_ind2_errors:c1}, \subref{fig:hess_2_ind2_errors:c2}, \subref{fig:hess_2_ind2_errors:c3}) and the solution at nodes $e_{G}^{u}$ (\subref{fig:hess_2_ind2_errors:d1}, \subref{fig:hess_2_ind2_errors:d2}, \subref{fig:hess_2_ind2_errors:d3}), $e_{G}^{v}$ (\subref{fig:hess_2_ind2_errors:e1}, \subref{fig:hess_2_ind2_errors:e2}, \subref{fig:hess_2_ind2_errors:e3}), $e_{G}^{g}$ (\subref{fig:hess_2_ind2_errors:f1}, \subref{fig:hess_2_ind2_errors:f2}, \subref{fig:hess_2_ind2_errors:f3}) on the discretization step $\mathrm{\Delta}t$, obtained in the norms $L_{1}$, $L_{2}$ and $L_{\infty}$, by numerical solution of the DAE system (\ref{eq:hess_dae_ind_2}) of index 2 obtained using polynomials with degrees $N = 1$ (\subref{fig:hess_2_ind2_errors:a1}, \subref{fig:hess_2_ind2_errors:b1}, \subref{fig:hess_2_ind2_errors:c1}, \subref{fig:hess_2_ind2_errors:d1}, \subref{fig:hess_2_ind2_errors:e1}, \subref{fig:hess_2_ind2_errors:f1}), $N = 8$ (\subref{fig:hess_2_ind2_errors:a2}, \subref{fig:hess_2_ind2_errors:b2}, \subref{fig:hess_2_ind2_errors:c2}, \subref{fig:hess_2_ind2_errors:d2}, \subref{fig:hess_2_ind2_errors:e2}, \subref{fig:hess_2_ind2_errors:f2}), $N = 40$ (\subref{fig:hess_2_ind2_errors:a3}, \subref{fig:hess_2_ind2_errors:b3}, \subref{fig:hess_2_ind2_errors:c3}, \subref{fig:hess_2_ind2_errors:d3}, \subref{fig:hess_2_ind2_errors:e3}, \subref{fig:hess_2_ind2_errors:f3}).
}
\label{fig:hess_2_ind2_errors}
\end{figure} 

\begin{table*}[h!]
\centering
\caption{%
Convergence orders $p_{L_{1}}^{n}$, $p_{L_{2}}^{n}$, $p_{L_{\infty}}^{n}$, calculated in norms $L_{1}$, $L_{2}$, $L_{\infty}$, of \textit{the numerical solution at the nodes} $(\mathbf{u}_{n}, \mathbf{v}_{n})$ of the ADER-DG method for the DAE system (\ref{eq:hess_dae_ind_2}) of index 2; $N$ is the degree of the basis polynomials $\varphi_{p}$. Orders $p^{n, u}$ are calculated for solution $\mathbf{u}_{n}$; orders $p^{n, v}$ --- for solution $\mathbf{v}_{n}$; orders $p^{n, g}$ --- for the conditions $\mathbf{g} = 0$ on the numerical solution at the nodes $(\mathbf{u}_{n}, \mathbf{v}_{n})$. The theoretical value of convergence order $p_{\rm th.}^{n} = 2N+1$ is applicable for the ADER-DG method for ODE problems and is presented for comparison.
}
\label{tab:conv_orders_nodes_hess_2_ind2}
\begin{tabular}{@{}|l|lll|lll|lll|c|@{}}
\toprule
$N$ & $p_{L_{1}}^{n, u}$ & $p_{L_{2}}^{n, u}$ & $p_{L_{\infty}}^{n, u}$ & $p_{L_{1}}^{n, v}$ & $p_{L_{2}}^{n, v}$ & $p_{L_{\infty}}^{n, v}$ & $p_{L_{1}}^{n, g}$ & $p_{L_{2}}^{n, g}$ & $p_{L_{\infty}}^{n, g}$ & $p_{\rm th.}^{n}$ \\
\midrule
$1$	&	$2.07$	&	$2.04$	&	$1.98$	&	$1.06$	&	$1.03$	&	$0.98$	&	$1.06$	&	$1.03$	&	$0.98$	&	$3$\\
$2$	&	$2.99$	&	$2.99$	&	$2.88$	&	$2.03$	&	$2.03$	&	$1.88$	&	$2.03$	&	$2.03$	&	$1.88$	&	$5$\\
$3$	&	$4.00$	&	$3.97$	&	$3.67$	&	$2.96$	&	$2.91$	&	$2.68$	&	$2.96$	&	$2.91$	&	$2.68$	&	$7$\\
$4$	&	$5.02$	&	$5.01$	&	$4.94$	&	$4.09$	&	$4.10$	&	$3.94$	&	$4.09$	&	$4.10$	&	$3.94$	&	$9$\\
$5$	&	$6.00$	&	$6.01$	&	$5.98$	&	$5.00$	&	$5.00$	&	$4.98$	&	$5.00$	&	$5.00$	&	$4.98$	&	$11$\\
$6$	&	$6.98$	&	$6.97$	&	$6.75$	&	$5.99$	&	$6.00$	&	$5.76$	&	$5.99$	&	$6.00$	&	$5.76$	&	$13$\\
$7$	&	$7.98$	&	$7.97$	&	$7.75$	&	$6.97$	&	$6.94$	&	$6.76$	&	$6.97$	&	$6.94$	&	$6.76$	&	$15$\\
$8$	&	$9.01$	&	$9.01$	&	$8.93$	&	$8.05$	&	$8.06$	&	$8.00$	&	$8.05$	&	$8.06$	&	$8.00$	&	$17$\\
$9$	&	$10.00$	&	$10.00$	&	$10.00$	&	$8.98$	&	$8.99$	&	$9.00$	&	$8.98$	&	$8.99$	&	$9.00$	&	$19$\\
$10$	&	$10.96$	&	$10.95$	&	$10.70$	&	$9.98$	&	$9.98$	&	$9.71$	&	$9.98$	&	$9.98$	&	$9.71$	&	$21$\\
$11$	&	$11.99$	&	$11.97$	&	$11.78$	&	$10.99$	&	$10.96$	&	$10.80$	&	$10.99$	&	$10.96$	&	$10.80$	&	$23$\\
$12$	&	$13.01$	&	$13.01$	&	$12.94$	&	$12.04$	&	$12.04$	&	$11.97$	&	$12.04$	&	$12.04$	&	$11.97$	&	$25$\\
$13$	&	$13.99$	&	$13.96$	&	$13.69$	&	$12.97$	&	$12.96$	&	$12.81$	&	$12.97$	&	$12.96$	&	$12.81$	&	$27$\\
$14$	&	$14.97$	&	$14.94$	&	$14.69$	&	$13.99$	&	$13.97$	&	$13.70$	&	$13.99$	&	$13.97$	&	$13.70$	&	$29$\\
$15$	&	$16.00$	&	$15.98$	&	$15.82$	&	$15.00$	&	$14.97$	&	$14.83$	&	$15.00$	&	$14.97$	&	$14.83$	&	$31$\\
$16$	&	$17.02$	&	$17.01$	&	$16.95$	&	$16.04$	&	$16.03$	&	$15.97$	&	$16.04$	&	$16.03$	&	$15.97$	&	$33$\\
$17$	&	$17.97$	&	$17.92$	&	$17.42$	&	$16.94$	&	$16.91$	&	$16.43$	&	$16.94$	&	$16.91$	&	$16.43$	&	$35$\\
$18$	&	$18.97$	&	$18.94$	&	$18.70$	&	$17.98$	&	$17.97$	&	$17.71$	&	$17.98$	&	$17.97$	&	$17.71$	&	$37$\\
$19$	&	$20.00$	&	$19.99$	&	$19.87$	&	$19.00$	&	$18.98$	&	$18.88$	&	$19.00$	&	$18.98$	&	$18.88$	&	$39$\\
$20$	&	$21.03$	&	$21.00$	&	$20.98$	&	$20.05$	&	$20.02$	&	$19.99$	&	$20.05$	&	$20.02$	&	$19.99$	&	$41$\\
\midrule
$25$	&	$25.93$	&	$25.87$	&	$25.51$	&	$24.92$	&	$24.86$	&	$24.52$	&	$24.92$	&	$24.86$	&	$24.52$	&	$51$\\
$30$	&	$30.98$	&	$30.96$	&	$30.76$	&	$30.00$	&	$29.98$	&	$29.77$	&	$30.00$	&	$29.98$	&	$29.77$	&	$61$\\
$35$	&	$36.03$	&	$36.02$	&	$35.86$	&	$35.03$	&	$35.02$	&	$34.86$	&	$35.03$	&	$35.02$	&	$34.86$	&	$71$\\
$40$	&	$40.96$	&	$40.75$	&	$40.23$	&	$39.99$	&	$39.77$	&	$39.24$	&	$39.99$	&	$39.77$	&	$39.24$	&	$81$\\
\bottomrule
\end{tabular}
\end{table*} 

\begin{table*}[h!]
\centering
\caption{%
Convergence orders $p_{L_{1}}^{l}$, $p_{L_{2}}^{l}$, $p_{L_{\infty}}^{l}$, calculated in norms $L_{1}$, $L_{2}$, $L_{\infty}$, of \textit{the local solution} $(\mathbf{u}_{L}, \mathbf{v}_{L})$ (represented between the nodes) of the ADER-DG method for the DAE system (\ref{eq:hess_dae_ind_2}) of index 2; $N$ is the degree of the basis polynomials $\varphi_{p}$. Orders $p^{l, u}$ are calculated for solution $\mathbf{u}_{L}$; orders $p^{l, v}$ --- for solution $\mathbf{v}_{L}$; orders $p^{l, g}$ --- for the conditions $\mathbf{g} = 0$ on the local solution $(\mathbf{u}_{L}, \mathbf{v}_{L})$. The theoretical value of convergence order $p_{\rm th.}^{l} = N+1$ is applicable for the ADER-DG method for ODE problems and is presented for comparison.
}
\label{tab:conv_orders_local_hess_2_ind2}
\begin{tabular}{@{}|l|lll|lll|lll|c|@{}}
\toprule
$N$ & $p_{L_{1}}^{l, u}$ & $p_{L_{2}}^{l, u}$ & $p_{L_{\infty}}^{l, u}$ & $p_{L_{1}}^{l, v}$ & $p_{L_{2}}^{l, v}$ & $p_{L_{\infty}}^{l, v}$ & $p_{L_{1}}^{l, g}$ & $p_{L_{2}}^{l, g}$ & $p_{L_{\infty}}^{l, g}$ & $p_{\rm th.}^{l}$ \\
\midrule
$1$	&	$2.01$	&	$2.01$	&	$2.00$	&	$0.99$	&	$0.98$	&	$0.98$	&	$0.99$	&	$0.98$	&	$0.98$	&	$2$\\
$2$	&	$3.00$	&	$3.00$	&	$2.74$	&	$1.97$	&	$1.98$	&	$1.88$	&	$1.97$	&	$1.98$	&	$1.88$	&	$3$\\
$3$	&	$4.00$	&	$3.99$	&	$3.85$	&	$2.97$	&	$2.88$	&	$2.54$	&	$2.97$	&	$2.88$	&	$2.54$	&	$4$\\
$4$	&	$5.00$	&	$5.00$	&	$4.96$	&	$4.00$	&	$3.98$	&	$3.98$	&	$4.00$	&	$3.98$	&	$3.98$	&	$5$\\
$5$	&	$6.01$	&	$6.00$	&	$6.00$	&	$4.95$	&	$4.94$	&	$4.93$	&	$4.95$	&	$4.94$	&	$4.93$	&	$6$\\
$6$	&	$6.98$	&	$6.96$	&	$6.72$	&	$5.98$	&	$5.98$	&	$5.79$	&	$5.98$	&	$5.98$	&	$5.79$	&	$7$\\
$7$	&	$8.00$	&	$7.98$	&	$7.84$	&	$6.97$	&	$6.89$	&	$6.70$	&	$6.97$	&	$6.89$	&	$6.70$	&	$8$\\
$8$	&	$9.01$	&	$9.00$	&	$9.00$	&	$7.97$	&	$7.93$	&	$8.01$	&	$7.97$	&	$7.93$	&	$8.01$	&	$9$\\
$9$	&	$10.01$	&	$10.01$	&	$9.99$	&	$8.92$	&	$8.91$	&	$8.98$	&	$8.92$	&	$8.91$	&	$8.98$	&	$10$\\
$10$	&	$10.97$	&	$10.95$	&	$10.72$	&	$10.02$	&	$9.98$	&	$9.74$	&	$10.02$	&	$9.98$	&	$9.74$	&	$11$\\
$11$	&	$12.00$	&	$11.99$	&	$11.86$	&	$10.97$	&	$10.89$	&	$10.76$	&	$10.97$	&	$10.89$	&	$10.76$	&	$12$\\
$12$	&	$13.02$	&	$13.01$	&	$12.97$	&	$11.96$	&	$11.90$	&	$11.92$	&	$11.96$	&	$11.90$	&	$11.92$	&	$13$\\
$13$	&	$13.97$	&	$13.96$	&	$13.79$	&	$12.92$	&	$12.89$	&	$12.76$	&	$12.92$	&	$12.89$	&	$12.76$	&	$14$\\
$14$	&	$14.98$	&	$14.95$	&	$14.73$	&	$14.05$	&	$13.97$	&	$13.72$	&	$14.05$	&	$13.97$	&	$13.72$	&	$15$\\
$15$	&	$16.00$	&	$15.99$	&	$15.89$	&	$14.95$	&	$14.90$	&	$14.80$	&	$14.95$	&	$14.90$	&	$14.80$	&	$16$\\
$16$	&	$17.03$	&	$17.02$	&	$16.98$	&	$15.94$	&	$15.89$	&	$15.80$	&	$15.94$	&	$15.89$	&	$15.80$	&	$17$\\
$17$	&	$17.96$	&	$17.90$	&	$17.38$	&	$16.92$	&	$16.87$	&	$16.39$	&	$16.92$	&	$16.87$	&	$16.39$	&	$18$\\
$18$	&	$18.98$	&	$18.95$	&	$18.75$	&	$18.00$	&	$17.95$	&	$17.73$	&	$18.00$	&	$17.95$	&	$17.73$	&	$19$\\
$19$	&	$20.01$	&	$20.00$	&	$19.93$	&	$18.95$	&	$18.91$	&	$18.86$	&	$18.95$	&	$18.91$	&	$18.86$	&	$20$\\
$20$	&	$21.04$	&	$21.01$	&	$20.99$	&	$19.94$	&	$19.87$	&	$19.83$	&	$19.94$	&	$19.87$	&	$19.83$	&	$21$\\
\midrule
$25$	&	$25.95$	&	$25.88$	&	$25.54$	&	$24.95$	&	$24.86$	&	$24.50$	&	$24.95$	&	$24.86$	&	$24.50$	&	$26$\\
$30$	&	$30.99$	&	$30.98$	&	$30.81$	&	$29.94$	&	$29.95$	&	$29.78$	&	$29.94$	&	$29.95$	&	$29.78$	&	$31$\\
$35$	&	$36.04$	&	$36.03$	&	$35.94$	&	$34.92$	&	$34.93$	&	$34.86$	&	$34.92$	&	$34.93$	&	$34.86$	&	$36$\\
$40$	&	$40.96$	&	$40.75$	&	$40.26$	&	$39.94$	&	$39.76$	&	$39.26$	&	$39.94$	&	$39.76$	&	$39.26$	&	$41$\\
\bottomrule
\end{tabular}
\end{table*}

\begin{figure}[h!]
\captionsetup[subfigure]{%
	position=bottom,
	font+=smaller,
	textfont=normalfont,
	singlelinecheck=off,
	justification=raggedright
}
\centering
\begin{subfigure}{0.320\textwidth}
\includegraphics[width=\textwidth]{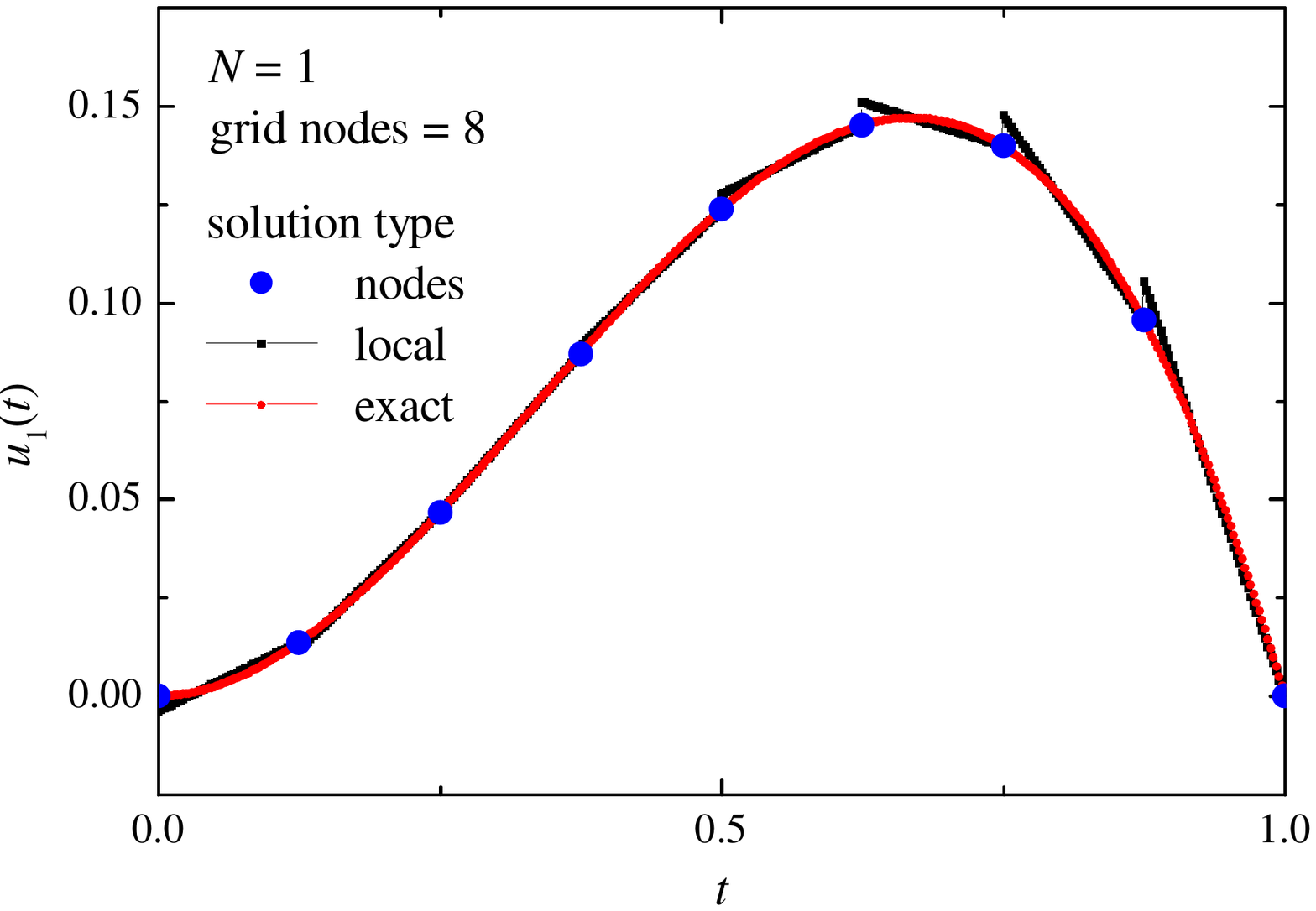}
\vspace{-8mm}\caption{\label{fig:hess_2_ind1_sol_uv:a1}}
\end{subfigure}
\begin{subfigure}{0.320\textwidth}
\includegraphics[width=\textwidth]{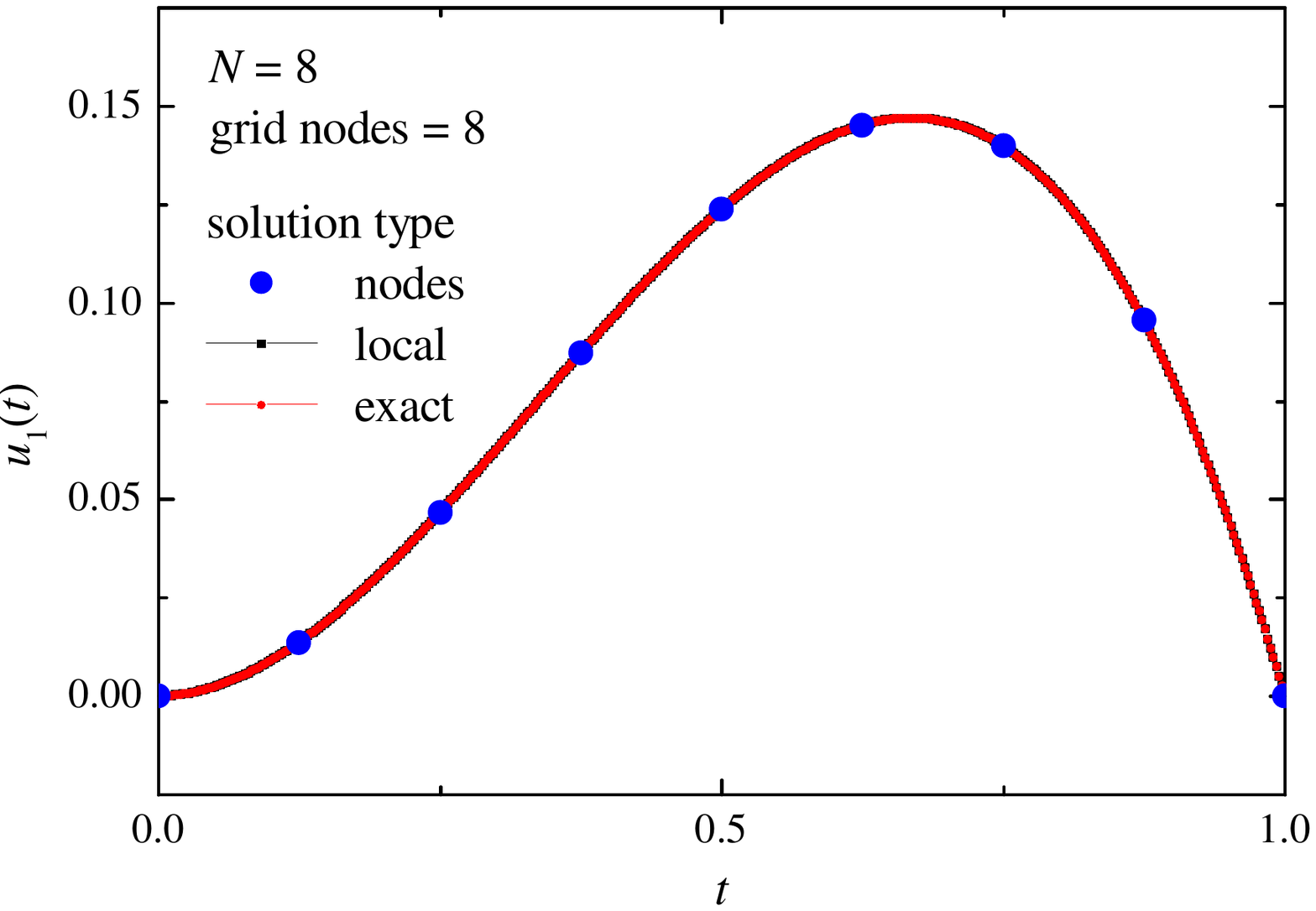}
\vspace{-8mm}\caption{\label{fig:hess_2_ind1_sol_uv:a2}}
\end{subfigure}
\begin{subfigure}{0.320\textwidth}
\includegraphics[width=\textwidth]{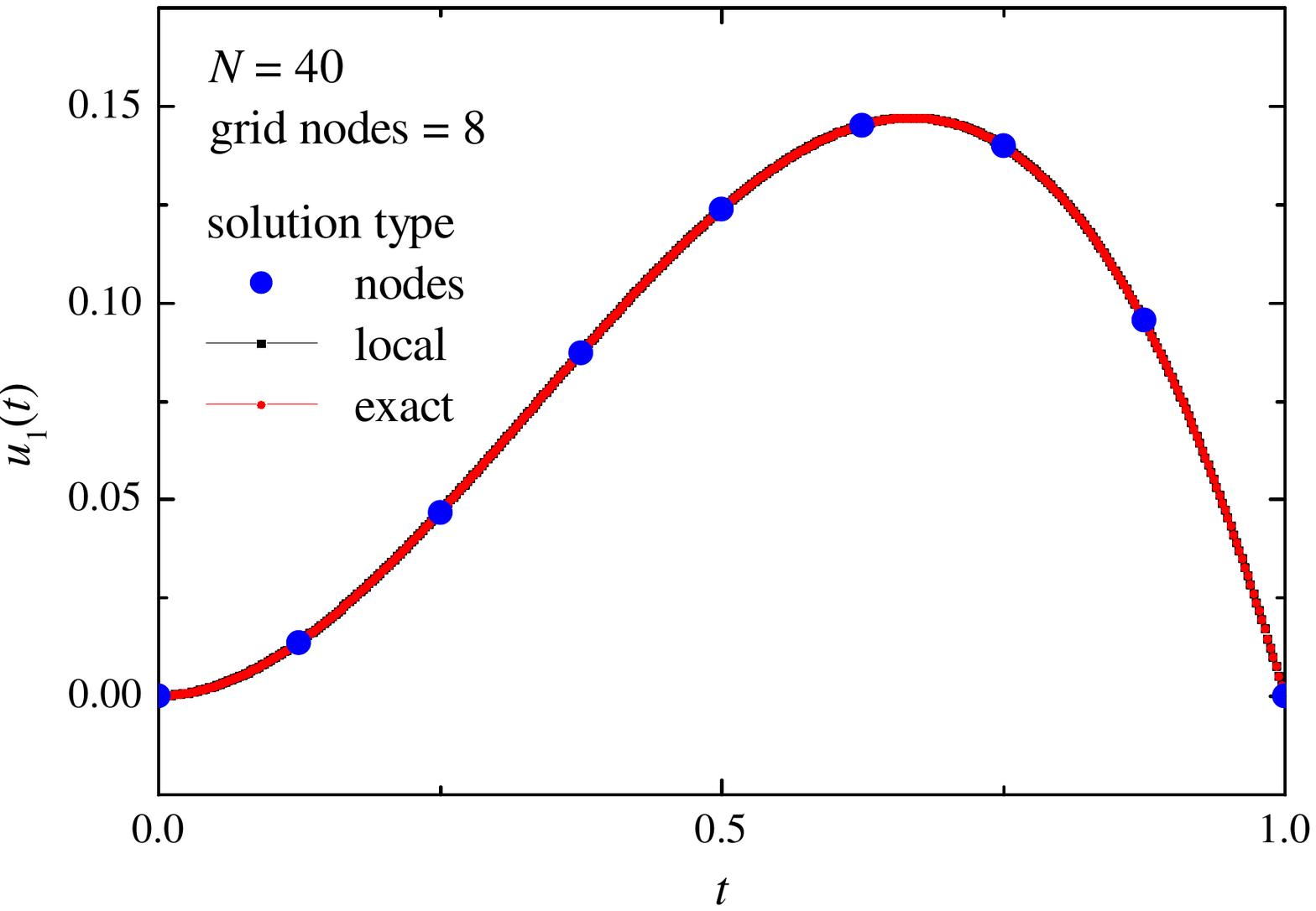}
\vspace{-8mm}\caption{\label{fig:hess_2_ind1_sol_uv:a3}}
\end{subfigure}\\
\begin{subfigure}{0.320\textwidth}
\includegraphics[width=\textwidth]{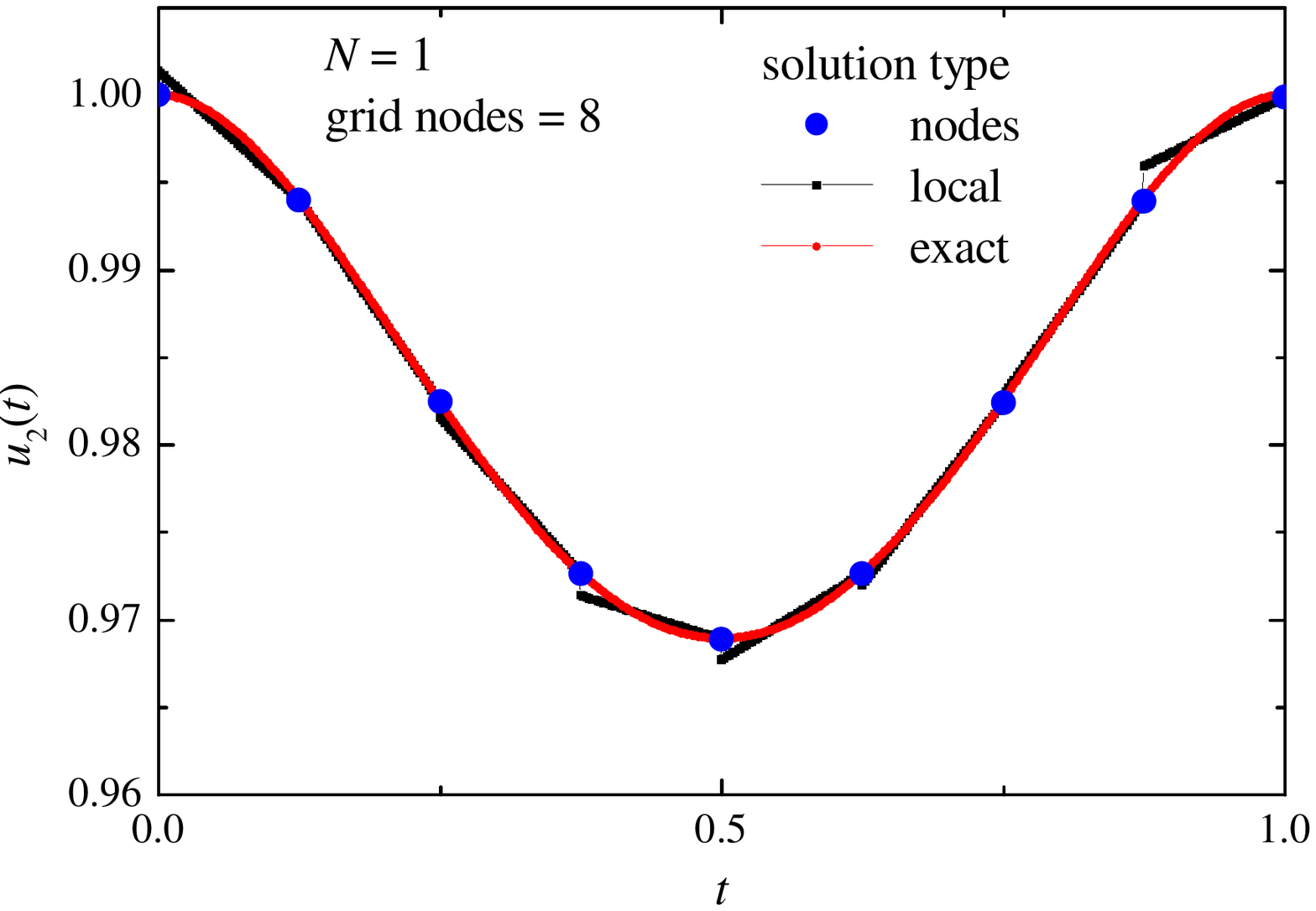}
\vspace{-8mm}\caption{\label{fig:hess_2_ind1_sol_uv:b1}}
\end{subfigure}
\begin{subfigure}{0.320\textwidth}
\includegraphics[width=\textwidth]{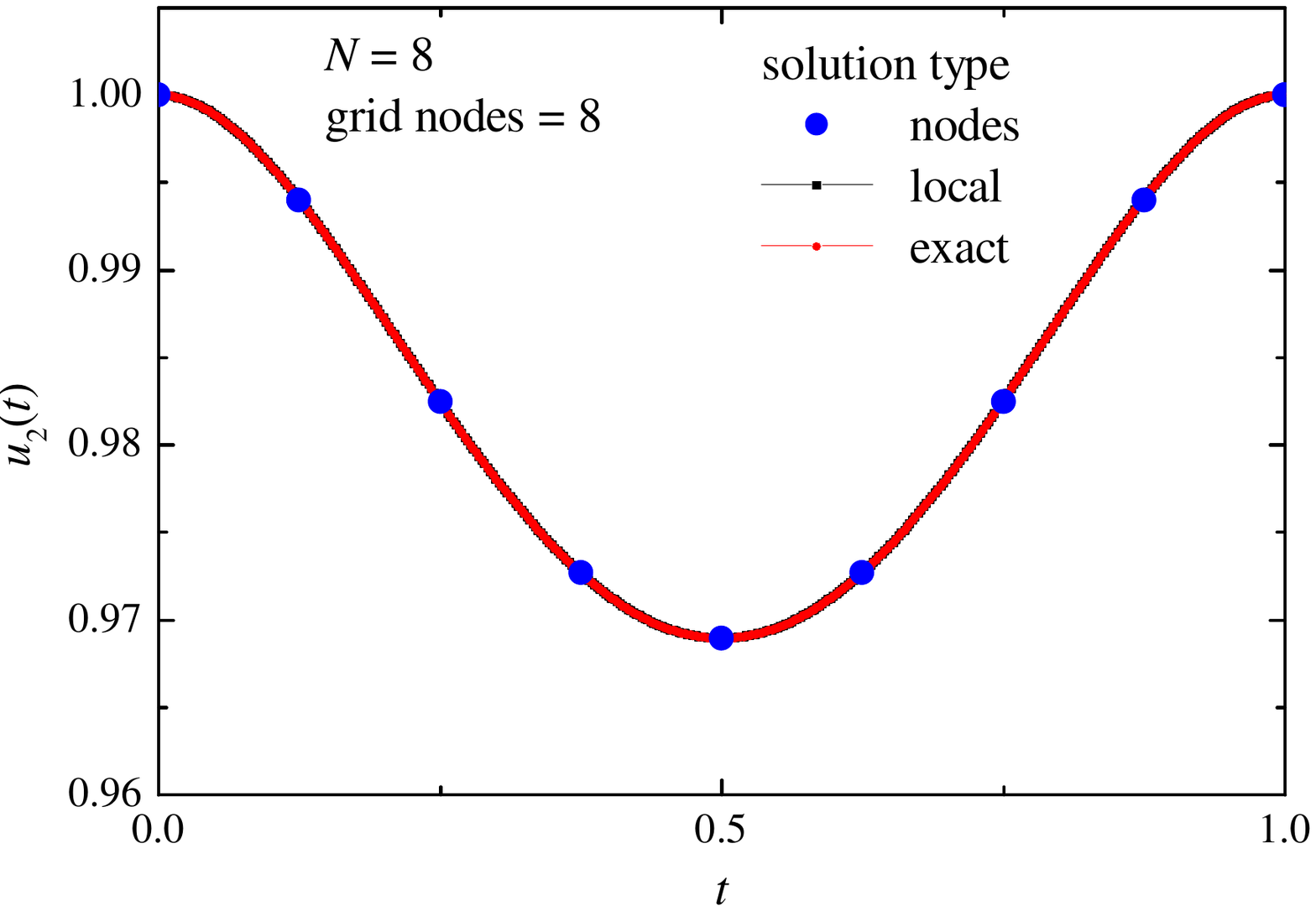}
\vspace{-8mm}\caption{\label{fig:hess_2_ind1_sol_uv:b2}}
\end{subfigure}
\begin{subfigure}{0.320\textwidth}
\includegraphics[width=\textwidth]{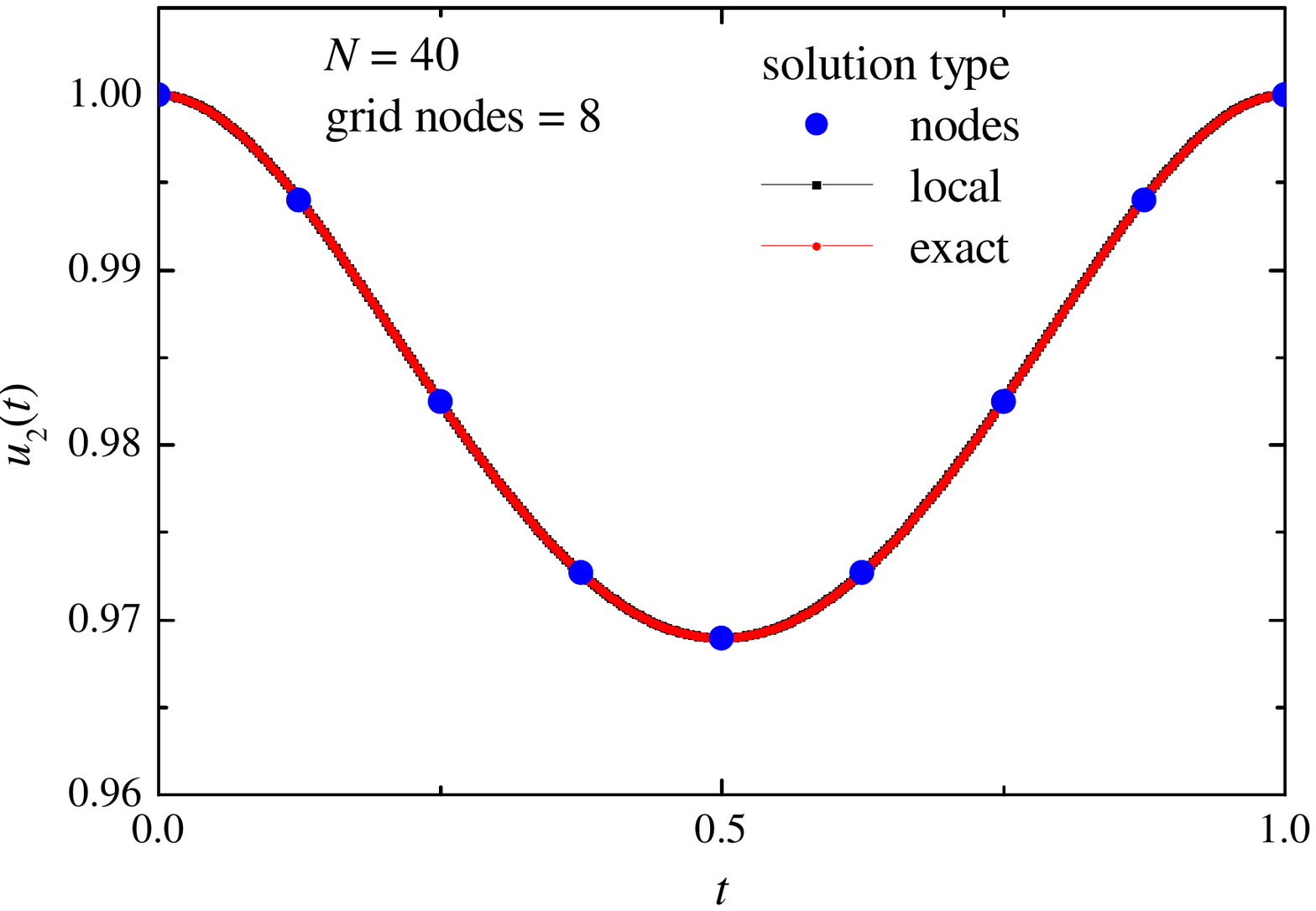}
\vspace{-8mm}\caption{\label{fig:hess_2_ind1_sol_uv:b3}}
\end{subfigure}\\
\begin{subfigure}{0.320\textwidth}
\includegraphics[width=\textwidth]{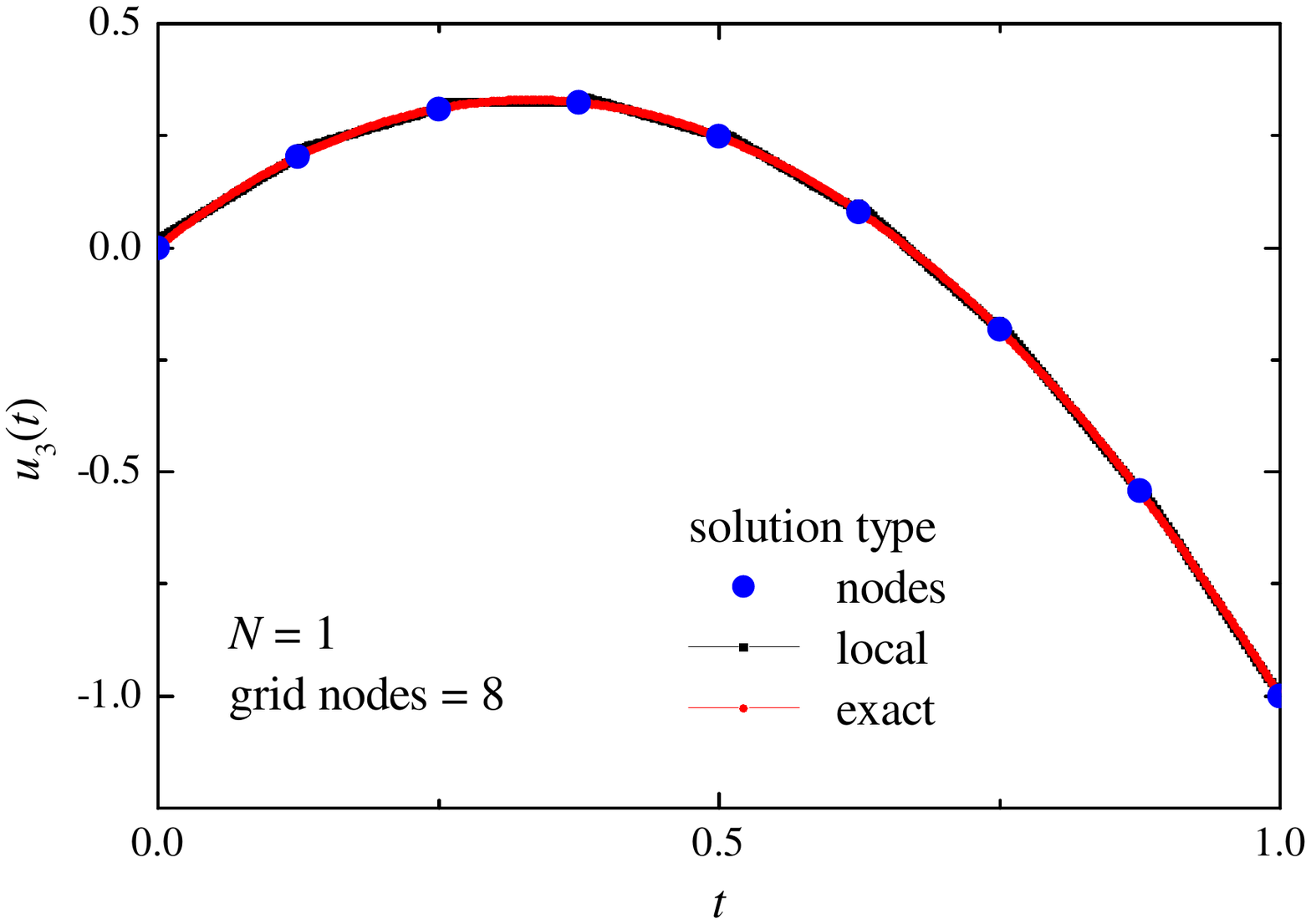}
\vspace{-8mm}\caption{\label{fig:hess_2_ind1_sol_uv:c1}}
\end{subfigure}
\begin{subfigure}{0.320\textwidth}
\includegraphics[width=\textwidth]{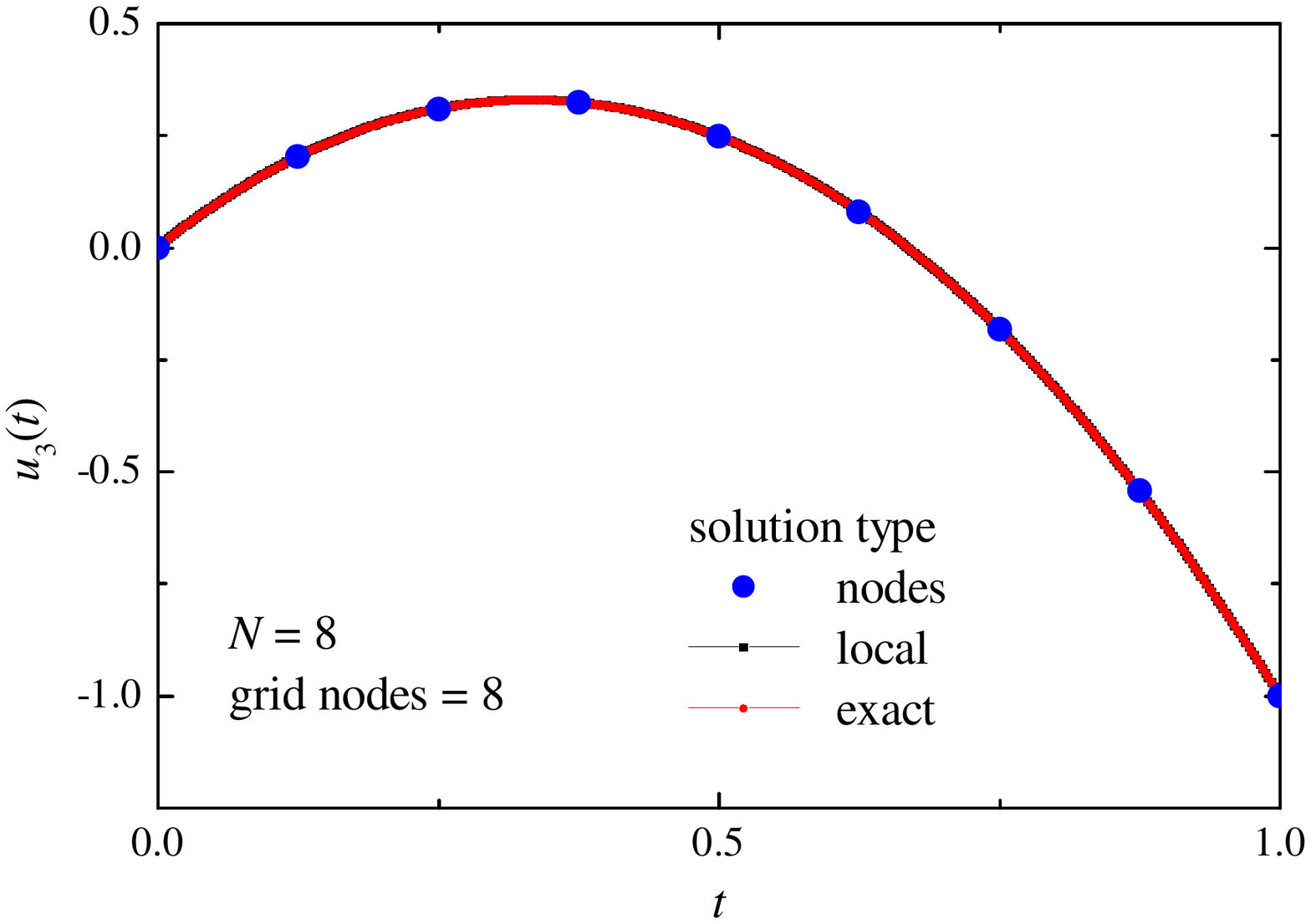}
\vspace{-8mm}\caption{\label{fig:hess_2_ind1_sol_uv:c2}}
\end{subfigure}
\begin{subfigure}{0.320\textwidth}
\includegraphics[width=\textwidth]{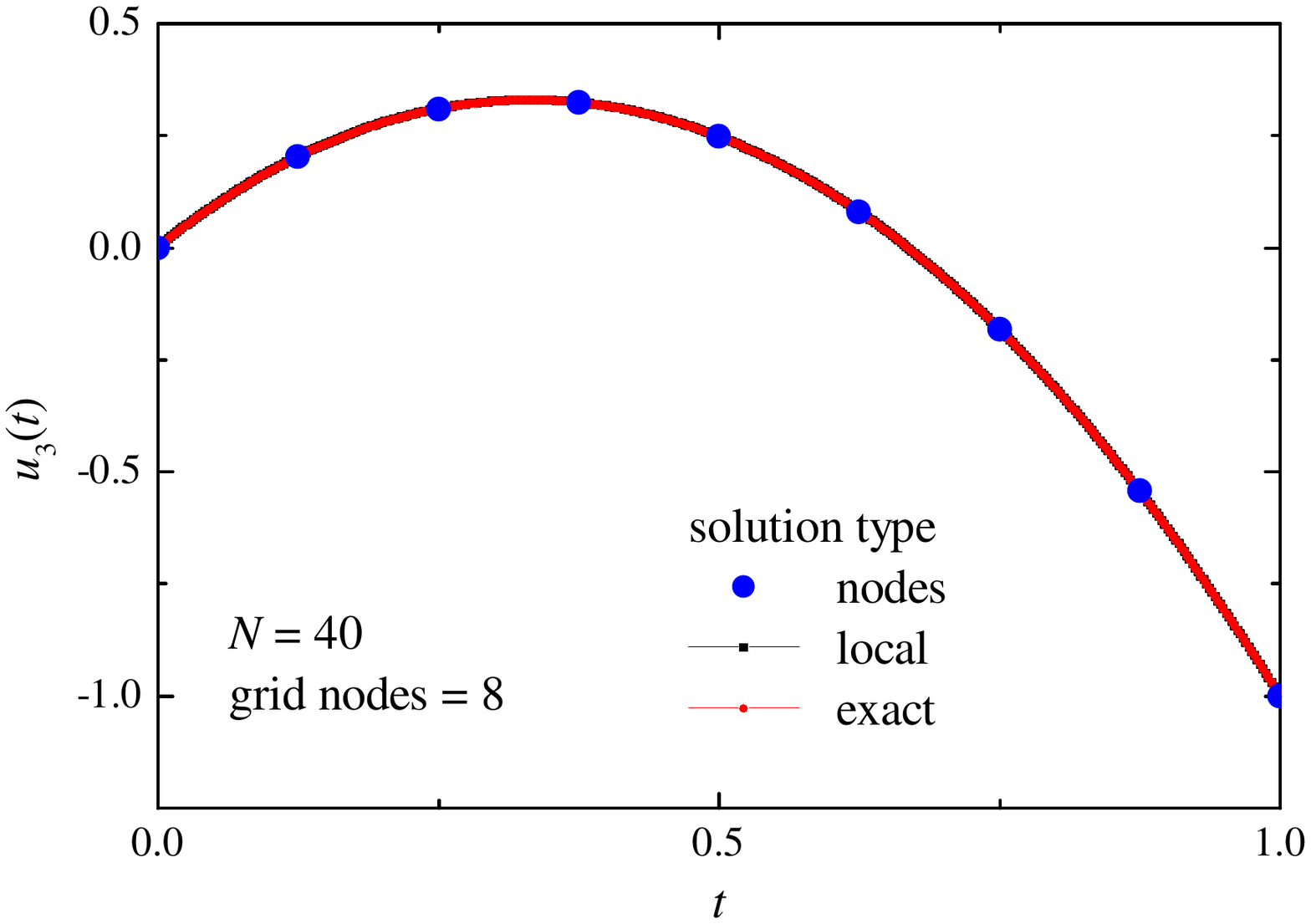}
\vspace{-8mm}\caption{\label{fig:hess_2_ind1_sol_uv:c3}}
\end{subfigure}\\
\begin{subfigure}{0.320\textwidth}
\includegraphics[width=\textwidth]{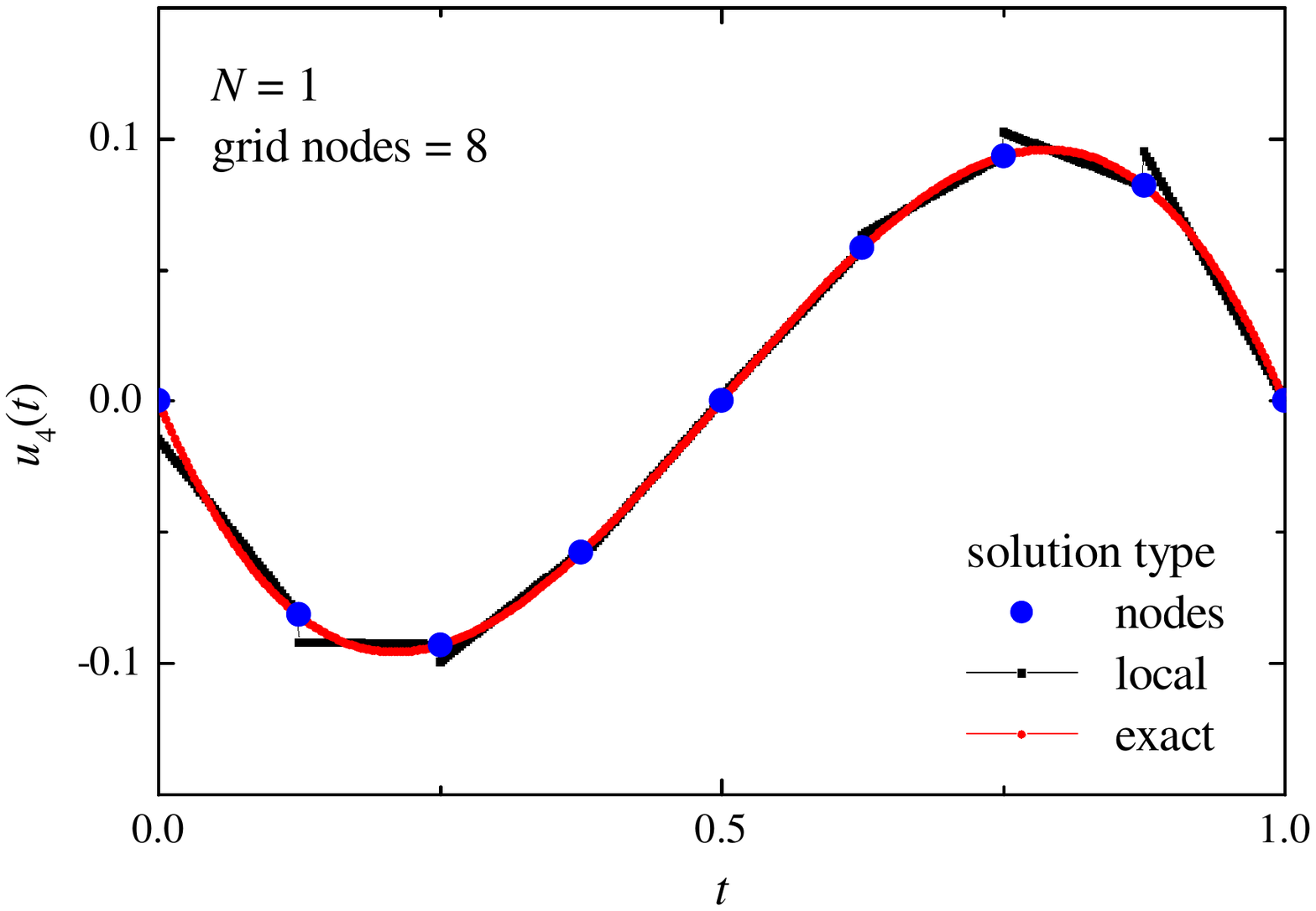}
\vspace{-8mm}\caption{\label{fig:hess_2_ind1_sol_uv:d1}}
\end{subfigure}
\begin{subfigure}{0.320\textwidth}
\includegraphics[width=\textwidth]{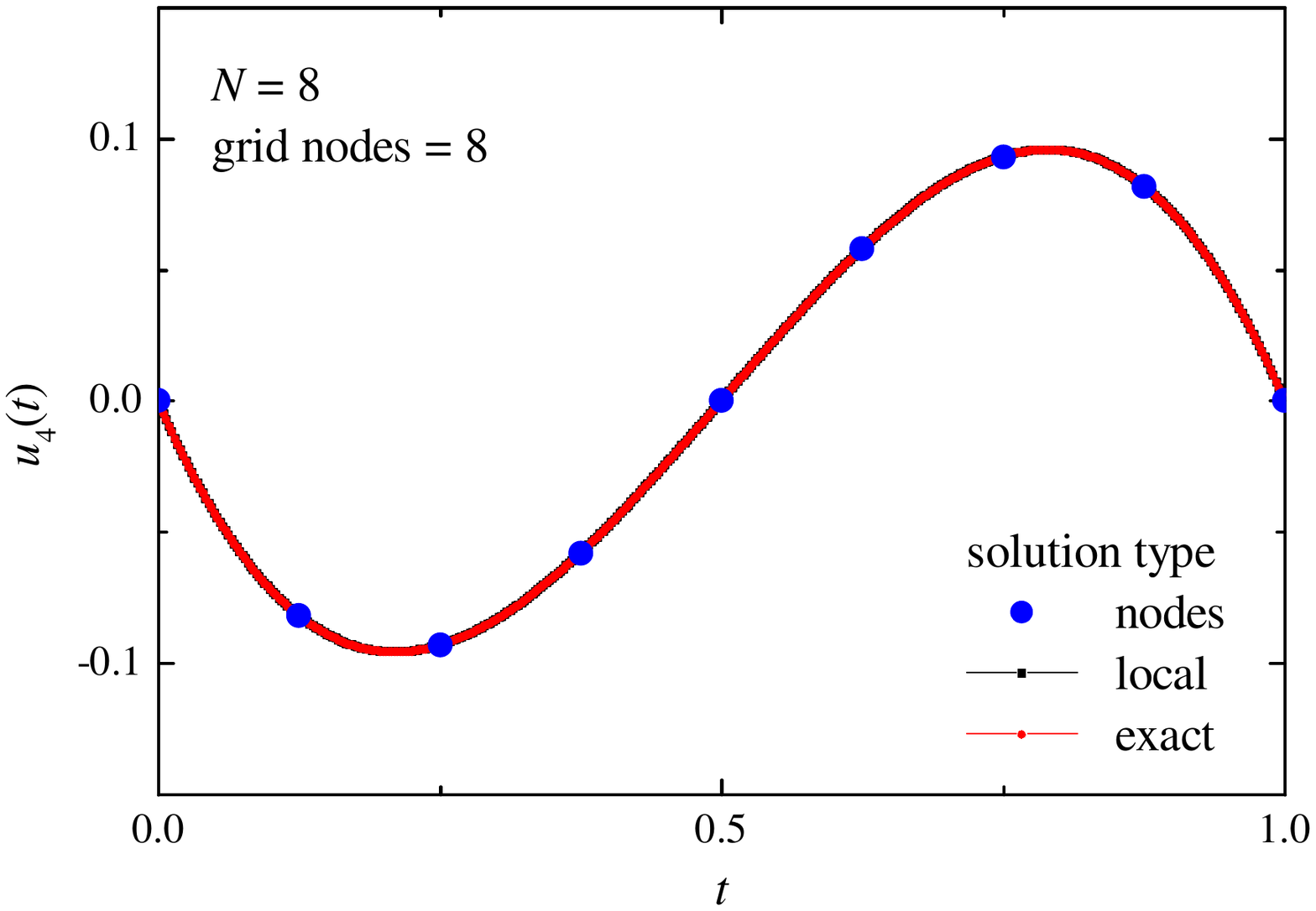}
\vspace{-8mm}\caption{\label{fig:hess_2_ind1_sol_uv:d2}}
\end{subfigure}
\begin{subfigure}{0.320\textwidth}
\includegraphics[width=\textwidth]{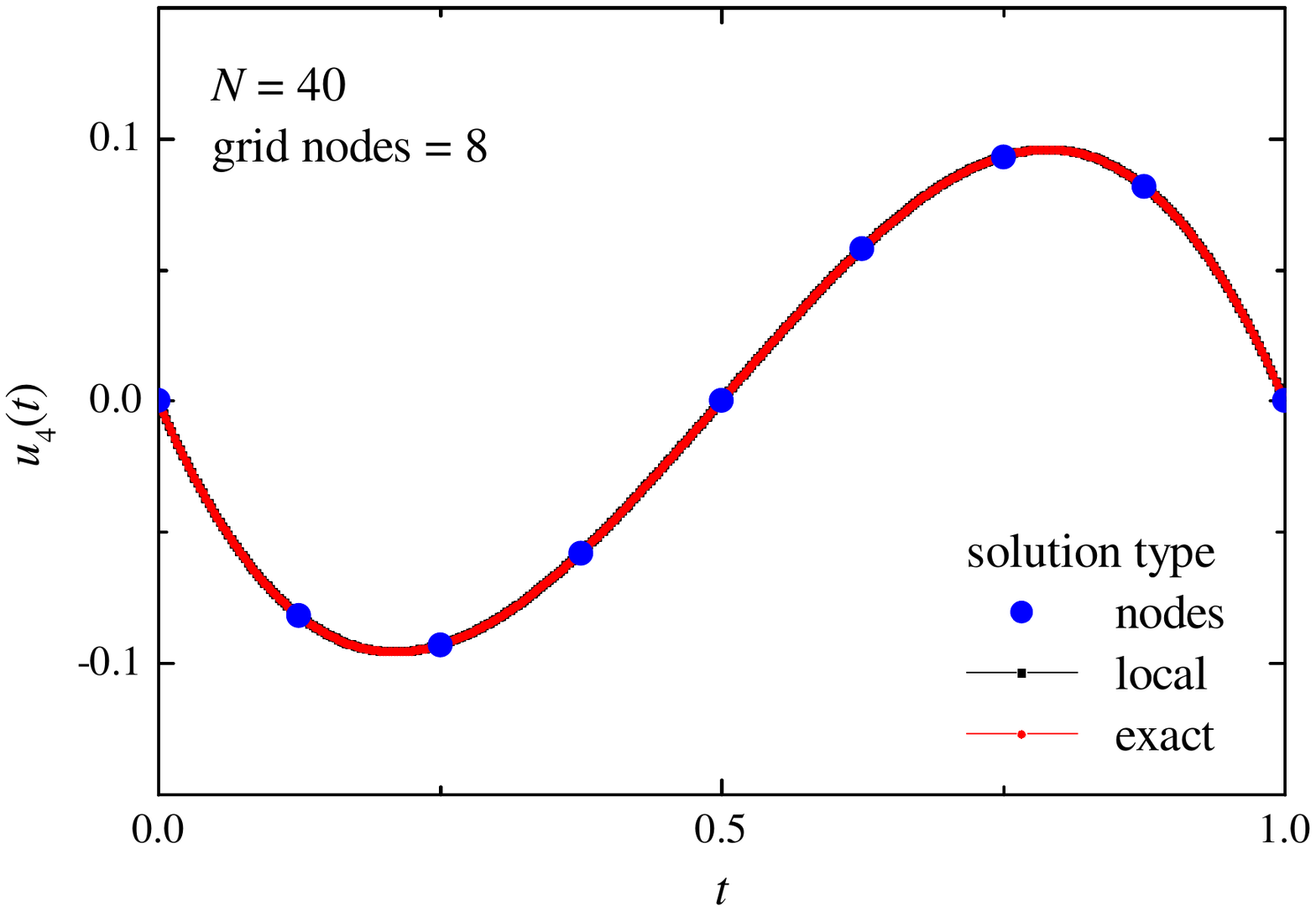}
\vspace{-8mm}\caption{\label{fig:hess_2_ind1_sol_uv:d3}}
\end{subfigure}\\
\begin{subfigure}{0.320\textwidth}
\includegraphics[width=\textwidth]{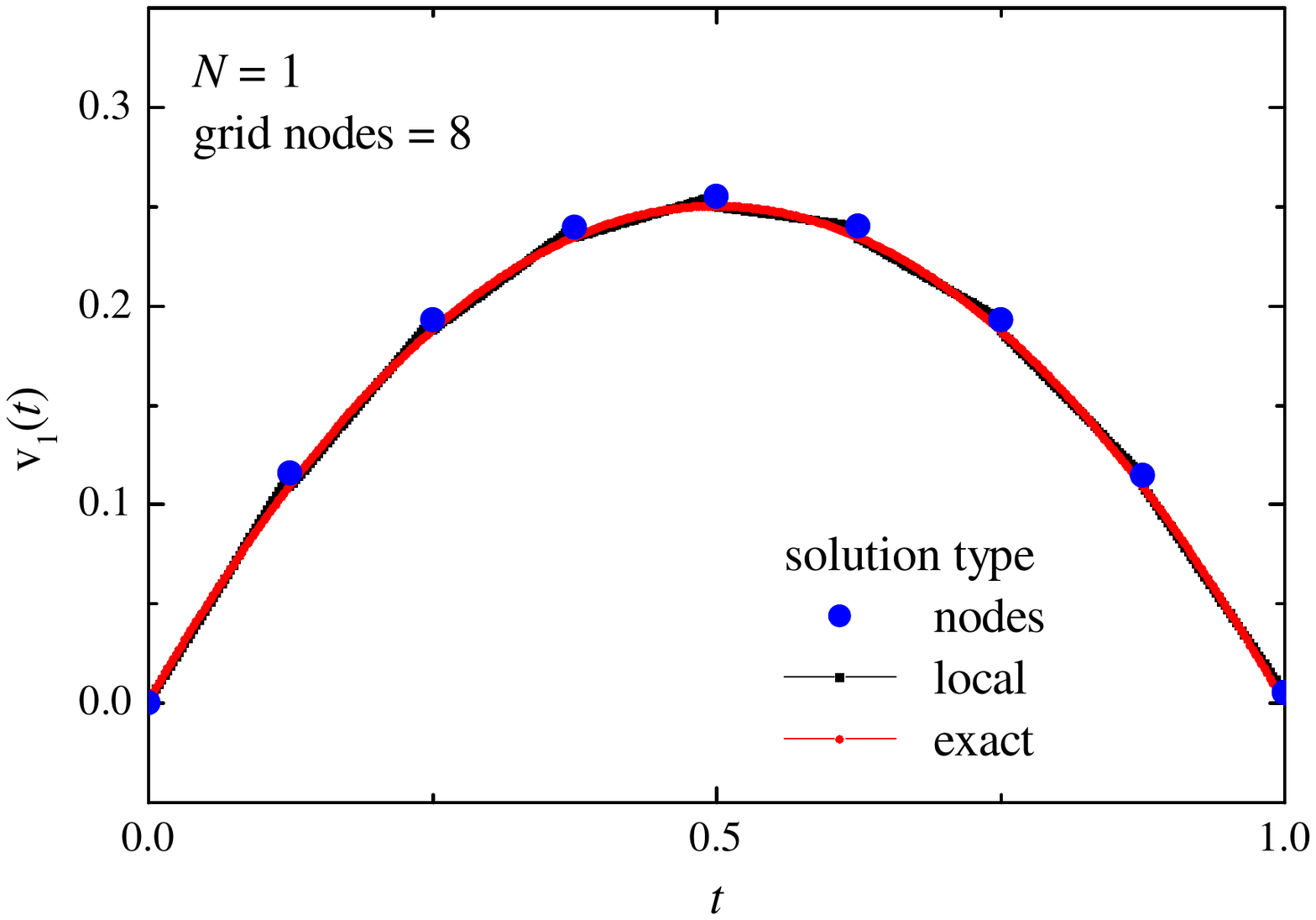}
\vspace{-8mm}\caption{\label{fig:hess_2_ind1_sol_uv:e1}}
\end{subfigure}
\begin{subfigure}{0.320\textwidth}
\includegraphics[width=\textwidth]{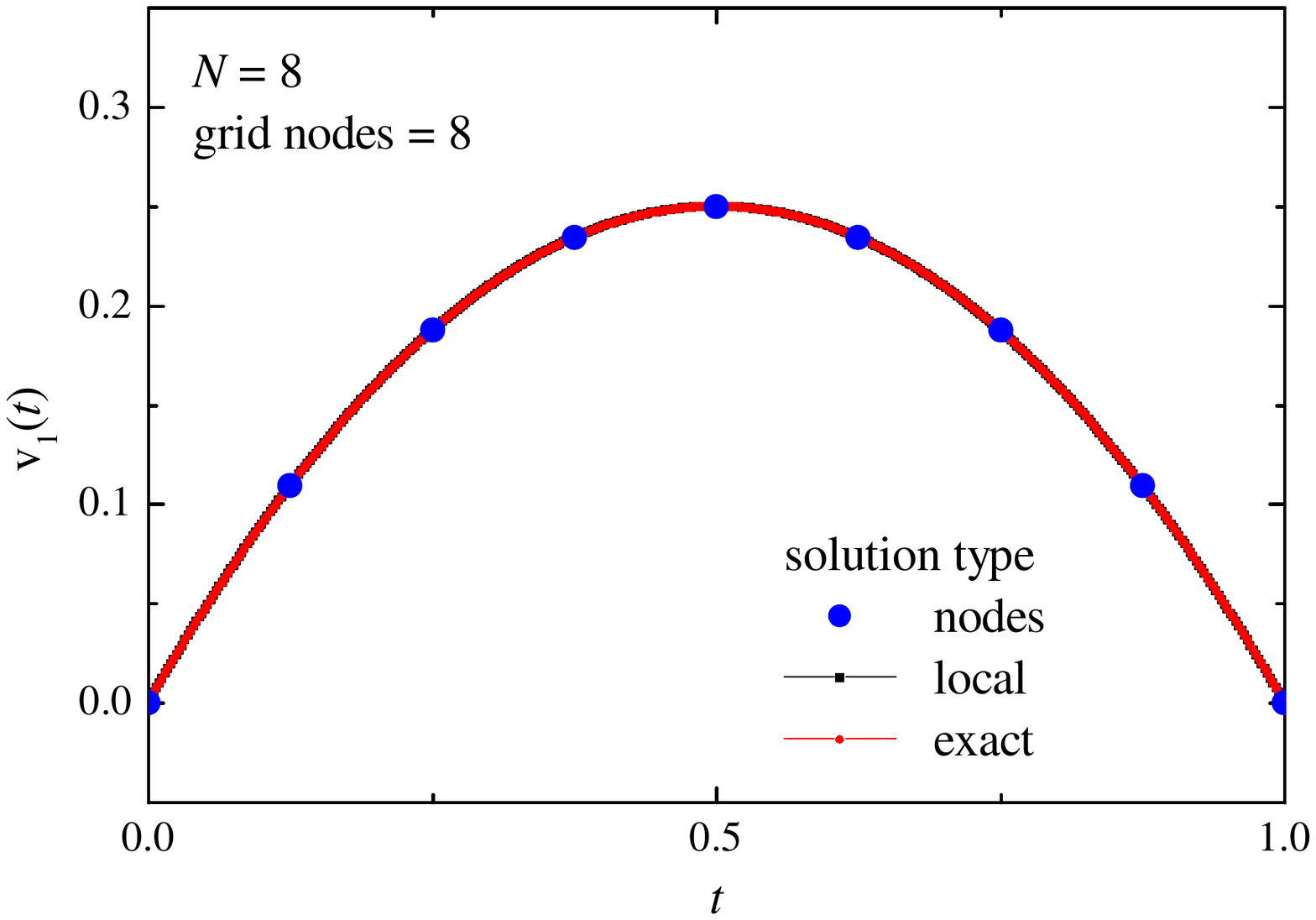}
\vspace{-8mm}\caption{\label{fig:hess_2_ind1_sol_uv:e2}}
\end{subfigure}
\begin{subfigure}{0.320\textwidth}
\includegraphics[width=\textwidth]{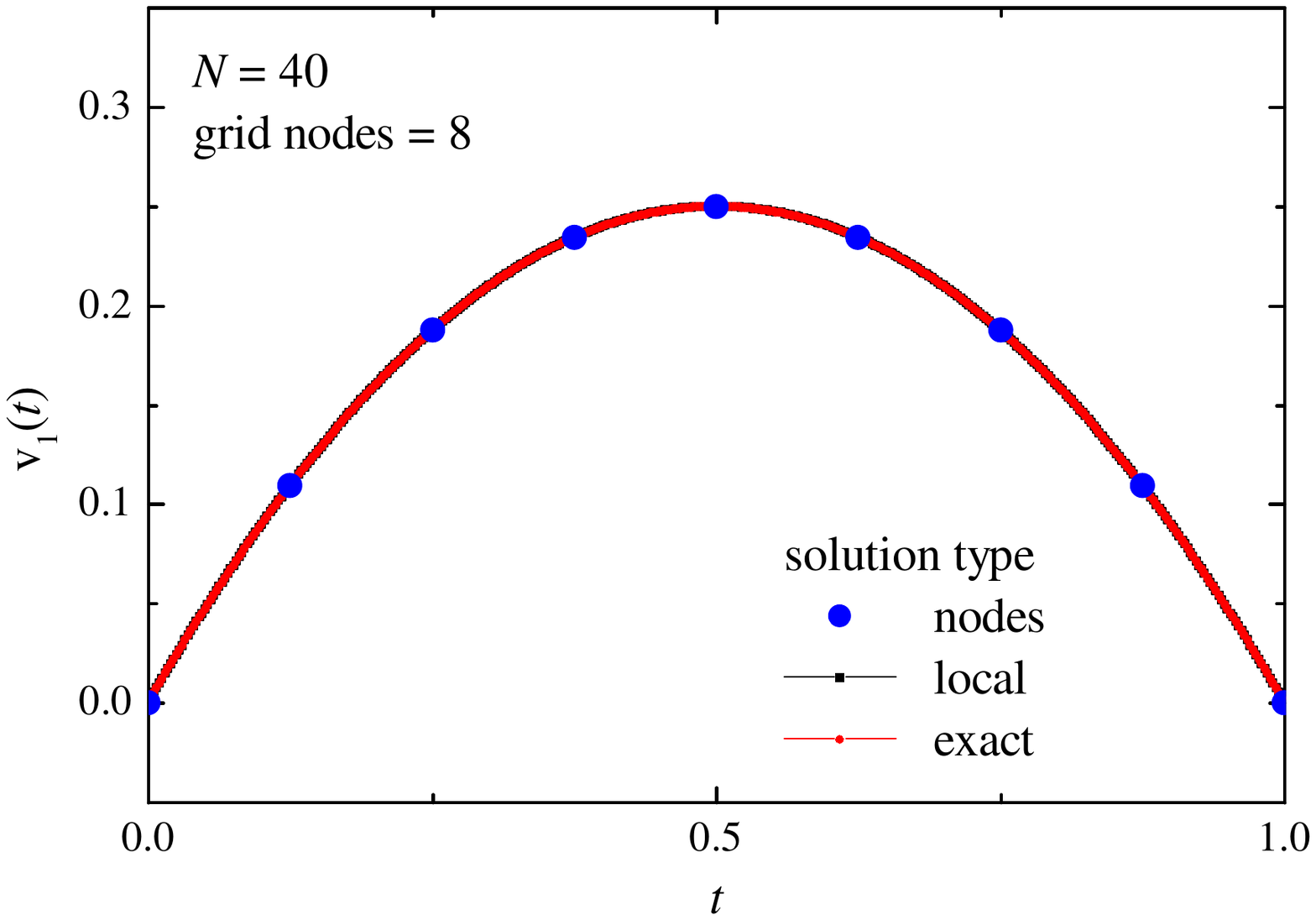}
\vspace{-8mm}\caption{\label{fig:hess_2_ind1_sol_uv:e3}}
\end{subfigure}\\
\caption{%
Numerical solution of the DAE system (\ref{eq:hess_dae_ind_2}) of index 1. Comparison of the solution at nodes $\mathbf{u}_{n}$, the local solution $\mathbf{u}_{L}(t)$ and the exact solution $\mathbf{u}^{\rm ex}(t)$ for components $u_{1}$ (\subref{fig:hess_2_ind1_sol_uv:a1}, \subref{fig:hess_2_ind1_sol_uv:a2}, \subref{fig:hess_2_ind1_sol_uv:a3}), $u_{2}$ (\subref{fig:hess_2_ind1_sol_uv:b1}, \subref{fig:hess_2_ind1_sol_uv:b2}, \subref{fig:hess_2_ind1_sol_uv:b3}), $u_{3}$ (\subref{fig:hess_2_ind1_sol_uv:c1}, \subref{fig:hess_2_ind1_sol_uv:c2}, \subref{fig:hess_2_ind1_sol_uv:c3}), $u_{4}$ (\subref{fig:hess_2_ind1_sol_uv:d1}, \subref{fig:hess_2_ind1_sol_uv:d2}, \subref{fig:hess_2_ind1_sol_uv:d3}) and $v_{1}$ (\subref{fig:hess_2_ind1_sol_uv:e1}, \subref{fig:hess_2_ind1_sol_uv:e2}, \subref{fig:hess_2_ind1_sol_uv:e3}), obtained using polynomials with degrees $N = 1$ (\subref{fig:hess_2_ind1_sol_uv:a1}, \subref{fig:hess_2_ind1_sol_uv:b1}, \subref{fig:hess_2_ind1_sol_uv:c1}, \subref{fig:hess_2_ind1_sol_uv:d1}, \subref{fig:hess_2_ind1_sol_uv:e1}), $N = 8$ (\subref{fig:hess_2_ind1_sol_uv:a2}, \subref{fig:hess_2_ind1_sol_uv:b2}, \subref{fig:hess_2_ind1_sol_uv:c2}, \subref{fig:hess_2_ind1_sol_uv:d2}, \subref{fig:hess_2_ind1_sol_uv:e2}), $N = 40$ (\subref{fig:hess_2_ind1_sol_uv:a3}, \subref{fig:hess_2_ind1_sol_uv:b3}, \subref{fig:hess_2_ind1_sol_uv:c3}, \subref{fig:hess_2_ind1_sol_uv:d3}, \subref{fig:hess_2_ind1_sol_uv:e3}).
}
\label{fig:hess_2_ind1_sol_uv}
\end{figure} 

\begin{figure}[h!]
\captionsetup[subfigure]{%
	position=bottom,
	font+=smaller,
	textfont=normalfont,
	singlelinecheck=off,
	justification=raggedright
}
\centering
\begin{subfigure}{0.320\textwidth}
\includegraphics[width=\textwidth]{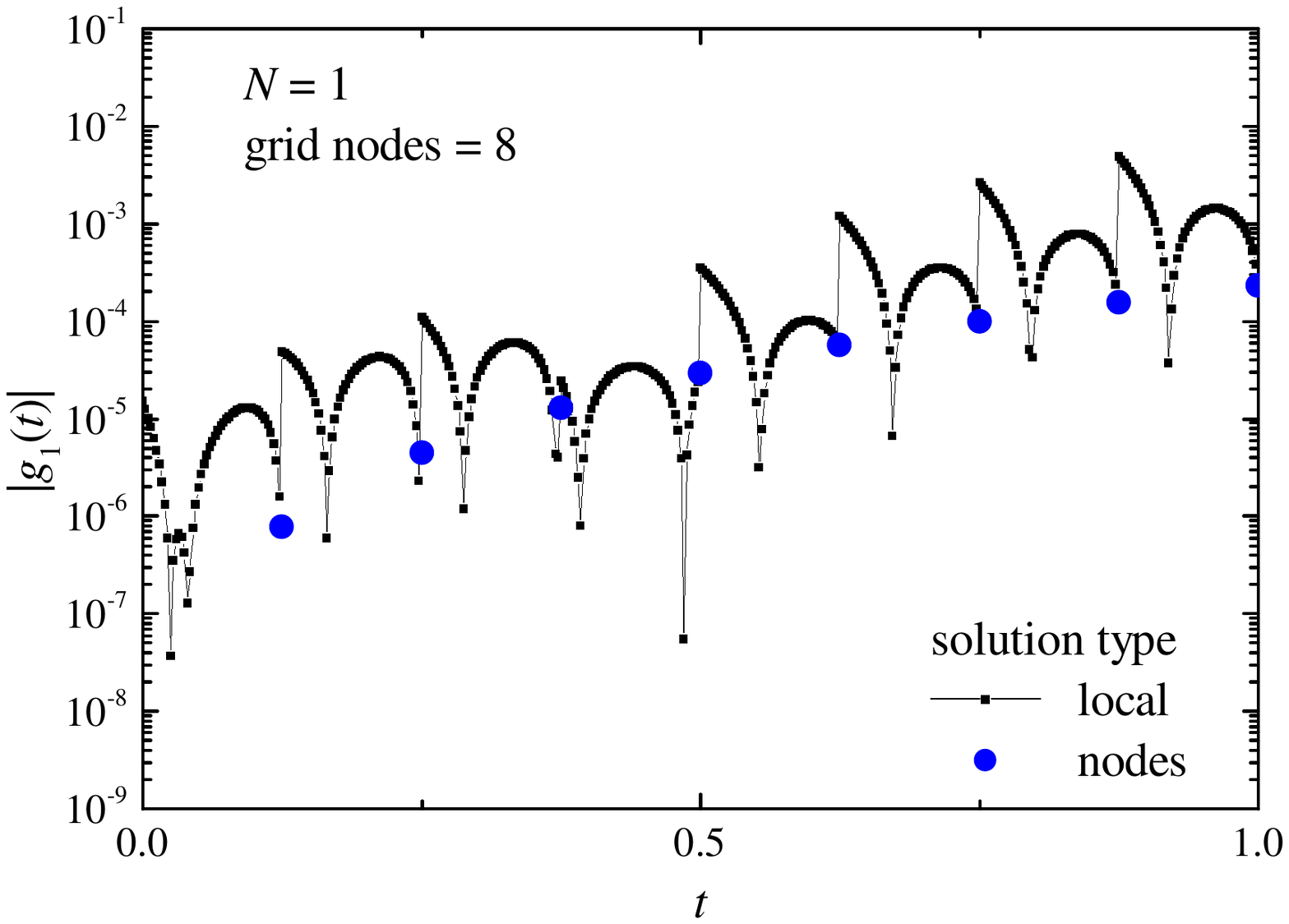}
\vspace{-8mm}\caption{\label{fig:hess_2_ind1_sol_g_eps:a1}}
\end{subfigure}
\begin{subfigure}{0.320\textwidth}
\includegraphics[width=\textwidth]{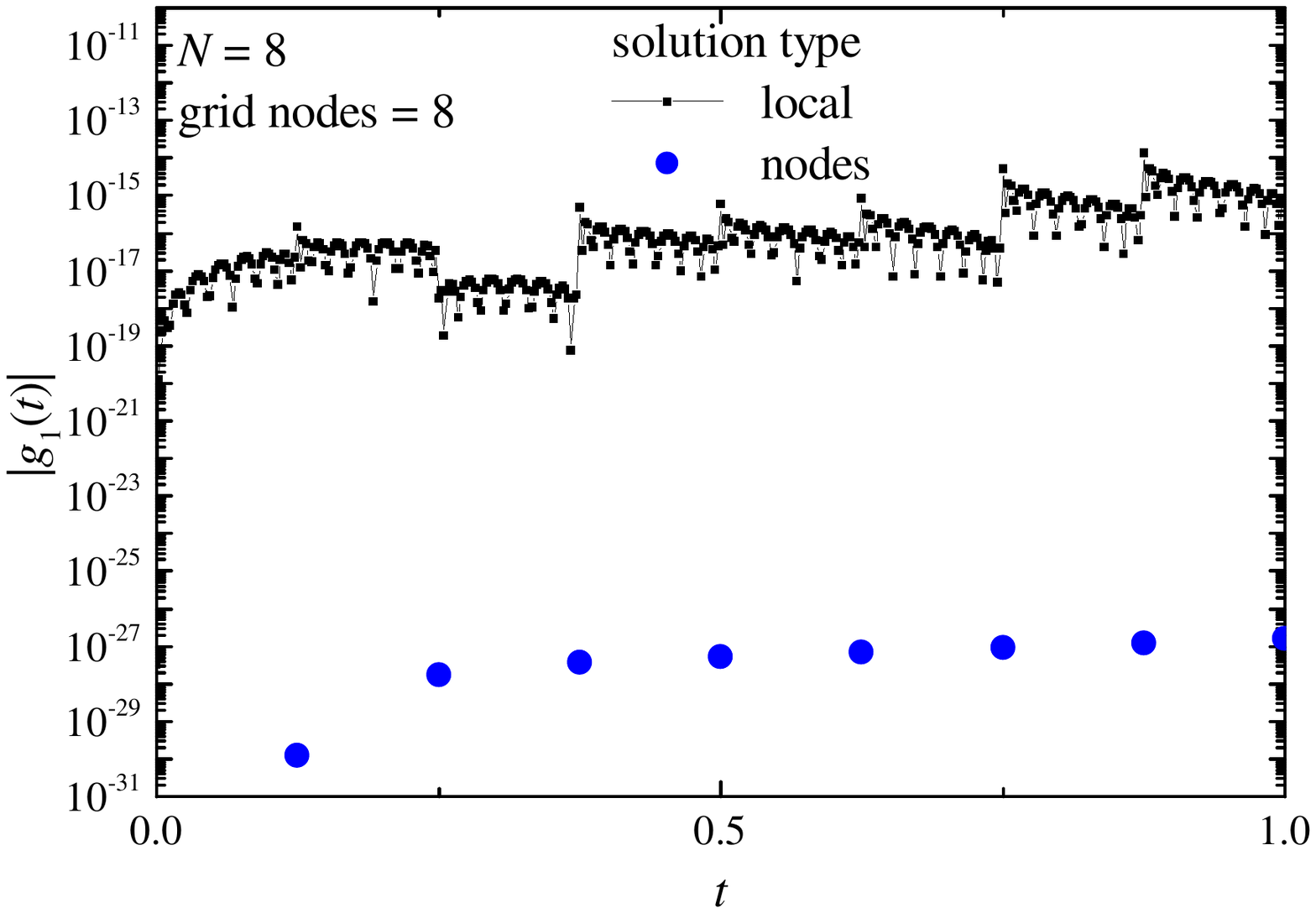}
\vspace{-8mm}\caption{\label{fig:hess_2_ind1_sol_g_eps:a2}}
\end{subfigure}
\begin{subfigure}{0.320\textwidth}
\includegraphics[width=\textwidth]{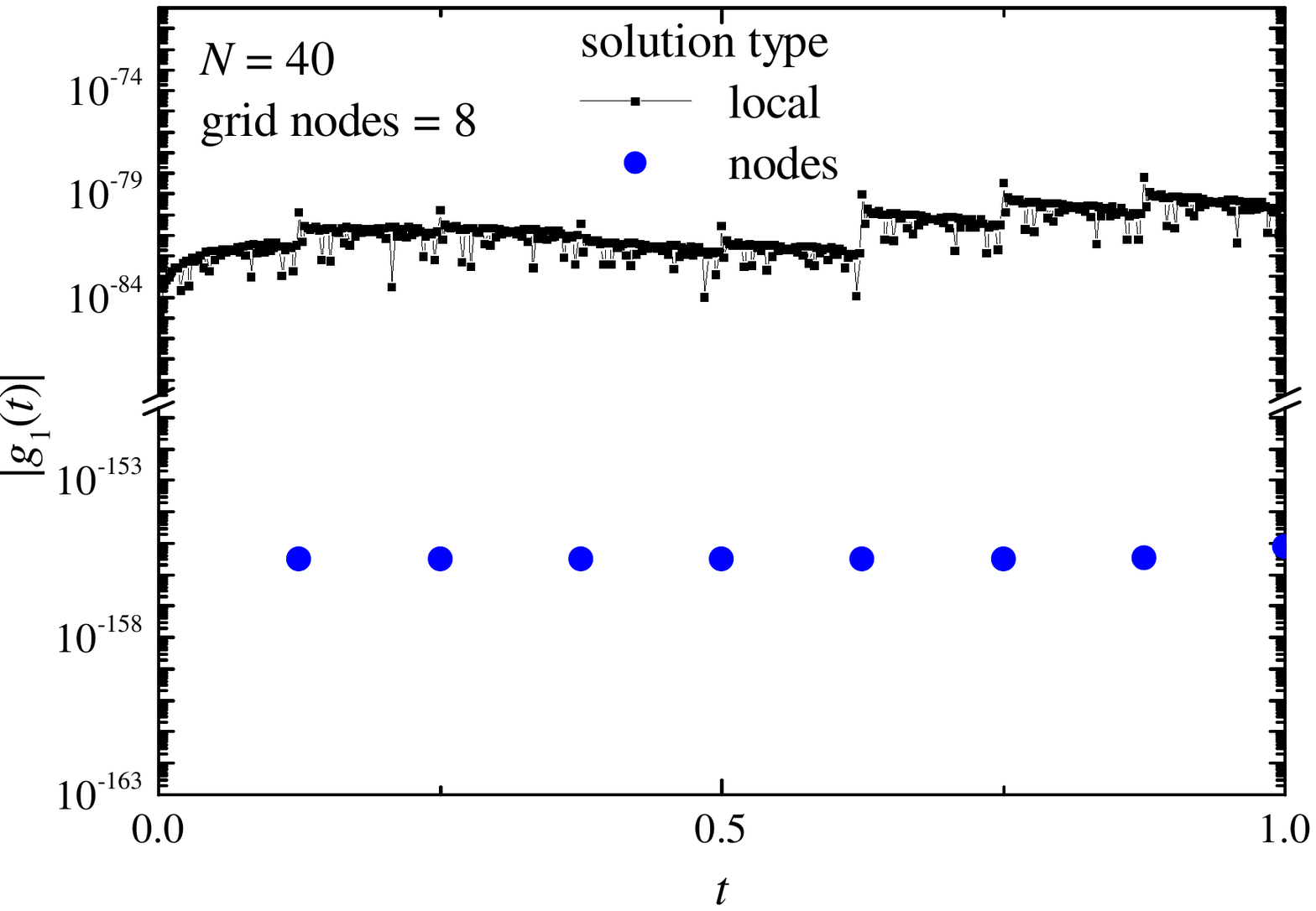}
\vspace{-8mm}\caption{\label{fig:hess_2_ind1_sol_g_eps:a3}}
\end{subfigure}\\
\begin{subfigure}{0.320\textwidth}
\includegraphics[width=\textwidth]{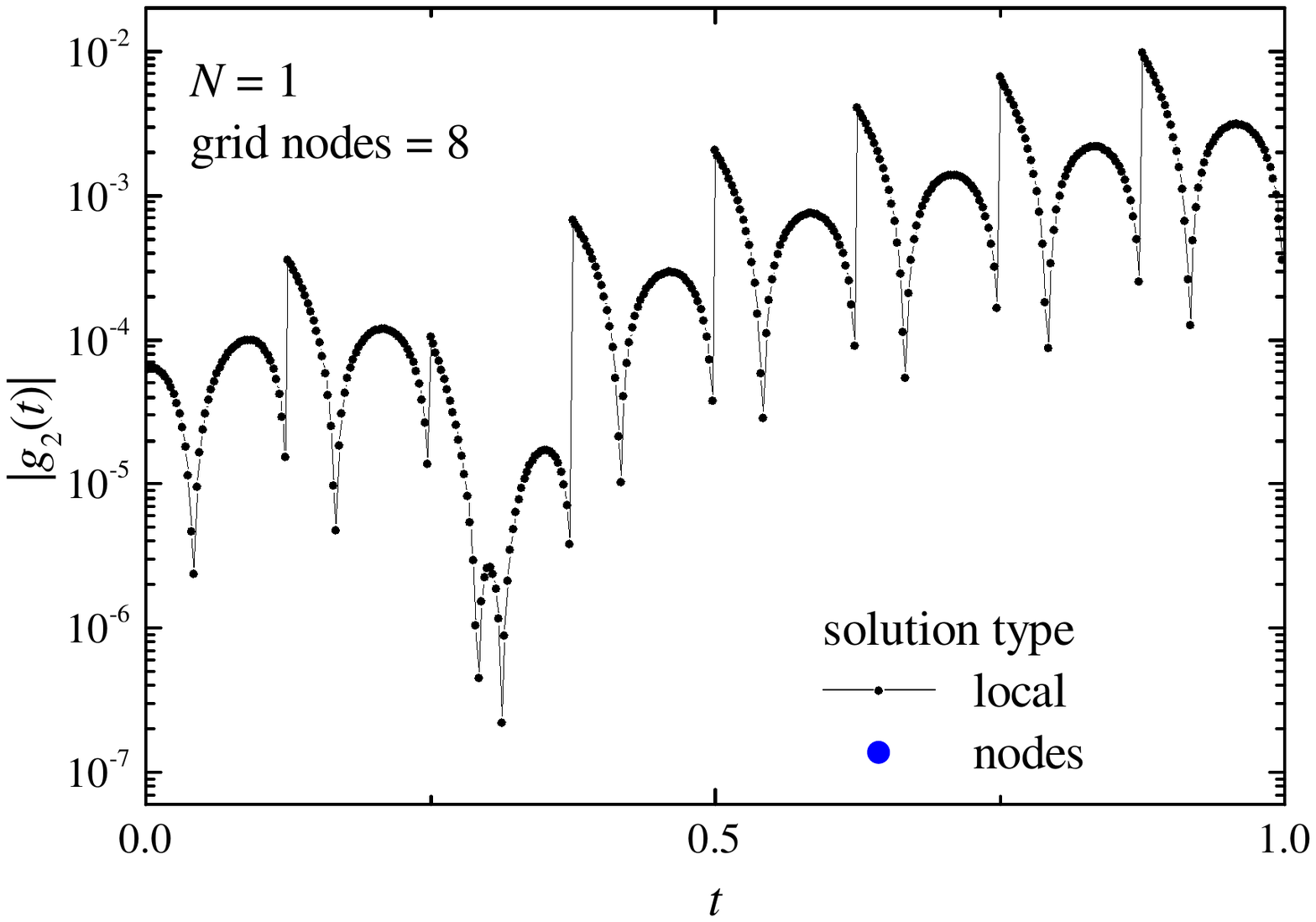}
\vspace{-8mm}\caption{\label{fig:hess_2_ind1_sol_g_eps:b1}}
\end{subfigure}
\begin{subfigure}{0.320\textwidth}
\includegraphics[width=\textwidth]{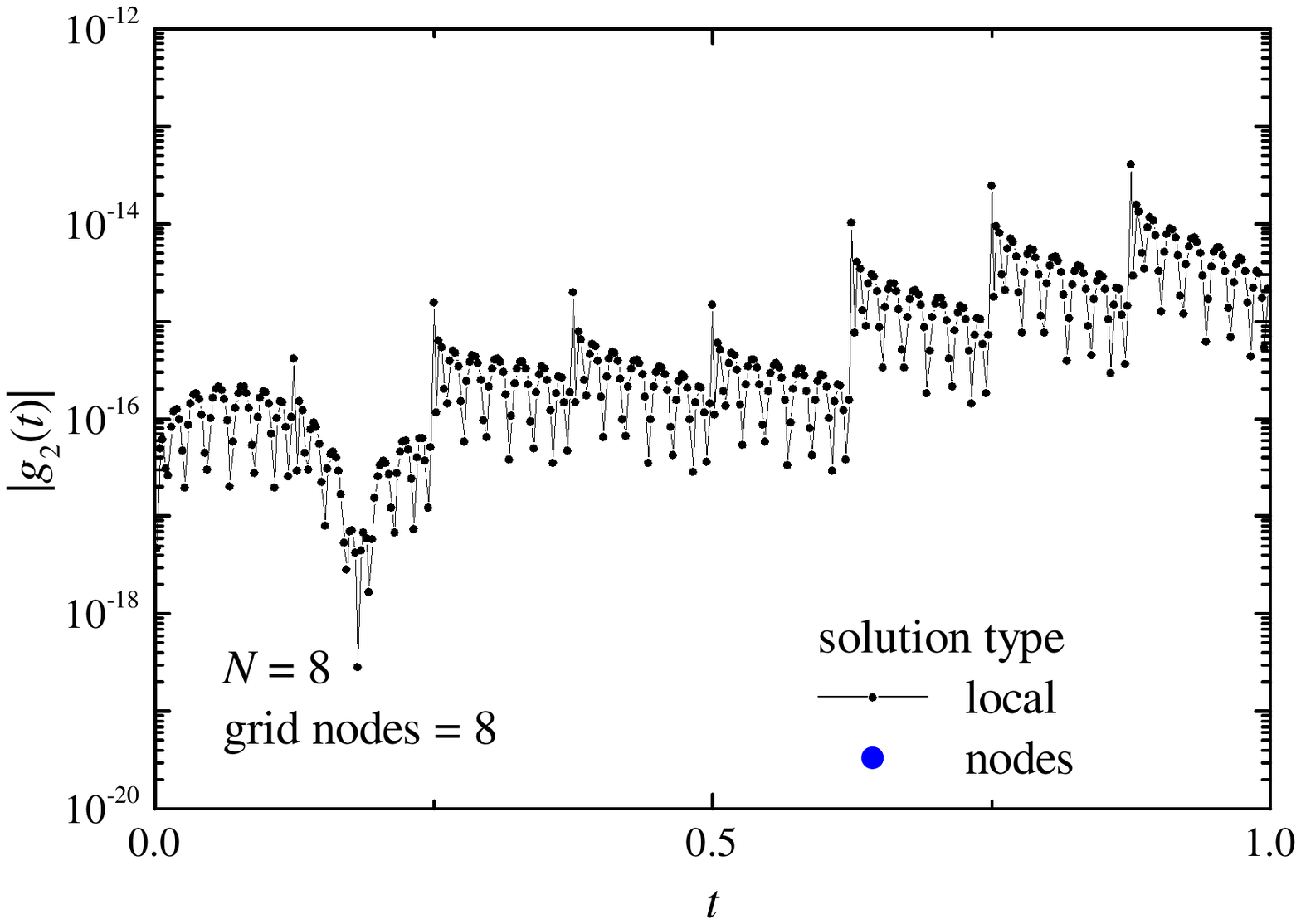}
\vspace{-8mm}\caption{\label{fig:hess_2_ind1_sol_g_eps:b2}}
\end{subfigure}
\begin{subfigure}{0.320\textwidth}
\includegraphics[width=\textwidth]{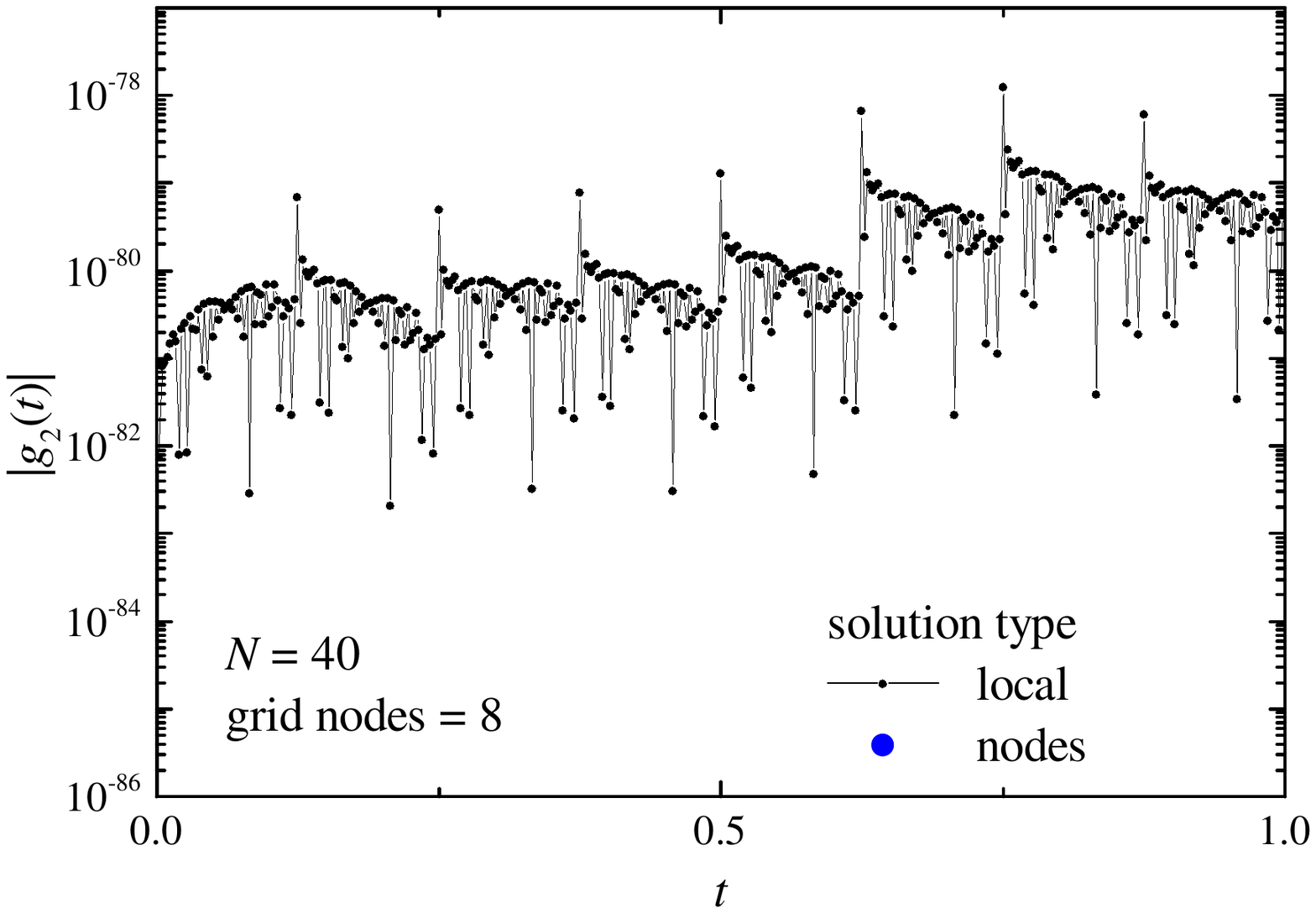}
\vspace{-8mm}\caption{\label{fig:hess_2_ind1_sol_g_eps:b3}}
\end{subfigure}\\
\begin{subfigure}{0.320\textwidth}
\includegraphics[width=\textwidth]{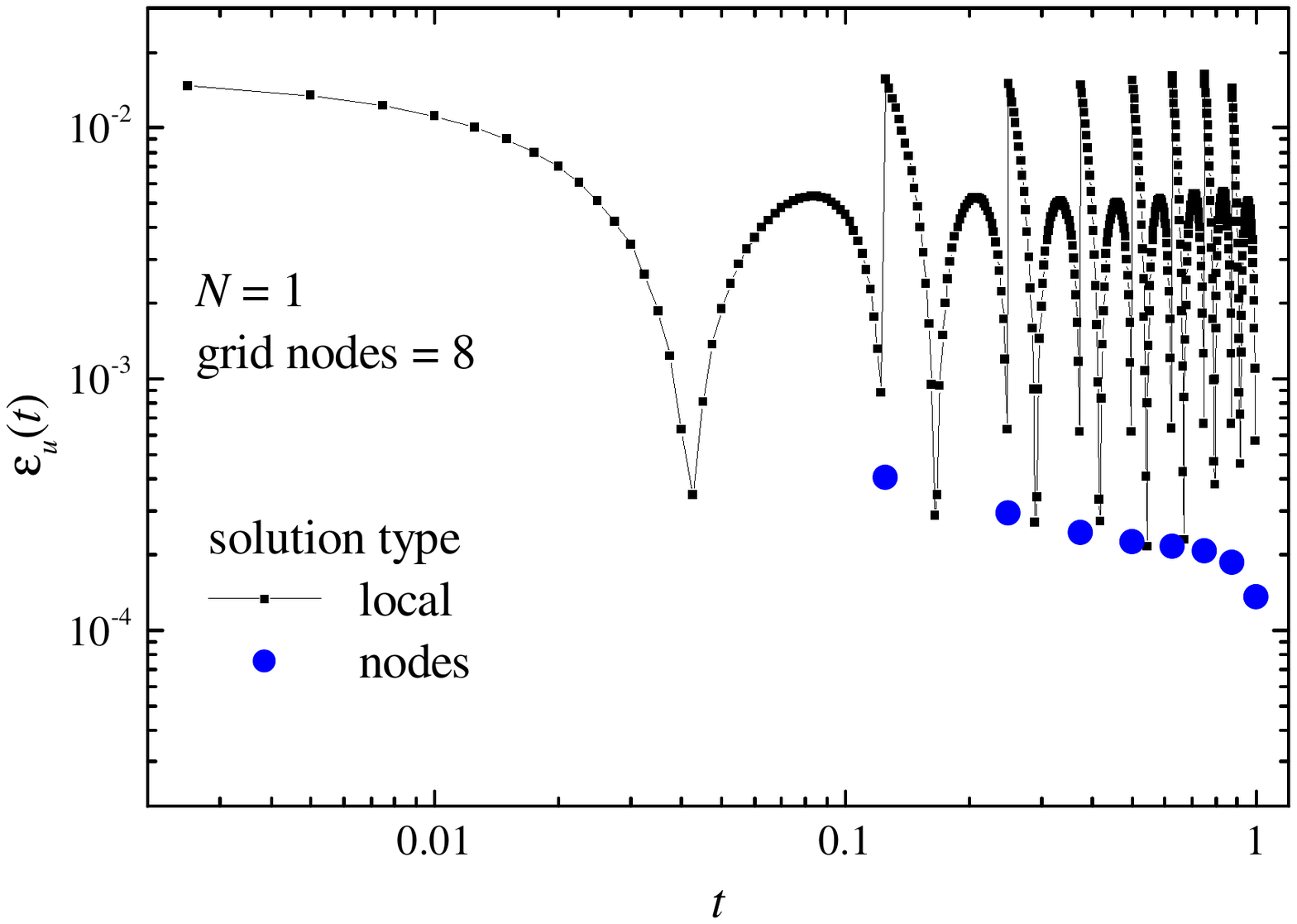}
\vspace{-8mm}\caption{\label{fig:hess_2_ind1_sol_g_eps:c1}}
\end{subfigure}
\begin{subfigure}{0.320\textwidth}
\includegraphics[width=\textwidth]{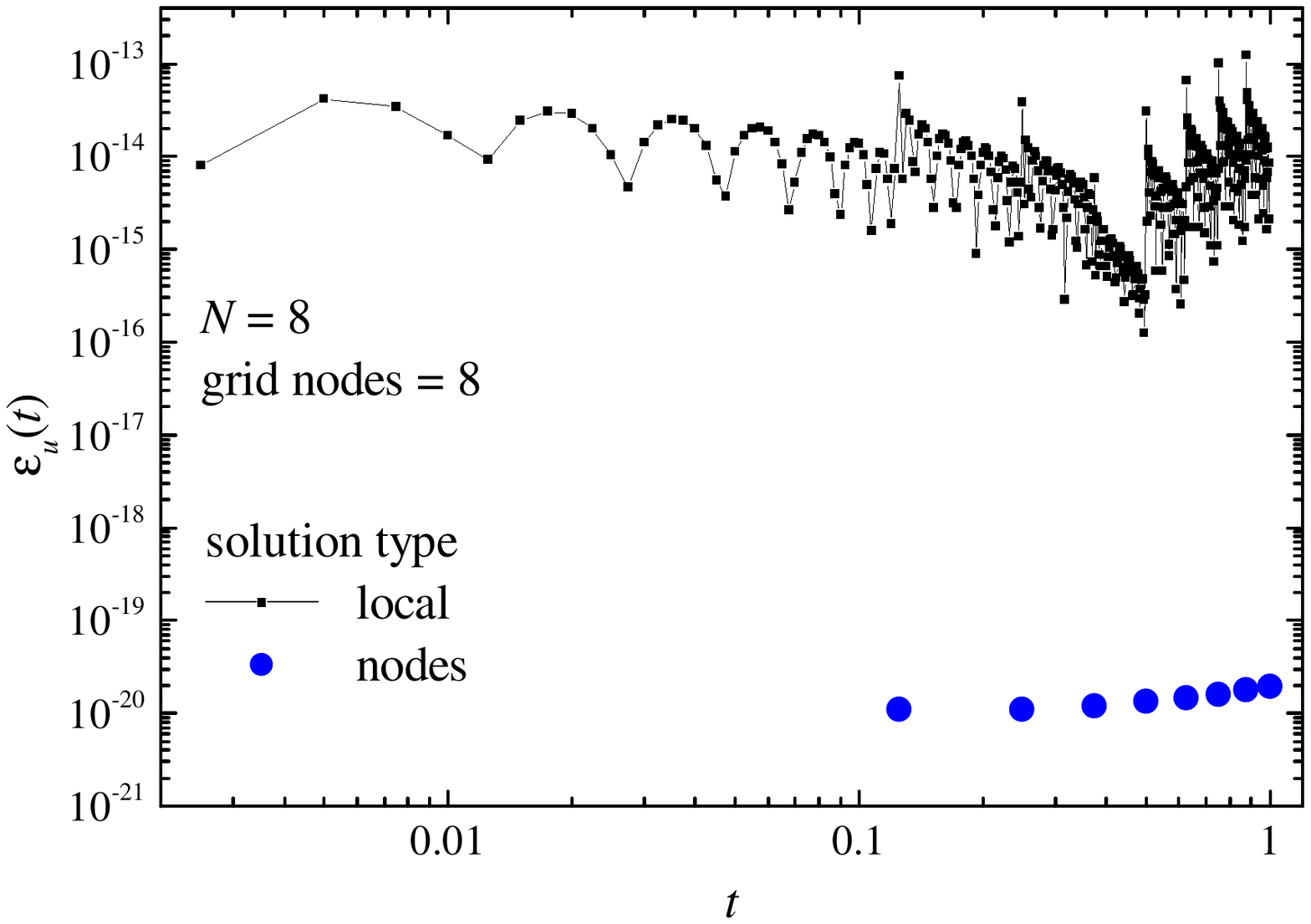}
\vspace{-8mm}\caption{\label{fig:hess_2_ind1_sol_g_eps:c2}}
\end{subfigure}
\begin{subfigure}{0.320\textwidth}
\includegraphics[width=\textwidth]{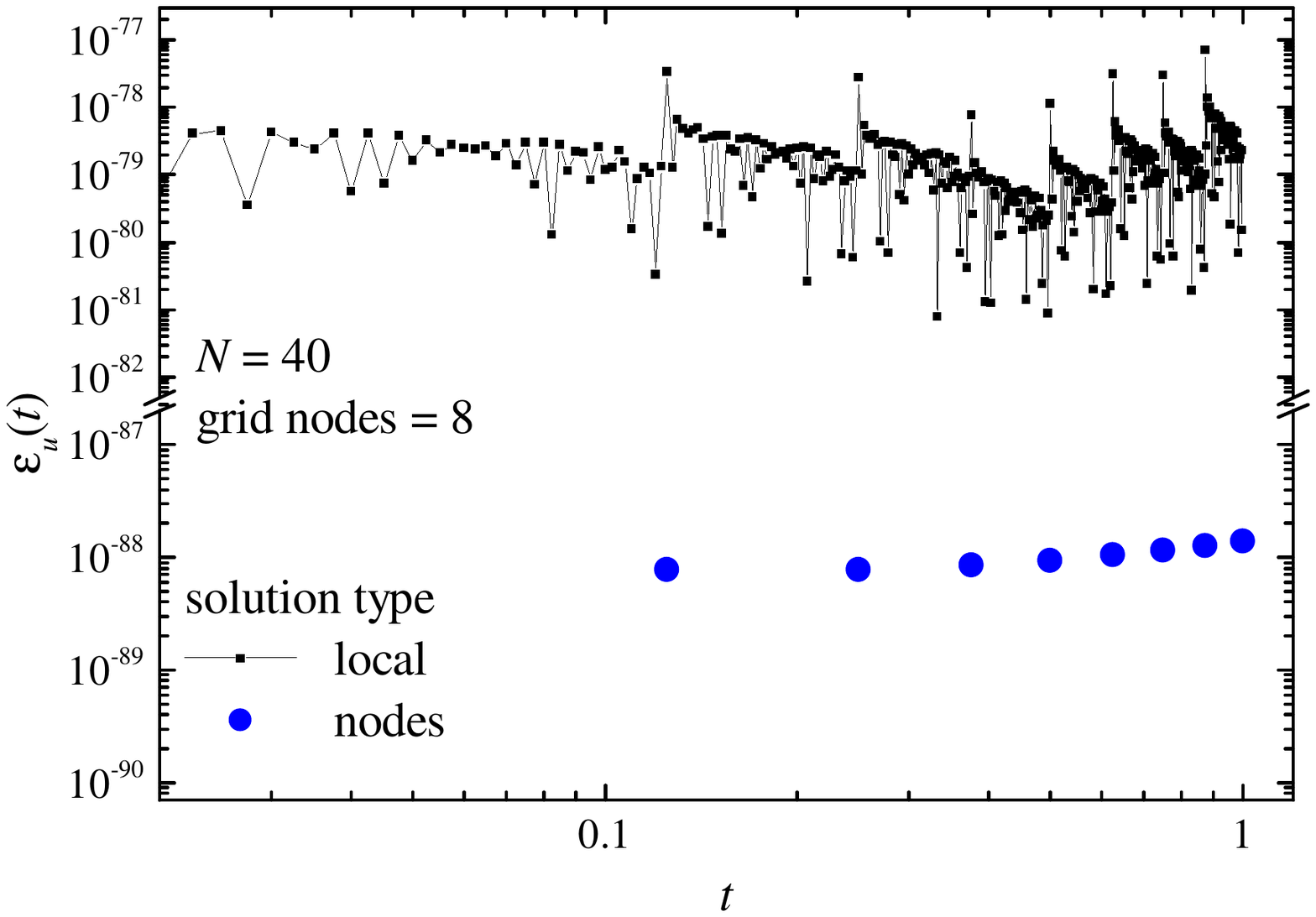}
\vspace{-8mm}\caption{\label{fig:hess_2_ind1_sol_g_eps:c3}}
\end{subfigure}\\
\begin{subfigure}{0.320\textwidth}
\includegraphics[width=\textwidth]{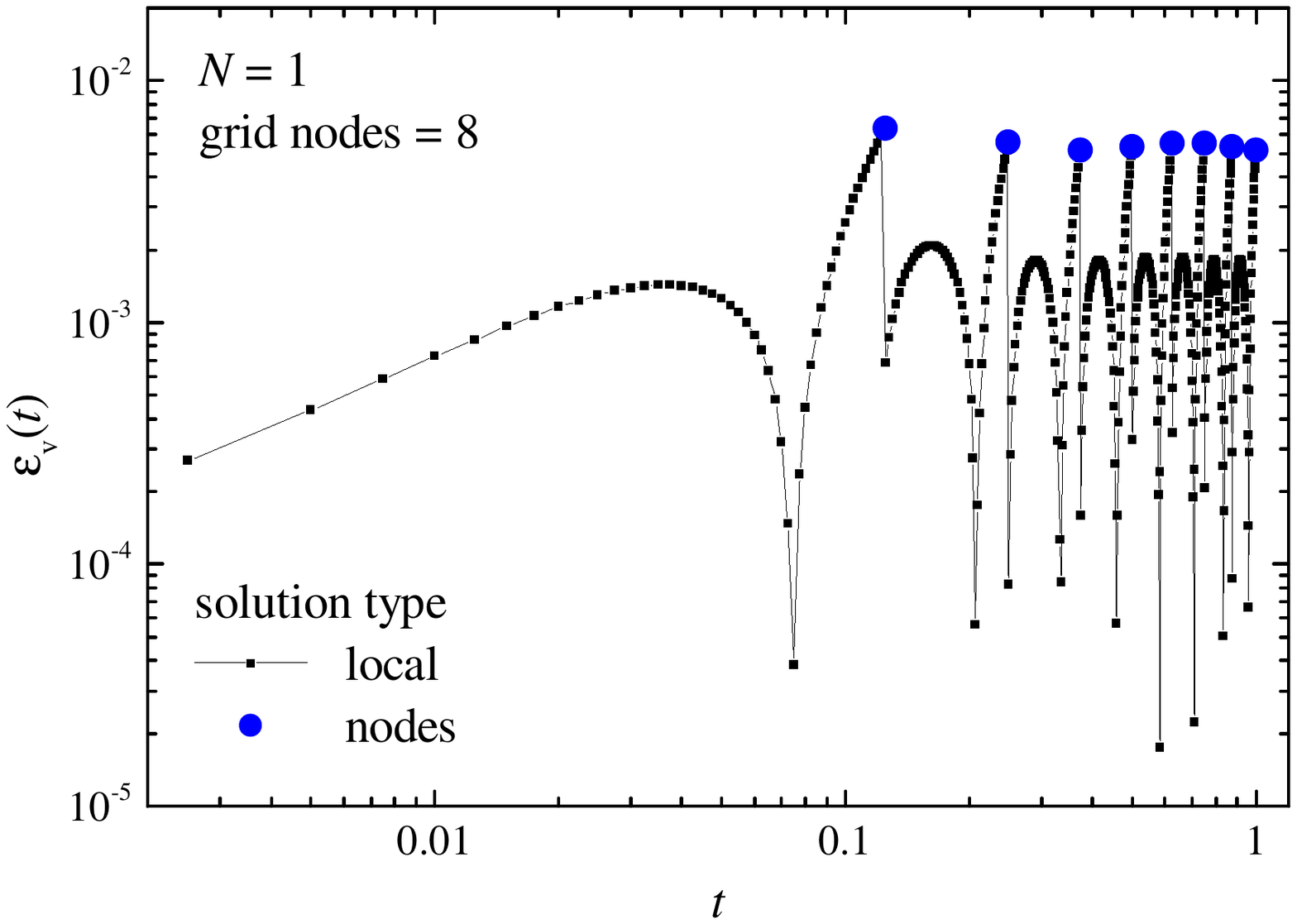}
\vspace{-8mm}\caption{\label{fig:hess_2_ind1_sol_g_eps:d1}}
\end{subfigure}
\begin{subfigure}{0.320\textwidth}
\includegraphics[width=\textwidth]{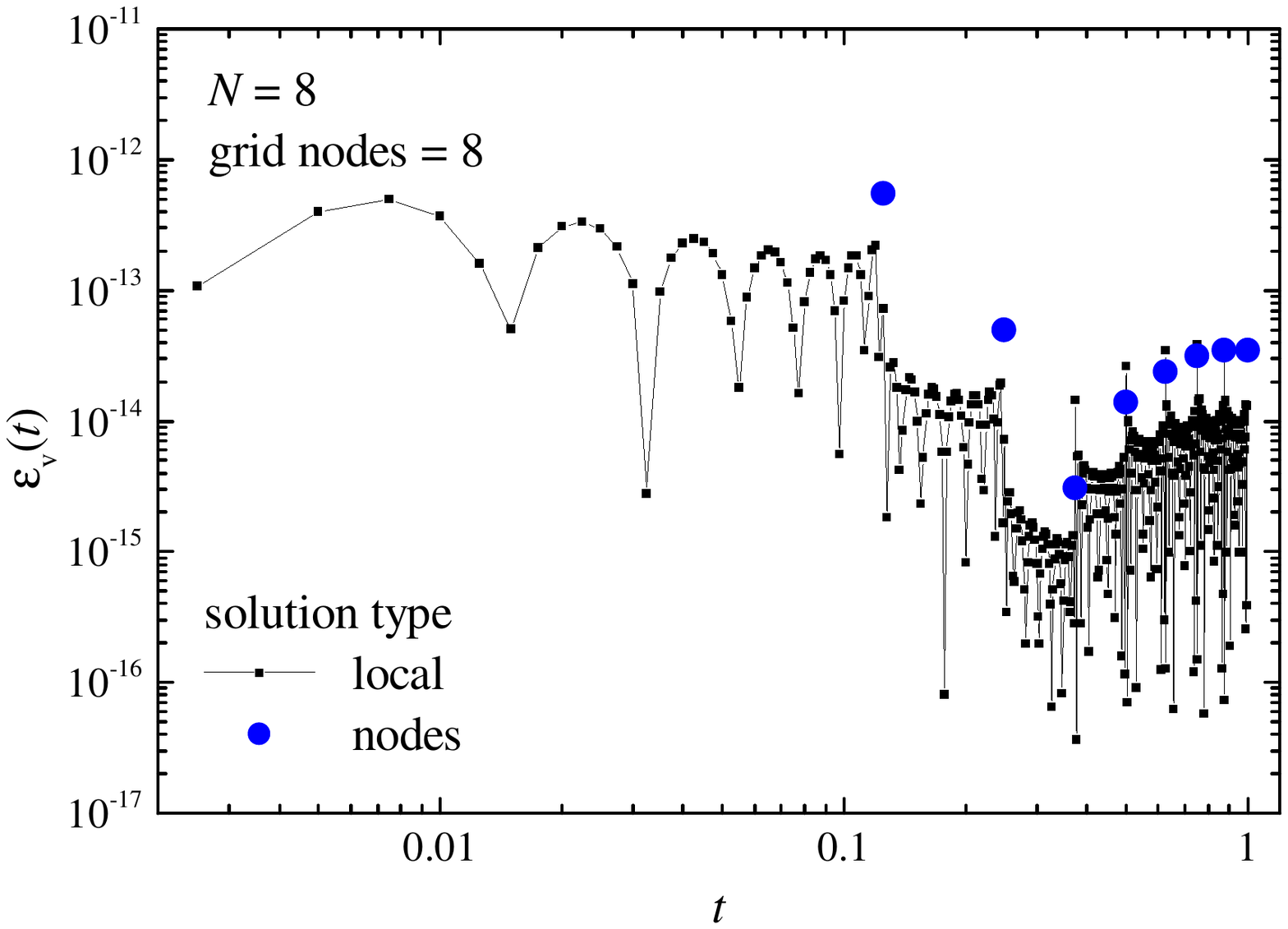}
\vspace{-8mm}\caption{\label{fig:hess_2_ind1_sol_g_eps:d2}}
\end{subfigure}
\begin{subfigure}{0.320\textwidth}
\includegraphics[width=\textwidth]{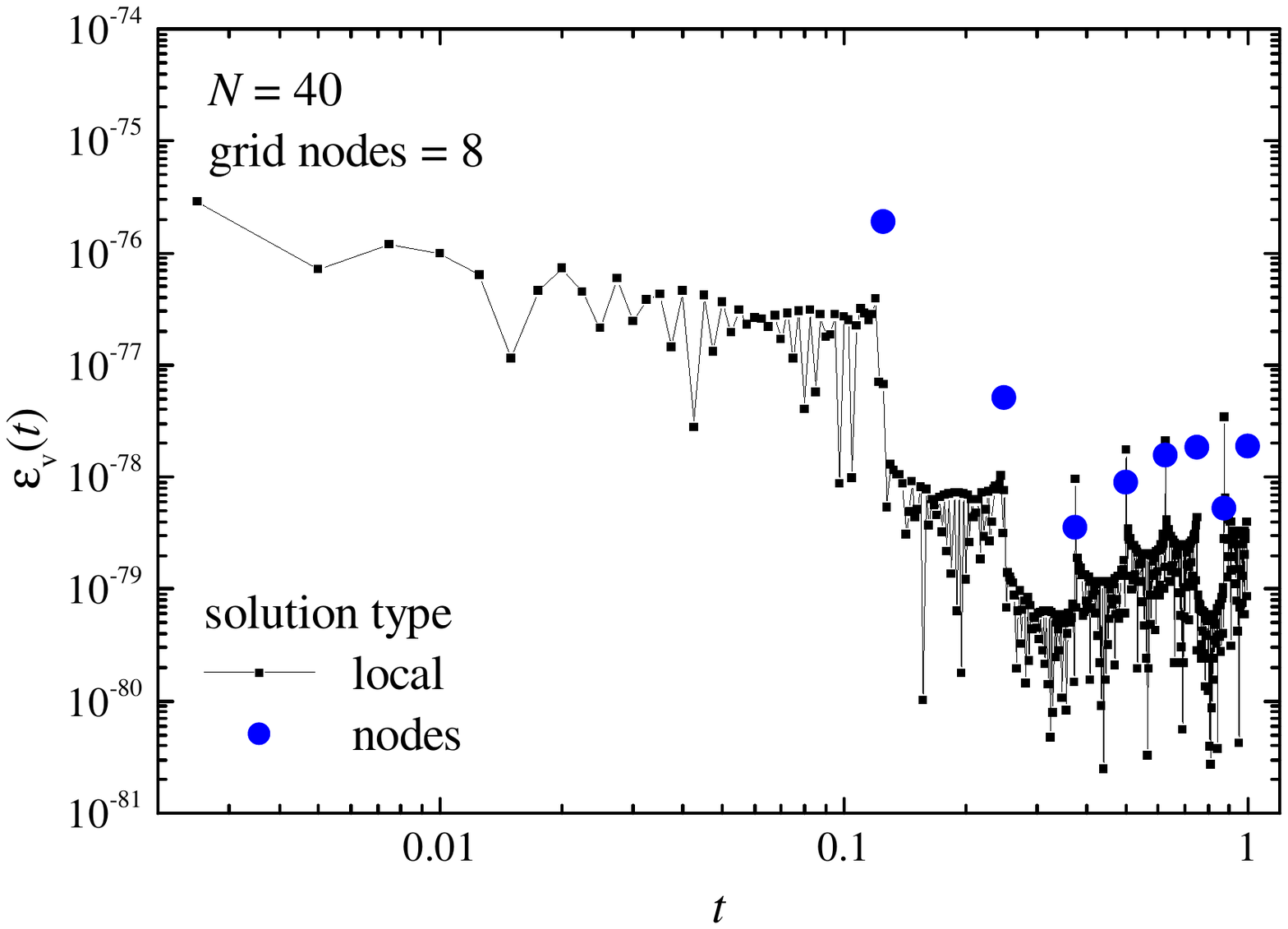}
\vspace{-8mm}\caption{\label{fig:hess_2_ind1_sol_g_eps:d3}}
\end{subfigure}\\
\begin{subfigure}{0.320\textwidth}
\includegraphics[width=\textwidth]{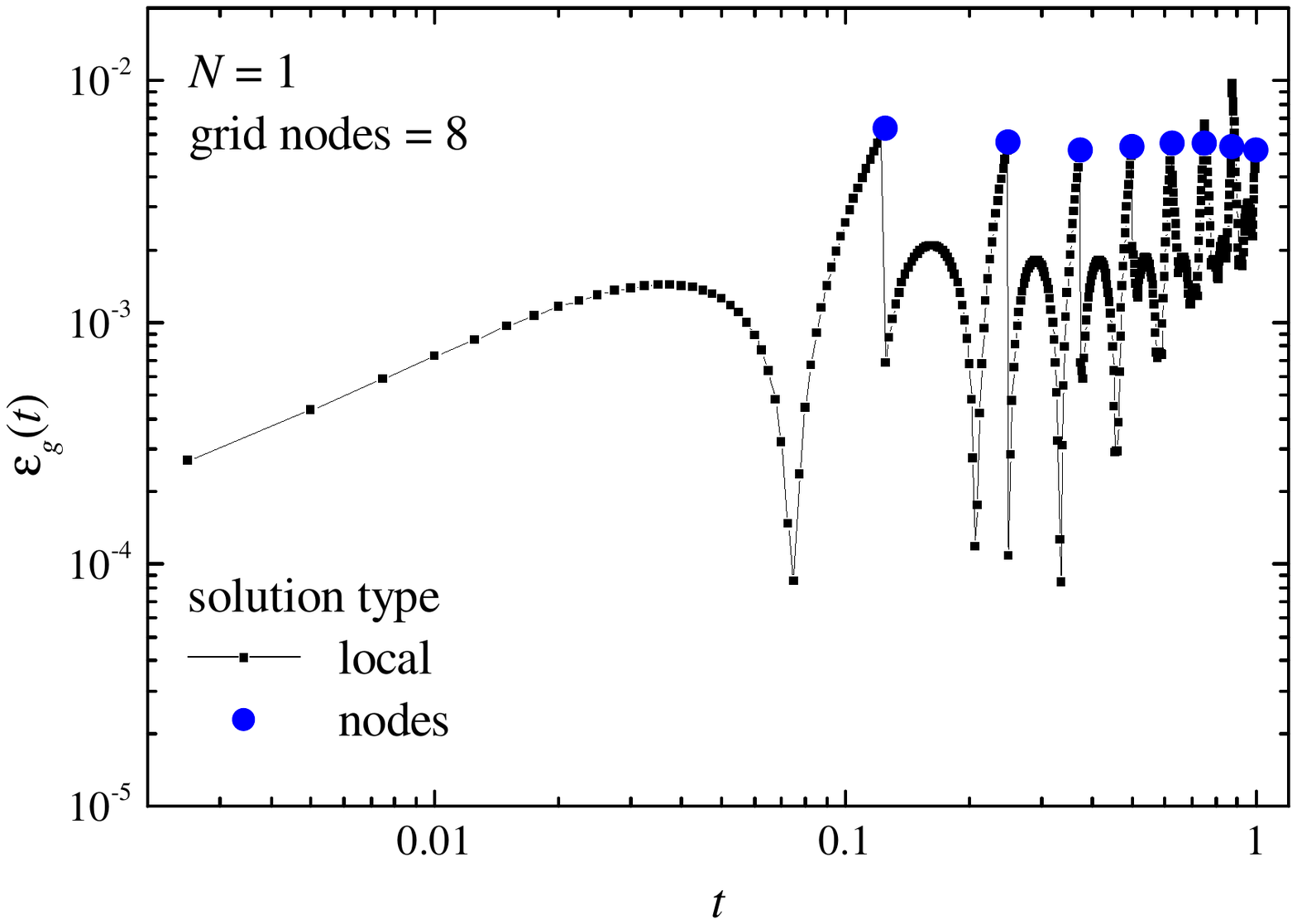}
\vspace{-8mm}\caption{\label{fig:hess_2_ind1_sol_g_eps:e1}}
\end{subfigure}
\begin{subfigure}{0.320\textwidth}
\includegraphics[width=\textwidth]{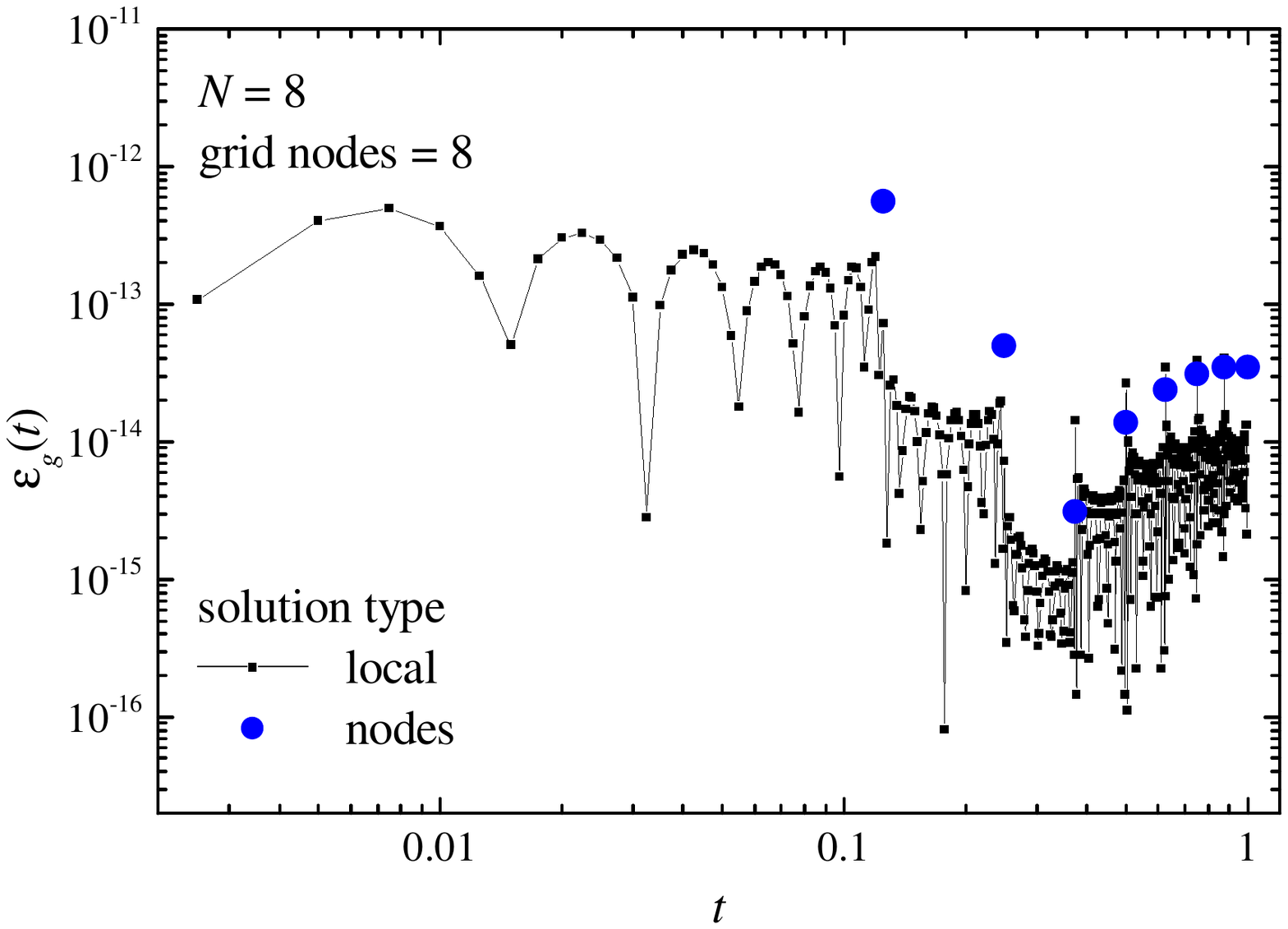}
\vspace{-8mm}\caption{\label{fig:hess_2_ind1_sol_g_eps:e2}}
\end{subfigure}
\begin{subfigure}{0.320\textwidth}
\includegraphics[width=\textwidth]{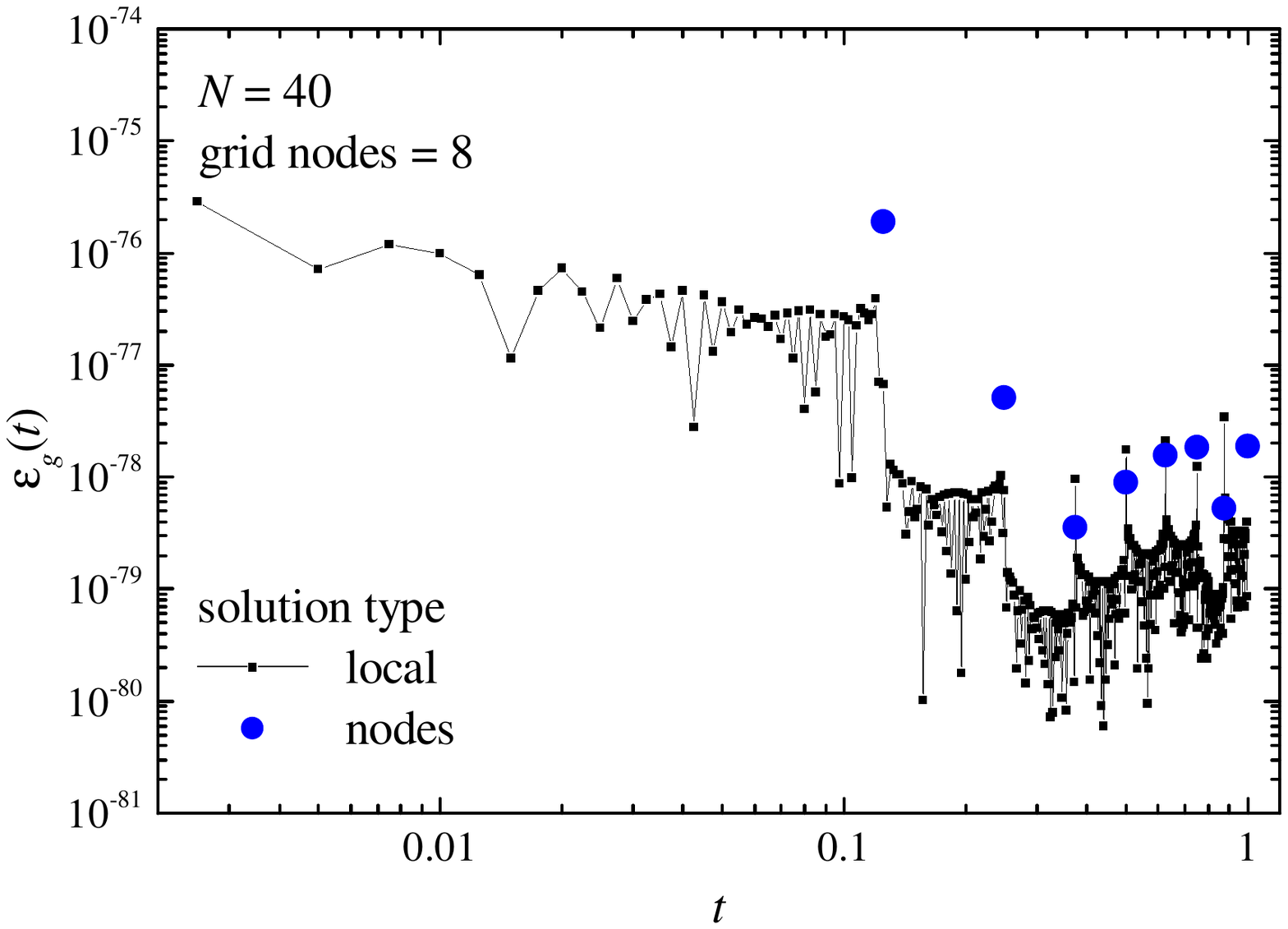}
\vspace{-8mm}\caption{\label{fig:hess_2_ind1_sol_g_eps:e3}}
\end{subfigure}\\
\caption{%
Numerical solution of the DAE system (\ref{eq:hess_dae_ind_2}) of index 1. Comparison of quantitative satisfiability of the conditions $g_{1} = 0$ (\subref{fig:hess_2_ind1_sol_g_eps:a1}, \subref{fig:hess_2_ind1_sol_g_eps:a2}, \subref{fig:hess_2_ind1_sol_g_eps:a3}) and $g_{2} = 0$ (\subref{fig:hess_2_ind1_sol_g_eps:b1}, \subref{fig:hess_2_ind1_sol_g_eps:b2}, \subref{fig:hess_2_ind1_sol_g_eps:b3}), the errors $\varepsilon_{u}(t)$ (\subref{fig:hess_2_ind1_sol_g_eps:c1}, \subref{fig:hess_2_ind1_sol_g_eps:c2}, \subref{fig:hess_2_ind1_sol_g_eps:c3}), $\varepsilon_{v}(t)$ (\subref{fig:hess_2_ind1_sol_g_eps:d1}, \subref{fig:hess_2_ind1_sol_g_eps:d2}, \subref{fig:hess_2_ind1_sol_g_eps:d3}), $\varepsilon_{g}(t)$ (\subref{fig:hess_2_ind1_sol_g_eps:e1}, \subref{fig:hess_2_ind1_sol_g_eps:e2}, \subref{fig:hess_2_ind1_sol_g_eps:e3}), for numerical solution obtained using polynomials with degrees $N = 1$ (\subref{fig:hess_2_ind1_sol_g_eps:a1}, \subref{fig:hess_2_ind1_sol_g_eps:b1}, \subref{fig:hess_2_ind1_sol_g_eps:c1}, \subref{fig:hess_2_ind1_sol_g_eps:d1}, \subref{fig:hess_2_ind1_sol_g_eps:e1}), $N = 8$ (\subref{fig:hess_2_ind1_sol_g_eps:a2}, \subref{fig:hess_2_ind1_sol_g_eps:b2}, \subref{fig:hess_2_ind1_sol_g_eps:c2}, \subref{fig:hess_2_ind1_sol_g_eps:d2}, \subref{fig:hess_2_ind1_sol_g_eps:e2}), $N = 40$ (\subref{fig:hess_2_ind1_sol_g_eps:a3}, \subref{fig:hess_2_ind1_sol_g_eps:b3}, \subref{fig:hess_2_ind1_sol_g_eps:c3}, \subref{fig:hess_2_ind1_sol_g_eps:d3}, \subref{fig:hess_2_ind1_sol_g_eps:e3}).
}
\label{fig:hess_2_ind1_sol_g_eps}
\end{figure} 

\begin{figure}[h!]
\captionsetup[subfigure]{%
	position=bottom,
	font+=smaller,
	textfont=normalfont,
	singlelinecheck=off,
	justification=raggedright
}
\centering
\begin{subfigure}{0.275\textwidth}
\includegraphics[width=\textwidth]{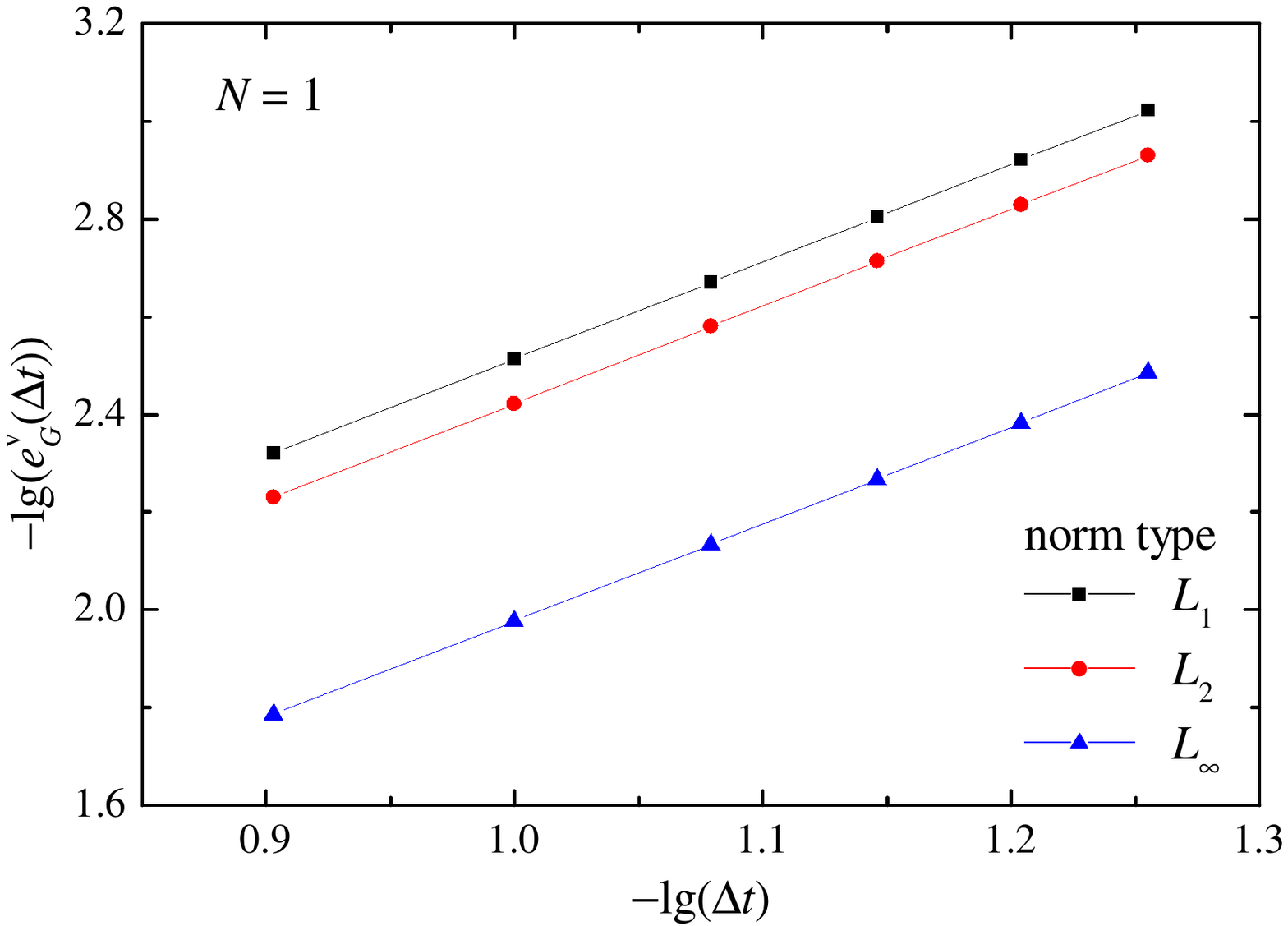}
\vspace{-8mm}\caption{\label{fig:hess_2_ind1_errors:a1}}
\end{subfigure}\hspace{6mm}
\begin{subfigure}{0.275\textwidth}
\includegraphics[width=\textwidth]{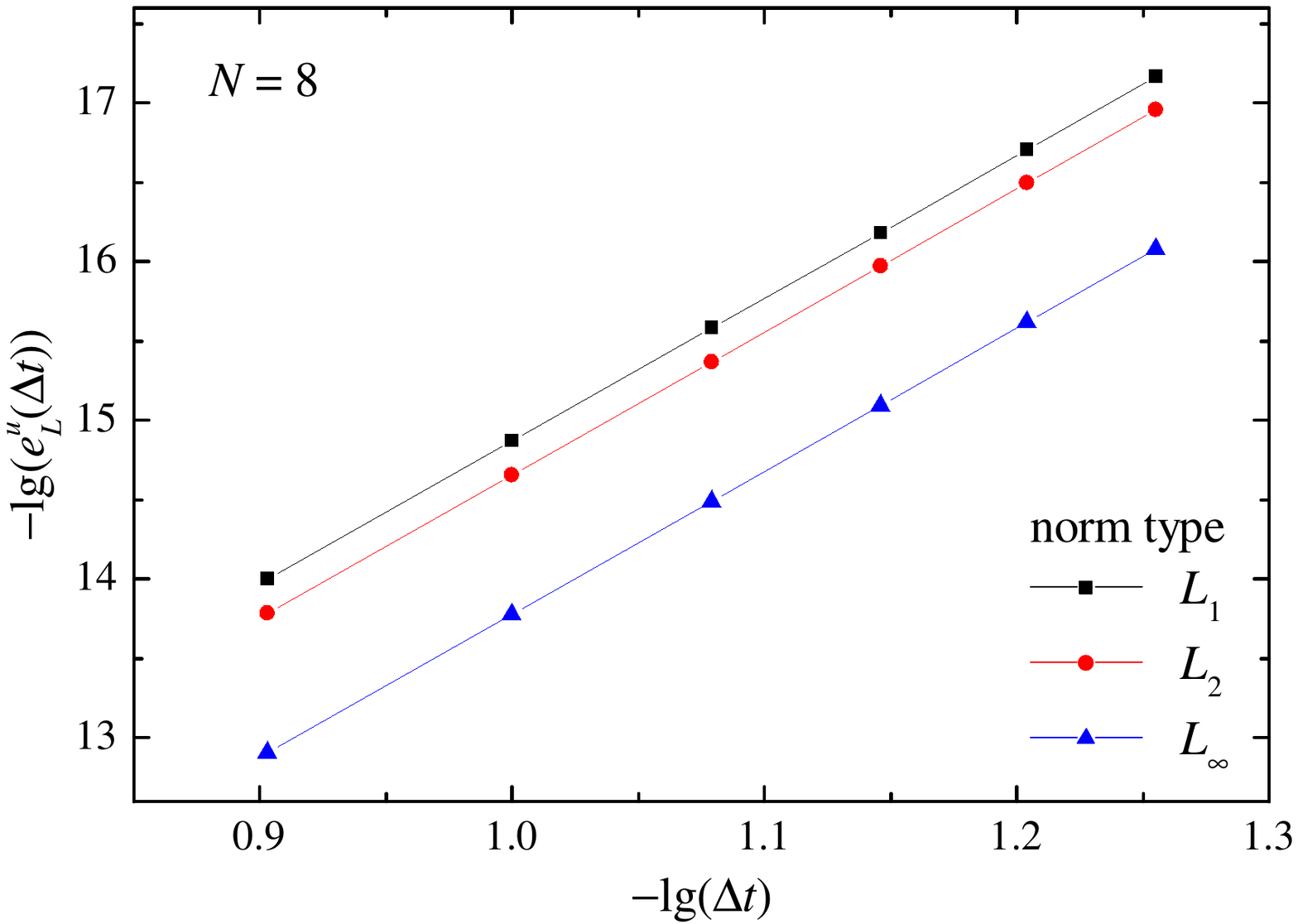}
\vspace{-8mm}\caption{\label{fig:hess_2_ind1_errors:a2}}
\end{subfigure}\hspace{6mm}
\begin{subfigure}{0.275\textwidth}
\includegraphics[width=\textwidth]{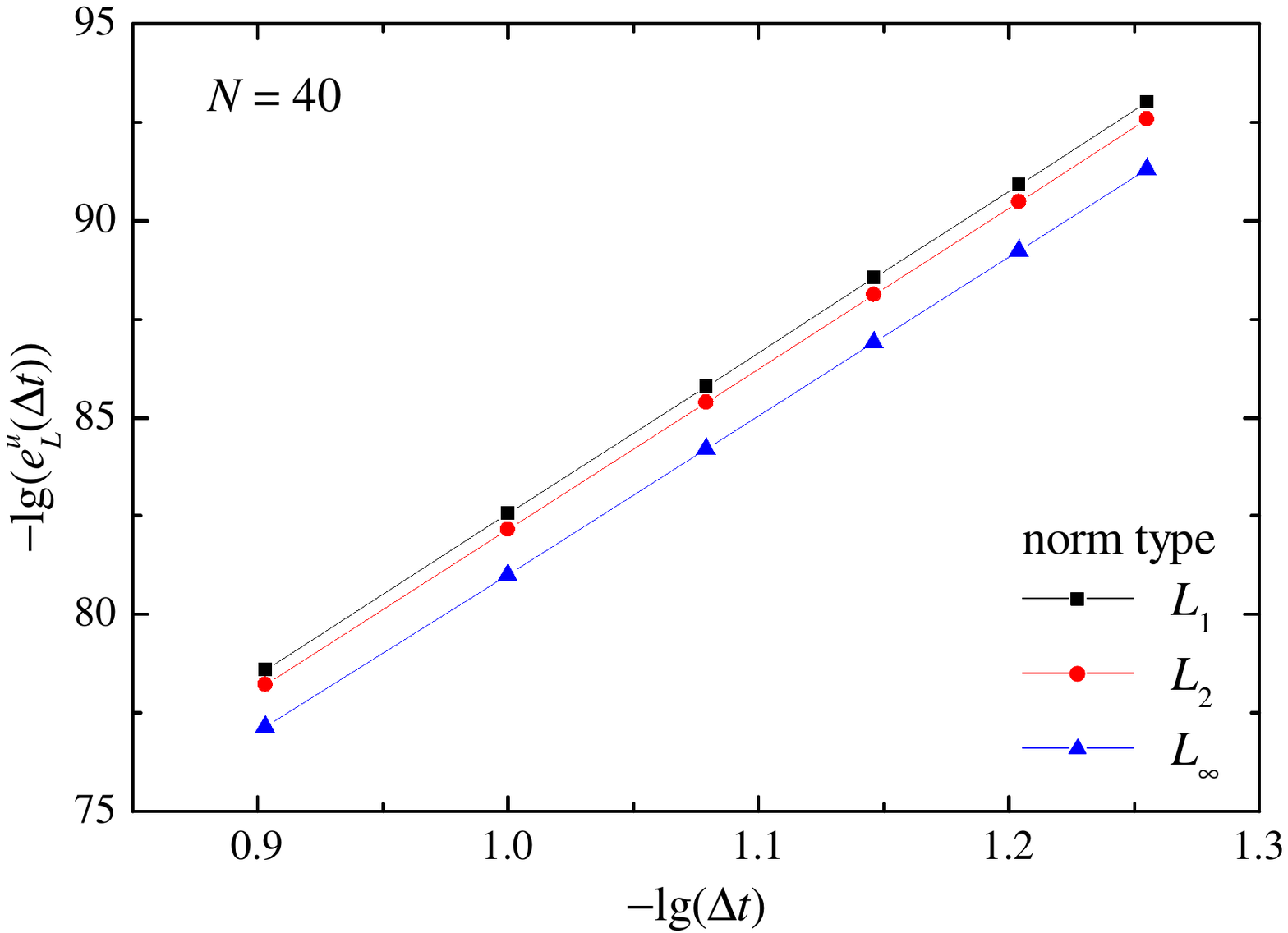}
\vspace{-8mm}\caption{\label{fig:hess_2_ind1_errors:a3}}
\end{subfigure}\\[-2mm]
\begin{subfigure}{0.275\textwidth}
\includegraphics[width=\textwidth]{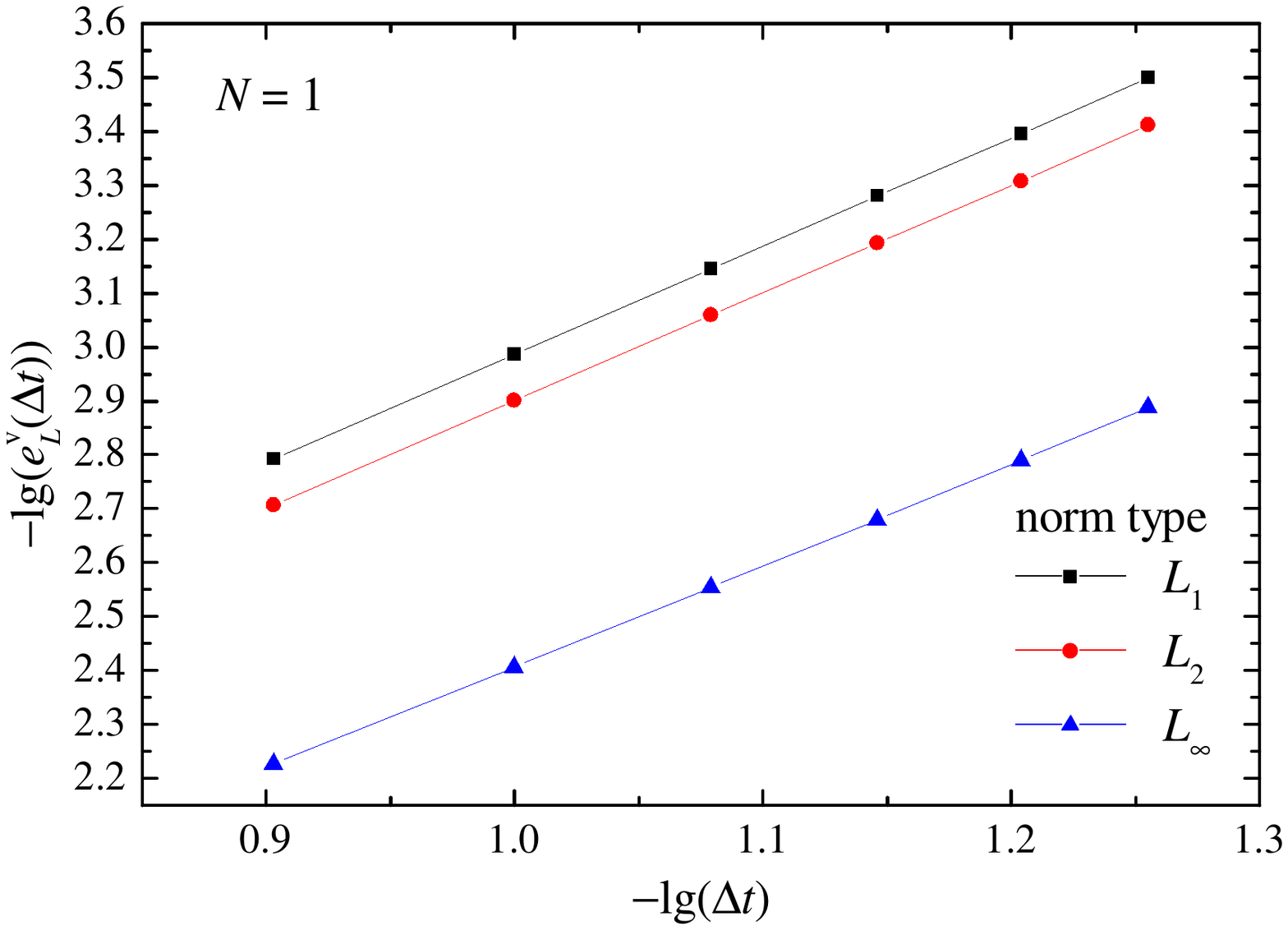}
\vspace{-8mm}\caption{\label{fig:hess_2_ind1_errors:b1}}
\end{subfigure}\hspace{6mm}
\begin{subfigure}{0.275\textwidth}
\includegraphics[width=\textwidth]{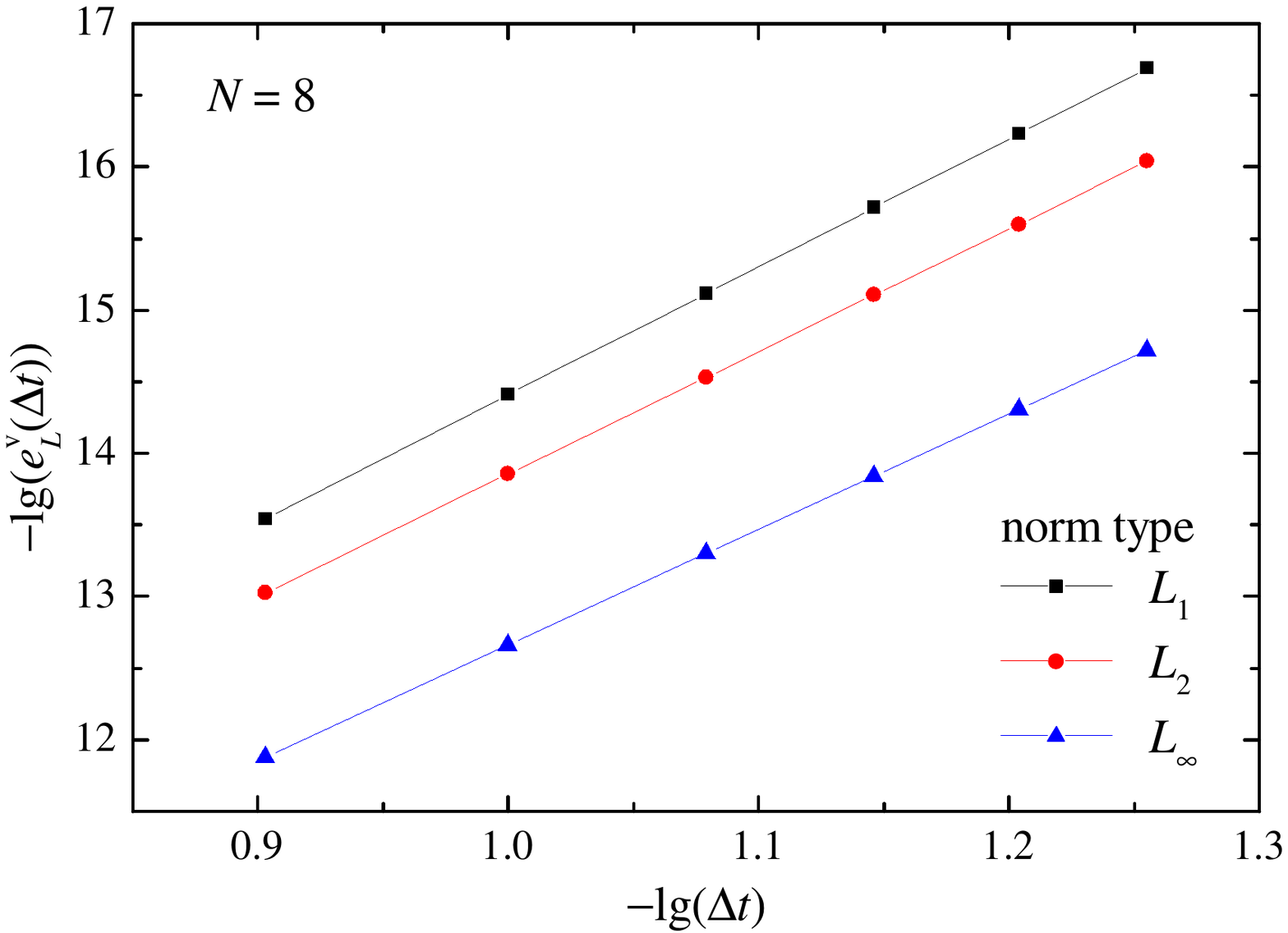}
\vspace{-8mm}\caption{\label{fig:hess_2_ind1_errors:b2}}
\end{subfigure}\hspace{6mm}
\begin{subfigure}{0.275\textwidth}
\includegraphics[width=\textwidth]{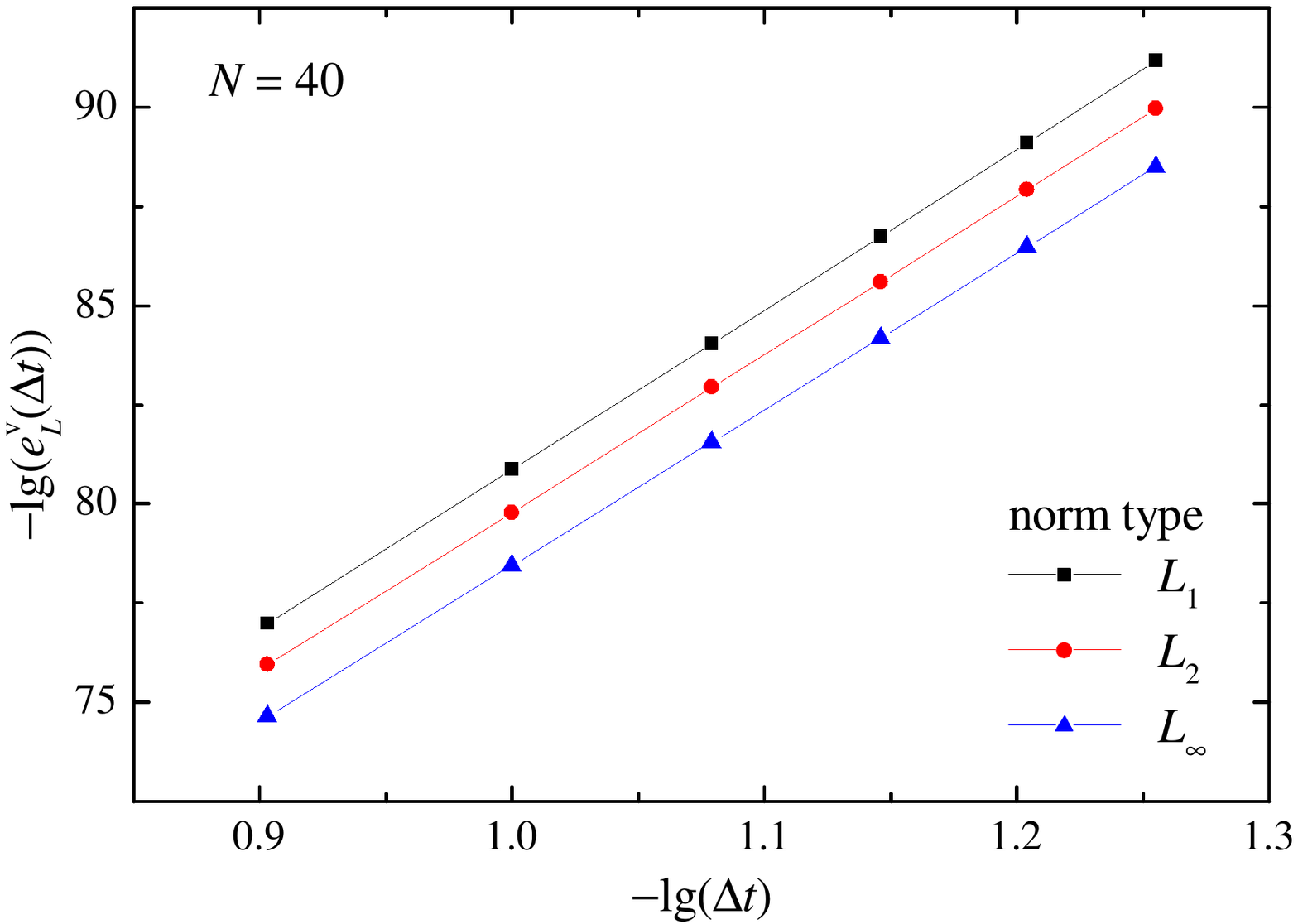}
\vspace{-8mm}\caption{\label{fig:hess_2_ind1_errors:b3}}
\end{subfigure}\\[-2mm]
\begin{subfigure}{0.275\textwidth}
\includegraphics[width=\textwidth]{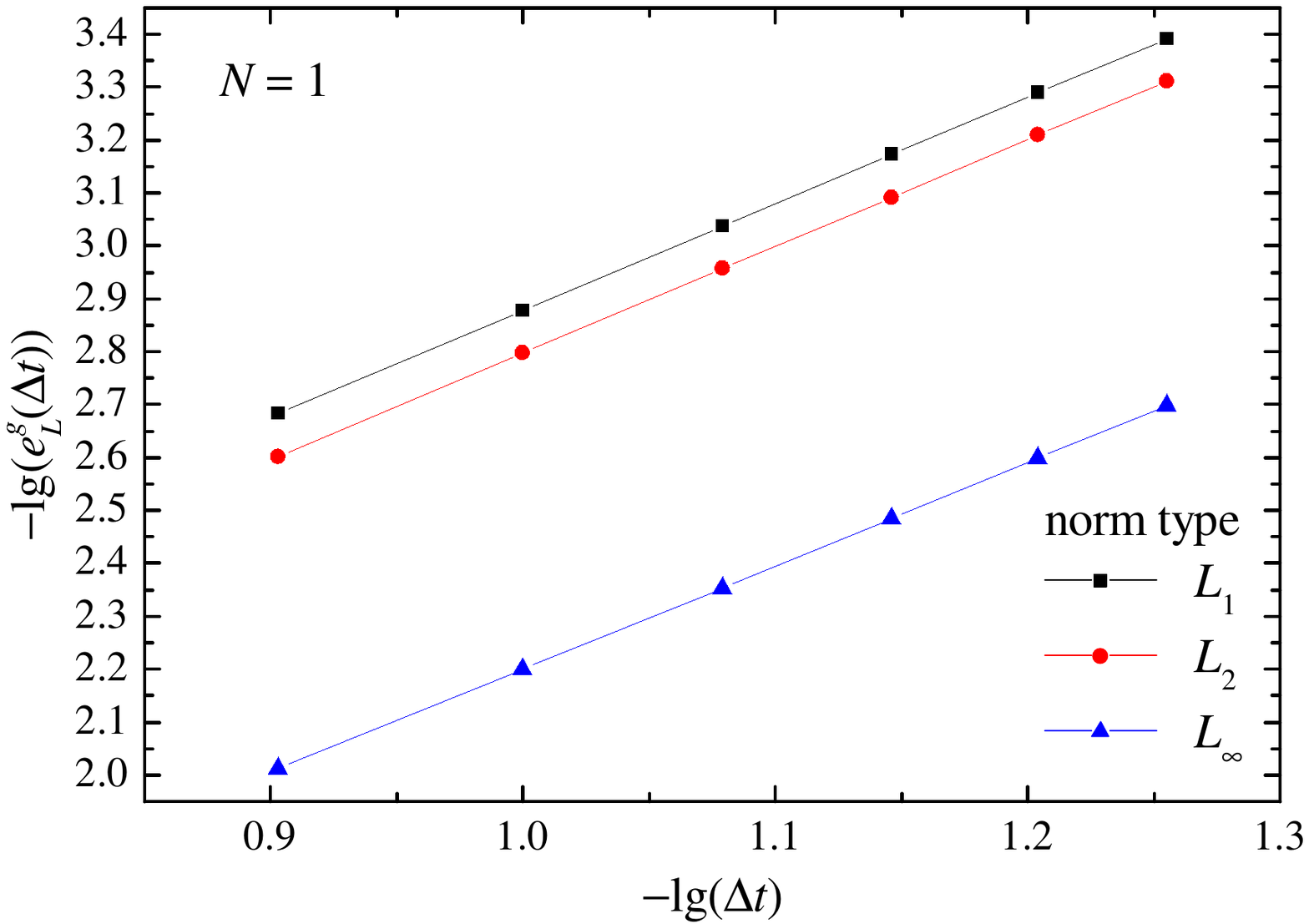}
\vspace{-8mm}\caption{\label{fig:hess_2_ind1_errors:c1}}
\end{subfigure}\hspace{6mm}
\begin{subfigure}{0.275\textwidth}
\includegraphics[width=\textwidth]{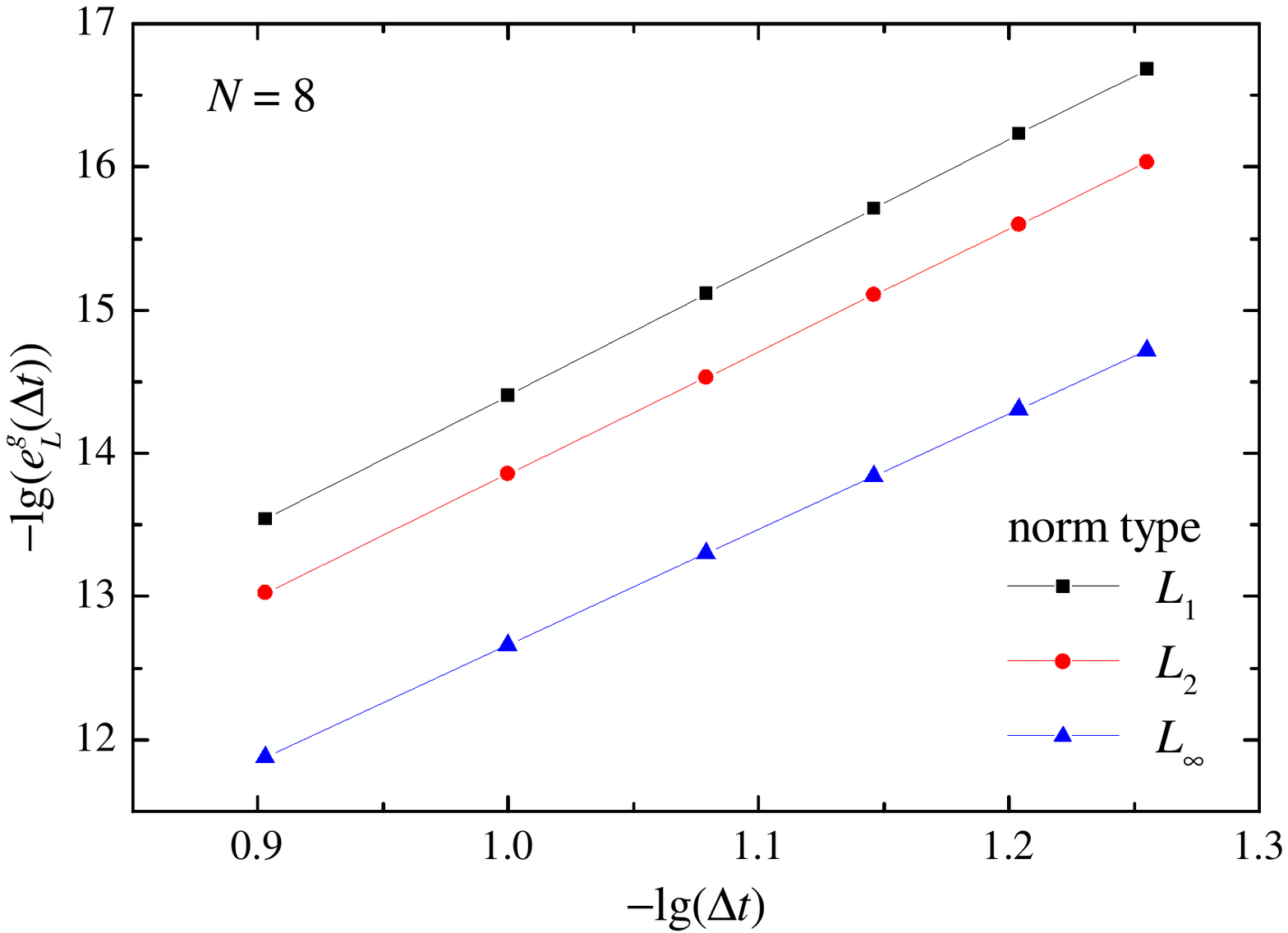}
\vspace{-8mm}\caption{\label{fig:hess_2_ind1_errors:c2}}
\end{subfigure}\hspace{6mm}
\begin{subfigure}{0.275\textwidth}
\includegraphics[width=\textwidth]{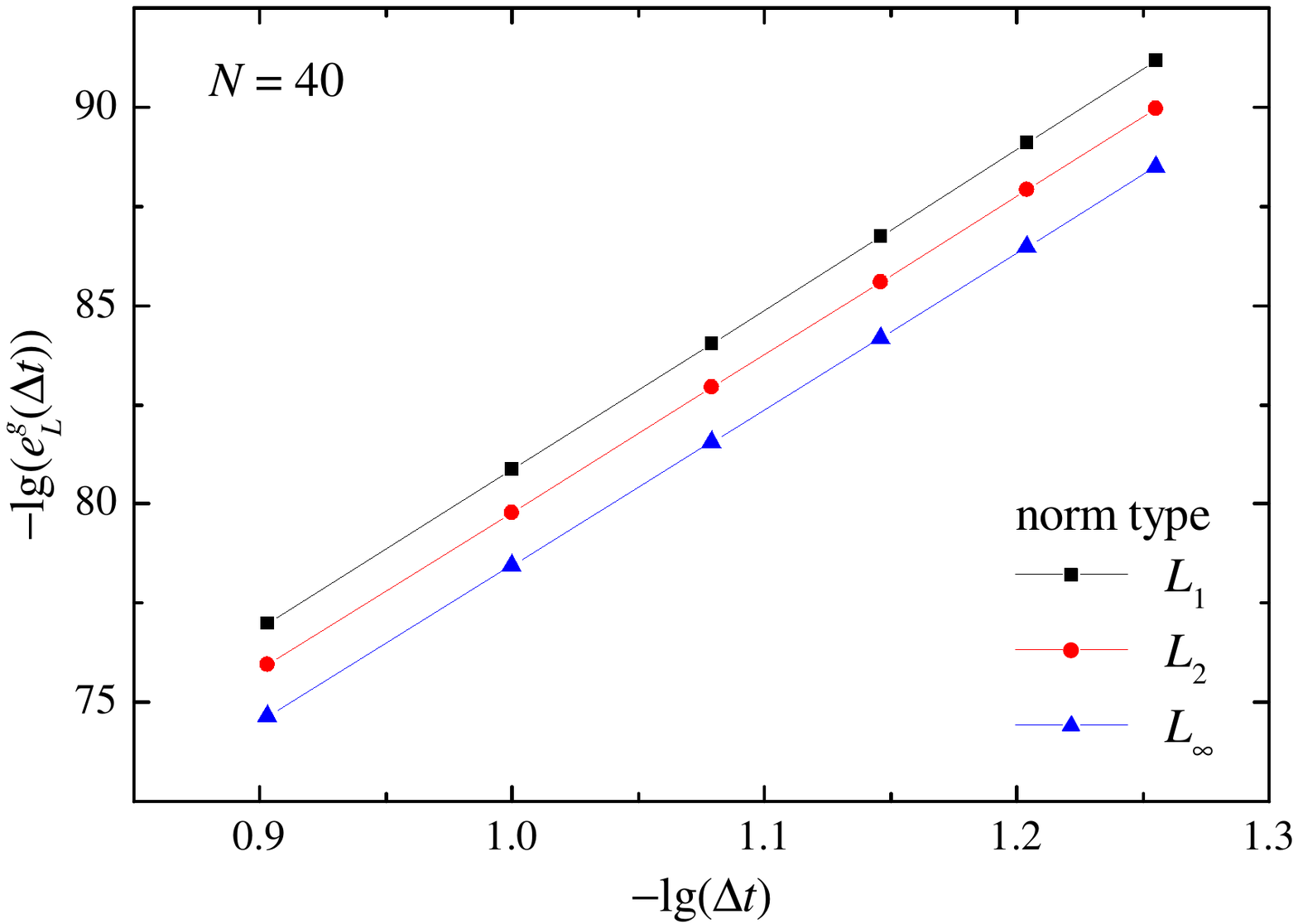}
\vspace{-8mm}\caption{\label{fig:hess_2_ind1_errors:c3}}
\end{subfigure}\\[-2mm]
\begin{subfigure}{0.275\textwidth}
\includegraphics[width=\textwidth]{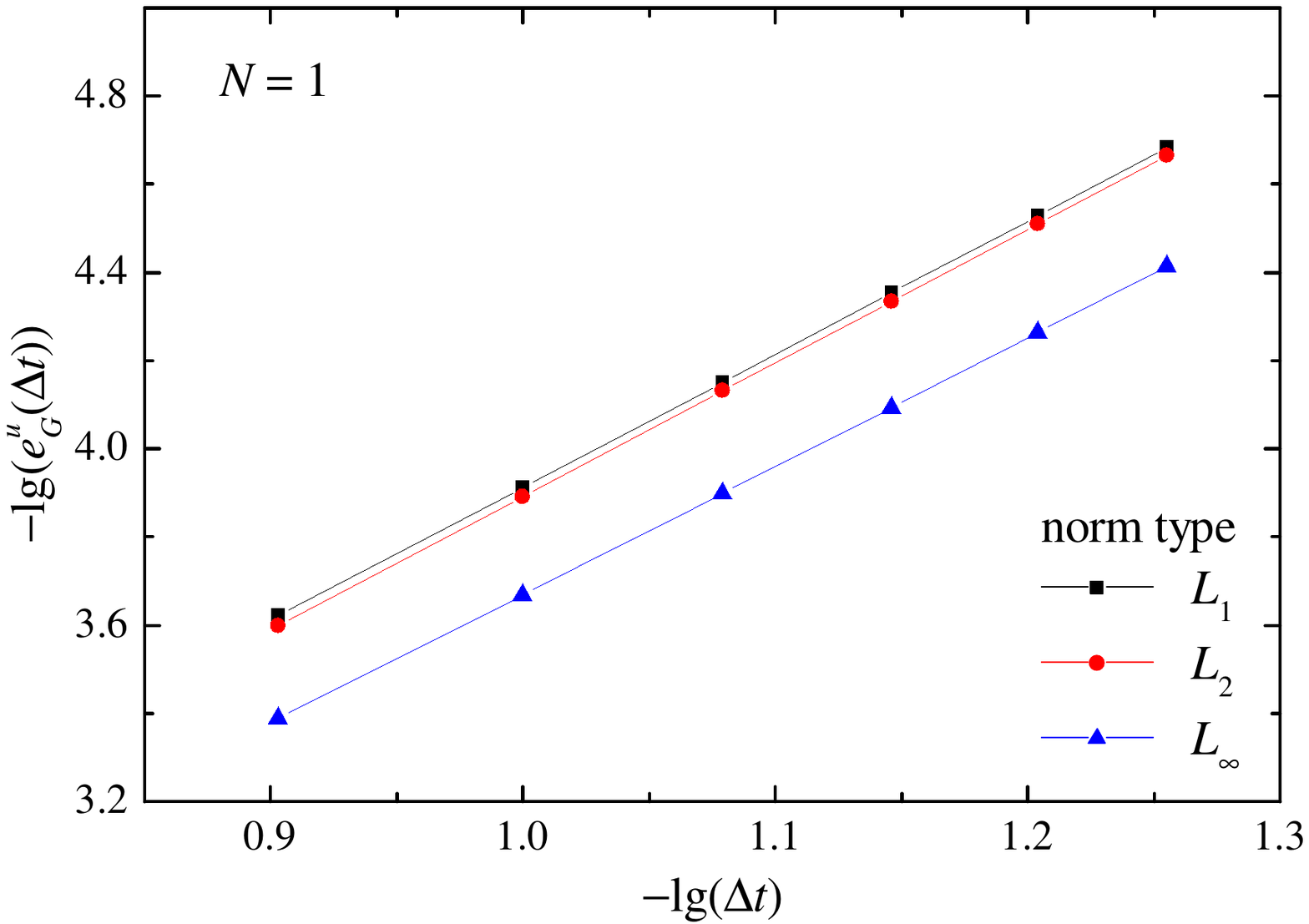}
\vspace{-8mm}\caption{\label{fig:hess_2_ind1_errors:d1}}
\end{subfigure}\hspace{6mm}
\begin{subfigure}{0.275\textwidth}
\includegraphics[width=\textwidth]{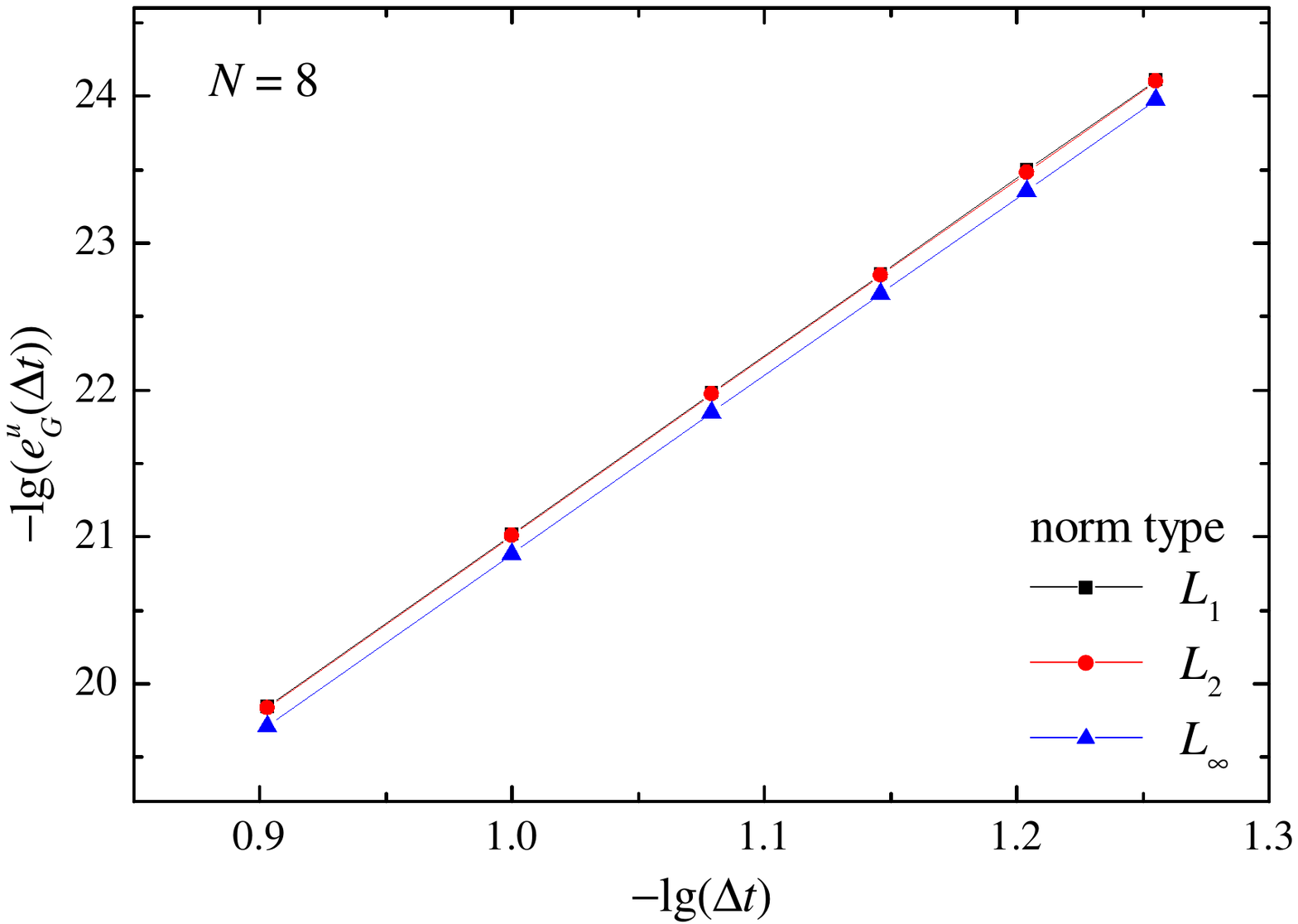}
\vspace{-8mm}\caption{\label{fig:hess_2_ind1_errors:d2}}
\end{subfigure}\hspace{6mm}
\begin{subfigure}{0.275\textwidth}
\includegraphics[width=\textwidth]{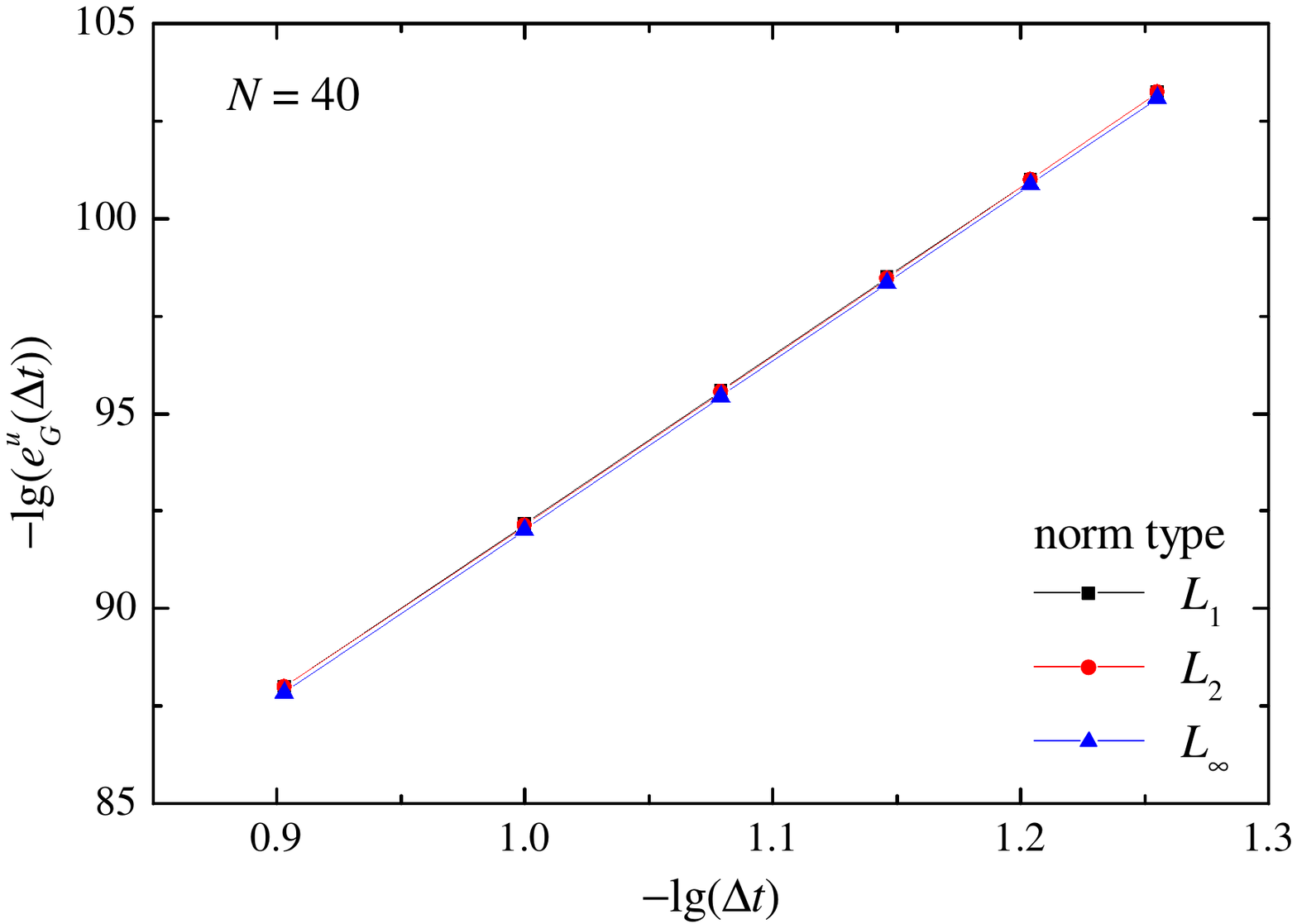}
\vspace{-8mm}\caption{\label{fig:hess_2_ind1_errors:d3}}
\end{subfigure}\\[-2mm]
\begin{subfigure}{0.275\textwidth}
\includegraphics[width=\textwidth]{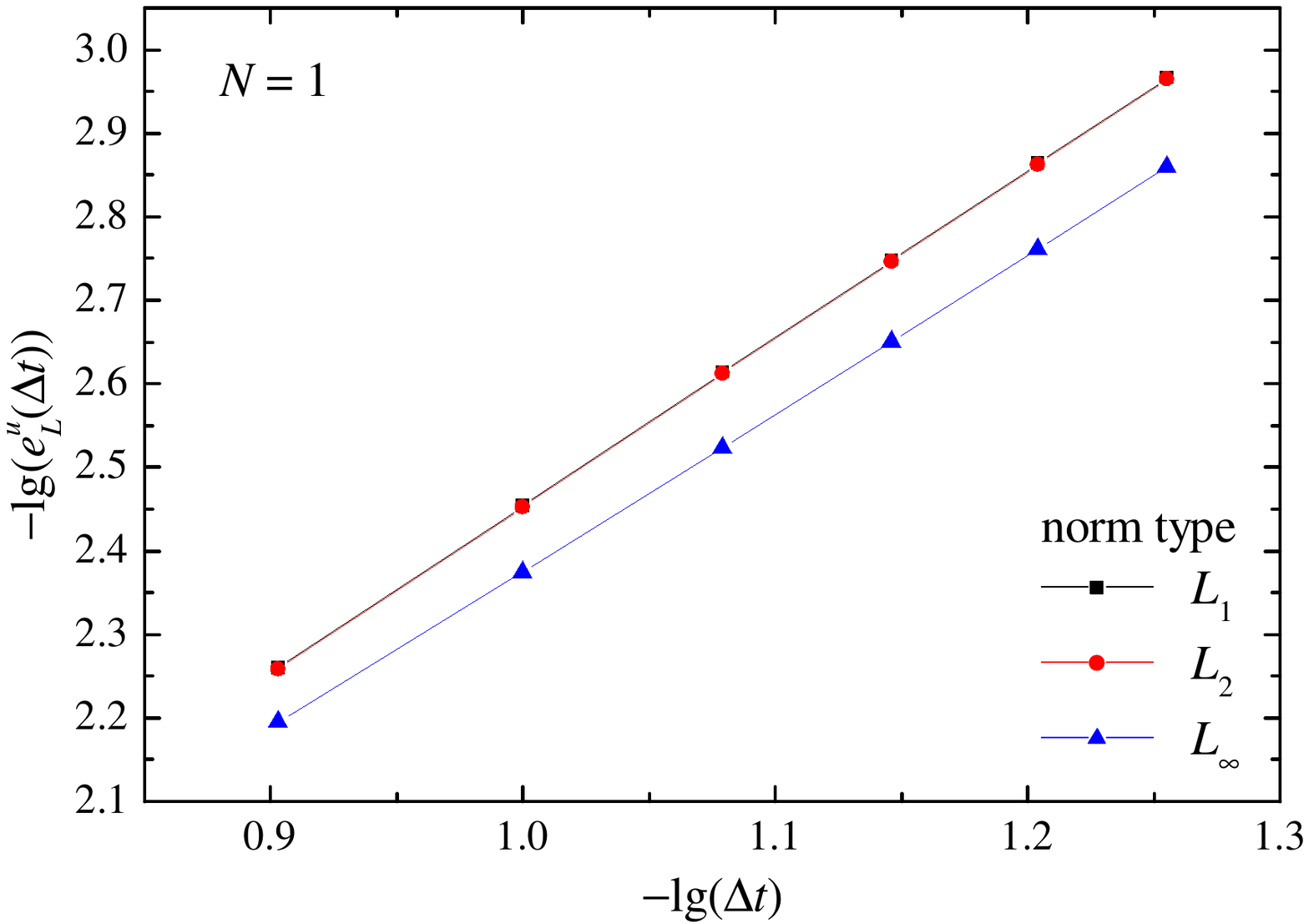}
\vspace{-8mm}\caption{\label{fig:hess_2_ind1_errors:e1}}
\end{subfigure}\hspace{6mm}
\begin{subfigure}{0.275\textwidth}
\includegraphics[width=\textwidth]{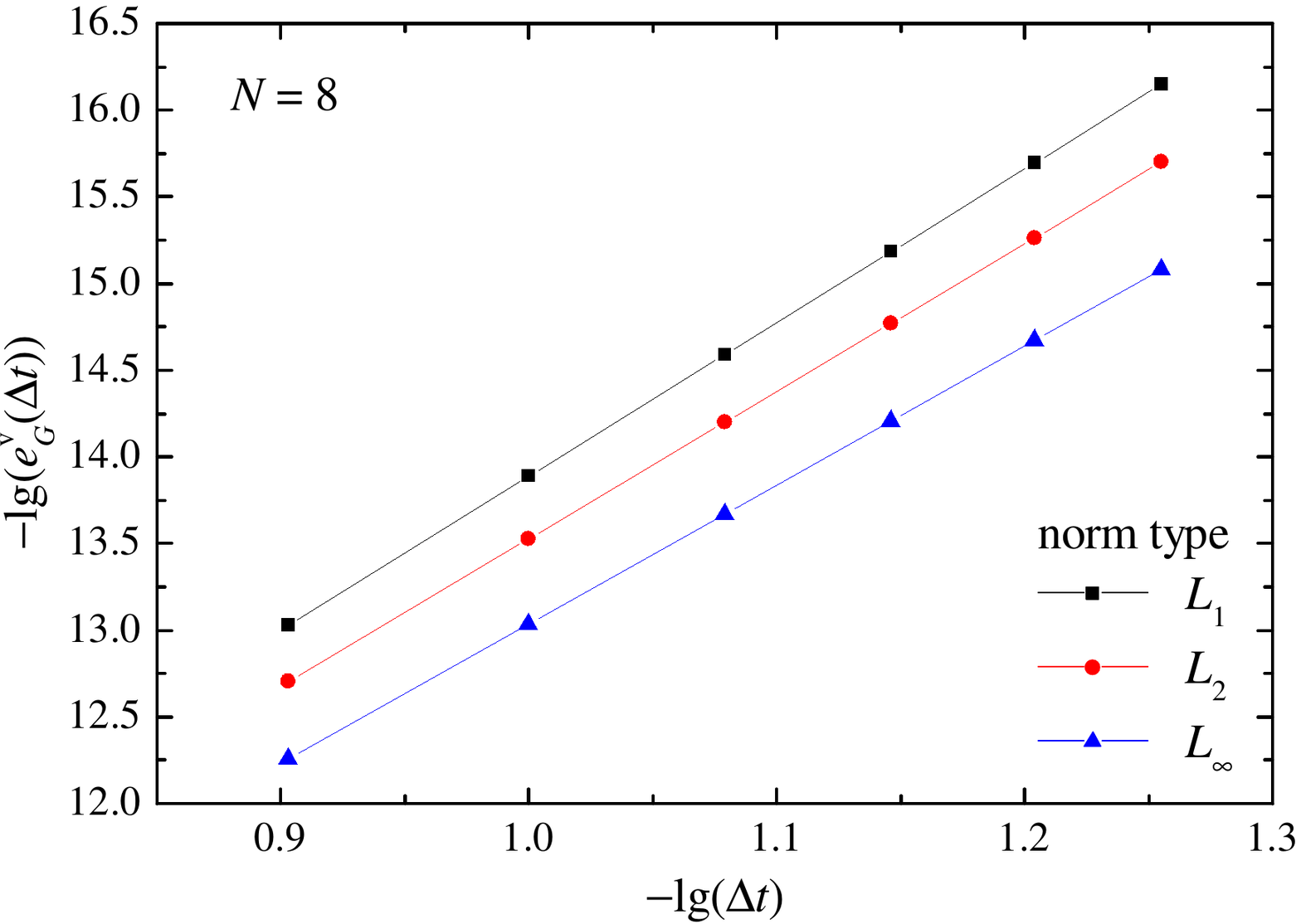}
\vspace{-8mm}\caption{\label{fig:hess_2_ind1_errors:e2}}
\end{subfigure}\hspace{6mm}
\begin{subfigure}{0.275\textwidth}
\includegraphics[width=\textwidth]{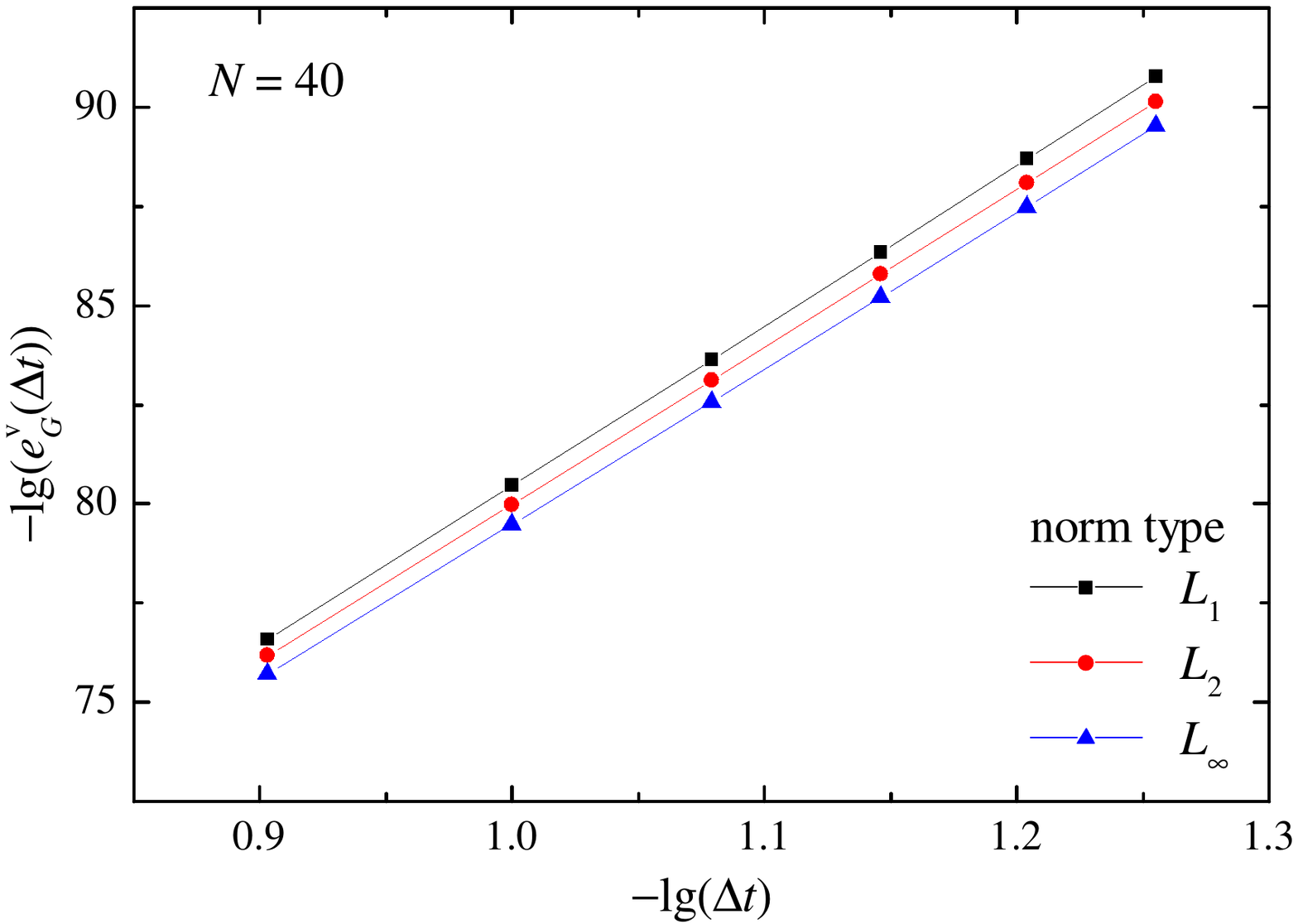}
\vspace{-8mm}\caption{\label{fig:hess_2_ind1_errors:e3}}
\end{subfigure}\\[-2mm]
\begin{subfigure}{0.275\textwidth}
\includegraphics[width=\textwidth]{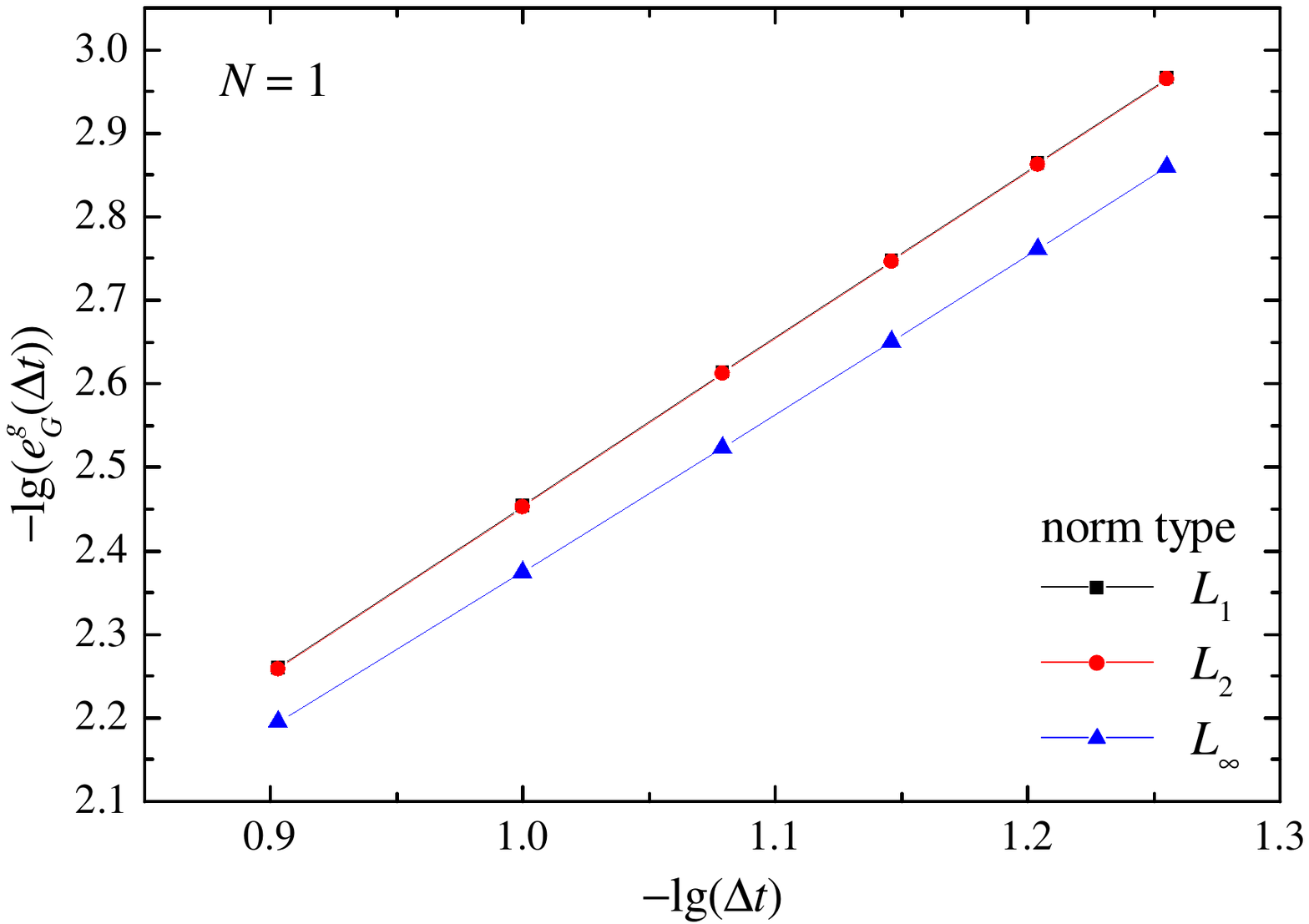}
\vspace{-8mm}\caption{\label{fig:hess_2_ind1_errors:f1}}
\end{subfigure}\hspace{6mm}
\begin{subfigure}{0.275\textwidth}
\includegraphics[width=\textwidth]{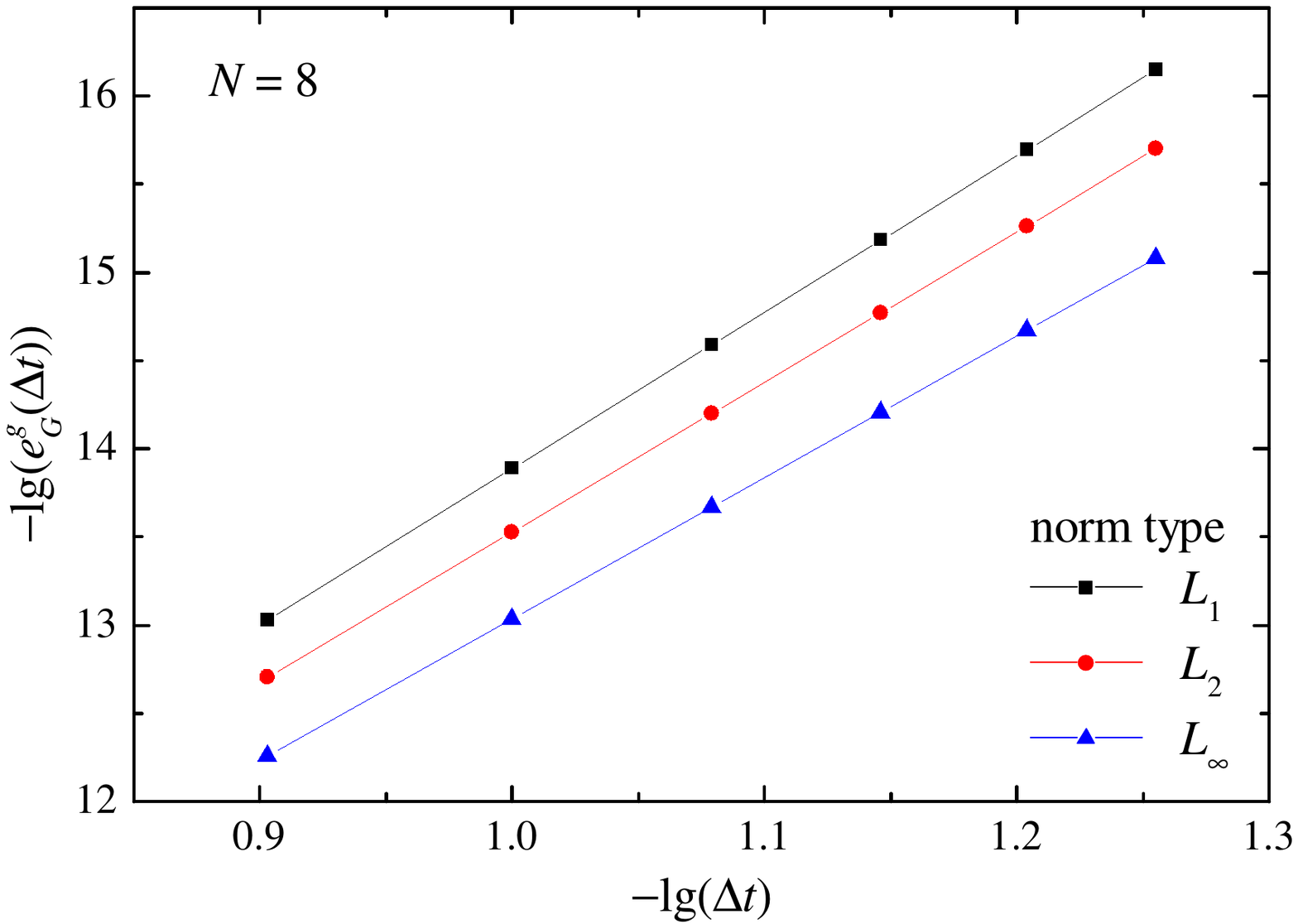}
\vspace{-8mm}\caption{\label{fig:hess_2_ind1_errors:f2}}
\end{subfigure}\hspace{6mm}
\begin{subfigure}{0.275\textwidth}
\includegraphics[width=\textwidth]{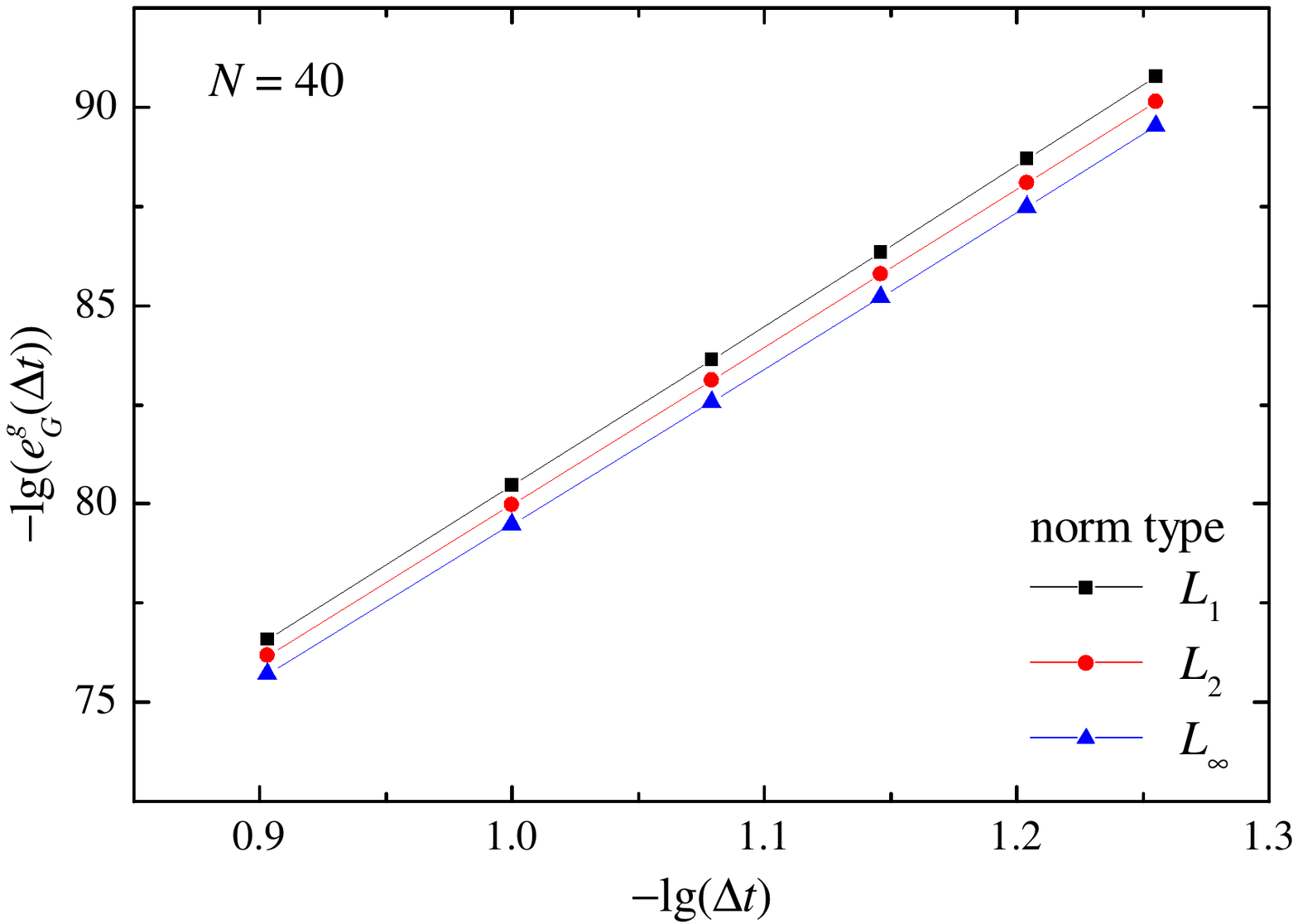}
\vspace{-8mm}\caption{\label{fig:hess_2_ind1_errors:f3}}
\end{subfigure}\\[-2mm]
\caption{%
Log-log plot of the dependence of the global errors for the local solution $e_{L}^{u}$ (\subref{fig:hess_2_ind1_errors:a1}, \subref{fig:hess_2_ind1_errors:a2}, \subref{fig:hess_2_ind1_errors:a3}), $e_{L}^{v}$ (\subref{fig:hess_2_ind1_errors:b1}, \subref{fig:hess_2_ind1_errors:b2}, \subref{fig:hess_2_ind1_errors:b3}), $e_{L}^{g}$ (\subref{fig:hess_2_ind1_errors:c1}, \subref{fig:hess_2_ind1_errors:c2}, \subref{fig:hess_2_ind1_errors:c3}) and the solution at nodes $e_{G}^{u}$ (\subref{fig:hess_2_ind1_errors:d1}, \subref{fig:hess_2_ind1_errors:d2}, \subref{fig:hess_2_ind1_errors:d3}), $e_{G}^{v}$ (\subref{fig:hess_2_ind1_errors:e1}, \subref{fig:hess_2_ind1_errors:e2}, \subref{fig:hess_2_ind1_errors:e3}), $e_{G}^{g}$ (\subref{fig:hess_2_ind1_errors:f1}, \subref{fig:hess_2_ind1_errors:f2}, \subref{fig:hess_2_ind1_errors:f3}) on the discretization step $\mathrm{\Delta}t$, obtained in the norms $L_{1}$, $L_{2}$ and $L_{\infty}$, by numerical solution of the DAE system (\ref{eq:hess_dae_ind_2}) of index 1 obtained using polynomials with degrees $N = 1$ (\subref{fig:hess_2_ind1_errors:a1}, \subref{fig:hess_2_ind1_errors:b1}, \subref{fig:hess_2_ind1_errors:c1}, \subref{fig:hess_2_ind1_errors:d1}, \subref{fig:hess_2_ind1_errors:e1}, \subref{fig:hess_2_ind1_errors:f1}), $N = 8$ (\subref{fig:hess_2_ind1_errors:a2}, \subref{fig:hess_2_ind1_errors:b2}, \subref{fig:hess_2_ind1_errors:c2}, \subref{fig:hess_2_ind1_errors:d2}, \subref{fig:hess_2_ind1_errors:e2}, \subref{fig:hess_2_ind1_errors:f2}), $N = 40$ (\subref{fig:hess_2_ind1_errors:a3}, \subref{fig:hess_2_ind1_errors:b3}, \subref{fig:hess_2_ind1_errors:c3}, \subref{fig:hess_2_ind1_errors:d3}, \subref{fig:hess_2_ind1_errors:e3}, \subref{fig:hess_2_ind1_errors:f3}).
}
\label{fig:hess_2_ind1_errors}
\end{figure} 

\begin{table*}[h!]
\centering
\caption{%
Convergence orders $p_{L_{1}}^{n}$, $p_{L_{2}}^{n}$, $p_{L_{\infty}}^{n}$, calculated in norms $L_{1}$, $L_{2}$, $L_{\infty}$, of \textit{the numerical solution at the nodes} $(\mathbf{u}_{n}, \mathbf{v}_{n})$ of the ADER-DG method for the DAE system (\ref{eq:hess_dae_ind_2}) of index 1; $N$ is the degree of the basis polynomials $\varphi_{p}$. Orders $p^{n, u}$ are calculated for solution $\mathbf{u}_{n}$; orders $p^{n, v}$ --- for solution $\mathbf{v}_{n}$; orders $p^{n, g}$ --- for the conditions $\mathbf{g} = 0$ on the numerical solution at the nodes $(\mathbf{u}_{n}, \mathbf{v}_{n})$. The theoretical value of convergence order $p_{\rm th.}^{n} = 2N+1$ is applicable for the ADER-DG method for ODE problems and is presented for comparison.
}
\label{tab:conv_orders_nodes_hess_2_ind1}
\begin{tabular}{@{}|l|lll|lll|lll|c|@{}}
\toprule
$N$ & $p_{L_{1}}^{n, u}$ & $p_{L_{2}}^{n, u}$ & $p_{L_{\infty}}^{n, u}$ & $p_{L_{1}}^{n, v}$ & $p_{L_{2}}^{n, v}$ & $p_{L_{\infty}}^{n, v}$ & $p_{L_{1}}^{n, g}$ & $p_{L_{2}}^{n, g}$ & $p_{L_{\infty}}^{n, g}$ & $p_{\rm th.}^{n}$ \\
\midrule
$1$	&	$3.02$	&	$3.03$	&	$2.91$	&	$2.01$	&	$2.00$	&	$1.89$	&	$2.01$	&	$2.00$	&	$1.89$	&	$3$\\
$2$	&	$3.84$	&	$3.34$	&	$2.83$	&	$2.71$	&	$2.37$	&	$1.85$	&	$2.71$	&	$2.37$	&	$1.85$	&	$5$\\
$3$	&	$5.63$	&	$5.15$	&	$4.66$	&	$3.58$	&	$3.08$	&	$2.57$	&	$3.58$	&	$3.08$	&	$2.57$	&	$7$\\
$4$	&	$8.07$	&	$8.07$	&	$7.98$	&	$4.77$	&	$4.42$	&	$3.93$	&	$4.77$	&	$4.42$	&	$3.93$	&	$9$\\
$5$	&	$9.89$	&	$9.92$	&	$9.96$	&	$6.05$	&	$6.05$	&	$5.98$	&	$6.05$	&	$6.05$	&	$5.98$	&	$11$\\
$6$	&	$9.86$	&	$9.87$	&	$9.84$	&	$6.64$	&	$6.22$	&	$5.72$	&	$6.64$	&	$6.22$	&	$5.72$	&	$13$\\
$7$	&	$10.90$	&	$10.90$	&	$10.88$	&	$7.65$	&	$7.23$	&	$6.73$	&	$7.65$	&	$7.23$	&	$6.73$	&	$15$\\
$8$	&	$12.12$	&	$12.12$	&	$12.09$	&	$8.86$	&	$8.51$	&	$8.02$	&	$8.86$	&	$8.51$	&	$8.02$	&	$17$\\
$9$	&	$11.98$	&	$11.98$	&	$11.96$	&	$9.34$	&	$8.39$	&	$7.47$	&	$9.34$	&	$8.39$	&	$7.47$	&	$19$\\
$10$	&	$13.81$	&	$13.81$	&	$13.78$	&	$10.63$	&	$10.18$	&	$9.69$	&	$10.63$	&	$10.18$	&	$9.69$	&	$21$\\
$11$	&	$14.91$	&	$14.91$	&	$14.89$	&	$11.72$	&	$11.30$	&	$10.80$	&	$11.72$	&	$11.30$	&	$10.80$	&	$23$\\
$12$	&	$16.22$	&	$16.22$	&	$16.20$	&	$12.97$	&	$12.64$	&	$12.15$	&	$12.97$	&	$12.64$	&	$12.15$	&	$25$\\
$13$	&	$16.39$	&	$16.39$	&	$16.37$	&	$13.38$	&	$12.70$	&	$12.18$	&	$13.38$	&	$12.70$	&	$12.18$	&	$27$\\
$14$	&	$17.79$	&	$17.79$	&	$17.77$	&	$14.64$	&	$14.19$	&	$13.69$	&	$14.64$	&	$14.19$	&	$13.69$	&	$29$\\
$15$	&	$18.93$	&	$18.94$	&	$18.91$	&	$15.78$	&	$15.35$	&	$14.85$	&	$15.78$	&	$15.35$	&	$14.85$	&	$31$\\
$16$	&	$20.49$	&	$20.49$	&	$20.46$	&	$17.17$	&	$16.91$	&	$16.43$	&	$17.17$	&	$16.91$	&	$16.43$	&	$33$\\
$17$	&	$20.51$	&	$20.51$	&	$20.49$	&	$17.43$	&	$16.87$	&	$16.37$	&	$17.43$	&	$16.87$	&	$16.37$	&	$35$\\
$18$	&	$21.80$	&	$21.80$	&	$21.77$	&	$18.67$	&	$18.21$	&	$17.71$	&	$18.67$	&	$18.21$	&	$17.71$	&	$37$\\
$19$	&	$22.98$	&	$22.98$	&	$22.95$	&	$19.84$	&	$19.41$	&	$18.91$	&	$19.84$	&	$19.41$	&	$18.91$	&	$39$\\
$20$	&	$25.57$	&	$25.57$	&	$25.55$	&	$21.64$	&	$21.75$	&	$21.36$	&	$21.64$	&	$21.75$	&	$21.36$	&	$41$\\
\midrule
$25$	&	$28.60$	&	$28.60$	&	$28.58$	&	$25.52$	&	$25.01$	&	$24.51$	&	$25.52$	&	$25.01$	&	$24.51$	&	$51$\\
$30$	&	$33.85$	&	$33.86$	&	$33.83$	&	$30.76$	&	$30.29$	&	$29.79$	&	$30.76$	&	$30.29$	&	$29.79$	&	$61$\\
$35$	&	$39.57$	&	$39.57$	&	$39.55$	&	$36.36$	&	$36.01$	&	$35.52$	&	$36.36$	&	$36.01$	&	$35.52$	&	$71$\\
$40$	&	$43.31$	&	$43.32$	&	$43.29$	&	$40.24$	&	$39.71$	&	$39.21$	&	$40.24$	&	$39.71$	&	$39.21$	&	$81$\\
\bottomrule
\end{tabular}
\end{table*} 

\begin{table*}[h!]
\centering
\caption{%
Convergence orders $p_{L_{1}}^{l}$, $p_{L_{2}}^{l}$, $p_{L_{\infty}}^{l}$, calculated in norms $L_{1}$, $L_{2}$, $L_{\infty}$, of \textit{the local solution} $(\mathbf{u}_{L}, \mathbf{v}_{L})$ (represented between the nodes) of the ADER-DG method for the DAE system (\ref{eq:hess_dae_ind_2}) of index 1; $N$ is the degree of the basis polynomials $\varphi_{p}$. Orders $p^{l, u}$ are calculated for solution $\mathbf{u}_{L}$; orders $p^{l, v}$ --- for solution $\mathbf{v}_{L}$; orders $p^{l, g}$ --- for the conditions $\mathbf{g} = 0$ on the local solution $(\mathbf{u}_{L}, \mathbf{v}_{L})$. The theoretical value of convergence order $p_{\rm th.}^{l} = N+1$ is applicable for the ADER-DG method for ODE problems and is presented for comparison.
}
\label{tab:conv_orders_local_hess_2_ind1}
\begin{tabular}{@{}|l|lll|lll|lll|c|@{}}
\toprule
$N$ & $p_{L_{1}}^{l, u}$ & $p_{L_{2}}^{l, u}$ & $p_{L_{\infty}}^{l, u}$ & $p_{L_{1}}^{l, v}$ & $p_{L_{2}}^{l, v}$ & $p_{L_{\infty}}^{l, v}$ & $p_{L_{1}}^{l, g}$ & $p_{L_{2}}^{l, g}$ & $p_{L_{\infty}}^{l, g}$ & $p_{\rm th.}^{l}$ \\
\midrule
$1$	&	$2.00$	&	$1.99$	&	$1.99$	&	$2.01$	&	$2.00$	&	$1.88$	&	$2.01$	&	$2.02$	&	$1.95$	&	$2$\\
$2$	&	$2.98$	&	$2.95$	&	$2.84$	&	$2.70$	&	$2.37$	&	$1.85$	&	$2.74$	&	$2.40$	&	$1.85$	&	$3$\\
$3$	&	$4.00$	&	$4.00$	&	$3.85$	&	$3.61$	&	$3.13$	&	$2.64$	&	$3.63$	&	$3.14$	&	$2.64$	&	$4$\\
$4$	&	$5.01$	&	$5.00$	&	$4.96$	&	$4.81$	&	$4.46$	&	$3.99$	&	$4.81$	&	$4.46$	&	$3.99$	&	$5$\\
$5$	&	$6.00$	&	$6.00$	&	$5.99$	&	$6.08$	&	$6.10$	&	$5.95$	&	$6.08$	&	$6.08$	&	$5.92$	&	$6$\\
$6$	&	$7.00$	&	$6.97$	&	$6.71$	&	$6.69$	&	$6.26$	&	$5.77$	&	$6.70$	&	$6.26$	&	$5.77$	&	$7$\\
$7$	&	$7.99$	&	$7.98$	&	$7.84$	&	$7.72$	&	$7.29$	&	$6.81$	&	$7.72$	&	$7.29$	&	$6.81$	&	$8$\\
$8$	&	$9.00$	&	$9.00$	&	$9.01$	&	$8.93$	&	$8.55$	&	$8.06$	&	$8.93$	&	$8.55$	&	$8.06$	&	$9$\\
$9$	&	$10.00$	&	$10.00$	&	$9.99$	&	$9.29$	&	$8.35$	&	$7.78$	&	$9.33$	&	$8.37$	&	$7.78$	&	$10$\\
$10$	&	$10.98$	&	$10.95$	&	$10.72$	&	$10.68$	&	$10.24$	&	$9.75$	&	$10.69$	&	$10.24$	&	$9.75$	&	$11$\\
$11$	&	$11.99$	&	$11.98$	&	$11.86$	&	$11.79$	&	$11.35$	&	$10.86$	&	$11.79$	&	$11.35$	&	$10.86$	&	$12$\\
$12$	&	$13.00$	&	$13.01$	&	$12.95$	&	$13.05$	&	$12.68$	&	$12.18$	&	$13.05$	&	$12.68$	&	$12.18$	&	$13$\\
$13$	&	$13.96$	&	$13.95$	&	$13.82$	&	$13.38$	&	$12.81$	&	$12.31$	&	$13.40$	&	$12.81$	&	$12.31$	&	$14$\\
$14$	&	$14.98$	&	$14.95$	&	$14.73$	&	$14.70$	&	$14.25$	&	$13.75$	&	$14.70$	&	$14.25$	&	$13.75$	&	$15$\\
$15$	&	$16.00$	&	$15.99$	&	$15.89$	&	$15.84$	&	$15.39$	&	$14.90$	&	$15.84$	&	$15.39$	&	$14.90$	&	$16$\\
$16$	&	$17.02$	&	$17.01$	&	$16.96$	&	$17.29$	&	$16.95$	&	$16.46$	&	$17.28$	&	$16.95$	&	$16.46$	&	$17$\\
$17$	&	$17.95$	&	$17.90$	&	$17.40$	&	$17.47$	&	$16.96$	&	$16.46$	&	$17.47$	&	$16.96$	&	$16.46$	&	$18$\\
$18$	&	$18.98$	&	$18.95$	&	$18.75$	&	$18.72$	&	$18.26$	&	$17.76$	&	$18.72$	&	$18.26$	&	$17.76$	&	$19$\\
$19$	&	$20.00$	&	$20.00$	&	$19.93$	&	$19.90$	&	$19.44$	&	$18.95$	&	$19.90$	&	$19.44$	&	$18.95$	&	$20$\\
$20$	&	$21.00$	&	$20.99$	&	$21.00$	&	$21.98$	&	$21.99$	&	$21.51$	&	$21.97$	&	$21.99$	&	$21.51$	&	$21$\\
\midrule
$25$	&	$25.96$	&	$25.89$	&	$25.55$	&	$25.56$	&	$25.07$	&	$24.57$	&	$25.56$	&	$25.07$	&	$24.57$	&	$26$\\
$30$	&	$30.98$	&	$30.97$	&	$30.81$	&	$30.81$	&	$30.33$	&	$29.83$	&	$30.81$	&	$30.33$	&	$29.83$	&	$31$\\
$35$	&	$36.03$	&	$36.03$	&	$35.94$	&	$36.47$	&	$36.04$	&	$35.54$	&	$36.47$	&	$36.04$	&	$35.54$	&	$36$\\
$40$	&	$40.95$	&	$40.76$	&	$40.26$	&	$40.28$	&	$39.78$	&	$39.28$	&	$40.28$	&	$39.78$	&	$39.28$	&	$41$\\
\bottomrule
\end{tabular}
\end{table*}

The obtained results of the numerical solution of DAE system (\ref{eq:hess_dae_ind_2}) of index 2 are presented in Figs.~\ref{fig:hess_2_ind2_sol_uv}, \ref{fig:hess_2_ind2_sol_g_eps}, \ref{fig:hess_2_ind2_errors} and in Tables~\ref{tab:conv_orders_nodes_hess_2_ind2}, \ref{tab:conv_orders_local_hess_2_ind2}. Fig.~\ref{fig:hess_2_ind2_sol_uv} shows a comparison of the numerical solution at the nodes $(\mathbf{u}_{n}, \mathbf{v}_{n})$, the numerical local solution $(\mathbf{u}_{L}, \mathbf{v}_{L})$ and the exact analytical solution $(\mathbf{u}^{\rm ex}, \mathbf{v}^{\rm ex})$ separately for each differential $\mathbf{u}$ and algebraic $\mathbf{v}$ variable. Fig.~\ref{fig:hess_2_ind2_sol_g_eps} shows the dependencies of the feasibility of conditions $g_{1} = 0$ and $g_{2} = 0$ on the coordinate $t$, as well as the dependencies of the local errors $\varepsilon_{u}$, $\varepsilon_{v}$, $\varepsilon_{g}$ on the coordinate $t$, which allows us to quantitatively estimate the accuracy of the numerical solution, especially taking into account the fact that the numerical solution obtained by the ADER-DG method with a local DG predictor with a high degree of polynomials $N$ very accurately corresponds to the exact analytical solution, and it is impossible to visually identify the error from the solution plots separately (in Fig.~\ref{fig:hess_2_ind2_sol_uv}). Fig.~\ref{fig:hess_2_ind2_errors} shows the dependencies of the global errors $e^{u}$, $e^{v}$, $e^{g}$ of the numerical solution at the nodes $(\mathbf{u}_{n}, \mathbf{v}_{n})$ and the local solution $(\mathbf{u}_{L}, \mathbf{v}_{L})$ on the discretization step ${\Delta t}$, separately for each differential $\mathbf{u}$ and algebraic variable $\mathbf{v}$ and the algebraic equations $\mathbf{g} = \mathbf{0}$, on the basis of which the empirical convergence orders $p$ were calculated.

The numerical solution of the DAE system (\ref{eq:hess_dae_ind_2}) of index 2, presented in Fig.~\ref{fig:hess_2_ind2_sol_uv}, clearly demonstrates the high accuracy achieved using the ADER-DG method with a local DG predictor. The most significant features of the numerical solution are demonstrated in the case of polynomial degree $N = 1$, the numerical solutions for which are shown in Figs.~\ref{fig:hess_2_ind2_sol_uv}~(\subref{fig:hess_2_ind2_sol_uv:a1}, \subref{fig:hess_2_ind2_sol_uv:b1}, \subref{fig:hess_2_ind2_sol_uv:c1}, \subref{fig:hess_2_ind2_sol_uv:d1}, \subref{fig:hess_2_ind2_sol_uv:e1}). In the case of polynomial degree $N = 1$, the dependence of the local solution for differential variables $\mathbf{u}_{L}$ on the coordinate $t$ shows the expected discontinuities of the solution at the grid nodes $t_{n}$, which was already observed in the previous examples considered, however, the behavior of the local solution for the algebraic variable $\mathbf{v}_{L}$ demonstrates a different type of feature --- in the first discretization domain on the grid $\Omega_{0}$, the local solution $\mathbf{v}_{L}$ agrees well with the exact analytical solution $\mathbf{v}^{\rm ex}$, however, in subsequent discretization domains $\Omega_{n}$, the point-wise correspondence between the local solution $\mathbf{v}_{L}$ and the exact analytical solution $\mathbf{v}^{\rm ex}$ significantly worsens. This is due to the condition that is set for obtaining a local solution in the original system of nonlinear algebraic equations of the predictor (\ref{eq:dae_weak_form}) --- despite the point-wise mismatch, the condition of orthogonality of the local error $\boldsymbol{\sigma}_{\rm v}$ to all basis functions $\varphi_{p}$, the expansion of which represents the numerical solution, is satisfied; however, it is necessary to use such a local solution $\mathbf{v}_{L}$ with great ``caution'' as a sufficiently accurate solution between grid nodes. In the previous examples, such behavior of the local solution $\mathbf{v}_{L}$ for algebraic variables in the case of polynomial degree $N = 1$ was not observed, which is due to the index 1 of the problems in the previous examples; this phenomenon begins to manifest itself significantly for DAE systems of index higher than 1. However, despite this feature of the behavior of the local solution $\mathbf{v}_{L}$ for algebraic variables, the convergence order $p_{j, v}$ for this case turned out to be positive --- with a decrease in the discretization step ${\Delta t}$, the amplitude of discontinuities in the local solution $\mathbf{v}_{L}$ at the grid nodes $t_{n}$ decreases. The presented features of the numerical solution practically disappear in the case of polynomial degrees $N > 1$, in particular, on the presented plots of the numerical solution for polynomials degree $N = 8$ (see Figs.~\ref{fig:hess_2_ind2_sol_uv}~(\subref{fig:hess_2_ind2_sol_uv:a2}, \subref{fig:hess_2_ind2_sol_uv:b2}, \subref{fig:hess_2_ind2_sol_uv:c2}, \subref{fig:hess_2_ind2_sol_uv:d2}, \subref{fig:hess_2_ind2_sol_uv:e2})) and $40$ (see Figs.~\ref{fig:hess_2_ind2_sol_uv}~(\subref{fig:hess_2_ind2_sol_uv:a3}, \subref{fig:hess_2_ind2_sol_uv:b3}, \subref{fig:hess_2_ind2_sol_uv:c3}, \subref{fig:hess_2_ind2_sol_uv:d3}, \subref{fig:hess_2_ind2_sol_uv:e3})). In these cases, the numerical solution at the nodes $(\mathbf{u}_{n}, \mathbf{v}_{n})$ and the local solution $(\mathbf{u}_{L}, \mathbf{v}_{L})$ coincide with the exact analytical solution $(\mathbf{u}^{\rm ex}, \mathbf{v}^{\rm ex})$ with very high accuracy, and the differences that arise on these plots are visually indistinguishable.

To quantitatively determine the accuracy of the local solution $(\mathbf{u}_{L}, \mathbf{v}_{L})$ and the solution at the nodes $(\mathbf{u}_{n}, \mathbf{v}_{n})$, plots of the local errors $|g_{1}|$, $|g_{2}|$ of satisfaction of the algebraic equations $g_{1} = 0$, $g_{2} = 0$ and local errors $\varepsilon_{u}$, $\varepsilon_{v}$, $\varepsilon_{g}$ were constructed, which were presented in Fig.~\ref{fig:hess_2_ind2_sol_g_eps}. From the presented dependencies of the error $|g_{1}|$, it was concluded that in the case of the solution at the nodes $(\mathbf{u}_{n}, \mathbf{v}_{n})$, the algebraic equation $g_{1} = 0$, which was explicitly included in the DAE system (\ref{eq:hess_dae_ind_2}) of index 2, is satisfied exactly (of course, within the accuracy of the representation of real numbers with a floating point). However, in the case of a local solution $(\mathbf{u}_{L}, \mathbf{v}_{L})$, the algebraic equation $g_{1} = 0$ is satisfied with a nonzero error $|g_{1}| \sim 10^{-7}$--$10^{-2}$ in the case of polynomial degree $N = 1$, $|g_{1}| \sim 10^{-20}$--$10^{-14}$ in the case of $N = 8$, $|g_{1}| \sim 10^{-85}$--$10^{-87}$ in the case of $N = 40$. The presented dependencies of the error $|g_{2}|$, which demonstrate the accuracy of satisfying the algebraic equation $g_{2} = 0$, which is a derivative with respect to $g_{1} = 0$ by $t$, show a nonzero accuracy in the condition $g_{2} = 0$ both for the solution at the nodes $(\mathbf{u}_{n}, \mathbf{v}_{n})$ and for the local solution $(\mathbf{u}_{L}, \mathbf{v}_{L})$. It is necessary to note an important feature of the behavior of the error $|g_{2}|$ --- for almost all values of the coordinate $t$, the error $|g_{2}|$ on the local solution $(\mathbf{u}_{L}, \mathbf{v}_{L})$ is lower than the error $|g_{2}|$ on the solution at the nodes $(\mathbf{u}_{n}, \mathbf{v}_{n})$. This phenomenon begins to manifest itself significantly for DAE systems of index higher than 1, and is associated with the commensurability of the convergence orders for the local solution $p^{l, u}$, $p^{l, v}$, $p^{l, g}$ and the solution at the nodes $p^{n, u}$, $p^{n, v}$, $p^{n, g}$ for DAE systems with a high index. It should also be noted that the error $|g_{2}| \sim 10^{2}$--$10^{3}$ times higher than the error $|g_{1}|$ for local solution $(\mathbf{u}_{L}, \mathbf{v}_{L})$, therefore it was found that the condition of index 1 is satisfied with less accuracy than the condition of index 2 in the case of the DAE system (\ref{eq:hess_dae_ind_2}) of index 2. The errors $\varepsilon_{u}$ calculated for the differential variables $\mathbf{u}_{L}$ demonstrate the classical expected behavior for the solution at the nodes and the local solution. However, for high degrees of the polynomials $N > 1$, the error $\varepsilon_{u}$ of the solution at the nodes $\mathbf{u}_{n}$ becomes higher than the error of the local solution $\mathbf{u}_{L}$ for almost all values of the coordinate $t$. The errors $\varepsilon_{v}$ and $\varepsilon_{g}$ demonstrate such behavior for all values of the polynomial degrees $N$. The characteristic values of the errors $\varepsilon_{u}$, $\varepsilon_{v}$, $\varepsilon_{g}$ for the polynomial degree $N = 1$ are in the range of $10^{-5}$--$10^{-2}$, for the polynomial degree $N = 8$ --- $10^{-17}$--$10^{-13}$, for the polynomial degree $N = 40$ --- $10^{-80}$--$10^{-74}$. Unlike the previous examples, in this case there are no significant differences in the characteristic values of local errors of the local solution $(\mathbf{u}_{L}, \mathbf{v}_{L})$ and the solution at the nodes $(\mathbf{u}_{n}, \mathbf{v}_{n})$. Fig.~\ref{fig:hess_2_ind2_errors} shows the log-log dependencies of global errors $e$ for the solution at nodes $(\mathbf{u}_{n}, \mathbf{v}_{n})$ and the local solution $(\mathbf{u}_{L}, \mathbf{v}_{L})$ separately for differential and algebraic variables and algebraic equations. The presented dependencies clearly demonstrate the power law $e(\Delta t) \sim (\Delta t)^{p}$, from which the empirical convergence orders $p$ were calculated. Unlike the previous examples, in this case there are no significant differences in the characteristic values of global errors $e$ of the local solution $(\mathbf{u}_{L}, \mathbf{v}_{L})$ and the solution at nodes $(\mathbf{u}_{n}, \mathbf{v}_{n})$.

The calculated empirical convergence orders $p^{n, u}$, $p^{n, v}$, $p^{n, g}$ for the solution at nodes $(\mathbf{u}_{n}, \mathbf{v}_{n})$ are presented in Table~\ref{tab:conv_orders_nodes_hess_2_ind2}, convergence orders $p^{l, u}$, $p^{l, v}$, $p^{l, g}$ for the local solution $(\mathbf{u}_{L}, \mathbf{v}_{L})$ --- in Table~\ref{tab:conv_orders_local_hess_2_ind2}. The convergence orders were calculated separately for the norms $L_{1}$, $L_{2}$, $L_{\infty}$ (\ref{eq:norms_def}). The obtained results differ significantly from similar results for the previous examples. The convergence orders of the local solution $p^{l, u}$ and the solution at nodes $p^{n, u}$ for differential variables approximately coincide with each other and with the expected convergence order for the local solution $p_{\rm local} = N+1$, therefore superconvergence does not occur in this case. The convergence orders of the local solution $p^{l, v}$, $p^{l, g}$ and the solution at the nodes $p^{l, v}$, $p^{l, g}$ for algebraic variables and algebraic equations approximately coincide with each other and are one unit less than the convergence orders for differential variables --- their convergence order is approximately equal to $N$. The worsening of the accuracy and convergence properties of the numerical ADER-DG method with a local DG predictor was expected for a DAE index higher than 1 --- this property is well known for other numerical methods for solving DAE systems~\cite{Hairer_book_2}, while in the case of the numerical ADER-DG method this phenomenon manifested itself in the form of a decrease in the empirical convergence orders $p$ compared to the expected convergence orders $p_{\rm nodes} = 2N+1$ and $p_{\rm local} = N+1$ (\ref{eq:expect_orders}) of the numerical method, which occurs when solving the initial value problem for ODE system.

The obtained results of the numerical solution of DAE system (\ref{eq:hess_dae_ind_2}) of index 1 are presented in Figs.~\ref{fig:hess_2_ind1_sol_uv}, \ref{fig:hess_2_ind1_sol_g_eps}, \ref{fig:hess_2_ind1_errors} and in Tables~\ref{tab:conv_orders_nodes_hess_2_ind1}, \ref{tab:conv_orders_local_hess_2_ind1}. Fig.~\ref{fig:hess_2_ind1_sol_uv} shows a comparison of the numerical solution at the nodes $(\mathbf{u}_{n}, \mathbf{v}_{n})$, the numerical local solution $(\mathbf{u}_{L}, \mathbf{v}_{L})$ and the exact analytical solution $(\mathbf{u}^{\rm ex}, \mathbf{v}^{\rm ex})$ separately for each differential $\mathbf{u}$ and algebraic $\mathbf{v}$ variable. Fig.~\ref{fig:hess_2_ind1_sol_g_eps} shows the dependencies of the feasibility of conditions $g_{1} = 0$ and $g_{2} = 0$ on the coordinate $t$, as well as the dependencies of the local errors $\varepsilon_{u}$, $\varepsilon_{v}$, $\varepsilon_{g}$ on the coordinate $t$, which allows us to quantitatively estimate the accuracy of the numerical solution, especially taking into account the fact that the numerical solution obtained by the ADER-DG method with a local DG predictor with a high degree of polynomials $N$ very accurately corresponds to the exact analytical solution, and it is impossible to visually identify the error from the solution plots separately (in Fig.~\ref{fig:hess_2_ind1_sol_uv}). Fig.~\ref{fig:hess_2_ind1_errors} shows the dependencies of the global errors $e^{u}$, $e^{v}$, $e^{g}$ of the numerical solution at the nodes $(\mathbf{u}_{n}, \mathbf{v}_{n})$ and the local solution $(\mathbf{u}_{L}, \mathbf{v}_{L})$ on the discretization step ${\Delta t}$, separately for each differential $\mathbf{u}$ and algebraic $\mathbf{v}$ variable and the algebraic equations $\mathbf{g} = \mathbf{0}$, on the basis of which the empirical convergence orders $p$ were calculated. 

The numerical solution of the DAE system (\ref{eq:hess_dae_ind_2}) of index 1, presented in Fig.~\ref{fig:hess_2_ind1_sol_uv}, clearly demonstrates the high accuracy achieved using the ADER-DG method with a local DG predictor. The obtained results of comparison of the numerical solution at the nodes and the local solution with the exact analytical solution, in general, correspond to the results from Example 1, presented in Subsection~\ref{sec:2:ct:ex1}, and Example 2, presented in Subsection~\ref{sec:2:ct:ex2}, where test examples of DAE of index 1 systems were solved. To quantitatively determine the accuracy of the local solution $(\mathbf{u}_{L}, \mathbf{v}_{L})$ and the solution at the nodes $(\mathbf{u}_{n}, \mathbf{v}_{n})$, plots of the local errors $|g_{1}|$, $|g_{2}|$ of satisfaction of the algebraic equations $g_{1} = 0$, $g_{2} = 0$ and local errors $\varepsilon_{u}$, $\varepsilon_{v}$, $\varepsilon_{g}$ were constructed, which were presented in Fig.~\ref{fig:hess_2_ind1_sol_g_eps}. From the presented dependencies of the error $|g_{2}|$, it was concluded that in the case of the solution at the nodes $(\mathbf{u}_{n}, \mathbf{v}_{n})$, the algebraic equation $g_{2} = 0$, which was explicitly included in the DAE system (\ref{eq:hess_dae_ind_2}) of index 1, is satisfied exactly (of course, within the accuracy of the representation of real numbers with a floating point). However, in the case of a local solution $(\mathbf{u}_{L}, \mathbf{v}_{L})$, the algebraic equation $g_{2} = 0$ is satisfied with a nonzero error $|g_{2}| \sim 10^{-7}$--$10^{-2}$ in the case of polynomial degree $N = 1$, $|g_{2}| \sim 10^{-19}$--$10^{-13}$ in the case of $N = 8$, $|g_{2}| \sim 10^{-83}$--$10^{-78}$ in the case of $N = 40$. In this case, the error $|g_{1}|$, determining the satisfaction of the algebraic equation $g_{1} = 0$, turned out to be non-zero, however, significantly less than the error $|g_{2}|$ in the case of DAE system (\ref{eq:hess_dae_ind_2}) of index 2. The obtained results confirm that the numerical ADER-DG method with a local LG predictor makes it possible to obtain a solution at the nodes $(\mathbf{u}_{n}, \mathbf{v}_{n})$ that exactly satisfies the algebraic equations explicitly included in the DAE system, while the algebraic equations derived from them (arising, for example, within the framework of decreasing the index of the DAE system) are satisfied with the accuracy characteristic of the local solution $(\mathbf{u}_{L}, \mathbf{v}_{L})$ and the solution at the nodes $(\mathbf{u}_{n}, \mathbf{v}_{n})$. Another important difference is the significant difference in the errors $\varepsilon_{u}$ of the local solution $\mathbf{u}_{n}$ and the solution at the nodes $\mathbf{u}_{L}$ for differential variables $\mathbf{u}$, which reaches $10$--$40$ times in the case of polynomial degrees of $N = 1$, $10^{4}$ times in the case of polynomial degrees of $N = 8$ and $10^{9}$ times in the case of polynomial degrees of $N = 40$, while there is no such difference for the errors $\varepsilon_{v}$, $\varepsilon_{g}$ of the local solution and the solution at the nodes for algebraic variables $\mathbf{v}$ and algebraic equations $\mathbf{g} = \mathbf{0}$.

Fig.~\ref{fig:hess_2_ind1_errors} shows the log-log dependencies of global errors $e$ for the solution at nodes $(\mathbf{u}_{n}, \mathbf{v}_{n})$ and the local solution $(\mathbf{u}_{L}, \mathbf{v}_{L})$ separately for differential and algebraic variables and algebraic equations. The presented dependencies clearly demonstrate the power law $e(\Delta t) \sim (\Delta t)^{p}$, from which the empirical convergence orders $p$ were calculated. The calculated empirical convergence orders $p^{n, u}$, $p^{n, v}$, $p^{n, g}$ for the solution at nodes $(\mathbf{u}_{n}, \mathbf{v}_{n})$ are presented in Table~\ref{tab:conv_orders_nodes_hess_2_ind1}, convergence orders $p^{l, u}$, $p^{l, v}$, $p^{l, g}$ for the local solution $(\mathbf{u}_{L}, \mathbf{v}_{L})$ --- in Table~\ref{tab:conv_orders_local_hess_2_ind1}. The convergence orders were calculated separately for the norms $L_{1}$, $L_{2}$, $L_{\infty}$ (\ref{eq:norms_def}). The empirical convergence orders in the case of the DAE system (\ref{eq:hess_dae_ind_2}) of index 1 practically coincide with the similar values in the case of the DAE system (\ref{eq:hess_dae_ind_2}) of index 2. Some difference is that the empirical convergence orders $p^{n, u}$ of the solution at the nodes $\mathbf{u}_{n}$ for differential variables are $1$--$3$ orders greater than $N+1$.

Therefore, it can be concluded that the ADER-DG method with a local DG predictor allows obtaining a highly accurate numerical solution of the DAE system of index 2, however, in this case, a significant decrease in the convergence order $p^{n}$ of the solution at the nodes $(\mathbf{u}_{n}, \mathbf{v}_{n})$ is observed --- the phenomenon of superconvergence disappears, the convergence orders were close to the values of $N+1$ expected (\ref{eq:expect_orders}) for the local solution  $(\mathbf{u}_{L}, \mathbf{v}_{L})$. An important feature of the numerical solution is the exact satisfaction of the algebraic equation of the corresponding index included in the DAE system. Decreasing the index of the DAE system to 1 in this specific example led to empirical convergence orders $p$ approximately equal to the convergence orders $p$ of the original DAE system of index 2.

\subsubsection{Example 4: mathematical pendulum}
\label{sec:2:ct:ex4}

The fourth example (which contains three sub-examples) of application of the numerical method ADER-DG with local DG predictor consisted in solving a well-known problem of mathematical pendulum, which is well known in the field of DAE systems, and which is a DAE system of index 3 in the original formulation of the problem:
\begin{equation}\label{eq:math_pend_dae_ind_3}
\begin{split}
&\ddot{x} = -\lambda x,\hspace{59.8mm} x(0) = \sin(\phi_{0}),\quad \dot{x}(0) = 0,\\
&\ddot{y} = -\lambda y - g,\hspace{54mm} y(0) = -\cos(\phi_{0}),\quad \dot{y}(0) = 0,\\
&g_{1} = x^{2} + y^{2} - 1 = 0\hspace{24mm} (\text{index}\ 3),\quad \lambda(0) = g\cos(\phi_{0}),\quad t\in[0, 10],\\
&g_{2} = x\dot{x} + y\dot{y} = 0\hspace{29.5mm} (\text{index}\ 2),\\
&g_{3} = \dot{x}^{2} + \dot{y}^{2} - \lambda(x^{2}+y^{2}) - gy = 0\hspace{2mm} (\text{index}\ 1),
\end{split}
\end{equation}
where the initial condition was given by the initial angle of the pendulum $\phi_{0}$. This test example is the first example in this paper for which the numerical solution of the DAE of index 3 was investigated. The presented DAE system of index 3 contains only one algebraic equation (constraint) $g_{1} = 0$. 

The exact analytical solution to the problem (\ref{eq:math_pend_dae_ind_3}) was determined through the dynamic dependence of the angle $\phi(t)$, which determines the deviation of the position of a mathematical pendulum from the position of stable mechanical equilibrium:
\begin{equation}
\begin{split}
&x = \sin(\phi(t)),\quad
 y = -\cos(\phi(t)),\\
&\dot{x} = \dot{\phi}(t)\cos(\phi(t)),\quad
 \dot{y} = \dot{\phi}(t)\sin(\phi(t)),\\
&\lambda = \dot{\phi}(t)^{2} + g\cos(\phi(t)),
\end{split}
\end{equation}
where the exact analytical solution for the function $\phi(t)$ was obtained by trivial integration of the ODE problem $\ddot{\phi} + \omega_{0}^{2}\sin(\phi) = 0$ with initial conditions $\phi(0) = \phi_{0}$, $\dot{\phi}(0) = 0$ (see, for example, interesting derivation and detailed analysis of the analytical solution in the work~\cite{math_pend_exact_sol_ref}):
\begin{equation}
\begin{split}
&\phi(t) = 2\cdot\arcsin\left\{\Upsilon\cdot
	\mathrm{sn}\left(K\left(\Upsilon\right) - \omega_{0}t,\, \Upsilon\right)\right\},\\
&\dot{\phi}(t) = \frac{\displaystyle
	-2\omega_{0}\Upsilon\cdot
	\mathrm{cn}\left(K\left(\Upsilon\right) - \omega_{0}t,\, \Upsilon\right)\cdot
	\mathrm{dn}\left(K\left(\Upsilon\right) - \omega_{0}t,\, \Upsilon\right)
}{\displaystyle
	\sqrt{1 - \Upsilon^{2}\cdot\mathrm{sn}^{2}\left(K\left(\Upsilon\right) - \omega_{0}t,\, \Upsilon\right)}
},\\
&\Upsilon = \sin\left(\frac{\phi_{0}}{2}\right),\qquad
\omega_{0} = \sqrt{g},
\end{split}
\end{equation}
where $K$ is the complete elliptical integral of the first kind, $(\mathrm{sn},\, \mathrm{cn},\, \mathrm{dn})$ are the set of Jacobi elliptic functions: sine, cosine, delta amplitude, respectively. It should be noted that the problem solved in this Subsection is nonlinear --- the presented solution defines the law of motion of a mathematical pendulum for an arbitrary amplitude of oscillations (limited only by the zero value of the initial velocity $\dot{\phi}(0) = 0$, therefore, it is impossible to obtain a full rotation in this formulation of the problem), and not only for small values of amplitude, when it would be possible to limit oneself only to the harmonic approximation of small oscillations. 

The problem (\ref{eq:math_pend_dae_ind_3}) was rewritten in a form consistent with the formulation of the original DAE system (\ref{eq:dae_chosen_form}):
\begin{equation}
\begin{split}
&\frac{du_{1}}{dt} = u_{3},\hspace{63mm} u_{1}(0) = \sin(\phi_{0}),\\
&\frac{du_{2}}{dt} = u_{4},\hspace{63mm} u_{2}(0) = -\cos(\phi_{0}),\\
&\frac{du_{3}}{dt} = -u_{1}v_{1},\hspace{57mm} u_{3}(0) = 0,\\
&\frac{du_{4}}{dt} = -u_{2}v_{1} - g,\hspace{51mm}  u_{4}(0) = 0,\\[1.9mm]
&g_{1} = u_{1}^{2} + u_{2}^{2} - 1 = 0\hspace{28.2mm} (\text{index}\ 3),\quad v_{1}(0) = g\cos(\phi_{0}),\\[2.4mm]
&g_{2} = u_{1}u_{3} + u_{2}u_{4} = 0\hspace{27.2mm} (\text{index}\ 2),\\[2.0mm]
&g_{3} = u_{3}^{2} + u_{4}^{2} - v_{1}(u_{1}^{2} + u_{2}^{2}) - gu_{2} = 0\hspace{3mm} (\text{index}\ 1),\\
\end{split}
\end{equation}
where sets of differential variables $\mathbf{u} = [x,\, y,\, \dot{x},\, \dot{y}]^{T}$, $\mathbf{v} = [\lambda]$ were defined. The full-component exact analytical solution of the problem was written in the following form:
\begin{equation}
\begin{split}
&\mathbf{u}^{\rm ex} = \left[
\begin{array}{l}
\sin(\phi(t))\\
-\cos(\phi(t))\\
\dot{\phi}(t)\cos(\phi(t))\\
\dot{\phi}(t)\sin(\phi(t))
\end{array}
\right],\\
&\mathbf{v}^{\rm ex} = \Big[\dot{\phi}(t)^{2} + g\cos(\phi(t))\Big].
\end{split}
\end{equation}
The value $g = 1$ was chosen in the calculations. The domain of definition $[0,\, 10]$ of the desired functions $\mathbf{u}$ and $\mathbf{v}$ was discretized into $L = 30$, $32$, $34$, $36$, $38$, $40$ discretization domains $\Omega_{n}$ in the case of polynomial degrees $N \leqslant 5$ and into $L = 10$, $12$, $14$, $16$, $18$, $20$ discretization domains $\Omega_{n}$ in the case of polynomial degrees $N > 5$. The discretization step ${\Delta t}_{n}$ was chosen equal ${\Delta t} = 10/(L-1)$ for all discretization domains $\Omega_{n}$ on the grid, which was done to be able to calculate the empirical convergence orders $p$. Therefore, the empirical convergence orders $p$ were calculated by least squares approximation of the dependence of the global error $e$ on the grid discretization step ${\Delta t}$ at 6 data points.

\begin{figure}[h!]
\captionsetup[subfigure]{%
	position=bottom,
	font+=smaller,
	textfont=normalfont,
	singlelinecheck=off,
	justification=raggedright
}
\centering
\begin{subfigure}{0.320\textwidth}
\includegraphics[width=\textwidth]{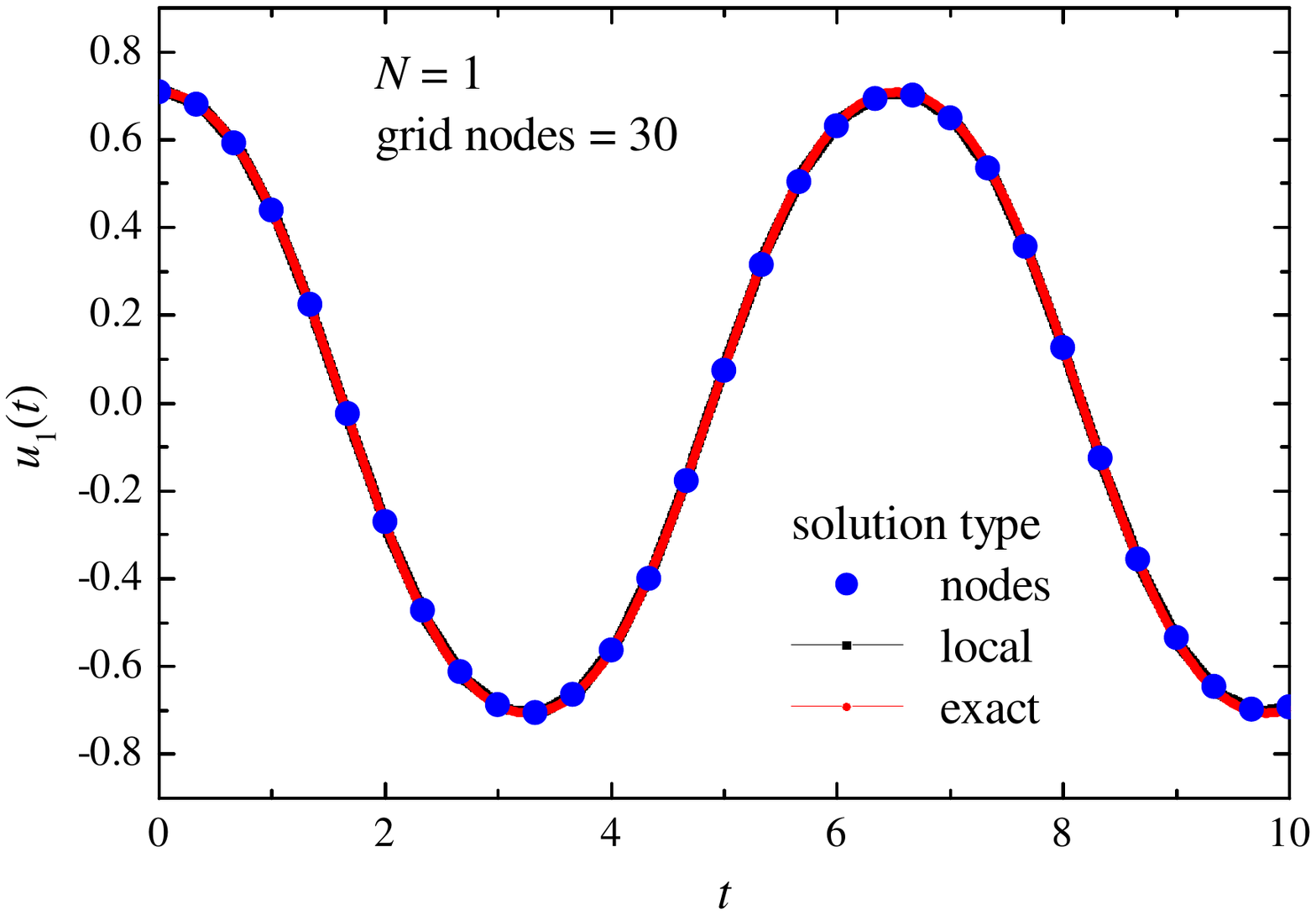}
\vspace{-8mm}\caption{\label{fig:pend_ind3_sol_uv:a1}}
\end{subfigure}
\begin{subfigure}{0.320\textwidth}
\includegraphics[width=\textwidth]{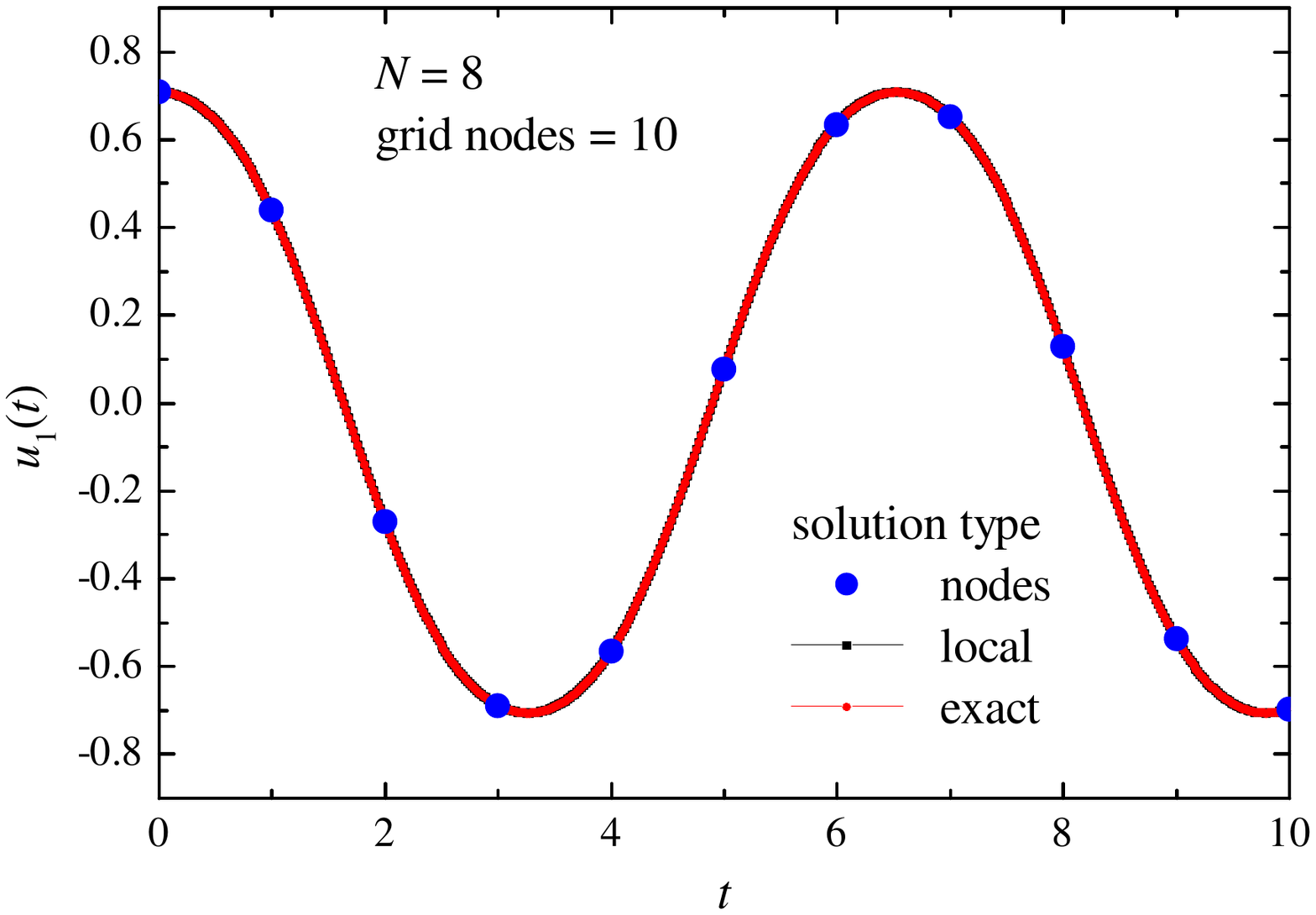}
\vspace{-8mm}\caption{\label{fig:pend_ind3_sol_uv:a2}}
\end{subfigure}
\begin{subfigure}{0.320\textwidth}
\includegraphics[width=\textwidth]{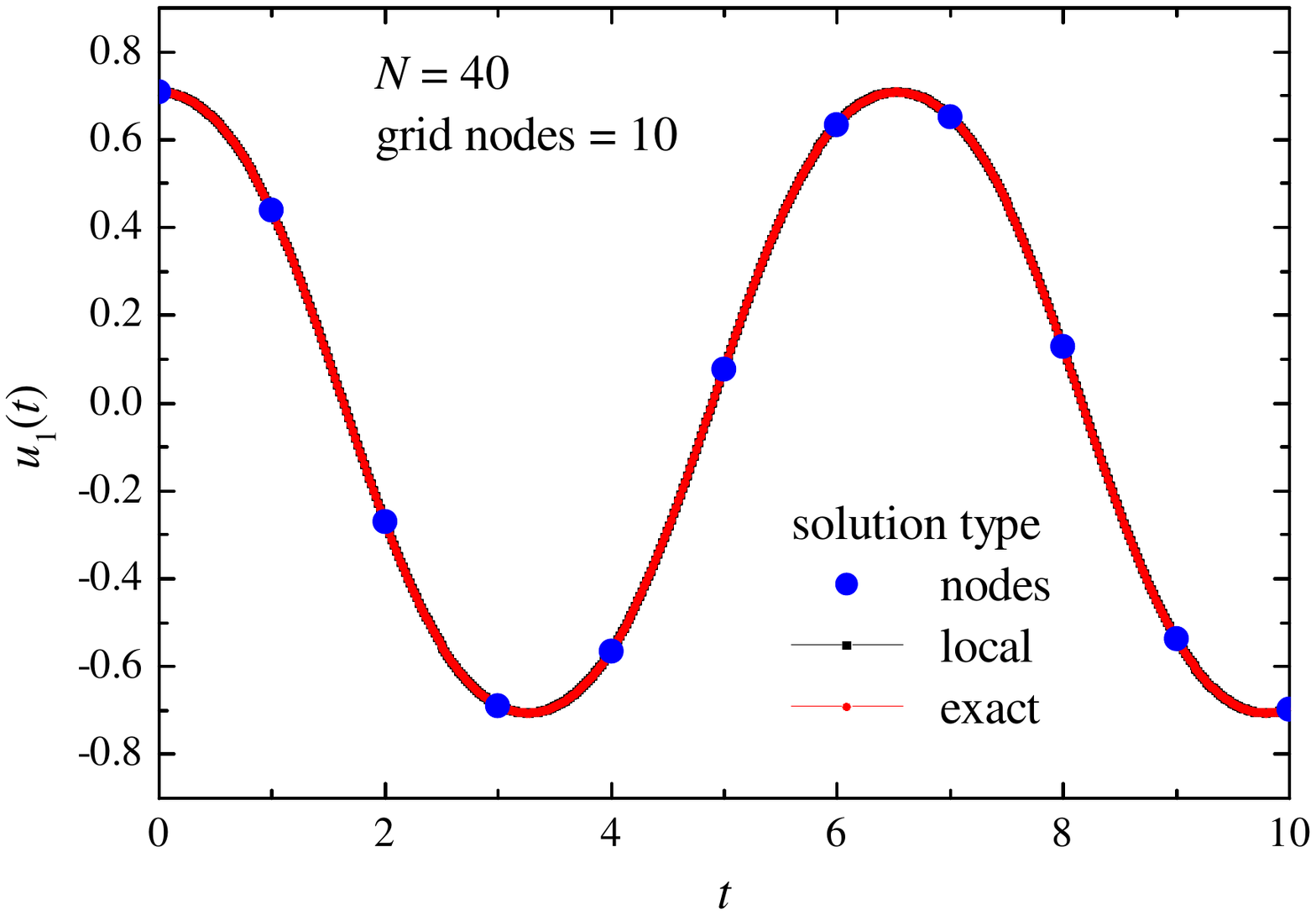}
\vspace{-8mm}\caption{\label{fig:pend_ind3_sol_uv:a3}}
\end{subfigure}\\
\begin{subfigure}{0.320\textwidth}
\includegraphics[width=\textwidth]{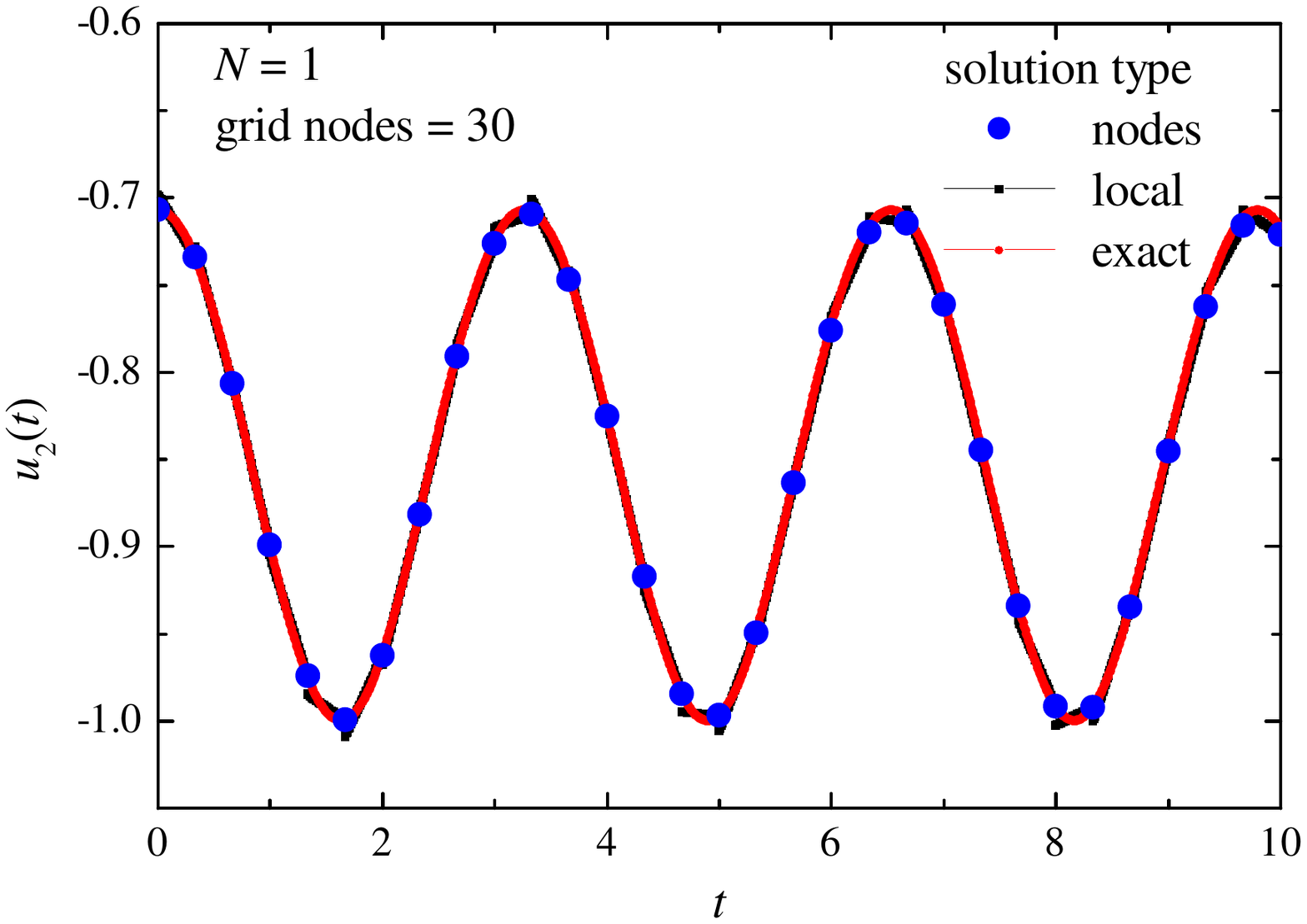}
\vspace{-8mm}\caption{\label{fig:pend_ind3_sol_uv:b1}}
\end{subfigure}
\begin{subfigure}{0.320\textwidth}
\includegraphics[width=\textwidth]{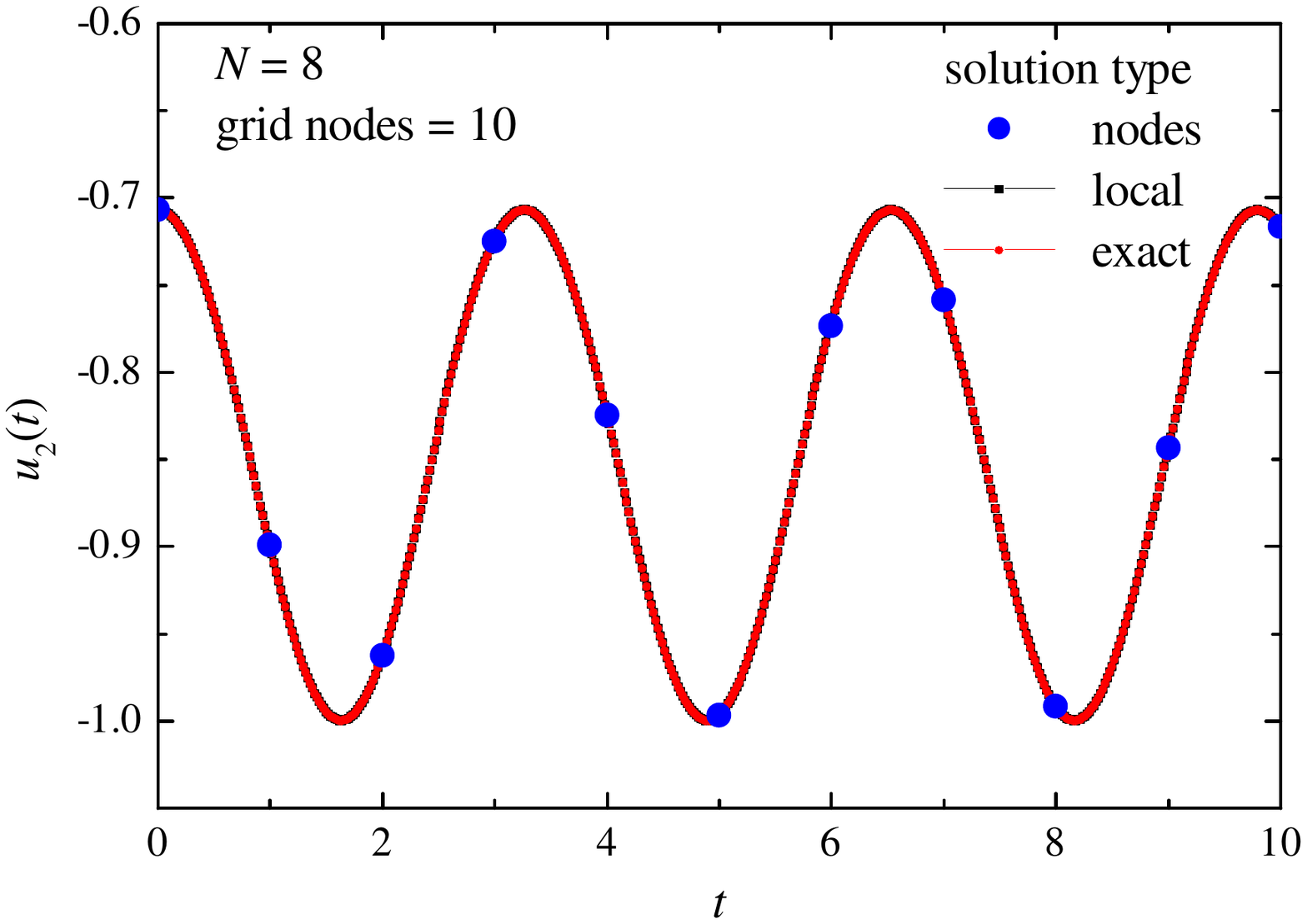}
\vspace{-8mm}\caption{\label{fig:pend_ind3_sol_uv:b2}}
\end{subfigure}
\begin{subfigure}{0.320\textwidth}
\includegraphics[width=\textwidth]{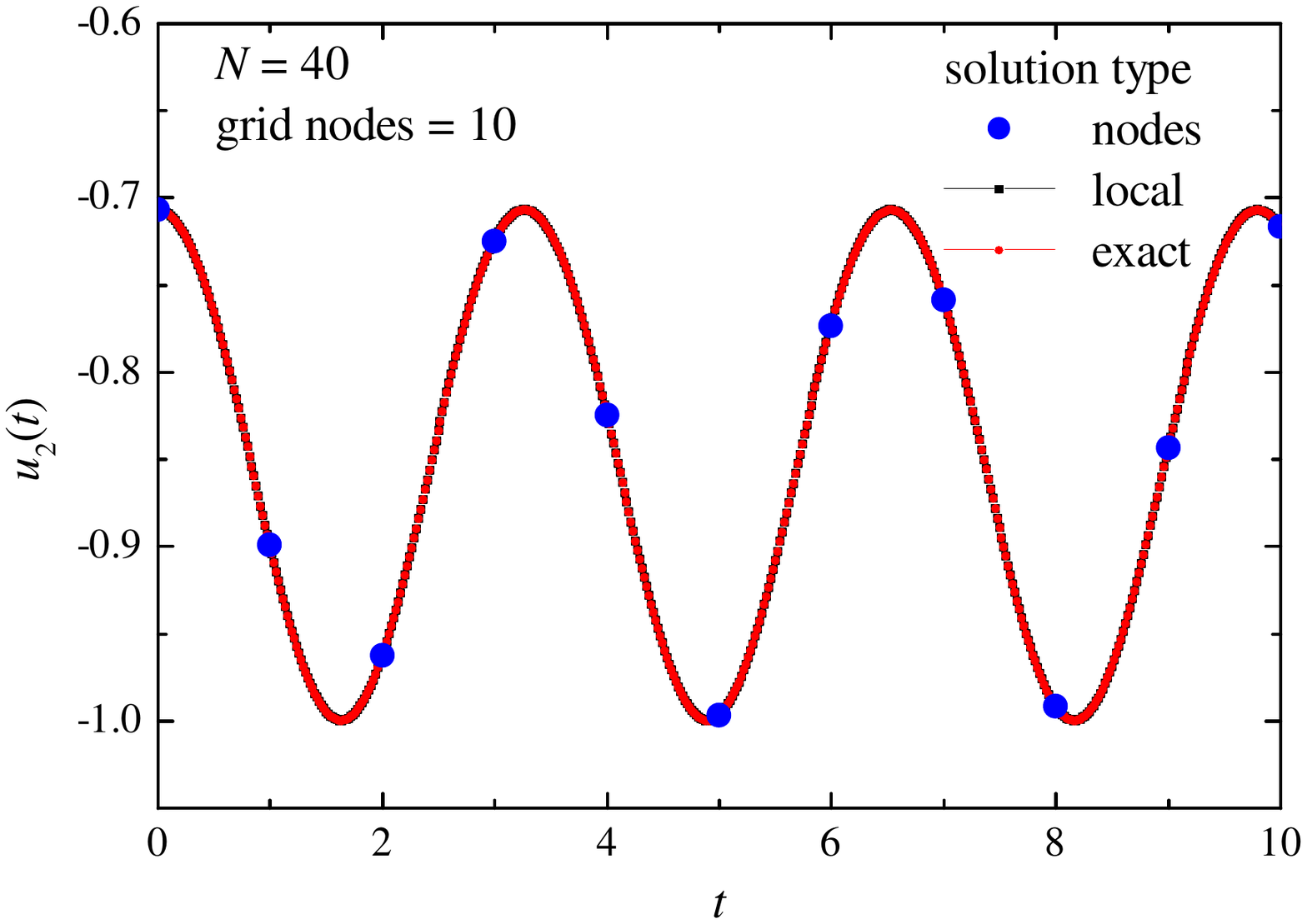}
\vspace{-8mm}\caption{\label{fig:pend_ind3_sol_uv:b3}}
\end{subfigure}\\
\begin{subfigure}{0.320\textwidth}
\includegraphics[width=\textwidth]{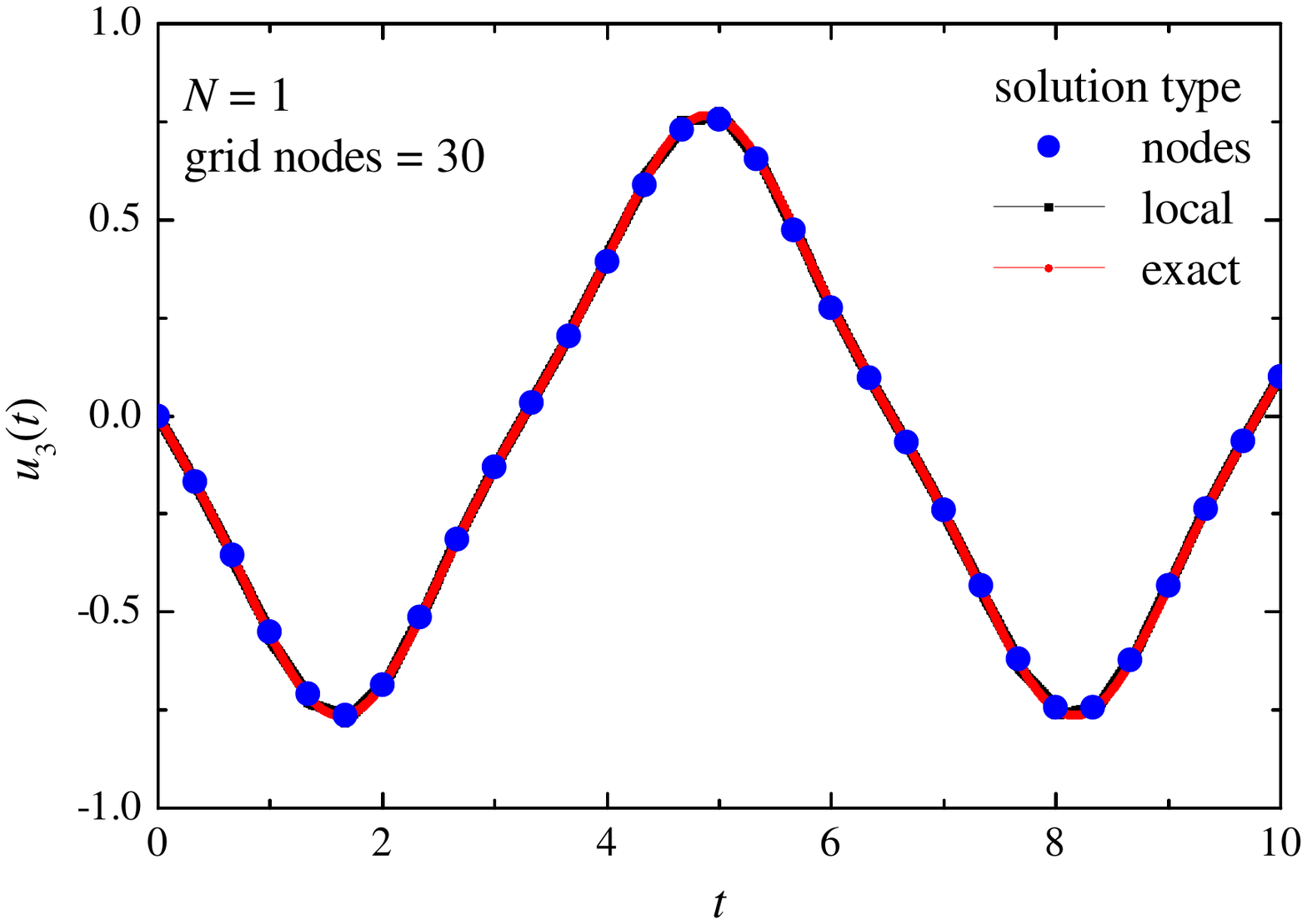}
\vspace{-8mm}\caption{\label{fig:pend_ind3_sol_uv:c1}}
\end{subfigure}
\begin{subfigure}{0.320\textwidth}
\includegraphics[width=\textwidth]{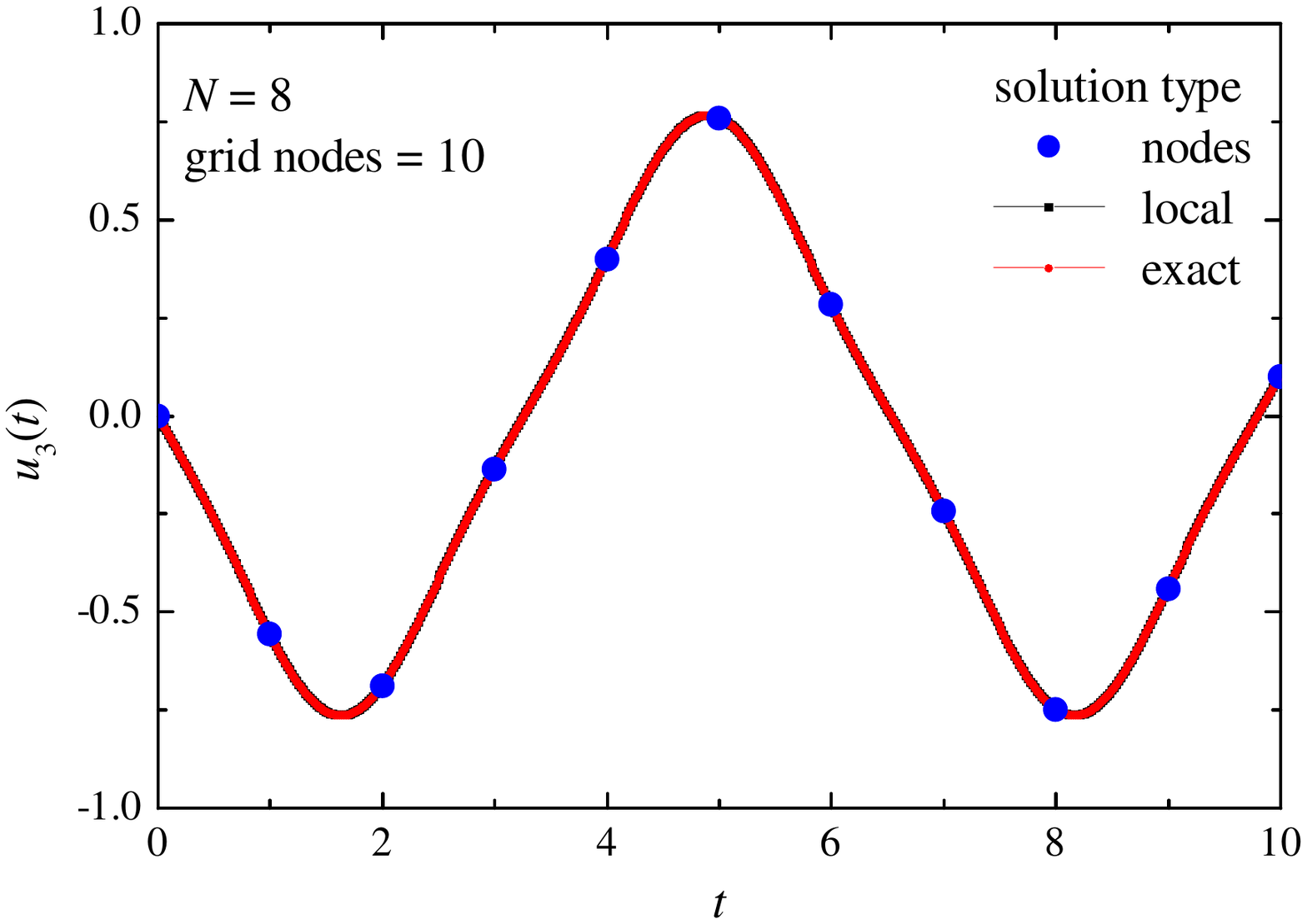}
\vspace{-8mm}\caption{\label{fig:pend_ind3_sol_uv:c2}}
\end{subfigure}
\begin{subfigure}{0.320\textwidth}
\includegraphics[width=\textwidth]{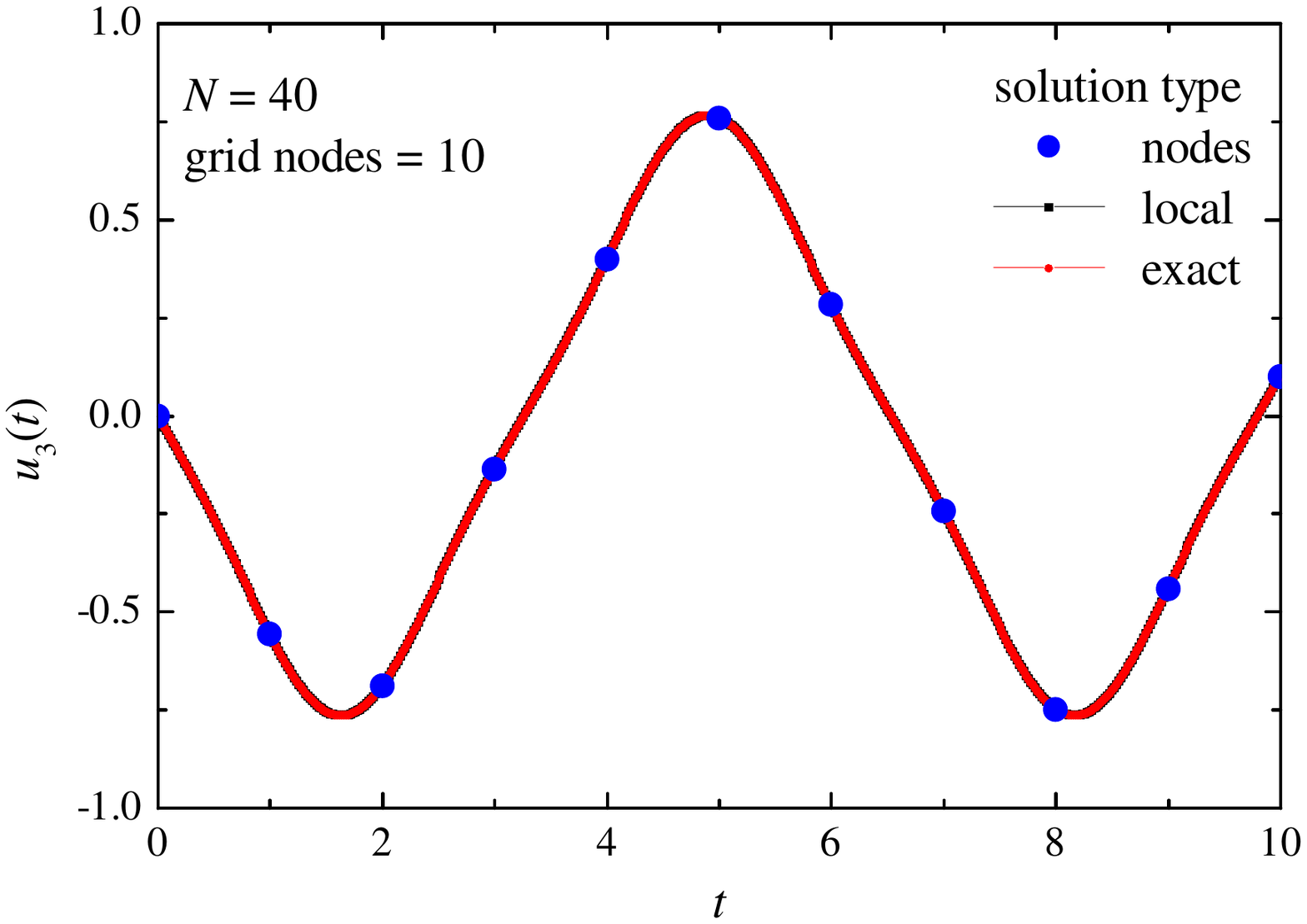}
\vspace{-8mm}\caption{\label{fig:pend_ind3_sol_uv:c3}}
\end{subfigure}\\
\begin{subfigure}{0.320\textwidth}
\includegraphics[width=\textwidth]{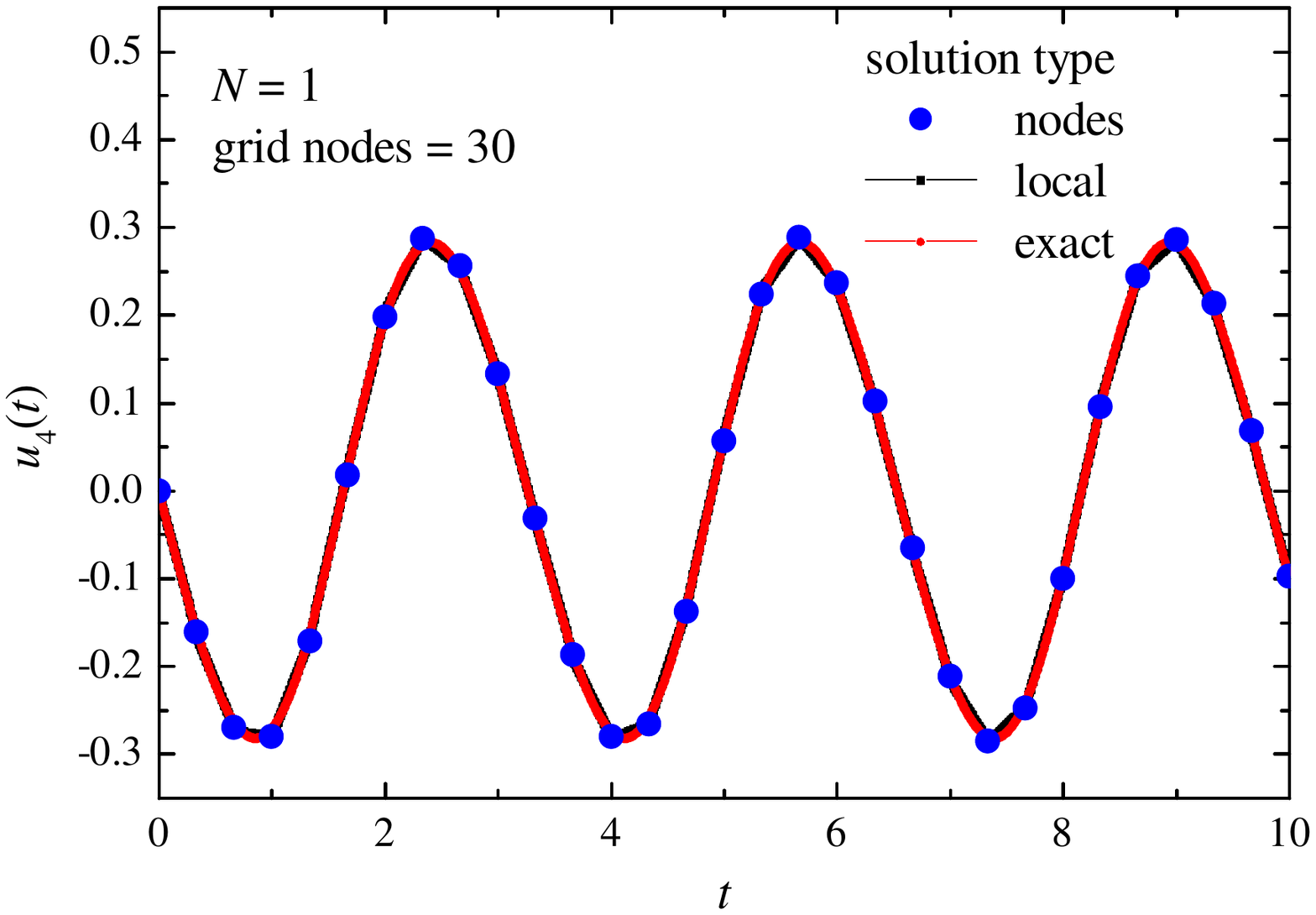}
\vspace{-8mm}\caption{\label{fig:pend_ind3_sol_uv:d1}}
\end{subfigure}
\begin{subfigure}{0.320\textwidth}
\includegraphics[width=\textwidth]{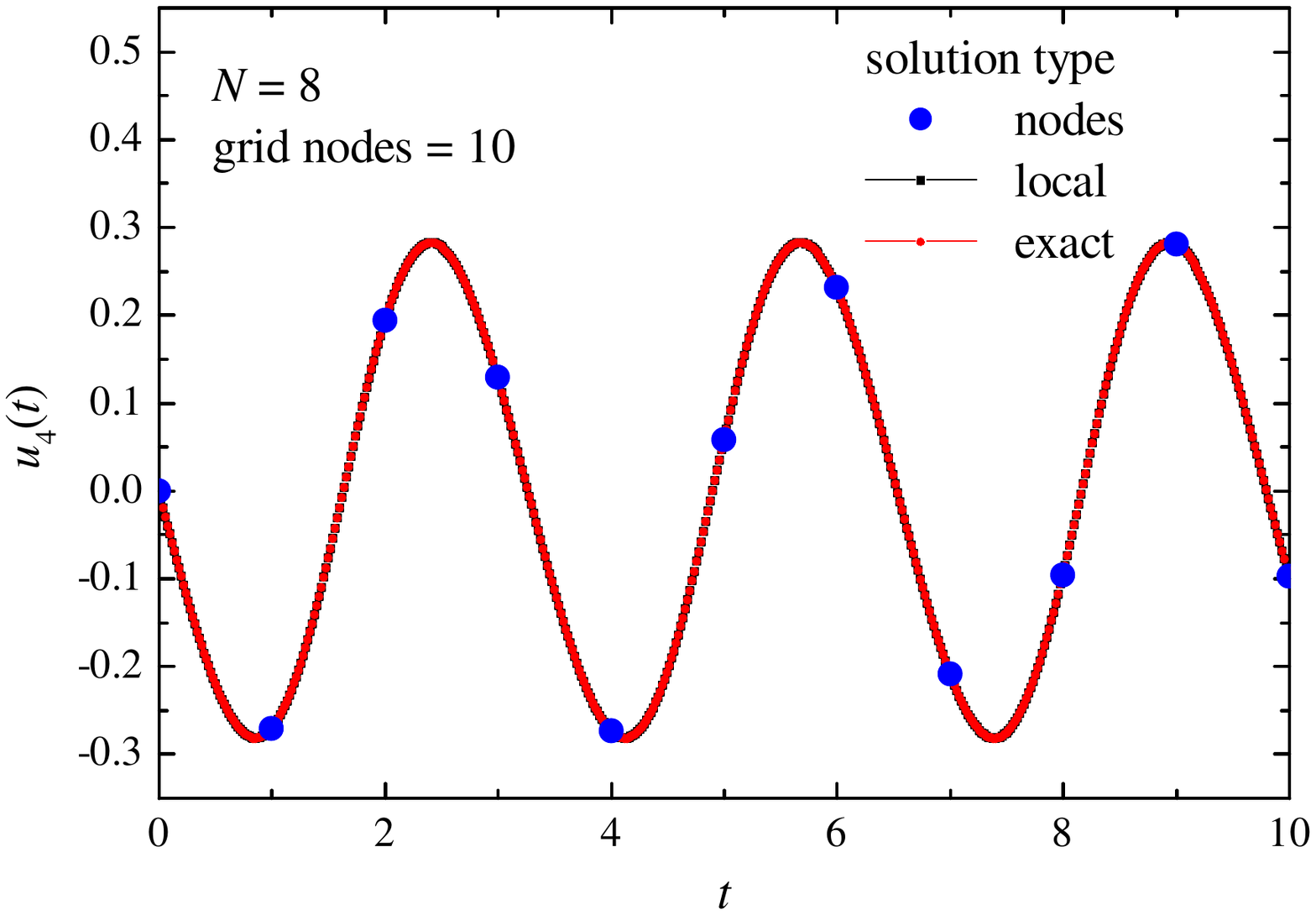}
\vspace{-8mm}\caption{\label{fig:pend_ind3_sol_uv:d2}}
\end{subfigure}
\begin{subfigure}{0.320\textwidth}
\includegraphics[width=\textwidth]{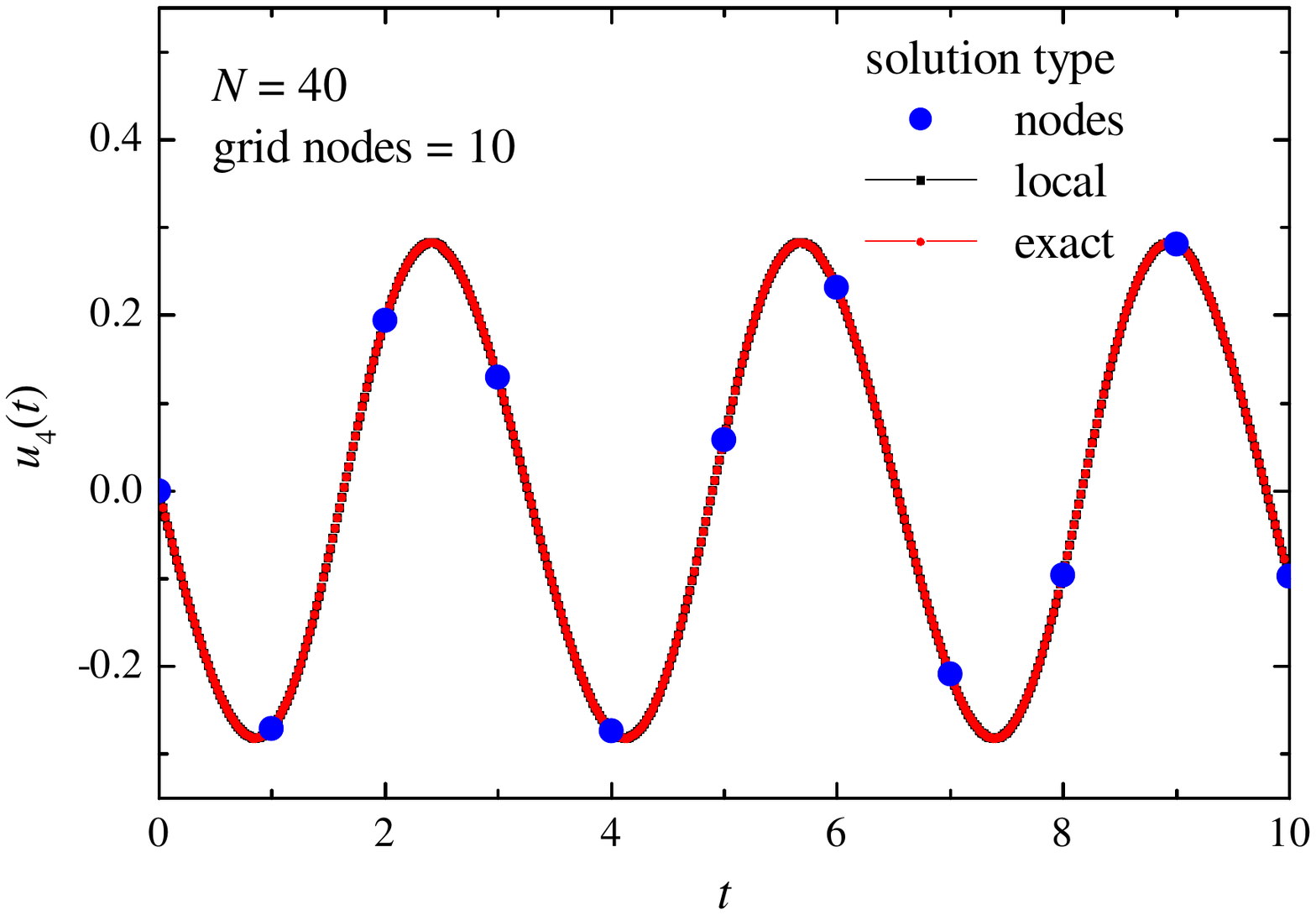}
\vspace{-8mm}\caption{\label{fig:pend_ind3_sol_uv:d3}}
\end{subfigure}\\
\begin{subfigure}{0.320\textwidth}
\includegraphics[width=\textwidth]{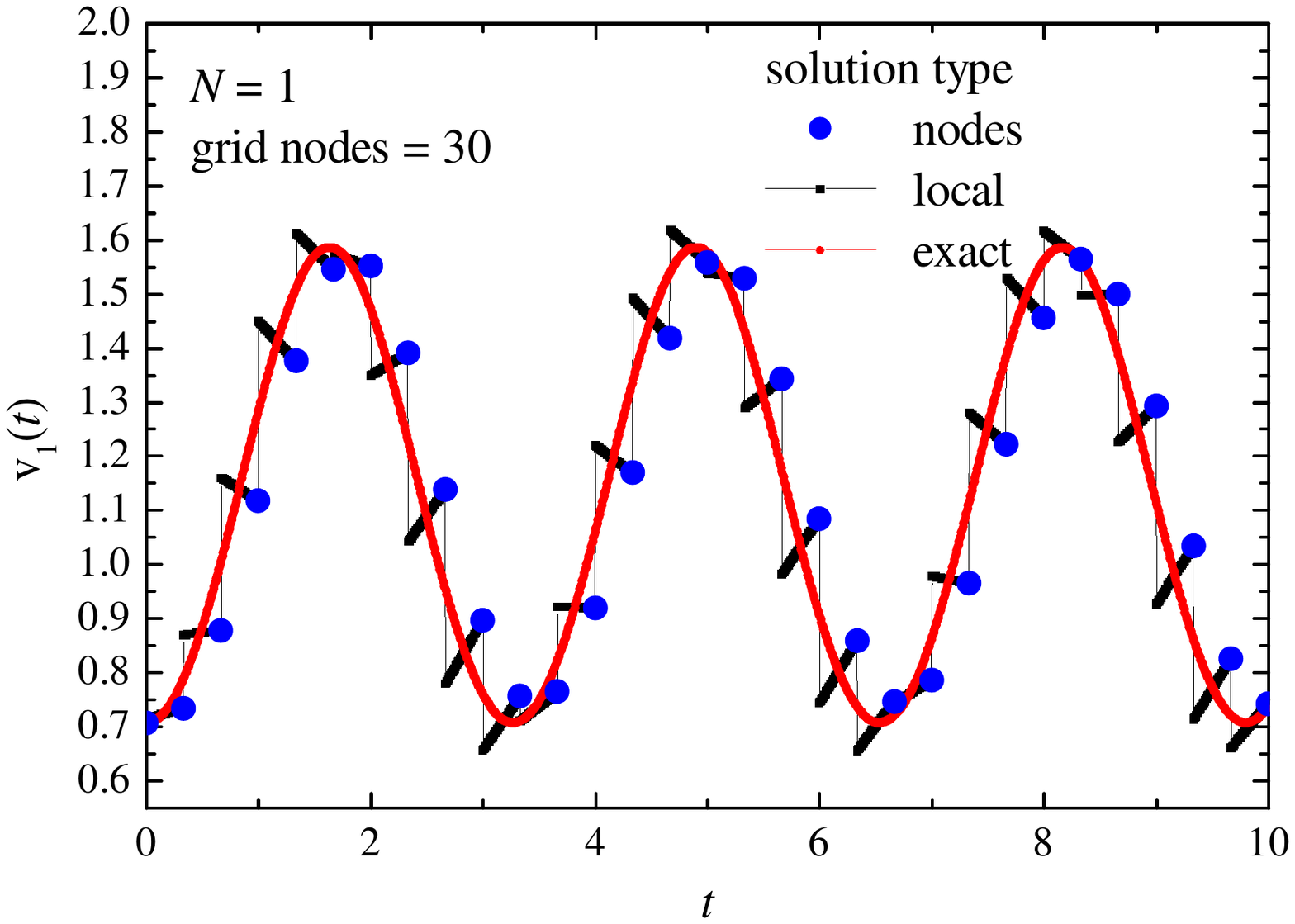}
\vspace{-8mm}\caption{\label{fig:pend_ind3_sol_uv:e1}}
\end{subfigure}
\begin{subfigure}{0.320\textwidth}
\includegraphics[width=\textwidth]{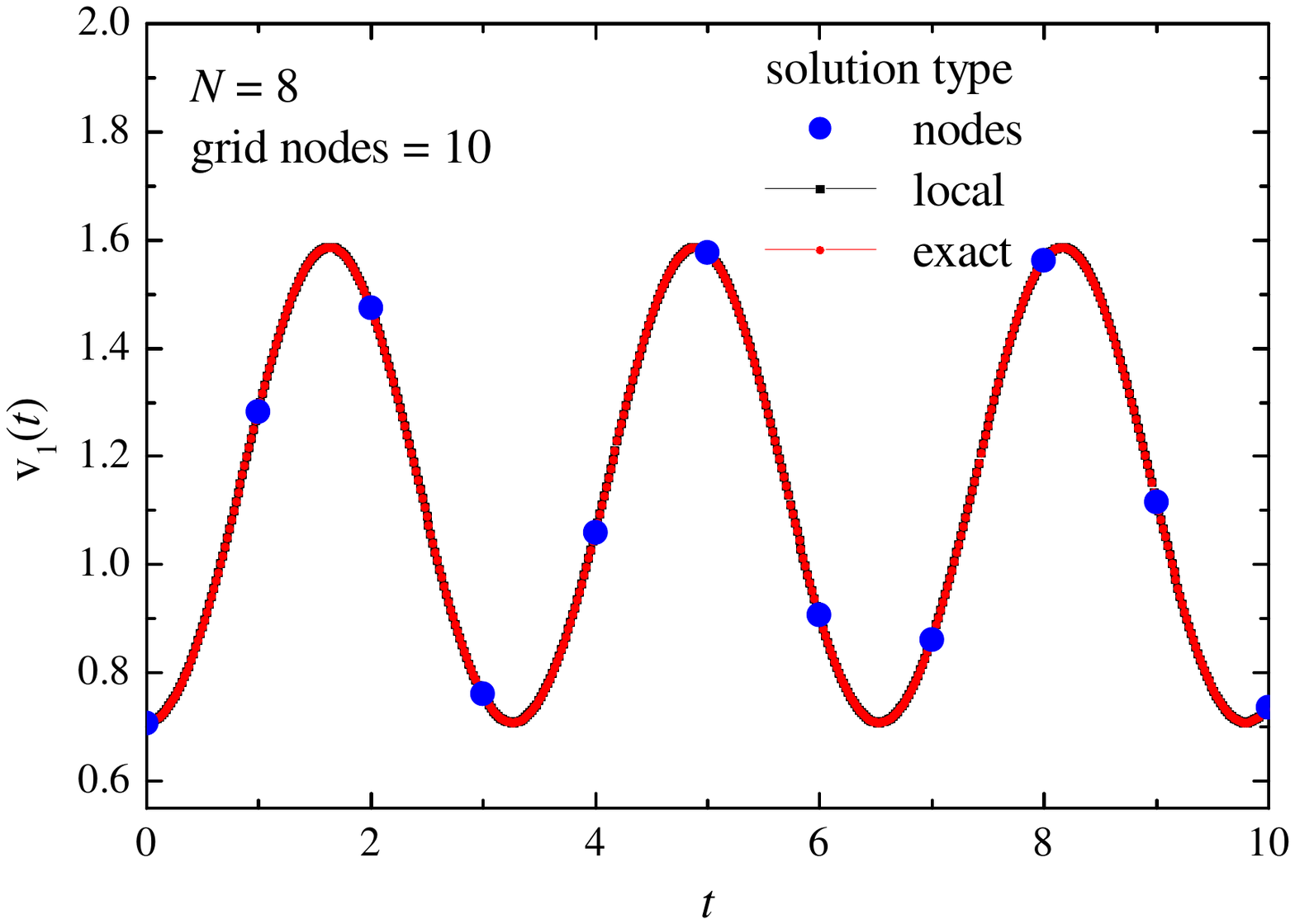}
\vspace{-8mm}\caption{\label{fig:pend_ind3_sol_uv:e2}}
\end{subfigure}
\begin{subfigure}{0.320\textwidth}
\includegraphics[width=\textwidth]{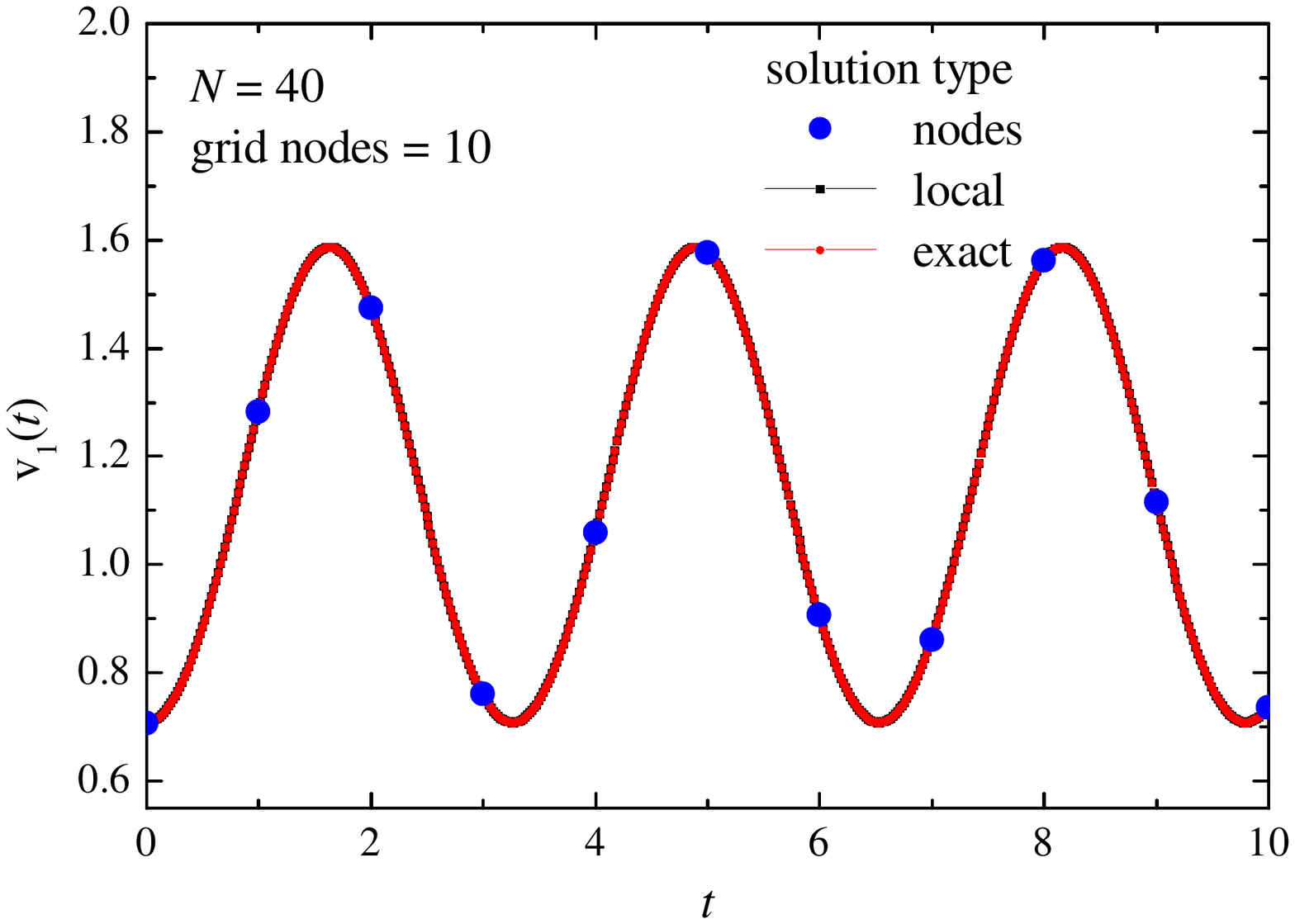}
\vspace{-8mm}\caption{\label{fig:pend_ind3_sol_uv:e3}}
\end{subfigure}\\
\caption{%
Numerical solution of the DAE system (\ref{eq:math_pend_dae_ind_3}) of index 3. Comparison of the solution at nodes $\mathbf{u}_{n}$, the local solution $\mathbf{u}_{L}(t)$ and the exact solution $\mathbf{u}^{\rm ex}(t)$ for components $u_{1}$ (\subref{fig:pend_ind3_sol_uv:a1}, \subref{fig:pend_ind3_sol_uv:a2}, \subref{fig:pend_ind3_sol_uv:a3}), $u_{2}$ (\subref{fig:pend_ind3_sol_uv:b1}, \subref{fig:pend_ind3_sol_uv:b2}, \subref{fig:pend_ind3_sol_uv:b3}), $u_{3}$ (\subref{fig:pend_ind3_sol_uv:c1}, \subref{fig:pend_ind3_sol_uv:c2}, \subref{fig:pend_ind3_sol_uv:c3}), $u_{4}$ (\subref{fig:pend_ind3_sol_uv:d1}, \subref{fig:pend_ind3_sol_uv:d2}, \subref{fig:pend_ind3_sol_uv:d3}) and $v_{1}$ (\subref{fig:pend_ind3_sol_uv:e1}, \subref{fig:pend_ind3_sol_uv:e2}, \subref{fig:pend_ind3_sol_uv:e3}), obtained using polynomials with degrees $N = 1$ (\subref{fig:pend_ind3_sol_uv:a1}, \subref{fig:pend_ind3_sol_uv:b1}, \subref{fig:pend_ind3_sol_uv:c1}, \subref{fig:pend_ind3_sol_uv:d1}, \subref{fig:pend_ind3_sol_uv:e1}), $N = 8$ (\subref{fig:pend_ind3_sol_uv:a2}, \subref{fig:pend_ind3_sol_uv:b2}, \subref{fig:pend_ind3_sol_uv:c2}, \subref{fig:pend_ind3_sol_uv:d2}, \subref{fig:pend_ind3_sol_uv:e2}) and $N = 40$ (\subref{fig:pend_ind3_sol_uv:a3}, \subref{fig:pend_ind3_sol_uv:b3}, \subref{fig:pend_ind3_sol_uv:c3}, \subref{fig:pend_ind3_sol_uv:d3}, \subref{fig:pend_ind3_sol_uv:e3}).
}
\label{fig:pend_ind3_sol_uv}
\end{figure} 

\begin{figure}[h!]
\captionsetup[subfigure]{%
	position=bottom,
	font+=smaller,
	textfont=normalfont,
	singlelinecheck=off,
	justification=raggedright
}
\centering
\begin{subfigure}{0.275\textwidth}
\includegraphics[width=\textwidth]{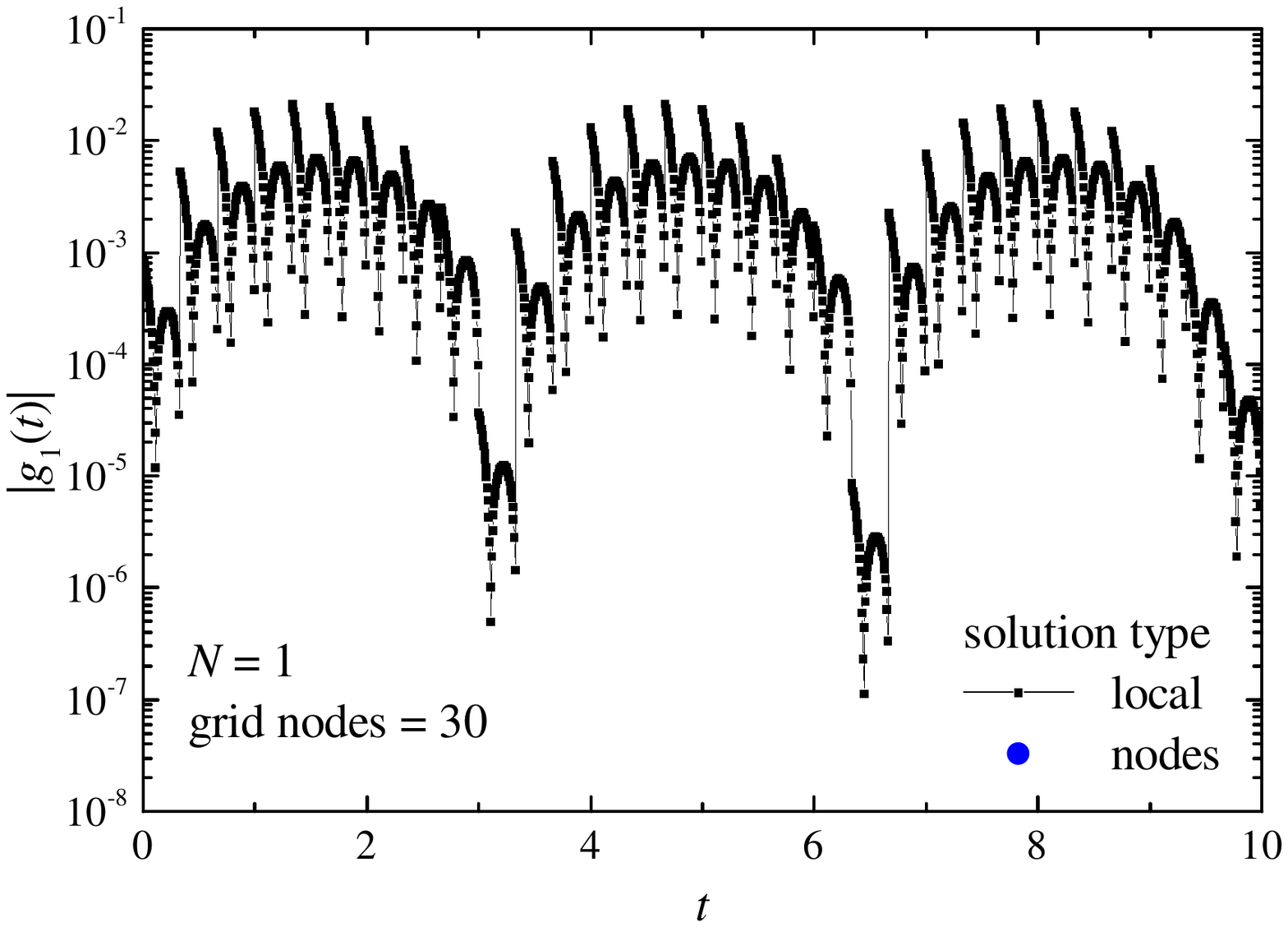}
\vspace{-8mm}\caption{\label{fig:pend_ind3_sol_g_eps:a1}}
\end{subfigure}\hspace{6mm}
\begin{subfigure}{0.275\textwidth}
\includegraphics[width=\textwidth]{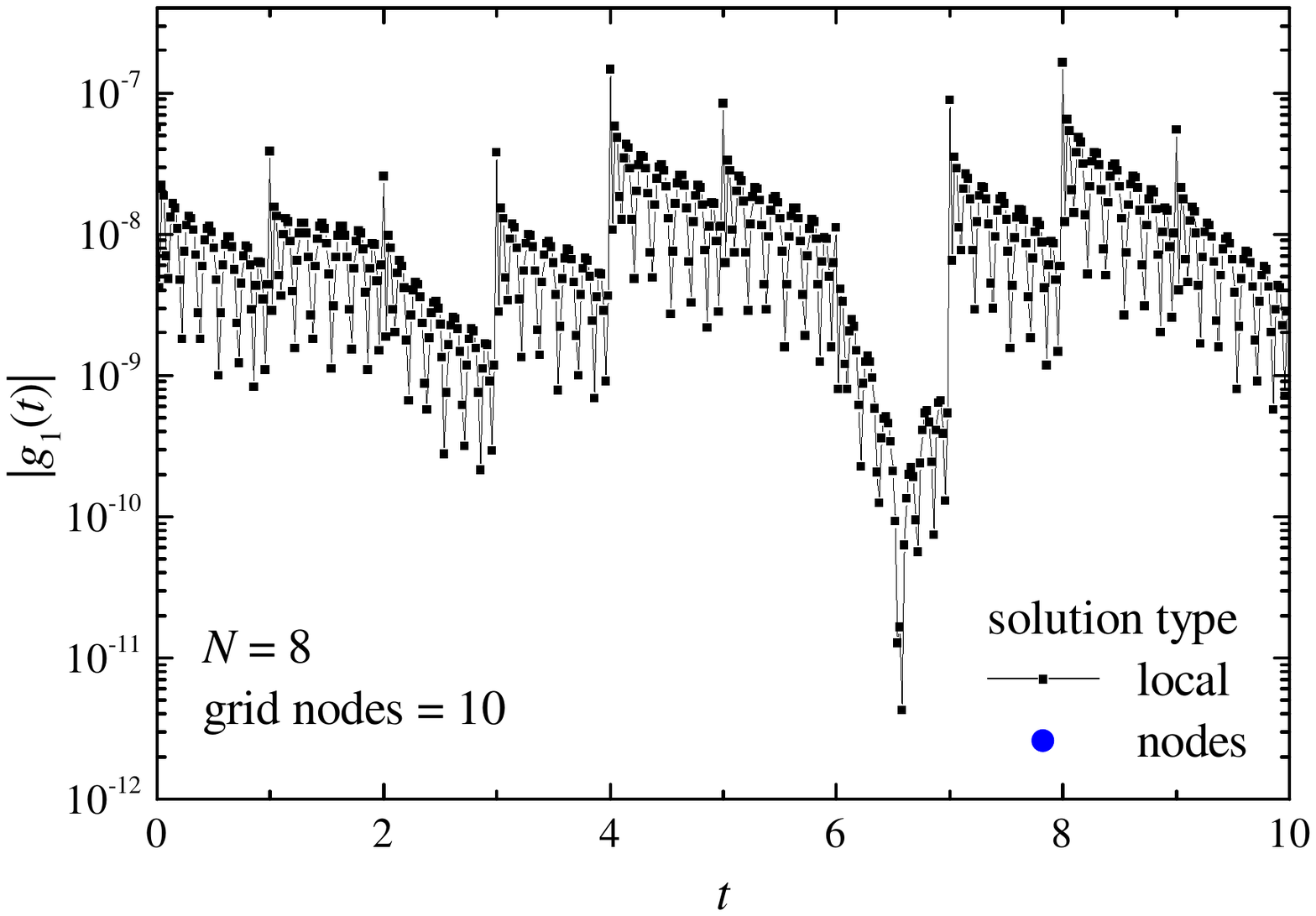}
\vspace{-8mm}\caption{\label{fig:pend_ind3_sol_g_eps:a2}}
\end{subfigure}\hspace{6mm}
\begin{subfigure}{0.275\textwidth}
\includegraphics[width=\textwidth]{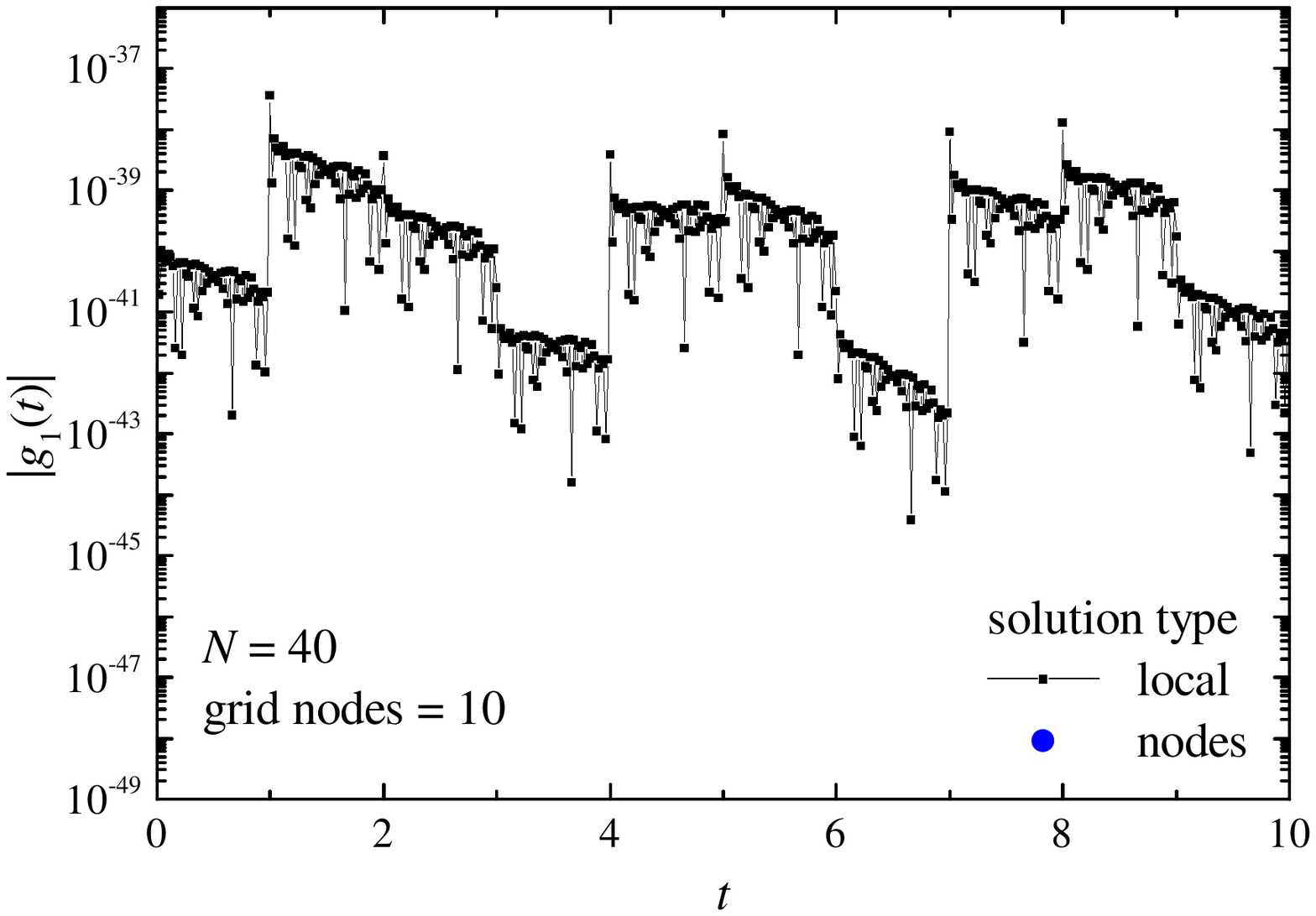}
\vspace{-8mm}\caption{\label{fig:pend_ind3_sol_g_eps:a3}}
\end{subfigure}\\[-2mm]
\begin{subfigure}{0.275\textwidth}
\includegraphics[width=\textwidth]{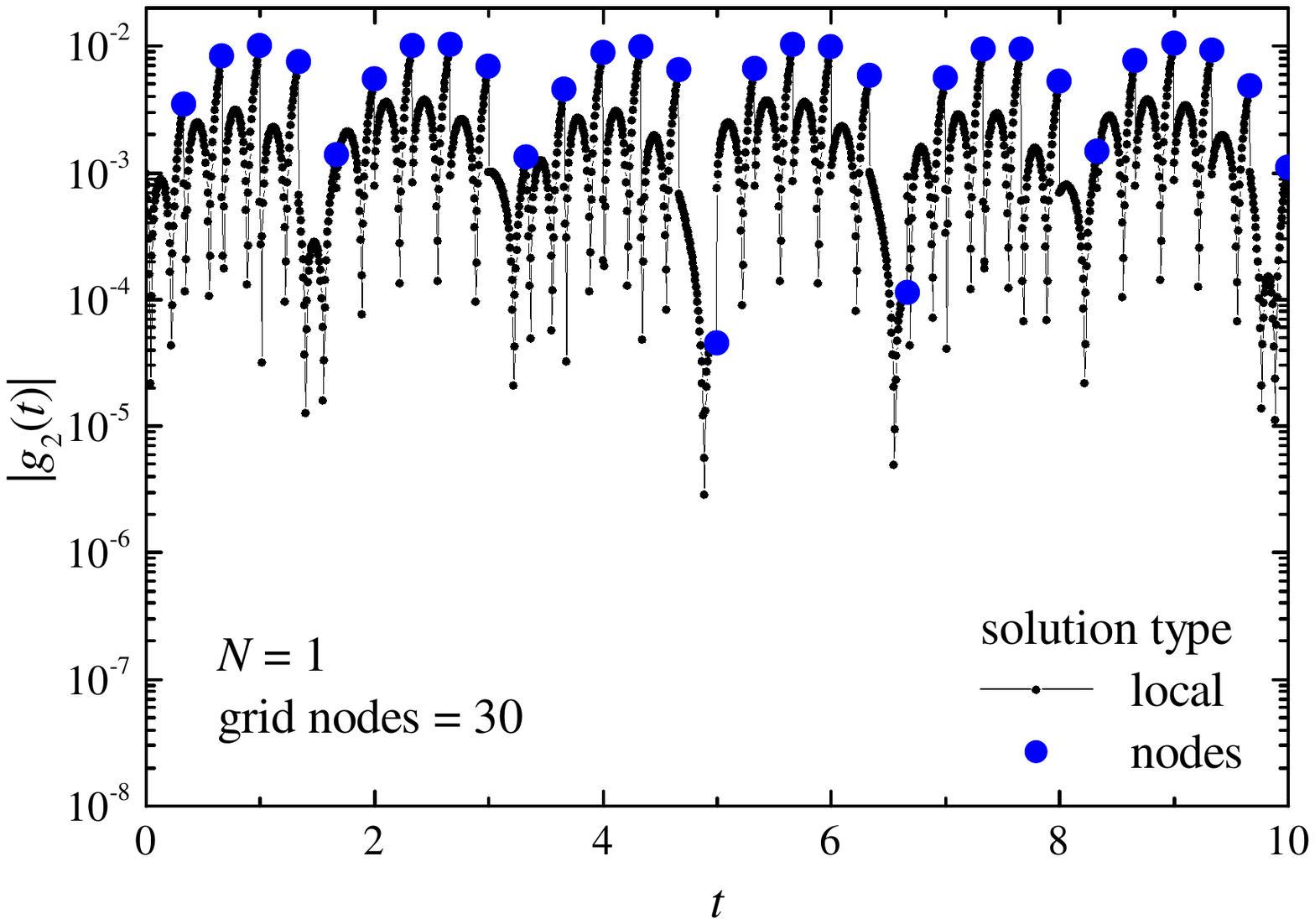}
\vspace{-8mm}\caption{\label{fig:pend_ind3_sol_g_eps:b1}}
\end{subfigure}\hspace{6mm}
\begin{subfigure}{0.275\textwidth}
\includegraphics[width=\textwidth]{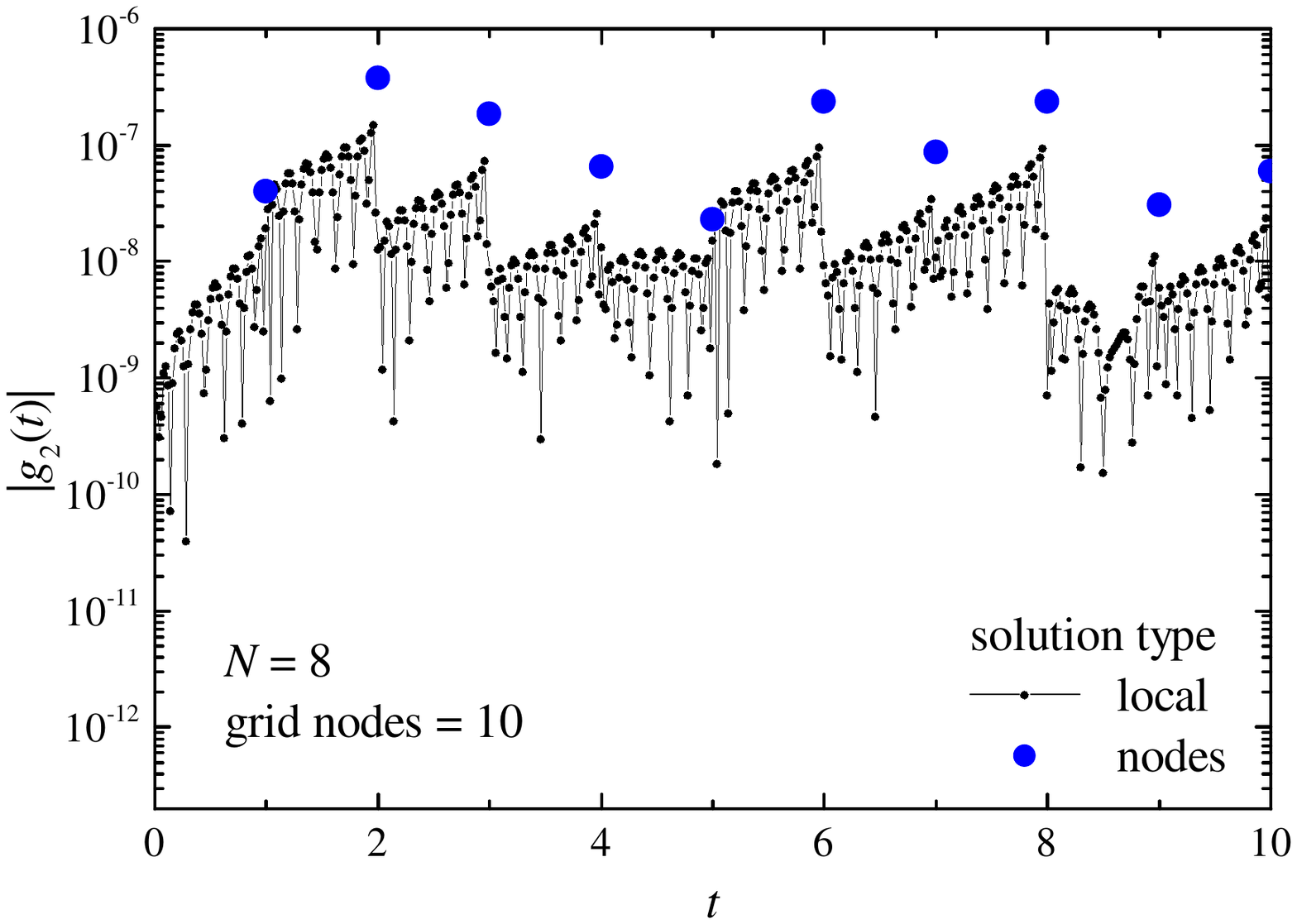}
\vspace{-8mm}\caption{\label{fig:pend_ind3_sol_g_eps:b2}}
\end{subfigure}\hspace{6mm}
\begin{subfigure}{0.275\textwidth}
\includegraphics[width=\textwidth]{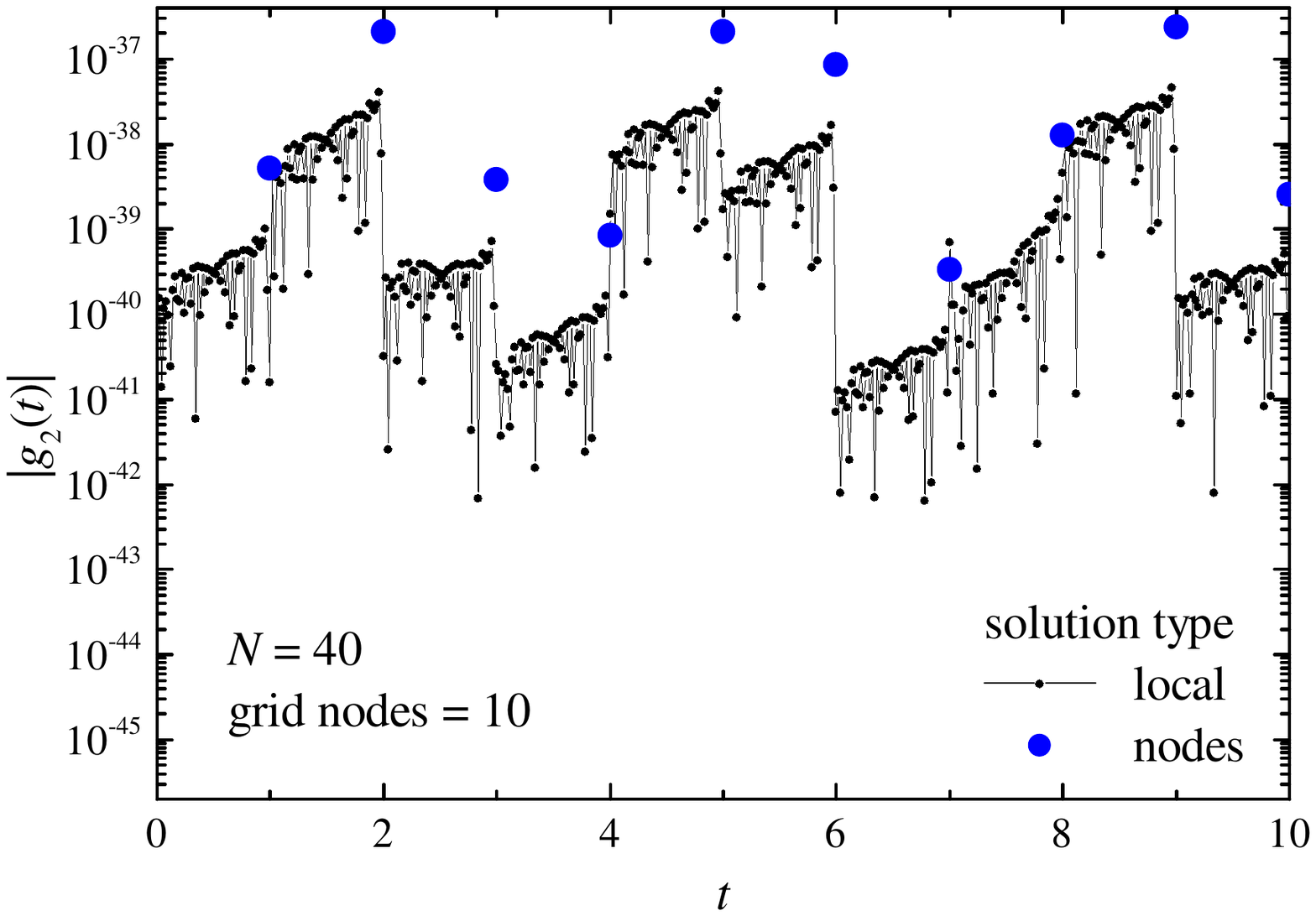}
\vspace{-8mm}\caption{\label{fig:pend_ind3_sol_g_eps:b3}}
\end{subfigure}\\[-2mm]
\begin{subfigure}{0.275\textwidth}
\includegraphics[width=\textwidth]{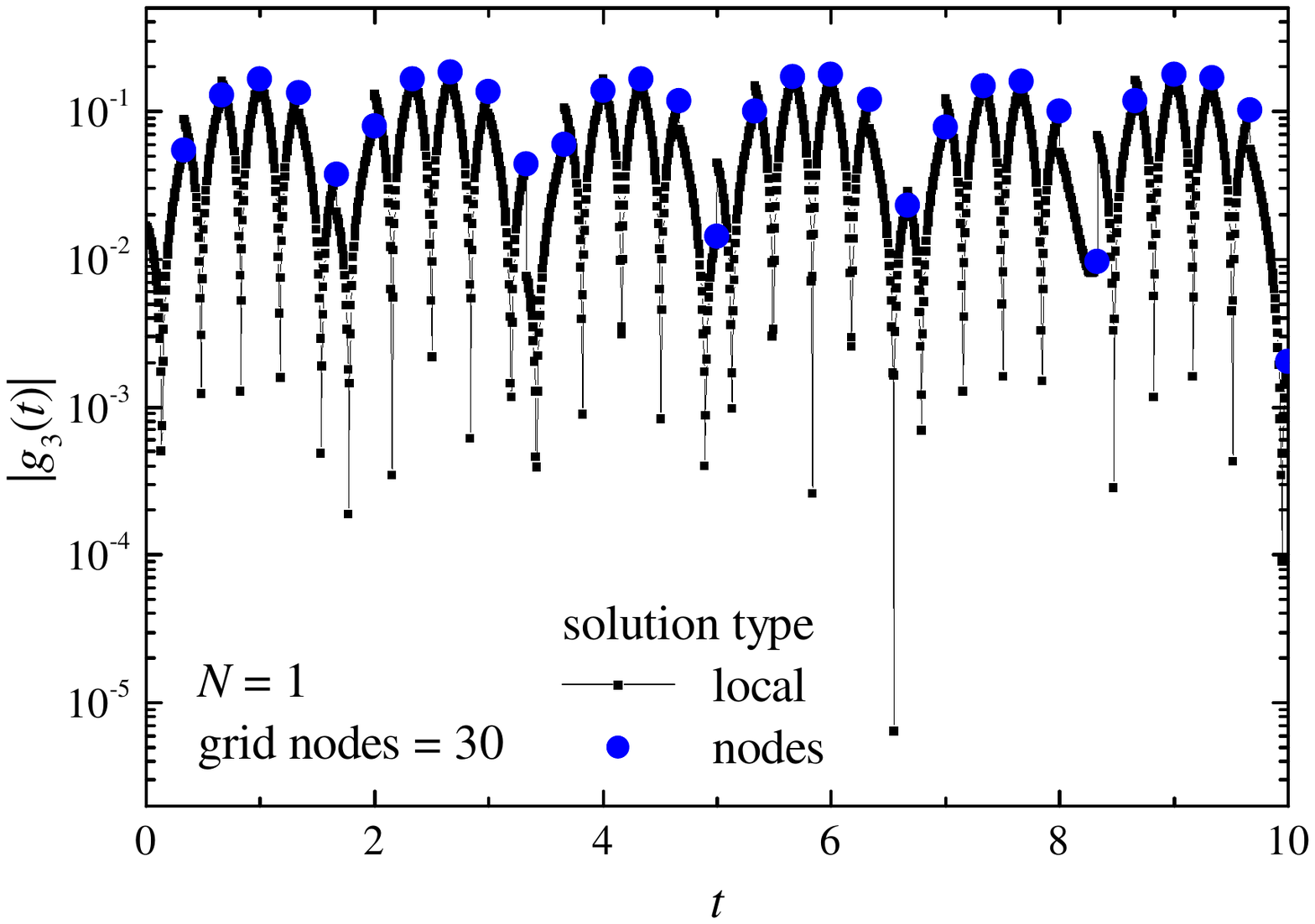}
\vspace{-8mm}\caption{\label{fig:pend_ind3_sol_g_eps:c1}}
\end{subfigure}\hspace{6mm}
\begin{subfigure}{0.275\textwidth}
\includegraphics[width=\textwidth]{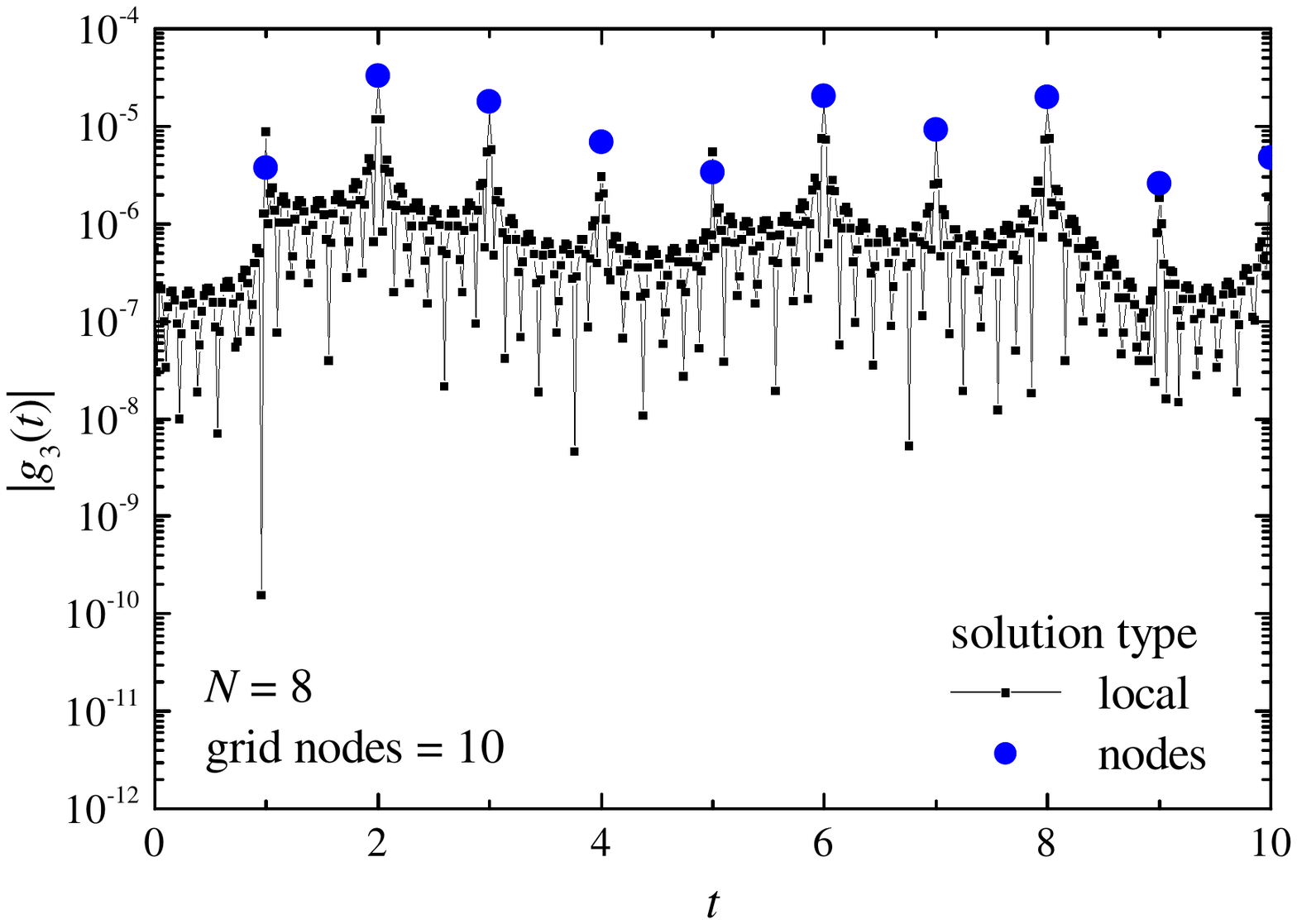}
\vspace{-8mm}\caption{\label{fig:pend_ind3_sol_g_eps:c2}}
\end{subfigure}\hspace{6mm}
\begin{subfigure}{0.275\textwidth}
\includegraphics[width=\textwidth]{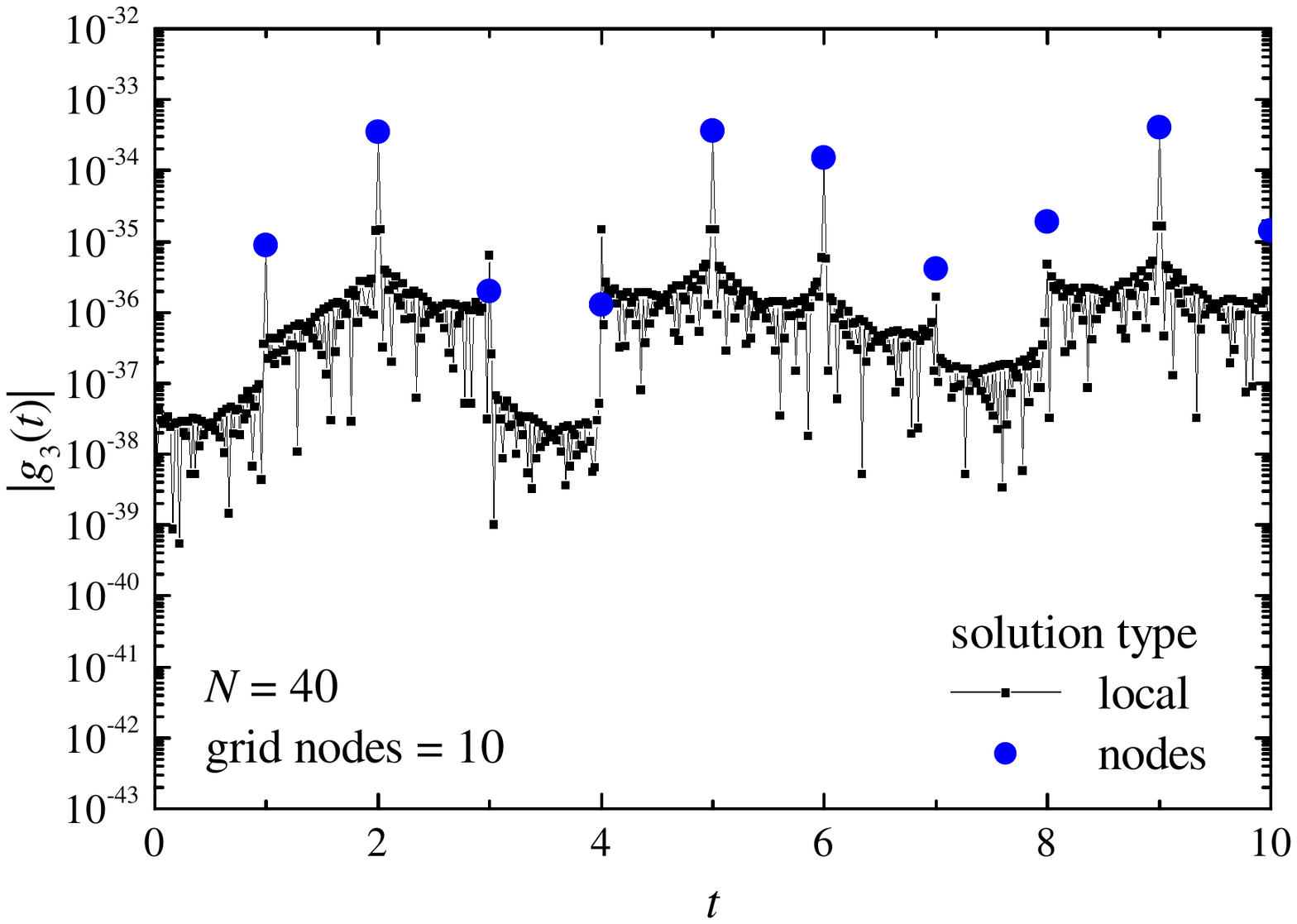}
\vspace{-8mm}\caption{\label{fig:pend_ind3_sol_g_eps:c3}}
\end{subfigure}\\[-2mm]
\begin{subfigure}{0.275\textwidth}
\includegraphics[width=\textwidth]{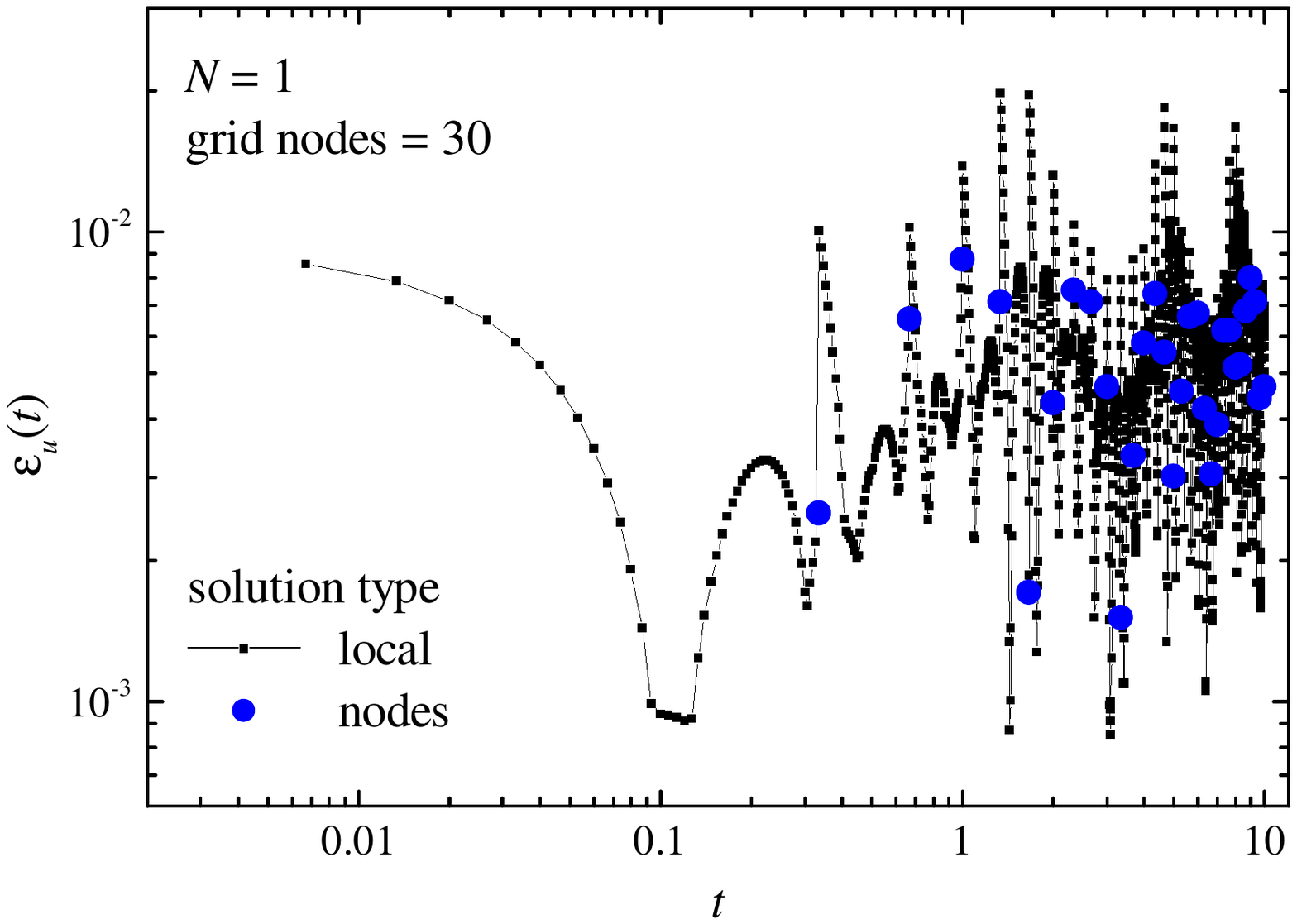}
\vspace{-8mm}\caption{\label{fig:pend_ind3_sol_g_eps:d1}}
\end{subfigure}\hspace{6mm}
\begin{subfigure}{0.275\textwidth}
\includegraphics[width=\textwidth]{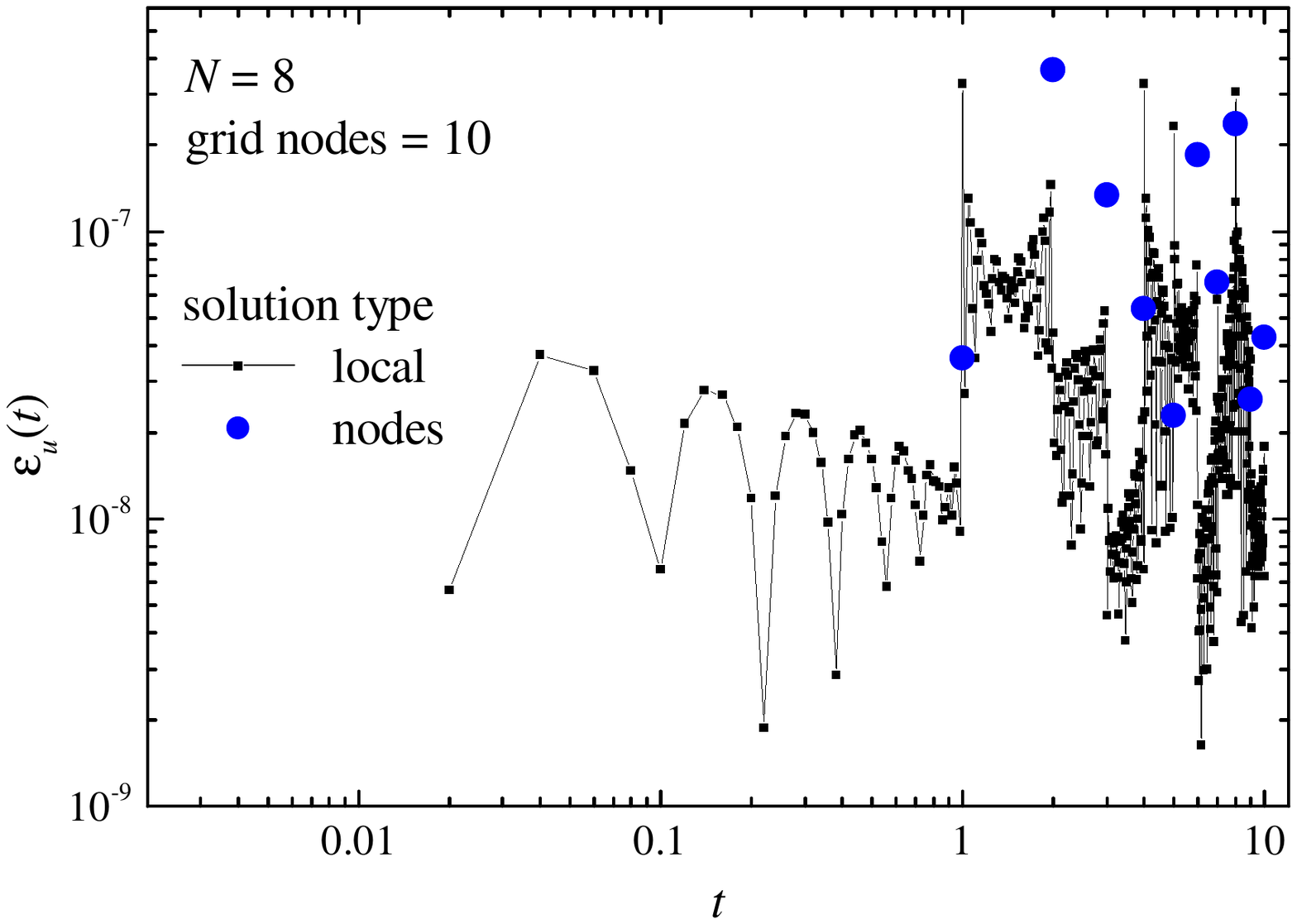}
\vspace{-8mm}\caption{\label{fig:pend_ind3_sol_g_eps:d2}}
\end{subfigure}\hspace{6mm}
\begin{subfigure}{0.275\textwidth}
\includegraphics[width=\textwidth]{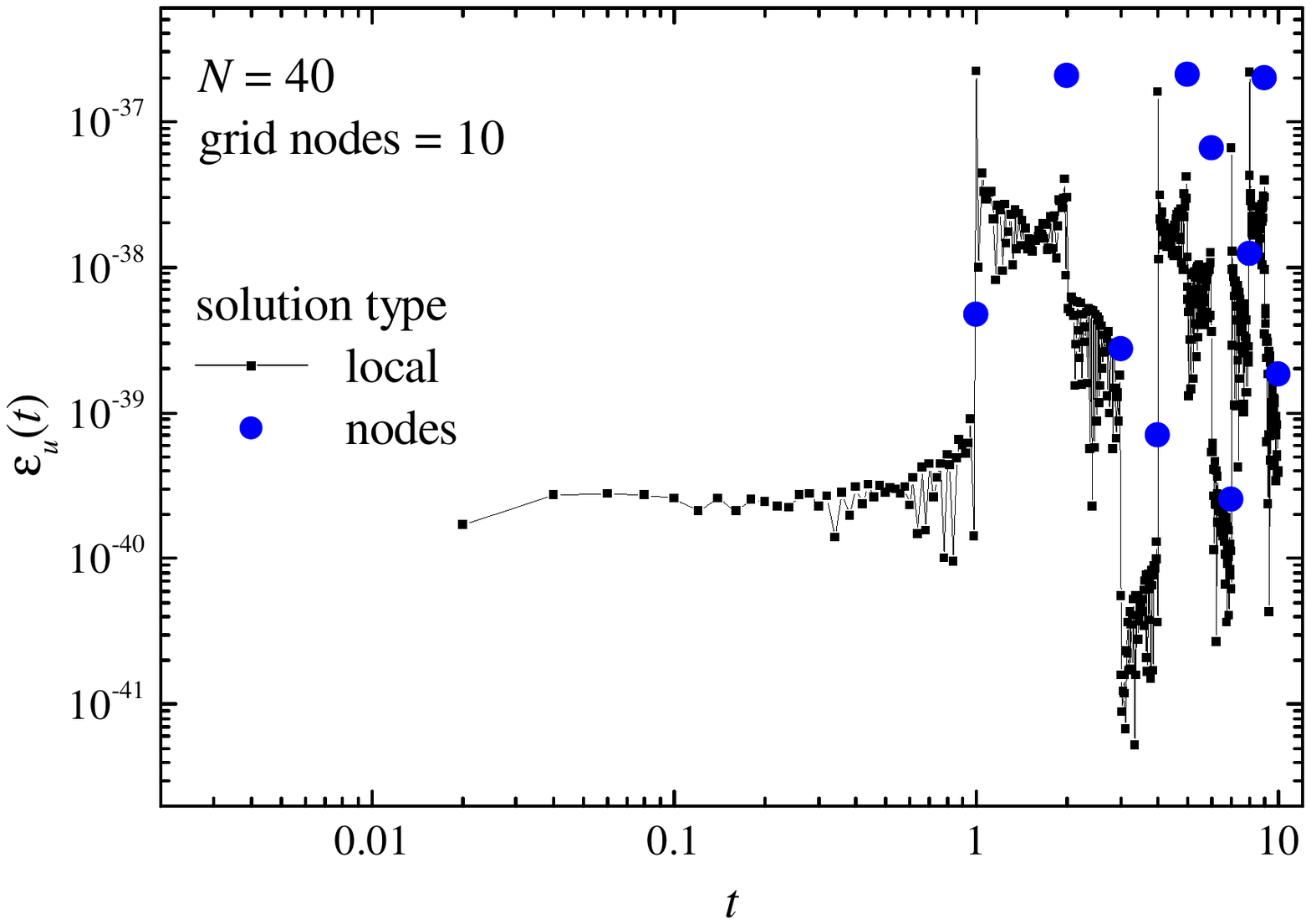}
\vspace{-8mm}\caption{\label{fig:pend_ind3_sol_g_eps:d3}}
\end{subfigure}\\[-2mm]
\begin{subfigure}{0.275\textwidth}
\includegraphics[width=\textwidth]{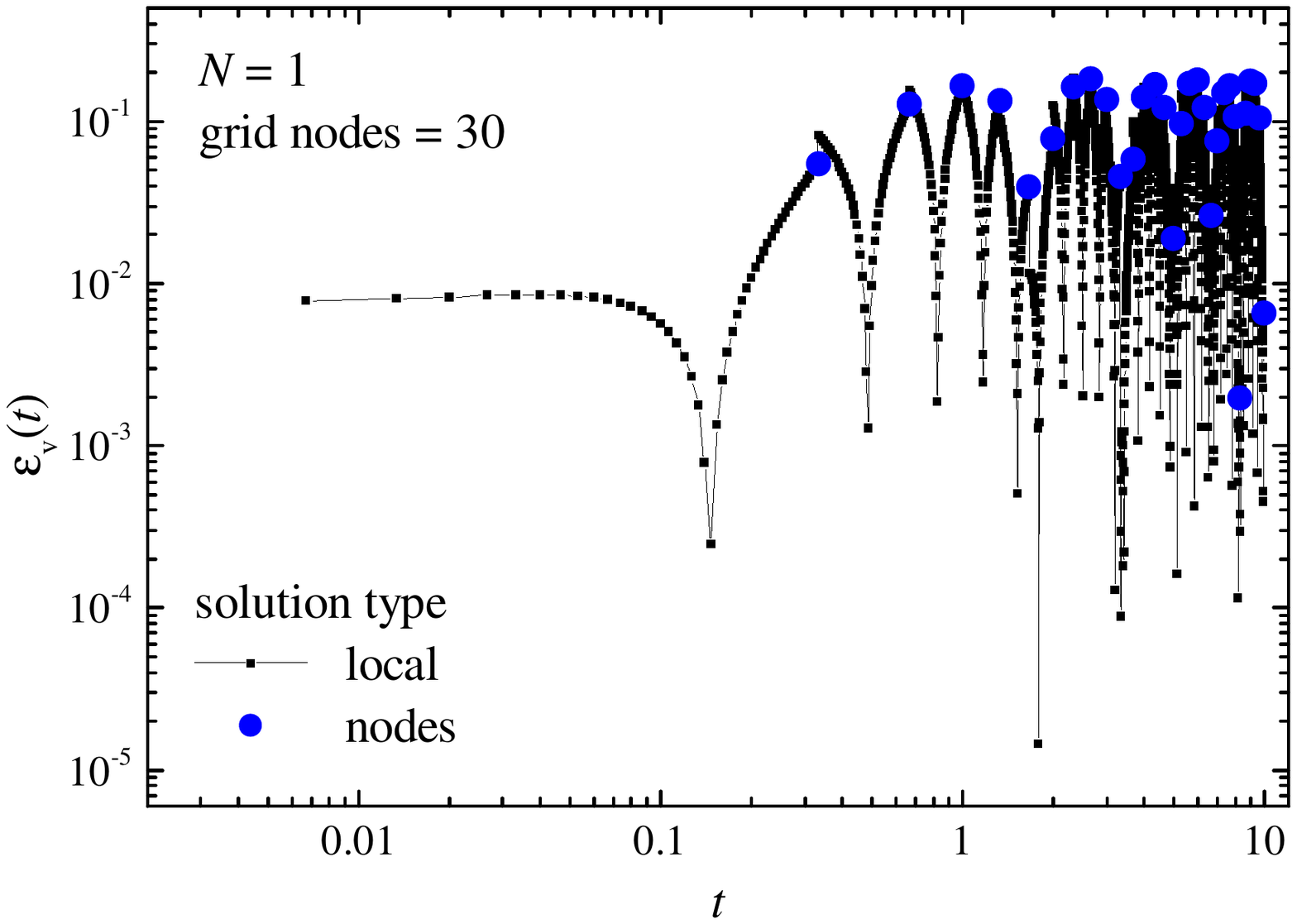}
\vspace{-8mm}\caption{\label{fig:pend_ind3_sol_g_eps:e1}}
\end{subfigure}\hspace{6mm}
\begin{subfigure}{0.275\textwidth}
\includegraphics[width=\textwidth]{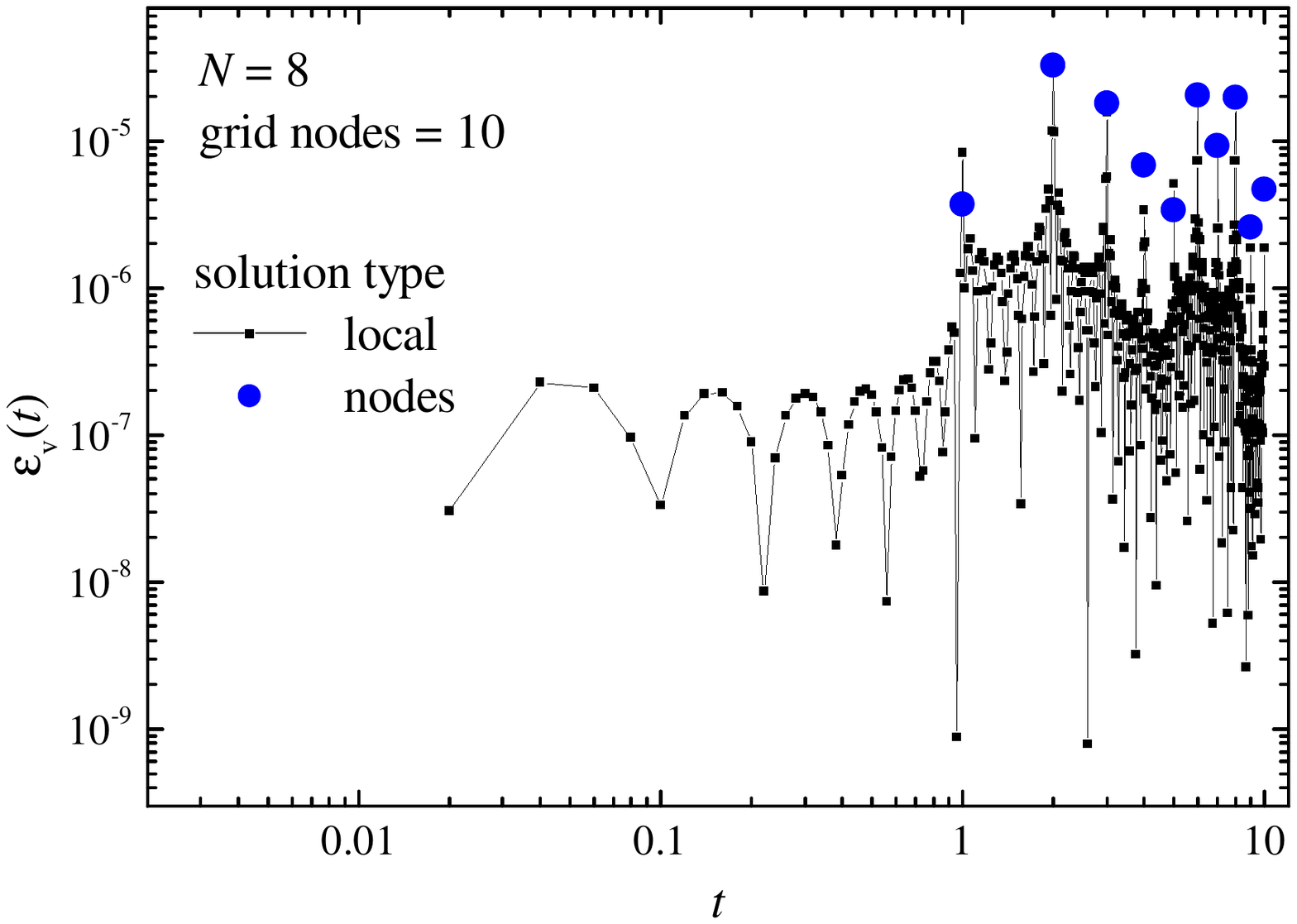}
\vspace{-8mm}\caption{\label{fig:pend_ind3_sol_g_eps:e2}}
\end{subfigure}\hspace{6mm}
\begin{subfigure}{0.275\textwidth}
\includegraphics[width=\textwidth]{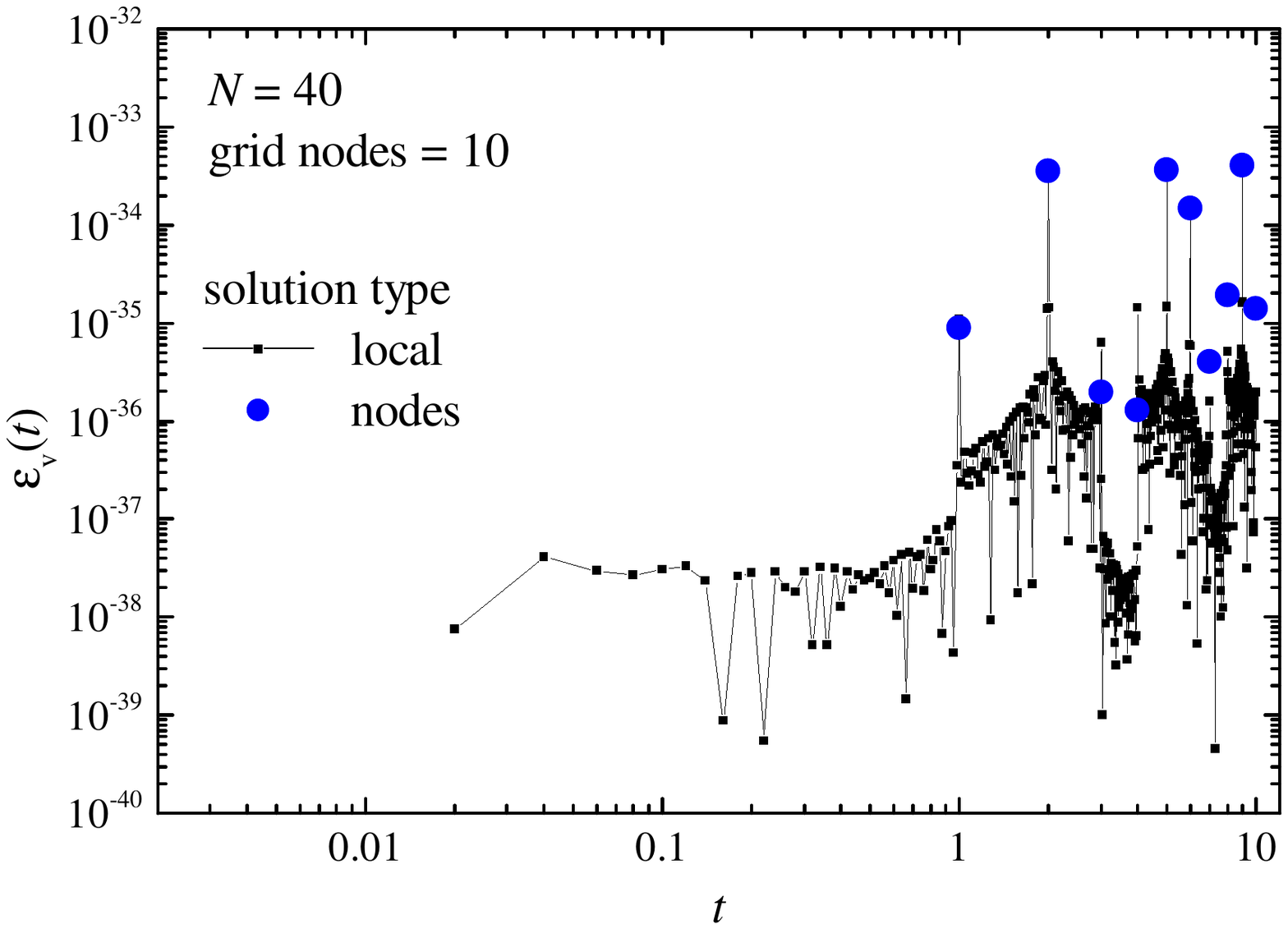}
\vspace{-8mm}\caption{\label{fig:pend_ind3_sol_g_eps:e3}}
\end{subfigure}\\[-2mm]
\begin{subfigure}{0.275\textwidth}
\includegraphics[width=\textwidth]{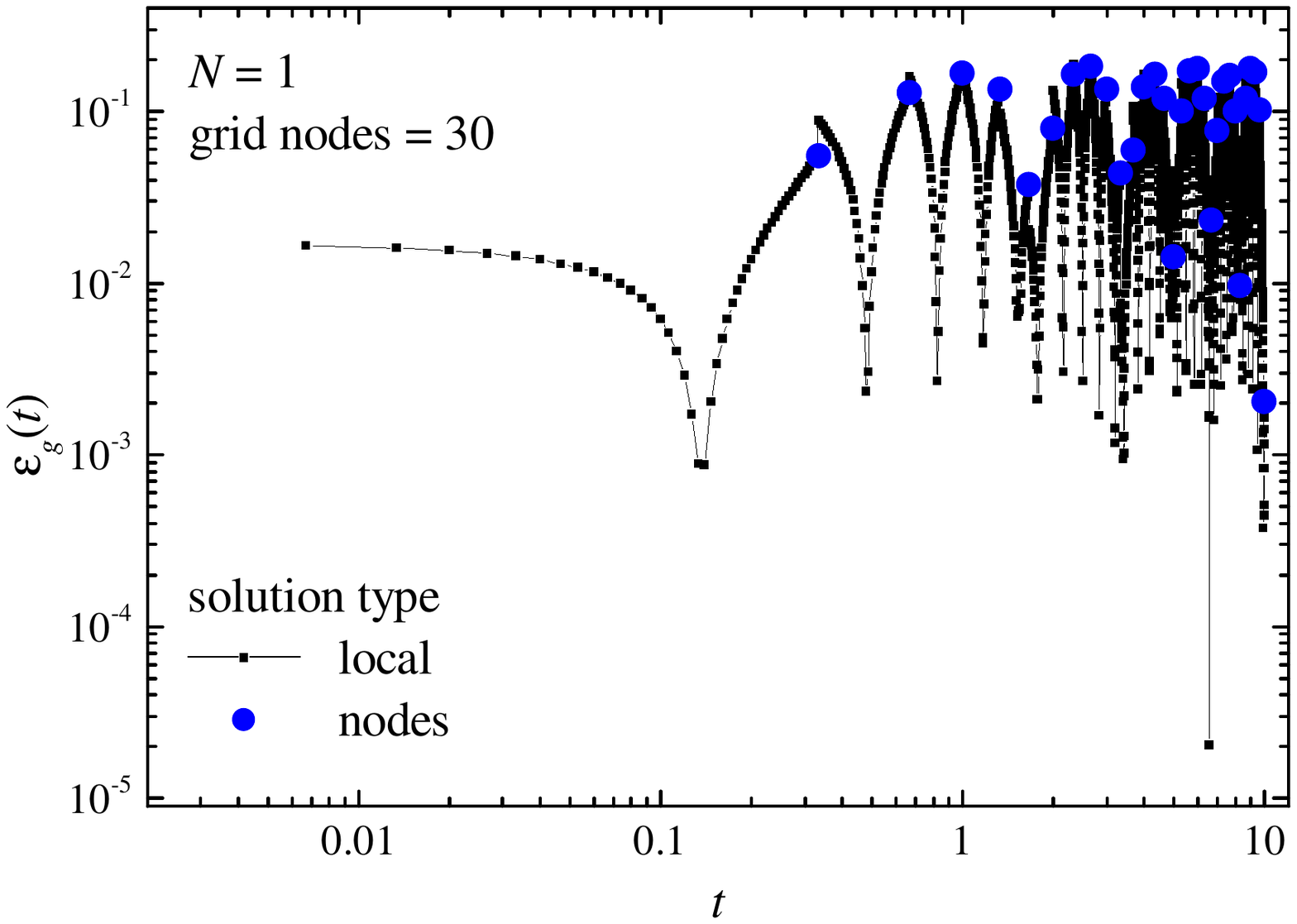}
\vspace{-8mm}\caption{\label{fig:pend_ind3_sol_g_eps:f1}}
\end{subfigure}\hspace{6mm}
\begin{subfigure}{0.275\textwidth}
\includegraphics[width=\textwidth]{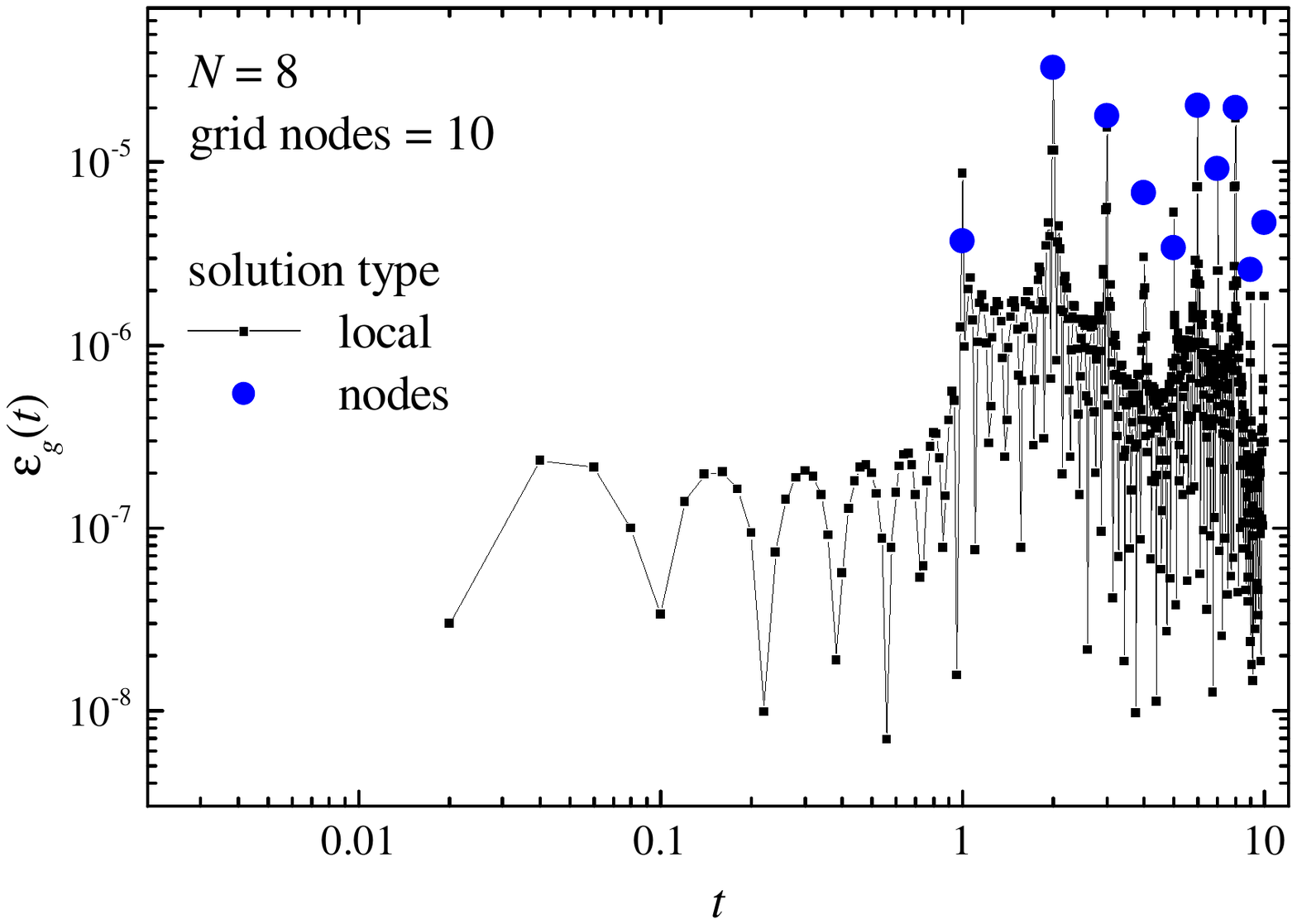}
\vspace{-8mm}\caption{\label{fig:pend_ind3_sol_g_eps:f2}}
\end{subfigure}\hspace{6mm}
\begin{subfigure}{0.275\textwidth}
\includegraphics[width=\textwidth]{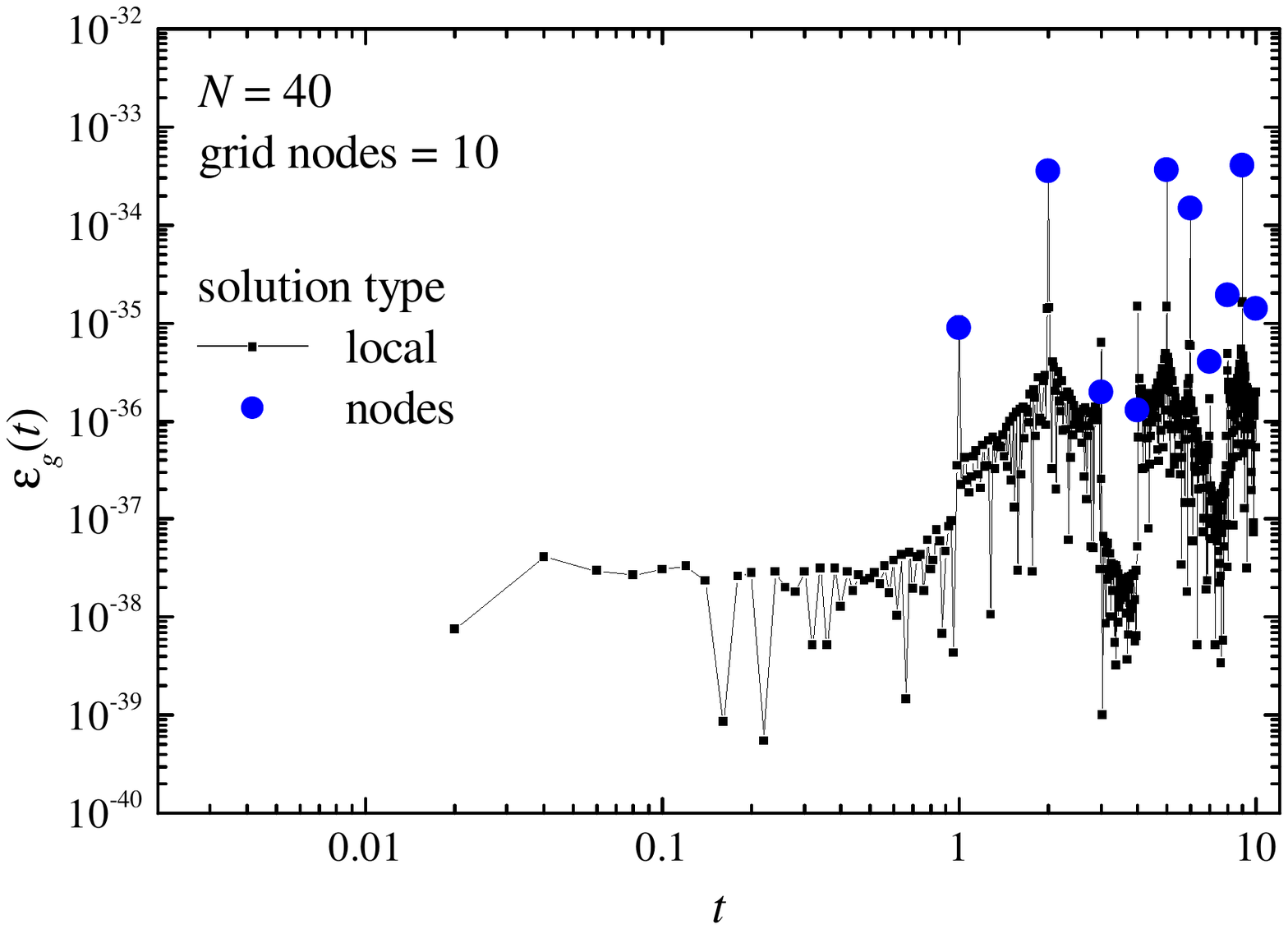}
\vspace{-8mm}\caption{\label{fig:pend_ind3_sol_g_eps:f3}}
\end{subfigure}\\[-2mm]
\caption{%
Numerical solution of the DAE system (\ref{eq:math_pend_dae_ind_3}) of index 3. Comparison of quantitative satisfiability of the conditions $g_{1} = 0$ (\subref{fig:pend_ind3_sol_g_eps:a1}, \subref{fig:pend_ind3_sol_g_eps:a2}, \subref{fig:pend_ind3_sol_g_eps:a3}), $g_{2} = 0$ (\subref{fig:pend_ind3_sol_g_eps:b1}, \subref{fig:pend_ind3_sol_g_eps:b2}, \subref{fig:pend_ind3_sol_g_eps:b3}) and $g_{3} = 0$ (\subref{fig:pend_ind3_sol_g_eps:c1}, \subref{fig:pend_ind3_sol_g_eps:c2}, \subref{fig:pend_ind3_sol_g_eps:c3}), the errors $\varepsilon_{u}(t)$ (\subref{fig:pend_ind3_sol_g_eps:d1}, \subref{fig:pend_ind3_sol_g_eps:d2}, \subref{fig:pend_ind3_sol_g_eps:d3}), $\varepsilon_{v}(t)$ (\subref{fig:pend_ind3_sol_g_eps:e1}, \subref{fig:pend_ind3_sol_g_eps:e2}, \subref{fig:pend_ind3_sol_g_eps:e3}), $\varepsilon_{g}(t)$ (\subref{fig:pend_ind3_sol_g_eps:f1}, \subref{fig:pend_ind3_sol_g_eps:f2}, \subref{fig:pend_ind3_sol_g_eps:f3}), for numerical solution obtained using polynomials with degrees $N = 1$ (\subref{fig:pend_ind3_sol_g_eps:a1}, \subref{fig:pend_ind3_sol_g_eps:b1}, \subref{fig:pend_ind3_sol_g_eps:c1}, \subref{fig:pend_ind3_sol_g_eps:d1}, \subref{fig:pend_ind3_sol_g_eps:e1}, \subref{fig:pend_ind3_sol_g_eps:f1}), $N = 8$ (\subref{fig:pend_ind3_sol_g_eps:a2}, \subref{fig:pend_ind3_sol_g_eps:b2}, \subref{fig:pend_ind3_sol_g_eps:c2}, \subref{fig:pend_ind3_sol_g_eps:d2}, \subref{fig:pend_ind3_sol_g_eps:e2}, \subref{fig:pend_ind3_sol_g_eps:f2}) and $N = 40$ (\subref{fig:pend_ind3_sol_g_eps:a3}, \subref{fig:pend_ind3_sol_g_eps:b3}, \subref{fig:pend_ind3_sol_g_eps:c3}, \subref{fig:pend_ind3_sol_g_eps:d3}, \subref{fig:pend_ind3_sol_g_eps:e3}, \subref{fig:pend_ind3_sol_g_eps:f3}).
}
\label{fig:pend_ind3_sol_g_eps}
\end{figure} 

\begin{figure}[h!]
\captionsetup[subfigure]{%
	position=bottom,
	font+=smaller,
	textfont=normalfont,
	singlelinecheck=off,
	justification=raggedright
}
\centering
\begin{subfigure}{0.275\textwidth}
\includegraphics[width=\textwidth]{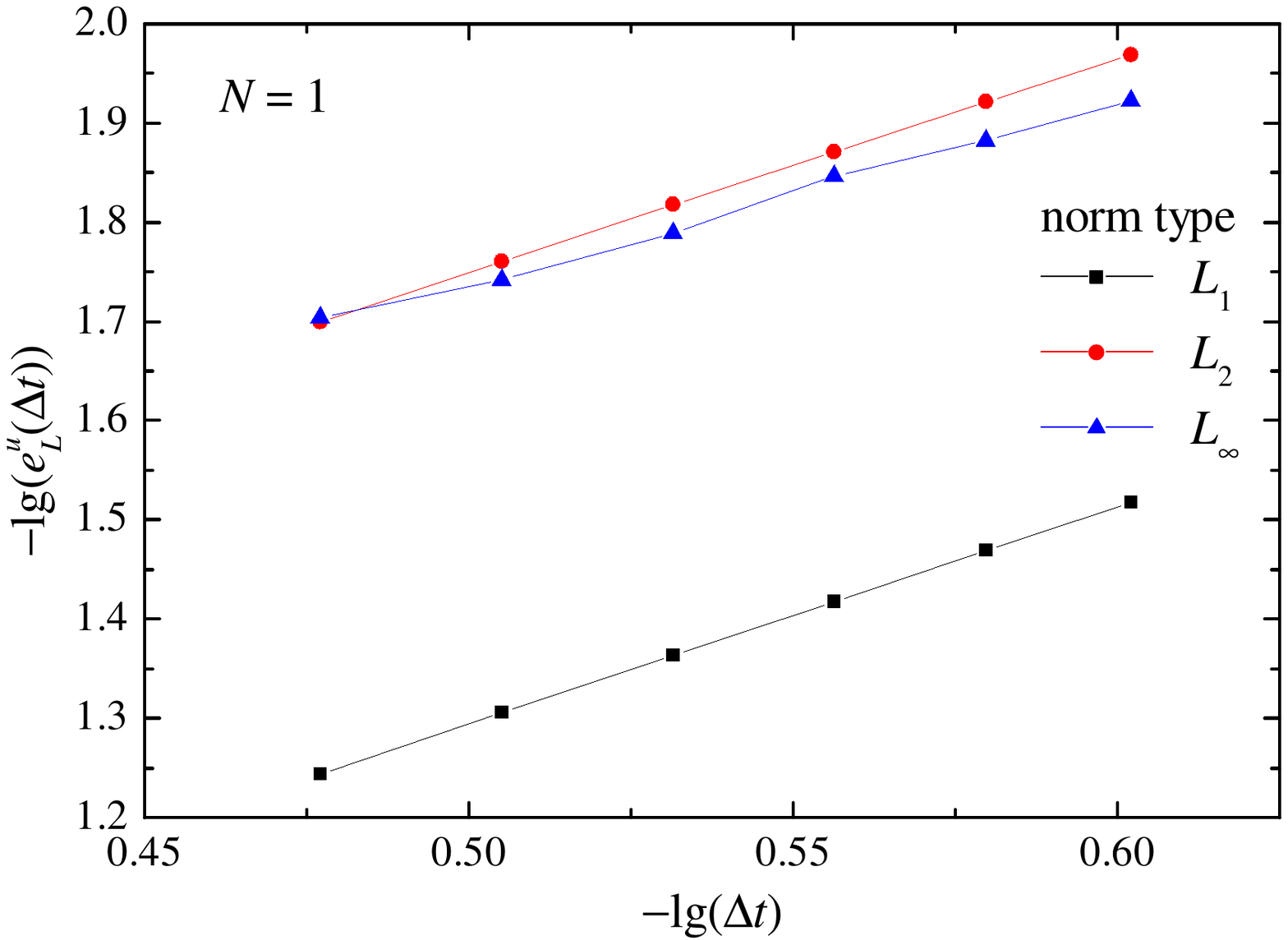}
\vspace{-8mm}\caption{\label{fig:pend_ind3_errors:a1}}
\end{subfigure}\hspace{6mm}
\begin{subfigure}{0.275\textwidth}
\includegraphics[width=\textwidth]{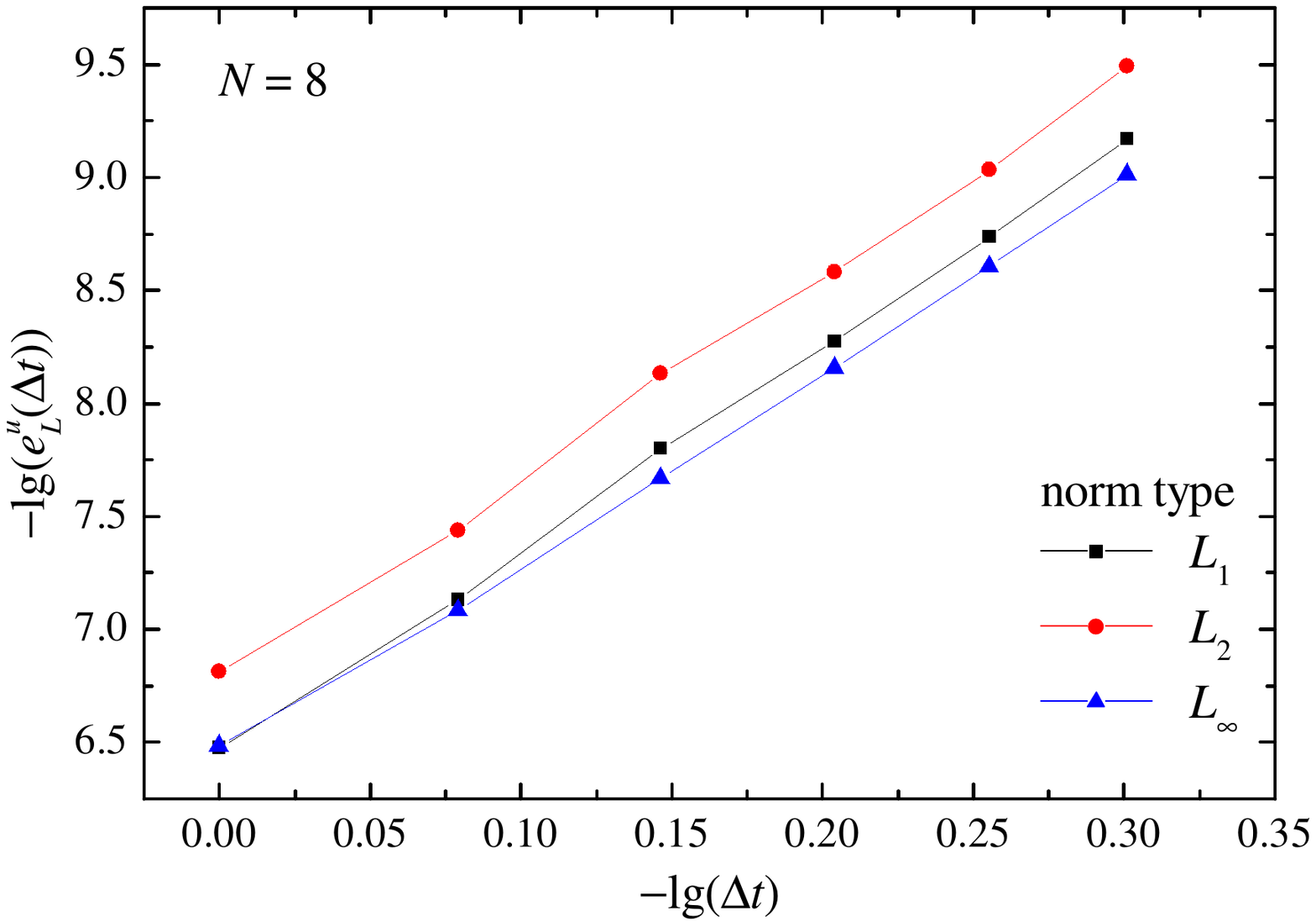}
\vspace{-8mm}\caption{\label{fig:pend_ind3_errors:a2}}
\end{subfigure}\hspace{6mm}
\begin{subfigure}{0.275\textwidth}
\includegraphics[width=\textwidth]{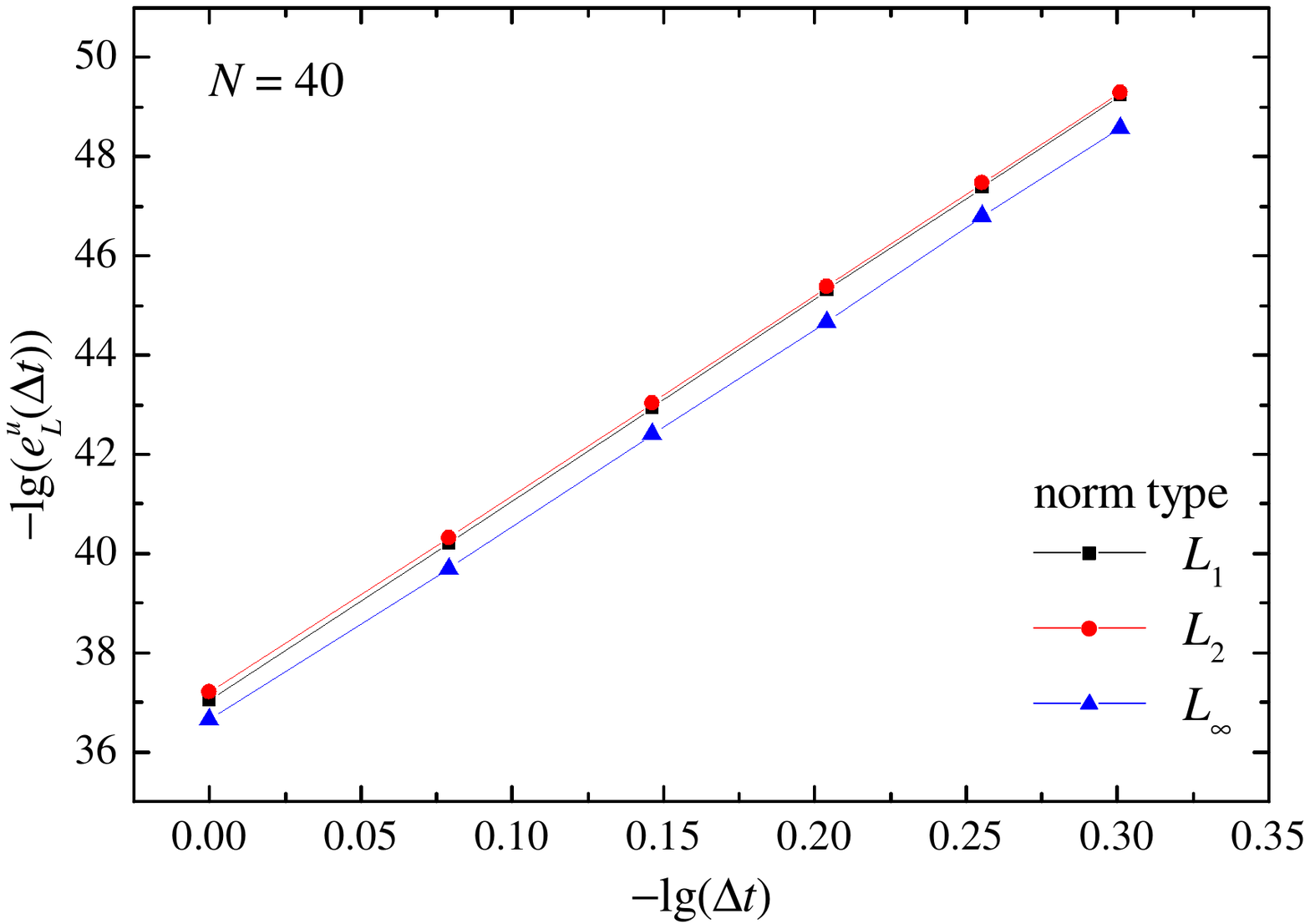}
\vspace{-8mm}\caption{\label{fig:pend_ind3_errors:a3}}
\end{subfigure}\\[-2mm]
\begin{subfigure}{0.275\textwidth}
\includegraphics[width=\textwidth]{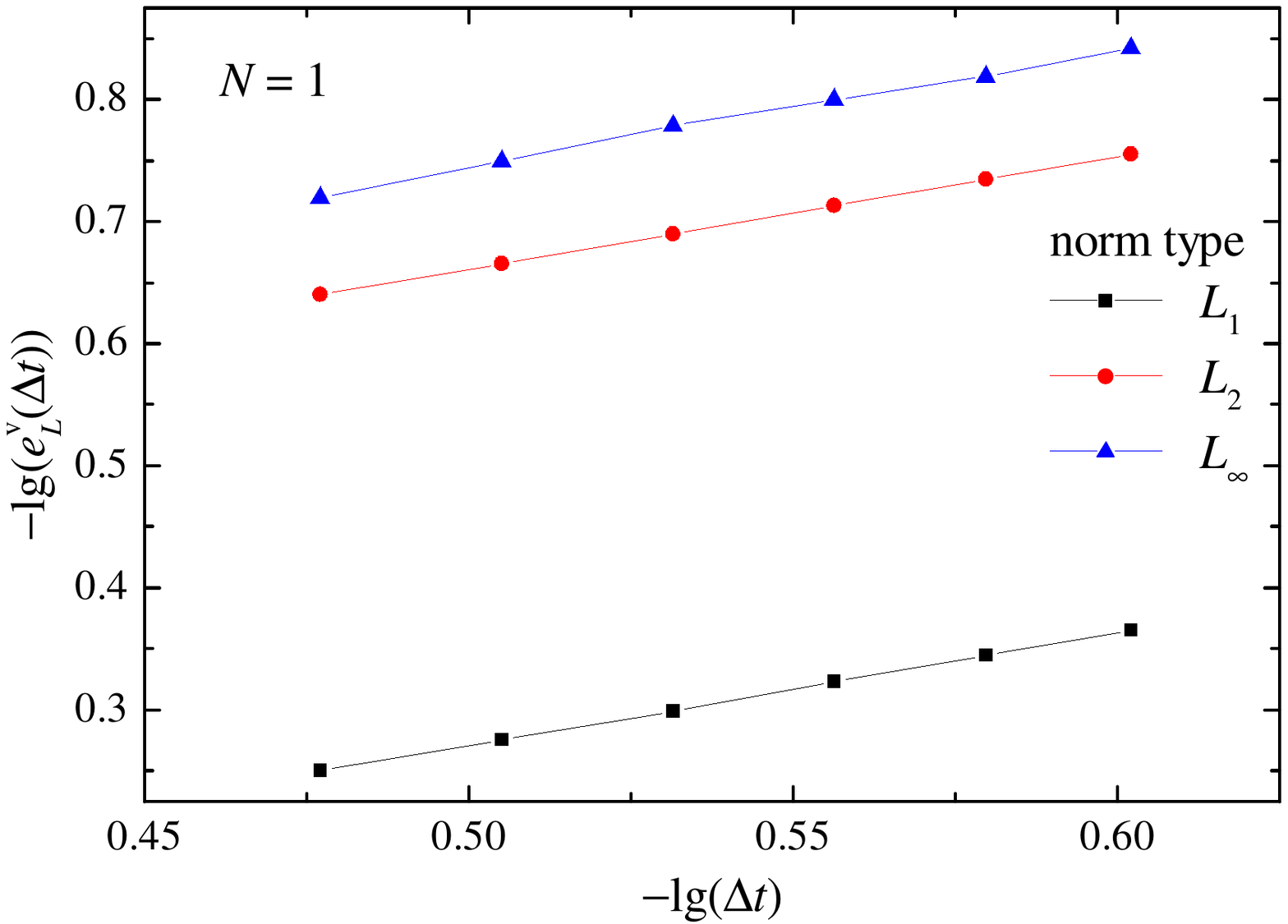}
\vspace{-8mm}\caption{\label{fig:pend_ind3_errors:b1}}
\end{subfigure}\hspace{6mm}
\begin{subfigure}{0.275\textwidth}
\includegraphics[width=\textwidth]{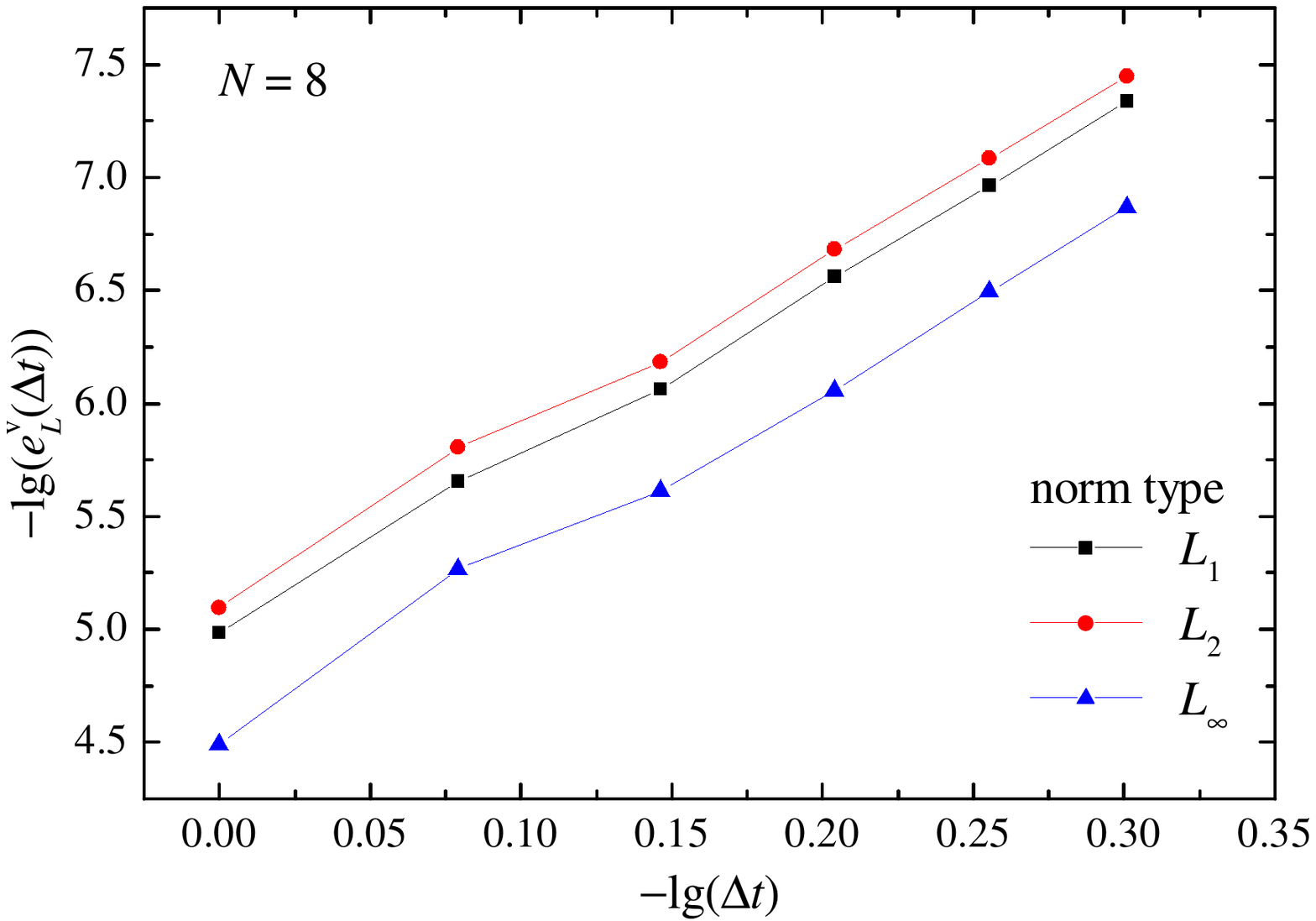}
\vspace{-8mm}\caption{\label{fig:pend_ind3_errors:b2}}
\end{subfigure}\hspace{6mm}
\begin{subfigure}{0.275\textwidth}
\includegraphics[width=\textwidth]{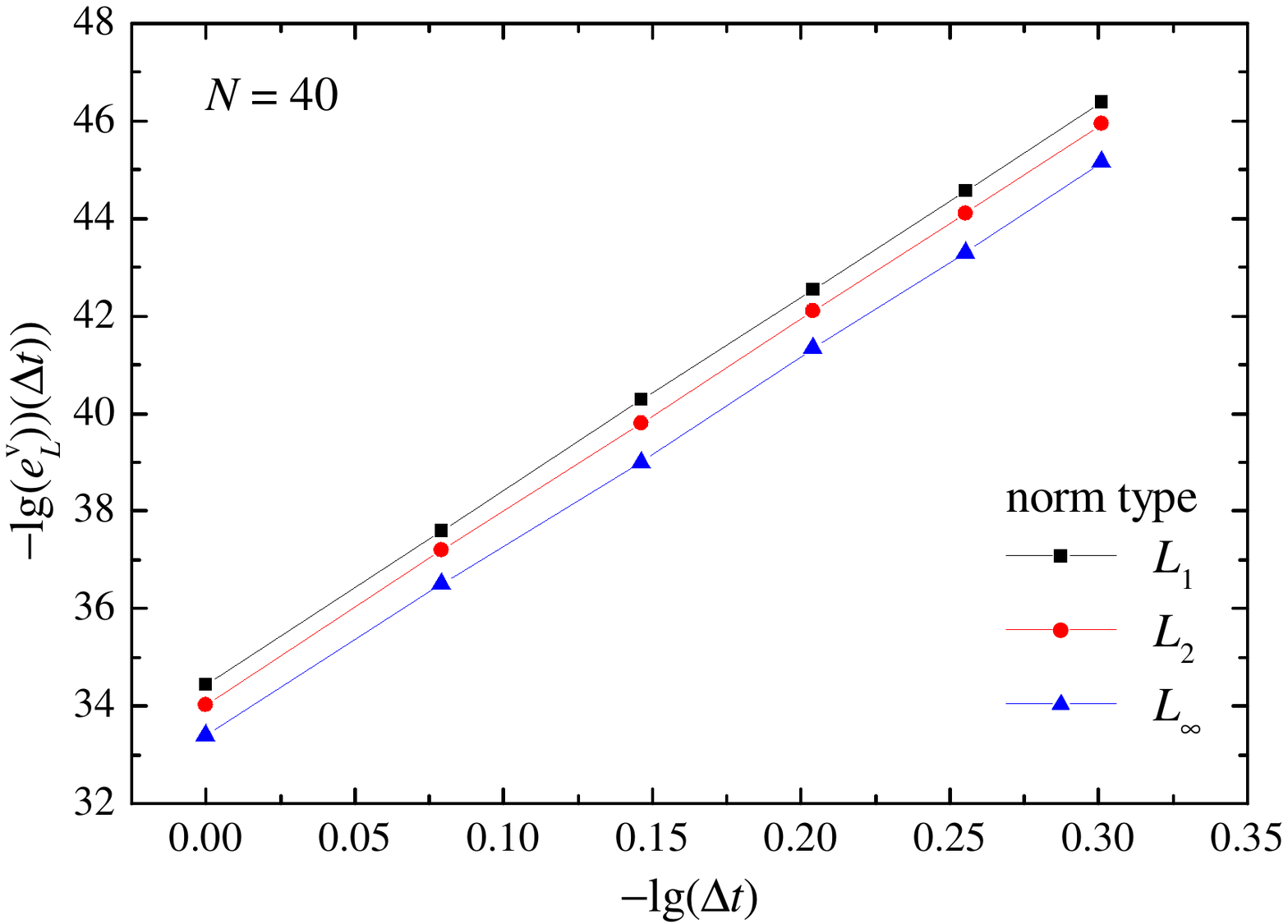}
\vspace{-8mm}\caption{\label{fig:pend_ind3_errors:b3}}
\end{subfigure}\\[-2mm]
\begin{subfigure}{0.275\textwidth}
\includegraphics[width=\textwidth]{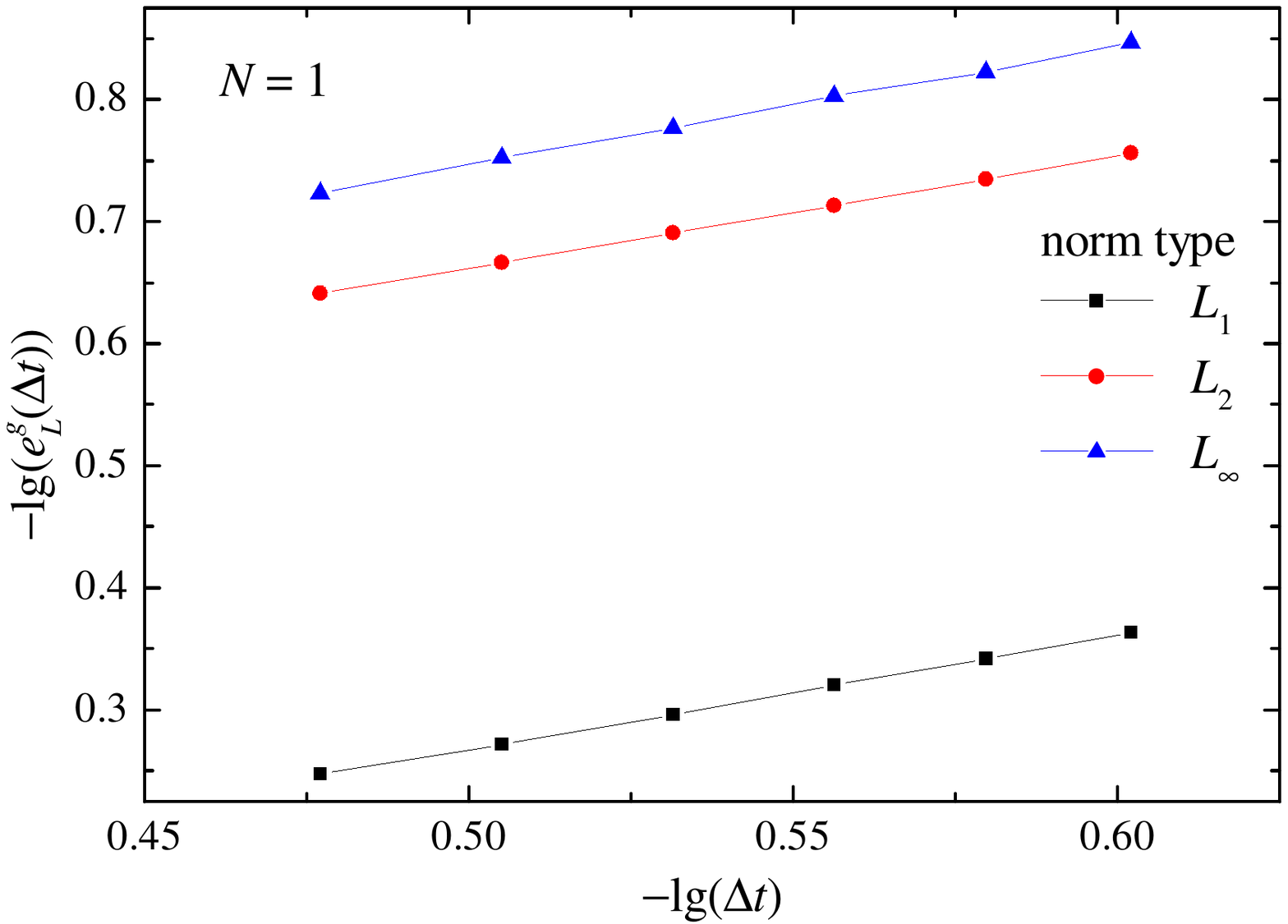}
\vspace{-8mm}\caption{\label{fig:pend_ind3_errors:c1}}
\end{subfigure}\hspace{6mm}
\begin{subfigure}{0.275\textwidth}
\includegraphics[width=\textwidth]{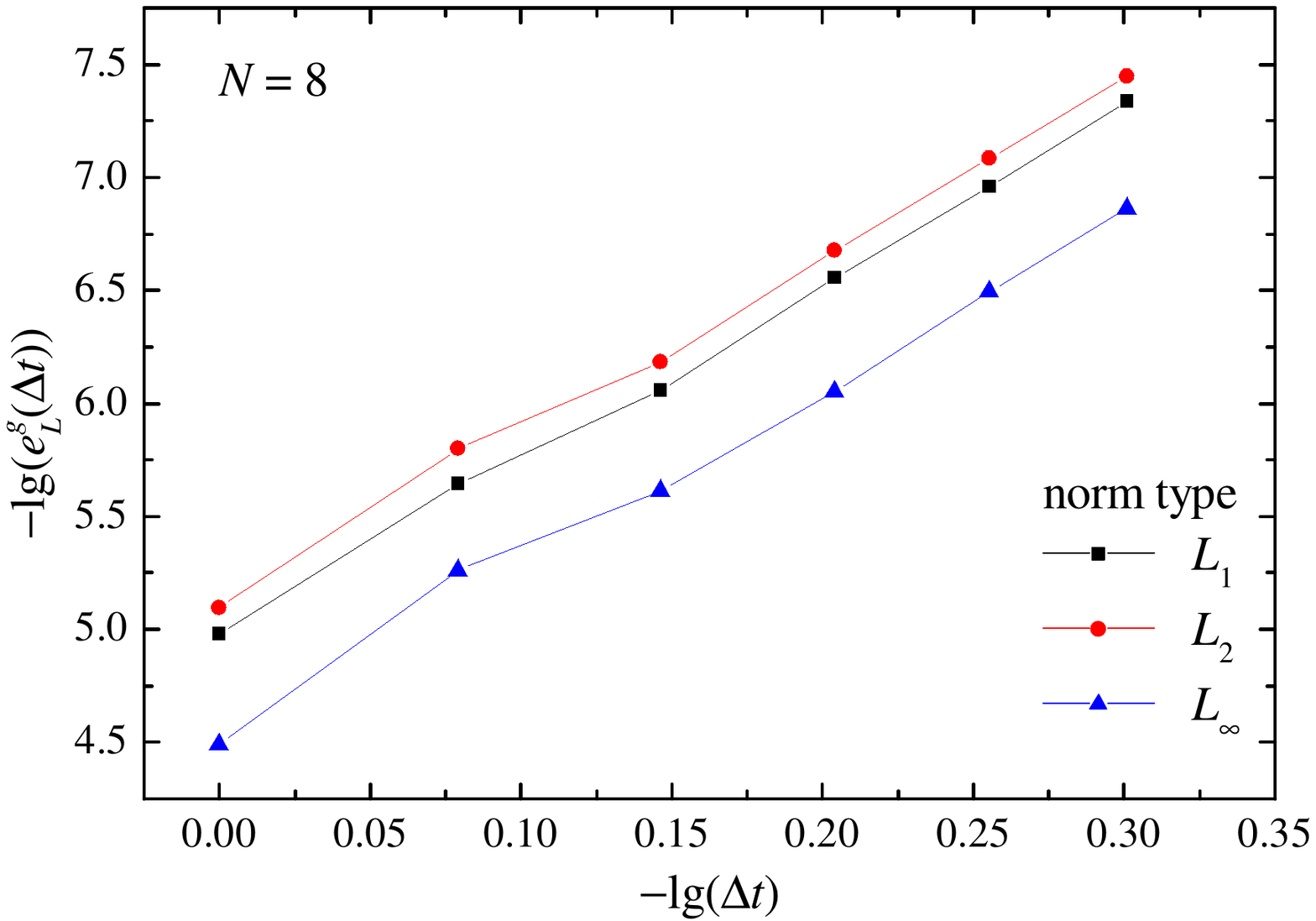}
\vspace{-8mm}\caption{\label{fig:pend_ind3_errors:c2}}
\end{subfigure}\hspace{6mm}
\begin{subfigure}{0.275\textwidth}
\includegraphics[width=\textwidth]{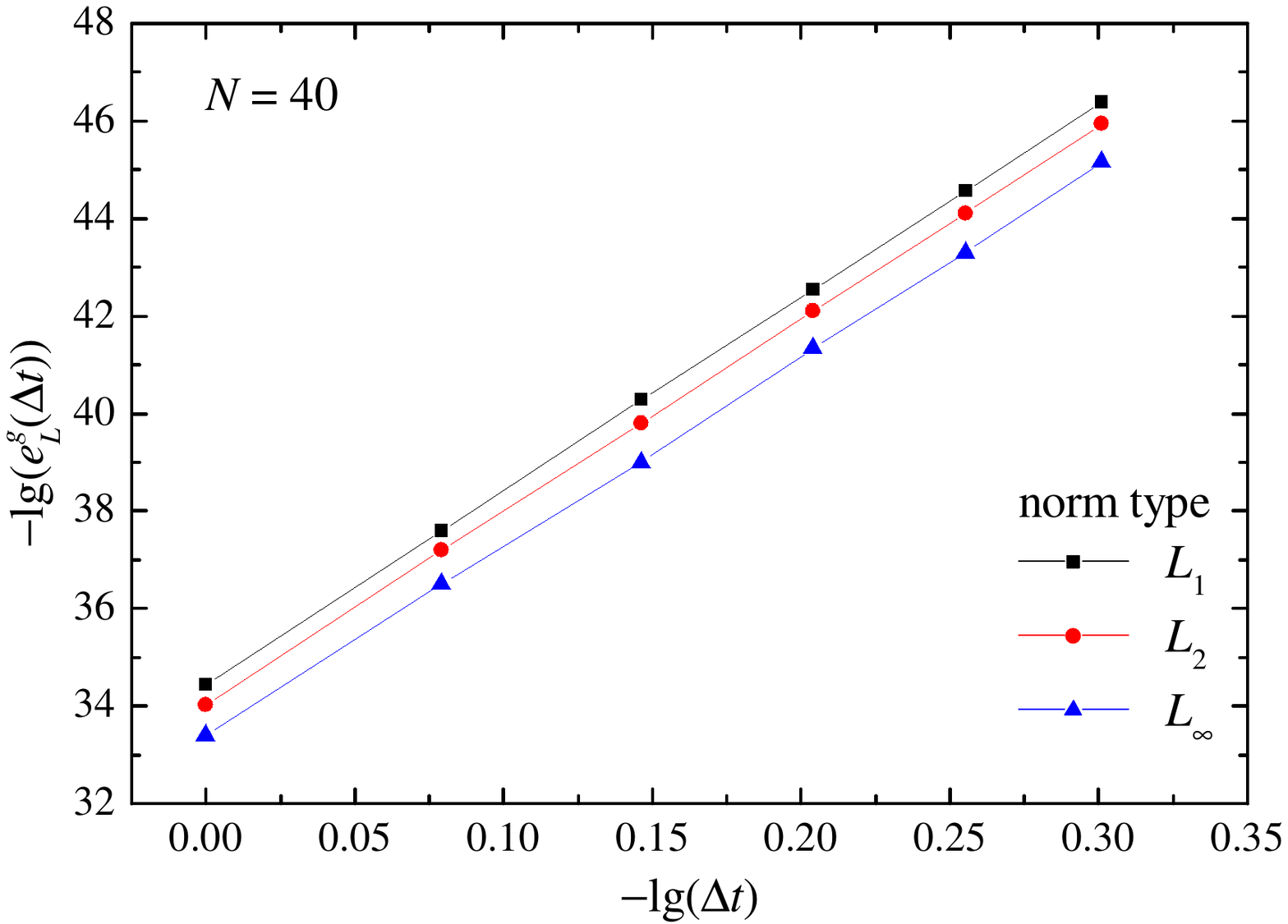}
\vspace{-8mm}\caption{\label{fig:pend_ind3_errors:c3}}
\end{subfigure}\\[-2mm]
\begin{subfigure}{0.275\textwidth}
\includegraphics[width=\textwidth]{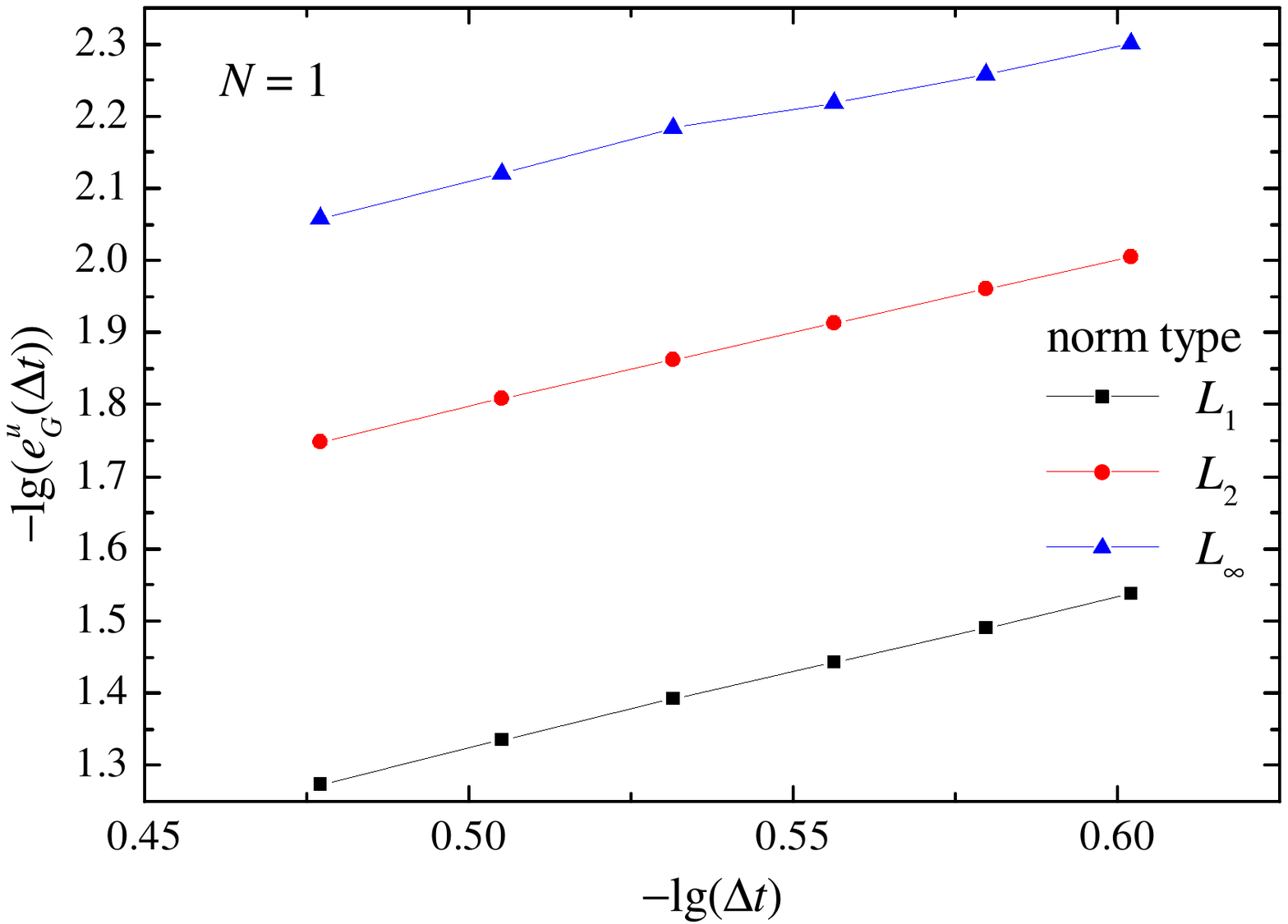}
\vspace{-8mm}\caption{\label{fig:pend_ind3_errors:d1}}
\end{subfigure}\hspace{6mm}
\begin{subfigure}{0.275\textwidth}
\includegraphics[width=\textwidth]{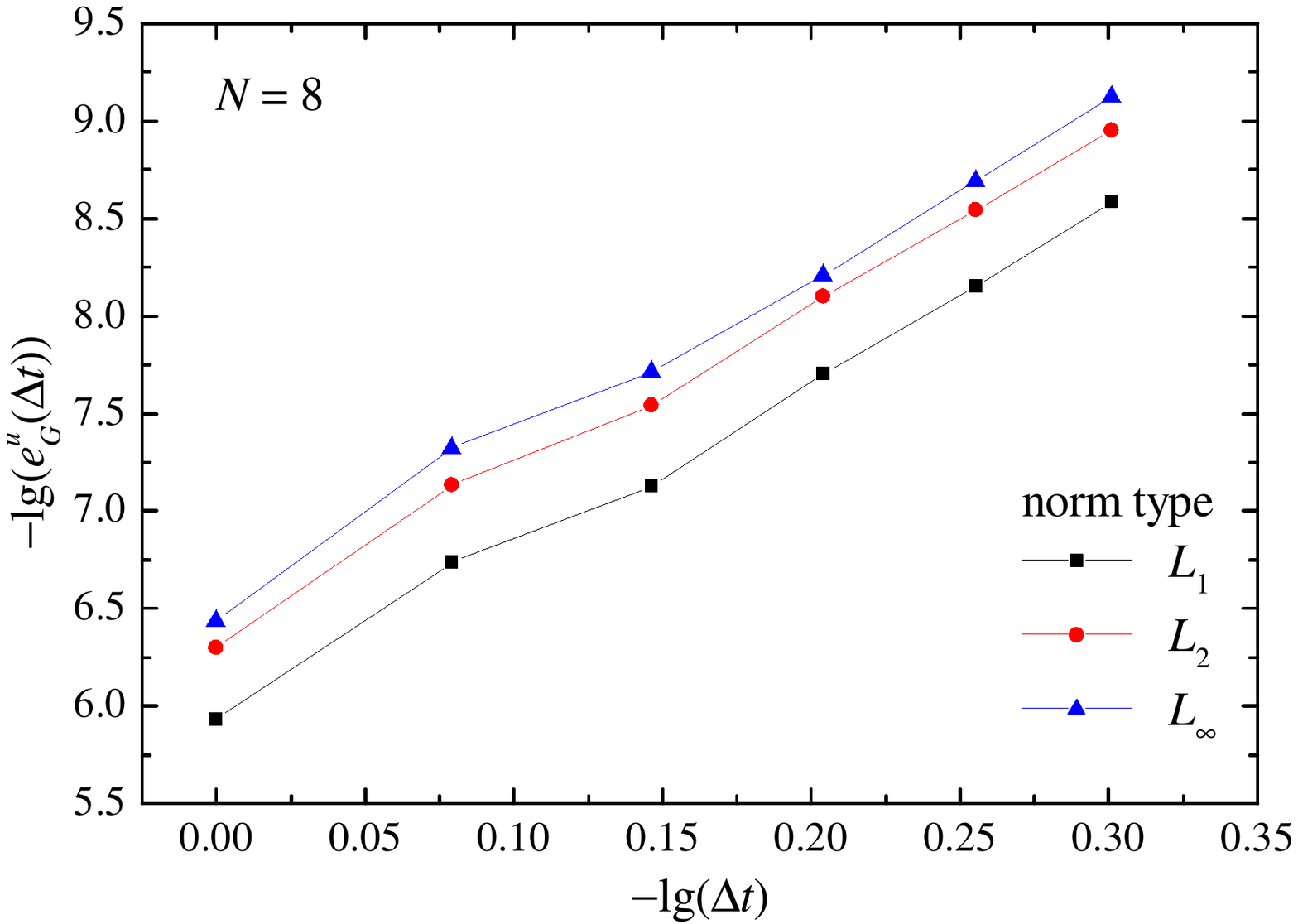}
\vspace{-8mm}\caption{\label{fig:pend_ind3_errors:d2}}
\end{subfigure}\hspace{6mm}
\begin{subfigure}{0.275\textwidth}
\includegraphics[width=\textwidth]{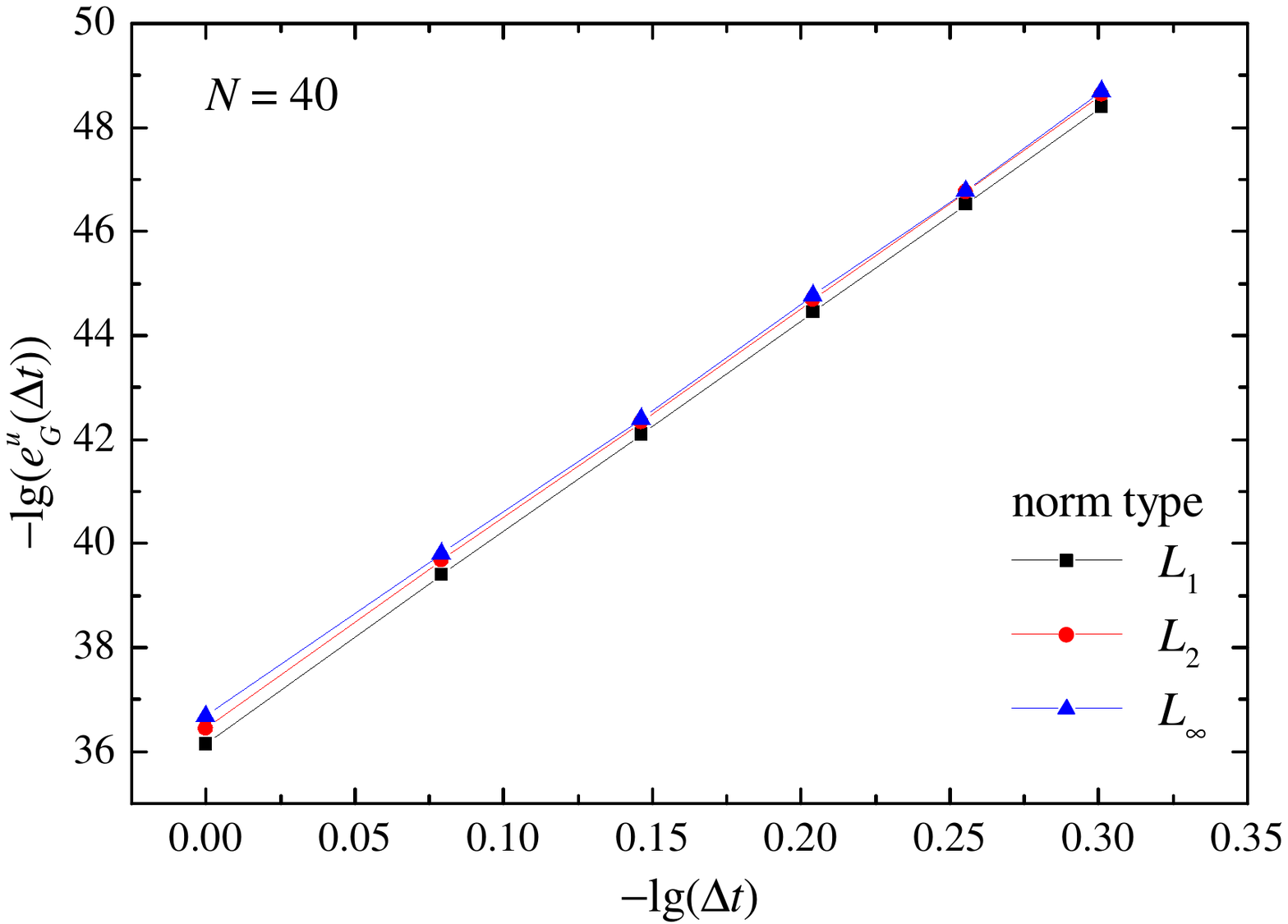}
\vspace{-8mm}\caption{\label{fig:pend_ind3_errors:d3}}
\end{subfigure}\\[-2mm]
\begin{subfigure}{0.275\textwidth}
\includegraphics[width=\textwidth]{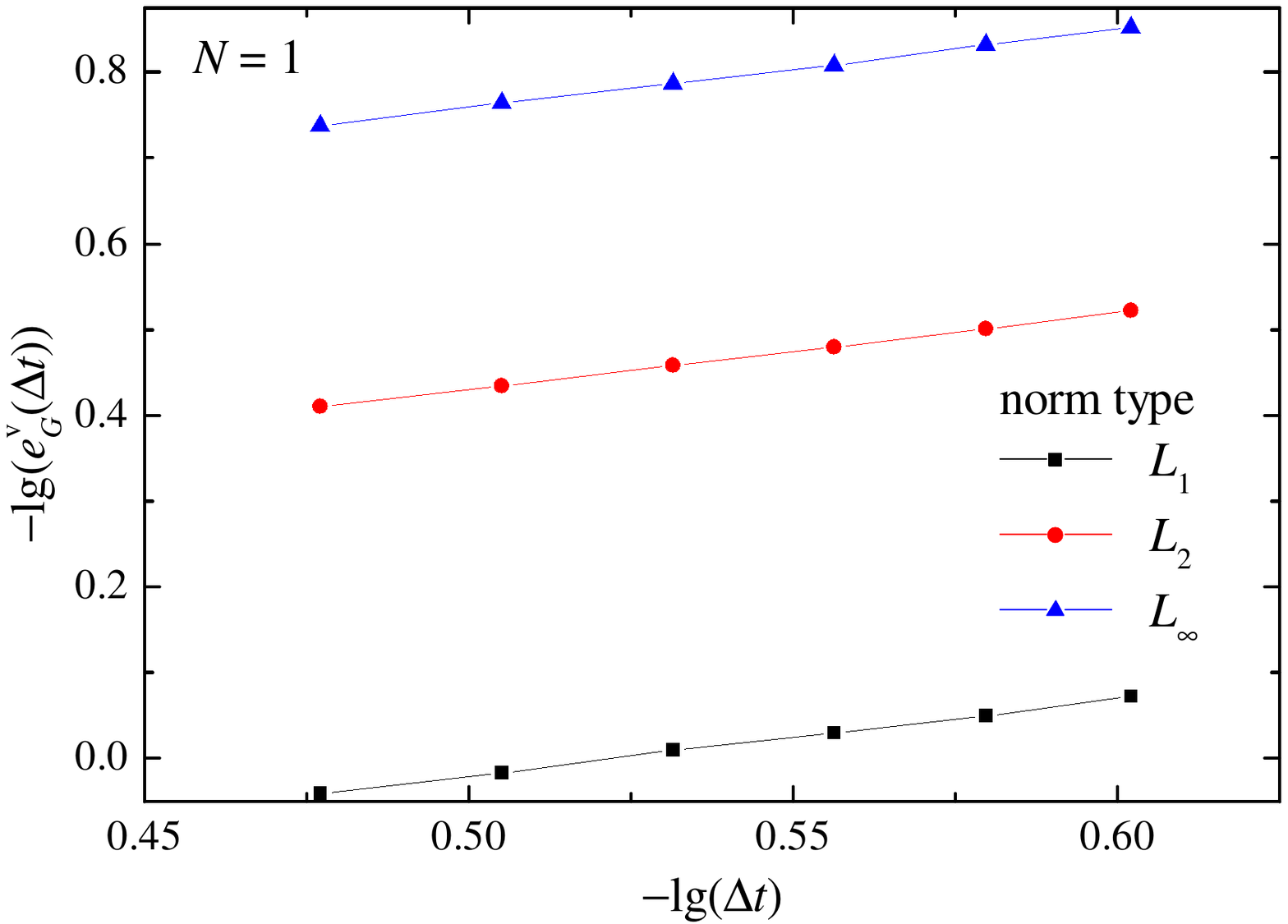}
\vspace{-8mm}\caption{\label{fig:pend_ind3_errors:e1}}
\end{subfigure}\hspace{6mm}
\begin{subfigure}{0.275\textwidth}
\includegraphics[width=\textwidth]{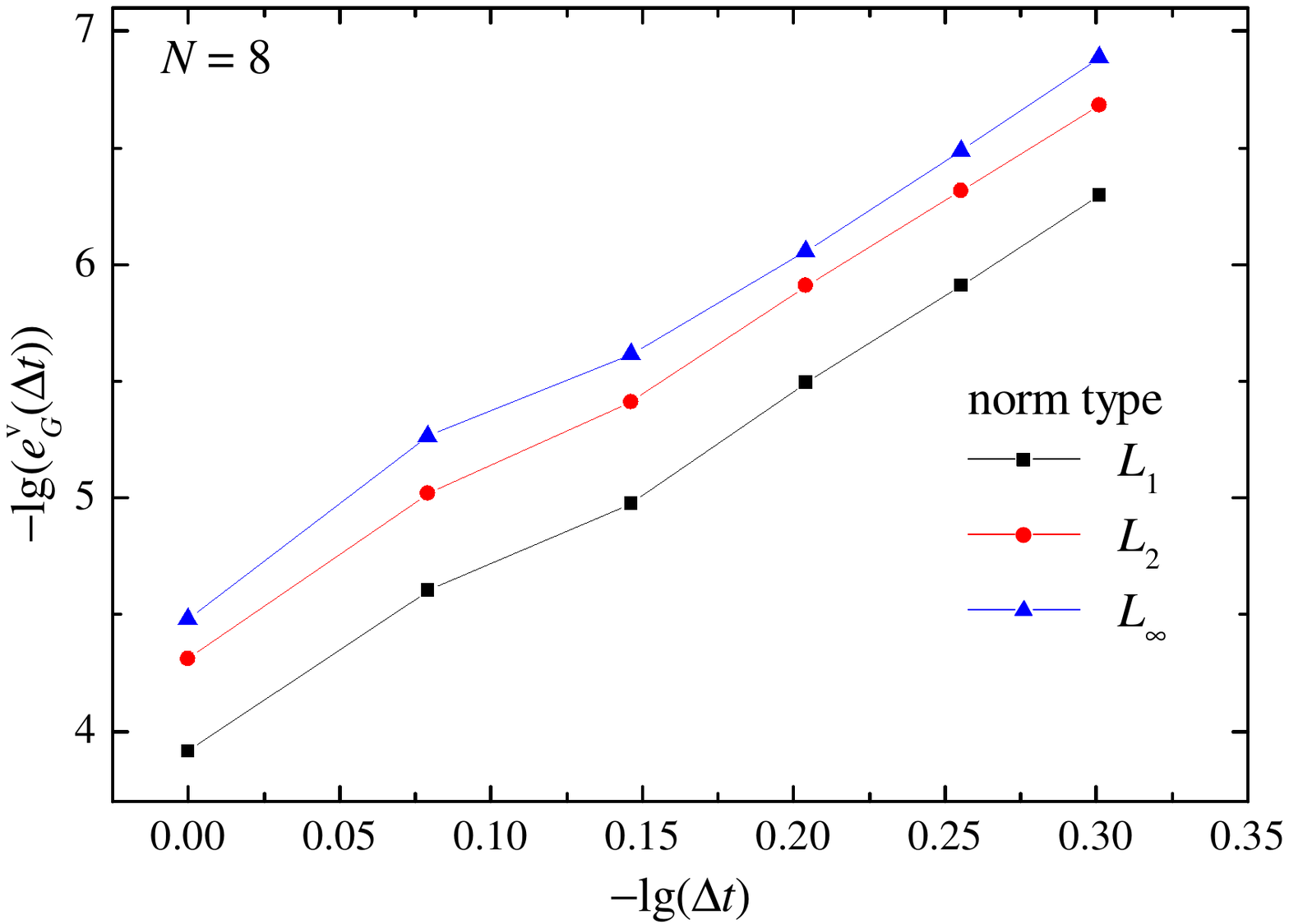}
\vspace{-8mm}\caption{\label{fig:pend_ind3_errors:e2}}
\end{subfigure}\hspace{6mm}
\begin{subfigure}{0.275\textwidth}
\includegraphics[width=\textwidth]{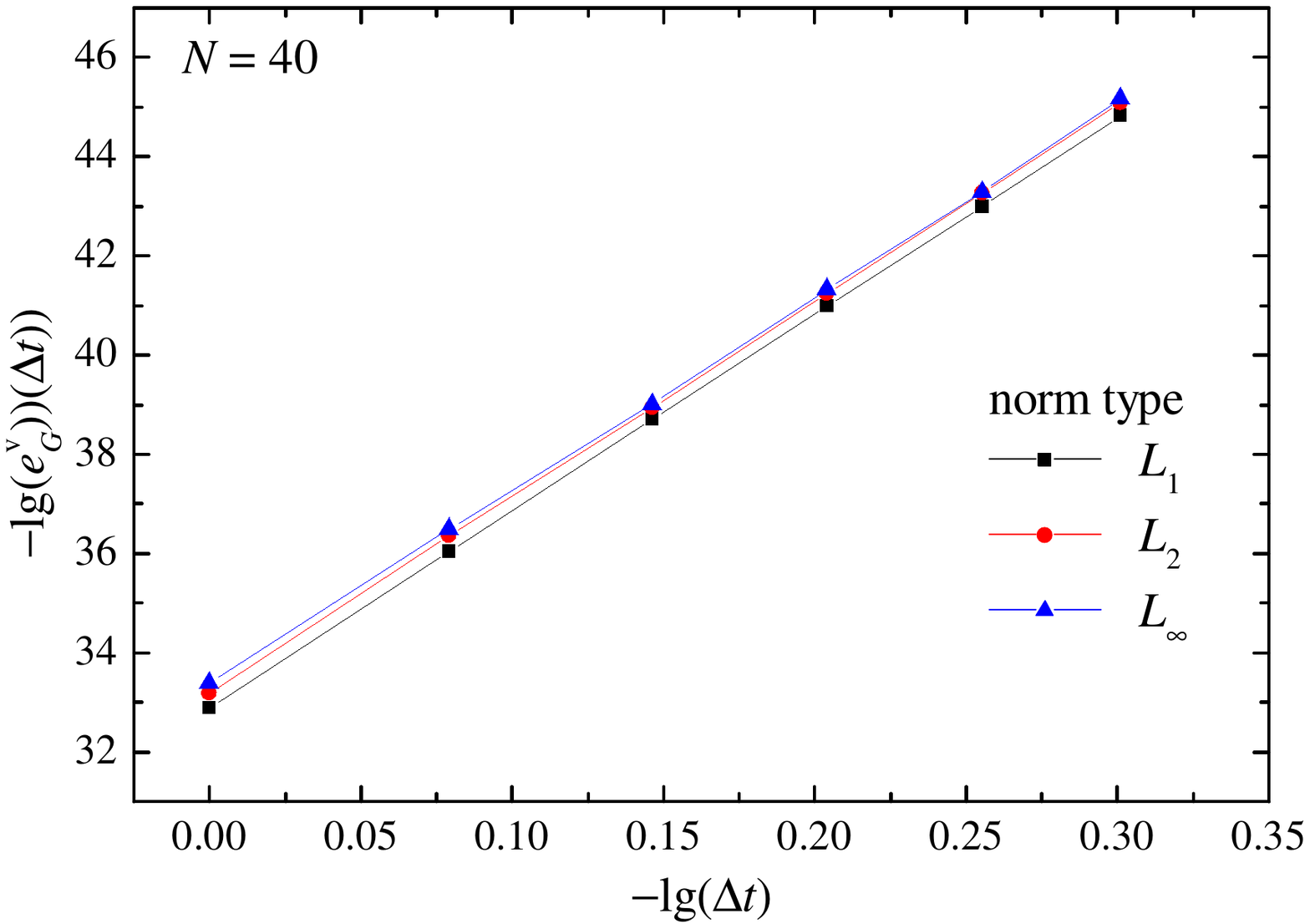}
\vspace{-8mm}\caption{\label{fig:pend_ind3_errors:e3}}
\end{subfigure}\\[-2mm]
\begin{subfigure}{0.275\textwidth}
\includegraphics[width=\textwidth]{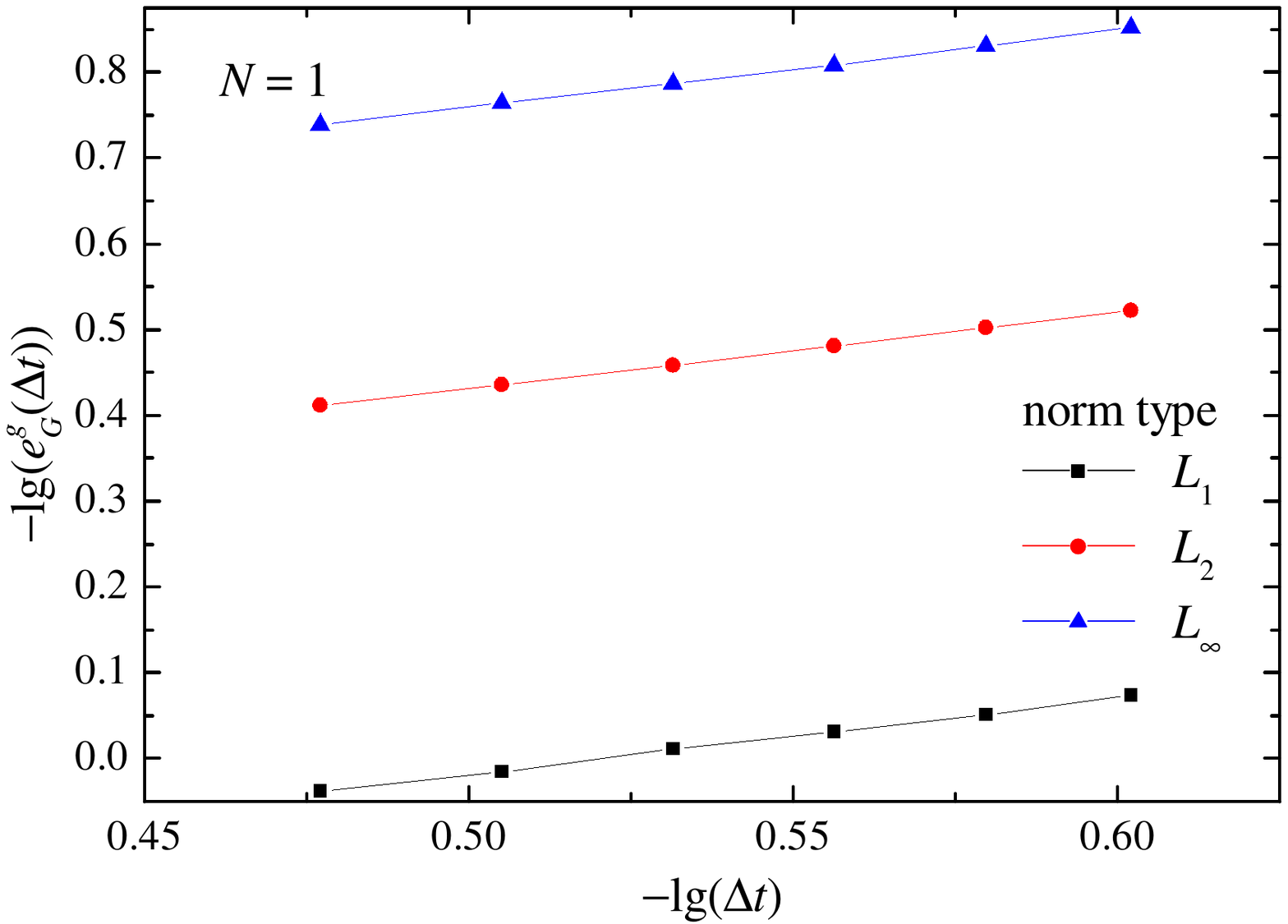}
\vspace{-8mm}\caption{\label{fig:pend_ind3_errors:f1}}
\end{subfigure}\hspace{6mm}
\begin{subfigure}{0.275\textwidth}
\includegraphics[width=\textwidth]{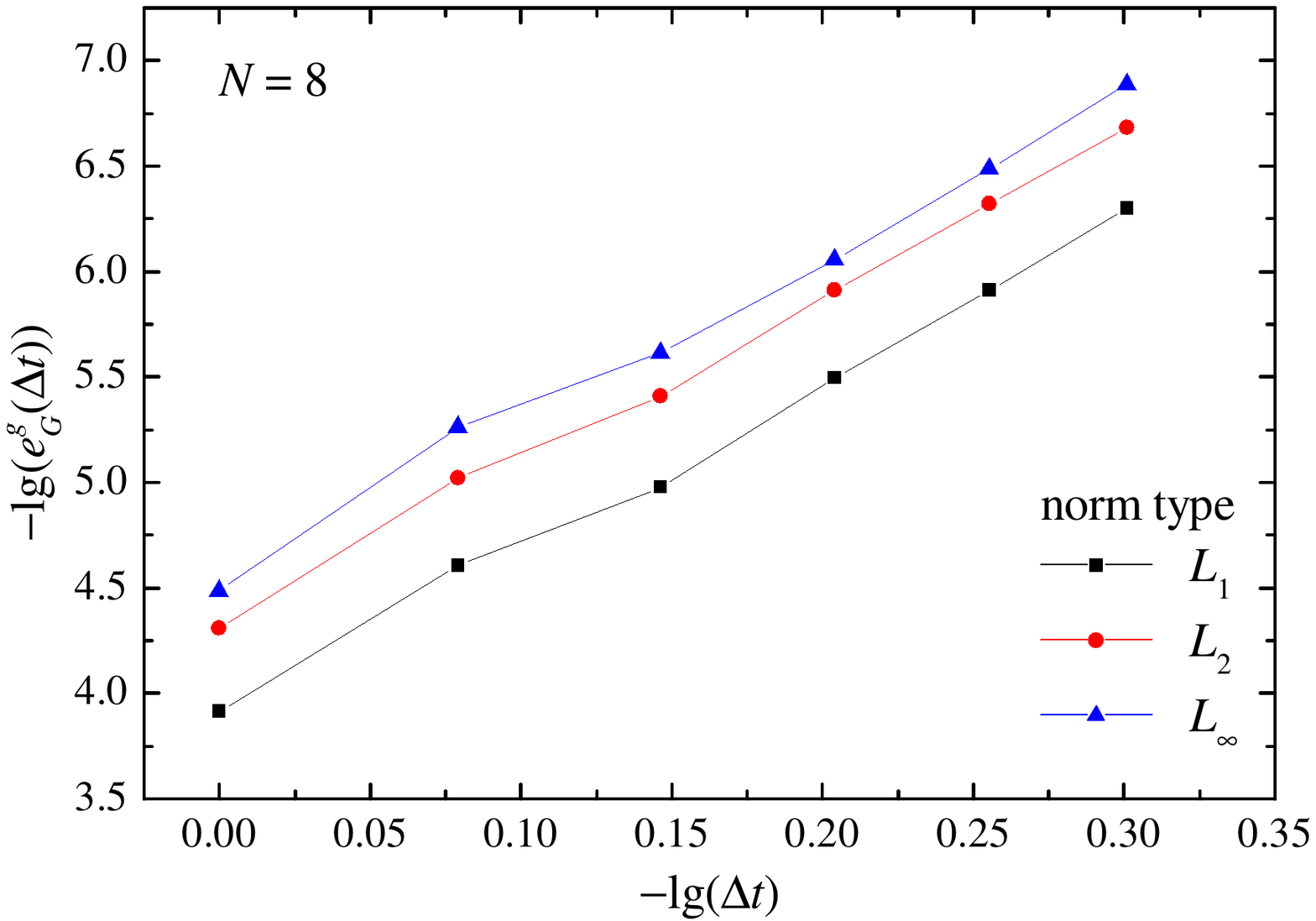}
\vspace{-8mm}\caption{\label{fig:pend_ind3_errors:f2}}
\end{subfigure}\hspace{6mm}
\begin{subfigure}{0.275\textwidth}
\includegraphics[width=\textwidth]{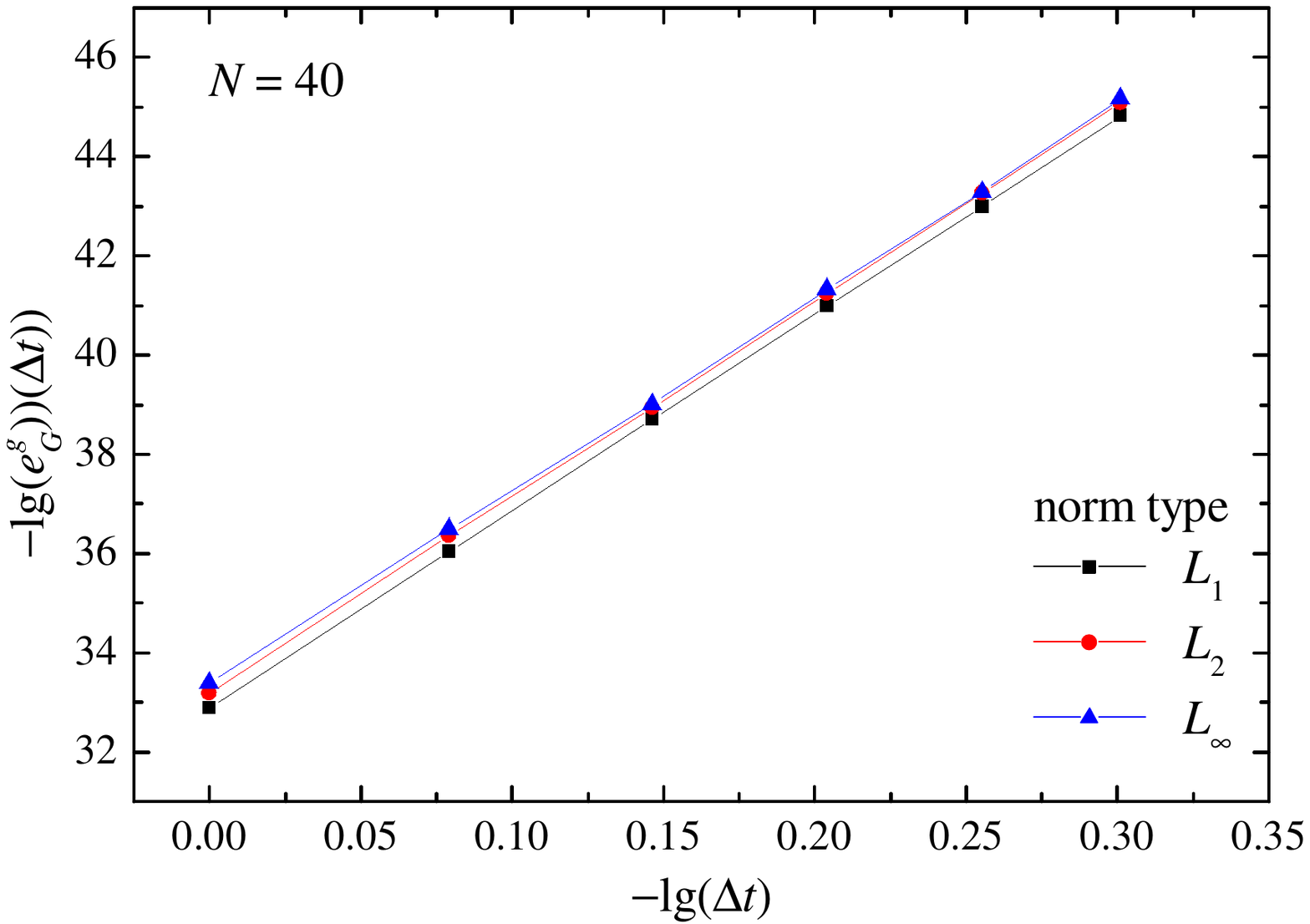}
\vspace{-8mm}\caption{\label{fig:pend_ind3_errors:f3}}
\end{subfigure}\\[-2mm]
\caption{%
Log-log plot of the dependence of the global errors for the local solution $e_{L}^{u}$ (\subref{fig:pend_ind3_errors:a1}, \subref{fig:pend_ind3_errors:a2}, \subref{fig:pend_ind3_errors:a3}), $e_{L}^{v}$ (\subref{fig:pend_ind3_errors:b1}, \subref{fig:pend_ind3_errors:b2}, \subref{fig:pend_ind3_errors:b3}), $e_{L}^{g}$ (\subref{fig:pend_ind3_errors:c1}, \subref{fig:pend_ind3_errors:c2}, \subref{fig:pend_ind3_errors:c3}) and the solution at nodes $e_{G}^{u}$ (\subref{fig:pend_ind3_errors:d1}, \subref{fig:pend_ind3_errors:d2}, \subref{fig:pend_ind3_errors:d3}), $e_{G}^{v}$ (\subref{fig:pend_ind3_errors:e1}, \subref{fig:pend_ind3_errors:e2}, \subref{fig:pend_ind3_errors:e3}), $e_{G}^{g}$ (\subref{fig:pend_ind3_errors:f1}, \subref{fig:pend_ind3_errors:f2}, \subref{fig:pend_ind3_errors:f3}) on the discretization step $\mathrm{\Delta}t$, obtained in the norms $L_{1}$, $L_{2}$ and $L_{\infty}$, by numerical solution of the DAE system (\ref{eq:math_pend_dae_ind_3}) of index 3 obtained using polynomials with degrees $N = 1$ (\subref{fig:pend_ind3_errors:a1}, \subref{fig:pend_ind3_errors:b1}, \subref{fig:pend_ind3_errors:c1}, \subref{fig:pend_ind3_errors:d1}, \subref{fig:pend_ind3_errors:e1}, \subref{fig:pend_ind3_errors:f1}), $N = 8$ (\subref{fig:pend_ind3_errors:a2}, \subref{fig:pend_ind3_errors:b2}, \subref{fig:pend_ind3_errors:c2}, \subref{fig:pend_ind3_errors:d2}, \subref{fig:pend_ind3_errors:e2}, \subref{fig:pend_ind3_errors:f2}) and $N = 40$ (\subref{fig:pend_ind3_errors:a3}, \subref{fig:pend_ind3_errors:b3}, \subref{fig:pend_ind3_errors:c3}, \subref{fig:pend_ind3_errors:d3}, \subref{fig:pend_ind3_errors:e3}, \subref{fig:pend_ind3_errors:f3}).
}
\label{fig:pend_ind3_errors}
\end{figure} 

\begin{table*}[h!]
\centering
\caption{%
Convergence orders $p_{L_{1}}^{n}$, $p_{L_{2}}^{n}$, $p_{L_{\infty}}^{n}$, calculated in norms $L_{1}$, $L_{2}$, $L_{\infty}$, of \textit{the numerical solution at the nodes} $(\mathbf{u}_{n}, \mathbf{v}_{n})$ of the ADER-DG method for the DAE system (\ref{eq:math_pend_dae_ind_3}) of index 3; $N$ is the degree of the basis polynomials $\varphi_{p}$. Orders $p^{n, u}$ are calculated for solution $\mathbf{u}_{n}$; orders $p^{n, v}$ --- for solution $\mathbf{v}_{n}$; orders $p^{n, g}$ --- for the conditions $\mathbf{g} = 0$ on the numerical solution at the nodes $(\mathbf{u}_{n}, \mathbf{v}_{n})$. The theoretical value of convergence order $p_{\rm th.}^{n} = 2N+1$ is applicable for the ADER-DG method for ODE problems and is presented for comparison.
}
\label{tab:conv_orders_nodes_pend_ind3}
\begin{tabular}{@{}|l|lll|lll|lll|c|@{}}
\toprule
$N$ & $p_{L_{1}}^{n, u}$ & $p_{L_{2}}^{n, u}$ & $p_{L_{\infty}}^{n, u}$ & $p_{L_{1}}^{n, v}$ & $p_{L_{2}}^{n, v}$ & $p_{L_{\infty}}^{n, v}$ & $p_{L_{1}}^{n, g}$ & $p_{L_{2}}^{n, g}$ & $p_{L_{\infty}}^{n, g}$ & $p_{\rm th.}^{n}$ \\
\midrule
$1$	&	$2.10$	&	$2.04$	&	$1.90$	&	$0.90$	&	$0.90$	&	$0.91$	&	$0.88$	&	$0.89$	&	$0.90$	&	$3$\\
$2$	&	$2.99$	&	$2.99$	&	$2.97$	&	$2.20$	&	$2.19$	&	$2.20$	&	$2.20$	&	$2.19$	&	$2.20$	&	$5$\\
$3$	&	$4.00$	&	$3.99$	&	$4.00$	&	$2.95$	&	$2.94$	&	$2.93$	&	$2.95$	&	$2.94$	&	$2.93$	&	$7$\\
$4$	&	$4.94$	&	$4.98$	&	$4.96$	&	$4.10$	&	$4.12$	&	$4.15$	&	$4.10$	&	$4.12$	&	$4.15$	&	$9$\\
$5$	&	$5.98$	&	$5.98$	&	$6.03$	&	$4.90$	&	$4.90$	&	$4.86$	&	$4.90$	&	$4.90$	&	$4.86$	&	$11$\\
$6$	&	$6.68$	&	$6.71$	&	$6.82$	&	$5.90$	&	$5.91$	&	$5.91$	&	$5.90$	&	$5.91$	&	$5.91$	&	$13$\\
$7$	&	$8.09$	&	$8.01$	&	$7.92$	&	$7.01$	&	$6.96$	&	$6.86$	&	$7.01$	&	$6.96$	&	$6.86$	&	$15$\\
$8$	&	$8.64$	&	$8.64$	&	$8.63$	&	$7.81$	&	$7.78$	&	$7.71$	&	$7.81$	&	$7.78$	&	$7.71$	&	$17$\\
$9$	&	$9.94$	&	$9.95$	&	$9.93$	&	$8.90$	&	$9.01$	&	$8.97$	&	$8.90$	&	$9.01$	&	$8.97$	&	$19$\\
$10$	&	$10.76$	&	$10.71$	&	$10.52$	&	$9.80$	&	$9.74$	&	$9.66$	&	$9.80$	&	$9.74$	&	$9.66$	&	$21$\\
$11$	&	$11.77$	&	$11.87$	&	$11.90$	&	$11.02$	&	$11.05$	&	$11.02$	&	$11.02$	&	$11.05$	&	$11.02$	&	$23$\\
$12$	&	$12.83$	&	$12.76$	&	$12.50$	&	$11.69$	&	$11.70$	&	$11.58$	&	$11.69$	&	$11.70$	&	$11.58$	&	$25$\\
$13$	&	$13.83$	&	$13.86$	&	$13.60$	&	$13.08$	&	$13.07$	&	$12.89$	&	$13.08$	&	$13.07$	&	$12.89$	&	$27$\\
$14$	&	$14.78$	&	$14.75$	&	$14.54$	&	$13.63$	&	$13.67$	&	$13.57$	&	$13.63$	&	$13.67$	&	$13.57$	&	$29$\\
$15$	&	$15.96$	&	$15.84$	&	$15.61$	&	$15.19$	&	$15.02$	&	$14.72$	&	$15.19$	&	$15.02$	&	$14.72$	&	$31$\\
$16$	&	$16.67$	&	$16.76$	&	$16.79$	&	$15.54$	&	$15.71$	&	$15.64$	&	$15.54$	&	$15.71$	&	$15.64$	&	$33$\\
$17$	&	$17.84$	&	$17.73$	&	$17.59$	&	$17.00$	&	$16.87$	&	$16.70$	&	$17.00$	&	$16.87$	&	$16.70$	&	$35$\\
$18$	&	$18.83$	&	$18.83$	&	$18.85$	&	$17.81$	&	$17.83$	&	$17.80$	&	$17.81$	&	$17.83$	&	$17.80$	&	$37$\\
$19$	&	$19.59$	&	$19.56$	&	$19.53$	&	$18.69$	&	$18.64$	&	$18.61$	&	$18.69$	&	$18.64$	&	$18.61$	&	$39$\\
$20$	&	$20.95$	&	$20.92$	&	$20.79$	&	$20.00$	&	$19.97$	&	$19.85$	&	$20.00$	&	$19.97$	&	$19.85$	&	$41$\\
\midrule
$25$	&	$25.47$	&	$25.43$	&	$25.51$	&	$24.48$	&	$24.42$	&	$24.50$	&	$24.48$	&	$24.42$	&	$24.50$	&	$51$\\
$30$	&	$30.18$	&	$30.11$	&	$29.91$	&	$29.20$	&	$29.12$	&	$28.88$	&	$29.20$	&	$29.12$	&	$28.88$	&	$61$\\
$35$	&	$35.66$	&	$35.53$	&	$35.39$	&	$34.77$	&	$34.61$	&	$34.42$	&	$34.77$	&	$34.61$	&	$34.42$	&	$71$\\
$40$	&	$40.62$	&	$40.46$	&	$39.86$	&	$39.66$	&	$39.50$	&	$39.01$	&	$39.66$	&	$39.50$	&	$39.01$	&	$81$\\
\bottomrule
\end{tabular}
\end{table*} 

\begin{table*}[h!]
\centering
\caption{%
Convergence orders $p_{L_{1}}^{l}$, $p_{L_{2}}^{l}$, $p_{L_{\infty}}^{l}$, calculated in norms $L_{1}$, $L_{2}$, $L_{\infty}$, of \textit{the local solution} $(\mathbf{u}_{L}, \mathbf{v}_{L})$ (represented between the nodes) of the ADER-DG method for the DAE system (\ref{eq:math_pend_dae_ind_3}) of index 3; $N$ is the degree of the basis polynomials $\varphi_{p}$. Orders $p^{l, u}$ are calculated for solution $\mathbf{u}_{L}$; orders $p^{l, v}$ --- for solution $\mathbf{v}_{L}$; orders $p^{l, g}$ --- for the conditions $\mathbf{g} = 0$ on the local solution $(\mathbf{u}_{L}, \mathbf{v}_{L})$. The theoretical value of convergence order $p_{\rm th.}^{l} = N+1$ is applicable for the ADER-DG method for ODE problems and is presented for comparison.
}
\label{tab:conv_orders_local_pend_ind3}
\begin{tabular}{@{}|l|lll|lll|lll|c|@{}}
\toprule
$N$ & $p_{L_{1}}^{l, u}$ & $p_{L_{2}}^{l, u}$ & $p_{L_{\infty}}^{l, u}$ & $p_{L_{1}}^{l, v}$ & $p_{L_{2}}^{l, v}$ & $p_{L_{\infty}}^{l, v}$ & $p_{L_{1}}^{l, g}$ & $p_{L_{2}}^{l, g}$ & $p_{L_{\infty}}^{l, g}$ & $p_{\rm th.}^{l}$ \\
\midrule
$1$	&	$2.19$	&	$2.16$	&	$1.80$	&	$0.92$	&	$0.92$	&	$0.97$	&	$0.93$	&	$0.92$	&	$0.98$	&	$2$\\
$2$	&	$3.04$	&	$3.04$	&	$3.15$	&	$2.06$	&	$2.05$	&	$2.14$	&	$2.08$	&	$2.07$	&	$2.17$	&	$3$\\
$3$	&	$3.98$	&	$3.97$	&	$3.89$	&	$2.98$	&	$2.95$	&	$2.97$	&	$2.99$	&	$2.95$	&	$2.97$	&	$4$\\
$4$	&	$5.01$	&	$5.03$	&	$5.06$	&	$4.05$	&	$4.06$	&	$4.14$	&	$4.06$	&	$4.06$	&	$4.15$	&	$5$\\
$5$	&	$5.97$	&	$5.95$	&	$5.87$	&	$4.97$	&	$4.92$	&	$4.84$	&	$4.97$	&	$4.92$	&	$4.84$	&	$6$\\
$6$	&	$7.07$	&	$7.07$	&	$6.71$	&	$5.83$	&	$5.83$	&	$5.90$	&	$5.86$	&	$5.84$	&	$5.91$	&	$7$\\
$7$	&	$7.97$	&	$7.98$	&	$7.81$	&	$7.16$	&	$7.04$	&	$6.84$	&	$7.18$	&	$7.06$	&	$6.87$	&	$8$\\
$8$	&	$8.96$	&	$8.90$	&	$8.44$	&	$7.72$	&	$7.70$	&	$7.65$	&	$7.75$	&	$7.70$	&	$7.65$	&	$9$\\
$9$	&	$10.01$	&	$10.10$	&	$10.01$	&	$9.14$	&	$9.06$	&	$8.99$	&	$9.16$	&	$9.07$	&	$9.01$	&	$10$\\
$10$	&	$10.83$	&	$10.76$	&	$10.46$	&	$9.75$	&	$9.66$	&	$9.54$	&	$9.76$	&	$9.65$	&	$9.53$	&	$11$\\
$11$	&	$12.02$	&	$12.12$	&	$12.11$	&	$11.03$	&	$11.07$	&	$11.05$	&	$11.05$	&	$11.09$	&	$11.06$	&	$12$\\
$12$	&	$12.79$	&	$12.73$	&	$12.59$	&	$11.82$	&	$11.64$	&	$11.49$	&	$11.82$	&	$11.63$	&	$11.49$	&	$13$\\
$13$	&	$13.95$	&	$14.00$	&	$13.90$	&	$13.10$	&	$13.08$	&	$12.90$	&	$13.11$	&	$13.08$	&	$12.90$	&	$14$\\
$14$	&	$14.86$	&	$14.84$	&	$14.68$	&	$13.75$	&	$13.65$	&	$13.50$	&	$13.75$	&	$13.65$	&	$13.50$	&	$15$\\
$15$	&	$15.90$	&	$15.81$	&	$15.69$	&	$15.18$	&	$15.01$	&	$14.72$	&	$15.18$	&	$15.01$	&	$14.72$	&	$16$\\
$16$	&	$16.97$	&	$17.03$	&	$16.78$	&	$15.84$	&	$15.74$	&	$15.69$	&	$15.85$	&	$15.74$	&	$15.70$	&	$17$\\
$17$	&	$17.81$	&	$17.64$	&	$17.42$	&	$16.98$	&	$16.84$	&	$16.63$	&	$16.99$	&	$16.84$	&	$16.63$	&	$18$\\
$18$	&	$19.05$	&	$19.17$	&	$18.86$	&	$18.00$	&	$17.87$	&	$17.79$	&	$18.01$	&	$17.88$	&	$17.79$	&	$19$\\
$19$	&	$19.73$	&	$19.57$	&	$19.29$	&	$18.64$	&	$18.62$	&	$18.55$	&	$18.65$	&	$18.61$	&	$18.55$	&	$20$\\
$20$	&	$21.07$	&	$21.19$	&	$21.00$	&	$20.15$	&	$20.00$	&	$19.84$	&	$20.15$	&	$20.01$	&	$19.84$	&	$21$\\
\midrule
$25$	&	$25.74$	&	$25.78$	&	$25.82$	&	$24.52$	&	$24.43$	&	$24.49$	&	$24.52$	&	$24.43$	&	$24.49$	&	$26$\\
$30$	&	$30.63$	&	$30.51$	&	$30.31$	&	$29.37$	&	$29.18$	&	$28.96$	&	$29.38$	&	$29.18$	&	$28.96$	&	$31$\\
$35$	&	$35.82$	&	$35.72$	&	$35.32$	&	$34.83$	&	$34.61$	&	$34.42$	&	$34.83$	&	$34.61$	&	$34.42$	&	$36$\\
$40$	&	$40.60$	&	$40.31$	&	$39.78$	&	$39.67$	&	$39.50$	&	$38.97$	&	$39.67$	&	$39.50$	&	$38.97$	&	$41$\\
\bottomrule
\end{tabular}
\end{table*}

\begin{figure}[h!]
\captionsetup[subfigure]{%
	position=bottom,
	font+=smaller,
	textfont=normalfont,
	singlelinecheck=off,
	justification=raggedright
}
\centering
\begin{subfigure}{0.320\textwidth}
\includegraphics[width=\textwidth]{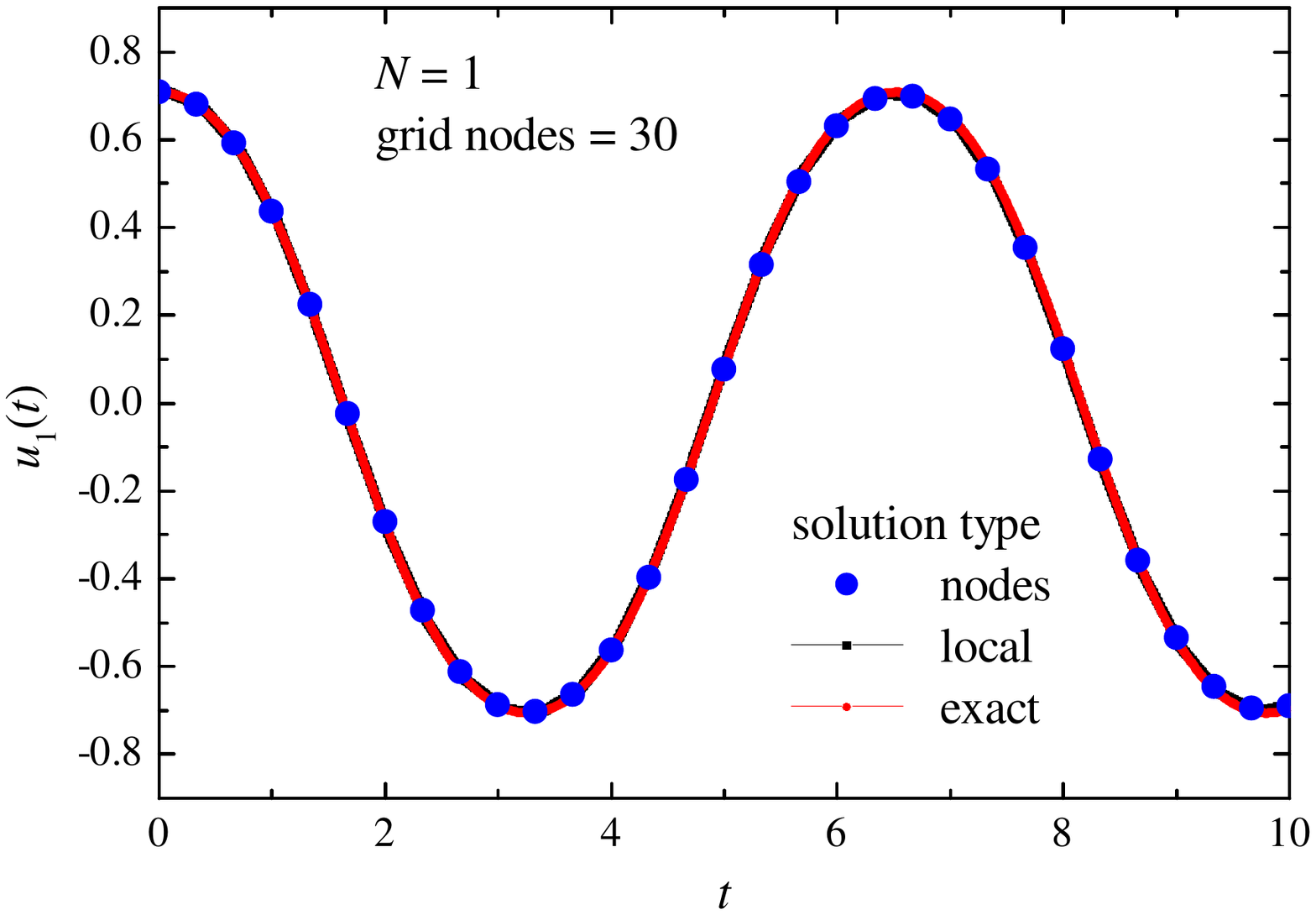}
\vspace{-8mm}\caption{\label{fig:pend_ind2_sol_uv:a1}}
\end{subfigure}
\begin{subfigure}{0.320\textwidth}
\includegraphics[width=\textwidth]{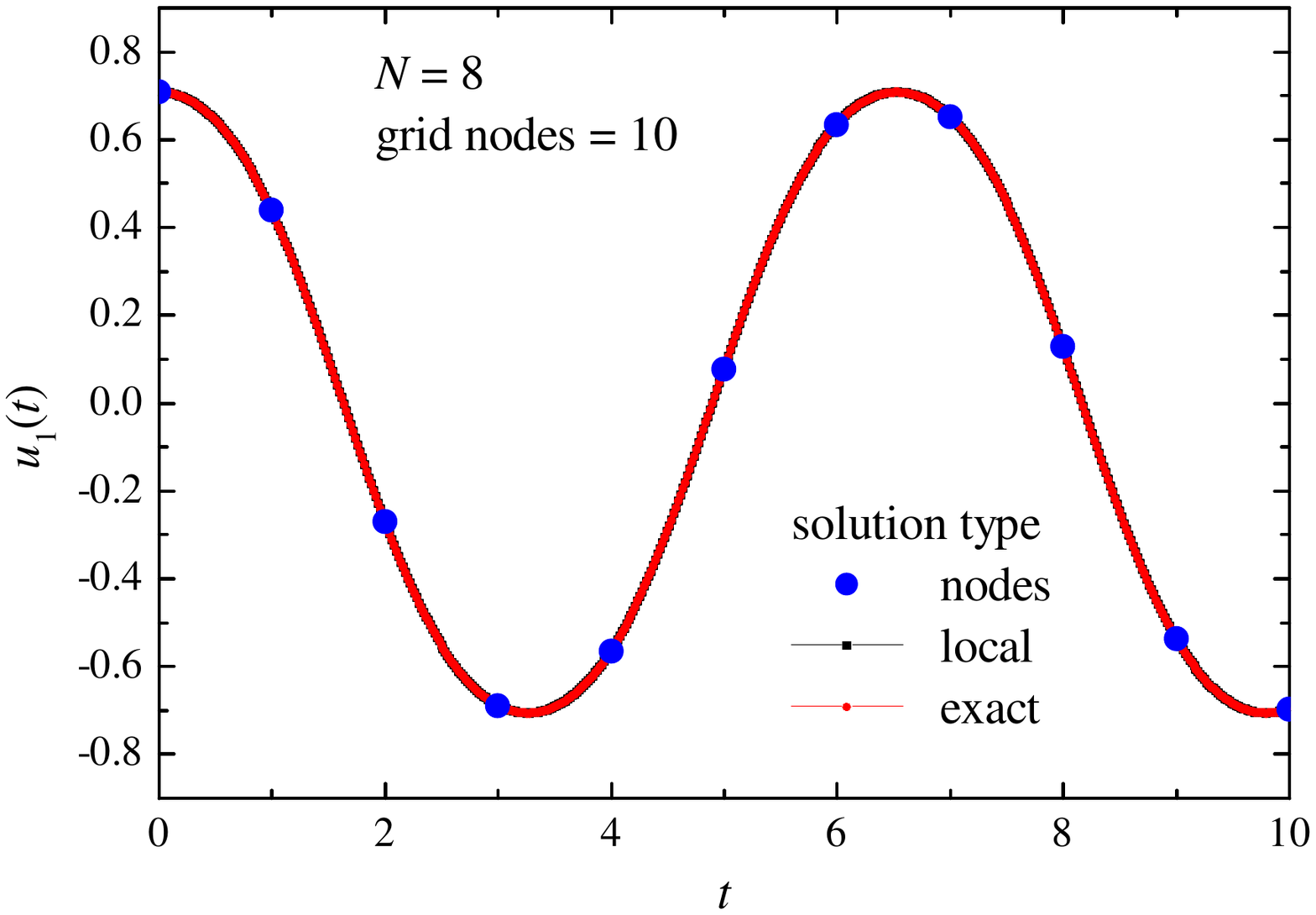}
\vspace{-8mm}\caption{\label{fig:pend_ind2_sol_uv:a2}}
\end{subfigure}
\begin{subfigure}{0.320\textwidth}
\includegraphics[width=\textwidth]{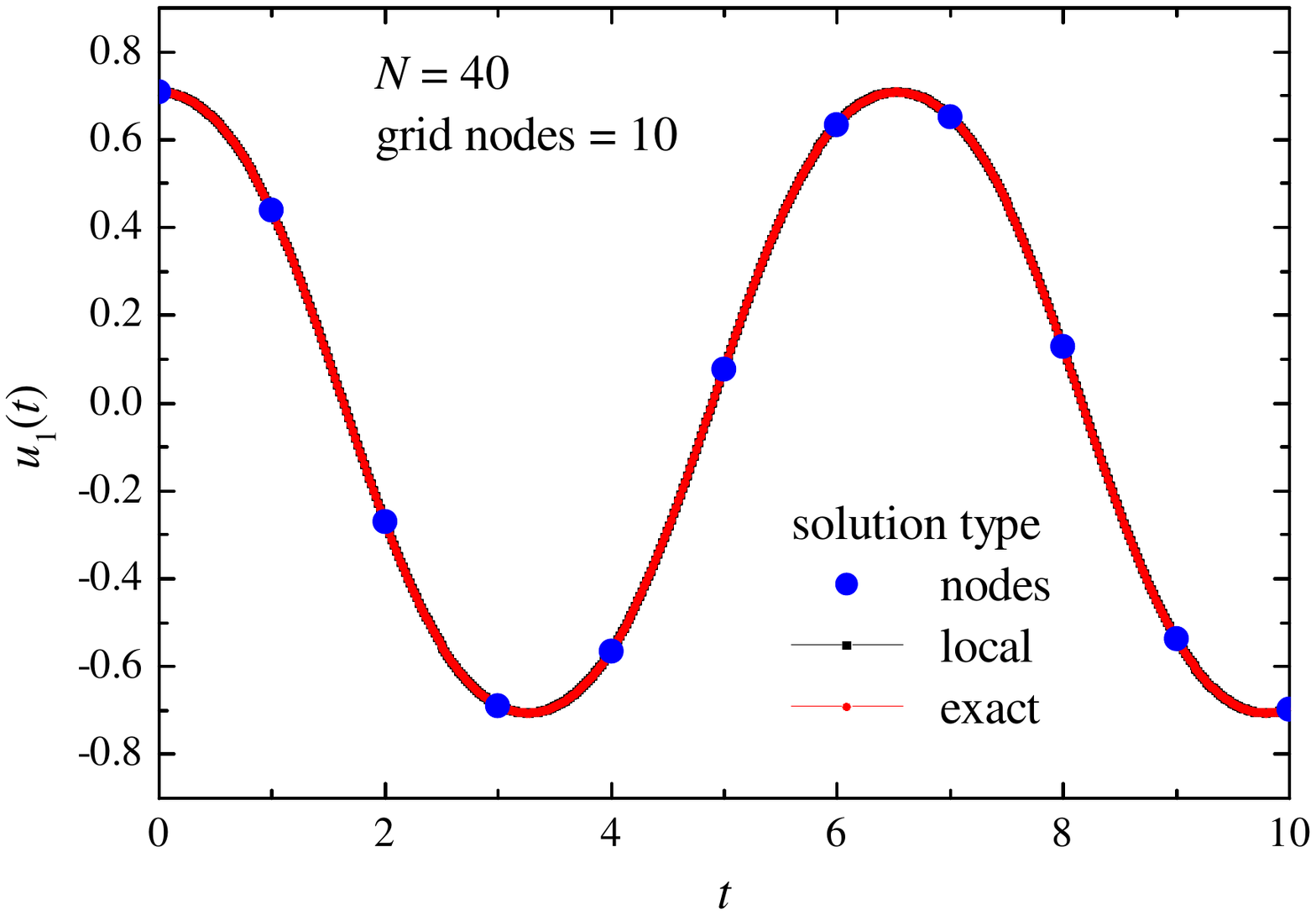}
\vspace{-8mm}\caption{\label{fig:pend_ind2_sol_uv:a3}}
\end{subfigure}\\
\begin{subfigure}{0.320\textwidth}
\includegraphics[width=\textwidth]{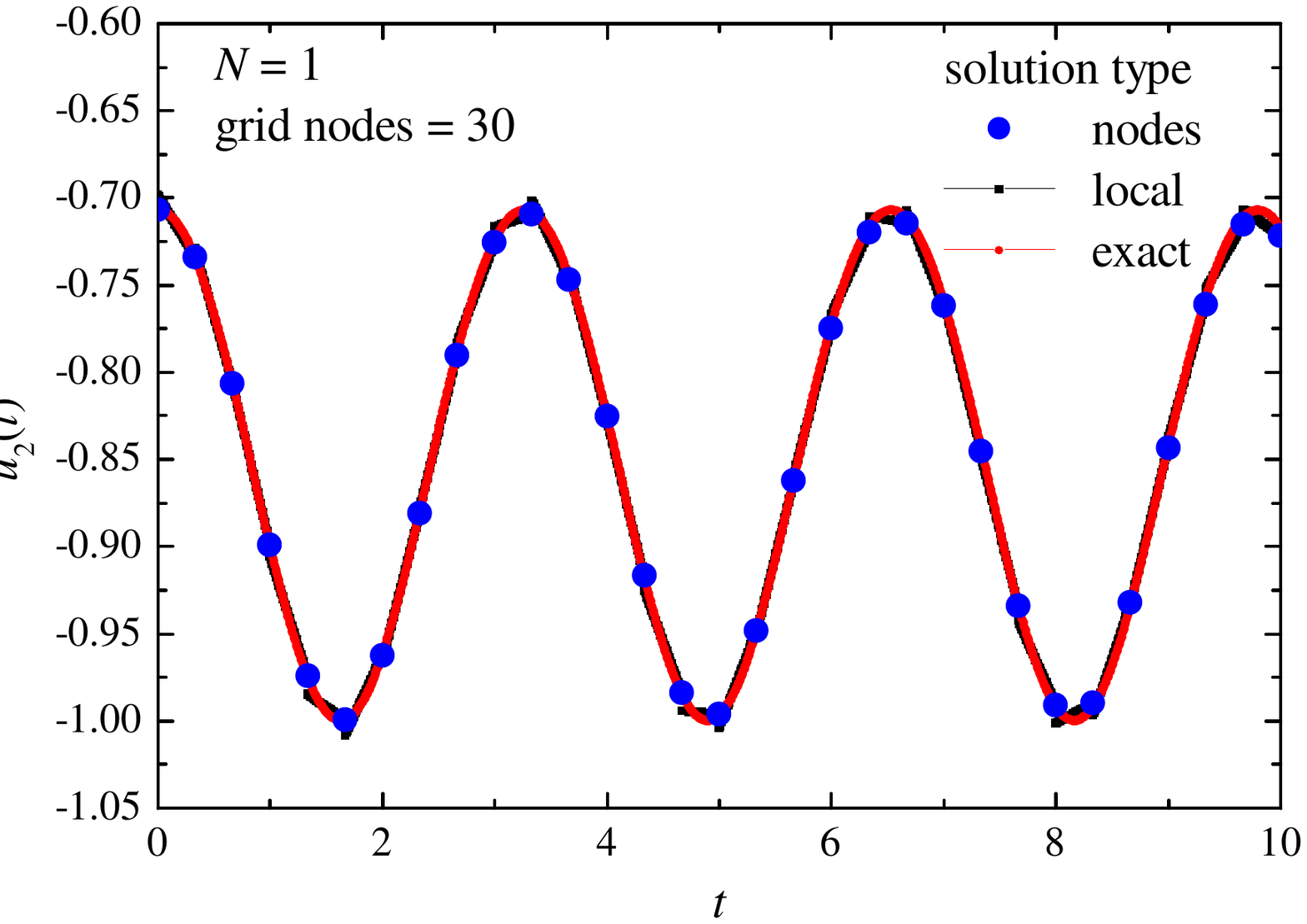}
\vspace{-8mm}\caption{\label{fig:pend_ind2_sol_uv:b1}}
\end{subfigure}
\begin{subfigure}{0.320\textwidth}
\includegraphics[width=\textwidth]{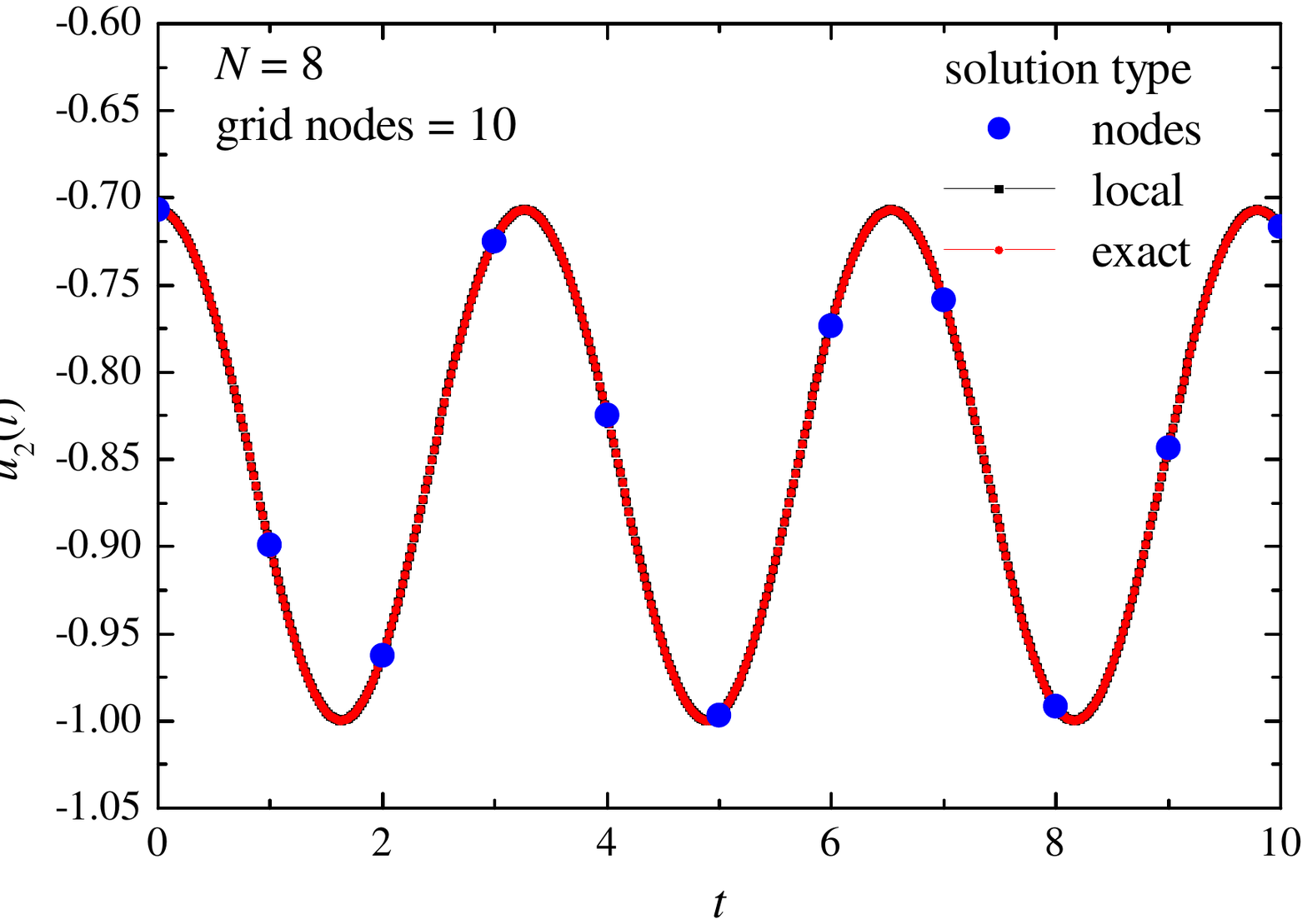}
\vspace{-8mm}\caption{\label{fig:pend_ind2_sol_uv:b2}}
\end{subfigure}
\begin{subfigure}{0.320\textwidth}
\includegraphics[width=\textwidth]{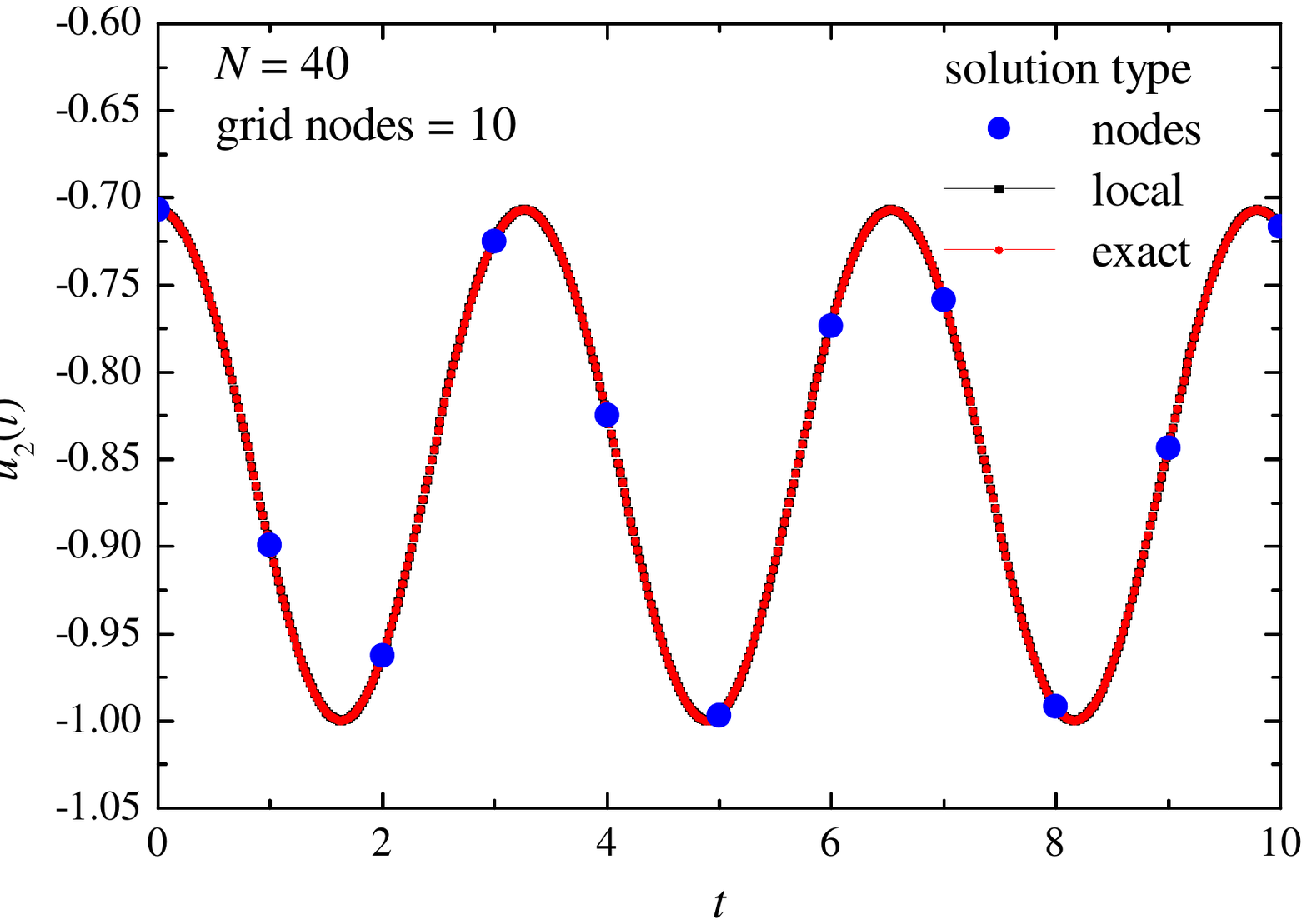}
\vspace{-8mm}\caption{\label{fig:pend_ind2_sol_uv:b3}}
\end{subfigure}\\
\begin{subfigure}{0.320\textwidth}
\includegraphics[width=\textwidth]{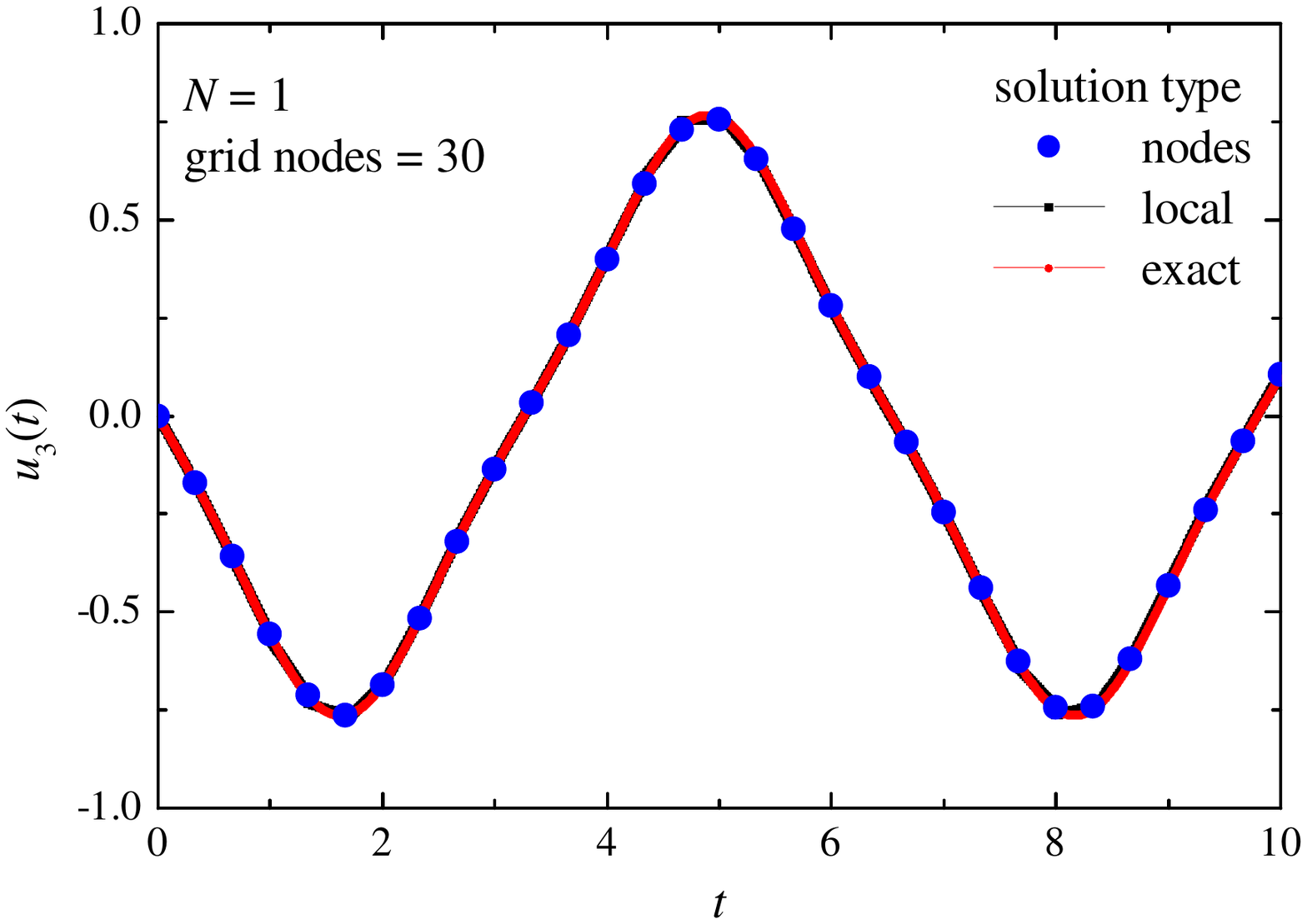}
\vspace{-8mm}\caption{\label{fig:pend_ind2_sol_uv:c1}}
\end{subfigure}
\begin{subfigure}{0.320\textwidth}
\includegraphics[width=\textwidth]{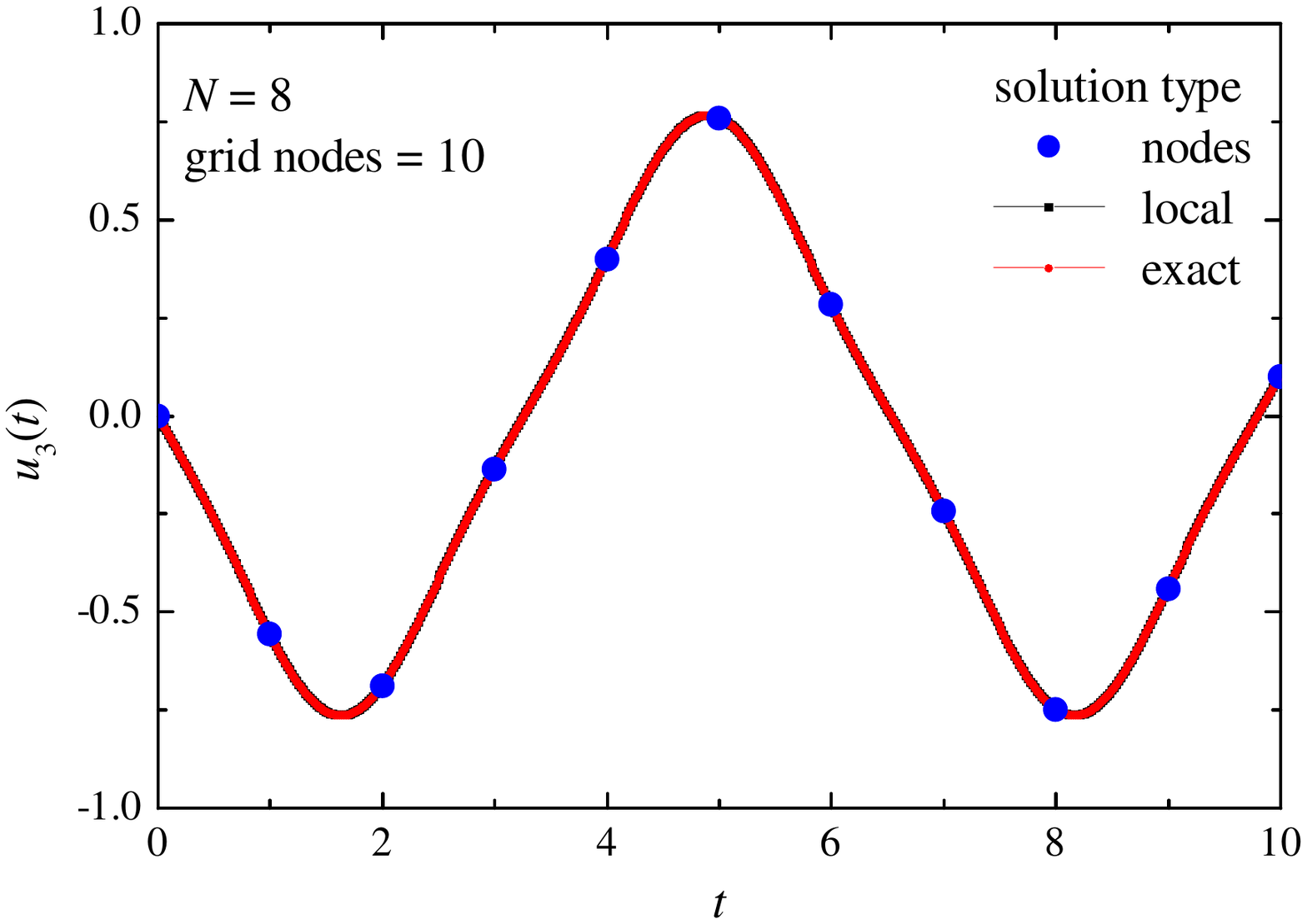}
\vspace{-8mm}\caption{\label{fig:pend_ind2_sol_uv:c2}}
\end{subfigure}
\begin{subfigure}{0.320\textwidth}
\includegraphics[width=\textwidth]{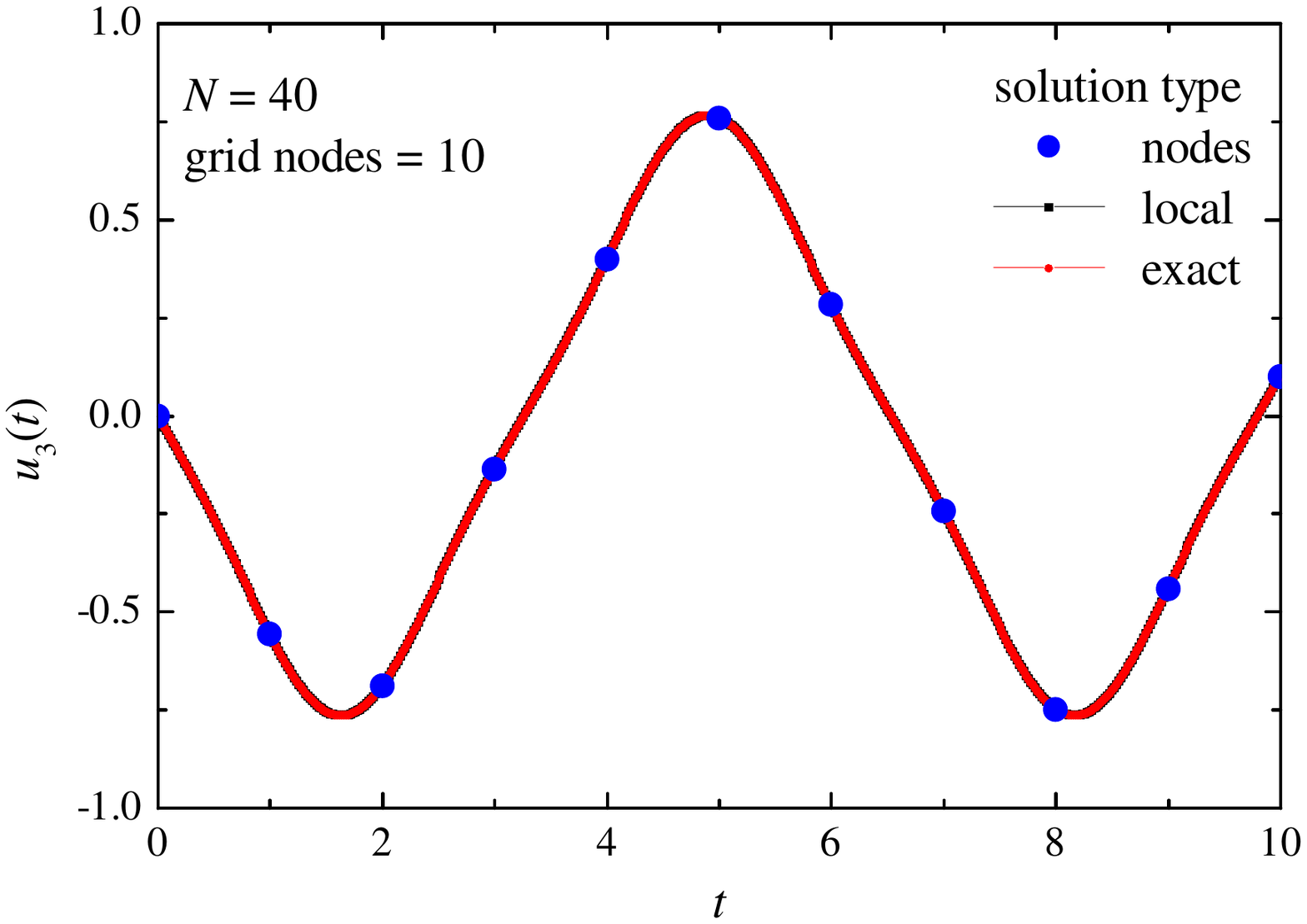}
\vspace{-8mm}\caption{\label{fig:pend_ind2_sol_uv:c3}}
\end{subfigure}\\
\begin{subfigure}{0.320\textwidth}
\includegraphics[width=\textwidth]{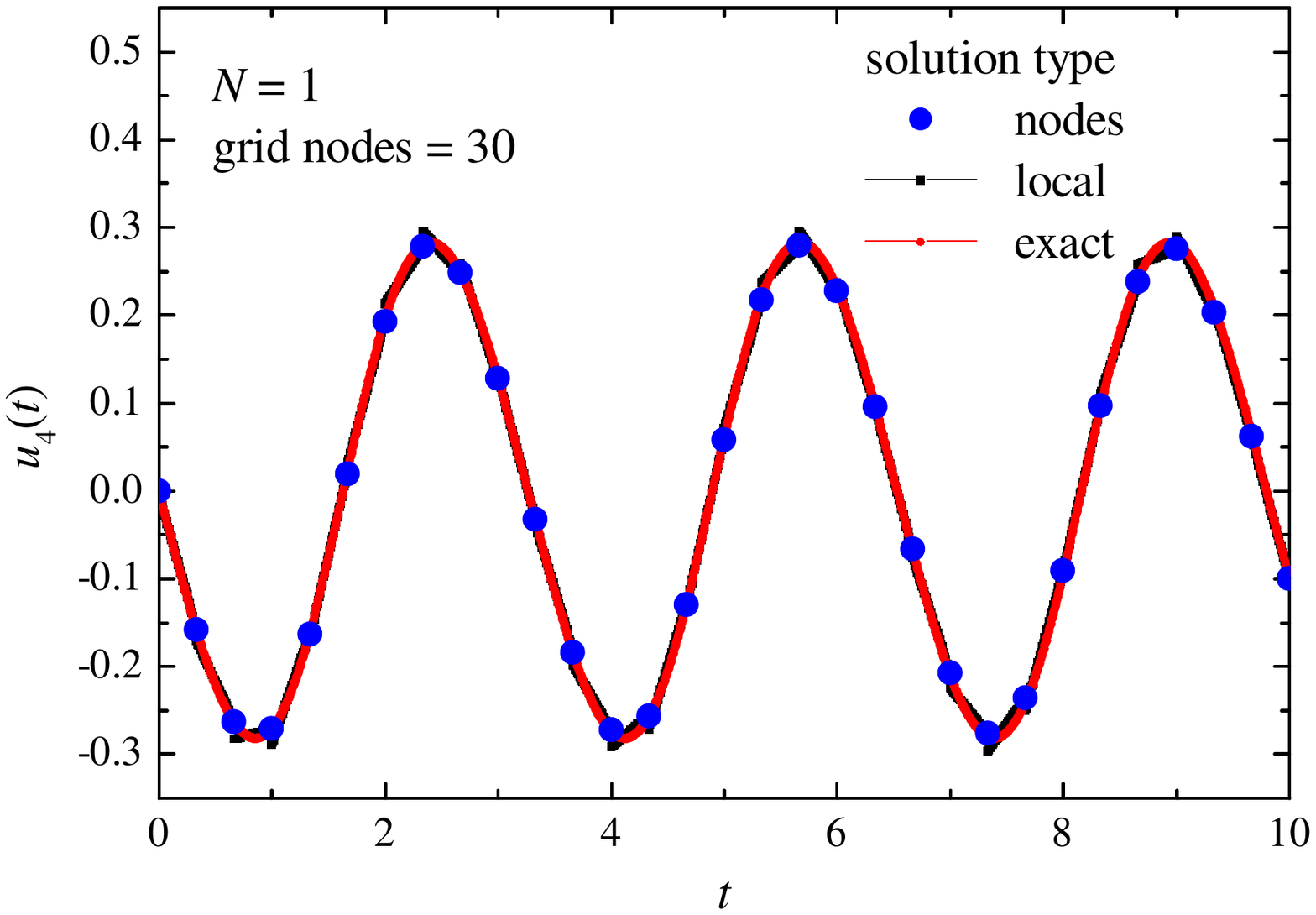}
\vspace{-8mm}\caption{\label{fig:pend_ind2_sol_uv:d1}}
\end{subfigure}
\begin{subfigure}{0.320\textwidth}
\includegraphics[width=\textwidth]{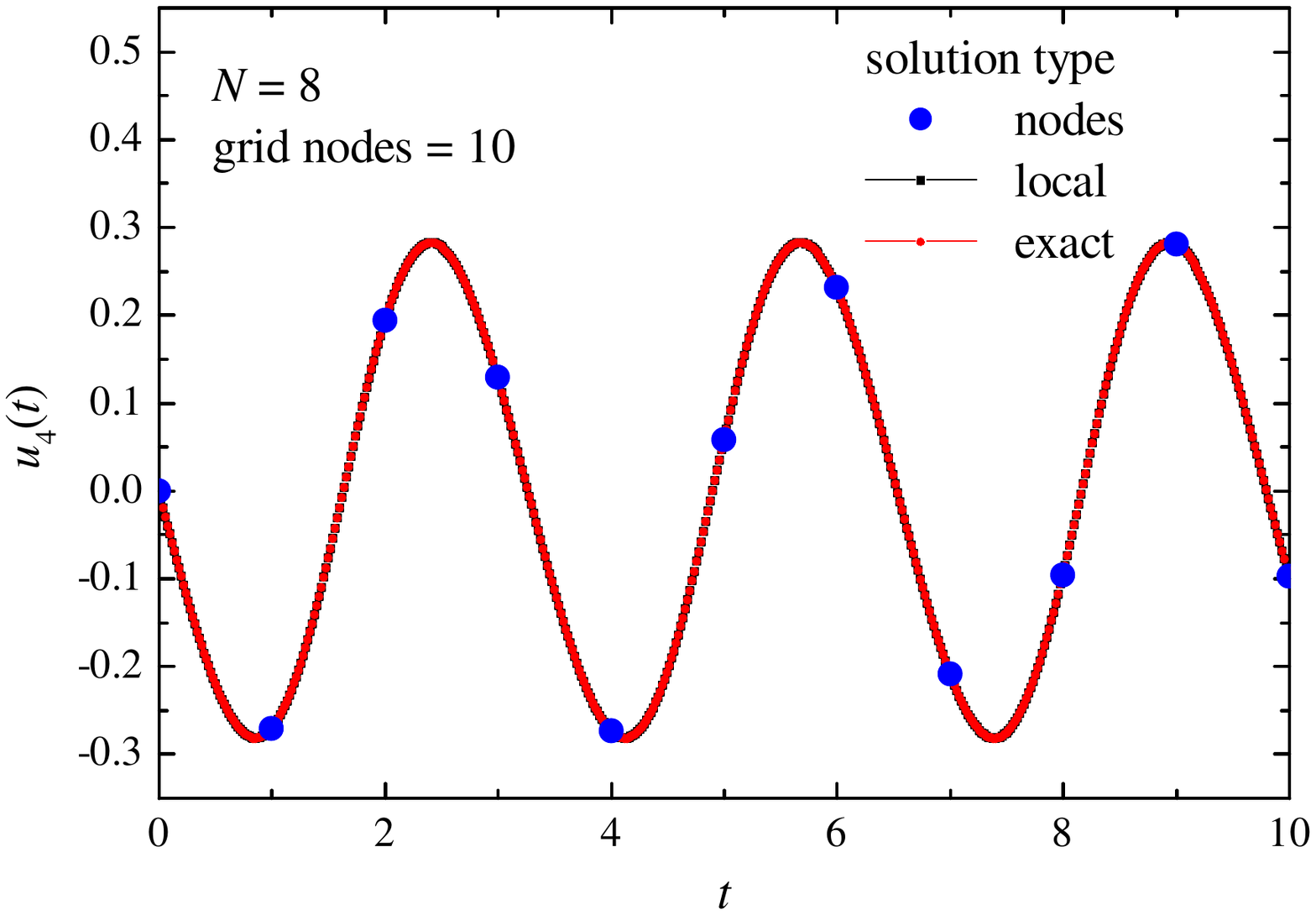}
\vspace{-8mm}\caption{\label{fig:pend_ind2_sol_uv:d2}}
\end{subfigure}
\begin{subfigure}{0.320\textwidth}
\includegraphics[width=\textwidth]{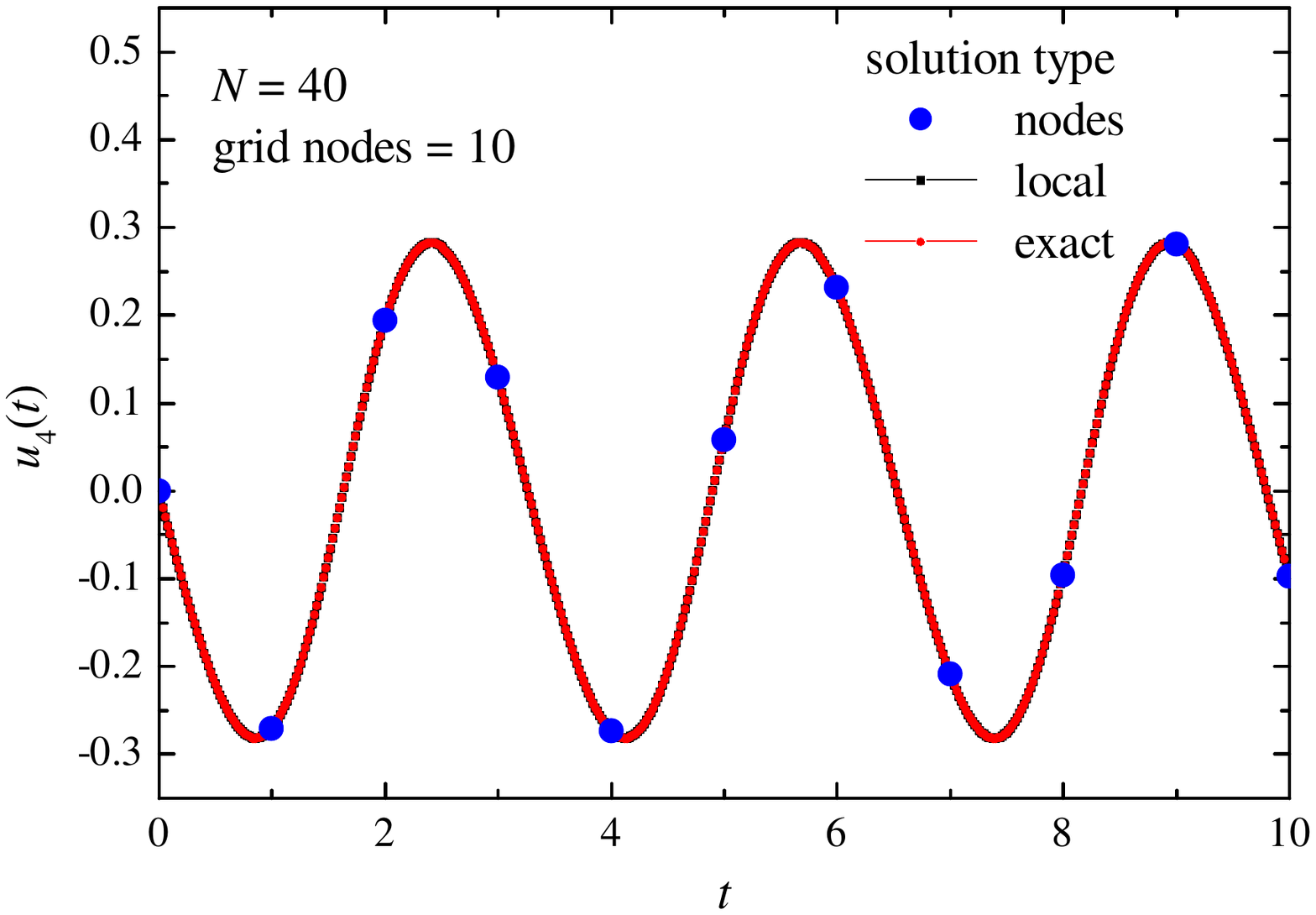}
\vspace{-8mm}\caption{\label{fig:pend_ind2_sol_uv:d3}}
\end{subfigure}\\
\begin{subfigure}{0.320\textwidth}
\includegraphics[width=\textwidth]{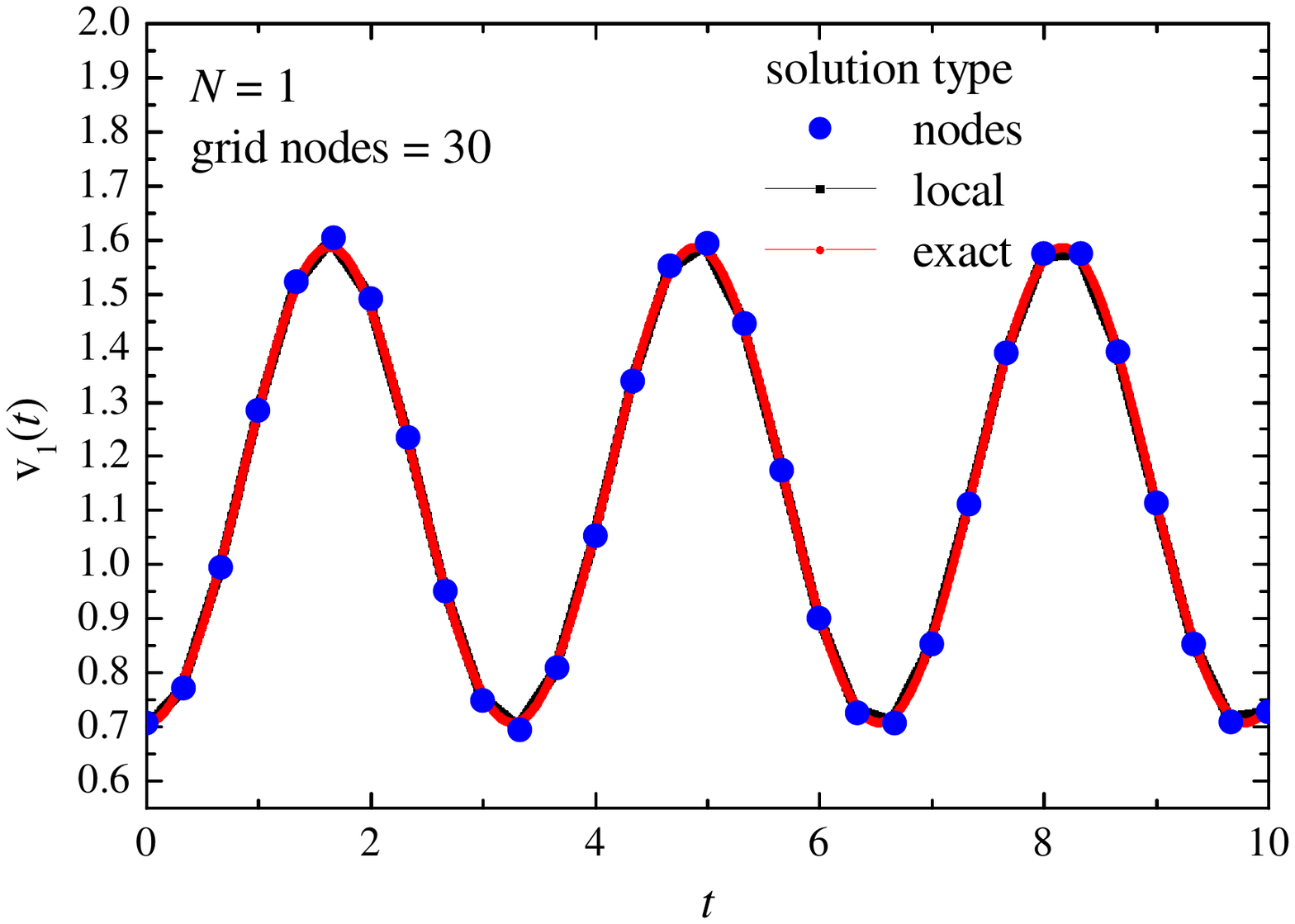}
\vspace{-8mm}\caption{\label{fig:pend_ind2_sol_uv:e1}}
\end{subfigure}
\begin{subfigure}{0.320\textwidth}
\includegraphics[width=\textwidth]{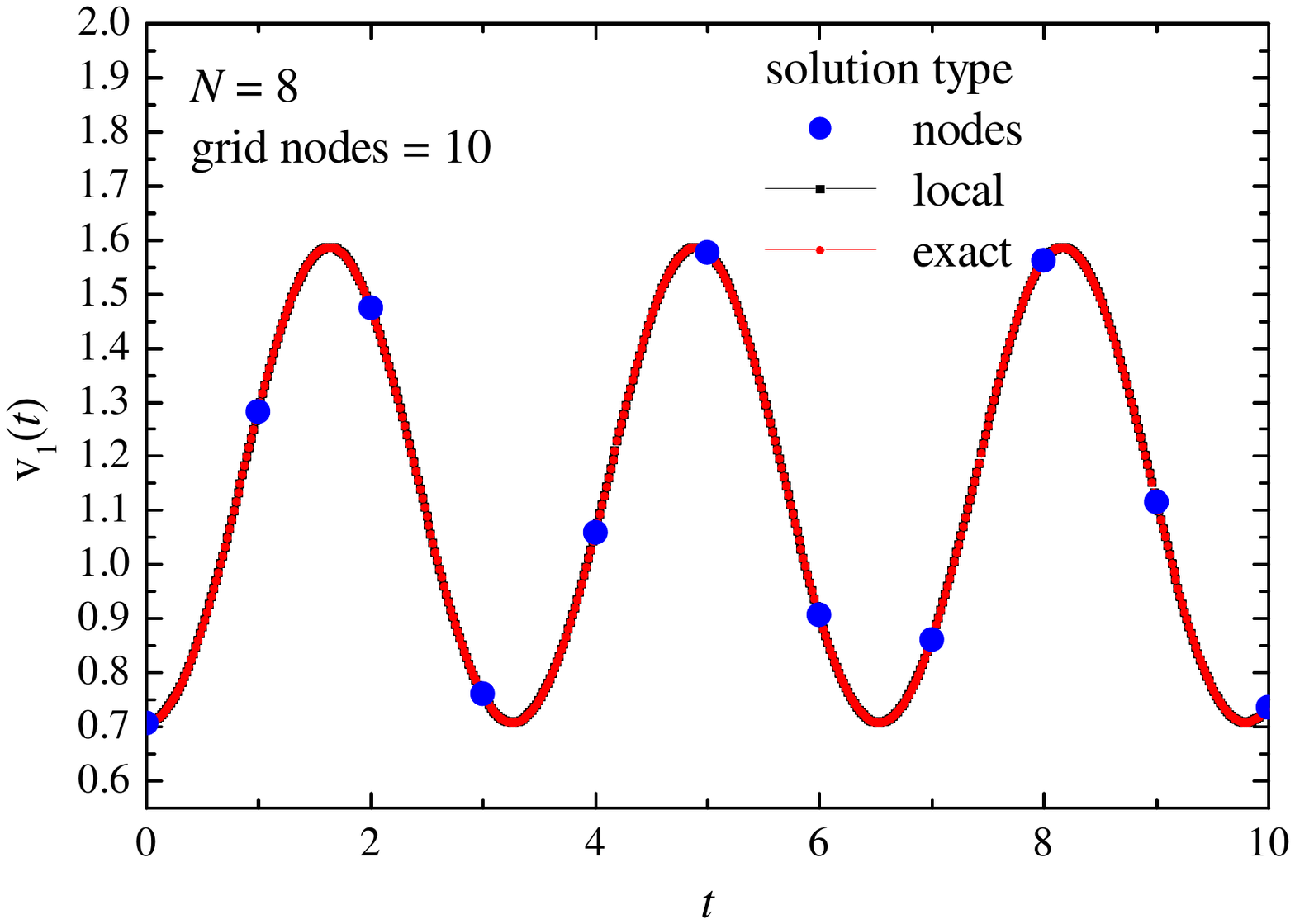}
\vspace{-8mm}\caption{\label{fig:pend_ind2_sol_uv:e2}}
\end{subfigure}
\begin{subfigure}{0.320\textwidth}
\includegraphics[width=\textwidth]{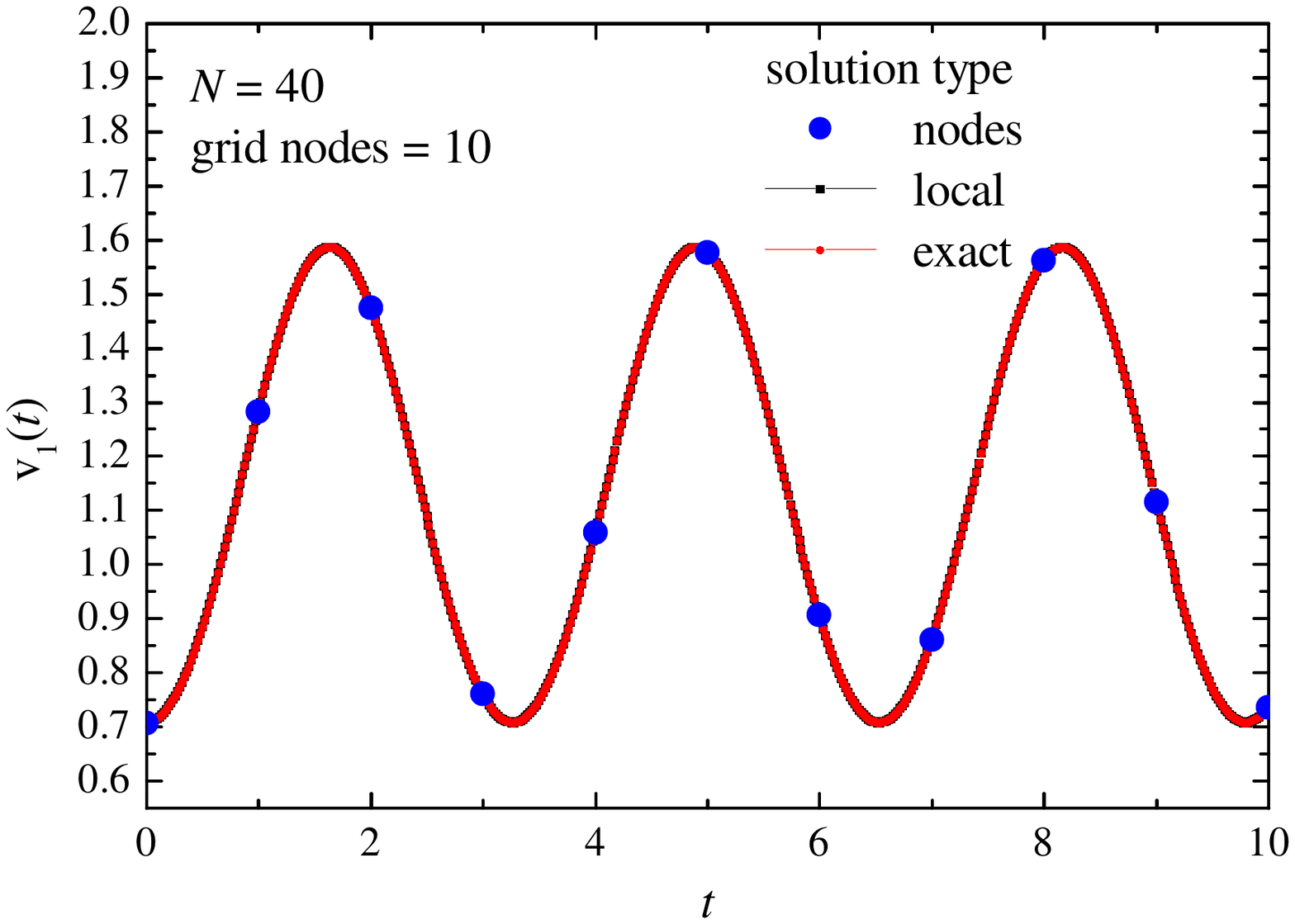}
\vspace{-8mm}\caption{\label{fig:pend_ind2_sol_uv:e3}}
\end{subfigure}\\
\caption{%
Numerical solution of the DAE system (\ref{eq:math_pend_dae_ind_3}) of index 2. Comparison of the solution at nodes $\mathbf{u}_{n}$, the local solution $\mathbf{u}_{L}(t)$ and the exact solution $\mathbf{u}^{\rm ex}(t)$ for components $u_{1}$ (\subref{fig:pend_ind2_sol_uv:a1}, \subref{fig:pend_ind2_sol_uv:a2}, \subref{fig:pend_ind2_sol_uv:a3}), $u_{2}$ (\subref{fig:pend_ind2_sol_uv:b1}, \subref{fig:pend_ind2_sol_uv:b2}, \subref{fig:pend_ind2_sol_uv:b3}), $u_{3}$ (\subref{fig:pend_ind2_sol_uv:c1}, \subref{fig:pend_ind2_sol_uv:c2}, \subref{fig:pend_ind2_sol_uv:c3}), $u_{4}$ (\subref{fig:pend_ind2_sol_uv:d1}, \subref{fig:pend_ind2_sol_uv:d2}, \subref{fig:pend_ind2_sol_uv:d3}) and $v_{1}$ (\subref{fig:pend_ind2_sol_uv:e1}, \subref{fig:pend_ind2_sol_uv:e2}, \subref{fig:pend_ind2_sol_uv:e3}), obtained using polynomials with degrees $N = 1$ (\subref{fig:pend_ind2_sol_uv:a1}, \subref{fig:pend_ind2_sol_uv:b1}, \subref{fig:pend_ind2_sol_uv:c1}, \subref{fig:pend_ind2_sol_uv:d1}, \subref{fig:pend_ind2_sol_uv:e1}), $N = 8$ (\subref{fig:pend_ind2_sol_uv:a2}, \subref{fig:pend_ind2_sol_uv:b2}, \subref{fig:pend_ind2_sol_uv:c2}, \subref{fig:pend_ind2_sol_uv:d2}, \subref{fig:pend_ind2_sol_uv:e2}) and $N = 40$ (\subref{fig:pend_ind2_sol_uv:a3}, \subref{fig:pend_ind2_sol_uv:b3}, \subref{fig:pend_ind2_sol_uv:c3}, \subref{fig:pend_ind2_sol_uv:d3}, \subref{fig:pend_ind2_sol_uv:e3}).
}
\label{fig:pend_ind2_sol_uv}
\end{figure} 

\begin{figure}[h!]
\captionsetup[subfigure]{%
	position=bottom,
	font+=smaller,
	textfont=normalfont,
	singlelinecheck=off,
	justification=raggedright
}
\centering
\begin{subfigure}{0.275\textwidth}
\includegraphics[width=\textwidth]{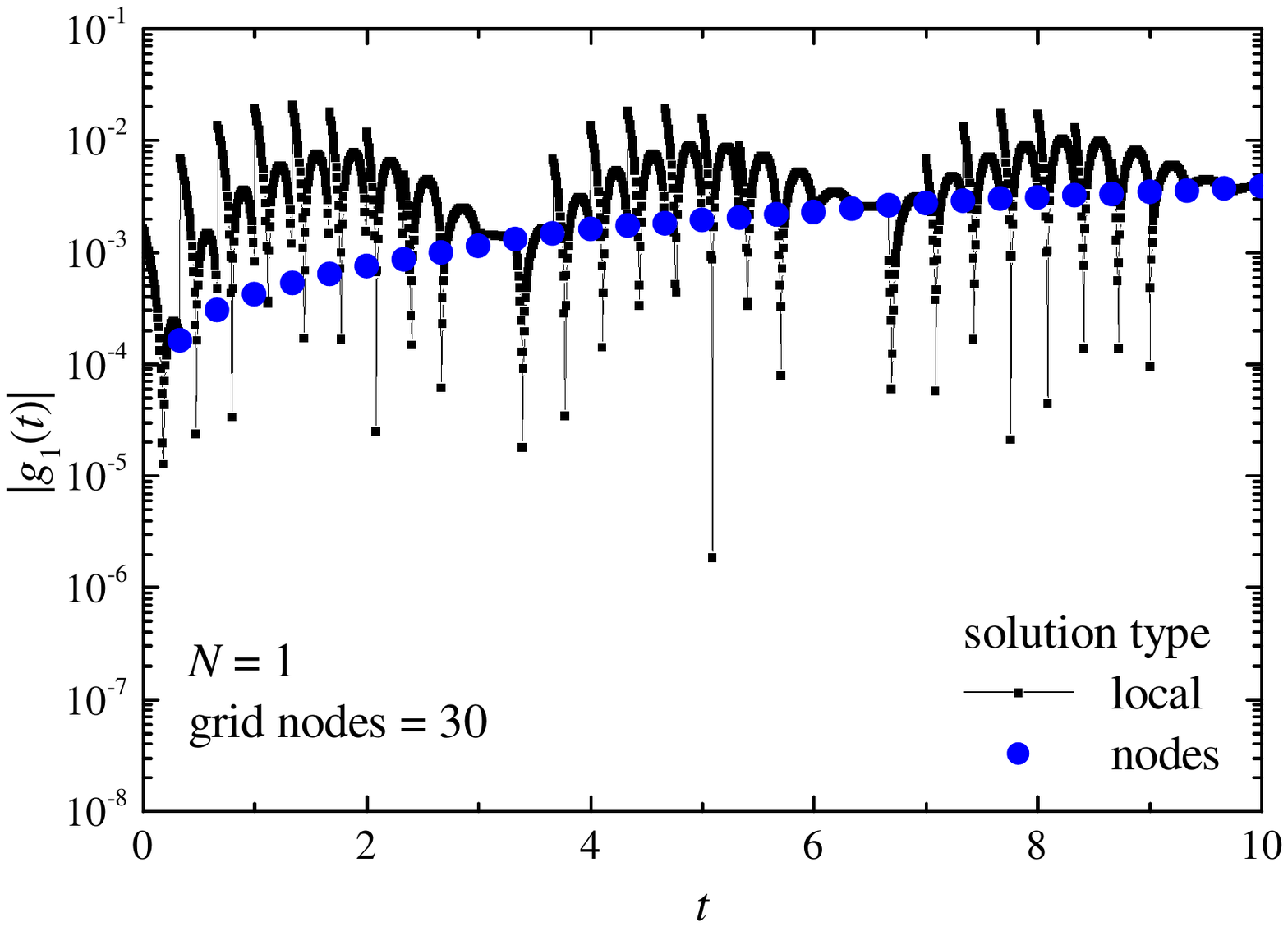}
\vspace{-8mm}\caption{\label{fig:pend_ind2_sol_g_eps:a1}}
\end{subfigure}\hspace{6mm}
\begin{subfigure}{0.275\textwidth}
\includegraphics[width=\textwidth]{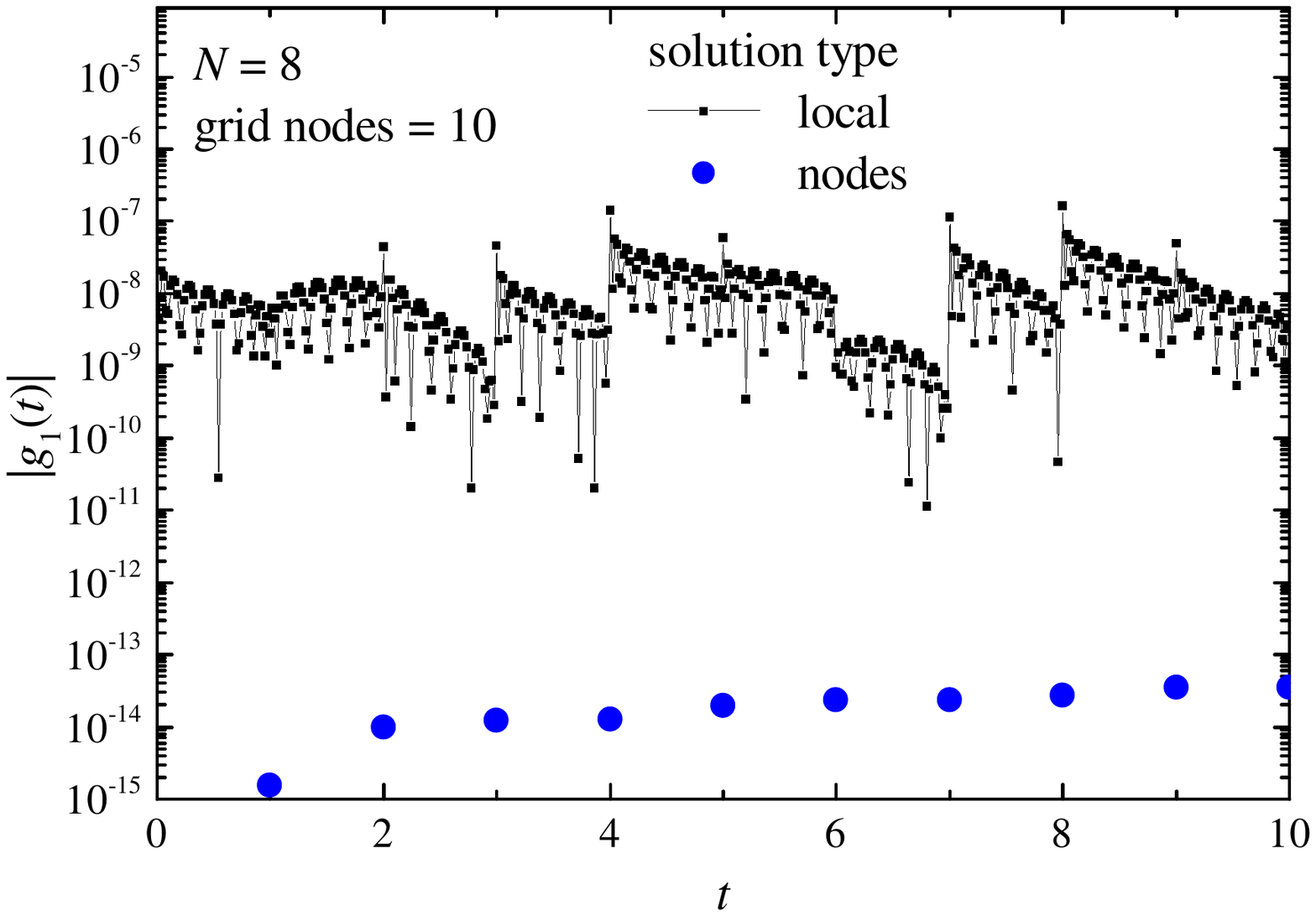}
\vspace{-8mm}\caption{\label{fig:pend_ind2_sol_g_eps:a2}}
\end{subfigure}\hspace{6mm}
\begin{subfigure}{0.275\textwidth}
\includegraphics[width=\textwidth]{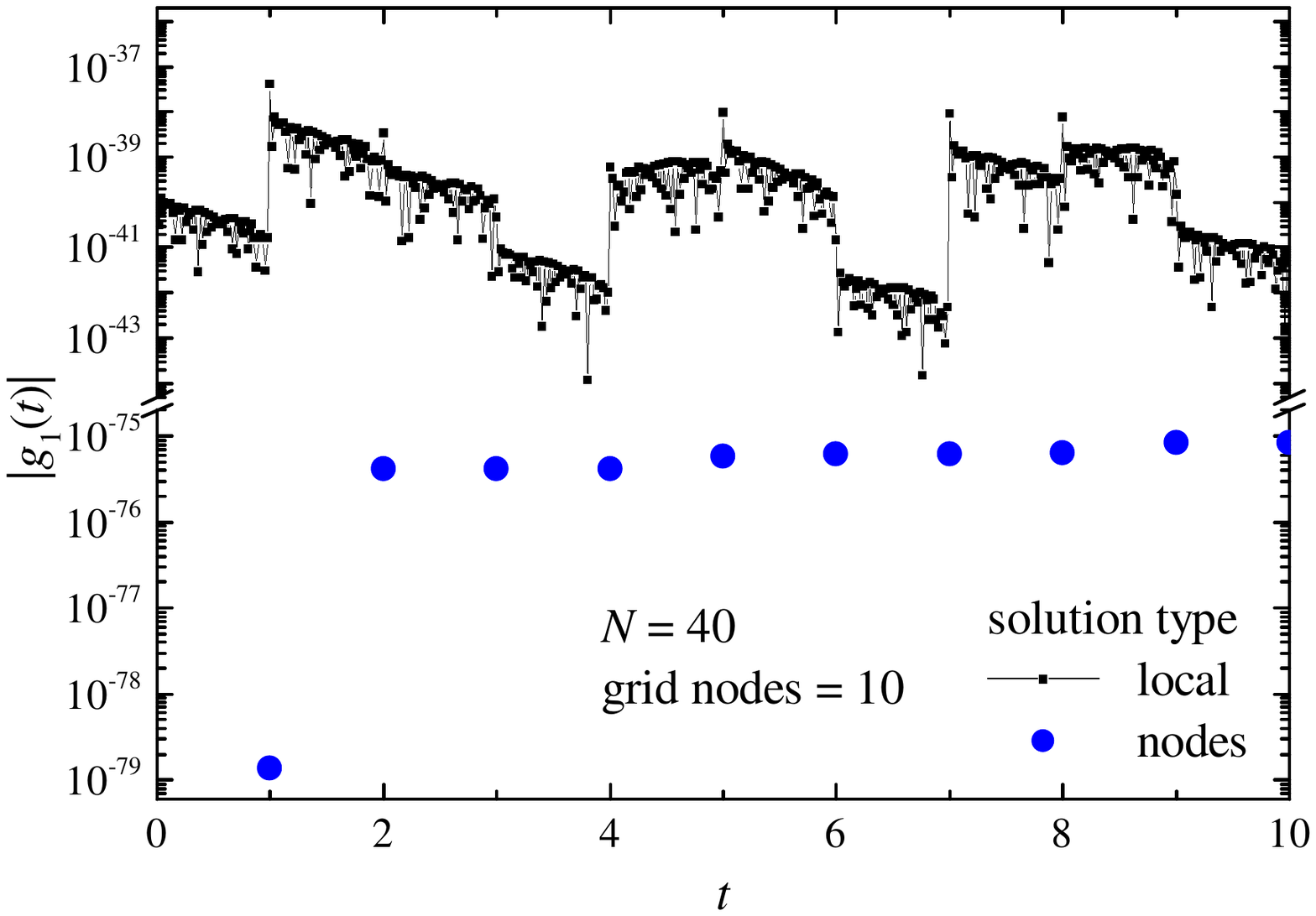}
\vspace{-8mm}\caption{\label{fig:pend_ind2_sol_g_eps:a3}}
\end{subfigure}\\[-2mm]
\begin{subfigure}{0.275\textwidth}
\includegraphics[width=\textwidth]{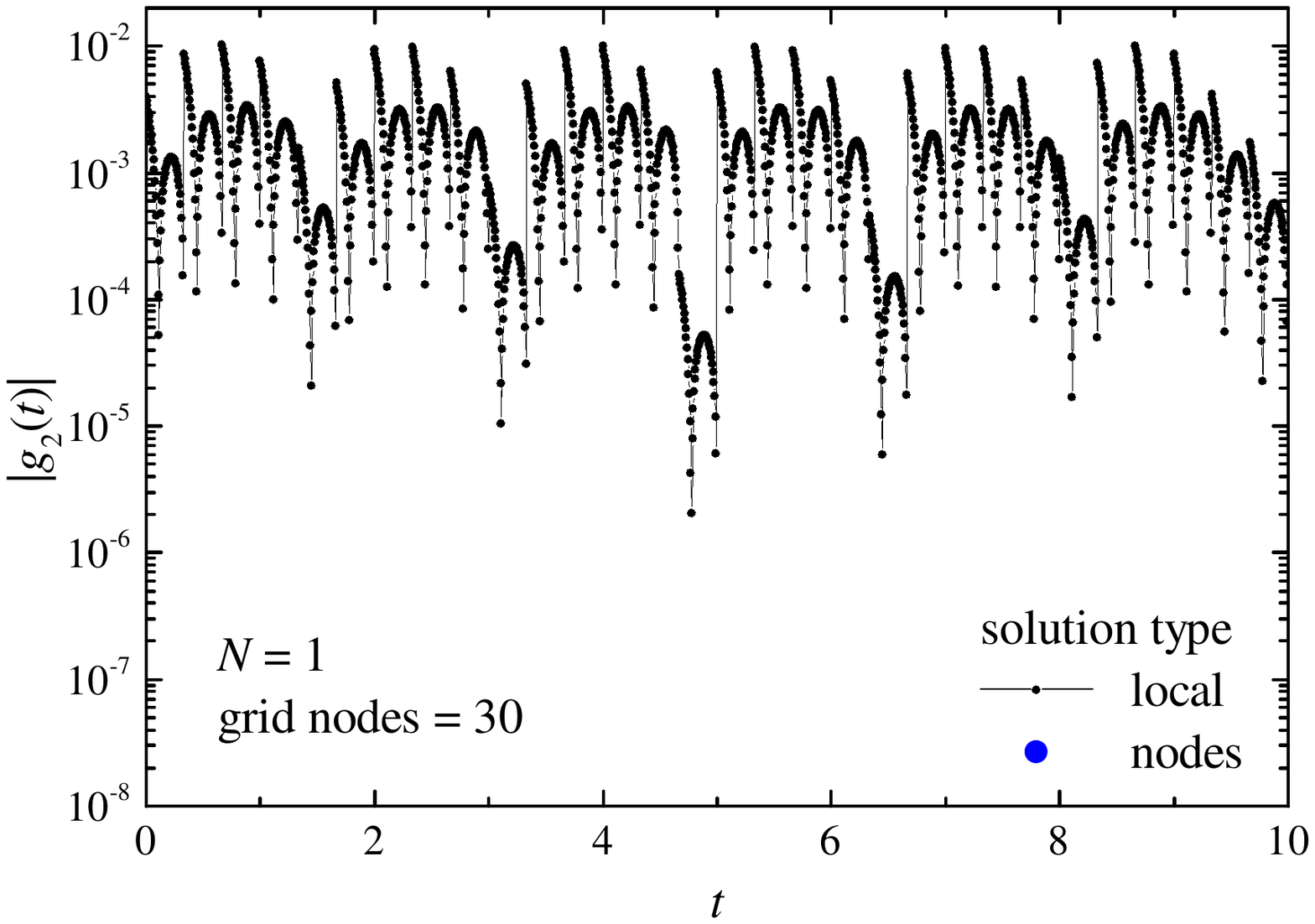}
\vspace{-8mm}\caption{\label{fig:pend_ind2_sol_g_eps:b1}}
\end{subfigure}\hspace{6mm}
\begin{subfigure}{0.275\textwidth}
\includegraphics[width=\textwidth]{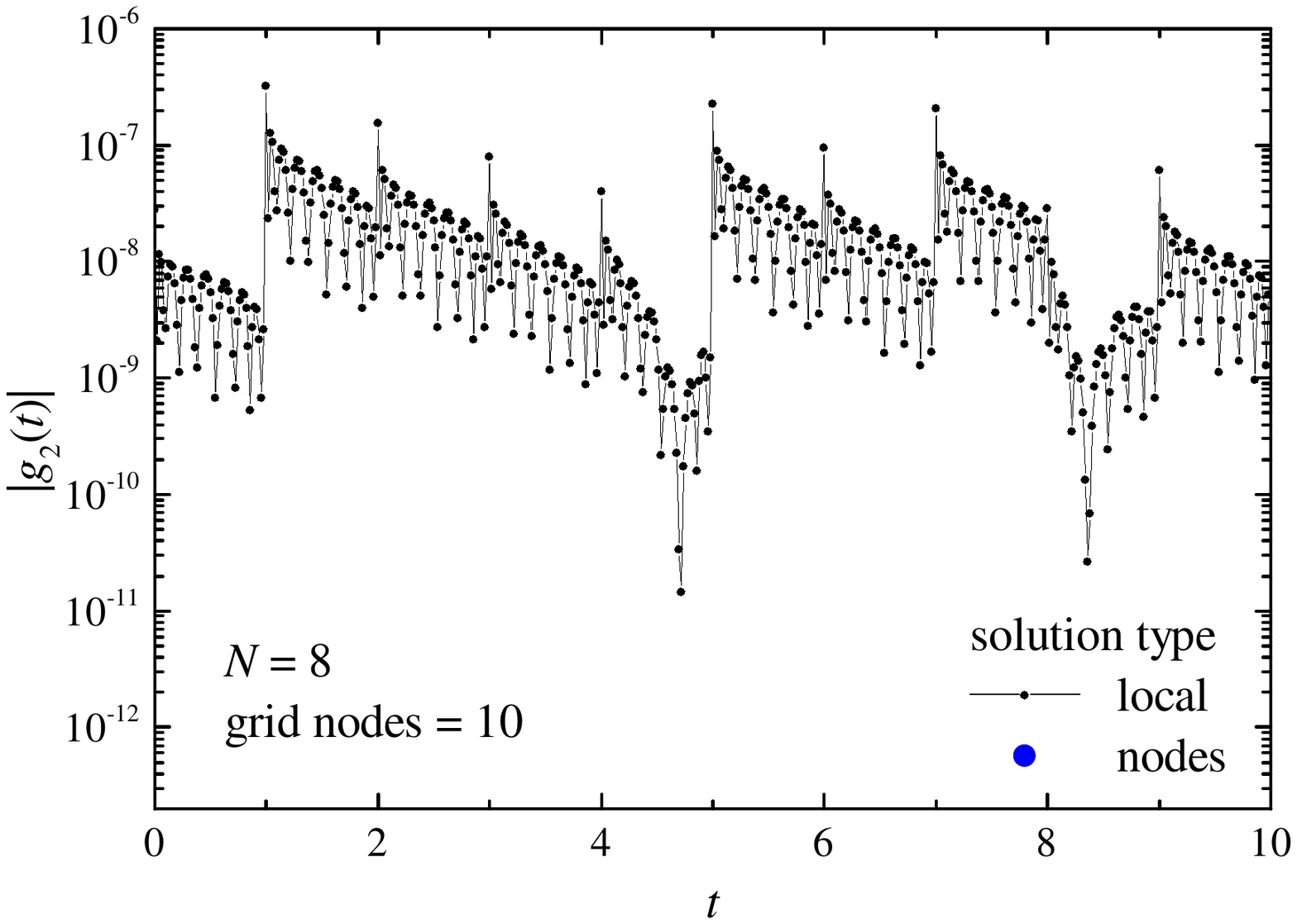}
\vspace{-8mm}\caption{\label{fig:pend_ind2_sol_g_eps:b2}}
\end{subfigure}\hspace{6mm}
\begin{subfigure}{0.275\textwidth}
\includegraphics[width=\textwidth]{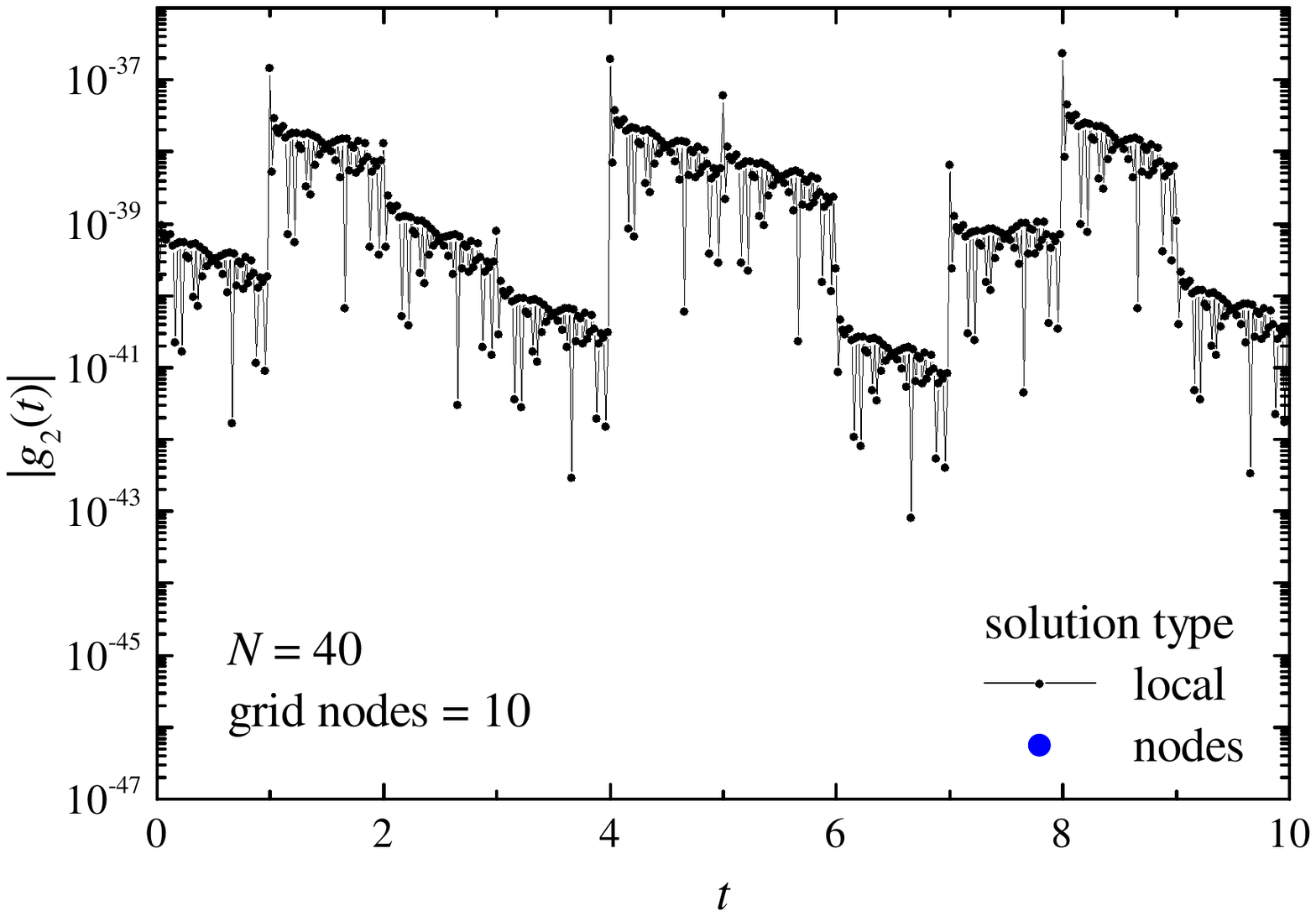}
\vspace{-8mm}\caption{\label{fig:pend_ind2_sol_g_eps:b3}}
\end{subfigure}\\[-2mm]
\begin{subfigure}{0.275\textwidth}
\includegraphics[width=\textwidth]{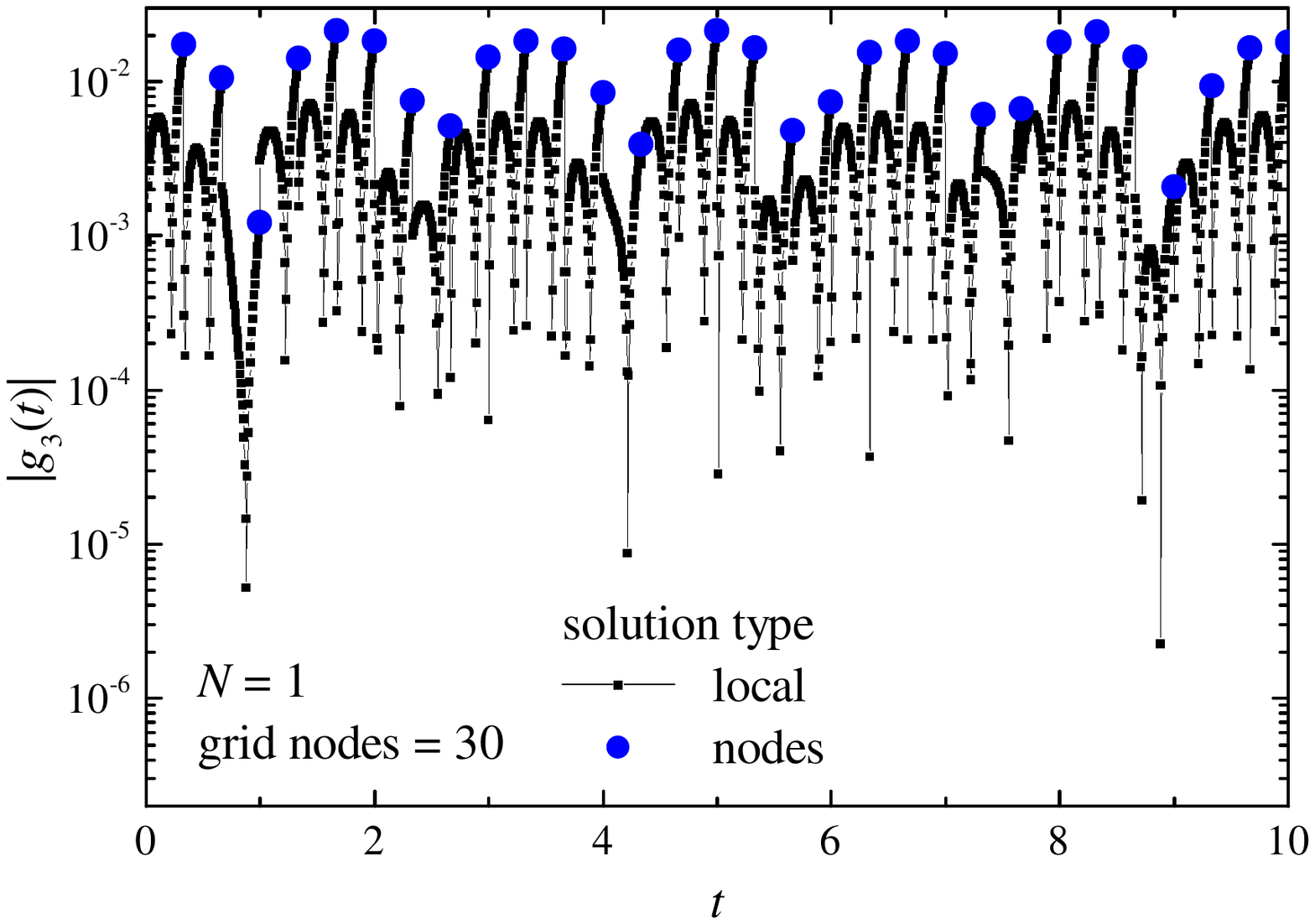}
\vspace{-8mm}\caption{\label{fig:pend_ind2_sol_g_eps:c1}}
\end{subfigure}\hspace{6mm}
\begin{subfigure}{0.275\textwidth}
\includegraphics[width=\textwidth]{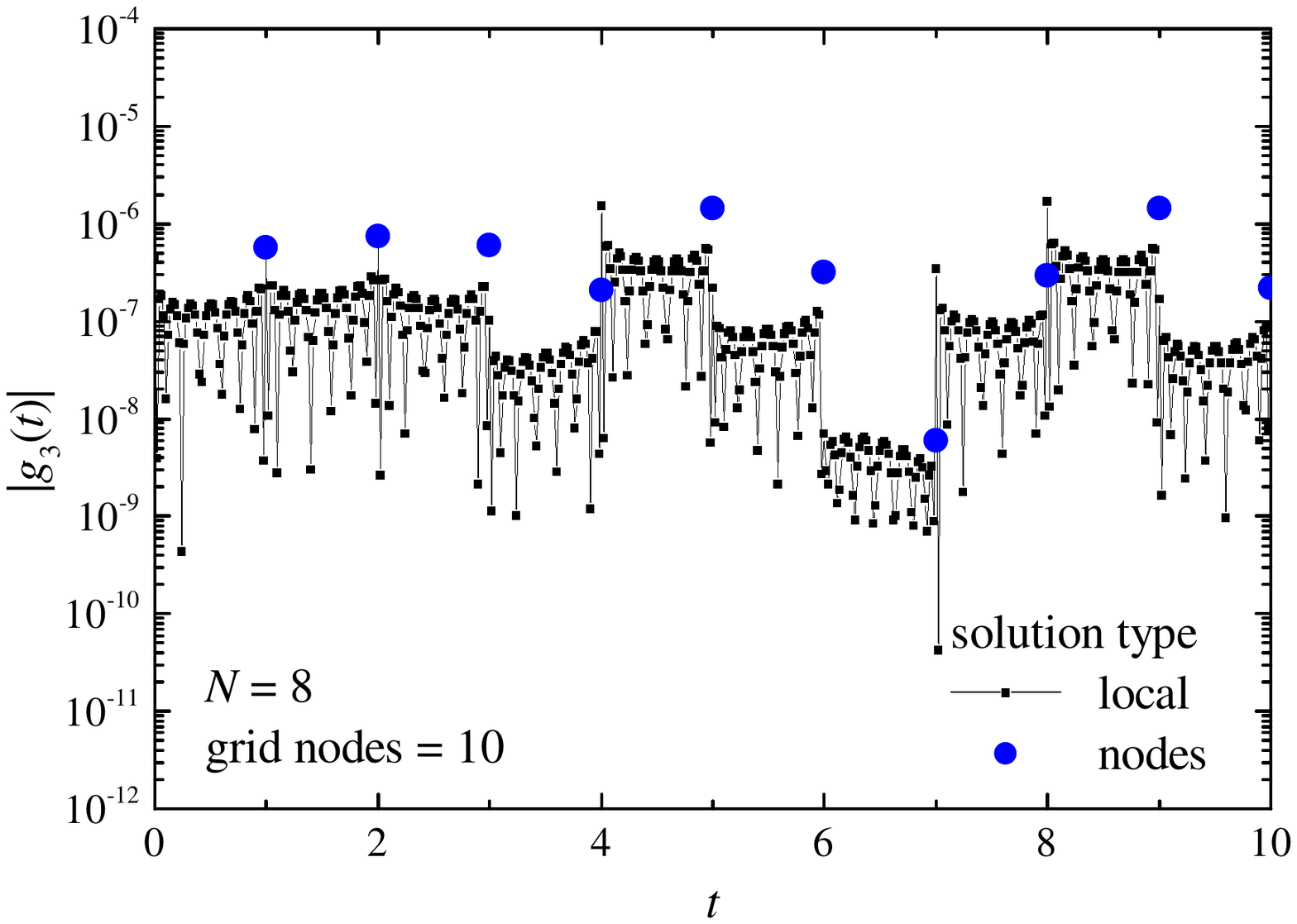}
\vspace{-8mm}\caption{\label{fig:pend_ind2_sol_g_eps:c2}}
\end{subfigure}\hspace{6mm}
\begin{subfigure}{0.275\textwidth}
\includegraphics[width=\textwidth]{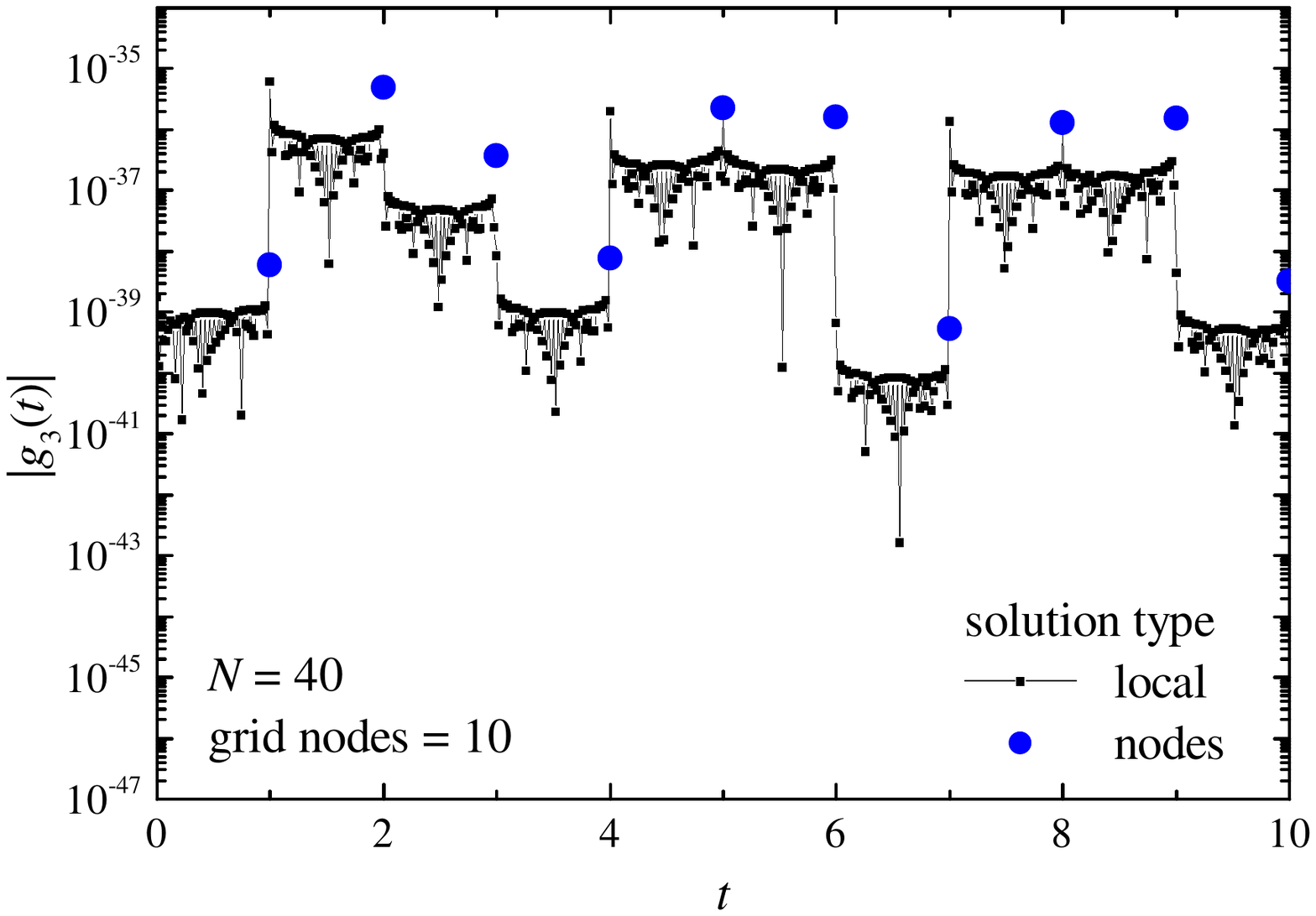}
\vspace{-8mm}\caption{\label{fig:pend_ind2_sol_g_eps:c3}}
\end{subfigure}\\[-2mm]
\begin{subfigure}{0.275\textwidth}
\includegraphics[width=\textwidth]{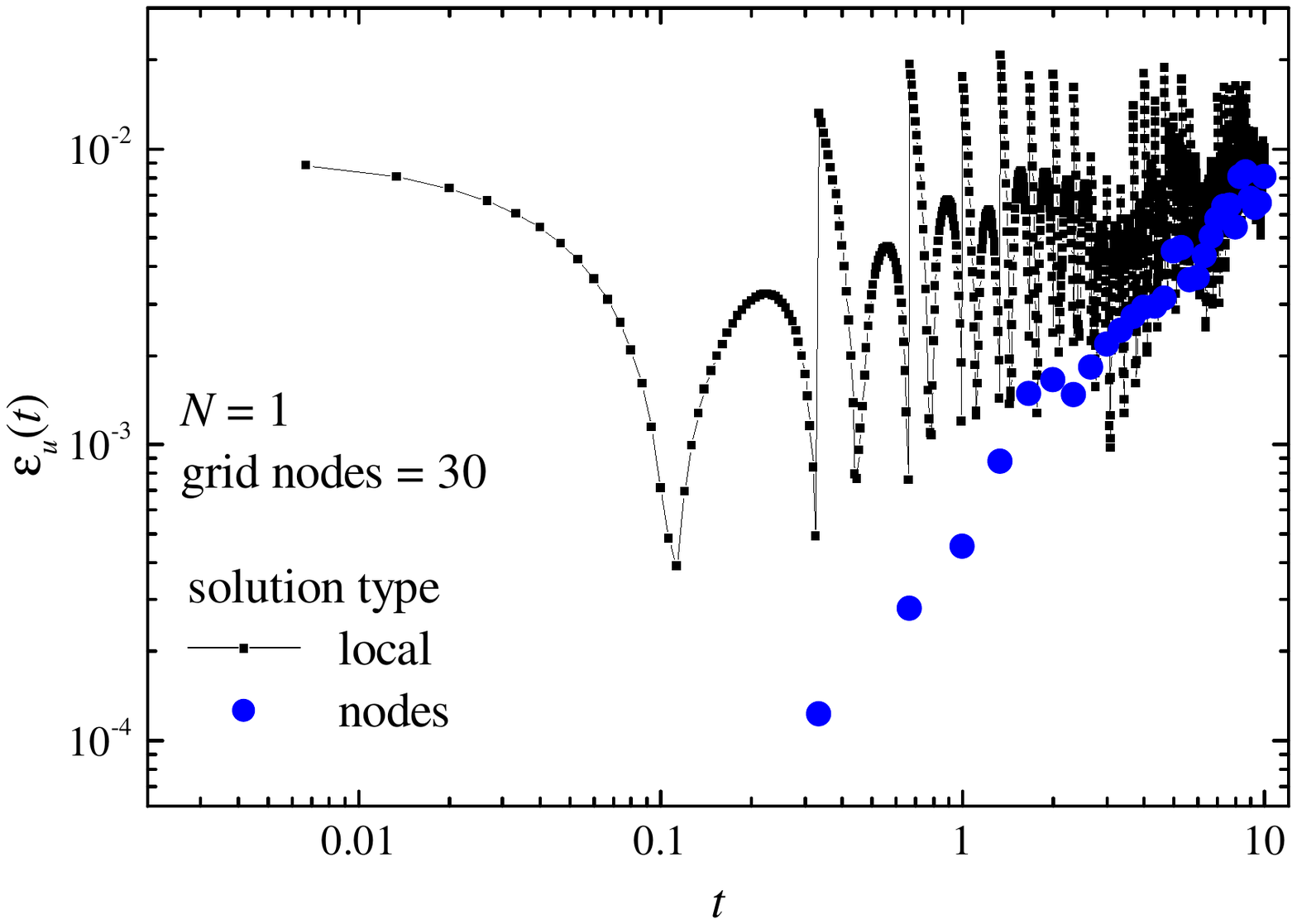}
\vspace{-8mm}\caption{\label{fig:pend_ind2_sol_g_eps:d1}}
\end{subfigure}\hspace{6mm}
\begin{subfigure}{0.275\textwidth}
\includegraphics[width=\textwidth]{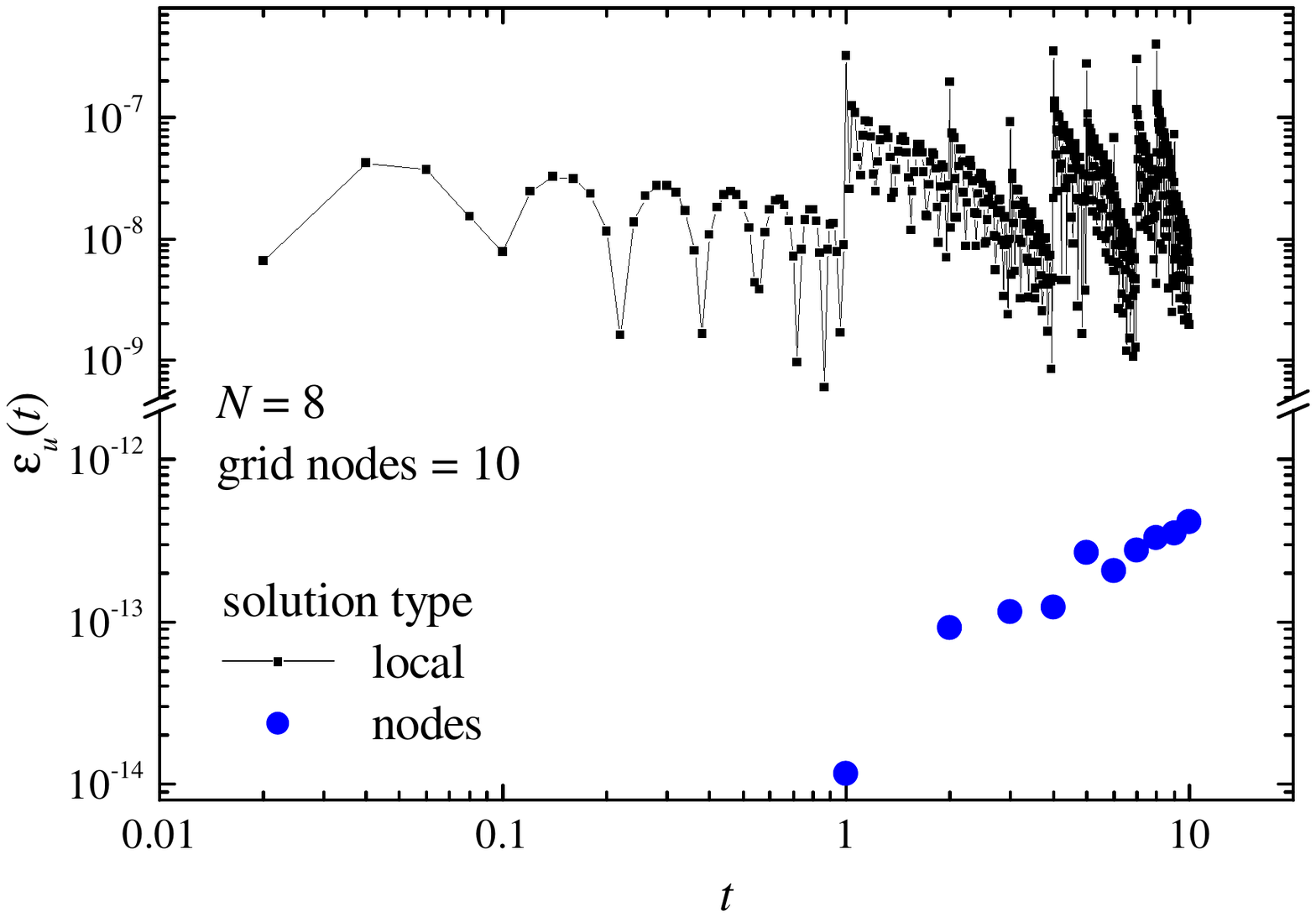}
\vspace{-8mm}\caption{\label{fig:pend_ind2_sol_g_eps:d2}}
\end{subfigure}\hspace{6mm}
\begin{subfigure}{0.275\textwidth}
\includegraphics[width=\textwidth]{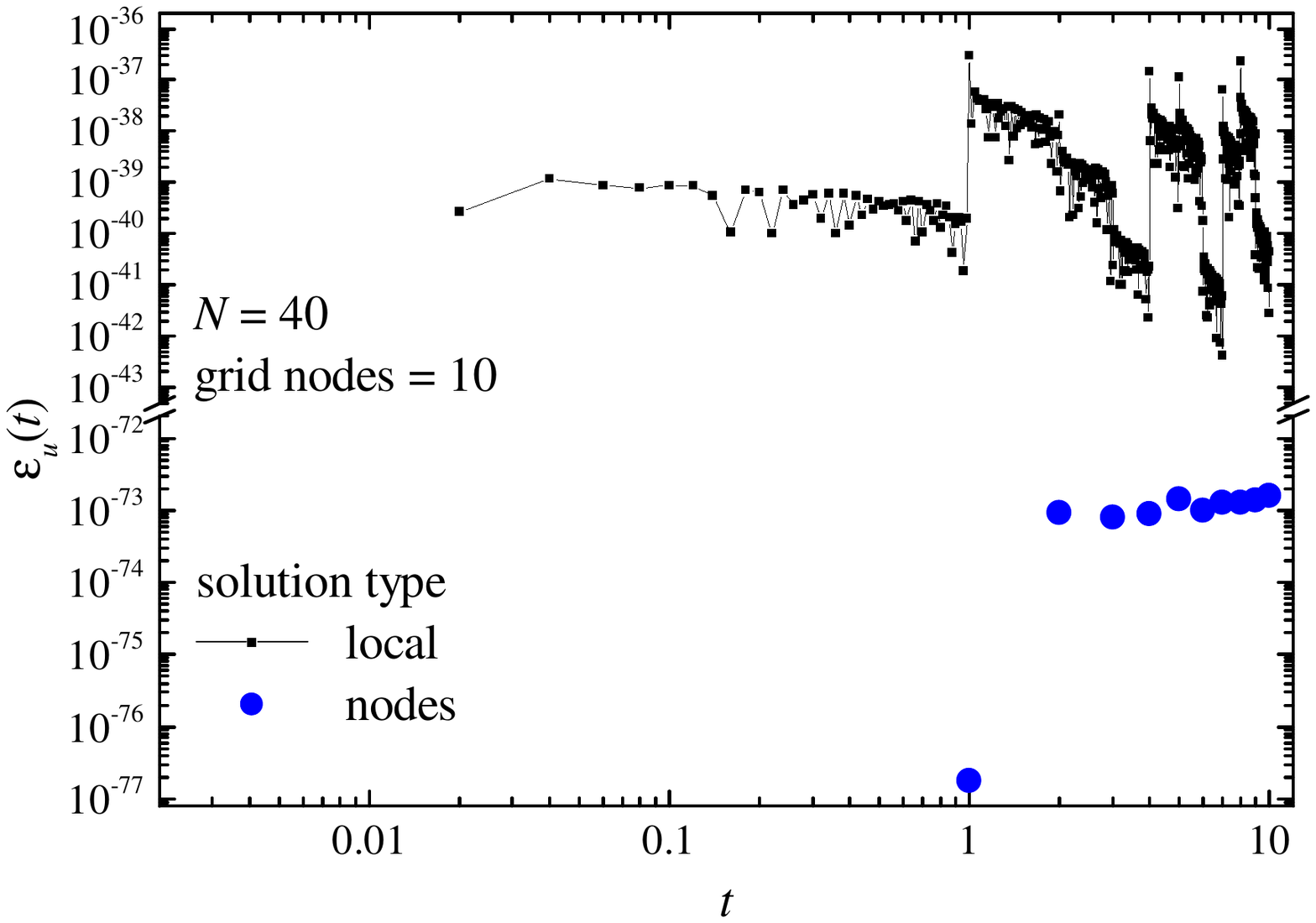}
\vspace{-8mm}\caption{\label{fig:pend_ind2_sol_g_eps:d3}}
\end{subfigure}\\[-2mm]
\begin{subfigure}{0.275\textwidth}
\includegraphics[width=\textwidth]{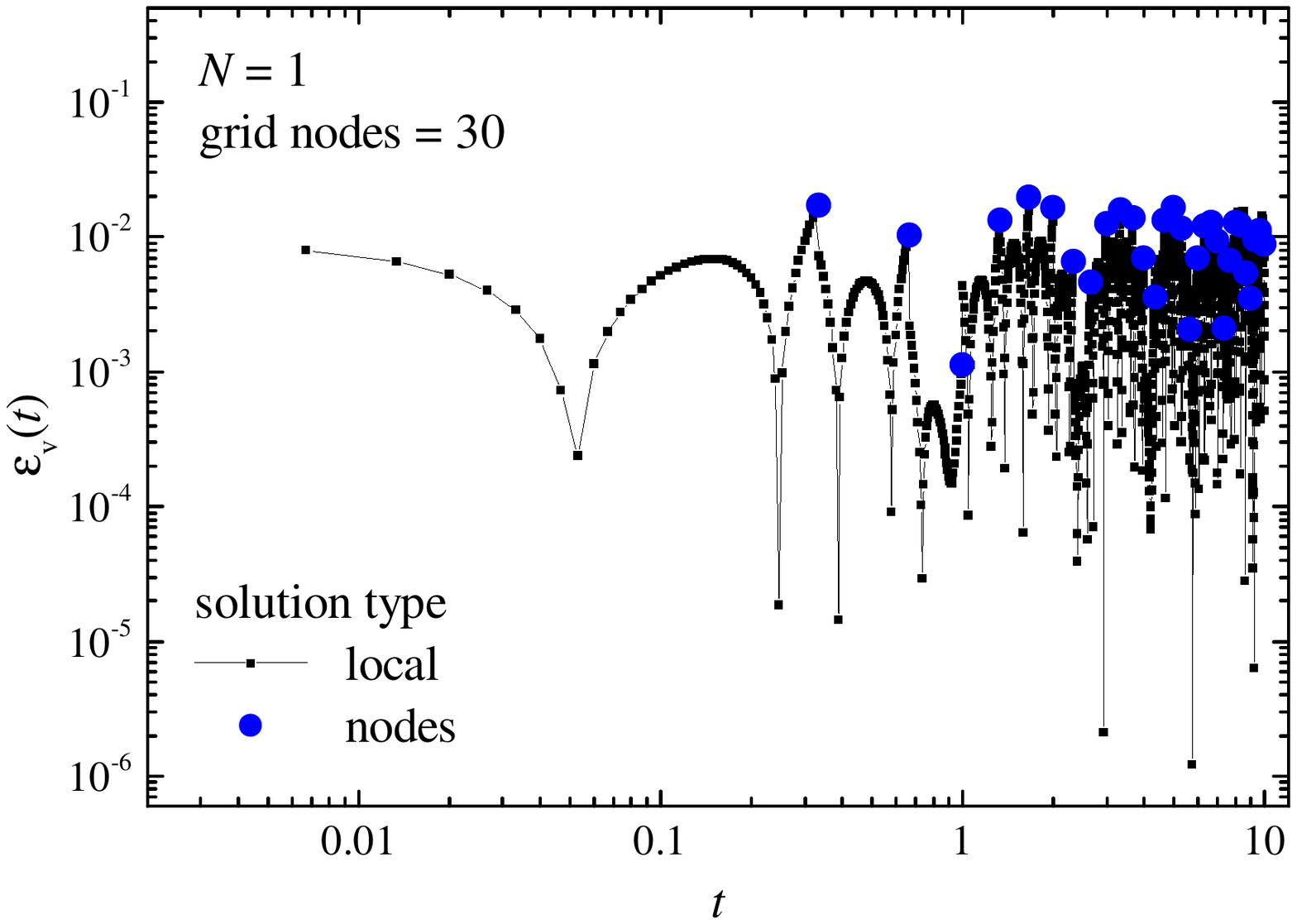}
\vspace{-8mm}\caption{\label{fig:pend_ind2_sol_g_eps:e1}}
\end{subfigure}\hspace{6mm}
\begin{subfigure}{0.275\textwidth}
\includegraphics[width=\textwidth]{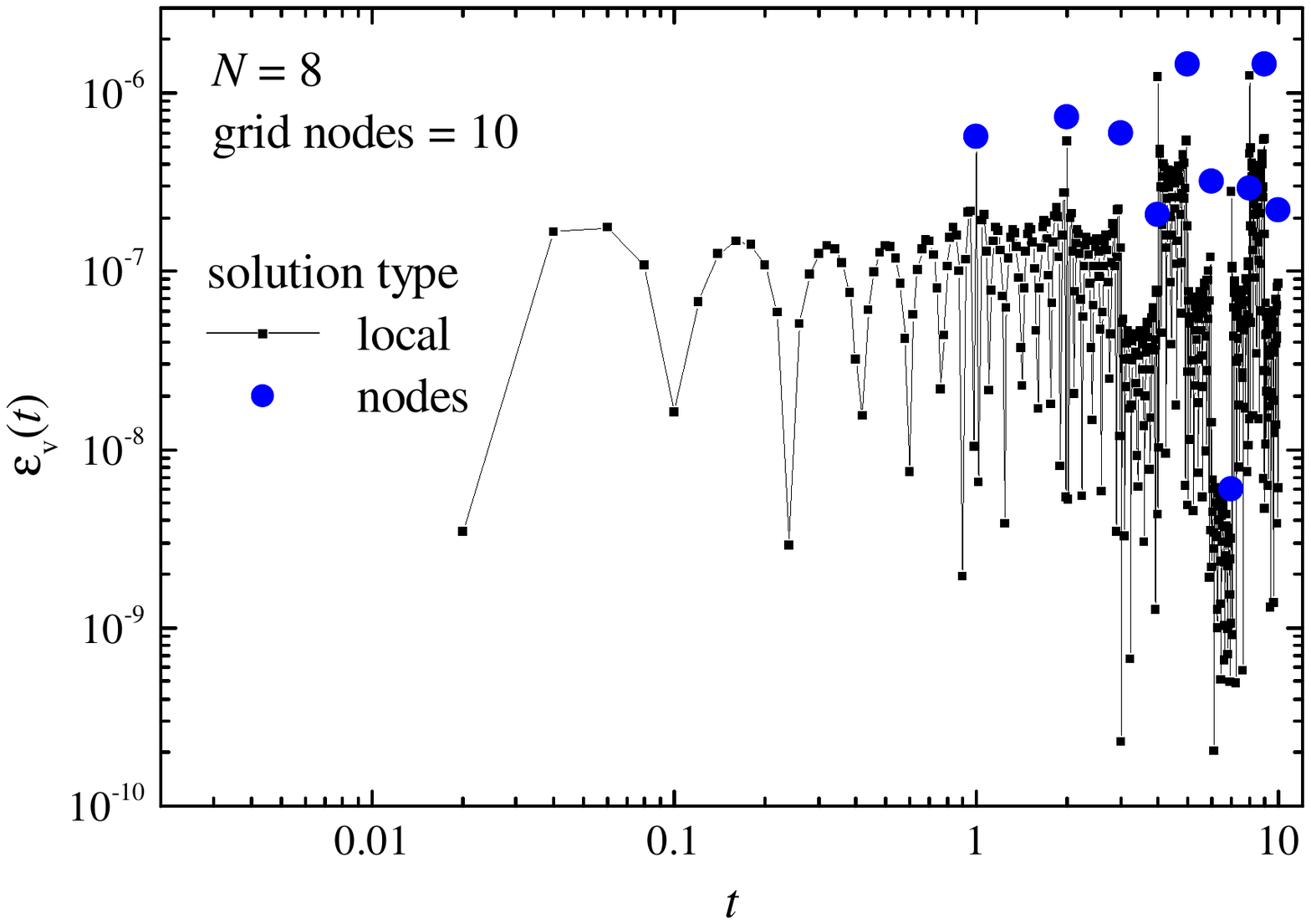}
\vspace{-8mm}\caption{\label{fig:pend_ind2_sol_g_eps:e2}}
\end{subfigure}\hspace{6mm}
\begin{subfigure}{0.275\textwidth}
\includegraphics[width=\textwidth]{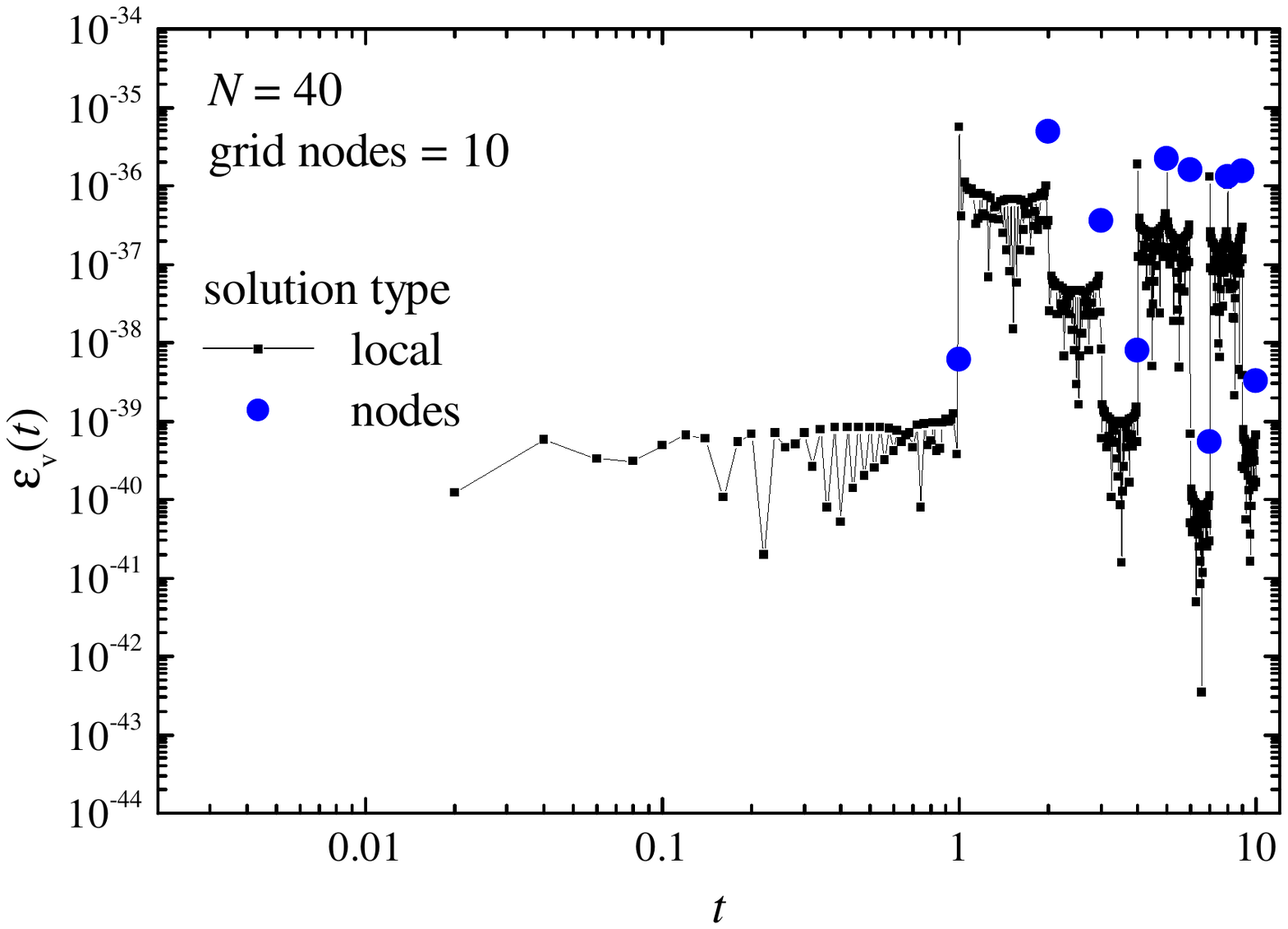}
\vspace{-8mm}\caption{\label{fig:pend_ind2_sol_g_eps:e3}}
\end{subfigure}\\[-2mm]
\begin{subfigure}{0.275\textwidth}
\includegraphics[width=\textwidth]{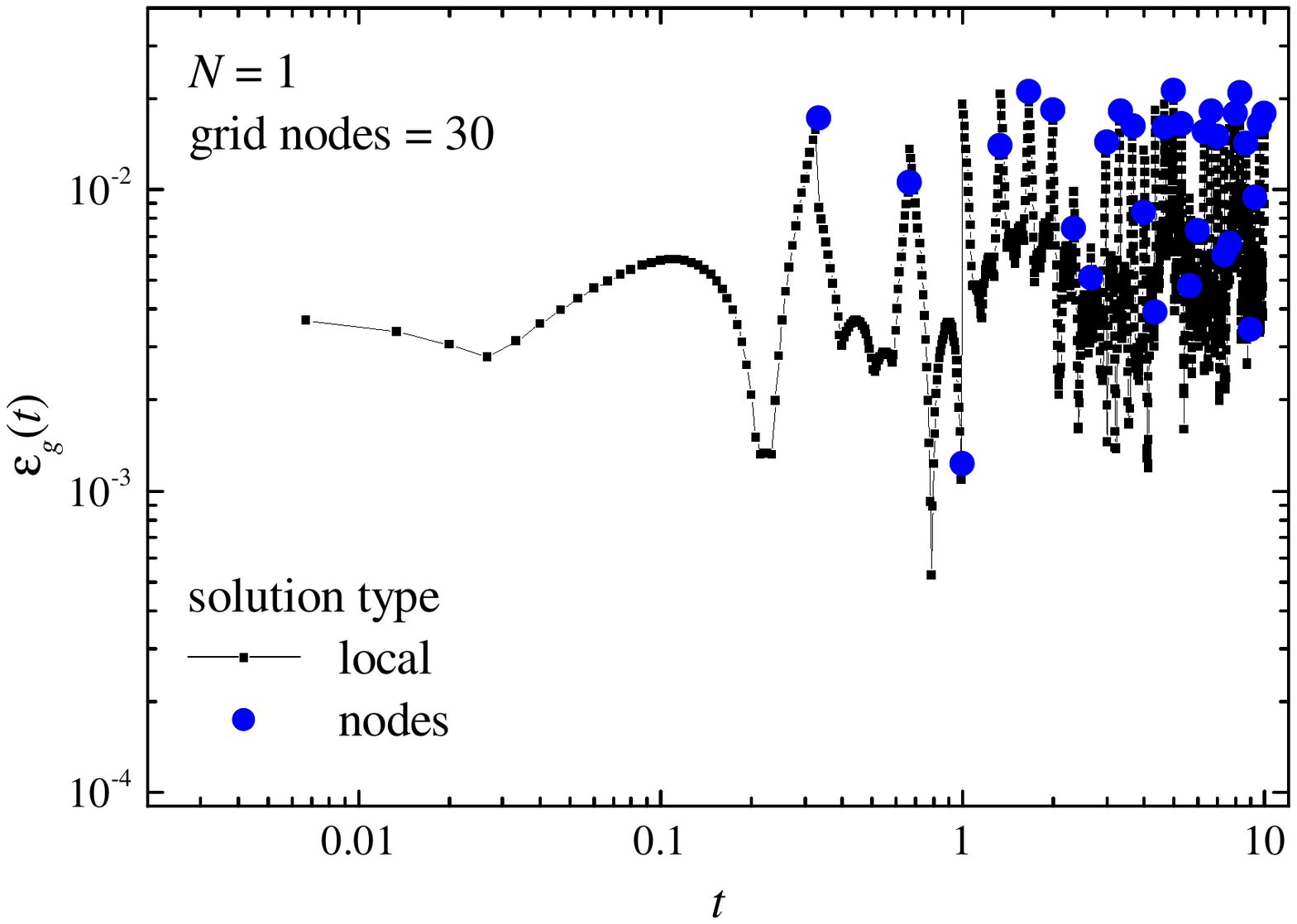}
\vspace{-8mm}\caption{\label{fig:pend_ind2_sol_g_eps:f1}}
\end{subfigure}\hspace{6mm}
\begin{subfigure}{0.275\textwidth}
\includegraphics[width=\textwidth]{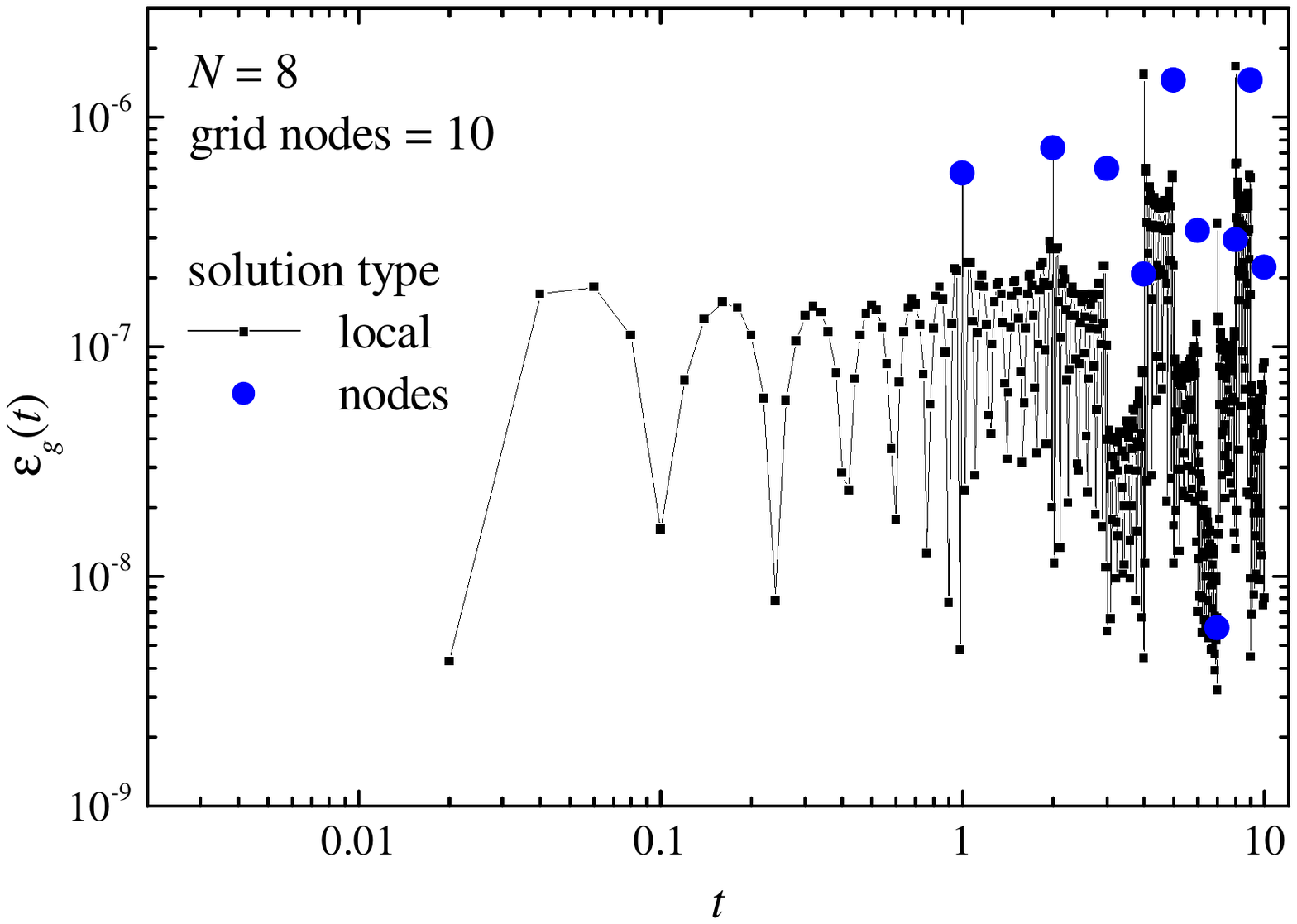}
\vspace{-8mm}\caption{\label{fig:pend_ind2_sol_g_eps:f2}}
\end{subfigure}\hspace{6mm}
\begin{subfigure}{0.275\textwidth}
\includegraphics[width=\textwidth]{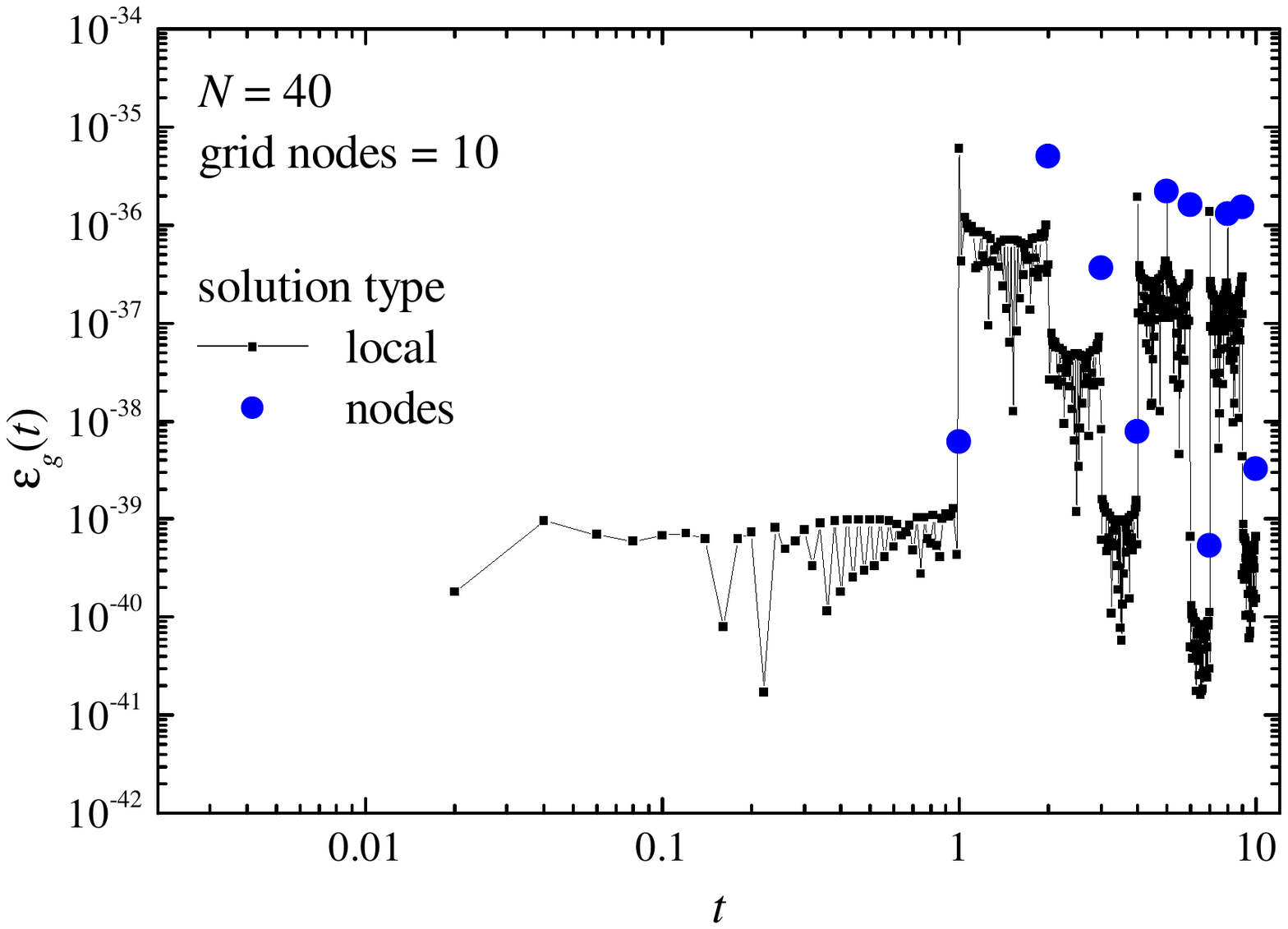}
\vspace{-8mm}\caption{\label{fig:pend_ind2_sol_g_eps:f3}}
\end{subfigure}\\[-2mm]
\caption{%
Numerical solution of the DAE system (\ref{eq:math_pend_dae_ind_3}) of index 2. Comparison of quantitative satisfiability of the conditions $g_{1} = 0$ (\subref{fig:pend_ind2_sol_g_eps:a1}, \subref{fig:pend_ind2_sol_g_eps:a2}, \subref{fig:pend_ind2_sol_g_eps:a3}), $g_{2} = 0$ (\subref{fig:pend_ind2_sol_g_eps:b1}, \subref{fig:pend_ind2_sol_g_eps:b2}, \subref{fig:pend_ind2_sol_g_eps:b3}) and $g_{3} = 0$ (\subref{fig:pend_ind2_sol_g_eps:c1}, \subref{fig:pend_ind2_sol_g_eps:c2}, \subref{fig:pend_ind2_sol_g_eps:c3}), the errors $\varepsilon_{u}(t)$ (\subref{fig:pend_ind2_sol_g_eps:d1}, \subref{fig:pend_ind2_sol_g_eps:d2}, \subref{fig:pend_ind2_sol_g_eps:d3}), $\varepsilon_{v}(t)$ (\subref{fig:pend_ind2_sol_g_eps:e1}, \subref{fig:pend_ind2_sol_g_eps:e2}, \subref{fig:pend_ind2_sol_g_eps:e3}), $\varepsilon_{g}(t)$ (\subref{fig:pend_ind2_sol_g_eps:f1}, \subref{fig:pend_ind2_sol_g_eps:f2}, \subref{fig:pend_ind2_sol_g_eps:f3}), for numerical solution obtained using polynomials with degrees $N = 1$ (\subref{fig:pend_ind2_sol_g_eps:a1}, \subref{fig:pend_ind2_sol_g_eps:b1}, \subref{fig:pend_ind2_sol_g_eps:c1}, \subref{fig:pend_ind2_sol_g_eps:d1}, \subref{fig:pend_ind2_sol_g_eps:e1}, \subref{fig:pend_ind2_sol_g_eps:f1}), $N = 8$ (\subref{fig:pend_ind2_sol_g_eps:a2}, \subref{fig:pend_ind2_sol_g_eps:b2}, \subref{fig:pend_ind2_sol_g_eps:c2}, \subref{fig:pend_ind2_sol_g_eps:d2}, \subref{fig:pend_ind2_sol_g_eps:e2}, \subref{fig:pend_ind2_sol_g_eps:f2}) and $N = 40$ (\subref{fig:pend_ind2_sol_g_eps:a3}, \subref{fig:pend_ind2_sol_g_eps:b3}, \subref{fig:pend_ind2_sol_g_eps:c3}, \subref{fig:pend_ind2_sol_g_eps:d3}, \subref{fig:pend_ind2_sol_g_eps:e3}, \subref{fig:pend_ind2_sol_g_eps:f3}).
}
\label{fig:pend_ind2_sol_g_eps}
\end{figure} 

\begin{figure}[h!]
\captionsetup[subfigure]{%
	position=bottom,
	font+=smaller,
	textfont=normalfont,
	singlelinecheck=off,
	justification=raggedright
}
\centering
\begin{subfigure}{0.275\textwidth}
\includegraphics[width=\textwidth]{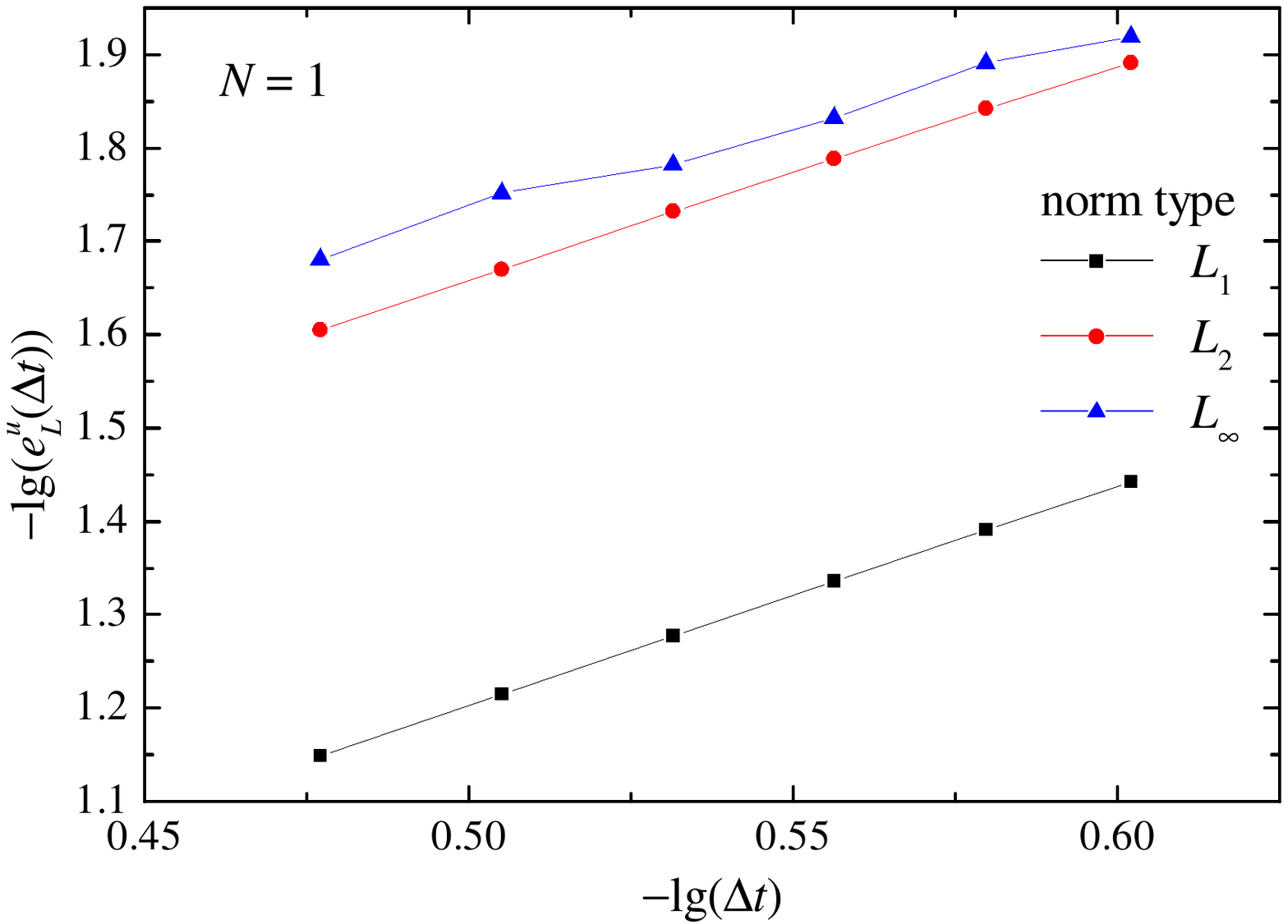}
\vspace{-8mm}\caption{\label{fig:pend_ind2_errors:a1}}
\end{subfigure}\hspace{6mm}
\begin{subfigure}{0.275\textwidth}
\includegraphics[width=\textwidth]{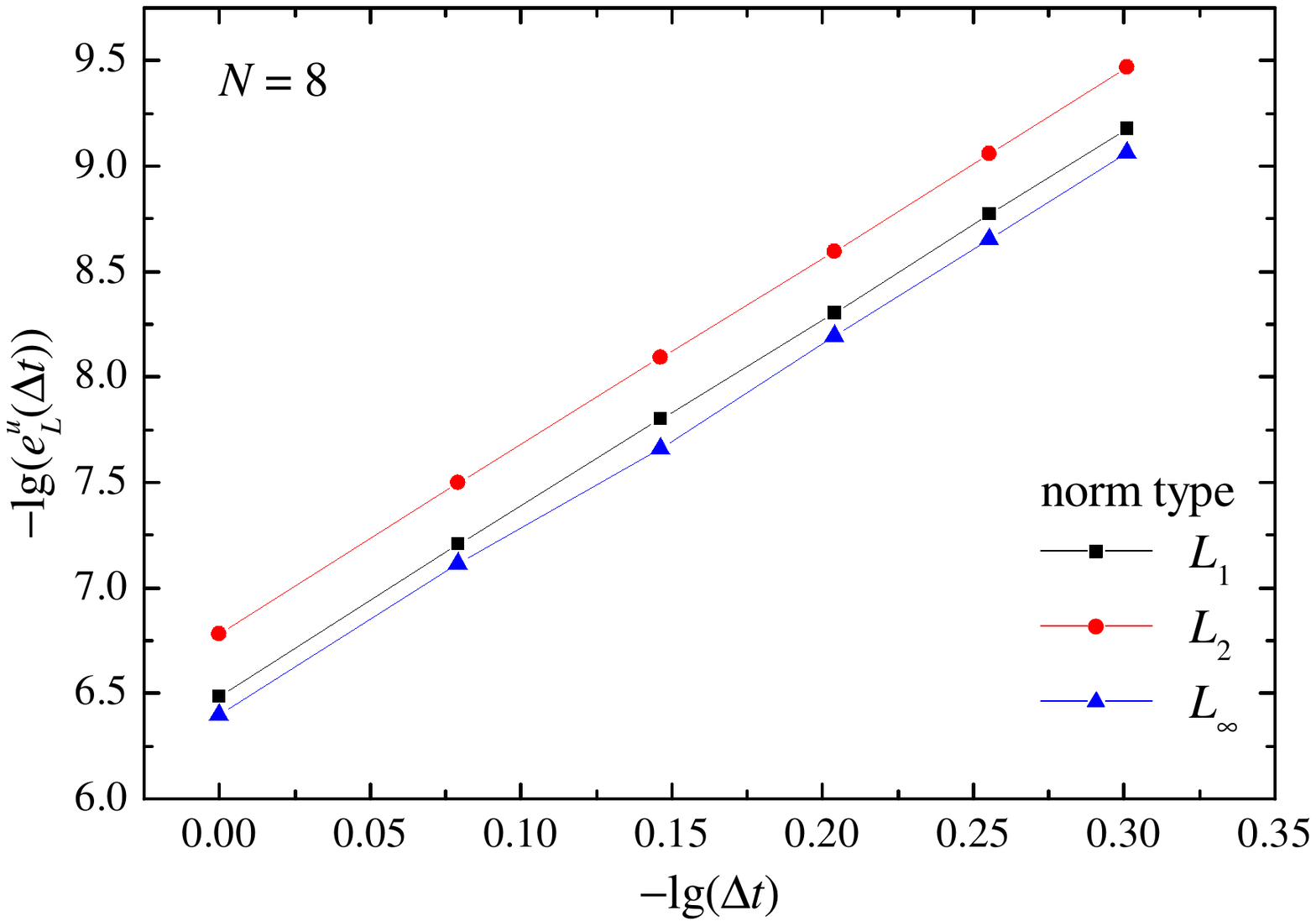}
\vspace{-8mm}\caption{\label{fig:pend_ind2_errors:a2}}
\end{subfigure}\hspace{6mm}
\begin{subfigure}{0.275\textwidth}
\includegraphics[width=\textwidth]{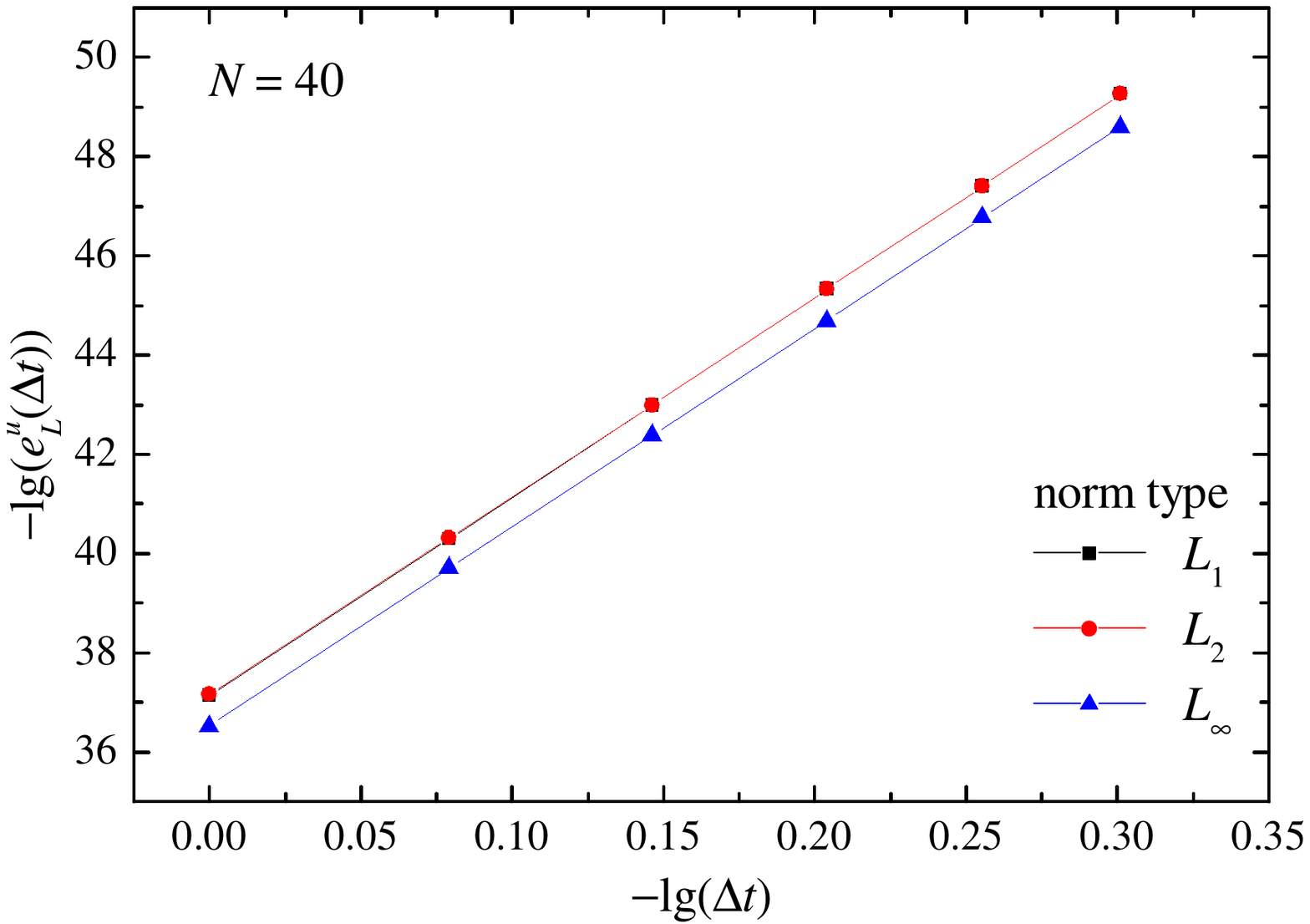}
\vspace{-8mm}\caption{\label{fig:pend_ind2_errors:a3}}
\end{subfigure}\\[-2mm]
\begin{subfigure}{0.275\textwidth}
\includegraphics[width=\textwidth]{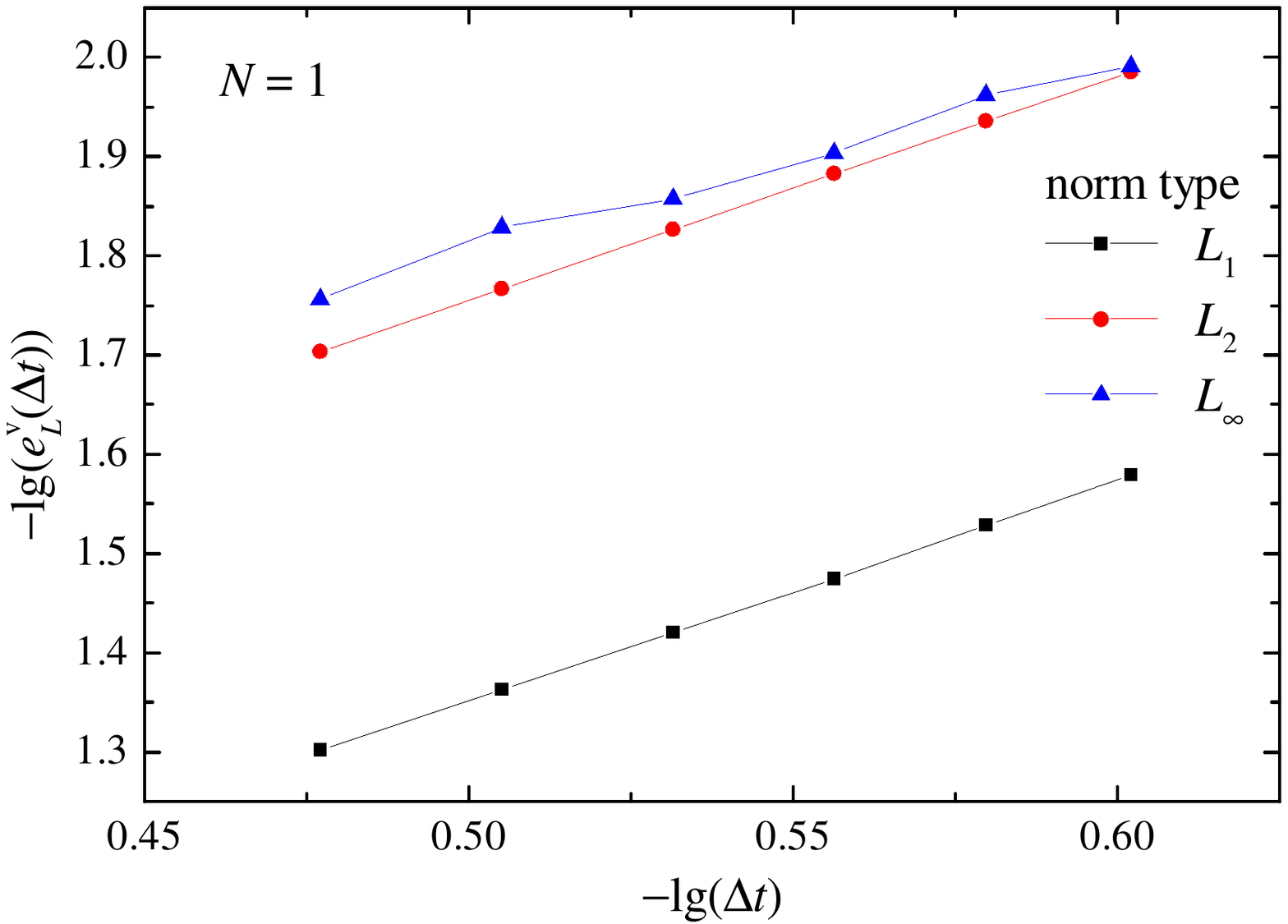}
\vspace{-8mm}\caption{\label{fig:pend_ind2_errors:b1}}
\end{subfigure}\hspace{6mm}
\begin{subfigure}{0.275\textwidth}
\includegraphics[width=\textwidth]{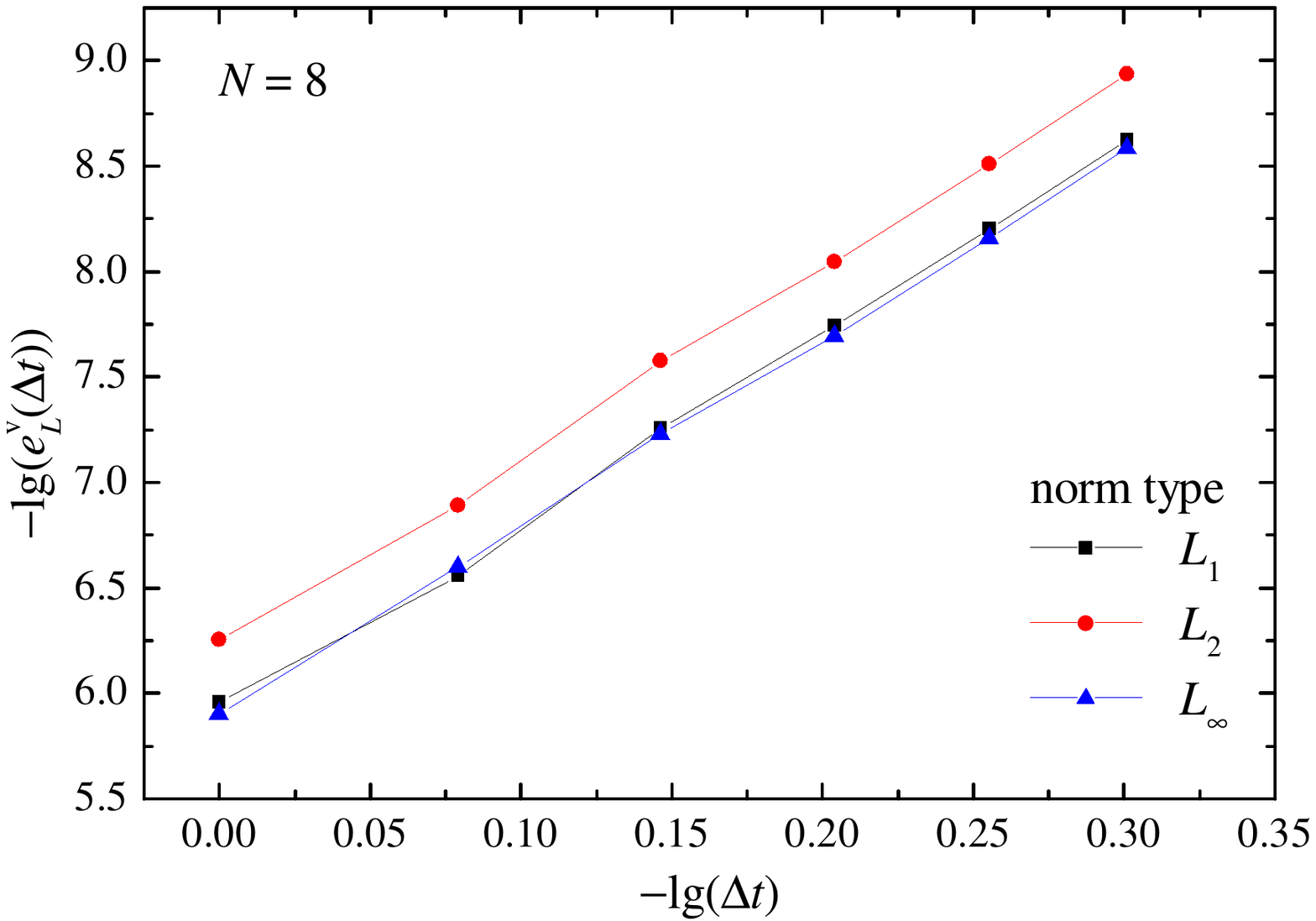}
\vspace{-8mm}\caption{\label{fig:pend_ind2_errors:b2}}
\end{subfigure}\hspace{6mm}
\begin{subfigure}{0.275\textwidth}
\includegraphics[width=\textwidth]{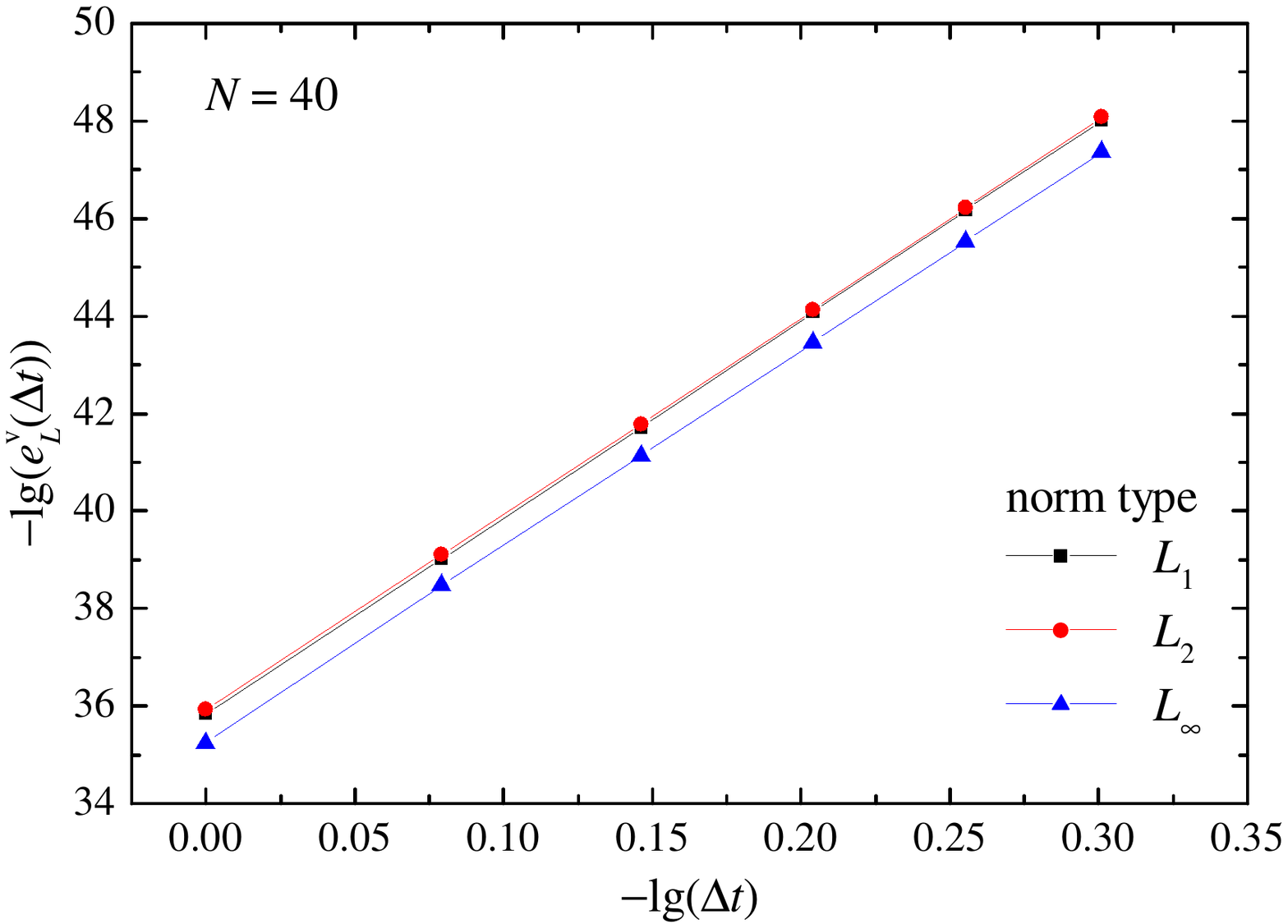}
\vspace{-8mm}\caption{\label{fig:pend_ind2_errors:b3}}
\end{subfigure}\\[-2mm]
\begin{subfigure}{0.275\textwidth}
\includegraphics[width=\textwidth]{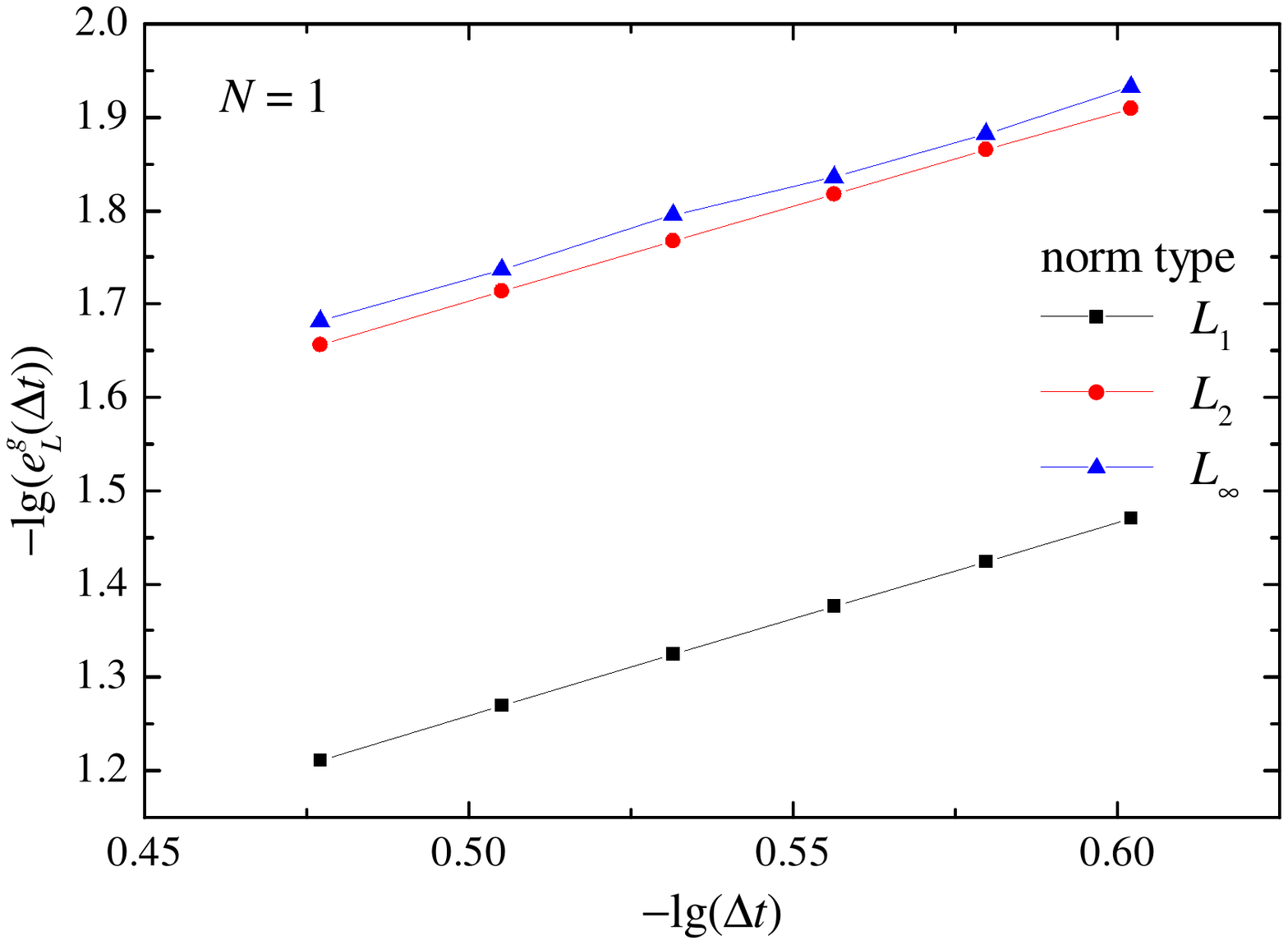}
\vspace{-8mm}\caption{\label{fig:pend_ind2_errors:c1}}
\end{subfigure}\hspace{6mm}
\begin{subfigure}{0.275\textwidth}
\includegraphics[width=\textwidth]{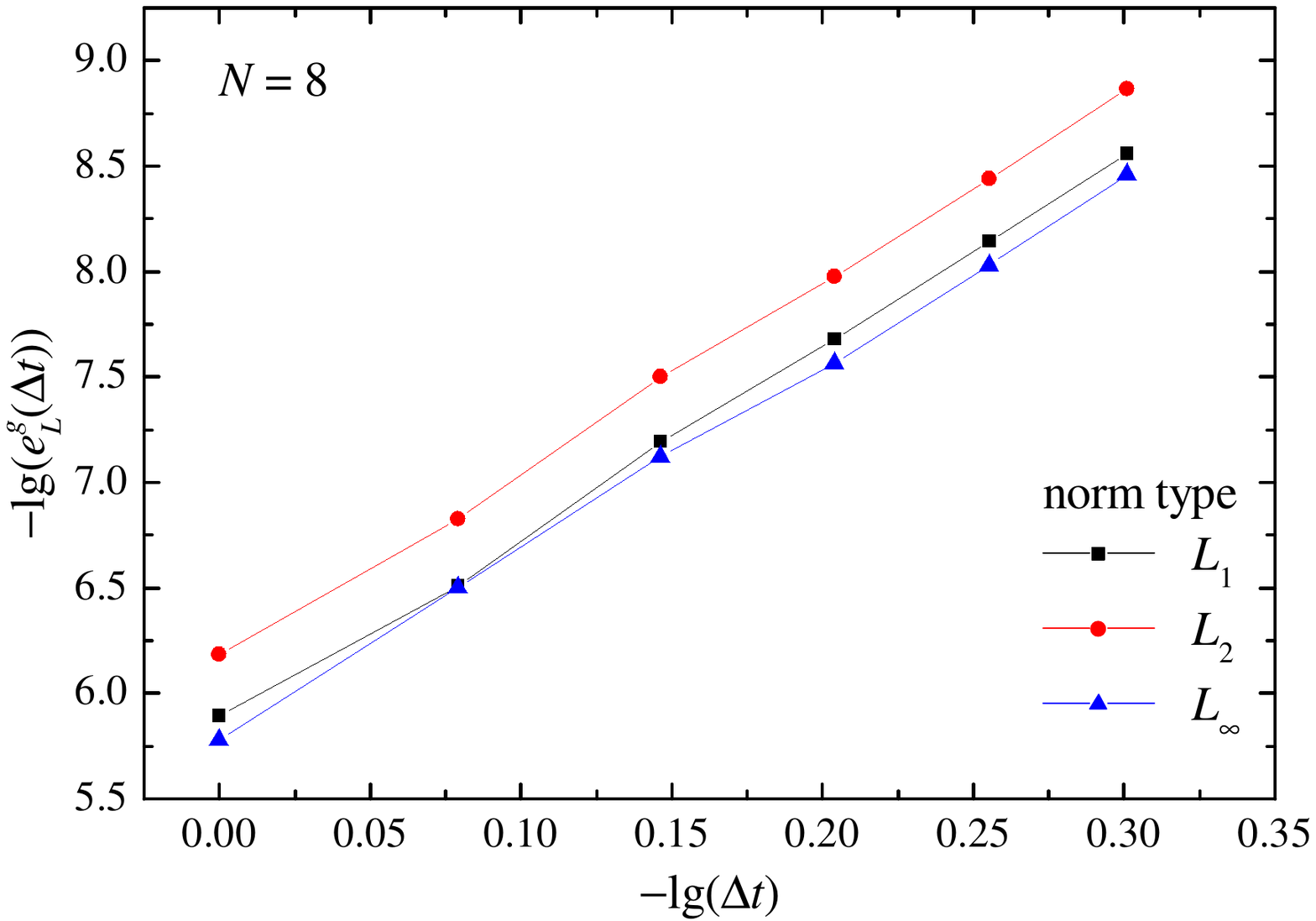}
\vspace{-8mm}\caption{\label{fig:pend_ind2_errors:c2}}
\end{subfigure}\hspace{6mm}
\begin{subfigure}{0.275\textwidth}
\includegraphics[width=\textwidth]{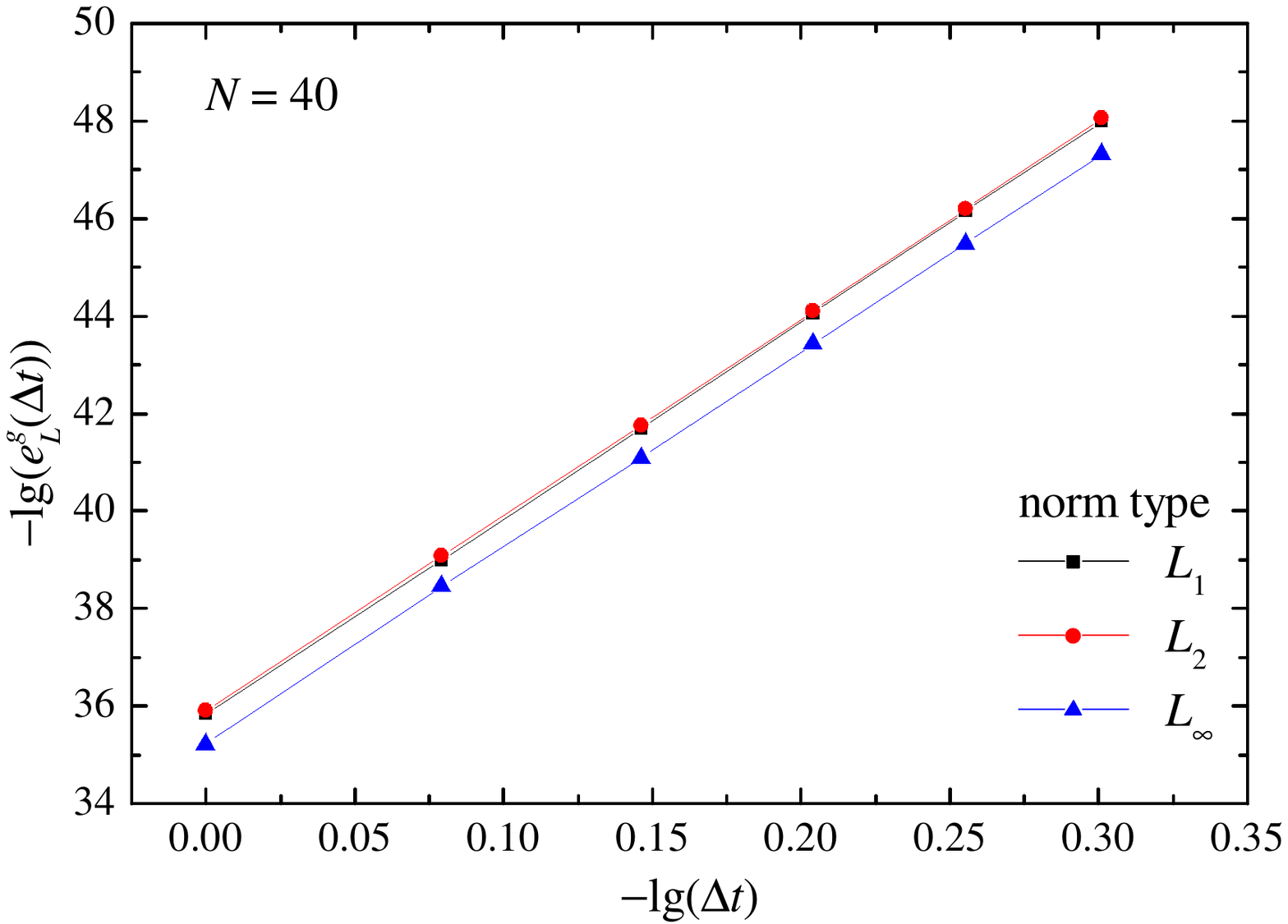}
\vspace{-8mm}\caption{\label{fig:pend_ind2_errors:c3}}
\end{subfigure}\\[-2mm]
\begin{subfigure}{0.275\textwidth}
\includegraphics[width=\textwidth]{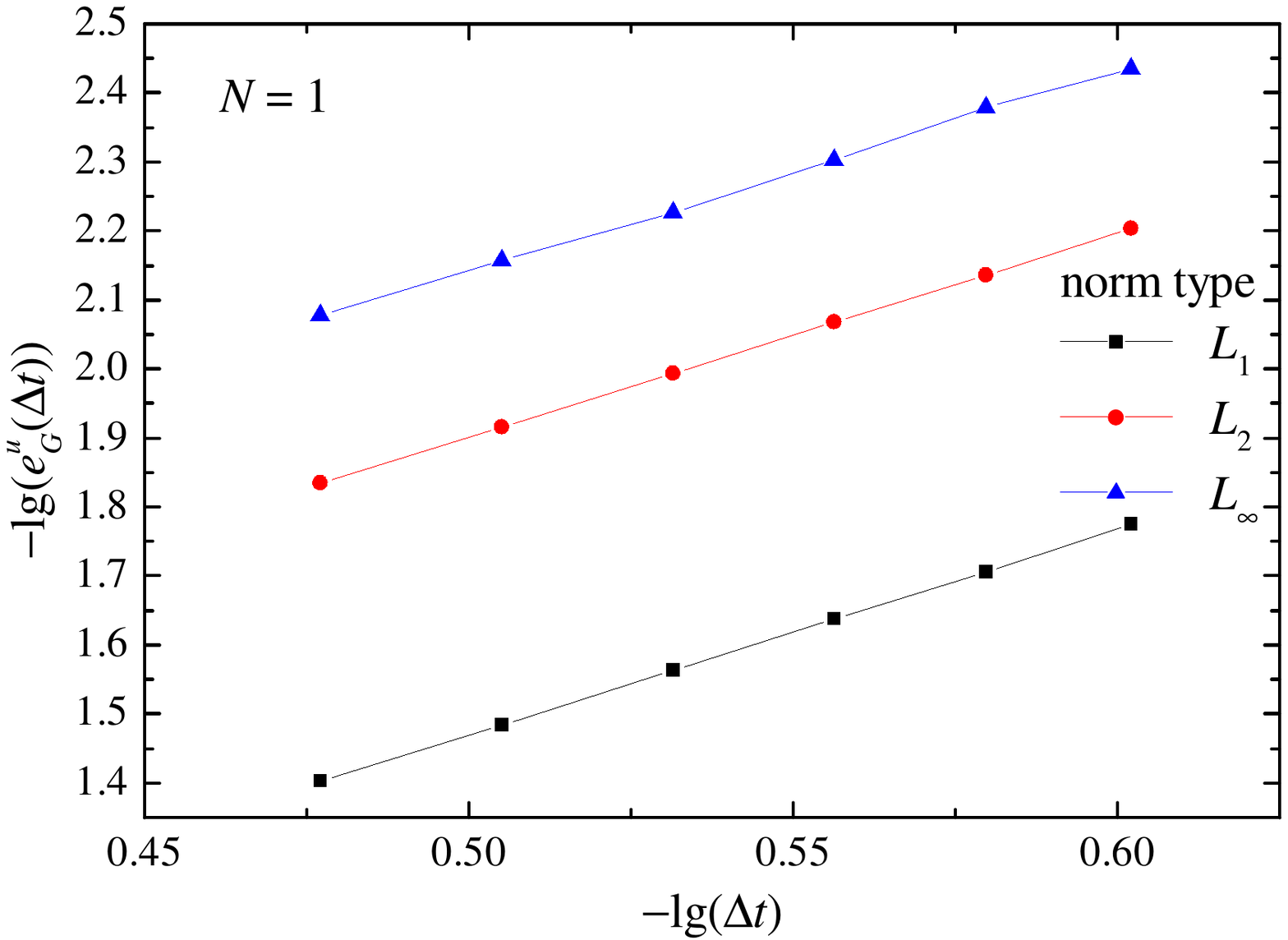}
\vspace{-8mm}\caption{\label{fig:pend_ind2_errors:d1}}
\end{subfigure}\hspace{6mm}
\begin{subfigure}{0.275\textwidth}
\includegraphics[width=\textwidth]{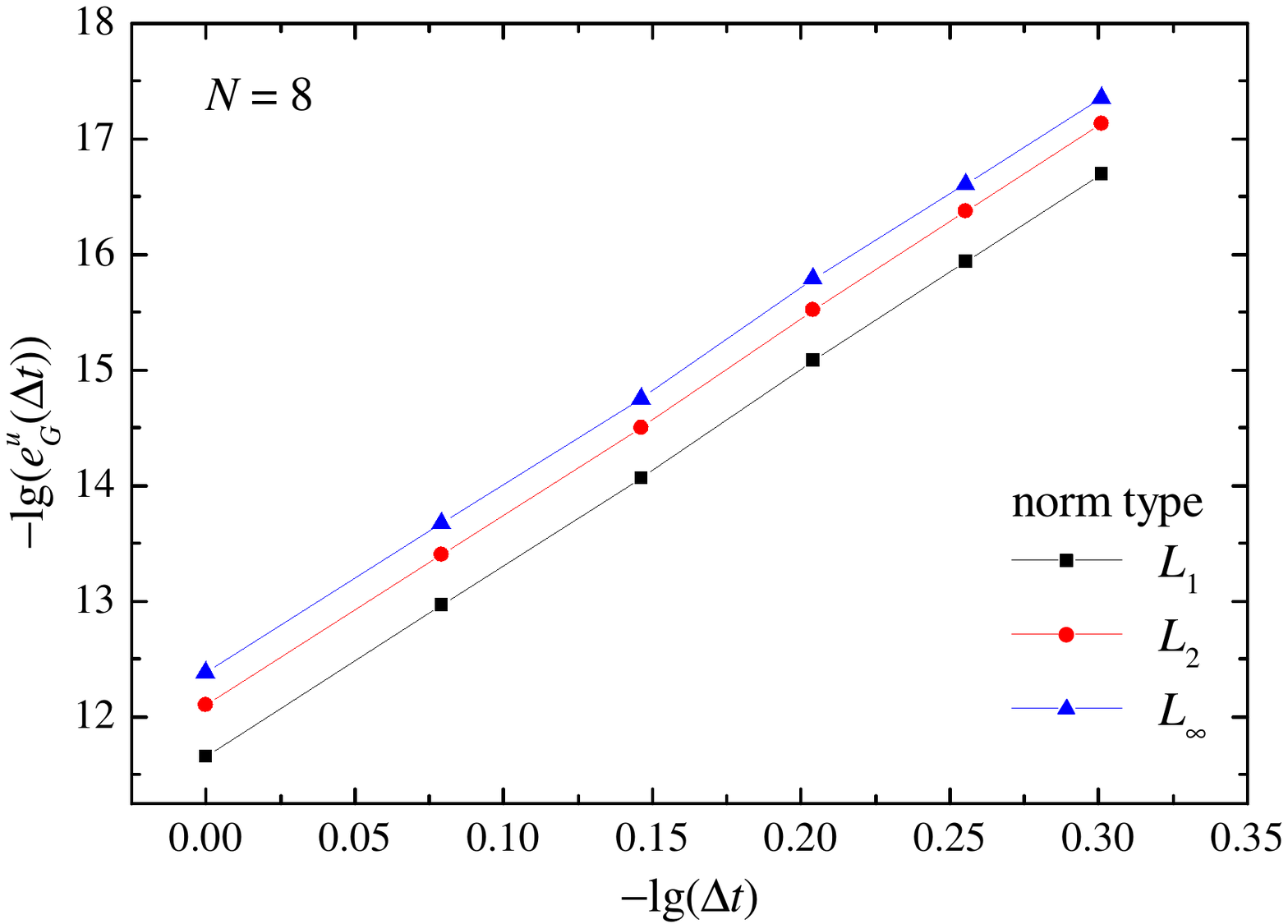}
\vspace{-8mm}\caption{\label{fig:pend_ind2_errors:d2}}
\end{subfigure}\hspace{6mm}
\begin{subfigure}{0.275\textwidth}
\includegraphics[width=\textwidth]{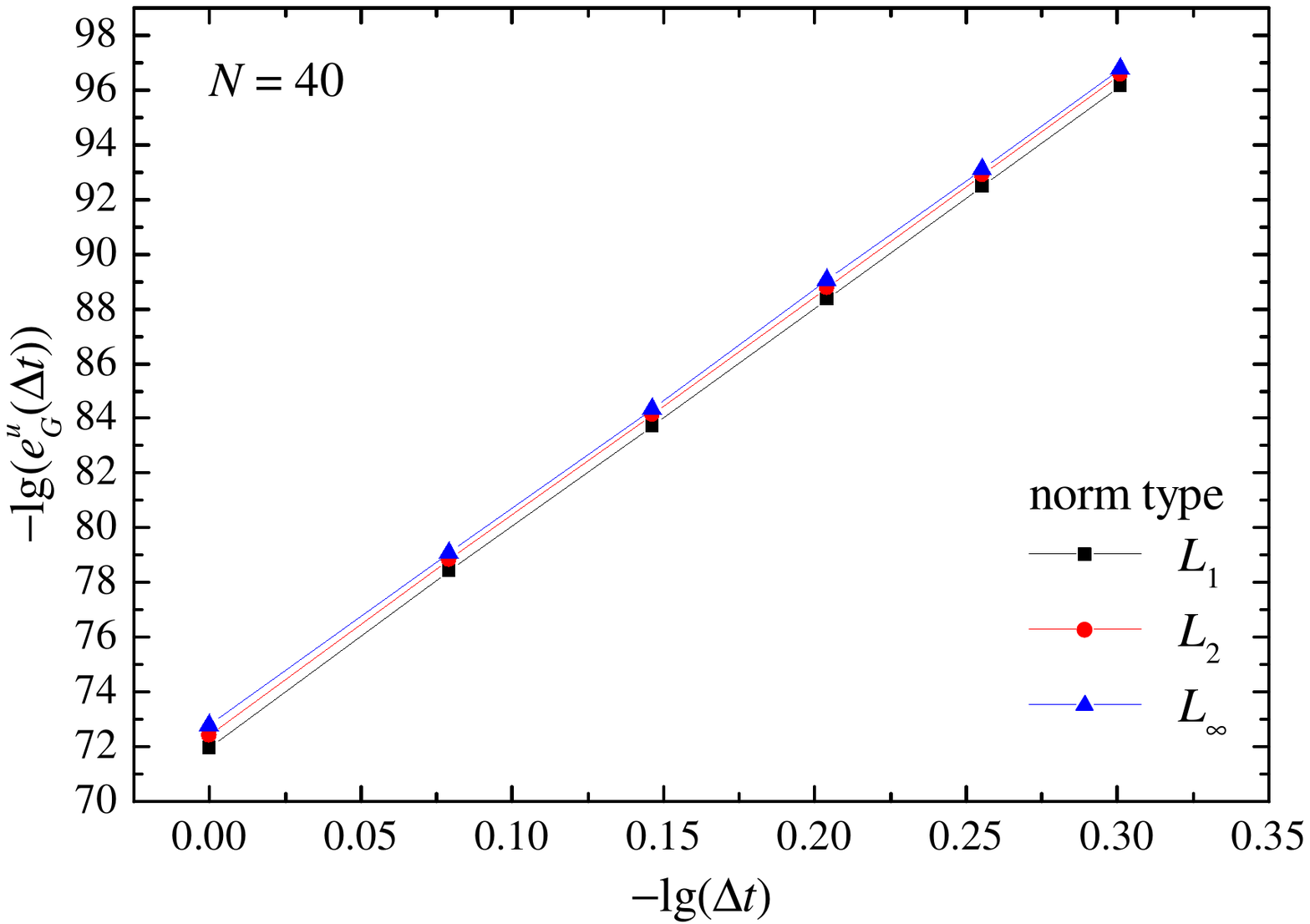}
\vspace{-8mm}\caption{\label{fig:pend_ind2_errors:d3}}
\end{subfigure}\\[-2mm]
\begin{subfigure}{0.275\textwidth}
\includegraphics[width=\textwidth]{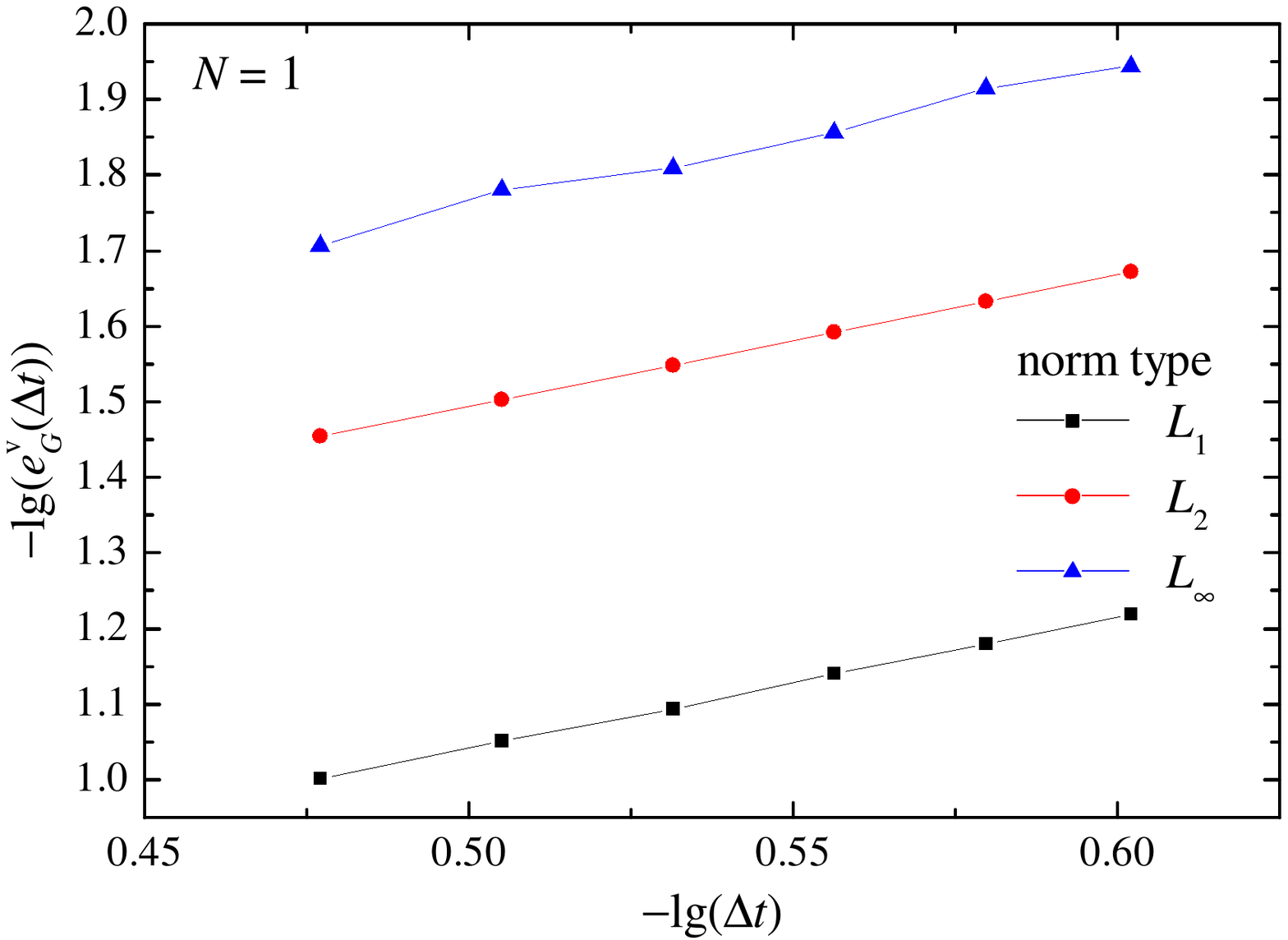}
\vspace{-8mm}\caption{\label{fig:pend_ind2_errors:e1}}
\end{subfigure}\hspace{6mm}
\begin{subfigure}{0.275\textwidth}
\includegraphics[width=\textwidth]{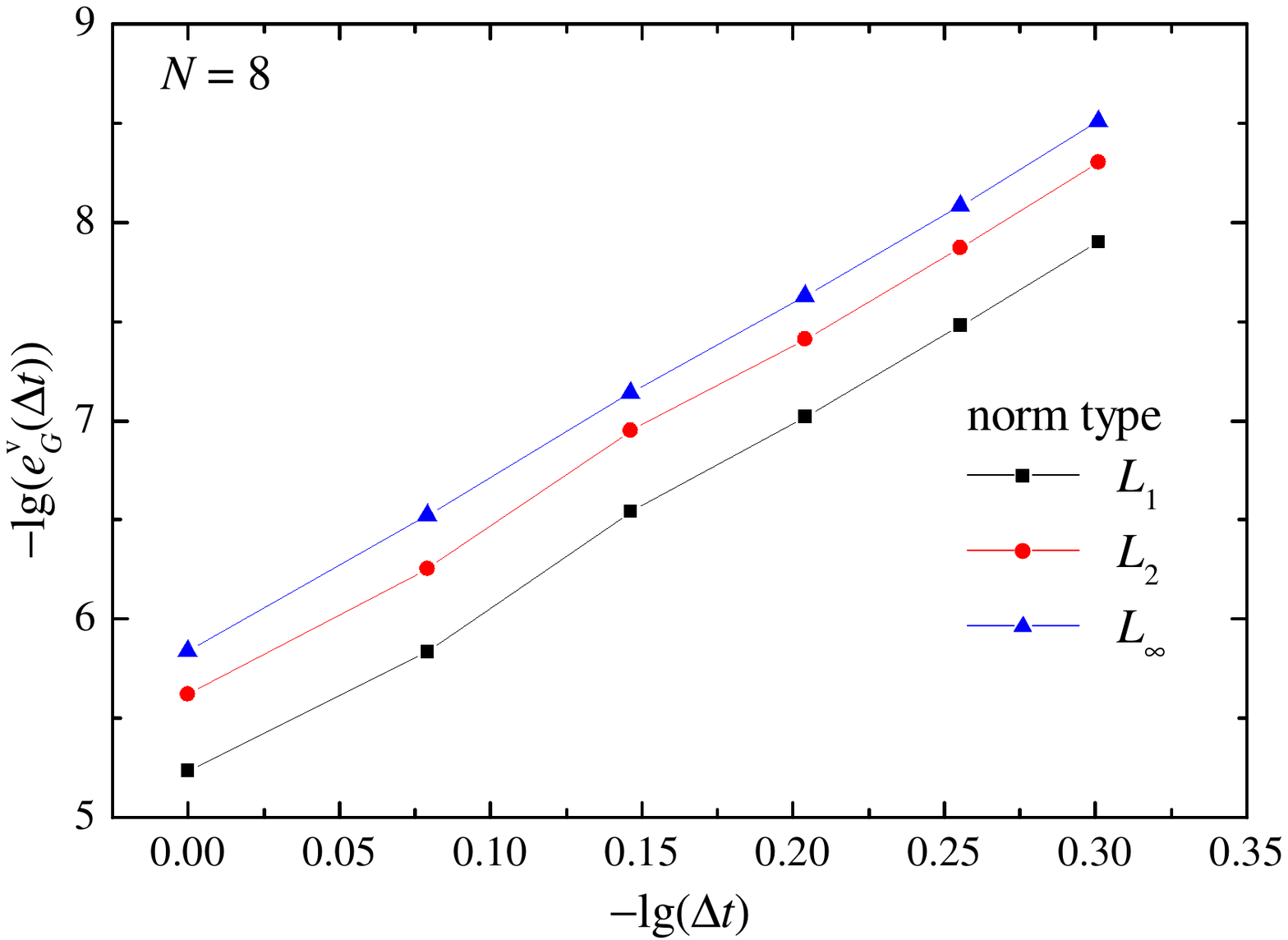}
\vspace{-8mm}\caption{\label{fig:pend_ind2_errors:e2}}
\end{subfigure}\hspace{6mm}
\begin{subfigure}{0.275\textwidth}
\includegraphics[width=\textwidth]{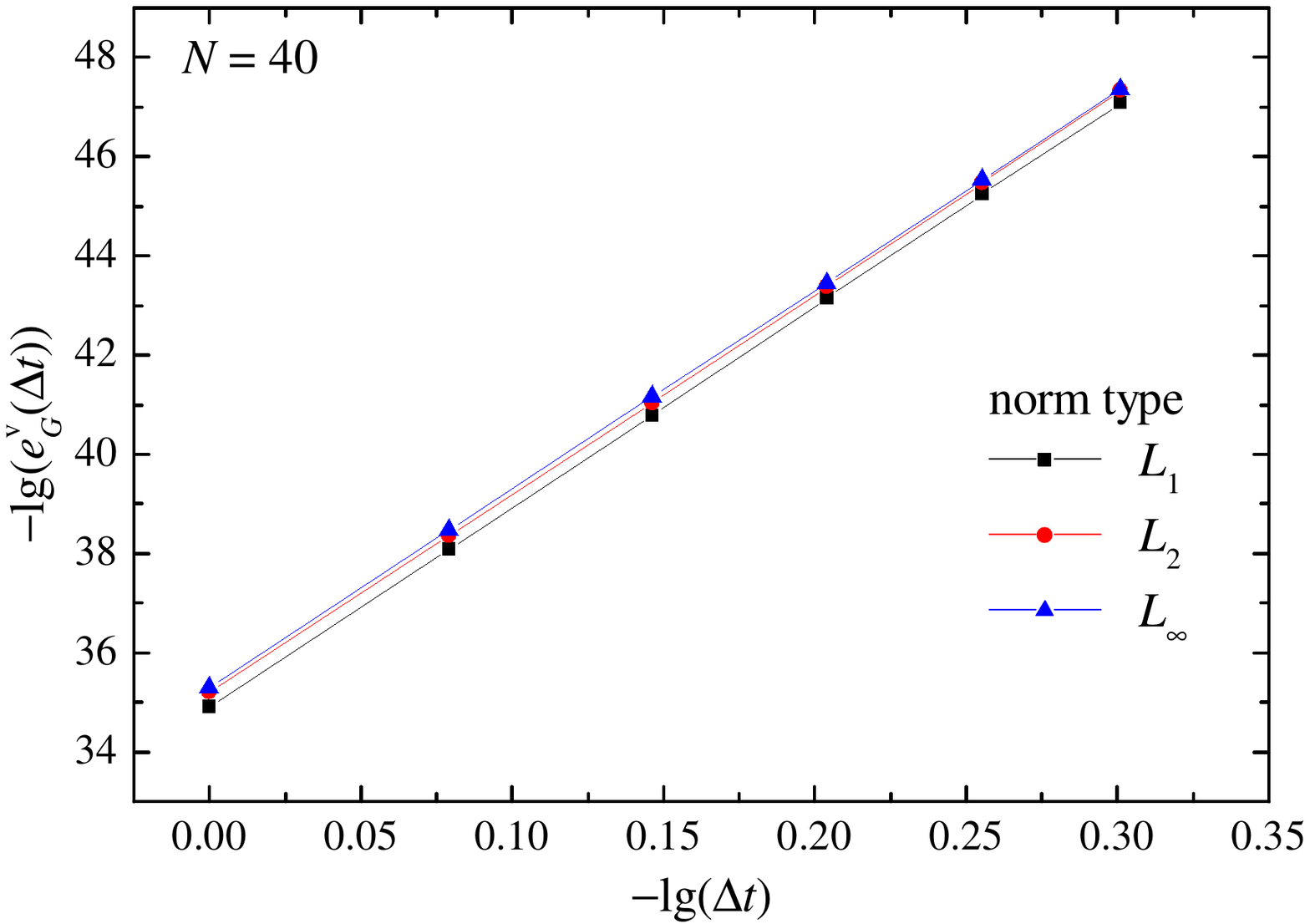}
\vspace{-8mm}\caption{\label{fig:pend_ind2_errors:e3}}
\end{subfigure}\\[-2mm]
\begin{subfigure}{0.275\textwidth}
\includegraphics[width=\textwidth]{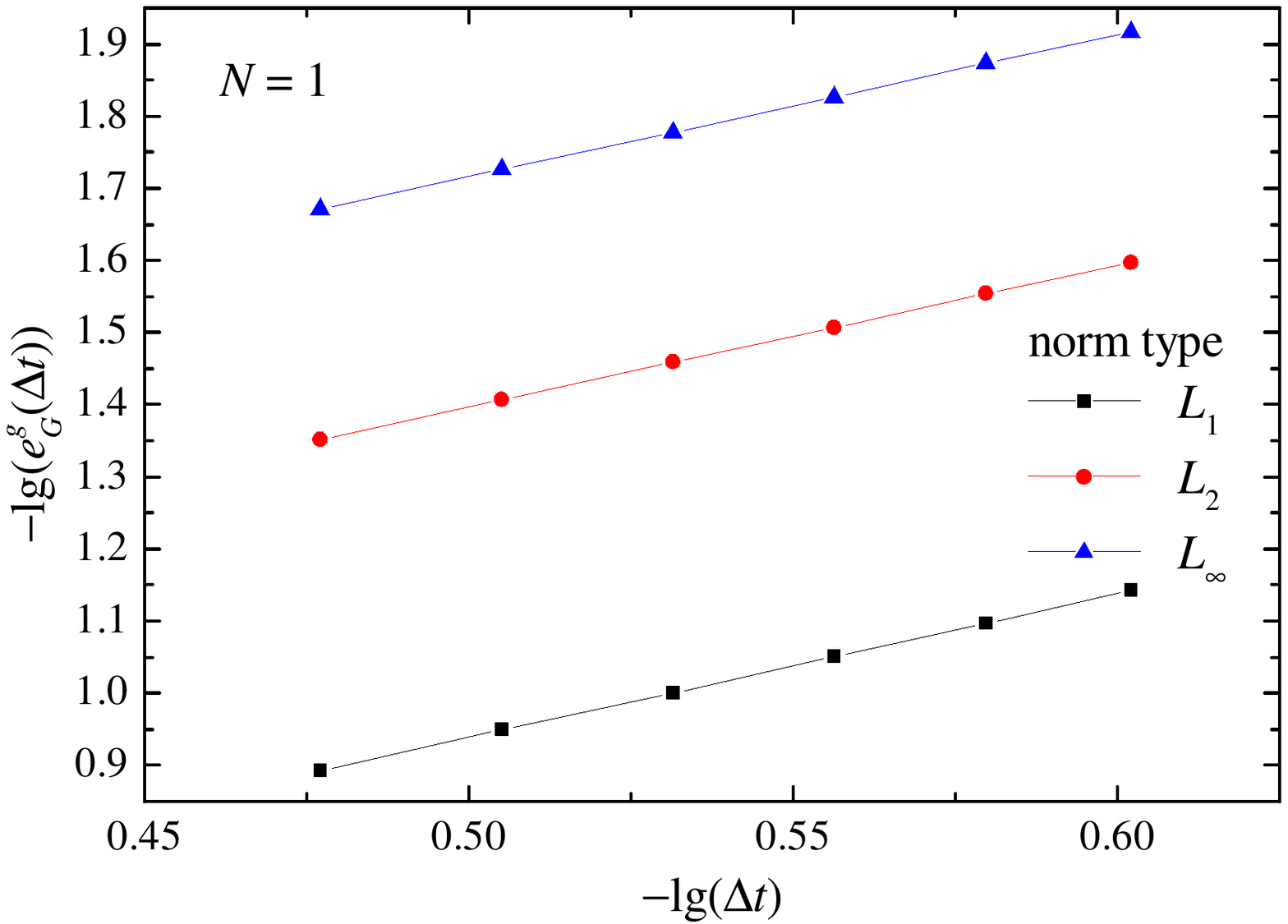}
\vspace{-8mm}\caption{\label{fig:pend_ind2_errors:f1}}
\end{subfigure}\hspace{6mm}
\begin{subfigure}{0.275\textwidth}
\includegraphics[width=\textwidth]{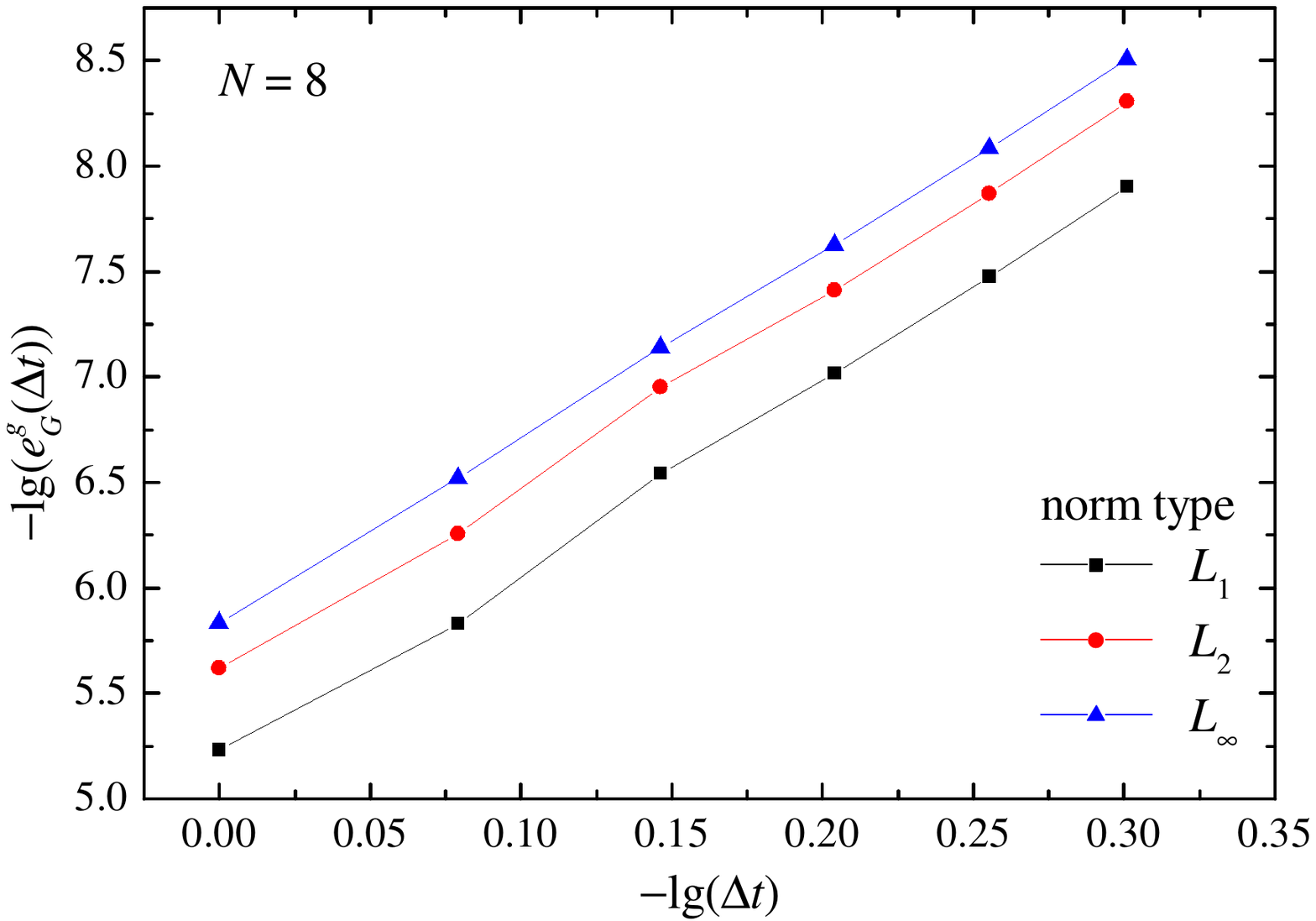}
\vspace{-8mm}\caption{\label{fig:pend_ind2_errors:f2}}
\end{subfigure}\hspace{6mm}
\begin{subfigure}{0.275\textwidth}
\includegraphics[width=\textwidth]{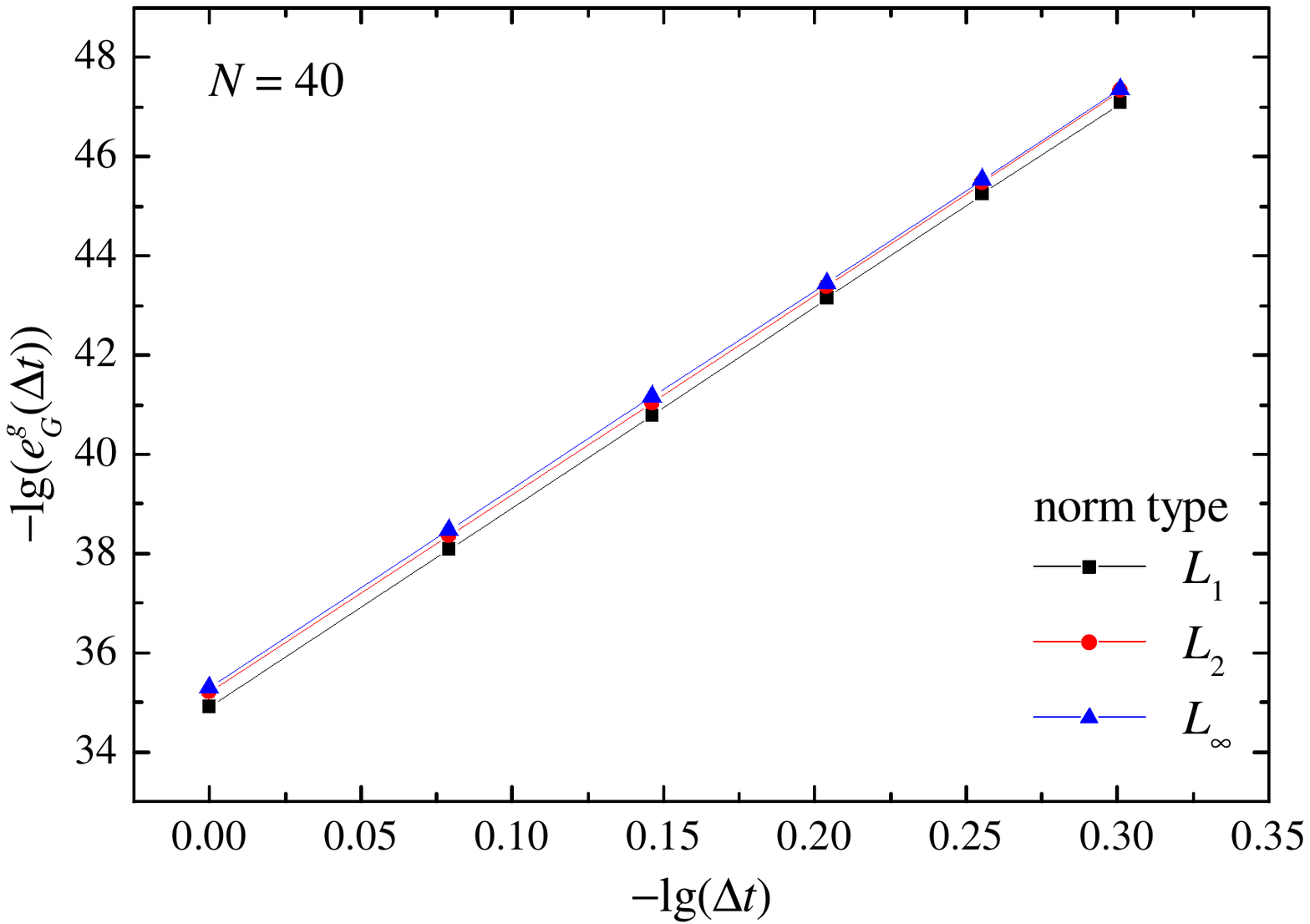}
\vspace{-8mm}\caption{\label{fig:pend_ind2_errors:f3}}
\end{subfigure}\\[-2mm]
\caption{%
Log-log plot of the dependence of the global errors for the local solution $e_{L}^{u}$ (\subref{fig:pend_ind2_errors:a1}, \subref{fig:pend_ind2_errors:a2}, \subref{fig:pend_ind2_errors:a3}), $e_{L}^{v}$ (\subref{fig:pend_ind2_errors:b1}, \subref{fig:pend_ind2_errors:b2}, \subref{fig:pend_ind2_errors:b3}), $e_{L}^{g}$ (\subref{fig:pend_ind2_errors:c1}, \subref{fig:pend_ind2_errors:c2}, \subref{fig:pend_ind2_errors:c3}) and the solution at nodes $e_{G}^{u}$ (\subref{fig:pend_ind2_errors:d1}, \subref{fig:pend_ind2_errors:d2}, \subref{fig:pend_ind2_errors:d3}), $e_{G}^{v}$ (\subref{fig:pend_ind2_errors:e1}, \subref{fig:pend_ind2_errors:e2}, \subref{fig:pend_ind2_errors:e3}), $e_{G}^{g}$ (\subref{fig:pend_ind2_errors:f1}, \subref{fig:pend_ind2_errors:f2}, \subref{fig:pend_ind2_errors:f3}) on the discretization step $\mathrm{\Delta}t$, obtained in the norms $L_{1}$, $L_{2}$ and $L_{\infty}$, by numerical solution of the DAE system (\ref{eq:math_pend_dae_ind_3}) of index 2 obtained using polynomials with degrees $N = 1$ (\subref{fig:pend_ind2_errors:a1}, \subref{fig:pend_ind2_errors:b1}, \subref{fig:pend_ind2_errors:c1}, \subref{fig:pend_ind2_errors:d1}, \subref{fig:pend_ind2_errors:e1}, \subref{fig:pend_ind2_errors:f1}), $N = 8$ (\subref{fig:pend_ind2_errors:a2}, \subref{fig:pend_ind2_errors:b2}, \subref{fig:pend_ind2_errors:c2}, \subref{fig:pend_ind2_errors:d2}, \subref{fig:pend_ind2_errors:e2}, \subref{fig:pend_ind2_errors:f2}) and $N = 40$ (\subref{fig:pend_ind2_errors:a3}, \subref{fig:pend_ind2_errors:b3}, \subref{fig:pend_ind2_errors:c3}, \subref{fig:pend_ind2_errors:d3}, \subref{fig:pend_ind2_errors:e3}, \subref{fig:pend_ind2_errors:f3}).
}
\label{fig:pend_ind2_errors}
\end{figure} 

\begin{table*}[h!]
\centering
\caption{%
Convergence orders $p_{L_{1}}^{n}$, $p_{L_{2}}^{n}$, $p_{L_{\infty}}^{n}$, calculated in norms $L_{1}$, $L_{2}$, $L_{\infty}$, of \textit{the numerical solution at the nodes} $(\mathbf{u}_{n}, \mathbf{v}_{n})$ of the ADER-DG method for the DAE system (\ref{eq:math_pend_dae_ind_3}) of index 2; $N$ is the degree of the basis polynomials $\varphi_{p}$. Orders $p^{n, u}$ are calculated for solution $\mathbf{u}_{n}$; orders $p^{n, v}$ --- for solution $\mathbf{v}_{n}$; orders $p^{n, g}$ --- for the conditions $\mathbf{g} = 0$ on the numerical solution at the nodes $(\mathbf{u}_{n}, \mathbf{v}_{n})$. The theoretical value of convergence order $p_{\rm th.}^{n} = 2N+1$ is applicable for the ADER-DG method for ODE problems and is presented for comparison.
}
\label{tab:conv_orders_nodes_pend_ind2}
\begin{tabular}{@{}|l|lll|lll|lll|c|@{}}
\toprule
$N$ & $p_{L_{1}}^{n, u}$ & $p_{L_{2}}^{n, u}$ & $p_{L_{\infty}}^{n, u}$ & $p_{L_{1}}^{n, v}$ & $p_{L_{2}}^{n, v}$ & $p_{L_{\infty}}^{n, v}$ & $p_{L_{1}}^{n, g}$ & $p_{L_{2}}^{n, g}$ & $p_{L_{\infty}}^{n, g}$ & $p_{\rm th.}^{n}$ \\
\midrule
$1$	&	$2.97$	&	$2.96$	&	$2.90$	&	$1.74$	&	$1.74$	&	$1.87$	&	$1.99$	&	$1.97$	&	$1.97$	&	$3$\\
$2$	&	$5.00$	&	$4.99$	&	$4.91$	&	$2.99$	&	$2.99$	&	$3.03$	&	$2.99$	&	$2.99$	&	$3.04$	&	$5$\\
$3$	&	$7.00$	&	$6.99$	&	$6.90$	&	$3.95$	&	$3.99$	&	$3.98$	&	$3.95$	&	$3.99$	&	$3.98$	&	$7$\\
$4$	&	$8.99$	&	$8.99$	&	$8.90$	&	$4.98$	&	$4.99$	&	$4.98$	&	$4.98$	&	$4.99$	&	$4.98$	&	$9$\\
$5$	&	$11.00$	&	$10.99$	&	$10.92$	&	$5.97$	&	$5.98$	&	$5.97$	&	$5.97$	&	$5.98$	&	$5.97$	&	$11$\\
$6$	&	$12.75$	&	$12.74$	&	$12.63$	&	$7.12$	&	$7.05$	&	$6.87$	&	$7.12$	&	$7.05$	&	$6.87$	&	$13$\\
$7$	&	$14.96$	&	$14.92$	&	$14.72$	&	$7.64$	&	$7.68$	&	$7.63$	&	$7.64$	&	$7.68$	&	$7.63$	&	$15$\\
$8$	&	$16.79$	&	$16.76$	&	$16.59$	&	$8.96$	&	$8.95$	&	$8.85$	&	$8.96$	&	$8.95$	&	$8.85$	&	$17$\\
$9$	&	$18.87$	&	$18.82$	&	$18.65$	&	$9.85$	&	$9.82$	&	$9.65$	&	$9.85$	&	$9.82$	&	$9.65$	&	$19$\\
$10$	&	$20.78$	&	$20.76$	&	$20.57$	&	$10.72$	&	$10.83$	&	$10.85$	&	$10.72$	&	$10.83$	&	$10.85$	&	$21$\\
$11$	&	$22.82$	&	$22.77$	&	$22.59$	&	$11.89$	&	$11.86$	&	$11.64$	&	$11.89$	&	$11.86$	&	$11.64$	&	$23$\\
$12$	&	$24.73$	&	$24.71$	&	$24.52$	&	$12.78$	&	$12.80$	&	$12.68$	&	$12.78$	&	$12.80$	&	$12.68$	&	$25$\\
$13$	&	$26.79$	&	$26.73$	&	$26.53$	&	$13.84$	&	$13.81$	&	$13.63$	&	$13.84$	&	$13.81$	&	$13.63$	&	$27$\\
$14$	&	$28.69$	&	$28.67$	&	$28.47$	&	$14.91$	&	$14.79$	&	$14.59$	&	$14.91$	&	$14.79$	&	$14.59$	&	$29$\\
$15$	&	$30.73$	&	$30.67$	&	$30.45$	&	$15.74$	&	$15.78$	&	$15.73$	&	$15.74$	&	$15.78$	&	$15.73$	&	$31$\\
$16$	&	$32.65$	&	$32.62$	&	$32.42$	&	$16.81$	&	$16.72$	&	$16.58$	&	$16.81$	&	$16.72$	&	$16.58$	&	$33$\\
$17$	&	$34.70$	&	$34.64$	&	$34.40$	&	$17.82$	&	$17.81$	&	$17.75$	&	$17.82$	&	$17.81$	&	$17.75$	&	$35$\\
$18$	&	$36.57$	&	$36.53$	&	$36.33$	&	$18.63$	&	$18.60$	&	$18.54$	&	$18.63$	&	$18.60$	&	$18.54$	&	$37$\\
$19$	&	$38.70$	&	$38.64$	&	$38.40$	&	$19.90$	&	$19.87$	&	$19.74$	&	$19.90$	&	$19.87$	&	$19.74$	&	$39$\\
$20$	&	$40.46$	&	$40.41$	&	$40.22$	&	$20.49$	&	$20.49$	&	$20.48$	&	$20.49$	&	$20.49$	&	$20.48$	&	$41$\\
\midrule
$25$	&	$50.58$	&	$50.52$	&	$50.25$	&	$25.91$	&	$25.83$	&	$25.63$	&	$25.91$	&	$25.83$	&	$25.63$	&	$51$\\
$30$	&	$60.40$	&	$60.32$	&	$60.03$	&	$30.67$	&	$30.57$	&	$30.41$	&	$30.67$	&	$30.57$	&	$30.41$	&	$61$\\
$35$	&	$70.00$	&	$69.97$	&	$69.79$	&	$35.57$	&	$35.43$	&	$35.23$	&	$35.57$	&	$35.43$	&	$35.23$	&	$71$\\
$40$	&	$80.22$	&	$80.12$	&	$79.73$	&	$40.49$	&	$40.27$	&	$40.07$	&	$40.49$	&	$40.27$	&	$40.07$	&	$81$\\
\bottomrule
\end{tabular}
\end{table*} 

\begin{table*}[h!]
\centering
\caption{%
Convergence orders $p_{L_{1}}^{l}$, $p_{L_{2}}^{l}$, $p_{L_{\infty}}^{l}$, calculated in norms $L_{1}$, $L_{2}$, $L_{\infty}$, of \textit{the local solution} $(\mathbf{u}_{L}, \mathbf{v}_{L})$ (represented between the nodes) of the ADER-DG method for the DAE system (\ref{eq:math_pend_dae_ind_3}) of index 2; $N$ is the degree of the basis polynomials $\varphi_{p}$. Orders $p^{l, u}$ are calculated for solution $\mathbf{u}_{L}$; orders $p^{l, v}$ --- for solution $\mathbf{v}_{L}$; orders $p^{l, g}$ --- for the conditions $\mathbf{g} = 0$ on the local solution $(\mathbf{u}_{L}, \mathbf{v}_{L})$. The theoretical value of convergence order $p_{\rm th.}^{l} = N+1$ is applicable for the ADER-DG method for ODE problems and is presented for comparison.
}
\label{tab:conv_orders_local_pend_ind2}
\begin{tabular}{@{}|l|lll|lll|lll|c|@{}}
\toprule
$N$ & $p_{L_{1}}^{l, u}$ & $p_{L_{2}}^{l, u}$ & $p_{L_{\infty}}^{l, u}$ & $p_{L_{1}}^{l, v}$ & $p_{L_{2}}^{l, v}$ & $p_{L_{\infty}}^{l, v}$ & $p_{L_{1}}^{l, g}$ & $p_{L_{2}}^{l, g}$ & $p_{L_{\infty}}^{l, g}$ & $p_{\rm th.}^{l}$ \\
\midrule
$1$	&	$2.35$	&	$2.29$	&	$1.90$	&	$2.22$	&	$2.26$	&	$1.86$	&	$2.07$	&	$2.03$	&	$1.98$	&	$2$\\
$2$	&	$2.99$	&	$2.97$	&	$2.95$	&	$3.00$	&	$2.99$	&	$3.04$	&	$2.99$	&	$2.99$	&	$3.05$	&	$3$\\
$3$	&	$4.03$	&	$4.02$	&	$3.98$	&	$3.97$	&	$3.99$	&	$3.98$	&	$3.99$	&	$4.00$	&	$3.98$	&	$4$\\
$4$	&	$5.05$	&	$5.04$	&	$4.96$	&	$4.99$	&	$4.99$	&	$4.98$	&	$4.98$	&	$4.99$	&	$5.03$	&	$5$\\
$5$	&	$6.00$	&	$5.99$	&	$5.97$	&	$5.97$	&	$5.99$	&	$5.97$	&	$5.99$	&	$5.99$	&	$5.97$	&	$6$\\
$6$	&	$6.87$	&	$6.87$	&	$6.87$	&	$7.12$	&	$7.05$	&	$6.89$	&	$7.10$	&	$7.04$	&	$6.85$	&	$7$\\
$7$	&	$7.96$	&	$7.94$	&	$7.86$	&	$7.68$	&	$7.73$	&	$7.75$	&	$7.72$	&	$7.76$	&	$7.75$	&	$8$\\
$8$	&	$8.91$	&	$8.90$	&	$8.82$	&	$8.96$	&	$8.95$	&	$8.87$	&	$8.94$	&	$8.94$	&	$8.83$	&	$9$\\
$9$	&	$9.94$	&	$9.92$	&	$9.82$	&	$9.88$	&	$9.86$	&	$9.77$	&	$9.87$	&	$9.86$	&	$9.77$	&	$10$\\
$10$	&	$10.91$	&	$10.88$	&	$10.74$	&	$10.74$	&	$10.85$	&	$10.87$	&	$10.74$	&	$10.84$	&	$10.86$	&	$11$\\
$11$	&	$11.89$	&	$11.85$	&	$11.74$	&	$11.91$	&	$11.89$	&	$11.78$	&	$11.90$	&	$11.89$	&	$11.79$	&	$12$\\
$12$	&	$12.82$	&	$12.79$	&	$12.65$	&	$12.81$	&	$12.83$	&	$12.66$	&	$12.81$	&	$12.82$	&	$12.62$	&	$13$\\
$13$	&	$13.81$	&	$13.78$	&	$13.73$	&	$13.85$	&	$13.84$	&	$13.77$	&	$13.85$	&	$13.84$	&	$13.78$	&	$14$\\
$14$	&	$14.79$	&	$14.77$	&	$14.74$	&	$14.93$	&	$14.81$	&	$14.55$	&	$14.92$	&	$14.80$	&	$14.47$	&	$15$\\
$15$	&	$15.83$	&	$15.78$	&	$15.66$	&	$15.75$	&	$15.81$	&	$15.73$	&	$15.76$	&	$15.81$	&	$15.70$	&	$16$\\
$16$	&	$16.81$	&	$16.79$	&	$16.58$	&	$16.84$	&	$16.73$	&	$16.57$	&	$16.84$	&	$16.72$	&	$16.52$	&	$17$\\
$17$	&	$17.79$	&	$17.78$	&	$17.61$	&	$17.85$	&	$17.84$	&	$17.79$	&	$17.84$	&	$17.83$	&	$17.77$	&	$18$\\
$18$	&	$18.82$	&	$18.79$	&	$18.61$	&	$18.62$	&	$18.60$	&	$18.55$	&	$18.62$	&	$18.60$	&	$18.50$	&	$19$\\
$19$	&	$19.80$	&	$19.78$	&	$19.67$	&	$19.93$	&	$19.91$	&	$19.75$	&	$19.92$	&	$19.90$	&	$19.72$	&	$20$\\
$20$	&	$20.78$	&	$20.74$	&	$20.66$	&	$20.52$	&	$20.49$	&	$20.50$	&	$20.52$	&	$20.48$	&	$20.47$	&	$21$\\
\midrule
$25$	&	$25.68$	&	$25.66$	&	$25.64$	&	$26.04$	&	$25.93$	&	$25.65$	&	$26.05$	&	$25.94$	&	$25.64$	&	$26$\\
$30$	&	$30.64$	&	$30.58$	&	$30.46$	&	$30.69$	&	$30.60$	&	$30.45$	&	$30.68$	&	$30.59$	&	$30.43$	&	$31$\\
$35$	&	$35.54$	&	$35.42$	&	$35.14$	&	$35.58$	&	$35.46$	&	$35.36$	&	$35.57$	&	$35.45$	&	$35.35$	&	$36$\\
$40$	&	$40.37$	&	$40.25$	&	$40.11$	&	$40.48$	&	$40.34$	&	$40.18$	&	$40.48$	&	$40.34$	&	$40.17$	&	$41$\\
\bottomrule
\end{tabular}
\end{table*}

\begin{figure}[h!]
\captionsetup[subfigure]{%
	position=bottom,
	font+=smaller,
	textfont=normalfont,
	singlelinecheck=off,
	justification=raggedright
}
\centering
\begin{subfigure}{0.320\textwidth}
\includegraphics[width=\textwidth]{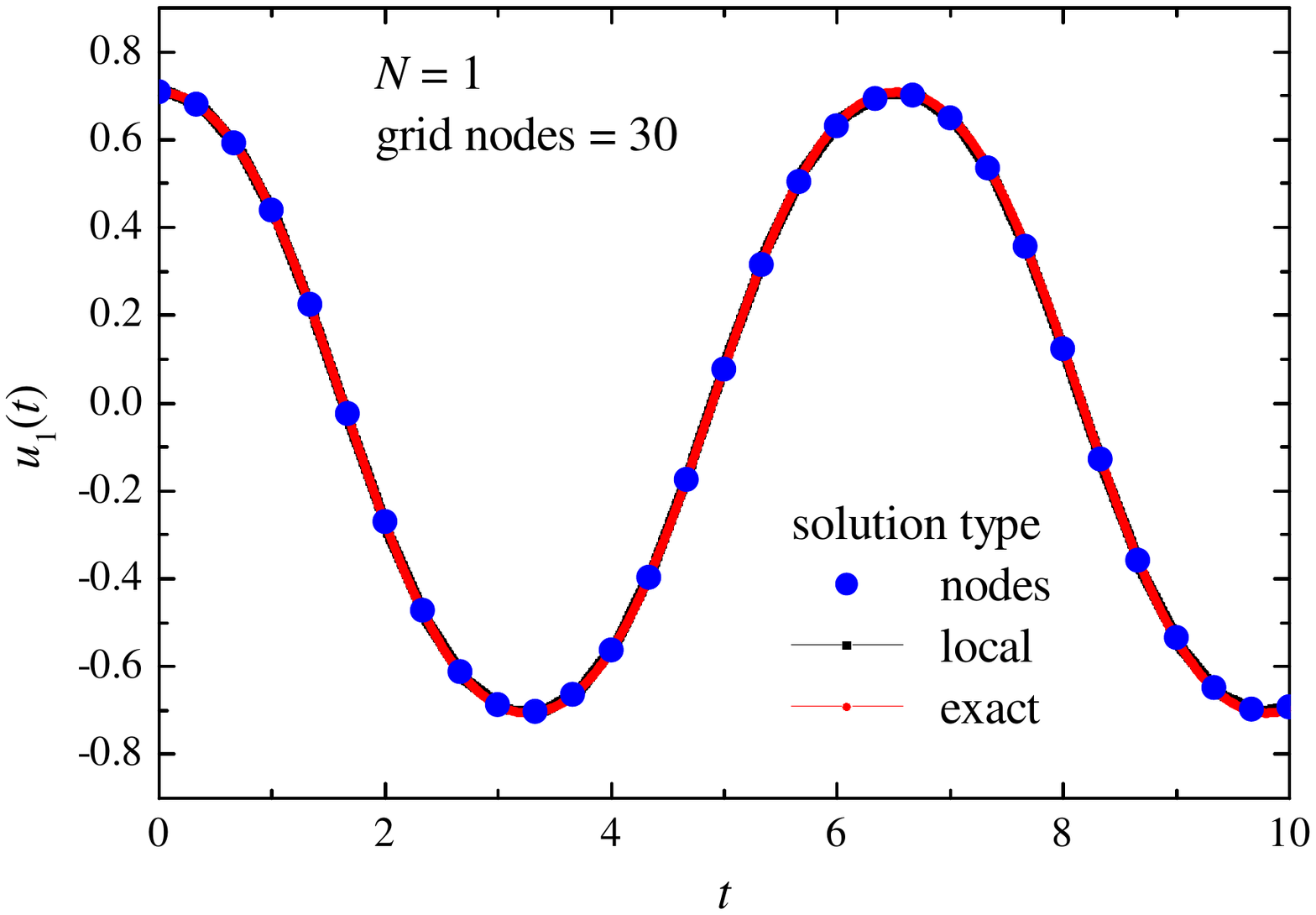}
\vspace{-8mm}\caption{\label{fig:pend_ind1_sol_uv:a1}}
\end{subfigure}
\begin{subfigure}{0.320\textwidth}
\includegraphics[width=\textwidth]{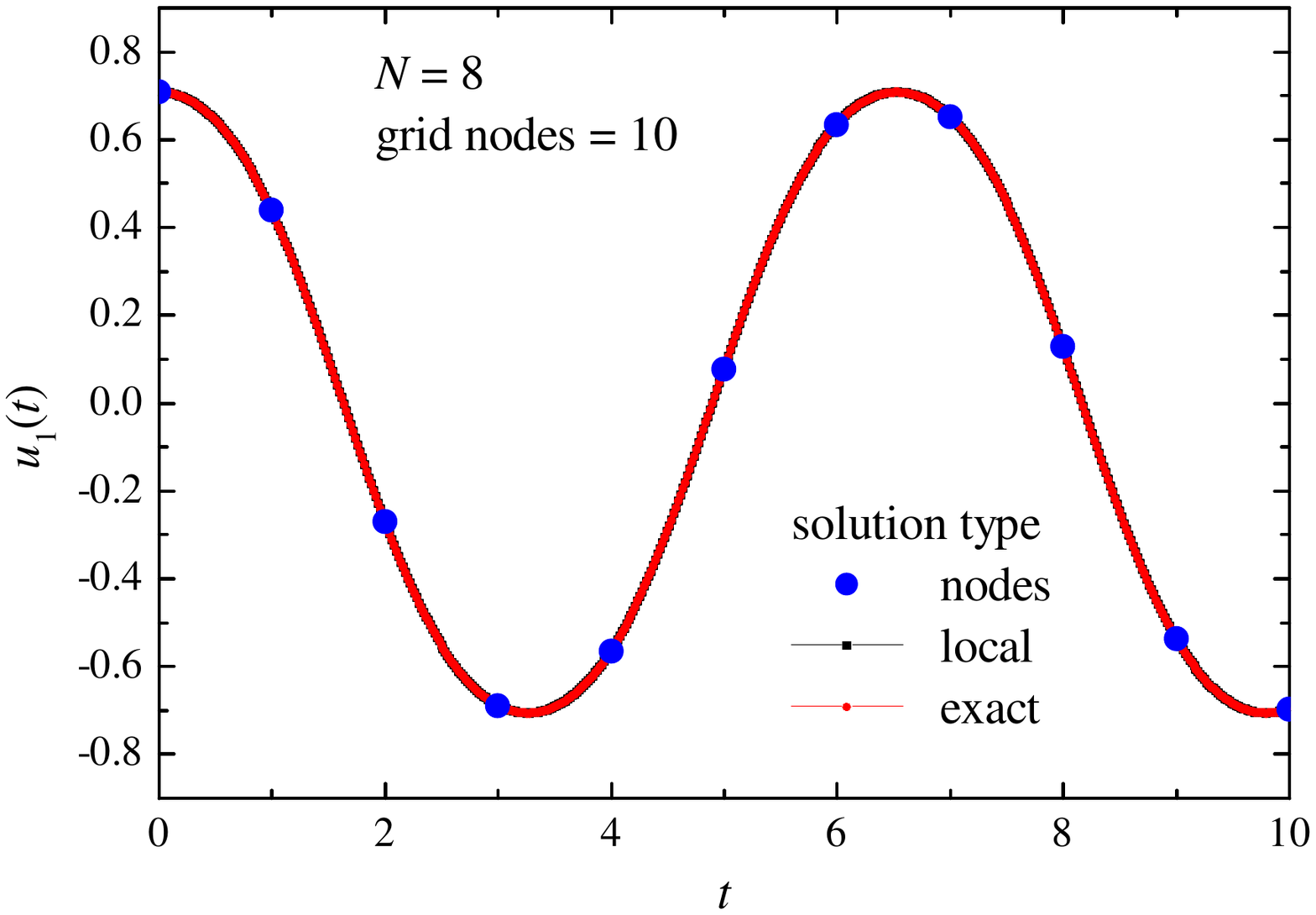}
\vspace{-8mm}\caption{\label{fig:pend_ind1_sol_uv:a2}}
\end{subfigure}
\begin{subfigure}{0.320\textwidth}
\includegraphics[width=\textwidth]{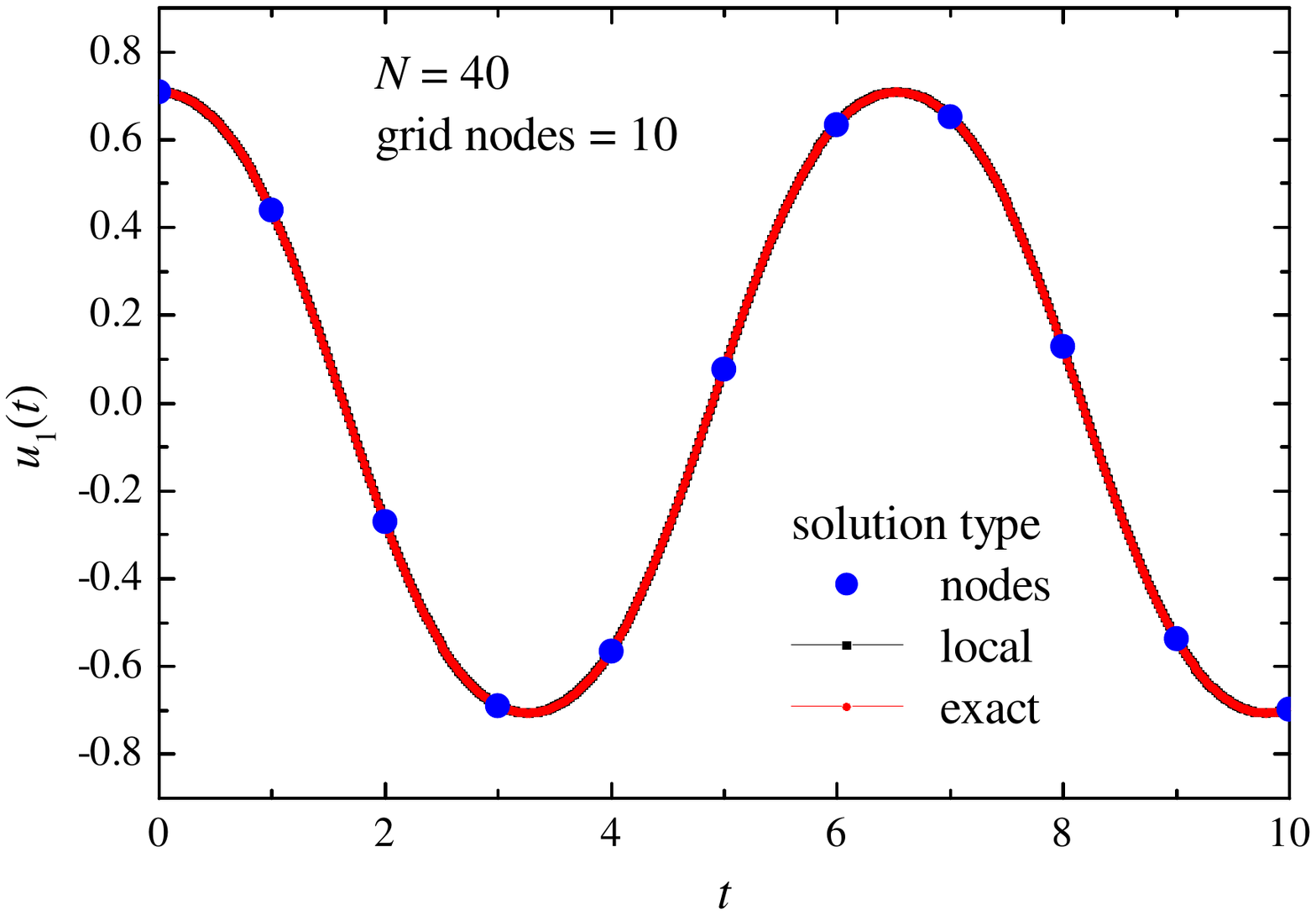}
\vspace{-8mm}\caption{\label{fig:pend_ind1_sol_uv:a3}}
\end{subfigure}\\
\begin{subfigure}{0.320\textwidth}
\includegraphics[width=\textwidth]{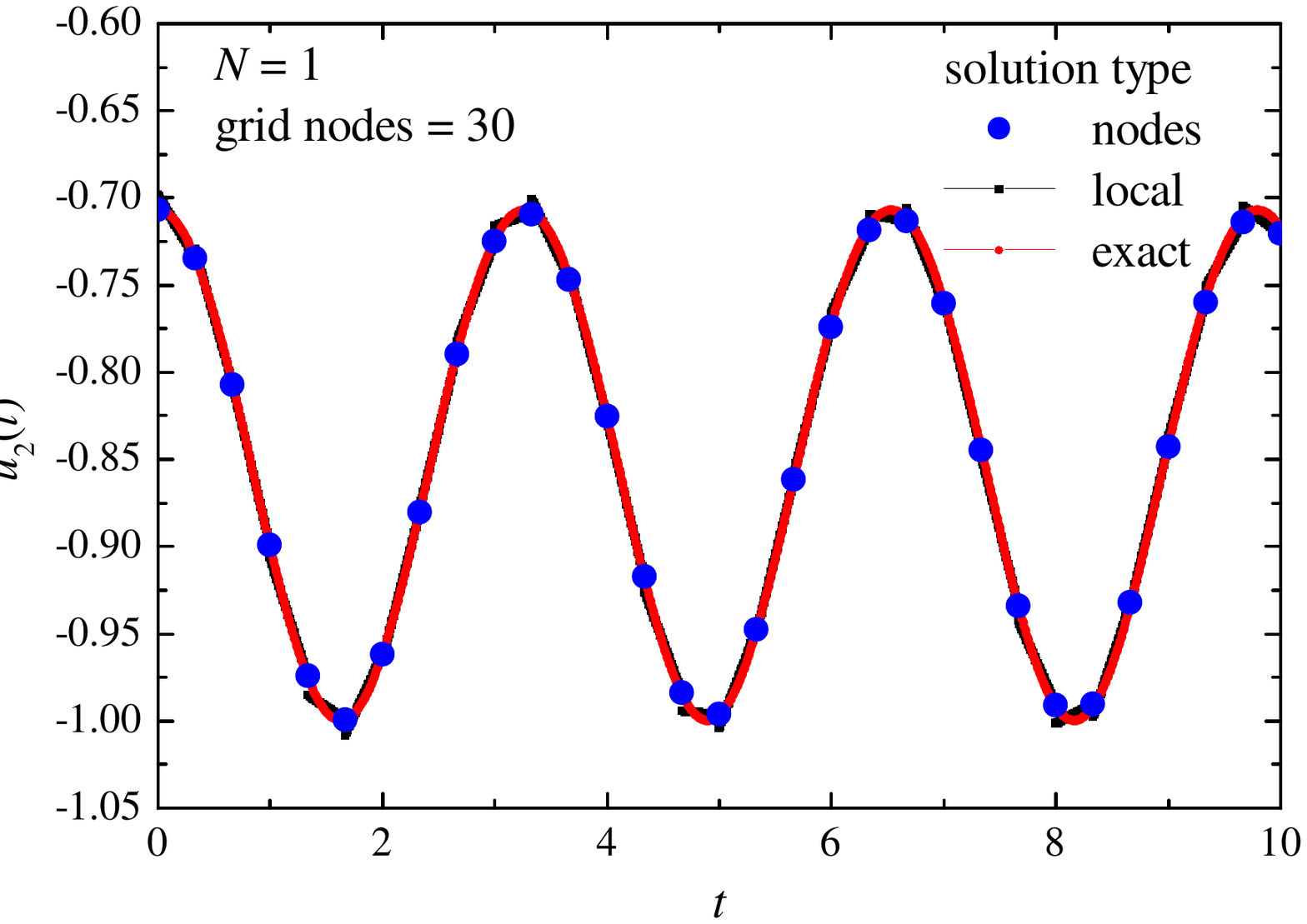}
\vspace{-8mm}\caption{\label{fig:pend_ind1_sol_uv:b1}}
\end{subfigure}
\begin{subfigure}{0.320\textwidth}
\includegraphics[width=\textwidth]{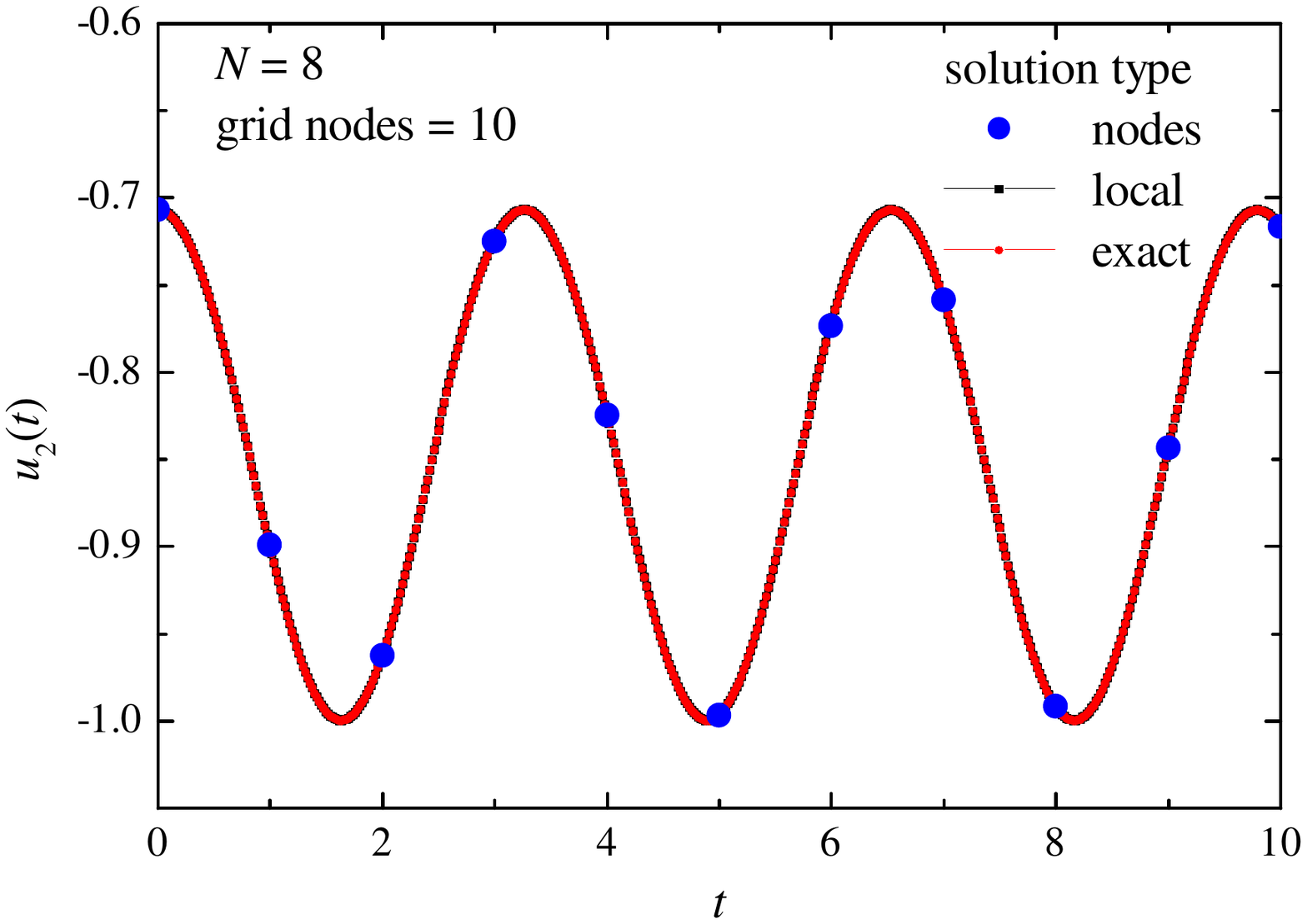}
\vspace{-8mm}\caption{\label{fig:pend_ind1_sol_uv:b2}}
\end{subfigure}
\begin{subfigure}{0.320\textwidth}
\includegraphics[width=\textwidth]{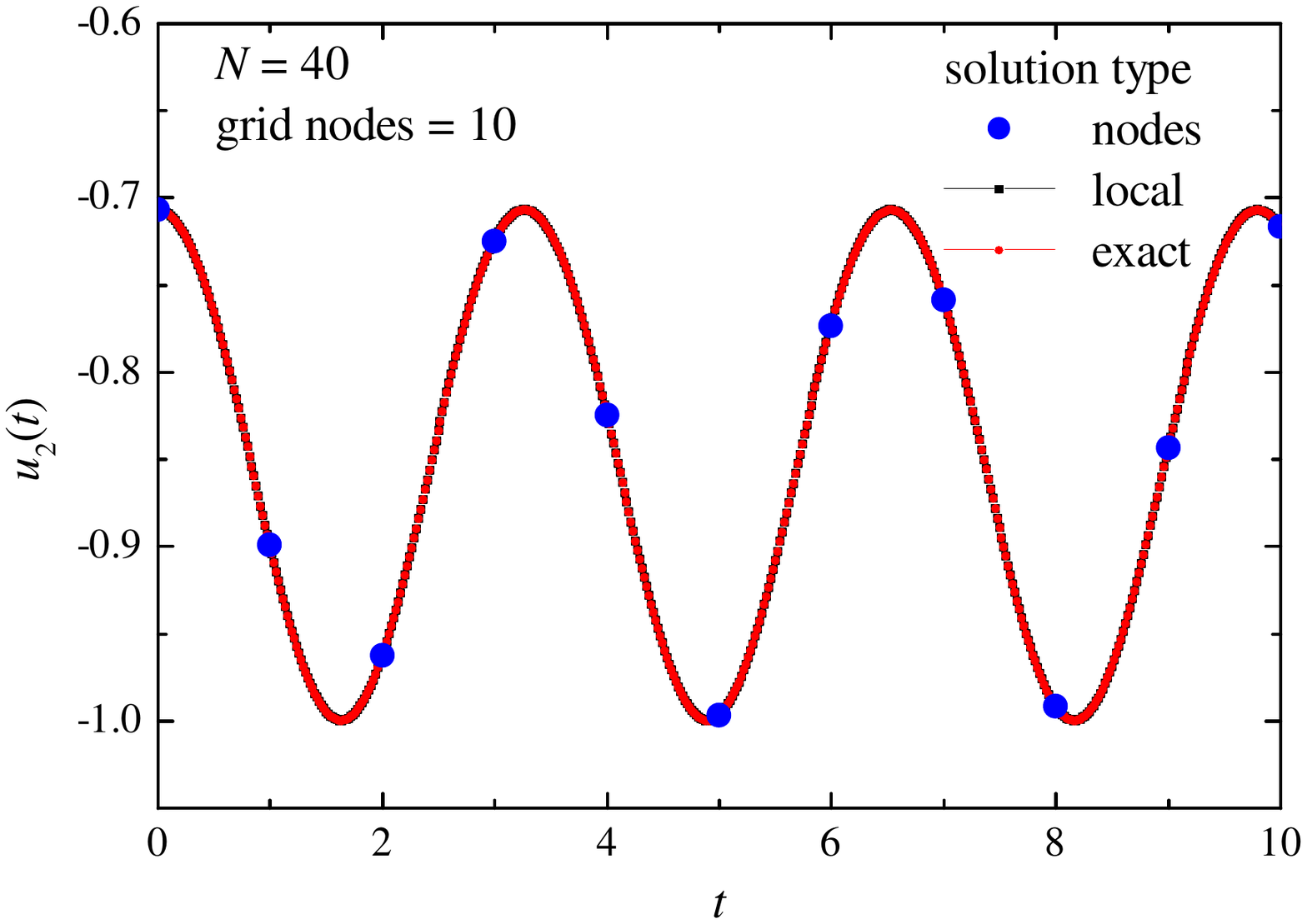}
\vspace{-8mm}\caption{\label{fig:pend_ind1_sol_uv:b3}}
\end{subfigure}\\
\begin{subfigure}{0.320\textwidth}
\includegraphics[width=\textwidth]{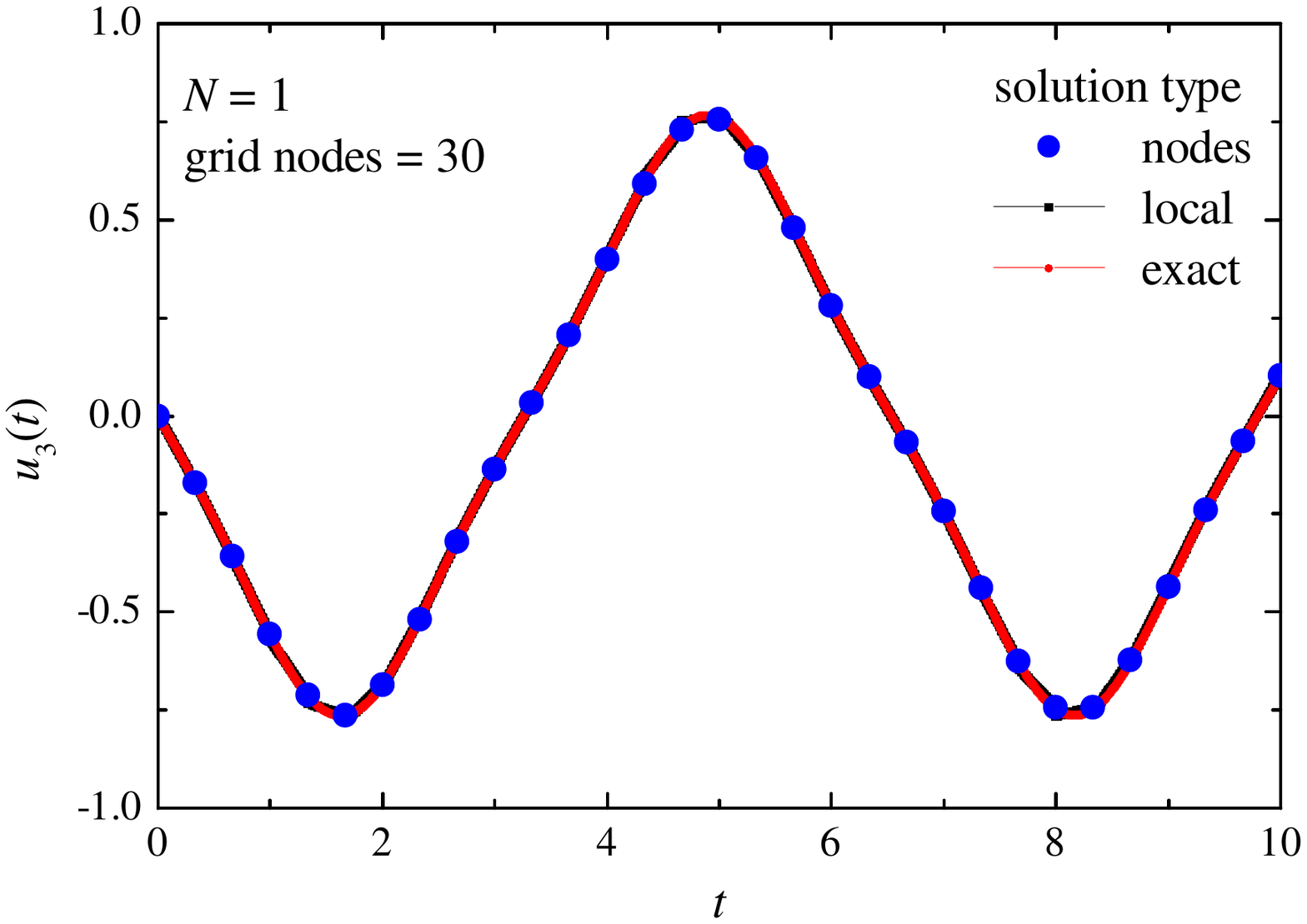}
\vspace{-8mm}\caption{\label{fig:pend_ind1_sol_uv:c1}}
\end{subfigure}
\begin{subfigure}{0.320\textwidth}
\includegraphics[width=\textwidth]{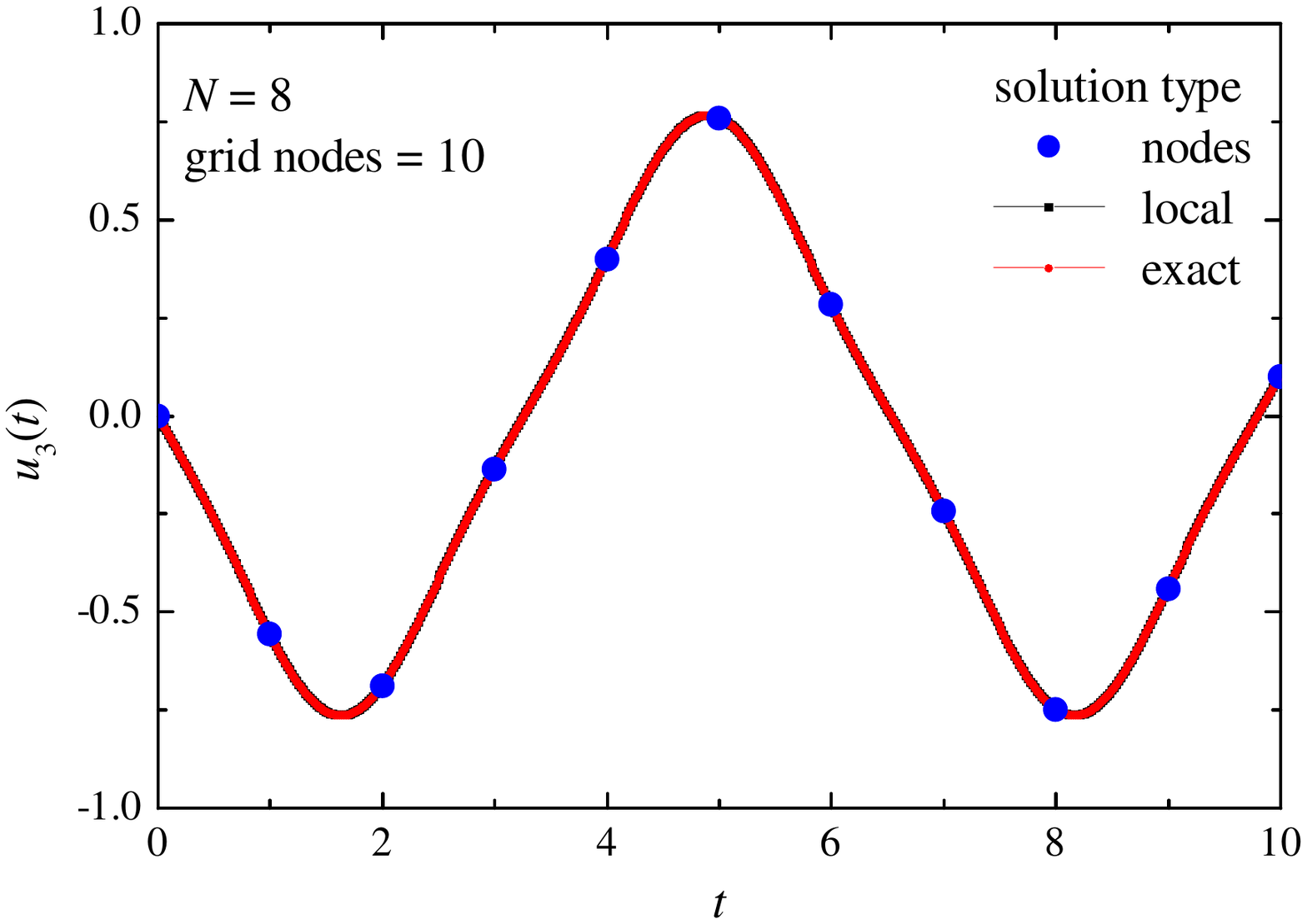}
\vspace{-8mm}\caption{\label{fig:pend_ind1_sol_uv:c2}}
\end{subfigure}
\begin{subfigure}{0.320\textwidth}
\includegraphics[width=\textwidth]{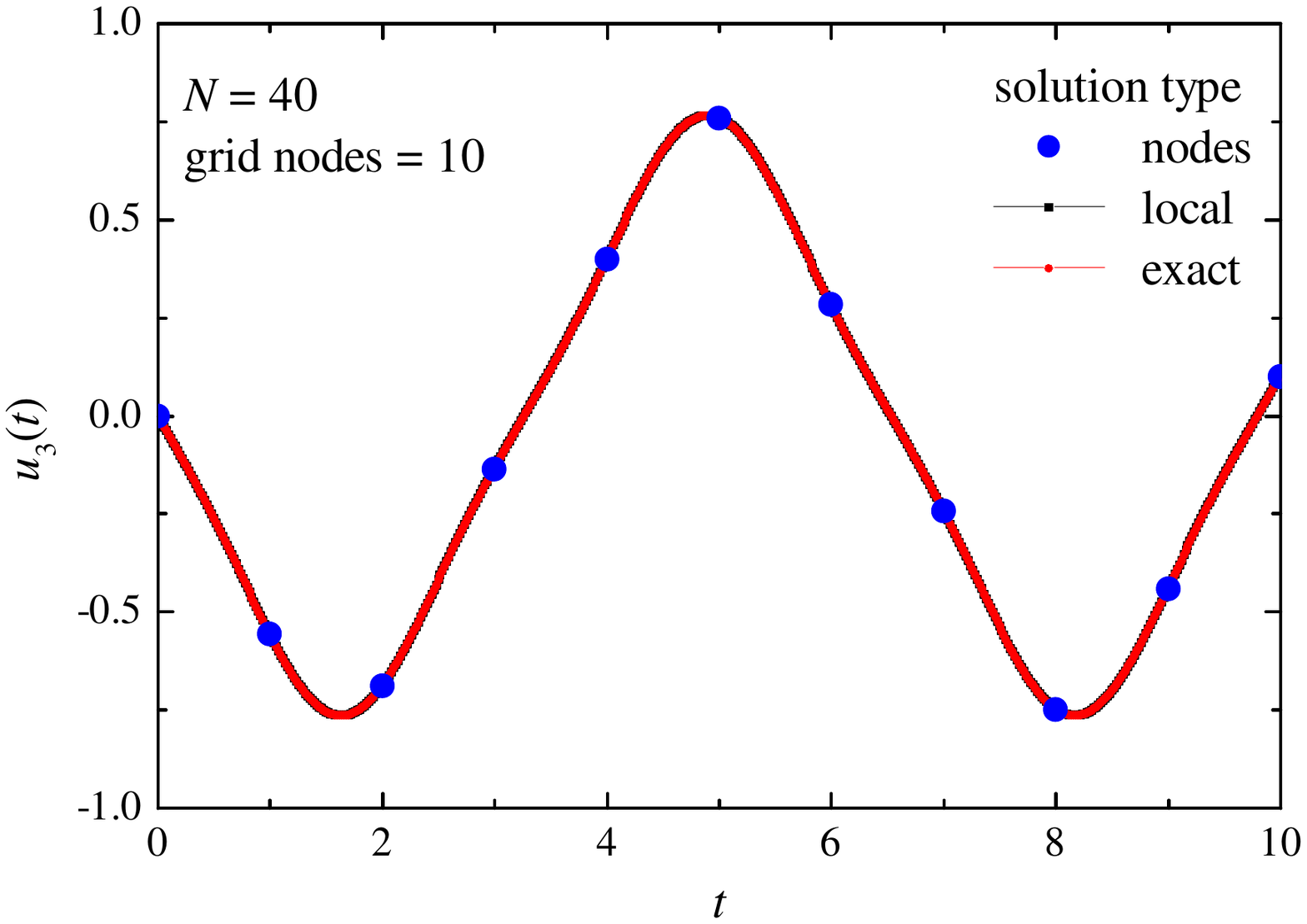}
\vspace{-8mm}\caption{\label{fig:pend_ind1_sol_uv:c3}}
\end{subfigure}\\
\begin{subfigure}{0.320\textwidth}
\includegraphics[width=\textwidth]{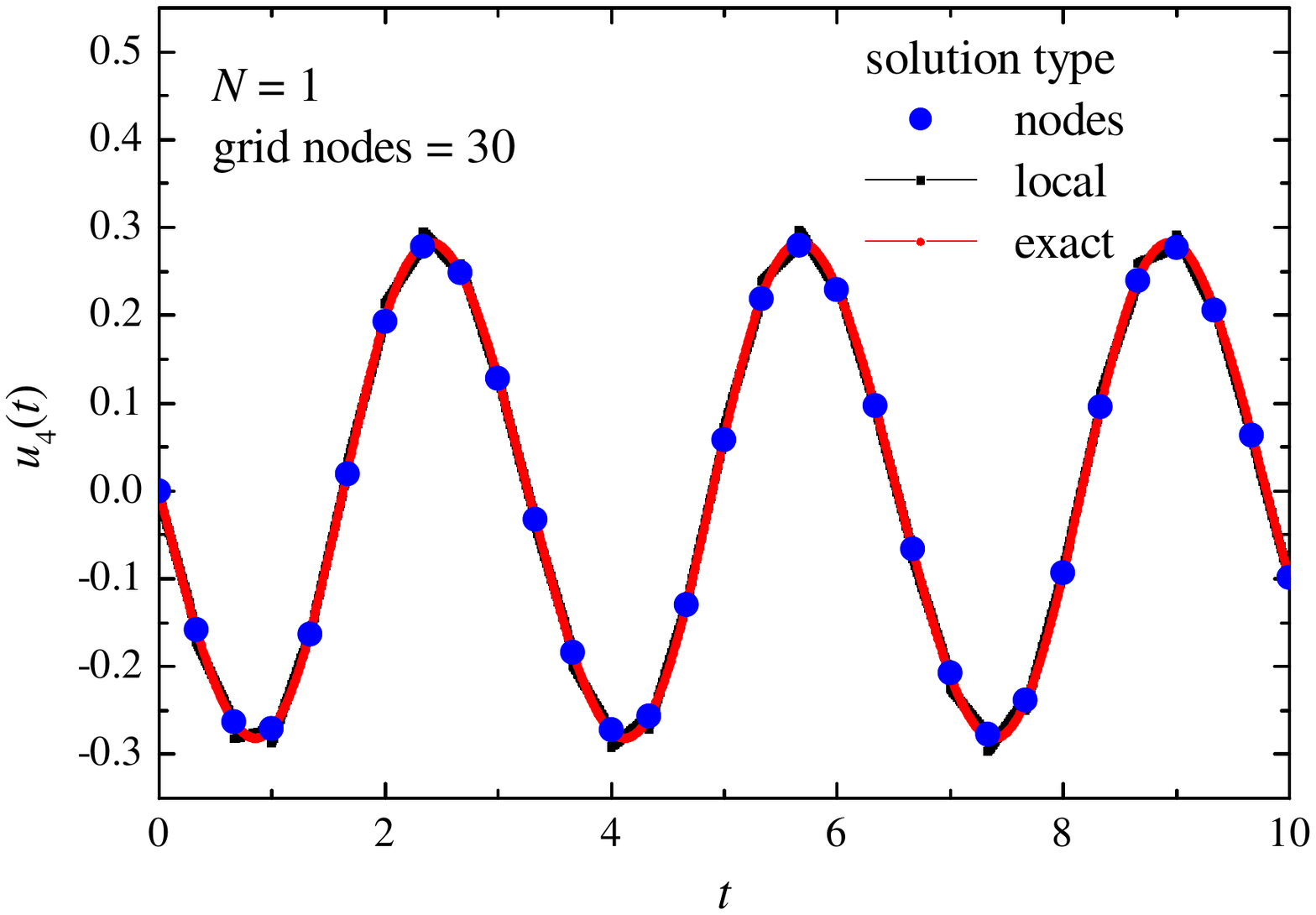}
\vspace{-8mm}\caption{\label{fig:pend_ind1_sol_uv:d1}}
\end{subfigure}
\begin{subfigure}{0.320\textwidth}
\includegraphics[width=\textwidth]{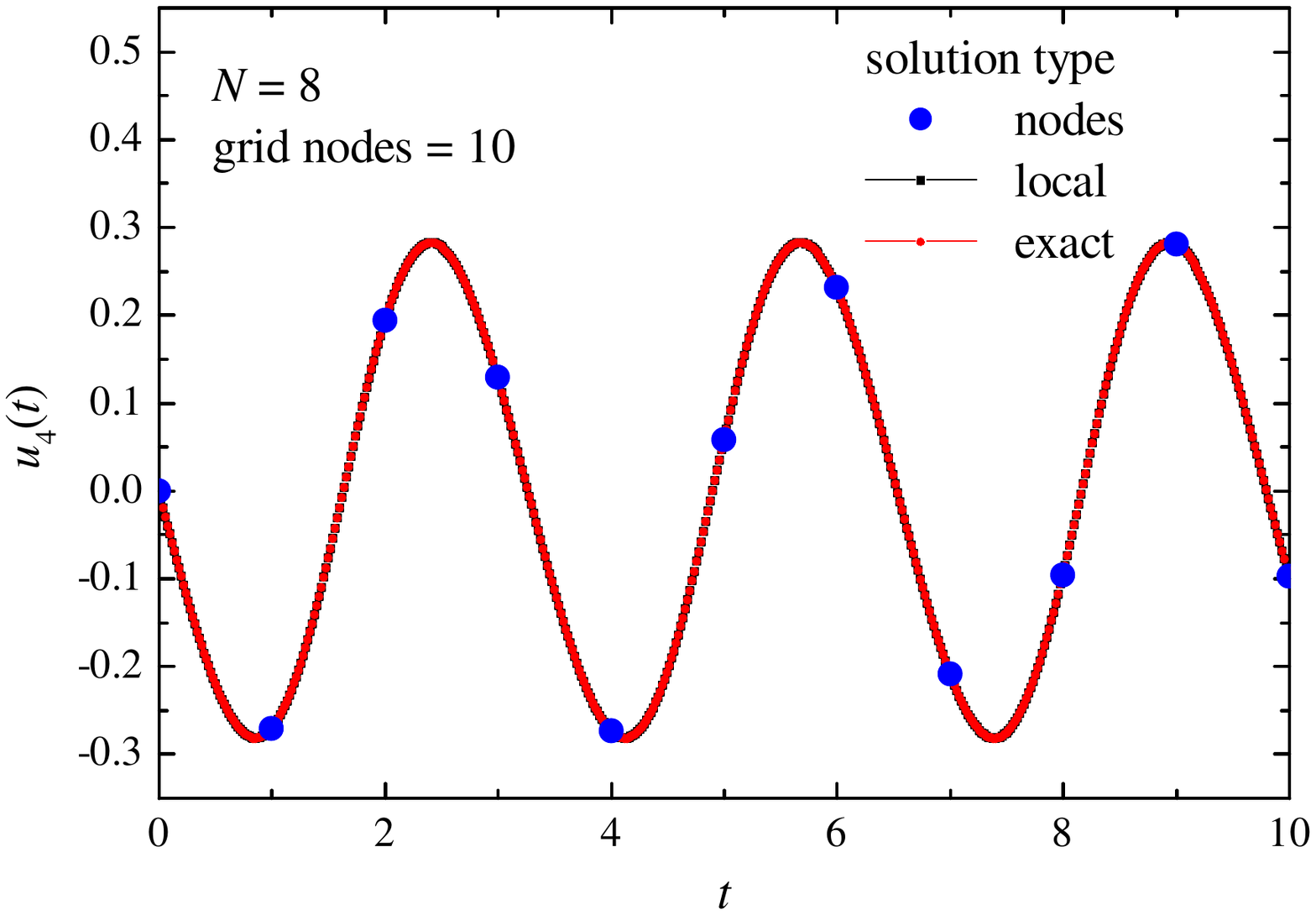}
\vspace{-8mm}\caption{\label{fig:pend_ind1_sol_uv:d2}}
\end{subfigure}
\begin{subfigure}{0.320\textwidth}
\includegraphics[width=\textwidth]{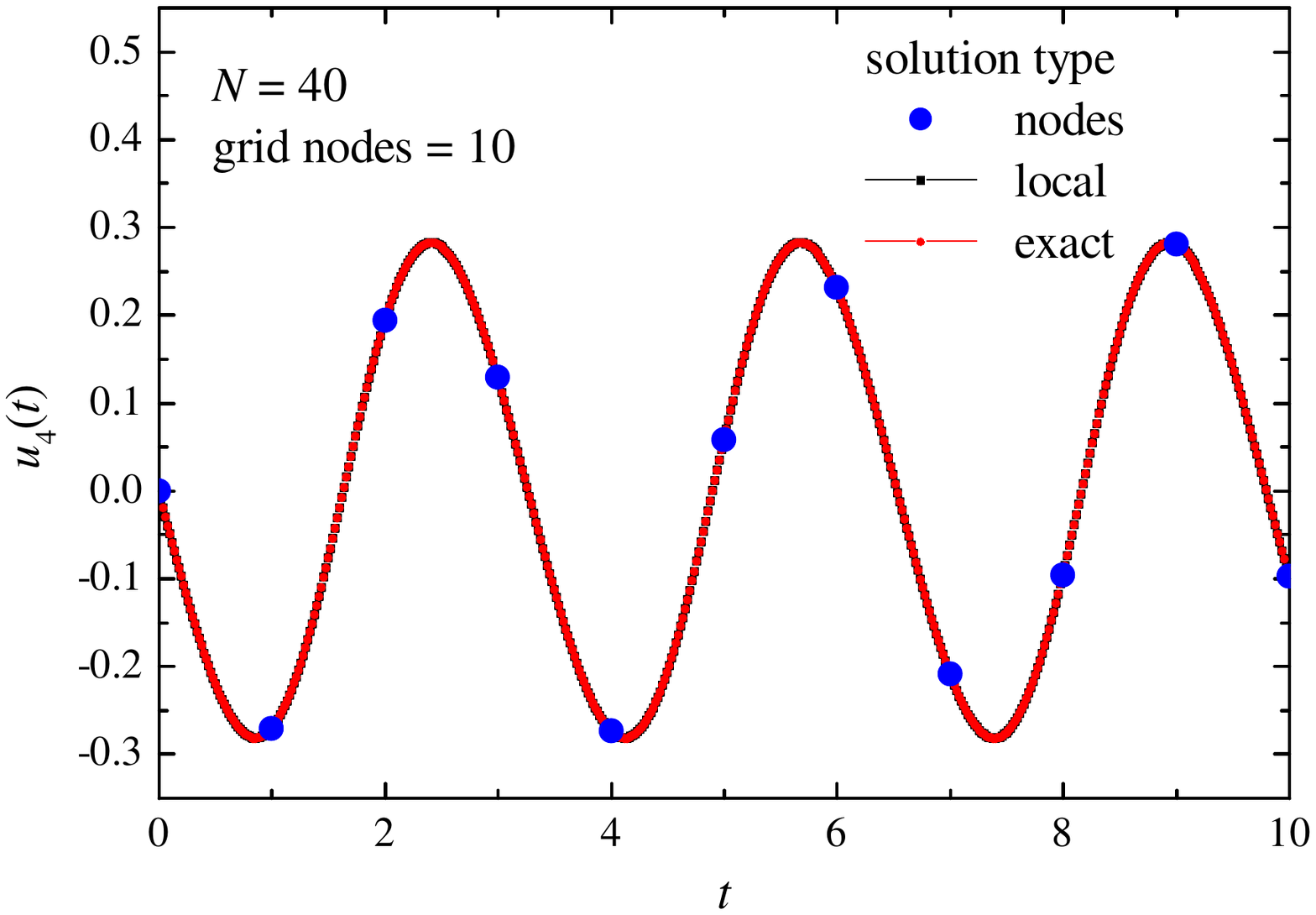}
\vspace{-8mm}\caption{\label{fig:pend_ind1_sol_uv:d3}}
\end{subfigure}\\
\begin{subfigure}{0.320\textwidth}
\includegraphics[width=\textwidth]{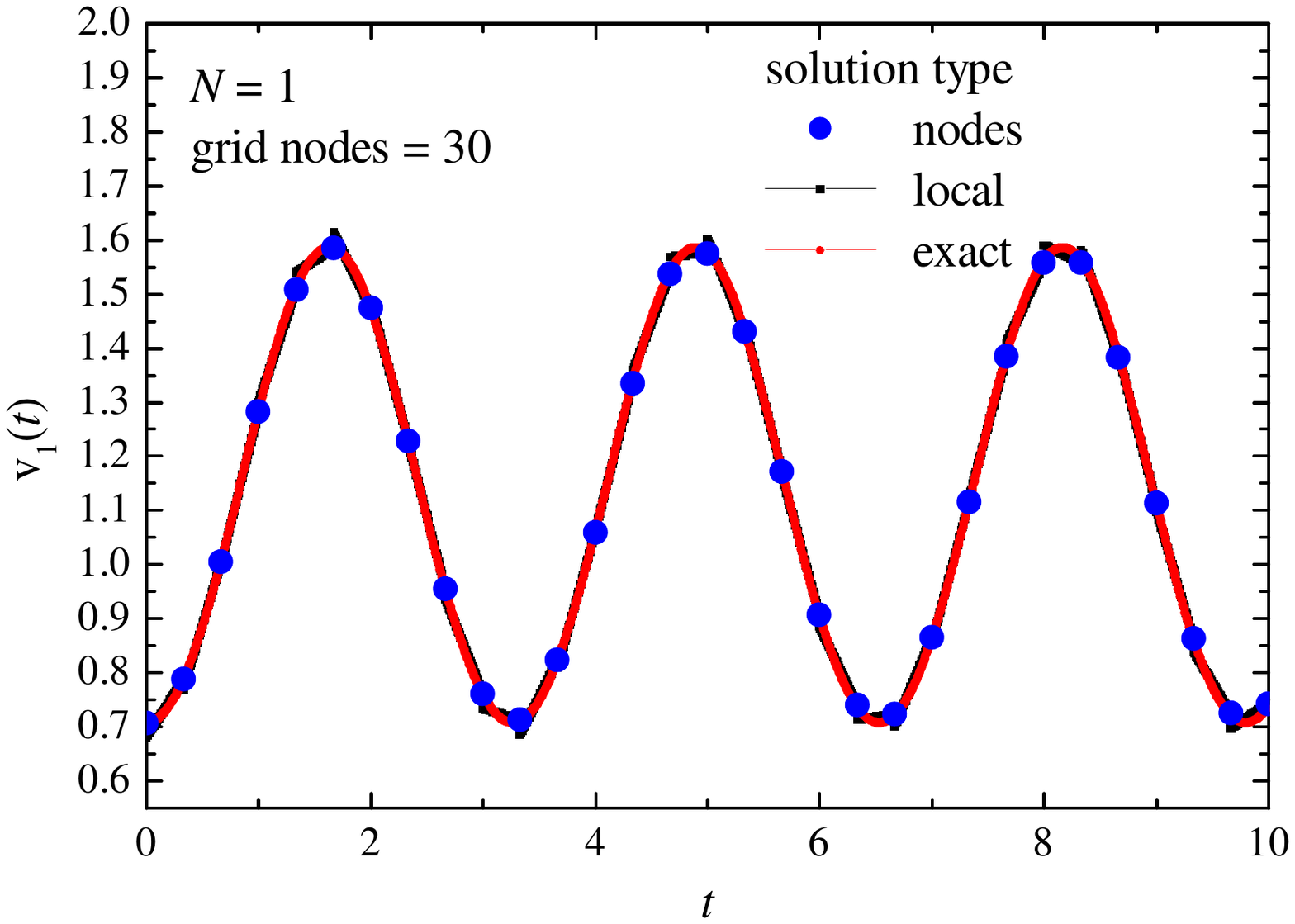}
\vspace{-8mm}\caption{\label{fig:pend_ind1_sol_uv:e1}}
\end{subfigure}
\begin{subfigure}{0.320\textwidth}
\includegraphics[width=\textwidth]{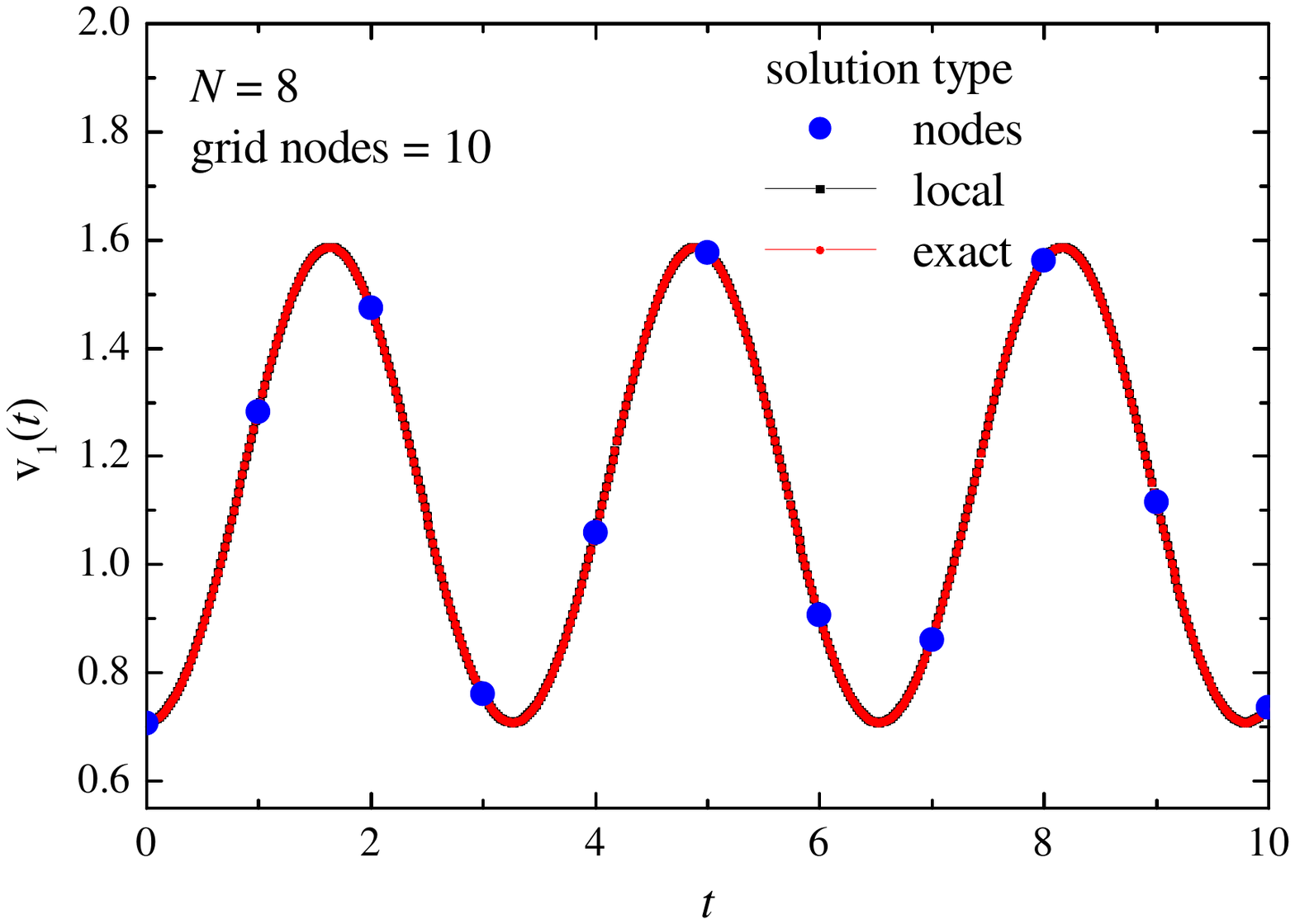}
\vspace{-8mm}\caption{\label{fig:pend_ind1_sol_uv:e2}}
\end{subfigure}
\begin{subfigure}{0.320\textwidth}
\includegraphics[width=\textwidth]{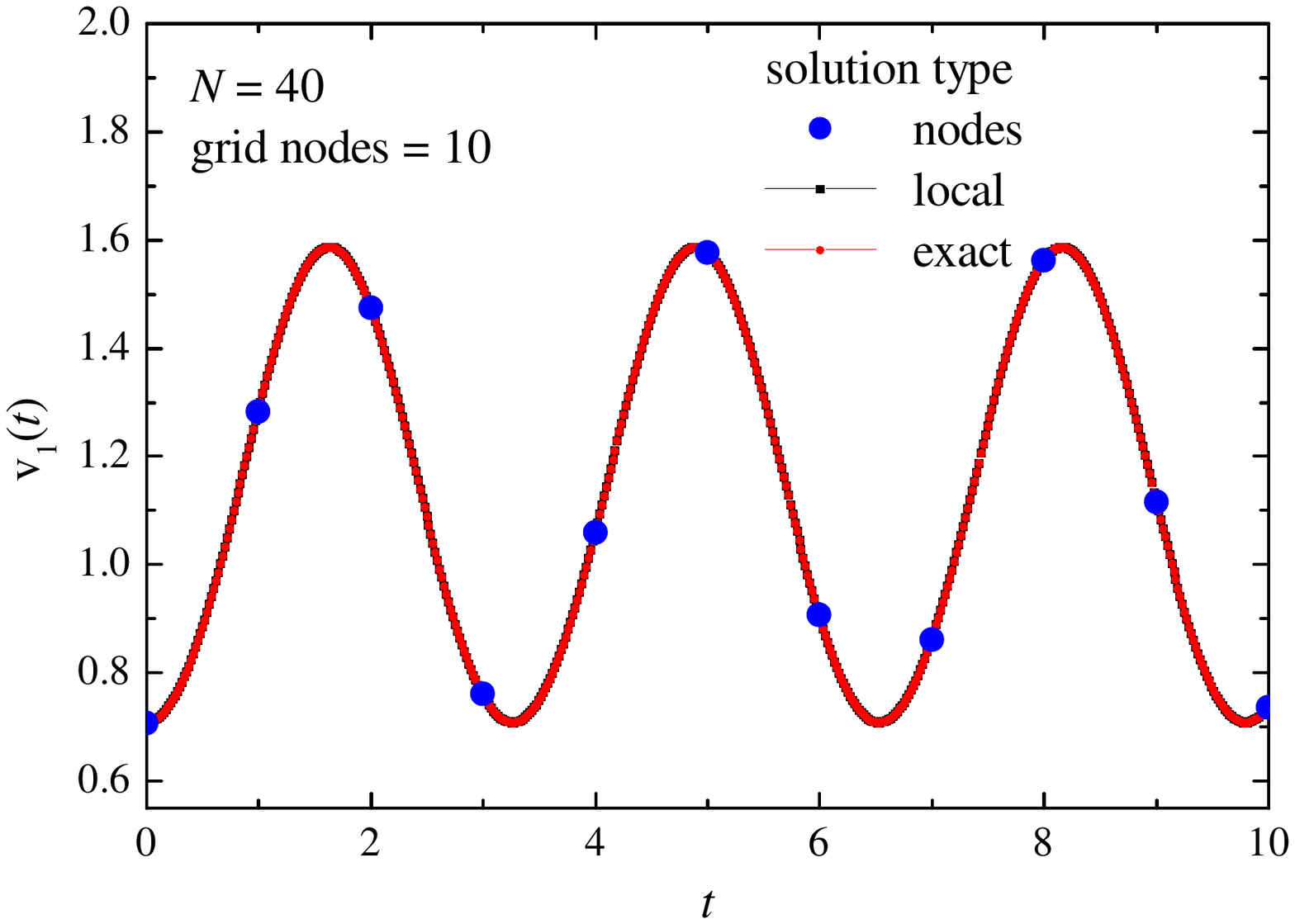}
\vspace{-8mm}\caption{\label{fig:pend_ind1_sol_uv:e3}}
\end{subfigure}\\
\caption{%
Numerical solution of the DAE system (\ref{eq:math_pend_dae_ind_3}) of index 1. Comparison of the solution at nodes $\mathbf{u}_{n}$, the local solution $\mathbf{u}_{L}(t)$ and the exact solution $\mathbf{u}^{\rm ex}(t)$ for components $u_{1}$ (\subref{fig:pend_ind1_sol_uv:a1}, \subref{fig:pend_ind1_sol_uv:a2}, \subref{fig:pend_ind1_sol_uv:a3}), $u_{2}$ (\subref{fig:pend_ind1_sol_uv:b1}, \subref{fig:pend_ind1_sol_uv:b2}, \subref{fig:pend_ind1_sol_uv:b3}), $u_{3}$ (\subref{fig:pend_ind1_sol_uv:c1}, \subref{fig:pend_ind1_sol_uv:c2}, \subref{fig:pend_ind1_sol_uv:c3}), $u_{4}$ (\subref{fig:pend_ind1_sol_uv:d1}, \subref{fig:pend_ind1_sol_uv:d2}, \subref{fig:pend_ind1_sol_uv:d3}) and $v_{1}$ (\subref{fig:pend_ind1_sol_uv:e1}, \subref{fig:pend_ind1_sol_uv:e2}, \subref{fig:pend_ind1_sol_uv:e3}), obtained using polynomials with degrees $N = 1$ (\subref{fig:pend_ind1_sol_uv:a1}, \subref{fig:pend_ind1_sol_uv:b1}, \subref{fig:pend_ind1_sol_uv:c1}, \subref{fig:pend_ind1_sol_uv:d1}, \subref{fig:pend_ind1_sol_uv:e1}), $N = 8$ (\subref{fig:pend_ind1_sol_uv:a2}, \subref{fig:pend_ind1_sol_uv:b2}, \subref{fig:pend_ind1_sol_uv:c2}, \subref{fig:pend_ind1_sol_uv:d2}, \subref{fig:pend_ind1_sol_uv:e2}) and $N = 40$ (\subref{fig:pend_ind1_sol_uv:a3}, \subref{fig:pend_ind1_sol_uv:b3}, \subref{fig:pend_ind1_sol_uv:c3}, \subref{fig:pend_ind1_sol_uv:d3}, \subref{fig:pend_ind1_sol_uv:e3}).
}
\label{fig:pend_ind1_sol_uv}
\end{figure} 

\begin{figure}[h!]
\captionsetup[subfigure]{%
	position=bottom,
	font+=smaller,
	textfont=normalfont,
	singlelinecheck=off,
	justification=raggedright
}
\centering
\begin{subfigure}{0.275\textwidth}
\includegraphics[width=\textwidth]{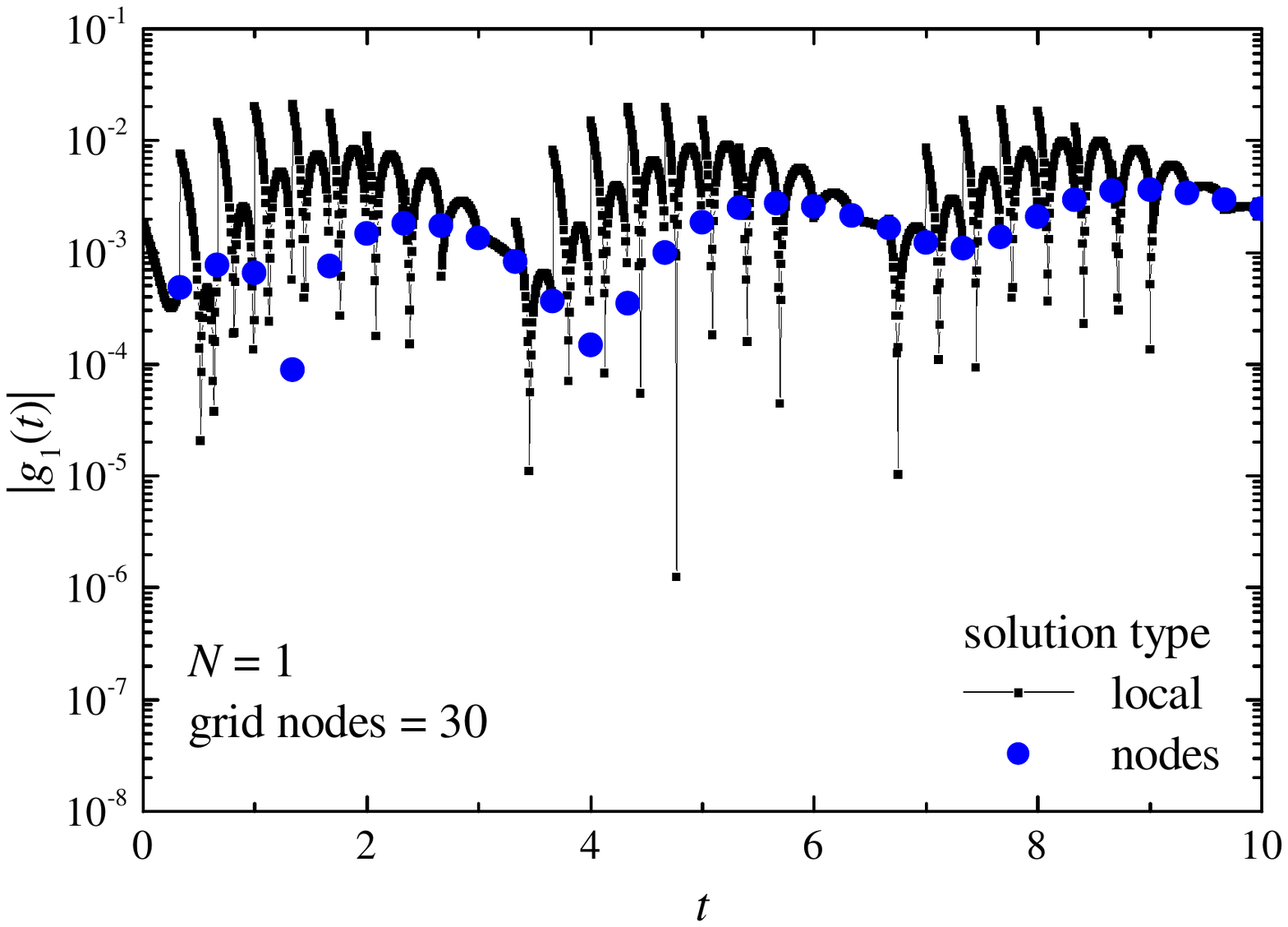}
\vspace{-8mm}\caption{\label{fig:pend_ind1_sol_g_eps:a1}}
\end{subfigure}\hspace{6mm}
\begin{subfigure}{0.275\textwidth}
\includegraphics[width=\textwidth]{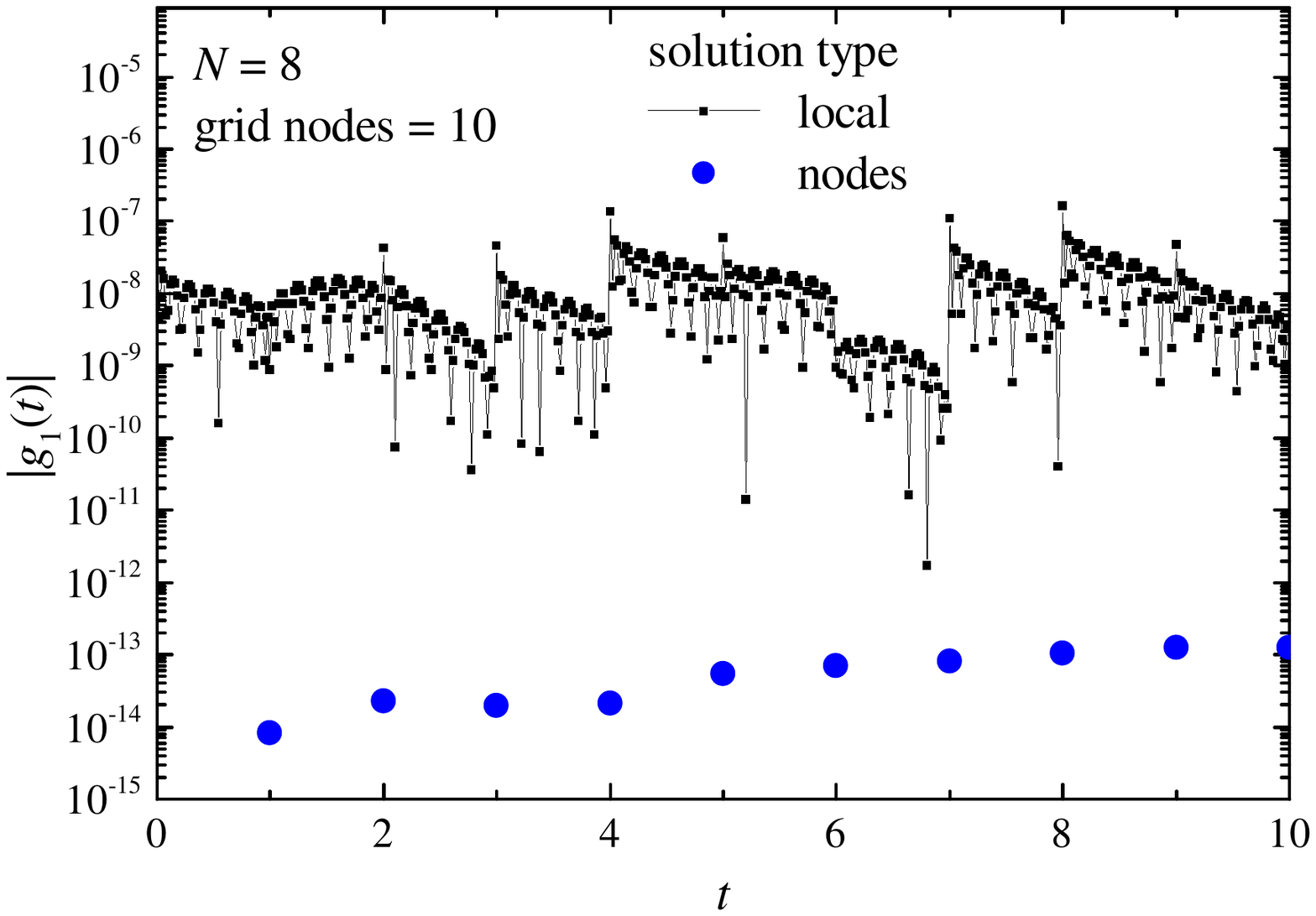}
\vspace{-8mm}\caption{\label{fig:pend_ind1_sol_g_eps:a2}}
\end{subfigure}\hspace{6mm}
\begin{subfigure}{0.275\textwidth}
\includegraphics[width=\textwidth]{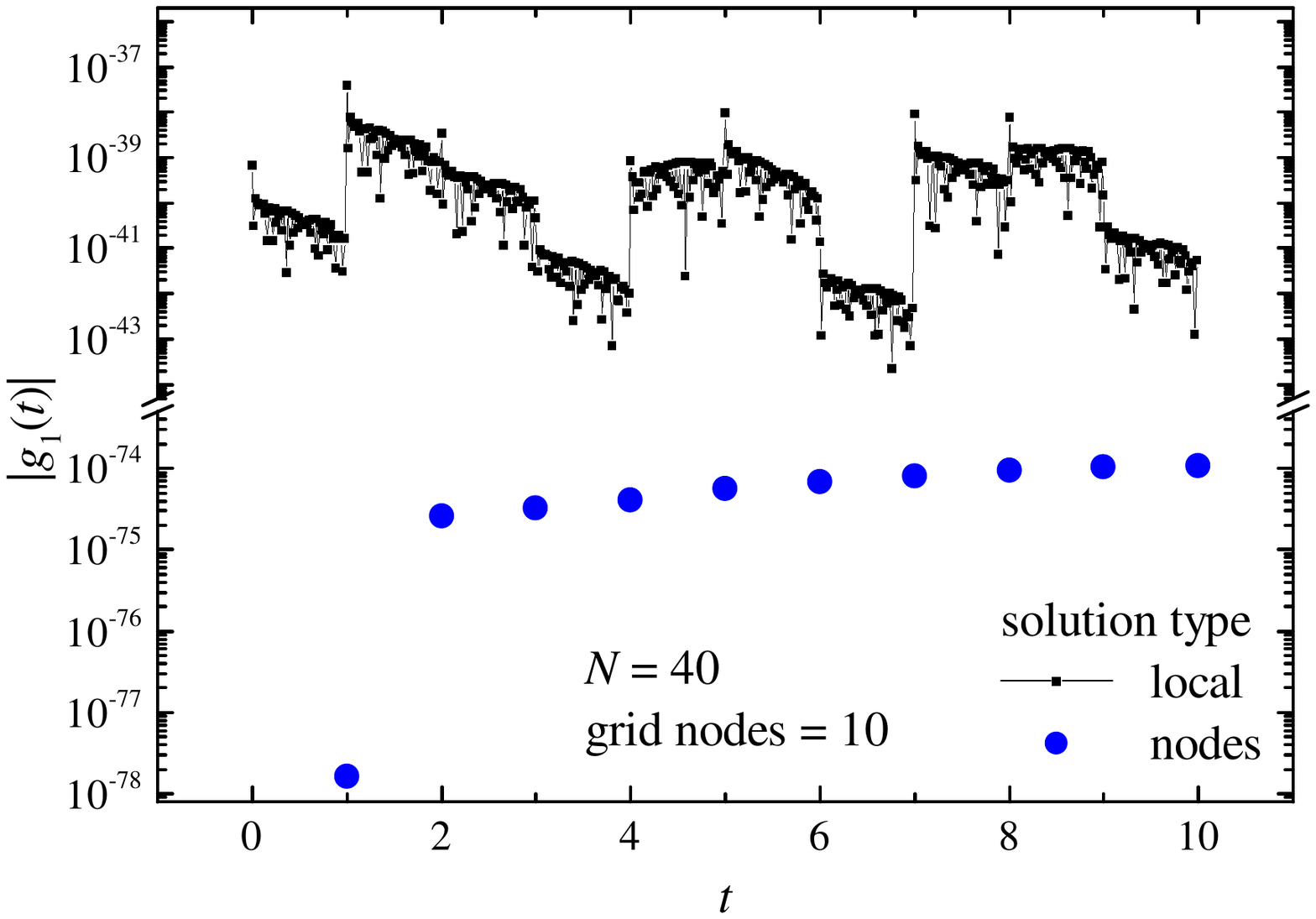}
\vspace{-8mm}\caption{\label{fig:pend_ind1_sol_g_eps:a3}}
\end{subfigure}\\[-2mm]
\begin{subfigure}{0.275\textwidth}
\includegraphics[width=\textwidth]{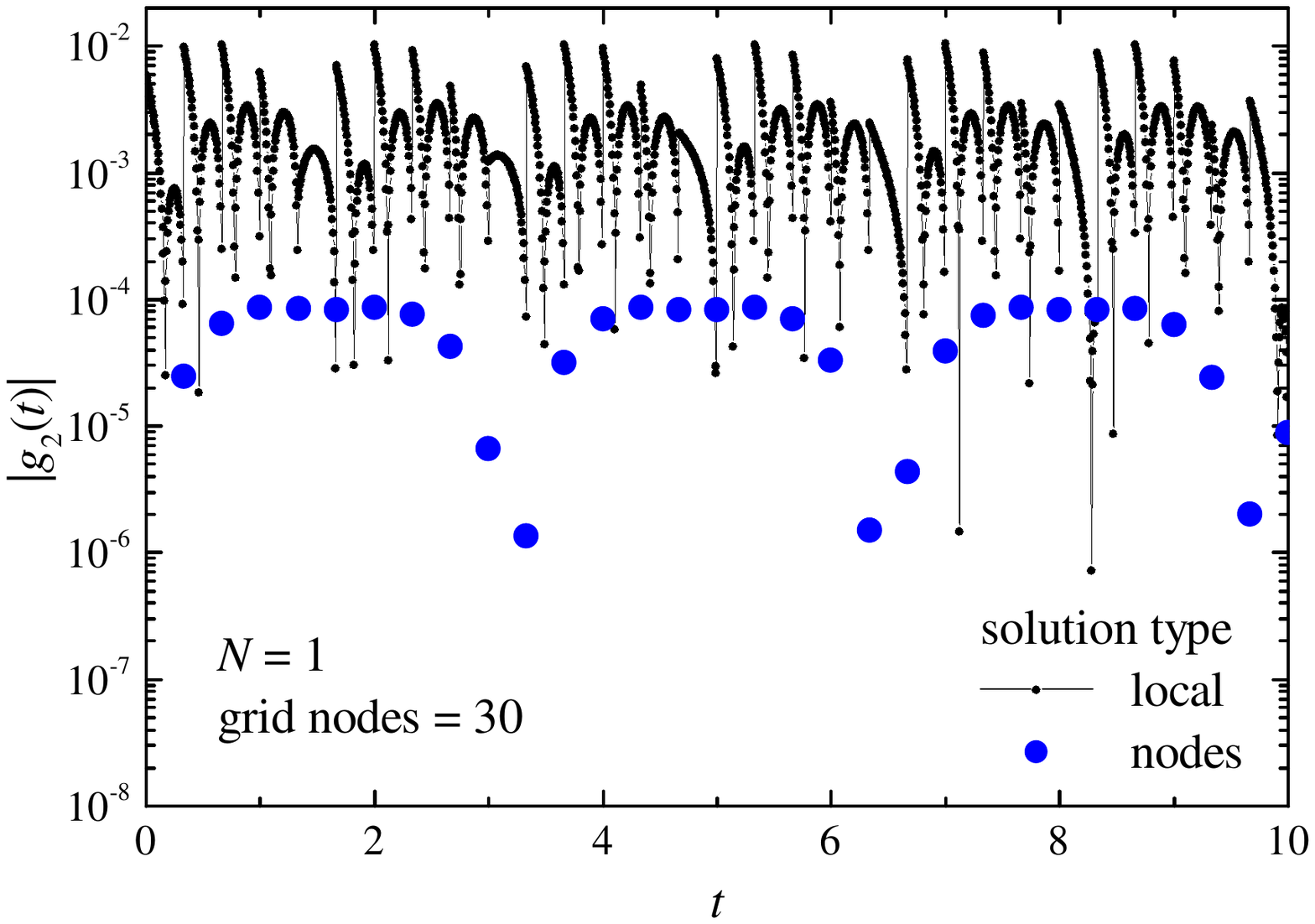}
\vspace{-8mm}\caption{\label{fig:pend_ind1_sol_g_eps:b1}}
\end{subfigure}\hspace{6mm}
\begin{subfigure}{0.275\textwidth}
\includegraphics[width=\textwidth]{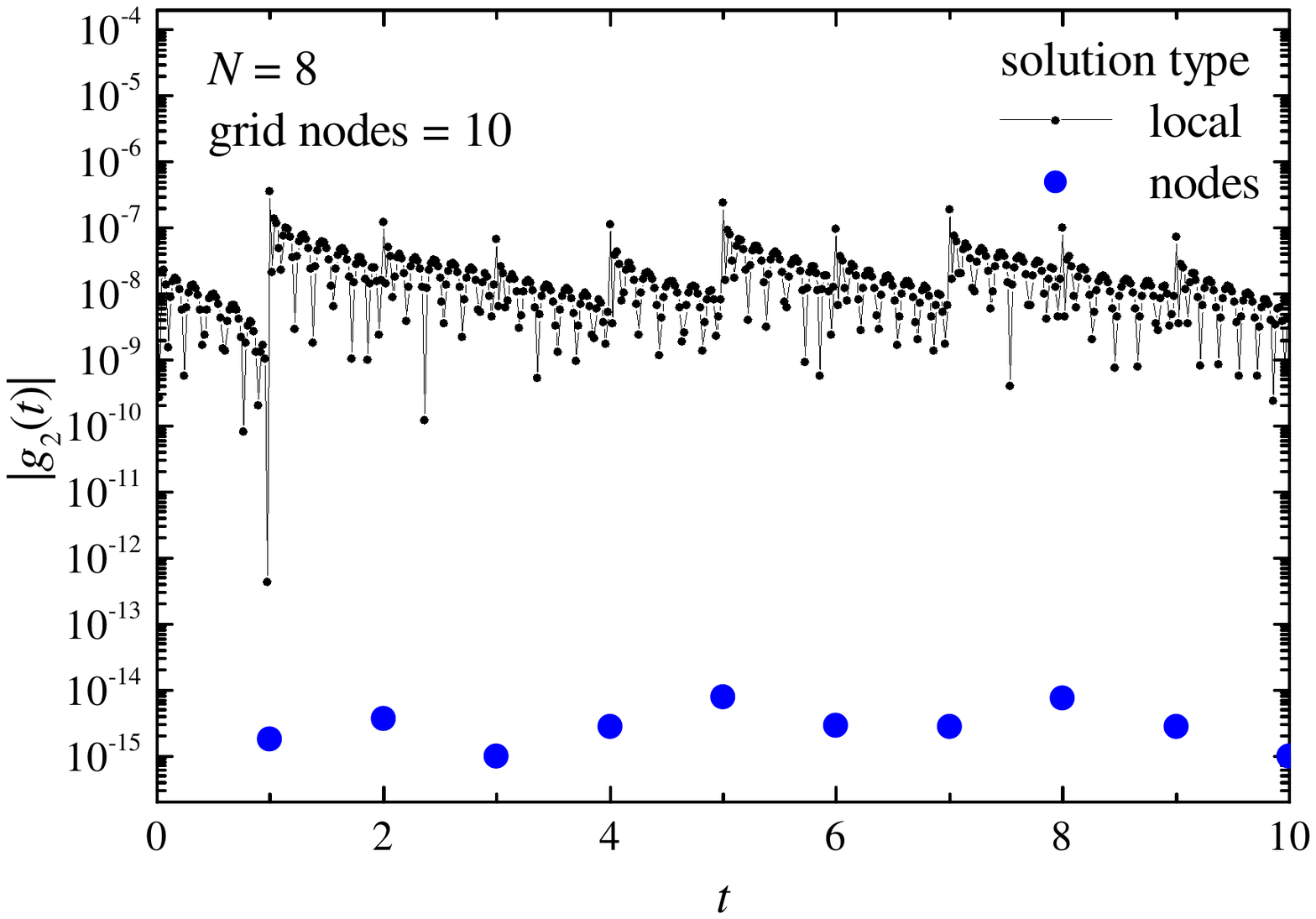}
\vspace{-8mm}\caption{\label{fig:pend_ind1_sol_g_eps:b2}}
\end{subfigure}\hspace{6mm}
\begin{subfigure}{0.275\textwidth}
\includegraphics[width=\textwidth]{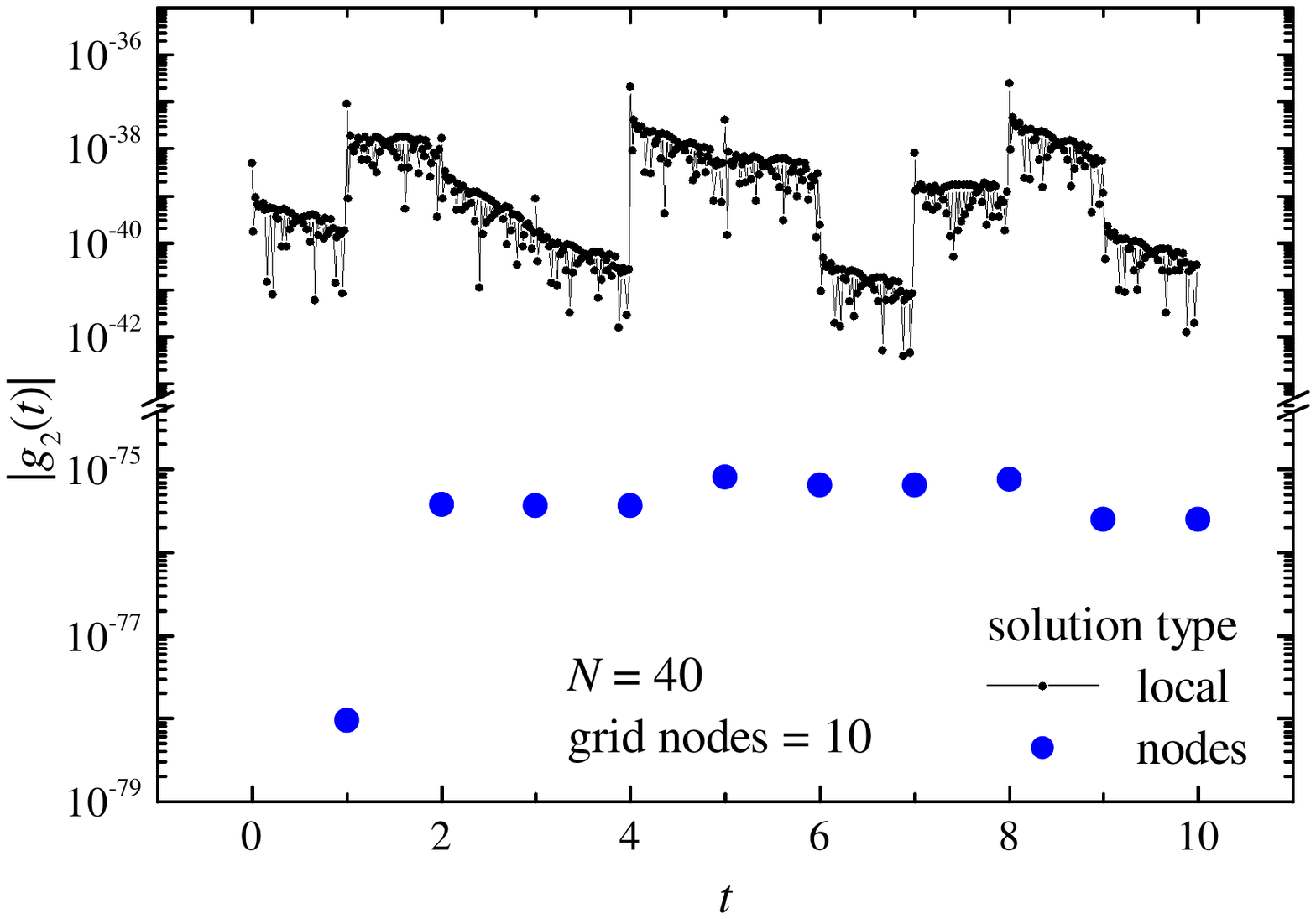}
\vspace{-8mm}\caption{\label{fig:pend_ind1_sol_g_eps:b3}}
\end{subfigure}\\[-2mm]
\begin{subfigure}{0.275\textwidth}
\includegraphics[width=\textwidth]{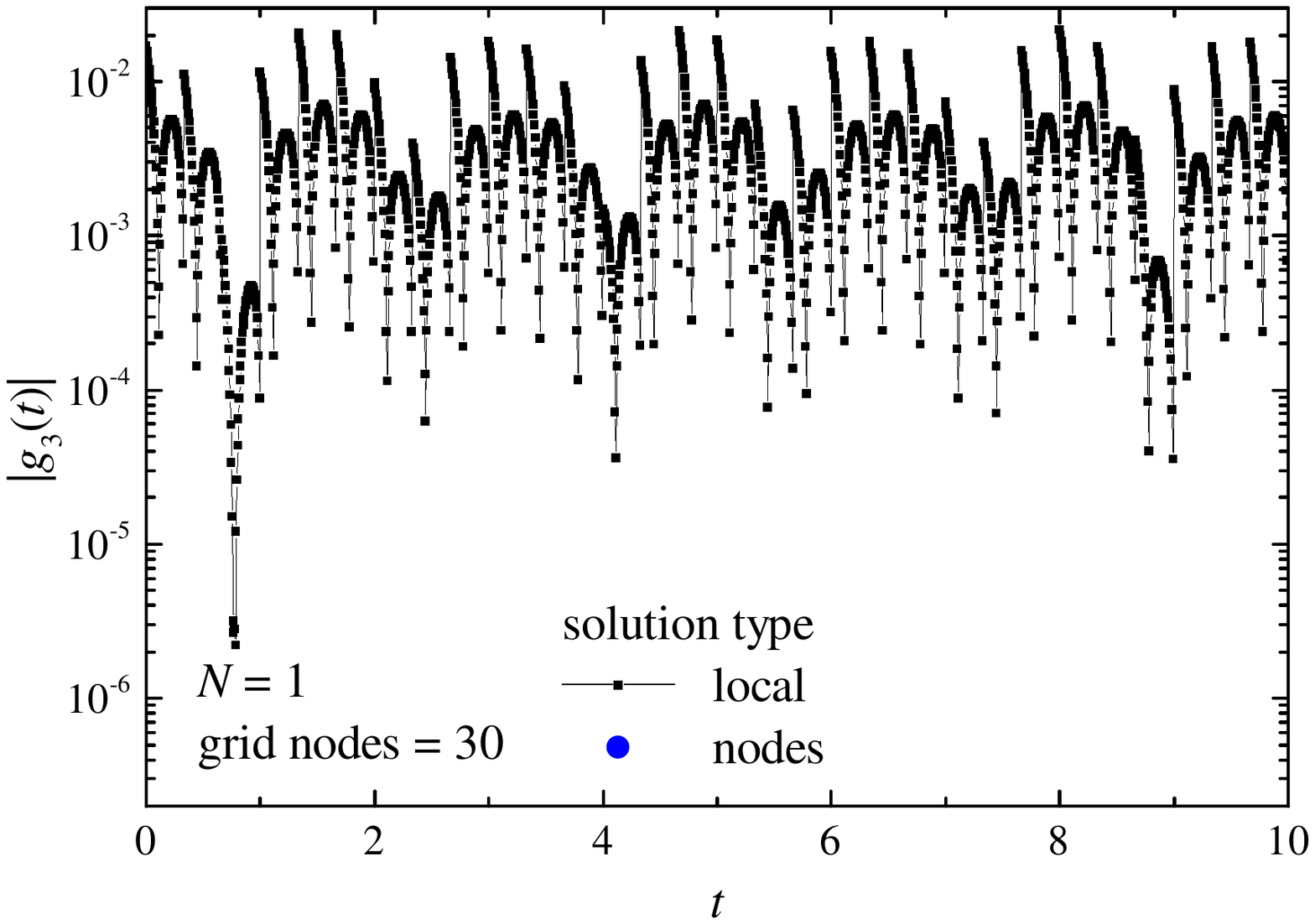}
\vspace{-8mm}\caption{\label{fig:pend_ind1_sol_g_eps:c1}}
\end{subfigure}\hspace{6mm}
\begin{subfigure}{0.275\textwidth}
\includegraphics[width=\textwidth]{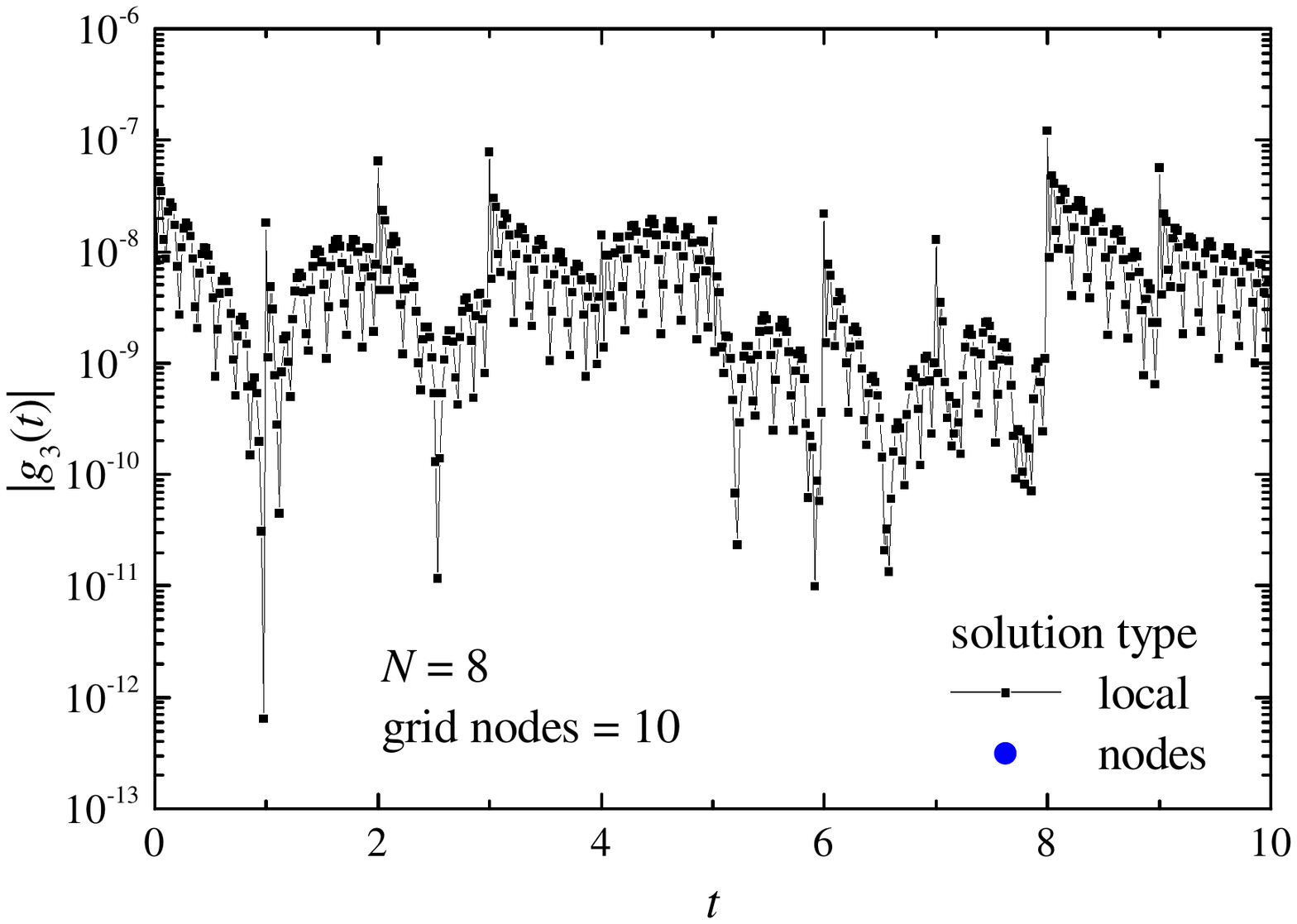}
\vspace{-8mm}\caption{\label{fig:pend_ind1_sol_g_eps:c2}}
\end{subfigure}\hspace{6mm}
\begin{subfigure}{0.275\textwidth}
\includegraphics[width=\textwidth]{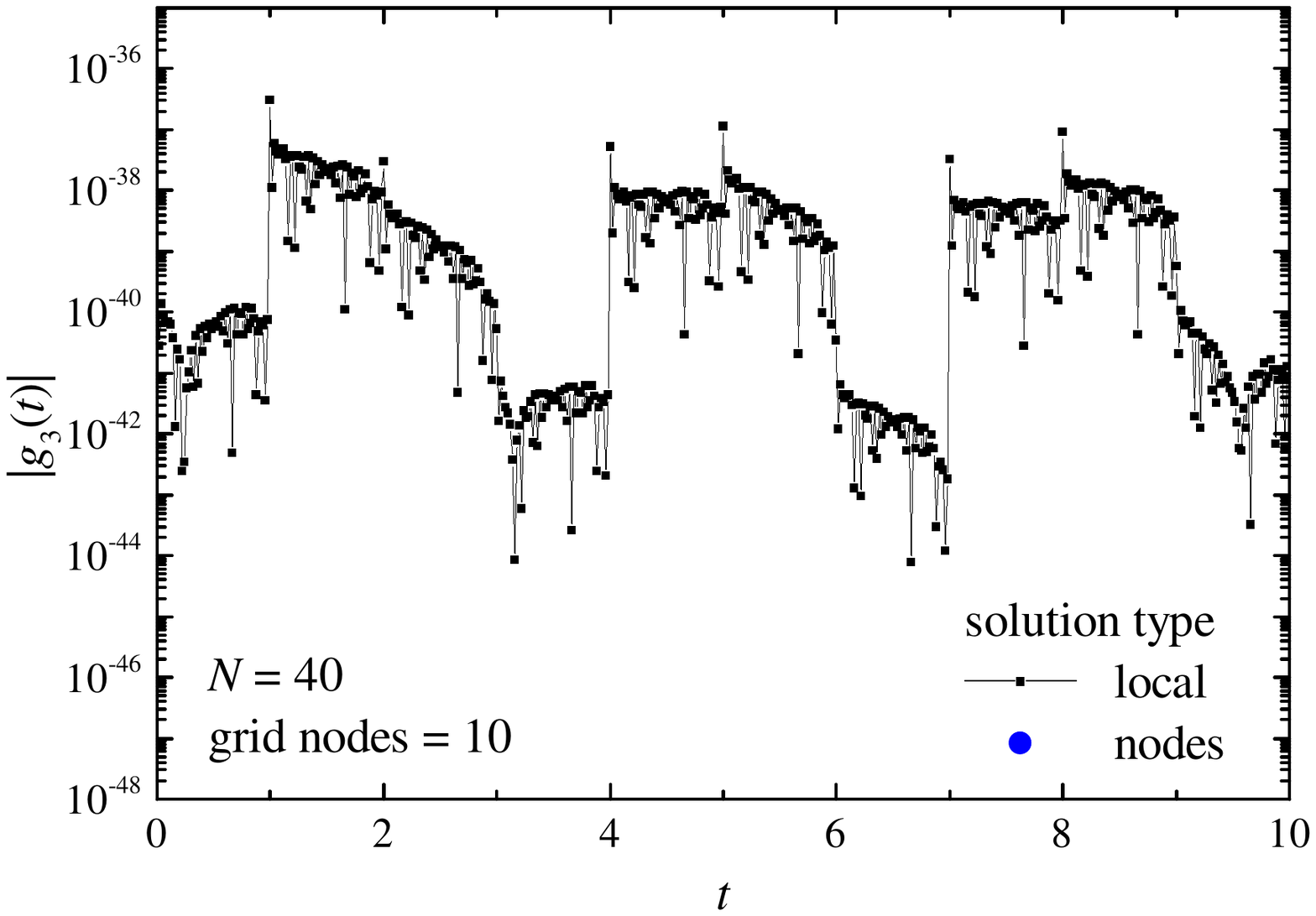}
\vspace{-8mm}\caption{\label{fig:pend_ind1_sol_g_eps:c3}}
\end{subfigure}\\[-2mm]
\begin{subfigure}{0.275\textwidth}
\includegraphics[width=\textwidth]{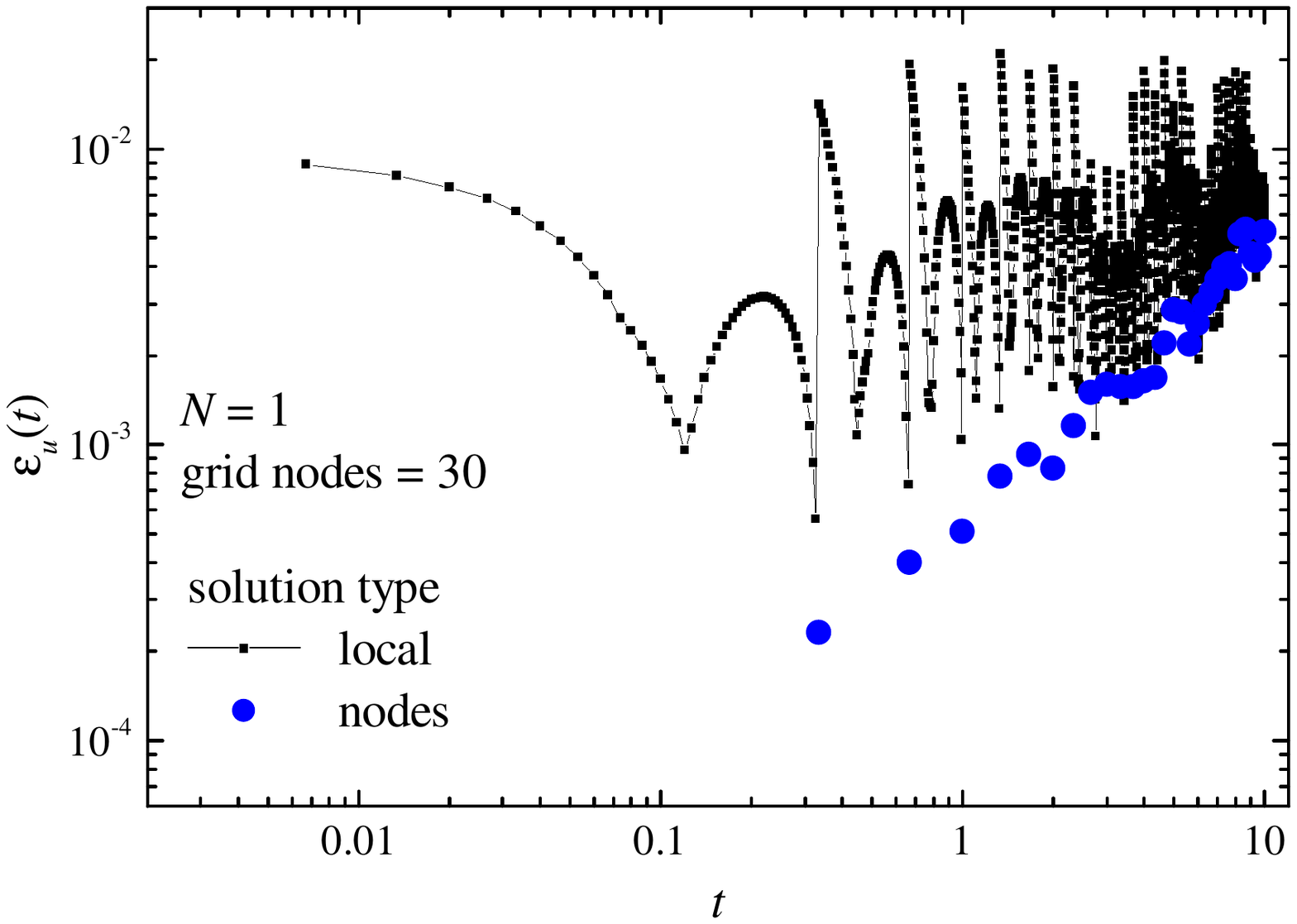}
\vspace{-8mm}\caption{\label{fig:pend_ind1_sol_g_eps:d1}}
\end{subfigure}\hspace{6mm}
\begin{subfigure}{0.275\textwidth}
\includegraphics[width=\textwidth]{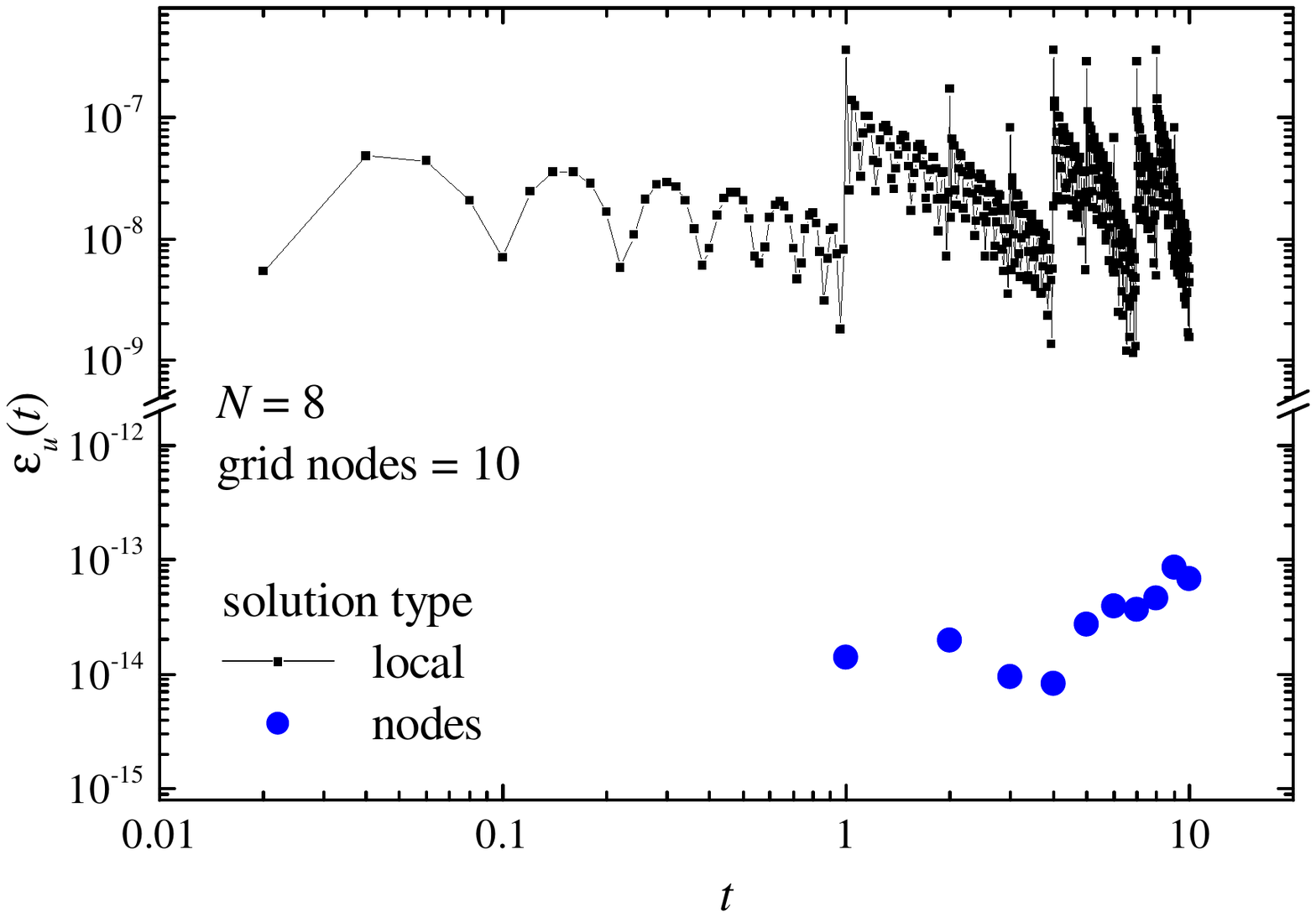}
\vspace{-8mm}\caption{\label{fig:pend_ind1_sol_g_eps:d2}}
\end{subfigure}\hspace{6mm}
\begin{subfigure}{0.275\textwidth}
\includegraphics[width=\textwidth]{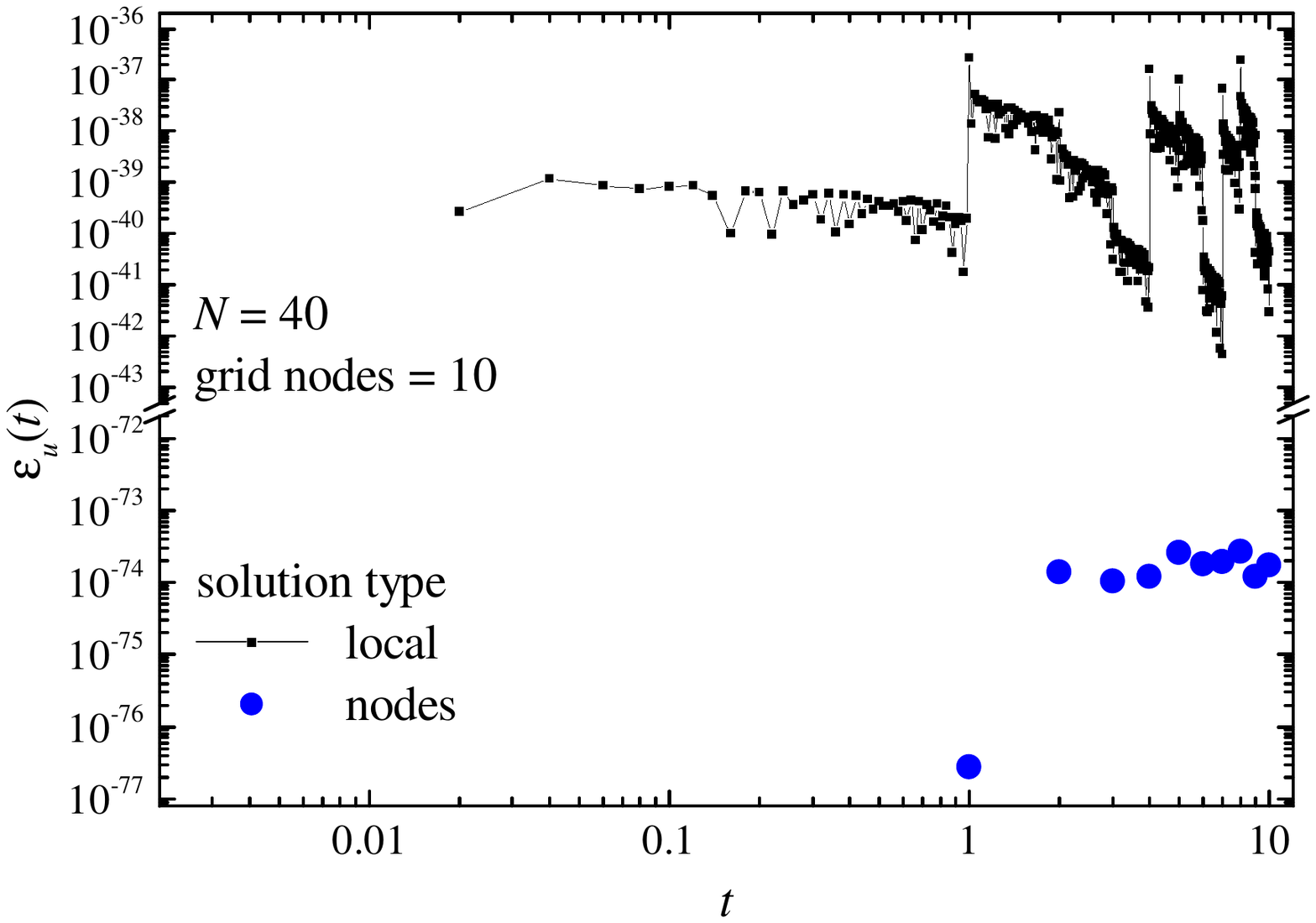}
\vspace{-8mm}\caption{\label{fig:pend_ind1_sol_g_eps:d3}}
\end{subfigure}\\[-2mm]
\begin{subfigure}{0.275\textwidth}
\includegraphics[width=\textwidth]{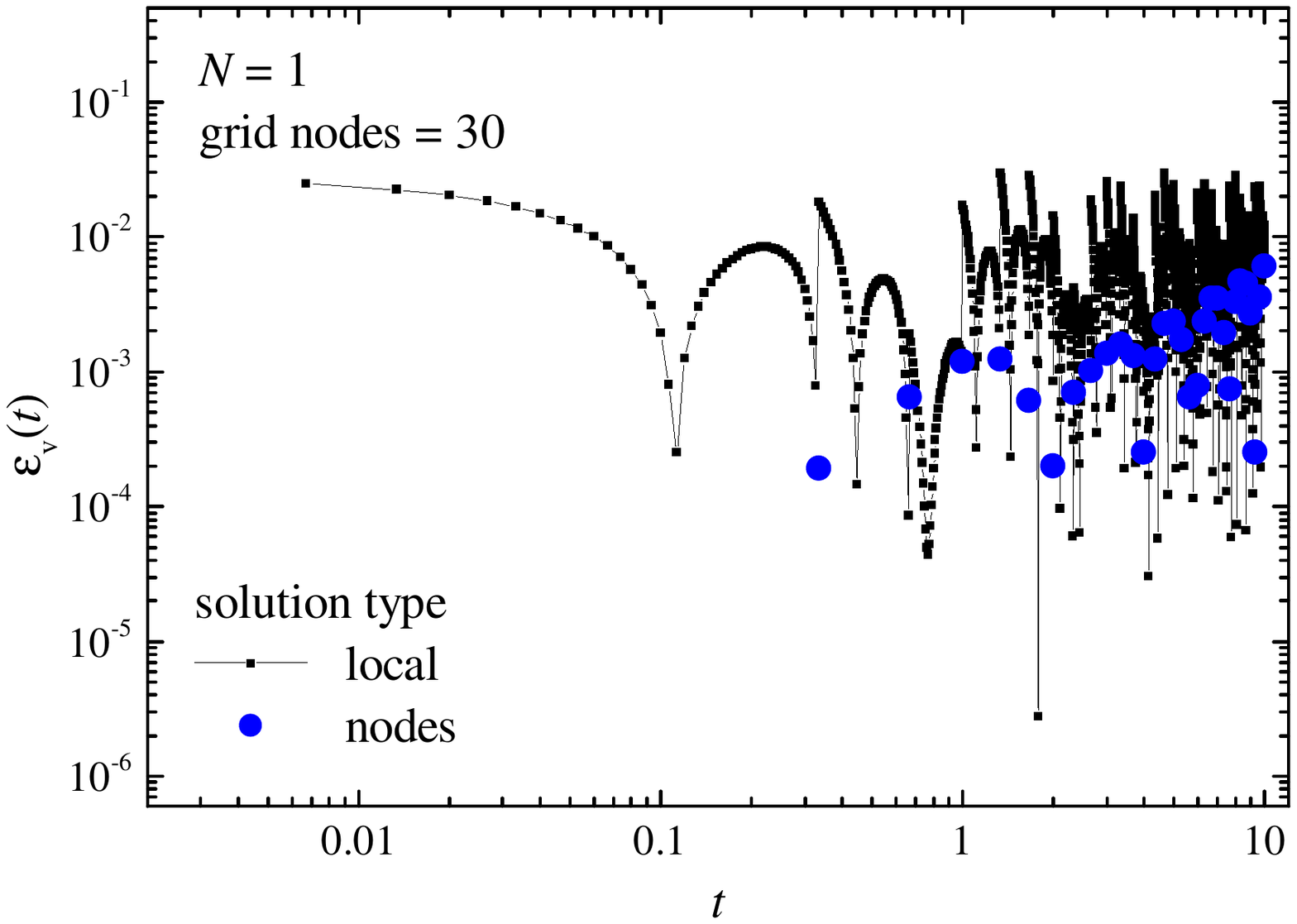}
\vspace{-8mm}\caption{\label{fig:pend_ind1_sol_g_eps:e1}}
\end{subfigure}\hspace{6mm}
\begin{subfigure}{0.275\textwidth}
\includegraphics[width=\textwidth]{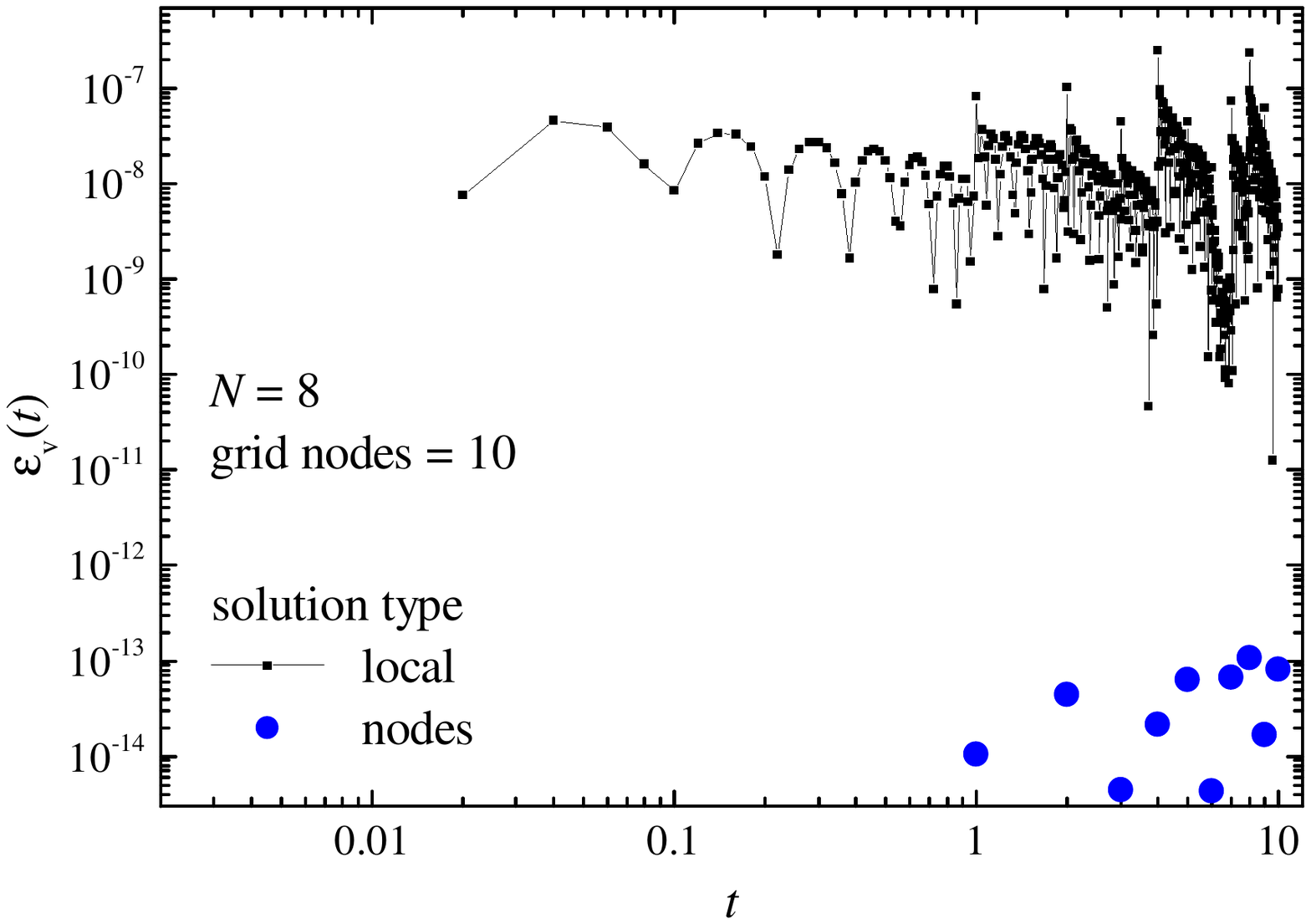}
\vspace{-8mm}\caption{\label{fig:pend_ind1_sol_g_eps:e2}}
\end{subfigure}\hspace{6mm}
\begin{subfigure}{0.275\textwidth}
\includegraphics[width=\textwidth]{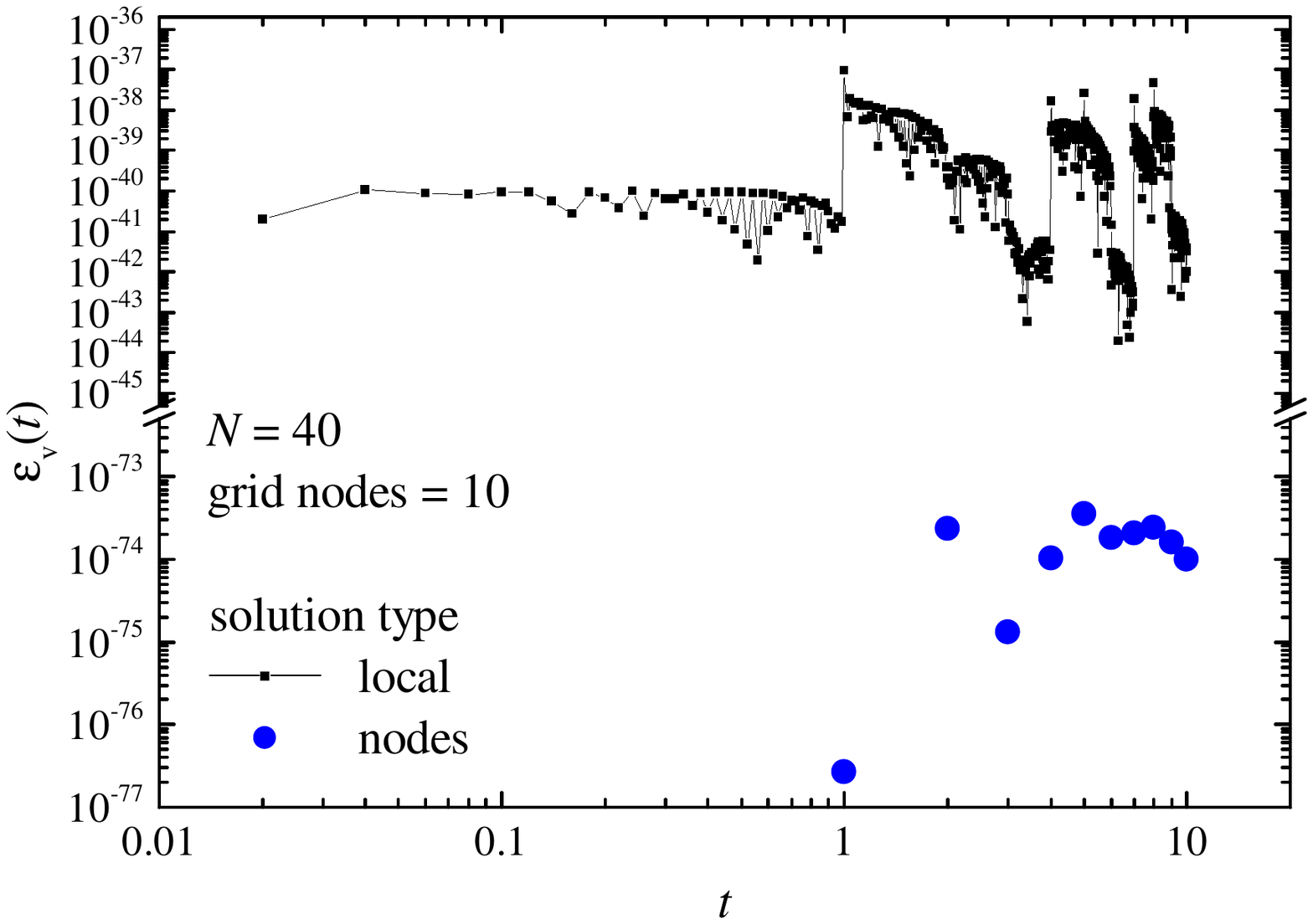}
\vspace{-8mm}\caption{\label{fig:pend_ind1_sol_g_eps:e3}}
\end{subfigure}\\[-2mm]
\begin{subfigure}{0.275\textwidth}
\includegraphics[width=\textwidth]{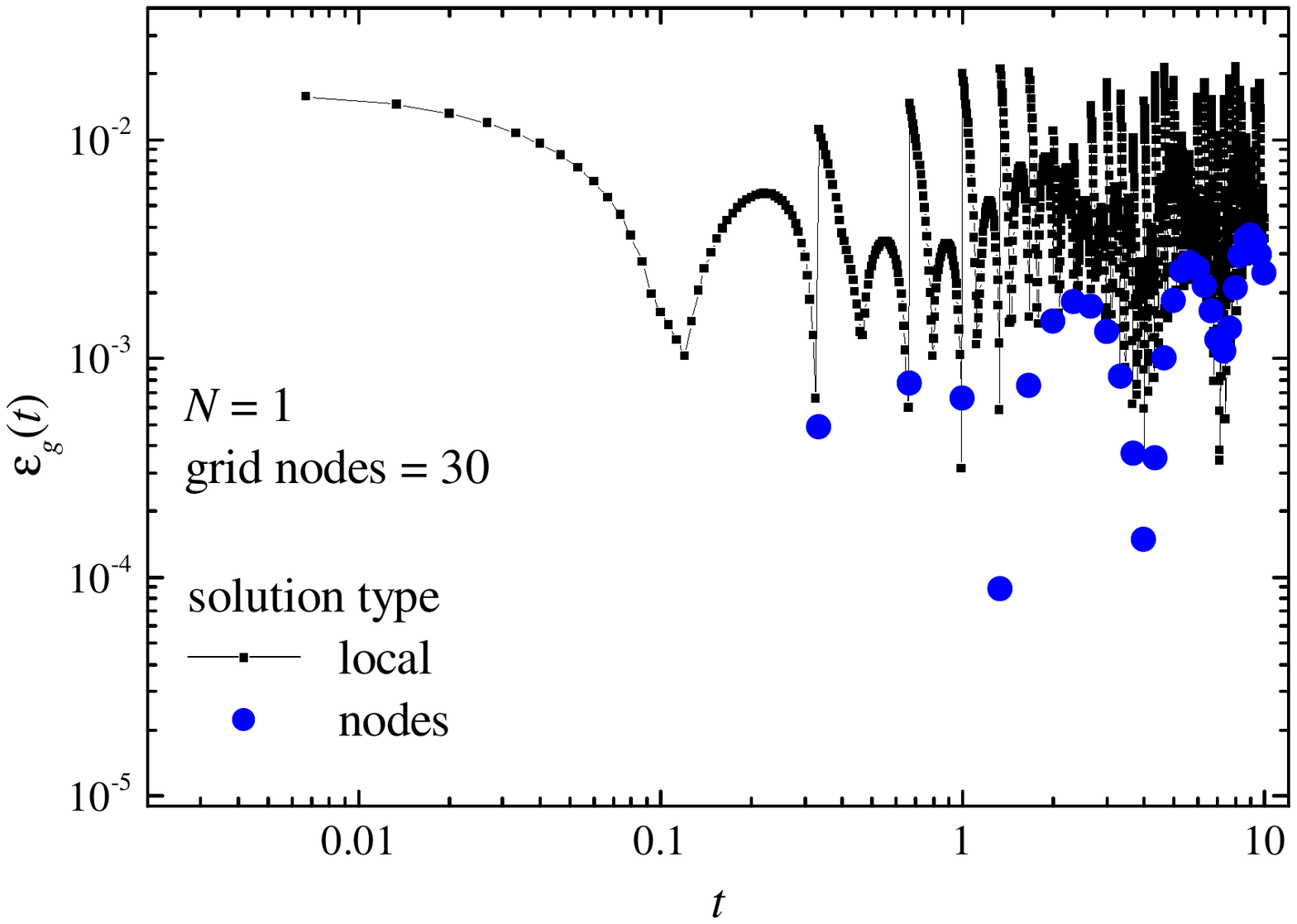}
\vspace{-8mm}\caption{\label{fig:pend_ind1_sol_g_eps:f1}}
\end{subfigure}\hspace{6mm}
\begin{subfigure}{0.275\textwidth}
\includegraphics[width=\textwidth]{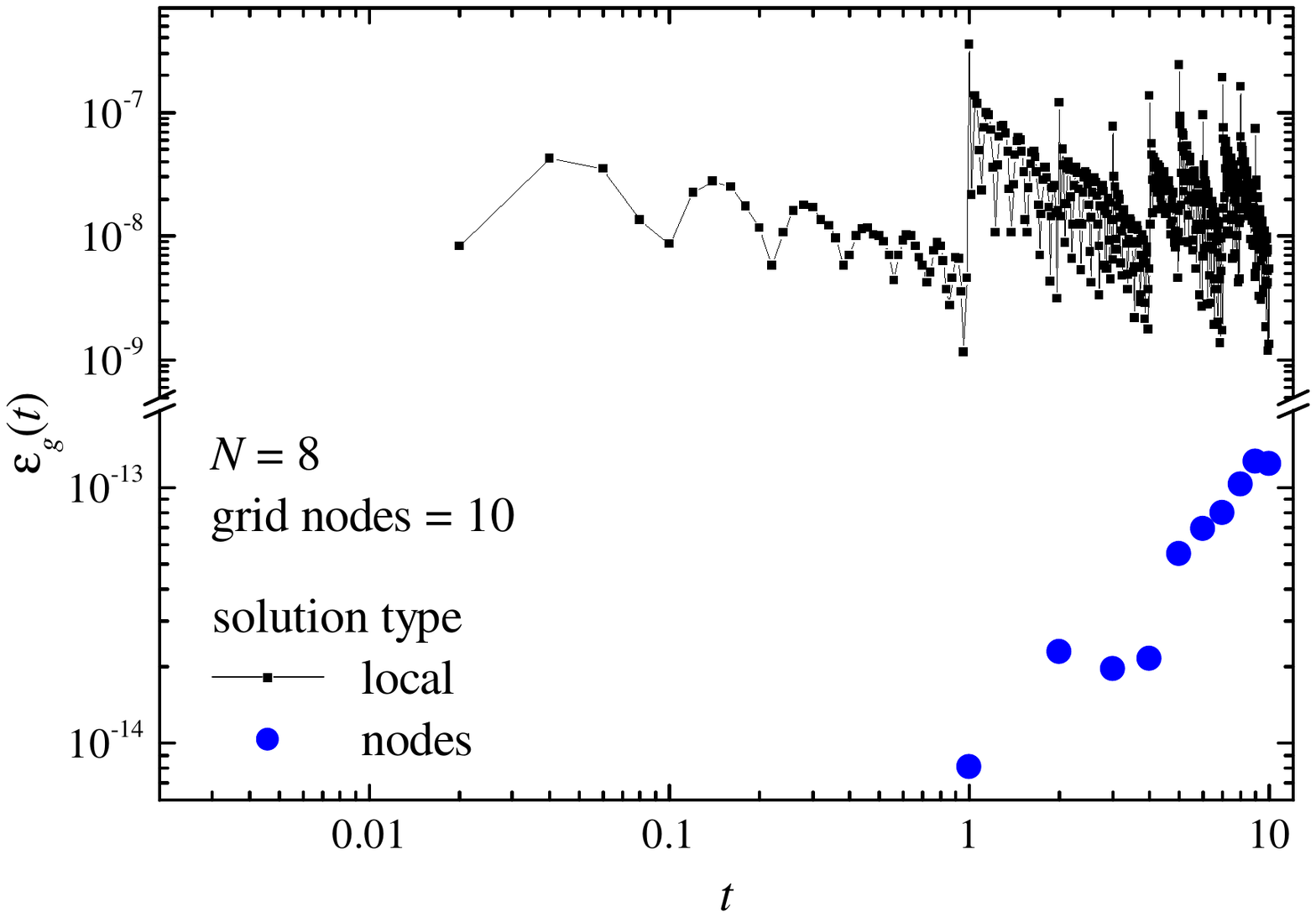}
\vspace{-8mm}\caption{\label{fig:pend_ind1_sol_g_eps:f2}}
\end{subfigure}\hspace{6mm}
\begin{subfigure}{0.275\textwidth}
\includegraphics[width=\textwidth]{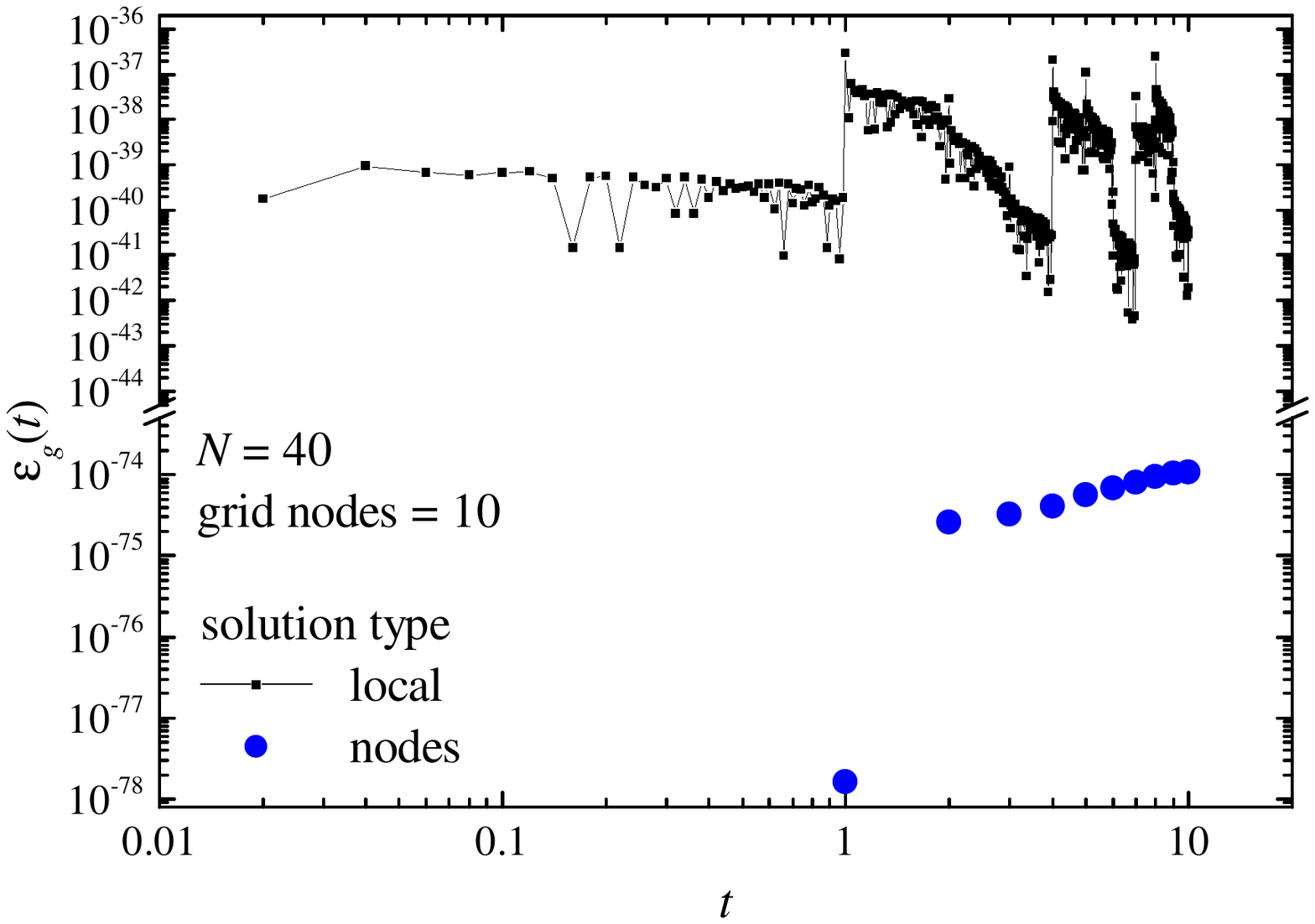}
\vspace{-8mm}\caption{\label{fig:pend_ind1_sol_g_eps:f3}}
\end{subfigure}\\[-2mm]
\caption{%
Numerical solution of the DAE system (\ref{eq:math_pend_dae_ind_3}) of index 1. Comparison of quantitative satisfiability of the conditions $g_{1} = 0$ (\subref{fig:pend_ind1_sol_g_eps:a1}, \subref{fig:pend_ind1_sol_g_eps:a2}, \subref{fig:pend_ind1_sol_g_eps:a3}), $g_{2} = 0$ (\subref{fig:pend_ind1_sol_g_eps:b1}, \subref{fig:pend_ind1_sol_g_eps:b2}, \subref{fig:pend_ind1_sol_g_eps:b3}) and $g_{3} = 0$ (\subref{fig:pend_ind1_sol_g_eps:c1}, \subref{fig:pend_ind1_sol_g_eps:c2}, \subref{fig:pend_ind1_sol_g_eps:c3}), the errors $\varepsilon_{u}(t)$ (\subref{fig:pend_ind1_sol_g_eps:d1}, \subref{fig:pend_ind1_sol_g_eps:d2}, \subref{fig:pend_ind1_sol_g_eps:d3}), $\varepsilon_{v}(t)$ (\subref{fig:pend_ind1_sol_g_eps:e1}, \subref{fig:pend_ind1_sol_g_eps:e2}, \subref{fig:pend_ind1_sol_g_eps:e3}), $\varepsilon_{g}(t)$ (\subref{fig:pend_ind1_sol_g_eps:f1}, \subref{fig:pend_ind1_sol_g_eps:f2}, \subref{fig:pend_ind1_sol_g_eps:f3}), for numerical solution obtained using polynomials with degrees $N = 1$ (\subref{fig:pend_ind1_sol_g_eps:a1}, \subref{fig:pend_ind1_sol_g_eps:b1}, \subref{fig:pend_ind1_sol_g_eps:c1}, \subref{fig:pend_ind1_sol_g_eps:d1}, \subref{fig:pend_ind1_sol_g_eps:e1}, \subref{fig:pend_ind1_sol_g_eps:f1}), $N = 8$ (\subref{fig:pend_ind1_sol_g_eps:a2}, \subref{fig:pend_ind1_sol_g_eps:b2}, \subref{fig:pend_ind1_sol_g_eps:c2}, \subref{fig:pend_ind1_sol_g_eps:d2}, \subref{fig:pend_ind1_sol_g_eps:e2}, \subref{fig:pend_ind1_sol_g_eps:f2}) and $N = 40$ (\subref{fig:pend_ind1_sol_g_eps:a3}, \subref{fig:pend_ind1_sol_g_eps:b3}, \subref{fig:pend_ind1_sol_g_eps:c3}, \subref{fig:pend_ind1_sol_g_eps:d3}, \subref{fig:pend_ind1_sol_g_eps:e3}, \subref{fig:pend_ind1_sol_g_eps:f3}).
}
\label{fig:pend_ind1_sol_g_eps}
\end{figure} 

\begin{figure}[h!]
\captionsetup[subfigure]{%
	position=bottom,
	font+=smaller,
	textfont=normalfont,
	singlelinecheck=off,
	justification=raggedright
}
\centering
\begin{subfigure}{0.275\textwidth}
\includegraphics[width=\textwidth]{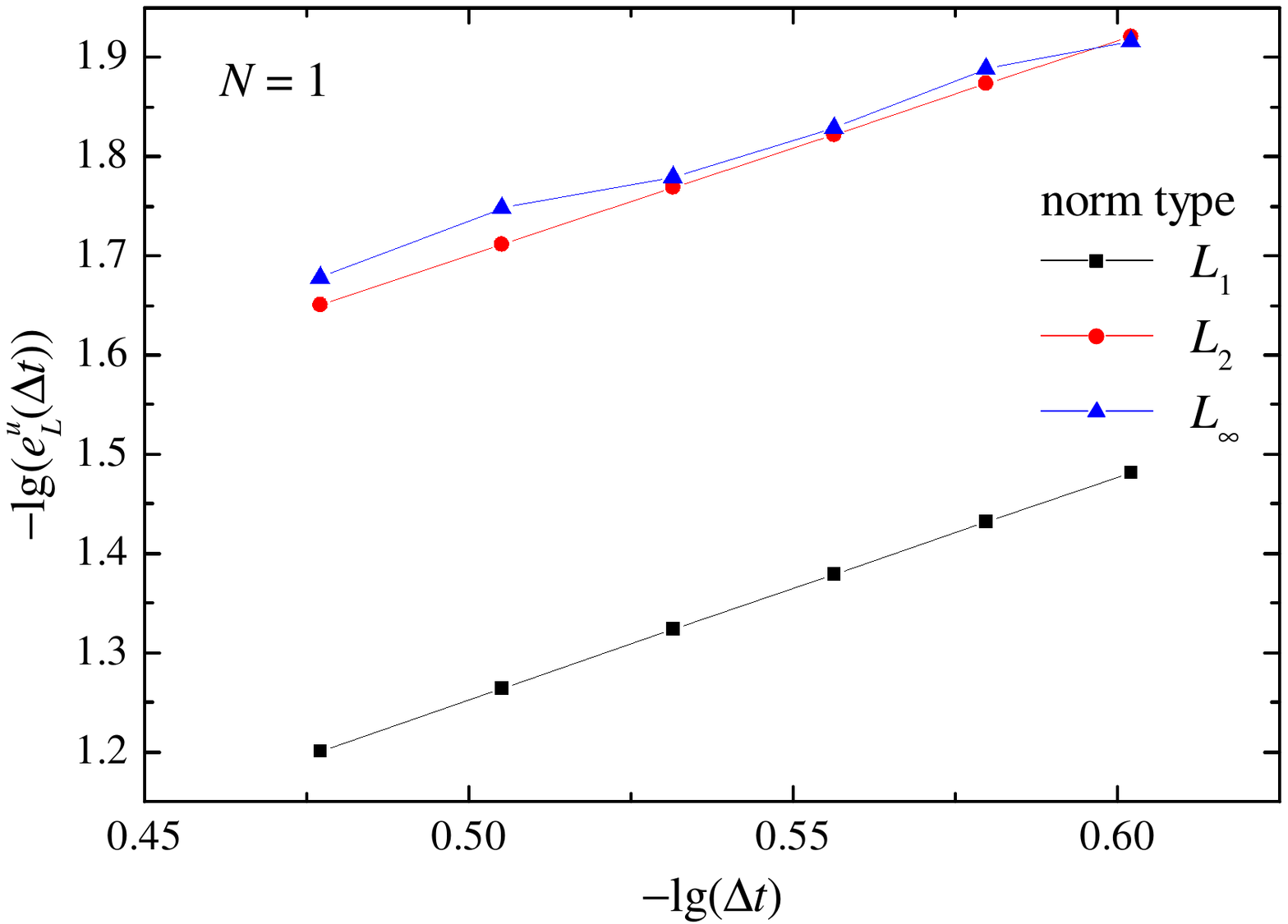}
\vspace{-8mm}\caption{\label{fig:pend_ind1_errors:a1}}
\end{subfigure}\hspace{6mm}
\begin{subfigure}{0.275\textwidth}
\includegraphics[width=\textwidth]{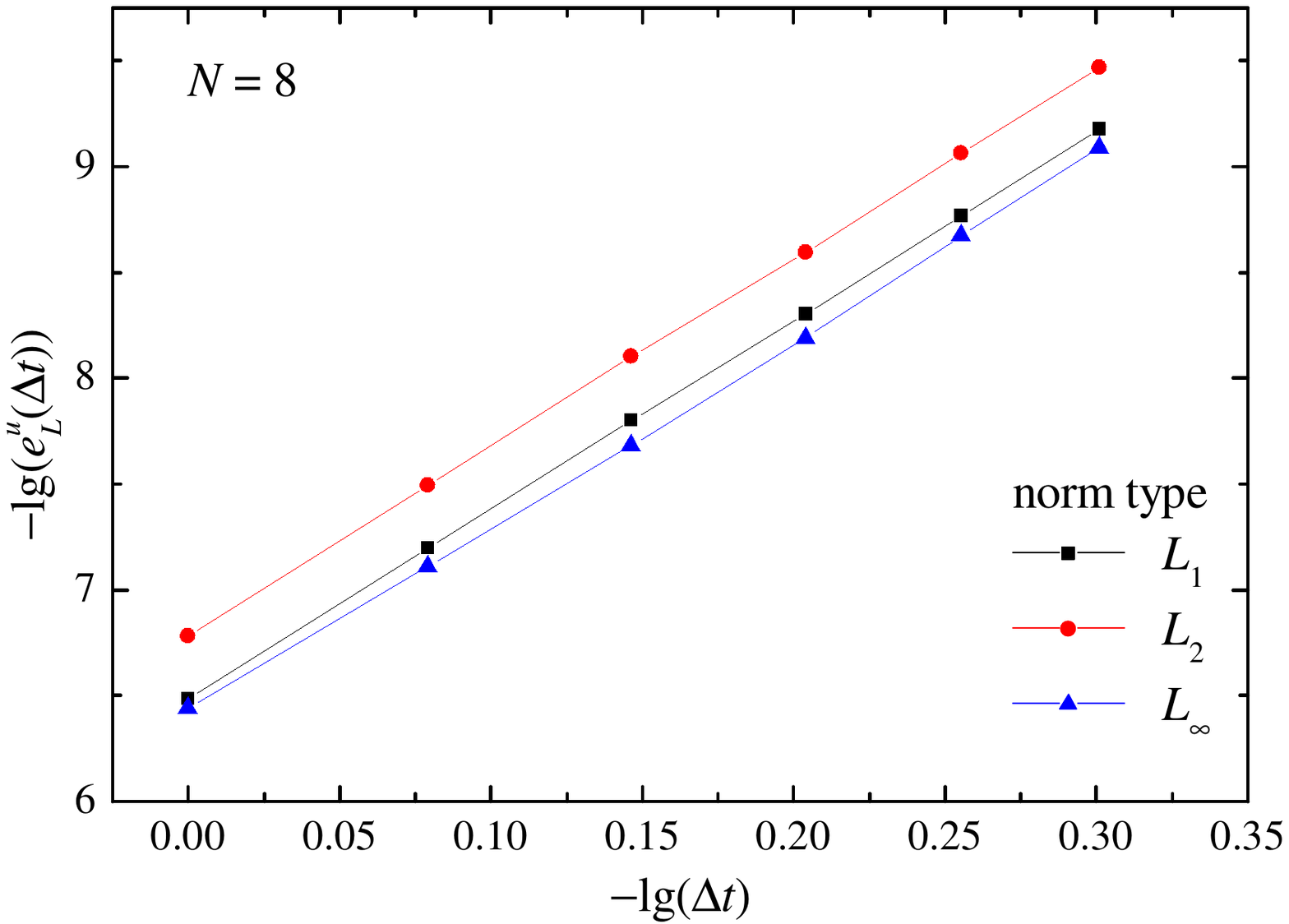}
\vspace{-8mm}\caption{\label{fig:pend_ind1_errors:a2}}
\end{subfigure}\hspace{6mm}
\begin{subfigure}{0.275\textwidth}
\includegraphics[width=\textwidth]{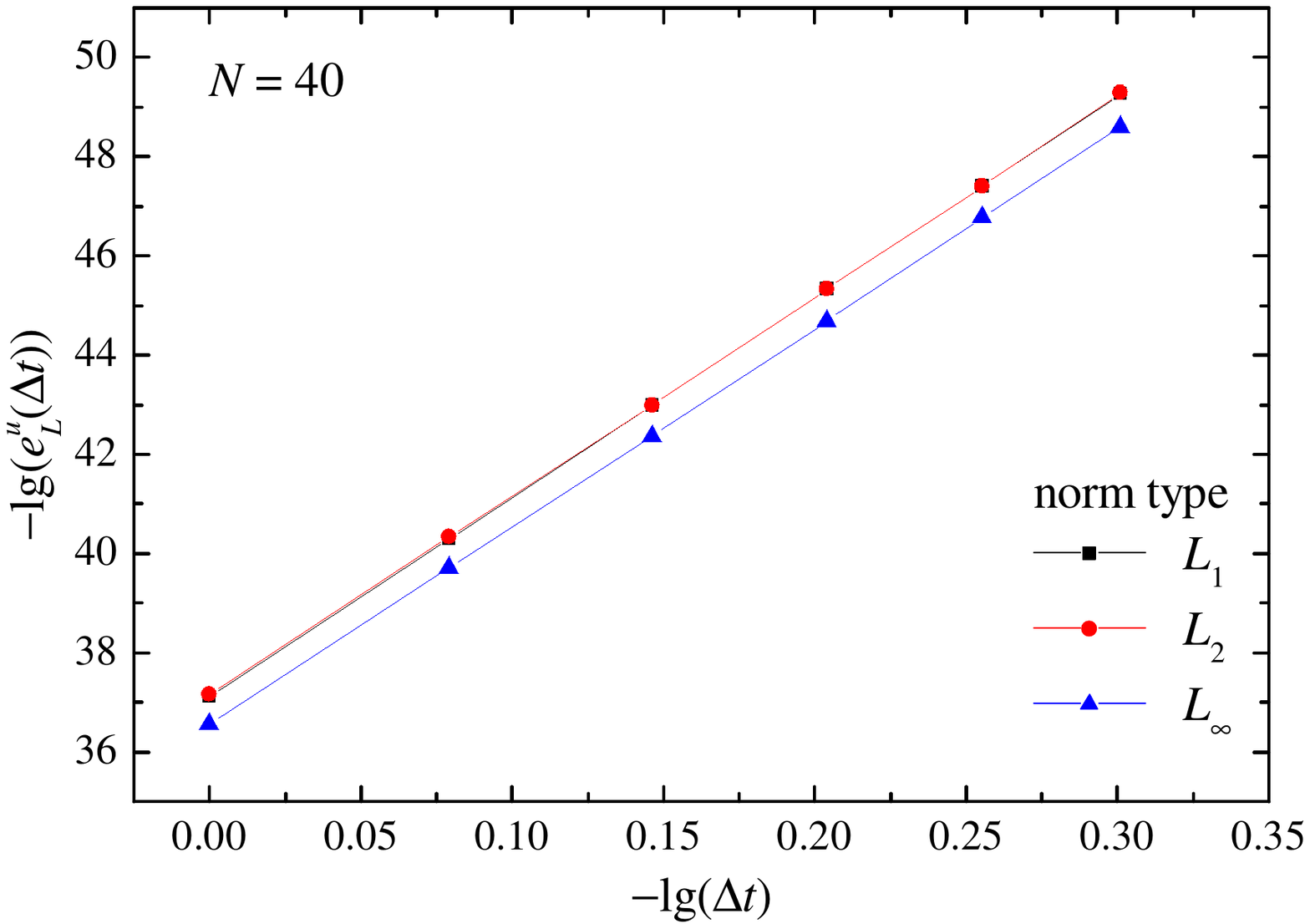}
\vspace{-8mm}\caption{\label{fig:pend_ind1_errors:a3}}
\end{subfigure}\\[-2mm]
\begin{subfigure}{0.275\textwidth}
\includegraphics[width=\textwidth]{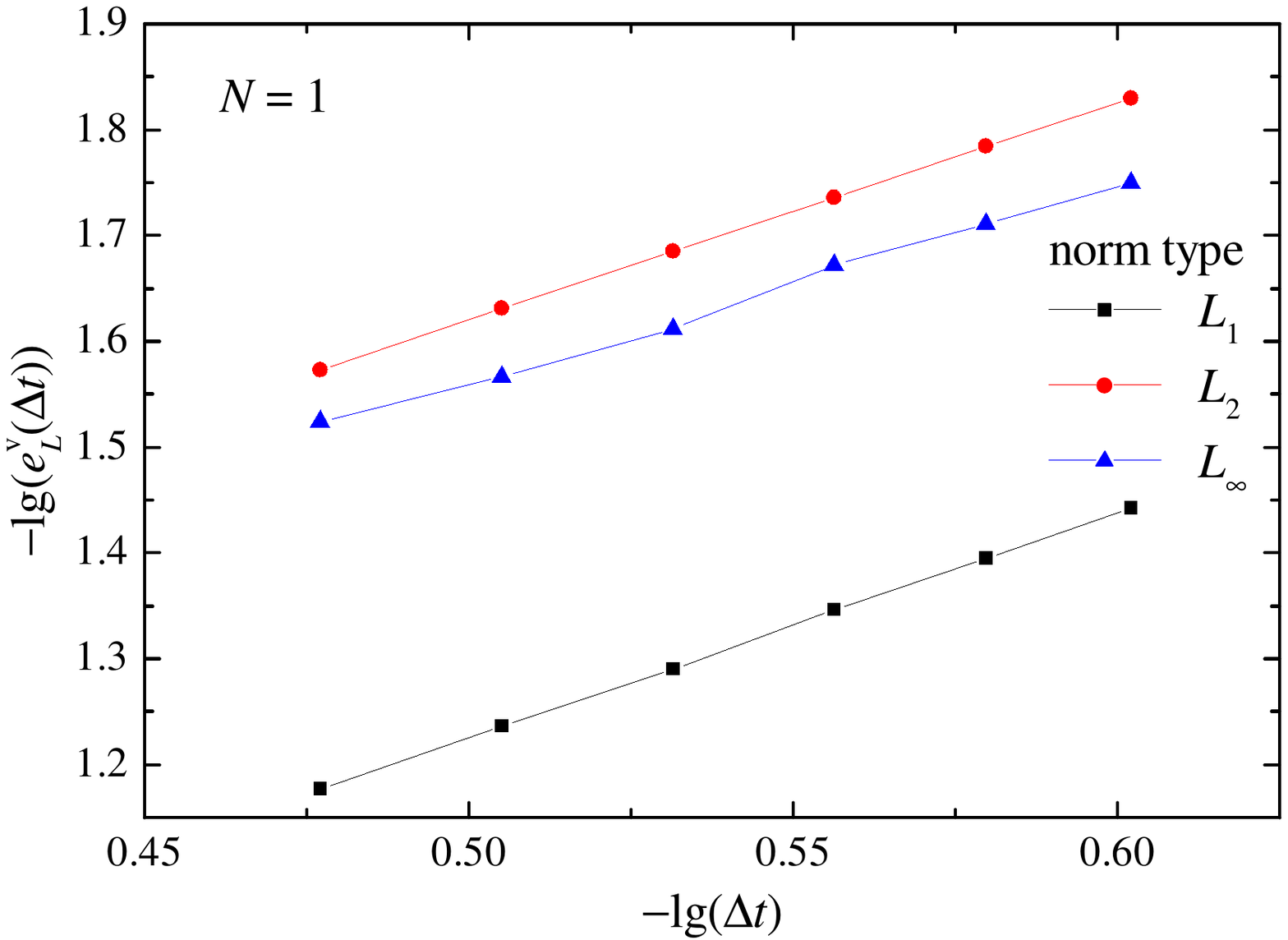}
\vspace{-8mm}\caption{\label{fig:pend_ind1_errors:b1}}
\end{subfigure}\hspace{6mm}
\begin{subfigure}{0.275\textwidth}
\includegraphics[width=\textwidth]{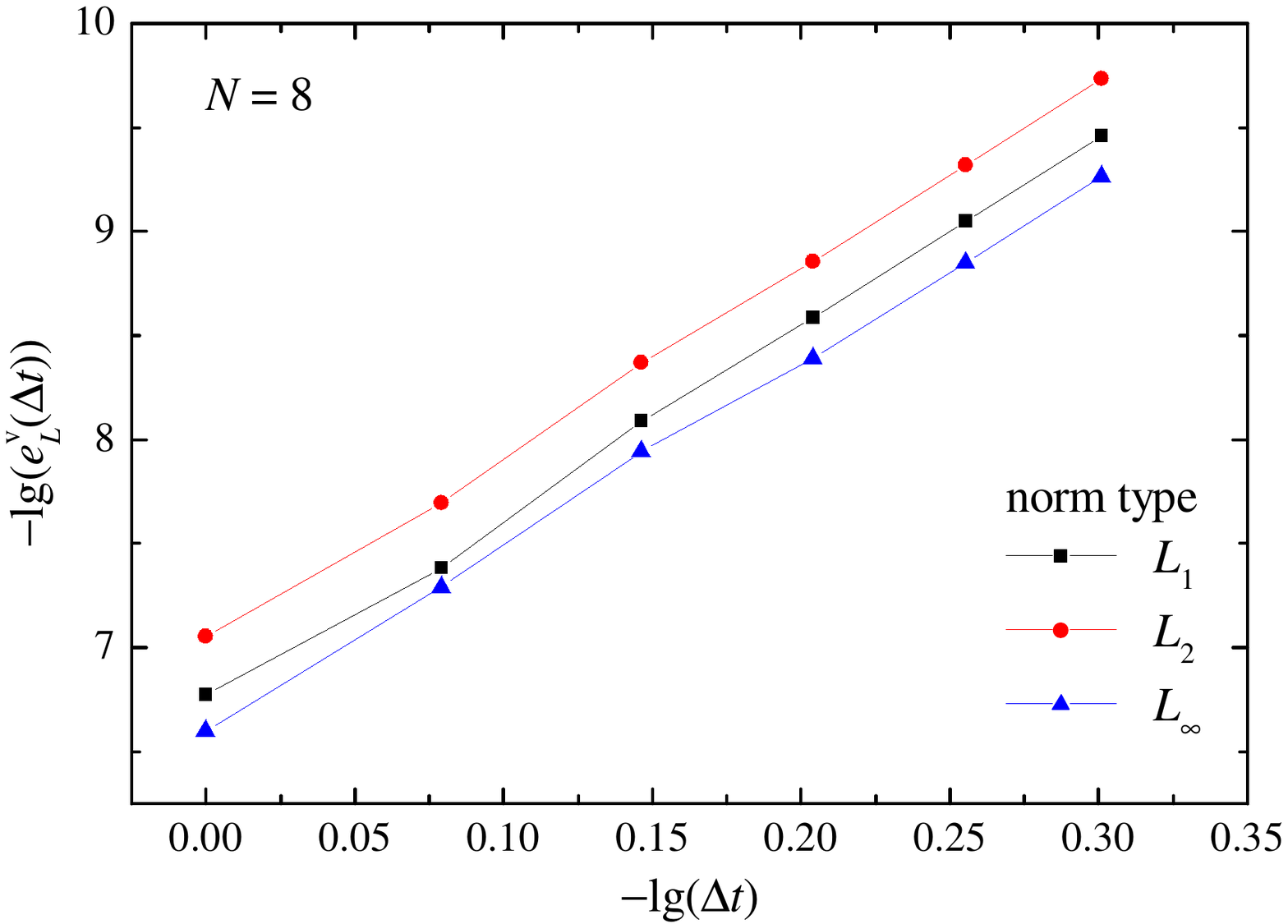}
\vspace{-8mm}\caption{\label{fig:pend_ind1_errors:b2}}
\end{subfigure}\hspace{6mm}
\begin{subfigure}{0.275\textwidth}
\includegraphics[width=\textwidth]{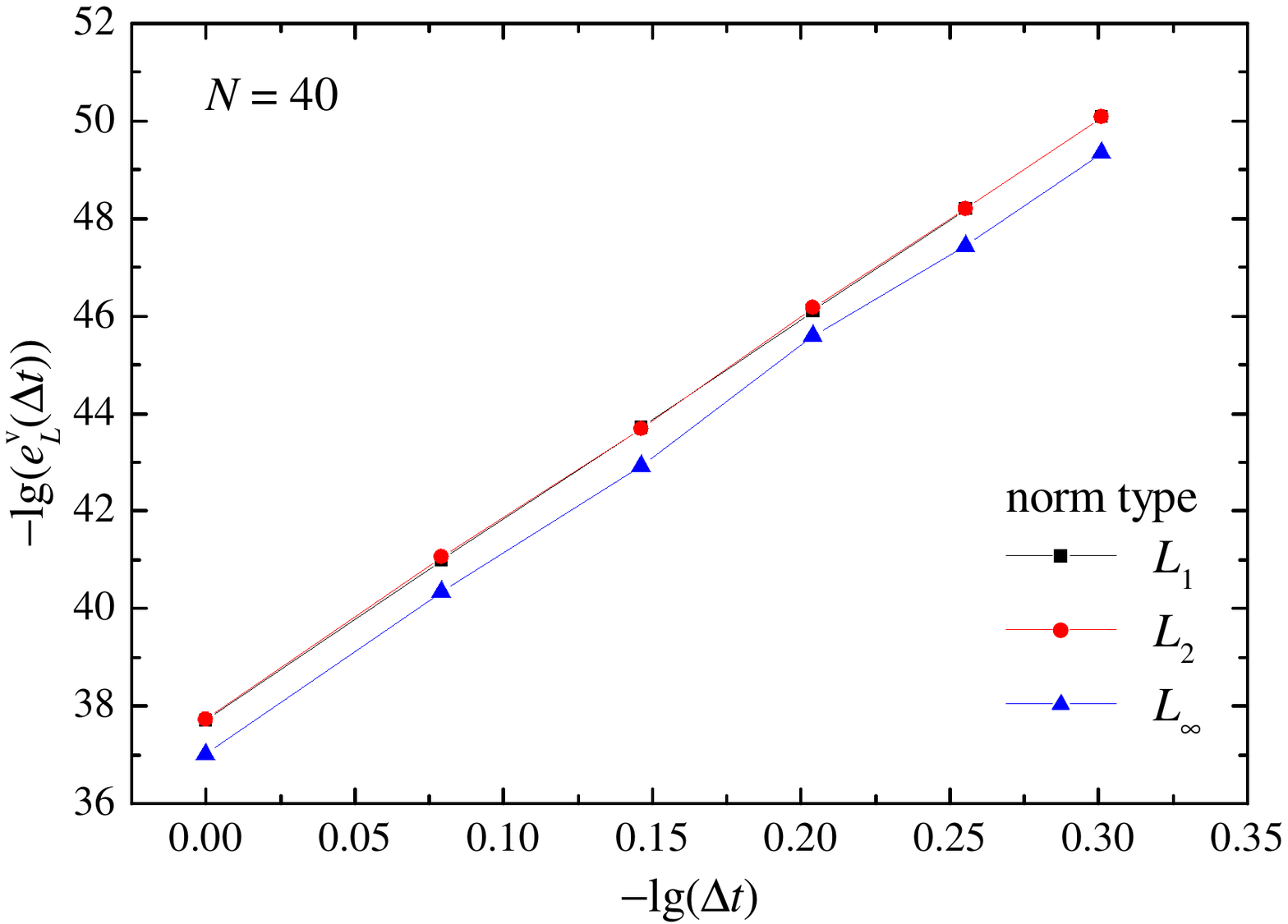}
\vspace{-8mm}\caption{\label{fig:pend_ind1_errors:b3}}
\end{subfigure}\\[-2mm]
\begin{subfigure}{0.275\textwidth}
\includegraphics[width=\textwidth]{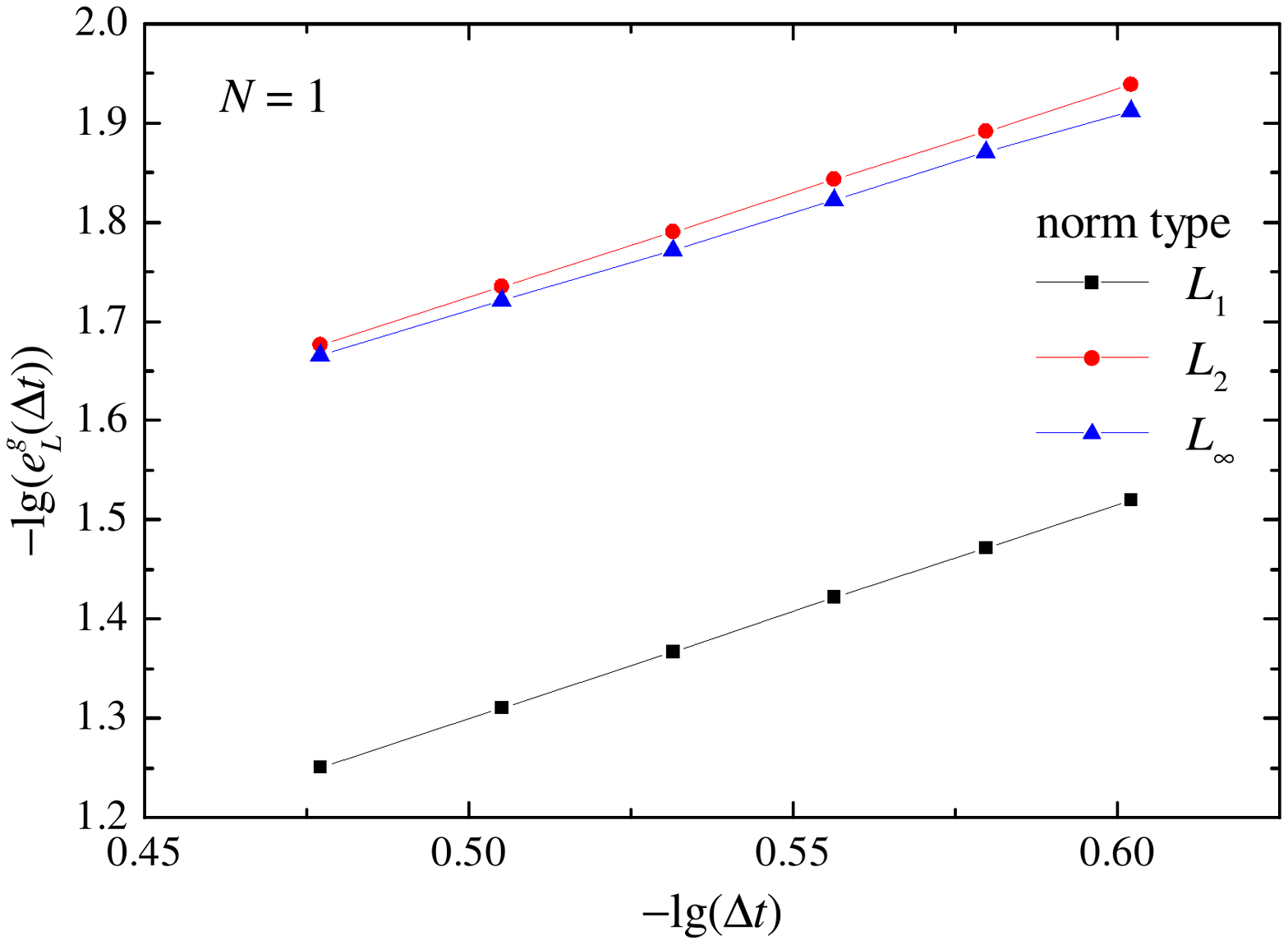}
\vspace{-8mm}\caption{\label{fig:pend_ind1_errors:c1}}
\end{subfigure}\hspace{6mm}
\begin{subfigure}{0.275\textwidth}
\includegraphics[width=\textwidth]{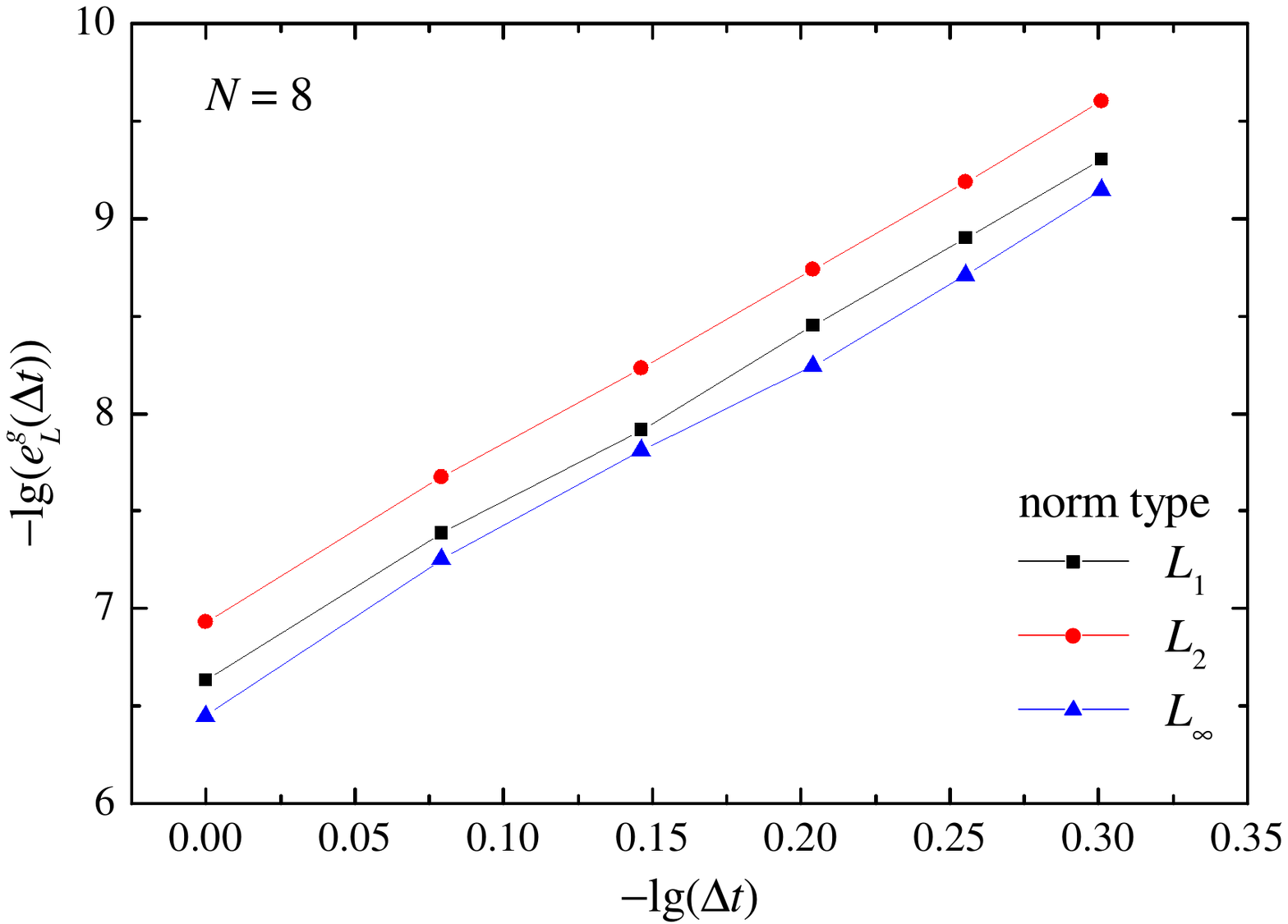}
\vspace{-8mm}\caption{\label{fig:pend_ind1_errors:c2}}
\end{subfigure}\hspace{6mm}
\begin{subfigure}{0.275\textwidth}
\includegraphics[width=\textwidth]{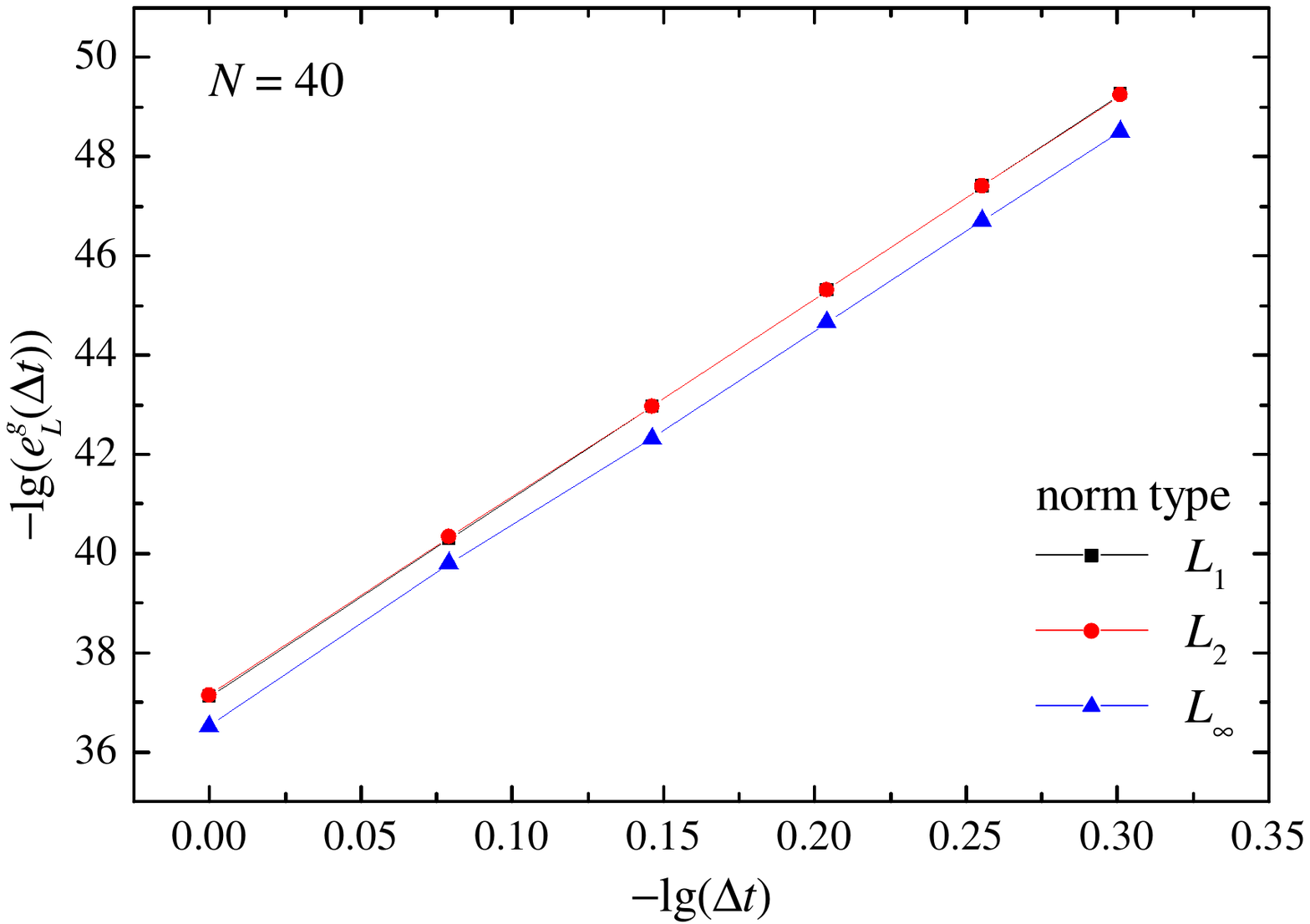}
\vspace{-8mm}\caption{\label{fig:pend_ind1_errors:c3}}
\end{subfigure}\\[-2mm]
\begin{subfigure}{0.275\textwidth}
\includegraphics[width=\textwidth]{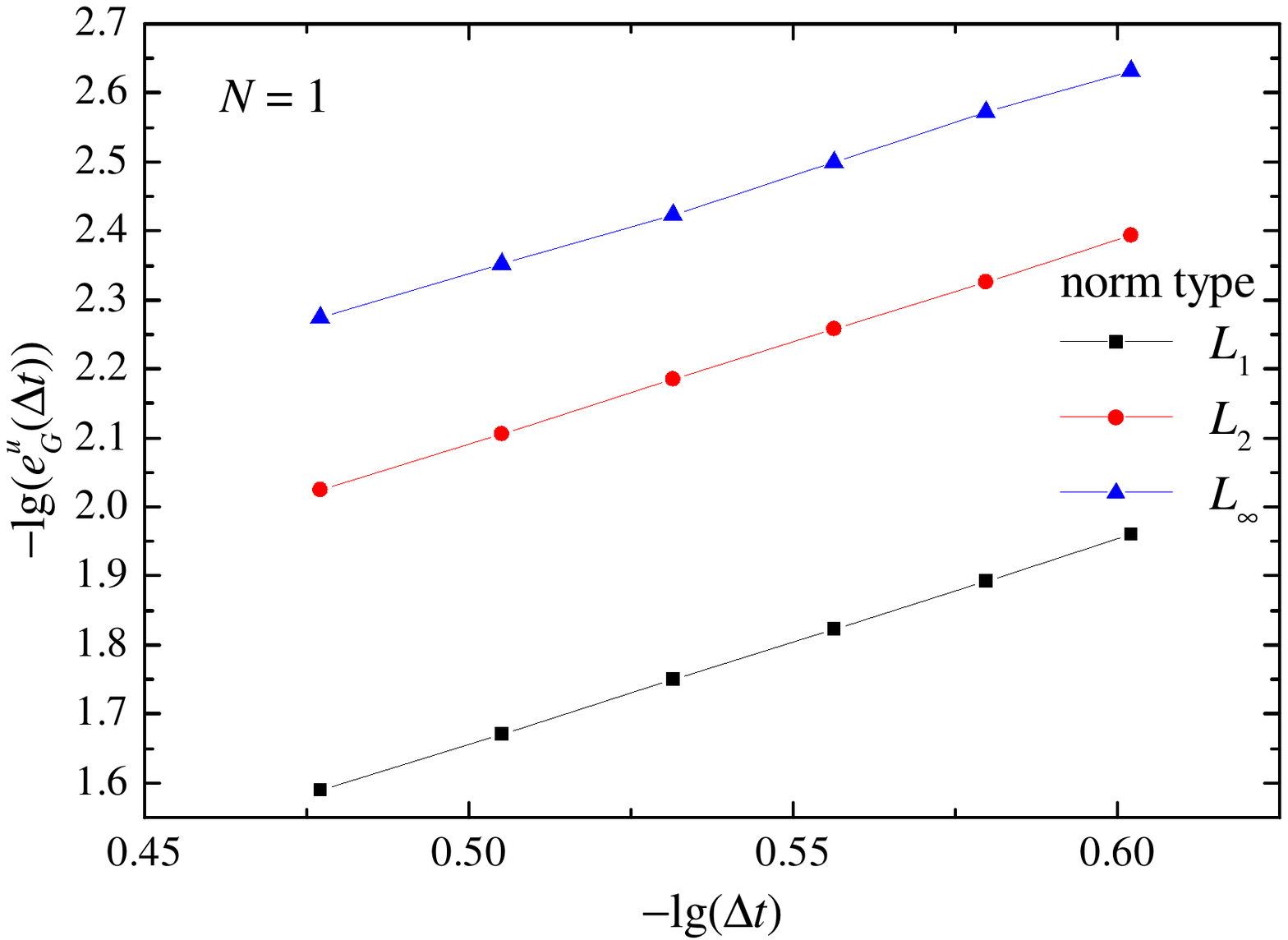}
\vspace{-8mm}\caption{\label{fig:pend_ind1_errors:d1}}
\end{subfigure}\hspace{6mm}
\begin{subfigure}{0.275\textwidth}
\includegraphics[width=\textwidth]{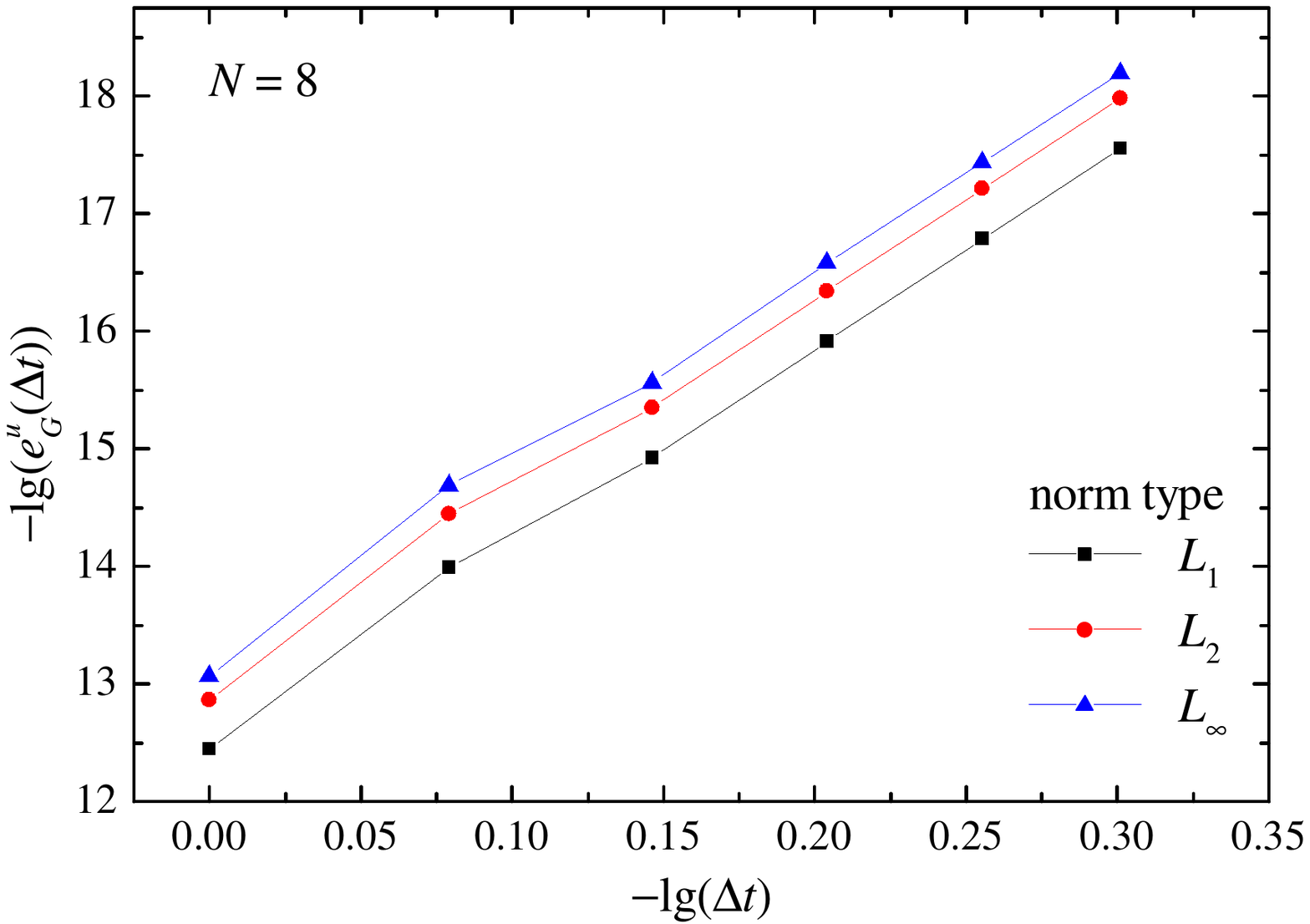}
\vspace{-8mm}\caption{\label{fig:pend_ind1_errors:d2}}
\end{subfigure}\hspace{6mm}
\begin{subfigure}{0.275\textwidth}
\includegraphics[width=\textwidth]{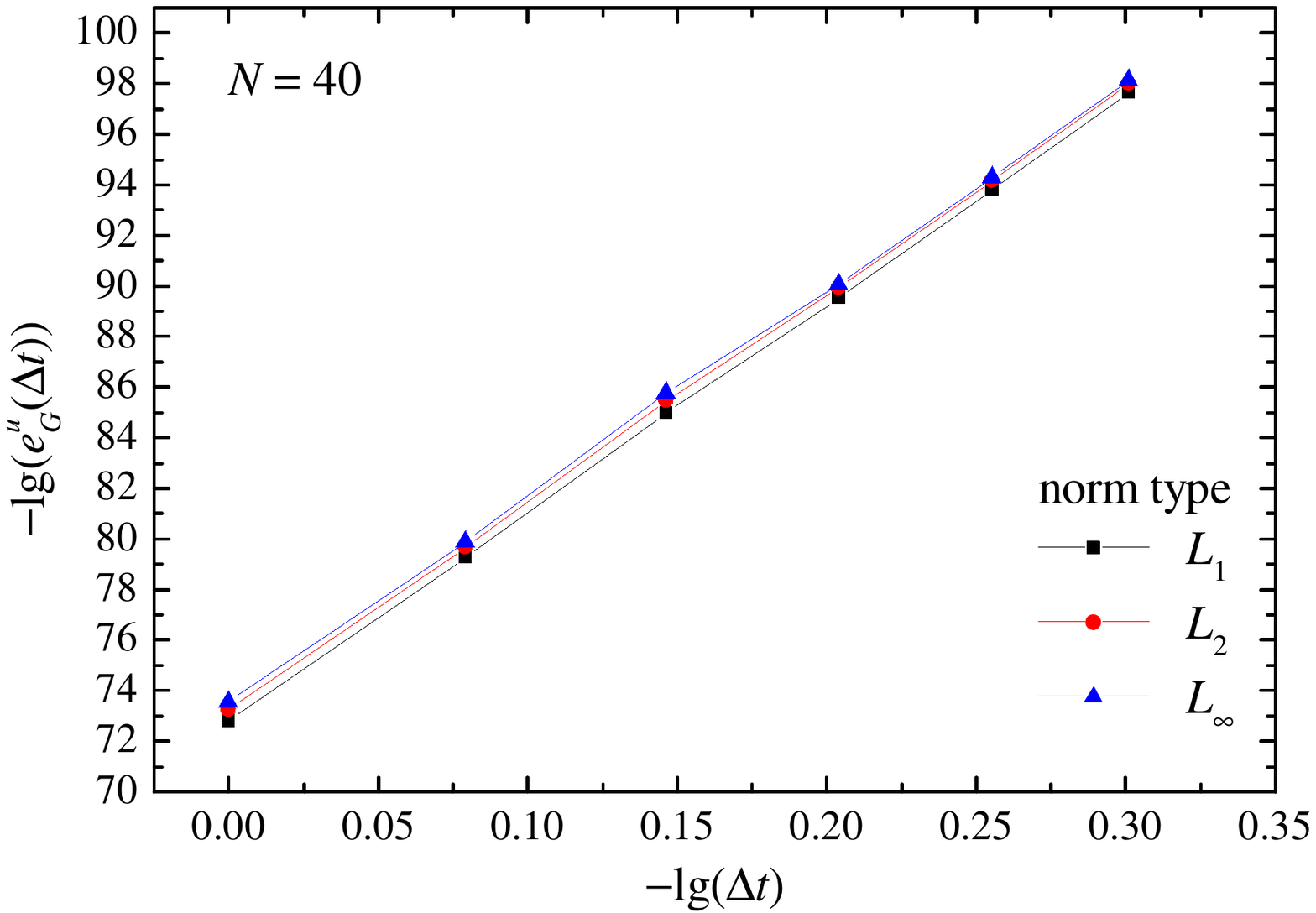}
\vspace{-8mm}\caption{\label{fig:pend_ind1_errors:d3}}
\end{subfigure}\\[-2mm]
\begin{subfigure}{0.275\textwidth}
\includegraphics[width=\textwidth]{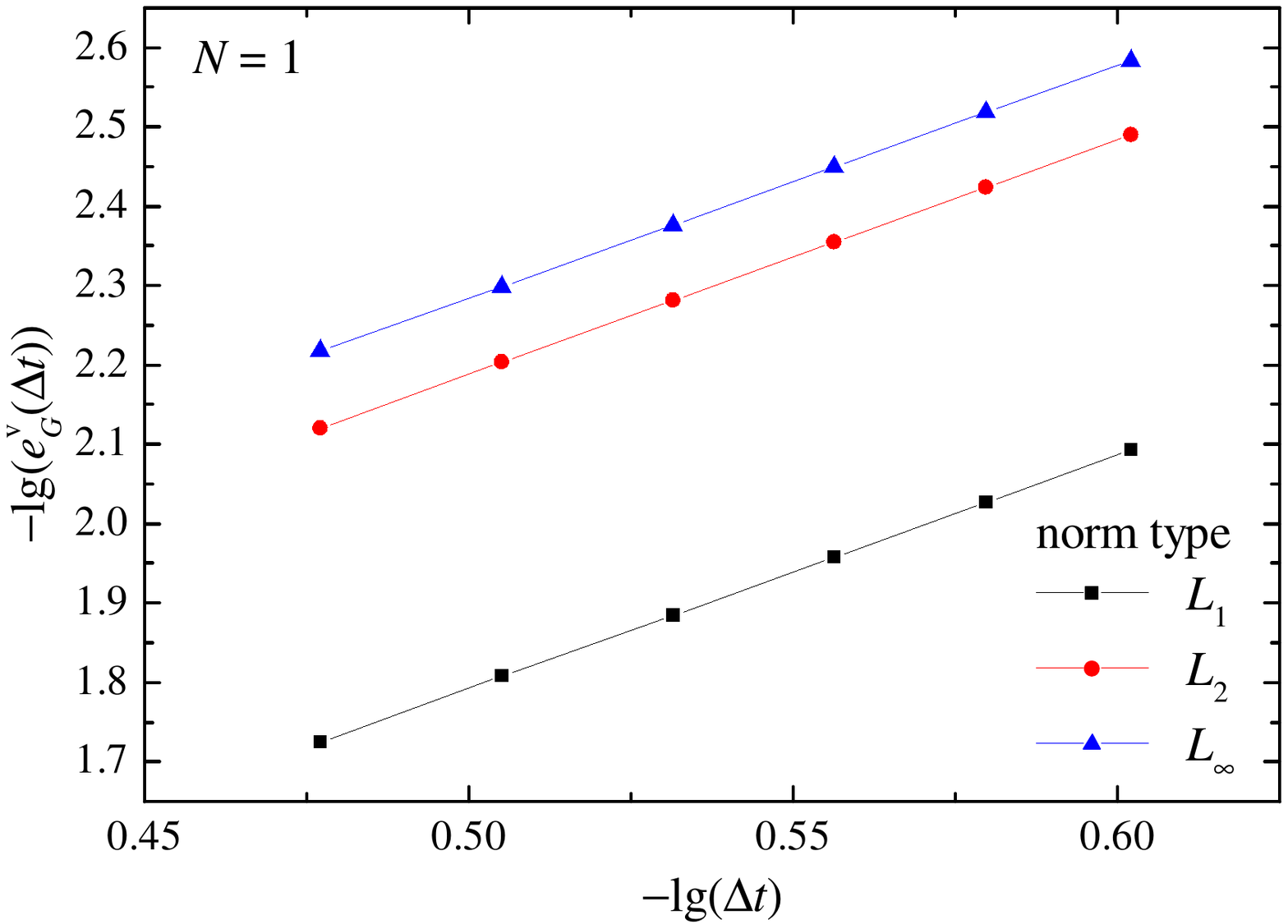}
\vspace{-8mm}\caption{\label{fig:pend_ind1_errors:e1}}
\end{subfigure}\hspace{6mm}
\begin{subfigure}{0.275\textwidth}
\includegraphics[width=\textwidth]{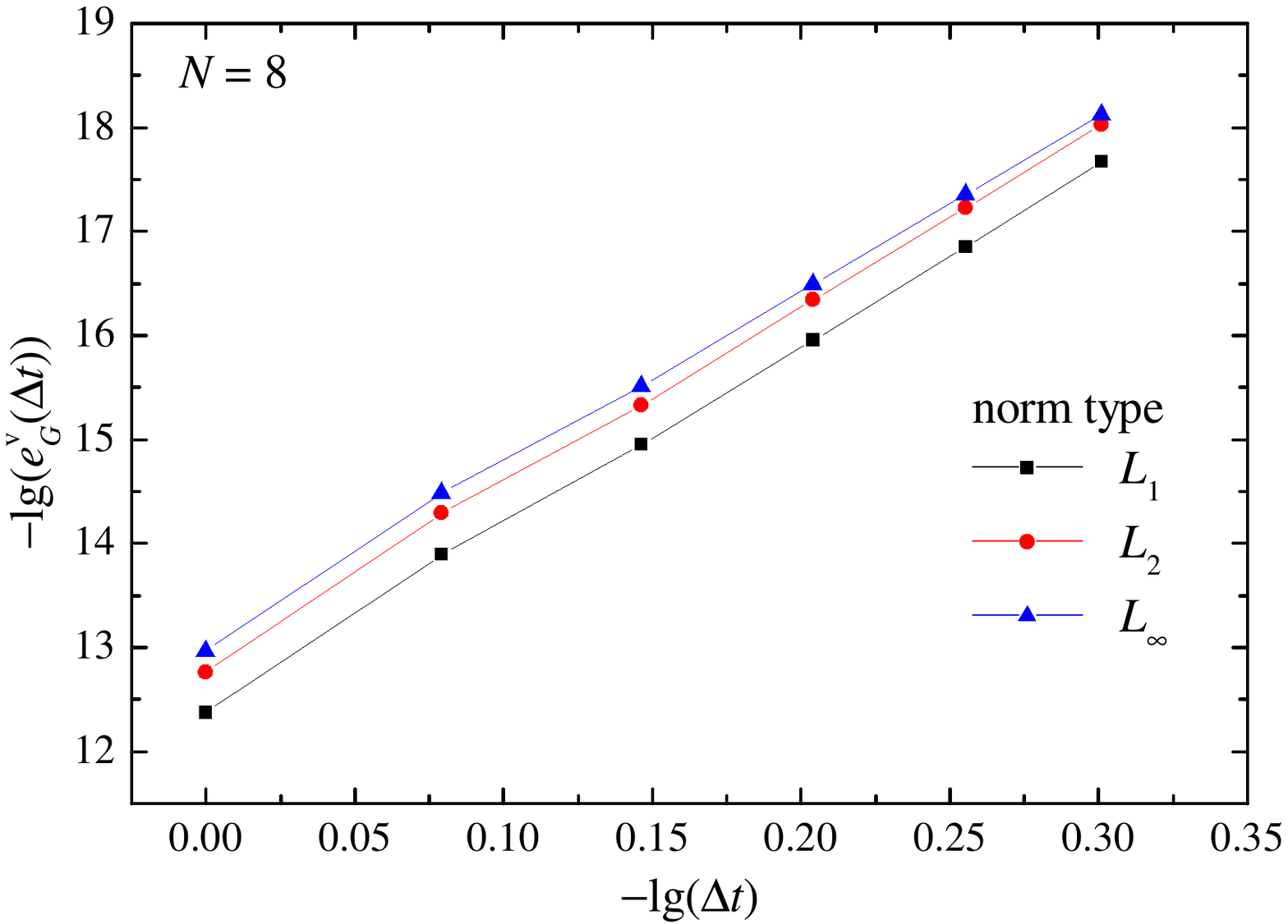}
\vspace{-8mm}\caption{\label{fig:pend_ind1_errors:e2}}
\end{subfigure}\hspace{6mm}
\begin{subfigure}{0.275\textwidth}
\includegraphics[width=\textwidth]{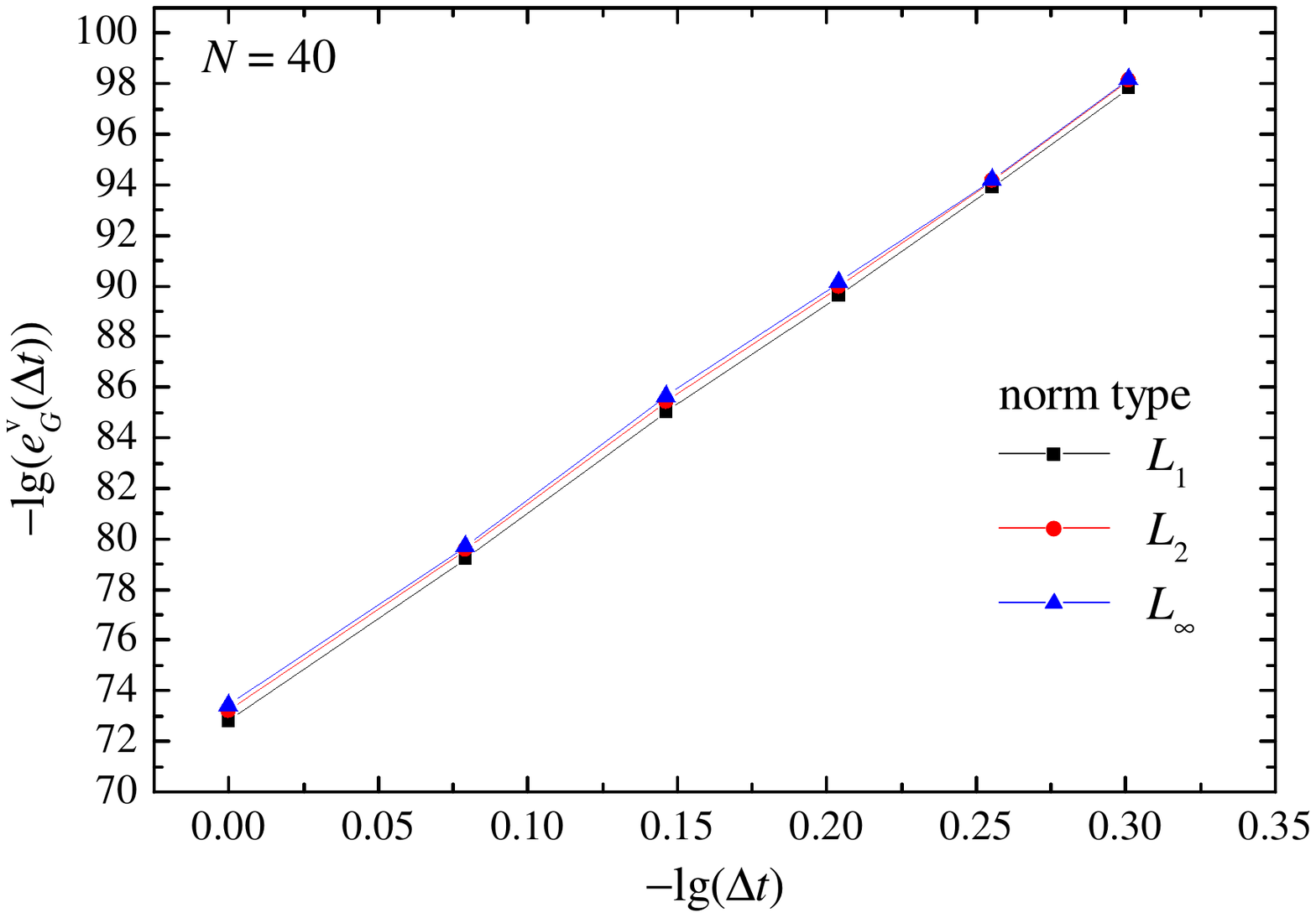}
\vspace{-8mm}\caption{\label{fig:pend_ind1_errors:e3}}
\end{subfigure}\\[-2mm]
\begin{subfigure}{0.275\textwidth}
\includegraphics[width=\textwidth]{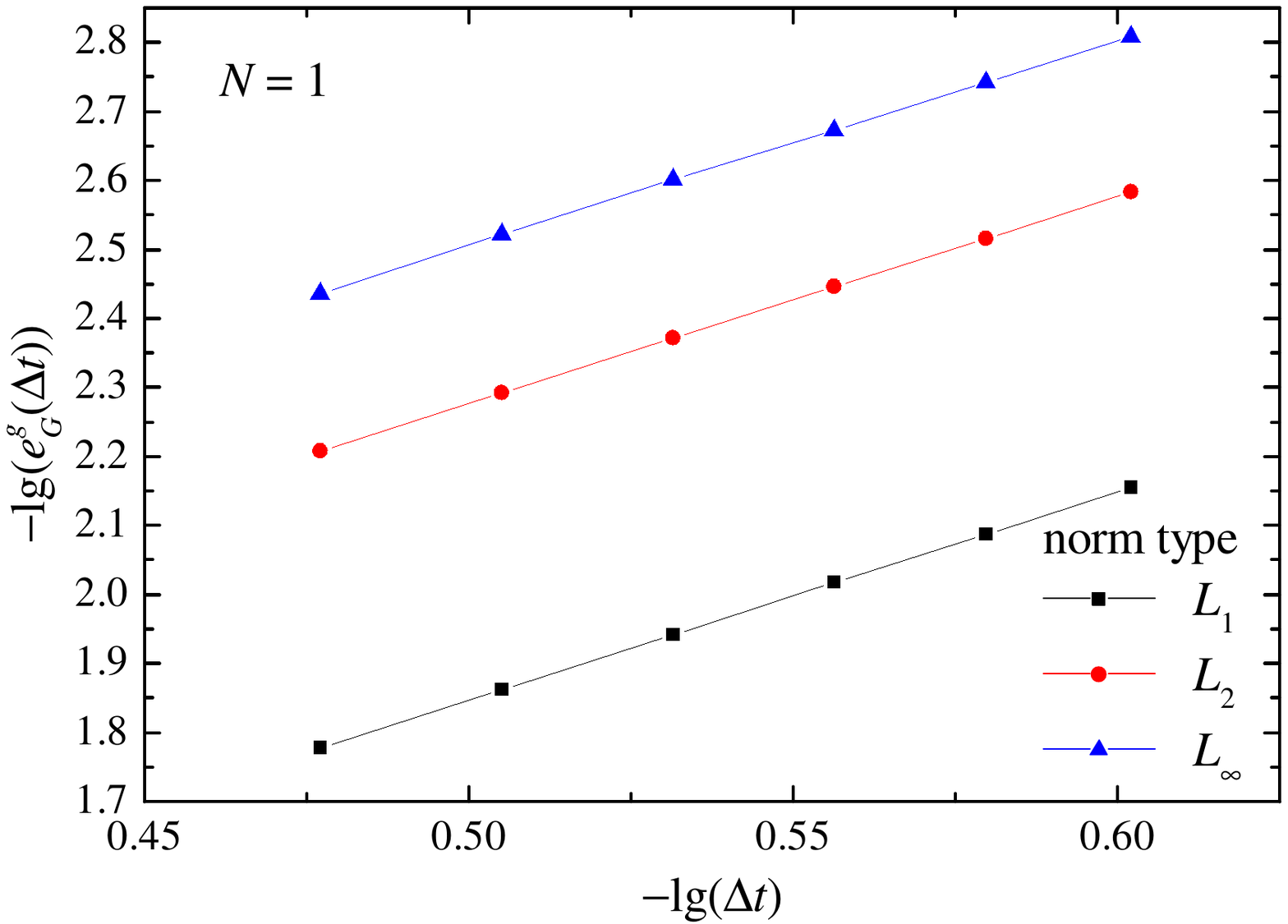}
\vspace{-8mm}\caption{\label{fig:pend_ind1_errors:f1}}
\end{subfigure}\hspace{6mm}
\begin{subfigure}{0.275\textwidth}
\includegraphics[width=\textwidth]{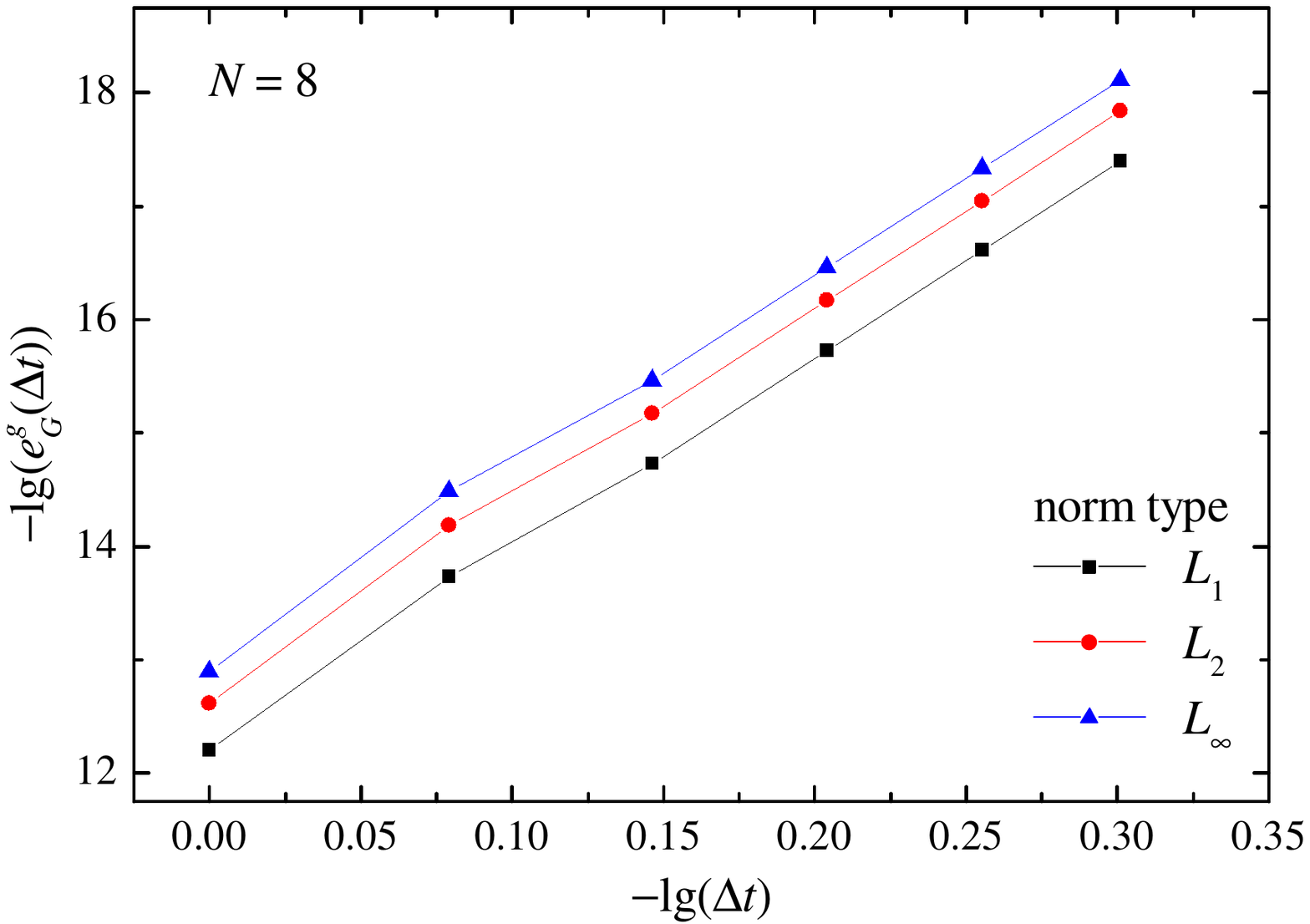}
\vspace{-8mm}\caption{\label{fig:pend_ind1_errors:f2}}
\end{subfigure}\hspace{6mm}
\begin{subfigure}{0.275\textwidth}
\includegraphics[width=\textwidth]{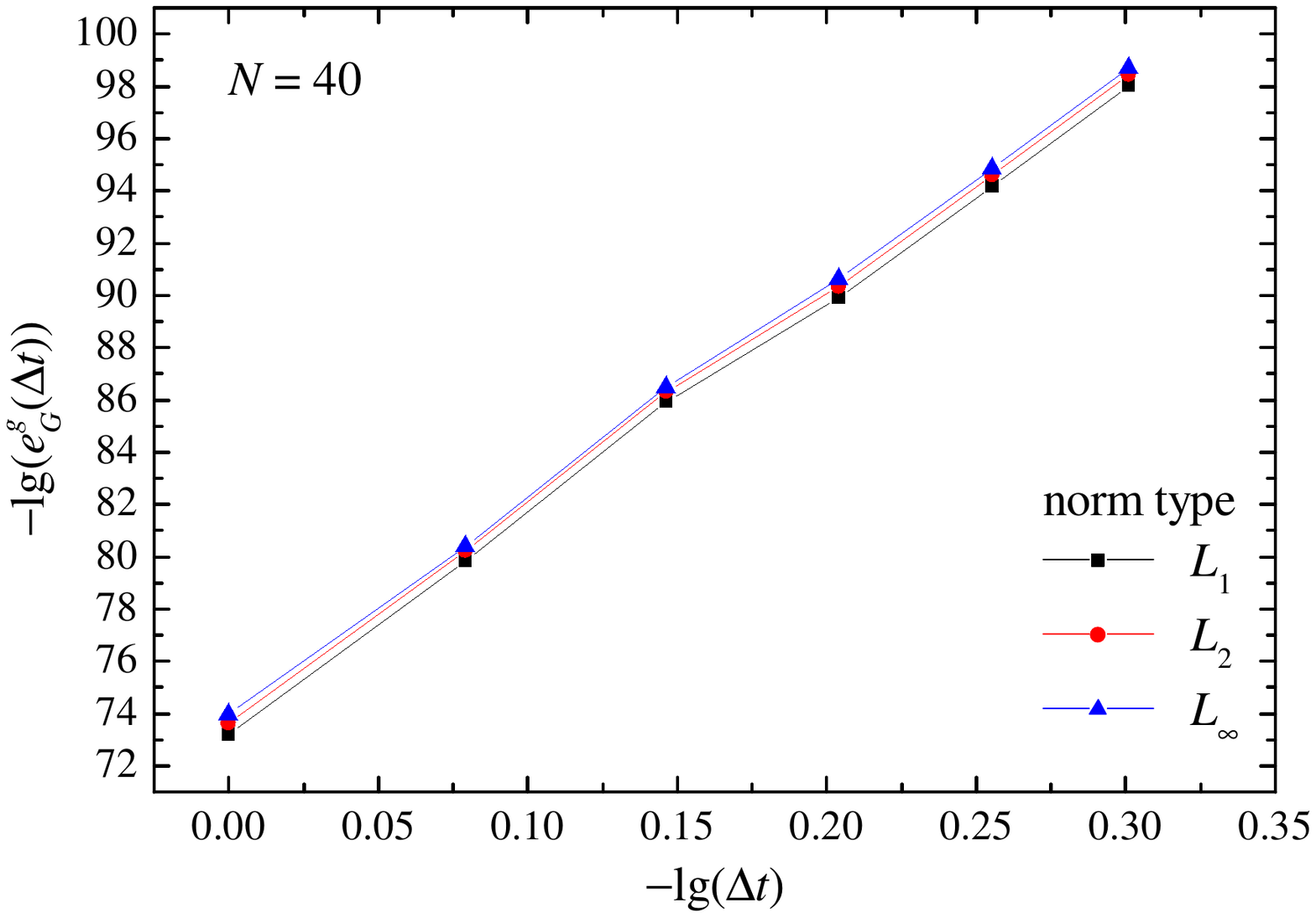}
\vspace{-8mm}\caption{\label{fig:pend_ind1_errors:f3}}
\end{subfigure}\\[-2mm]
\caption{%
Log-log plot of the dependence of the global errors for the local solution $e_{L}^{u}$ (\subref{fig:pend_ind1_errors:a1}, \subref{fig:pend_ind1_errors:a2}, \subref{fig:pend_ind1_errors:a3}), $e_{L}^{v}$ (\subref{fig:pend_ind1_errors:b1}, \subref{fig:pend_ind1_errors:b2}, \subref{fig:pend_ind1_errors:b3}), $e_{L}^{g}$ (\subref{fig:pend_ind1_errors:c1}, \subref{fig:pend_ind1_errors:c2}, \subref{fig:pend_ind1_errors:c3}) and the solution at nodes $e_{G}^{u}$ (\subref{fig:pend_ind1_errors:d1}, \subref{fig:pend_ind1_errors:d2}, \subref{fig:pend_ind1_errors:d3}), $e_{G}^{v}$ (\subref{fig:pend_ind1_errors:e1}, \subref{fig:pend_ind1_errors:e2}, \subref{fig:pend_ind1_errors:e3}), $e_{G}^{g}$ (\subref{fig:pend_ind1_errors:f1}, \subref{fig:pend_ind1_errors:f2}, \subref{fig:pend_ind1_errors:f3}) on the discretization step $\mathrm{\Delta}t$, obtained in the norms $L_{1}$, $L_{2}$ and $L_{\infty}$, by numerical solution of the DAE system (\ref{eq:math_pend_dae_ind_3}) of index 1 obtained using polynomials with degrees $N = 1$ (\subref{fig:pend_ind1_errors:a1}, \subref{fig:pend_ind1_errors:b1}, \subref{fig:pend_ind1_errors:c1}, \subref{fig:pend_ind1_errors:d1}, \subref{fig:pend_ind1_errors:e1}, \subref{fig:pend_ind1_errors:f1}), $N = 8$ (\subref{fig:pend_ind1_errors:a2}, \subref{fig:pend_ind1_errors:b2}, \subref{fig:pend_ind1_errors:c2}, \subref{fig:pend_ind1_errors:d2}, \subref{fig:pend_ind1_errors:e2}, \subref{fig:pend_ind1_errors:f2}) and $N = 40$ (\subref{fig:pend_ind1_errors:a3}, \subref{fig:pend_ind1_errors:b3}, \subref{fig:pend_ind1_errors:c3}, \subref{fig:pend_ind1_errors:d3}, \subref{fig:pend_ind1_errors:e3}, \subref{fig:pend_ind1_errors:f3}).
}
\label{fig:pend_ind1_errors}
\end{figure} 

\begin{table*}[h!]
\centering
\caption{%
Convergence orders $p_{L_{1}}^{n}$, $p_{L_{2}}^{n}$, $p_{L_{\infty}}^{n}$, calculated in norms $L_{1}$, $L_{2}$, $L_{\infty}$, of \textit{the numerical solution at the nodes} $(\mathbf{u}_{n}, \mathbf{v}_{n})$ of the ADER-DG method for the DAE system (\ref{eq:math_pend_dae_ind_3}) of index 1; $N$ is the degree of the basis polynomials $\varphi_{p}$. Orders $p^{n, u}$ are calculated for solution $\mathbf{u}_{n}$; orders $p^{n, v}$ --- for solution $\mathbf{v}_{n}$; orders $p^{n, g}$ --- for the conditions $\mathbf{g} = 0$ on the numerical solution at the nodes $(\mathbf{u}_{n}, \mathbf{v}_{n})$. The theoretical value of convergence order $p_{\rm th.}^{n} = 2N+1$ is applicable for the ADER-DG method for ODE problems and is presented for comparison.
}
\label{tab:conv_orders_nodes_pend_ind1}
\begin{tabular}{@{}|l|lll|lll|lll|c|@{}}
\toprule
$N$ & $p_{L_{1}}^{n, u}$ & $p_{L_{2}}^{n, u}$ & $p_{L_{\infty}}^{n, u}$ & $p_{L_{1}}^{n, v}$ & $p_{L_{2}}^{n, v}$ & $p_{L_{\infty}}^{n, v}$ & $p_{L_{1}}^{n, g}$ & $p_{L_{2}}^{n, g}$ & $p_{L_{\infty}}^{n, g}$ & $p_{\rm th.}^{n}$ \\
\midrule
$1$	&	$2.97$	&	$2.96$	&	$2.89$	&	$2.94$	&	$2.96$	&	$2.93$	&	$3.02$	&	$3.00$	&	$2.97$	&	$3$\\
$2$	&	$5.02$	&	$5.00$	&	$4.93$	&	$5.10$	&	$5.05$	&	$4.97$	&	$5.05$	&	$5.04$	&	$5.01$	&	$5$\\
$3$	&	$6.97$	&	$6.96$	&	$6.94$	&	$6.94$	&	$7.01$	&	$6.95$	&	$7.02$	&	$7.01$	&	$6.96$	&	$7$\\
$4$	&	$9.01$	&	$9.00$	&	$8.94$	&	$9.13$	&	$9.07$	&	$8.97$	&	$9.05$	&	$9.04$	&	$9.02$	&	$9$\\
$5$	&	$10.96$	&	$10.95$	&	$10.93$	&	$10.99$	&	$11.02$	&	$10.94$	&	$11.02$	&	$11.01$	&	$11.00$	&	$11$\\
$6$	&	$13.23$	&	$13.24$	&	$13.20$	&	$13.82$	&	$13.67$	&	$13.21$	&	$13.49$	&	$13.56$	&	$13.57$	&	$13$\\
$7$	&	$14.81$	&	$14.74$	&	$14.55$	&	$15.30$	&	$15.11$	&	$14.35$	&	$14.82$	&	$14.72$	&	$14.50$	&	$15$\\
$8$	&	$16.71$	&	$16.70$	&	$16.71$	&	$17.39$	&	$17.30$	&	$16.96$	&	$17.06$	&	$17.10$	&	$17.08$	&	$17$\\
$9$	&	$19.25$	&	$19.29$	&	$19.19$	&	$19.86$	&	$19.78$	&	$19.37$	&	$19.40$	&	$19.33$	&	$19.17$	&	$19$\\
$10$	&	$20.68$	&	$20.64$	&	$20.62$	&	$21.23$	&	$21.25$	&	$21.11$	&	$20.78$	&	$20.75$	&	$20.70$	&	$21$\\
$11$	&	$23.27$	&	$23.33$	&	$23.41$	&	$23.87$	&	$23.91$	&	$23.76$	&	$23.49$	&	$23.47$	&	$23.37$	&	$23$\\
$12$	&	$25.02$	&	$24.99$	&	$25.02$	&	$25.43$	&	$25.55$	&	$25.59$	&	$25.01$	&	$24.98$	&	$24.91$	&	$25$\\
$13$	&	$27.29$	&	$27.30$	&	$27.48$	&	$27.77$	&	$27.93$	&	$28.04$	&	$27.42$	&	$27.40$	&	$27.32$	&	$27$\\
$14$	&	$29.37$	&	$29.33$	&	$29.32$	&	$29.62$	&	$29.86$	&	$30.07$	&	$29.26$	&	$29.23$	&	$29.14$	&	$29$\\
$15$	&	$31.50$	&	$31.45$	&	$31.29$	&	$31.71$	&	$32.02$	&	$32.34$	&	$31.40$	&	$31.38$	&	$31.31$	&	$31$\\
$16$	&	$33.67$	&	$33.59$	&	$33.23$	&	$33.77$	&	$34.08$	&	$34.36$	&	$33.42$	&	$33.40$	&	$33.32$	&	$33$\\
$17$	&	$35.75$	&	$35.65$	&	$35.17$	&	$35.80$	&	$36.14$	&	$36.42$	&	$35.45$	&	$35.43$	&	$35.37$	&	$35$\\
$18$	&	$37.91$	&	$37.78$	&	$37.12$	&	$37.91$	&	$38.24$	&	$38.38$	&	$37.53$	&	$37.52$	&	$37.46$	&	$37$\\
$19$	&	$39.96$	&	$39.82$	&	$39.06$	&	$39.98$	&	$40.27$	&	$40.36$	&	$39.52$	&	$39.51$	&	$39.46$	&	$39$\\
$20$	&	$42.07$	&	$41.90$	&	$41.14$	&	$42.16$	&	$42.36$	&	$42.38$	&	$41.62$	&	$41.61$	&	$41.57$	&	$41$\\
\midrule
$25$	&	$52.30$	&	$52.05$	&	$51.26$	&	$52.68$	&	$52.54$	&	$52.33$	&	$51.74$	&	$51.75$	&	$51.75$	&	$51$\\
$30$	&	$62.45$	&	$62.15$	&	$61.39$	&	$62.96$	&	$62.66$	&	$62.31$	&	$61.88$	&	$61.90$	&	$61.91$	&	$61$\\
$35$	&	$72.53$	&	$72.19$	&	$71.45$	&	$73.11$	&	$72.70$	&	$72.26$	&	$71.93$	&	$71.96$	&	$71.98$	&	$71$\\
$40$	&	$82.47$	&	$82.14$	&	$81.46$	&	$83.03$	&	$82.65$	&	$82.14$	&	$81.89$	&	$81.94$	&	$81.95$	&	$81$\\
\bottomrule
\end{tabular}
\end{table*} 

\begin{table*}[h!]
\centering
\caption{%
Convergence orders $p_{L_{1}}^{l}$, $p_{L_{2}}^{l}$, $p_{L_{\infty}}^{l}$, calculated in norms $L_{1}$, $L_{2}$, $L_{\infty}$, of \textit{the local solution} $(\mathbf{u}_{L}, \mathbf{v}_{L})$ (represented between the nodes) of the ADER-DG method for the DAE system (\ref{eq:math_pend_dae_ind_3}) of index 1; $N$ is the degree of the basis polynomials $\varphi_{p}$. Orders $p^{l, u}$ are calculated for solution $\mathbf{u}_{L}$; orders $p^{l, v}$ --- for solution $\mathbf{v}_{L}$; orders $p^{l, g}$ --- for the conditions $\mathbf{g} = 0$ on the local solution $(\mathbf{u}_{L}, \mathbf{v}_{L})$. The theoretical value of convergence order $p_{\rm th.}^{l} = N+1$ is applicable for the ADER-DG method for ODE problems and is presented for comparison.
}
\label{tab:conv_orders_local_pend_ind1}
\begin{tabular}{@{}|l|lll|lll|lll|c|@{}}
\toprule
$N$ & $p_{L_{1}}^{l, u}$ & $p_{L_{2}}^{l, u}$ & $p_{L_{\infty}}^{l, u}$ & $p_{L_{1}}^{l, v}$ & $p_{L_{2}}^{l, v}$ & $p_{L_{\infty}}^{l, v}$ & $p_{L_{1}}^{l, g}$ & $p_{L_{2}}^{l, g}$ & $p_{L_{\infty}}^{l, g}$ & $p_{\rm th.}^{l}$ \\
\midrule
$1$	&	$2.25$	&	$2.17$	&	$1.90$	&	$2.13$	&	$2.05$	&	$1.86$	&	$2.16$	&	$2.10$	&	$1.98$	&	$2$\\
$2$	&	$3.00$	&	$2.97$	&	$2.96$	&	$2.98$	&	$2.99$	&	$3.02$	&	$3.00$	&	$2.99$	&	$3.02$	&	$3$\\
$3$	&	$4.03$	&	$4.02$	&	$3.96$	&	$3.97$	&	$3.99$	&	$3.98$	&	$4.00$	&	$4.00$	&	$3.97$	&	$4$\\
$4$	&	$5.05$	&	$5.04$	&	$4.96$	&	$5.01$	&	$4.99$	&	$4.98$	&	$4.98$	&	$4.97$	&	$4.99$	&	$5$\\
$5$	&	$6.00$	&	$5.99$	&	$5.96$	&	$6.09$	&	$5.99$	&	$5.97$	&	$6.02$	&	$6.02$	&	$5.94$	&	$6$\\
$6$	&	$6.90$	&	$6.87$	&	$6.83$	&	$7.16$	&	$7.06$	&	$6.86$	&	$6.89$	&	$6.88$	&	$6.83$	&	$7$\\
$7$	&	$7.94$	&	$7.89$	&	$7.81$	&	$7.75$	&	$7.71$	&	$7.70$	&	$8.00$	&	$7.93$	&	$7.74$	&	$8$\\
$8$	&	$8.93$	&	$8.91$	&	$8.80$	&	$9.04$	&	$8.96$	&	$8.83$	&	$8.82$	&	$8.79$	&	$8.76$	&	$9$\\
$9$	&	$9.95$	&	$9.92$	&	$9.77$	&	$9.91$	&	$9.86$	&	$9.76$	&	$9.96$	&	$9.92$	&	$9.83$	&	$10$\\
$10$	&	$10.92$	&	$10.90$	&	$10.71$	&	$10.83$	&	$10.84$	&	$10.88$	&	$10.91$	&	$10.85$	&	$10.74$	&	$11$\\
$11$	&	$11.91$	&	$11.85$	&	$11.68$	&	$12.00$	&	$11.92$	&	$11.84$	&	$11.83$	&	$11.83$	&	$11.90$	&	$12$\\
$12$	&	$12.83$	&	$12.76$	&	$12.63$	&	$12.84$	&	$12.83$	&	$12.75$	&	$12.88$	&	$12.84$	&	$12.81$	&	$13$\\
$13$	&	$13.81$	&	$13.74$	&	$13.59$	&	$13.92$	&	$13.87$	&	$13.88$	&	$13.81$	&	$13.78$	&	$13.71$	&	$14$\\
$14$	&	$14.80$	&	$14.75$	&	$14.61$	&	$15.00$	&	$14.85$	&	$14.40$	&	$14.82$	&	$14.81$	&	$14.83$	&	$15$\\
$15$	&	$15.82$	&	$15.75$	&	$15.58$	&	$15.91$	&	$15.82$	&	$15.68$	&	$15.82$	&	$15.74$	&	$15.47$	&	$16$\\
$16$	&	$16.81$	&	$16.76$	&	$16.54$	&	$17.02$	&	$16.81$	&	$16.18$	&	$16.73$	&	$16.71$	&	$16.62$	&	$17$\\
$17$	&	$17.80$	&	$17.75$	&	$17.52$	&	$17.93$	&	$17.86$	&	$17.75$	&	$17.74$	&	$17.61$	&	$17.56$	&	$18$\\
$18$	&	$18.80$	&	$18.75$	&	$18.52$	&	$18.97$	&	$18.69$	&	$18.05$	&	$18.72$	&	$18.68$	&	$18.54$	&	$19$\\
$19$	&	$19.79$	&	$19.74$	&	$19.56$	&	$20.04$	&	$19.98$	&	$19.80$	&	$19.71$	&	$19.56$	&	$19.53$	&	$20$\\
$20$	&	$20.78$	&	$20.72$	&	$20.59$	&	$20.88$	&	$20.58$	&	$20.12$	&	$20.73$	&	$20.65$	&	$20.38$	&	$21$\\
\midrule
$25$	&	$25.67$	&	$25.61$	&	$25.61$	&	$26.34$	&	$26.42$	&	$26.27$	&	$25.60$	&	$25.43$	&	$25.20$	&	$26$\\
$30$	&	$30.63$	&	$30.56$	&	$30.48$	&	$30.84$	&	$30.66$	&	$30.42$	&	$30.52$	&	$30.39$	&	$30.21$	&	$31$\\
$35$	&	$35.54$	&	$35.41$	&	$35.23$	&	$36.01$	&	$35.91$	&	$35.77$	&	$35.45$	&	$35.34$	&	$35.35$	&	$36$\\
$40$	&	$40.40$	&	$40.24$	&	$40.04$	&	$41.04$	&	$40.98$	&	$40.96$	&	$40.35$	&	$40.21$	&	$39.72$	&	$41$\\
\bottomrule
\end{tabular}
\end{table*}

As a result of differentiation of the algebraic equation $g_{1} = 0$, the DAE index of the system decreases by 1 and becomes equal to 2, while the algebraic equation $g_{1} = 0$ in the DAE system is replaced by the algebraic equation $g_{2} = 0$, and this procedure is a decrease in the index of the DAE system. As a result of repeated differentiation of the algebraic equation $g_{1} = 0$, i.e. differentiation of the algebraic equation $g_{2} = 0$, the DAE index of the system decreases by 1 more and becomes equal to 1, while the algebraic equation $g_{1} = 0$ in the DAE system is replaced by the algebraic equation $g_{3} = 0$, and this procedure is also a decrease in the DAE index of the system. In this example in the present work, three separate DAE systems were investigated: the DAE system of equations with the algebraic equation $g_{1} = 0$, which was called the DAE system (\ref{eq:math_pend_dae_ind_3}) of index 3; system of equations of DAE of equations with algebraic equation $g_{2} = 0$, which was called DAE system (\ref{eq:math_pend_dae_ind_3}) of index 2; system of equations of DAE of equations with algebraic equation $g_{3} = 0$, which was called DAE system (\ref{eq:math_pend_dae_ind_3}) of index 1. With an increase in the index of the DAE system, a deterioration in the properties of accuracy and convergence of the numerical method ADER-DG with a local DG predictor was expected, this property is well known for other numerical methods for solving DAE systems of equations --- usually a decrease in the empirical convergence order is observed compared to the expected convergence order of the numerical method, which occurs when solving the initial value problem for the ODE system~\cite{Hairer_book_2}. Therefore, interest arose in a quantitative study of the convergence of the numerical solution for the original DAE system and the DAE system with a decreased index.

It is necessary to note an interesting property of the numerical solution --- an algebraic equation $g_{1} = 0$, which is solved numerically explicitly, is satisfied exactly (within the limits of the accuracy of representing real numbers by floating-point numbers) for the solution at the nodes $(\mathbf{u}_{n}, \mathbf{v}_{n})$. Therefore, when separately solving the DAE system (\ref{eq:math_pend_dae_ind_3}) of index 3, index 2 and index 1, the errors were calculated that determine the satisfaction of both algebraic equation $g_{1} = 0$, $g_{2} = 0$ and $g_{3} = 0$ --- it was expected that the numerical solution at the nodes $(\mathbf{u}_{n}, \mathbf{v}_{n})$ would exactly satisfy only that algebraic equation, which is explicitly included in the solved DAE system.

The obtained results of the numerical solution of DAE system (\ref{eq:math_pend_dae_ind_3}) of index 3 are presented in Figs.~\ref{fig:pend_ind3_sol_uv}, \ref{fig:pend_ind3_sol_g_eps}, \ref{fig:pend_ind3_errors} and in Tables~\ref{tab:conv_orders_nodes_pend_ind3}, \ref{tab:conv_orders_local_pend_ind3}. Fig.~\ref{fig:pend_ind3_sol_uv} shows a comparison of the numerical solution at the nodes $(\mathbf{u}_{n}, \mathbf{v}_{n})$, the numerical local solution $(\mathbf{u}_{L}, \mathbf{v}_{L})$ and the exact analytical solution $(\mathbf{u}^{\rm ex}, \mathbf{v}^{\rm ex})$ separately for each differential $\mathbf{u}$ and algebraic $\mathbf{v}$ variable. Fig.~\ref{fig:pend_ind3_sol_g_eps} shows the dependencies of the feasibility of conditions $g_{1} = 0$, $g_{2} = 0$ and  $g_{3} = 0$ on the coordinate $t$, as well as the dependencies of the local errors $\varepsilon_{u}$, $\varepsilon_{v}$, $\varepsilon_{g}$ on the coordinate $t$, which allows one to quantitatively estimate the accuracy of the numerical solution, especially taking into account that the numerical solution obtained by the ADER-DG method with a local DG predictor with a high degree of polynomials $N$ very accurately corresponds to the exact analytical solution, and it is impossible to visually identify the error from the solution plots separately (in Fig.~\ref{fig:pend_ind3_sol_uv}). Fig.~\ref{fig:pend_ind3_errors} shows the dependencies of the global errors $e^{u}$, $e^{v}$, $e^{g}$ of the numerical solution at the nodes $(\mathbf{u}_{n}, \mathbf{v}_{n})$ and the local solution $(\mathbf{u}_{L}, \mathbf{v}_{L})$ on the discretization step ${\Delta t}$, separately for each differential $\mathbf{u}$ and algebraic $\mathbf{v}$ variable and the algebraic equations $\mathbf{g} = \mathbf{0}$, on the basis of which the empirical convergence orders $p$ were calculated.

The numerical solution of the DAE system (\ref{eq:math_pend_dae_ind_3}) of index 3, presented in Fig.~\ref{fig:pend_ind3_sol_uv}, clearly demonstrates the high accuracy achieved using the ADER-DG method with a local DG predictor. The most significant features of the numerical solution are demonstrated in the case of polynomial degree $N = 1$, the numerical solutions for which are shown in Figs.~\ref{fig:pend_ind3_sol_uv}~(\subref{fig:pend_ind3_sol_uv:a1}, \subref{fig:pend_ind3_sol_uv:b1}, \subref{fig:pend_ind3_sol_uv:c1}, \subref{fig:pend_ind3_sol_uv:d1}, \subref{fig:pend_ind3_sol_uv:e1}). In the case of polynomial degree $N = 1$, the dependence of the local solution for differential variables $\mathbf{u}_{L}$ on the coordinate $t$ shows the expected discontinuities of the solution at the grid nodes $t_{n}$, which was already observed in the previous examples considered, however, the behavior of the local solution for the algebraic variable $\mathbf{v}_{L}$ demonstrates a different type of feature --- in the first discretization domain on the grid $\Omega_{0}$, the local solution $\mathbf{v}_{L}$ agrees well with the exact analytical solution $\mathbf{v}^{\rm ex}$, however, in subsequent discretization domains $\Omega_{n}$, the point-wise correspondence between the local solution $\mathbf{v}_{L}$ and the exact analytical solution $\mathbf{v}^{\rm ex}$ significantly worsens. This is due to the condition that is set for obtaining a local solution in the original system of nonlinear algebraic equations of the predictor (\ref{eq:dae_weak_form}) --- despite the point-wise mismatch, the condition of orthogonality of the local error $\boldsymbol{\sigma}_{\rm v}$ to all basis functions $\varphi_{p}$, the expansion of which represents the numerical solution, is satisfied; however, it is necessary to use such a local solution $\mathbf{v}_{L}$ with great ``caution'' as a sufficiently accurate solution between grid nodes. In the previous examples of DAE systems of index 1 systems, such behavior of the local solution $\mathbf{v}_{L}$ for algebraic variables in the case of polynomial degree $N = 1$ was not observed, which is due to the index 1 of the problems in the previous examples; this phenomenon begins to manifest itself significantly for DAE systems of index higher than 1, and it was already observed in the numerical solution of the DAE system (\ref{eq:hess_dae_ind_2}) of index 2 in Example 3, presented in Subsection~\ref{sec:2:ct:ex3}. However, despite this feature of the behavior of the local solution $\mathbf{v}_{L}$ for algebraic variables, the convergence order $p_{j, v}$ for this case turned out to be positive --- with a decrease in the discretization step ${\Delta t}$, the amplitude of discontinuities in the local solution $\mathbf{v}_{L}$ at the grid nodes $t_{n}$ decreases. The presented features of the numerical solution practically disappear in the case of polynomial degrees $N > 1$, in particular, on the presented plots of the numerical solution for polynomials degree $N = 8$ (see Figs.~\ref{fig:pend_ind3_sol_uv}~(\subref{fig:pend_ind3_sol_uv:a2}, \subref{fig:pend_ind3_sol_uv:b2}, \subref{fig:pend_ind3_sol_uv:c2}, \subref{fig:pend_ind3_sol_uv:d2}, \subref{fig:pend_ind3_sol_uv:e2})) and $40$ (see Figs.~\ref{fig:pend_ind3_sol_uv}~(\subref{fig:pend_ind3_sol_uv:a3}, \subref{fig:pend_ind3_sol_uv:b3}, \subref{fig:pend_ind3_sol_uv:c3}, \subref{fig:pend_ind3_sol_uv:d3}, \subref{fig:pend_ind3_sol_uv:e3})). In these cases, the numerical solution at the nodes $(\mathbf{u}_{n}, \mathbf{v}_{n})$ and the local solution $(\mathbf{u}_{L}, \mathbf{v}_{L})$ coincide with the exact analytical solution $(\mathbf{u}^{\rm ex}, \mathbf{v}^{\rm ex})$ with very high accuracy, and the differences that arise on these plots are visually indistinguishable.

To quantitatively determine the accuracy of the local solution $(\mathbf{u}_{L}, \mathbf{v}_{L})$ and the solution at the nodes $(\mathbf{u}_{n}, \mathbf{v}_{n})$, plots of the local errors $|g_{1}|$, $|g_{2}|$, $|g_{3}|$ of satisfaction of the algebraic equations $g_{1} = 0$, $g_{2} = 0$, $g_{3} = 0$ and local errors $\varepsilon_{u}$, $\varepsilon_{v}$, $\varepsilon_{g}$ were constructed, which were presented in Fig.~\ref{fig:pend_ind3_sol_g_eps}. From the presented error dependencies $|g_{1}|$ it was concluded that in the case of a solution at the nodes $(\mathbf{u}_{n}, \mathbf{v}_{n})$, the algebraic equation $g_{1} = 0$, which was explicitly included in the DAE system (\ref{eq:math_pend_dae_ind_3}) of index 3, is satisfied exactly (of course, within the accuracy of representation of real numbers with a floating point). However, in the case of a local solution $(\mathbf{u}_{L}, \mathbf{v}_{L})$, the algebraic equation $g_{1} = 0$ is satisfied with a non-zero error of $|g_{1}| \sim 10^{-7}$--$10^{-1}$ in the case of polynomial degree $N = 1$, $|g_{1}| \sim 10^{-11}$--$10^{-7}$ in the case of $N = 8$, $|g_{1}| \sim 10^{-44}$--$10^{-37}$ in the case of $N = 40$. The presented error dependencies $|g_{2}|$ and $|g_{3}|$, which demonstrate the accuracy of satisfaction of the algebraic equations $g_{2} = 0$ and $g_{3} = 0$, respectively, which are derivatives with respect to $g_{1} = 0$, show non-zero accuracy of the conditions $g_{2} = 0$ and $g_{3} = 0$ both for the solution at the nodes $(\mathbf{u}_{n}, \mathbf{v}_{n})$ and the local solution $(\mathbf{u}_{L}, \mathbf{v}_{L})$. It is necessary to note an important feature of the behavior of the errors $|g_{2}|$ and $|g_{3}|$ --- for almost all values of the coordinate $t$, the errors $|g_{2}|$ and $|g_{3}|$ on the local solution $(\mathbf{u}_{L}, \mathbf{v}_{L})$ are lower than the errors $|g_{2}|$ and $|g_{3}|$ on the solution at the nodes $(\mathbf{u}_{n}, \mathbf{v}_{n})$, respectively. This phenomenon begins to manifest itself significantly for DAE systems of index higher than 1, and is associated with the commensurability of the convergence of orders for the local solution $p^{l, u}$, $p^{l, v}$, $p^{l, g}$ and the solution at the nodes $p^{n, u}$, $p^{n, v}$, $p^{n, g}$ for DAE systems with a high index. It should also be noted that the error $|g_{2}|$ corresponds in order of magnitude to the error $|g_{1}|$ for a local solution ($|g_{1}|$ for the solution at the nodes turned out to be zero), however, the errors $|g_{3}|$ are $10^{2}$--$10^{5}$ times higher than the error $|g_{1}|$, therefore it was revealed that the condition of index 2 is satisfied with less accuracy than the condition of index 3 in the case of the DAE system (\ref{eq:math_pend_dae_ind_3}) of index 3. In general, this result can be considered expected --- with ``removal'' from the original algebraic equation (by the number of differentiations of the original constraint) the errors $|g|$ increases.

The characteristic values of the errors $\varepsilon_{u}$, $\varepsilon_{v}$, $\varepsilon_{g}$ (see Fig.~\ref{fig:pend_ind3_sol_g_eps}) for the polynomial degree $N = 1$ are in the range of $10^{-5}$--$10^{-1}$, for the polynomial degree $N = 8$ --- $10^{-9}$--$10^{-4}$, for the polynomial degree $N = 40$ --- $10^{-39}$--$10^{-33}$. Unlike the previous examples, in this case there are no significant differences in the characteristic values of local errors of the local solution $(\mathbf{u}_{L}, \mathbf{v}_{L})$ and the solution at the nodes $(\mathbf{u}_{n}, \mathbf{v}_{n})$. Fig.~\ref{fig:hess_2_ind2_errors} shows the log-log dependencies of global errors $e$ for the solution at nodes $(\mathbf{u}_{n}, \mathbf{v}_{n})$ and the local solution $(\mathbf{u}_{L}, \mathbf{v}_{L})$ separately for differential and algebraic variables and algebraic equations. The presented dependencies clearly demonstrate the power law $e(\Delta t) \sim (\Delta t)^{p}$, from which the empirical convergence orders $p$ were calculated. Unlike the previous examples, in this case there are no significant differences in the characteristic values of global errors $e$ of the local solution $(\mathbf{u}_{L}, \mathbf{v}_{L})$ and the solution at nodes $(\mathbf{u}_{n}, \mathbf{v}_{n})$.

The calculated empirical convergence orders $p^{n, u}$, $p^{n, v}$, $p^{n, g}$ for the solution at nodes $(\mathbf{u}_{n}, \mathbf{v}_{n})$ are presented in Table~\ref{tab:conv_orders_nodes_pend_ind3}, convergence orders $p^{l, u}$, $p^{l, v}$, $p^{l, g}$ for the local solution $(\mathbf{u}_{L}, \mathbf{v}_{L})$ --- in Table~\ref{tab:conv_orders_local_pend_ind3}. The convergence orders were calculated separately for the norms $L_{1}$, $L_{2}$, $L_{\infty}$ (\ref{eq:norms_def}). The obtained results differ significantly from similar results for previous examples of DAE systems of index 1, and are in good agreement with the results of the DAE system (\ref{eq:hess_dae_ind_2}) of index 2 in Example 3, presented in Subsection~\ref{sec:2:ct:ex3}. The convergence orders of the local solution $p^{l, u}$ and the solution at nodes $p^{n, u}$ for differential variables approximately coincide with each other and with the expected convergence order for the local solution $p_{\rm local} = N+1$, therefore superconvergence does not occur in this case. The convergence orders of the local solution $p^{l, v}$, $p^{l, g}$ and the solution at the nodes $p^{l, v}$, $p^{l, g}$ for algebraic variables and algebraic equations approximately coincide with each other and are one unit less than the convergence orders for differential variables --- their convergence order is approximately equal to $N$. The worsening of the accuracy and convergence properties of the numerical ADER-DG method with a local DG predictor was expected for a DAE index higher than 1 --- this property is well known for other numerical methods for solving DAE systems~\cite{Hairer_book_2}, while in the case of the numerical ADER-DG method this phenomenon manifested itself in the form of a decrease in the empirical convergence orders $p$ compared to the expected convergence orders $p_{\rm nodes} = 2N+1$ and $p_{\rm local} = N+1$ (\ref{eq:expect_orders}) of the numerical method, which occurs when solving the initial value problem for ODE system.

The obtained results of the numerical solution of DAE system (\ref{eq:math_pend_dae_ind_3}) of index 2 are presented in Figs.~\ref{fig:pend_ind2_sol_uv}, \ref{fig:pend_ind2_sol_g_eps}, \ref{fig:pend_ind2_errors} and in Tables~\ref{tab:conv_orders_nodes_pend_ind2}, \ref{tab:conv_orders_local_pend_ind2}. Fig.~\ref{fig:pend_ind2_sol_uv} shows a comparison of the numerical solution at the nodes $(\mathbf{u}_{n}, \mathbf{v}_{n})$, the numerical local solution $(\mathbf{u}_{L}, \mathbf{v}_{L})$ and the exact analytical solution $(\mathbf{u}^{\rm ex}, \mathbf{v}^{\rm ex})$ separately for each differential $\mathbf{u}$ and algebraic $\mathbf{v}$ variable. Fig.~\ref{fig:pend_ind2_sol_g_eps} shows the dependencies of the feasibility of conditions $g_{1} = 0$, $g_{2} = 0$ and  $g_{3} = 0$ on the coordinate $t$, as well as the dependencies of the local errors $\varepsilon_{u}$, $\varepsilon_{v}$, $\varepsilon_{g}$ on the coordinate $t$, which allows us to quantitatively estimate the accuracy of the numerical solution, especially taking into account that the numerical solution obtained by the ADER-DG method with a local DG predictor with a high degree of polynomials $N$ very accurately corresponds to the exact analytical solution, and it is impossible to visually identify the error from the solution plots separately (in Fig.~\ref{fig:pend_ind2_sol_uv}). Fig.~\ref{fig:pend_ind2_errors} shows the dependencies of the global errors $e^{u}$, $e^{v}$, $e^{g}$ of the numerical solution at the nodes $(\mathbf{u}_{n}, \mathbf{v}_{n})$ and the local solution $(\mathbf{u}_{L}, \mathbf{v}_{L})$ on the discretization step ${\Delta t}$, separately for each differential $\mathbf{u}$ and algebraic $\mathbf{v}$ variable and the algebraic equations $\mathbf{g} = \mathbf{0}$, on the basis of which the empirical convergence orders $p$ were calculated.

The numerical solution of the DAE system (\ref{eq:math_pend_dae_ind_3}) of index 2, presented in Fig.~\ref{fig:pend_ind2_sol_uv}, clearly demonstrates the high accuracy achieved using the ADER-DG method with a local DG predictor. The obtained results of comparison of the numerical solution at the nodes and the local solution with the exact analytical solution, in general, correspond to the results from Example 1, presented in Subsection~\ref{sec:2:ct:ex1}, and Example 2, presented in Subsection~\ref{sec:2:ct:ex2}, where test examples of DAE of index 1 systems were solved, but this result differs significantly from the DAE system (\ref{eq:hess_dae_ind_2}) of index 2. To quantitatively determine the accuracy of the local solution $(\mathbf{u}_{L}, \mathbf{v}_{L})$ and the solution at the nodes $(\mathbf{u}_{n}, \mathbf{v}_{n})$, plots of the local errors $|g_{1}|$, $|g_{2}|$, $|g_{3}|$ of satisfaction of the algebraic equations $g_{1} = 0$, $g_{2} = 0$, $g_{3} = 0$ and local errors $\varepsilon_{u}$, $\varepsilon_{v}$, $\varepsilon_{g}$ were constructed, which were presented in Fig.~\ref{fig:pend_ind2_sol_g_eps}. From the presented error dependencies $|g_{1}|$ it was concluded that in the case of a solution at the nodes $(\mathbf{u}_{n}, \mathbf{v}_{n})$, the algebraic equation $g_{2} = 0$, which was explicitly included in the DAE system (\ref{eq:math_pend_dae_ind_3}) of index 2, is satisfied exactly (of course, within the accuracy of representation of real numbers with a floating point). However, in the case of a local solution $(\mathbf{u}_{L}, \mathbf{v}_{L})$, the algebraic equation $g_{2} = 0$ is satisfied with a non-zero error of $|g_{2}| \sim 10^{-6}$--$10^{-2}$ in the case of polynomial degree $N = 1$, $|g_{2}| \sim 10^{-10}$--$10^{-6}$ in the case of $N = 8$, $|g_{2}| \sim 10^{-43}$--$10^{-37}$ in the case of $N = 40$. The presented error dependencies $|g_{1}|$ and $|g_{3}|$, which demonstrate the accuracy of satisfaction of the algebraic equations $g_{1} = 0$ and $g_{3} = 0$, respectively, show non-zero accuracy of the conditions $g_{2} = 0$ and $g_{3} = 0$ both for the solution at the nodes $(\mathbf{u}_{n}, \mathbf{v}_{n})$ and the local solution $(\mathbf{u}_{L}, \mathbf{v}_{L})$. It is necessary to note an important feature of the behavior of the errors $|g_{3}|$ --- for almost all values of the coordinate $t$, the errors and $|g_{3}|$ on the local solution $(\mathbf{u}_{L}, \mathbf{v}_{L})$ are lower than the errors $|g_{2}|$ and $|g_{3}|$ on the solution at the nodes $(\mathbf{u}_{n}, \mathbf{v}_{n})$, respectively. This phenomenon begins to manifest itself significantly for DAE systems of index higher than 1, and is associated with the commensurability of the convergence of orders for the local solution $p^{l, v}$, $p^{l, g}$ and the solution at the nodes $p^{n, v}$, $p^{n, g}$ for DAE systems with a high index. However, in the case of $g_{2} = 0$, the error $|g_{2}|$ for the solution at the nodes $(\mathbf{u}_{n}, \mathbf{v}_{n})$ is significantly smaller than the error $|g_{2}|$ for the local solution $(\mathbf{u}_{L}, \mathbf{v}_{L})$: by $10^{0}$--$10^{2}$ times in the case of polynomial degree $N = 1$, by $10^{6}$--$10^{7}$ times in the case of polynomial degree $N = 8$, by $10^{32}$--$10^{36}$ times in the case of polynomial degree $N = 40$ --- this result differs significantly from the DAE system (\ref{eq:math_pend_dae_ind_3}) of index 3 and the DAE system (\ref{eq:hess_dae_ind_2}) of index 2, but corresponds to the results of Examples 1 and 2, presented in Subsections~\ref{sec:2:ct:ex1} and \ref{sec:2:ct:ex2}, respectively, where the DAE systems of index 1 were studied. It should also be noted that the error $|g_{2}|$ corresponds in order of magnitude to the errors $|g_{1}|$ and $|g_{3}|$ for a local solution ($|g_{2}|$ for the solution at the nodes turned out to be zero). In general, this result can be considered expected --- with ``removal'' from the original algebraic equation (by the number of differentiations or integrations of the original constraint) the errors $|g|$ increases. Another important difference is the difference in the errors $\varepsilon_{u}$ of the local solution $\mathbf{u}_{n}$ and the solution at the nodes $\mathbf{u}_{L}$ for differential variables $\mathbf{u}$, which reaches $10^{1}$--$10^{2}$ times in the case of polynomial degrees of $N = 1$, $10^{2}$--$10^{7}$ times in the case of polynomial degrees of $N = 8$ and $10^{34}$--$10^{37}$ times in the case of polynomial degrees of $N = 40$, while there is no such difference for the errors $\varepsilon_{v}$, $\varepsilon_{g}$ of the local solution and the solution at the nodes for algebraic variables $\mathbf{v}$ and algebraic equations $\mathbf{g} = \mathbf{0}$. Fig.~\ref{fig:pend_ind2_errors} shows the log-log dependencies of global errors $e$ for the solution at nodes $(\mathbf{u}_{n}, \mathbf{v}_{n})$ and the local solution $(\mathbf{u}_{L}, \mathbf{v}_{L})$ separately for differential and algebraic variables and algebraic equations. The presented dependencies clearly demonstrate the power law $e(\Delta t) \sim (\Delta t)^{p}$, from which the empirical convergence orders $p$ were calculated. 

The calculated empirical convergence orders $p^{n, u}$, $p^{n, v}$, $p^{n, g}$ for the solution at nodes $(\mathbf{u}_{n}, \mathbf{v}_{n})$ are presented in Table~\ref{tab:conv_orders_nodes_pend_ind2}, convergence orders $p^{l, u}$, $p^{l, v}$, $p^{l, g}$ for the local solution $(\mathbf{u}_{L}, \mathbf{v}_{L})$ --- in Table~\ref{tab:conv_orders_local_pend_ind2}. The convergence orders were calculated separately for the norms $L_{1}$, $L_{2}$, $L_{\infty}$ (\ref{eq:norms_def}). The obtained empirical convergence orders $p^{l, u}$, $p^{l, v}$, $p^{l, g}$ for the local solution $(\mathbf{u}_{L}, \mathbf{v}_{L})$ approximately correspond to the expected values for the local solution $p_{\rm local} = N+1$ (\ref{eq:expect_orders}). However, the solution at the nodes  $\mathbf{u}_{n}$ for differential variables demonstrates unexpectedly high empirical convergence orders $p^{n, u}$ and demonstrates superconvergence characteristic of the numerical method ADER-DG with a local DG predictor --- from $2N+1$ for polynomial degrees $N \leqslant 25$ to $2N$ in the case of polynomial degrees $N > 25$. At the same time, the empirical convergence orders $p^{n, v}$, $p^{n, g}$ of the solution at the nodes for algebraic variables $\mathbf{v}_{n}$ and algebraic equations $\mathbf{g} = \mathbf{0}$ correspond to the expected values for the local solution --- superconvergence is not observed.

The obtained results of the numerical solution of DAE system (\ref{eq:math_pend_dae_ind_3}) of index 1 are presented in Figs.~\ref{fig:pend_ind1_sol_uv}, \ref{fig:pend_ind1_sol_g_eps}, \ref{fig:pend_ind1_errors} and in Tables~\ref{tab:conv_orders_nodes_pend_ind1}, \ref{tab:conv_orders_local_pend_ind1}. Fig.~\ref{fig:pend_ind1_sol_uv} shows a comparison of the numerical solution at the nodes $(\mathbf{u}_{n}, \mathbf{v}_{n})$, the numerical local solution $(\mathbf{u}_{L}, \mathbf{v}_{L})$ and the exact analytical solution $(\mathbf{u}^{\rm ex}, \mathbf{v}^{\rm ex})$ separately for each differential $\mathbf{u}$ and algebraic $\mathbf{v}$ variable. Fig.~\ref{fig:pend_ind1_sol_g_eps} shows the dependencies of the feasibility of conditions $g_{1} = 0$, $g_{2} = 0$ and  $g_{3} = 0$ on the coordinate $t$, as well as the dependencies of the local errors $\varepsilon_{u}$, $\varepsilon_{v}$, $\varepsilon_{g}$ on the coordinate $t$, which allows us to quantitatively estimate the accuracy of the numerical solution, especially taking into account that the numerical solution obtained by the ADER-DG method with a local DG predictor with a high degree of polynomials $N$ very accurately corresponds to the exact analytical solution, and it is impossible to visually identify the error from the solution plots separately (in Fig.~\ref{fig:pend_ind1_sol_uv}). Fig.~\ref{fig:pend_ind1_errors} shows the dependencies of the global errors $e^{u}$, $e^{v}$, $e^{g}$ of the numerical solution at the nodes $(\mathbf{u}_{n}, \mathbf{v}_{n})$ and the local solution $(\mathbf{u}_{L}, \mathbf{v}_{L})$ on the discretization step ${\Delta t}$, separately for each differential $\mathbf{u}$ and algebraic $\mathbf{v}$ variable and the algebraic equations $\mathbf{g} = \mathbf{0}$, on the basis of which the empirical convergence orders $p$ were calculated.

The numerical solution of the DAE system (\ref{eq:math_pend_dae_ind_3}) of index 1, presented in Fig.~\ref{fig:pend_ind1_sol_uv}, clearly demonstrates the high accuracy achieved using the ADER-DG method with a local DG predictor. The obtained numerical solution at the nodes $(\mathbf{u}_{n}, \mathbf{v}_{n})$ and the local solution $(\mathbf{u}_{L}, \mathbf{v}_{L})$, in its qualitative and quantitative properties when compared with the exact analytical solution $(\mathbf{u}^{\rm ex}, \mathbf{v}^{\rm ex})$, is similar to the numerical solution of the DAE system of index 1, which was studied in Example 1, presented in Subsection~\ref{sec:2:ct:ex1}. Artifacts of the local solution $\mathbf{v}_{L}$, which arose in the solution of the DAE system (\ref{eq:math_pend_dae_ind_3}) of index 3, presented in this Subsection above, and the solution of the DAE system (\ref{eq:hess_dae_ind_2}) of index 2, which was studied in Example 3, presented in Subsection~\ref{sec:2:ct:ex3}, in the case of polynomials of degree $N = 1$, do not arise in this case. 

To quantitatively determine the accuracy of the local solution $(\mathbf{u}_{L}, \mathbf{v}_{L})$ and the solution at the nodes $(\mathbf{u}_{n}, \mathbf{v}_{n})$, plots of the local errors $|g_{1}|$, $|g_{2}|$, $|g_{3}|$ of satisfaction of the algebraic equations $g_{1} = 0$, $g_{2} = 0$, $g_{3} = 0$ and local errors $\varepsilon_{u}$, $\varepsilon_{v}$, $\varepsilon_{g}$ were constructed, which were presented in Fig.~\ref{fig:pend_ind1_sol_g_eps}. The obtained dependencies of the errors of the numerical solution on time $t$, in terms of qualitative and quantitative properties when compared with the exact analytical solution $(\mathbf{u}^{\rm ex}, \mathbf{v}^{\rm ex})$, also correspond to the properties of the dependencies of the error of the numerical solution of the DAE system of index 1, which was studied in Example 1, presented in Subsection~\ref{sec:2:ct:ex1}. Fig.~\ref{fig:pend_ind1_errors} shows the log-log dependencies of global errors $e$ for the solution at nodes $(\mathbf{u}_{n}, \mathbf{v}_{n})$ and the local solution $(\mathbf{u}_{L}, \mathbf{v}_{L})$, separately for each differential $\mathbf{u}$ and algebraic $\mathbf{v}$ variable and the algebraic equations $\mathbf{g} = \mathbf{0}$. The presented dependencies clearly demonstrate the power law $e(\Delta t) \sim (\Delta t)^{p}$, from which the empirical convergence orders $p$ were calculated.

The calculated empirical convergence orders $p^{n, u}$, $p^{n, v}$, $p^{n, g}$ for the solution at nodes $(\mathbf{u}_{n}, \mathbf{v}_{n})$ are presented in Table~\ref{tab:conv_orders_nodes_pend_ind1}, convergence orders $p^{l, u}$, $p^{l, v}$, $p^{l, g}$ for the local solution $(\mathbf{u}_{L}, \mathbf{v}_{L})$ --- in Table~\ref{tab:conv_orders_local_pend_ind1}. The convergence orders were calculated separately for the norms $L_{1}$, $L_{2}$, $L_{\infty}$ (\ref{eq:norms_def}). Also for comparison, the expected values of the convergence orders $p$, which the ADER-DG numerical method with the local DG predictor provides for solving the initial value problem for the ODE system, are presented in Table~\ref{tab:conv_orders_nodes_pend_ind1}: $p_{\rm nodes} = 2N+1$, $p_{\rm local} = N+1$ (\ref{eq:expect_orders}). Comparison of the obtained convergence orders $p^{n, u}$, $p^{n, v}$, $p^{n, g}$ for the solution at the nodes clearly demonstrates not only the expected superconvergence $2N+1$, but also higher values of the convergence orders --- the obtained empirical values are often greater than the expected values by $1.0$--$2.0$. Comparison of the obtained convergence orders $p^{l, u}$, $p^{l, v}$, $p^{l, g}$ for the local solution demonstrates similar behavior --- the obtained empirical values of the convergence orders are in good agreement or exceed the expected values $p_{\rm local} = N+1$ (\ref{eq:expect_orders}). These results differ significantly from the results for the DAE system (\ref{eq:math_pend_dae_ind_3}) of indices 3 and 2, presented in this Subsection above, as well as from the results for the DAE system (\ref{eq:hess_dae_ind_2}) of indices 2 and 1, which was studied in Example 3, presented in Subsection~\ref{sec:2:ct:ex3}.

Therefore, it can be concluded that the presented numerical solution of the DAE system (\ref{eq:math_pend_dae_ind_3}) of index 3 demonstrates a significant decrease in the empirical convergence orders $p$ compared to the expected values (\ref{eq:expect_orders}). Decreasing the DAE system (\ref{eq:math_pend_dae_ind_3}) index to 2 demonstrates significantly higher empirical convergence orders $p^{n, u}$ of the solution at the nodes $\mathbf{u}_{n}$ for differential variables $\mathbf{u}$, however, the empirical convergence orders of the solution at the nodes $p^{n, v}$, $p^{n, g}$ for algebraic variables $\mathbf{v}$ and algebraic equations $\mathbf{g} = \mathbf{0}$ are significantly decreased, to values of $N+1$. This result differs from the result for the system (\ref{eq:hess_dae_ind_2}) of index 1, which was studied in Example 3, presented in Subsection~\ref{sec:2:ct:ex3}. Decreasing the DAE system (\ref{eq:math_pend_dae_ind_3}) index to 1 demonstrates all the main convergence properties of the ADER-DG method with a local DG predictor, which were observed for the DAE system of index 1, which was studied in Example 1, presented in Subsection~\ref{sec:2:ct:ex1}, and for the initial value problem for the ODE system.

\subsubsection{Example 5: double pendulum}
\label{sec:2:ct:ex5}

The fifth example (which contains three sub-examples) of application of the numerical method ADER-DG with local DG predictor consisted in solving a well-known problem of double pendulum, which is well known in the field of DAE systems, and which in the original formulation of the problem is a DAE system of index 3:
\begin{equation}\label{eq:math_dpend_dae_ind_3}
\begin{split}
&\ddot{x}_{1} = -\lambda_{1} x_{1} - \lambda_{2} (x_{1} - x_{2}),\hspace{10mm} \ddot{x}_{2} = \lambda_{2} (x_{1} - x_{2}),\\
&\ddot{y}_{1} = -\lambda_{1} y_{1} - \lambda_{2} (y_{1} - y_{2}) - g,\hspace{5.1mm} \ddot{y}_{2} = \lambda_{2} (y_{1} - y_{2}) - g,\\
&g_{11} = x_{1}^{2} + y_{1}^{2} - 1 = 0,\\
&g_{12} = (x_{1} - x_{2})^{2} + (y_{1} - y_{2})^{2} - 1 = 0,\\
&g_{21} = x_{1}\dot{x}_{1} + y_{1}\dot{y}_{1} = 0,\\
&g_{22} = (x_{1} - x_{2})(\dot{x}_{1} - \dot{x}_{2}) + (y_{1} - y_{2})(\dot{y}_{1} - \dot{y}_{2}) = 0,\\
&g_{31} = \dot{x}_{1}^{2} + \dot{y}_{1}^{2} - gy_{1} - \lambda_{1}(x_{1}^{2} + y_{1}^{2})\\
&\phantom{g_{31} }- \lambda_{2}\left[x_{1}(x_{1} - x_{2}) + y_{1}(y_{1} - y_{2})\right] = 0,\\
&g_{32} = (\dot{x}_{1} - \dot{x}_{2})^{2} + (\dot{y}_{1} - \dot{y}_{2})^{2}\\
&\phantom{g_{32} }- (x_{1} - x_{2})\left[\lambda_{1}x_{1} + 2\lambda_{2}(x_{1} - x_{2})\right]\\
&\phantom{g_{32} }- (y_{1} - y_{2})\left[\lambda_{1}y_{1} + 2\lambda_{2}(y_{1} - y_{2})\right] = 0,\\
&x_{1}(0) = \sin(\phi_{1,0}),\quad y_{1}(0) = -\cos(\phi_{1,0}),\\
&x_{2}(0) = \sin(\phi_{1,0}) + \sin(\phi_{2,0}),\\
&y_{2}(0) = -(\cos(\phi_{1,0}) + \cos(\phi_{2,0})),\\
&\dot{x}_{1}(0) = 0,\quad \dot{y}_{1}(0) = 0,\quad \dot{x}_{2}(0) = 0,\quad \dot{y}_{2}(0) = 0,\\
&\lambda_{1}(0) = \frac{-2gy_{1}(0)}{2 - \left[1 - x_{1}(0)x_{2}(0) - y_{1}(0)y_{2}(0)\right]^{2}},\\
&\lambda_{2}(0) = \frac{gy_{1}(0)\left[1 - x_{1}(0)x_{2}(0) - y_{1}(0)y_{2}(0)\right]}{2 - \left[1 - x_{1}(0)x_{2}(0) - y_{1}(0)y_{2}(0)\right]^{2}},\\
\end{split}
\end{equation}
where the initial condition was given by the initial angles of the pendulum $\phi_{1, 0} = 0.25\pi$ and $\phi_{2, 0} = 0.30\pi$. The domain of definition was chosen as $[0,\, 20]$. This test example is the second example in this paper for which the numerical solution of the DAE of index 3 was investigated. It should be noted that in this case, in the positions $(\phi_{1}, \phi_{2}) = (\pi, \pi)$, $(\pi, 0)$, $(0, \pi)$ and $(\pi, \pi)$, the DAE system (\ref{eq:math_dpend_dae_ind_3}) becomes degenerate, and in the strict sense, in the general case, the index of the DAE system is undefined for the entire phase space, however, outside these points of the phase space, the index of the DAE system (\ref{eq:math_dpend_dae_ind_3}) is always equal to 3. The presented DAE system of index 3 contains only two algebraic equations (constraints) $g_{11} = 0$ and $g_{12} = 0$.

The equations of motion (\ref{eq:math_dpend_dae_ind_3}) of a double pendulum do not have an exact analytical solution expressed in elementary functions or classical special functions. Therefore, to analyze the accuracy and convergence of the numerical solution obtained using the ADER-DG method with a local DG predictor, a reference solution was used, for which a high-precision numerical solution was chosen. The reference solution was obtained for the dynamic dependencies of the angles $\phi_{1}(t)$ and $\phi_{2}(t)$, which determine the deviation of the double pendulum position from the position of stable mechanical equilibrium:
\begin{equation}\label{eq:math_dpend_xy_by_phis}
\begin{split}
&x_{1} = \sin(\phi_{1}(t)),\quad
 y_{1} = -\cos(\phi_{1}(t)),\\
&\dot{x}_{1} = \dot{\phi}_{1}(t)\cos(\phi_{1}(t)),\quad
 \dot{y}_{1} = \dot{\phi}_{1}(t)\sin(\phi_{1}(t)),\\
&x_{2} = \sin(\phi_{1}(t)) + \sin(\phi_{2}(t)),\\
&y_{2} = -(\cos(\phi_{1}(t)) + \cos(\phi_{2}(t))),\\
&\dot{x}_{2} = \dot{\phi}_{1}(t)\cos(\phi_{1}(t)) + \dot{\phi}_{2}(t)\cos(\phi_{2}(t)),\\
&\dot{y}_{2} = \dot{\phi}_{1}(t)\sin(\phi_{1}(t)) + \dot{\phi}_{2}(t)\sin(\phi_{2}(t)),\\
\end{split}
\end{equation}
and the expressions for the reaction forces $\lambda_{1}$ and $\lambda_{2}$ of the constraints were obtained in the following form:
\begin{equation}
\begin{split}
&\lambda_{1} = \frac{
	2(\dot{x}_{1}^{2} + \dot{y}_{1}^{2} - gy_{1}) - 
	\left[(\dot{x}_{1} - \dot{x}_{2})^{2} + (\dot{y}_{1} - \dot{y}_{2})^{2}\right]\left[1 - x_{1}x_{2} - y_{1}y_{2}\right]
}{2 - \left[1 - x_{1}x_{2} - y_{1}y_{2}\right]^{2}},\\
&\lambda_{2} = \frac{
	\left[(\dot{x}_{1} - \dot{x}_{2})^{2} + (\dot{y}_{1} - \dot{y}_{2})^{2}\right] -
	\left[1 - x_{1}x_{2} - y_{1}y_{2}\right](\dot{x}_{1}^{2} + \dot{y}_{1}^{2} - gy_{1})
}{2 - \left[1 - x_{1}x_{2} - y_{1}y_{2}\right]^{2}},\\
\end{split}
\end{equation}
where the expressions for the coordinates and velocities were then substituted from the previous relations (\ref{eq:math_dpend_xy_by_phis}). It should be noted that the problem solved in this Subsection is nonlinear --- the presented solution defines the law of motion of a mathematical pendulum for an arbitrary amplitude of oscillations (limited only by the zero value of the initial velocity $\dot{\phi}(0) = 0$), and not only for small values of amplitude, when it would be possible to limit oneself only to the harmonic approximation of small oscillations. The equations of motion for the angles were obtained within the Lagrangian formalism of classical mechanics:
\begin{equation}\label{eq:math_dpend_ode}
\begin{split}
&\ddot{\phi}_{1} = -\frac{
	\sin(\phi_{1}\hspace{-1.5mm} - \hspace{-1mm}\phi_{2})\dot{\phi}_{2}^{2} + 2\sin(\phi_{1}) +
	\cos(\phi_{1}\hspace{-1.5mm} - \hspace{-1mm}\phi_{2})\left(\sin(\phi_{1}\hspace{-1.5mm} - \hspace{-1mm}\phi_{2})\dot{\phi}_{1}^{2} - \sin(\phi_{2})\right)
}{2 - \cos^{2}(\phi_{1}\hspace{-1.5mm} - \hspace{-1mm}\phi_{2})},\\
&\ddot{\phi}_{2} = \frac{
	\cos(\phi_{1}\hspace{-1.5mm} - \hspace{-1mm}\phi_{2})\left(\sin(\phi_{1}\hspace{-1.5mm} - \hspace{-1mm}\phi_{2})\dot{\phi}_{2}^{2} + 2\sin(\phi_{1})\right) +
	2\left(\sin(\phi_{1}\hspace{-1.5mm} - \hspace{-1mm}\phi_{2})\dot{\phi}_{1}^{2} - \sin(\phi_{2})\right)
}{2 - \cos^{2}(\phi_{1}\hspace{-1.5mm} - \hspace{-1mm}\phi_{2})},\\[2mm]
&\phi_{1}(0) = \phi_{1,0},\quad \phi_{2}(0) = \phi_{2,0},\quad \dot{\phi}_{1}(0) = 0,\quad \dot{\phi}_{2}(0) = 0.
\end{split}
\end{equation}
The numerical solution for the initial value problem for the ODE system was obtained using the ADER-DG method with a local DG predictor~\cite{ader_dg_ivp_ode} with a polynomial degree of $N = 60$ and a number of discretization steps of $200$. Based on the analysis of the convergence of the numerical solution over a set of discretization steps of $50$, $100$, $150$ and $200$, it was found that the resulting final numerical solution, called the reference solution $(\mathbf{u}^{\rm ref},\, \mathbf{v}^{\rm ref})$, has an error at the final point $t_{f} = 20$ of no more than $10^{-400}$, which was sufficient for the analysis of the numerical solution of the DAE system studied in this Subsection.

The problem (\ref{eq:math_pend_dae_ind_3}) was rewritten in a form consistent with the formulation of the original DAE system (\ref{eq:dae_chosen_form}):
\begin{equation}
\begin{split}
&\frac{du_{1}}{dt} = u_{5},\qquad
 \frac{du_{2}}{dt} = u_{6},\qquad\quad\,\,
 \frac{du_{3}}{dt} = u_{7},\qquad
 \frac{du_{4}}{dt} = u_{8},\\
&\frac{du_{5}}{dt} = -v_{1} u_{1} - v_{2} (u_{1} - u_{3}),\qquad
 \frac{du_{6}}{dt} = -v_{1} u_{2} - v_{2} (u_{2} - u_{4}) - g,\\
&\frac{du_{7}}{dt} = v_{2} (u_{1} - u_{3}),\qquad\qquad\qquad
 \frac{du_{8}}{dt} = v_{2} (u_{2} - u_{4}) - g,\\
&g_{11} = u_{1}^{2} + u_{2}^{2} - 1 = 0,\\
&g_{12} = (u_{1} - u_{3})^{2} + (u_{2} - u_{4})^{2} - 1 = 0,\\
&g_{21} = u_{1}u_{5} + u_{2}u_{6} = 0,\\
&g_{22} = (u_{1} - u_{3})(u_{5} - u_{7}) + (u_{2} - u_{4})(u_{6} - u_{8}) = 0,\\
&g_{31} = u_{5}^{2} + u_{6}^{2} - gu_{2} - v_{1}(u_{1}^{2} + u_{2}^{2}) - v_{2}\left[u_{1}(u_{1} - u_{3}) + u_{2}(u_{2} - u_{4})\right] = 0,\\
&g_{32} = (u_{5} - u_{7})^{2} + (u_{6} - u_{8})^{2} - (u_{1} - u_{3})\left[v_{1}u_{1} + 2v_{2}(u_{1} - u_{3})\right]\\
&\phantom{g_{32} }- (u_{2} - u_{4})\left[v_{1}u_{2} + 2v_{2}(u_{2} - u_{4})\right] = 0,\\
\end{split}
\end{equation}
where $\mathbf{u} = [x_{1},\, y_{1},\, x_{2},\, y_{2},\, \dot{x}_{1},\, \dot{y}_{1},\, \dot{x}_{2},\, \dot{y}_{2}]^{T}$, $\mathbf{v} = [\lambda_{1},\, \lambda_{2}]$. The value $g = 1$ was chosen in the calculations. The domain of definition $[0,\, 20]$ of the desired functions $\mathbf{u}$ and $\mathbf{v}$ was discretized into $L = 30$, $32$, $34$, $36$, $38$, $40$ discretization domains $\Omega_{n}$ in the case of polynomial degrees $N \leqslant 5$ and into $L = 10$, $12$, $14$, $16$, $18$, $20$ discretization domains $\Omega_{n}$ in the case of polynomial degrees $N > 5$. The discretization step ${\Delta t}_{n}$ was chosen equal ${\Delta t} = 10/(L-1)$ for all discretization domains $\Omega_{n}$ on the grid, which was done to be able to calculate the empirical convergence orders $p$. Therefore, the empirical convergence orders $p$ were calculated by least squares approximation of the dependence of the global error $e$ on the grid discretization step ${\Delta t}$ at 6 data points.

As a result of differentiation of the algebraic equations $g_{11} = 0$ and $g_{12} = 0$, the DAE index of the system decreases by 1 and becomes equal to 2, while the algebraic equations $g_{11} = 0$ and $g_{12} = 0$ in the DAE system are replaced by the algebraic equations $g_{21} = 0$ and $g_{22} = 0$, and this procedure is a decrease in the index of the DAE system. As a result of repeated differentiation of the algebraic equations $g_{11} = 0$ and $g_{12} = 0$, i.e. differentiation of the algebraic equations $g_{21} = 0$ and $g_{22} = 0$, the DAE index of the system decreases by 1 more and becomes equal to 1, while the algebraic equations $g_{11} = 0$ and $g_{12} = 0$ in the DAE system are replaced by the algebraic equations $g_{31} = 0$ and $g_{32} = 0$, and this procedure is also a decrease in the DAE index of the system. In this example in the present work, three separate DAE systems were investigated: the DAE system of equations with the algebraic equations $g_{11} = 0$ and $g_{12} = 0$, which was called the DAE system (\ref{eq:math_dpend_dae_ind_3}) of index 3; system of equations of DAE of equations with algebraic equations $g_{21} = 0$ and $g_{22} = 0$, which was called DAE system (\ref{eq:math_dpend_dae_ind_3}) of index 2; system of equations of DAE of equations with algebraic equations $g_{31} = 0$ and $g_{32} = 0$, which was called DAE system (\ref{eq:math_dpend_dae_ind_3}) of index 1. With an increase in the index of the DAE system, a deterioration in the properties of accuracy and convergence of the numerical method ADER-DG with a local DG predictor was expected, this property is well known for other numerical methods for solving DAE systems of equations --- usually a decrease in the empirical convergence order is observed compared to the expected convergence order of the numerical method, which occurs when solving the initial value problem for the ODE system~\cite{Hairer_book_2}. Therefore, interest arose in a quantitative study of the convergence of the numerical solution for the original DAE system and the DAE system with a decreased index.

\begin{figure}[h!]
\captionsetup[subfigure]{%
	position=bottom,
	font+=smaller,
	textfont=normalfont,
	singlelinecheck=off,
	justification=raggedright
}
\centering
\begin{subfigure}{0.240\textwidth}
\includegraphics[width=\textwidth]{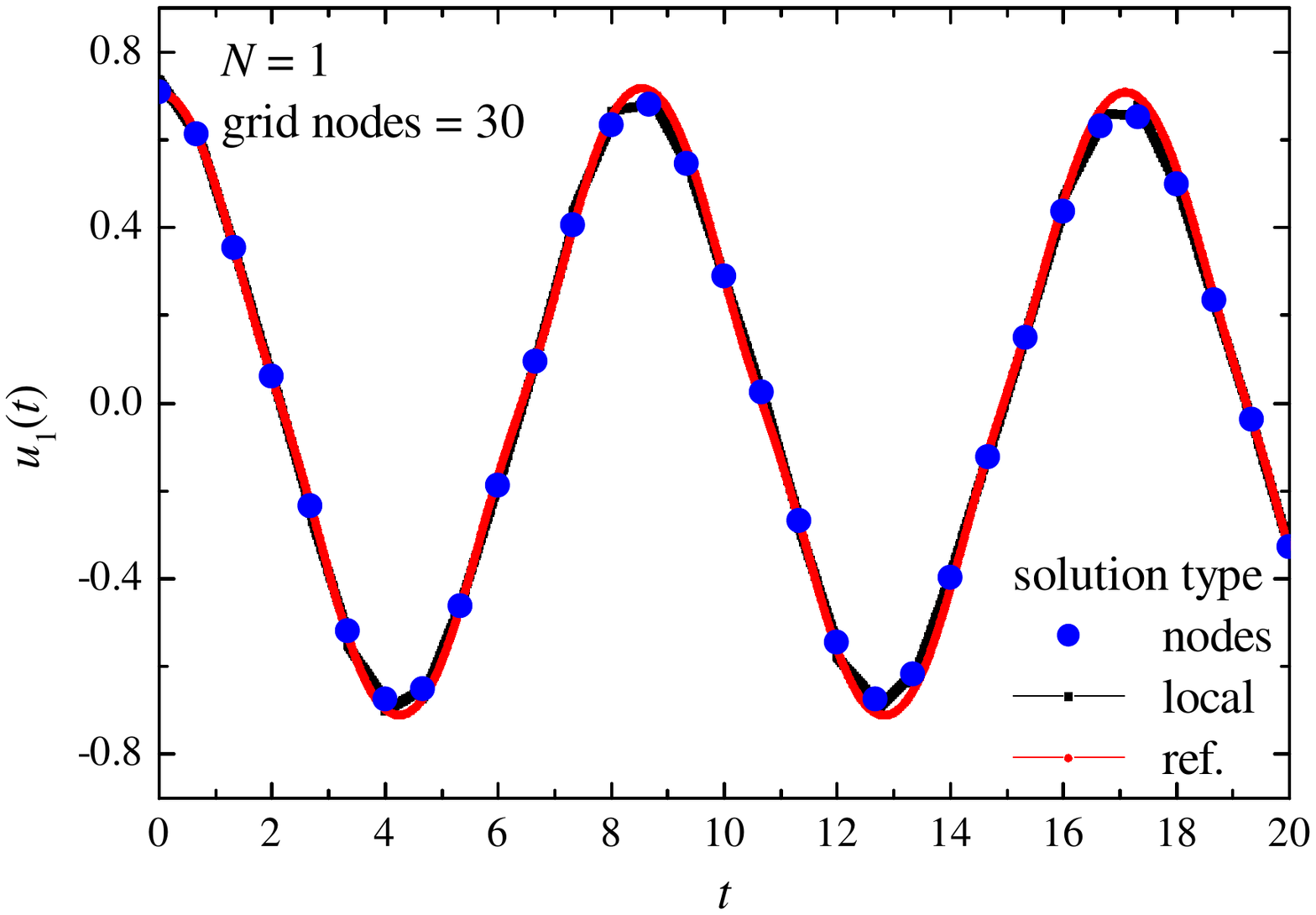}
\vspace{-8mm}\caption{\label{fig:dpend_ind3_sols_u:a1}}
\end{subfigure}
\begin{subfigure}{0.240\textwidth}
\includegraphics[width=\textwidth]{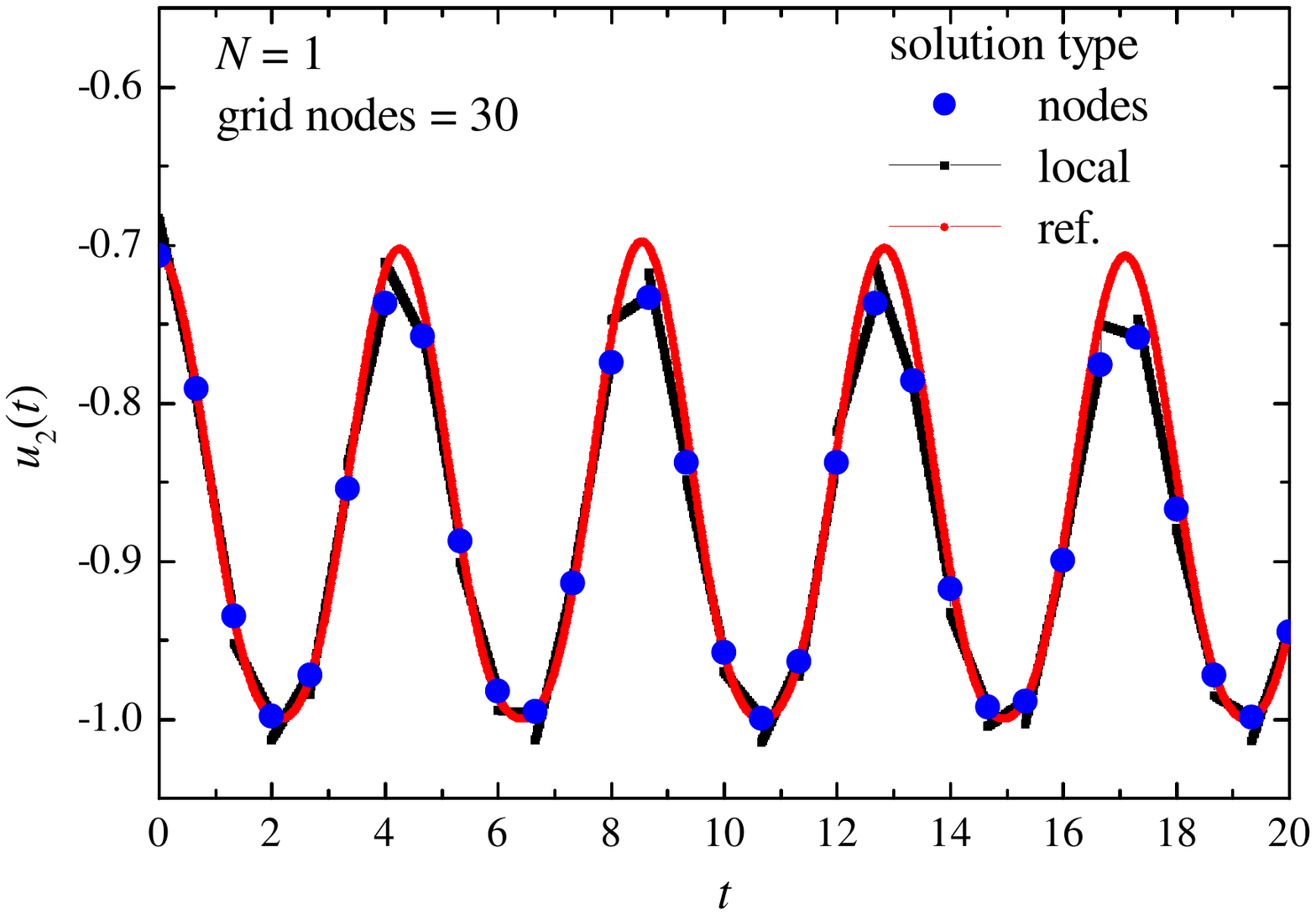}
\vspace{-8mm}\caption{\label{fig:dpend_ind3_sols_u:a2}}
\end{subfigure}
\begin{subfigure}{0.240\textwidth}
\includegraphics[width=\textwidth]{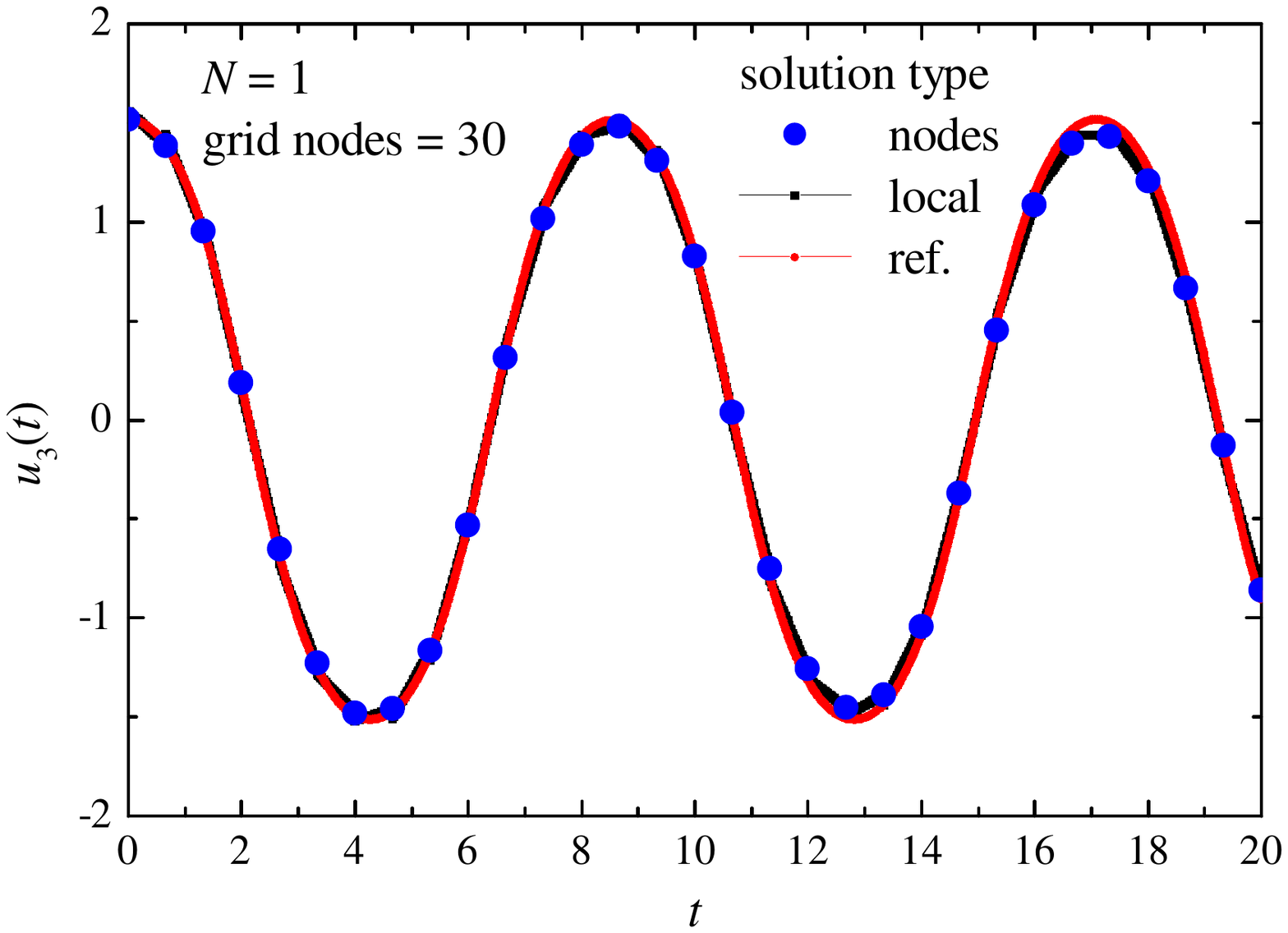}
\vspace{-8mm}\caption{\label{fig:dpend_ind3_sols_u:a3}}
\end{subfigure}
\begin{subfigure}{0.240\textwidth}
\includegraphics[width=\textwidth]{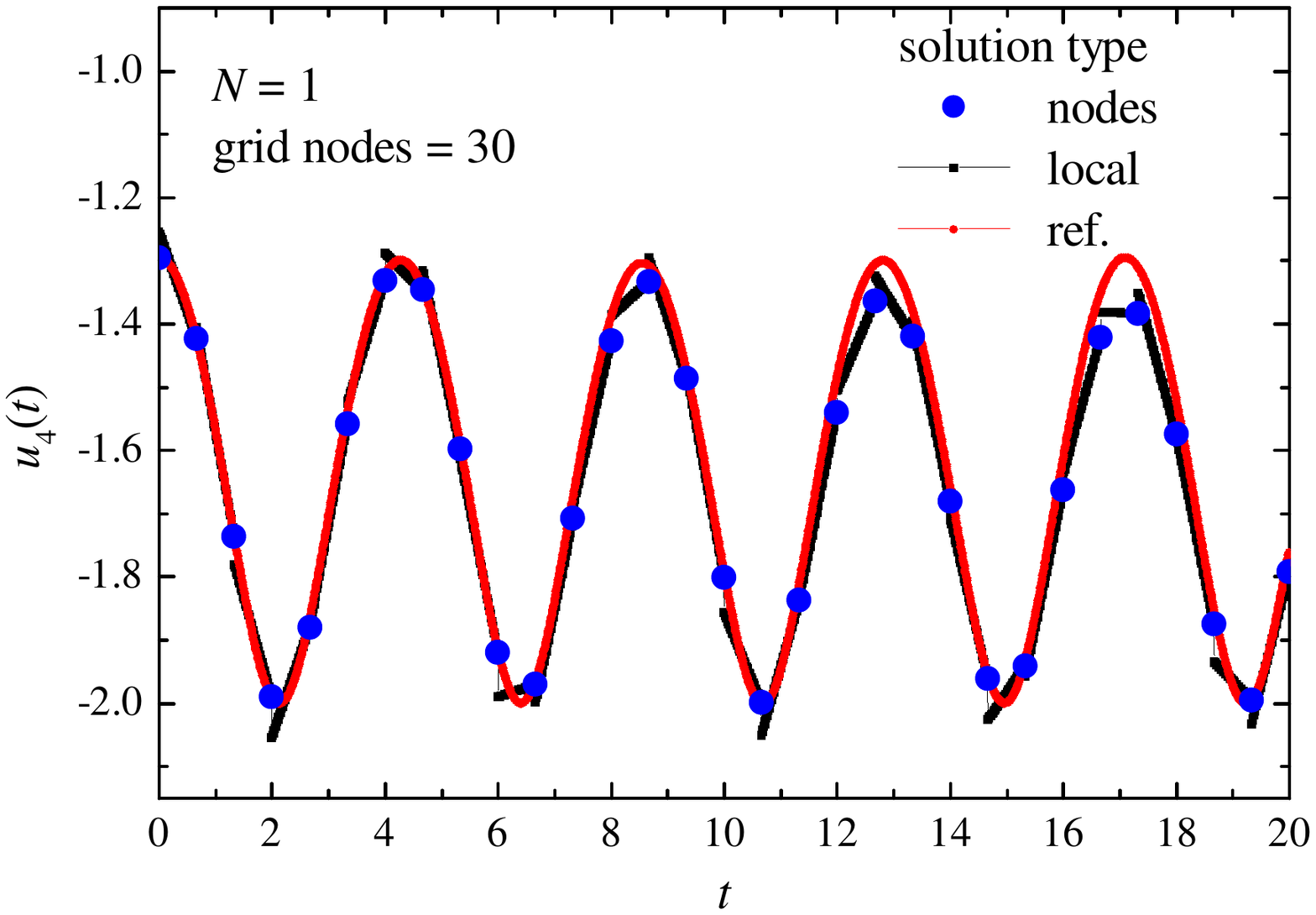}
\vspace{-8mm}\caption{\label{fig:dpend_ind3_sols_u:a4}}
\end{subfigure}\\[2mm]
\begin{subfigure}{0.240\textwidth}
\includegraphics[width=\textwidth]{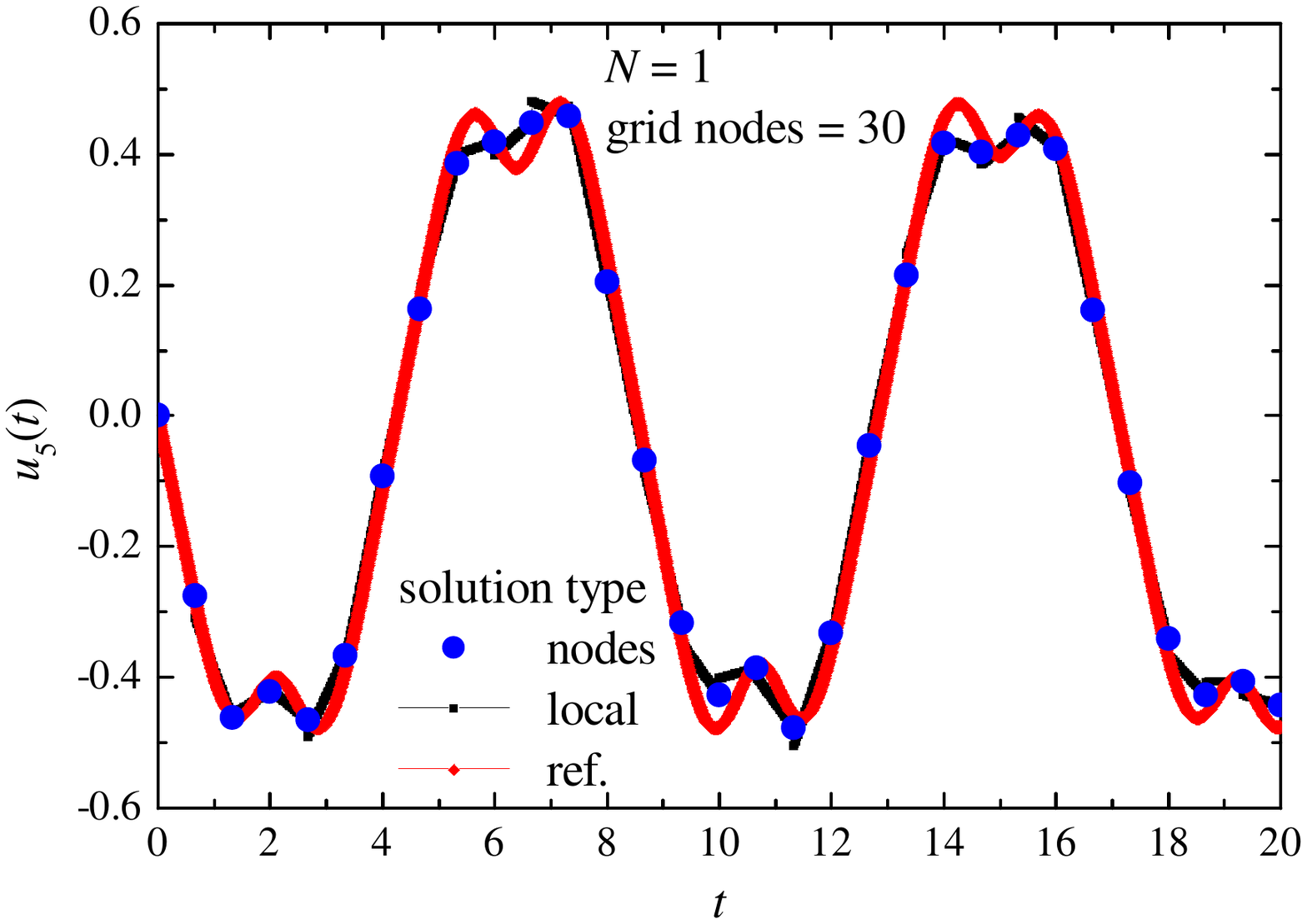}
\vspace{-8mm}\caption{\label{fig:dpend_ind3_sols_u:b1}}
\end{subfigure}
\begin{subfigure}{0.240\textwidth}
\includegraphics[width=\textwidth]{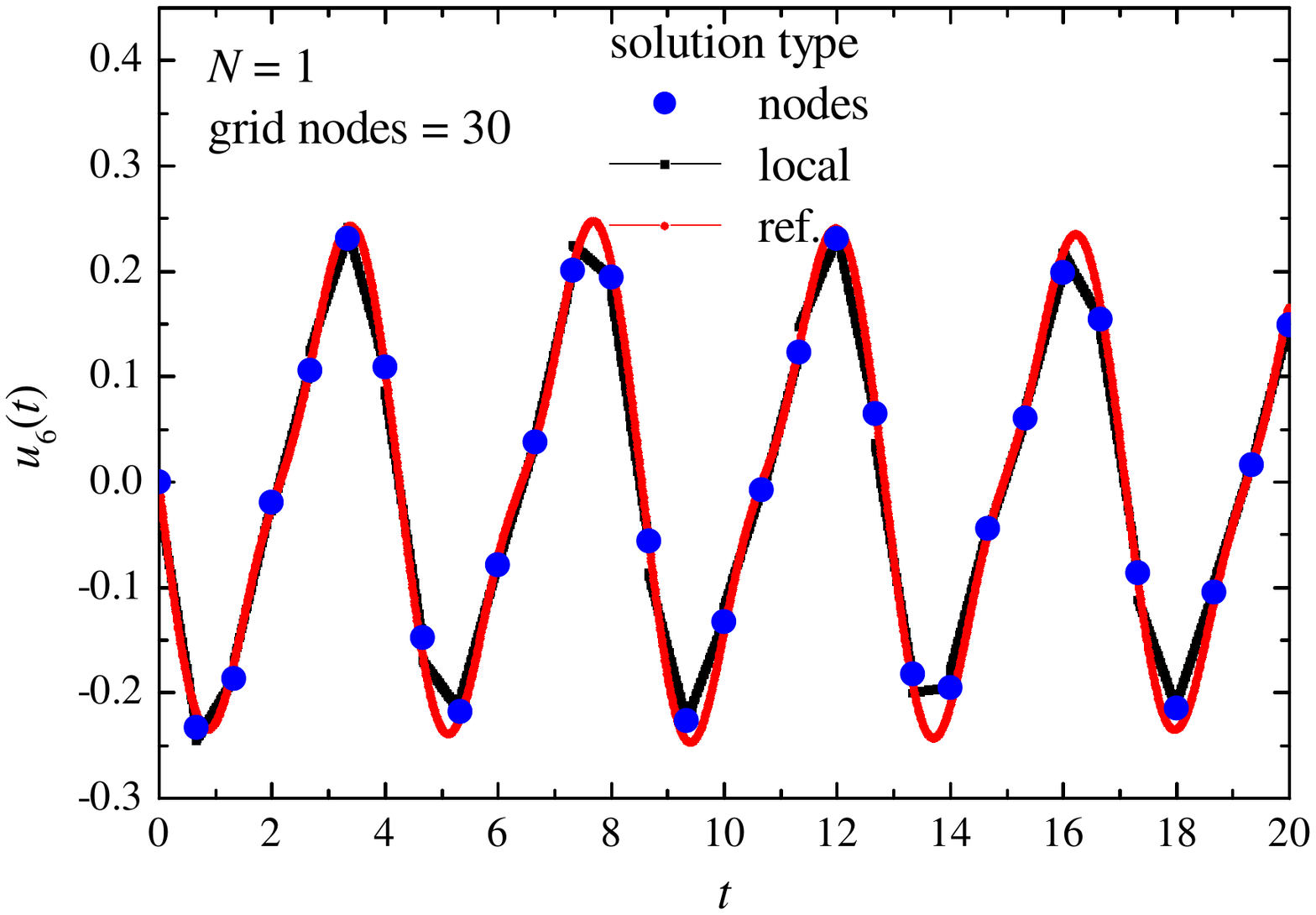}
\vspace{-8mm}\caption{\label{fig:dpend_ind3_sols_u:b2}}
\end{subfigure}
\begin{subfigure}{0.240\textwidth}
\includegraphics[width=\textwidth]{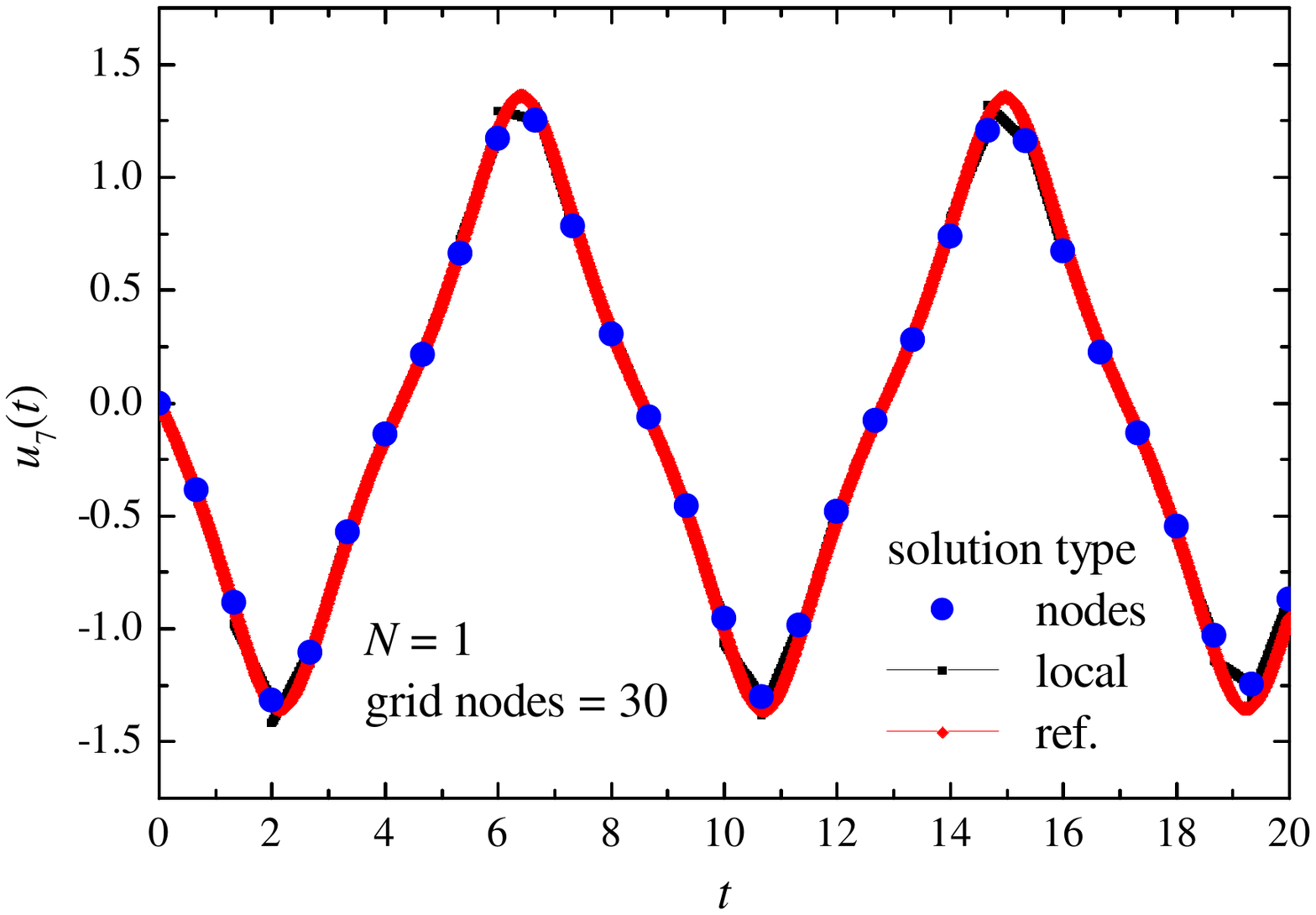}
\vspace{-8mm}\caption{\label{fig:dpend_ind3_sols_u:b3}}
\end{subfigure}
\begin{subfigure}{0.240\textwidth}
\includegraphics[width=\textwidth]{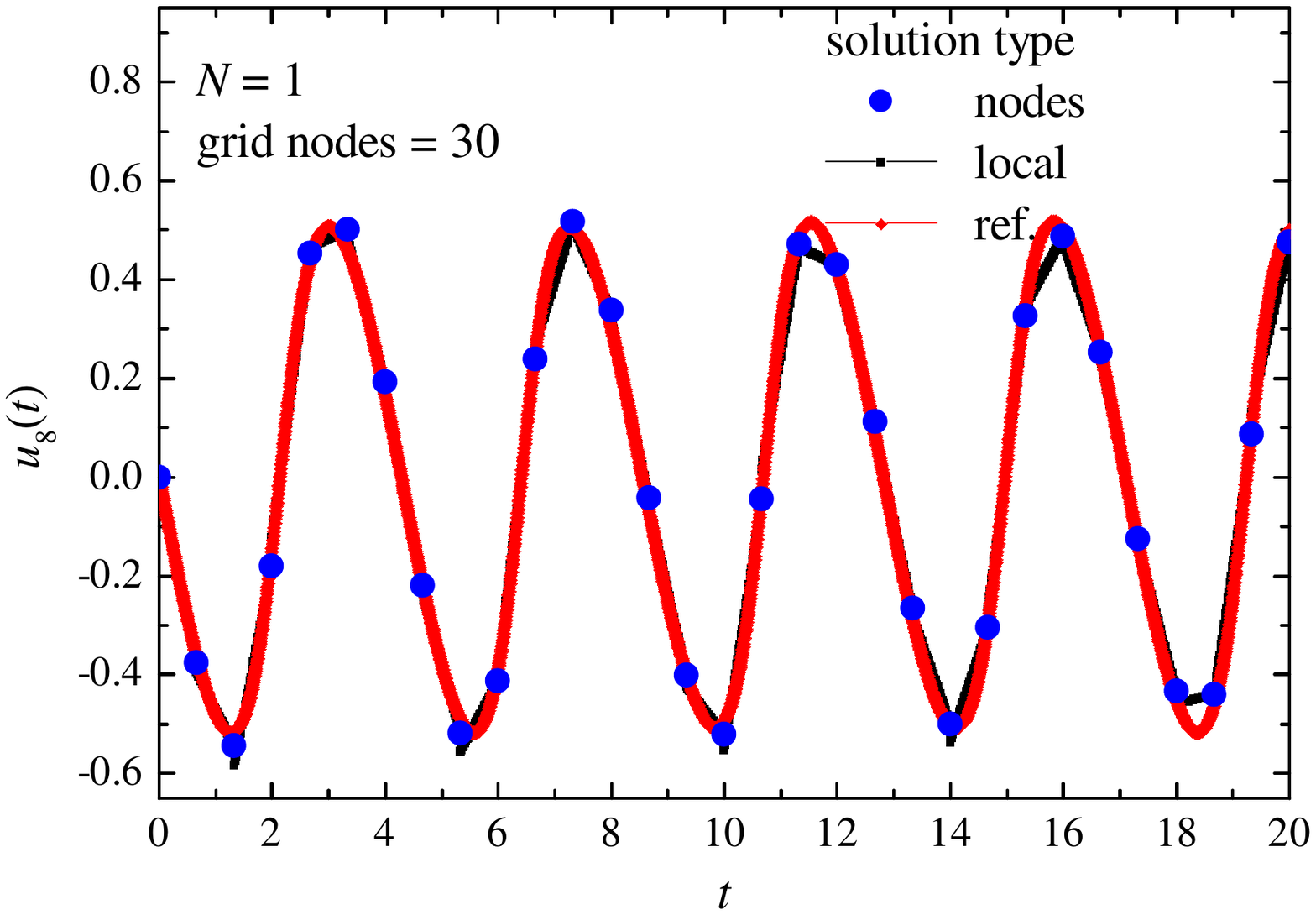}
\vspace{-8mm}\caption{\label{fig:dpend_ind3_sols_u:b4}}
\end{subfigure}\\[2mm]
\begin{subfigure}{0.240\textwidth}
\includegraphics[width=\textwidth]{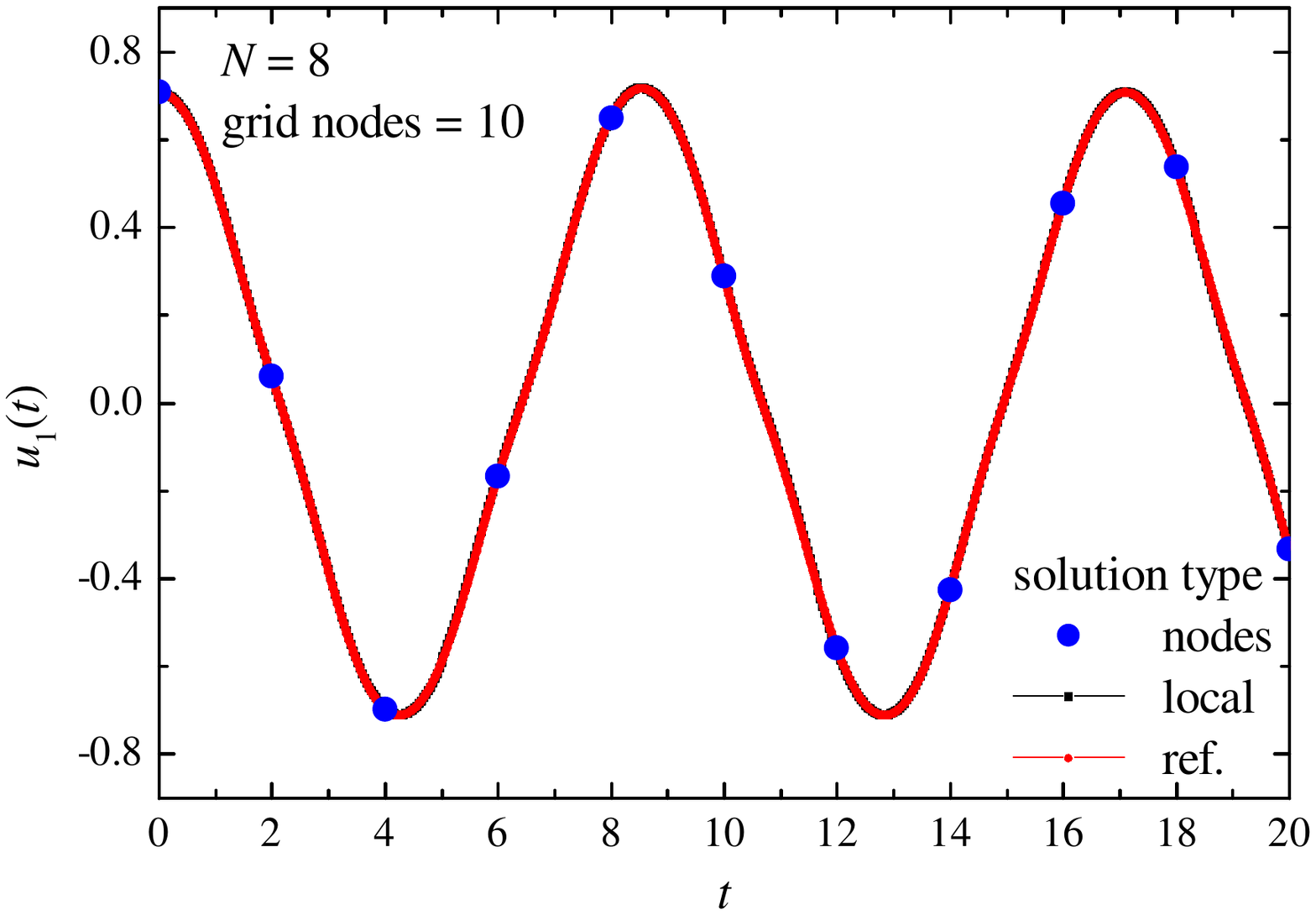}
\vspace{-8mm}\caption{\label{fig:dpend_ind3_sols_u:c1}}
\end{subfigure}
\begin{subfigure}{0.240\textwidth}
\includegraphics[width=\textwidth]{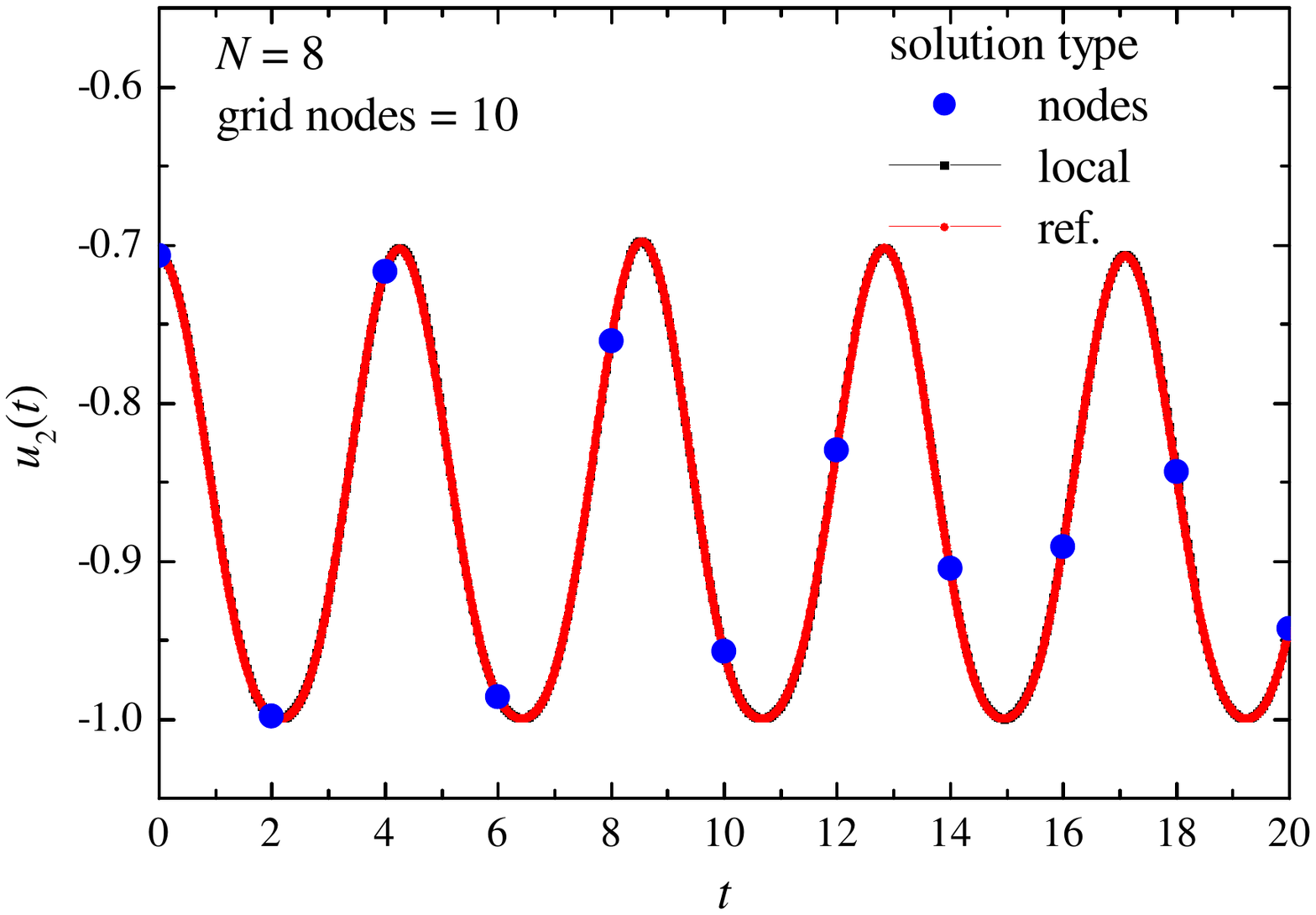}
\vspace{-8mm}\caption{\label{fig:dpend_ind3_sols_u:c2}}
\end{subfigure}
\begin{subfigure}{0.240\textwidth}
\includegraphics[width=\textwidth]{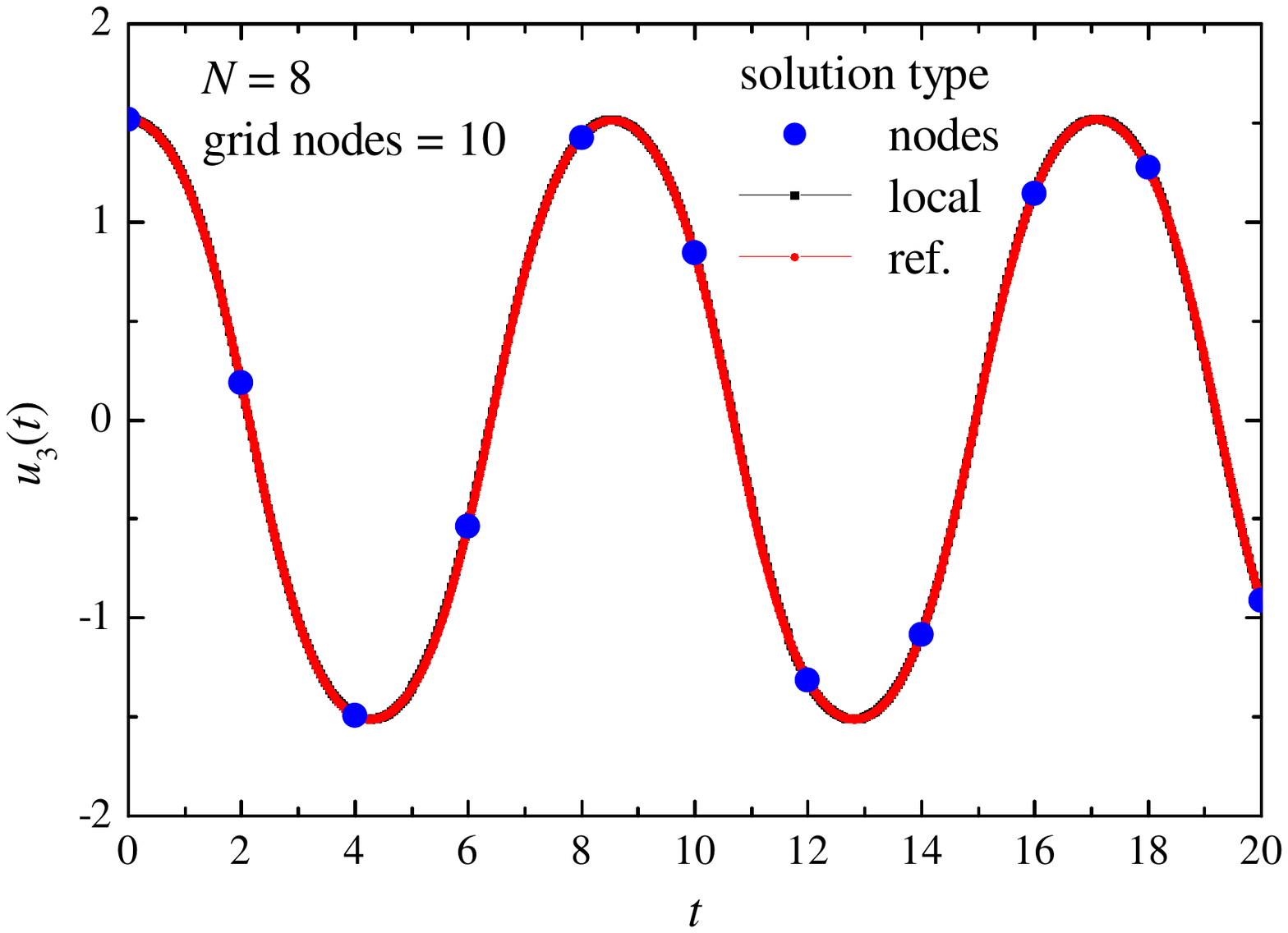}
\vspace{-8mm}\caption{\label{fig:dpend_ind3_sols_u:c3}}
\end{subfigure}
\begin{subfigure}{0.240\textwidth}
\includegraphics[width=\textwidth]{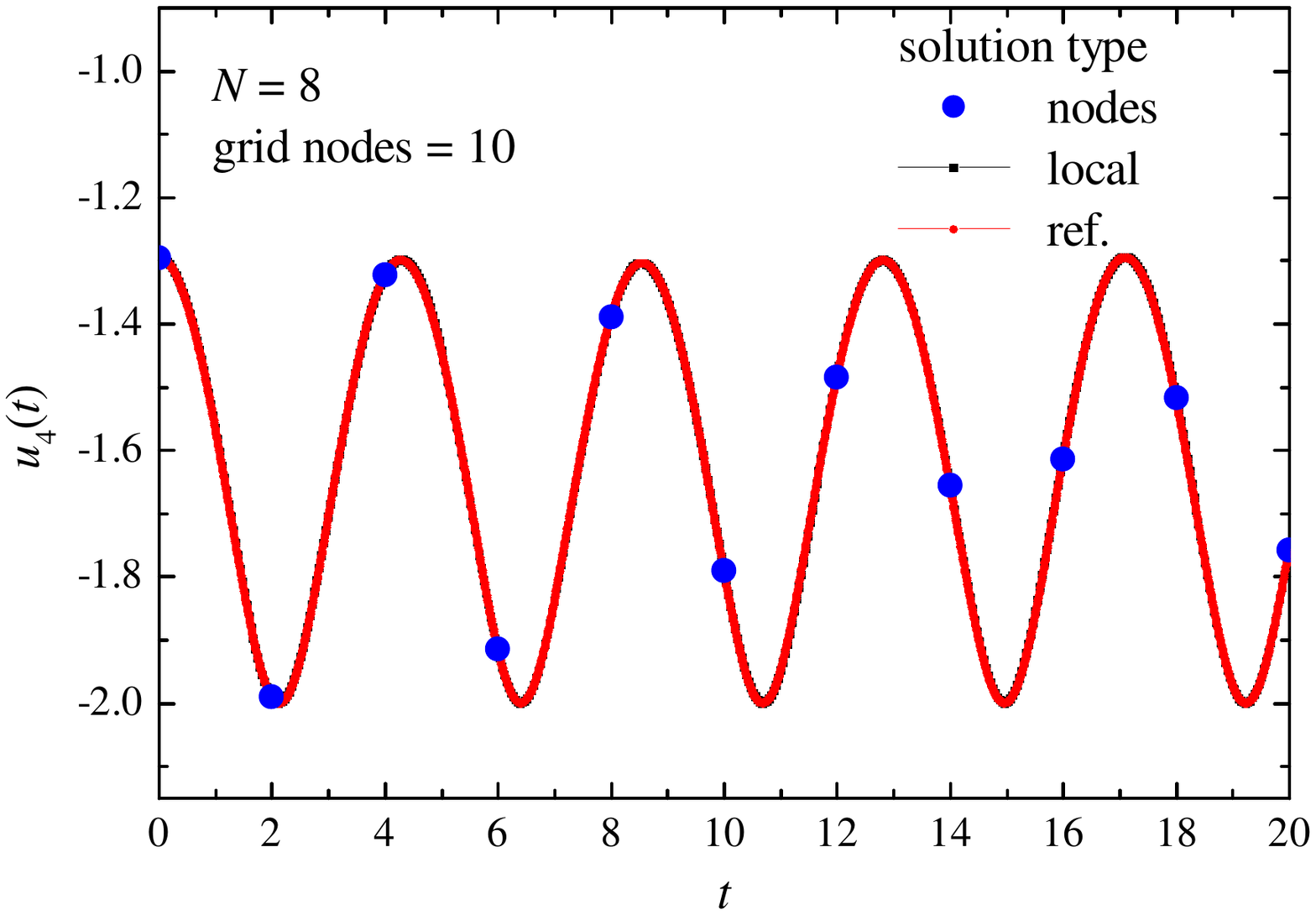}
\vspace{-8mm}\caption{\label{fig:dpend_ind3_sols_u:c4}}
\end{subfigure}\\[2mm]
\begin{subfigure}{0.240\textwidth}
\includegraphics[width=\textwidth]{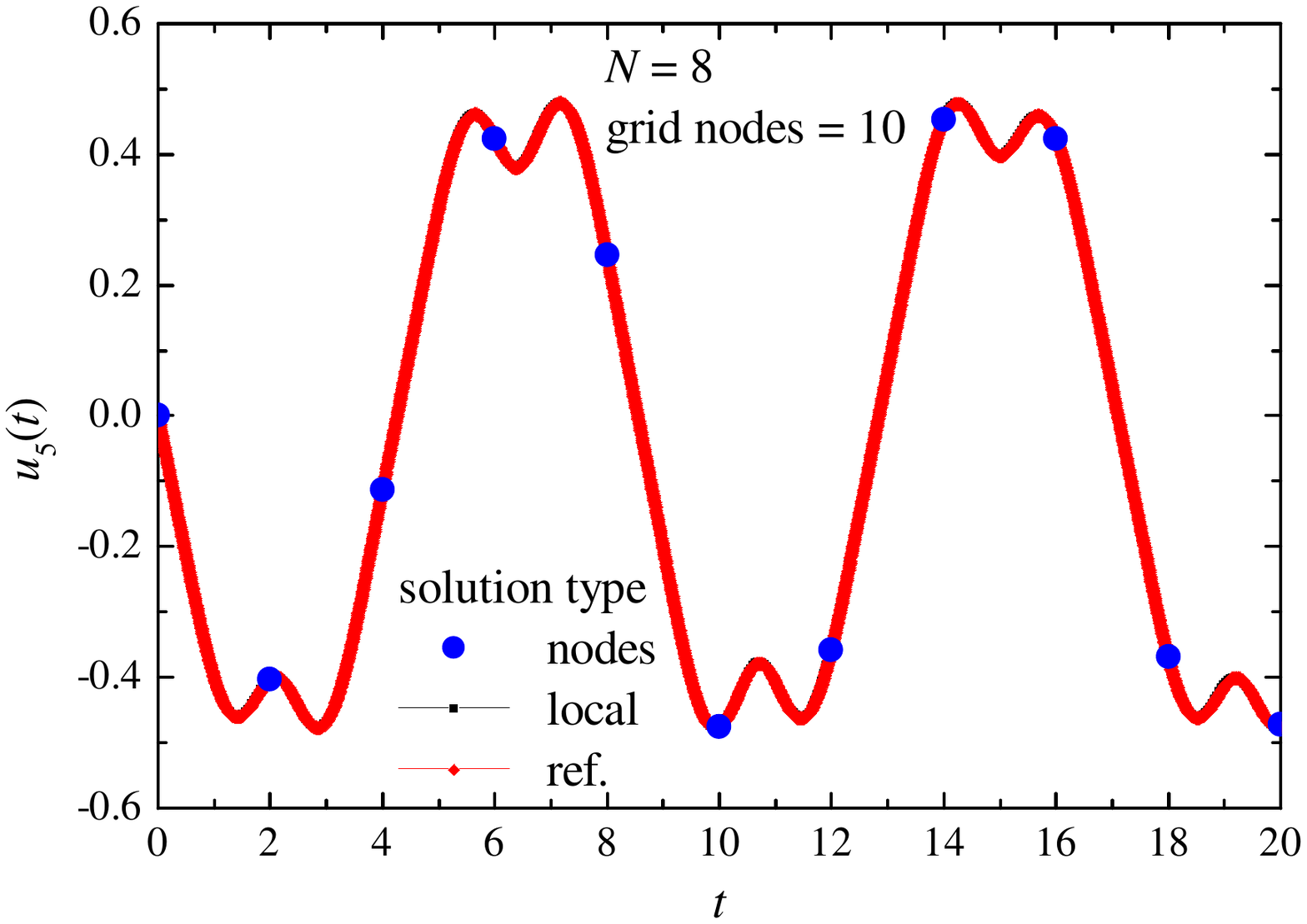}
\vspace{-8mm}\caption{\label{fig:dpend_ind3_sols_u:d1}}
\end{subfigure}
\begin{subfigure}{0.240\textwidth}
\includegraphics[width=\textwidth]{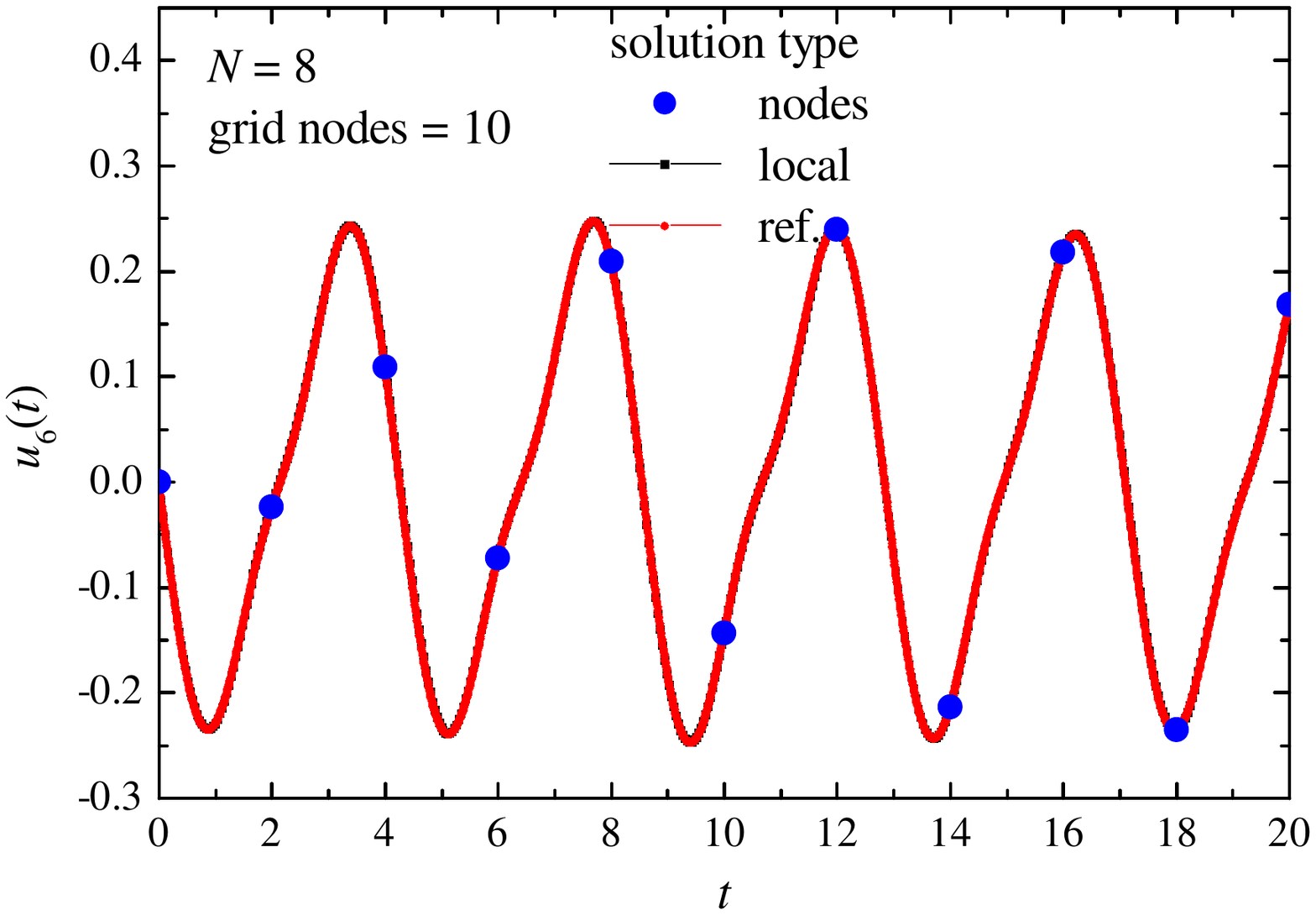}
\vspace{-8mm}\caption{\label{fig:dpend_ind3_sols_u:d2}}
\end{subfigure}
\begin{subfigure}{0.240\textwidth}
\includegraphics[width=\textwidth]{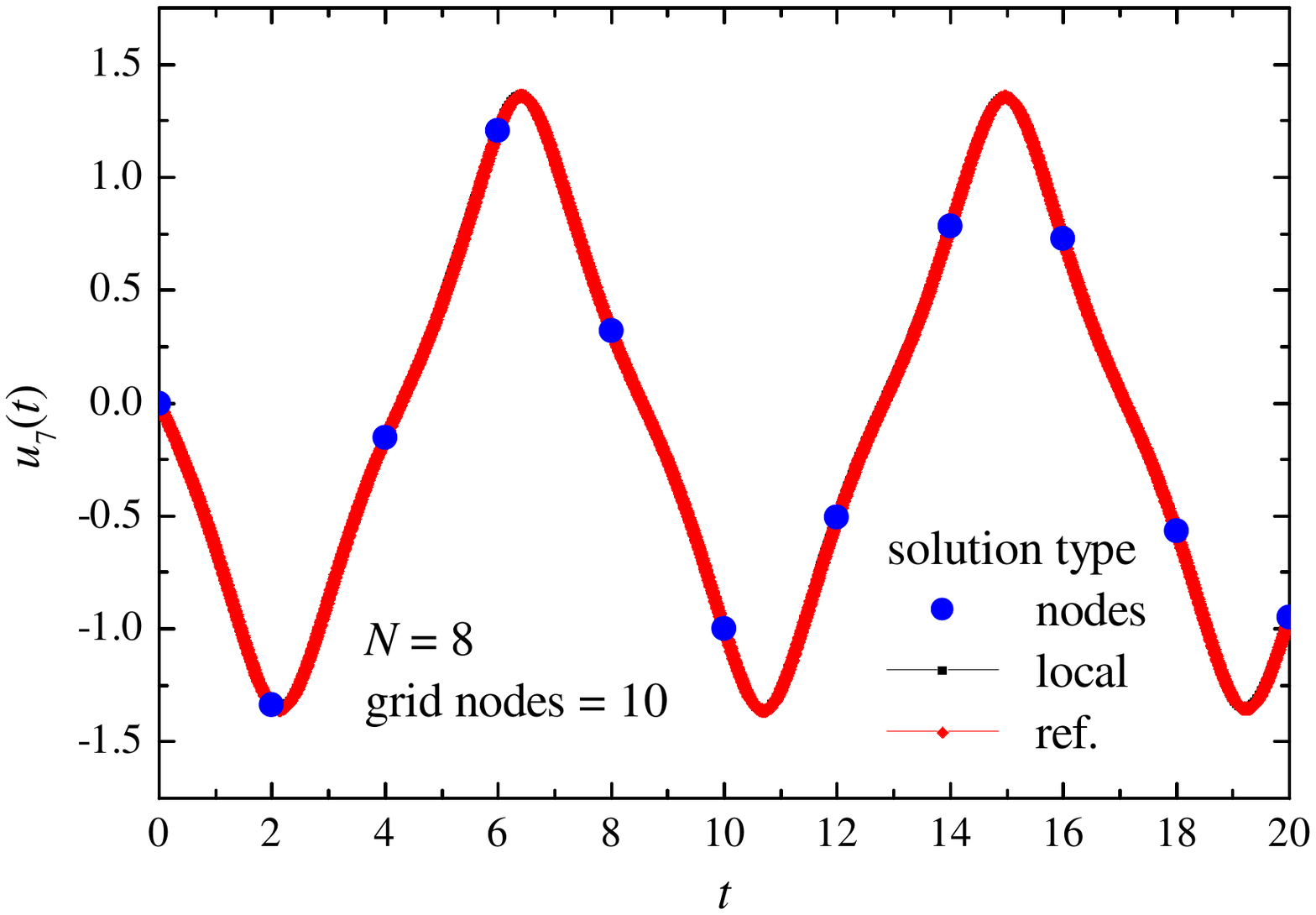}
\vspace{-8mm}\caption{\label{fig:dpend_ind3_sols_u:d3}}
\end{subfigure}
\begin{subfigure}{0.240\textwidth}
\includegraphics[width=\textwidth]{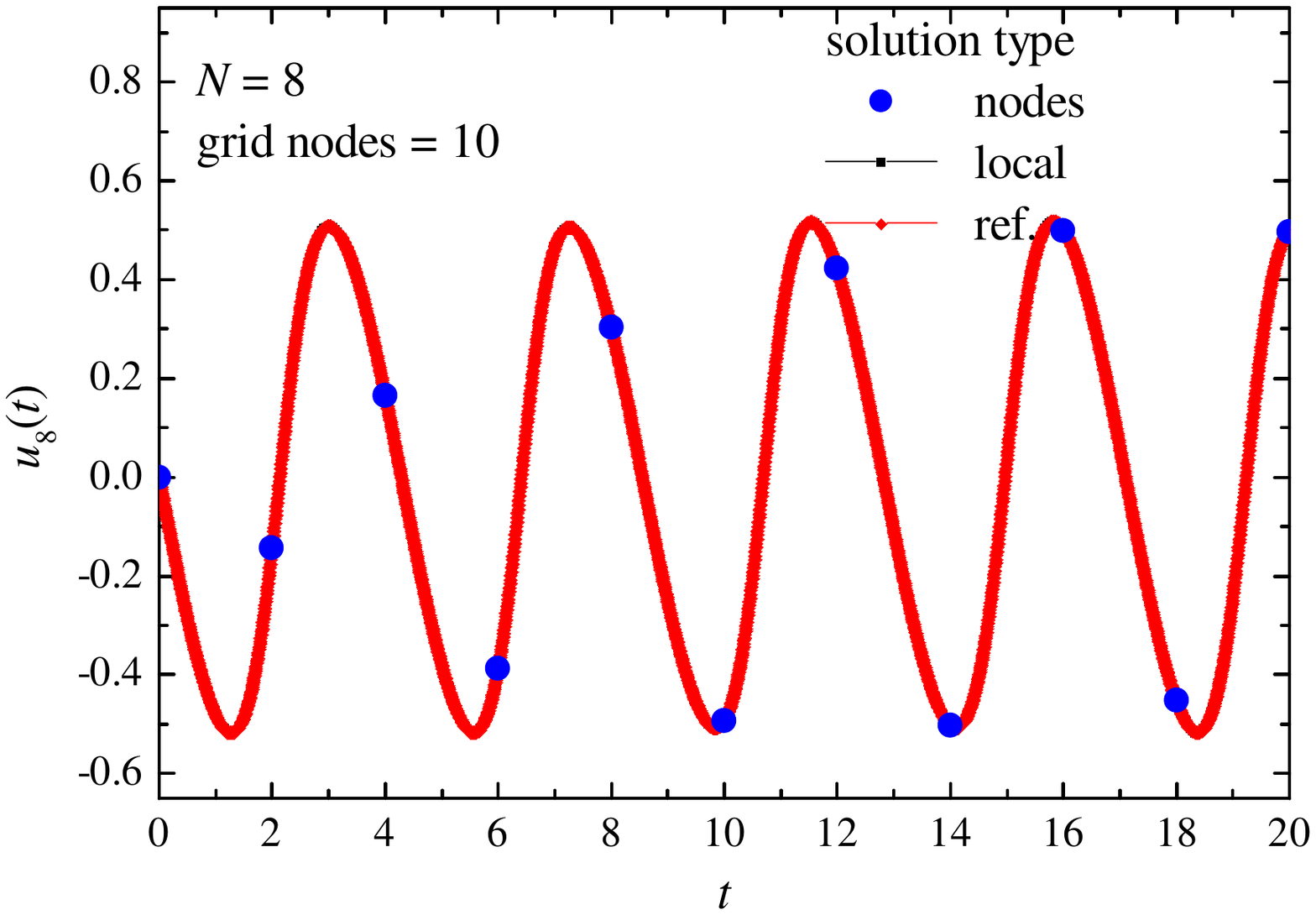}
\vspace{-8mm}\caption{\label{fig:dpend_ind3_sols_u:d4}}
\end{subfigure}\\[2mm]
\begin{subfigure}{0.240\textwidth}
\includegraphics[width=\textwidth]{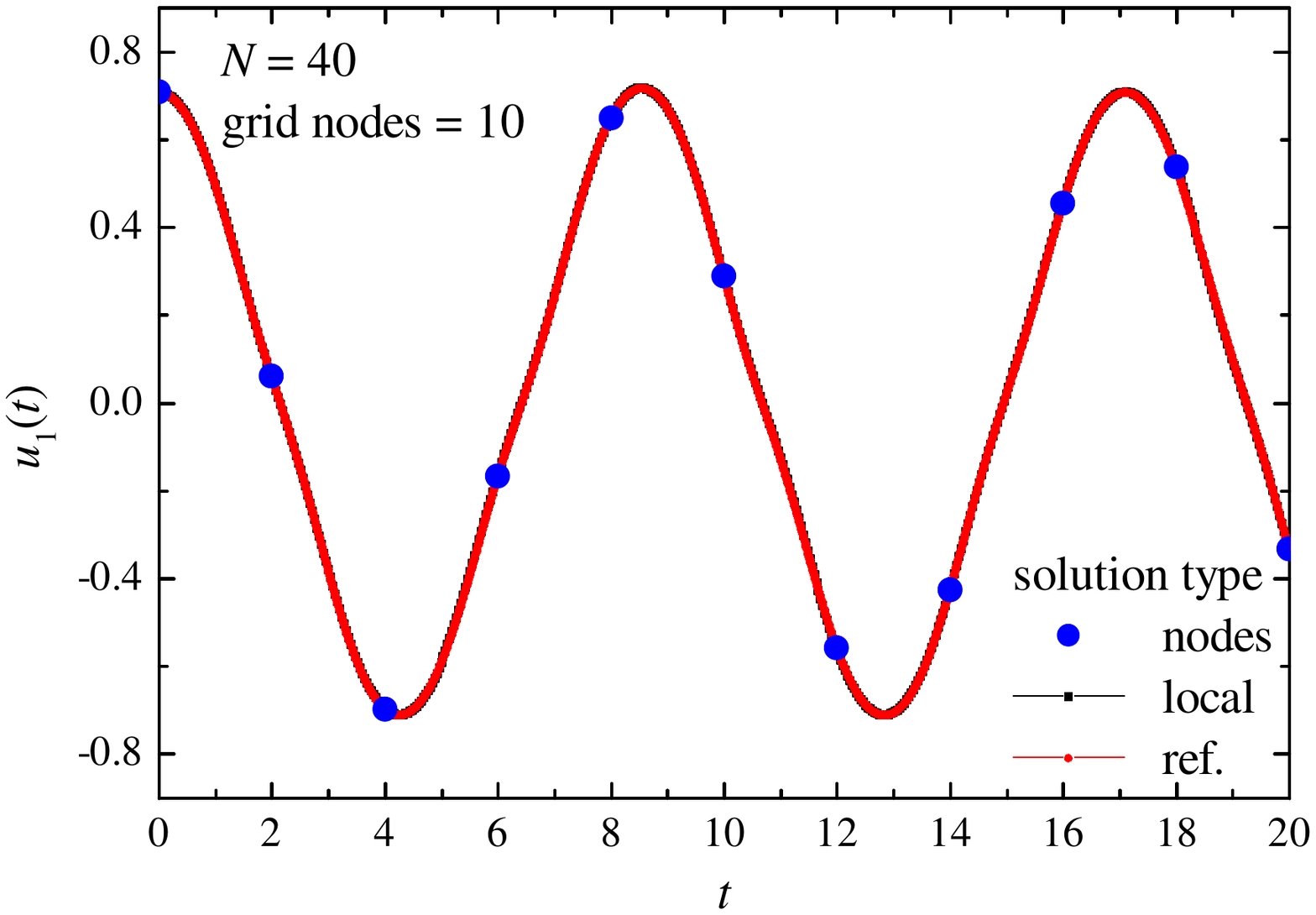}
\vspace{-8mm}\caption{\label{fig:dpend_ind3_sols_u:e1}}
\end{subfigure}
\begin{subfigure}{0.240\textwidth}
\includegraphics[width=\textwidth]{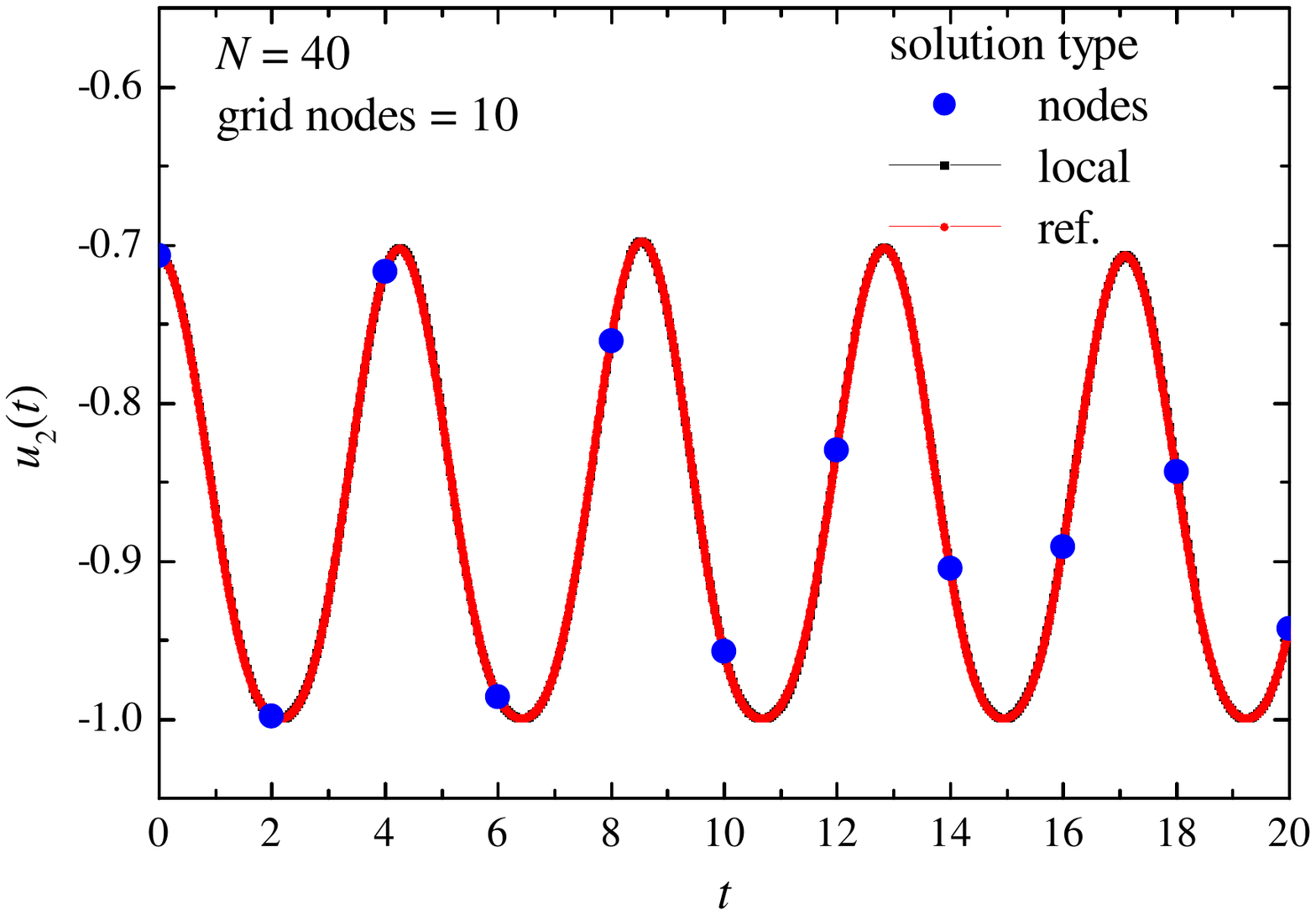}
\vspace{-8mm}\caption{\label{fig:dpend_ind3_sols_u:e2}}
\end{subfigure}
\begin{subfigure}{0.240\textwidth}
\includegraphics[width=\textwidth]{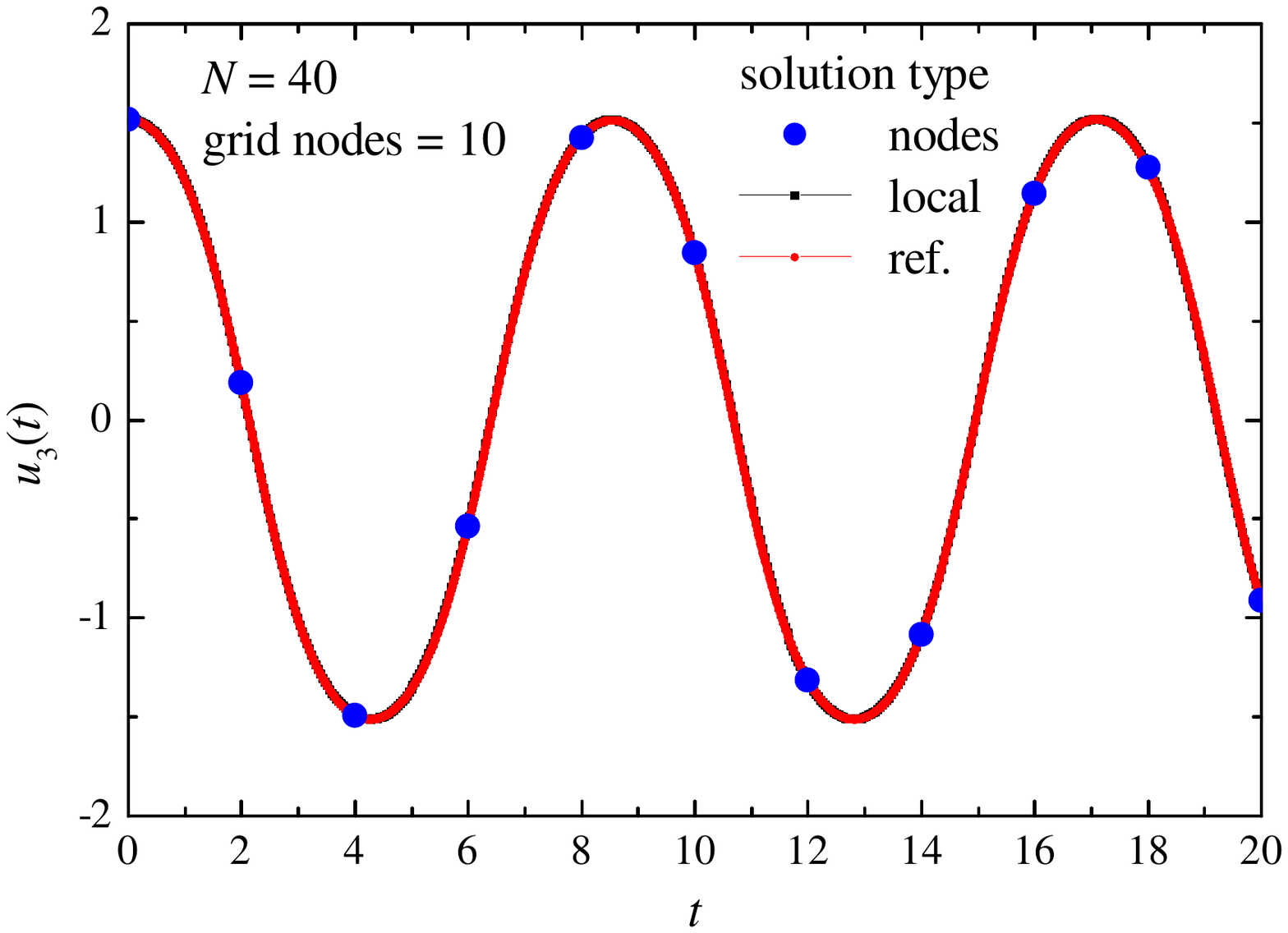}
\vspace{-8mm}\caption{\label{fig:dpend_ind3_sols_u:e3}}
\end{subfigure}
\begin{subfigure}{0.240\textwidth}
\includegraphics[width=\textwidth]{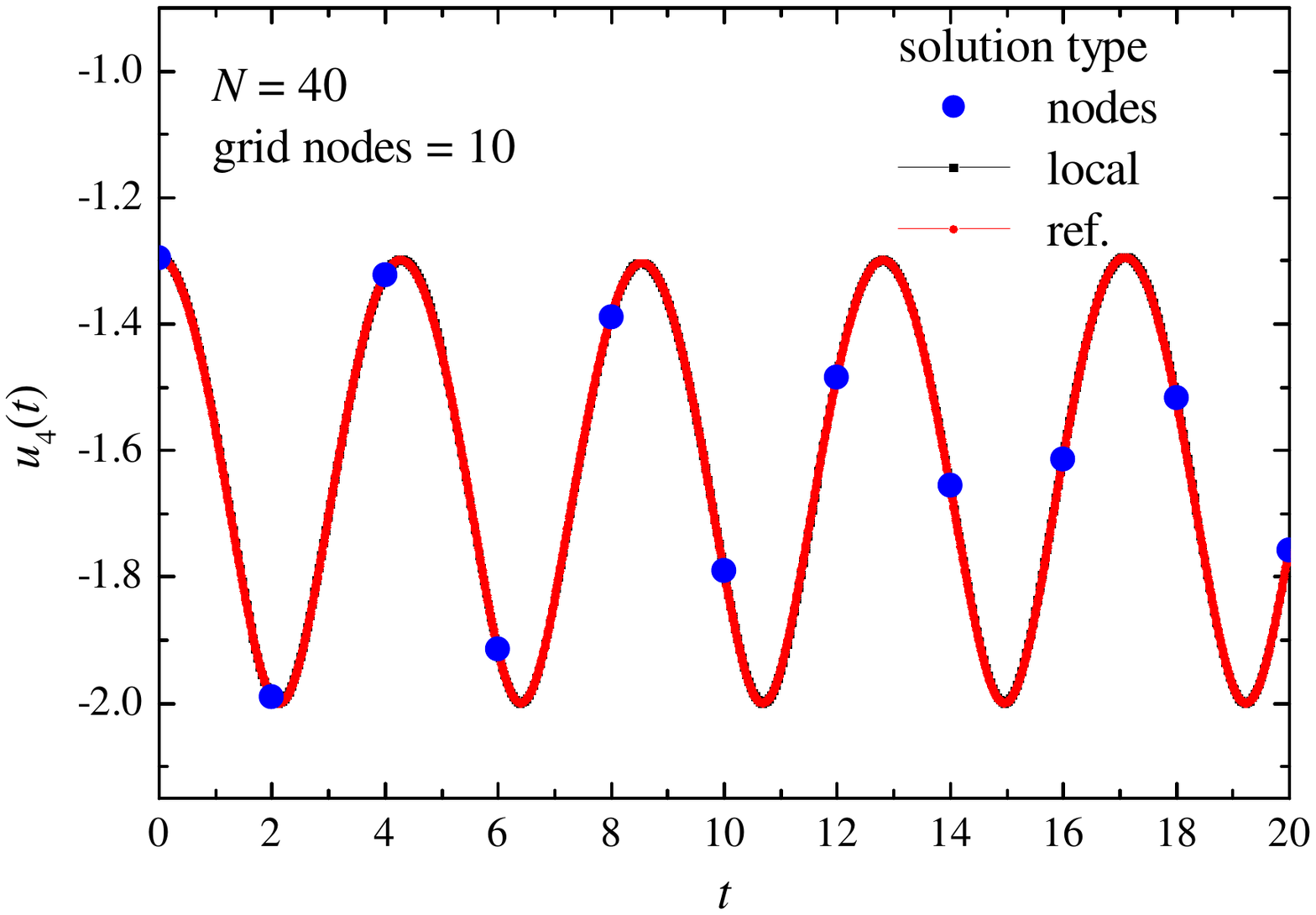}
\vspace{-8mm}\caption{\label{fig:dpend_ind3_sols_u:e4}}
\end{subfigure}\\[2mm]
\begin{subfigure}{0.240\textwidth}
\includegraphics[width=\textwidth]{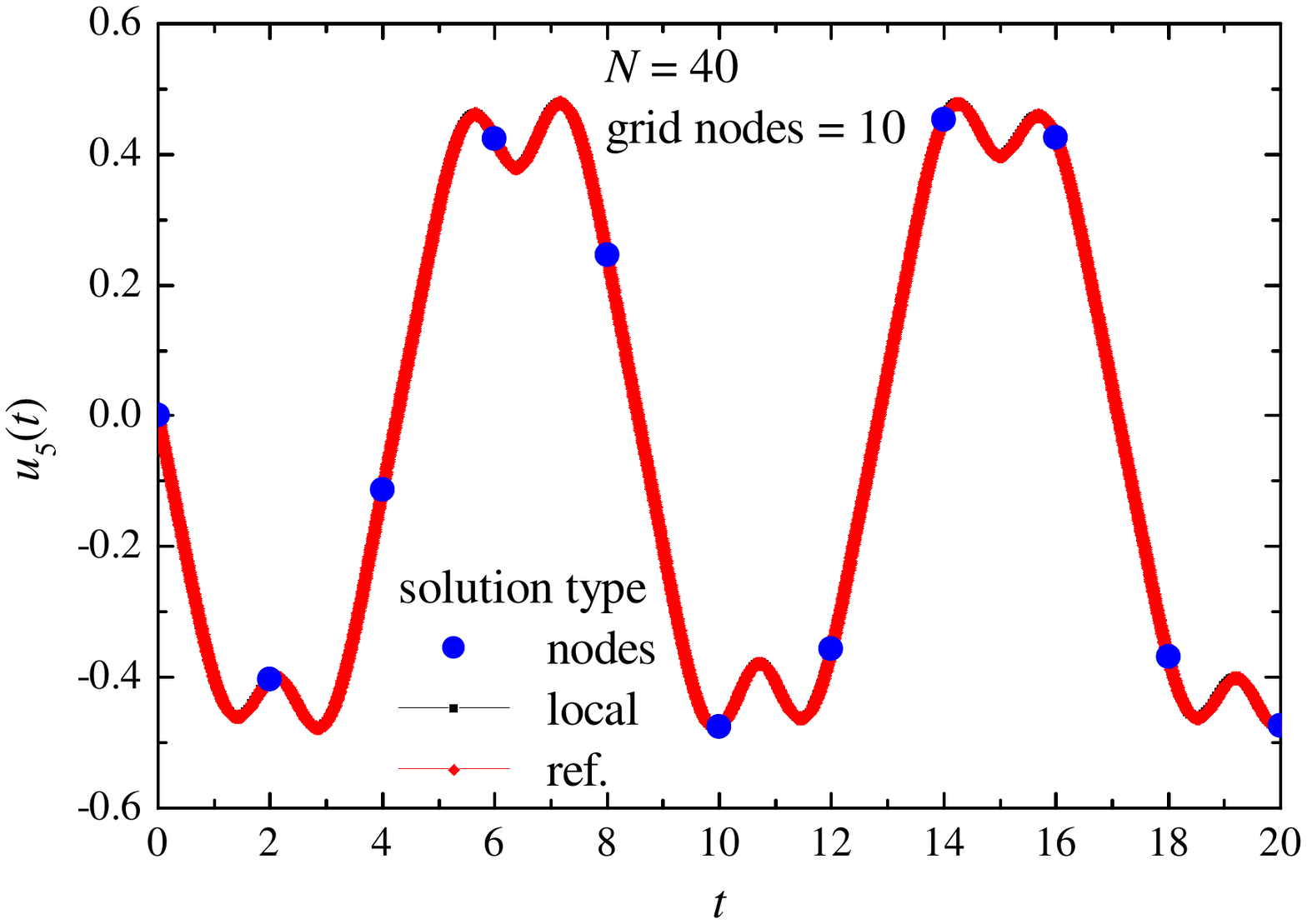}
\vspace{-8mm}\caption{\label{fig:dpend_ind3_sols_u:f1}}
\end{subfigure}
\begin{subfigure}{0.240\textwidth}
\includegraphics[width=\textwidth]{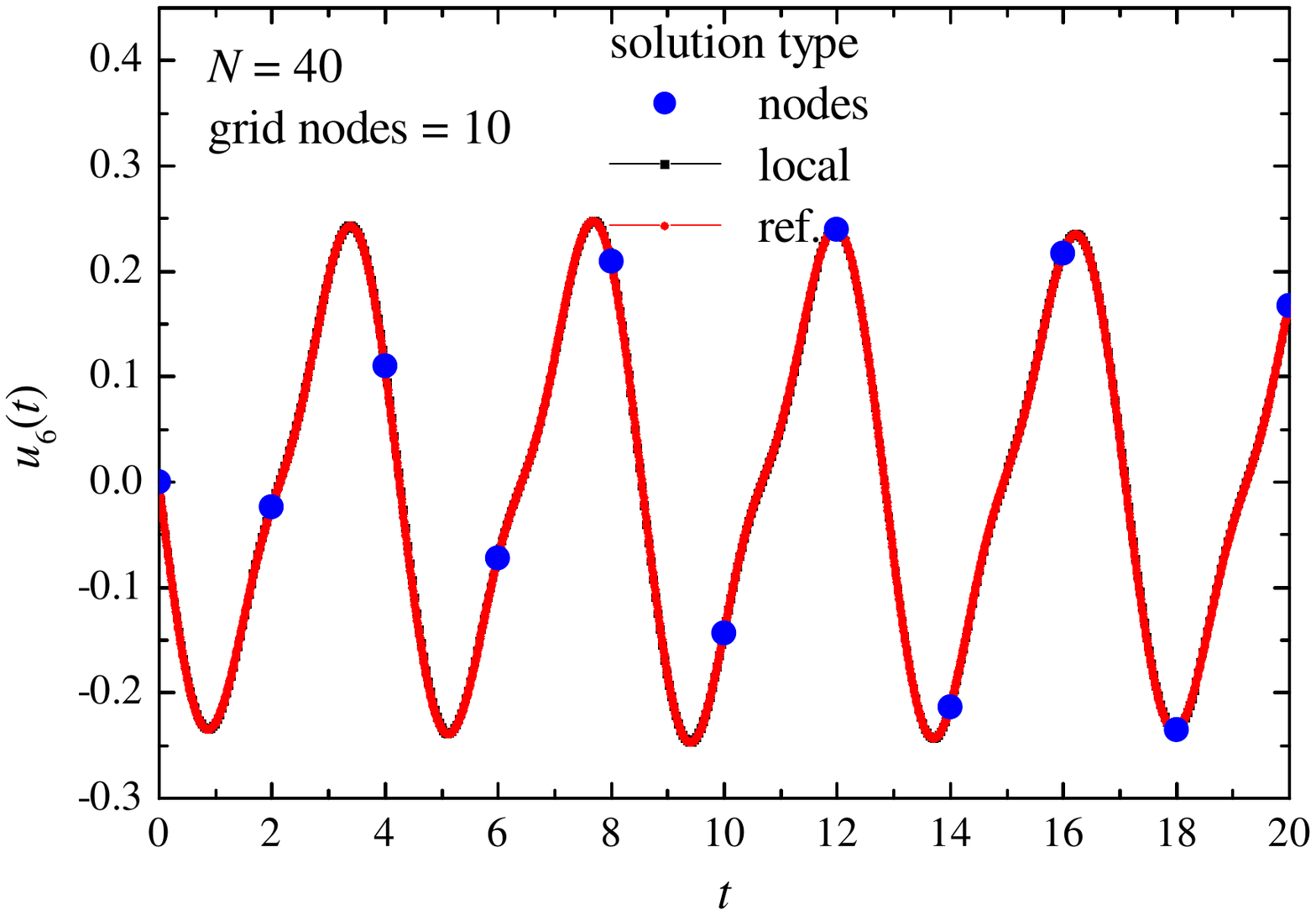}
\vspace{-8mm}\caption{\label{fig:dpend_ind3_sols_u:f2}}
\end{subfigure}
\begin{subfigure}{0.240\textwidth}
\includegraphics[width=\textwidth]{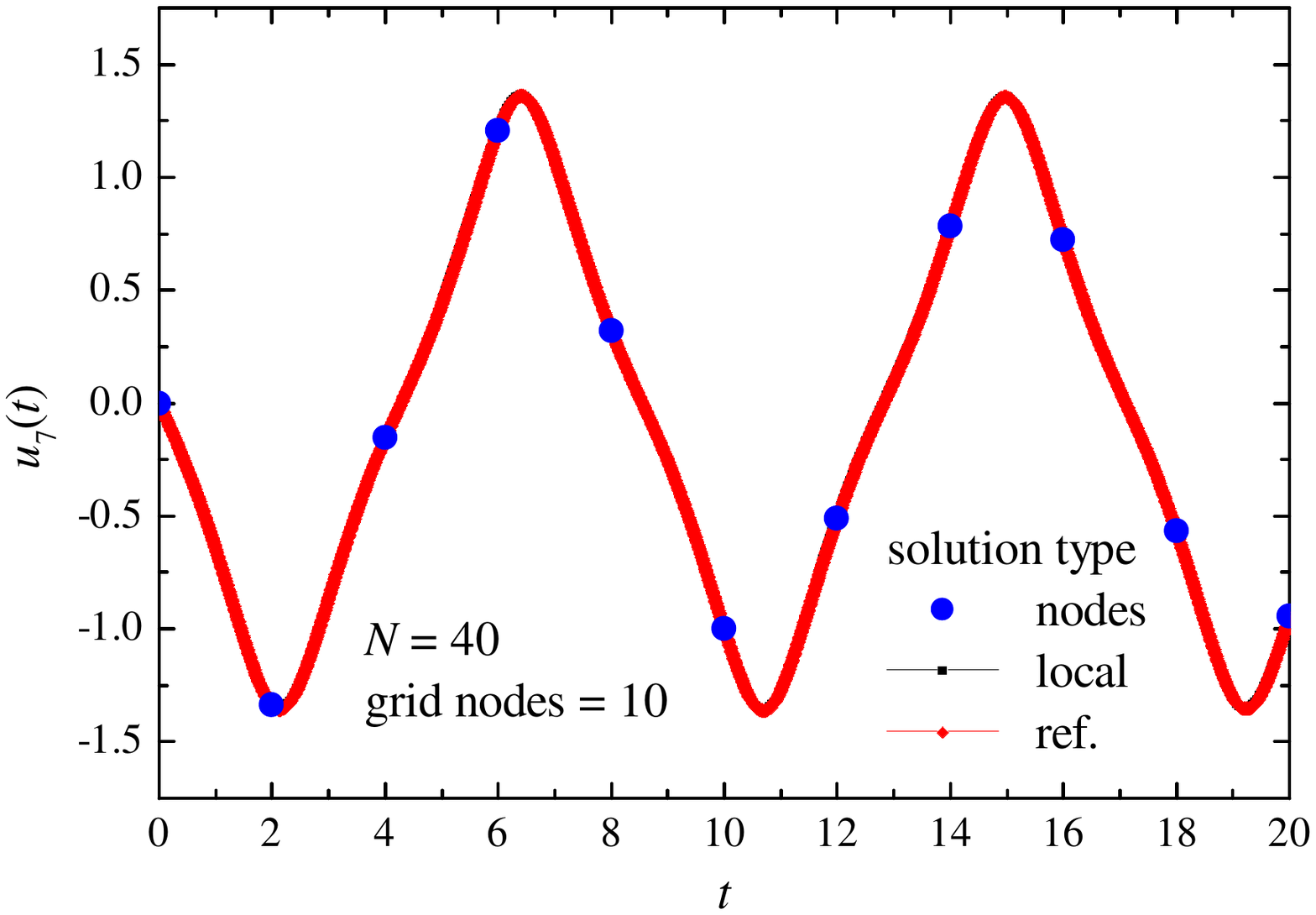}
\vspace{-8mm}\caption{\label{fig:dpend_ind3_sols_u:f3}}
\end{subfigure}
\begin{subfigure}{0.240\textwidth}
\includegraphics[width=\textwidth]{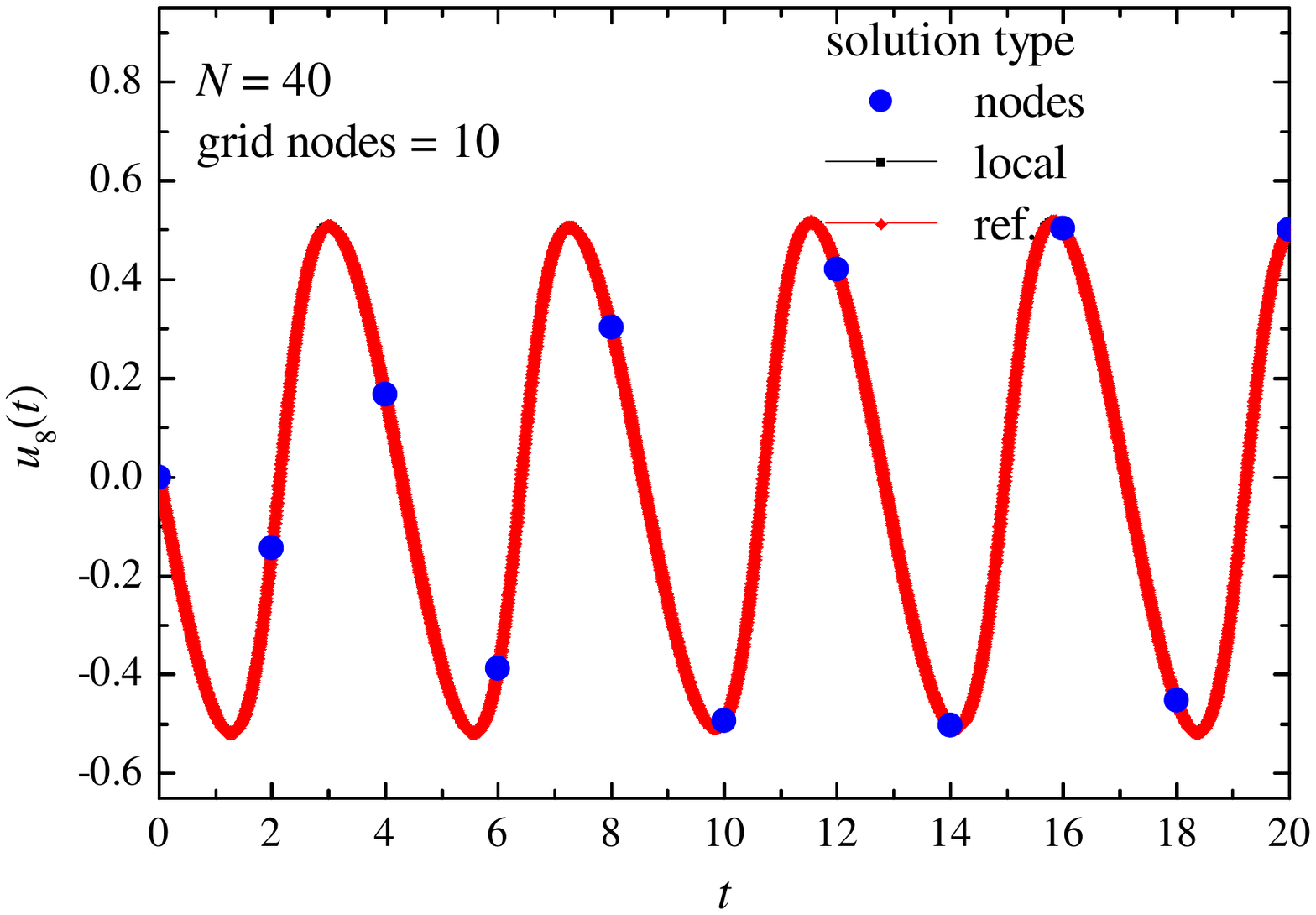}
\vspace{-8mm}\caption{\label{fig:dpend_ind3_sols_u:f4}}
\end{subfigure}\\[2mm]
\caption{%
Numerical solution of the DAE system (\ref{eq:math_dpend_dae_ind_3}) of index 3. Comparison of the solution at nodes $\mathbf{u}_{n}$, the local solution $\mathbf{u}_{L}(t)$ and the reference solution $\mathbf{u}^{\rm ref}(t)$ for components $u_{1}$ (\subref{fig:dpend_ind3_sols_u:a1}, \subref{fig:dpend_ind3_sols_u:c1}, \subref{fig:dpend_ind3_sols_u:e1}), $u_{2}$ (\subref{fig:dpend_ind3_sols_u:a2}, \subref{fig:dpend_ind3_sols_u:c2}, \subref{fig:dpend_ind3_sols_u:e2}), $u_{3}$ (\subref{fig:dpend_ind3_sols_u:a3}, \subref{fig:dpend_ind3_sols_u:c3}, \subref{fig:dpend_ind3_sols_u:e3}), $u_{4}$ (\subref{fig:dpend_ind3_sols_u:a4}, \subref{fig:dpend_ind3_sols_u:c4}, \subref{fig:dpend_ind3_sols_u:e4}), $u_{5}$ (\subref{fig:dpend_ind3_sols_u:b1}, \subref{fig:dpend_ind3_sols_u:d1}, \subref{fig:dpend_ind3_sols_u:f1}), $u_{6}$ (\subref{fig:dpend_ind3_sols_u:b2}, \subref{fig:dpend_ind3_sols_u:d2}, \subref{fig:dpend_ind3_sols_u:f2}), $u_{7}$ (\subref{fig:dpend_ind3_sols_u:b3}, \subref{fig:dpend_ind3_sols_u:d3}, \subref{fig:dpend_ind3_sols_u:f3}), $u_{8}$ (\subref{fig:dpend_ind3_sols_u:b4}, \subref{fig:dpend_ind3_sols_u:d4}, \subref{fig:dpend_ind3_sols_u:f4}), obtained using polynomials with degrees $N = 1$ (\subref{fig:dpend_ind3_sols_u:a1}, \subref{fig:dpend_ind3_sols_u:a2}, \subref{fig:dpend_ind3_sols_u:a3}, \subref{fig:dpend_ind3_sols_u:a4}, \subref{fig:dpend_ind3_sols_u:b1}, \subref{fig:dpend_ind3_sols_u:b2}, \subref{fig:dpend_ind3_sols_u:b3}, \subref{fig:dpend_ind3_sols_u:b4}), $N = 8$ (\subref{fig:dpend_ind3_sols_u:c1}, \subref{fig:dpend_ind3_sols_u:c2}, \subref{fig:dpend_ind3_sols_u:c3}, \subref{fig:dpend_ind3_sols_u:c4}, \subref{fig:dpend_ind3_sols_u:d1}, \subref{fig:dpend_ind3_sols_u:d2}, \subref{fig:dpend_ind3_sols_u:d3}, \subref{fig:dpend_ind3_sols_u:d4}) and $N = 40$ (\subref{fig:dpend_ind3_sols_u:e1}, \subref{fig:dpend_ind3_sols_u:e2}, \subref{fig:dpend_ind3_sols_u:e3}, \subref{fig:dpend_ind3_sols_u:e4}, \subref{fig:dpend_ind3_sols_u:f1}, \subref{fig:dpend_ind3_sols_u:f2}, \subref{fig:dpend_ind3_sols_u:f3}, \subref{fig:dpend_ind3_sols_u:f4}).
}
\label{fig:dpend_ind3_sols_u}
\end{figure} 

\begin{figure}[h!]
\captionsetup[subfigure]{%
	position=bottom,
	font+=smaller,
	textfont=normalfont,
	singlelinecheck=off,
	justification=raggedright
}
\centering
\begin{subfigure}{0.240\textwidth}
\includegraphics[width=\textwidth]{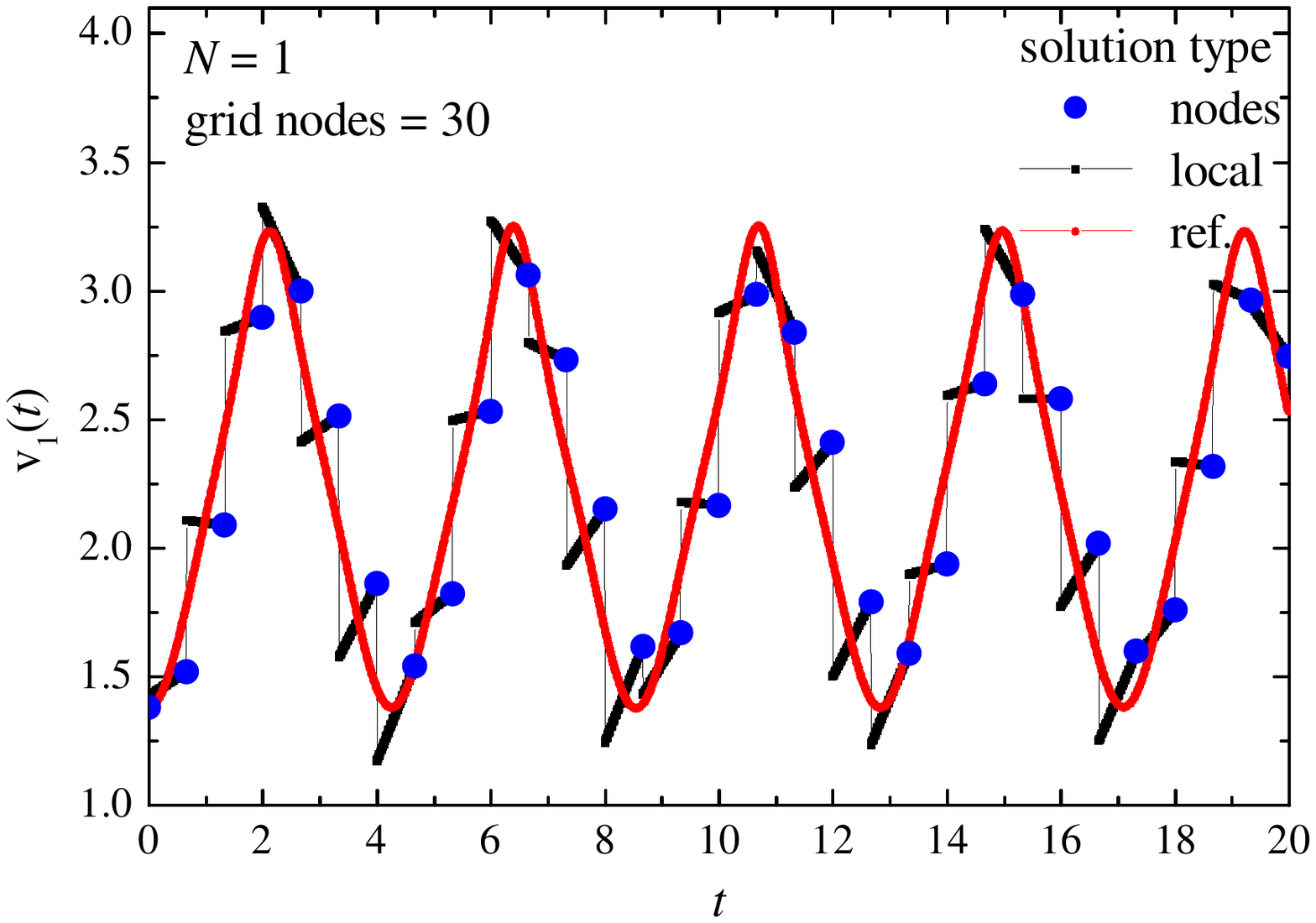}
\vspace{-8mm}\caption{\label{fig:dpend_ind3_sols_vg:a1}}
\end{subfigure}
\begin{subfigure}{0.240\textwidth}
\includegraphics[width=\textwidth]{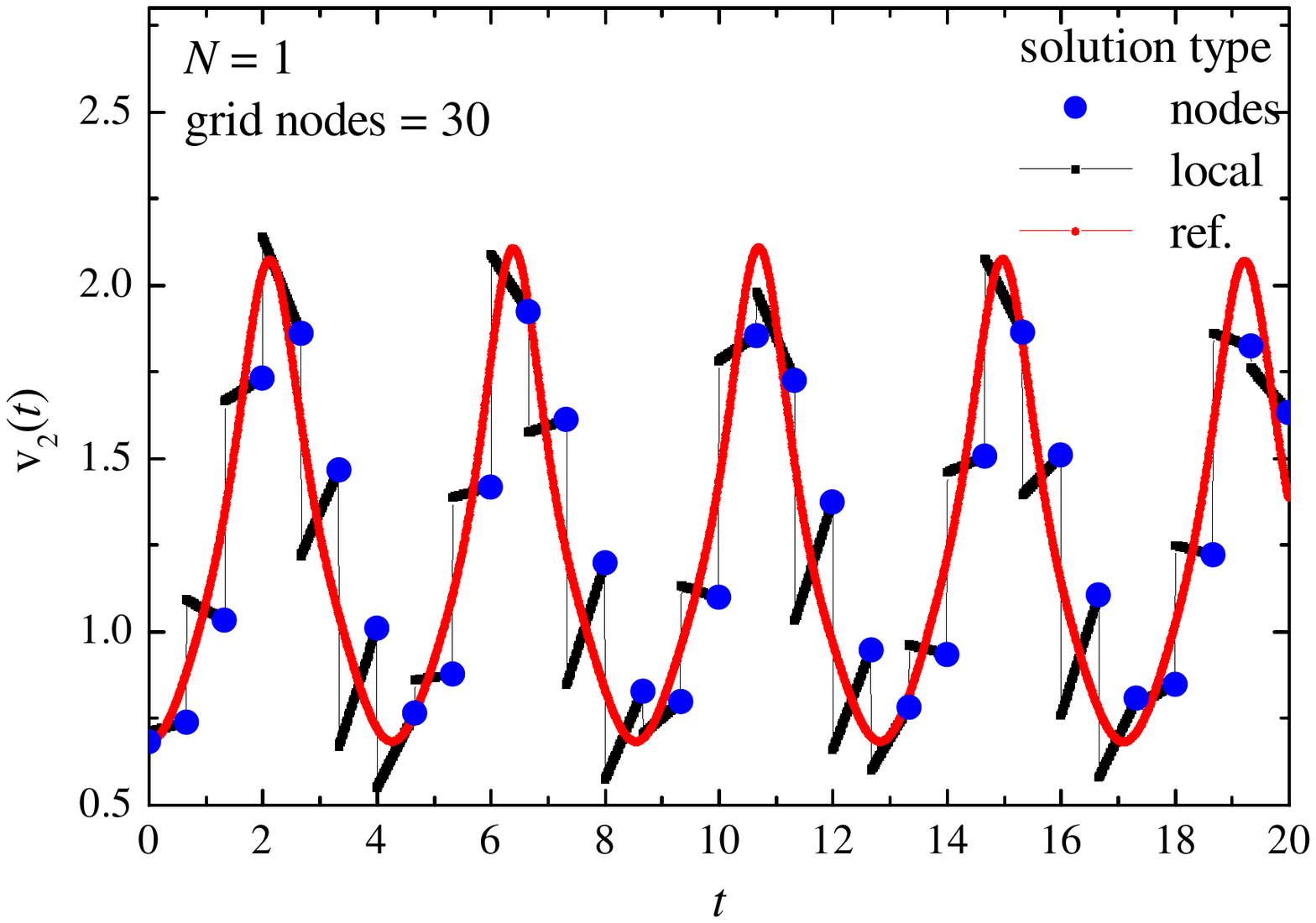}
\vspace{-8mm}\caption{\label{fig:dpend_ind3_sols_vg:a2}}
\end{subfigure}
\begin{subfigure}{0.240\textwidth}
\includegraphics[width=\textwidth]{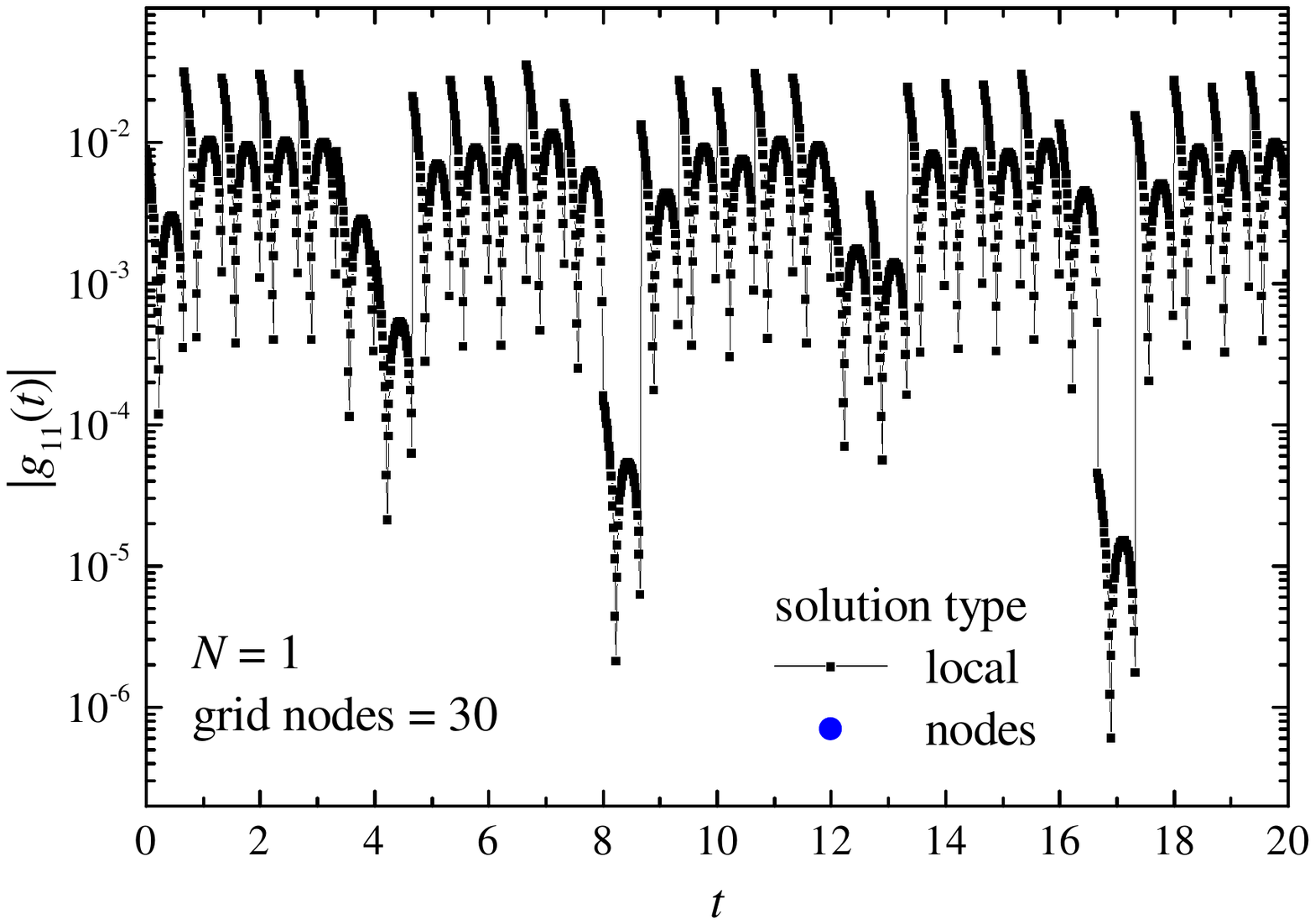}
\vspace{-8mm}\caption{\label{fig:dpend_ind3_sols_vg:a3}}
\end{subfigure}
\begin{subfigure}{0.240\textwidth}
\includegraphics[width=\textwidth]{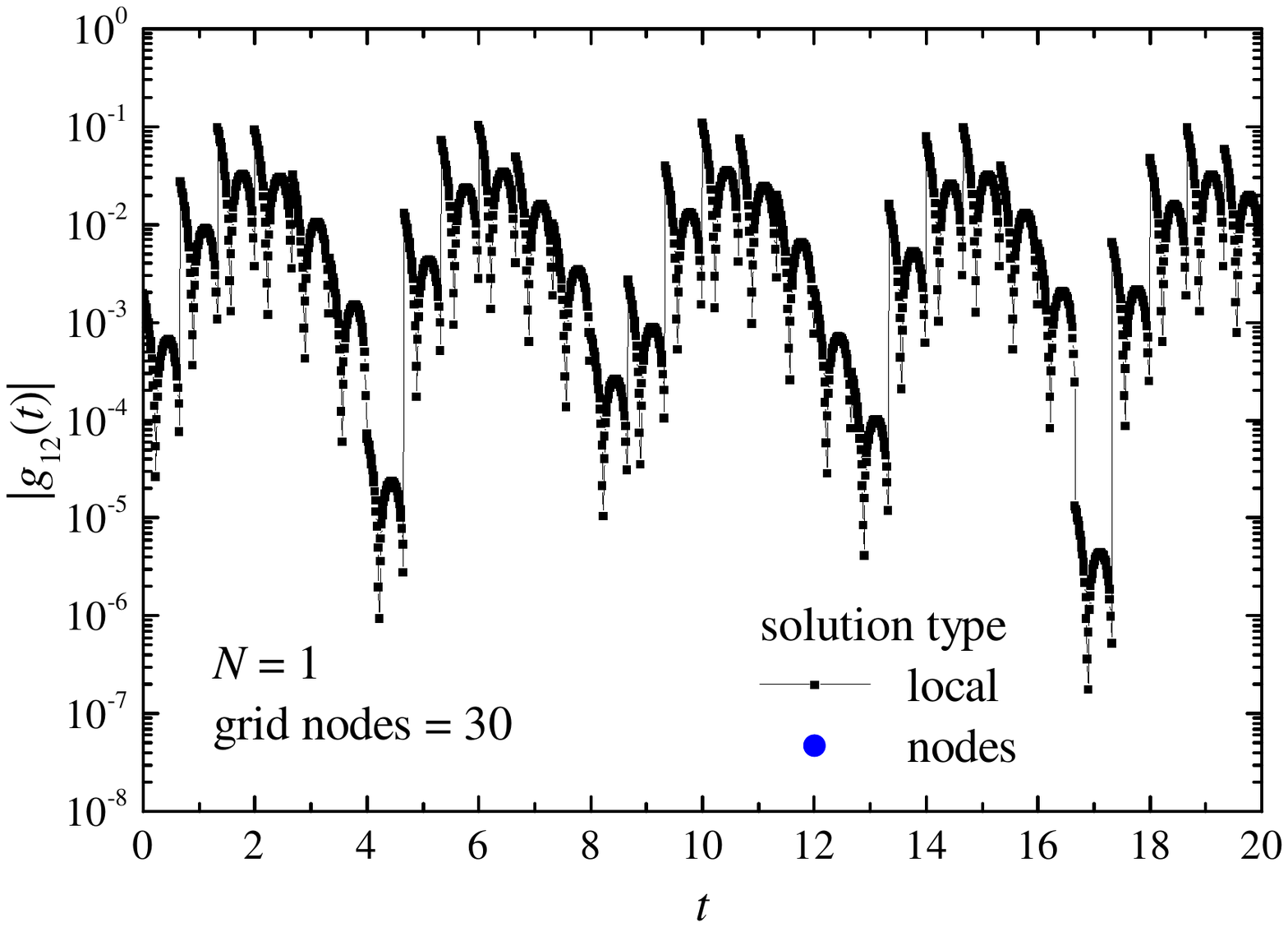}
\vspace{-8mm}\caption{\label{fig:dpend_ind3_sols_vg:a4}}
\end{subfigure}\\[2mm]
\begin{subfigure}{0.240\textwidth}
\includegraphics[width=\textwidth]{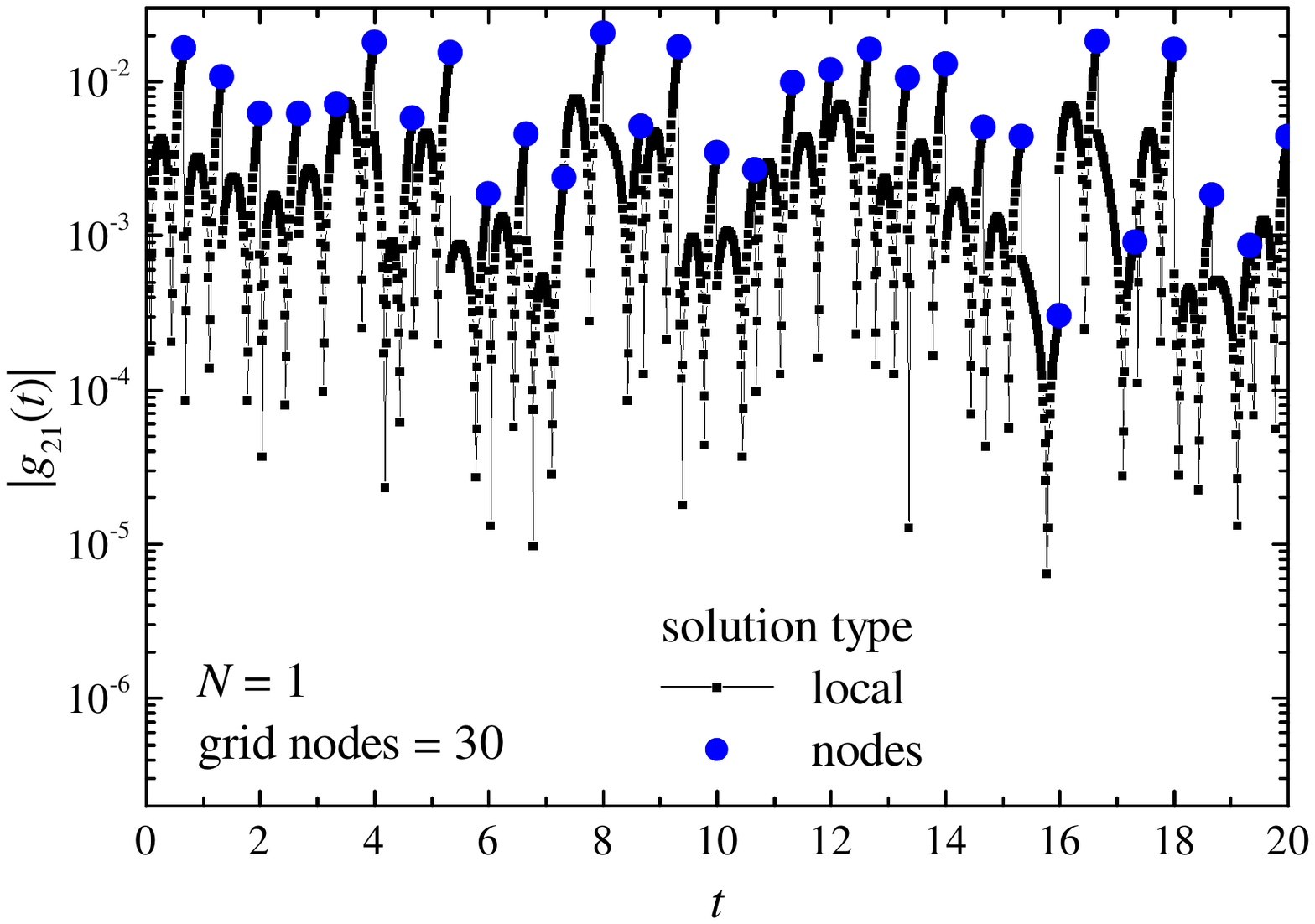}
\vspace{-8mm}\caption{\label{fig:dpend_ind3_sols_vg:b1}}
\end{subfigure}
\begin{subfigure}{0.240\textwidth}
\includegraphics[width=\textwidth]{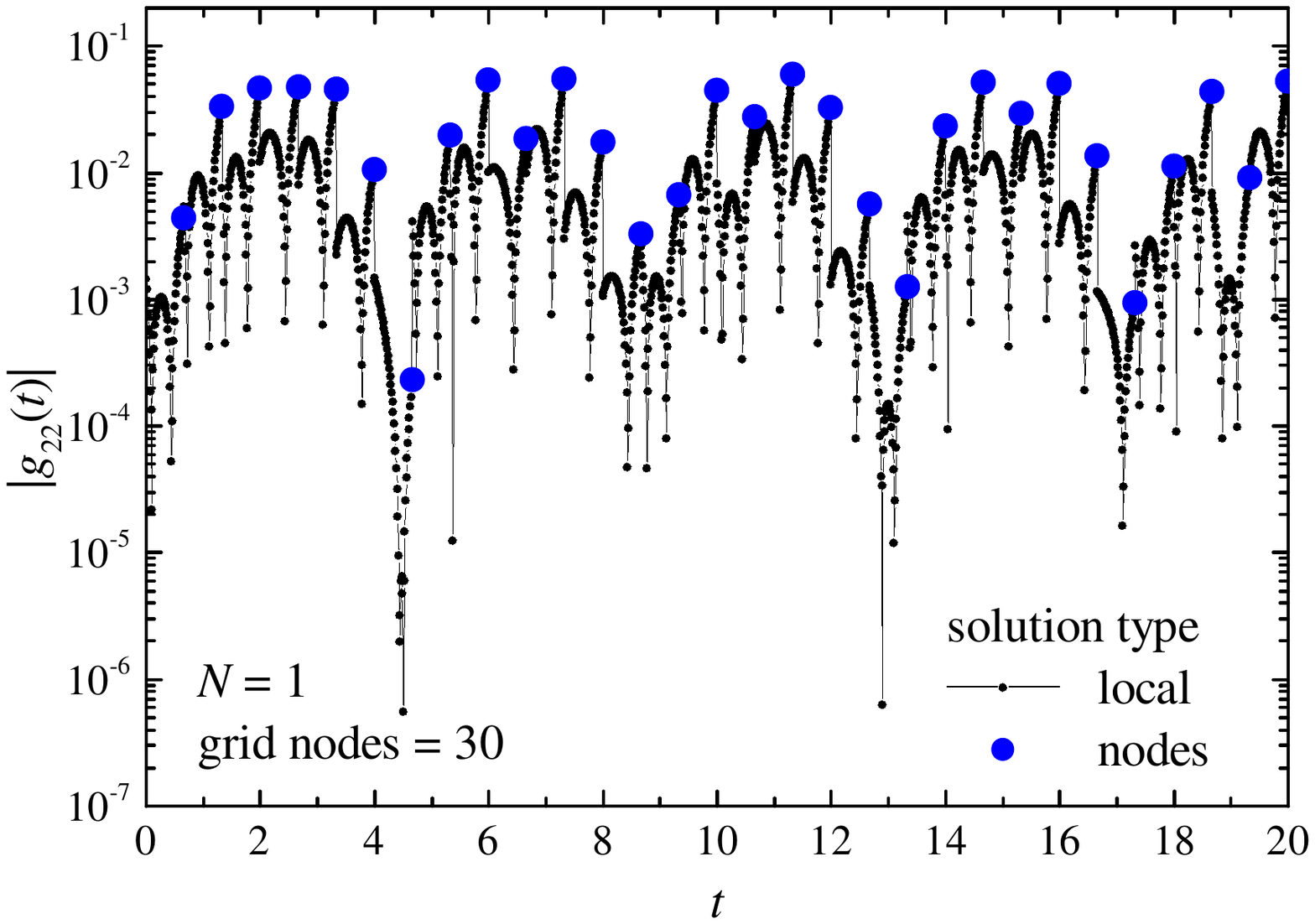}
\vspace{-8mm}\caption{\label{fig:dpend_ind3_sols_vg:b2}}
\end{subfigure}
\begin{subfigure}{0.240\textwidth}
\includegraphics[width=\textwidth]{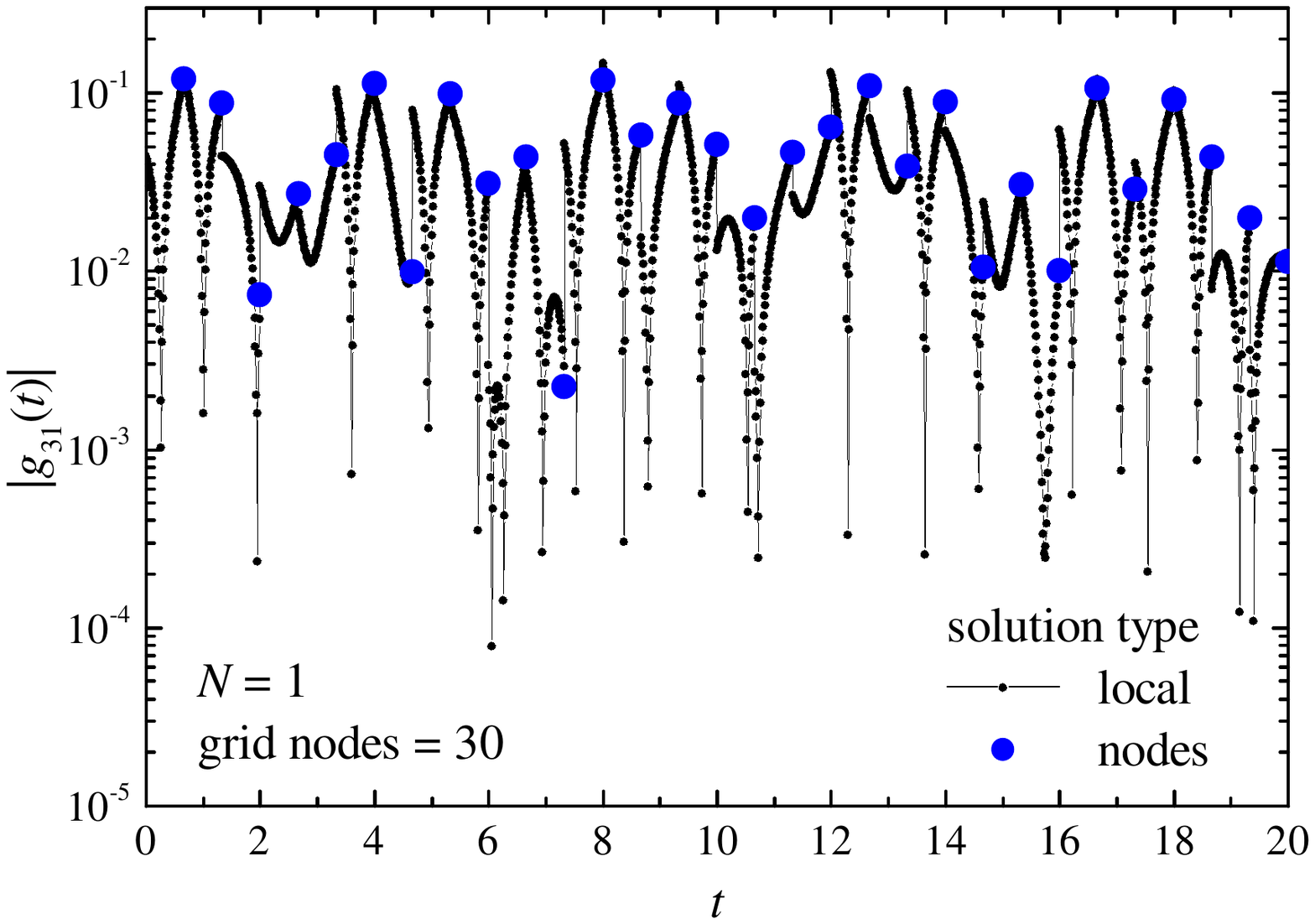}
\vspace{-8mm}\caption{\label{fig:dpend_ind3_sols_vg:b3}}
\end{subfigure}
\begin{subfigure}{0.240\textwidth}
\includegraphics[width=\textwidth]{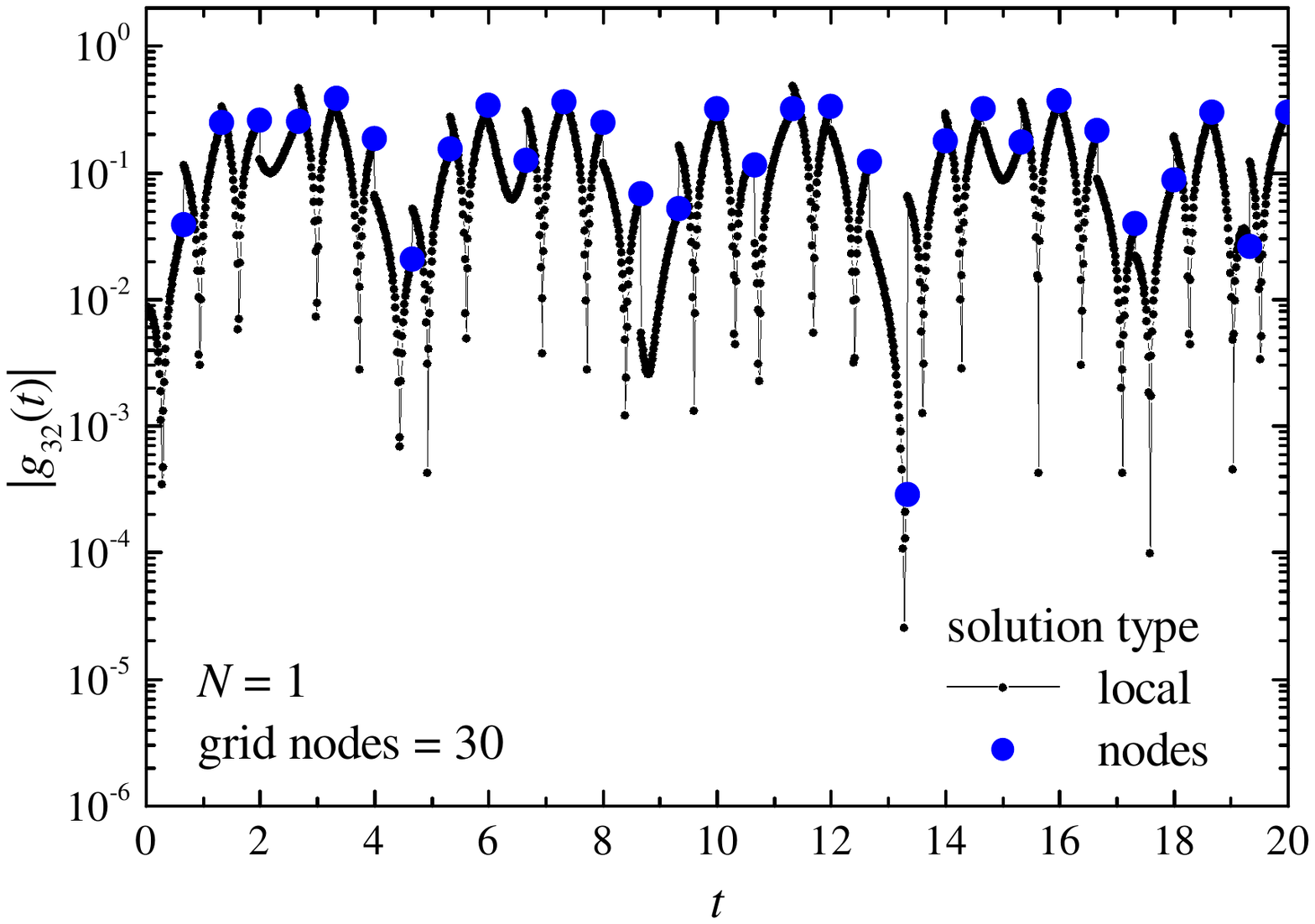}
\vspace{-8mm}\caption{\label{fig:dpend_ind3_sols_vg:b4}}
\end{subfigure}\\[2mm]
\begin{subfigure}{0.240\textwidth}
\includegraphics[width=\textwidth]{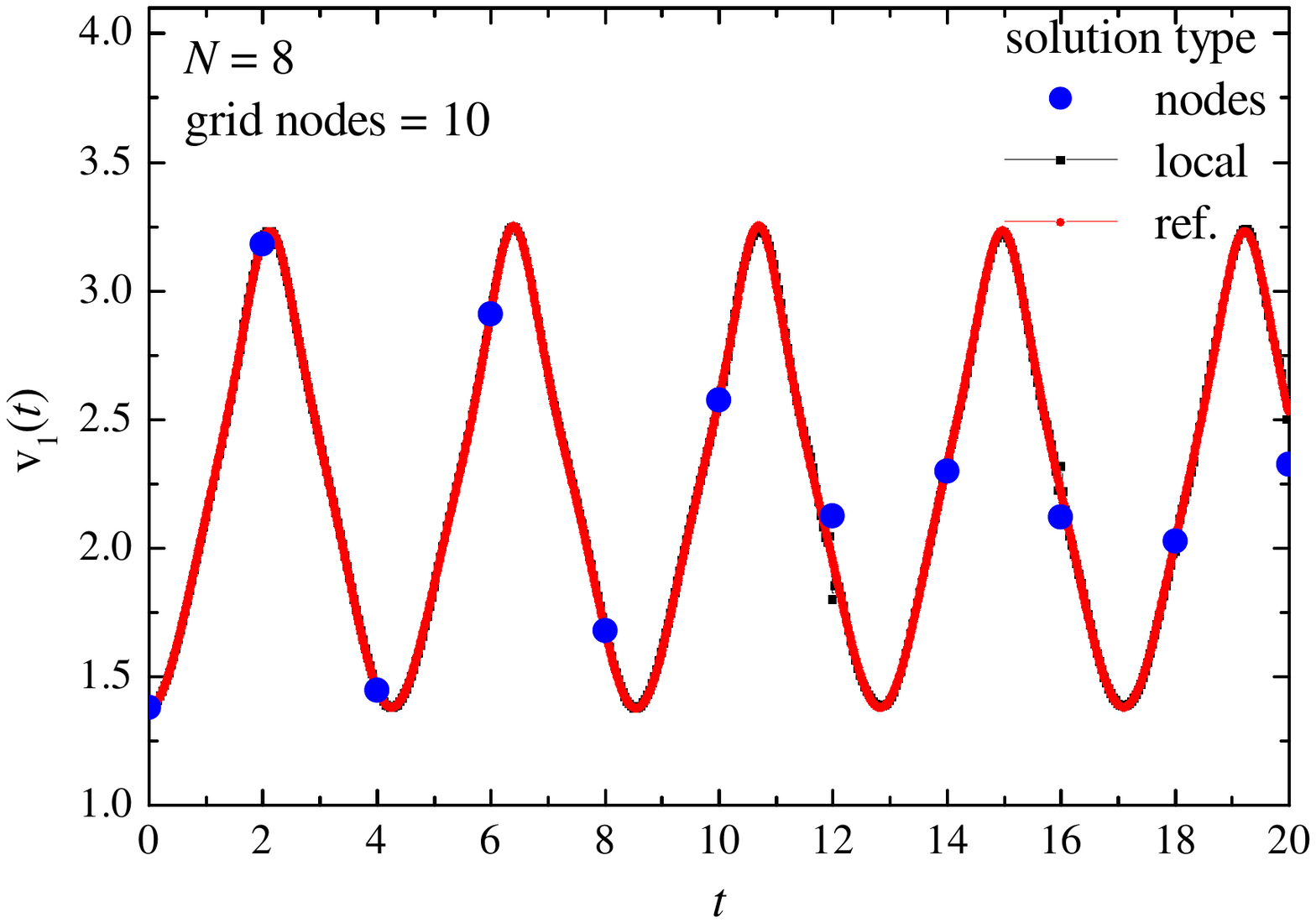}
\vspace{-8mm}\caption{\label{fig:dpend_ind3_sols_vg:c1}}
\end{subfigure}
\begin{subfigure}{0.240\textwidth}
\includegraphics[width=\textwidth]{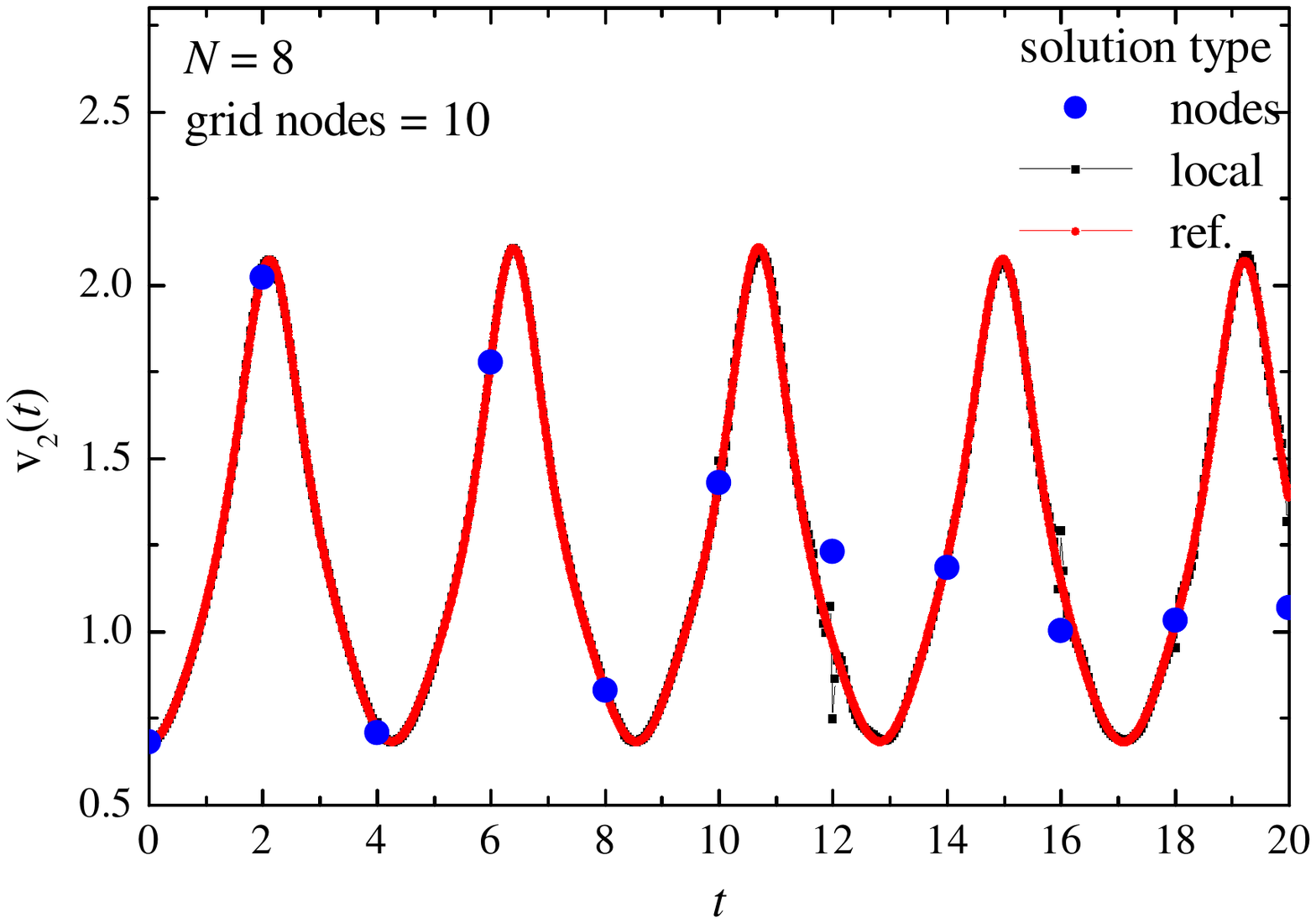}
\vspace{-8mm}\caption{\label{fig:dpend_ind3_sols_vg:c2}}
\end{subfigure}
\begin{subfigure}{0.240\textwidth}
\includegraphics[width=\textwidth]{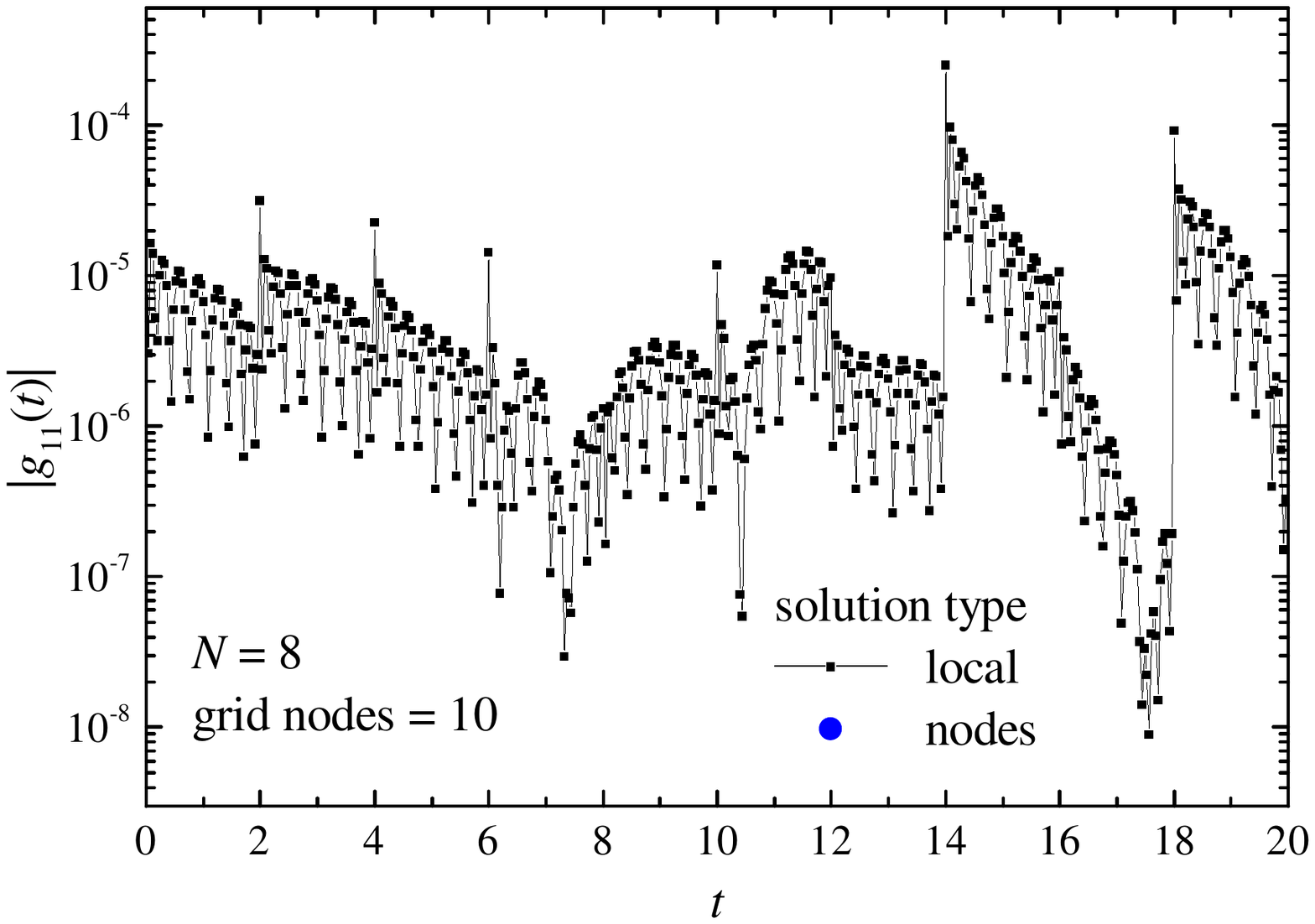}
\vspace{-8mm}\caption{\label{fig:dpend_ind3_sols_vg:c3}}
\end{subfigure}
\begin{subfigure}{0.240\textwidth}
\includegraphics[width=\textwidth]{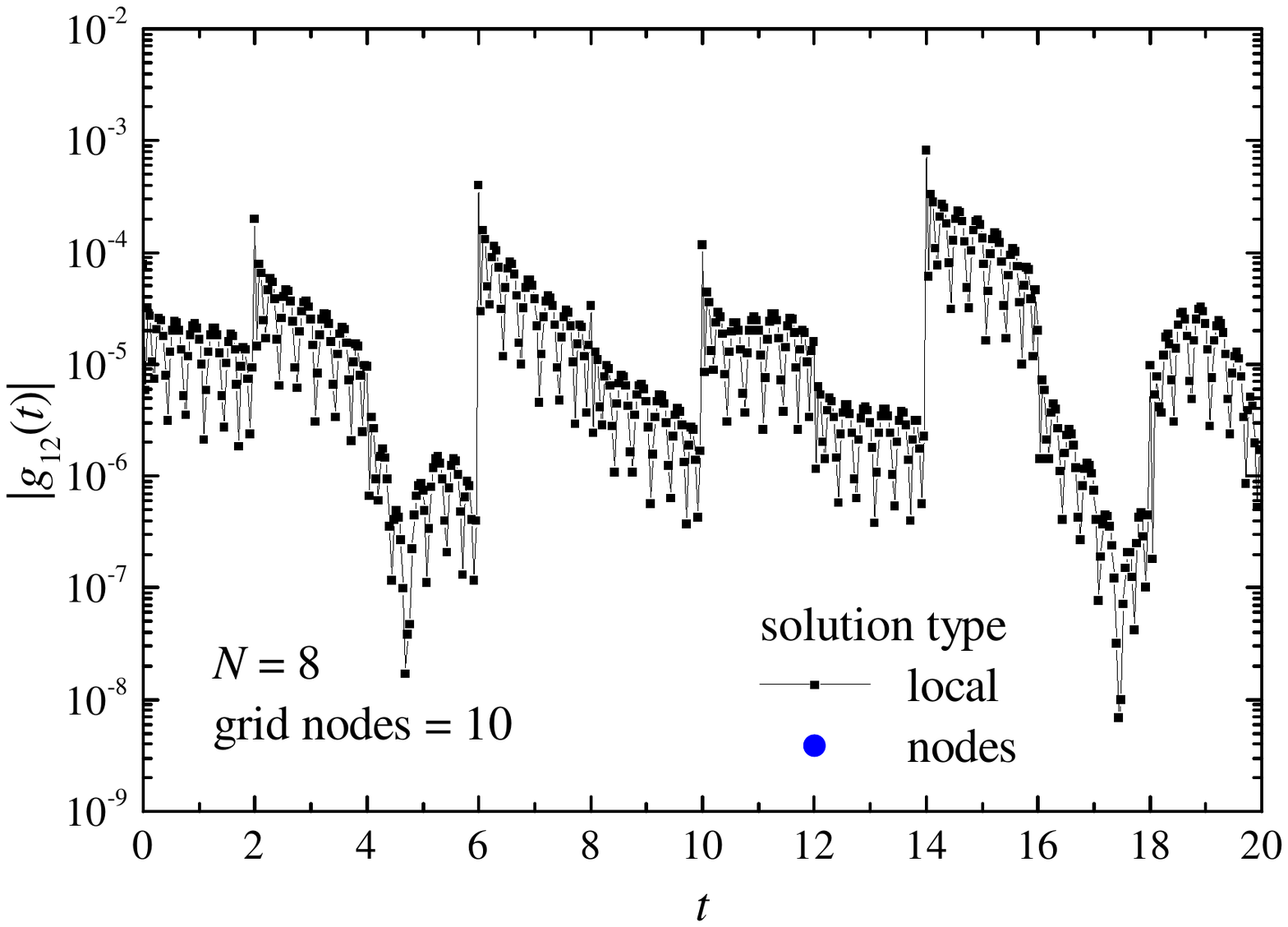}
\vspace{-8mm}\caption{\label{fig:dpend_ind3_sols_vg:c4}}
\end{subfigure}\\[2mm]
\begin{subfigure}{0.240\textwidth}
\includegraphics[width=\textwidth]{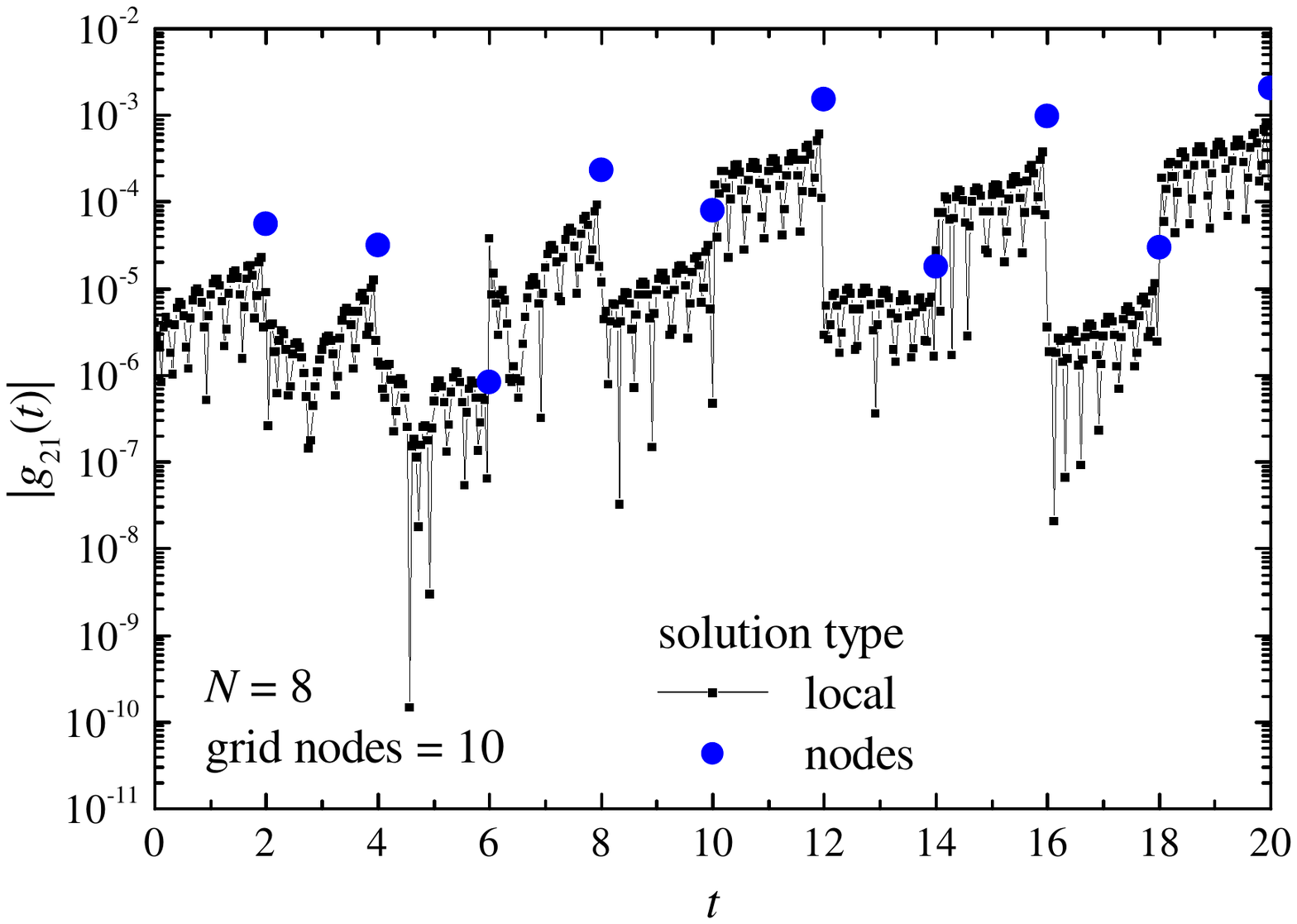}
\vspace{-8mm}\caption{\label{fig:dpend_ind3_sols_vg:d1}}
\end{subfigure}
\begin{subfigure}{0.240\textwidth}
\includegraphics[width=\textwidth]{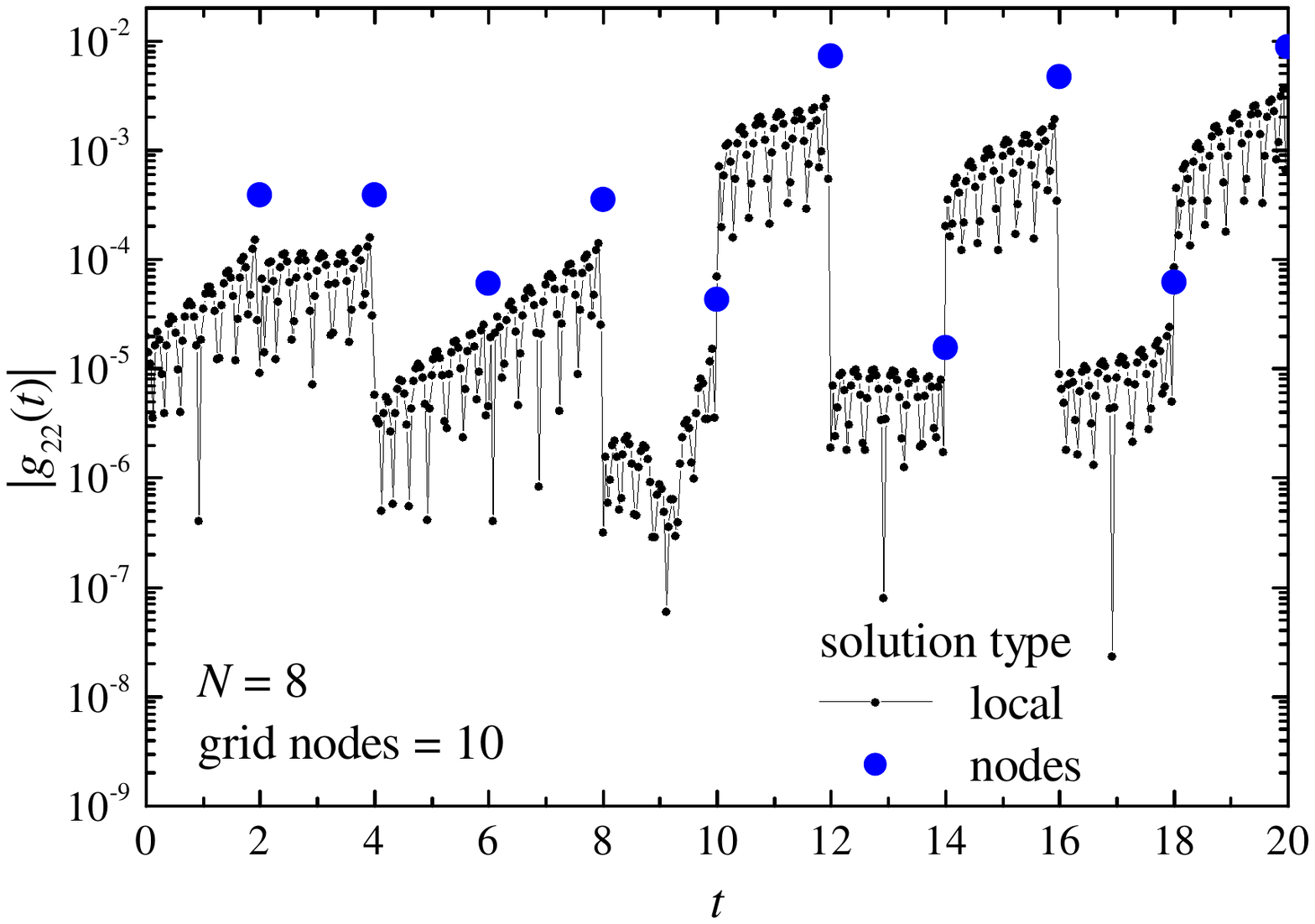}
\vspace{-8mm}\caption{\label{fig:dpend_ind3_sols_vg:d2}}
\end{subfigure}
\begin{subfigure}{0.240\textwidth}
\includegraphics[width=\textwidth]{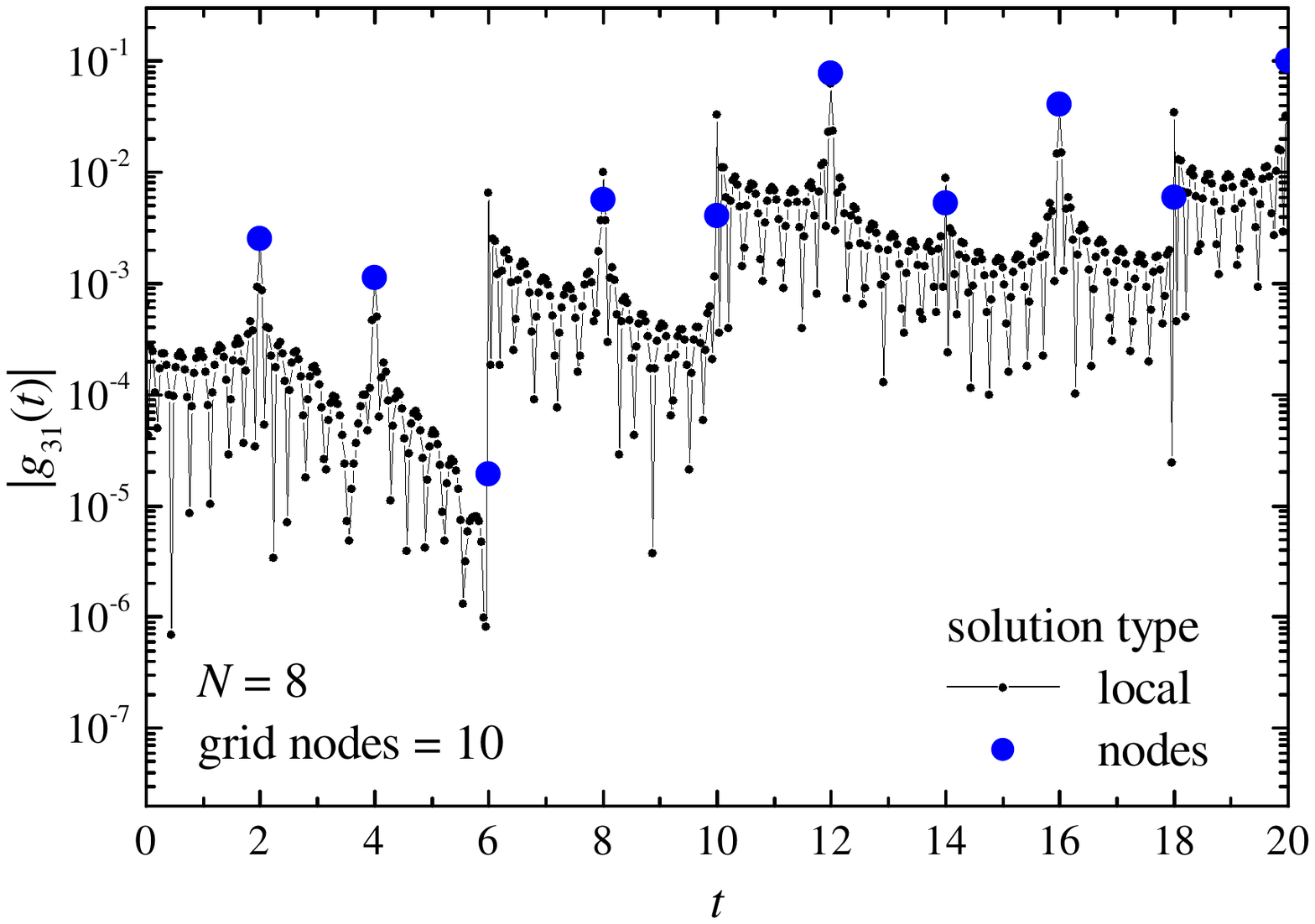}
\vspace{-8mm}\caption{\label{fig:dpend_ind3_sols_vg:d3}}
\end{subfigure}
\begin{subfigure}{0.240\textwidth}
\includegraphics[width=\textwidth]{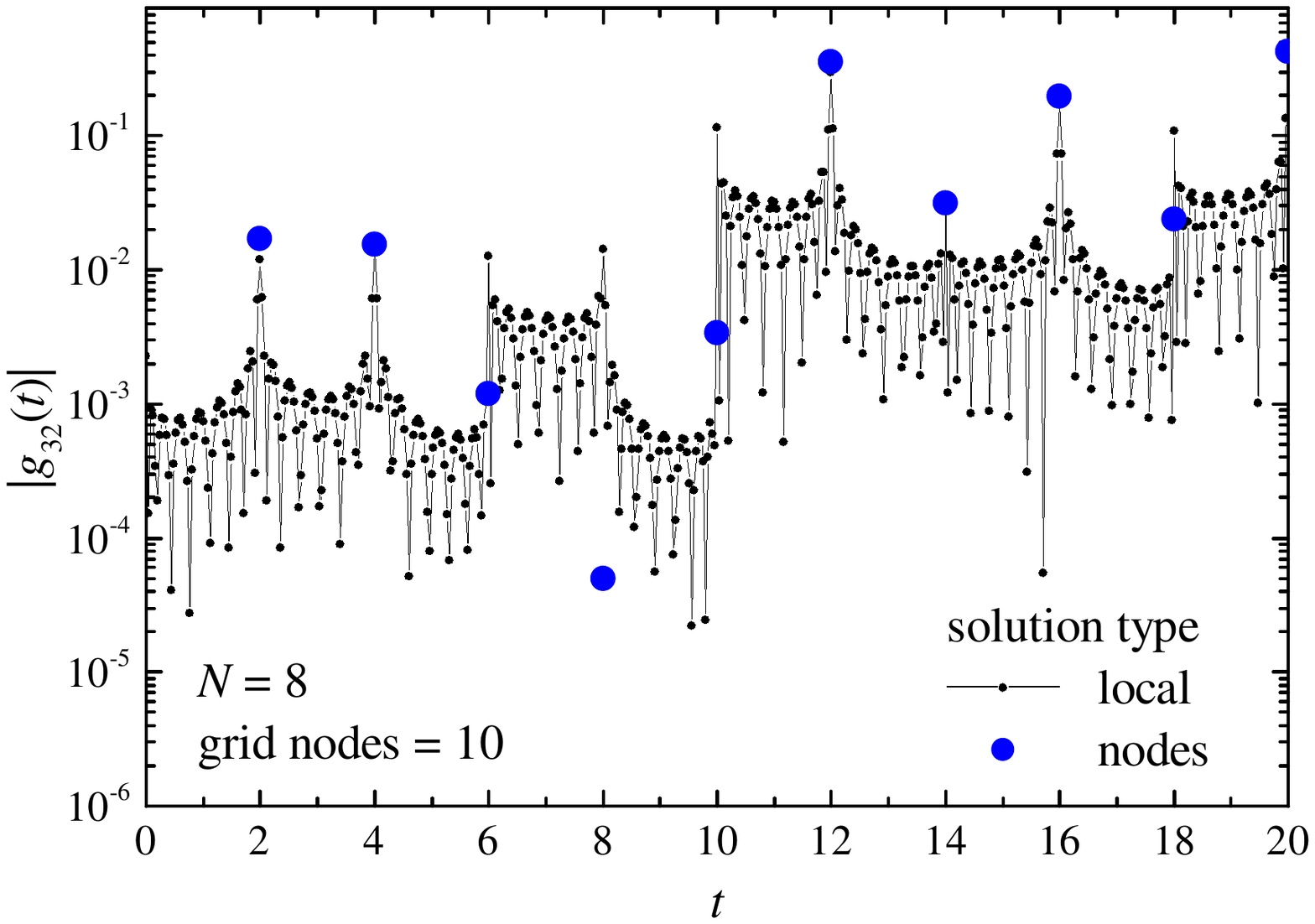}
\vspace{-8mm}\caption{\label{fig:dpend_ind3_sols_vg:d4}}
\end{subfigure}\\[2mm]
\begin{subfigure}{0.240\textwidth}
\includegraphics[width=\textwidth]{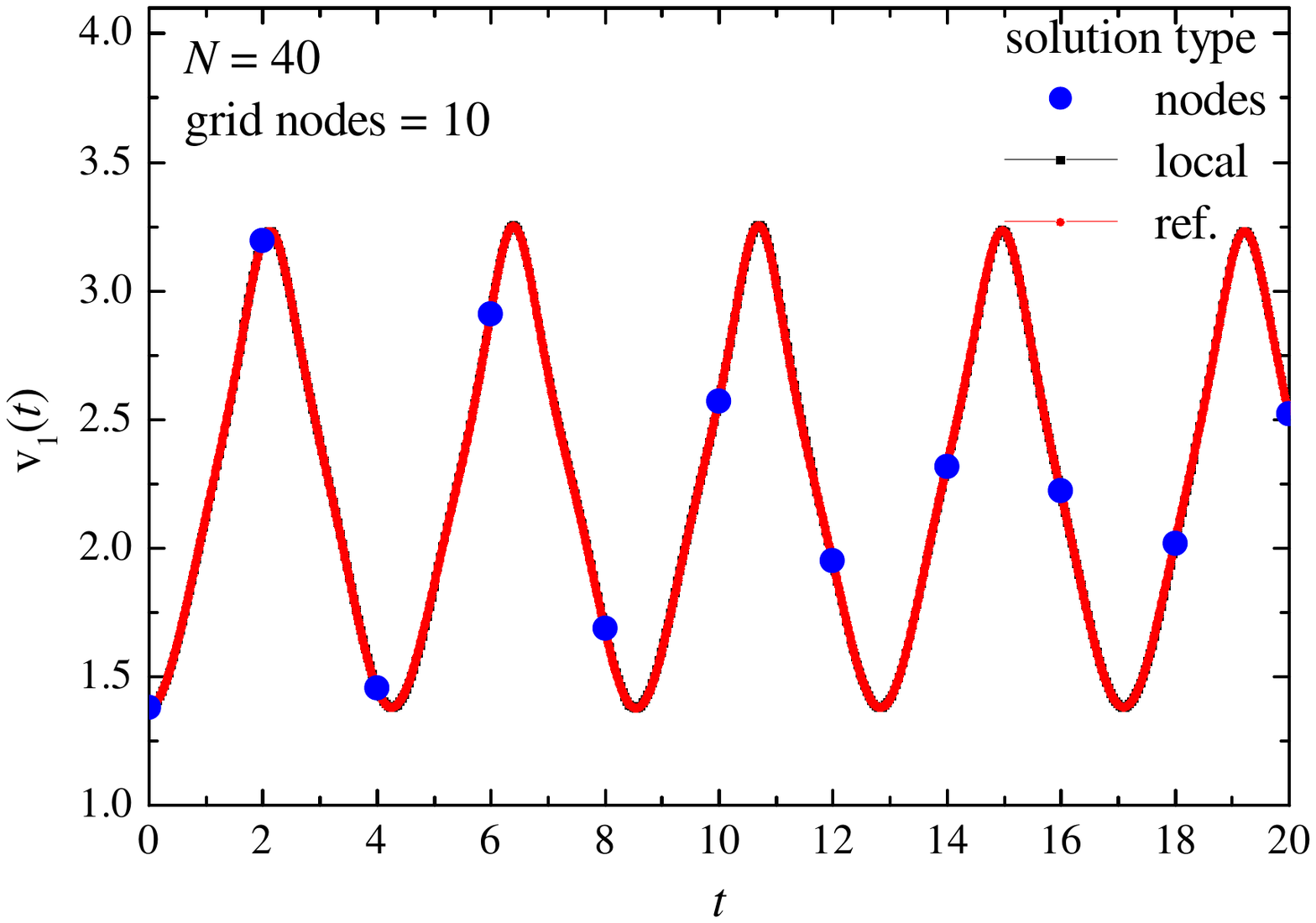}
\vspace{-8mm}\caption{\label{fig:dpend_ind3_sols_vg:e1}}
\end{subfigure}
\begin{subfigure}{0.240\textwidth}
\includegraphics[width=\textwidth]{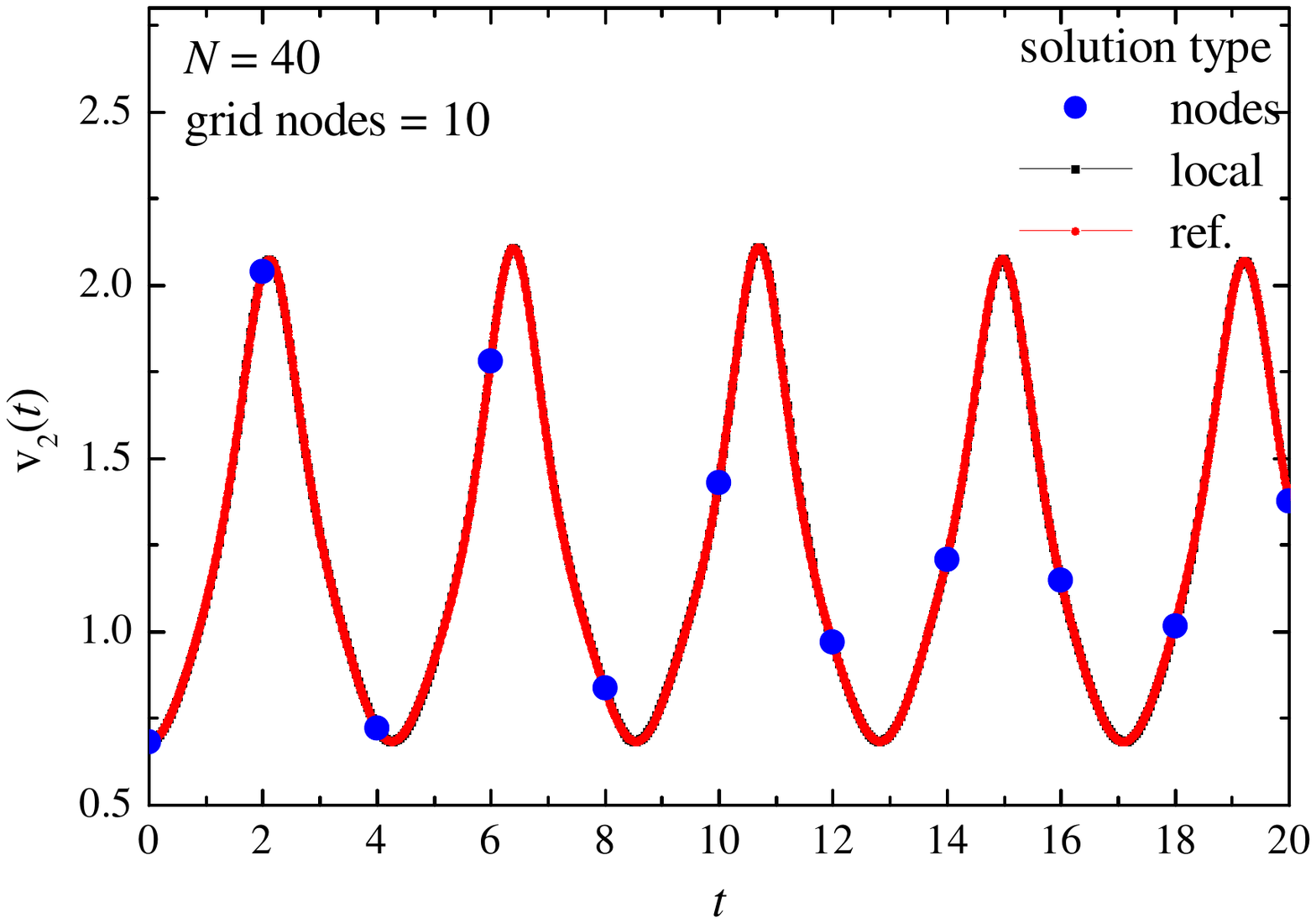}
\vspace{-8mm}\caption{\label{fig:dpend_ind3_sols_vg:e2}}
\end{subfigure}
\begin{subfigure}{0.240\textwidth}
\includegraphics[width=\textwidth]{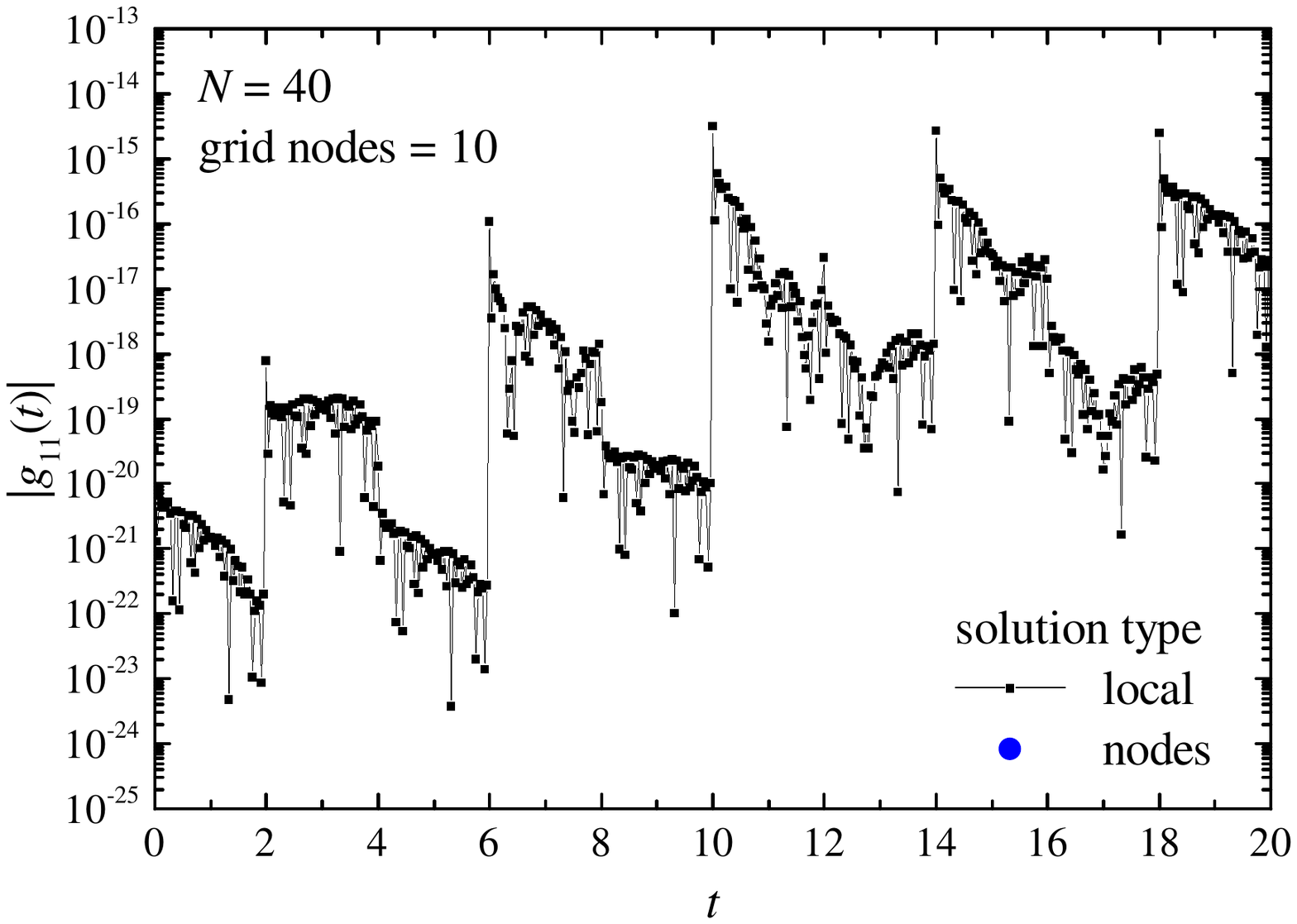}
\vspace{-8mm}\caption{\label{fig:dpend_ind3_sols_vg:e3}}
\end{subfigure}
\begin{subfigure}{0.240\textwidth}
\includegraphics[width=\textwidth]{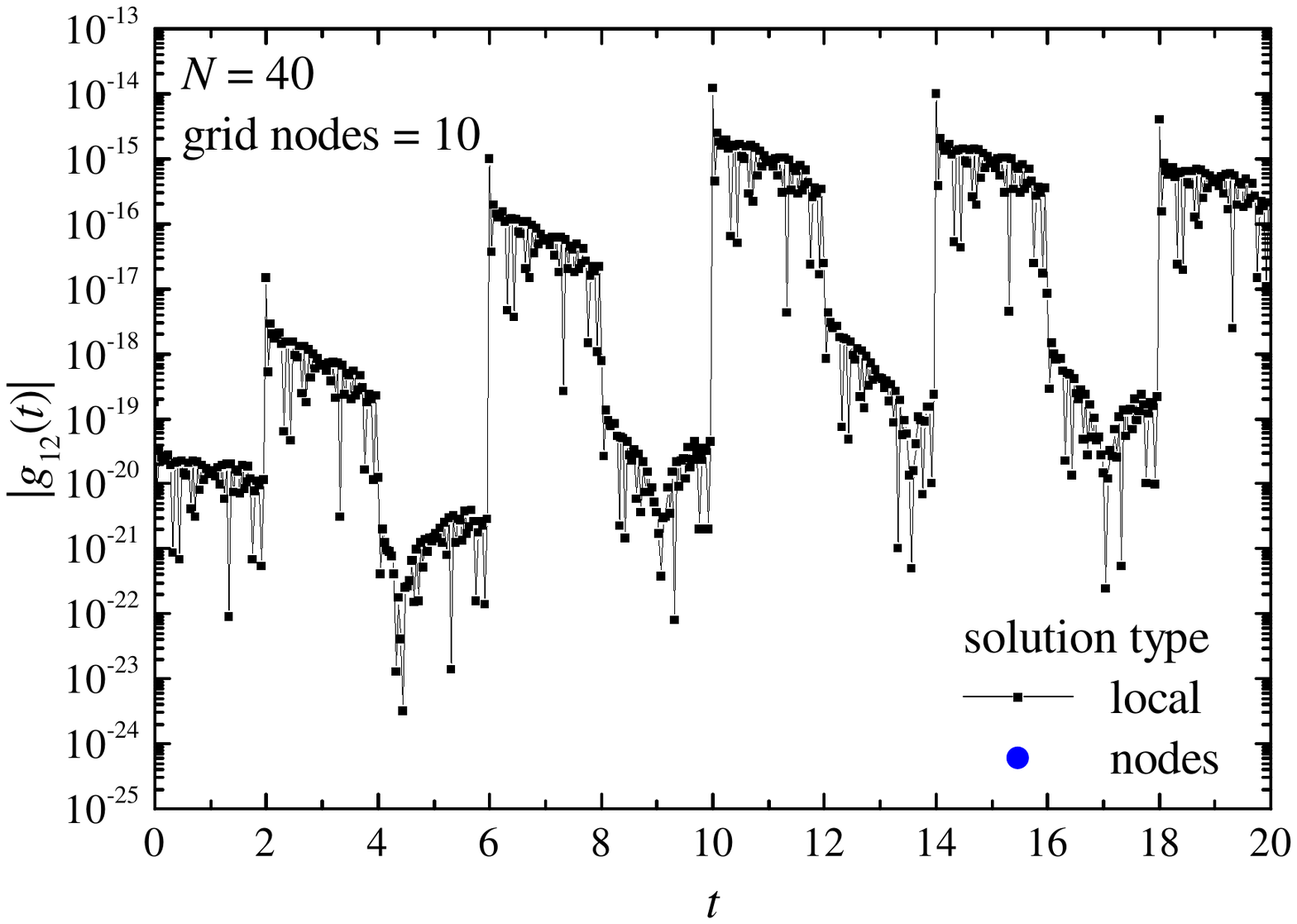}
\vspace{-8mm}\caption{\label{fig:dpend_ind3_sols_vg:e4}}
\end{subfigure}\\[2mm]
\begin{subfigure}{0.240\textwidth}
\includegraphics[width=\textwidth]{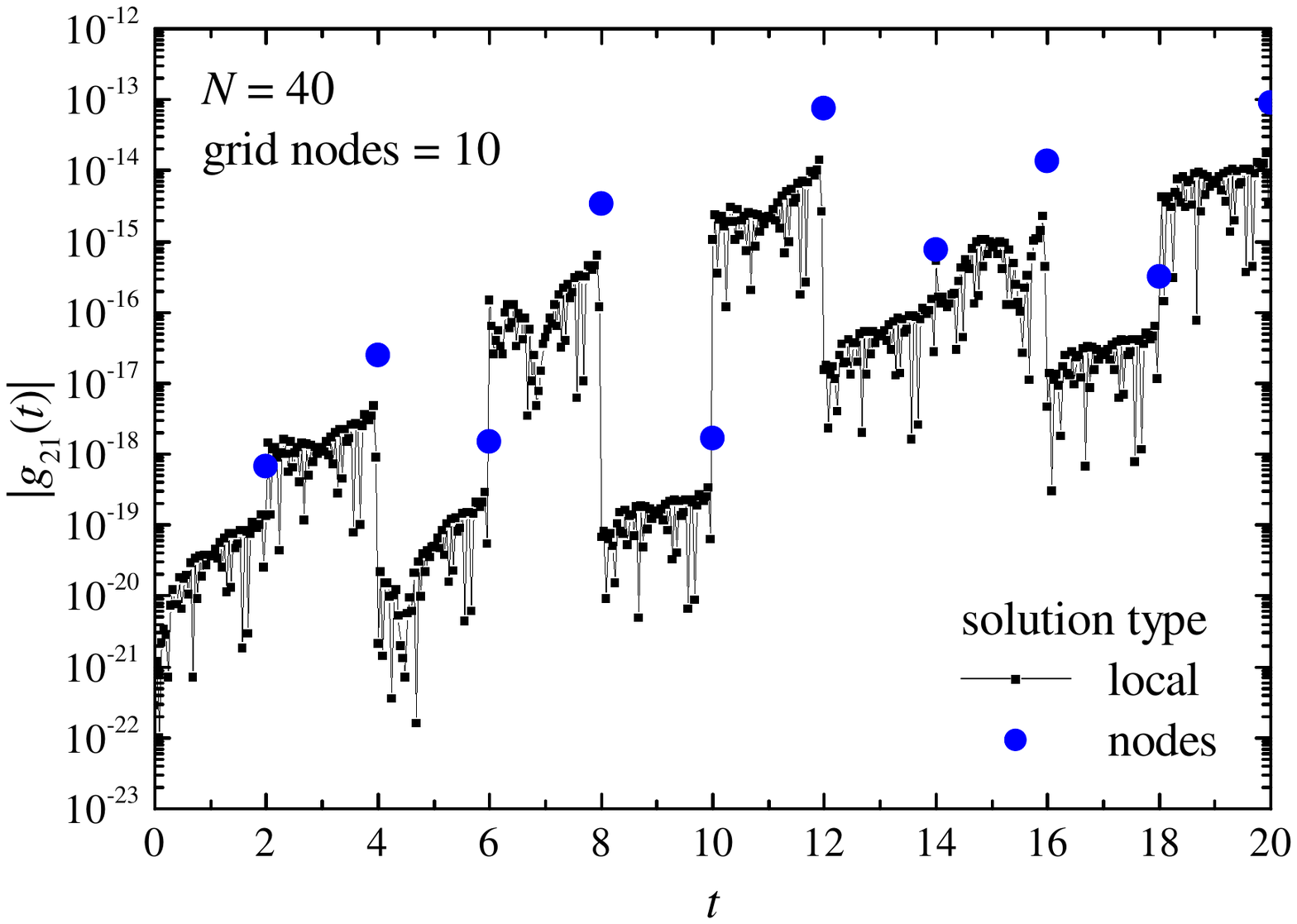}
\vspace{-8mm}\caption{\label{fig:dpend_ind3_sols_vg:f1}}
\end{subfigure}
\begin{subfigure}{0.240\textwidth}
\includegraphics[width=\textwidth]{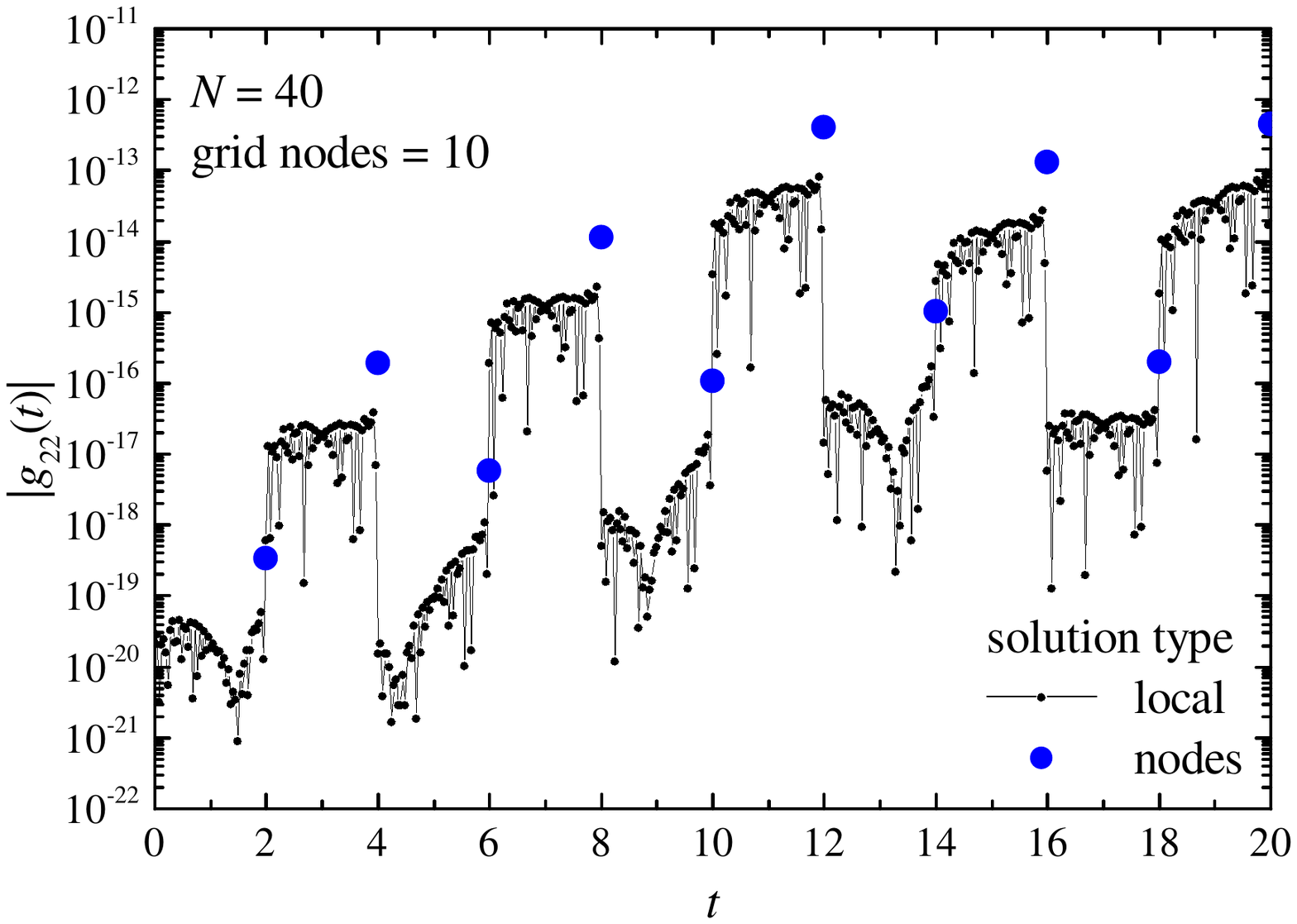}
\vspace{-8mm}\caption{\label{fig:dpend_ind3_sols_vg:f2}}
\end{subfigure}
\begin{subfigure}{0.240\textwidth}
\includegraphics[width=\textwidth]{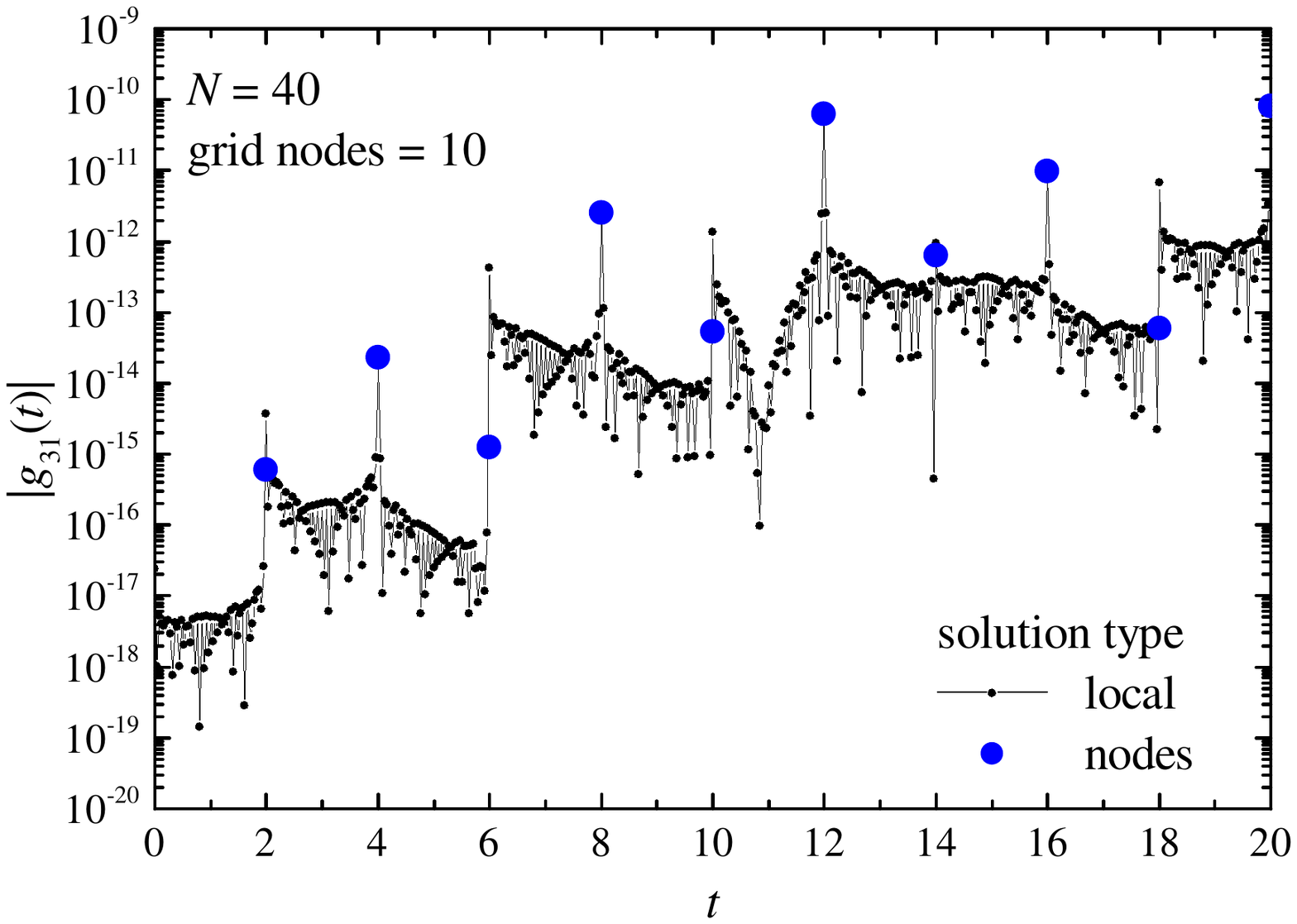}
\vspace{-8mm}\caption{\label{fig:dpend_ind3_sols_vg:f3}}
\end{subfigure}
\begin{subfigure}{0.240\textwidth}
\includegraphics[width=\textwidth]{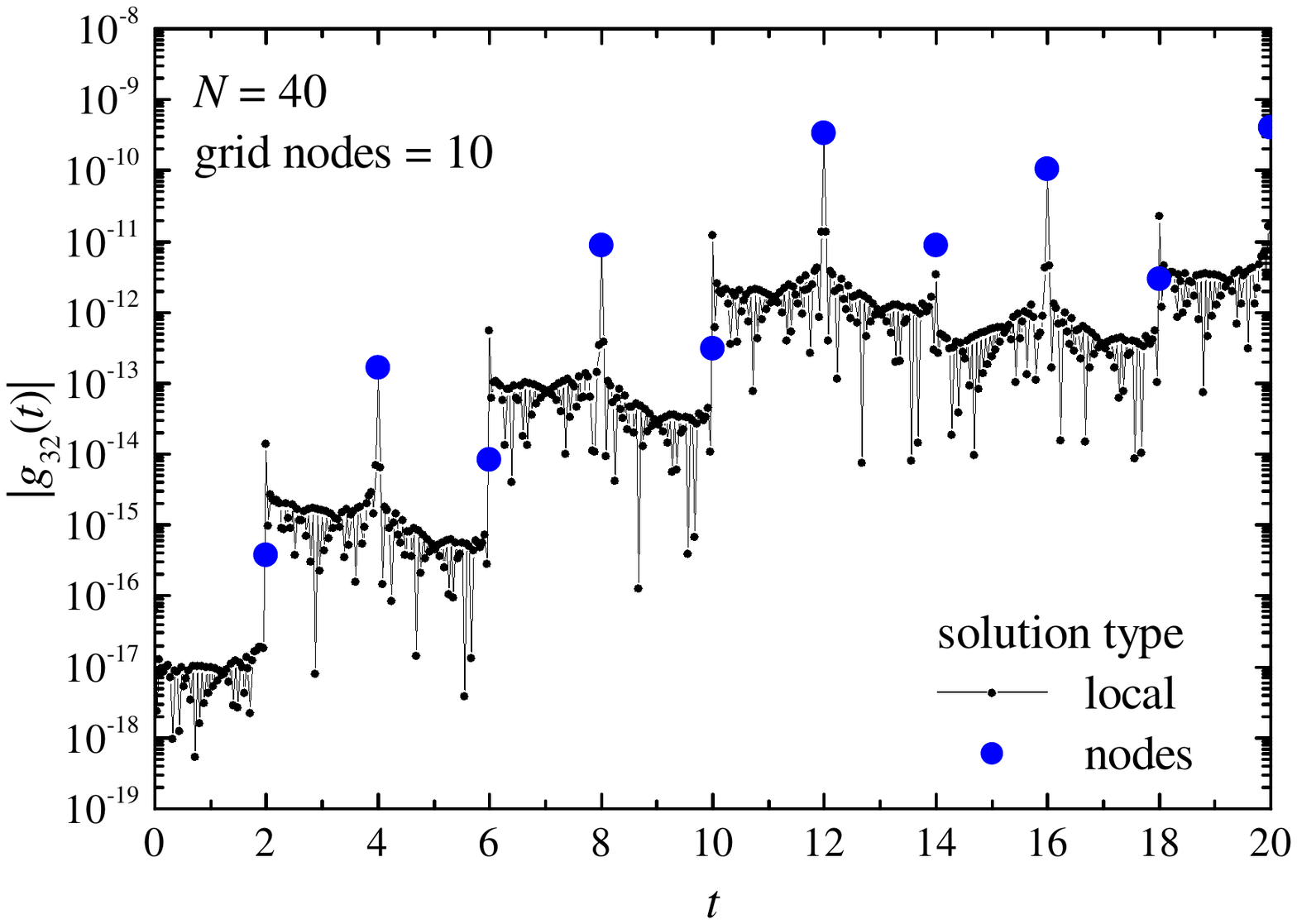}
\vspace{-8mm}\caption{\label{fig:dpend_ind3_sols_vg:f4}}
\end{subfigure}\\[2mm]
\caption{%
Numerical solution of the DAE system (\ref{eq:math_dpend_dae_ind_3}) of index 3. Comparison of the solution at nodes $\mathbf{v}_{n}$, the local solution $\mathbf{v}_{L}(t)$ and the reference solution $\mathbf{v}^{\rm ref}(t)$ for components $v_{1}$ (\subref{fig:dpend_ind3_sols_vg:a1}, \subref{fig:dpend_ind3_sols_vg:c1}, \subref{fig:dpend_ind3_sols_vg:e1}), $v_{2}$ (\subref{fig:dpend_ind3_sols_vg:a2}, \subref{fig:dpend_ind3_sols_vg:c2}, \subref{fig:dpend_ind3_sols_vg:e2}), quantitative satisfiability of the conditions $g_{11} = 0$ (\subref{fig:dpend_ind3_sols_vg:a3}, \subref{fig:dpend_ind3_sols_vg:c3}, \subref{fig:dpend_ind3_sols_vg:e3}), $g_{12} = 0$ (\subref{fig:dpend_ind3_sols_vg:a4}, \subref{fig:dpend_ind3_sols_vg:c4}, \subref{fig:dpend_ind3_sols_vg:e4}), $g_{21} = 0$ (\subref{fig:dpend_ind3_sols_vg:b1}, \subref{fig:dpend_ind3_sols_vg:d1}, \subref{fig:dpend_ind3_sols_vg:f1}), $g_{22} = 0$ (\subref{fig:dpend_ind3_sols_vg:b2}, \subref{fig:dpend_ind3_sols_vg:d2}, \subref{fig:dpend_ind3_sols_vg:f2}), $g_{31} = 0$ (\subref{fig:dpend_ind3_sols_vg:b3}, \subref{fig:dpend_ind3_sols_vg:d3}, \subref{fig:dpend_ind3_sols_vg:f3}), $g_{32} = 0$ (\subref{fig:dpend_ind3_sols_vg:b4}, \subref{fig:dpend_ind3_sols_vg:d4}, \subref{fig:dpend_ind3_sols_vg:f4}), obtained using polynomials with degrees $N = 1$ (\subref{fig:dpend_ind3_sols_vg:a1}, \subref{fig:dpend_ind3_sols_vg:a2}, \subref{fig:dpend_ind3_sols_vg:a3}, \subref{fig:dpend_ind3_sols_vg:a4}, \subref{fig:dpend_ind3_sols_vg:b1}, \subref{fig:dpend_ind3_sols_vg:b2}, \subref{fig:dpend_ind3_sols_vg:b3}, \subref{fig:dpend_ind3_sols_vg:b4}), $N = 8$ (\subref{fig:dpend_ind3_sols_vg:c1}, \subref{fig:dpend_ind3_sols_vg:c2}, \subref{fig:dpend_ind3_sols_vg:c3}, \subref{fig:dpend_ind3_sols_vg:c4}, \subref{fig:dpend_ind3_sols_vg:d1}, \subref{fig:dpend_ind3_sols_vg:d2}, \subref{fig:dpend_ind3_sols_vg:d3}, \subref{fig:dpend_ind3_sols_vg:d4}) and $N = 40$ (\subref{fig:dpend_ind3_sols_vg:e1}, \subref{fig:dpend_ind3_sols_vg:e2}, \subref{fig:dpend_ind3_sols_vg:e3}, \subref{fig:dpend_ind3_sols_vg:e4}, \subref{fig:dpend_ind3_sols_vg:f1}, \subref{fig:dpend_ind3_sols_vg:f2}, \subref{fig:dpend_ind3_sols_vg:f3}, \subref{fig:dpend_ind3_sols_vg:f4}).
}
\label{fig:dpend_ind3_sols_vg}
\end{figure} 

\begin{figure}[h!]
\captionsetup[subfigure]{%
	position=bottom,
	font+=smaller,
	textfont=normalfont,
	singlelinecheck=off,
	justification=raggedright
}
\centering
\begin{subfigure}{0.320\textwidth}
\includegraphics[width=\textwidth]{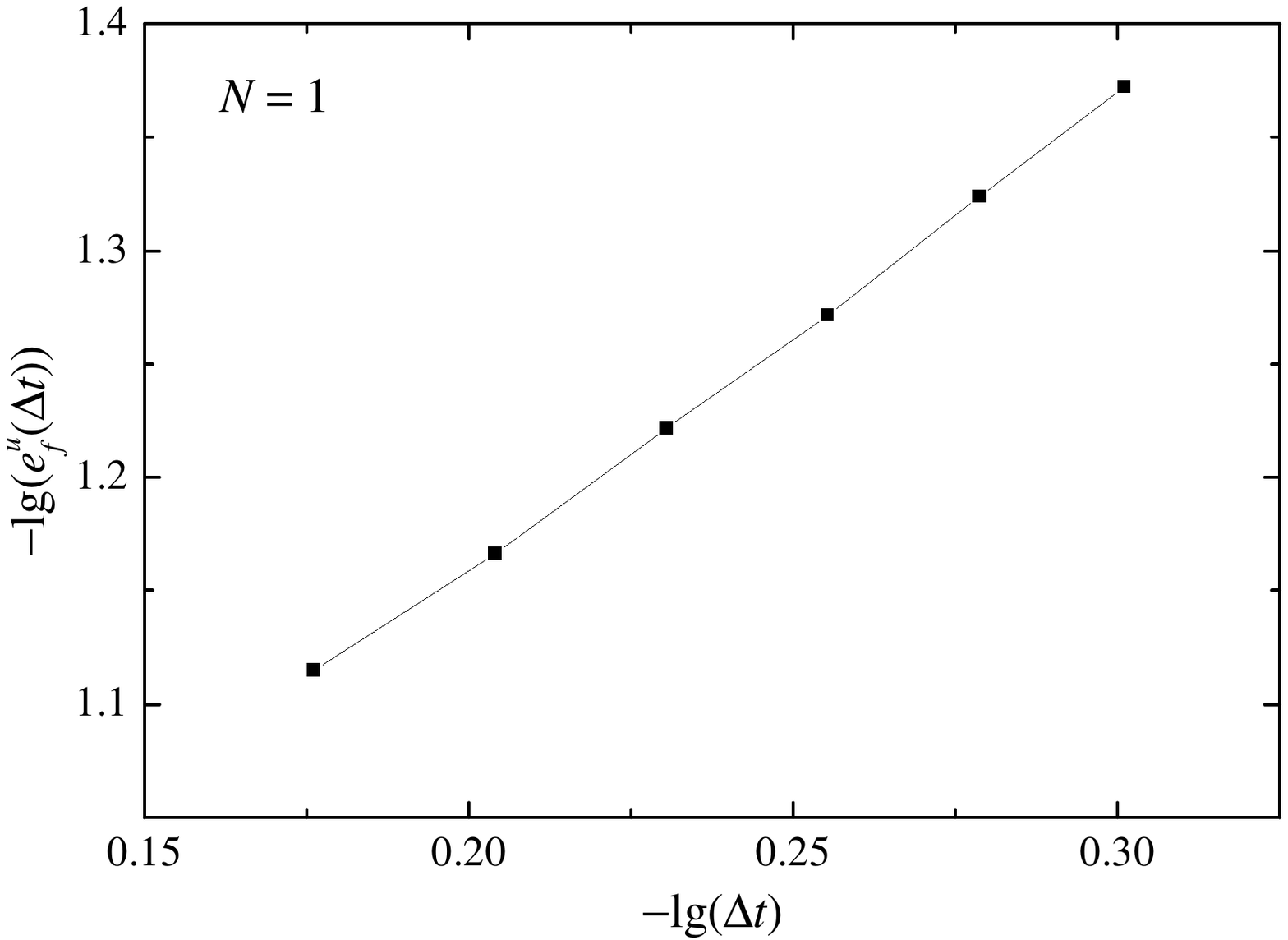}
\vspace{-8mm}\caption{\label{fig:dpend_ind3_errors:a1}}
\end{subfigure}
\begin{subfigure}{0.320\textwidth}
\includegraphics[width=\textwidth]{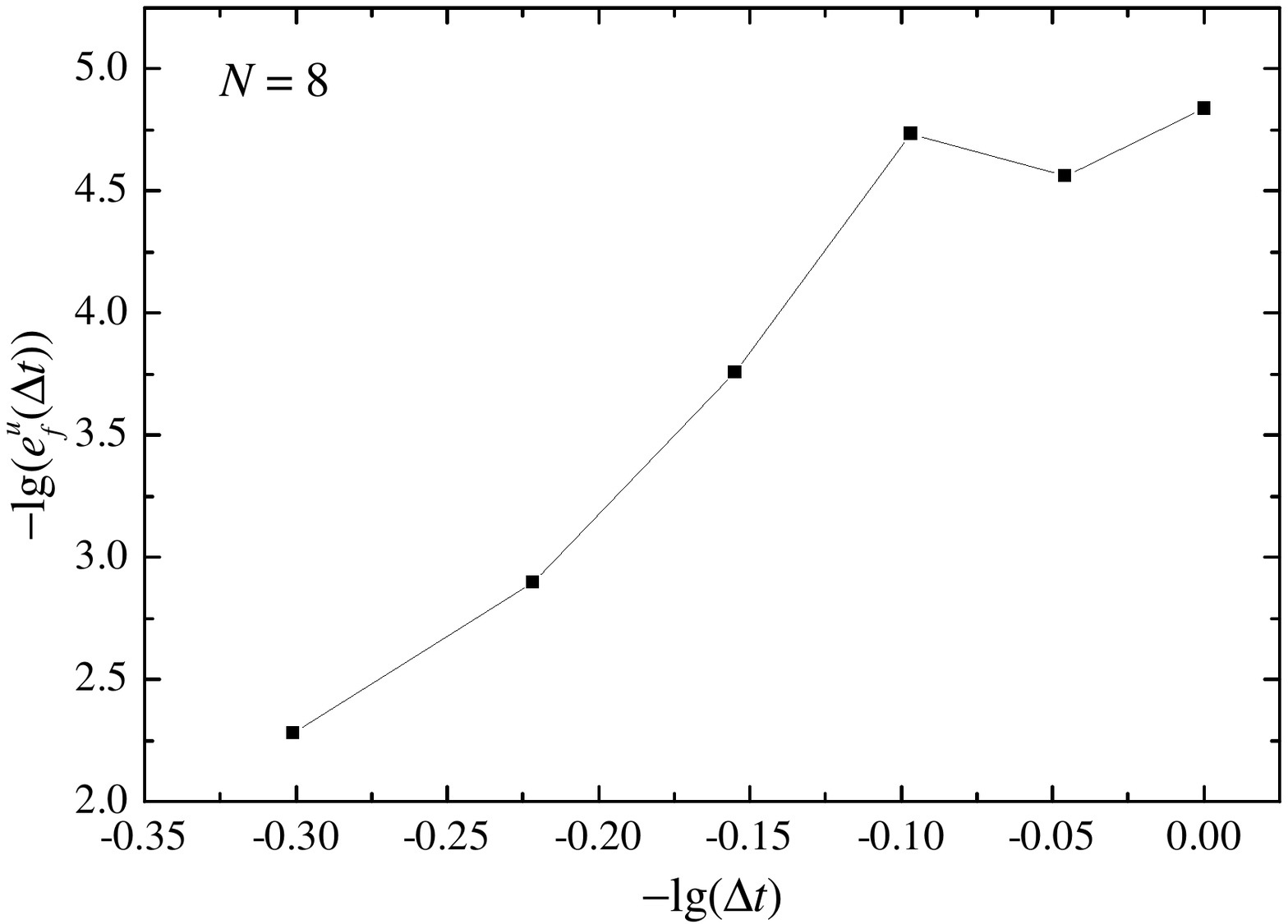}
\vspace{-8mm}\caption{\label{fig:dpend_ind3_errors:a2}}
\end{subfigure}
\begin{subfigure}{0.320\textwidth}
\includegraphics[width=\textwidth]{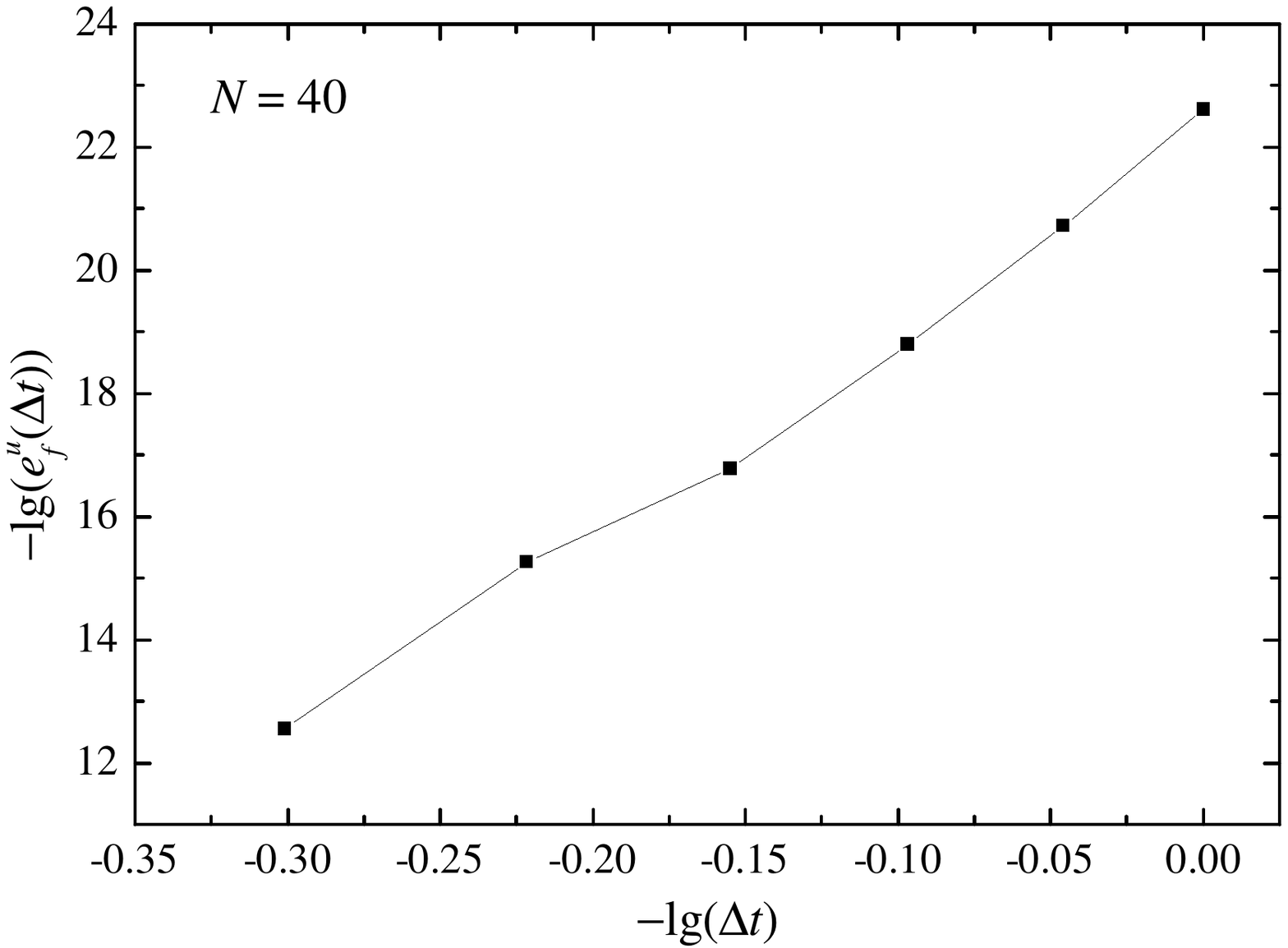}
\vspace{-8mm}\caption{\label{fig:dpend_ind3_errors:a3}}
\end{subfigure}\\[2mm]
\begin{subfigure}{0.320\textwidth}
\includegraphics[width=\textwidth]{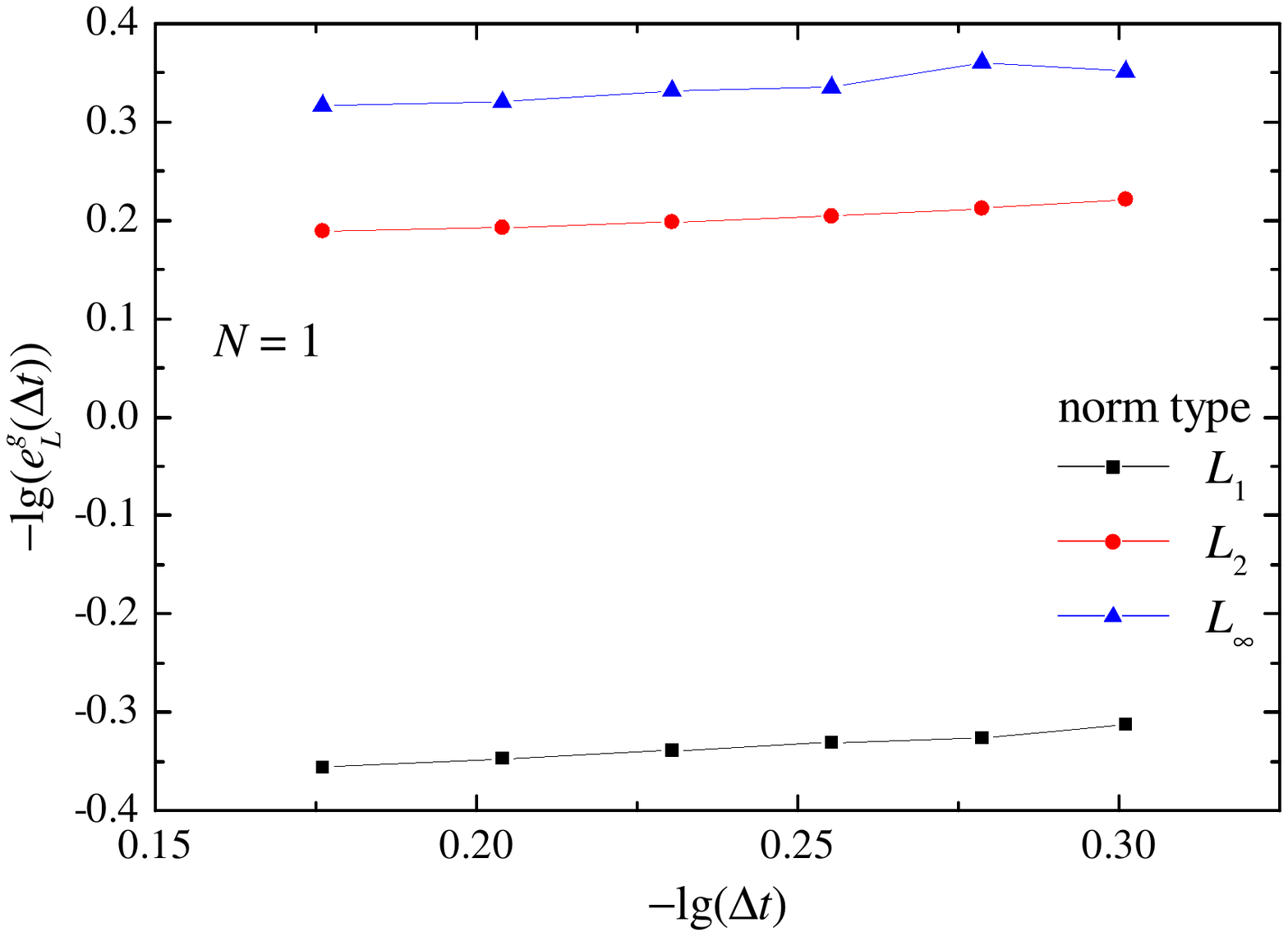}
\vspace{-8mm}\caption{\label{fig:dpend_ind3_errors:b1}}
\end{subfigure}
\begin{subfigure}{0.320\textwidth}
\includegraphics[width=\textwidth]{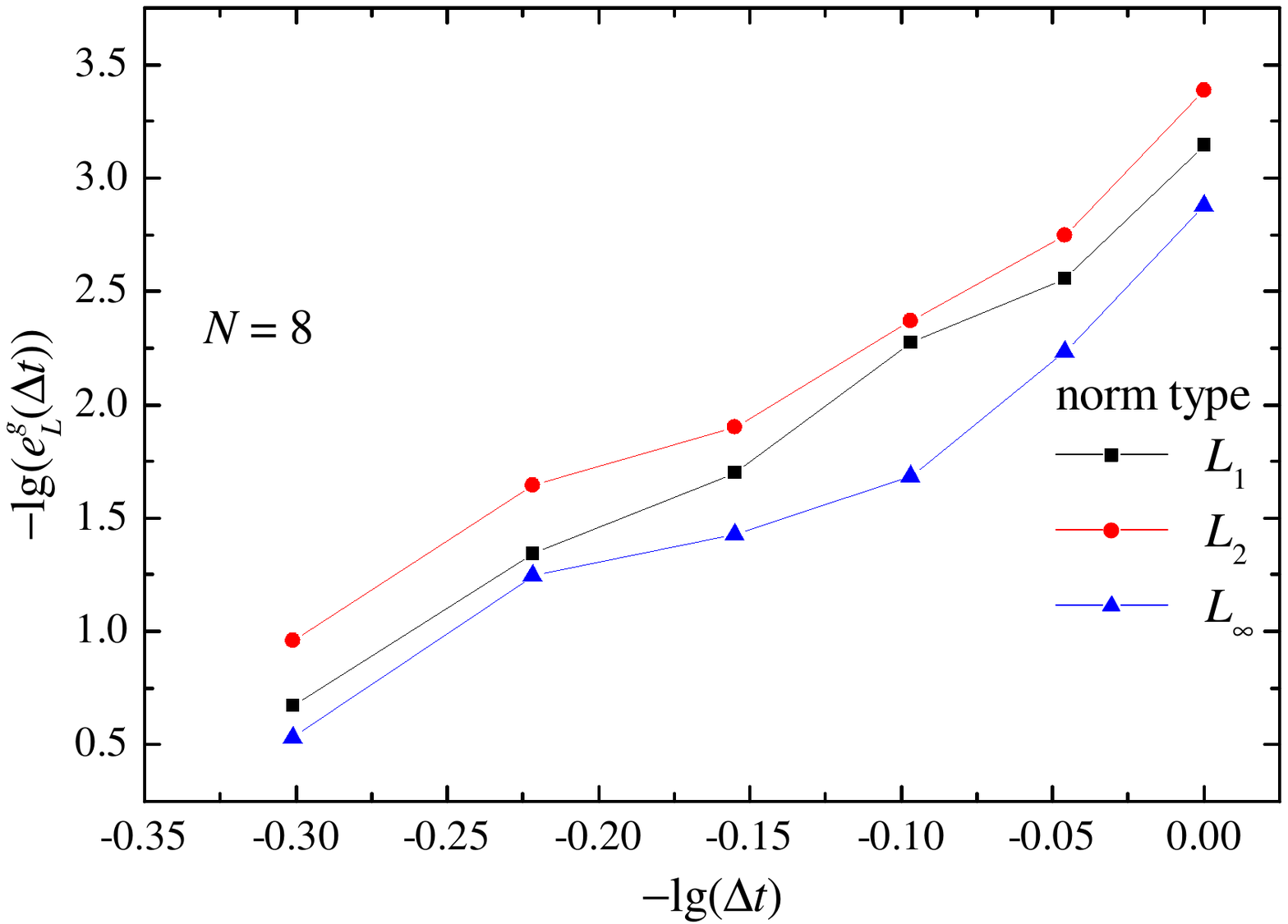}
\vspace{-8mm}\caption{\label{fig:dpend_ind3_errors:b2}}
\end{subfigure}
\begin{subfigure}{0.320\textwidth}
\includegraphics[width=\textwidth]{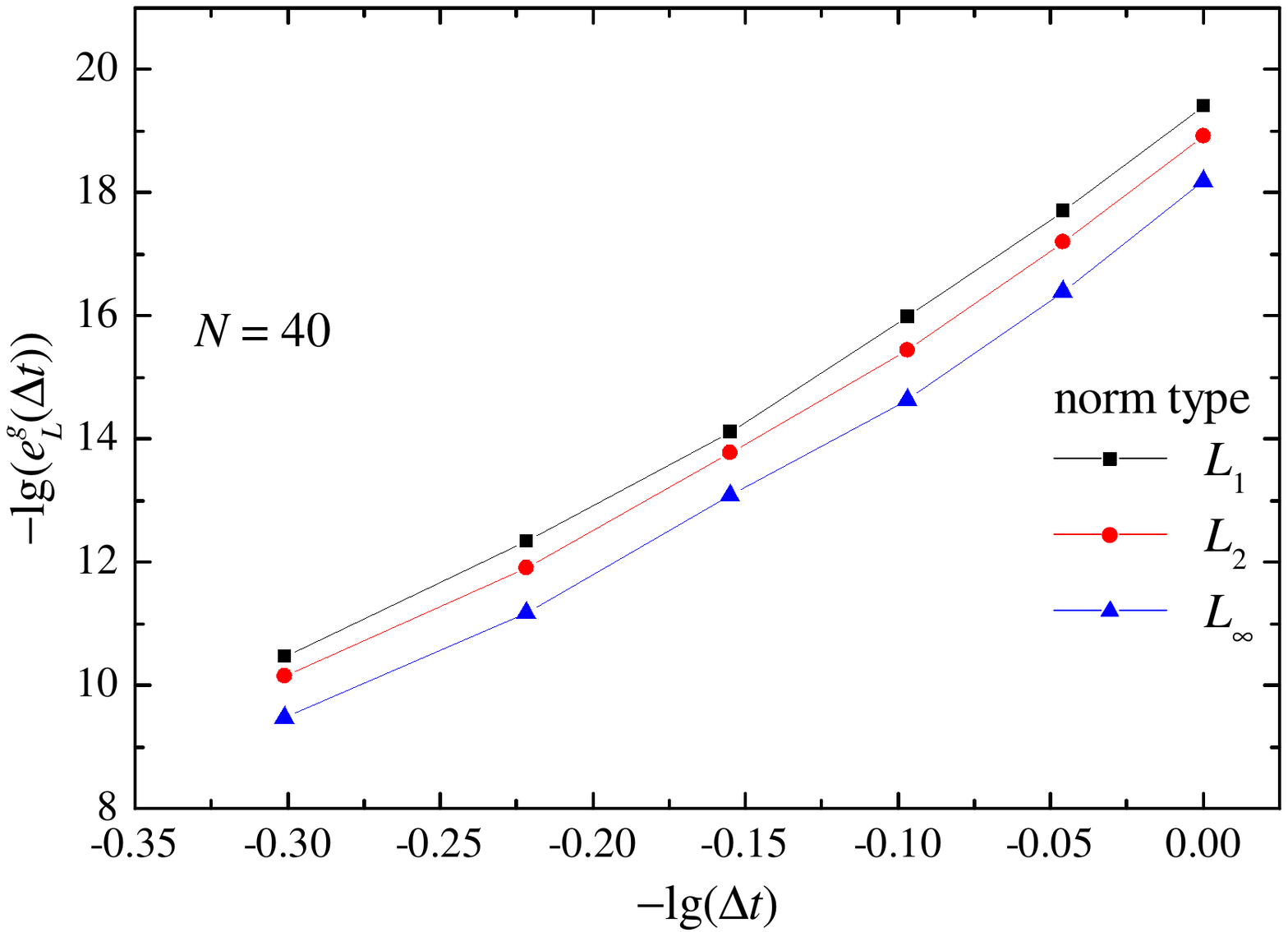}
\vspace{-8mm}\caption{\label{fig:dpend_ind3_errors:b3}}
\end{subfigure}\\[2mm]
\begin{subfigure}{0.320\textwidth}
\includegraphics[width=\textwidth]{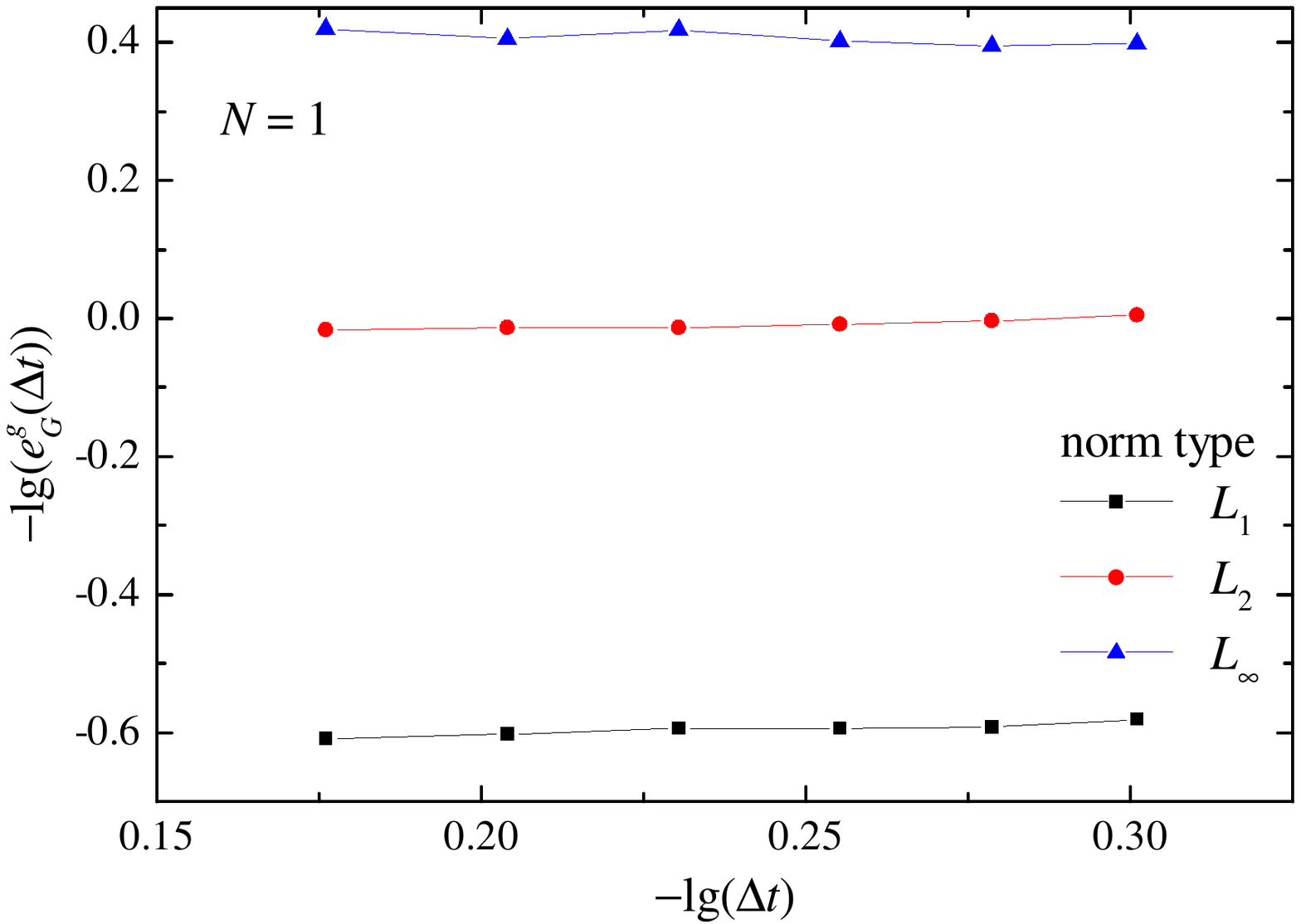}
\vspace{-8mm}\caption{\label{fig:dpend_ind3_errors:c1}}
\end{subfigure}
\begin{subfigure}{0.320\textwidth}
\includegraphics[width=\textwidth]{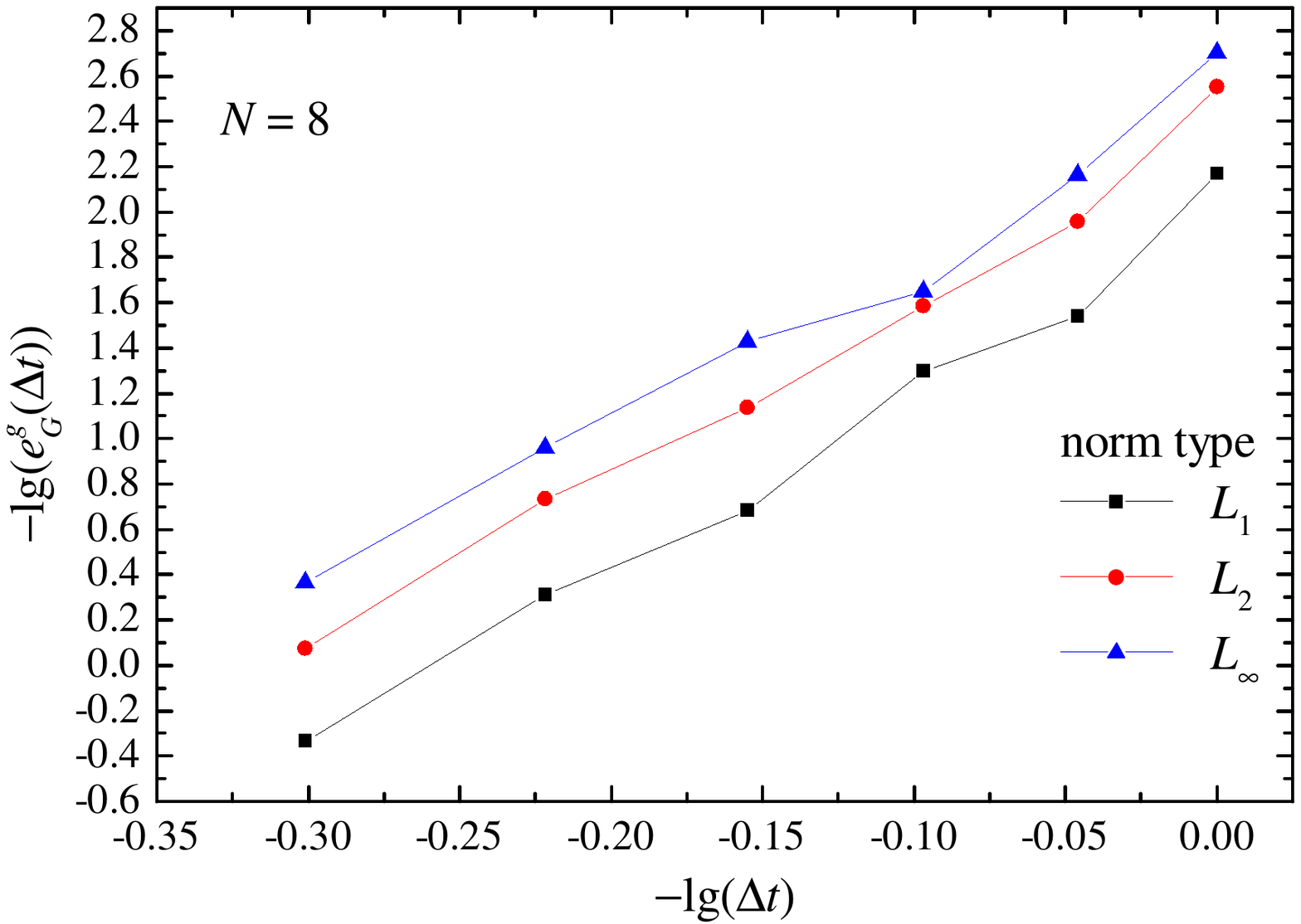}
\vspace{-8mm}\caption{\label{fig:dpend_ind3_errors:c2}}
\end{subfigure}
\begin{subfigure}{0.320\textwidth}
\includegraphics[width=\textwidth]{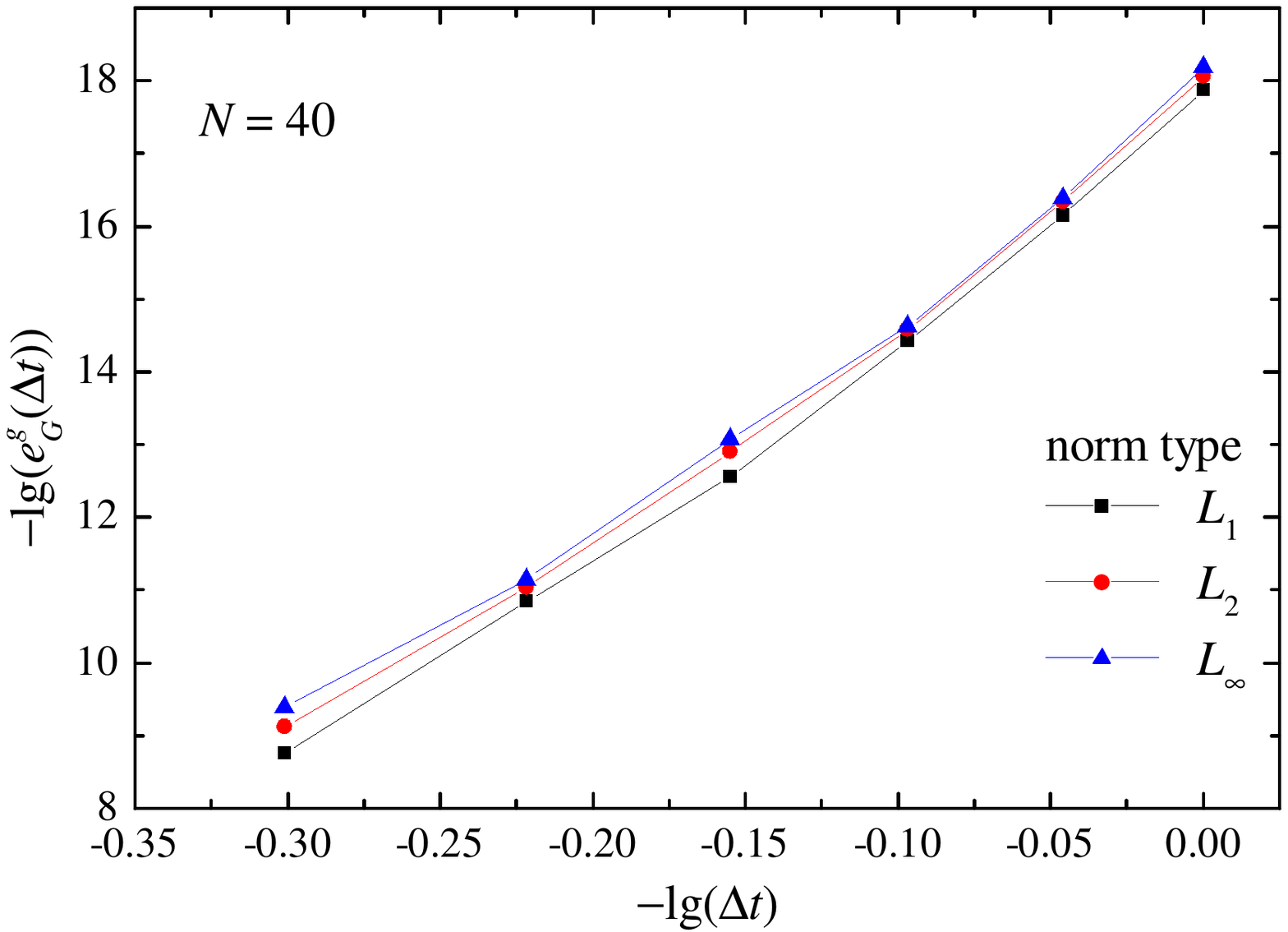}
\vspace{-8mm}\caption{\label{fig:dpend_ind3_errors:c3}}
\end{subfigure}\\[2mm]
\caption{%
Log-log plot of the dependence of the error for the solution at nodes at final time $t_{f}$ $e_{f}^{u}$ (\subref{fig:dpend_ind3_errors:a1}, \subref{fig:dpend_ind3_errors:a2}, \subref{fig:dpend_ind3_errors:a3}), global errors for the local solution $e_{L}^{g}$ (\subref{fig:dpend_ind3_errors:b1}, \subref{fig:dpend_ind3_errors:b2}, \subref{fig:dpend_ind3_errors:b3}) and the solution at nodes $e_{G}^{u}$ (\subref{fig:dpend_ind3_errors:c1}, \subref{fig:dpend_ind3_errors:c2}, \subref{fig:dpend_ind3_errors:c3}), by numerical solution of the DAE system (\ref{eq:math_dpend_dae_ind_3}) of index 3 obtained using polynomials with degrees $N = 1$ (\subref{fig:dpend_ind3_errors:a1}, \subref{fig:dpend_ind3_errors:b1}, \subref{fig:dpend_ind3_errors:c1}), $N = 8$ (\subref{fig:dpend_ind3_errors:a2}, \subref{fig:dpend_ind3_errors:b2}, \subref{fig:dpend_ind3_errors:c2}) and $N = 40$ (\subref{fig:dpend_ind3_errors:a3}, \subref{fig:dpend_ind3_errors:b3}, \subref{fig:dpend_ind3_errors:c3}).
}
\label{fig:dpend_ind3_errors}
\end{figure} 

\begin{table*}[h!]
\centering
\caption{%
Convergence orders $p_{f}$, $p_{L_{1}}$, $p_{L_{2}}$, $p_{L_{\infty}}$, calculated on final time step $t_{f}$ and in norms $L_{1}$, $L_{2}$, $L_{\infty}$ of the ADER-DG method for the DAE system (\ref{eq:math_dpend_dae_ind_3}) of index 3; $N$ is the degree of the basis polynomials $\varphi_{p}$. Orders $p^{n, u}$ are calculated for solution $\mathbf{u}_{n}$; orders $p^{n, g}$ and $p^{l, g}$ --- for the conditions $\mathbf{g} = 0$ on the numerical solution at the nodes $(\mathbf{u}_{n}, \mathbf{v}_{n})$ and on the local solution $(\mathbf{u}_{L}, \mathbf{v}_{L})$. The theoretical values of convergence order $p_{\rm th.}^{n} = 2N+1$ and $p_{\rm th.}^{l} = N+1$ are applicable for the ADER-DG method for ODE problems, and are presented for comparison.
}
\label{tab:conv_orders_dpend_ind3}
\begin{tabular}{@{}|l|l|lll|c|lll|c|@{}}
\toprule
$N$ & $p_{f}^{n, u}$ & $p_{L_{1}}^{n, g}$ & $p_{L_{2}}^{n, g}$ & $p_{L_{\infty}}^{n, g}$ & $p_{\rm th.}^{n}$ & $p_{L_{1}}^{l, g}$ & $p_{L_{2}}^{l, g}$ & $p_{L_{\infty}}^{l, g}$ & $p_{\rm th.}^{l}$ \\
\midrule
$1$	&	$2.07$	&	$0.20$	&	$0.16$	&	$-0.17$	&	$3$	&	$0.33$	&	$0.26$	&	$0.34$	&	$2$\\
$2$	&	$4.70$	&	$2.27$	&	$2.33$	&	$2.21$	&	$5$	&	$2.22$	&	$2.19$	&	$2.10$	&	$3$\\
$3$	&	$4.56$	&	$2.69$	&	$2.65$	&	$2.86$	&	$7$	&	$2.79$	&	$2.78$	&	$2.90$	&	$4$\\
$4$	&	$4.69$	&	$3.89$	&	$3.77$	&	$4.04$	&	$9$	&	$3.74$	&	$3.70$	&	$3.74$	&	$5$\\
$5$	&	$7.11$	&	$4.81$	&	$4.97$	&	$4.82$	&	$11$	&	$5.25$	&	$5.07$	&	$4.92$	&	$6$\\
$6$	&	$8.99$	&	$6.14$	&	$6.06$	&	$5.71$	&	$13$	&	$6.04$	&	$5.68$	&	$5.01$	&	$7$\\
$7$	&	$6.23$	&	$6.43$	&	$6.37$	&	$6.35$	&	$15$	&	$6.68$	&	$6.38$	&	$6.17$	&	$8$\\
$8$	&	$9.12$	&	$7.99$	&	$7.83$	&	$7.36$	&	$17$	&	$7.88$	&	$7.54$	&	$7.05$	&	$9$\\
$9$	&	$9.69$	&	$8.04$	&	$7.86$	&	$7.72$	&	$19$	&	$8.13$	&	$7.73$	&	$7.51$	&	$10$\\
$10$	&	$10.63$	&	$9.45$	&	$9.24$	&	$8.73$	&	$21$	&	$9.52$	&	$9.05$	&	$8.47$	&	$11$\\
$11$	&	$11.52$	&	$9.81$	&	$9.49$	&	$9.00$	&	$23$	&	$9.89$	&	$9.27$	&	$8.85$	&	$12$\\
$12$	&	$10.49$	&	$10.33$	&	$10.02$	&	$9.38$	&	$25$	&	$10.43$	&	$9.82$	&	$9.14$	&	$13$\\
$13$	&	$11.82$	&	$11.43$	&	$11.06$	&	$10.29$	&	$27$	&	$11.59$	&	$10.83$	&	$10.13$	&	$14$\\
$14$	&	$11.69$	&	$10.96$	&	$10.76$	&	$10.40$	&	$29$	&	$11.09$	&	$10.71$	&	$10.33$	&	$15$\\
$15$	&	$12.49$	&	$13.01$	&	$12.48$	&	$11.61$	&	$31$	&	$13.16$	&	$12.30$	&	$11.49$	&	$16$\\
$16$	&	$14.50$	&	$12.66$	&	$12.28$	&	$11.85$	&	$33$	&	$12.65$	&	$11.99$	&	$11.79$	&	$17$\\
$17$	&	$15.20$	&	$14.10$	&	$13.73$	&	$13.01$	&	$35$	&	$14.26$	&	$13.57$	&	$12.93$	&	$18$\\
$18$	&	$17.00$	&	$14.71$	&	$14.15$	&	$13.47$	&	$37$	&	$14.59$	&	$13.78$	&	$13.21$	&	$19$\\
$19$	&	$17.11$	&	$15.21$	&	$14.71$	&	$14.17$	&	$39$	&	$15.38$	&	$14.57$	&	$14.12$	&	$20$\\
$20$	&	$18.19$	&	$16.30$	&	$15.85$	&	$15.11$	&	$41$	&	$16.18$	&	$15.54$	&	$15.03$	&	$21$\\
\midrule
$25$	&	$21.88$	&	$20.08$	&	$19.68$	&	$19.28$	&	$51$	&	$20.13$	&	$19.54$	&	$19.14$	&	$26$\\
$30$	&	$24.31$	&	$22.43$	&	$21.95$	&	$21.44$	&	$61$	&	$22.36$	&	$21.69$	&	$21.48$	&	$31$\\
$35$	&	$29.38$	&	$26.71$	&	$26.32$	&	$25.84$	&	$71$	&	$26.59$	&	$26.14$	&	$25.81$	&	$36$\\
$40$	&	$32.65$	&	$30.05$	&	$29.58$	&	$29.01$	&	$81$	&	$29.71$	&	$29.11$	&	$28.77$	&	$41$\\
\bottomrule
\end{tabular}
\end{table*}

\begin{figure}[h!]
\captionsetup[subfigure]{%
	position=bottom,
	font+=smaller,
	textfont=normalfont,
	singlelinecheck=off,
	justification=raggedright
}
\centering
\begin{subfigure}{0.240\textwidth}
\includegraphics[width=\textwidth]{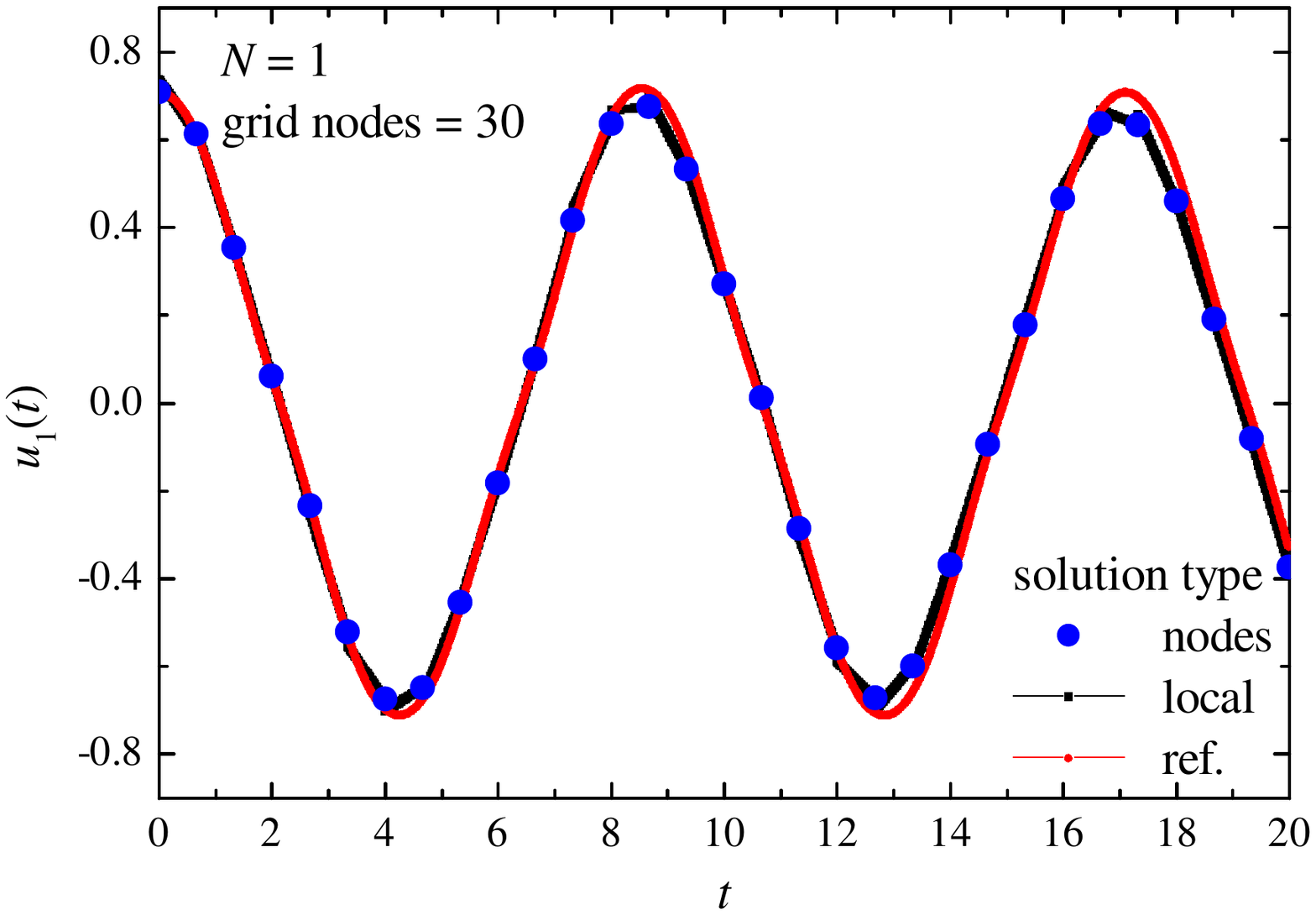}
\vspace{-8mm}\caption{\label{fig:dpend_ind2_sols_u:a1}}
\end{subfigure}
\begin{subfigure}{0.240\textwidth}
\includegraphics[width=\textwidth]{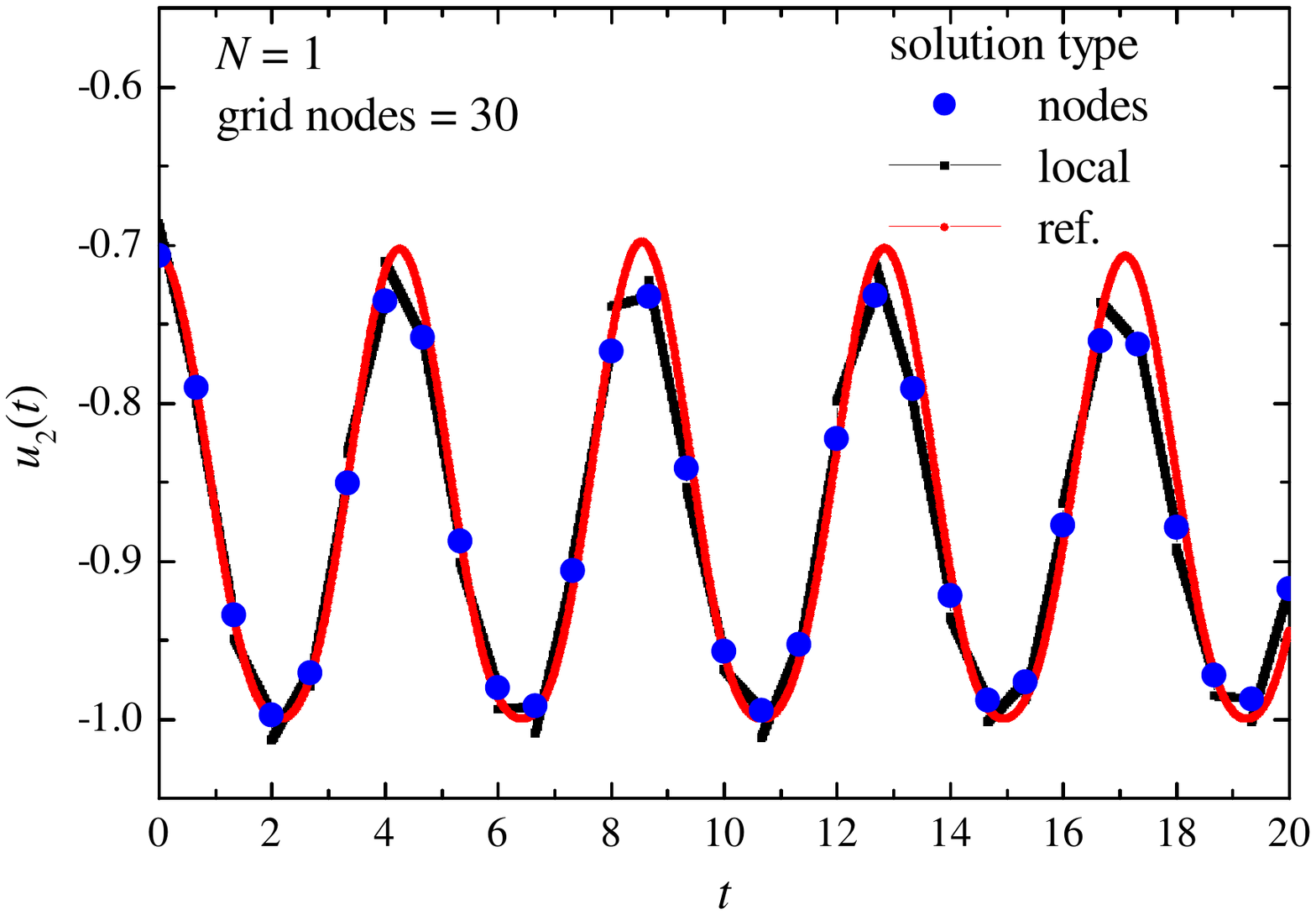}
\vspace{-8mm}\caption{\label{fig:dpend_ind2_sols_u:a2}}
\end{subfigure}
\begin{subfigure}{0.240\textwidth}
\includegraphics[width=\textwidth]{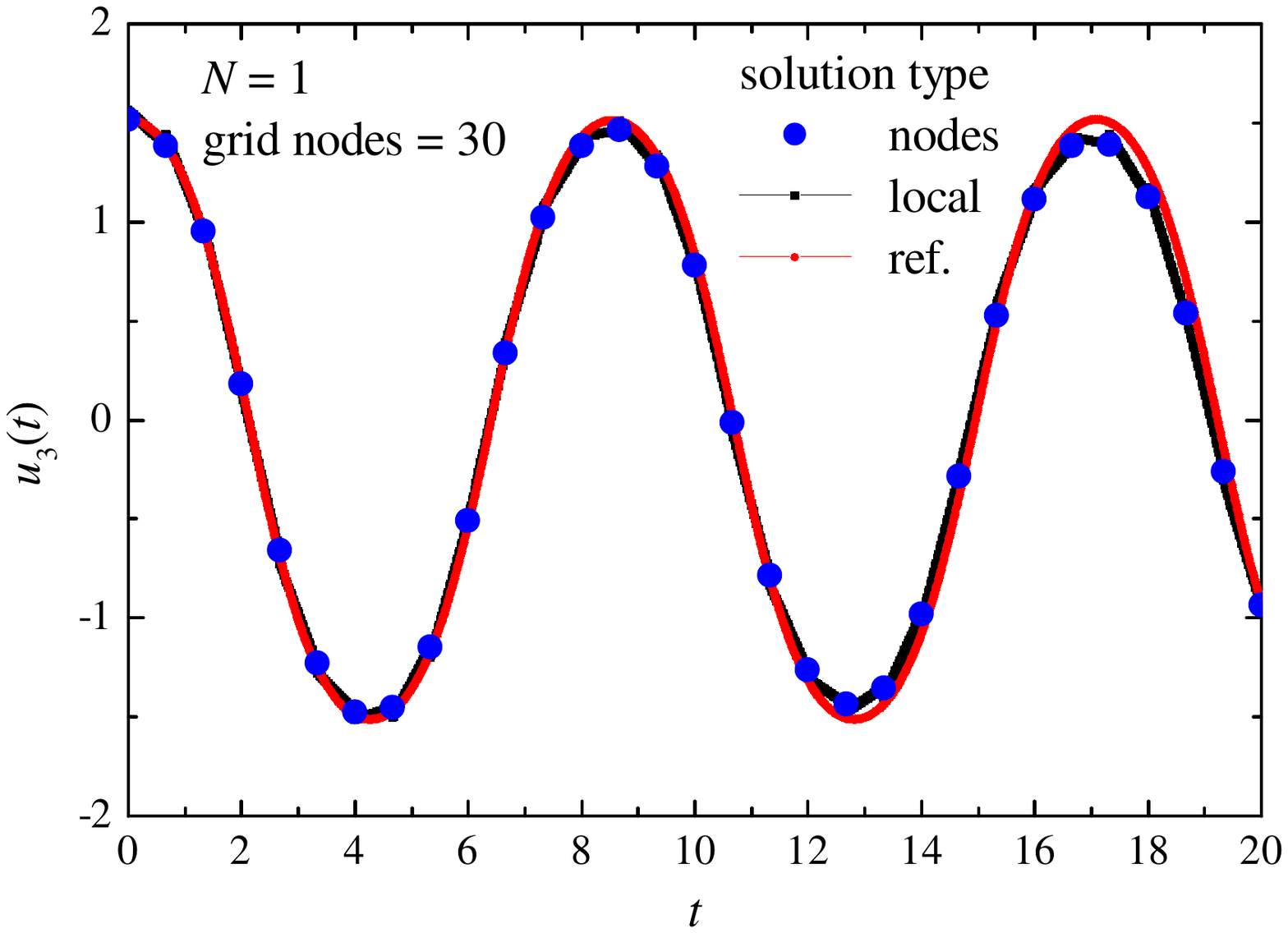}
\vspace{-8mm}\caption{\label{fig:dpend_ind2_sols_u:a3}}
\end{subfigure}
\begin{subfigure}{0.240\textwidth}
\includegraphics[width=\textwidth]{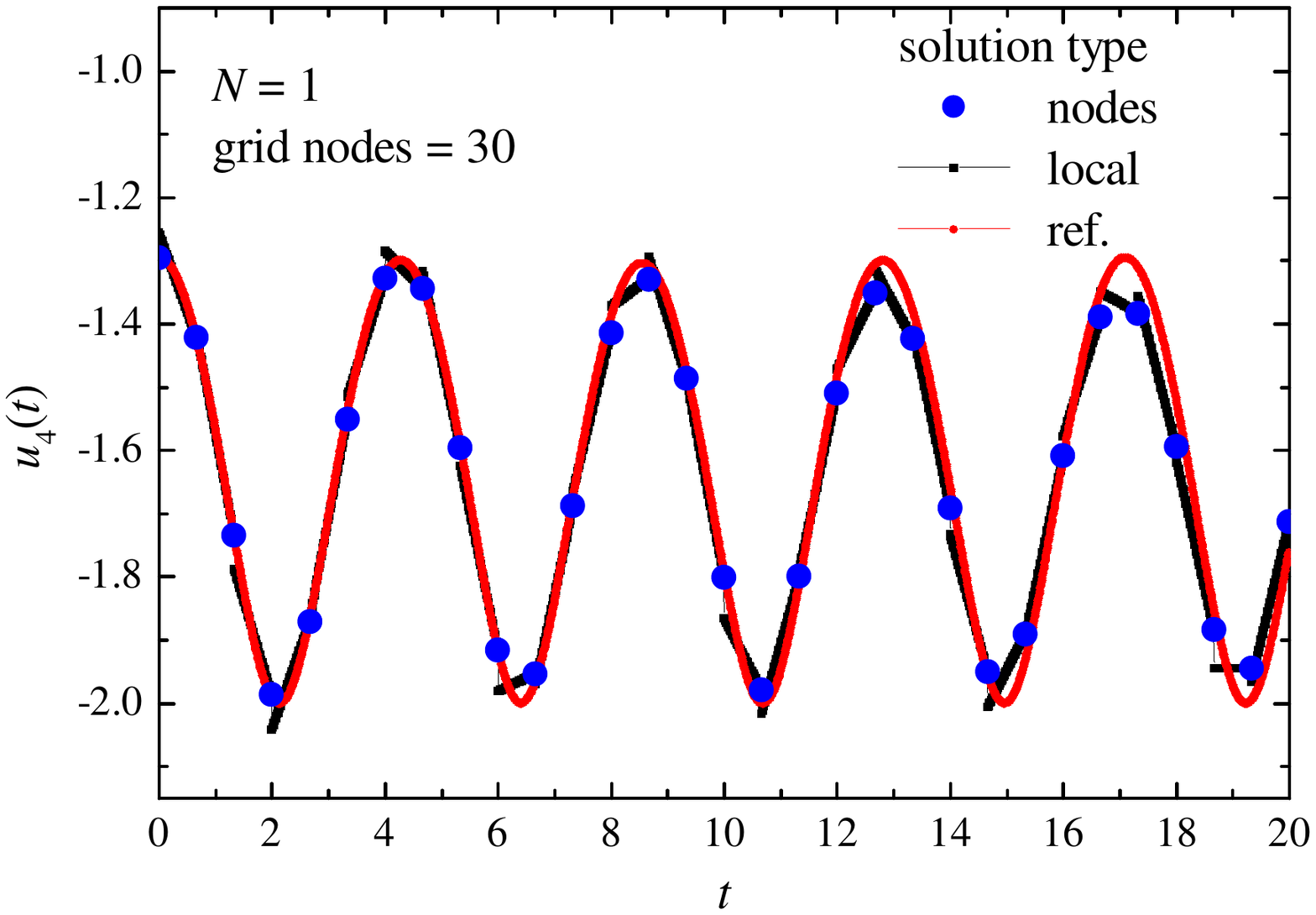}
\vspace{-8mm}\caption{\label{fig:dpend_ind2_sols_u:a4}}
\end{subfigure}\\[2mm]
\begin{subfigure}{0.240\textwidth}
\includegraphics[width=\textwidth]{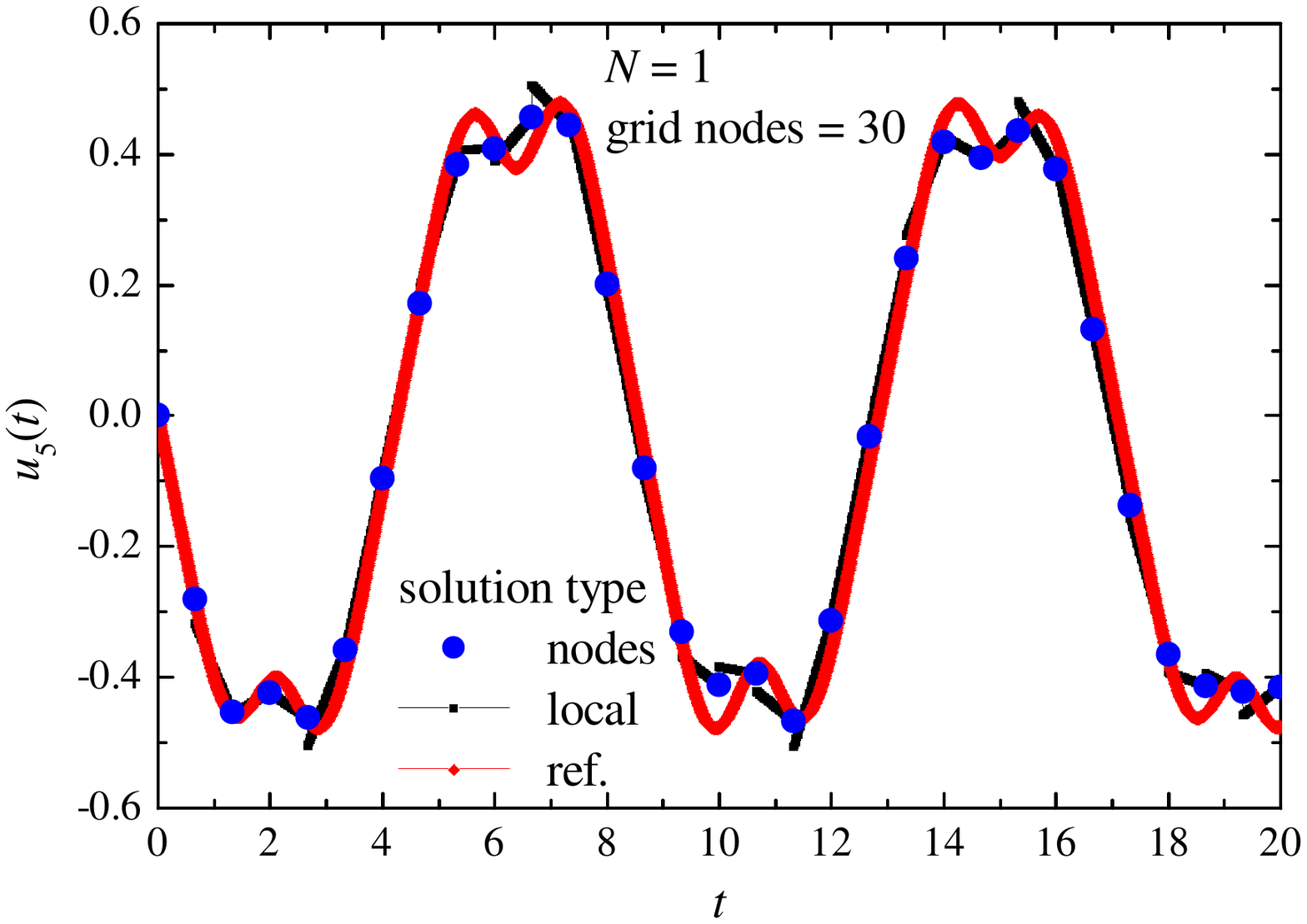}
\vspace{-8mm}\caption{\label{fig:dpend_ind2_sols_u:b1}}
\end{subfigure}
\begin{subfigure}{0.240\textwidth}
\includegraphics[width=\textwidth]{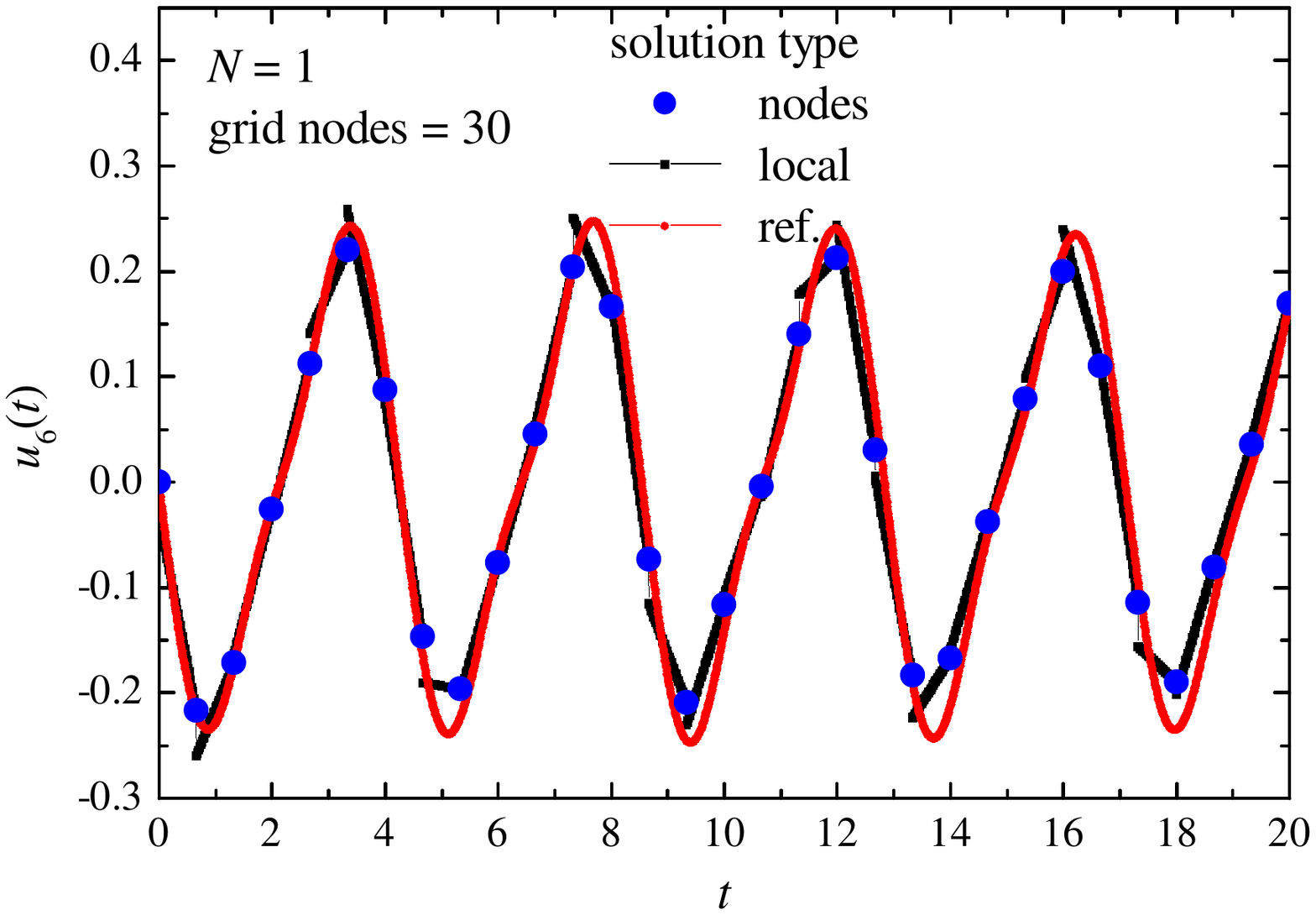}
\vspace{-8mm}\caption{\label{fig:dpend_ind2_sols_u:b2}}
\end{subfigure}
\begin{subfigure}{0.240\textwidth}
\includegraphics[width=\textwidth]{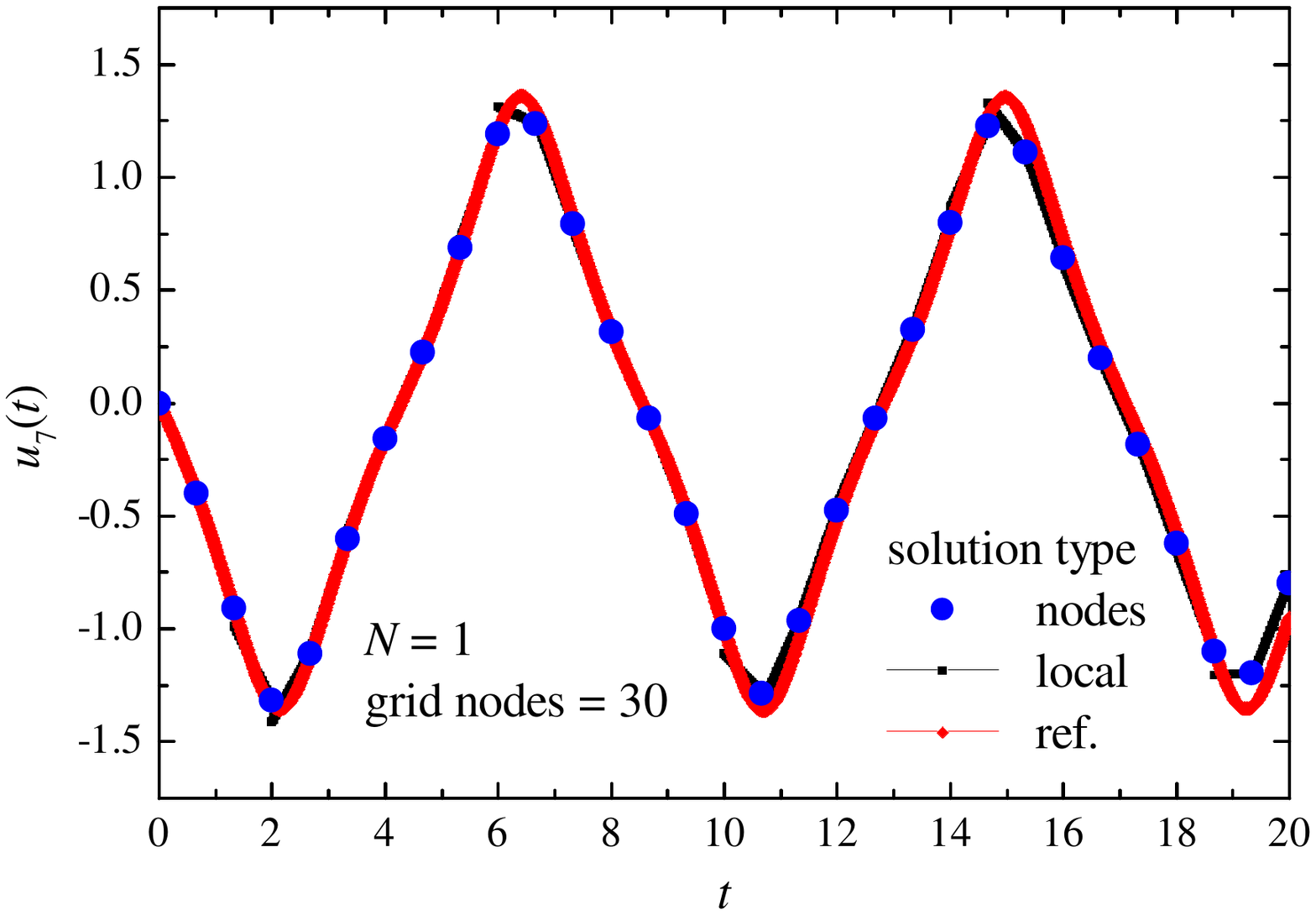}
\vspace{-8mm}\caption{\label{fig:dpend_ind2_sols_u:b3}}
\end{subfigure}
\begin{subfigure}{0.240\textwidth}
\includegraphics[width=\textwidth]{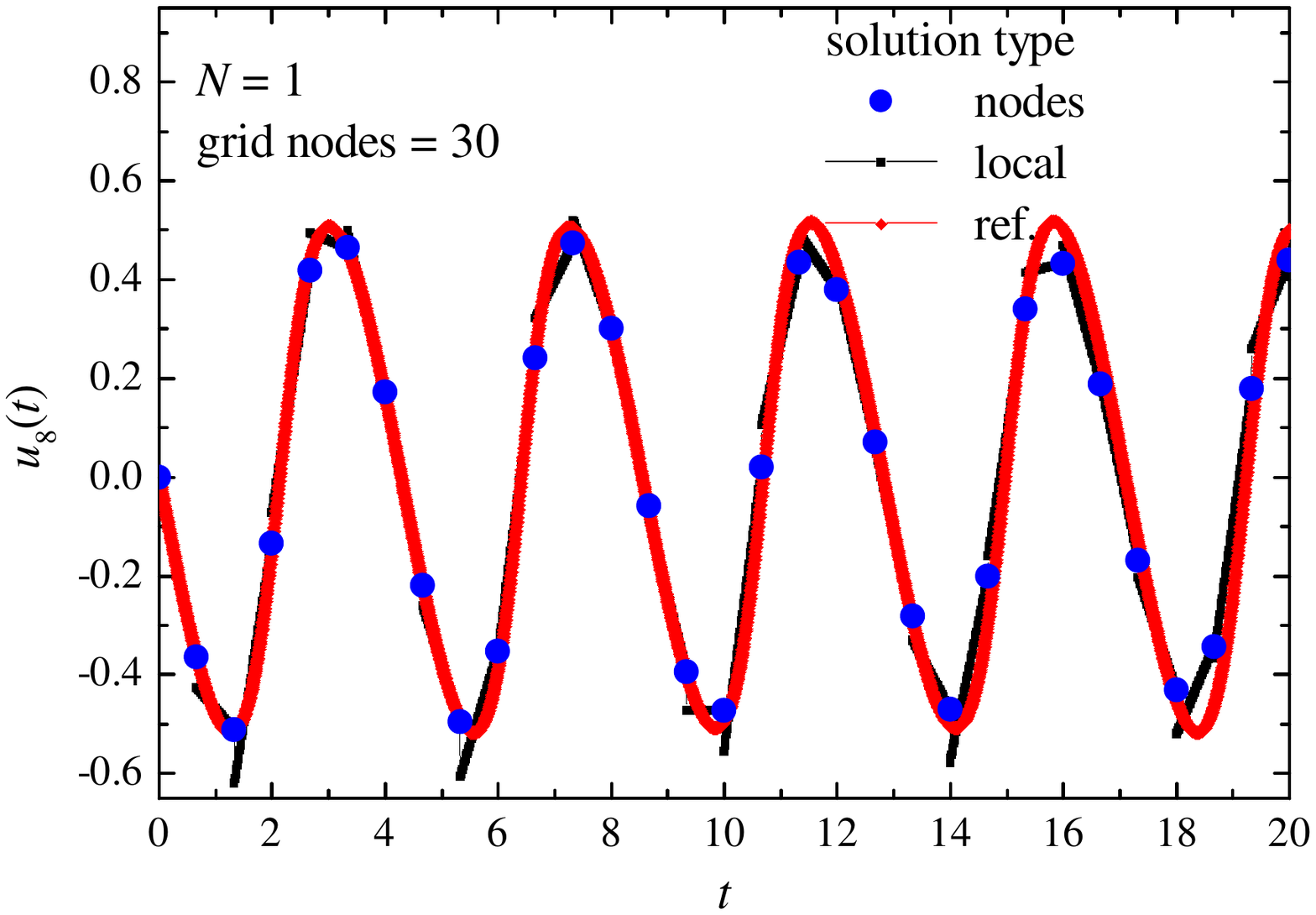}
\vspace{-8mm}\caption{\label{fig:dpend_ind2_sols_u:b4}}
\end{subfigure}\\[2mm]
\begin{subfigure}{0.240\textwidth}
\includegraphics[width=\textwidth]{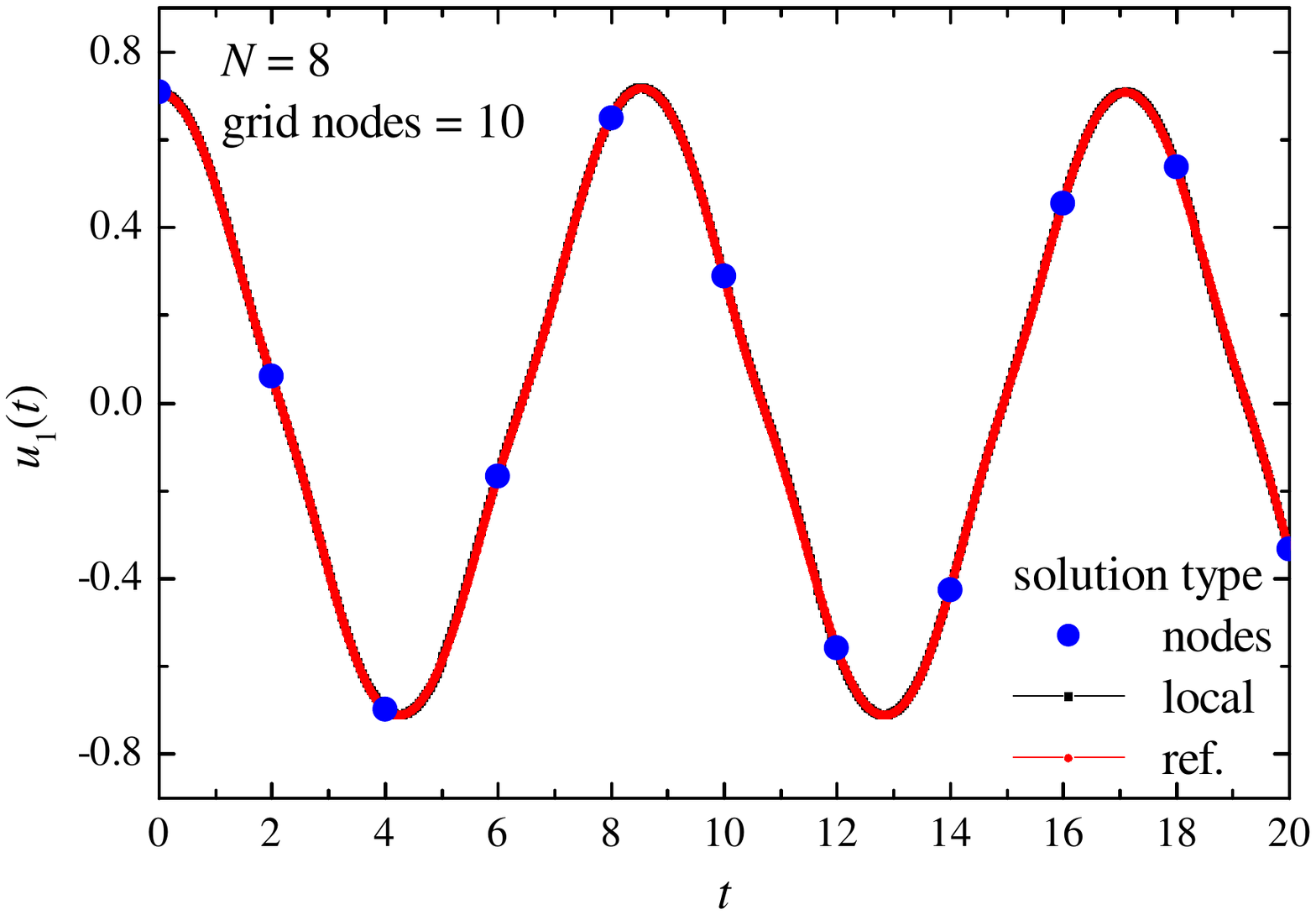}
\vspace{-8mm}\caption{\label{fig:dpend_ind2_sols_u:c1}}
\end{subfigure}
\begin{subfigure}{0.240\textwidth}
\includegraphics[width=\textwidth]{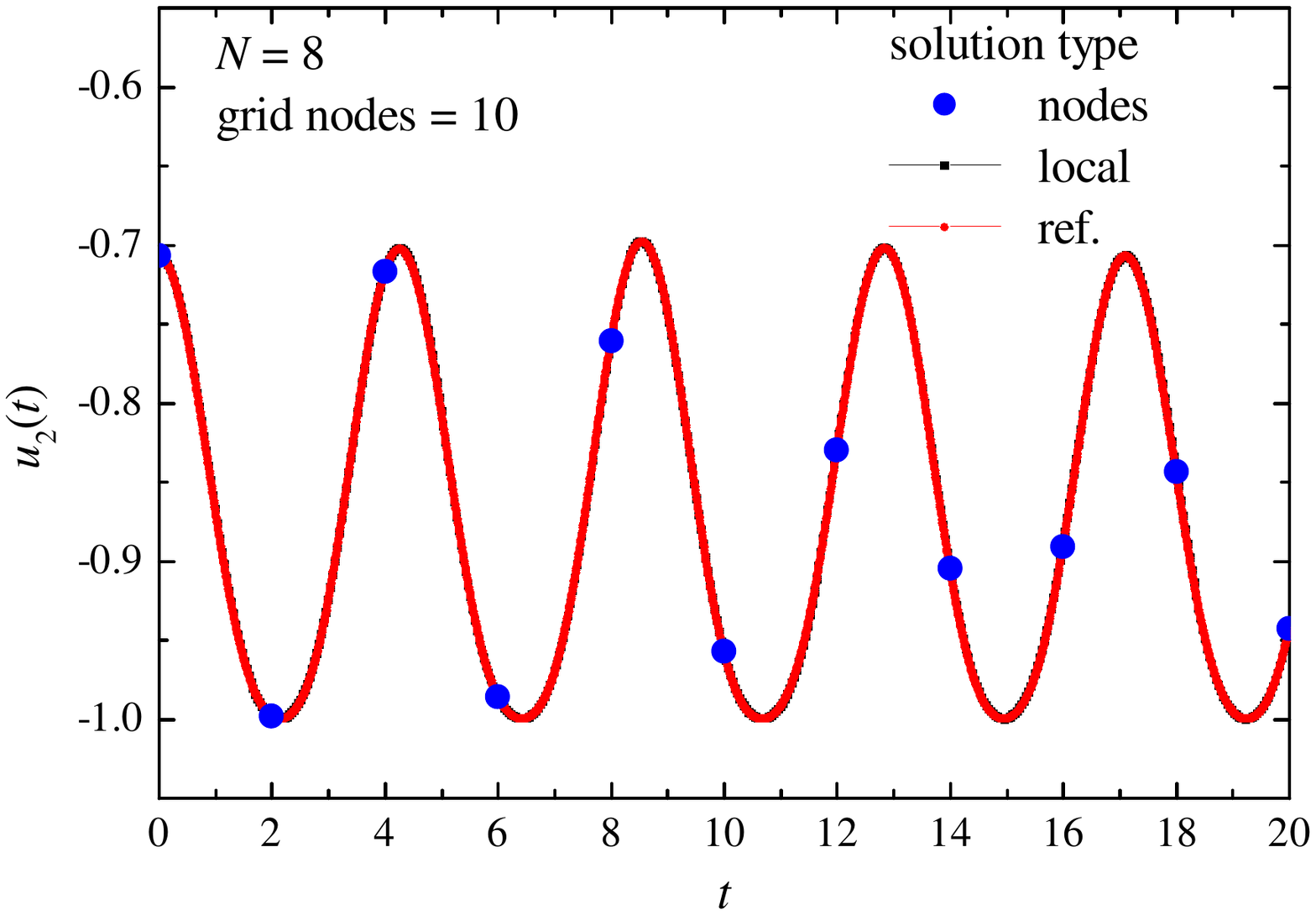}
\vspace{-8mm}\caption{\label{fig:dpend_ind2_sols_u:c2}}
\end{subfigure}
\begin{subfigure}{0.240\textwidth}
\includegraphics[width=\textwidth]{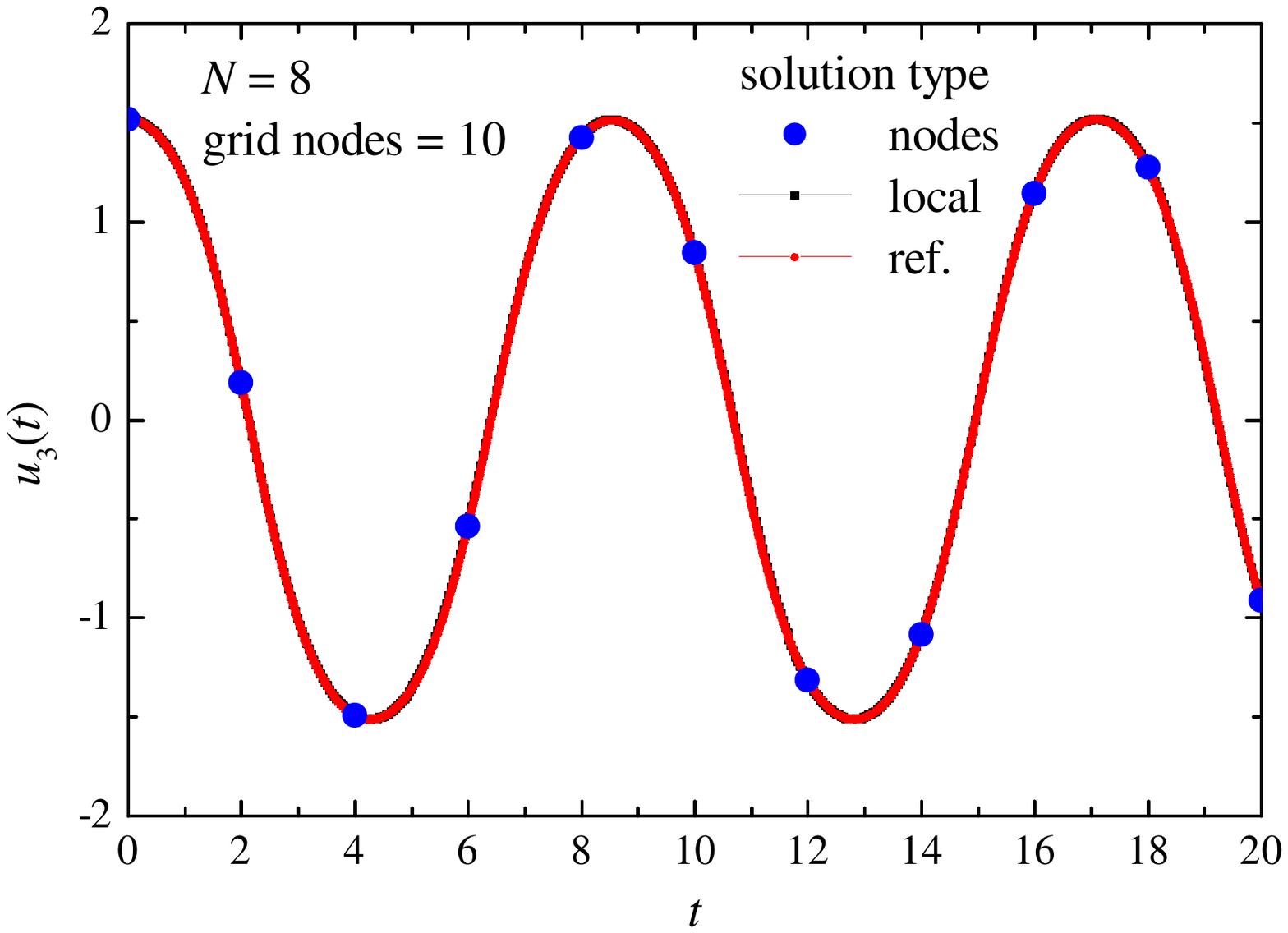}
\vspace{-8mm}\caption{\label{fig:dpend_ind2_sols_u:c3}}
\end{subfigure}
\begin{subfigure}{0.240\textwidth}
\includegraphics[width=\textwidth]{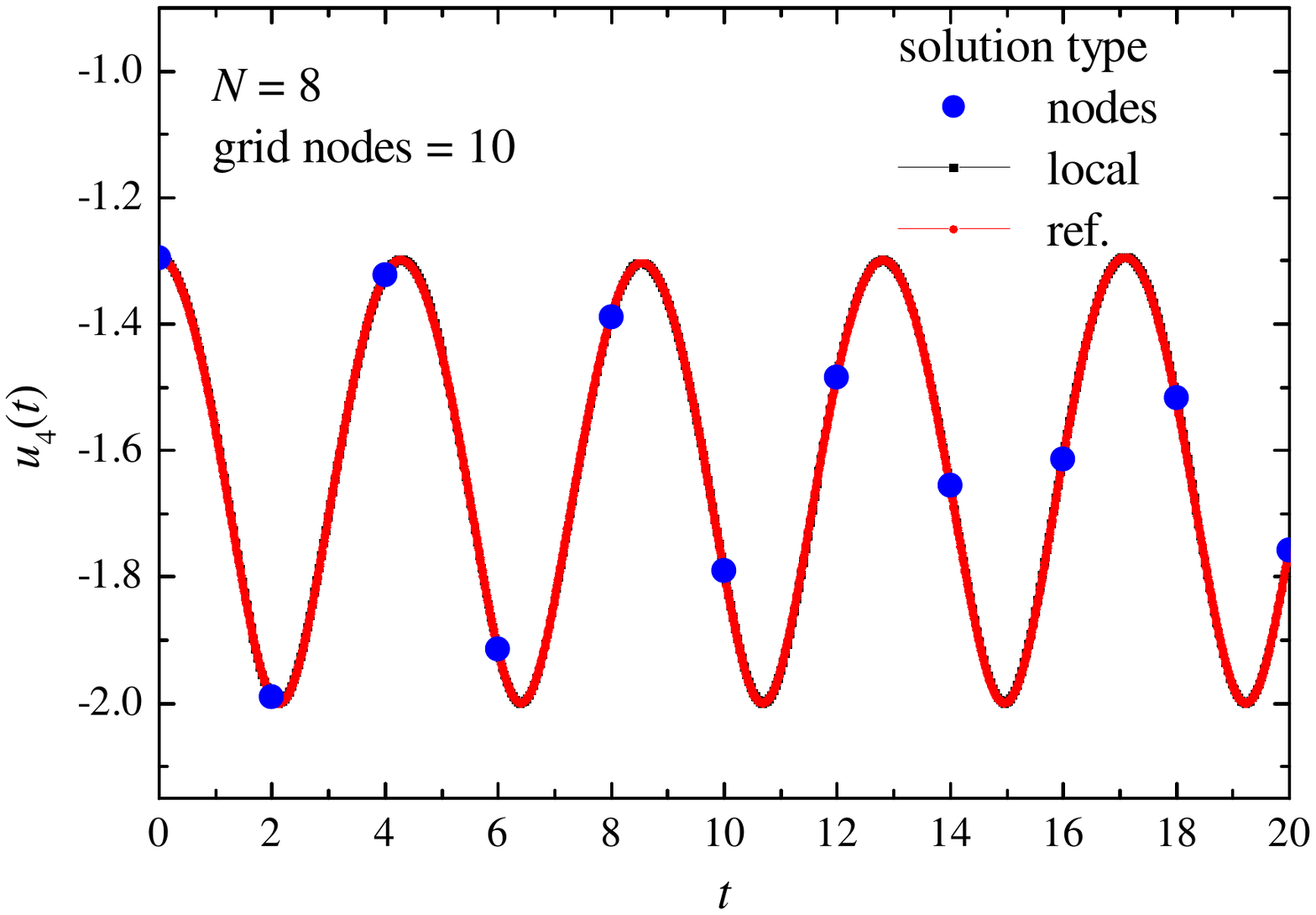}
\vspace{-8mm}\caption{\label{fig:dpend_ind2_sols_u:c4}}
\end{subfigure}\\[2mm]
\begin{subfigure}{0.240\textwidth}
\includegraphics[width=\textwidth]{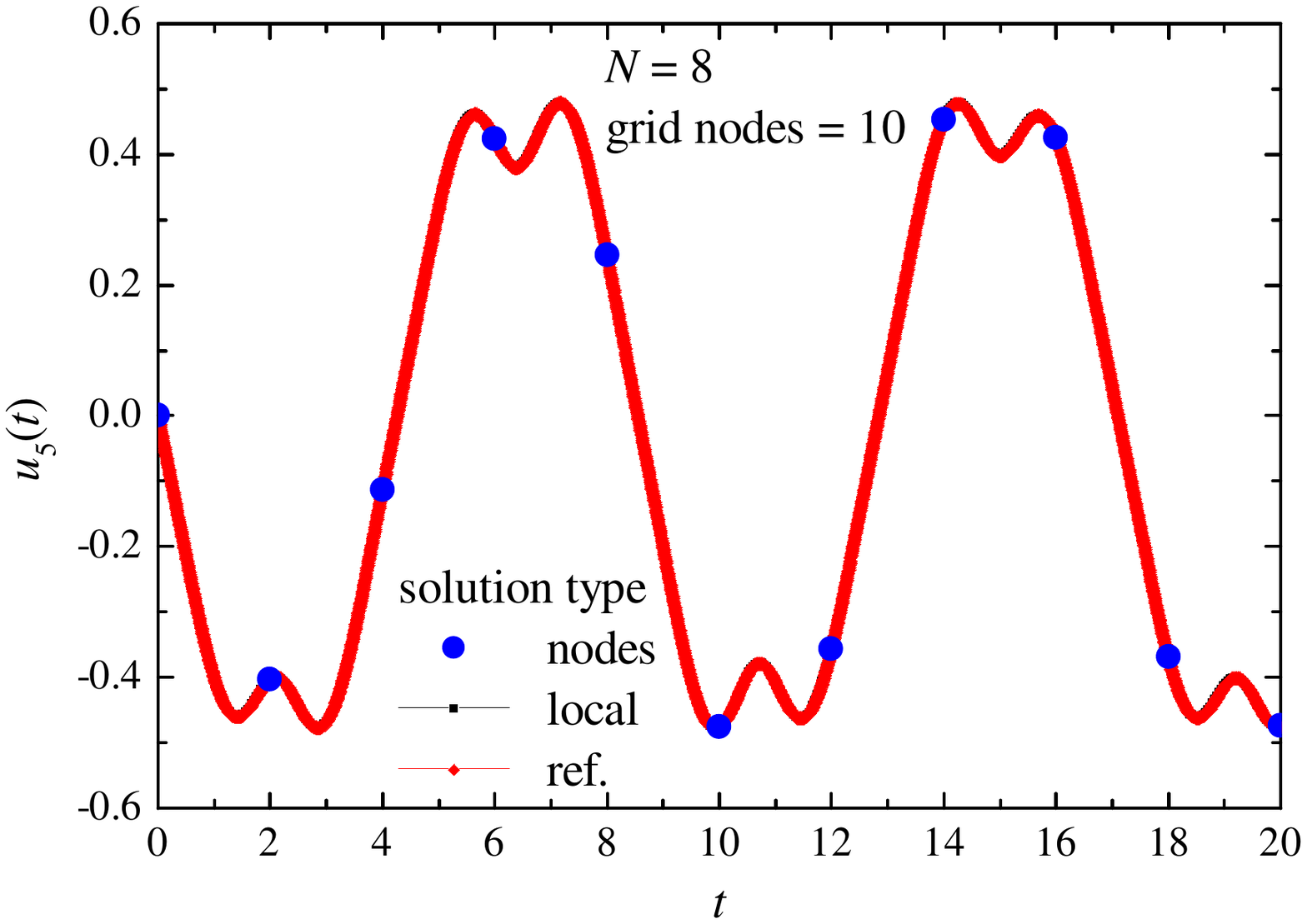}
\vspace{-8mm}\caption{\label{fig:dpend_ind2_sols_u:d1}}
\end{subfigure}
\begin{subfigure}{0.240\textwidth}
\includegraphics[width=\textwidth]{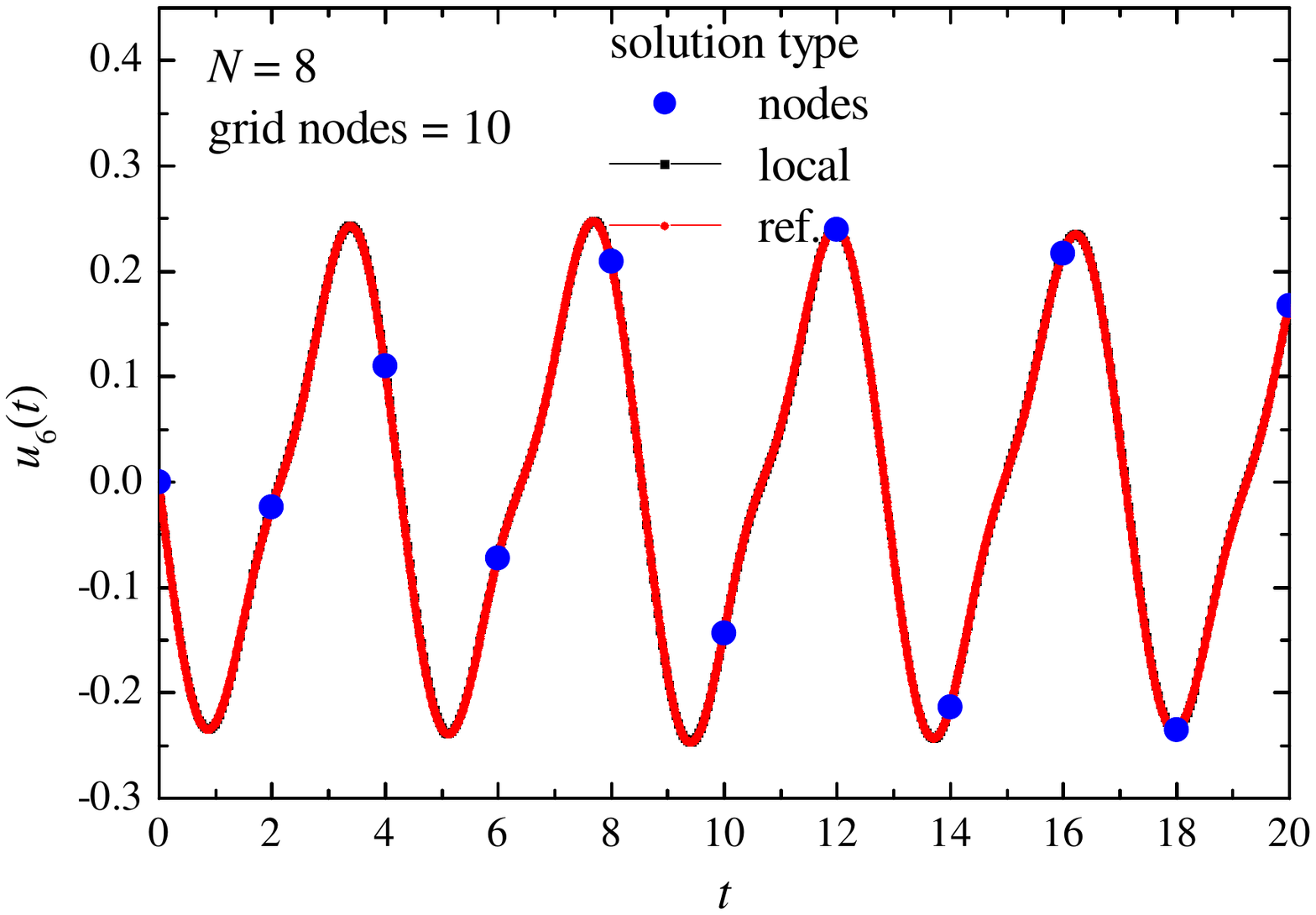}
\vspace{-8mm}\caption{\label{fig:dpend_ind2_sols_u:d2}}
\end{subfigure}
\begin{subfigure}{0.240\textwidth}
\includegraphics[width=\textwidth]{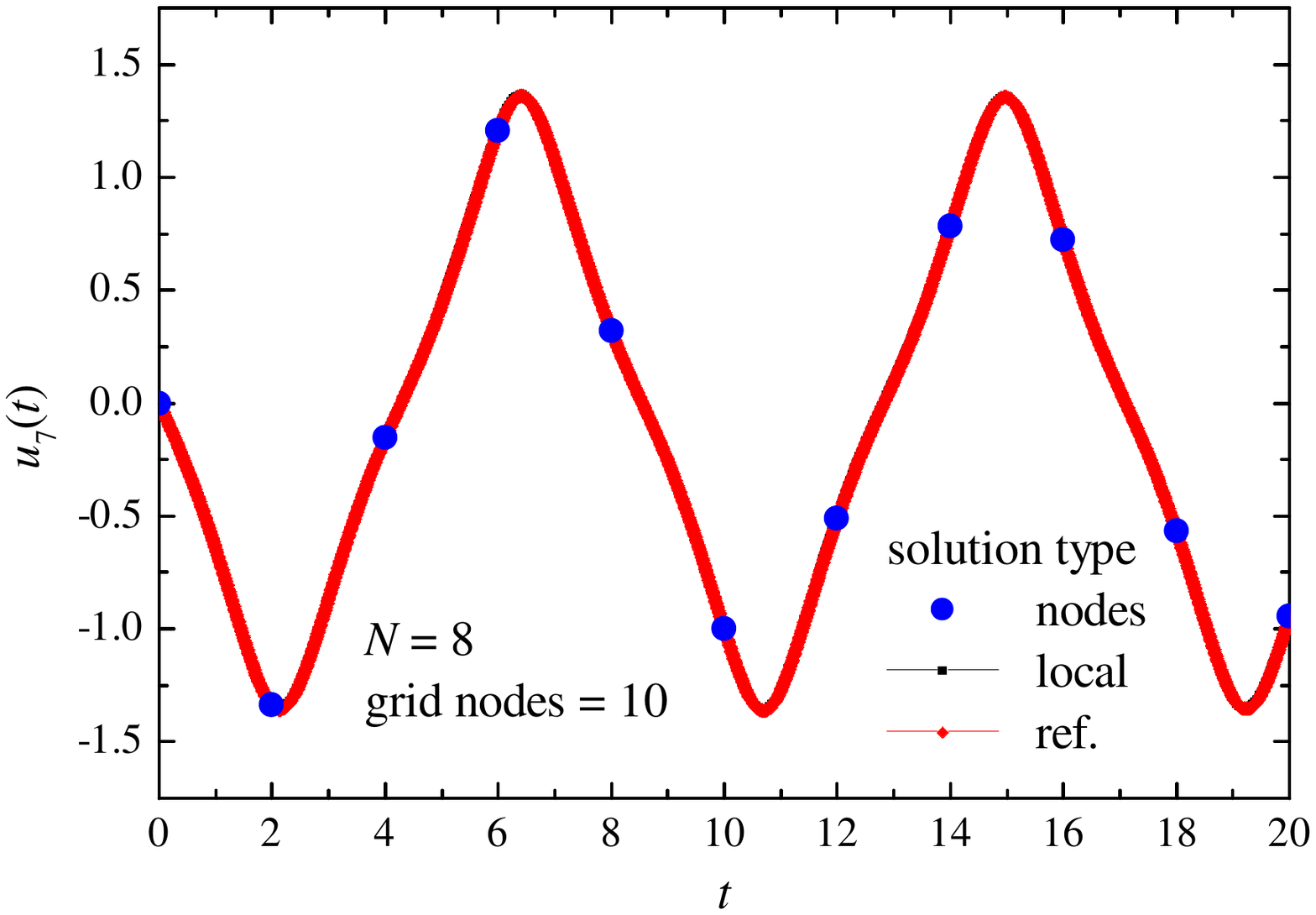}
\vspace{-8mm}\caption{\label{fig:dpend_ind2_sols_u:d3}}
\end{subfigure}
\begin{subfigure}{0.240\textwidth}
\includegraphics[width=\textwidth]{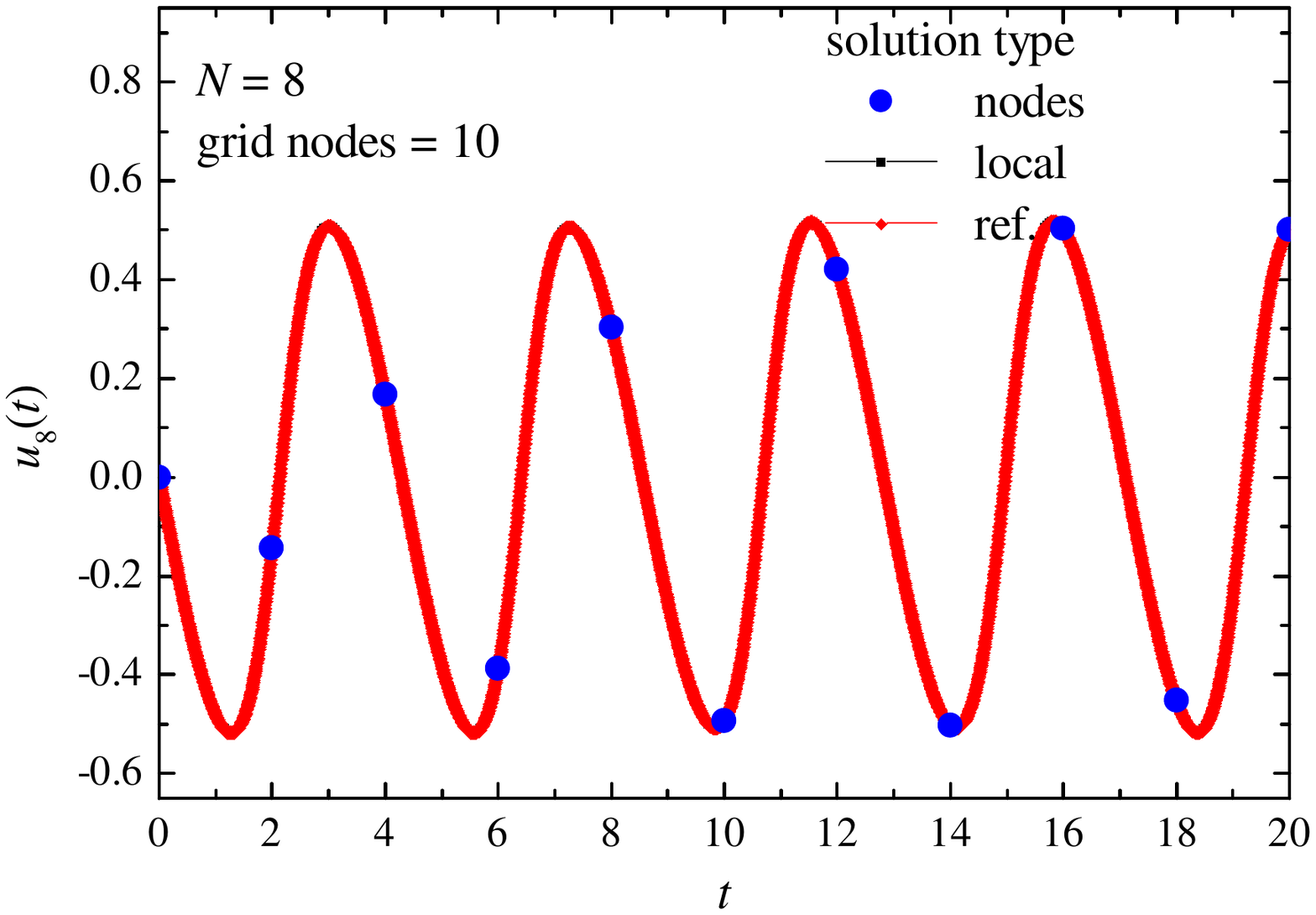}
\vspace{-8mm}\caption{\label{fig:dpend_ind2_sols_u:d4}}
\end{subfigure}\\[2mm]
\begin{subfigure}{0.240\textwidth}
\includegraphics[width=\textwidth]{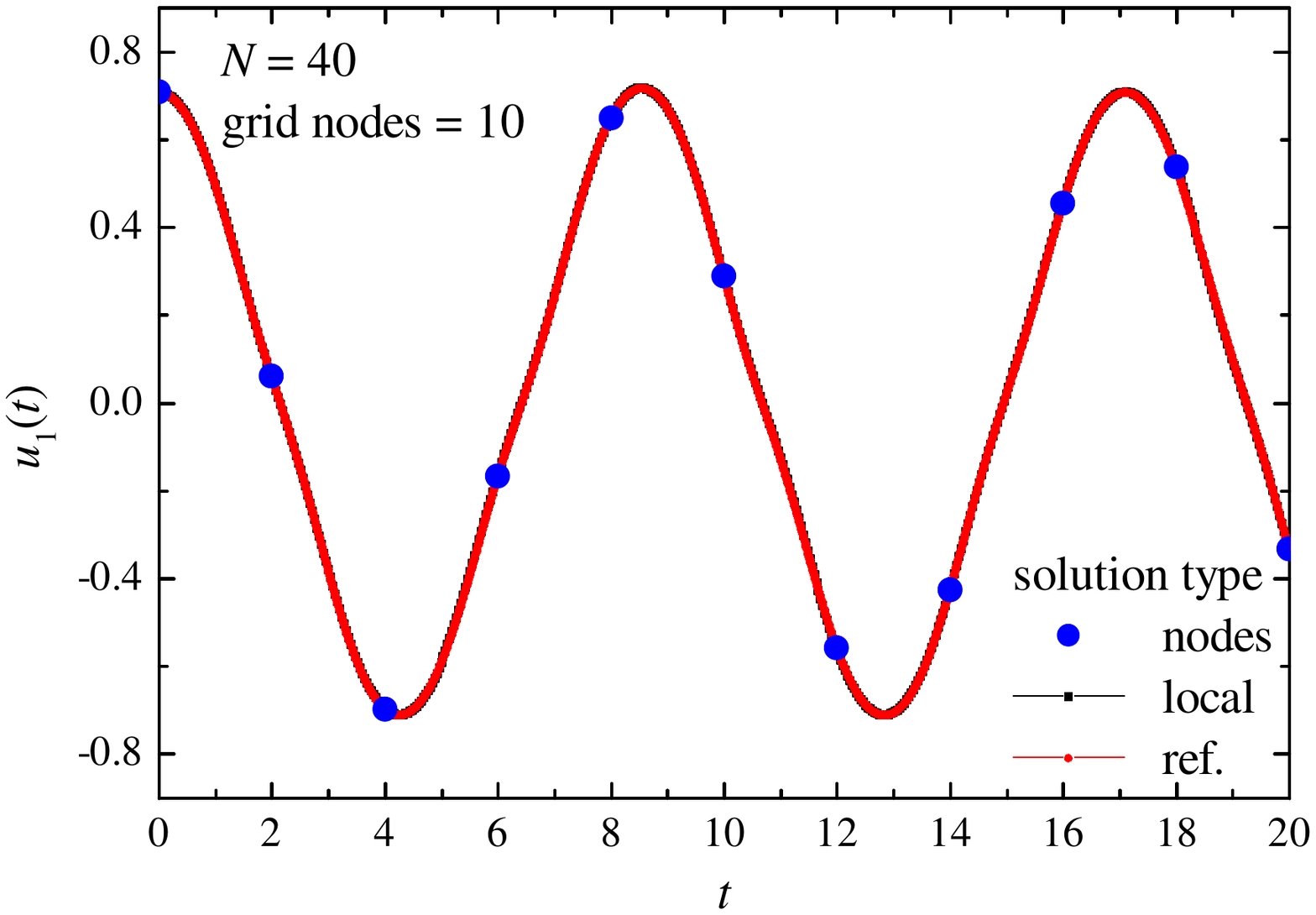}
\vspace{-8mm}\caption{\label{fig:dpend_ind2_sols_u:e1}}
\end{subfigure}
\begin{subfigure}{0.240\textwidth}
\includegraphics[width=\textwidth]{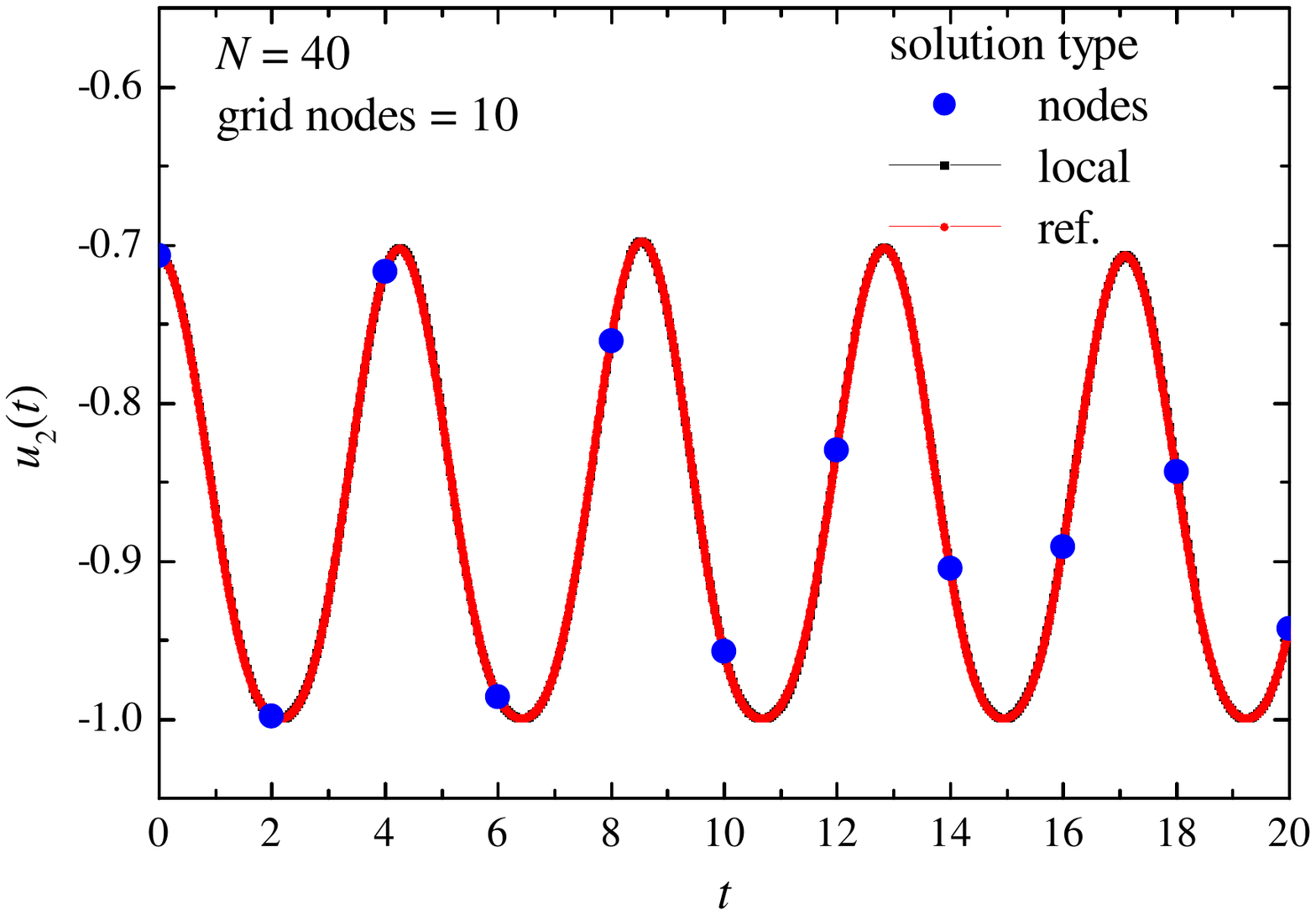}
\vspace{-8mm}\caption{\label{fig:dpend_ind2_sols_u:e2}}
\end{subfigure}
\begin{subfigure}{0.240\textwidth}
\includegraphics[width=\textwidth]{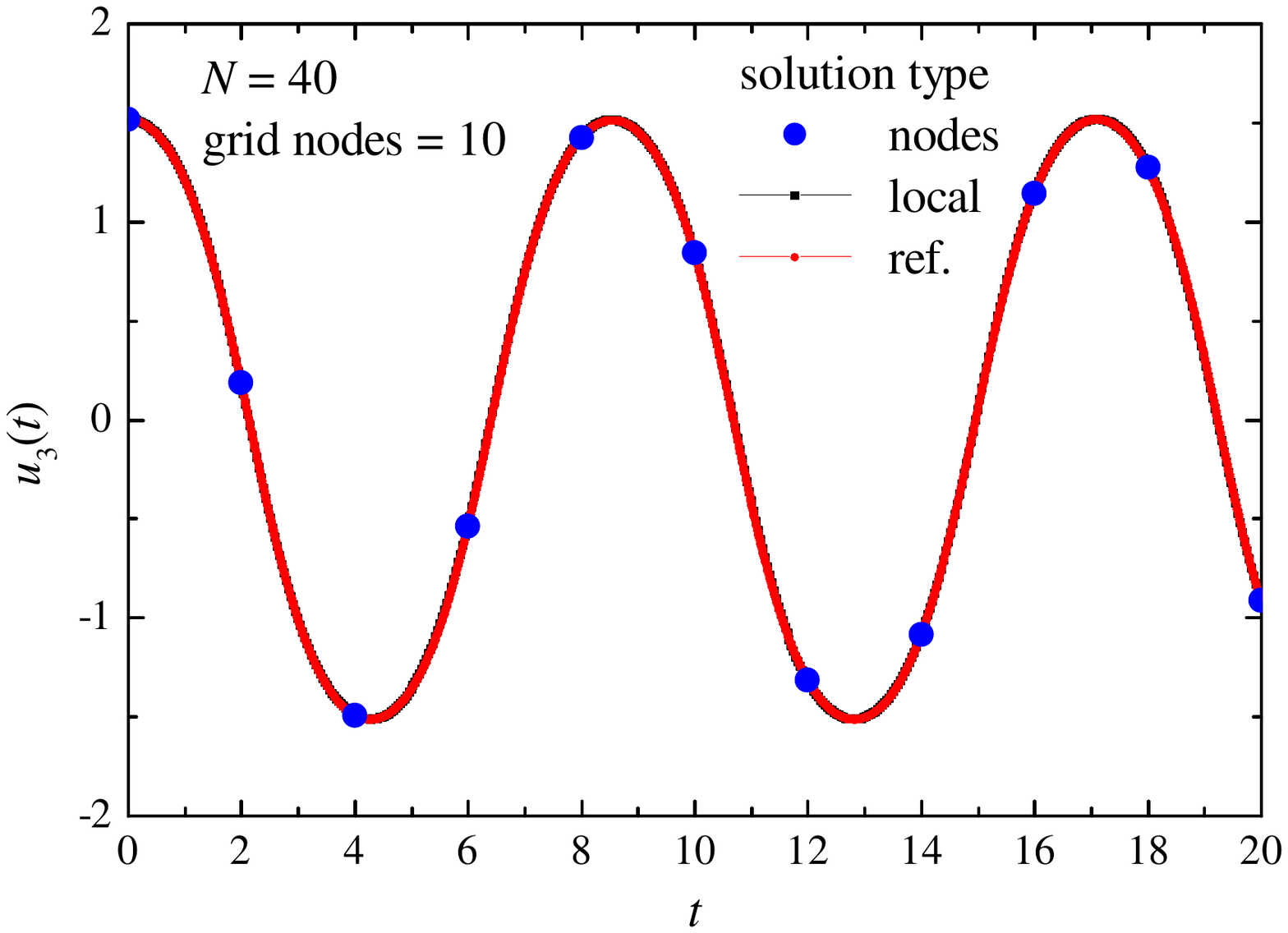}
\vspace{-8mm}\caption{\label{fig:dpend_ind2_sols_u:e3}}
\end{subfigure}
\begin{subfigure}{0.240\textwidth}
\includegraphics[width=\textwidth]{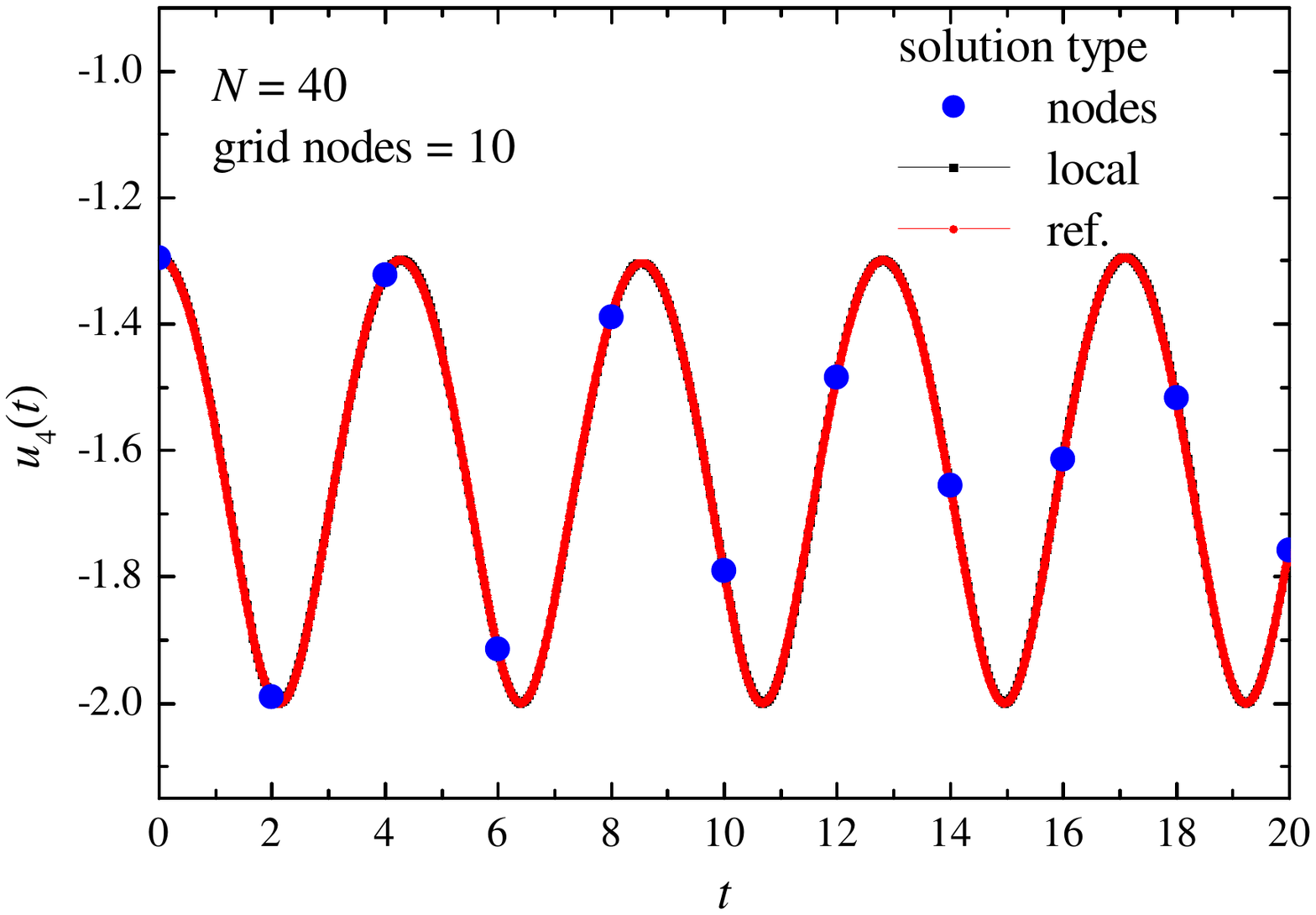}
\vspace{-8mm}\caption{\label{fig:dpend_ind2_sols_u:e4}}
\end{subfigure}\\[2mm]
\begin{subfigure}{0.240\textwidth}
\includegraphics[width=\textwidth]{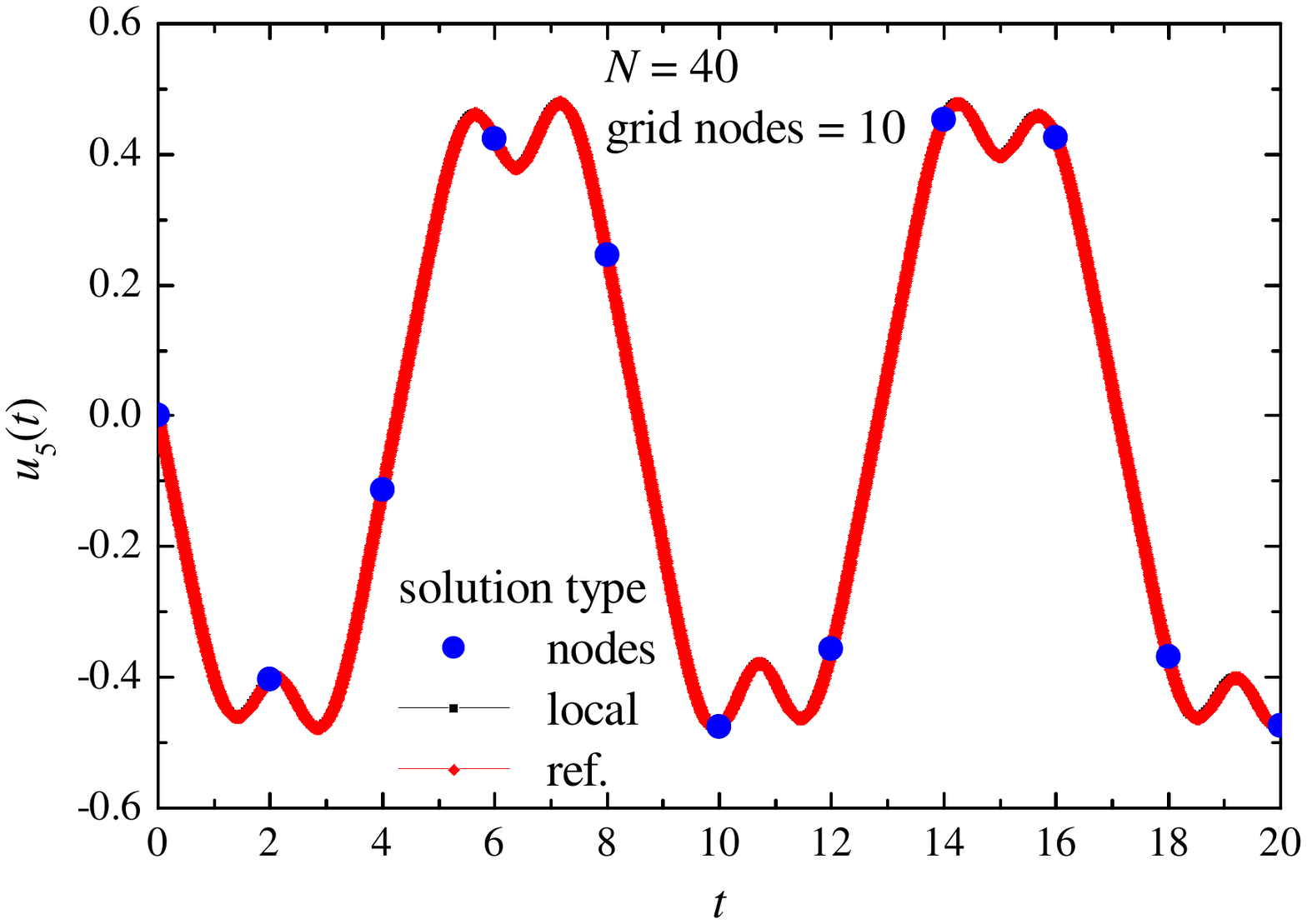}
\vspace{-8mm}\caption{\label{fig:dpend_ind2_sols_u:f1}}
\end{subfigure}
\begin{subfigure}{0.240\textwidth}
\includegraphics[width=\textwidth]{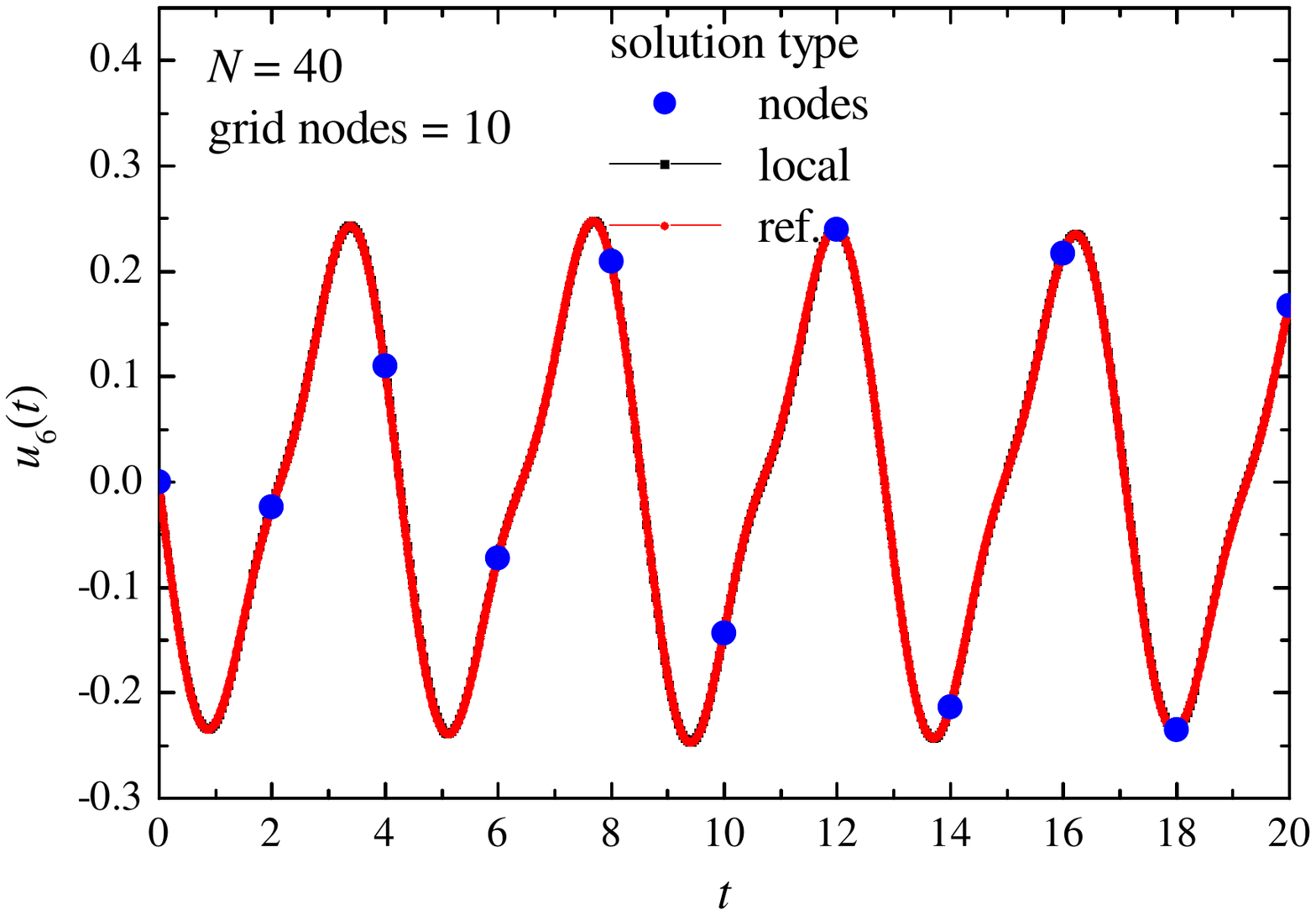}
\vspace{-8mm}\caption{\label{fig:dpend_ind2_sols_u:f2}}
\end{subfigure}
\begin{subfigure}{0.240\textwidth}
\includegraphics[width=\textwidth]{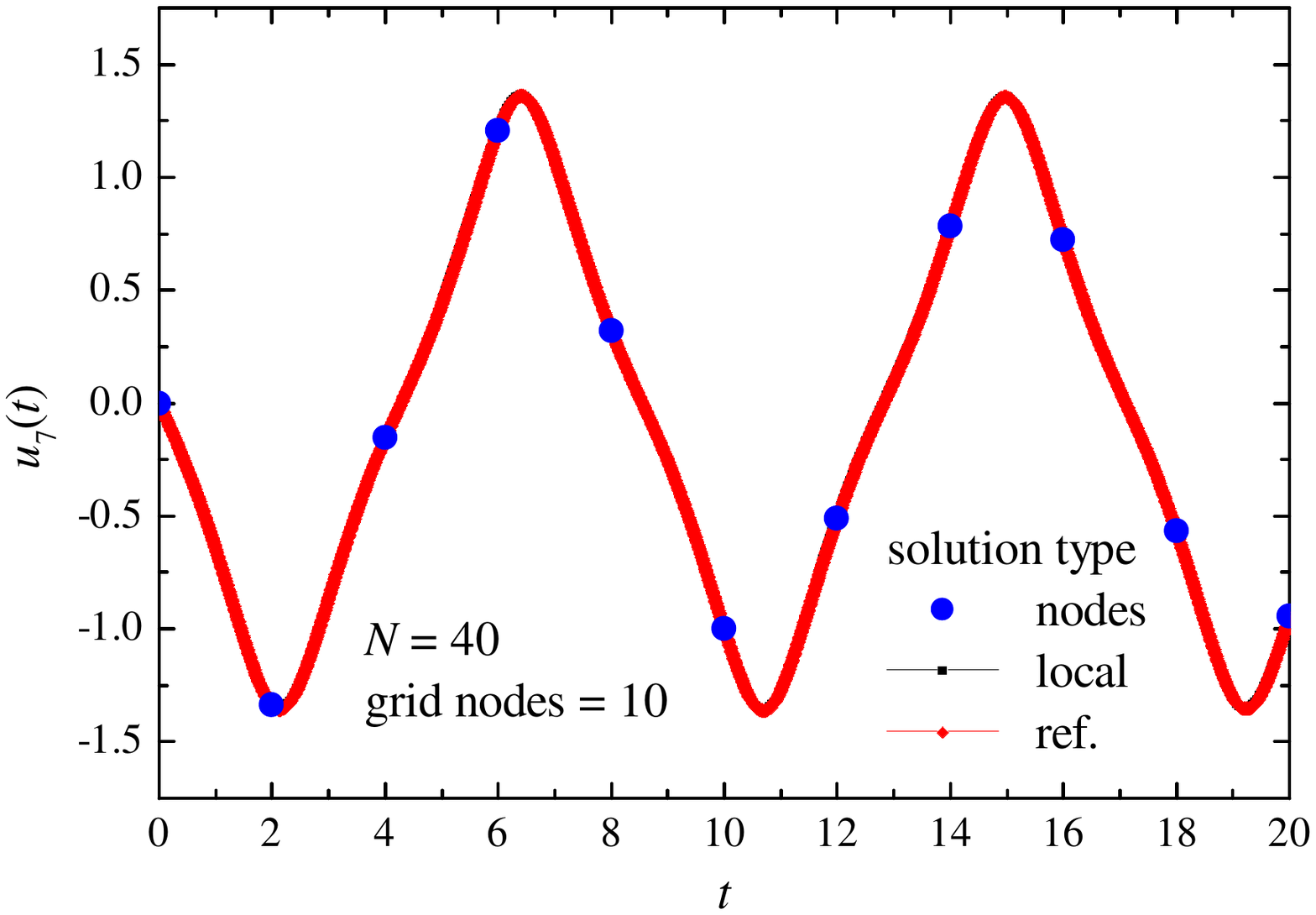}
\vspace{-8mm}\caption{\label{fig:dpend_ind2_sols_u:f3}}
\end{subfigure}
\begin{subfigure}{0.240\textwidth}
\includegraphics[width=\textwidth]{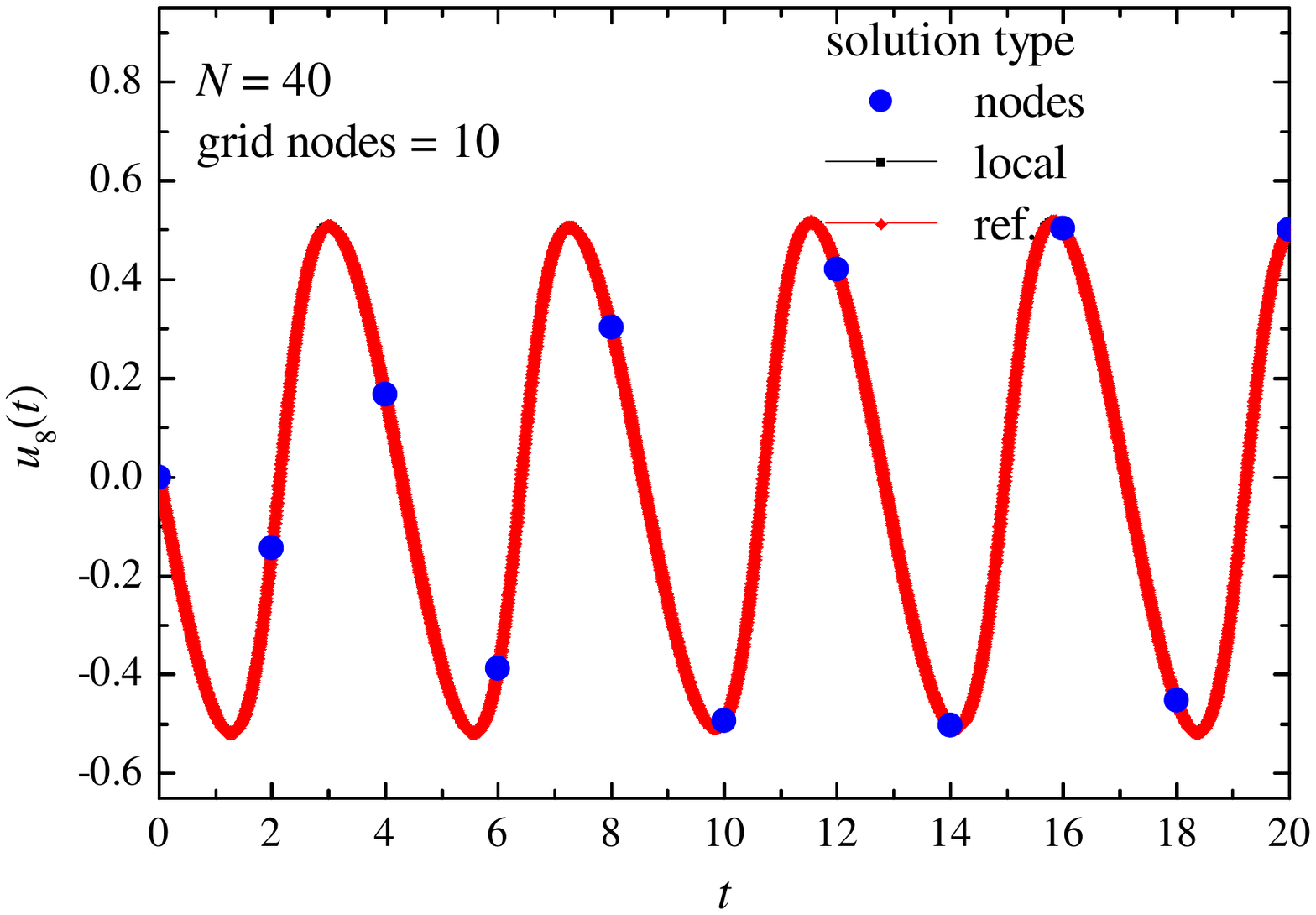}
\vspace{-8mm}\caption{\label{fig:dpend_ind2_sols_u:f4}}
\end{subfigure}\\[2mm]
\caption{%
Numerical solution of the DAE system (\ref{eq:math_dpend_dae_ind_3}) of index 2. Comparison of the solution at nodes $\mathbf{u}_{n}$, the local solution $\mathbf{u}_{L}(t)$ and the reference solution $\mathbf{u}^{\rm ref}(t)$ for components $u_{1}$ (\subref{fig:dpend_ind2_sols_u:a1}, \subref{fig:dpend_ind2_sols_u:c1}, \subref{fig:dpend_ind2_sols_u:e1}), $u_{2}$ (\subref{fig:dpend_ind2_sols_u:a2}, \subref{fig:dpend_ind2_sols_u:c2}, \subref{fig:dpend_ind2_sols_u:e2}), $u_{3}$ (\subref{fig:dpend_ind2_sols_u:a3}, \subref{fig:dpend_ind2_sols_u:c3}, \subref{fig:dpend_ind2_sols_u:e3}), $u_{4}$ (\subref{fig:dpend_ind2_sols_u:a4}, \subref{fig:dpend_ind2_sols_u:c4}, \subref{fig:dpend_ind2_sols_u:e4}), $u_{5}$ (\subref{fig:dpend_ind2_sols_u:b1}, \subref{fig:dpend_ind2_sols_u:d1}, \subref{fig:dpend_ind2_sols_u:f1}), $u_{6}$ (\subref{fig:dpend_ind2_sols_u:b2}, \subref{fig:dpend_ind2_sols_u:d2}, \subref{fig:dpend_ind2_sols_u:f2}), $u_{7}$ (\subref{fig:dpend_ind2_sols_u:b3}, \subref{fig:dpend_ind2_sols_u:d3}, \subref{fig:dpend_ind2_sols_u:f3}), $u_{8}$ (\subref{fig:dpend_ind2_sols_u:b4}, \subref{fig:dpend_ind2_sols_u:d4}, \subref{fig:dpend_ind2_sols_u:f4}), obtained using polynomials with degrees $N = 1$ (\subref{fig:dpend_ind2_sols_u:a1}, \subref{fig:dpend_ind2_sols_u:a2}, \subref{fig:dpend_ind2_sols_u:a3}, \subref{fig:dpend_ind2_sols_u:a4}, \subref{fig:dpend_ind2_sols_u:b1}, \subref{fig:dpend_ind2_sols_u:b2}, \subref{fig:dpend_ind2_sols_u:b3}, \subref{fig:dpend_ind2_sols_u:b4}), $N = 8$ (\subref{fig:dpend_ind2_sols_u:c1}, \subref{fig:dpend_ind2_sols_u:c2}, \subref{fig:dpend_ind2_sols_u:c3}, \subref{fig:dpend_ind2_sols_u:c4}, \subref{fig:dpend_ind2_sols_u:d1}, \subref{fig:dpend_ind2_sols_u:d2}, \subref{fig:dpend_ind2_sols_u:d3}, \subref{fig:dpend_ind2_sols_u:d4}) and $N = 40$ (\subref{fig:dpend_ind2_sols_u:e1}, \subref{fig:dpend_ind2_sols_u:e2}, \subref{fig:dpend_ind2_sols_u:e3}, \subref{fig:dpend_ind2_sols_u:e4}, \subref{fig:dpend_ind2_sols_u:f1}, \subref{fig:dpend_ind2_sols_u:f2}, \subref{fig:dpend_ind2_sols_u:f3}, \subref{fig:dpend_ind2_sols_u:f4}).
}
\label{fig:dpend_ind2_sols_u}
\end{figure} 

\begin{figure}[h!]
\captionsetup[subfigure]{%
	position=bottom,
	font+=smaller,
	textfont=normalfont,
	singlelinecheck=off,
	justification=raggedright
}
\centering
\begin{subfigure}{0.240\textwidth}
\includegraphics[width=\textwidth]{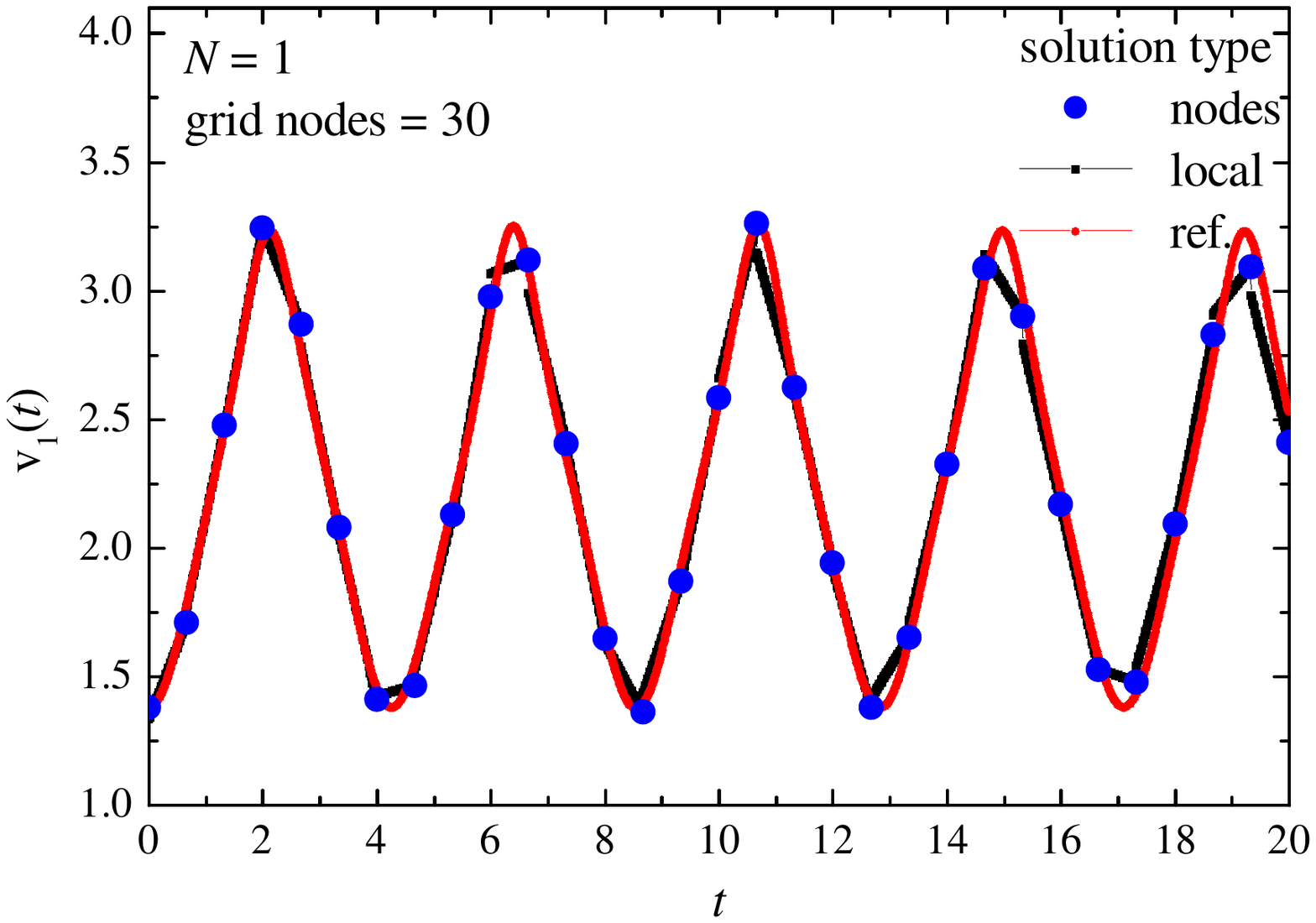}
\vspace{-8mm}\caption{\label{fig:dpend_ind2_sols_vg:a1}}
\end{subfigure}
\begin{subfigure}{0.240\textwidth}
\includegraphics[width=\textwidth]{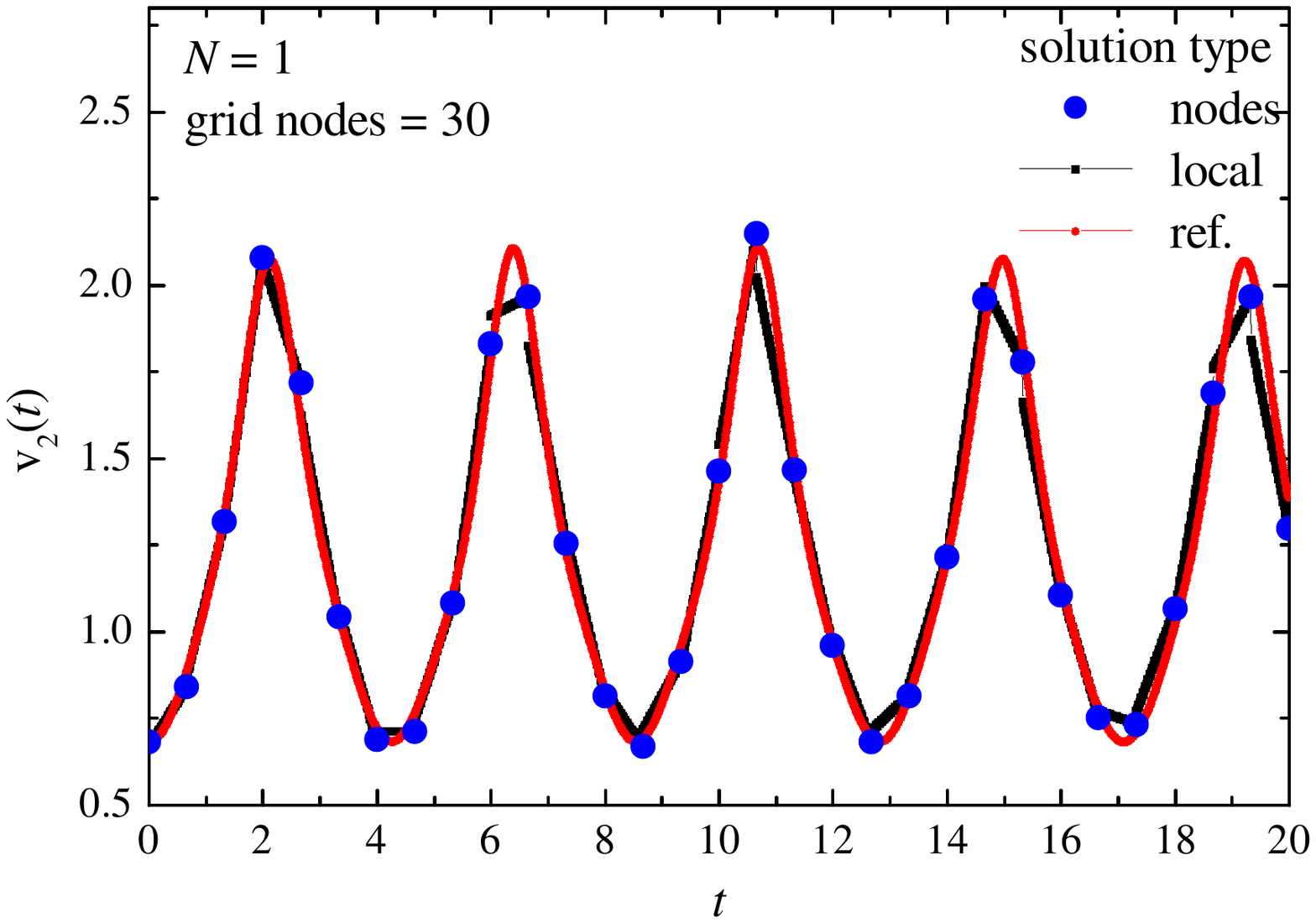}
\vspace{-8mm}\caption{\label{fig:dpend_ind2_sols_vg:a2}}
\end{subfigure}
\begin{subfigure}{0.240\textwidth}
\includegraphics[width=\textwidth]{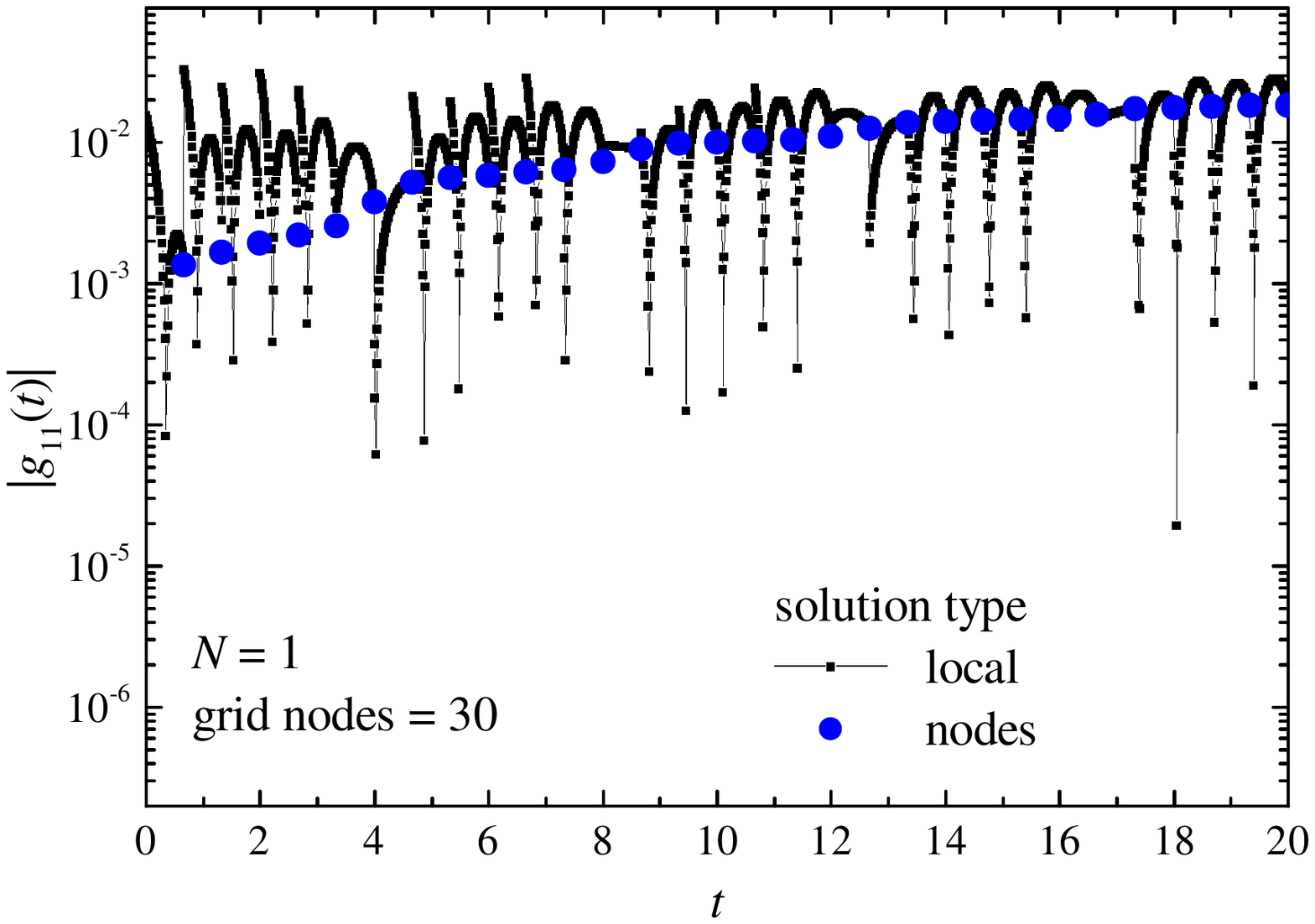}
\vspace{-8mm}\caption{\label{fig:dpend_ind2_sols_vg:a3}}
\end{subfigure}
\begin{subfigure}{0.240\textwidth}
\includegraphics[width=\textwidth]{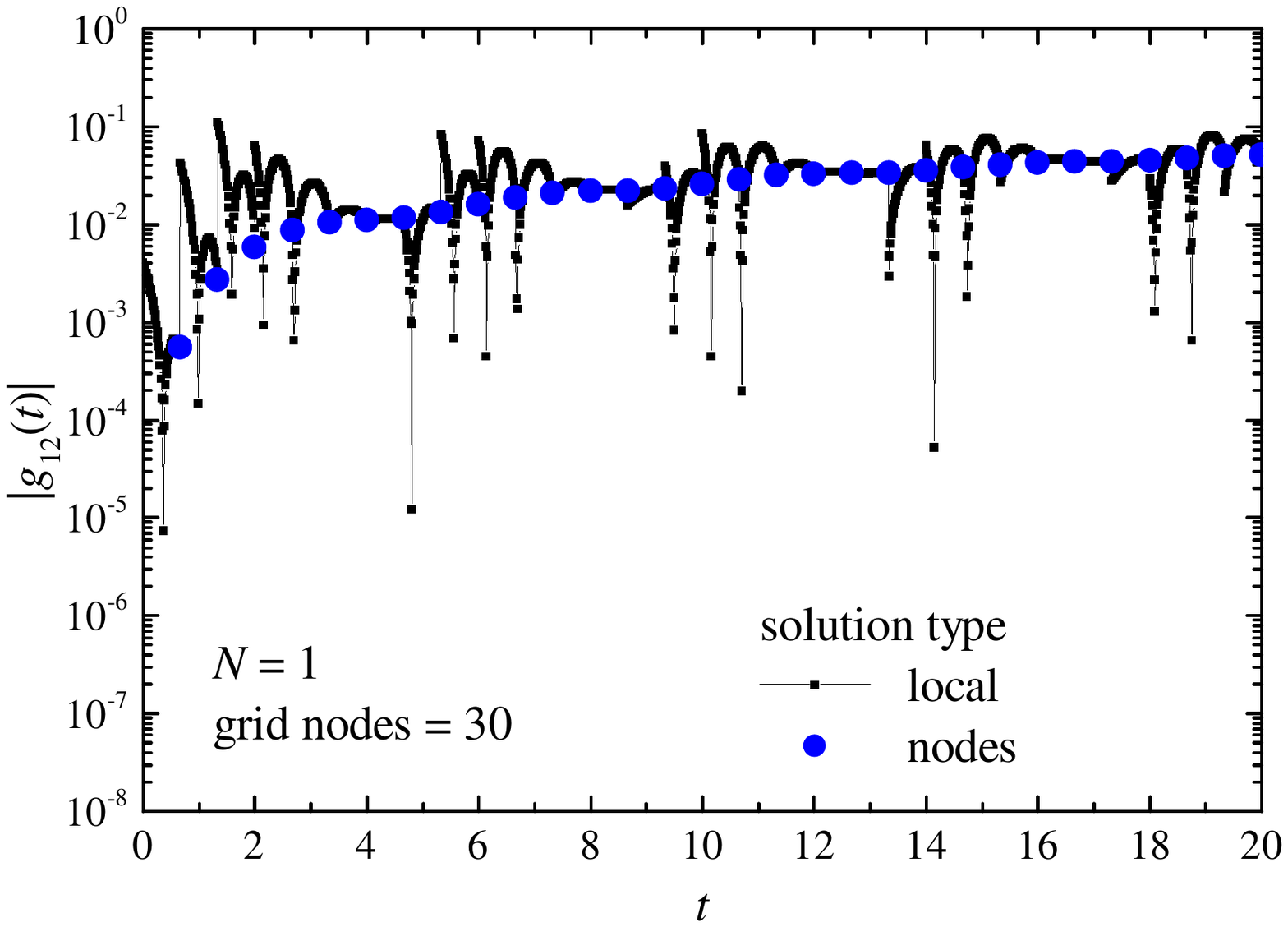}
\vspace{-8mm}\caption{\label{fig:dpend_ind2_sols_vg:a4}}
\end{subfigure}\\[2mm]
\begin{subfigure}{0.240\textwidth}
\includegraphics[width=\textwidth]{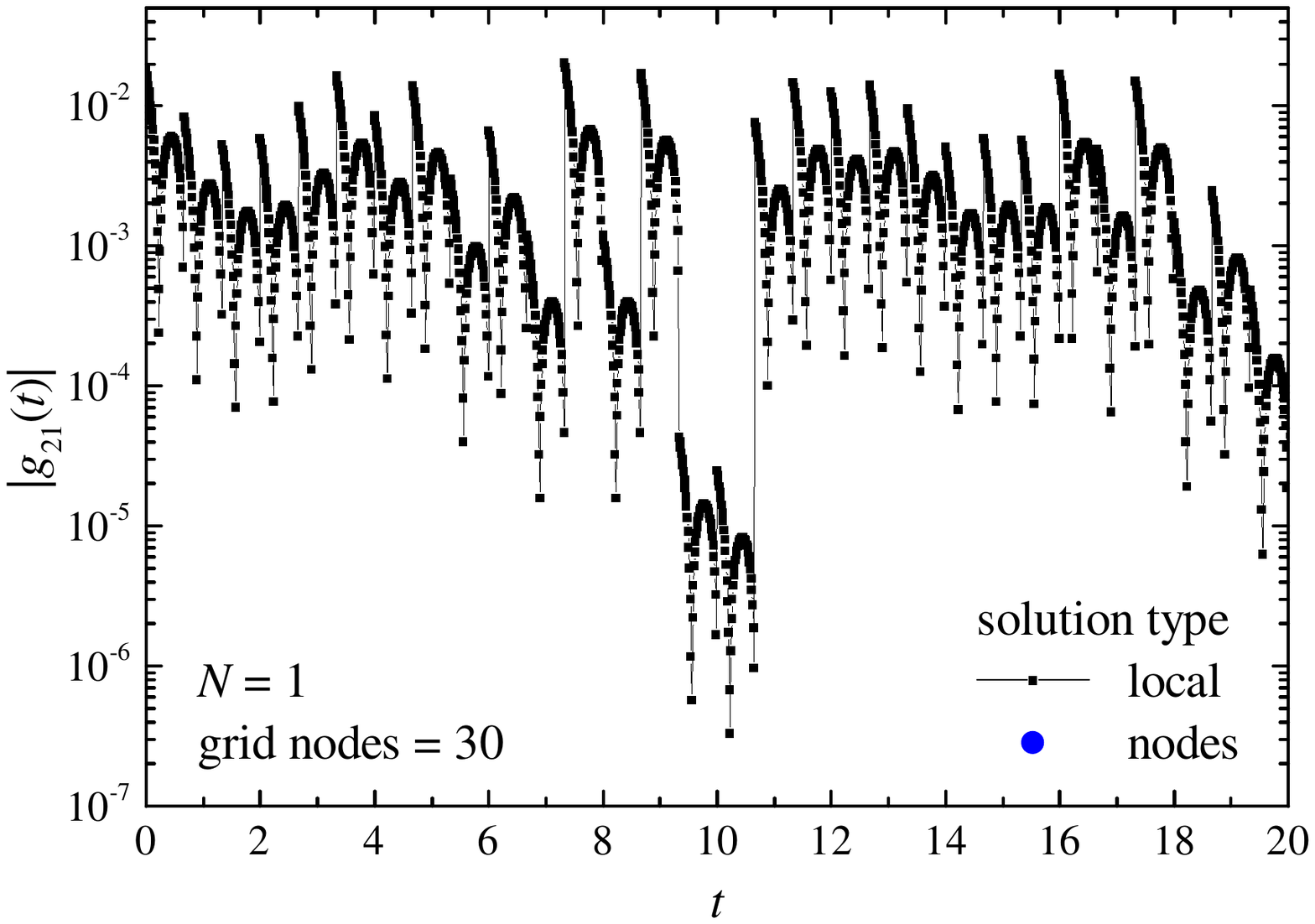}
\vspace{-8mm}\caption{\label{fig:dpend_ind2_sols_vg:b1}}
\end{subfigure}
\begin{subfigure}{0.240\textwidth}
\includegraphics[width=\textwidth]{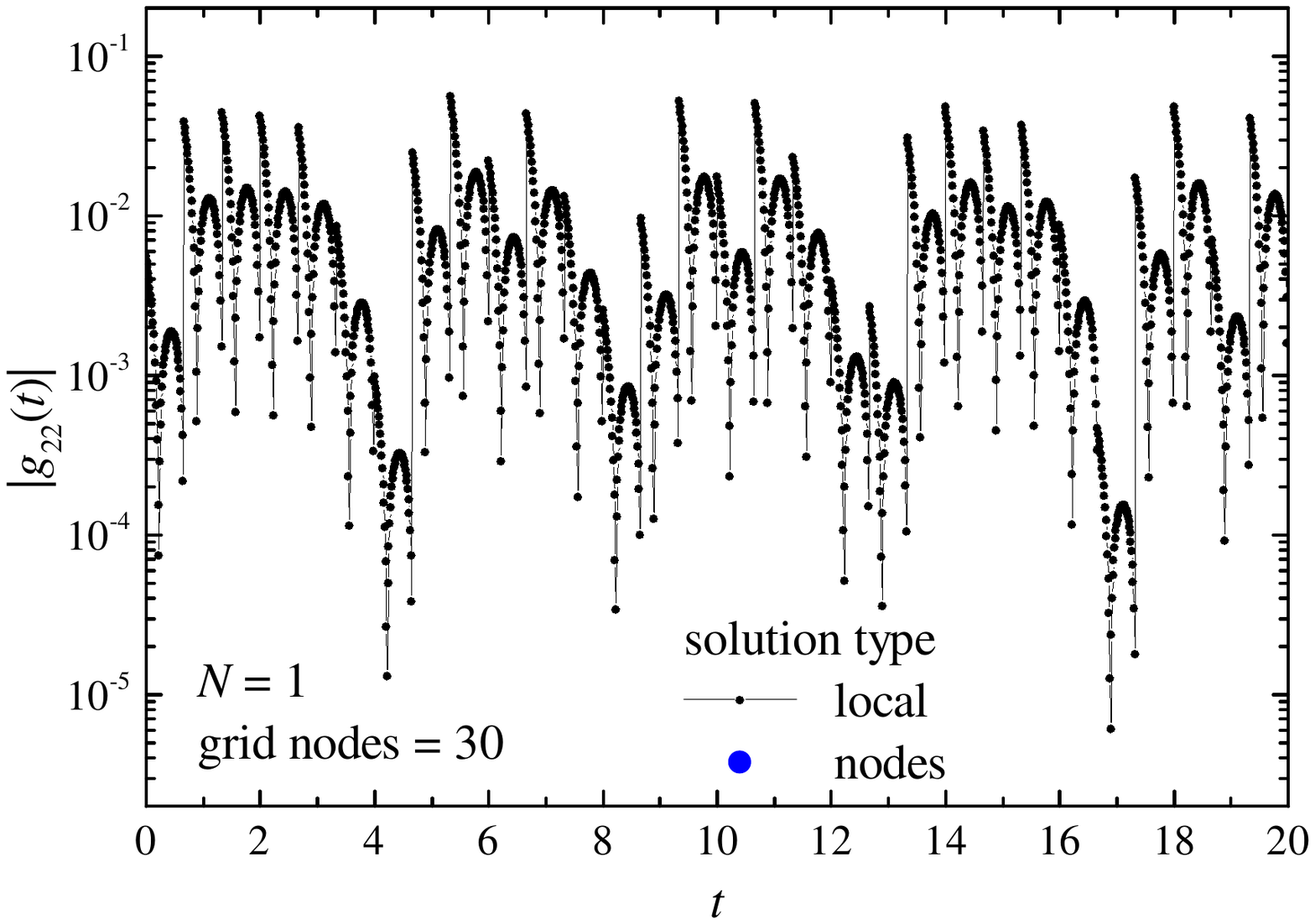}
\vspace{-8mm}\caption{\label{fig:dpend_ind2_sols_vg:b2}}
\end{subfigure}
\begin{subfigure}{0.240\textwidth}
\includegraphics[width=\textwidth]{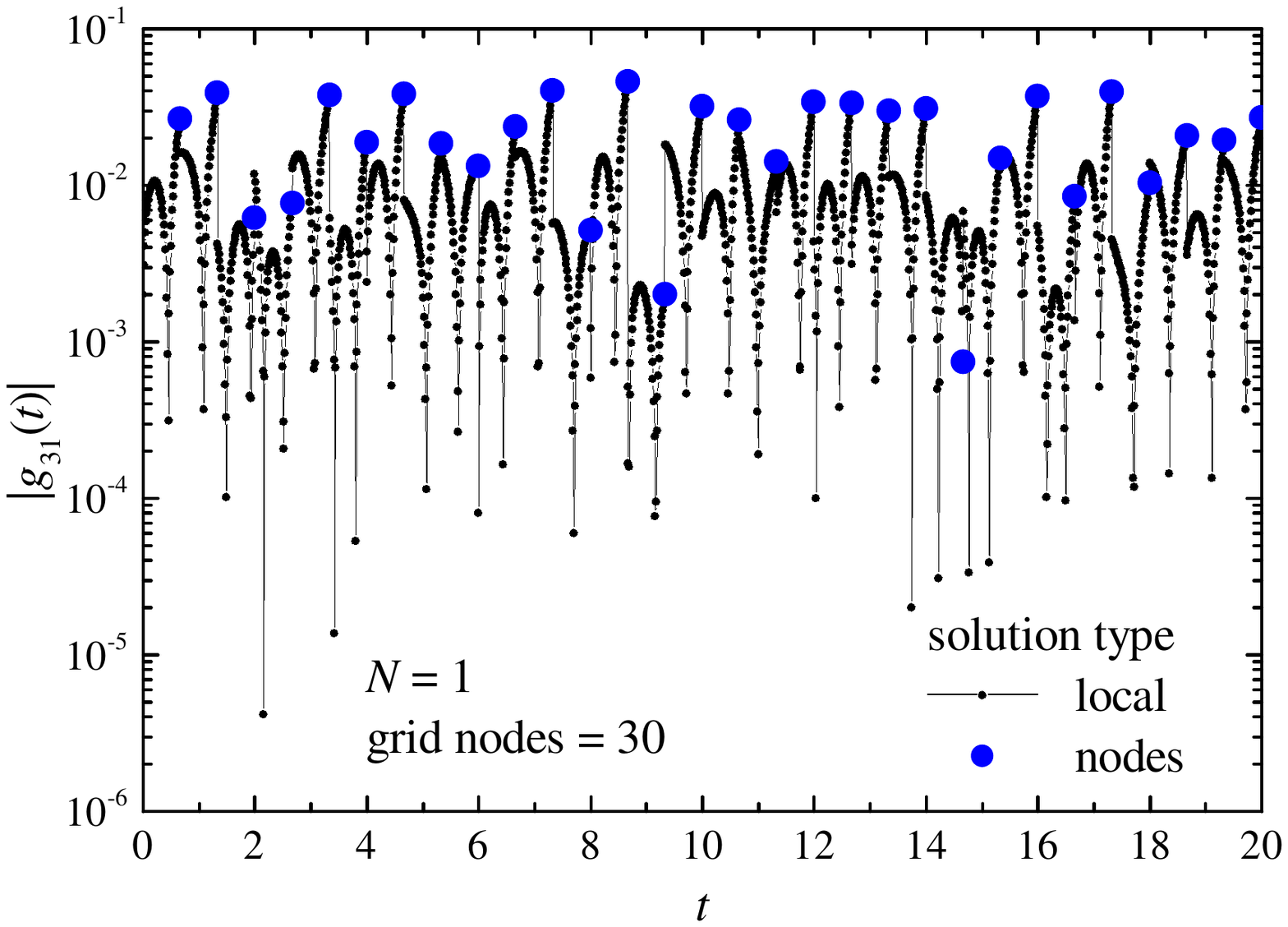}
\vspace{-8mm}\caption{\label{fig:dpend_ind2_sols_vg:b3}}
\end{subfigure}
\begin{subfigure}{0.240\textwidth}
\includegraphics[width=\textwidth]{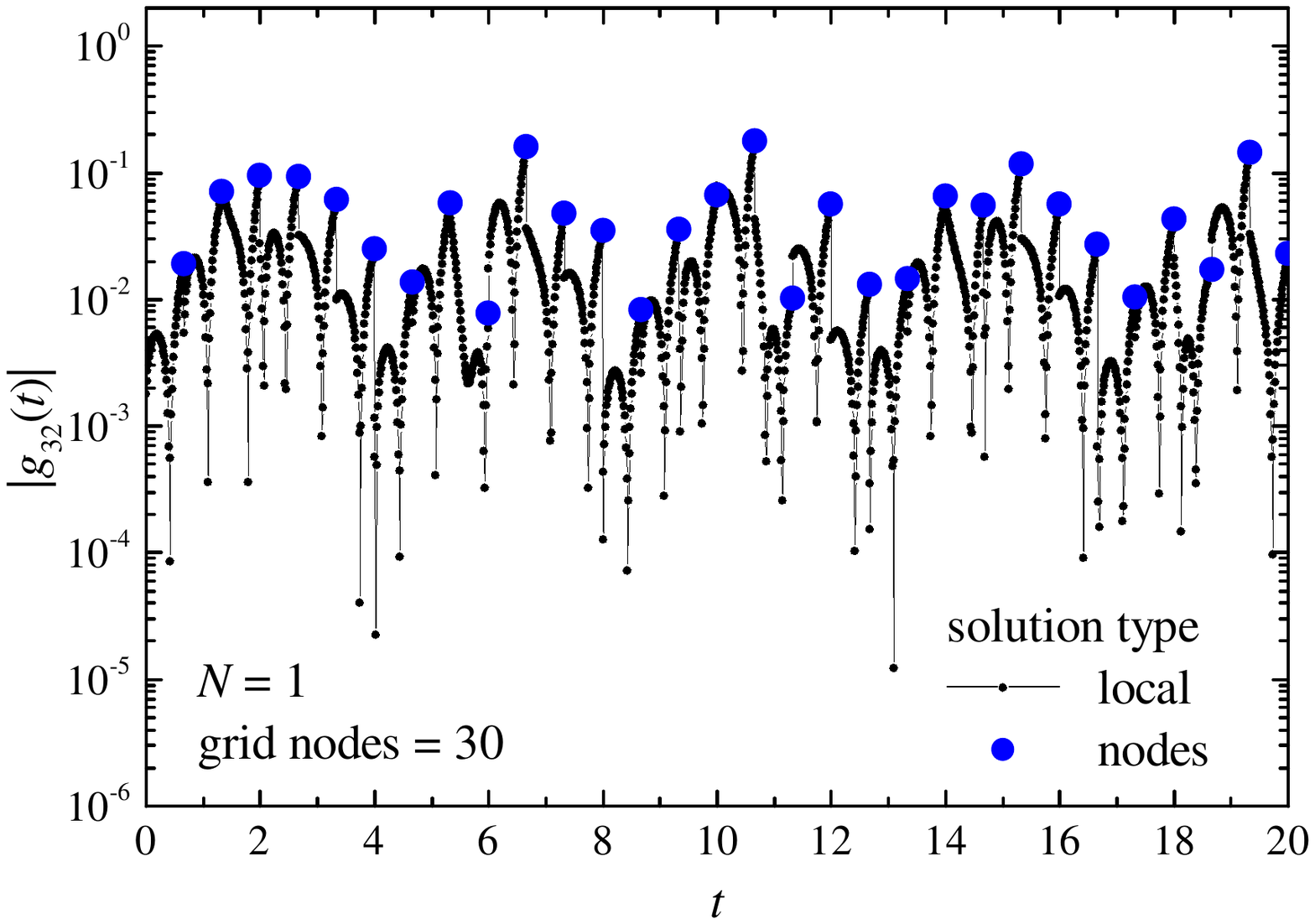}
\vspace{-8mm}\caption{\label{fig:dpend_ind2_sols_vg:b4}}
\end{subfigure}\\[2mm]
\begin{subfigure}{0.240\textwidth}
\includegraphics[width=\textwidth]{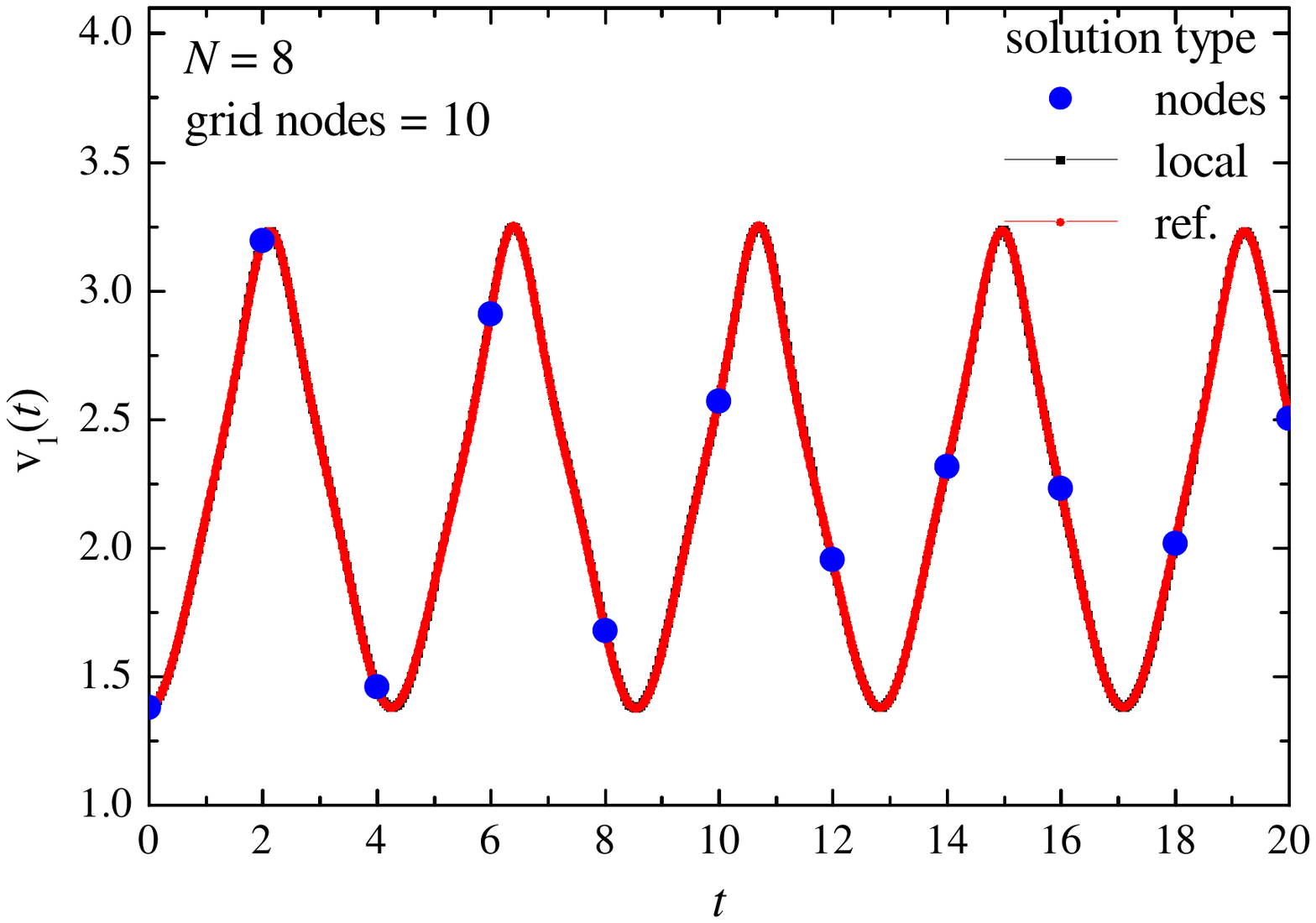}
\vspace{-8mm}\caption{\label{fig:dpend_ind2_sols_vg:c1}}
\end{subfigure}
\begin{subfigure}{0.240\textwidth}
\includegraphics[width=\textwidth]{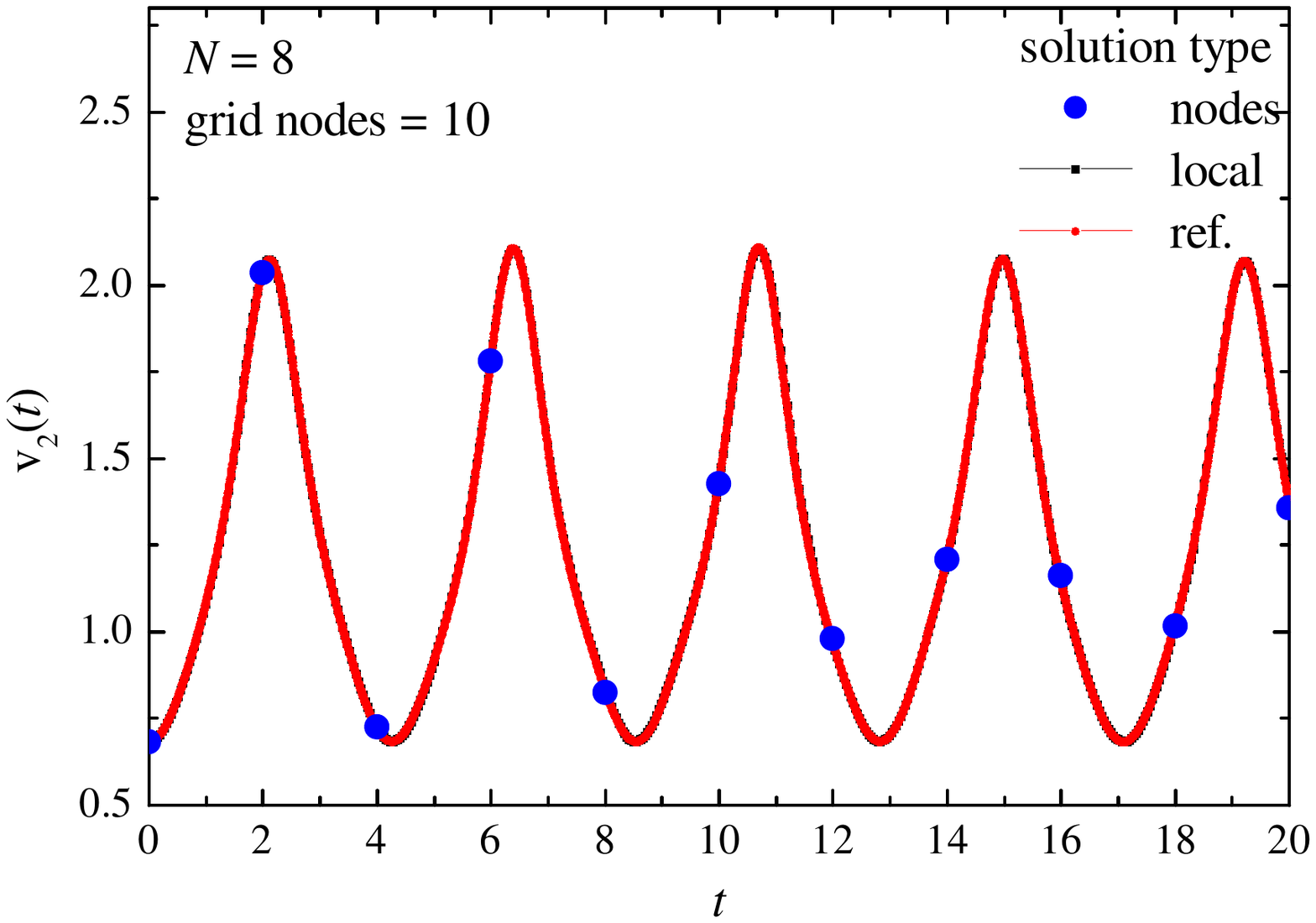}
\vspace{-8mm}\caption{\label{fig:dpend_ind2_sols_vg:c2}}
\end{subfigure}
\begin{subfigure}{0.240\textwidth}
\includegraphics[width=\textwidth]{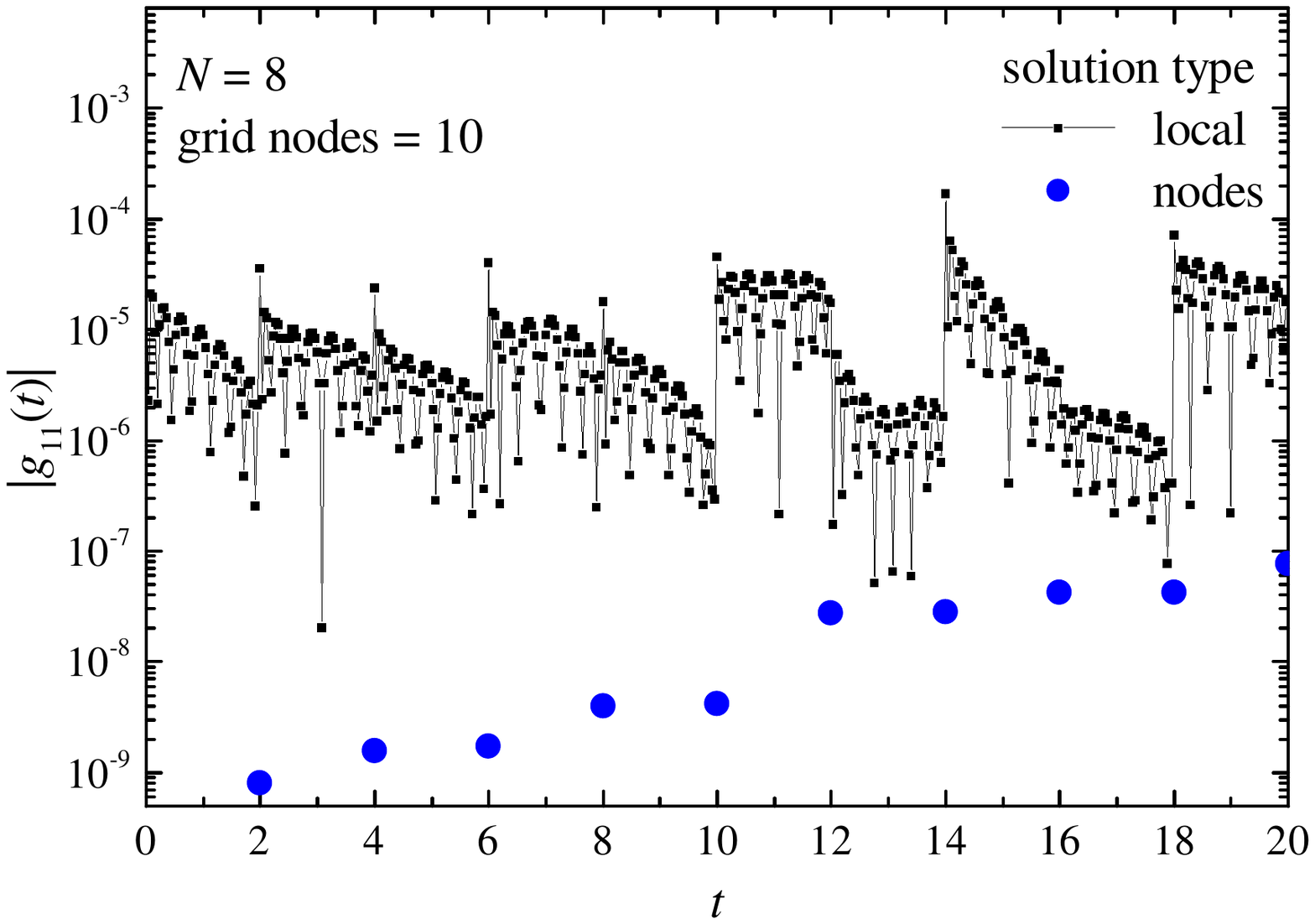}
\vspace{-8mm}\caption{\label{fig:dpend_ind2_sols_vg:c3}}
\end{subfigure}
\begin{subfigure}{0.240\textwidth}
\includegraphics[width=\textwidth]{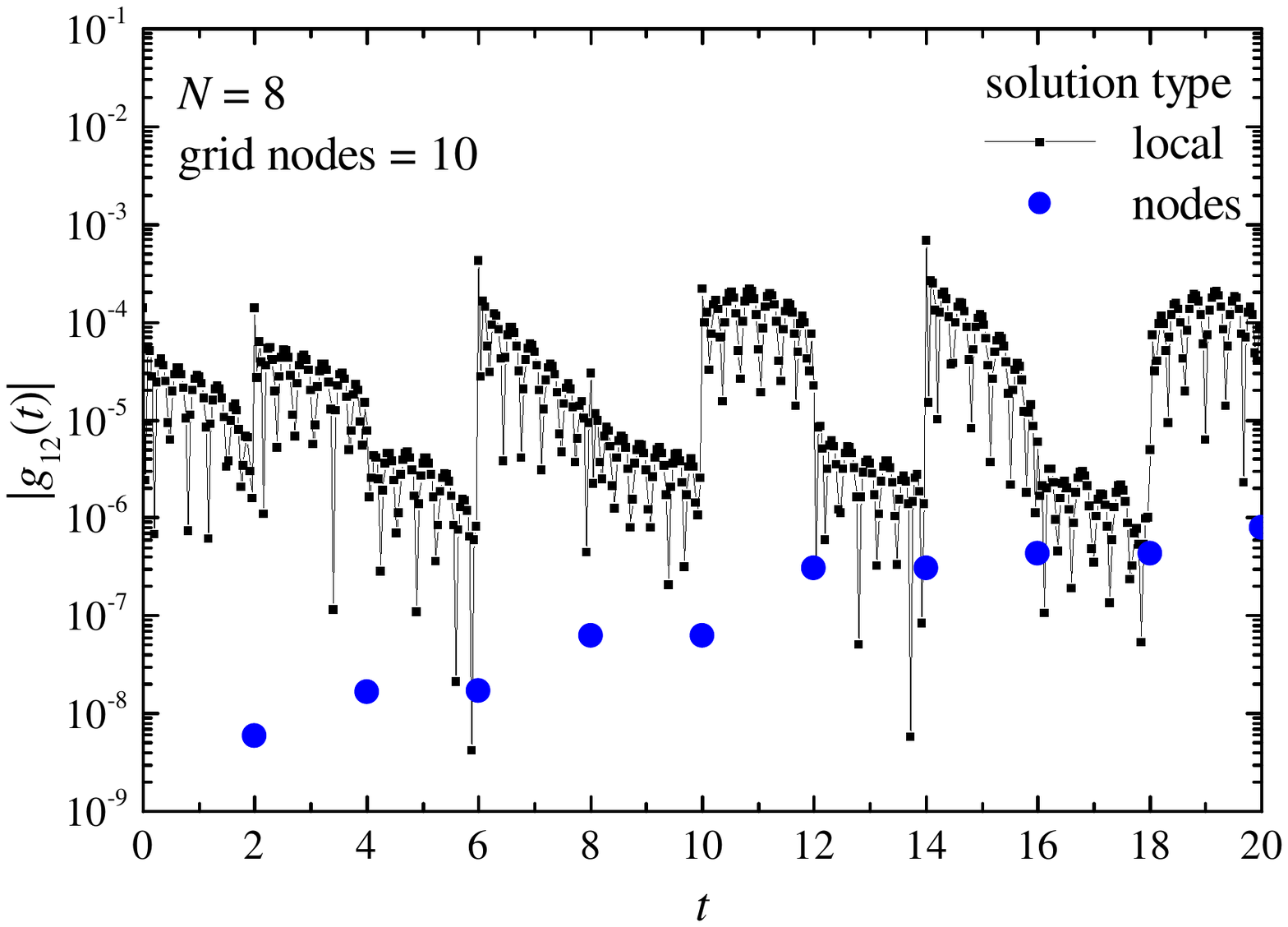}
\vspace{-8mm}\caption{\label{fig:dpend_ind2_sols_vg:c4}}
\end{subfigure}\\[2mm]
\begin{subfigure}{0.240\textwidth}
\includegraphics[width=\textwidth]{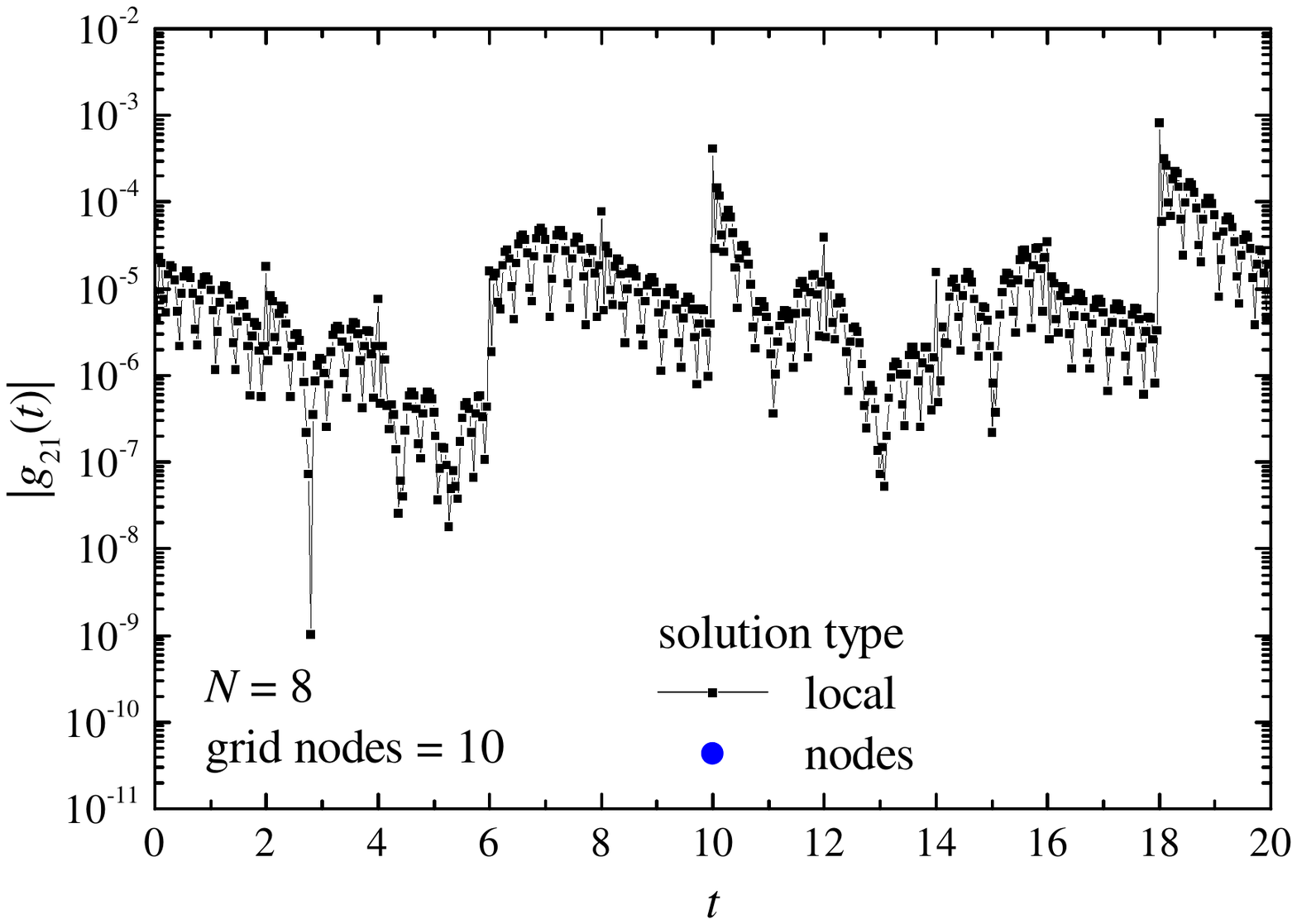}
\vspace{-8mm}\caption{\label{fig:dpend_ind2_sols_vg:d1}}
\end{subfigure}
\begin{subfigure}{0.240\textwidth}
\includegraphics[width=\textwidth]{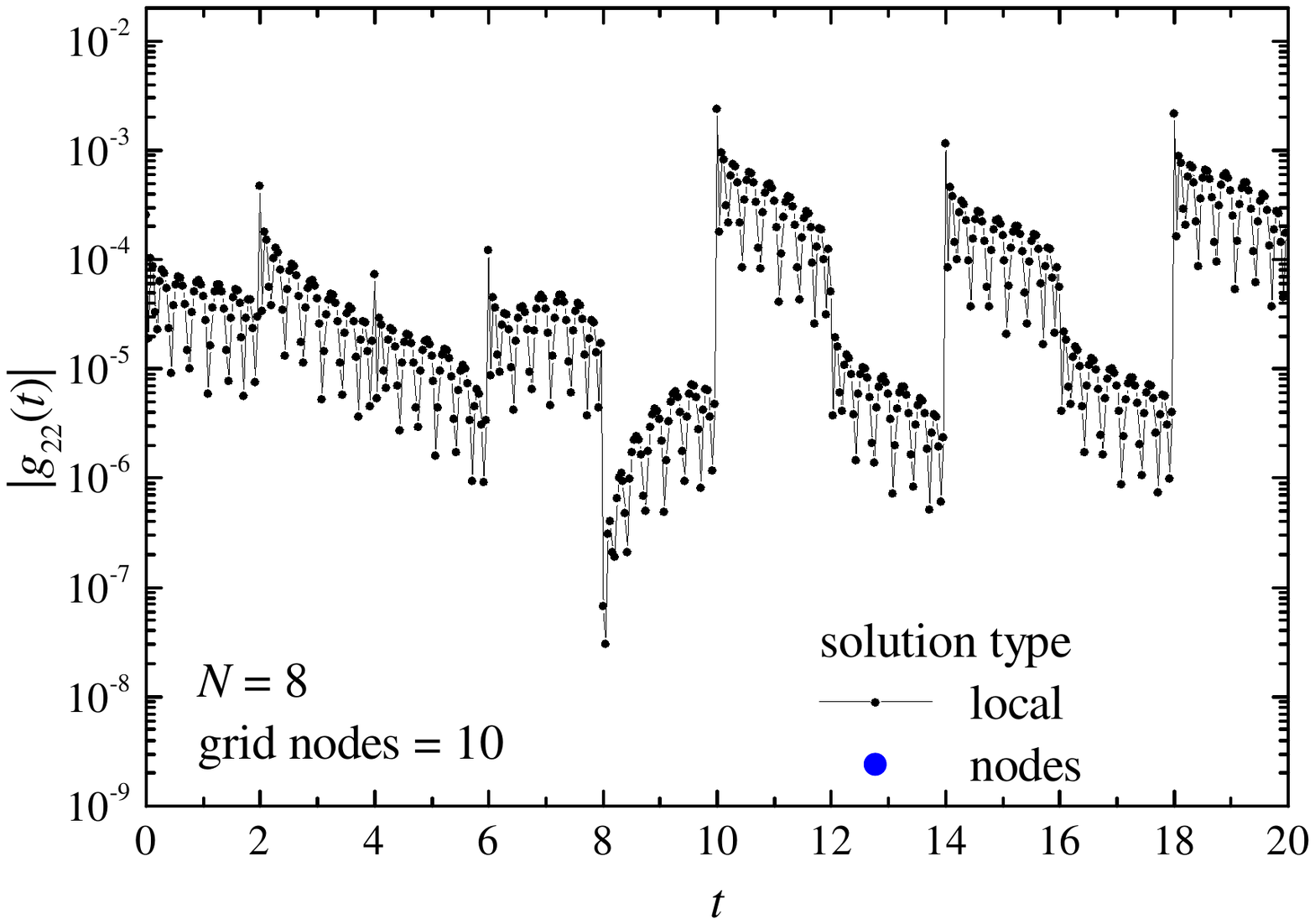}
\vspace{-8mm}\caption{\label{fig:dpend_ind2_sols_vg:d2}}
\end{subfigure}
\begin{subfigure}{0.240\textwidth}
\includegraphics[width=\textwidth]{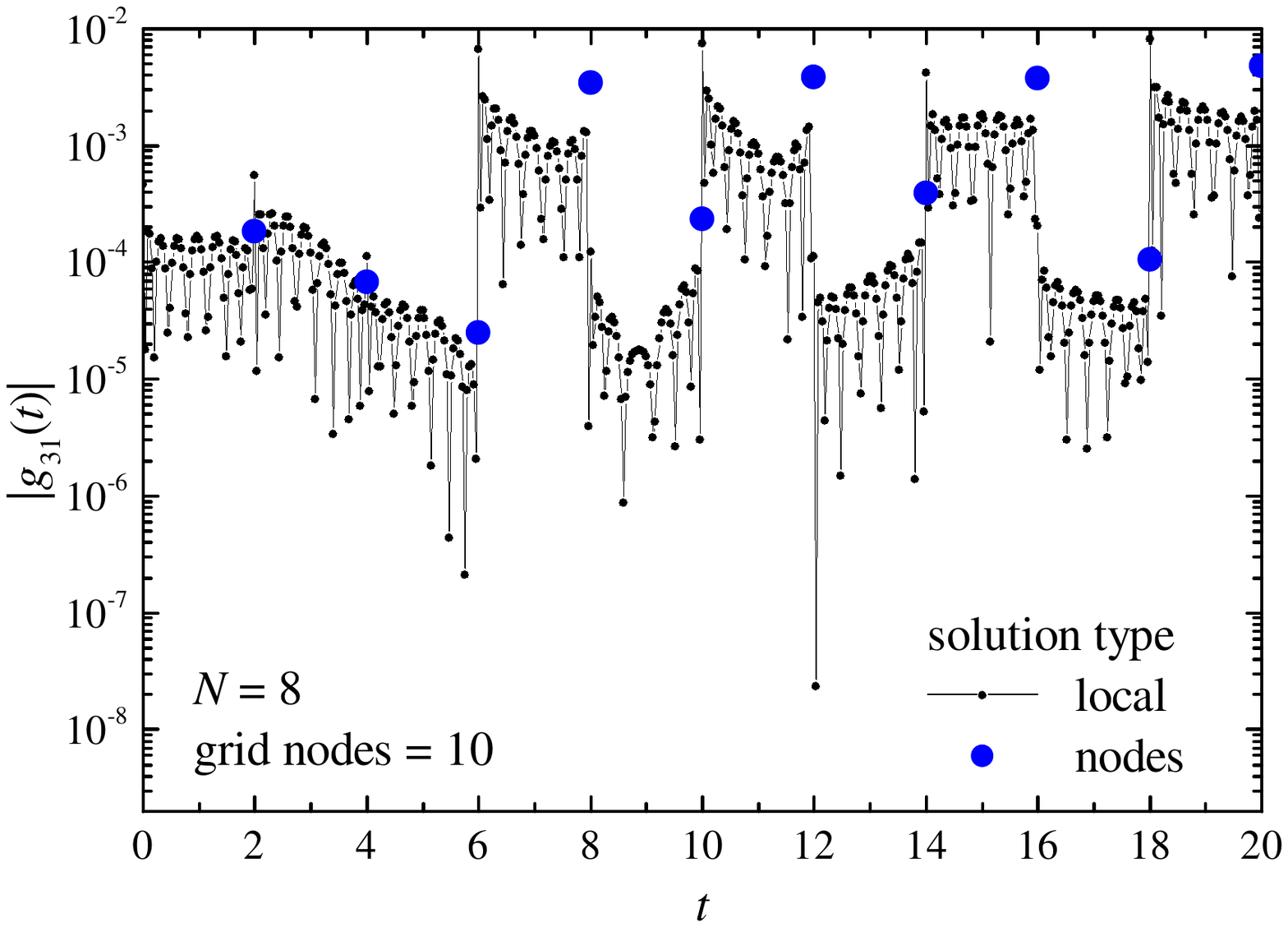}
\vspace{-8mm}\caption{\label{fig:dpend_ind2_sols_vg:d3}}
\end{subfigure}
\begin{subfigure}{0.240\textwidth}
\includegraphics[width=\textwidth]{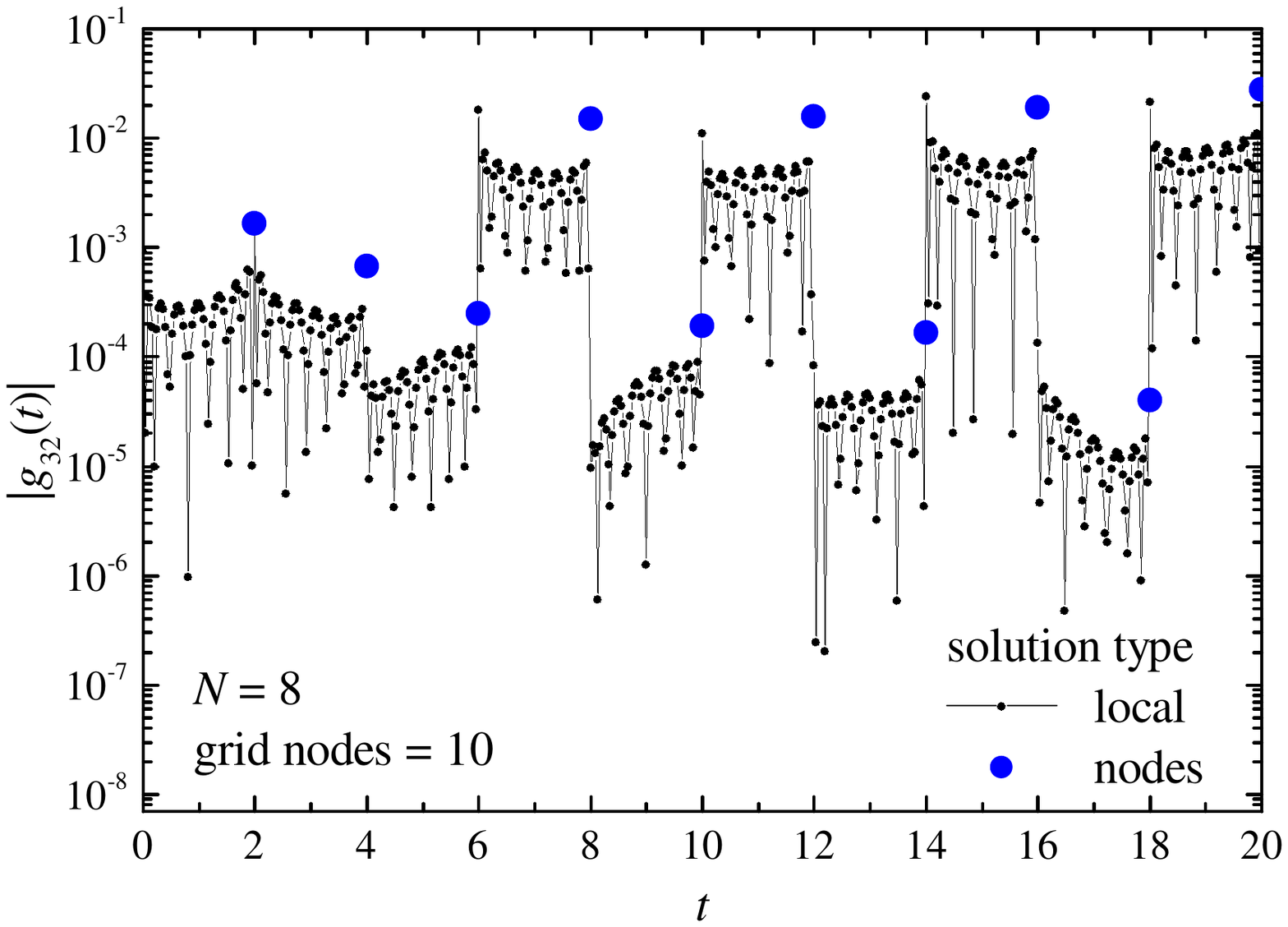}
\vspace{-8mm}\caption{\label{fig:dpend_ind2_sols_vg:d4}}
\end{subfigure}\\[2mm]
\begin{subfigure}{0.240\textwidth}
\includegraphics[width=\textwidth]{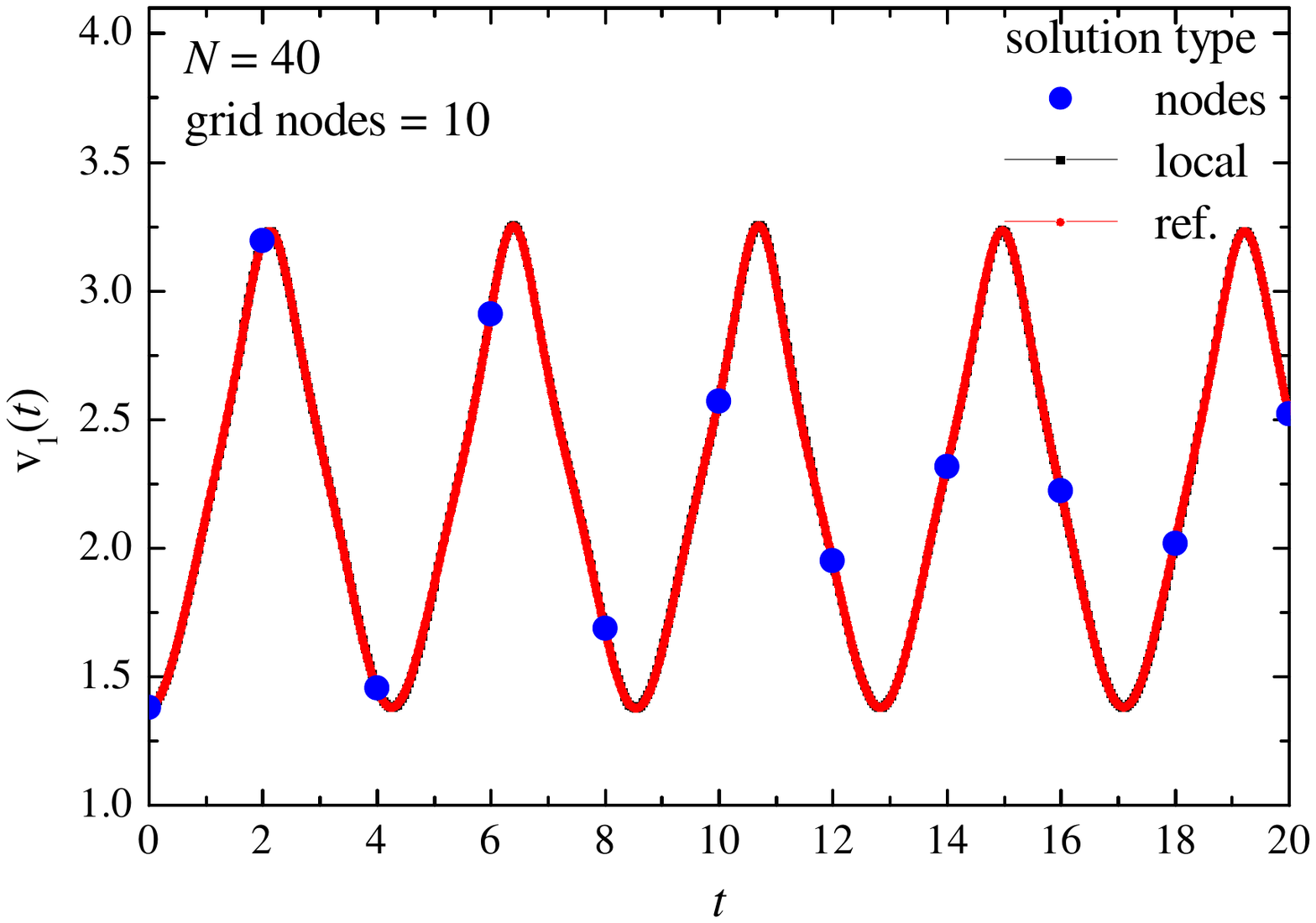}
\vspace{-8mm}\caption{\label{fig:dpend_ind2_sols_vg:e1}}
\end{subfigure}
\begin{subfigure}{0.240\textwidth}
\includegraphics[width=\textwidth]{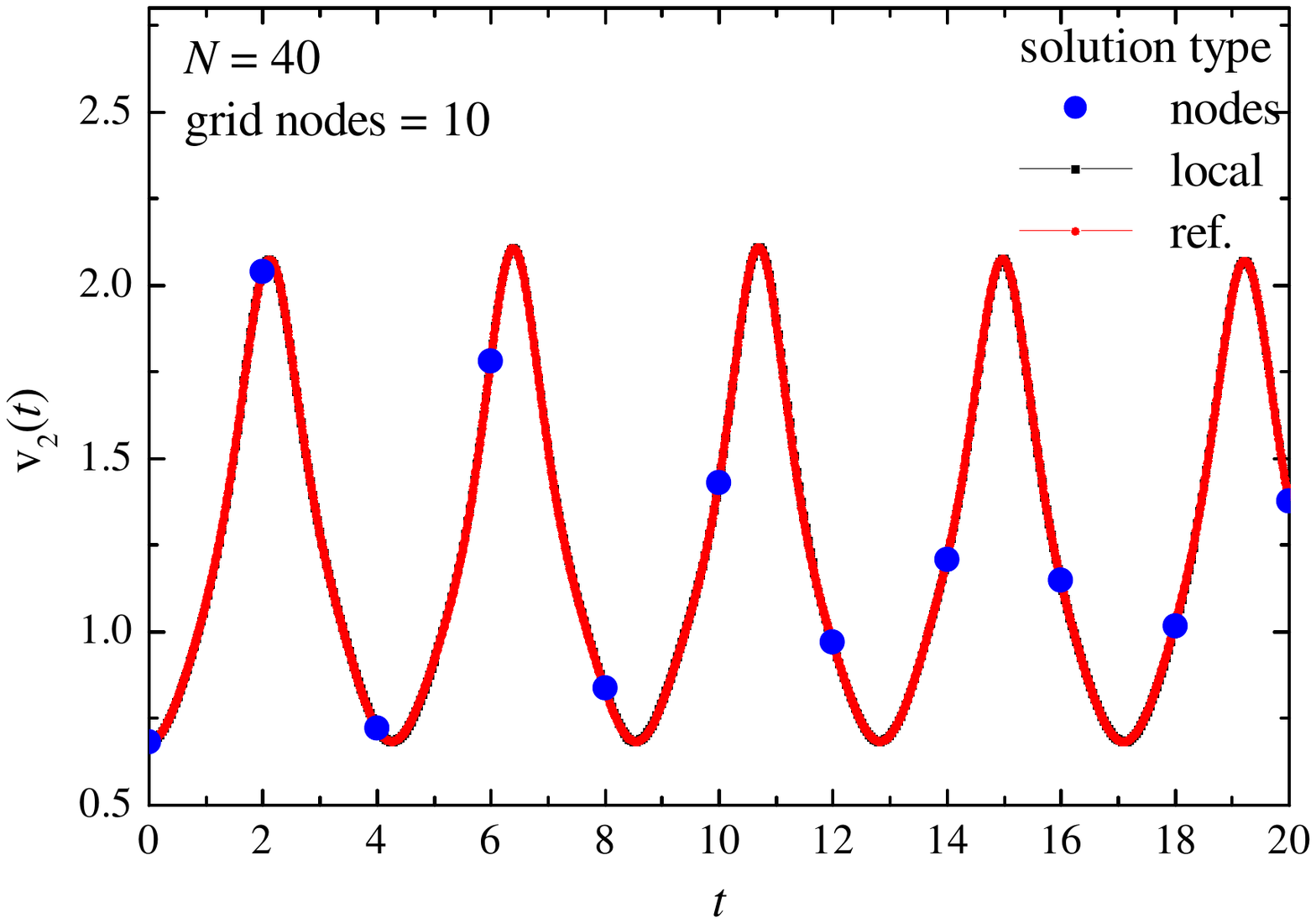}
\vspace{-8mm}\caption{\label{fig:dpend_ind2_sols_vg:e2}}
\end{subfigure}
\begin{subfigure}{0.240\textwidth}
\includegraphics[width=\textwidth]{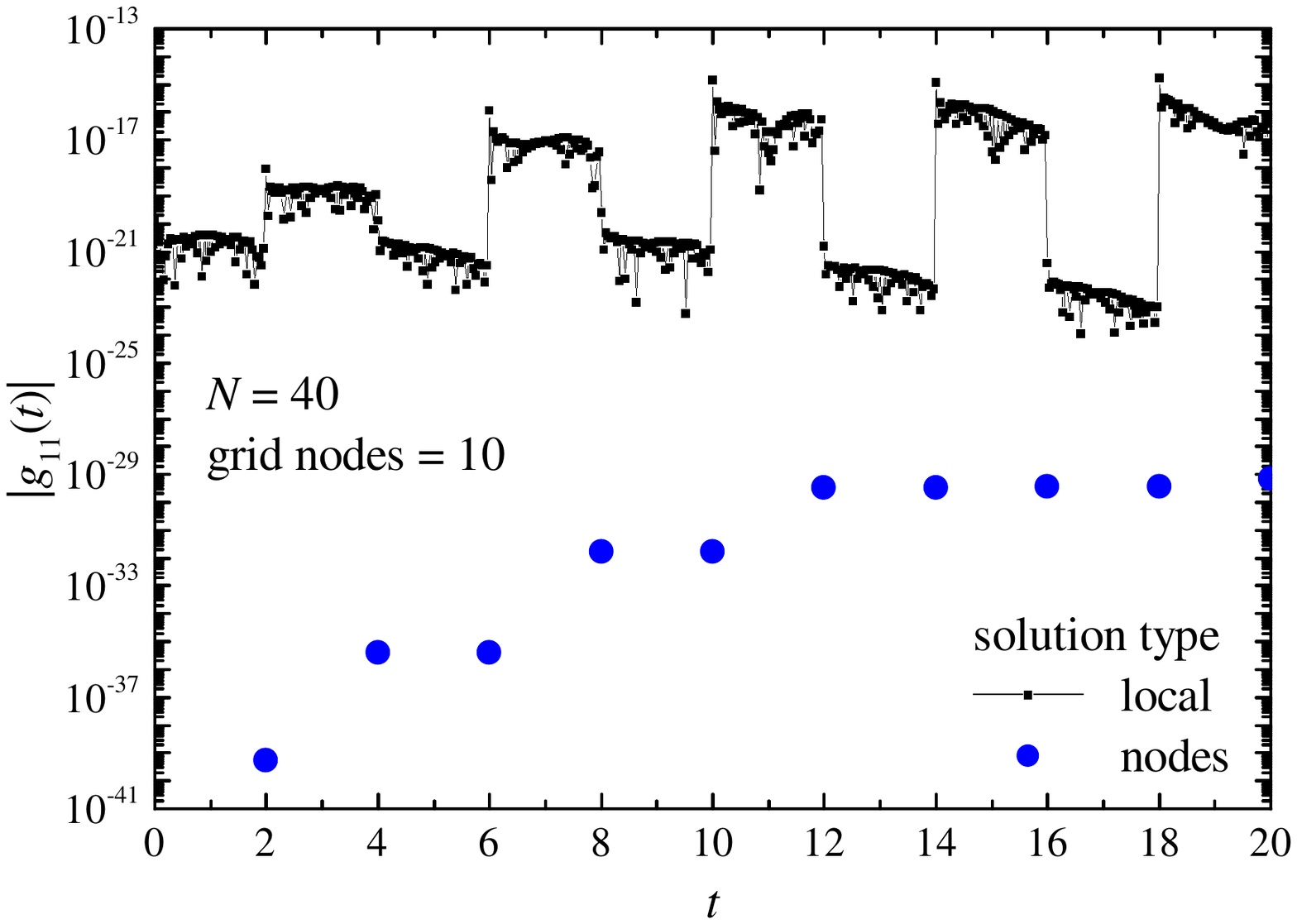}
\vspace{-8mm}\caption{\label{fig:dpend_ind2_sols_vg:e3}}
\end{subfigure}
\begin{subfigure}{0.240\textwidth}
\includegraphics[width=\textwidth]{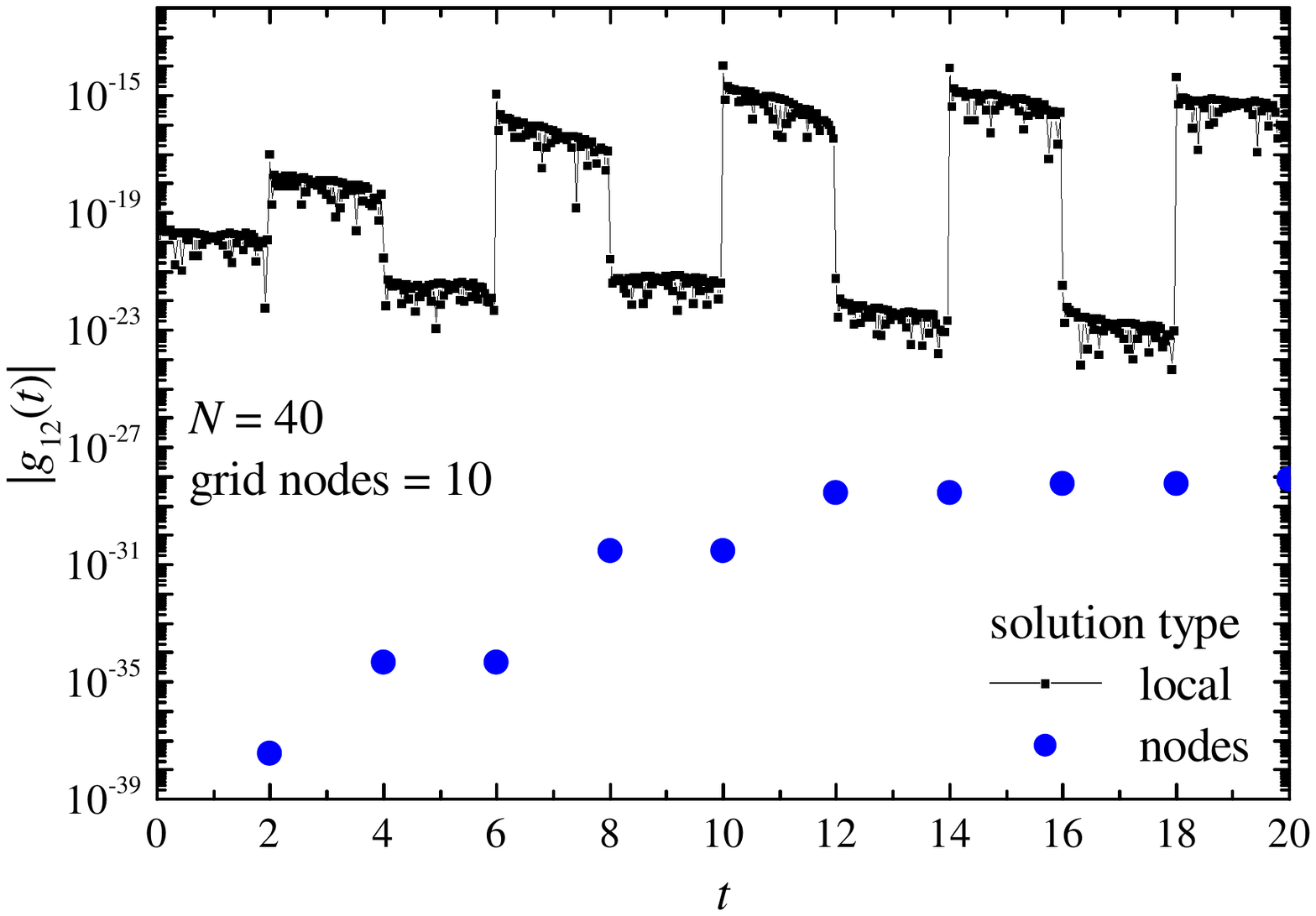}
\vspace{-8mm}\caption{\label{fig:dpend_ind2_sols_vg:e4}}
\end{subfigure}\\[2mm]
\begin{subfigure}{0.240\textwidth}
\includegraphics[width=\textwidth]{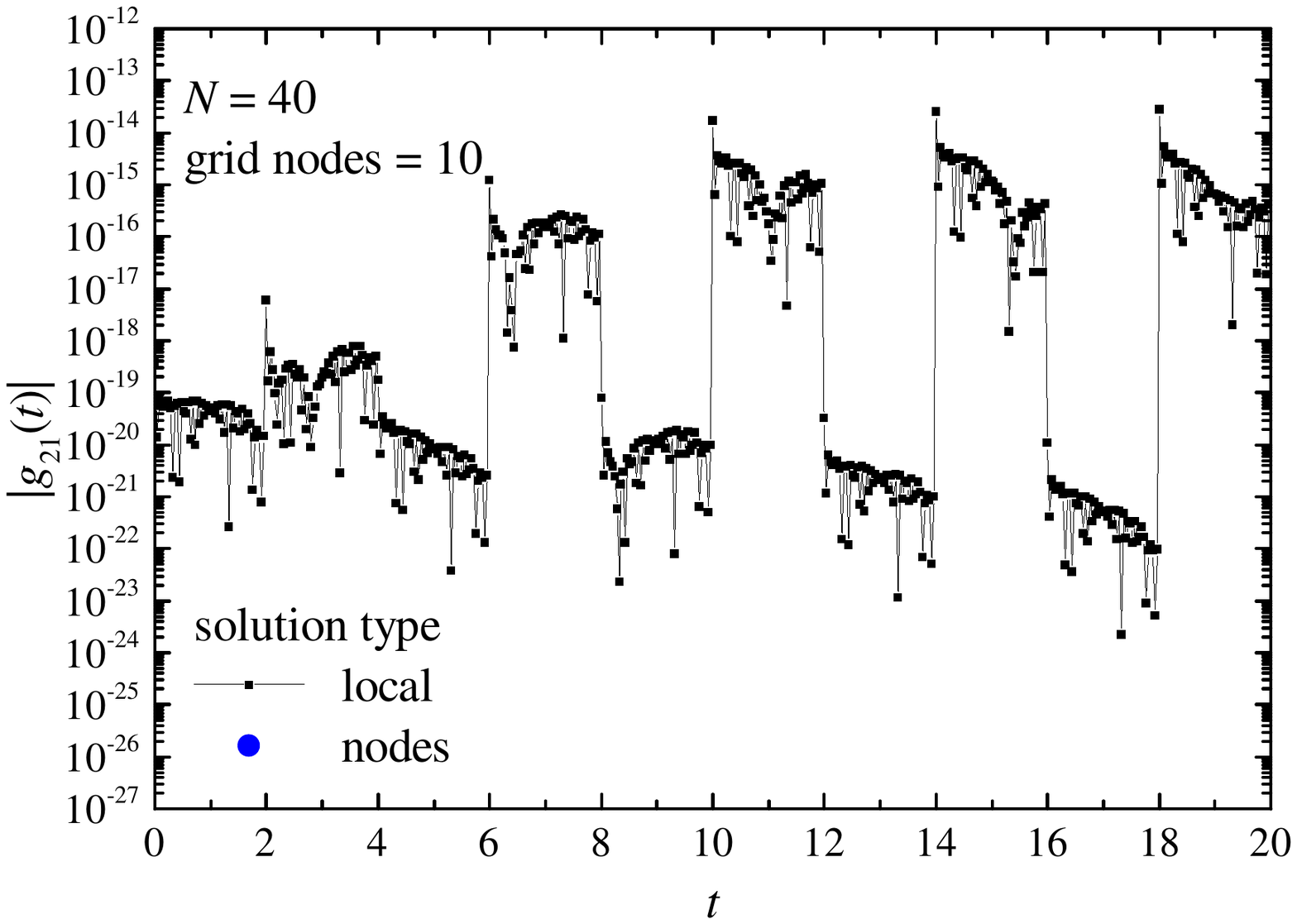}
\vspace{-8mm}\caption{\label{fig:dpend_ind2_sols_vg:f1}}
\end{subfigure}
\begin{subfigure}{0.240\textwidth}
\includegraphics[width=\textwidth]{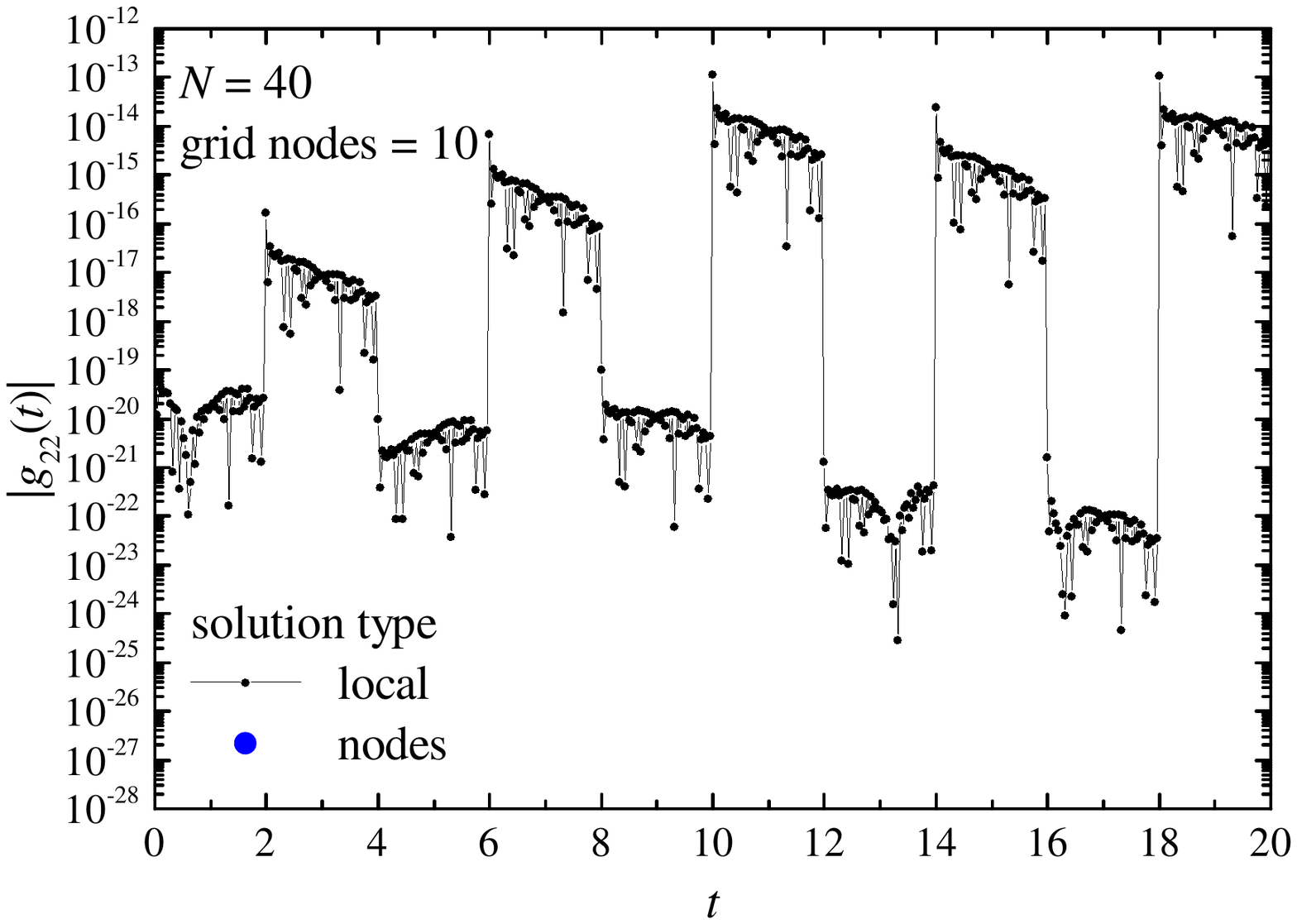}
\vspace{-8mm}\caption{\label{fig:dpend_ind2_sols_vg:f2}}
\end{subfigure}
\begin{subfigure}{0.240\textwidth}
\includegraphics[width=\textwidth]{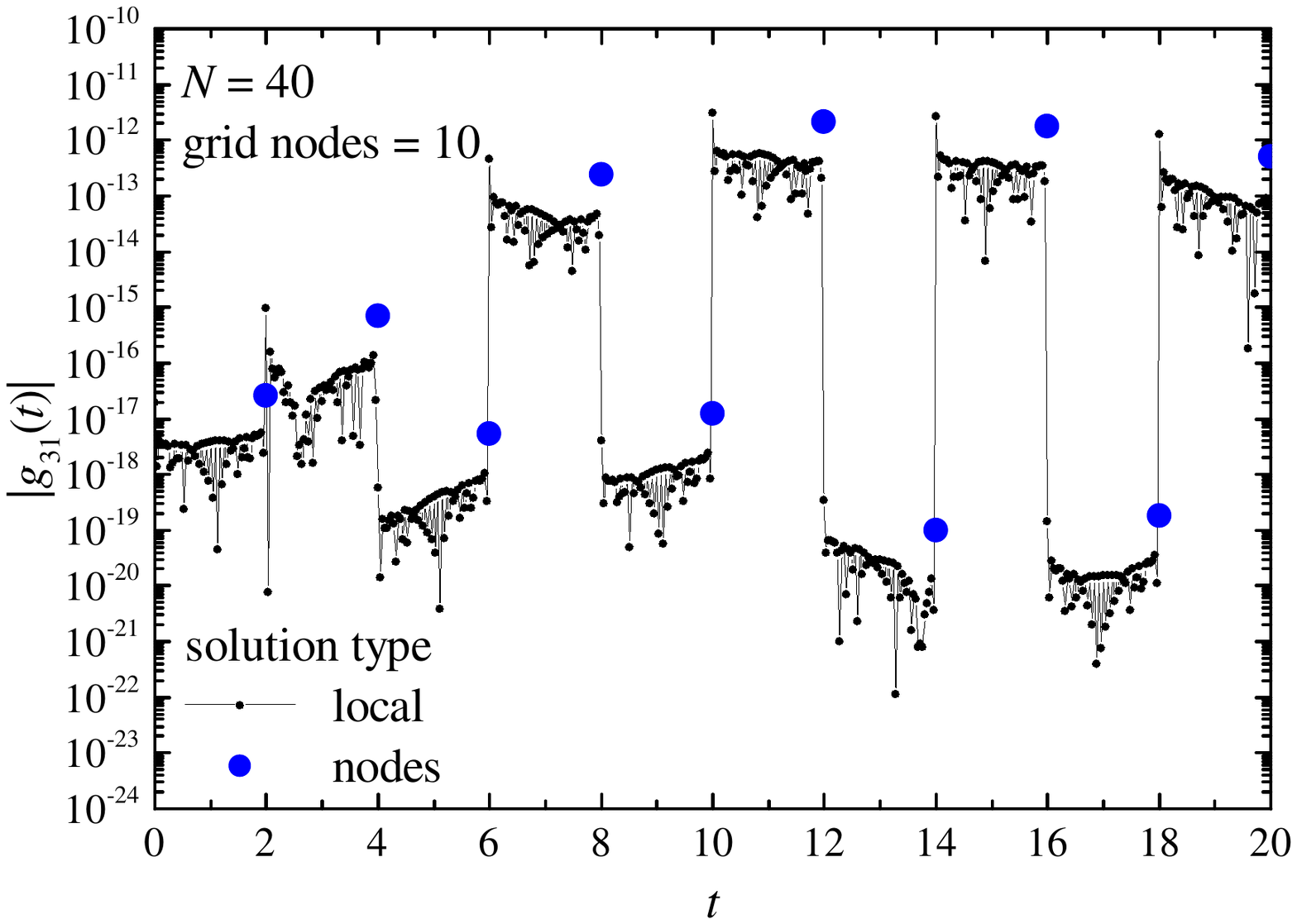}
\vspace{-8mm}\caption{\label{fig:dpend_ind2_sols_vg:f3}}
\end{subfigure}
\begin{subfigure}{0.240\textwidth}
\includegraphics[width=\textwidth]{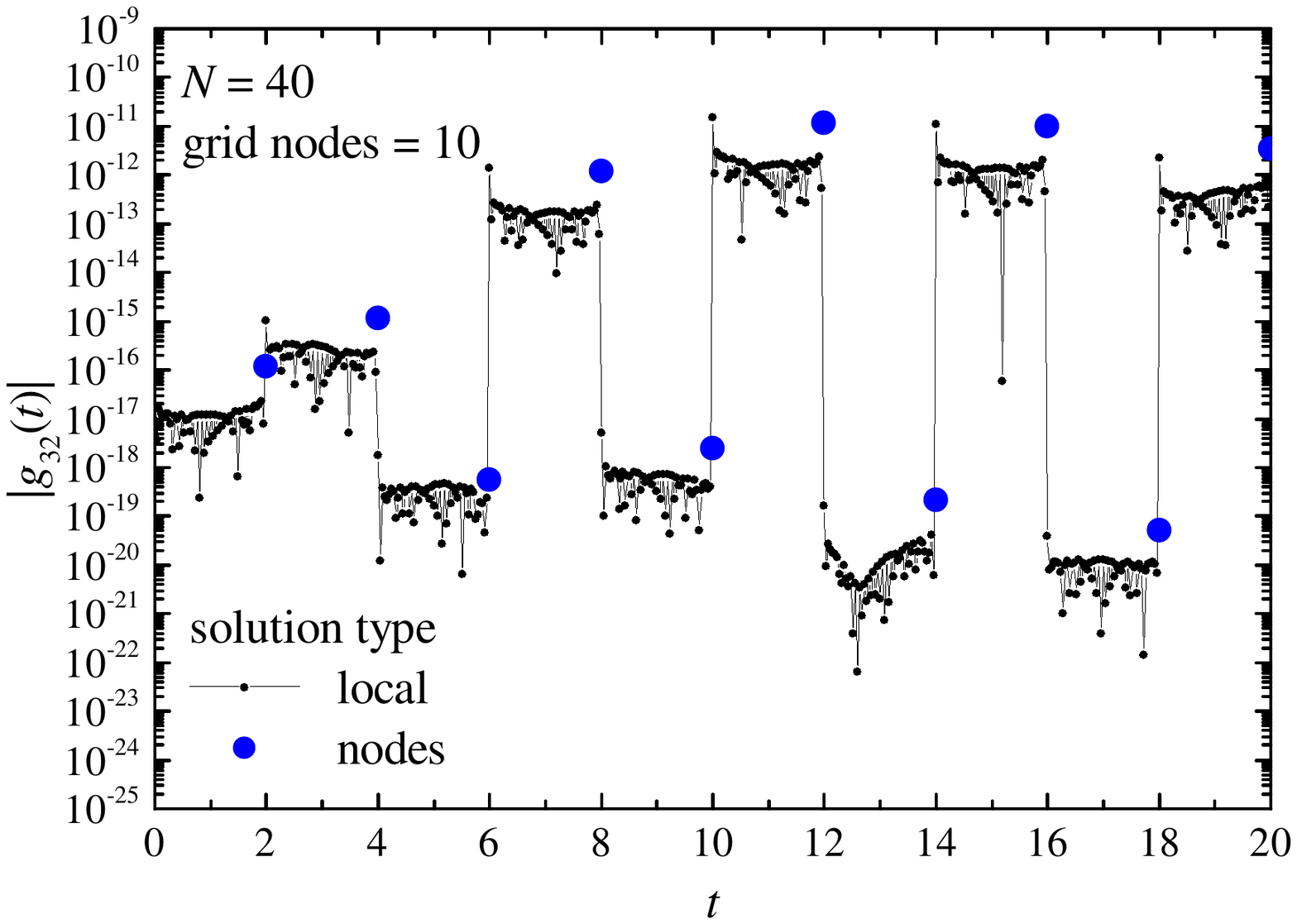}
\vspace{-8mm}\caption{\label{fig:dpend_ind2_sols_vg:f4}}
\end{subfigure}\\[2mm]
\caption{%
Numerical solution of the DAE system (\ref{eq:math_dpend_dae_ind_3}) of index 2. Comparison of the solution at nodes $\mathbf{v}_{n}$, the local solution $\mathbf{v}_{L}(t)$ and the reference solution $\mathbf{v}^{\rm ref}(t)$ for components $v_{1}$ (\subref{fig:dpend_ind2_sols_vg:a1}, \subref{fig:dpend_ind2_sols_vg:c1}, \subref{fig:dpend_ind2_sols_vg:e1}), $v_{2}$ (\subref{fig:dpend_ind2_sols_vg:a2}, \subref{fig:dpend_ind2_sols_vg:c2}, \subref{fig:dpend_ind2_sols_vg:e2}), quantitative satisfiability of the conditions $g_{11} = 0$ (\subref{fig:dpend_ind2_sols_vg:a3}, \subref{fig:dpend_ind2_sols_vg:c3}, \subref{fig:dpend_ind2_sols_vg:e3}), $g_{12} = 0$ (\subref{fig:dpend_ind2_sols_vg:a4}, \subref{fig:dpend_ind2_sols_vg:c4}, \subref{fig:dpend_ind2_sols_vg:e4}), $g_{21} = 0$ (\subref{fig:dpend_ind2_sols_vg:b1}, \subref{fig:dpend_ind2_sols_vg:d1}, \subref{fig:dpend_ind2_sols_vg:f1}), $g_{22} = 0$ (\subref{fig:dpend_ind2_sols_vg:b2}, \subref{fig:dpend_ind2_sols_vg:d2}, \subref{fig:dpend_ind2_sols_vg:f2}), $g_{31} = 0$ (\subref{fig:dpend_ind2_sols_vg:b3}, \subref{fig:dpend_ind2_sols_vg:d3}, \subref{fig:dpend_ind2_sols_vg:f3}), $g_{32} = 0$ (\subref{fig:dpend_ind2_sols_vg:b4}, \subref{fig:dpend_ind2_sols_vg:d4}, \subref{fig:dpend_ind2_sols_vg:f4}), obtained using polynomials with degrees $N = 1$ (\subref{fig:dpend_ind2_sols_vg:a1}, \subref{fig:dpend_ind2_sols_vg:a2}, \subref{fig:dpend_ind2_sols_vg:a3}, \subref{fig:dpend_ind2_sols_vg:a4}, \subref{fig:dpend_ind2_sols_vg:b1}, \subref{fig:dpend_ind2_sols_vg:b2}, \subref{fig:dpend_ind2_sols_vg:b3}, \subref{fig:dpend_ind2_sols_vg:b4}), $N = 8$ (\subref{fig:dpend_ind2_sols_vg:c1}, \subref{fig:dpend_ind2_sols_vg:c2}, \subref{fig:dpend_ind2_sols_vg:c3}, \subref{fig:dpend_ind2_sols_vg:c4}, \subref{fig:dpend_ind2_sols_vg:d1}, \subref{fig:dpend_ind2_sols_vg:d2}, \subref{fig:dpend_ind2_sols_vg:d3}, \subref{fig:dpend_ind2_sols_vg:d4}) and $N = 40$ (\subref{fig:dpend_ind2_sols_vg:e1}, \subref{fig:dpend_ind2_sols_vg:e2}, \subref{fig:dpend_ind2_sols_vg:e3}, \subref{fig:dpend_ind2_sols_vg:e4}, \subref{fig:dpend_ind2_sols_vg:f1}, \subref{fig:dpend_ind2_sols_vg:f2}, \subref{fig:dpend_ind2_sols_vg:f3}, \subref{fig:dpend_ind2_sols_vg:f4}).
}
\label{fig:dpend_ind2_sols_vg}
\end{figure} 

\begin{figure}[h!]
\captionsetup[subfigure]{%
	position=bottom,
	font+=smaller,
	textfont=normalfont,
	singlelinecheck=off,
	justification=raggedright
}
\centering
\begin{subfigure}{0.320\textwidth}
\includegraphics[width=\textwidth]{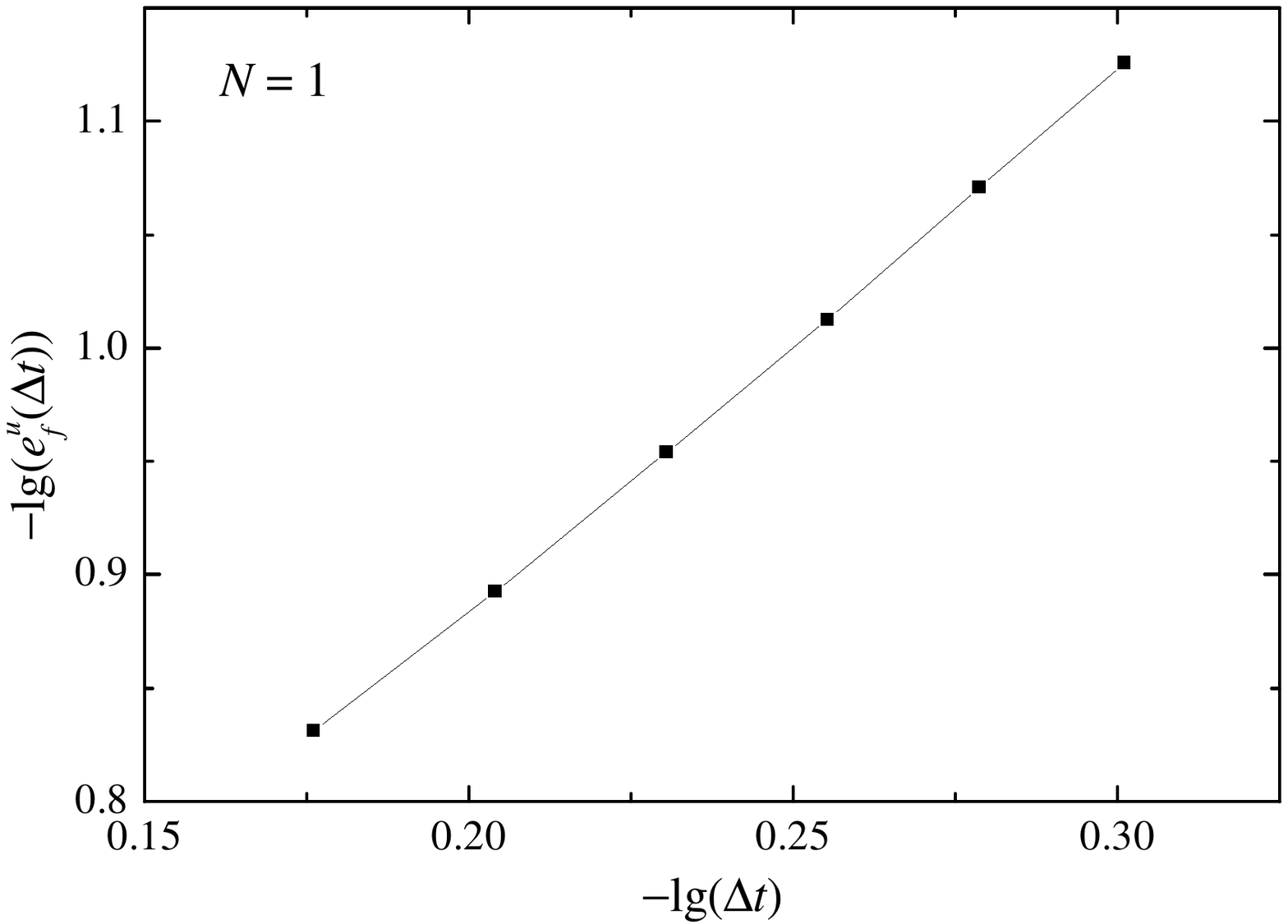}
\vspace{-8mm}\caption{\label{fig:dpend_ind2_errors:a1}}
\end{subfigure}
\begin{subfigure}{0.320\textwidth}
\includegraphics[width=\textwidth]{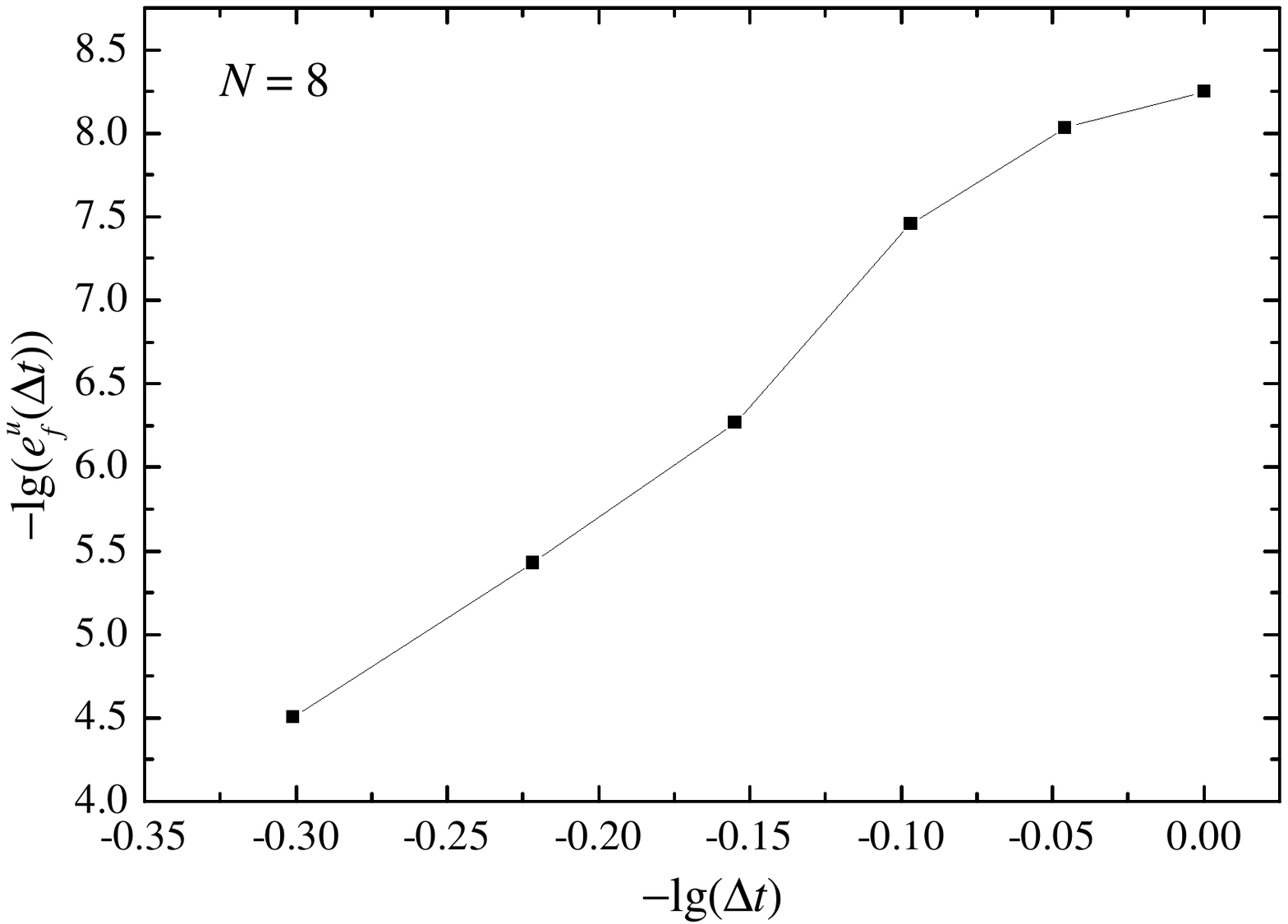}
\vspace{-8mm}\caption{\label{fig:dpend_ind2_errors:a2}}
\end{subfigure}
\begin{subfigure}{0.320\textwidth}
\includegraphics[width=\textwidth]{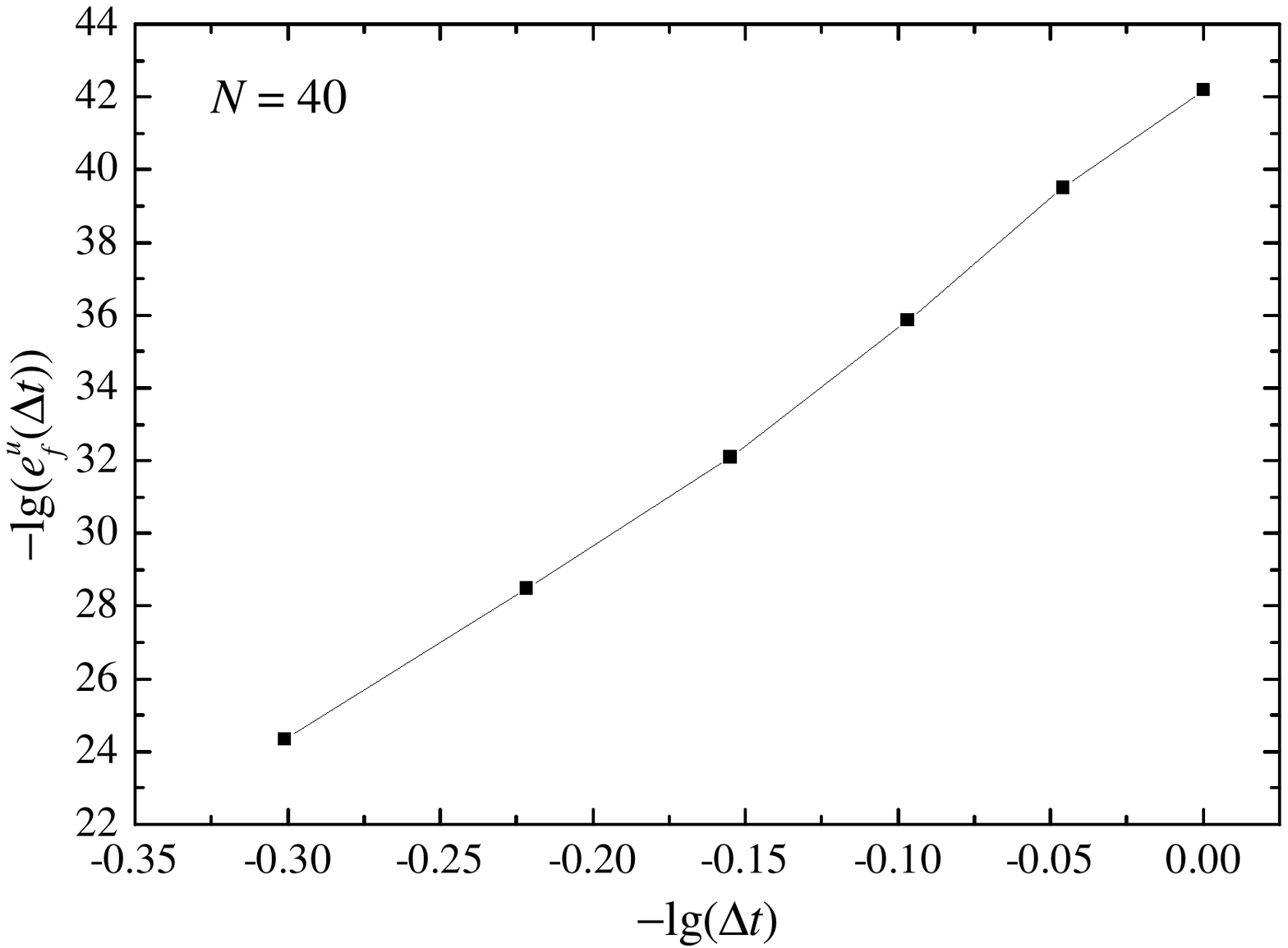}
\vspace{-8mm}\caption{\label{fig:dpend_ind2_errors:a3}}
\end{subfigure}\\[2mm]
\begin{subfigure}{0.320\textwidth}
\includegraphics[width=\textwidth]{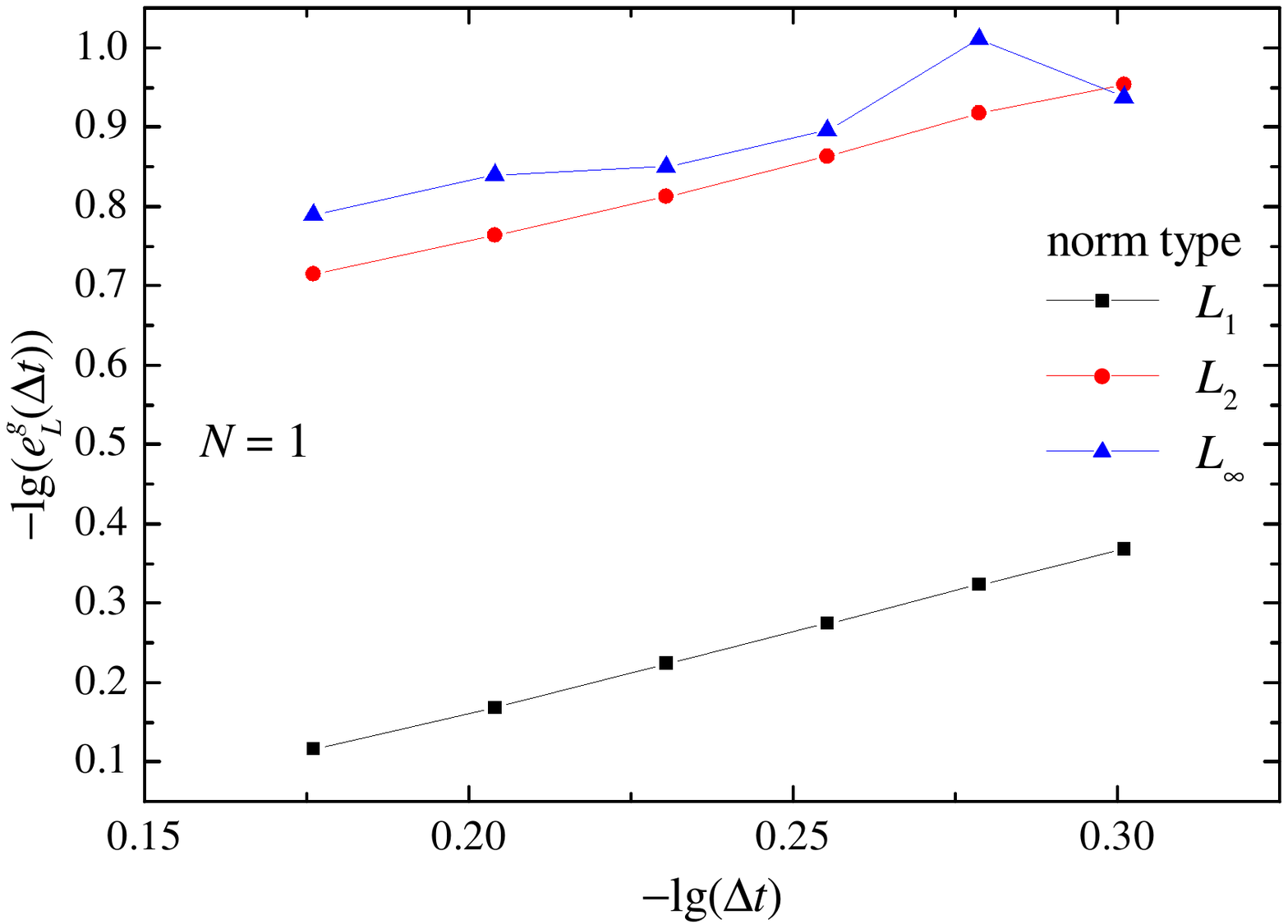}
\vspace{-8mm}\caption{\label{fig:dpend_ind2_errors:b1}}
\end{subfigure}
\begin{subfigure}{0.320\textwidth}
\includegraphics[width=\textwidth]{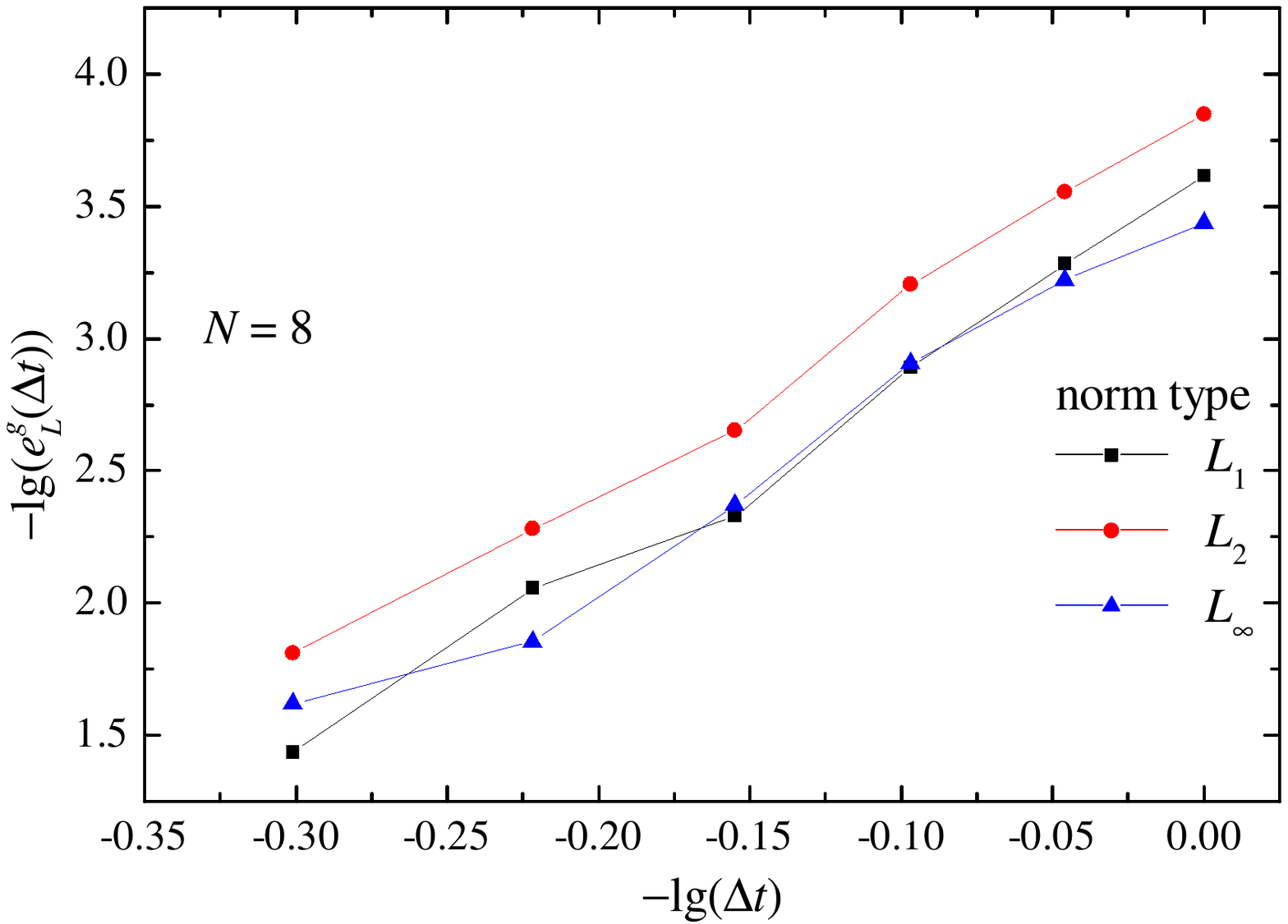}
\vspace{-8mm}\caption{\label{fig:dpend_ind2_errors:b2}}
\end{subfigure}
\begin{subfigure}{0.320\textwidth}
\includegraphics[width=\textwidth]{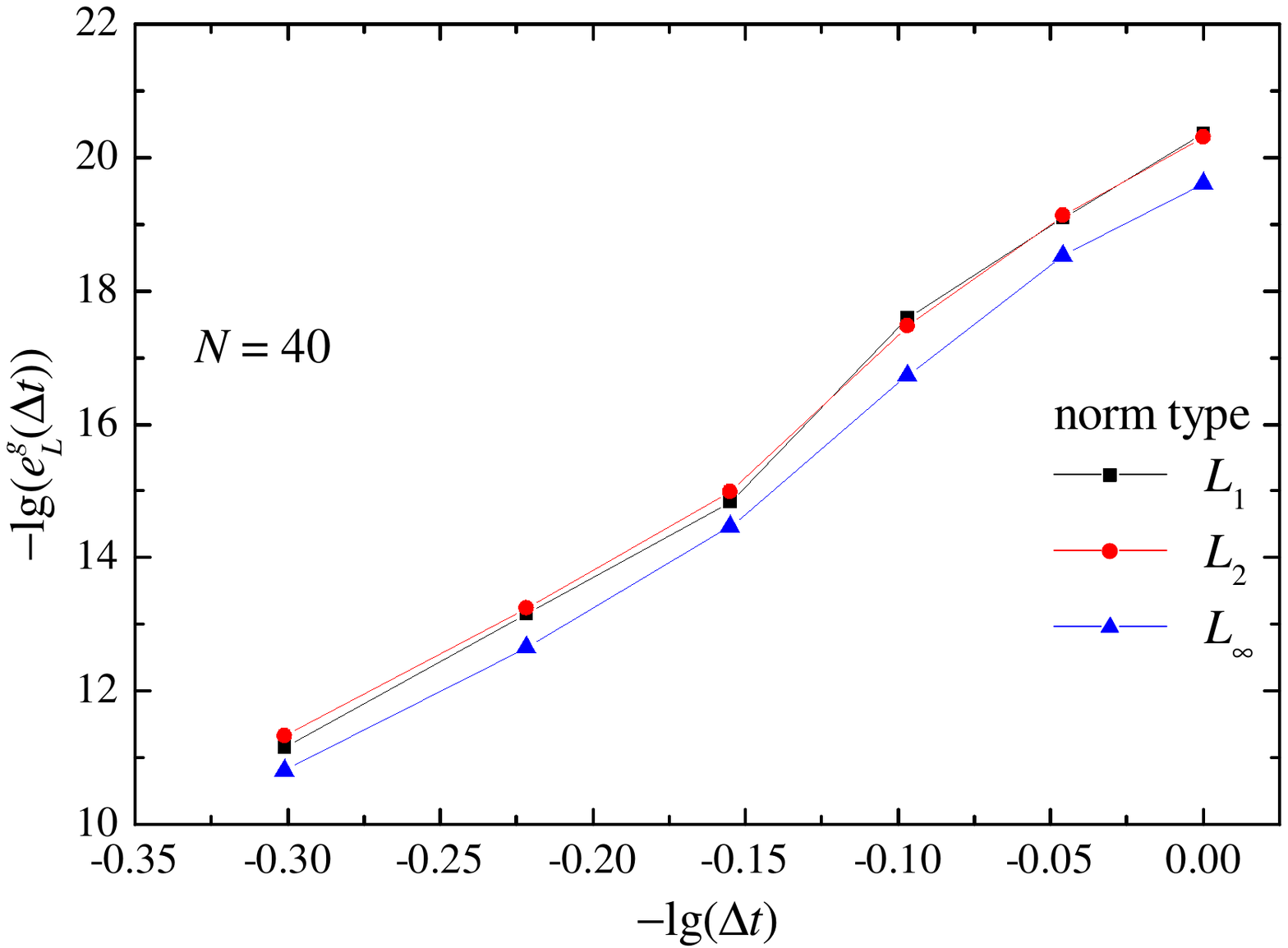}
\vspace{-8mm}\caption{\label{fig:dpend_ind2_errors:b3}}
\end{subfigure}\\[2mm]
\begin{subfigure}{0.320\textwidth}
\includegraphics[width=\textwidth]{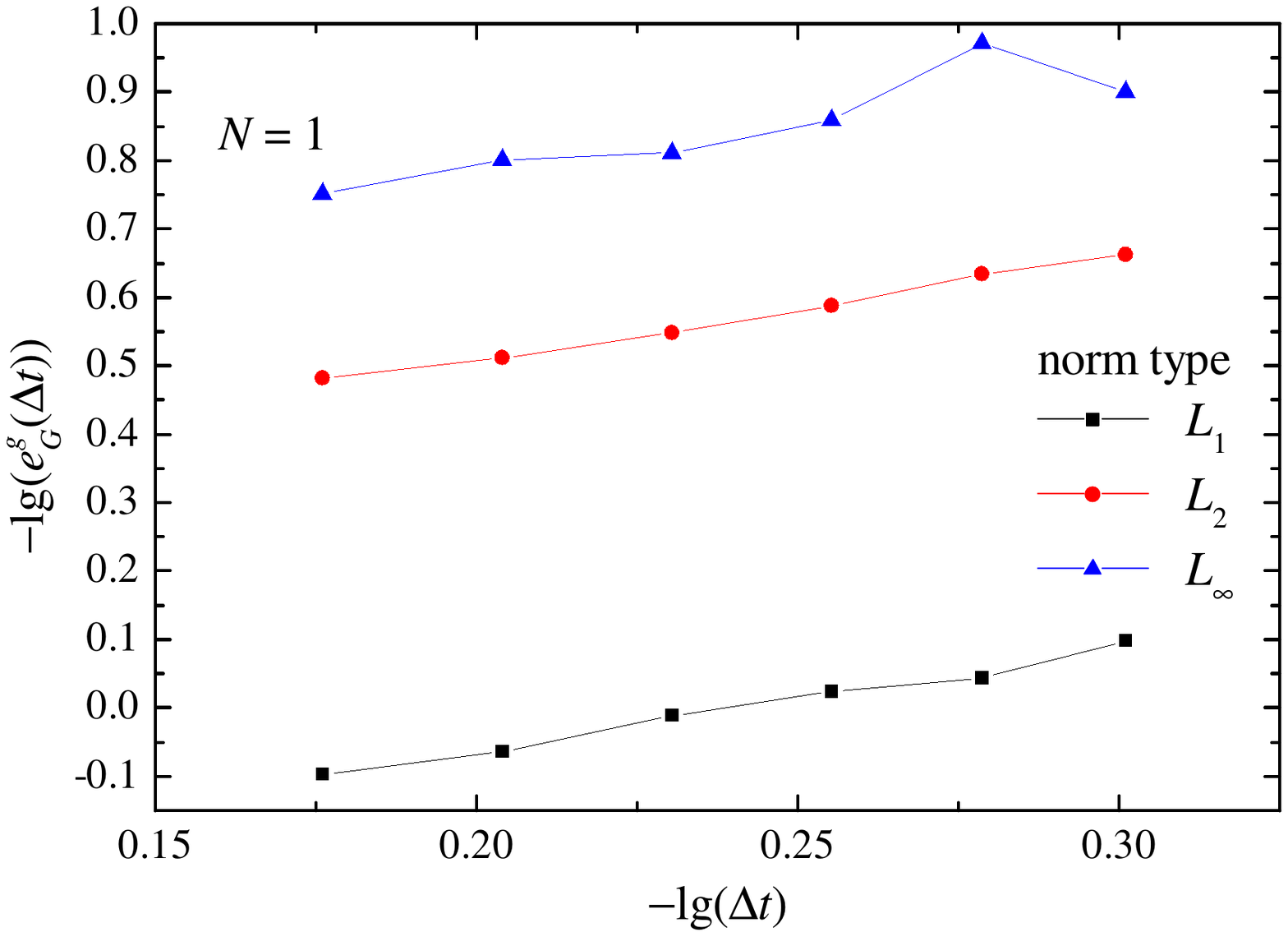}
\vspace{-8mm}\caption{\label{fig:dpend_ind2_errors:c1}}
\end{subfigure}
\begin{subfigure}{0.320\textwidth}
\includegraphics[width=\textwidth]{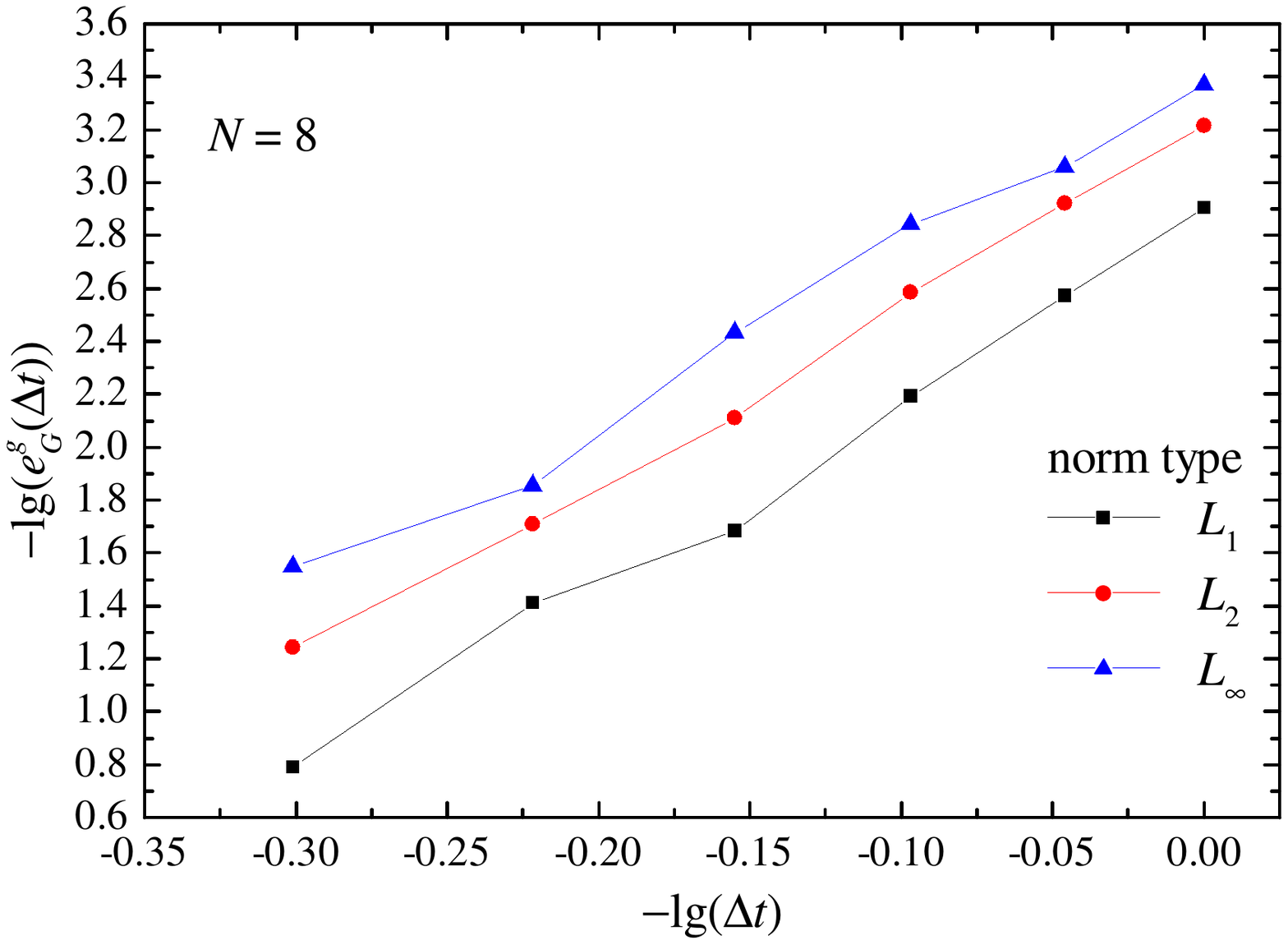}
\vspace{-8mm}\caption{\label{fig:dpend_ind2_errors:c2}}
\end{subfigure}
\begin{subfigure}{0.320\textwidth}
\includegraphics[width=\textwidth]{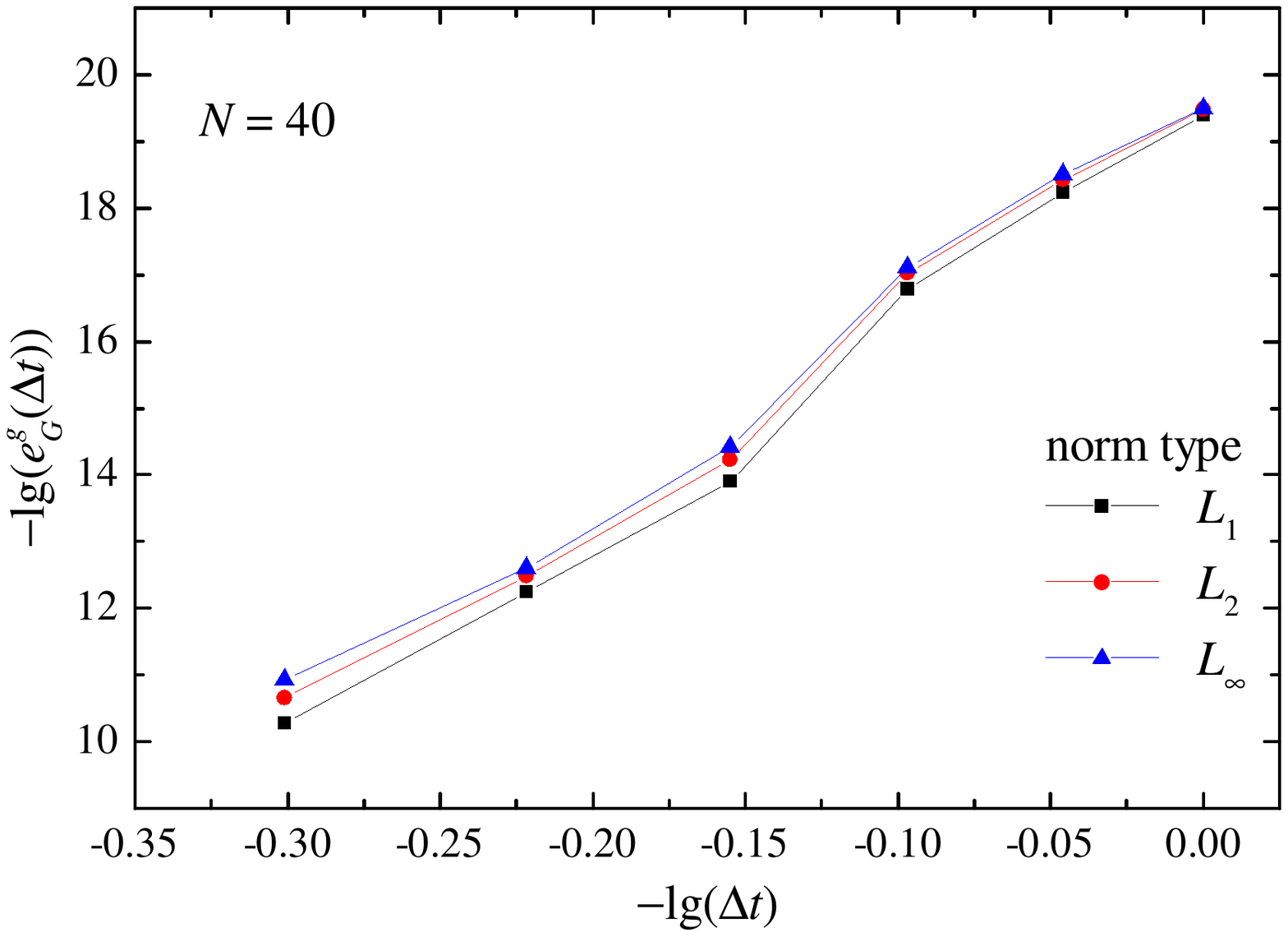}
\vspace{-8mm}\caption{\label{fig:dpend_ind2_errors:c3}}
\end{subfigure}\\[2mm]
\caption{%
Log-log plot of the dependence of the error for the solution at nodes at final time $t_{f}$ $e_{f}^{u}$ (\subref{fig:dpend_ind2_errors:a1}, \subref{fig:dpend_ind2_errors:a2}, \subref{fig:dpend_ind2_errors:a3}), global errors for the local solution $e_{L}^{g}$ (\subref{fig:dpend_ind2_errors:b1}, \subref{fig:dpend_ind2_errors:b2}, \subref{fig:dpend_ind2_errors:b3}) and the solution at nodes $e_{G}^{u}$ (\subref{fig:dpend_ind2_errors:c1}, \subref{fig:dpend_ind2_errors:c2}, \subref{fig:dpend_ind2_errors:c3}), by numerical solution of the DAE system (\ref{eq:math_dpend_dae_ind_3}) of index 2 obtained using polynomials with degrees $N = 1$ (\subref{fig:dpend_ind2_errors:a1}, \subref{fig:dpend_ind2_errors:b1}, \subref{fig:dpend_ind2_errors:c1}), $N = 8$ (\subref{fig:dpend_ind2_errors:a2}, \subref{fig:dpend_ind2_errors:b2}, \subref{fig:dpend_ind2_errors:c2}) and $N = 40$ (\subref{fig:dpend_ind2_errors:a3}, \subref{fig:dpend_ind2_errors:b3}, \subref{fig:dpend_ind2_errors:c3}).
}
\label{fig:dpend_ind2_errors}
\end{figure} 

\begin{table*}[h!]
\centering
\caption{%
Convergence orders $p_{f}$, $p_{L_{1}}$, $p_{L_{2}}$, $p_{L_{\infty}}$, calculated on final time step $t_{f}$ and in norms $L_{1}$, $L_{2}$, $L_{\infty}$ of the ADER-DG method for the DAE system (\ref{eq:math_dpend_dae_ind_3}) of index 2; $N$ is the degree of the basis polynomials $\varphi_{p}$. Orders $p^{n, u}$ are calculated for solution $\mathbf{u}_{n}$; orders $p^{n, g}$ and $p^{l, g}$ --- for the conditions $\mathbf{g} = 0$ on the numerical solution at the nodes $(\mathbf{u}_{n}, \mathbf{v}_{n})$ and on the local solution $(\mathbf{u}_{L}, \mathbf{v}_{L})$. The theoretical values of convergence order $p_{\rm th.}^{n} = 2N+1$ and $p_{\rm th.}^{l} = N+1$ are applicable for the ADER-DG method for ODE problems, and are presented for comparison.
}
\label{tab:conv_orders_dpend_ind2}
\begin{tabular}{@{}|l|l|lll|c|lll|c|@{}}
\toprule
$N$ & $p_{f}^{n, u}$ & $p_{L_{1}}^{n, g}$ & $p_{L_{2}}^{n, g}$ & $p_{L_{\infty}}^{n, g}$ & $p_{\rm th.}^{n}$ & $p_{L_{1}}^{l, g}$ & $p_{L_{2}}^{l, g}$ & $p_{L_{\infty}}^{l, g}$ & $p_{\rm th.}^{l}$ \\
\midrule
$1$	&	$2.36$	&	$1.52$	&	$1.49$	&	$1.50$	&	$3$	&	$2.03$	&	$1.95$	&	$1.49$	&	$2$\\
$2$	&	$4.66$	&	$2.83$	&	$2.95$	&	$3.05$	&	$5$	&	$2.93$	&	$2.98$	&	$3.06$	&	$3$\\
$3$	&	$6.63$	&	$3.57$	&	$3.52$	&	$3.67$	&	$7$	&	$3.66$	&	$3.62$	&	$3.66$	&	$4$\\
$4$	&	$8.89$	&	$5.04$	&	$4.92$	&	$4.73$	&	$9$	&	$5.03$	&	$4.96$	&	$4.77$	&	$5$\\
$5$	&	$11.13$	&	$5.64$	&	$5.60$	&	$5.38$	&	$11$	&	$5.72$	&	$5.70$	&	$6.01$	&	$6$\\
$6$	&	$10.03$	&	$5.55$	&	$5.19$	&	$4.86$	&	$13$	&	$5.81$	&	$5.51$	&	$5.10$	&	$7$\\
$7$	&	$12.62$	&	$7.10$	&	$6.73$	&	$6.22$	&	$15$	&	$7.21$	&	$6.92$	&	$6.42$	&	$8$\\
$8$	&	$13.30$	&	$6.94$	&	$6.66$	&	$6.27$	&	$17$	&	$7.21$	&	$6.92$	&	$6.52$	&	$9$\\
$9$	&	$15.52$	&	$8.32$	&	$8.13$	&	$7.78$	&	$19$	&	$8.66$	&	$8.47$	&	$8.30$	&	$10$\\
$10$	&	$17.14$	&	$8.49$	&	$8.13$	&	$7.71$	&	$21$	&	$8.74$	&	$8.39$	&	$7.94$	&	$11$\\
$11$	&	$17.61$	&	$9.63$	&	$9.31$	&	$8.79$	&	$23$	&	$9.89$	&	$9.66$	&	$9.47$	&	$12$\\
$12$	&	$20.05$	&	$10.14$	&	$9.68$	&	$8.98$	&	$25$	&	$10.42$	&	$10.01$	&	$9.34$	&	$13$\\
$13$	&	$19.96$	&	$11.01$	&	$10.55$	&	$9.86$	&	$27$	&	$11.24$	&	$10.80$	&	$10.25$	&	$14$\\
$14$	&	$22.67$	&	$12.01$	&	$11.27$	&	$10.41$	&	$29$	&	$12.25$	&	$11.64$	&	$11.06$	&	$15$\\
$15$	&	$23.02$	&	$12.38$	&	$11.87$	&	$11.27$	&	$31$	&	$12.57$	&	$12.06$	&	$11.68$	&	$16$\\
$16$	&	$25.02$	&	$13.14$	&	$12.71$	&	$12.26$	&	$33$	&	$13.46$	&	$13.10$	&	$12.74$	&	$17$\\
$17$	&	$26.28$	&	$13.84$	&	$13.34$	&	$12.74$	&	$35$	&	$14.09$	&	$13.53$	&	$12.86$	&	$18$\\
$18$	&	$26.67$	&	$14.42$	&	$14.05$	&	$13.85$	&	$37$	&	$14.75$	&	$14.41$	&	$14.15$	&	$19$\\
$19$	&	$29.17$	&	$15.53$	&	$14.91$	&	$14.07$	&	$39$	&	$15.75$	&	$15.13$	&	$14.34$	&	$20$\\
$20$	&	$29.33$	&	$16.18$	&	$15.48$	&	$14.66$	&	$41$	&	$16.38$	&	$15.72$	&	$14.77$	&	$21$\\
\midrule
$25$	&	$37.03$	&	$19.91$	&	$19.32$	&	$18.86$	&	$51$	&	$20.22$	&	$19.65$	&	$19.21$	&	$26$\\
$30$	&	$44.45$	&	$23.64$	&	$22.87$	&	$22.13$	&	$61$	&	$23.97$	&	$23.25$	&	$22.56$	&	$31$\\
$35$	&	$52.57$	&	$27.78$	&	$26.94$	&	$26.11$	&	$71$	&	$28.01$	&	$27.42$	&	$27.04$	&	$36$\\
$40$	&	$59.99$	&	$31.64$	&	$30.88$	&	$30.10$	&	$81$	&	$31.66$	&	$30.95$	&	$30.40$	&	$41$\\
\bottomrule
\end{tabular}
\end{table*}

\begin{figure}[h!]
\captionsetup[subfigure]{%
	position=bottom,
	font+=smaller,
	textfont=normalfont,
	singlelinecheck=off,
	justification=raggedright
}
\centering
\begin{subfigure}{0.240\textwidth}
\includegraphics[width=\textwidth]{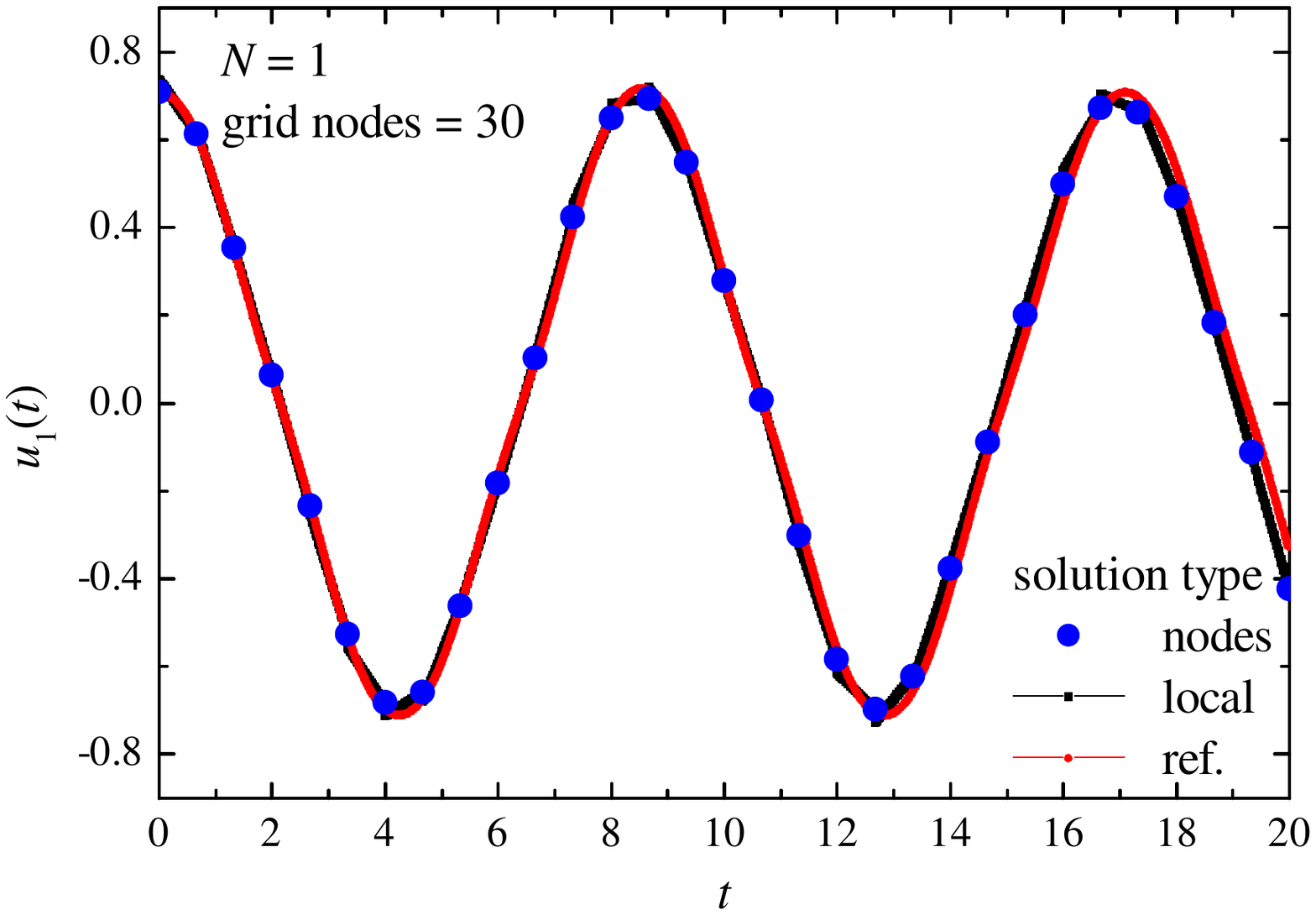}
\vspace{-8mm}\caption{\label{fig:dpend_ind1_sols_u:a1}}
\end{subfigure}
\begin{subfigure}{0.240\textwidth}
\includegraphics[width=\textwidth]{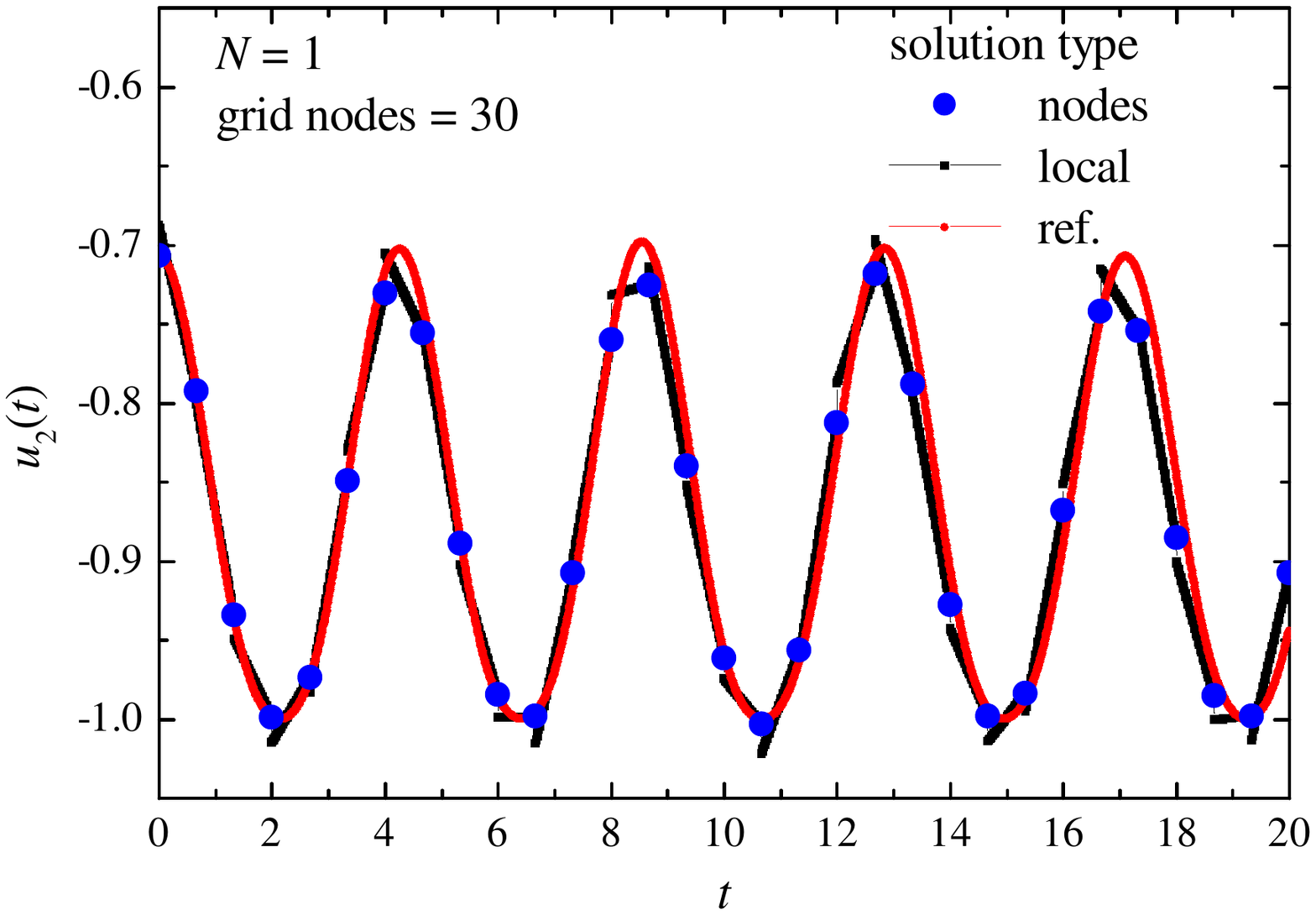}
\vspace{-8mm}\caption{\label{fig:dpend_ind1_sols_u:a2}}
\end{subfigure}
\begin{subfigure}{0.240\textwidth}
\includegraphics[width=\textwidth]{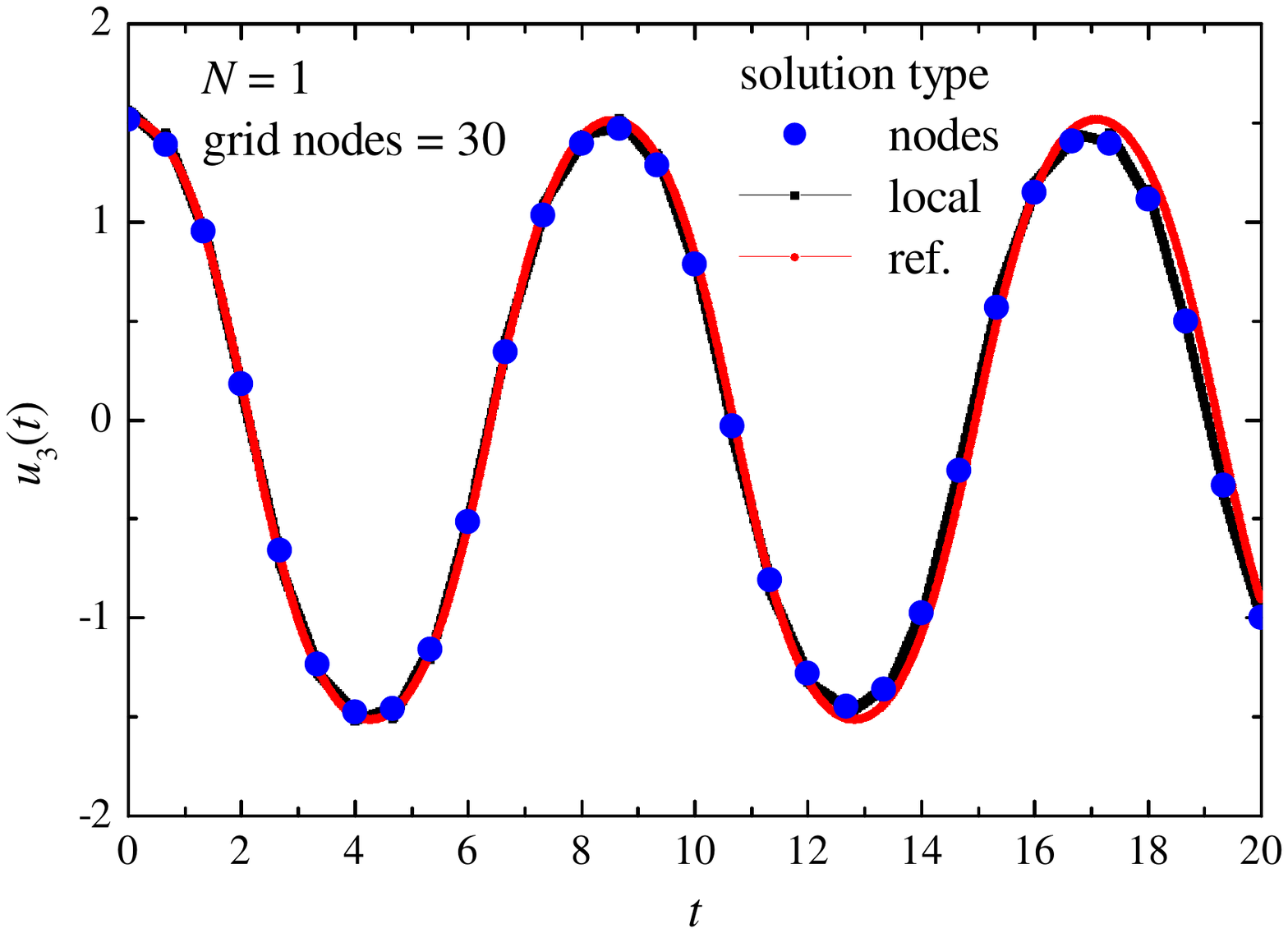}
\vspace{-8mm}\caption{\label{fig:dpend_ind1_sols_u:a3}}
\end{subfigure}
\begin{subfigure}{0.240\textwidth}
\includegraphics[width=\textwidth]{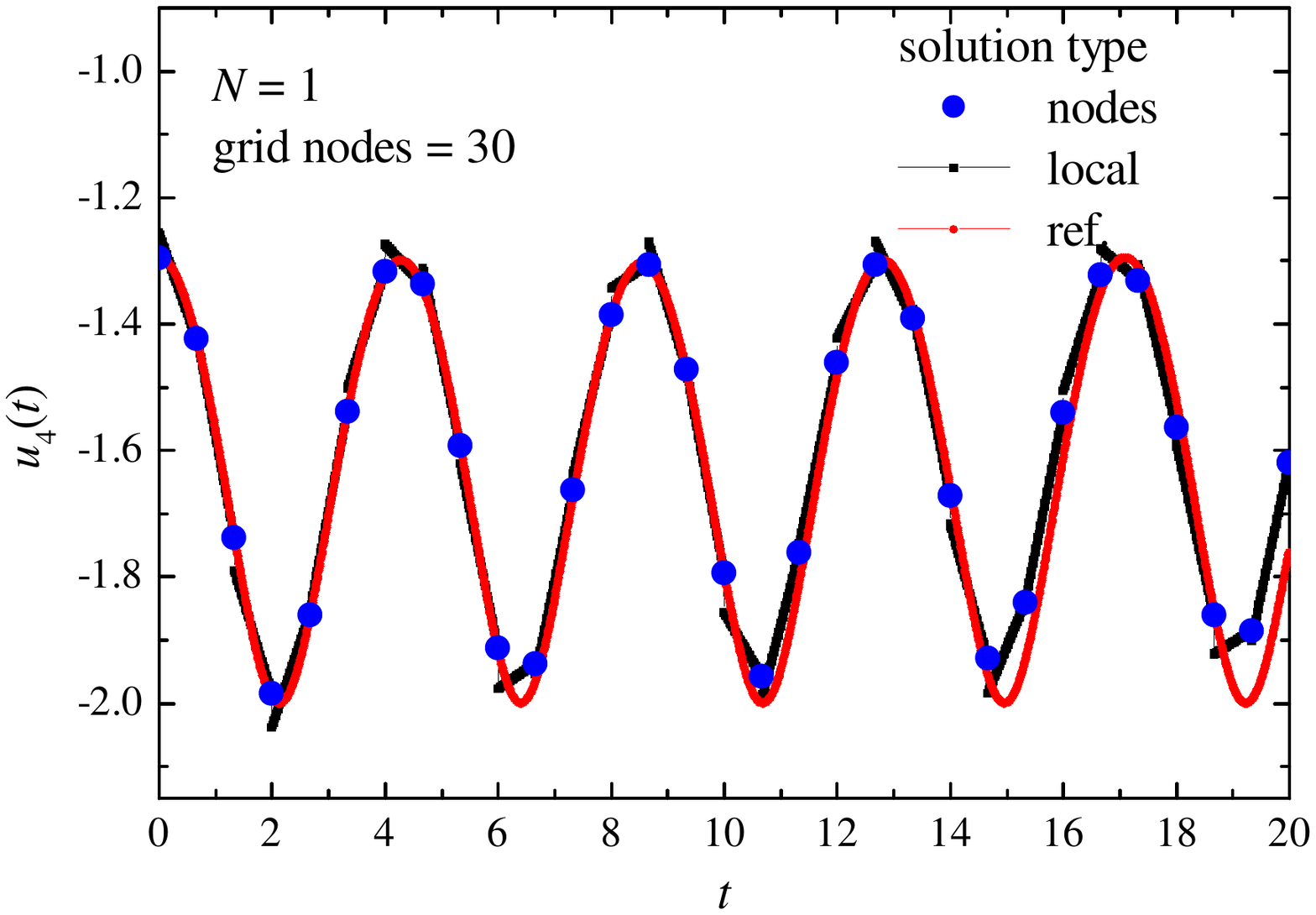}
\vspace{-8mm}\caption{\label{fig:dpend_ind1_sols_u:a4}}
\end{subfigure}\\[2mm]
\begin{subfigure}{0.240\textwidth}
\includegraphics[width=\textwidth]{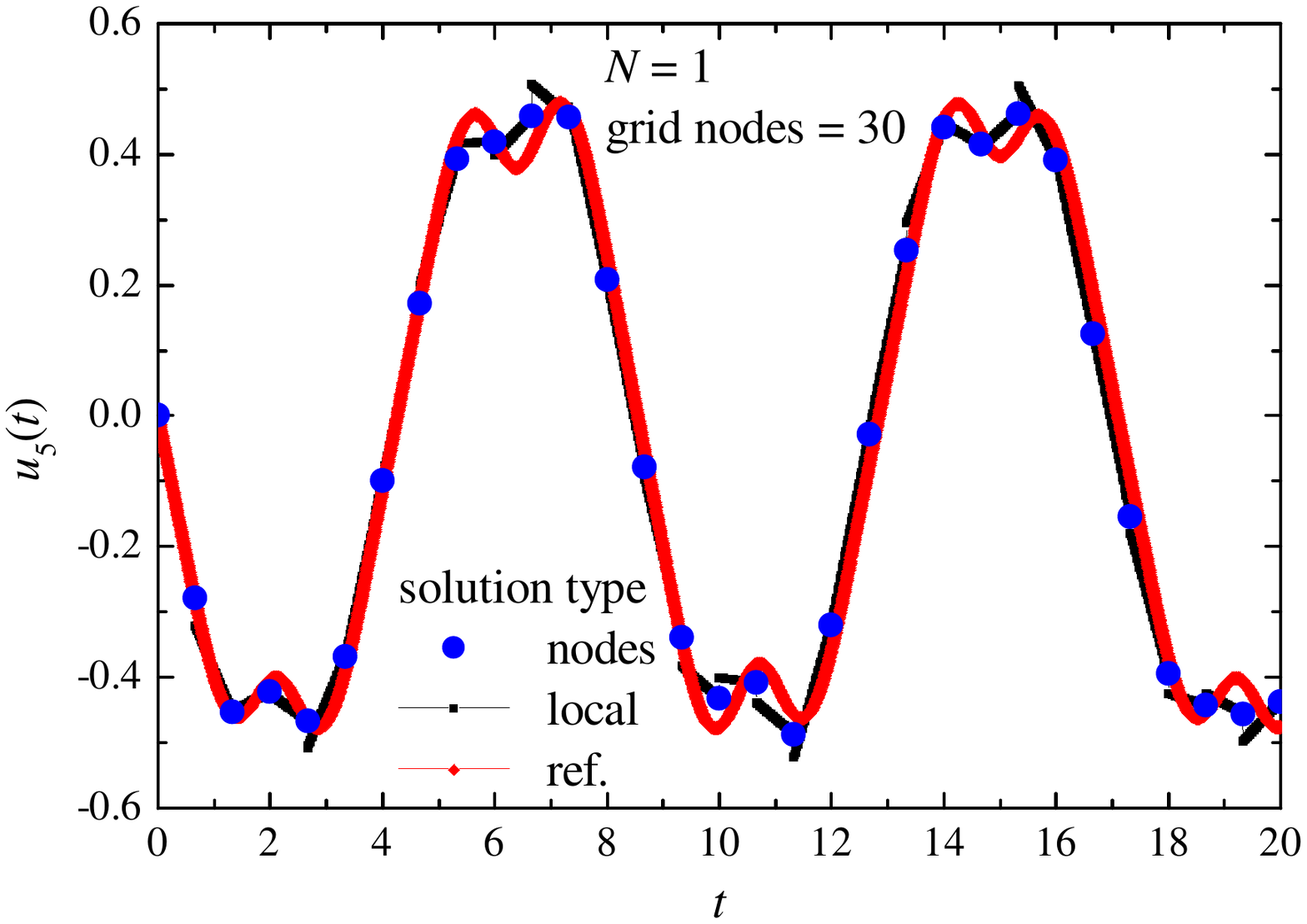}
\vspace{-8mm}\caption{\label{fig:dpend_ind1_sols_u:b1}}
\end{subfigure}
\begin{subfigure}{0.240\textwidth}
\includegraphics[width=\textwidth]{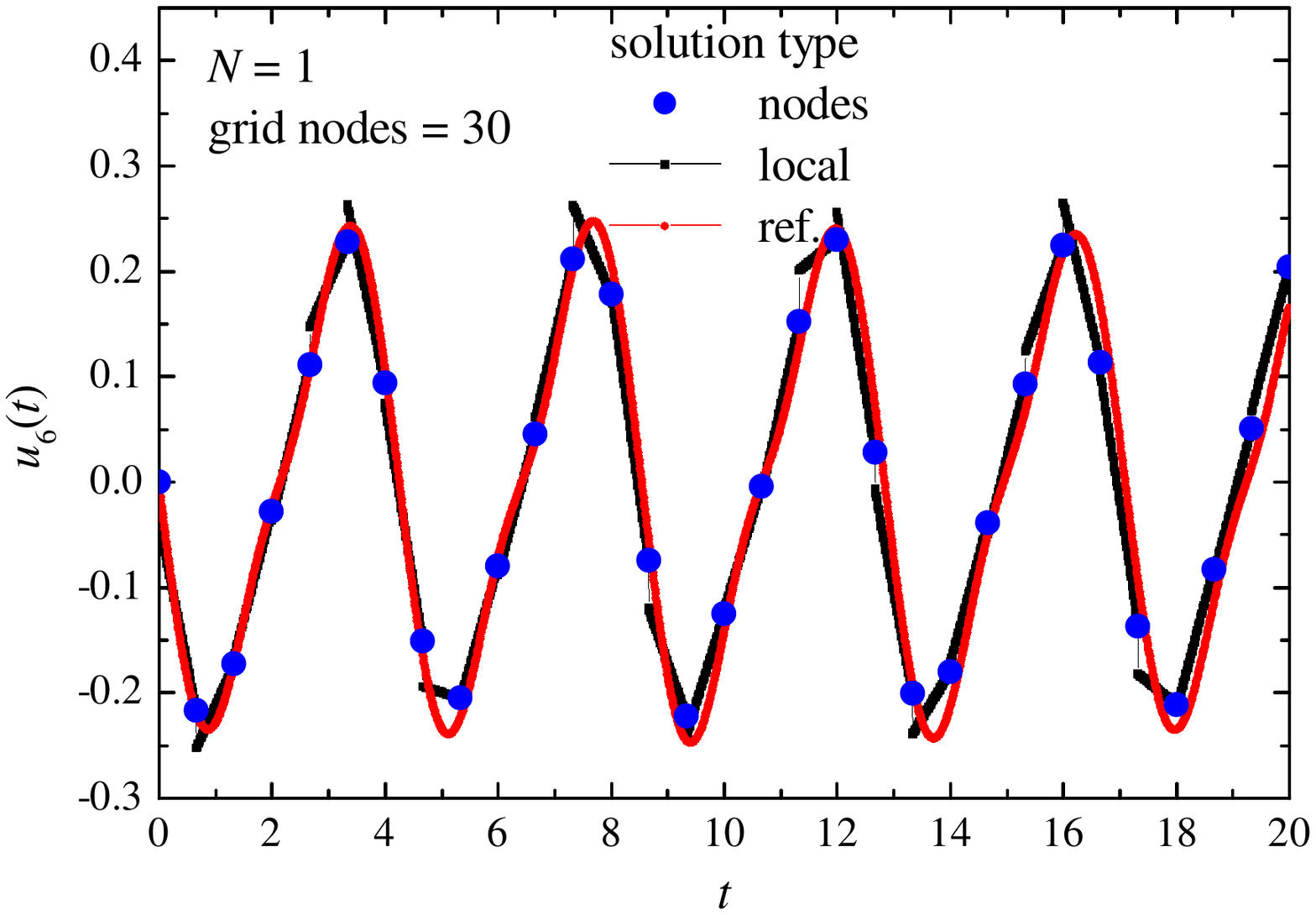}
\vspace{-8mm}\caption{\label{fig:dpend_ind1_sols_u:b2}}
\end{subfigure}
\begin{subfigure}{0.240\textwidth}
\includegraphics[width=\textwidth]{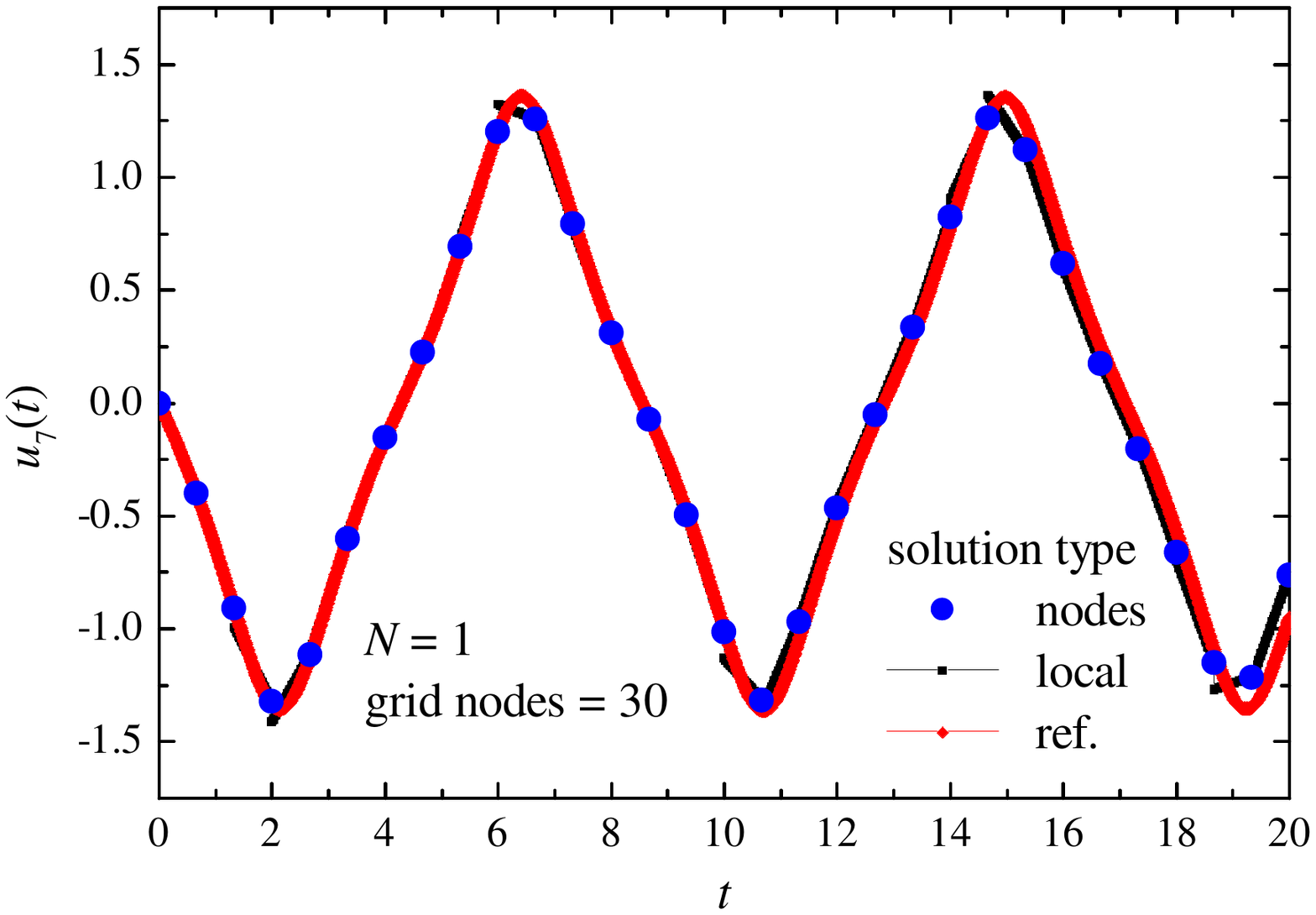}
\vspace{-8mm}\caption{\label{fig:dpend_ind1_sols_u:b3}}
\end{subfigure}
\begin{subfigure}{0.240\textwidth}
\includegraphics[width=\textwidth]{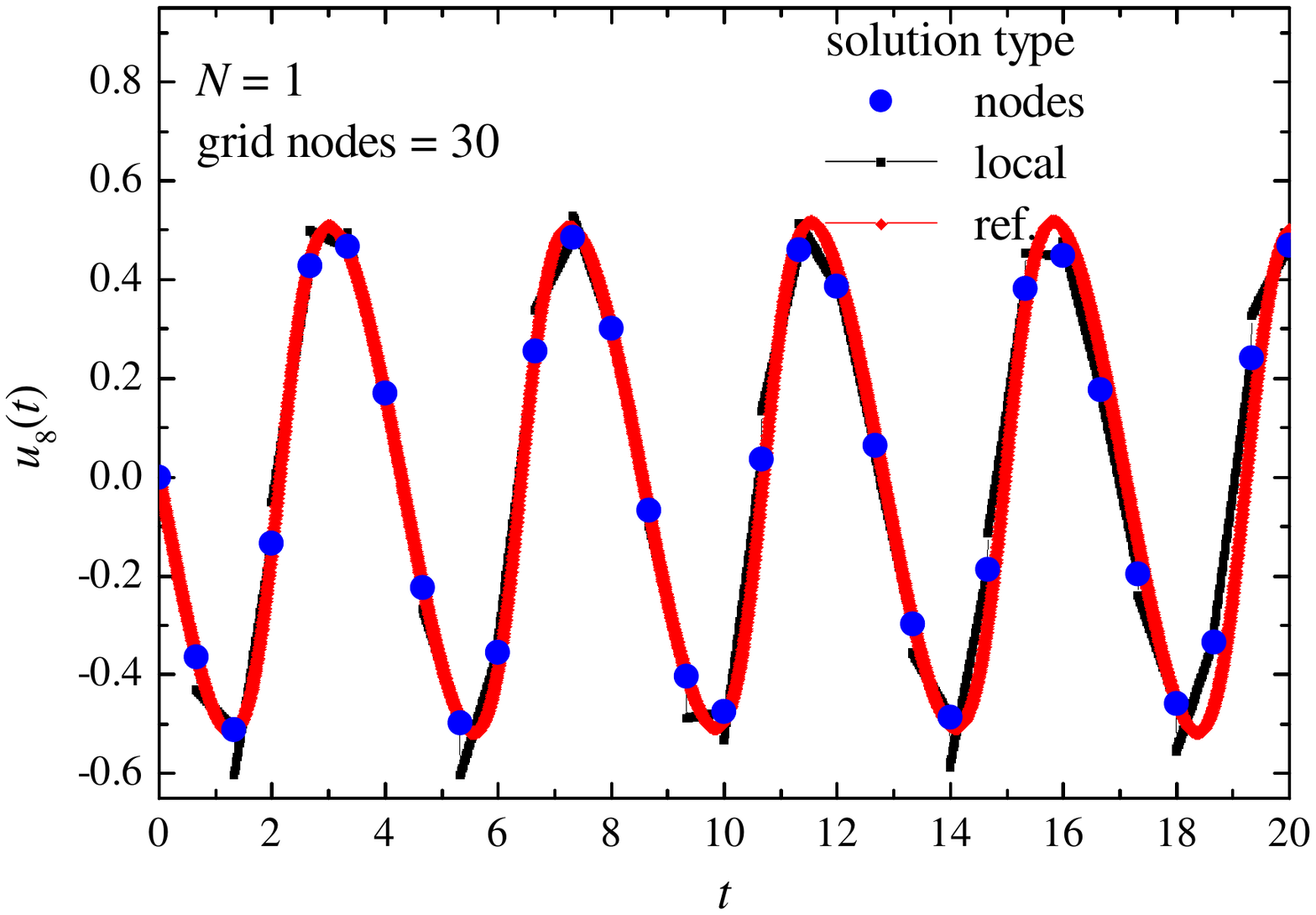}
\vspace{-8mm}\caption{\label{fig:dpend_ind1_sols_u:b4}}
\end{subfigure}\\[2mm]
\begin{subfigure}{0.240\textwidth}
\includegraphics[width=\textwidth]{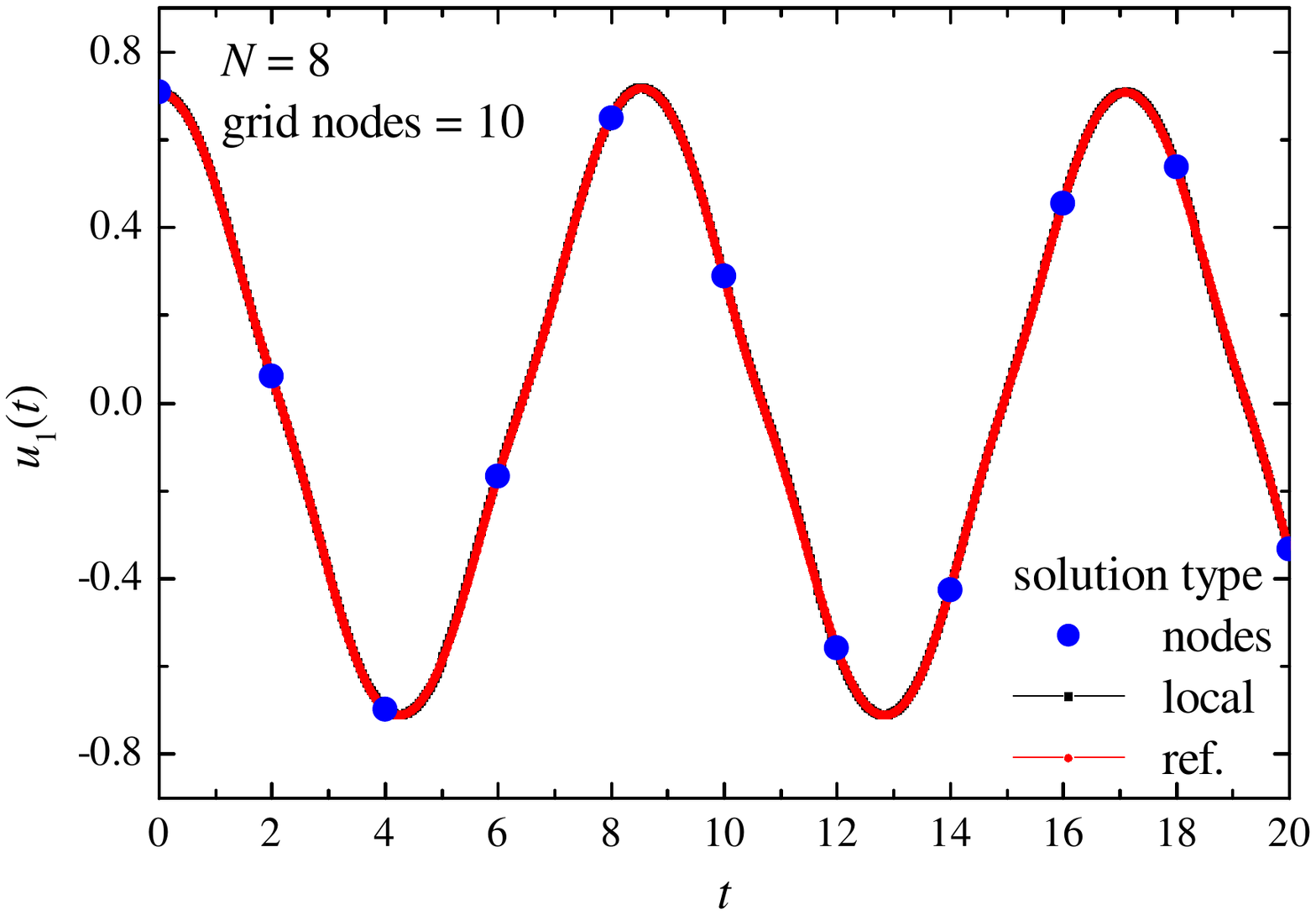}
\vspace{-8mm}\caption{\label{fig:dpend_ind1_sols_u:c1}}
\end{subfigure}
\begin{subfigure}{0.240\textwidth}
\includegraphics[width=\textwidth]{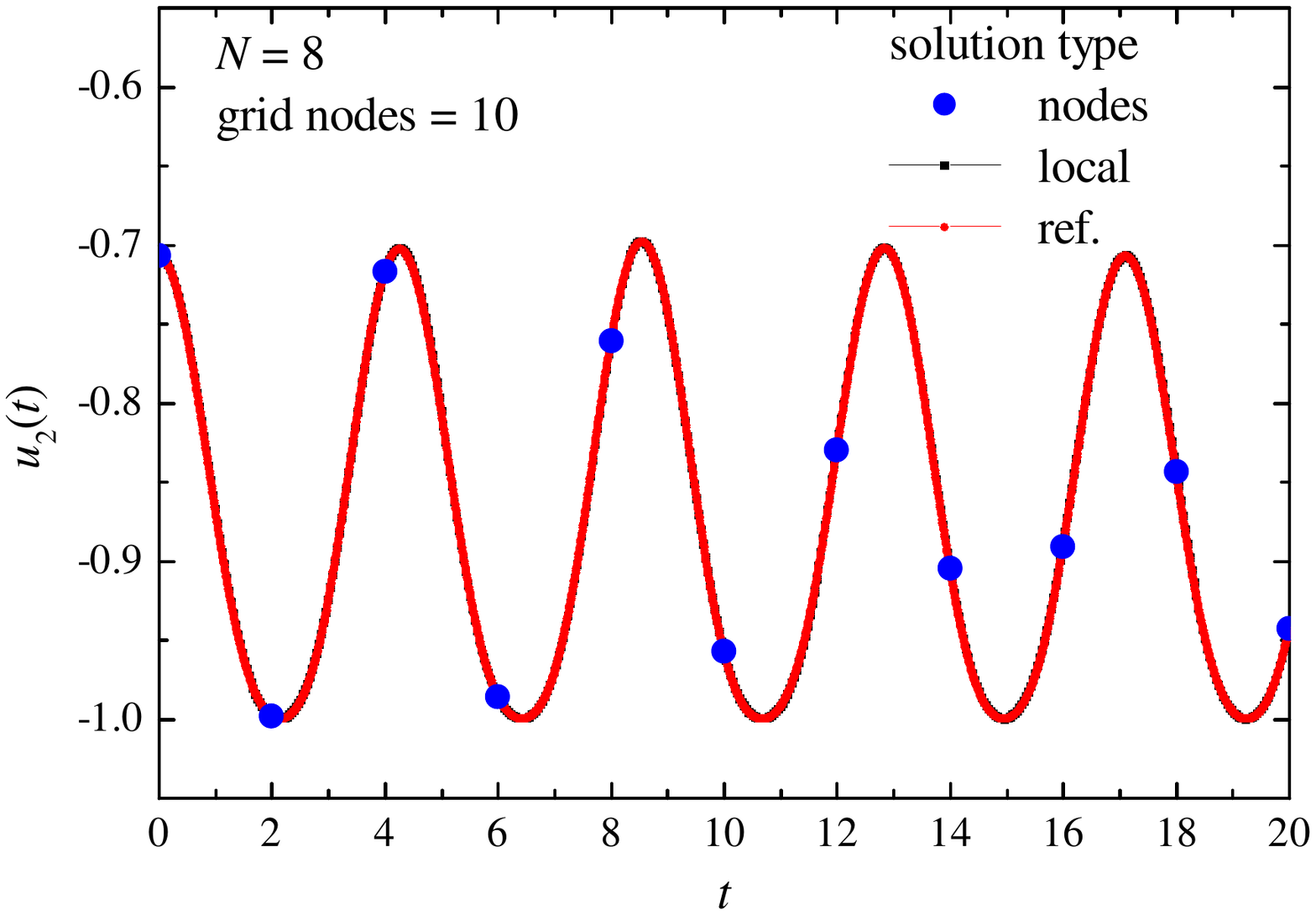}
\vspace{-8mm}\caption{\label{fig:dpend_ind1_sols_u:c2}}
\end{subfigure}
\begin{subfigure}{0.240\textwidth}
\includegraphics[width=\textwidth]{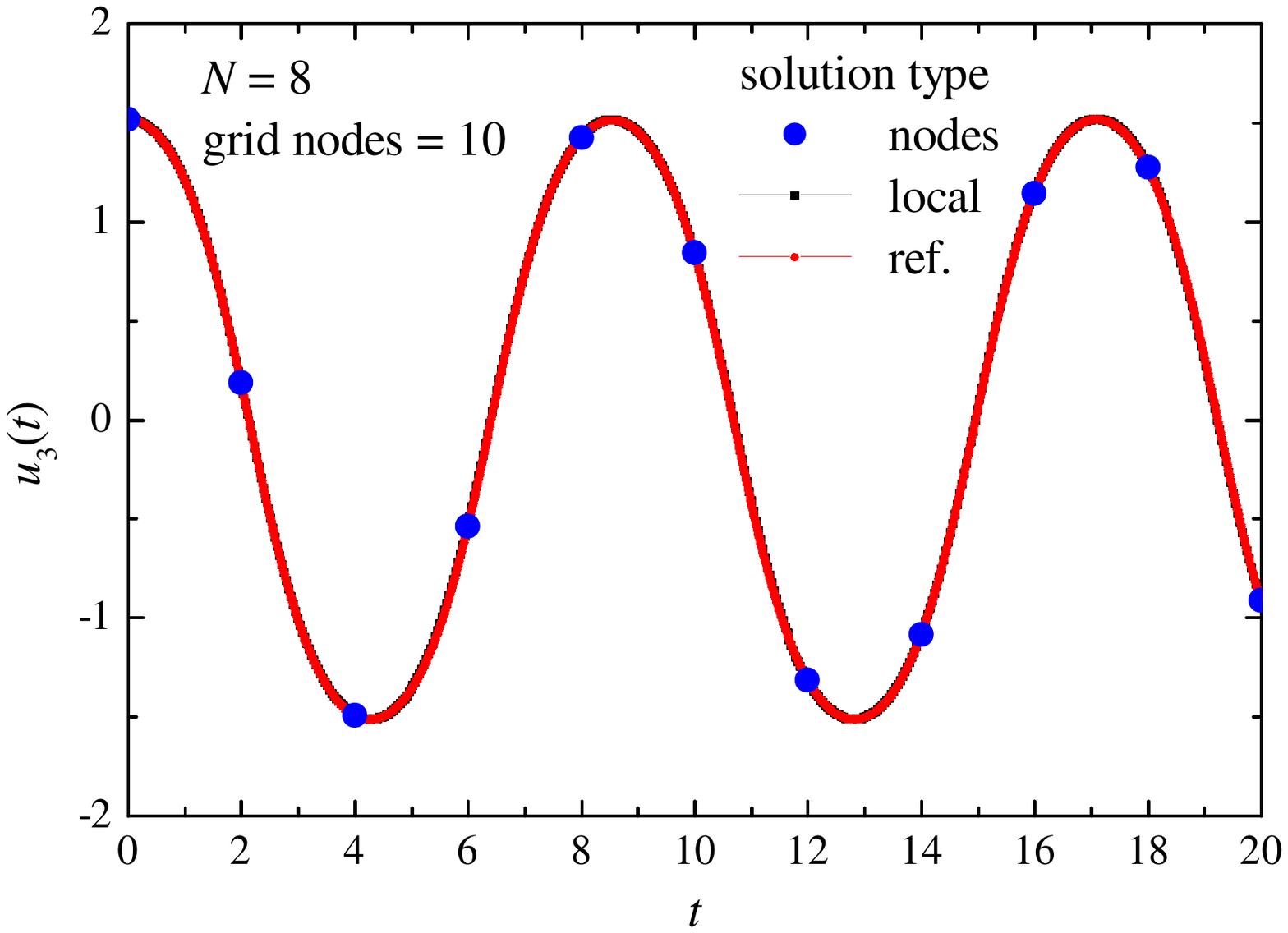}
\vspace{-8mm}\caption{\label{fig:dpend_ind1_sols_u:c3}}
\end{subfigure}
\begin{subfigure}{0.240\textwidth}
\includegraphics[width=\textwidth]{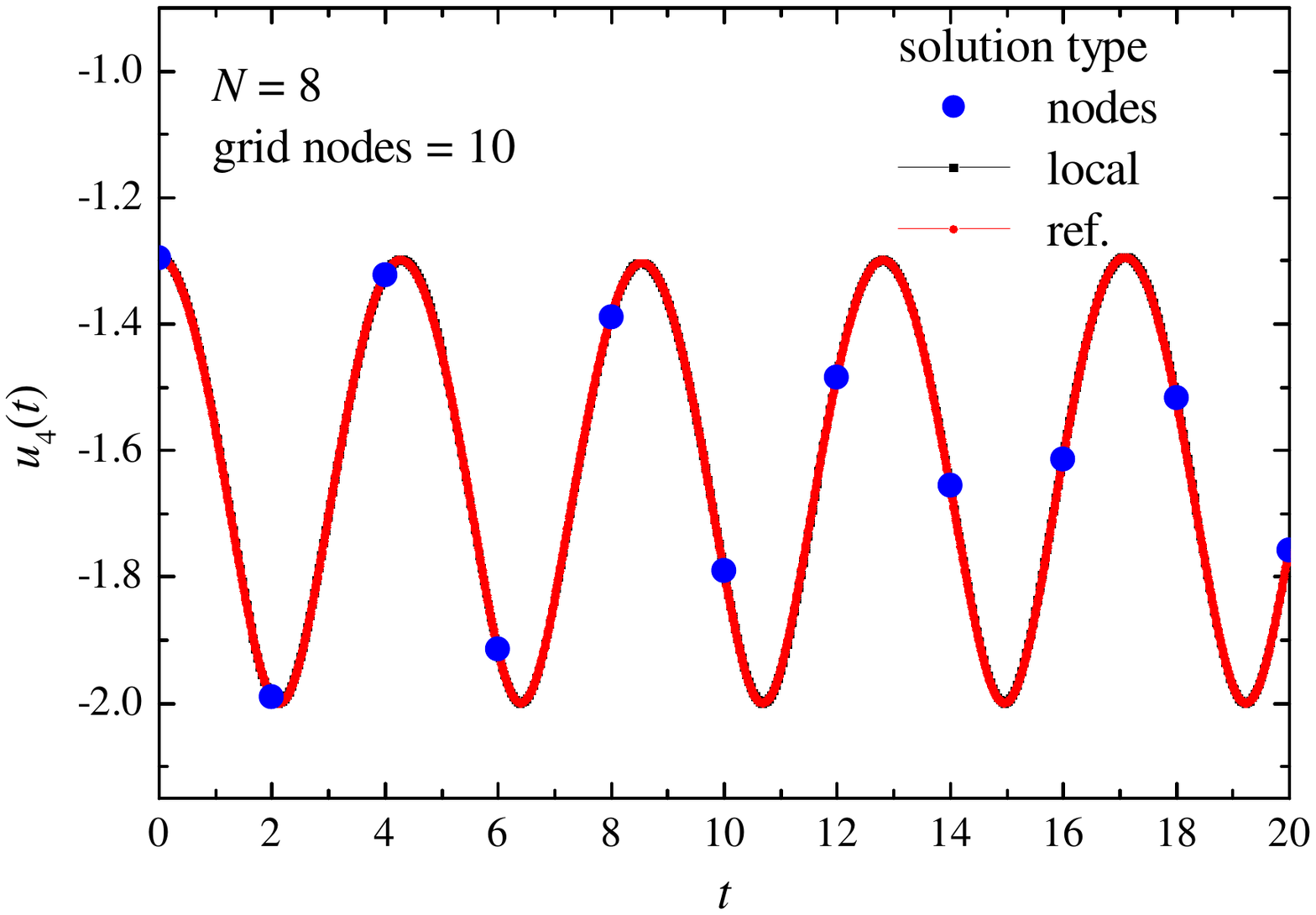}
\vspace{-8mm}\caption{\label{fig:dpend_ind1_sols_u:c4}}
\end{subfigure}\\[2mm]
\begin{subfigure}{0.240\textwidth}
\includegraphics[width=\textwidth]{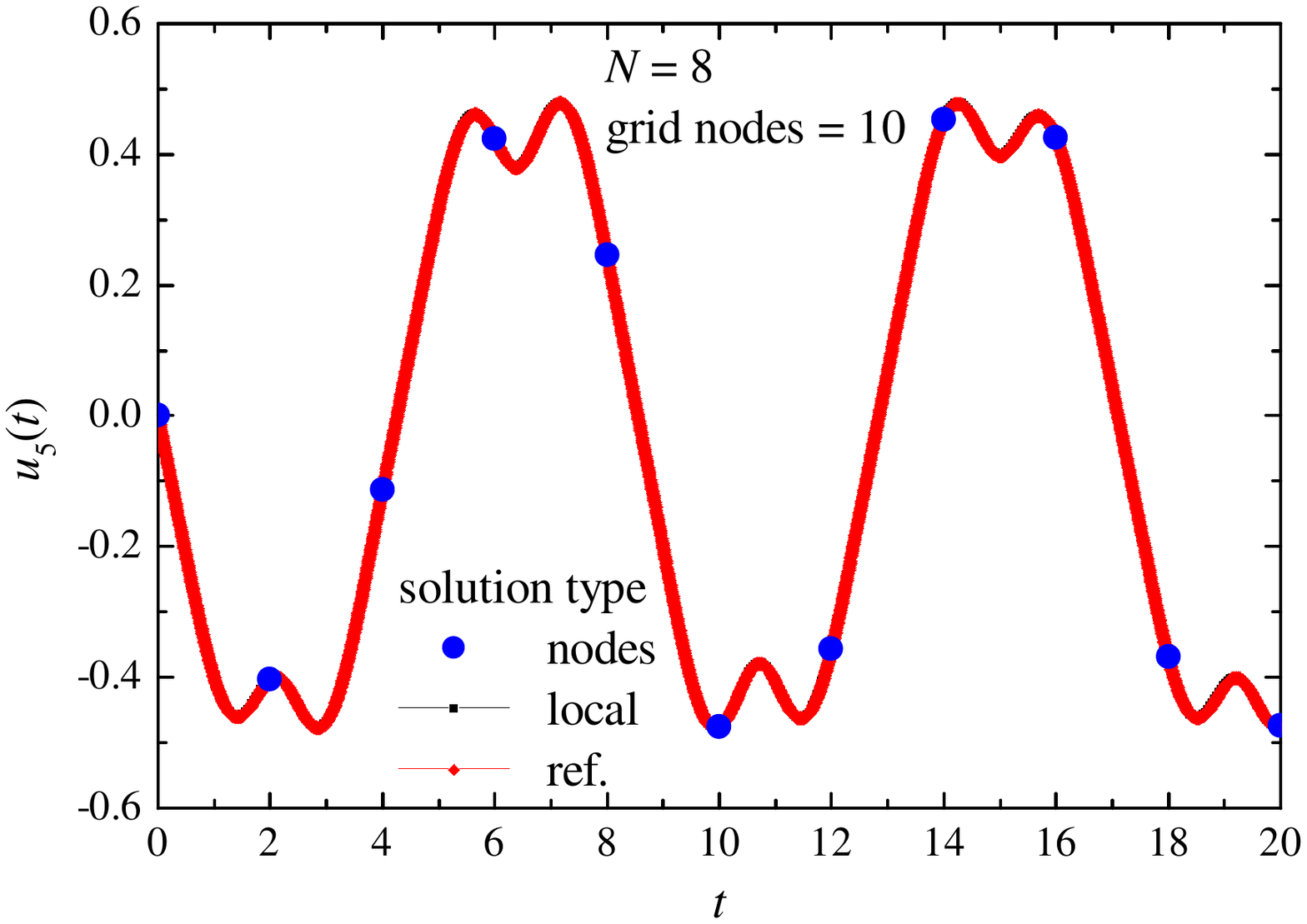}
\vspace{-8mm}\caption{\label{fig:dpend_ind1_sols_u:d1}}
\end{subfigure}
\begin{subfigure}{0.240\textwidth}
\includegraphics[width=\textwidth]{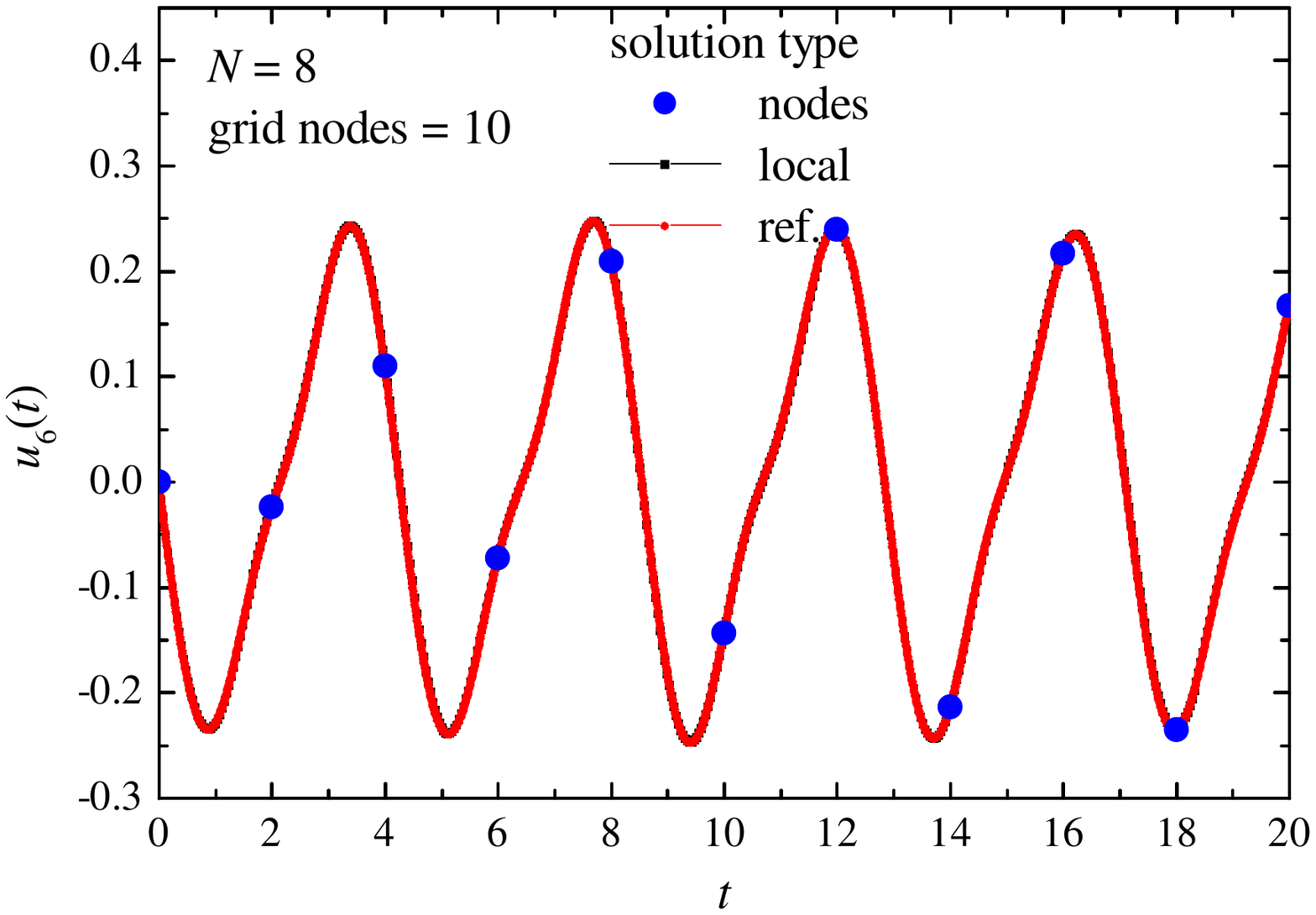}
\vspace{-8mm}\caption{\label{fig:dpend_ind1_sols_u:d2}}
\end{subfigure}
\begin{subfigure}{0.240\textwidth}
\includegraphics[width=\textwidth]{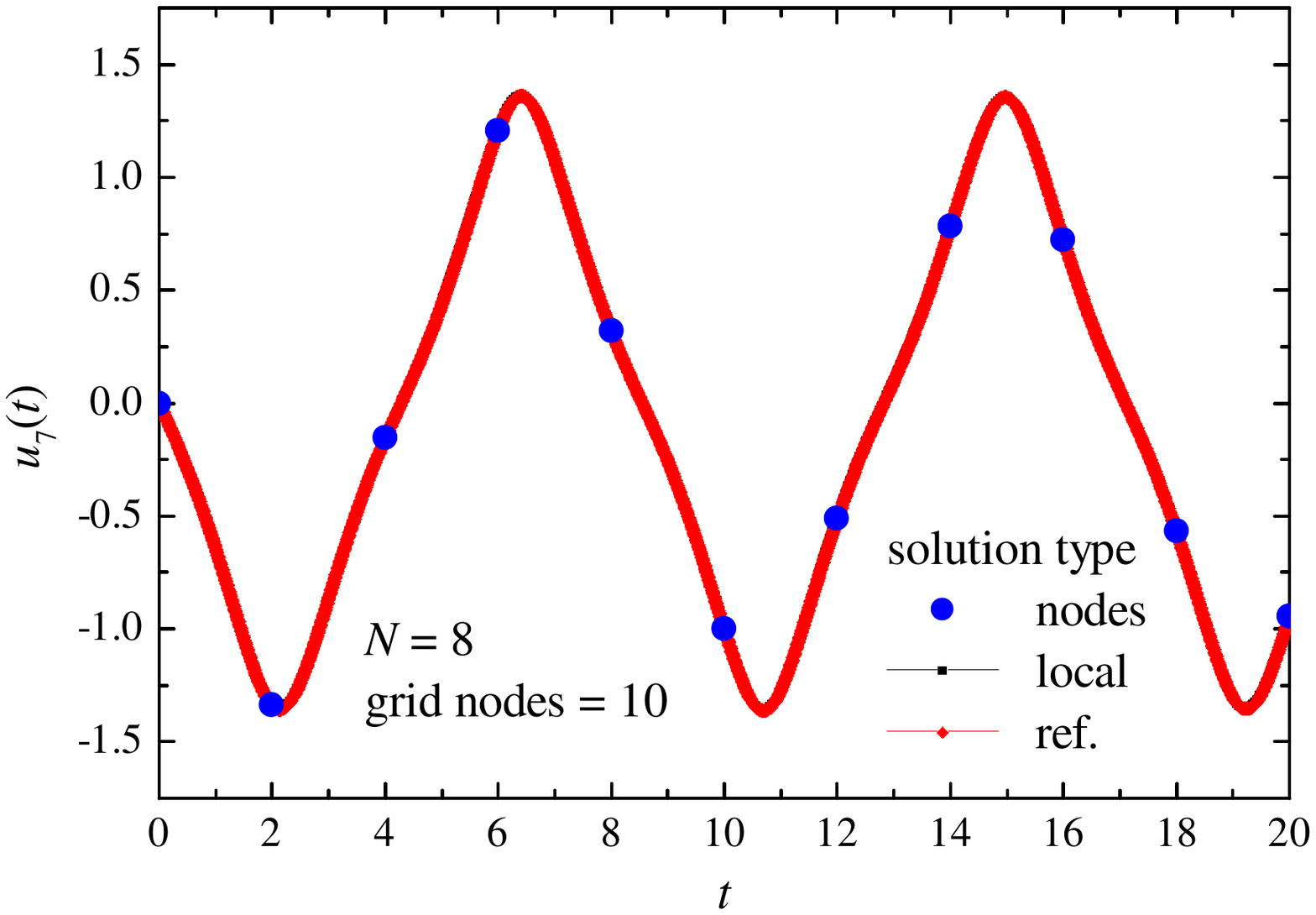}
\vspace{-8mm}\caption{\label{fig:dpend_ind1_sols_u:d3}}
\end{subfigure}
\begin{subfigure}{0.240\textwidth}
\includegraphics[width=\textwidth]{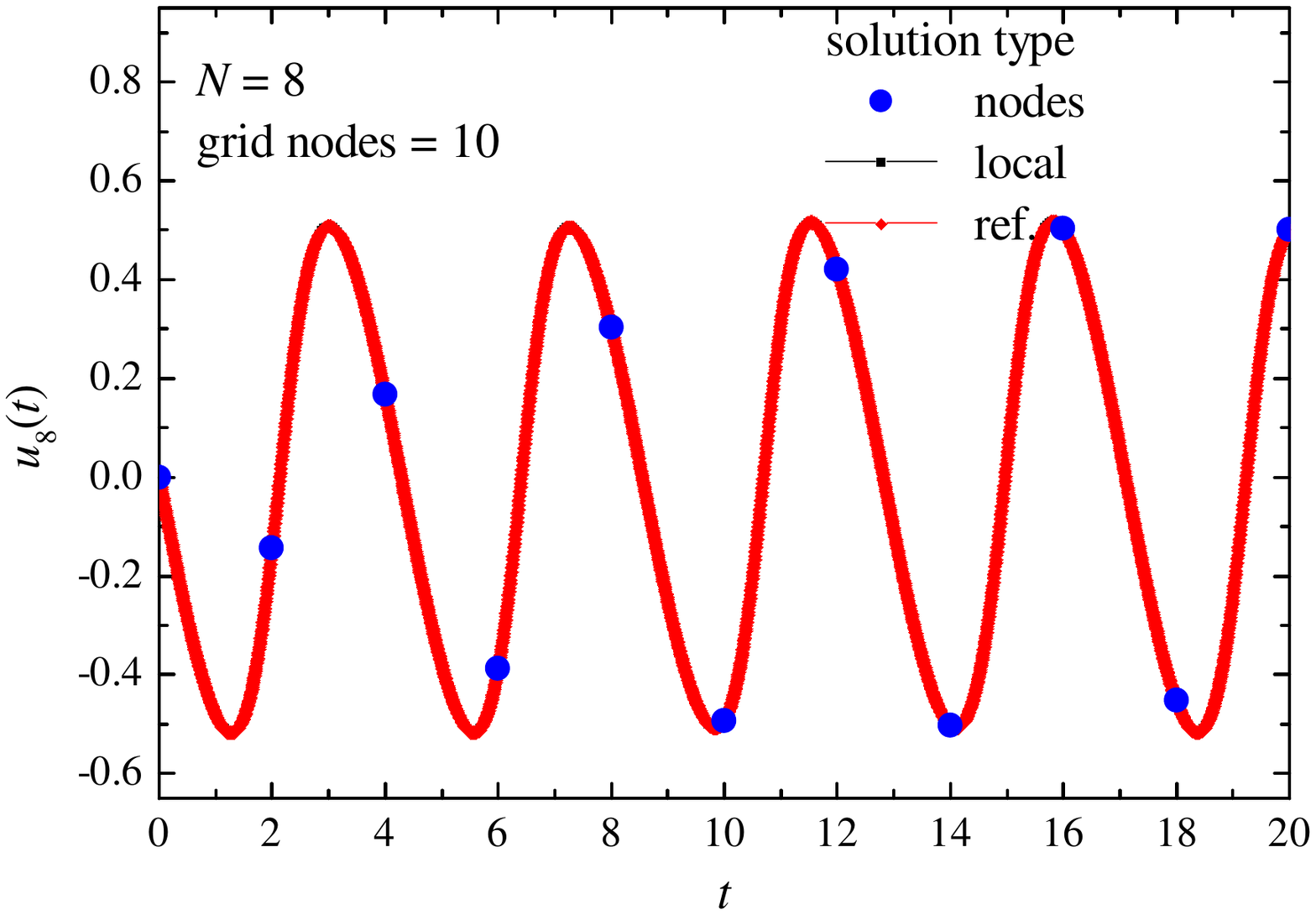}
\vspace{-8mm}\caption{\label{fig:dpend_ind1_sols_u:d4}}
\end{subfigure}\\[2mm]
\begin{subfigure}{0.240\textwidth}
\includegraphics[width=\textwidth]{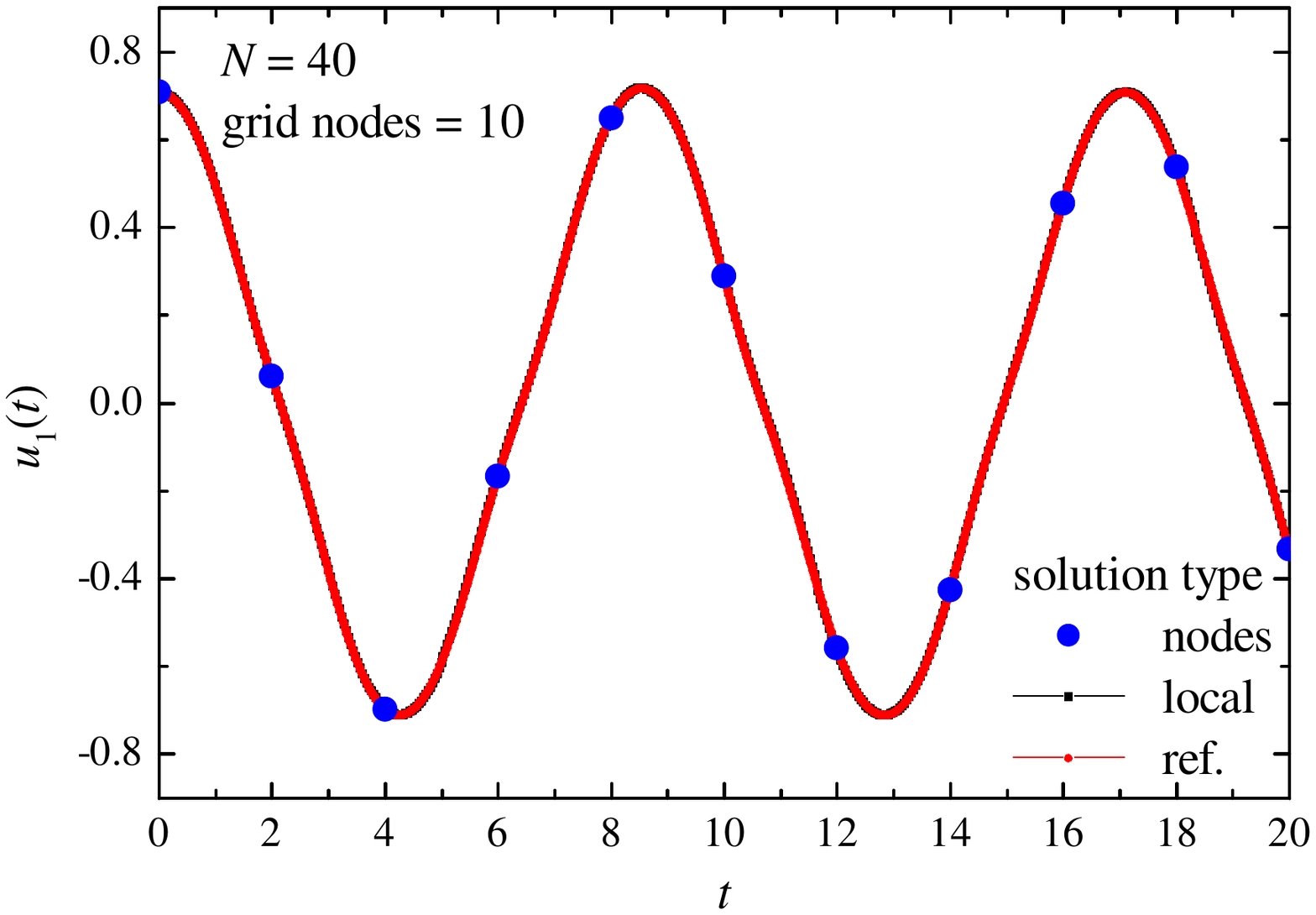}
\vspace{-8mm}\caption{\label{fig:dpend_ind1_sols_u:e1}}
\end{subfigure}
\begin{subfigure}{0.240\textwidth}
\includegraphics[width=\textwidth]{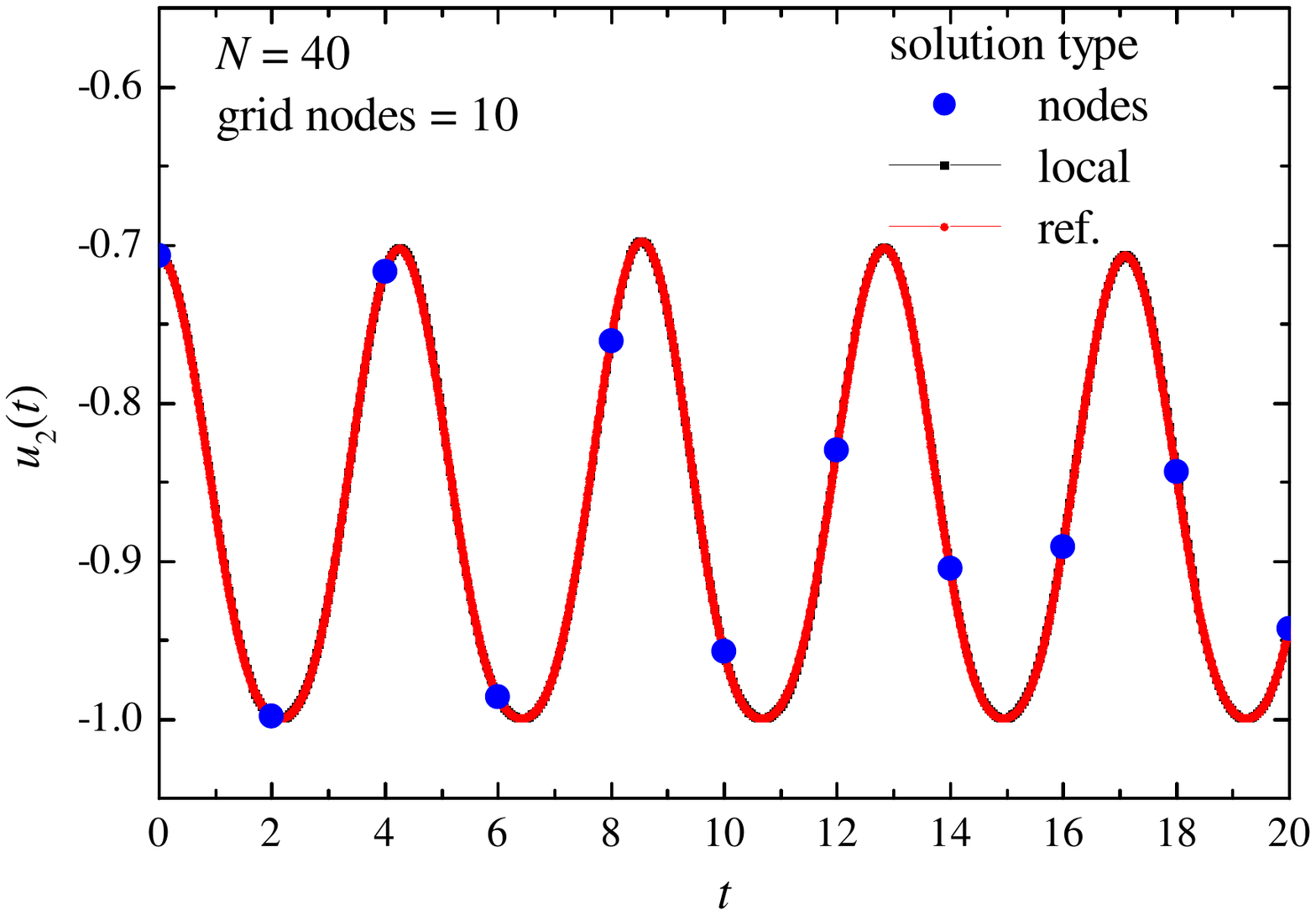}
\vspace{-8mm}\caption{\label{fig:dpend_ind1_sols_u:e2}}
\end{subfigure}
\begin{subfigure}{0.240\textwidth}
\includegraphics[width=\textwidth]{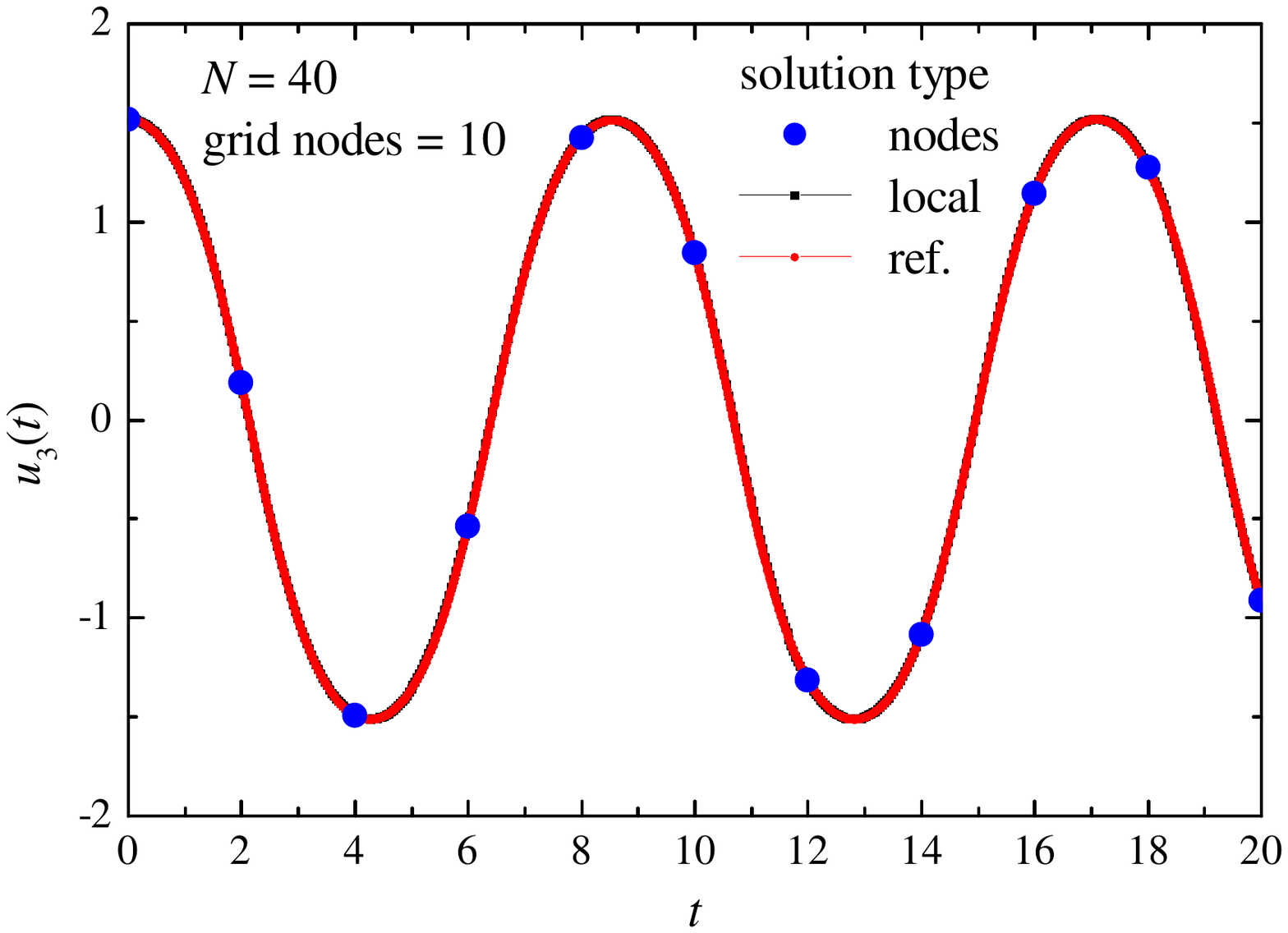}
\vspace{-8mm}\caption{\label{fig:dpend_ind1_sols_u:e3}}
\end{subfigure}
\begin{subfigure}{0.240\textwidth}
\includegraphics[width=\textwidth]{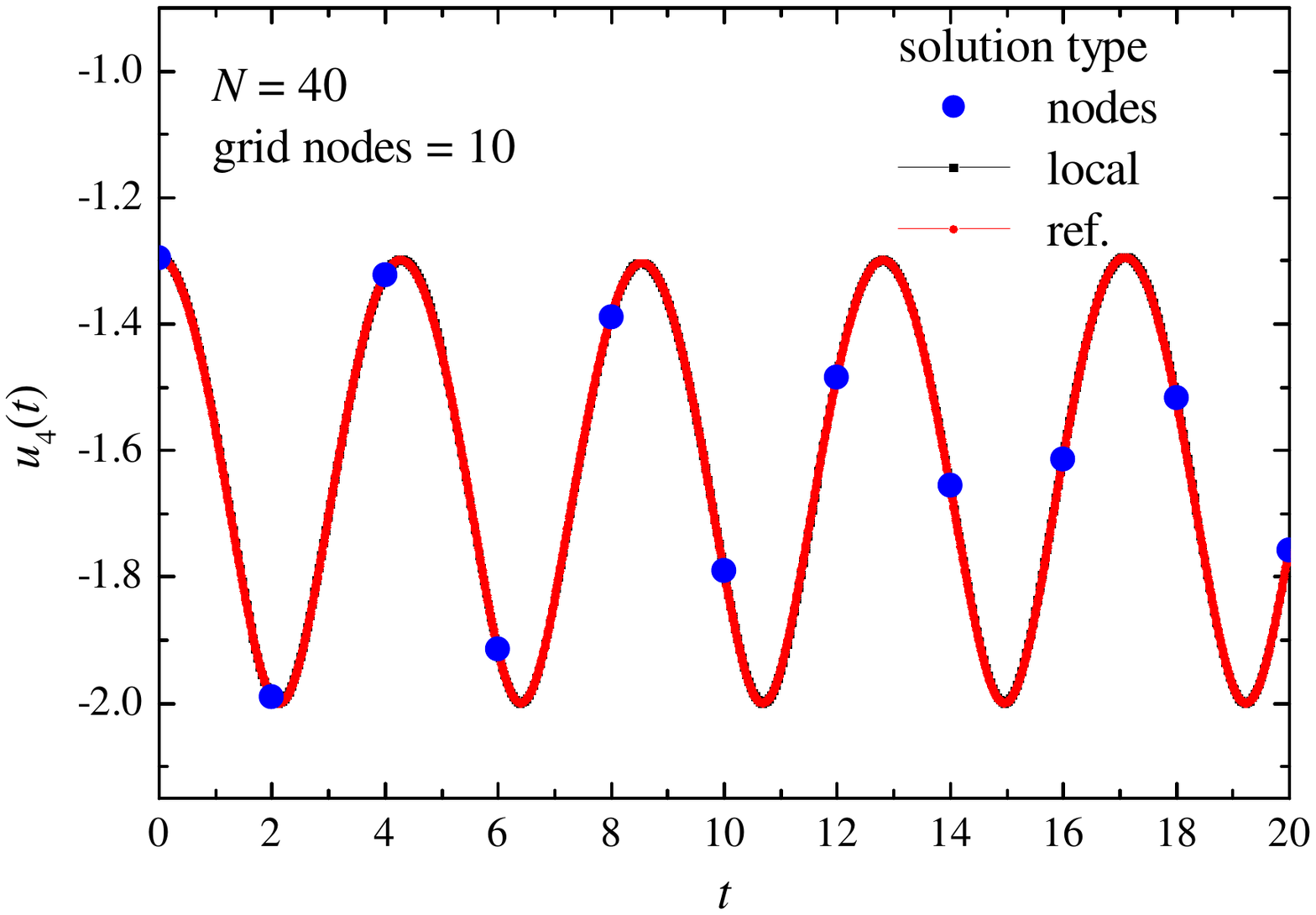}
\vspace{-8mm}\caption{\label{fig:dpend_ind1_sols_u:e4}}
\end{subfigure}\\[2mm]
\begin{subfigure}{0.240\textwidth}
\includegraphics[width=\textwidth]{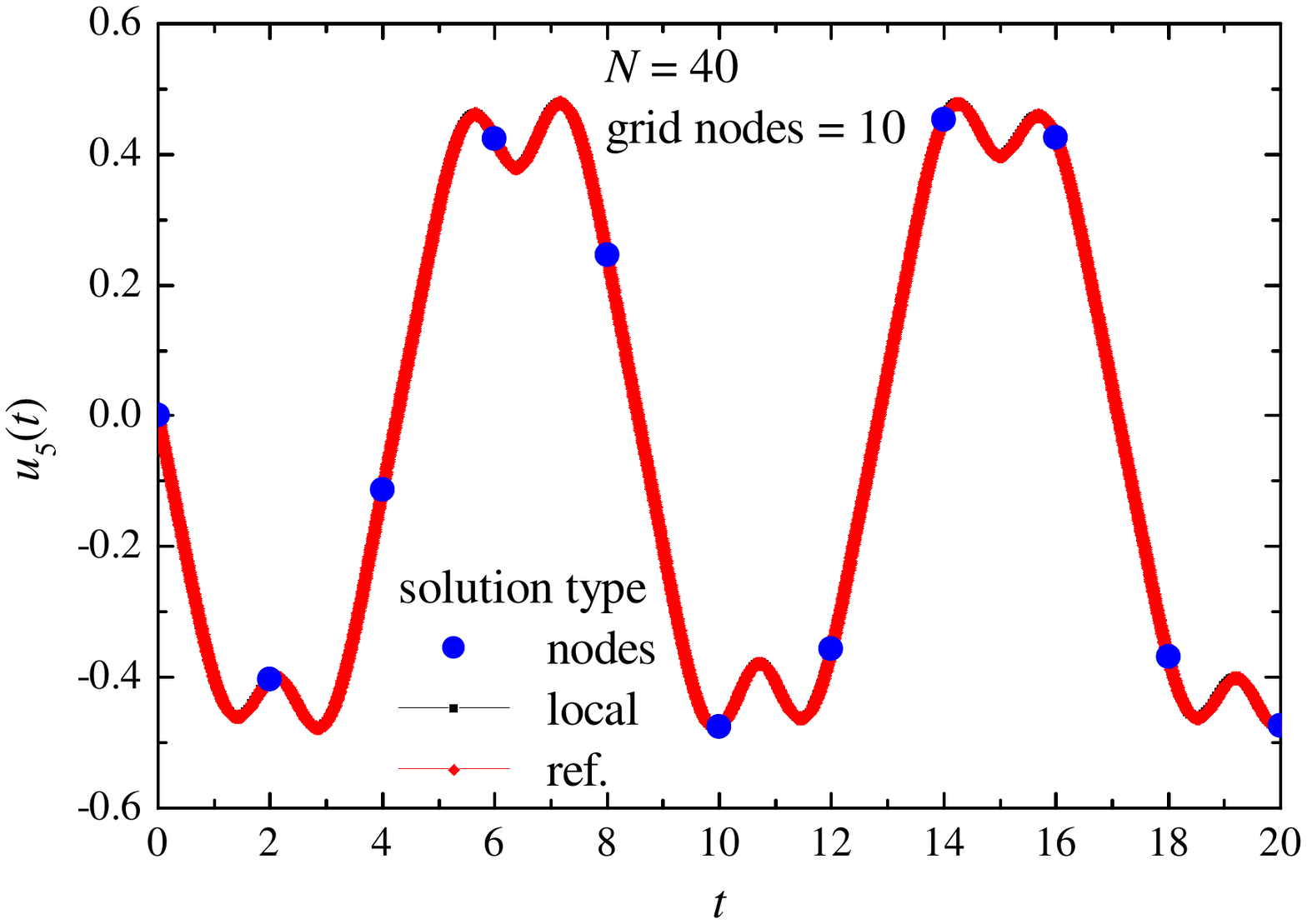}
\vspace{-8mm}\caption{\label{fig:dpend_ind1_sols_u:f1}}
\end{subfigure}
\begin{subfigure}{0.240\textwidth}
\includegraphics[width=\textwidth]{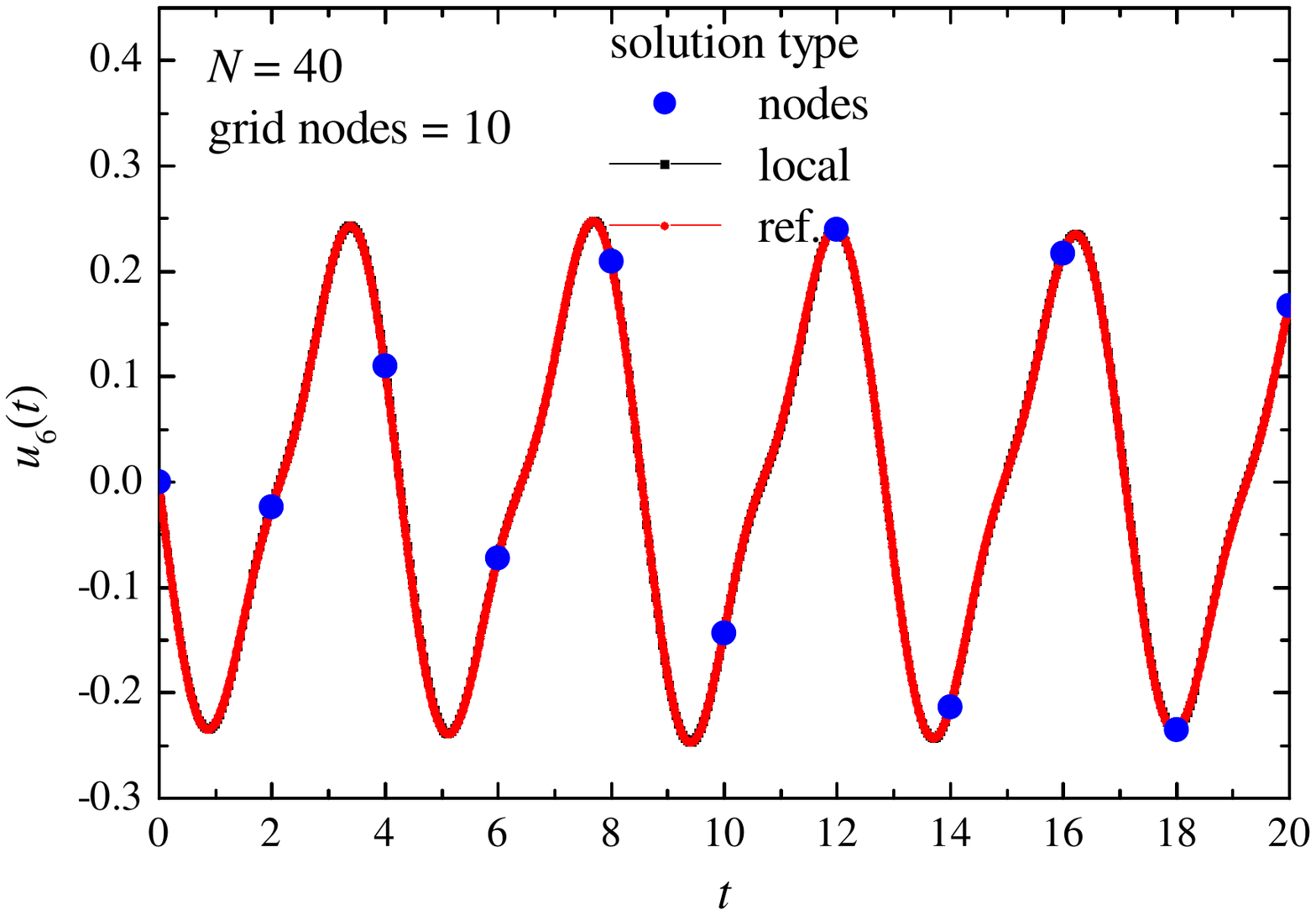}
\vspace{-8mm}\caption{\label{fig:dpend_ind1_sols_u:f2}}
\end{subfigure}
\begin{subfigure}{0.240\textwidth}
\includegraphics[width=\textwidth]{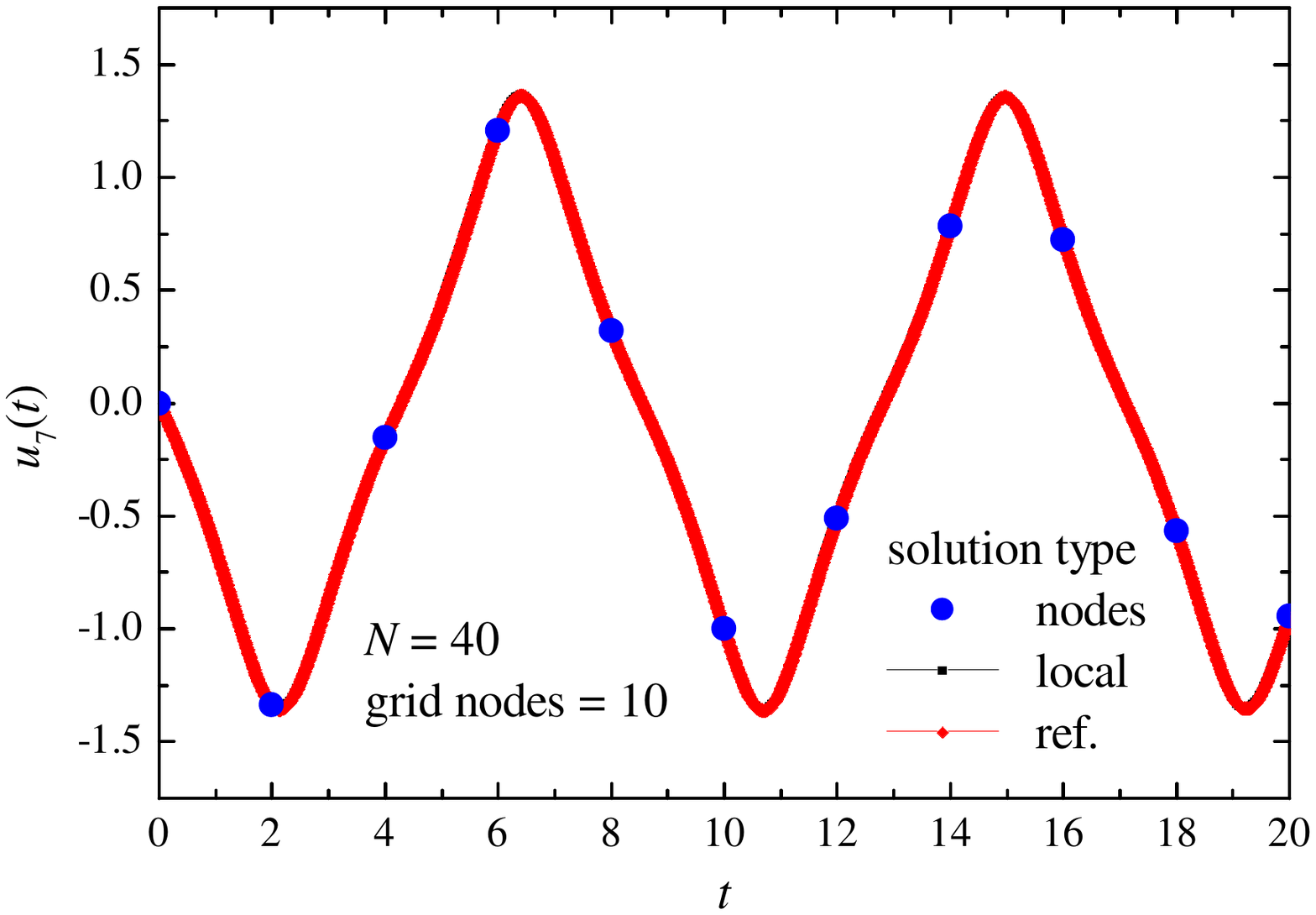}
\vspace{-8mm}\caption{\label{fig:dpend_ind1_sols_u:f3}}
\end{subfigure}
\begin{subfigure}{0.240\textwidth}
\includegraphics[width=\textwidth]{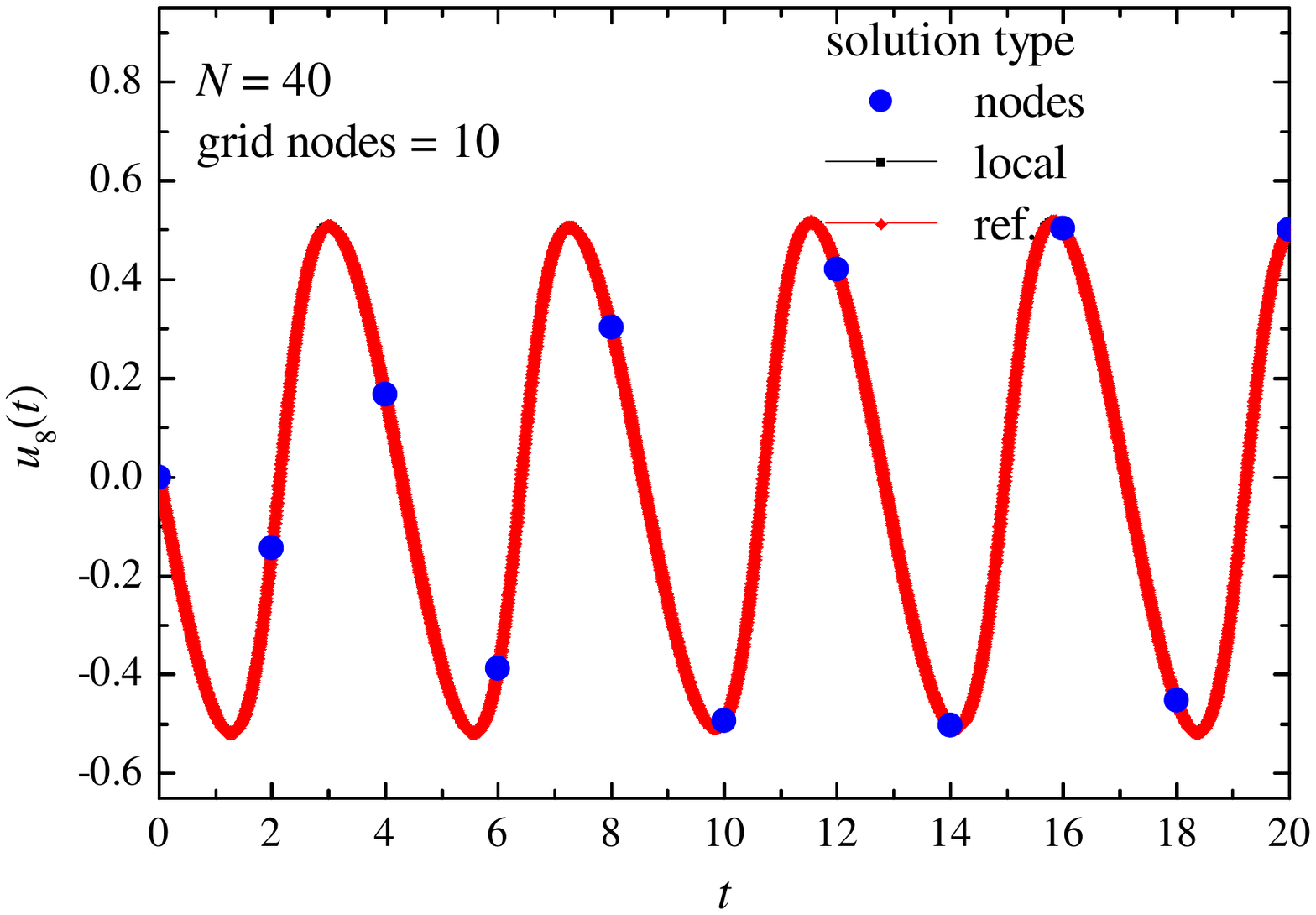}
\vspace{-8mm}\caption{\label{fig:dpend_ind1_sols_u:f4}}
\end{subfigure}\\[2mm]
\caption{%
Numerical solution of the DAE system (\ref{eq:math_dpend_dae_ind_3}) of index 1. Comparison of the solution at nodes $\mathbf{u}_{n}$, the local solution $\mathbf{u}_{L}(t)$ and the reference solution $\mathbf{u}^{\rm ref}(t)$ for components $u_{1}$ (\subref{fig:dpend_ind1_sols_u:a1}, \subref{fig:dpend_ind1_sols_u:c1}, \subref{fig:dpend_ind1_sols_u:e1}), $u_{2}$ (\subref{fig:dpend_ind1_sols_u:a2}, \subref{fig:dpend_ind1_sols_u:c2}, \subref{fig:dpend_ind1_sols_u:e2}), $u_{3}$ (\subref{fig:dpend_ind1_sols_u:a3}, \subref{fig:dpend_ind1_sols_u:c3}, \subref{fig:dpend_ind1_sols_u:e3}), $u_{4}$ (\subref{fig:dpend_ind1_sols_u:a4}, \subref{fig:dpend_ind1_sols_u:c4}, \subref{fig:dpend_ind1_sols_u:e4}), $u_{5}$ (\subref{fig:dpend_ind1_sols_u:b1}, \subref{fig:dpend_ind1_sols_u:d1}, \subref{fig:dpend_ind1_sols_u:f1}), $u_{6}$ (\subref{fig:dpend_ind1_sols_u:b2}, \subref{fig:dpend_ind1_sols_u:d2}, \subref{fig:dpend_ind1_sols_u:f2}), $u_{7}$ (\subref{fig:dpend_ind1_sols_u:b3}, \subref{fig:dpend_ind1_sols_u:d3}, \subref{fig:dpend_ind1_sols_u:f3}), $u_{8}$ (\subref{fig:dpend_ind1_sols_u:b4}, \subref{fig:dpend_ind1_sols_u:d4}, \subref{fig:dpend_ind1_sols_u:f4}), obtained using polynomials with degrees $N = 1$ (\subref{fig:dpend_ind1_sols_u:a1}, \subref{fig:dpend_ind1_sols_u:a2}, \subref{fig:dpend_ind1_sols_u:a3}, \subref{fig:dpend_ind1_sols_u:a4}, \subref{fig:dpend_ind1_sols_u:b1}, \subref{fig:dpend_ind1_sols_u:b2}, \subref{fig:dpend_ind1_sols_u:b3}, \subref{fig:dpend_ind1_sols_u:b4}), $N = 8$ (\subref{fig:dpend_ind1_sols_u:c1}, \subref{fig:dpend_ind1_sols_u:c2}, \subref{fig:dpend_ind1_sols_u:c3}, \subref{fig:dpend_ind1_sols_u:c4}, \subref{fig:dpend_ind1_sols_u:d1}, \subref{fig:dpend_ind1_sols_u:d2}, \subref{fig:dpend_ind1_sols_u:d3}, \subref{fig:dpend_ind1_sols_u:d4}) and $N = 40$ (\subref{fig:dpend_ind1_sols_u:e1}, \subref{fig:dpend_ind1_sols_u:e2}, \subref{fig:dpend_ind1_sols_u:e3}, \subref{fig:dpend_ind1_sols_u:e4}, \subref{fig:dpend_ind1_sols_u:f1}, \subref{fig:dpend_ind1_sols_u:f2}, \subref{fig:dpend_ind1_sols_u:f3}, \subref{fig:dpend_ind1_sols_u:f4}).
}
\label{fig:dpend_ind1_sols_u}
\end{figure} 

\begin{figure}[h!]
\captionsetup[subfigure]{%
	position=bottom,
	font+=smaller,
	textfont=normalfont,
	singlelinecheck=off,
	justification=raggedright
}
\centering
\begin{subfigure}{0.240\textwidth}
\includegraphics[width=\textwidth]{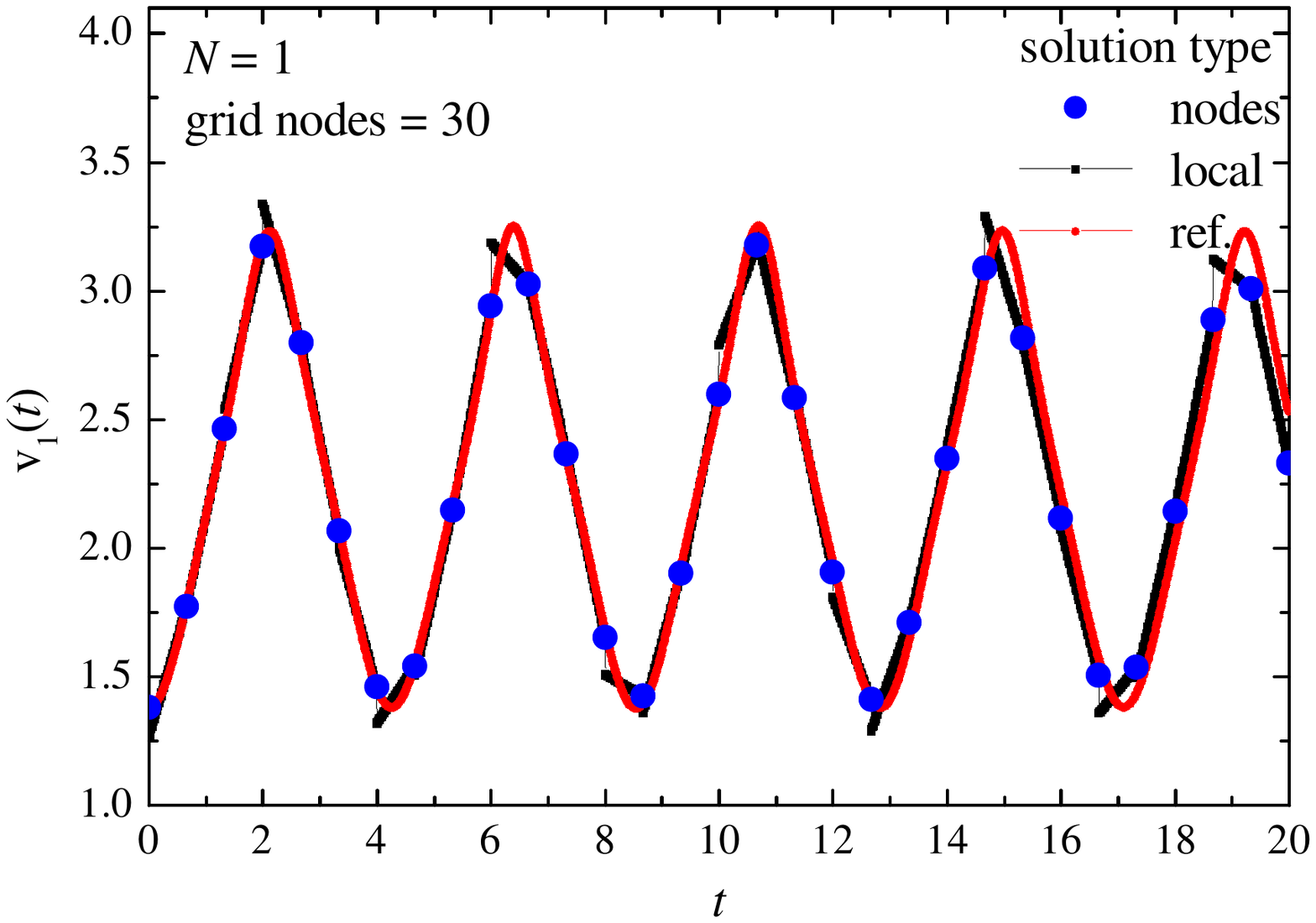}
\vspace{-8mm}\caption{\label{fig:dpend_ind1_sols_vg:a1}}
\end{subfigure}
\begin{subfigure}{0.240\textwidth}
\includegraphics[width=\textwidth]{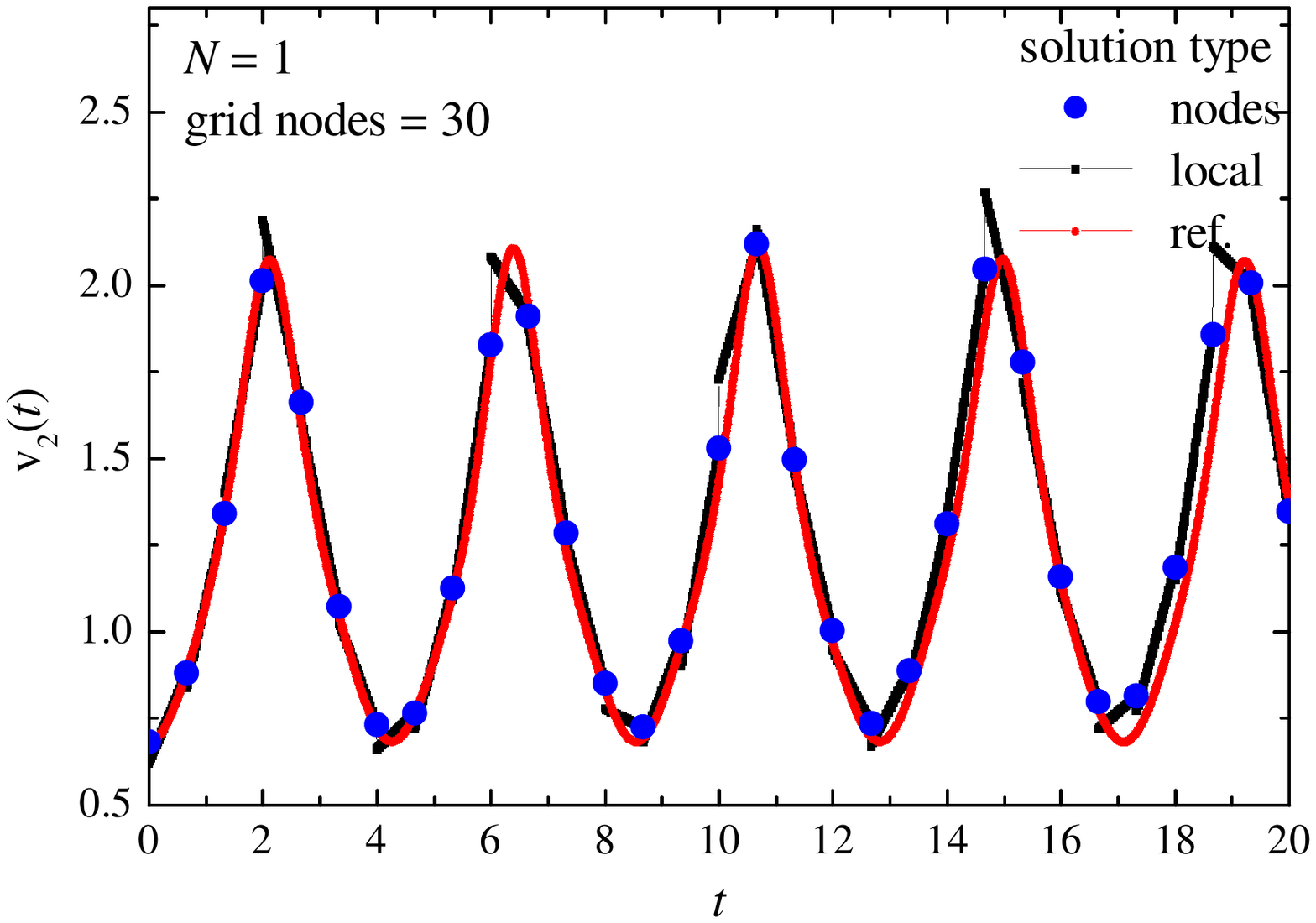}
\vspace{-8mm}\caption{\label{fig:dpend_ind1_sols_vg:a2}}
\end{subfigure}
\begin{subfigure}{0.240\textwidth}
\includegraphics[width=\textwidth]{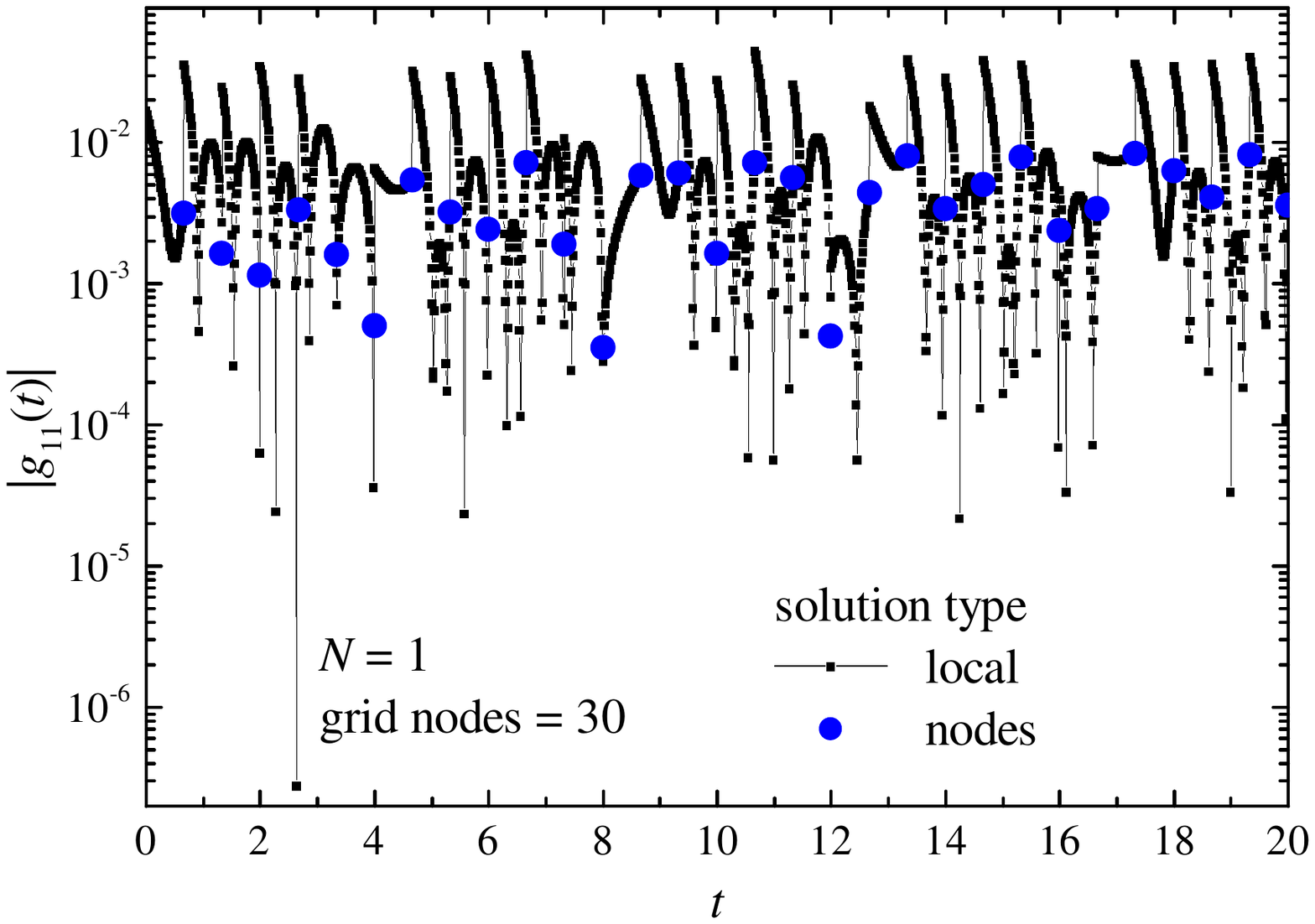}
\vspace{-8mm}\caption{\label{fig:dpend_ind1_sols_vg:a3}}
\end{subfigure}
\begin{subfigure}{0.240\textwidth}
\includegraphics[width=\textwidth]{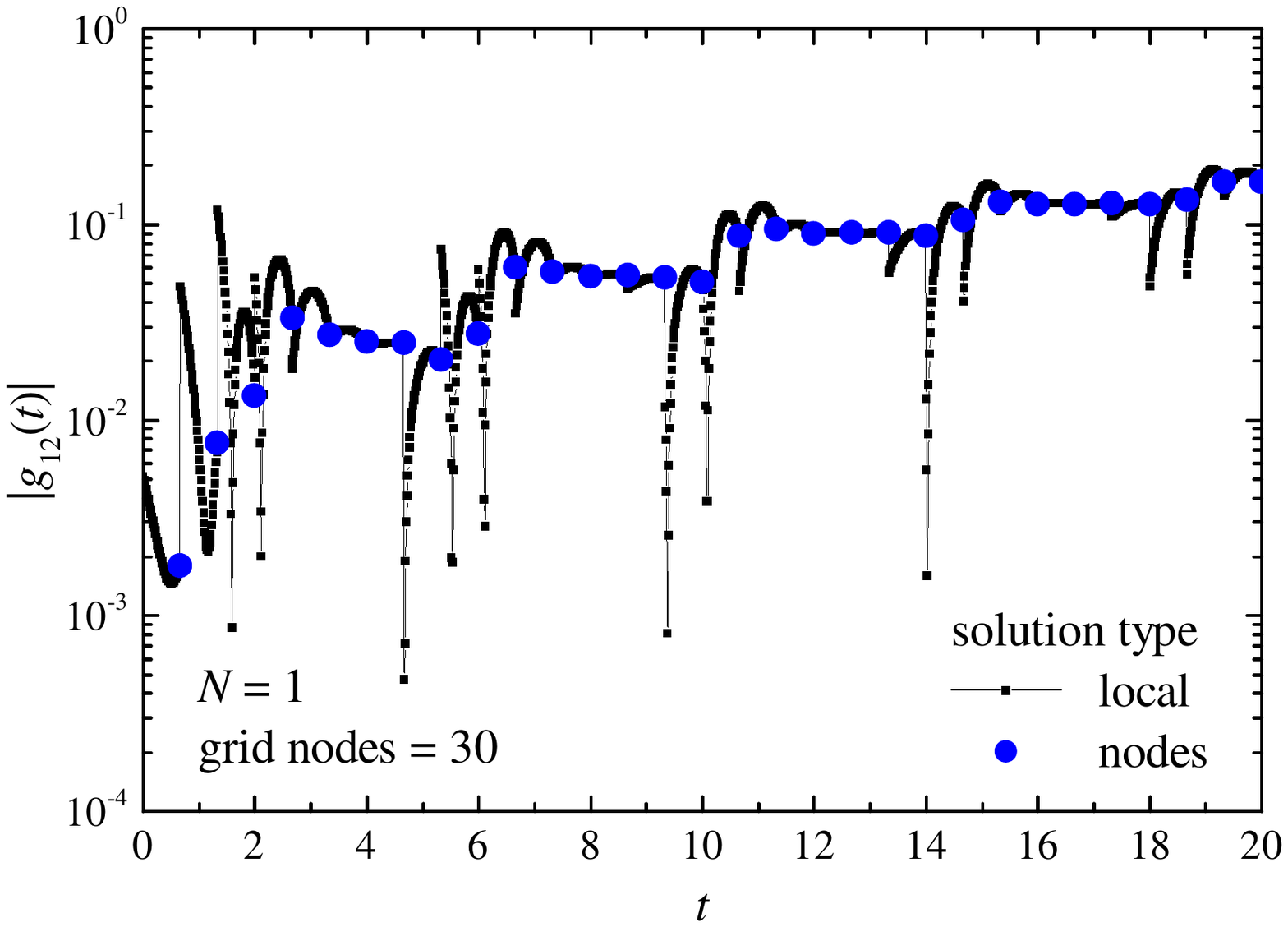}
\vspace{-8mm}\caption{\label{fig:dpend_ind1_sols_vg:a4}}
\end{subfigure}\\[2mm]
\begin{subfigure}{0.240\textwidth}
\includegraphics[width=\textwidth]{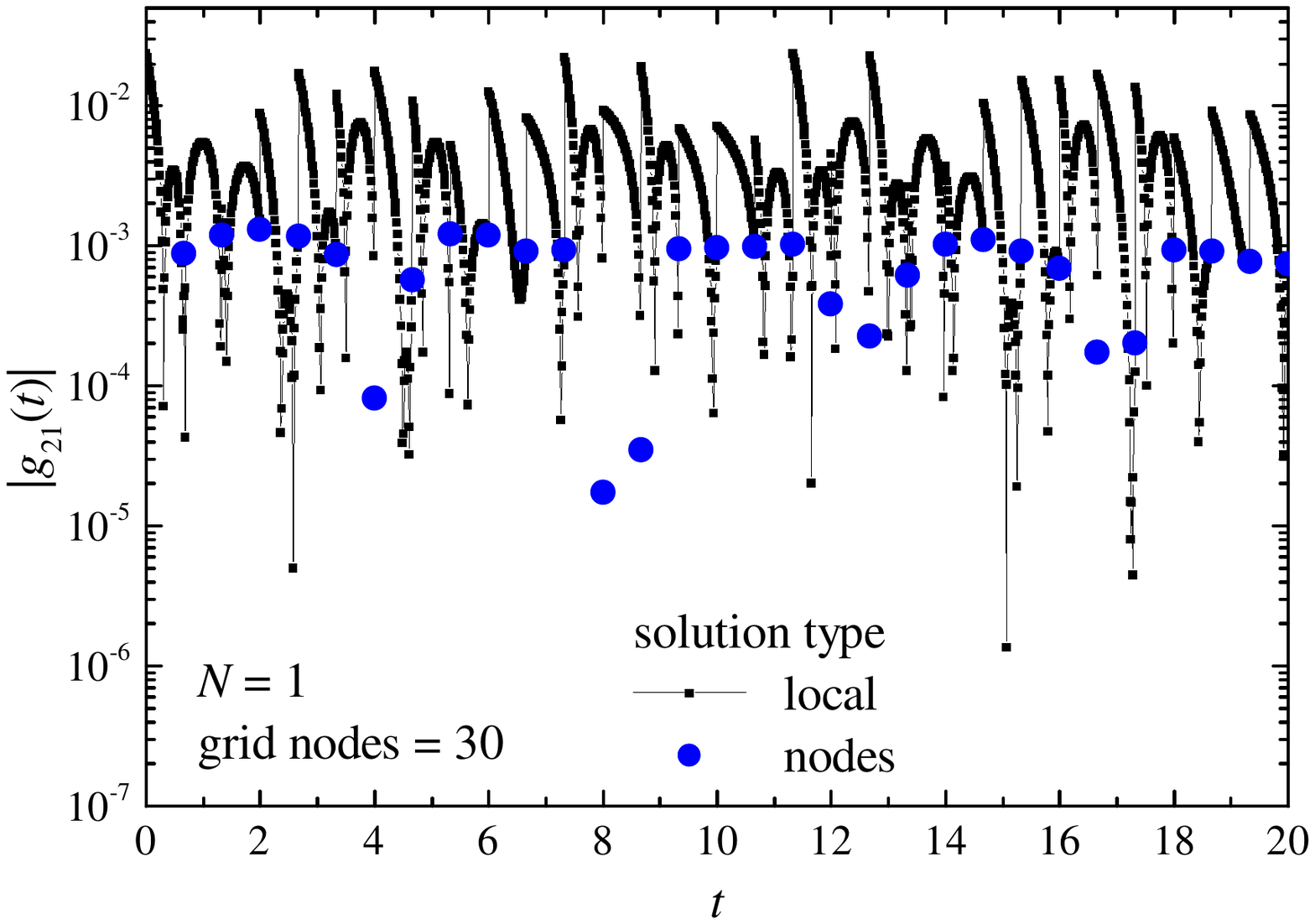}
\vspace{-8mm}\caption{\label{fig:dpend_ind1_sols_vg:b1}}
\end{subfigure}
\begin{subfigure}{0.240\textwidth}
\includegraphics[width=\textwidth]{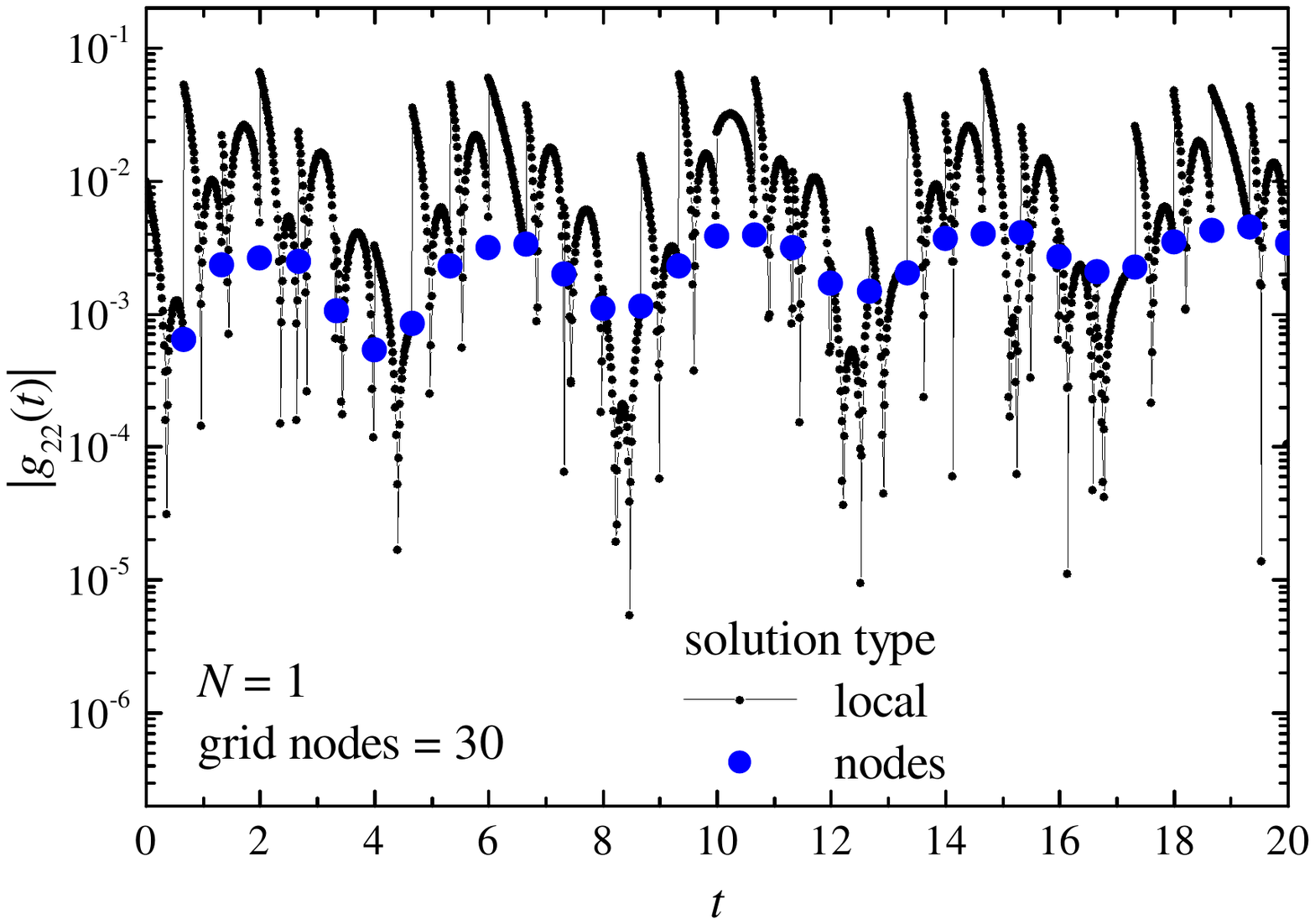}
\vspace{-8mm}\caption{\label{fig:dpend_ind1_sols_vg:b2}}
\end{subfigure}
\begin{subfigure}{0.240\textwidth}
\includegraphics[width=\textwidth]{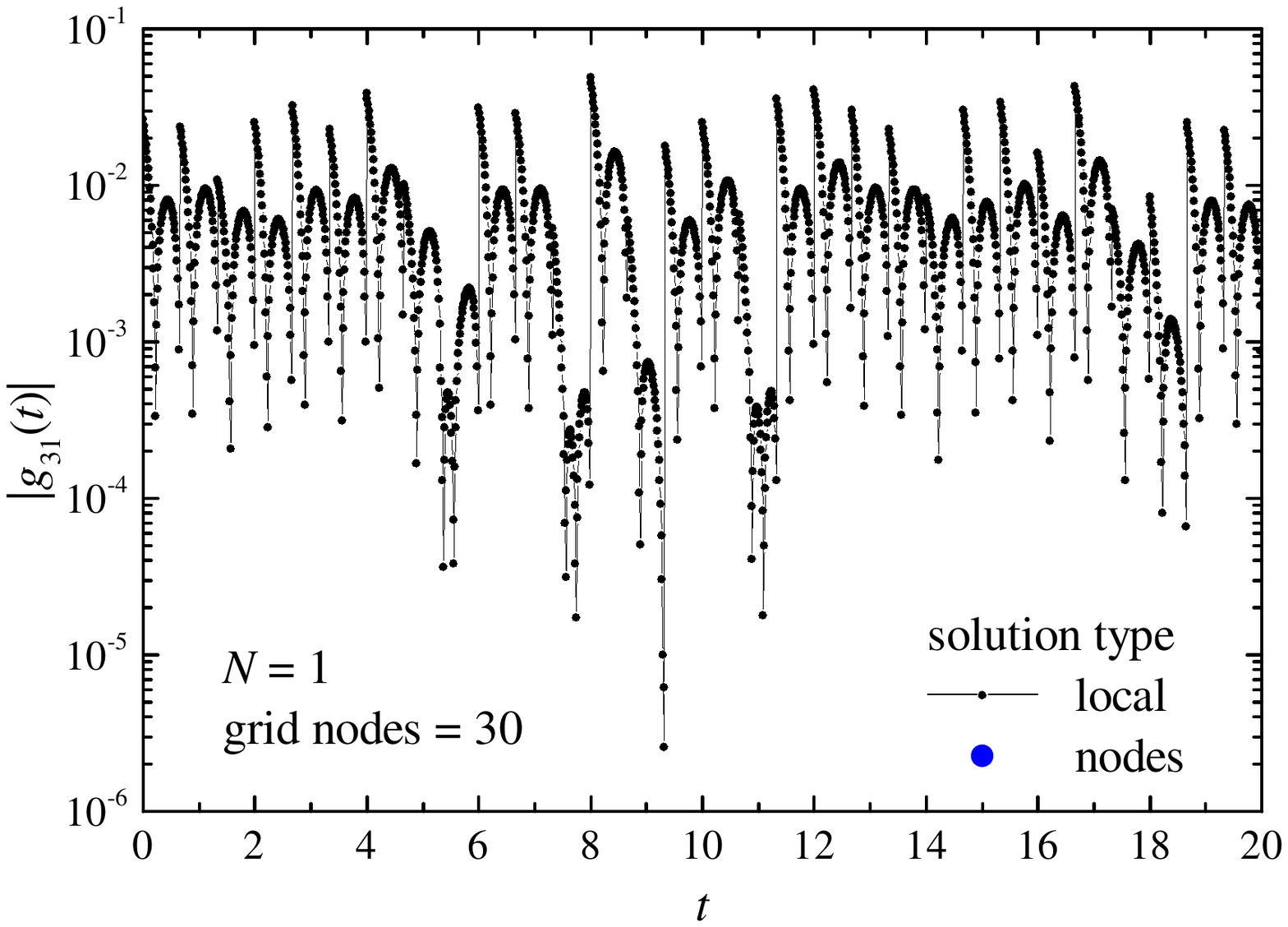}
\vspace{-8mm}\caption{\label{fig:dpend_ind1_sols_vg:b3}}
\end{subfigure}
\begin{subfigure}{0.240\textwidth}
\includegraphics[width=\textwidth]{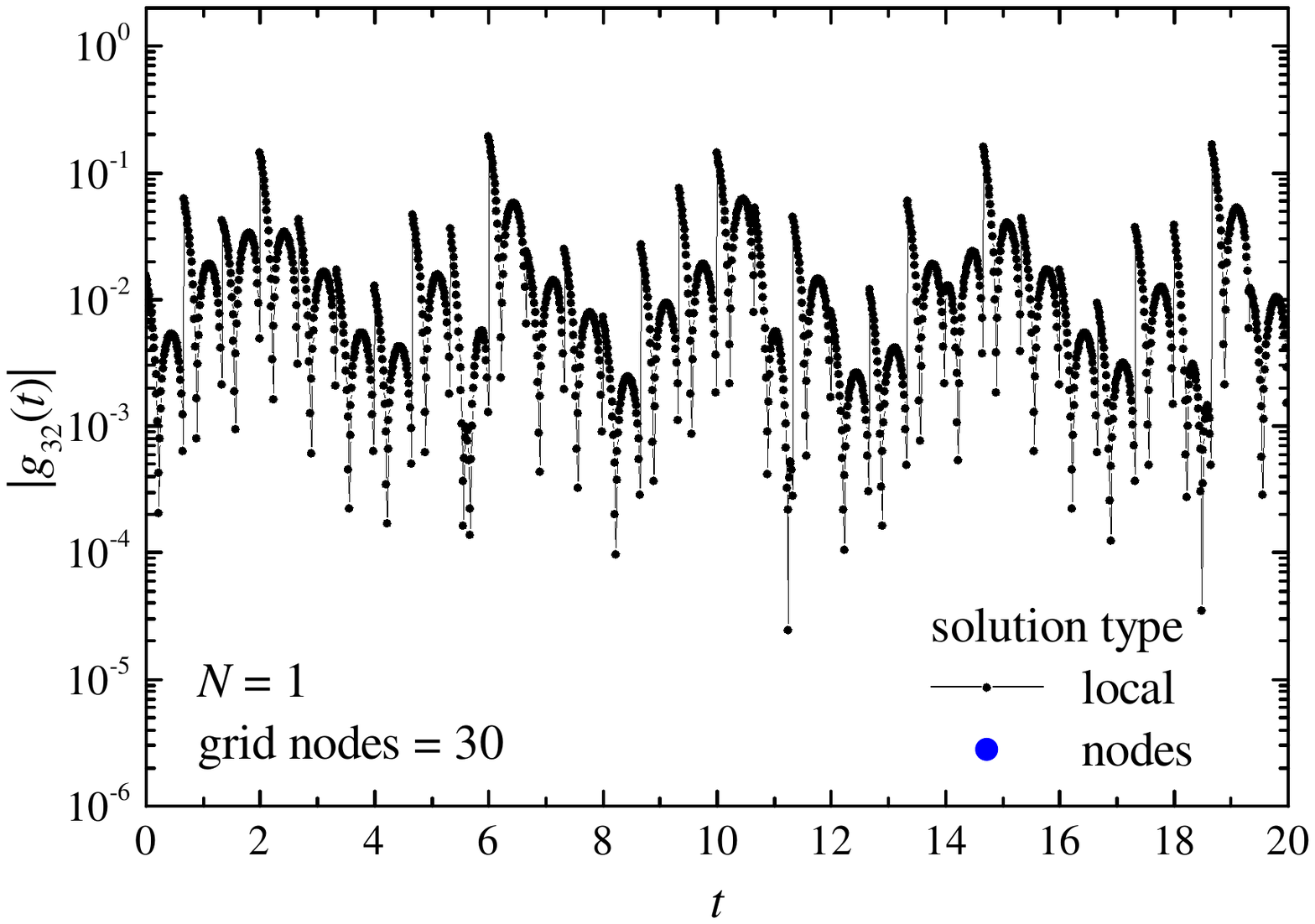}
\vspace{-8mm}\caption{\label{fig:dpend_ind1_sols_vg:b4}}
\end{subfigure}\\[2mm]
\begin{subfigure}{0.240\textwidth}
\includegraphics[width=\textwidth]{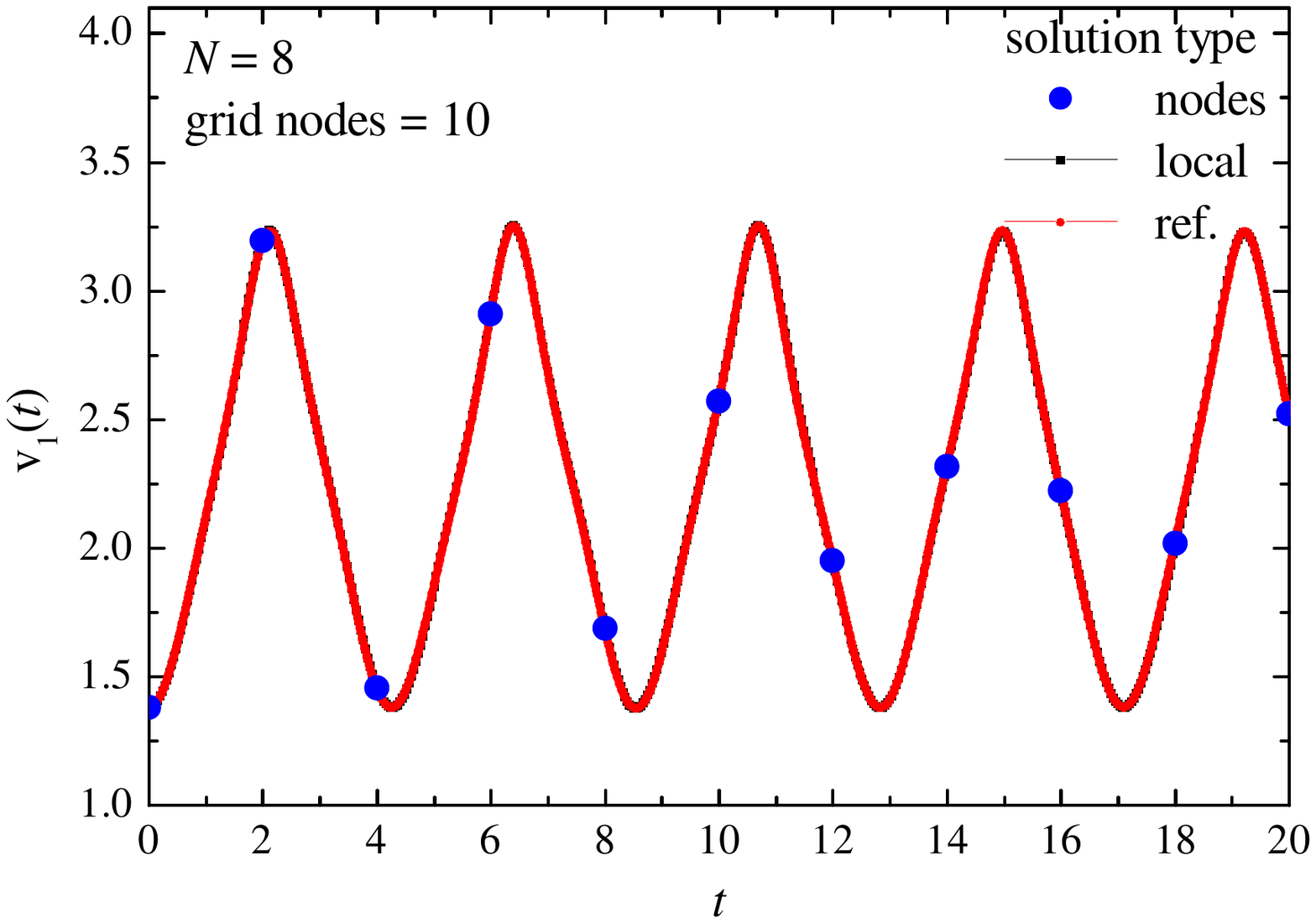}
\vspace{-8mm}\caption{\label{fig:dpend_ind1_sols_vg:c1}}
\end{subfigure}
\begin{subfigure}{0.240\textwidth}
\includegraphics[width=\textwidth]{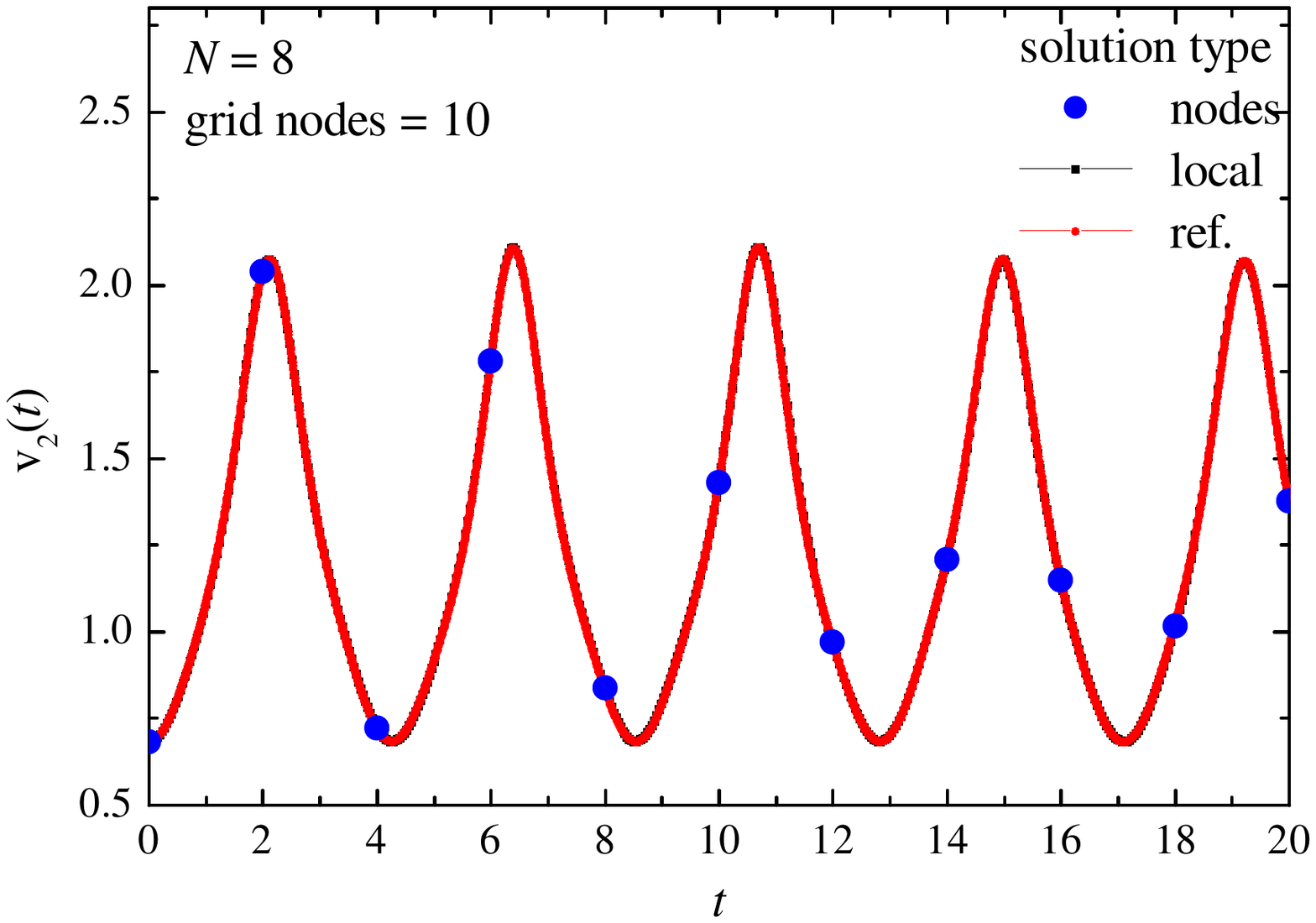}
\vspace{-8mm}\caption{\label{fig:dpend_ind1_sols_vg:c2}}
\end{subfigure}
\begin{subfigure}{0.240\textwidth}
\includegraphics[width=\textwidth]{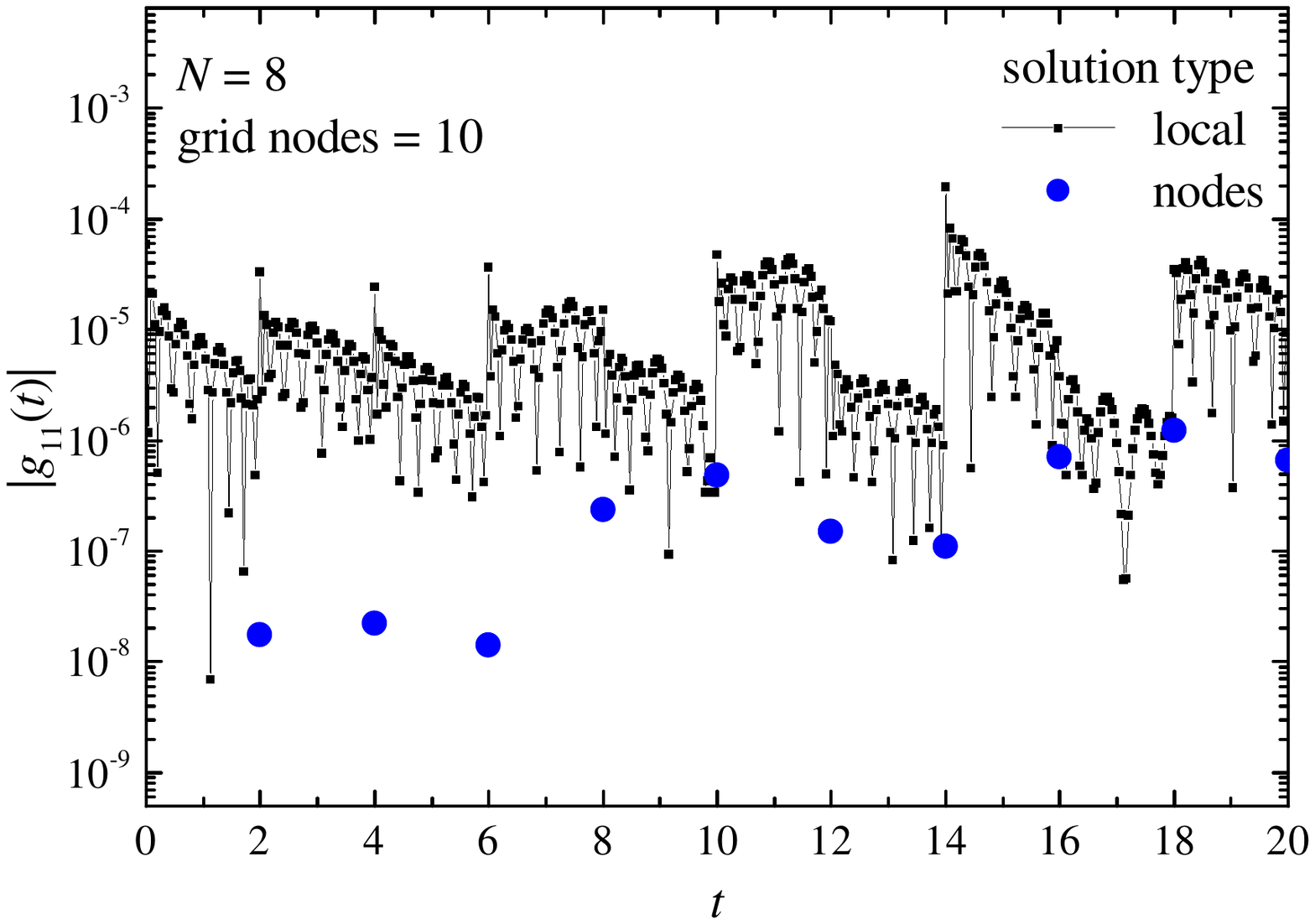}
\vspace{-8mm}\caption{\label{fig:dpend_ind1_sols_vg:c3}}
\end{subfigure}
\begin{subfigure}{0.240\textwidth}
\includegraphics[width=\textwidth]{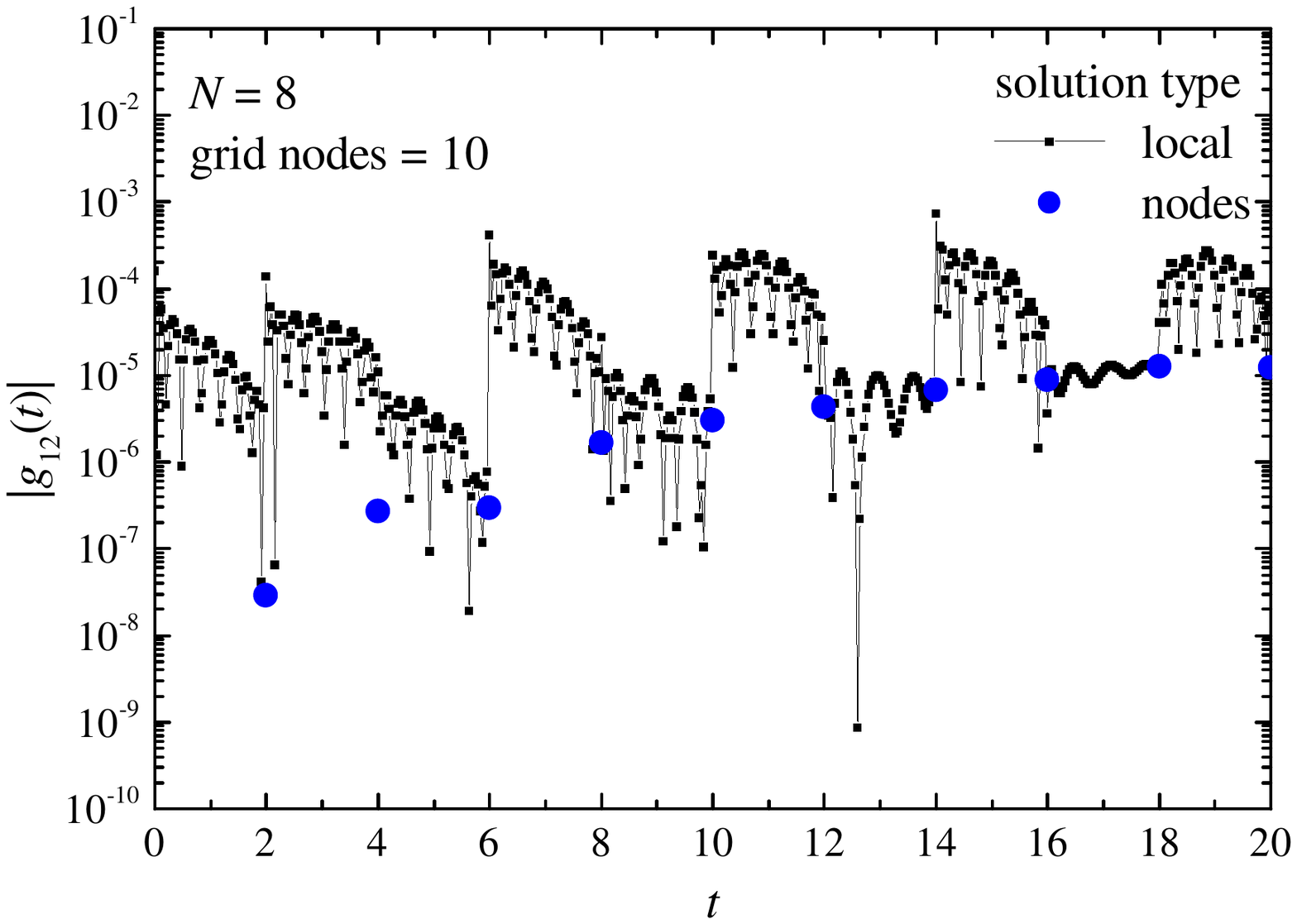}
\vspace{-8mm}\caption{\label{fig:dpend_ind1_sols_vg:c4}}
\end{subfigure}\\[2mm]
\begin{subfigure}{0.240\textwidth}
\includegraphics[width=\textwidth]{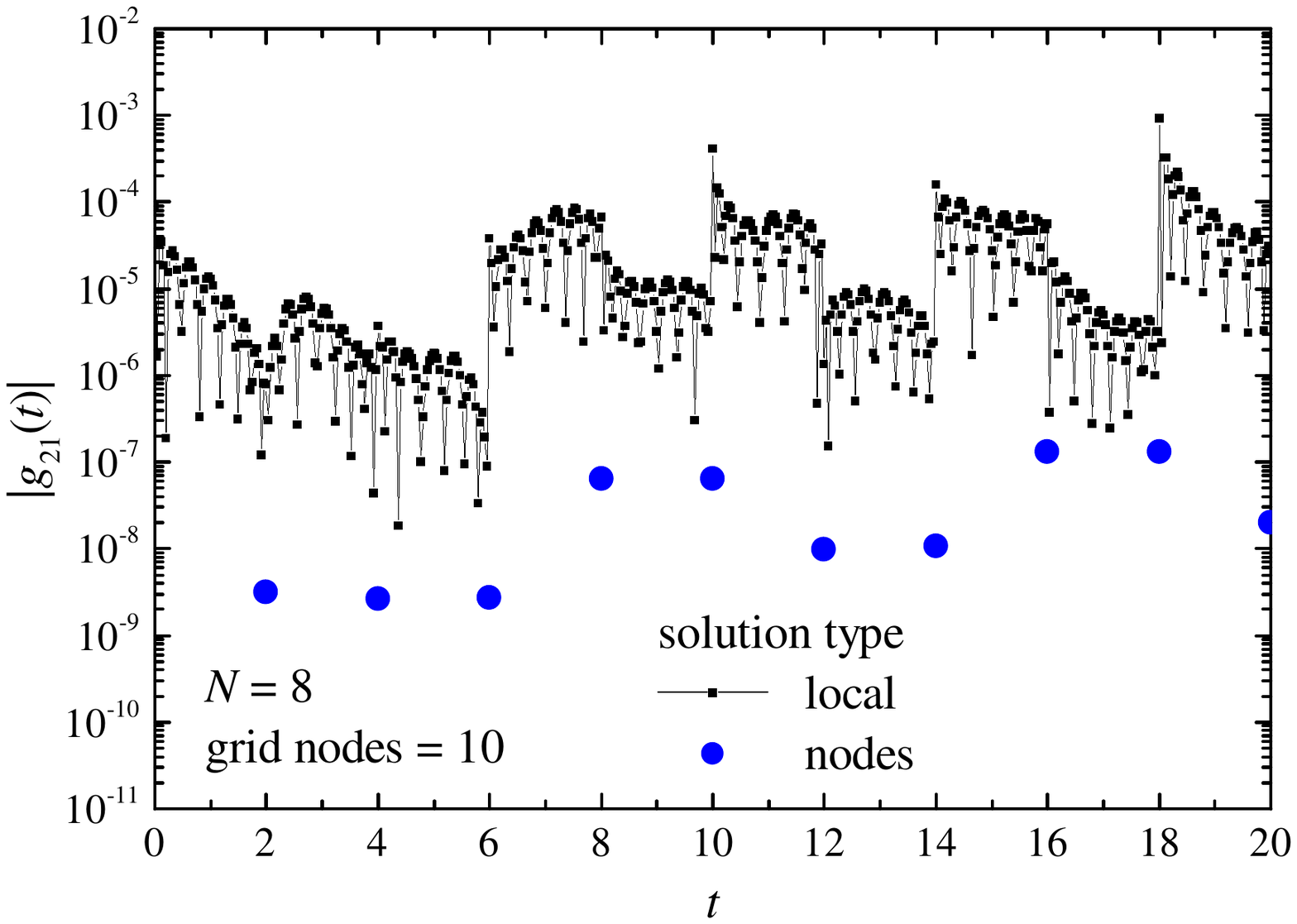}
\vspace{-8mm}\caption{\label{fig:dpend_ind1_sols_vg:d1}}
\end{subfigure}
\begin{subfigure}{0.240\textwidth}
\includegraphics[width=\textwidth]{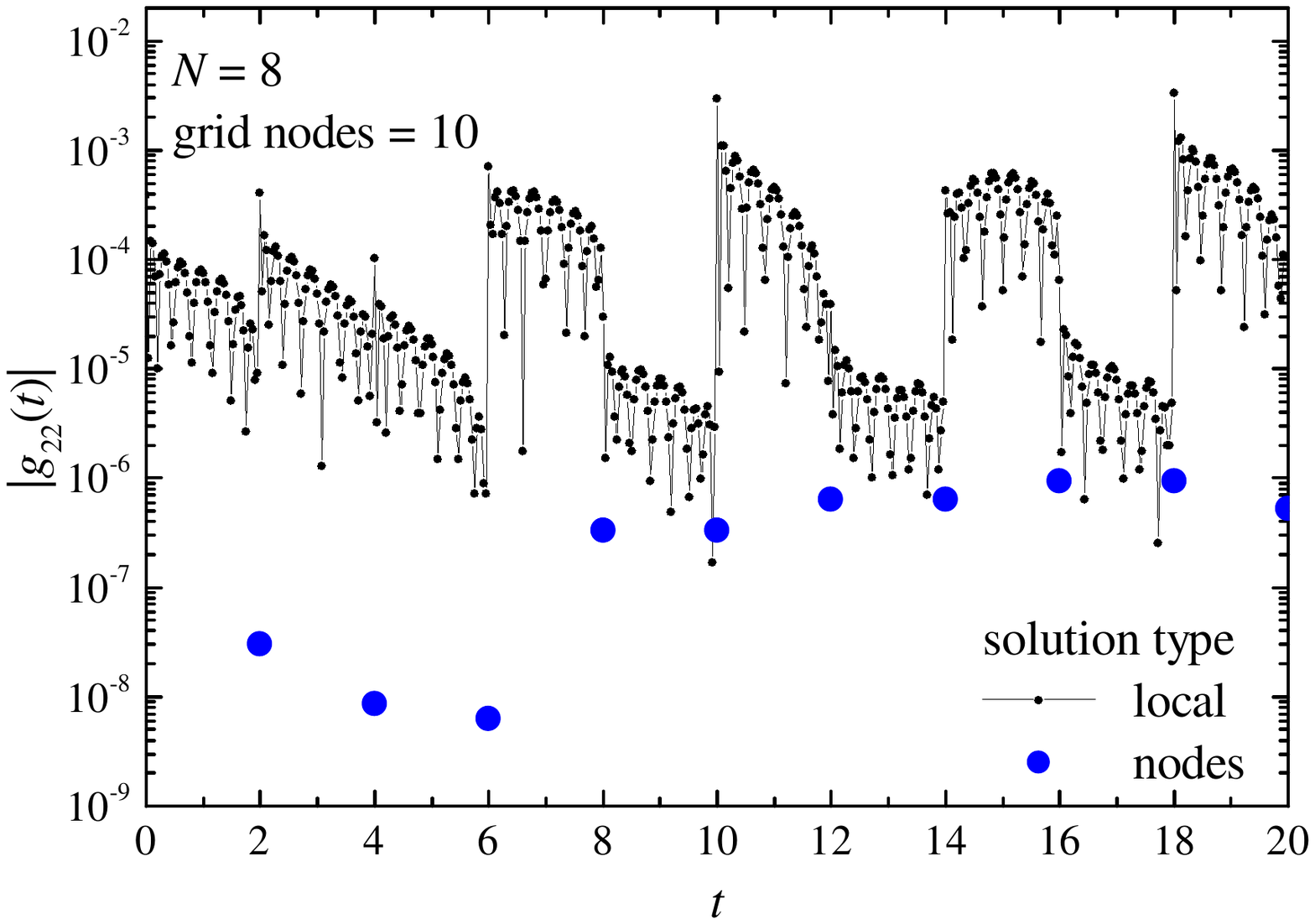}
\vspace{-8mm}\caption{\label{fig:dpend_ind1_sols_vg:d2}}
\end{subfigure}
\begin{subfigure}{0.240\textwidth}
\includegraphics[width=\textwidth]{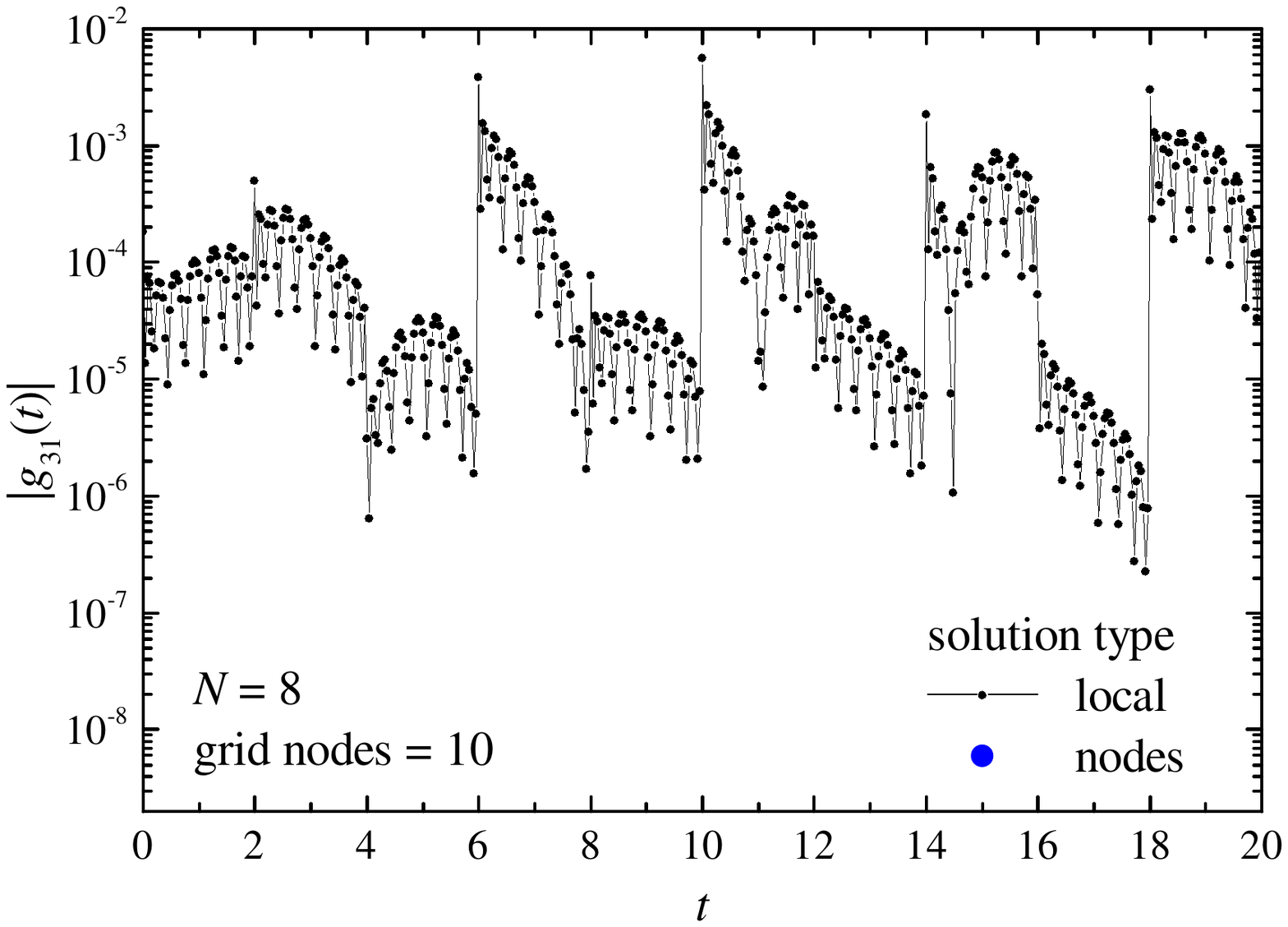}
\vspace{-8mm}\caption{\label{fig:dpend_ind1_sols_vg:d3}}
\end{subfigure}
\begin{subfigure}{0.240\textwidth}
\includegraphics[width=\textwidth]{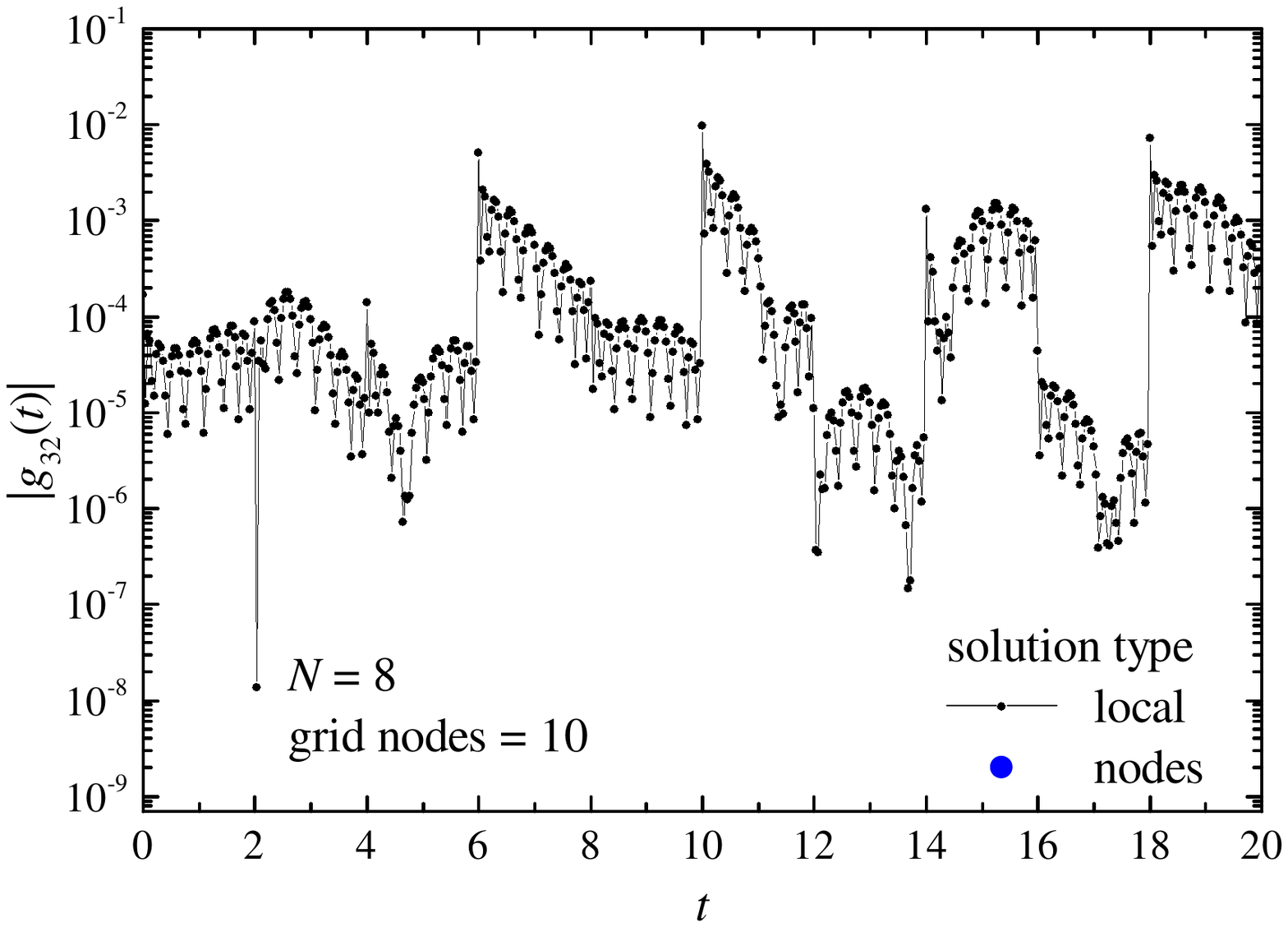}
\vspace{-8mm}\caption{\label{fig:dpend_ind1_sols_vg:d4}}
\end{subfigure}\\[2mm]
\begin{subfigure}{0.240\textwidth}
\includegraphics[width=\textwidth]{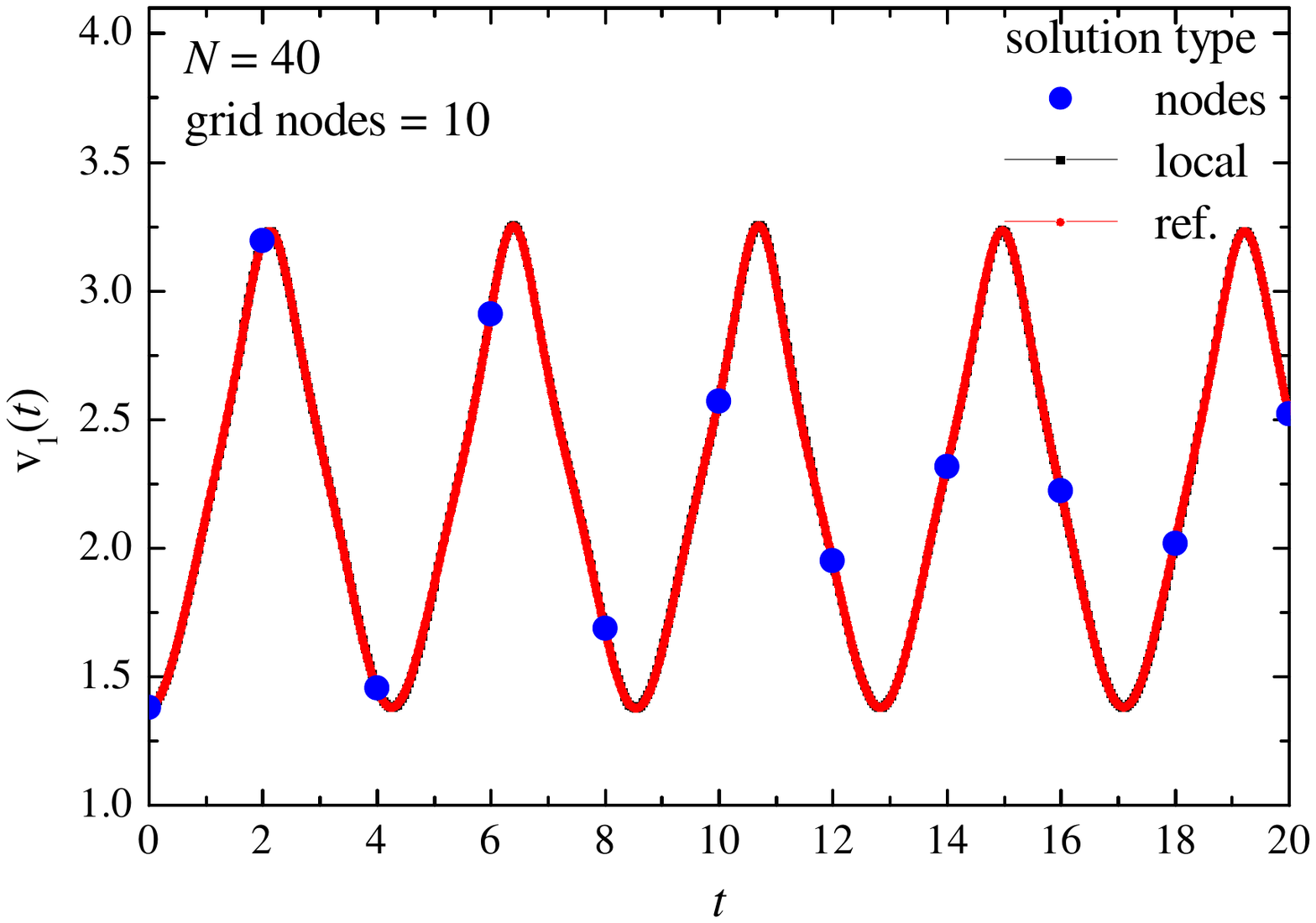}
\vspace{-8mm}\caption{\label{fig:dpend_ind1_sols_vg:e1}}
\end{subfigure}
\begin{subfigure}{0.240\textwidth}
\includegraphics[width=\textwidth]{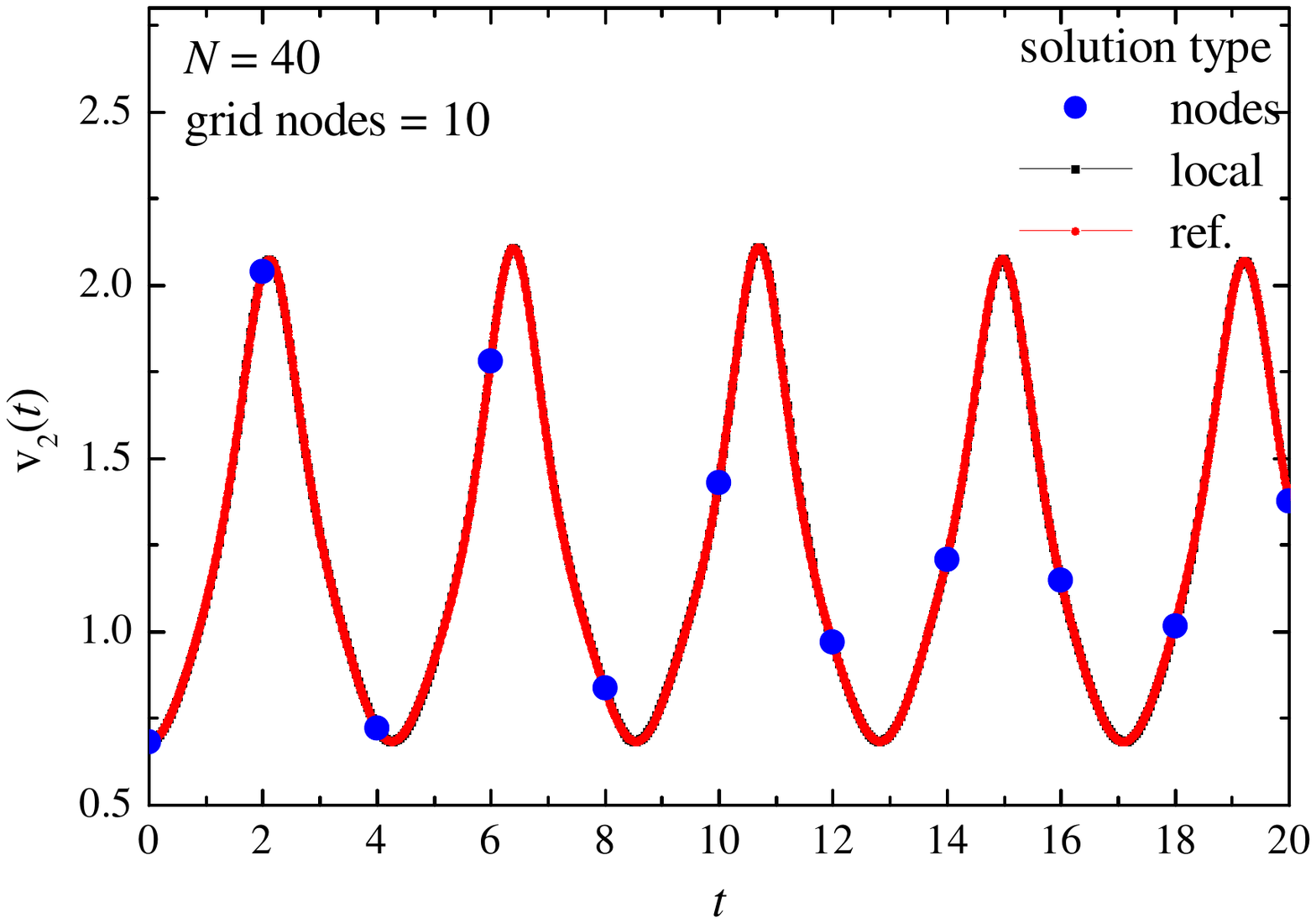}
\vspace{-8mm}\caption{\label{fig:dpend_ind1_sols_vg:e2}}
\end{subfigure}
\begin{subfigure}{0.240\textwidth}
\includegraphics[width=\textwidth]{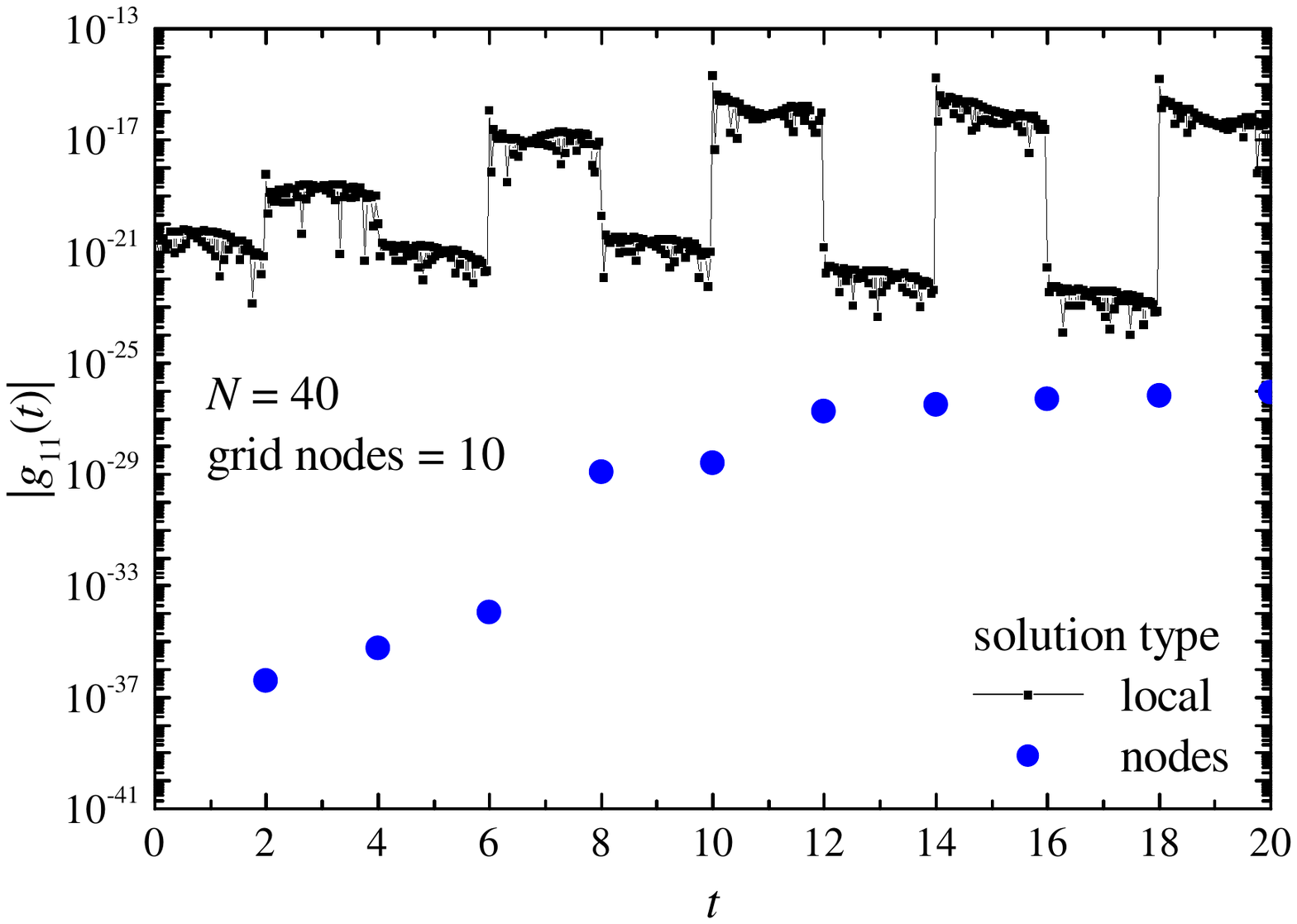}
\vspace{-8mm}\caption{\label{fig:dpend_ind1_sols_vg:e3}}
\end{subfigure}
\begin{subfigure}{0.240\textwidth}
\includegraphics[width=\textwidth]{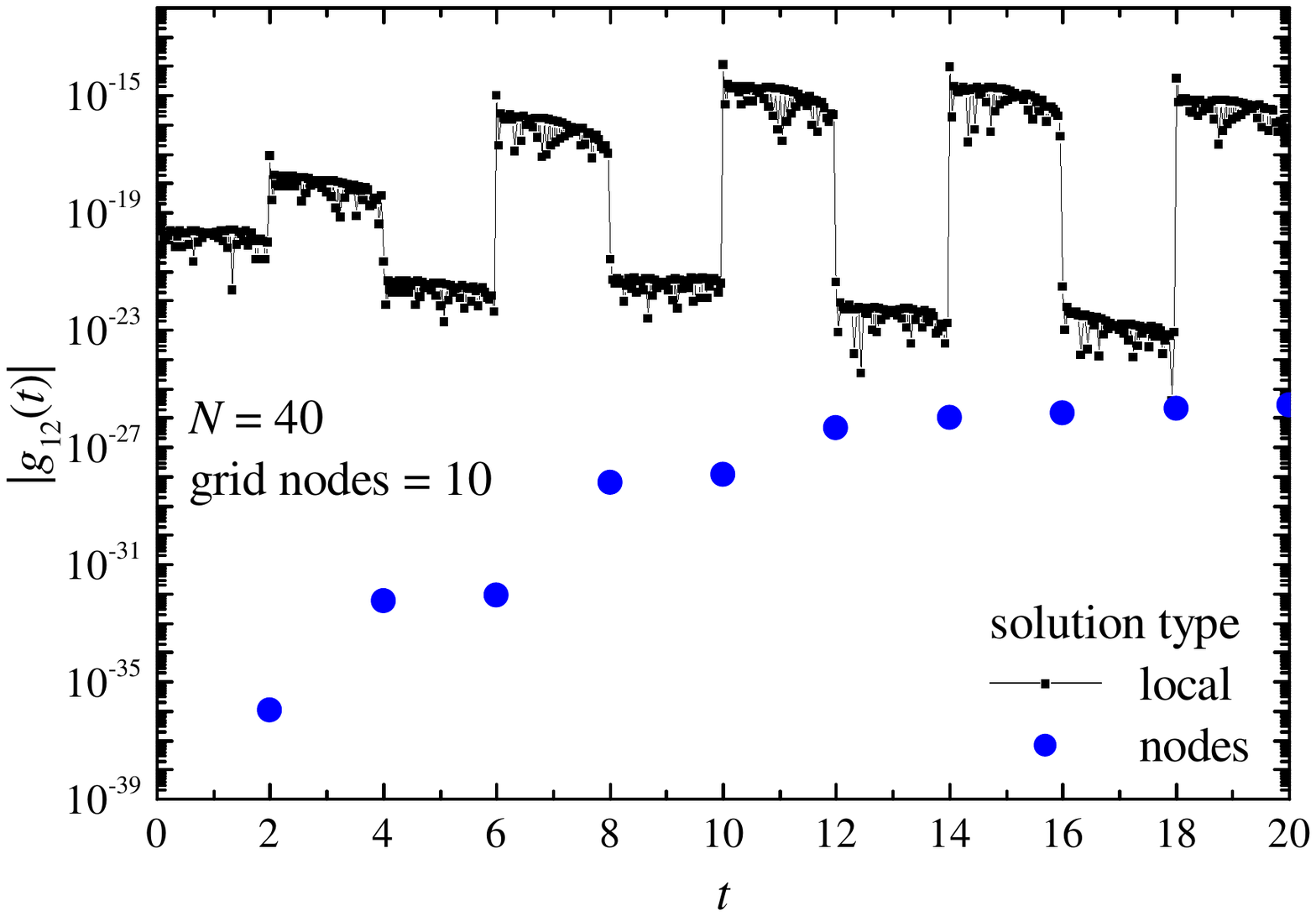}
\vspace{-8mm}\caption{\label{fig:dpend_ind1_sols_vg:e4}}
\end{subfigure}\\[2mm]
\begin{subfigure}{0.240\textwidth}
\includegraphics[width=\textwidth]{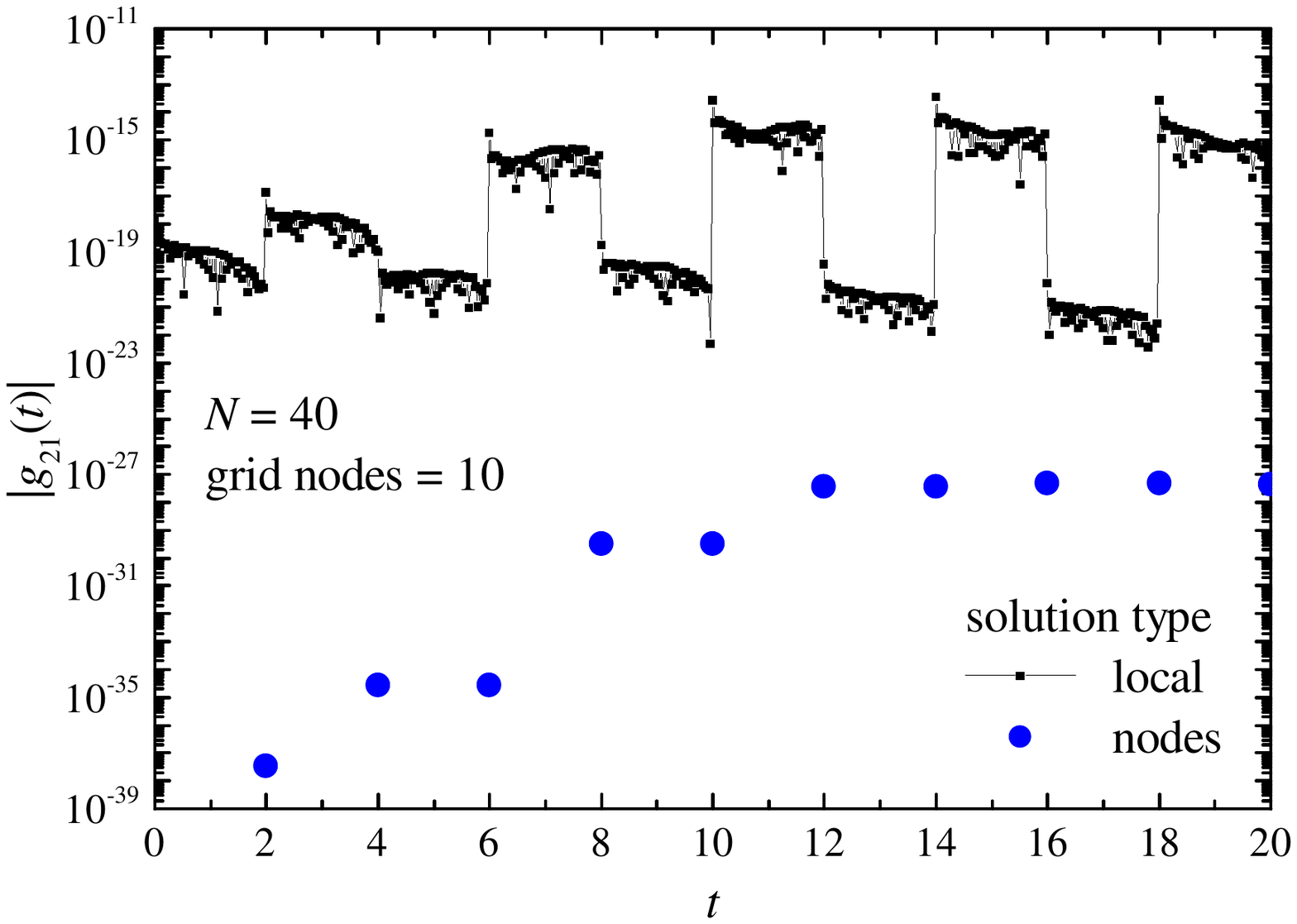}
\vspace{-8mm}\caption{\label{fig:dpend_ind1_sols_vg:f1}}
\end{subfigure}
\begin{subfigure}{0.240\textwidth}
\includegraphics[width=\textwidth]{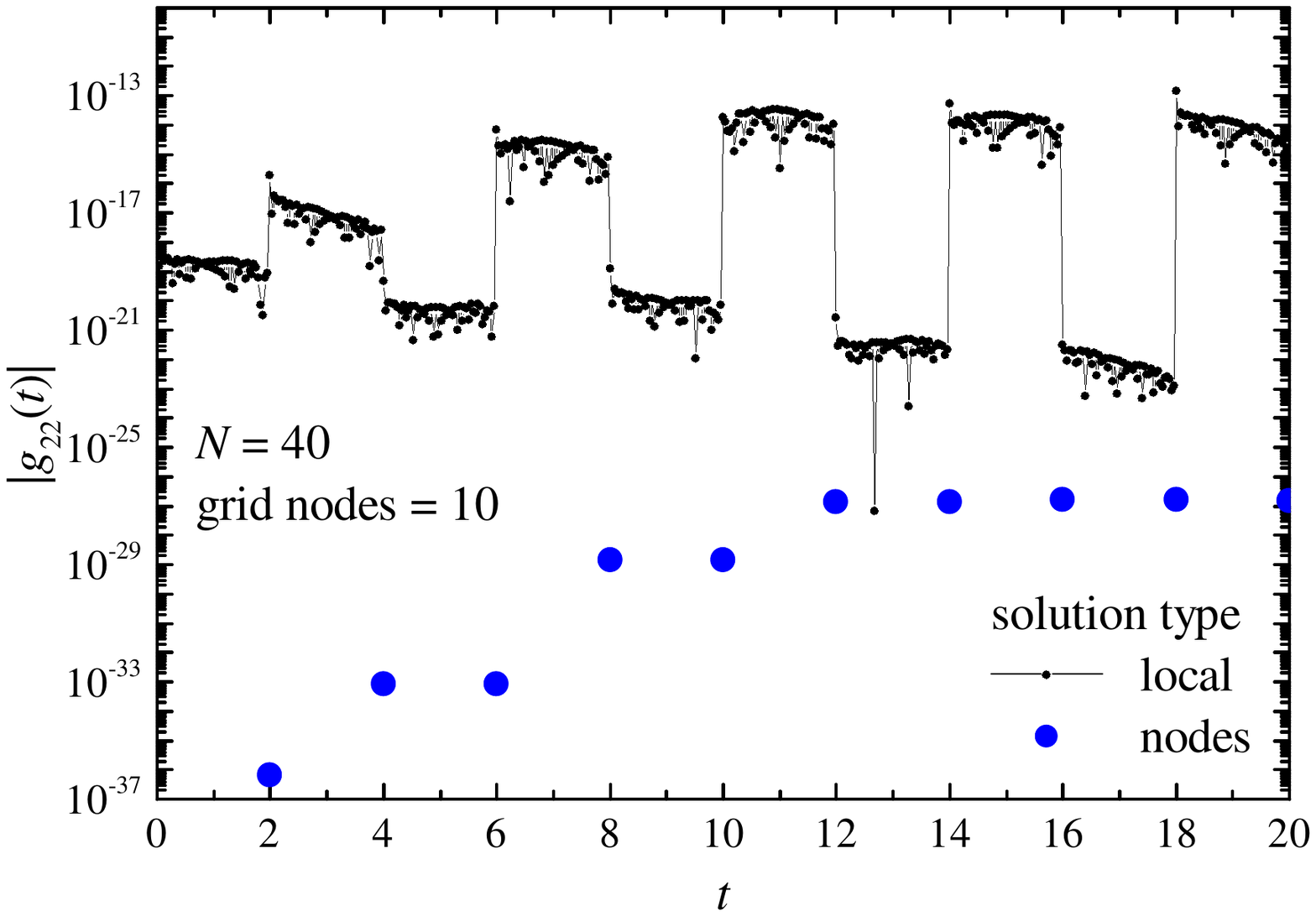}
\vspace{-8mm}\caption{\label{fig:dpend_ind1_sols_vg:f2}}
\end{subfigure}
\begin{subfigure}{0.240\textwidth}
\includegraphics[width=\textwidth]{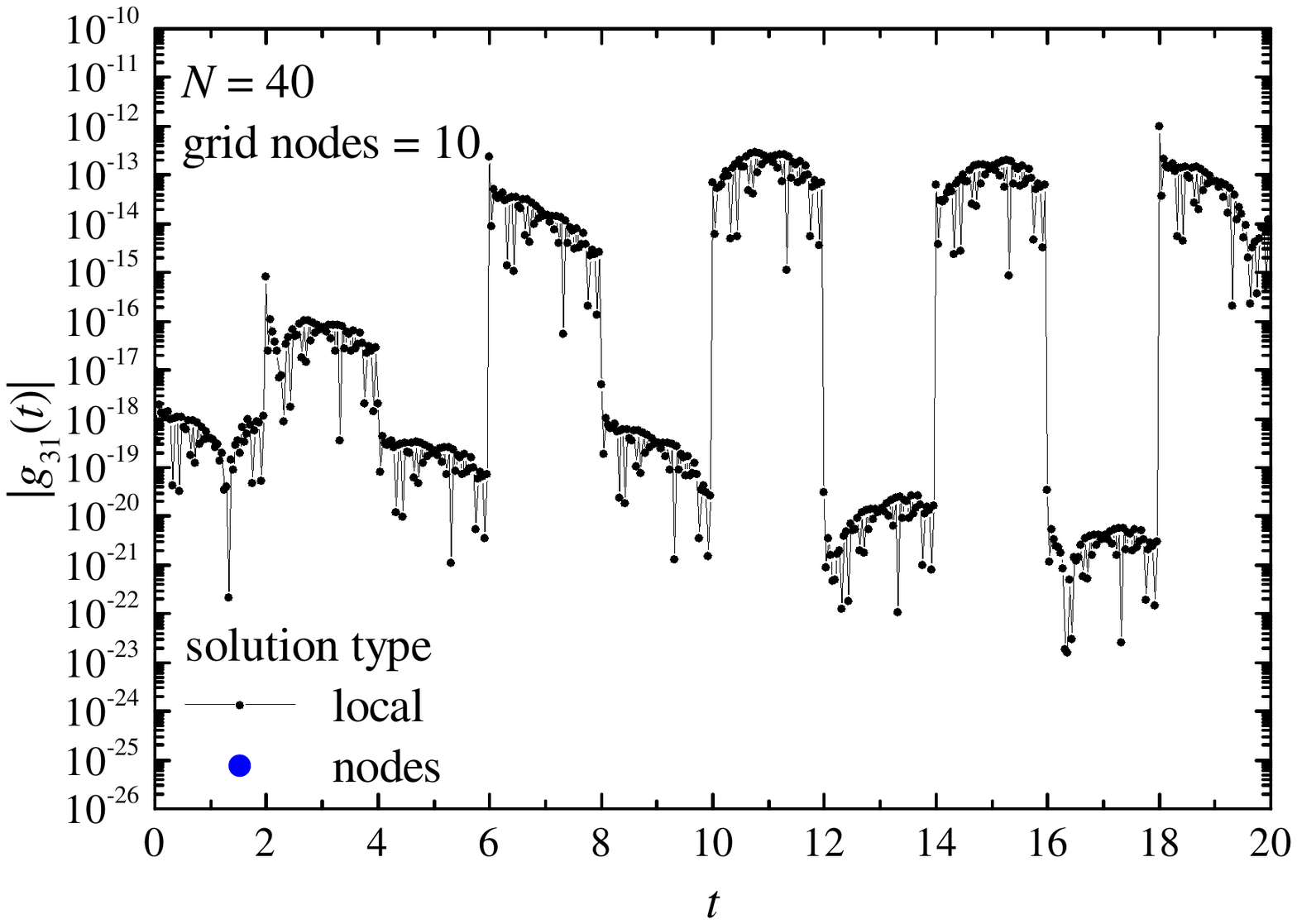}
\vspace{-8mm}\caption{\label{fig:dpend_ind1_sols_vg:f3}}
\end{subfigure}
\begin{subfigure}{0.240\textwidth}
\includegraphics[width=\textwidth]{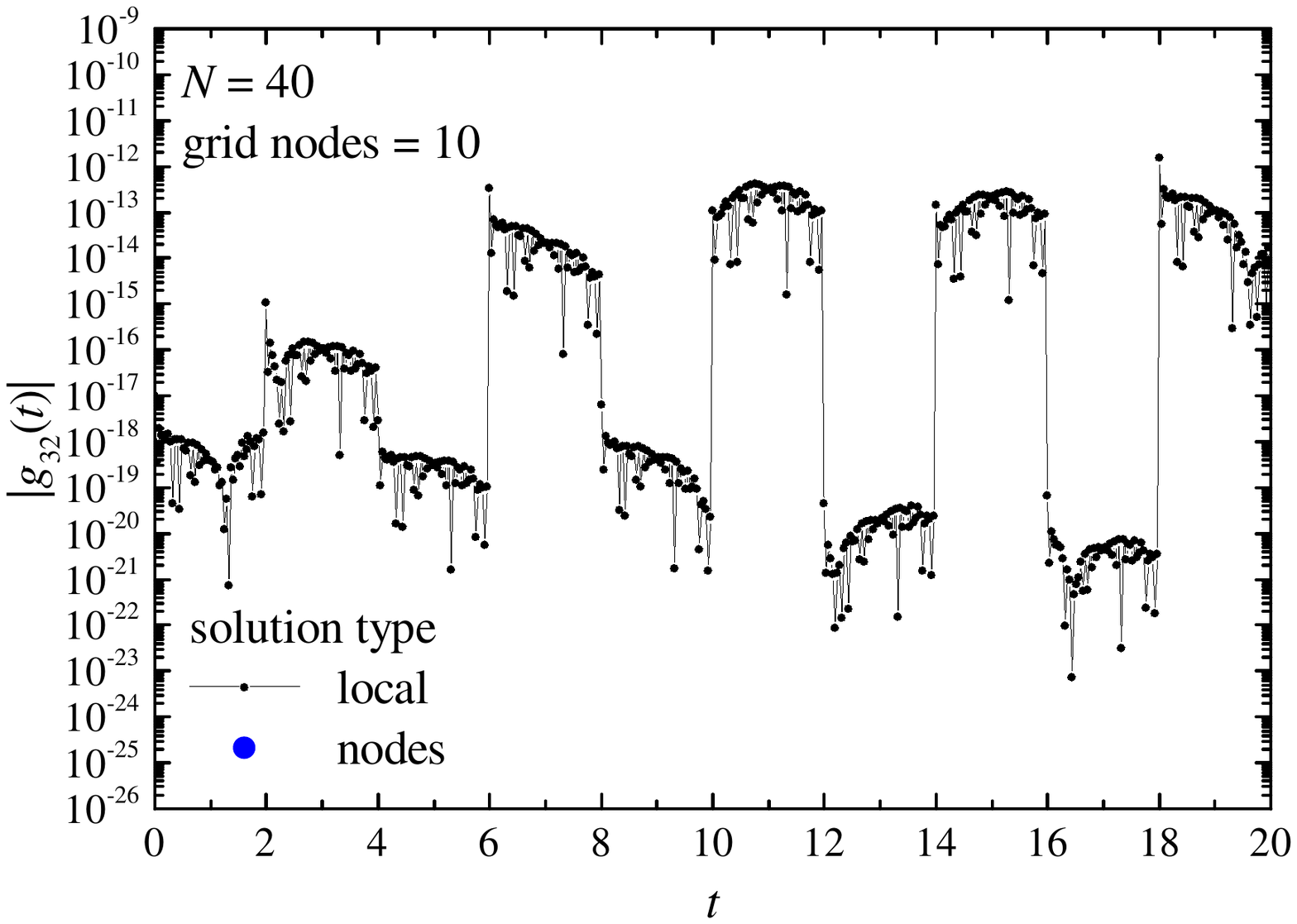}
\vspace{-8mm}\caption{\label{fig:dpend_ind1_sols_vg:f4}}
\end{subfigure}\\[2mm]
\caption{%
Numerical solution of the DAE system (\ref{eq:math_dpend_dae_ind_3}) of index 1. Comparison of the solution at nodes $\mathbf{v}_{n}$, the local solution $\mathbf{v}_{L}(t)$ and the reference solution $\mathbf{v}^{\rm ref}(t)$ for components $v_{1}$ (\subref{fig:dpend_ind1_sols_vg:a1}, \subref{fig:dpend_ind1_sols_vg:c1}, \subref{fig:dpend_ind1_sols_vg:e1}), $v_{2}$ (\subref{fig:dpend_ind1_sols_vg:a2}, \subref{fig:dpend_ind1_sols_vg:c2}, \subref{fig:dpend_ind1_sols_vg:e2}), quantitative satisfiability of the conditions $g_{11} = 0$ (\subref{fig:dpend_ind1_sols_vg:a3}, \subref{fig:dpend_ind1_sols_vg:c3}, \subref{fig:dpend_ind1_sols_vg:e3}), $g_{12} = 0$ (\subref{fig:dpend_ind1_sols_vg:a4}, \subref{fig:dpend_ind1_sols_vg:c4}, \subref{fig:dpend_ind1_sols_vg:e4}), $g_{21} = 0$ (\subref{fig:dpend_ind1_sols_vg:b1}, \subref{fig:dpend_ind1_sols_vg:d1}, \subref{fig:dpend_ind1_sols_vg:f1}), $g_{22} = 0$ (\subref{fig:dpend_ind1_sols_vg:b2}, \subref{fig:dpend_ind1_sols_vg:d2}, \subref{fig:dpend_ind1_sols_vg:f2}), $g_{31} = 0$ (\subref{fig:dpend_ind1_sols_vg:b3}, \subref{fig:dpend_ind1_sols_vg:d3}, \subref{fig:dpend_ind1_sols_vg:f3}), $g_{32} = 0$ (\subref{fig:dpend_ind1_sols_vg:b4}, \subref{fig:dpend_ind1_sols_vg:d4}, \subref{fig:dpend_ind1_sols_vg:f4}), obtained using polynomials with degrees $N = 1$ (\subref{fig:dpend_ind1_sols_vg:a1}, \subref{fig:dpend_ind1_sols_vg:a2}, \subref{fig:dpend_ind1_sols_vg:a3}, \subref{fig:dpend_ind1_sols_vg:a4}, \subref{fig:dpend_ind1_sols_vg:b1}, \subref{fig:dpend_ind1_sols_vg:b2}, \subref{fig:dpend_ind1_sols_vg:b3}, \subref{fig:dpend_ind1_sols_vg:b4}), $N = 8$ (\subref{fig:dpend_ind1_sols_vg:c1}, \subref{fig:dpend_ind1_sols_vg:c2}, \subref{fig:dpend_ind1_sols_vg:c3}, \subref{fig:dpend_ind1_sols_vg:c4}, \subref{fig:dpend_ind1_sols_vg:d1}, \subref{fig:dpend_ind1_sols_vg:d2}, \subref{fig:dpend_ind1_sols_vg:d3}, \subref{fig:dpend_ind1_sols_vg:d4}) and $N = 40$ (\subref{fig:dpend_ind1_sols_vg:e1}, \subref{fig:dpend_ind1_sols_vg:e2}, \subref{fig:dpend_ind1_sols_vg:e3}, \subref{fig:dpend_ind1_sols_vg:e4}, \subref{fig:dpend_ind1_sols_vg:f1}, \subref{fig:dpend_ind1_sols_vg:f2}, \subref{fig:dpend_ind1_sols_vg:f3}, \subref{fig:dpend_ind1_sols_vg:f4}).
}
\label{fig:dpend_ind1_sols_vg}
\end{figure} 

\begin{figure}[h!]
\captionsetup[subfigure]{%
	position=bottom,
	font+=smaller,
	textfont=normalfont,
	singlelinecheck=off,
	justification=raggedright
}
\centering
\begin{subfigure}{0.320\textwidth}
\includegraphics[width=\textwidth]{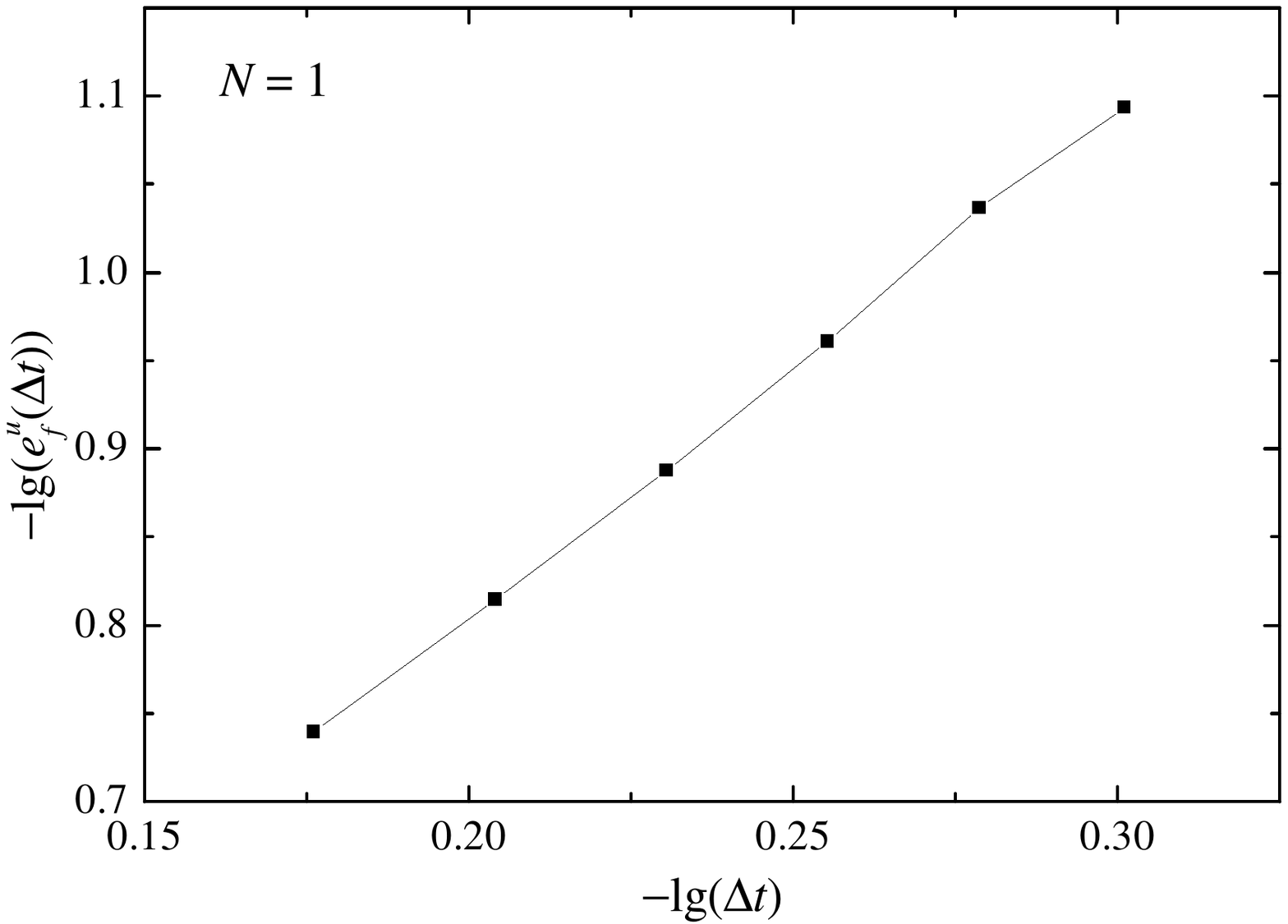}
\vspace{-8mm}\caption{\label{fig:dpend_ind1_errors:a1}}
\end{subfigure}
\begin{subfigure}{0.320\textwidth}
\includegraphics[width=\textwidth]{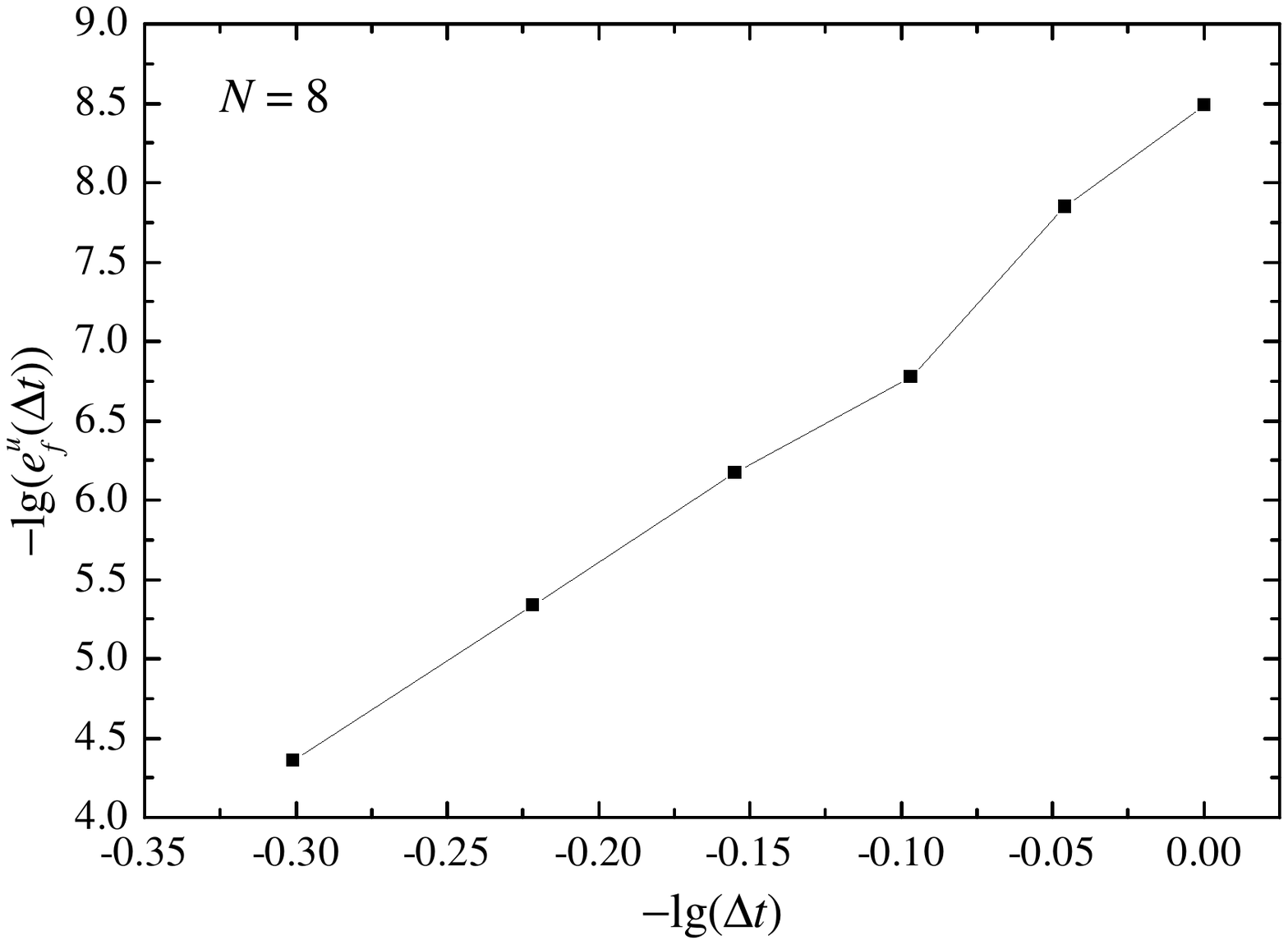}
\vspace{-8mm}\caption{\label{fig:dpend_ind1_errors:a2}}
\end{subfigure}
\begin{subfigure}{0.320\textwidth}
\includegraphics[width=\textwidth]{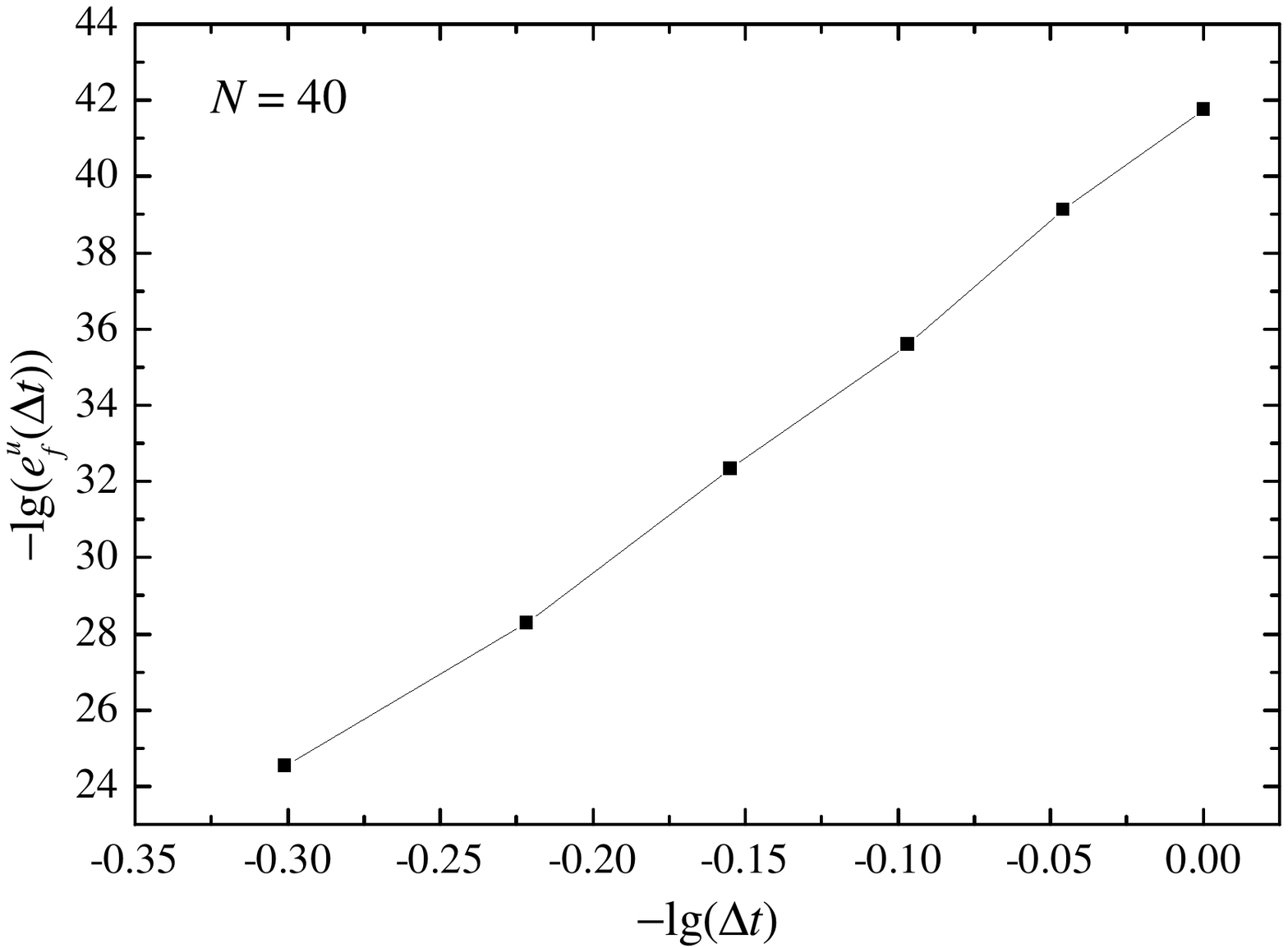}
\vspace{-8mm}\caption{\label{fig:dpend_ind1_errors:a3}}
\end{subfigure}\\[2mm]
\begin{subfigure}{0.320\textwidth}
\includegraphics[width=\textwidth]{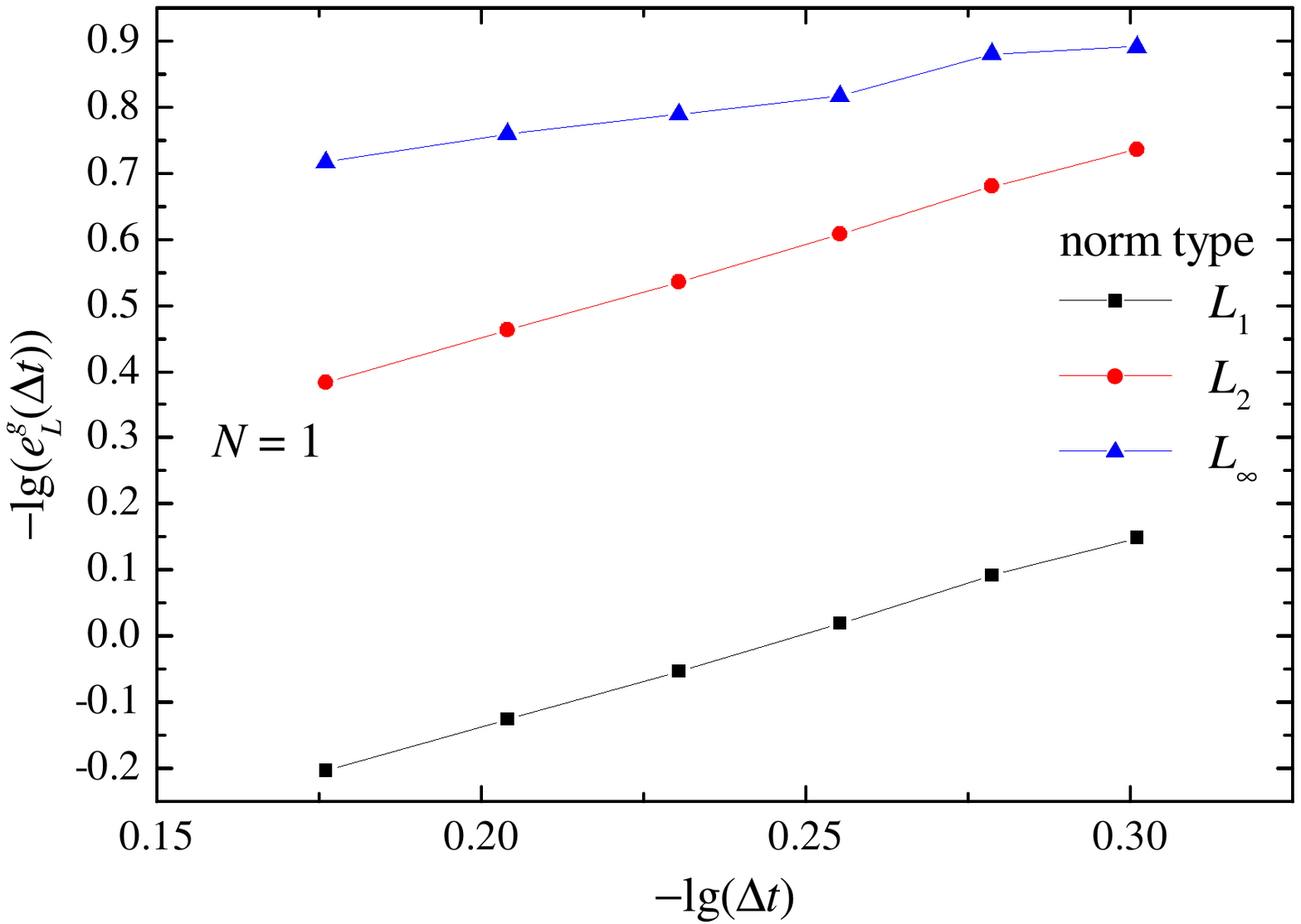}
\vspace{-8mm}\caption{\label{fig:dpend_ind1_errors:b1}}
\end{subfigure}
\begin{subfigure}{0.320\textwidth}
\includegraphics[width=\textwidth]{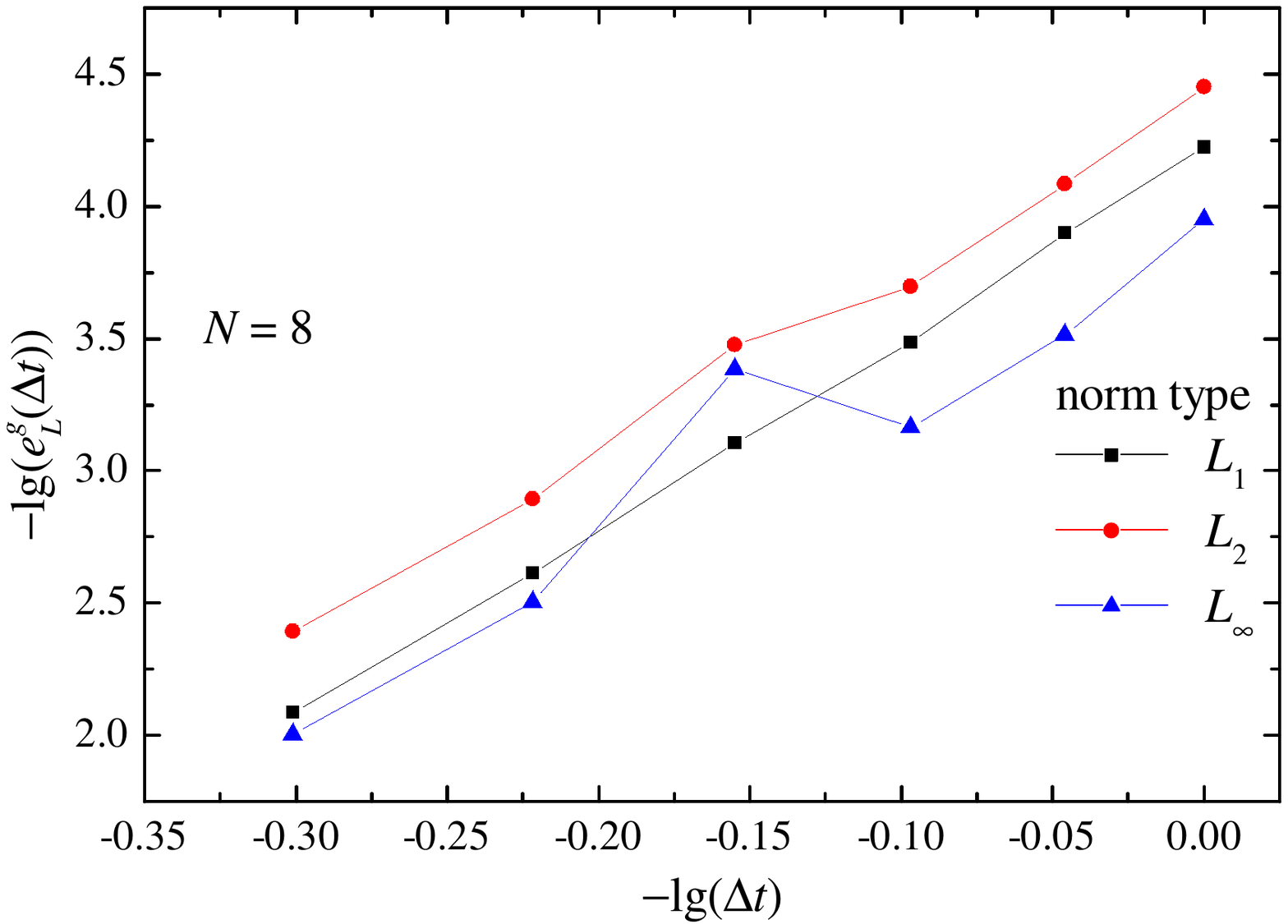}
\vspace{-8mm}\caption{\label{fig:dpend_ind1_errors:b2}}
\end{subfigure}
\begin{subfigure}{0.320\textwidth}
\includegraphics[width=\textwidth]{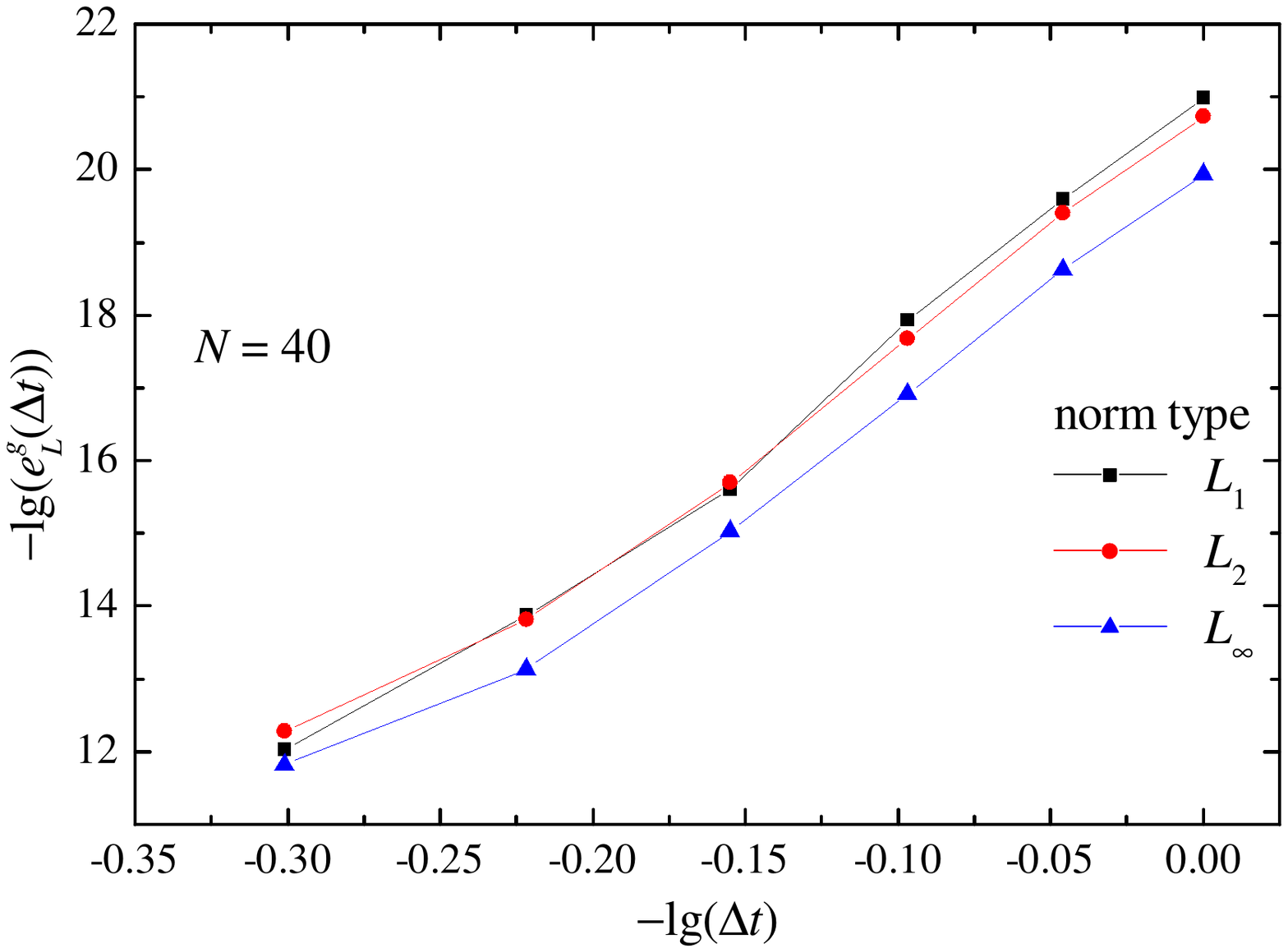}
\vspace{-8mm}\caption{\label{fig:dpend_ind1_errors:b3}}
\end{subfigure}\\[2mm]
\begin{subfigure}{0.320\textwidth}
\includegraphics[width=\textwidth]{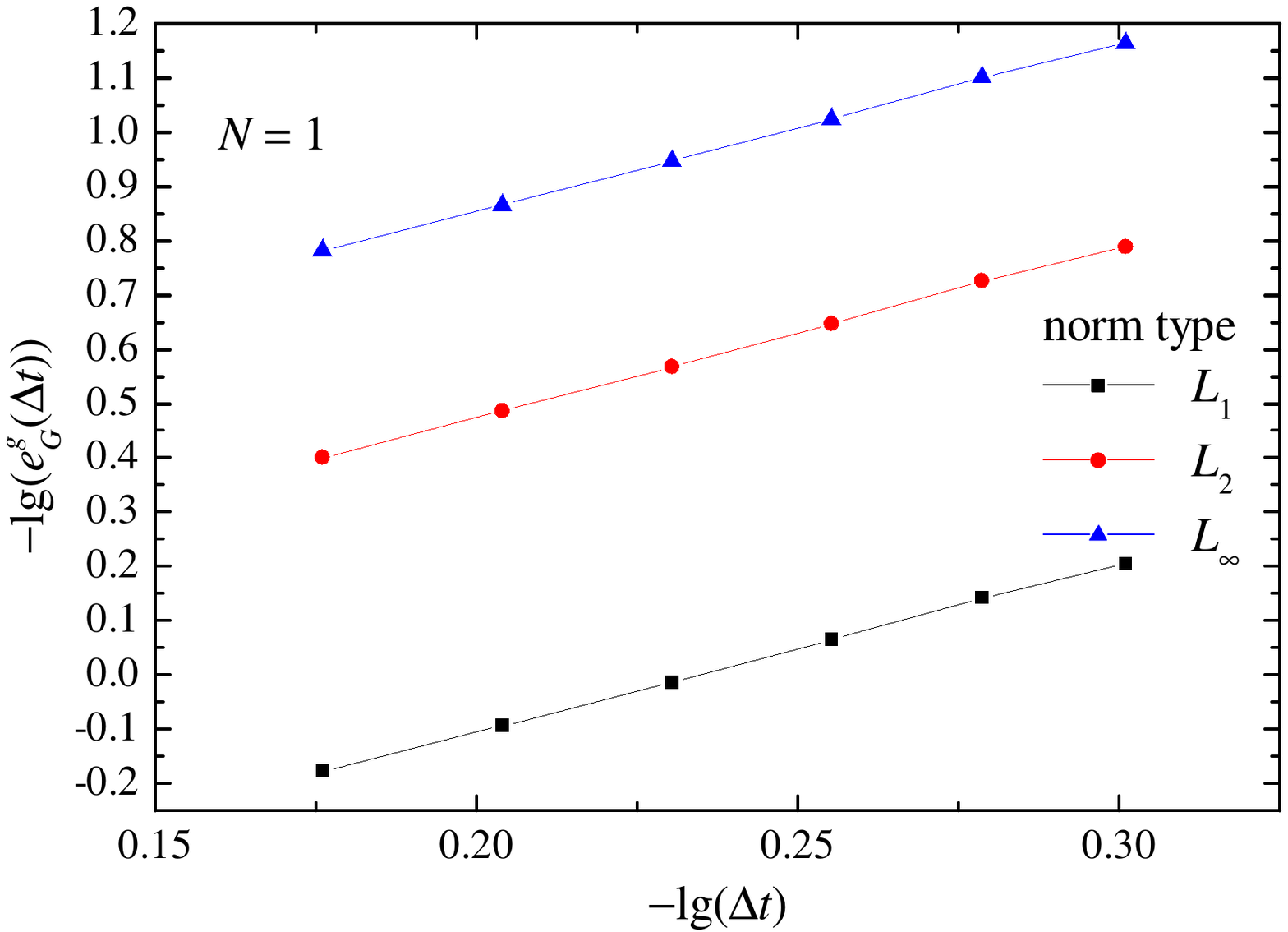}
\vspace{-8mm}\caption{\label{fig:dpend_ind1_errors:c1}}
\end{subfigure}
\begin{subfigure}{0.320\textwidth}
\includegraphics[width=\textwidth]{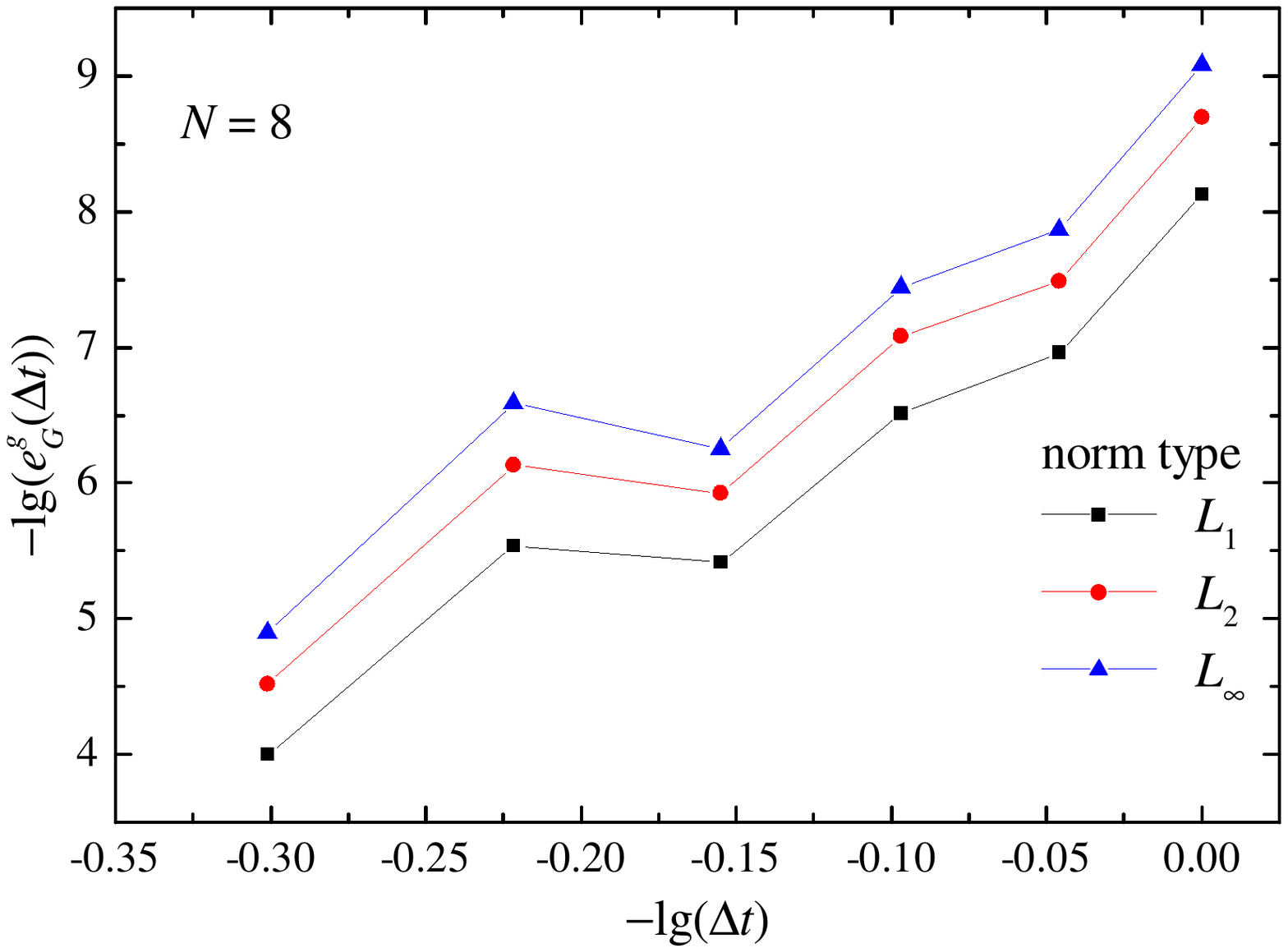}
\vspace{-8mm}\caption{\label{fig:dpend_ind1_errors:c2}}
\end{subfigure}
\begin{subfigure}{0.320\textwidth}
\includegraphics[width=\textwidth]{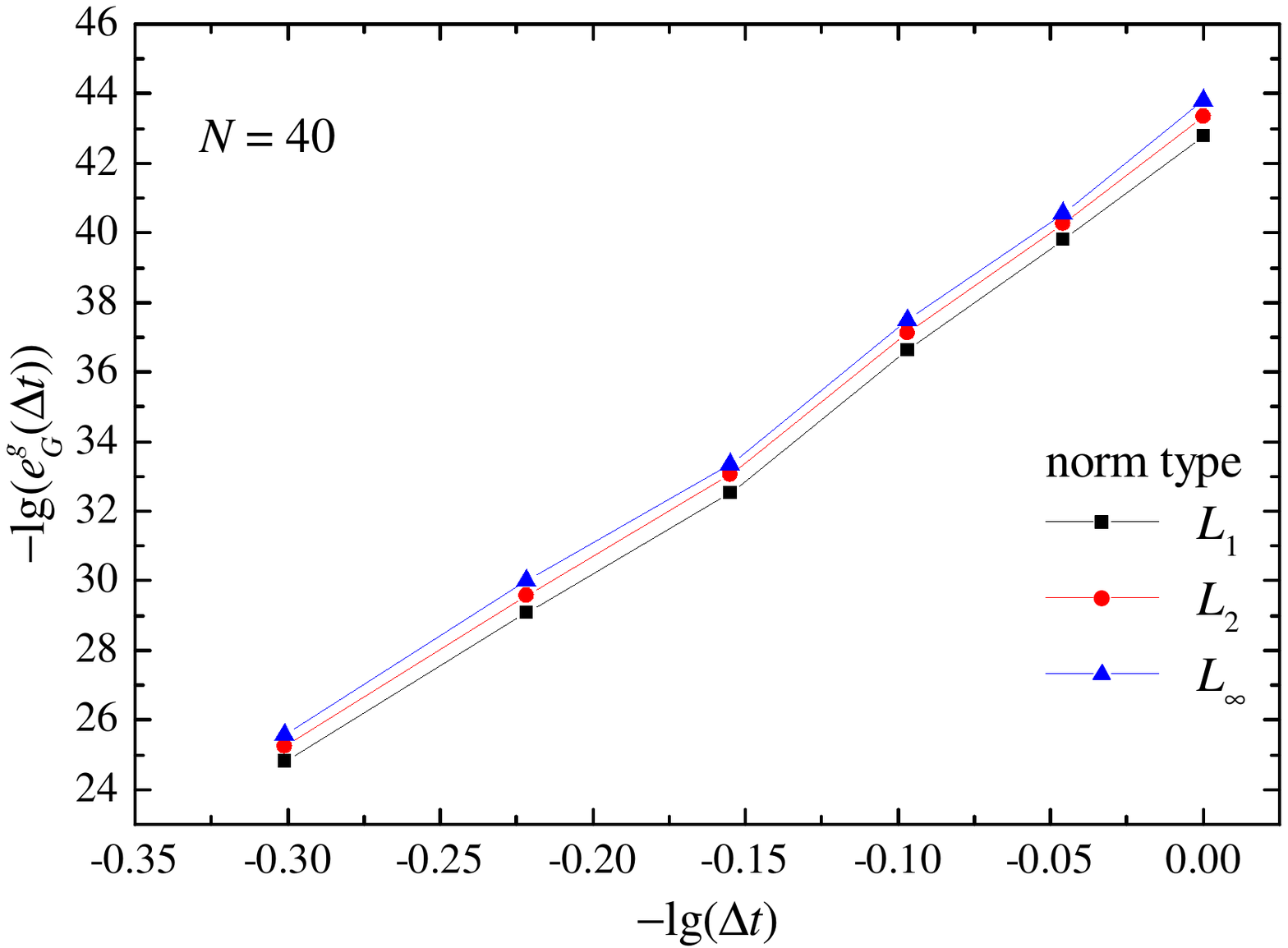}
\vspace{-8mm}\caption{\label{fig:dpend_ind1_errors:c3}}
\end{subfigure}\\[2mm]
\caption{%
Log-log plot of the dependence of the error for the solution at nodes at final time $t_{f}$ $e_{f}^{u}$ (\subref{fig:dpend_ind1_errors:a1}, \subref{fig:dpend_ind1_errors:a2}, \subref{fig:dpend_ind1_errors:a3}), global errors for the local solution $e_{L}^{g}$ (\subref{fig:dpend_ind1_errors:b1}, \subref{fig:dpend_ind1_errors:b2}, \subref{fig:dpend_ind1_errors:b3}) and the solution at nodes $e_{G}^{u}$ (\subref{fig:dpend_ind1_errors:c1}, \subref{fig:dpend_ind1_errors:c2}, \subref{fig:dpend_ind1_errors:c3}), by numerical solution of the DAE system (\ref{eq:math_dpend_dae_ind_3}) of index 1 obtained using polynomials with degrees $N = 1$ (\subref{fig:dpend_ind1_errors:a1}, \subref{fig:dpend_ind1_errors:b1}, \subref{fig:dpend_ind1_errors:c1}), $N = 8$ (\subref{fig:dpend_ind1_errors:a2}, \subref{fig:dpend_ind1_errors:b2}, \subref{fig:dpend_ind1_errors:c2}) and $N = 40$ (\subref{fig:dpend_ind1_errors:a3}, \subref{fig:dpend_ind1_errors:b3}, \subref{fig:dpend_ind1_errors:c3}).
}
\label{fig:dpend_ind1_errors}
\end{figure} 

\begin{table*}[h!]
\centering
\caption{%
Convergence orders $p_{f}$, $p_{L_{1}}$, $p_{L_{2}}$, $p_{L_{\infty}}$, calculated on final time step $t_{f}$ and in norms $L_{1}$, $L_{2}$, $L_{\infty}$ of the ADER-DG method for the DAE system (\ref{eq:math_dpend_dae_ind_3}) of index 1; $N$ is the degree of the basis polynomials $\varphi_{p}$. Orders $p^{n, u}$ are calculated for solution $\mathbf{u}_{n}$; orders $p^{n, g}$ and $p^{l, g}$ --- for the conditions $\mathbf{g} = 0$ on the numerical solution at the nodes $(\mathbf{u}_{n}, \mathbf{v}_{n})$ and on the local solution $(\mathbf{u}_{L}, \mathbf{v}_{L})$. The theoretical values of convergence order $p_{\rm th.}^{n} = 2N+1$ and $p_{\rm th.}^{l} = N+1$ are applicable for the ADER-DG method for ODE problems, and are presented for comparison.
}
\label{tab:conv_orders_dpend_ind1}
\begin{tabular}{@{}|l|l|lll|c|lll|c|@{}}
\toprule
$N$ & $p_{f}^{n, u}$ & $p_{L_{1}}^{n, g}$ & $p_{L_{2}}^{n, g}$ & $p_{L_{\infty}}^{n, g}$ & $p_{\rm th.}^{n}$ & $p_{L_{1}}^{l, g}$ & $p_{L_{2}}^{l, g}$ & $p_{L_{\infty}}^{l, g}$ & $p_{\rm th.}^{l}$ \\
\midrule
$1$	&	$2.87$	&	$3.10$	&	$3.12$	&	$3.09$	&	$3$	&	$2.84$	&	$2.85$	&	$1.44$	&	$2$\\
$2$	&	$3.83$	&	$4.28$	&	$4.27$	&	$4.36$	&	$5$	&	$3.12$	&	$3.12$	&	$3.19$	&	$3$\\
$3$	&	$7.91$	&	$7.65$	&	$7.55$	&	$7.30$	&	$7$	&	$3.75$	&	$3.60$	&	$3.19$	&	$4$\\
$4$	&	$7.93$	&	$8.81$	&	$8.92$	&	$8.65$	&	$9$	&	$4.95$	&	$4.99$	&	$4.44$	&	$5$\\
$5$	&	$9.86$	&	$7.98$	&	$7.70$	&	$6.95$	&	$11$	&	$5.86$	&	$5.52$	&	$4.76$	&	$6$\\
$6$	&	$10.40$	&	$10.27$	&	$10.39$	&	$10.37$	&	$13$	&	$5.77$	&	$5.35$	&	$4.63$	&	$7$\\
$7$	&	$11.75$	&	$12.44$	&	$12.57$	&	$12.68$	&	$15$	&	$6.83$	&	$6.36$	&	$5.86$	&	$8$\\
$8$	&	$13.67$	&	$12.32$	&	$12.37$	&	$12.26$	&	$17$	&	$7.12$	&	$6.74$	&	$6.06$	&	$9$\\
$9$	&	$15.81$	&	$14.90$	&	$15.03$	&	$15.08$	&	$19$	&	$7.98$	&	$7.55$	&	$6.93$	&	$10$\\
$10$	&	$17.62$	&	$15.80$	&	$15.78$	&	$15.60$	&	$21$	&	$8.59$	&	$8.00$	&	$7.10$	&	$11$\\
$11$	&	$18.68$	&	$16.92$	&	$17.04$	&	$17.10$	&	$23$	&	$9.56$	&	$9.07$	&	$8.16$	&	$12$\\
$12$	&	$19.14$	&	$19.68$	&	$19.72$	&	$19.70$	&	$25$	&	$9.90$	&	$9.05$	&	$7.60$	&	$13$\\
$13$	&	$20.40$	&	$19.08$	&	$19.26$	&	$19.50$	&	$27$	&	$10.97$	&	$10.53$	&	$9.62$	&	$14$\\
$14$	&	$20.71$	&	$22.01$	&	$22.16$	&	$22.18$	&	$29$	&	$11.38$	&	$10.54$	&	$9.31$	&	$15$\\
$15$	&	$22.30$	&	$21.87$	&	$21.84$	&	$21.63$	&	$31$	&	$12.42$	&	$12.08$	&	$11.21$	&	$16$\\
$16$	&	$22.20$	&	$23.93$	&	$24.11$	&	$24.29$	&	$33$	&	$12.80$	&	$12.32$	&	$11.81$	&	$17$\\
$17$	&	$24.79$	&	$24.71$	&	$24.74$	&	$24.58$	&	$35$	&	$13.96$	&	$13.49$	&	$12.76$	&	$18$\\
$18$	&	$26.31$	&	$26.20$	&	$26.29$	&	$26.41$	&	$37$	&	$14.45$	&	$14.00$	&	$13.68$	&	$19$\\
$19$	&	$28.93$	&	$27.73$	&	$27.88$	&	$27.82$	&	$39$	&	$15.20$	&	$14.57$	&	$13.94$	&	$20$\\
$20$	&	$31.02$	&	$28.78$	&	$28.76$	&	$28.80$	&	$41$	&	$15.99$	&	$15.43$	&	$15.00$	&	$21$\\
\midrule
$25$	&	$36.26$	&	$35.62$	&	$35.62$	&	$35.55$	&	$51$	&	$19.55$	&	$19.09$	&	$18.58$	&	$26$\\
$30$	&	$44.05$	&	$44.00$	&	$44.35$	&	$44.67$	&	$61$	&	$23.04$	&	$22.28$	&	$21.72$	&	$31$\\
$35$	&	$52.42$	&	$51.81$	&	$52.01$	&	$52.14$	&	$71$	&	$26.90$	&	$25.85$	&	$25.14$	&	$36$\\
$40$	&	$58.07$	&	$60.03$	&	$60.24$	&	$60.46$	&	$81$	&	$30.49$	&	$28.96$	&	$27.94$	&	$41$\\
\bottomrule
\end{tabular}
\end{table*}

It is necessary to note an interesting property of the numerical solution --- an algebraic equations $g_{11} = 0$ and $g_{12} = 0$, which is solved numerically explicitly, is satisfied exactly (within the limits of the accuracy of representing real numbers by floating-point numbers) for the solution at the nodes $(\mathbf{u}_{n}, \mathbf{v}_{n})$. Therefore, when separately solving the DAE system (\ref{eq:math_dpend_dae_ind_3}) of index 3, index 2 and index 1, the errors were calculated that determine the satisfaction of both algebraic equations $g_{11} = 0$ and $g_{12} = 0$, $g_{21} = 0$ and $g_{22} = 0$, $g_{31} = 0$ and $g_{32} = 0$ were calculated --- it was expected that the numerical solution at the nodes $(\mathbf{u}_{n}, \mathbf{v}_{n})$ would exactly satisfy only that algebraic equation, which is explicitly included in the solved DAE system.

The obtained results of the numerical solution of DAE system (\ref{eq:math_dpend_dae_ind_3}) of index 3 are presented in Figs.~\ref{fig:dpend_ind3_sols_u}, \ref{fig:dpend_ind3_sols_vg}, \ref{fig:dpend_ind3_errors} and in Table~\ref{tab:conv_orders_dpend_ind3}. Fig.~\ref{fig:dpend_ind3_sols_u} shows a comparison of the numerical solution at the nodes for differential variables $\mathbf{u}_{n}$, the numerical local solution for differential variables $\mathbf{u}_{L}$ and the reference solution for differential variables $\mathbf{u}^{\rm ref}$. Fig.~\ref{fig:dpend_ind3_sols_vg} shows a comparison of the numerical solution at the nodes for algebraic variables $\mathbf{v}_{n}$, the numerical local solution for algebraic variables $\mathbf{v}_{L}$ and the reference solution for algebraic variables $\mathbf{v}^{\rm ref}$, the dependencies of the feasibility of conditions $g_{11} = 0$, $g_{12} = 0$, $g_{21} = 0$, $g_{22} = 0$, $g_{31} = 0$, $g_{32} = 0$ on the coordinate $t$, which allows us to quantitatively estimate the accuracy of the numerical solution, especially taking into account that the numerical solution obtained by the ADER-DG method with a local DG predictor with a high degree of polynomials $N$ very accurately corresponds to the exact analytical solution, and it is impossible to visually identify the error from the solution plots separately (in Fig.~\ref{fig:dpend_ind3_sols_u}, Fig.~\ref{fig:dpend_ind3_sols_vg} (\subref{fig:dpend_ind3_sols_vg:a1}, \subref{fig:dpend_ind3_sols_vg:c1}, \subref{fig:dpend_ind3_sols_vg:e1}) and Fig.~\ref{fig:dpend_ind3_sols_vg} (\subref{fig:dpend_ind3_sols_vg:a2}, \subref{fig:dpend_ind3_sols_vg:c2}, \subref{fig:dpend_ind3_sols_vg:e2})). Fig.~\ref{fig:dpend_ind3_errors} shows the dependencies of the global errors $e^{u}$ and $e^{g}$ of the numerical solution at the nodes $(\mathbf{u}_{n}, \mathbf{v}_{n})$ on the discretization step ${\Delta t}$, separately for differential variables $\mathbf{u}$ and the algebraic equations $\mathbf{g} = \mathbf{0}$, satisfying separately on the numerical solution at the nodes $(\mathbf{u}_{n}, \mathbf{v}_{n})$ and the local solution $(\mathbf{u}_{L}, \mathbf{v}_{L})$, on the basis of which the empirical convergence orders $p$ were calculated. For the double pendulum test case presented in this Subsection, there is no exact analytical solution, so the numerical solution was compared with the reference numerical solution (\ref{eq:math_dpend_ode}), and the convergence orders $p_{f}^{n, u}$ for the solution at nodes $\mathbf{u}_{n}$ were calculated for the differential variables $\mathbf{u}$ in the ``final'' norm (\ref{eq:norms_def}), as well as the convergence orders $p^{n, g}$ and $p^{l, g}$ for the errors of satisfying the algebraic equations $\mathbf{g} = \mathbf{0}$ for the solution at the nodes $(\mathbf{u}_{n}, \mathbf{v}_{n})$ and the local solution $(\mathbf{u}_{L}, \mathbf{v}_{L})$ separately for the norms $L_{1}$, $L_{2}$, $L_{\infty}$ (\ref{eq:norms_def}) were calculated.

The dynamic dependencies of the numerical solution at the nodes $(\mathbf{u}_{n}, \mathbf{v}_{n})$ and the local solution $(\mathbf{u}_{L}, \mathbf{v}_{L})$ shown in Fig.~\ref{fig:dpend_ind3_sols_u} for differential variables $\mathbf{u}$ and in Fig.~\ref{fig:dpend_ind3_sols_vg} (\subref{fig:dpend_ind3_sols_vg:a1}, \subref{fig:dpend_ind3_sols_vg:c1}, \subref{fig:dpend_ind3_sols_vg:e1}) and Fig.~\ref{fig:dpend_ind3_sols_vg} (\subref{fig:dpend_ind3_sols_vg:a2}, \subref{fig:dpend_ind3_sols_vg:c2}, \subref{fig:dpend_ind3_sols_vg:e2}) for algebraic variables $\mathbf{v}$ demonstrate the high accuracy of the ADER-DG numerical method with a local DG predictor for solving the DAE system (\ref{eq:math_dpend_dae_ind_3}) of index 3 describing the motion of a double pendulum. This is clearly seen from a comparison of the dynamic dependencies of the solution at the nodes $(\mathbf{u}_{n}, \mathbf{v}_{n})$, the local solution $(\mathbf{u}_{L}, \mathbf{v}_{L})$, and the reference solution $(\mathbf{u}^{\rm ref}, \mathbf{v}^{\rm ref})$. In the case of polynomials of degree $N = 1$, the numerical solution for which is shown in Fig.~\ref{fig:dpend_ind3_sols_u} (\subref{fig:dpend_ind3_sols_u:a1}, \subref{fig:dpend_ind3_sols_u:a2}, \subref{fig:dpend_ind3_sols_u:a3}, \subref{fig:dpend_ind3_sols_u:a4}, \subref{fig:dpend_ind3_sols_u:b1}, \subref{fig:dpend_ind3_sols_u:b2}, \subref{fig:dpend_ind3_sols_u:b3}, \subref{fig:dpend_ind3_sols_u:b4}) for differential variables $\mathbf{u}$ and in Fig.~\ref{fig:dpend_ind3_sols_vg} (\subref{fig:dpend_ind3_sols_vg:a1}, \subref{fig:dpend_ind3_sols_vg:a2}) for algebraic variables $\mathbf{v}$, the artifacts expected for a local solution are observed --- sufficiently large discontinuities of the local solution at the grid nodes $t_{n}$, which was already observed in previous cases. In this case, in the case of algebraic variables $\mathbf{v}$ (see Fig.~\ref{fig:dpend_ind3_sols_vg} (\subref{fig:dpend_ind3_sols_vg:a1}, \subref{fig:dpend_ind3_sols_vg:a2})), the local solution $\mathbf{v}_{L}$ exhibits significant deficiencies in accuracy when the dynamic curves $\mathbf{v}_{L}(t)$ of the numerical solution are often orthogonal to the dynamic curves of the reference solution $\mathbf{v}^{\rm ref}(t)$, as was observed in Example 3 for the DAE system (\ref{eq:hess_dae_ind_2}) of index 2, presented in Subsection~\ref{sec:2:ct:ex3}, and in Example 4 for the DAE system (\ref{eq:math_dpend_dae_ind_3}) of index 3, presented in Subsection~\ref{sec:2:ct:ex4}. Therefore, in the case of using polynomials with degree $N = 1$, it is necessary to use the local solution $\mathbf{v}_{L}$ with ``caution'' as a sufficiently accurate solution between the grid nodes $\Omega_{n}$. In the case of polynomial degrees $N > 1$, these artifacts of the numerical solution arising in the dependencies of the local solution for algebraic variables $\mathbf{v}_{L}$ on the argument $t$ do not appear. In the presented dynamic dependencies of the numerical solution at the nodes $(\mathbf{u}_{n}, \mathbf{v}_{n})$ and the local solution $(\mathbf{u}_{L}, \mathbf{v}_{L})$ for the cases of polynomial degrees of $N = 8$ and $N = 40$, even on a very coarse grid (only $10$ discretization domains), the numerical solution demonstrates excellent point-wise agreement with the reference solution of the problem $(\mathbf{u}^{\rm ref}, \mathbf{v}^{\rm ref})$. To determine the accuracy of fulfillment of algebraic equations $g_{11} = 0$, $g_{12} = 0$, $g_{21} = 0$, $g_{22} = 0$, $g_{31} = 0$, $g_{32} = 0$ the dynamic dependencies of errors $|g_{11}|$, $|g_{12}|$, $|g_{21}|$, $|g_{22}|$, $|g_{31}|$, $|g_{32}|$ were calculated, which were presented in Fig.~\ref{fig:dpend_ind3_sols_vg} (\subref{fig:dpend_ind3_sols_vg:a3}, \subref{fig:dpend_ind3_sols_vg:a4}, \subref{fig:dpend_ind3_sols_vg:b1}, \subref{fig:dpend_ind3_sols_vg:b2}, \subref{fig:dpend_ind3_sols_vg:b3}, \subref{fig:dpend_ind3_sols_vg:b4}) for polynomial degree $N = 1$, in Fig.~\ref{fig:dpend_ind3_sols_vg} (\subref{fig:dpend_ind3_sols_vg:c3}, \subref{fig:dpend_ind3_sols_vg:c4}, \subref{fig:dpend_ind3_sols_vg:d1}, \subref{fig:dpend_ind3_sols_vg:d2}, \subref{fig:dpend_ind3_sols_vg:d3}, \subref{fig:dpend_ind3_sols_vg:d4}) for polynomial degree $N = 8$ and in Fig.~\ref{fig:dpend_ind3_sols_vg} (\subref{fig:dpend_ind3_sols_vg:e3}, \subref{fig:dpend_ind3_sols_vg:e4}, \subref{fig:dpend_ind3_sols_vg:f1}, \subref{fig:dpend_ind3_sols_vg:f2}, \subref{fig:dpend_ind3_sols_vg:f3}, \subref{fig:dpend_ind3_sols_vg:f4}) for polynomial degree $N = 40$. The presented dependencies $|g_{11}|$ and $|g_{12}|$ demonstrate that algebraic equations $g_{11} = 0$ and $g_{12} = 0$ are satisfied exactly for the solution at the nodes (of course, within the accuracy of representation of real numbers with a floating point), therefore, they are not presented in Fig.~\ref{fig:dpend_ind3_sols_vg}, including in the case of polynomial degree $N = 1$, when the numerical solution $(\mathbf{u}_{n}, \mathbf{v}_{n})$, $(\mathbf{u}_{L}, \mathbf{v}_{L})$ and the reference solution $(\mathbf{u}^{\rm ref}, \mathbf{v}^{\rm ref})$ demonstrate not very high point-wise agreement. Errors $|g_{11}|$ and $|g_{12}|$ for the local solution $(\mathbf{u}_{L}, \mathbf{v}_{L})$ are non-zero and are in the characteristic range of $10^{-1}$--$10^{-7}$ in the case of polynomial degree $N = 1$, $10^{-4}$--$10^{-8}$ in the case of polynomial degree $N = 8$ and $10^{-14}$--$10^{-24}$ in the case of polynomial degree $N = 40$. Errors $|g_{21}|$, $|g_{22}|$ and $|g_{31}|$, $|g_{32}|$ in satisfying algebraic equations $g_{21} = 0$, $g_{22} = 0$ and $g_{31} = 0$, $g_{32} = 0$, which are derivatives of algebraic equations $g_{11} = 0$ and $g_{12} = 0$ (with single and double differentiation, respectively) have characteristic values of $10^{-2}$--$10^{-6}$ and $10^{0}$--$10^{-4}$, respectively, in the case of polynomial degree $N = 1$, $10^{-3}$--$10^{-8}$ and $10^{-1}$--$10^{-5}$ --- in the case of polynomial degree $N = 8$, $10^{-13}$--$10^{-22}$ and $10^{-10}$--$10^{-19}$ --- in the case of polynomial degree $N = 40$. In this case, errors $|g_{21}|$ and $|g_{22}|$, $|g_{31}|$ and $|g_{32}|$ for the solution at the nodes $(\mathbf{u}_{n}, \mathbf{v}_{n})$ are larger than the corresponding errors for the local solution $(\mathbf{u}_{L}, \mathbf{v}_{L})$ at almost all points $t$. This feature of the numerical solution was already observed in Example 3 for the DAE system (\ref{eq:hess_dae_ind_2}) of index 2, presented in Subsection~\ref{sec:2:ct:ex3}, and was especially clearly observed in Example 4 for the DAE system (\ref{eq:math_dpend_dae_ind_3}) of index 3, presented in Subsection~\ref{sec:2:ct:ex4}. In Fig.~\ref{fig:dpend_ind3_errors} the dependencies of global errors $e_{f}^{u}$, $e_{L}^{g}$, $e_{G}^{u}$ on the discretization step ${\Delta t}$ are presented.  The presented dependencies demonstrate the power law $e(\Delta t) \sim (\Delta t)^{p}$, from which the empirical convergence orders $p$ were calculated. However, in comparison with the previously studied test examples, in this case the approximations are not so good. The orders of convergence calculated on the basis of the presented dependencies were presented in Table~\ref{tab:conv_orders_dpend_ind3}, where the expected values $p_{\rm nodes} = 2N+1$, $p_{\rm local} = N+1$ (\ref{eq:expect_orders}) were also presented for comparison. The obtained values of convergence of orders demonstrate significantly lower values compared to the expected values (\ref{eq:expect_orders}), which may be a consequence of the high index of the studied DAE system --- usually a decrease in the empirical convergence order is observed compared to the expected convergence order of the numerical method, which occurs when solving the initial value problem for the ODE system~\cite{Hairer_book_2}.

The obtained results of the numerical solution of DAE system (\ref{eq:math_dpend_dae_ind_3}) of index 2 are presented in Figs.~\ref{fig:dpend_ind2_sols_u}, \ref{fig:dpend_ind2_sols_vg}, \ref{fig:dpend_ind2_errors} and in Table~\ref{tab:conv_orders_dpend_ind2}. Fig.~\ref{fig:dpend_ind2_sols_u} shows a comparison of the numerical solution at the nodes for differential variables $\mathbf{u}_{n}$, the numerical local solution for differential variables $\mathbf{u}_{L}$ and the reference solution for differential variables $\mathbf{u}^{\rm ref}$. Fig.~\ref{fig:dpend_ind2_sols_vg} shows a comparison of the numerical solution at the nodes for algebraic variables $\mathbf{v}_{n}$, the numerical local solution for algebraic variables $\mathbf{v}_{L}$ and the reference solution for algebraic variables $\mathbf{v}^{\rm ref}$, the dependencies of the feasibility of conditions $g_{11} = 0$, $g_{12} = 0$, $g_{21} = 0$, $g_{22} = 0$, $g_{31} = 0$, $g_{32} = 0$ on the coordinate $t$, which allows us to quantitatively estimate the accuracy of the numerical solution, especially taking into account that the numerical solution obtained by the ADER-DG method with a local DG predictor with a high degree of polynomials $N$ very accurately corresponds to the exact analytical solution, and it is impossible to visually identify the error from the solution plots separately (in Fig.~\ref{fig:dpend_ind2_sols_u}, Fig.~\ref{fig:dpend_ind2_sols_vg} (\subref{fig:dpend_ind2_sols_vg:a1}, \subref{fig:dpend_ind2_sols_vg:c1}, \subref{fig:dpend_ind2_sols_vg:e1}) and Fig.~\ref{fig:dpend_ind2_sols_vg} (\subref{fig:dpend_ind2_sols_vg:a2}, \subref{fig:dpend_ind2_sols_vg:c2}, \subref{fig:dpend_ind2_sols_vg:e2})). Fig.~\ref{fig:dpend_ind2_errors} shows the dependencies of the global errors $e^{u}$ and $e^{g}$ of the numerical solution at the nodes $(\mathbf{u}_{n}, \mathbf{v}_{n})$ on the discretization step ${\Delta t}$, separately for differential variables $\mathbf{u}$ and the algebraic equations $\mathbf{g} = \mathbf{0}$, satisfying separately on the numerical solution at the nodes $(\mathbf{u}_{n}, \mathbf{v}_{n})$ and the local solution $(\mathbf{u}_{L}, \mathbf{v}_{L})$, on the basis of which the empirical convergence orders $p$ were calculated. For the double pendulum test case presented in this Subsection, there is no exact analytical solution, so the numerical solution was compared with the reference numerical solution (\ref{eq:math_dpend_ode}), and the convergence orders $p_{f}^{n, u}$ for the solution at nodes $\mathbf{u}_{n}$ for the differential variables $\mathbf{u}$ in the ``final'' norm (\ref{eq:norms_def}), as well the convergence orders $p^{n, g}$ and $p^{l, g}$ for the errors of satisfying the algebraic equations $\mathbf{g} = \mathbf{0}$ for the solution at the nodes $(\mathbf{u}_{n}, \mathbf{v}_{n})$ and the local solution $(\mathbf{u}_{L}, \mathbf{v}_{L})$ separately for the norms $L_{1}$, $L_{2}$, $L_{\infty}$ (\ref{eq:norms_def}).

The dynamic dependencies of the numerical solution at the nodes $(\mathbf{u}_{n}, \mathbf{v}_{n})$ and the local solution $(\mathbf{u}_{L}, \mathbf{v}_{L})$ shown in Fig.~\ref{fig:dpend_ind2_sols_u} for differential variables $\mathbf{u}$ and in Fig.~\ref{fig:dpend_ind2_sols_vg} (\subref{fig:dpend_ind2_sols_vg:a1}, \subref{fig:dpend_ind2_sols_vg:c1}, \subref{fig:dpend_ind2_sols_vg:e1}) and Fig.~\ref{fig:dpend_ind2_sols_vg} (\subref{fig:dpend_ind2_sols_vg:a2}, \subref{fig:dpend_ind2_sols_vg:c2}, \subref{fig:dpend_ind2_sols_vg:e2}) for algebraic variables $\mathbf{v}$ demonstrate the high accuracy of the ADER-DG numerical method with a local DG predictor for solving the DAE system (\ref{eq:math_dpend_dae_ind_3}) of index 2 describing the motion of a double pendulum. This is clearly seen from a comparison of the dynamic dependencies of the solution at the nodes $(\mathbf{u}_{n}, \mathbf{v}_{n})$, the local solution $(\mathbf{u}_{L}, \mathbf{v}_{L})$, and the reference solution $(\mathbf{u}^{\rm ref}, \mathbf{v}^{\rm ref})$. The obtained results of comparison of the numerical solution at the nodes and the local solution with the reference solution, in general, correspond to the results from Example 1, presented in Subsection~\ref{sec:2:ct:ex1}, Example 4 for the DAE system (\ref{eq:math_dpend_dae_ind_3}) of index 2, presented in Subsection~\ref{sec:2:ct:ex4}, and Example 2, presented in Subsection~\ref{sec:2:ct:ex2}, where test examples of DAE of index 1 systems were solved, but this result differs significantly from the DAE system (\ref{eq:hess_dae_ind_2}) of index 2. In the case of polynomials of degree $N = 1$, the numerical solution for which is shown in Fig.~\ref{fig:dpend_ind2_sols_u} (\subref{fig:dpend_ind2_sols_u:a1}, \subref{fig:dpend_ind2_sols_u:a2}, \subref{fig:dpend_ind2_sols_u:a3}, \subref{fig:dpend_ind2_sols_u:a4}, \subref{fig:dpend_ind2_sols_u:b1}, \subref{fig:dpend_ind2_sols_u:b2}, \subref{fig:dpend_ind2_sols_u:b3}, \subref{fig:dpend_ind2_sols_u:b4}) for differential variables $\mathbf{u}$ and in Fig.~\ref{fig:dpend_ind2_sols_vg} (\subref{fig:dpend_ind2_sols_vg:a1}, \subref{fig:dpend_ind2_sols_vg:a2}) for algebraic variables $\mathbf{v}$, artifacts expected for a local solution observed for DAE system (\ref{eq:math_dpend_dae_ind_3}) of index 3 do not occur. In the presented dynamic dependencies of the numerical solution at the nodes $(\mathbf{u}_{n}, \mathbf{v}_{n})$ and the local solution $(\mathbf{u}_{L}, \mathbf{v}_{L})$ for the cases of polynomial degrees of $N = 8$ and $N = 40$, even on a very coarse grid (only $10$ discretization domains), the numerical solution demonstrates excellent point-wise agreement with the reference solution of the problem $(\mathbf{u}^{\rm ref}, \mathbf{v}^{\rm ref})$. To determine the accuracy of fulfillment of algebraic equations $g_{11} = 0$, $g_{12} = 0$, $g_{21} = 0$, $g_{22} = 0$, $g_{31} = 0$, $g_{32} = 0$ the dynamic dependencies of errors $|g_{11}|$, $|g_{12}|$, $|g_{21}|$, $|g_{22}|$, $|g_{31}|$, $|g_{32}|$ were calculated, which were presented in Fig.~\ref{fig:dpend_ind2_sols_vg} (\subref{fig:dpend_ind2_sols_vg:a3}, \subref{fig:dpend_ind2_sols_vg:a4}, \subref{fig:dpend_ind2_sols_vg:b1}, \subref{fig:dpend_ind2_sols_vg:b2}, \subref{fig:dpend_ind2_sols_vg:b3}, \subref{fig:dpend_ind2_sols_vg:b4}) for polynomial degree $N = 1$, in Fig.~\ref{fig:dpend_ind2_sols_vg} (\subref{fig:dpend_ind2_sols_vg:c3}, \subref{fig:dpend_ind2_sols_vg:c4}, \subref{fig:dpend_ind2_sols_vg:d1}, \subref{fig:dpend_ind2_sols_vg:d2}, \subref{fig:dpend_ind2_sols_vg:d3}, \subref{fig:dpend_ind2_sols_vg:d4}) for polynomial degree $N = 8$ and in Fig.~\ref{fig:dpend_ind2_sols_vg} (\subref{fig:dpend_ind2_sols_vg:e3}, \subref{fig:dpend_ind2_sols_vg:e4}, \subref{fig:dpend_ind2_sols_vg:f1}, \subref{fig:dpend_ind2_sols_vg:f2}, \subref{fig:dpend_ind2_sols_vg:f3}, \subref{fig:dpend_ind2_sols_vg:f4}) for polynomial degree $N = 40$. The presented dependencies $|g_{21}|$ and $|g_{22}|$ demonstrate that algebraic equations $g_{21} = 0$ and $g_{22} = 0$ are satisfied exactly for the solution at the nodes (of course, within the accuracy of representation of real numbers with a floating point), therefore, they are not presented in Fig.~\ref{fig:dpend_ind2_sols_vg}. Errors $|g_{21}|$ and $|g_{22}|$ for the local solution $(\mathbf{u}_{L}, \mathbf{v}_{L})$ are non-zero and are in the characteristic range of $10^{-2}$--$10^{-6}$ in the case of polynomial degree $N = 1$, $10^{-3}$--$10^{-9}$ in the case of polynomial degree $N = 8$ and $10^{-13}$--$10^{-23}$ in the case of polynomial degree $N = 40$. Errors $|g_{31}|$, $|g_{32}|$ in satisfying algebraic equations $g_{31} = 0$, $g_{32} = 0$, which are associated with the algebraic equations $g_{21} = 0$, $g_{22} = 0$ single differentiation, respectively, have characteristic values of $10^{-1}$--$10^{-5}$ in the case of polynomial degree $N = 1$, $10^{-2}$--$10^{-7}$ --- in the case of polynomial degree $N = 8$, $10^{-11}$--$10^{-22}$ --- in the case of polynomial degree $N = 40$. Errors $|g_{11}|$ and $|g_{12}|$ in satisfying algebraic equations $g_{11} = 0$, $g_{12} = 0$ and $g_{31} = 0$, $g_{32} = 0$, which are associated with the algebraic equations $g_{21} = 0$ and $g_{22} = 0$ single integration, for the solution at the nodes $(\mathbf{u}_{n}, \mathbf{v}_{n})$ is significantly smaller than the errors $|g_{11}|$ and $|g_{12}|$ for the local solution $(\mathbf{u}_{L}, \mathbf{v}_{L})$: by $10^{0}$--$10^{2}$ times in the case of polynomial degree $N = 1$, by $10^{2}$--$10^{4}$ times in the case of polynomial degree $N = 8$, by $10^{10}$--$10^{22}$ times in the case of polynomial degree $N = 40$ --- this result differs significantly from the DAE system (\ref{eq:hess_dae_ind_2}) of index 2, but corresponds to the results of Example 4, presented in Subsections~\ref{sec:2:ct:ex4}, where the DAE systems of index 2 were studied. However, errors $|g_{31}|$ and $|g_{32}|$ for the solution at the nodes $(\mathbf{u}_{n}, \mathbf{v}_{n})$ are larger than the corresponding errors for the local solution $(\mathbf{u}_{L}, \mathbf{v}_{L})$ at almost all points $t$. This feature of the numerical solution was already observed in Example 4 for the DAE system (\ref{eq:math_dpend_dae_ind_3}) of index 2, presented in Subsection~\ref{sec:2:ct:ex4}. In Fig.~\ref{fig:dpend_ind2_errors} the dependencies of global errors $e_{f}^{u}$, $e_{L}^{g}$, $e_{G}^{u}$ on the discretization step ${\Delta t}$ are presented. The presented dependencies demonstrate the power law $e(\Delta t) \sim (\Delta t)^{p}$, from which the empirical convergence orders $p$ were calculated. However, in comparison with the previously studied test examples, in this case the approximations are not so good. The orders of convergence calculated on the basis of the presented dependencies were presented in Table~\ref{tab:conv_orders_dpend_ind3}, where the expected values $p_{\rm nodes} = 2N+1$, $p_{\rm local} = N+1$ (\ref{eq:expect_orders}) were also presented for comparison. The obtained values of convergence of orders demonstrate significantly lower values compared to the expected values (\ref{eq:expect_orders}), which may be a consequence of the high index of the studied DAE system --- usually a decrease in the empirical convergence order is observed compared to the expected convergence order of the numerical method, which occurs when solving the initial value problem for the ODE system~\cite{Hairer_book_2}. However, the obtained convergence orders are significantly greater than the values obtained for the DAE system of index 3, considered in this Subsection above.

The obtained results of the numerical solution of DAE system (\ref{eq:math_dpend_dae_ind_3}) of index 1 are presented in Figs.~\ref{fig:dpend_ind1_sols_u}, \ref{fig:dpend_ind1_sols_vg}, \ref{fig:dpend_ind1_errors} and in Table~\ref{tab:conv_orders_dpend_ind1}. Fig.~\ref{fig:pend_ind1_sol_uv} shows a comparison of the numerical solution at the nodes for differential variables $\mathbf{u}_{n}$, the numerical local solution for differential variables $\mathbf{u}_{L}$ and the reference solution for differential variables $\mathbf{u}^{\rm ref}$. Fig.~\ref{fig:dpend_ind1_sols_vg} shows a comparison of the numerical solution at the nodes for algebraic variables $\mathbf{v}_{n}$, the numerical local solution for algebraic variables $\mathbf{v}_{L}$ and the reference solution for algebraic variables $\mathbf{v}^{\rm ref}$, the dependencies of the feasibility of conditions $g_{11} = 0$, $g_{12} = 0$, $g_{21} = 0$, $g_{22} = 0$, $g_{31} = 0$, $g_{32} = 0$ on the coordinate $t$, which allows us to quantitatively estimate the accuracy of the numerical solution, especially taking into account that the numerical solution obtained by the ADER-DG method with a local DG predictor with a high degree of polynomials $N$ very accurately corresponds to the exact analytical solution, and it is impossible to visually identify the error from the solution plots separately (in Fig.~\ref{fig:dpend_ind1_sols_u}, Fig.~\ref{fig:dpend_ind1_sols_vg} (\subref{fig:dpend_ind1_sols_vg:a1}, \subref{fig:dpend_ind1_sols_vg:c1}, \subref{fig:dpend_ind1_sols_vg:e1}) and Fig.~\ref{fig:dpend_ind1_sols_vg} (\subref{fig:dpend_ind1_sols_vg:a2}, \subref{fig:dpend_ind1_sols_vg:c2}, \subref{fig:dpend_ind1_sols_vg:e2})). Fig.~\ref{fig:dpend_ind1_errors} shows the dependencies of the global errors $e^{u}$ and $e^{g}$ of the numerical solution at the nodes $(\mathbf{u}_{n}, \mathbf{v}_{n})$ on the discretization step ${\Delta t}$, separately for differential variables $\mathbf{u}$ and the algebraic equations $\mathbf{g} = \mathbf{0}$, satisfying separately on the numerical solution at the nodes $(\mathbf{u}_{n}, \mathbf{v}_{n})$ and the local solution $(\mathbf{u}_{L}, \mathbf{v}_{L})$, on the basis of which the empirical convergence orders $p$ were calculated. For the double pendulum test case presented in this Subsection, there is no exact analytical solution, so the numerical solution was compared with the reference numerical solution (\ref{eq:math_dpend_ode}), and the convergence orders $p_{f}^{n, u}$ were calculated for the solution at nodes $\mathbf{u}_{n}$ for the differential variables $\mathbf{u}$ in the ``final'' norm (\ref{eq:norms_def}), as well as the convergence orders $p^{n, g}$ and $p^{l, g}$ for the errors of satisfying the algebraic equations $\mathbf{g} = \mathbf{0}$ for the solution at the nodes $(\mathbf{u}_{n}, \mathbf{v}_{n})$ and the local solution $(\mathbf{u}_{L}, \mathbf{v}_{L})$ separately for the norms $L_{1}$, $L_{2}$, $L_{\infty}$ (\ref{eq:norms_def}) were calculated.

The dynamic dependencies of the numerical solution at the nodes $(\mathbf{u}_{n}, \mathbf{v}_{n})$ and the local solution $(\mathbf{u}_{L}, \mathbf{v}_{L})$ shown in Fig.~\ref{fig:dpend_ind1_sols_u} for differential variables $\mathbf{u}$ and in Fig.~\ref{fig:dpend_ind1_sols_vg} (\subref{fig:dpend_ind1_sols_vg:a1}, \subref{fig:dpend_ind1_sols_vg:c1}, \subref{fig:dpend_ind1_sols_vg:e1}) and Fig.~\ref{fig:dpend_ind1_sols_vg} (\subref{fig:dpend_ind1_sols_vg:a2}, \subref{fig:dpend_ind1_sols_vg:c2}, \subref{fig:dpend_ind1_sols_vg:e2}) for algebraic variables $\mathbf{v}$ demonstrate the high accuracy of the ADER-DG numerical method with a local DG predictor for solving the DAE system (\ref{eq:math_dpend_dae_ind_3}) of index 1 describing the motion of a double pendulum. This is clearly seen from a comparison of the dynamic dependencies of the solution at the nodes $(\mathbf{u}_{n}, \mathbf{v}_{n})$, the local solution $(\mathbf{u}_{L}, \mathbf{v}_{L})$, and the reference solution $(\mathbf{u}^{\rm ref}, \mathbf{v}^{\rm ref})$. The obtained results of comparison of the numerical solution at the nodes and the local solution with the reference solution, in general, correspond to the results from Example 1, presented in Subsection~\ref{sec:2:ct:ex1}, Example 4 for the DAE system (\ref{eq:math_dpend_dae_ind_3}) of index 1, presented in Subsection~\ref{sec:2:ct:ex4}, and Example 2, presented in Subsection~\ref{sec:2:ct:ex2}, where test examples of DAE of index 1 systems were solved. In the presented dynamic dependencies of the numerical solution at the nodes $(\mathbf{u}_{n}, \mathbf{v}_{n})$ and the local solution $(\mathbf{u}_{L}, \mathbf{v}_{L})$ for the cases of polynomial degrees of $N = 1$, $N = 8$ and $N = 40$, even on a very coarse grid (only $10$ discretization domains), the numerical solution demonstrates excellent point-wise agreement with the reference solution of the problem $(\mathbf{u}^{\rm ref}, \mathbf{v}^{\rm ref})$. To determine the accuracy of fulfillment of algebraic equations $g_{11} = 0$, $g_{12} = 0$, $g_{21} = 0$, $g_{22} = 0$, $g_{31} = 0$, $g_{32} = 0$ the dynamic dependencies of errors $|g_{11}|$, $|g_{12}|$, $|g_{21}|$, $|g_{22}|$, $|g_{31}|$, $|g_{32}|$ were calculated, which were presented in Fig.~\ref{fig:dpend_ind1_sols_vg} (\subref{fig:dpend_ind1_sols_vg:a3}, \subref{fig:dpend_ind1_sols_vg:a4}, \subref{fig:dpend_ind1_sols_vg:b1}, \subref{fig:dpend_ind1_sols_vg:b2}, \subref{fig:dpend_ind1_sols_vg:b3}, \subref{fig:dpend_ind1_sols_vg:b4}) for polynomial degree $N = 1$, in Fig.~\ref{fig:dpend_ind1_sols_vg} (\subref{fig:dpend_ind1_sols_vg:c3}, \subref{fig:dpend_ind1_sols_vg:c4}, \subref{fig:dpend_ind1_sols_vg:d1}, \subref{fig:dpend_ind1_sols_vg:d2}, \subref{fig:dpend_ind1_sols_vg:d3}, \subref{fig:dpend_ind1_sols_vg:d4}) for polynomial degree $N = 8$ and in Fig.~\ref{fig:dpend_ind1_sols_vg} (\subref{fig:dpend_ind1_sols_vg:e3}, \subref{fig:dpend_ind1_sols_vg:e4}, \subref{fig:dpend_ind1_sols_vg:f1}, \subref{fig:dpend_ind1_sols_vg:f2}, \subref{fig:dpend_ind1_sols_vg:f3}, \subref{fig:dpend_ind1_sols_vg:f4}) for polynomial degree $N = 40$. The presented dependencies $|g_{31}|$ and $|g_{32}|$ demonstrate that algebraic equations $g_{31} = 0$ and $g_{32} = 0$ are satisfied exactly for the solution at the nodes (of course, within the accuracy of representation of real numbers with a floating point), therefore, they are not presented in Fig.~\ref{fig:dpend_ind1_sols_vg}. Errors $|g_{31}|$ and $|g_{32}|$ for the local solution $(\mathbf{u}_{L}, \mathbf{v}_{L})$ are non-zero and are in the characteristic range of $10^{-2}$--$10^{-5}$ in the case of polynomial degree $N = 1$, $10^{-3}$--$10^{-6}$ in the case of polynomial degree $N = 8$ and $10^{-12}$--$10^{-22}$ in the case of polynomial degree $N = 40$. Errors $|g_{11}|$, $|g_{12}|$, $|g_{21}|$, $|g_{22}|$ for the solution at the nodes $(\mathbf{u}_{n}, \mathbf{v}_{n})$ is significantly smaller than the errors $|g_{11}|$ and $|g_{12}|$ for the local solution $(\mathbf{u}_{L}, \mathbf{v}_{L})$: by $10^{0}$--$10^{2}$ times in the case of polynomial degree $N = 1$, by $10^{2}$--$10^{3}$ times in the case of polynomial degree $N = 8$, by $10^{8}$--$10^{20}$ times in the case of polynomial degree $N = 40$. In Fig.~\ref{fig:dpend_ind1_errors} the dependencies of global errors $e_{f}^{u}$, $e_{L}^{g}$, $e_{G}^{u}$ on the discretization step ${\Delta t}$ are presented. The presented dependencies demonstrate the power law $e(\Delta t) \sim (\Delta t)^{p}$, from which the empirical convergence orders $p$ were calculated. However, in comparison with the previously studied test examples, in this case the approximations are not so good. The orders of convergence calculated on the basis of the presented dependencies were presented in Table~\ref{tab:conv_orders_dpend_ind3}, where the expected values $p_{\rm nodes} = 2N+1$, $p_{\rm local} = N+1$ (\ref{eq:expect_orders}) were also presented for comparison. The obtained values of convergence of orders demonstrate significantly lower values compared to the expected values (\ref{eq:expect_orders}), which may be a consequence of the high index of the studied DAE system --- usually a decrease in the empirical convergence order is observed compared to the expected convergence order of the numerical method, which occurs when solving the initial value problem for the ODE system~\cite{Hairer_book_2}. The most important difference from the case of the DAE system (\ref{eq:math_dpend_dae_ind_3}) of index 2 considered in this Subsection is the significantly higher convergence orders $p^{n, g}$, while the convergence orders $p_{f}$ takes approximately the same values as in the case of index 2. However, the obtained convergence orders are significantly greater than the values obtained for the DAE system of index 3, considered in this Subsection above.

Therefore, it can be concluded that the presented numerical solution of the DAE system (\ref{eq:math_dpend_dae_ind_3}) demonstrates a significant decrease in the empirical convergence orders $p$ compared to the expected values (\ref{eq:expect_orders}). A decrease in the DAE index of the system leads to a significant increase in convergence order, however, their values are significantly less than the expected values.

\subsection{Stiff tests}
\label{sec:2:st}

The examples of application of the numerical method ADER-DG with a local DG predictor presented in the Section~\ref{sec:2:ct} were related to the numerical solution of non-stiff problems for DAE systems: the simple DAE system of index 1 (\ref{eq:simple_dae_ind_1}) is presented in Subsection~\ref{sec:2:ct:ex1}, Hessenberg DAE systems of index 1 (\ref{eq:hess_dae_ind_1}) and index 2 (\ref{eq:hess_dae_ind_2}) are presented in Subsections~\ref{sec:2:ct:ex2} and~\ref{sec:2:ct:ex3}, a solution of the mathematical pendulum problem representing a DAE system of index 3 (\ref{eq:math_pend_dae_ind_3}) is presented in Subsection~\ref{sec:2:ct:ex4}, a solution of the double pendulum problem also representing a DAE system of index 3 (\ref{eq:math_dpend_dae_ind_3}) is presented in Subsection~\ref{sec:2:ct:ex5}. However, it is well known that the numerical method ADER-DG with a local DG predictor is well suited for solving stiff problems~\cite{ader_stiff_1, ader_stiff_2, ader_stiff_3, ader_dg_ivp_ode}, including those characterized by extreme stiffness~\cite{ader_dg_ivp_ode}. In the works~\cite{ader_stiff_1, ader_stiff_2}, numerical methods based on the ADER paradigm were used to solve stiff problems associated in particular with combustion and detonation. In the works~\cite{ader_stiff_3, ader_stiff_4}, the finite element numerical method ADER-DG with a posteriori correction of the solution by the finite-volume ADER-WENO method was also successfully used to simulate the development and propagation of a one-dimensional and multidimensional detonation waves.

This Section presents a numerical solution of a well-stiff problem known from the field of solving the initial value problem for ODE systems: $\dot{u} = u^{2} - u^{3}$, $u(0) = \delta$, $t\in[0,\, 2/\delta]$, where $\delta$ is a small constant. This equation is a very simple model of the fireball balances in a flammable medium~\cite{stiff_flame_1}: the quadratic term $u^{2}$ is related to the area of the combustion domain boundary to which oxygen is supplied, and the cubic term $u^{3}$ determines a value proportional to the volume of the fireball. This problem is well known among stiff ODE systems~\cite{stiff_flame_1, stiff_flame_2, stiff_flame_3, stiff_flame_4, stiff_flame_5}. This initial value problem for ODE systems was reformulated as a DAE system of the following form:
\begin{equation}\label{eq:shampine_flame}
\begin{split}
&\dot{u}_{1} = u_{1}^{2} - v_{1},\hspace{15mm} u_{1}(0) = \delta,\\
&g_{1} = u_{1}^{3} - v_{1} = 0,\hspace{8.85mm} v_{1}(0) = \delta^{3},\qquad t\in[0,\, 2/\delta],
\end{split}
\end{equation}
which corresponds to the formulation of the original DAE system (\ref{eq:dae_chosen_form}). The presented DAE system of index 1 contains one algebraic equation (constraint) $g_{1} = 0$. The full-component exact analytical solution of the problem was obtained by trivial integration of the ODE problem $\dot{u} = u^{2} - u^{3}$ with initial condition $u(0) = \delta$ and was written in the following form:
\begin{equation}\label{eq:shampine_flame:exact_solution}
\begin{split}
&u(t) = \frac{1}{W(a\cdot\exp(a - t)) + 1},\\
&v(t) = \frac{1}{\left[W(a\cdot\exp(a - t)) + 1\right]^{3}},\quad a = \frac{1}{\delta} - 1,
\end{split}
\end{equation}
where $W=W(x)$ is the Lambert $W$ function that determines the solution to the transcendental algebraic equation $W\exp(W)=x$. The investigations of the solution to the problem (\ref{eq:shampine_flame}) presented in the works~\cite{stiff_flame_1, stiff_flame_3, stiff_flame_4} and the work~\cite{stiff_flame_2} shows that in the case of small values of parameter $\delta \ll 1$, the equation demonstrate high stiffness in the domain near point $1/\delta$, while case $\delta = 10^{-4}$ can be called stiff, and it is usually used to demonstrate~\cite{stiff_flame_1, stiff_flame_2, stiff_flame_3, stiff_flame_4} the capabilities of the numerical method for solving stiff ODE systems, and cases $\delta = 10^{-5}$ and $10^{-6}$ is characterized by extreme stiffness. It is known that the solution $u(t)$ has asymptotic behavior $u(t \gg 1/\delta) \simeq 1$, and in the domain near point $1/\delta$ there is a very sharp change in the function from a small value $u \simeq \delta$ to an asymptotic one $u(t \rightarrow \infty) \rightarrow 1$. All dynamic dependencies of the solution were constructed and analysed in units $\delta^{-1}$ that have a dimension $[t]^{-1}$, i.e. depending on the argument $\delta \cdot t$, using which the domain of definition $\Omega = [0, 2]$ did not depend on the constant $\delta$.

A numerical solution to the DAE system problem (\ref{eq:shampine_flame}) was obtained using the numerical method ADER-DG with a local DG predictor using degrees of polynomials $N = 1,\, 2,\, \ldots,\, 40$, which is similar to the numerical solutions of non-stiff problems presented above, and for three parameter $\delta$ values: $\delta = 10^{-4}$, $\delta = 10^{-5}$, $\delta = 10^{-6}$. Due to the fact that the areas of occurrence of significant stiffness of the solution is localized in the vicinity of point $1/\delta$, a variable grid step $\Delta t$ was chosen. In case $\delta = 10^{-4}$, the definition domain $\Omega = [0, 2/\delta]$ was divided into three subdomains: $\Omega_{1} = [0, 0.4/\delta]$, $\Omega_{2} = [0.4/\delta, 0.6/\delta]$, $\Omega_{3} = [0.6/\delta, 2/\delta]$; cases of $L_{1} = 10$, $12$, $15$, $20$ grid nodes were selected in the $\Omega_{1}$ and $\Omega_{3}$ subdomains, cases of $L_{2} = 1000$, $1200$, $1500$, $2000$ grid nodes were selected in the $\Omega_{2}$ subdomain. In case $\delta = 10^{-5}$, the definition domain $\Omega = [0, 2/\delta]$ was divided into three subdomains: $\Omega_{1} = [0, 0.495/\delta]$, $\Omega_{2} = [0.495/\delta, 0.505/\delta]$, $\Omega_{3} = [0.505/\delta, 2/\delta]$; cases of $L_{1} = 40$, $48$, $60$, $80$ grid nodes were selected in the $\Omega_{1}$ and $\Omega_{3}$ subdomains, cases of $L_{2} = 1000$, $1200$, $1500$, $2000$ grid nodes were selected in the $\Omega_{2}$ subdomain. In case $\delta = 10^{-6}$, the definition domain $\Omega = [0, 2/\delta]$ was divided into three subdomains: $\Omega_{1} = [0, 0.4975/\delta]$, $\Omega_{2} = [0.4975/\delta, 0.5025/\delta]$, $\Omega_{3} = [0.5025/\delta, 2/\delta]$; cases of $L_{1} = 160$, $192$, $240$, $320$ grid nodes were selected in the $\Omega_{1}$ and $\Omega_{3}$ subdomains, cases of $L_{2} = 1000$, $1200$, $1500$, $2000$ grid nodes were selected in the $\Omega_{2}$ subdomain. The choice of the number of grid nodes in these cases was also associated with a proportional increase in the number of nodes in each of the subdomains $(\Omega_{1},\, \Omega_{2},\, \Omega_{3})$ with a decrease in the discretization step ${\Delta t}$, which was associated with the need to calculate empirical convergence orders $p$. The total number of discretization domains in case $\delta = 10^{-4}$ was $L = 1020$, $1224$, $1530$, $2040$, in case $\delta = 10^{-5}$ --- $1080$, $1296$, $1620$, $2180$, in case $\delta = 10^{-6}$ --- $1320$, $1584$, $1980$, $2640$, while the number of discretization domains differed only in subdomains $\Omega_{1}$ and $\Omega_{2}$, and in the subdomain $\Omega_{3}$ in all cases the number of discretization domains was the same.

The cases $\delta = 10^{-4}$, $\delta = 10^{-5}$ were used in the work~\cite{ader_dg_ivp_ode} to study the possibilities of applying the numerical method ADER-DG with a local DG predictor to solving the initial value problem for stiff ODE systems. It was found that in this case, a significant decrease in the empirical convergence orders is observed. The parameters of the division of the definition domain and the number of discretization domains for each subdomain in cases $\delta = 10^{-4}$ and $\delta = 10^{-5}$ were chosen in accordance with the work~\cite{ader_dg_ivp_ode}. Case $\delta = 10^{-6}$ was not studied in the work~\cite{ader_dg_ivp_ode}. It should be noted that the work~\cite{ader_dg_ivp_ode} used the numerical method ADER-DG with a local DG predictor based on Lagrange interpolation polynomials with nodal points at the roots of shifted Legendre polynomials as basis functions for representing a local discrete time solution (\ref{eq:qr_def_exp}), however, in this work, set of Lagrange interpolation polynomials with nodal points at the roots of the right Radau polynomials were used to represent a local discrete time solution (\ref{eq:qr_def_exp}).

\begin{figure}[h!]
\captionsetup[subfigure]{%
	position=bottom,
	font+=smaller,
	textfont=normalfont,
	singlelinecheck=off,
	justification=raggedright
}
\centering
\begin{subfigure}{0.320\textwidth}
\includegraphics[width=\textwidth]{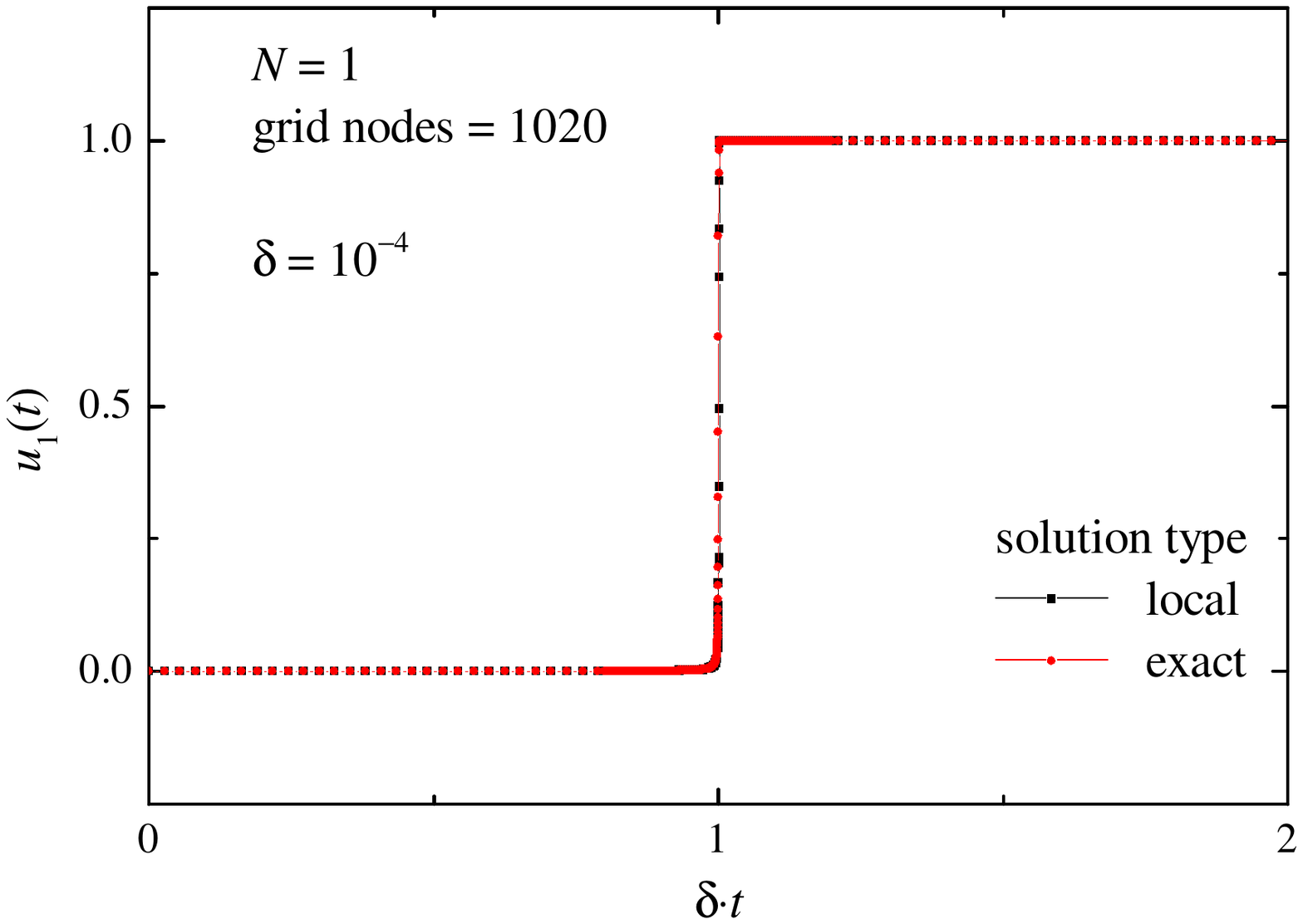}
\vspace{-8mm}\caption{\label{fig:shampine_flame_delta_10m4_sol_qug:a1}}
\end{subfigure}
\begin{subfigure}{0.320\textwidth}
\includegraphics[width=\textwidth]{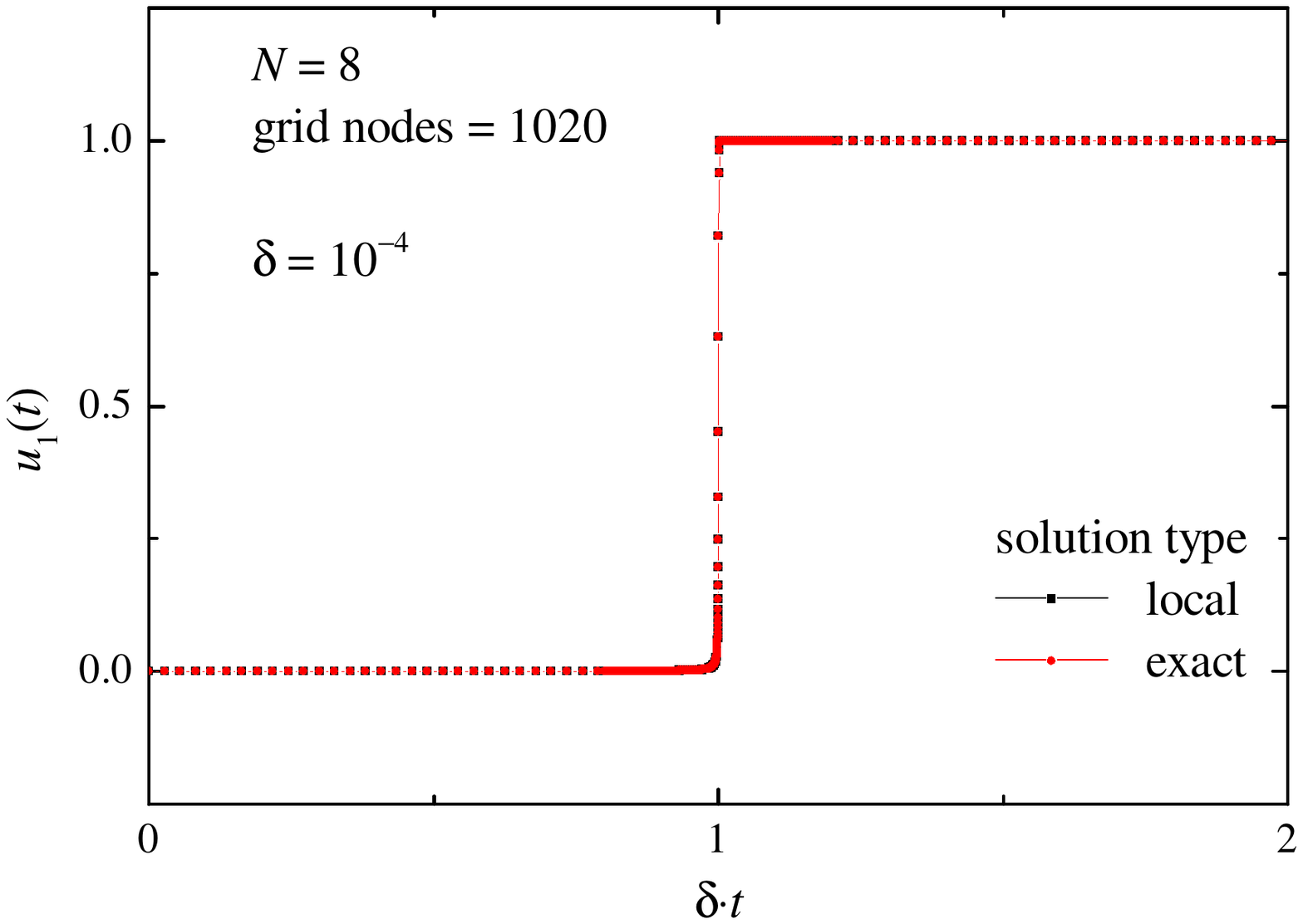}
\vspace{-8mm}\caption{\label{fig:shampine_flame_delta_10m4_sol_qug:a2}}
\end{subfigure}
\begin{subfigure}{0.320\textwidth}
\includegraphics[width=\textwidth]{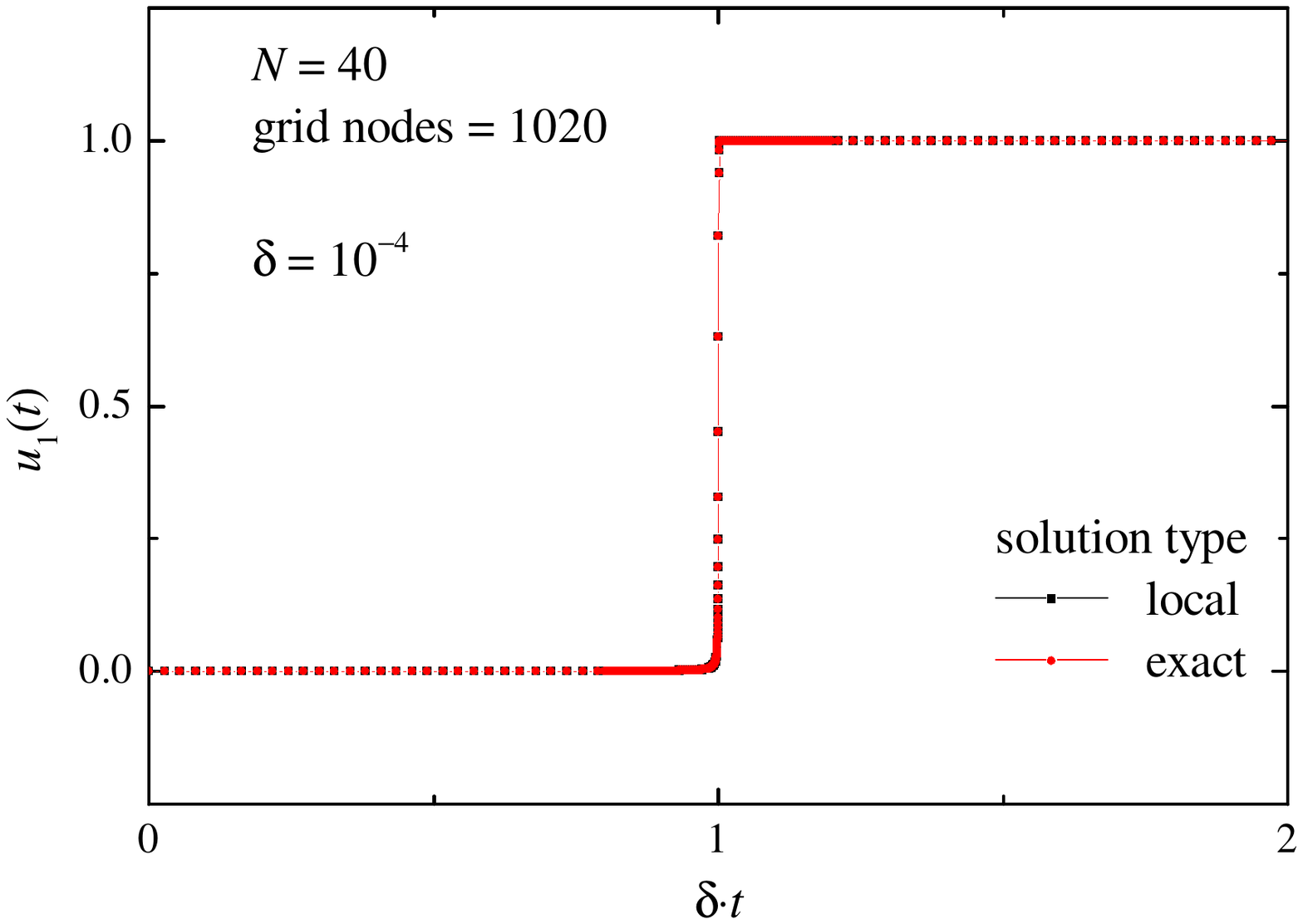}
\vspace{-8mm}\caption{\label{fig:shampine_flame_delta_10m4_sol_qug:a3}}
\end{subfigure}\\
\begin{subfigure}{0.320\textwidth}
\includegraphics[width=\textwidth]{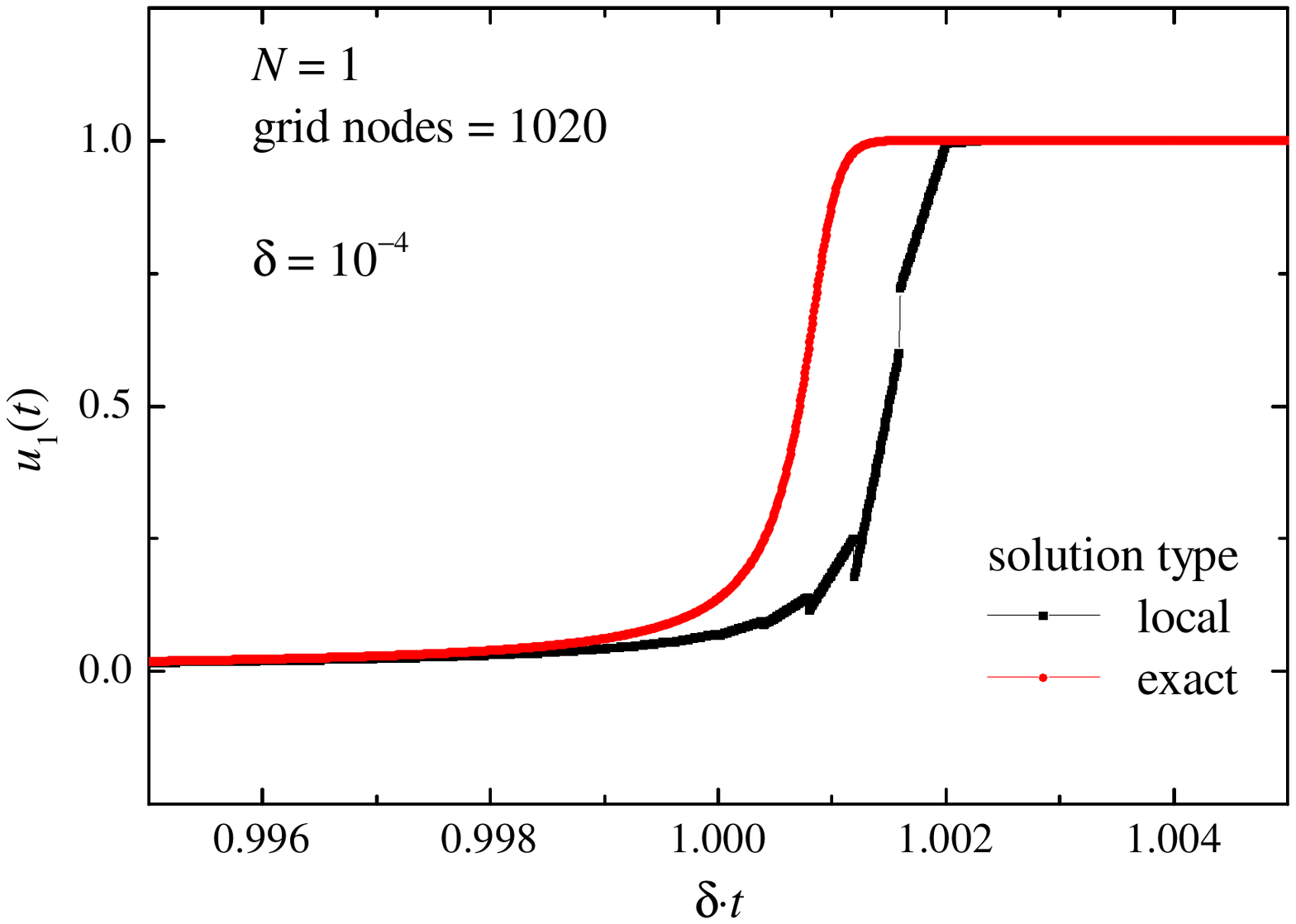}
\vspace{-8mm}\caption{\label{fig:shampine_flame_delta_10m4_sol_qug:b1}}
\end{subfigure}
\begin{subfigure}{0.320\textwidth}
\includegraphics[width=\textwidth]{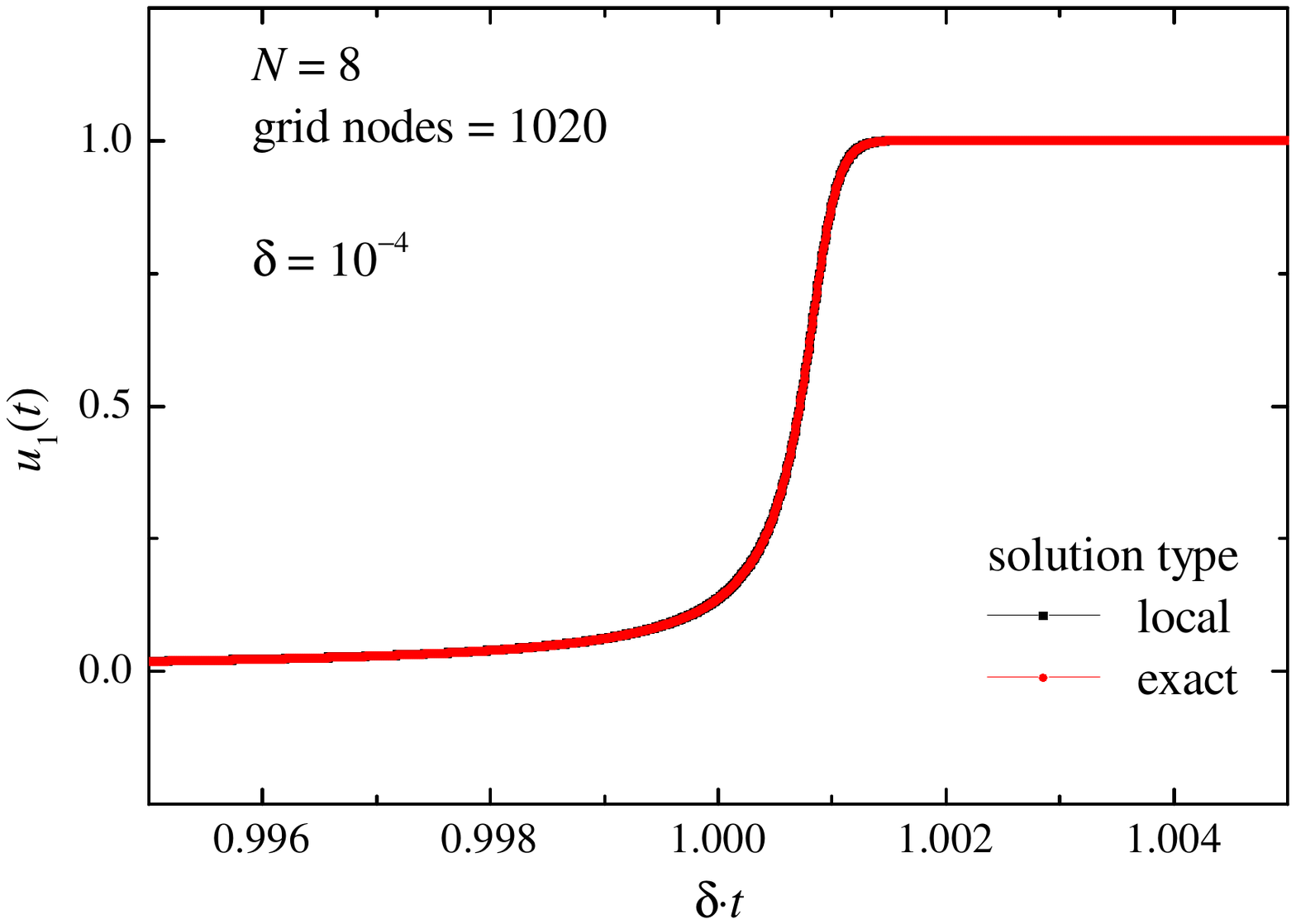}
\vspace{-8mm}\caption{\label{fig:shampine_flame_delta_10m4_sol_qug:b2}}
\end{subfigure}
\begin{subfigure}{0.320\textwidth}
\includegraphics[width=\textwidth]{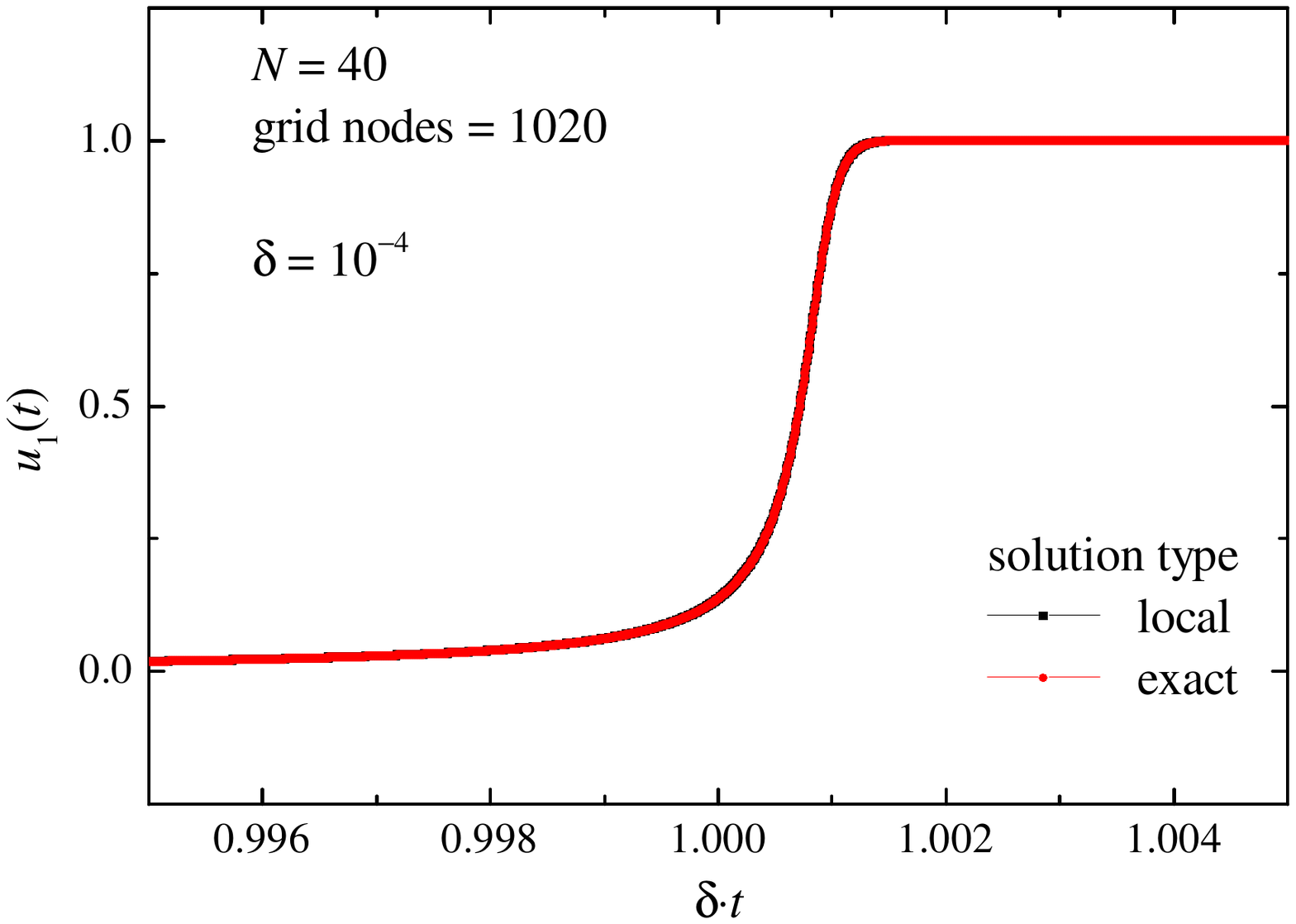}
\vspace{-8mm}\caption{\label{fig:shampine_flame_delta_10m4_sol_qug:b3}}
\end{subfigure}\\
\begin{subfigure}{0.320\textwidth}
\includegraphics[width=\textwidth]{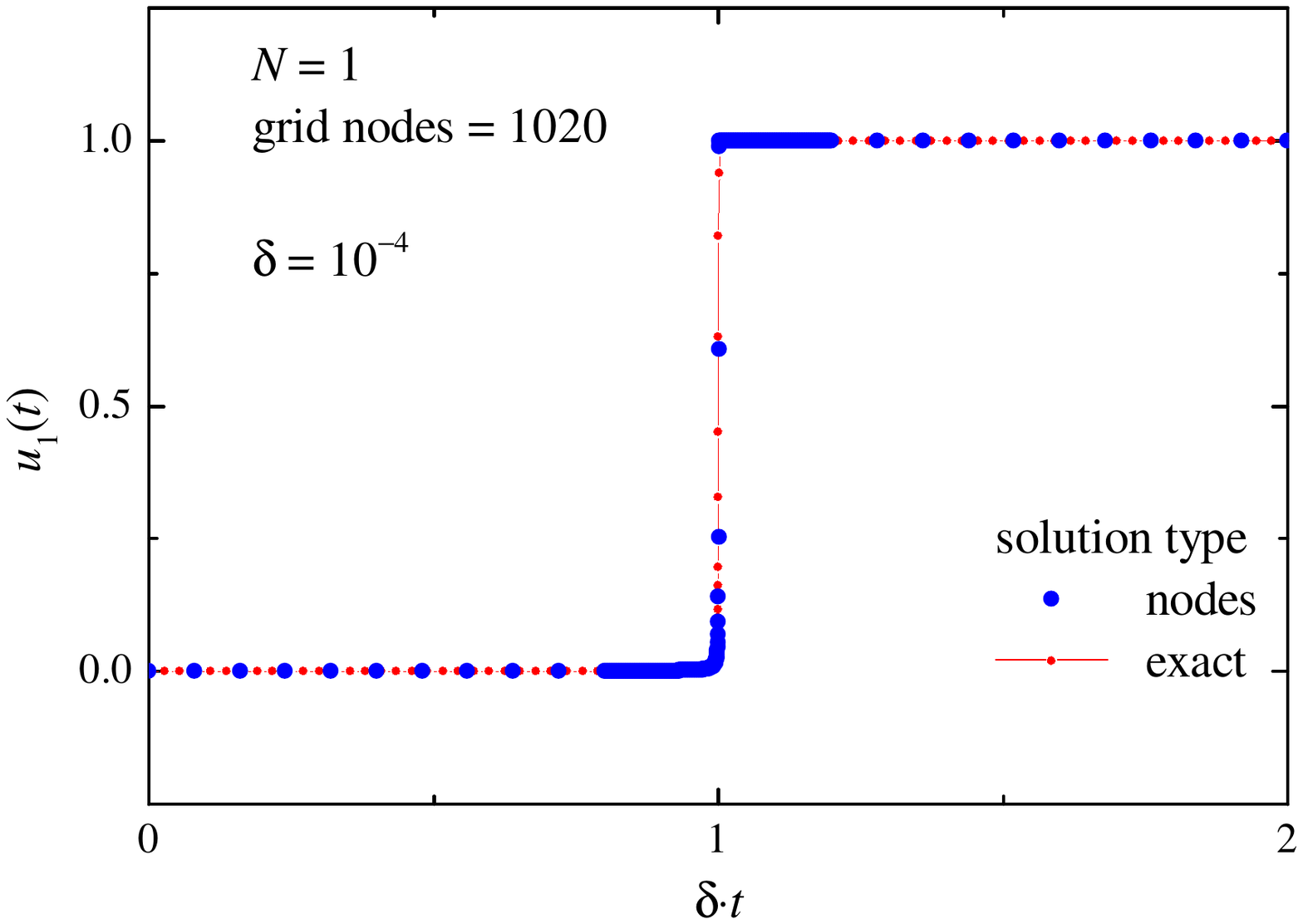}
\vspace{-8mm}\caption{\label{fig:shampine_flame_delta_10m4_sol_qug:c1}}
\end{subfigure}
\begin{subfigure}{0.320\textwidth}
\includegraphics[width=\textwidth]{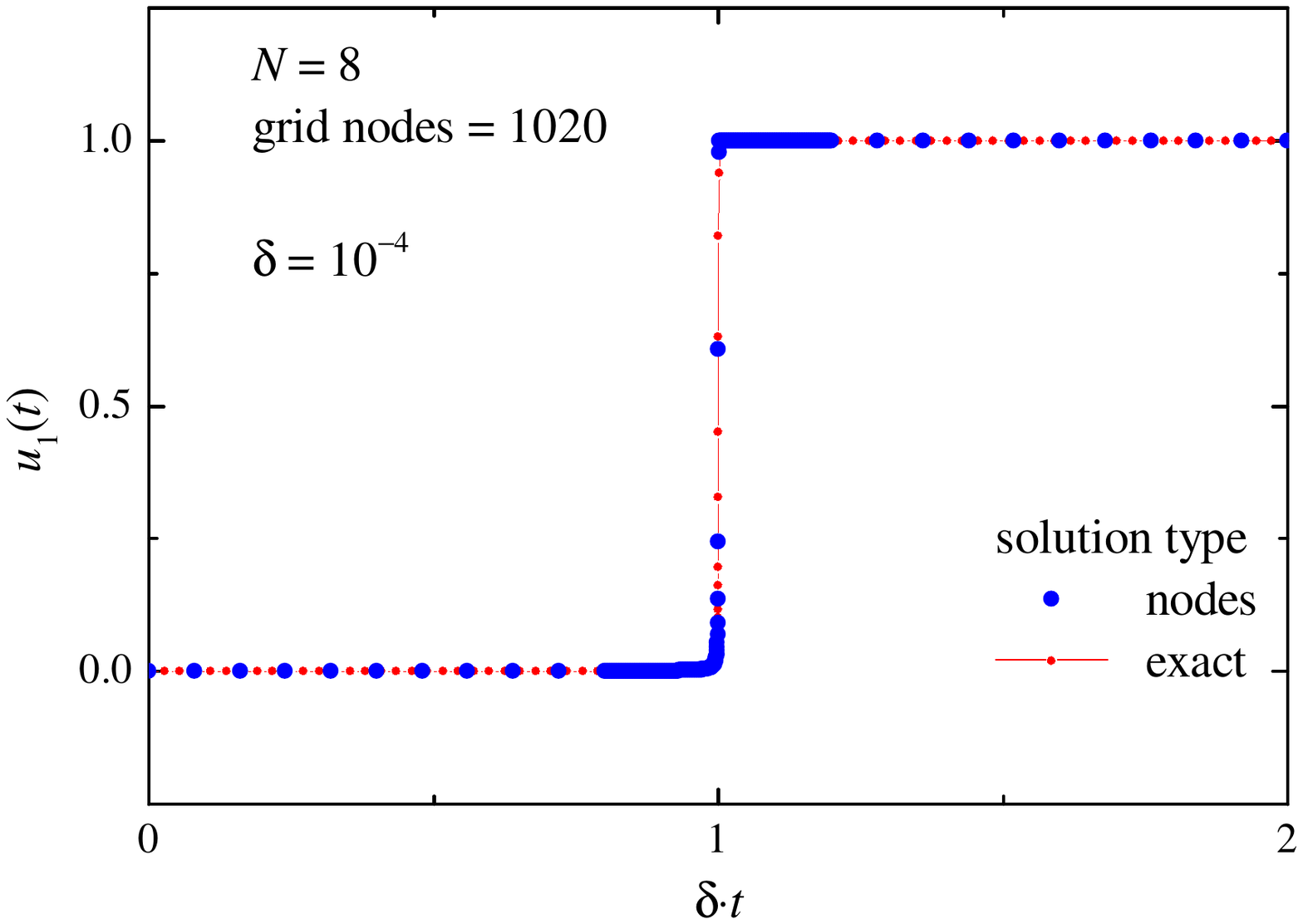}
\vspace{-8mm}\caption{\label{fig:shampine_flame_delta_10m4_sol_qug:c2}}
\end{subfigure}
\begin{subfigure}{0.320\textwidth}
\includegraphics[width=\textwidth]{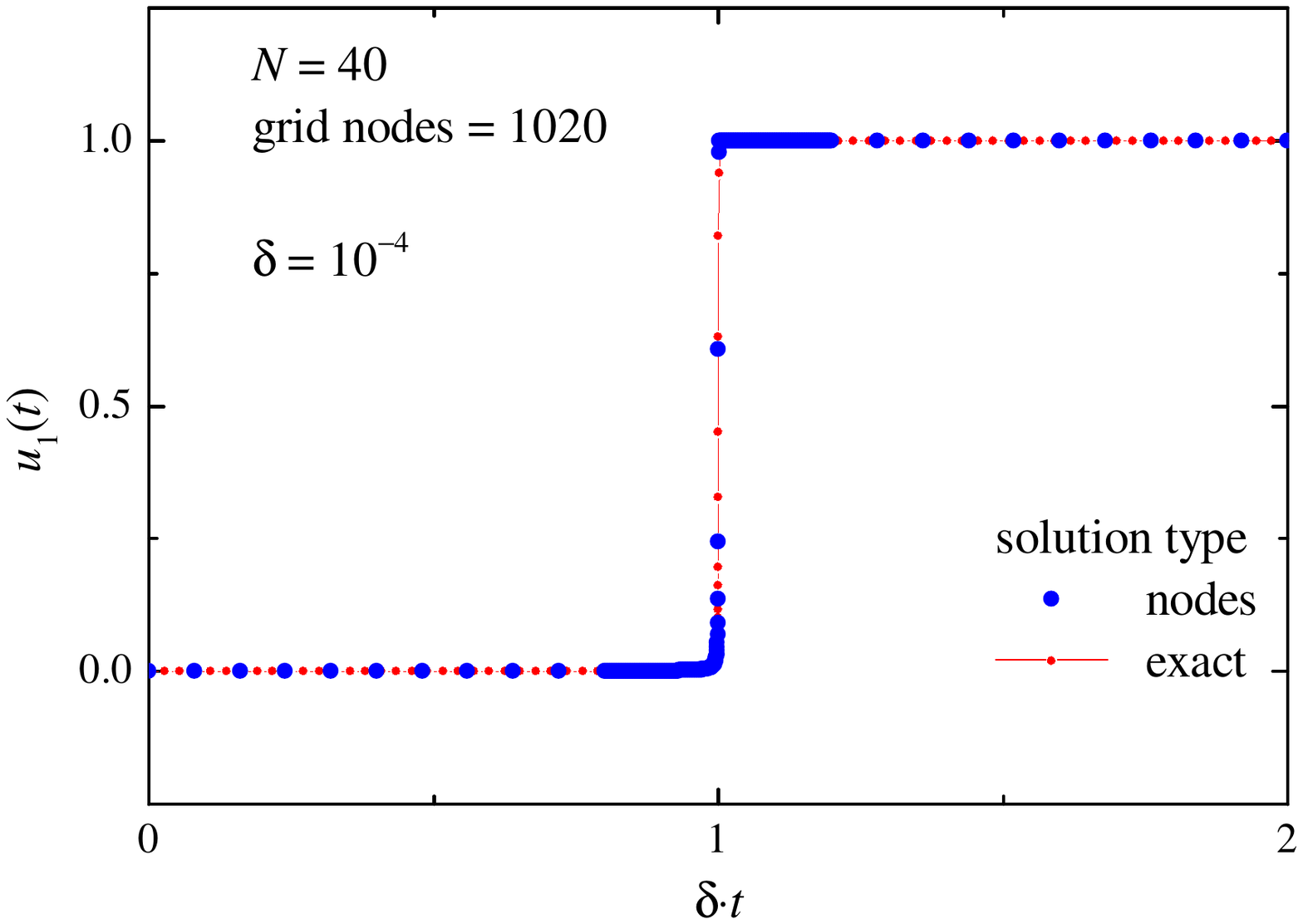}
\vspace{-8mm}\caption{\label{fig:shampine_flame_delta_10m4_sol_qug:c3}}
\end{subfigure}\\
\begin{subfigure}{0.320\textwidth}
\includegraphics[width=\textwidth]{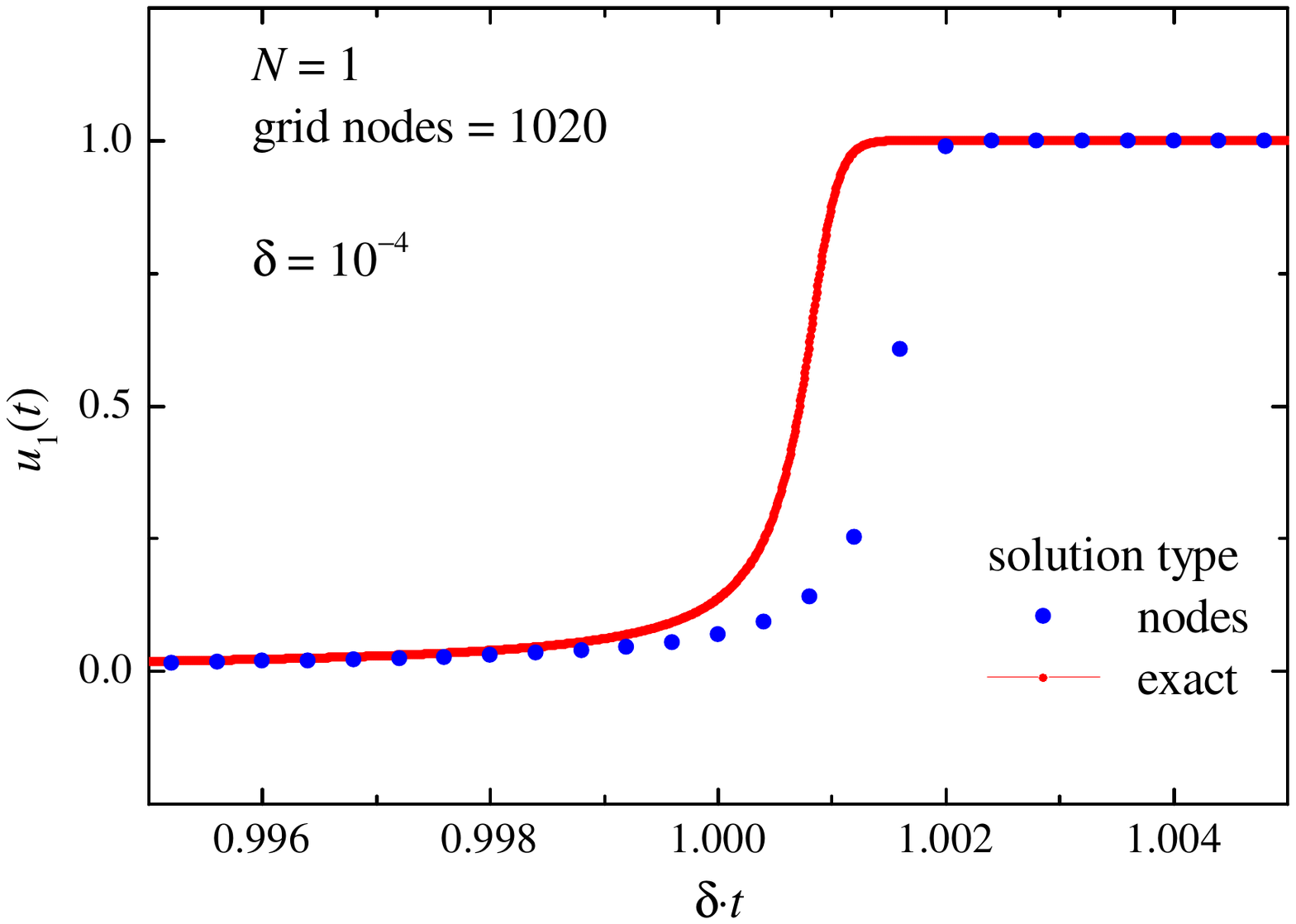}
\vspace{-8mm}\caption{\label{fig:shampine_flame_delta_10m4_sol_qug:d1}}
\end{subfigure}
\begin{subfigure}{0.320\textwidth}
\includegraphics[width=\textwidth]{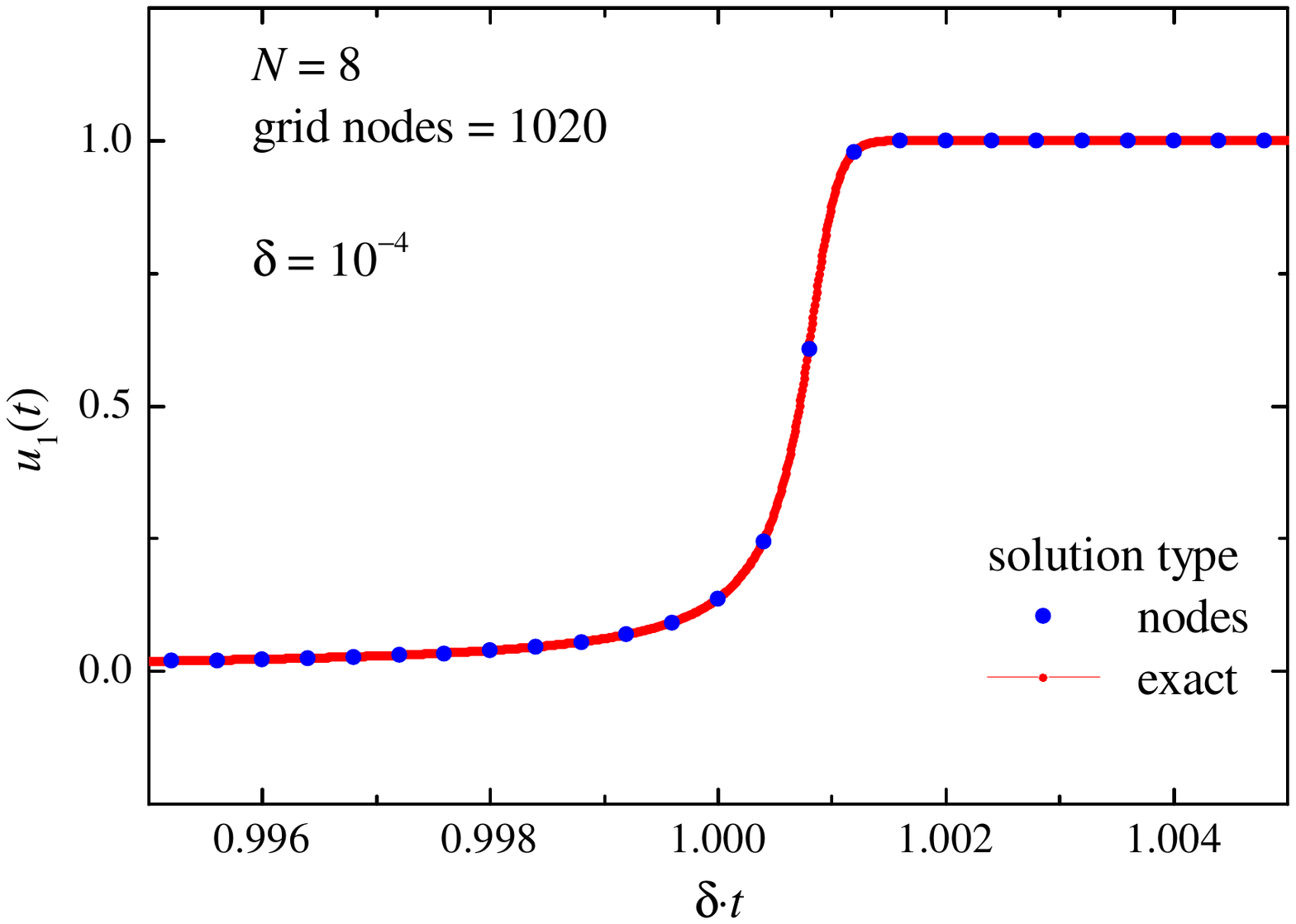}
\vspace{-8mm}\caption{\label{fig:shampine_flame_delta_10m4_sol_qug:d2}}
\end{subfigure}
\begin{subfigure}{0.320\textwidth}
\includegraphics[width=\textwidth]{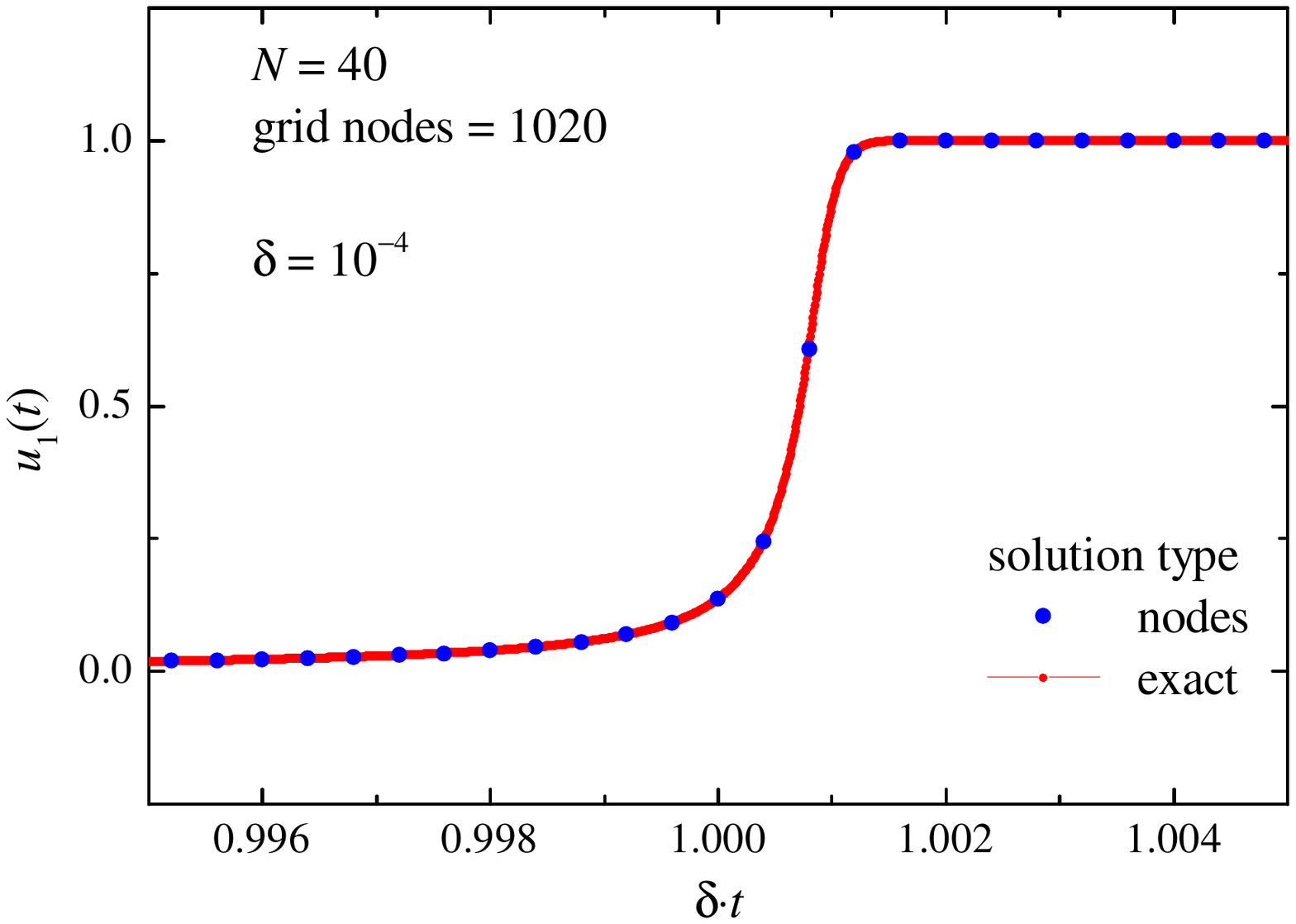}
\vspace{-8mm}\caption{\label{fig:shampine_flame_delta_10m4_sol_qug:d3}}
\end{subfigure}\\
\begin{subfigure}{0.320\textwidth}
\includegraphics[width=\textwidth]{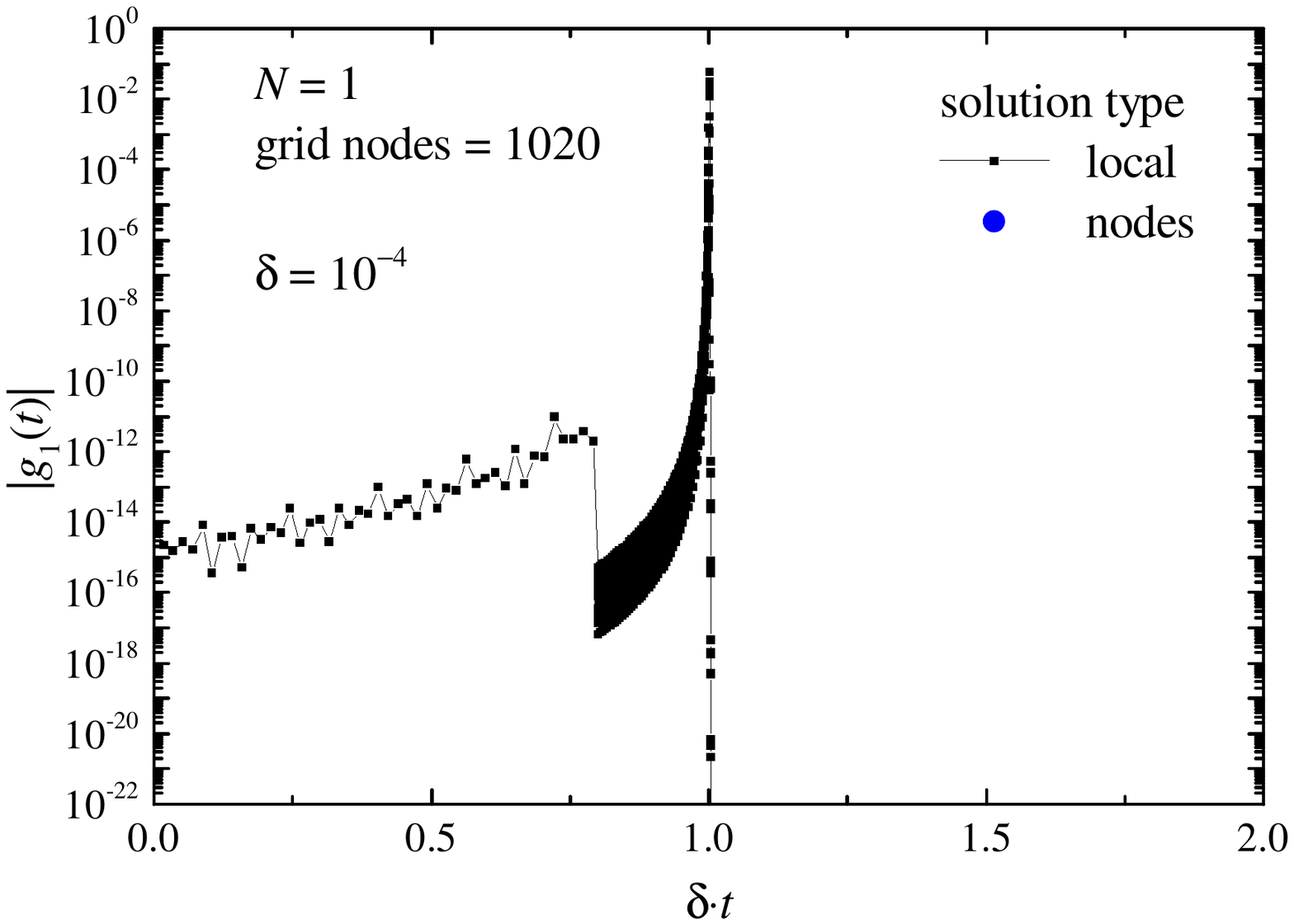}
\vspace{-8mm}\caption{\label{fig:shampine_flame_delta_10m4_sol_qug:e1}}
\end{subfigure}
\begin{subfigure}{0.320\textwidth}
\includegraphics[width=\textwidth]{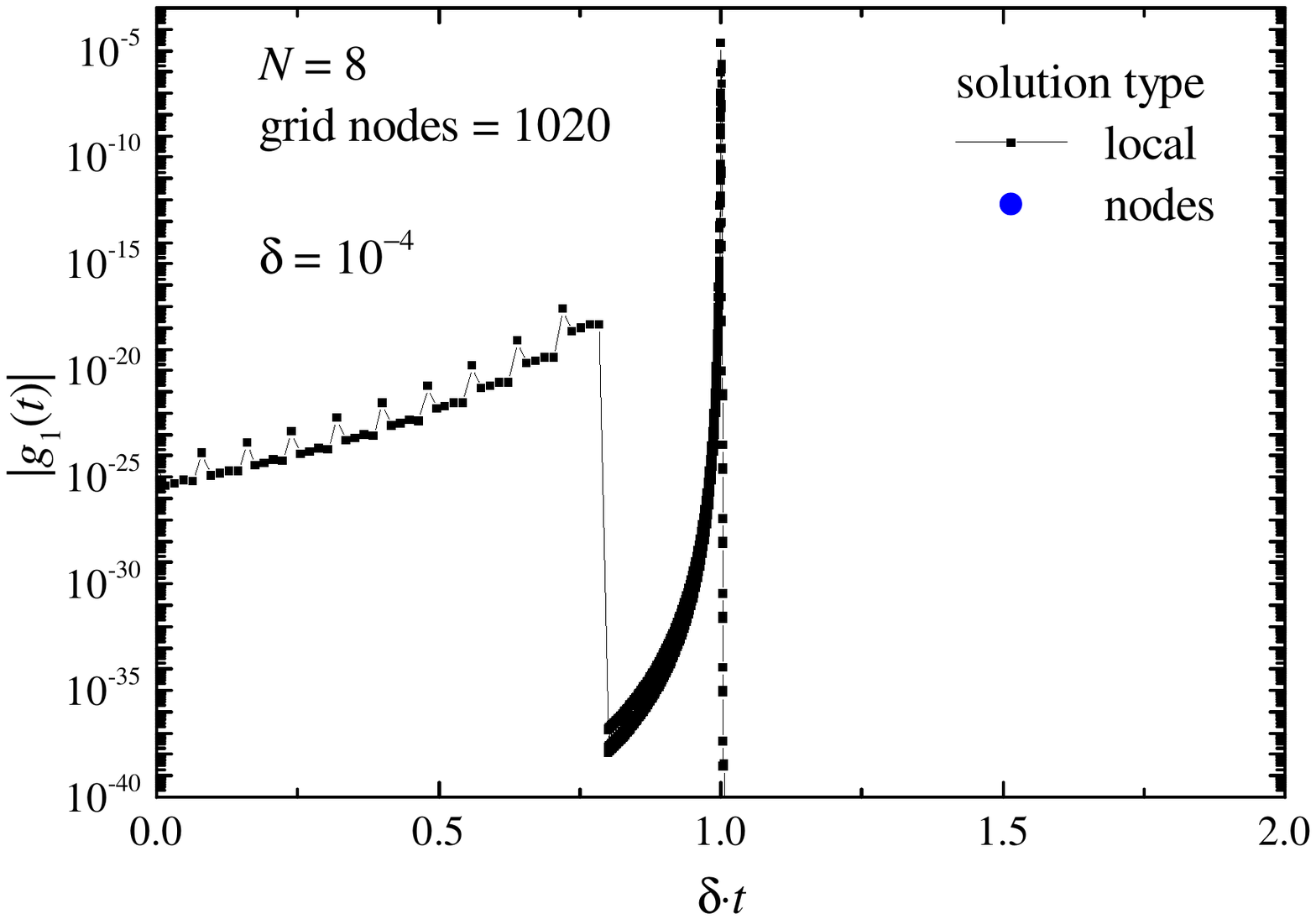}
\vspace{-8mm}\caption{\label{fig:shampine_flame_delta_10m4_sol_qug:e2}}
\end{subfigure}
\begin{subfigure}{0.320\textwidth}
\includegraphics[width=\textwidth]{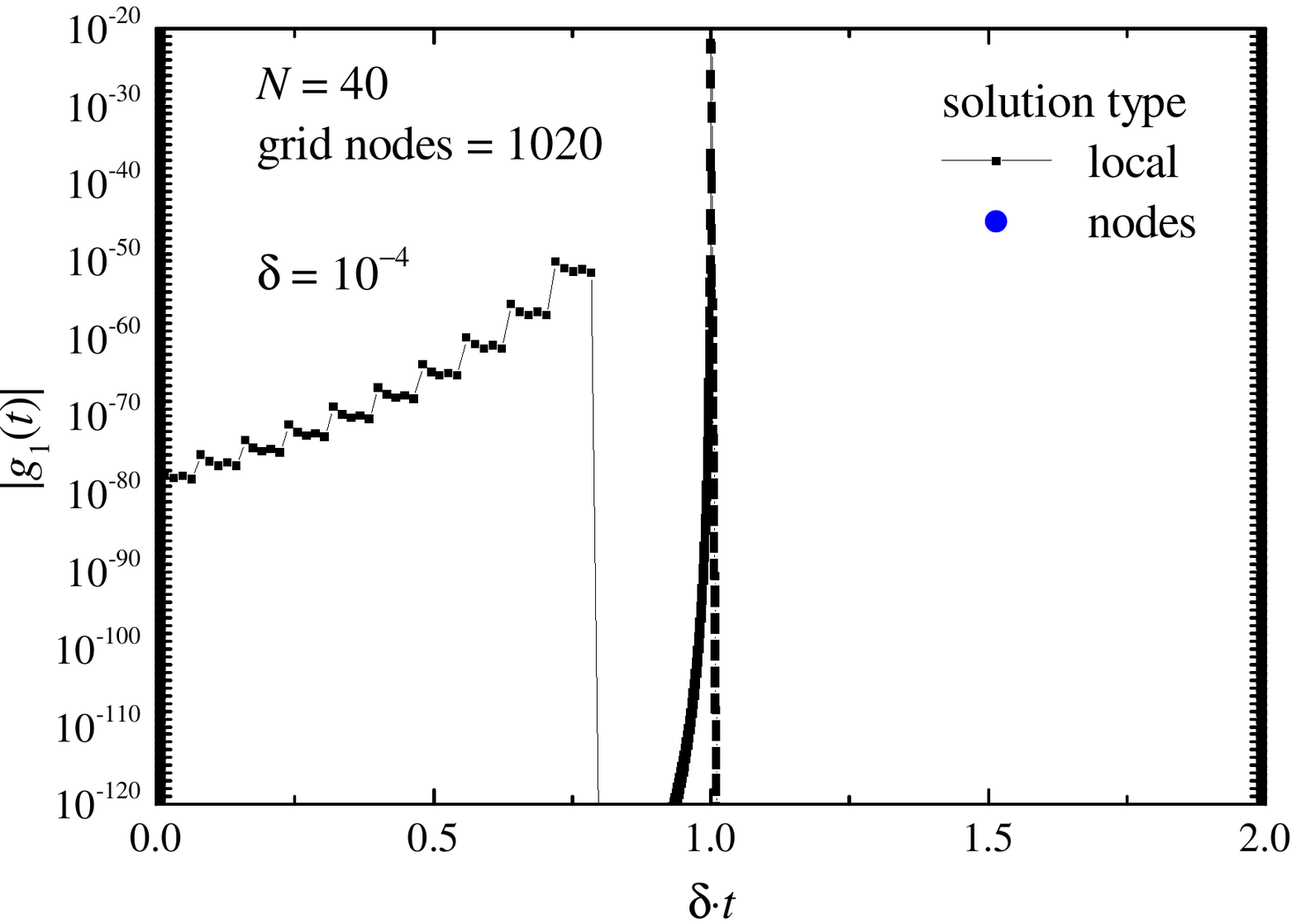}
\vspace{-8mm}\caption{\label{fig:shampine_flame_delta_10m4_sol_qug:e3}}
\end{subfigure}\\
\caption{%
Numerical solution of the stiff DAE system (\ref{eq:shampine_flame}) of index 1 with $\delta = 10^{-4}$. Comparison of the solution at nodes $\mathbf{u}_{n}$ (\subref{fig:shampine_flame_delta_10m4_sol_qug:c1}, \subref{fig:shampine_flame_delta_10m4_sol_qug:c2}, \subref{fig:shampine_flame_delta_10m4_sol_qug:c3}, \subref{fig:shampine_flame_delta_10m4_sol_qug:d1}, \subref{fig:shampine_flame_delta_10m4_sol_qug:d2}, \subref{fig:shampine_flame_delta_10m4_sol_qug:d3}), the local solution $\mathbf{u}_{L}(t)$ (\subref{fig:shampine_flame_delta_10m4_sol_qug:a1}, \subref{fig:shampine_flame_delta_10m4_sol_qug:a2}, \subref{fig:shampine_flame_delta_10m4_sol_qug:a3}, \subref{fig:shampine_flame_delta_10m4_sol_qug:b1}, \subref{fig:shampine_flame_delta_10m4_sol_qug:b2}, \subref{fig:shampine_flame_delta_10m4_sol_qug:b3}) and the exact solution $\mathbf{u}^{\rm ex}(t)$ for component $u_{1}$ (\subref{fig:shampine_flame_delta_10m4_sol_qug:a1}, \subref{fig:shampine_flame_delta_10m4_sol_qug:a2}, \subref{fig:shampine_flame_delta_10m4_sol_qug:a3}, \subref{fig:shampine_flame_delta_10m4_sol_qug:b1}, \subref{fig:shampine_flame_delta_10m4_sol_qug:b2}, \subref{fig:shampine_flame_delta_10m4_sol_qug:b3}, \subref{fig:shampine_flame_delta_10m4_sol_qug:c1}, \subref{fig:shampine_flame_delta_10m4_sol_qug:c2}, \subref{fig:shampine_flame_delta_10m4_sol_qug:c3}, \subref{fig:shampine_flame_delta_10m4_sol_qug:d1}, \subref{fig:shampine_flame_delta_10m4_sol_qug:d2}, \subref{fig:shampine_flame_delta_10m4_sol_qug:d3}), quantitative satisfiability of the conditions $g_{1} = 0$ (\subref{fig:shampine_flame_delta_10m4_sol_qug:e1}, \subref{fig:shampine_flame_delta_10m4_sol_qug:e2}, \subref{fig:shampine_flame_delta_10m4_sol_qug:e3}), obtained using polynomials with degrees $N = 1$ (\subref{fig:shampine_flame_delta_10m4_sol_qug:a1}, \subref{fig:shampine_flame_delta_10m4_sol_qug:b1}, \subref{fig:shampine_flame_delta_10m4_sol_qug:c1}, \subref{fig:shampine_flame_delta_10m4_sol_qug:d1}, \subref{fig:shampine_flame_delta_10m4_sol_qug:e1}), $N = 8$ (\subref{fig:shampine_flame_delta_10m4_sol_qug:a2}, \subref{fig:shampine_flame_delta_10m4_sol_qug:b2}, \subref{fig:shampine_flame_delta_10m4_sol_qug:c2}, \subref{fig:shampine_flame_delta_10m4_sol_qug:d2}, \subref{fig:shampine_flame_delta_10m4_sol_qug:e2}) and $N = 40$ (\subref{fig:shampine_flame_delta_10m4_sol_qug:a3}, \subref{fig:shampine_flame_delta_10m4_sol_qug:b3}, \subref{fig:shampine_flame_delta_10m4_sol_qug:c3}, \subref{fig:shampine_flame_delta_10m4_sol_qug:d3}, \subref{fig:shampine_flame_delta_10m4_sol_qug:e3}).
}
\label{fig:shampine_flame_delta_10m4_sol_qug}
\end{figure}

\begin{figure}[h!]
\captionsetup[subfigure]{%
	position=bottom,
	font+=smaller,
	textfont=normalfont,
	singlelinecheck=off,
	justification=raggedright
}
\centering
\begin{subfigure}{0.320\textwidth}
\includegraphics[width=\textwidth]{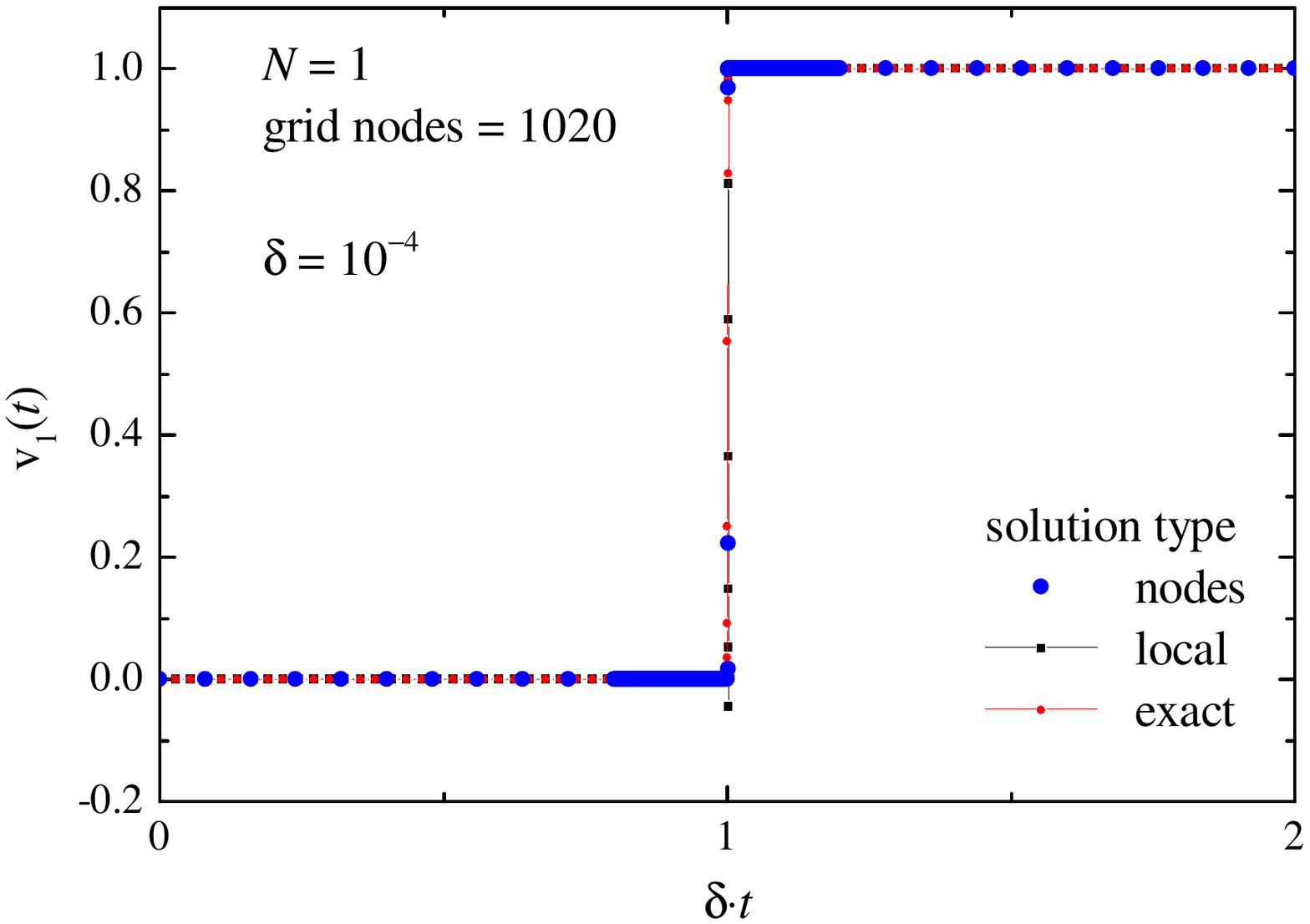}
\vspace{-8mm}\caption{\label{fig:shampine_flame_delta_10m4_sol_v_epss:a1}}
\end{subfigure}
\begin{subfigure}{0.320\textwidth}
\includegraphics[width=\textwidth]{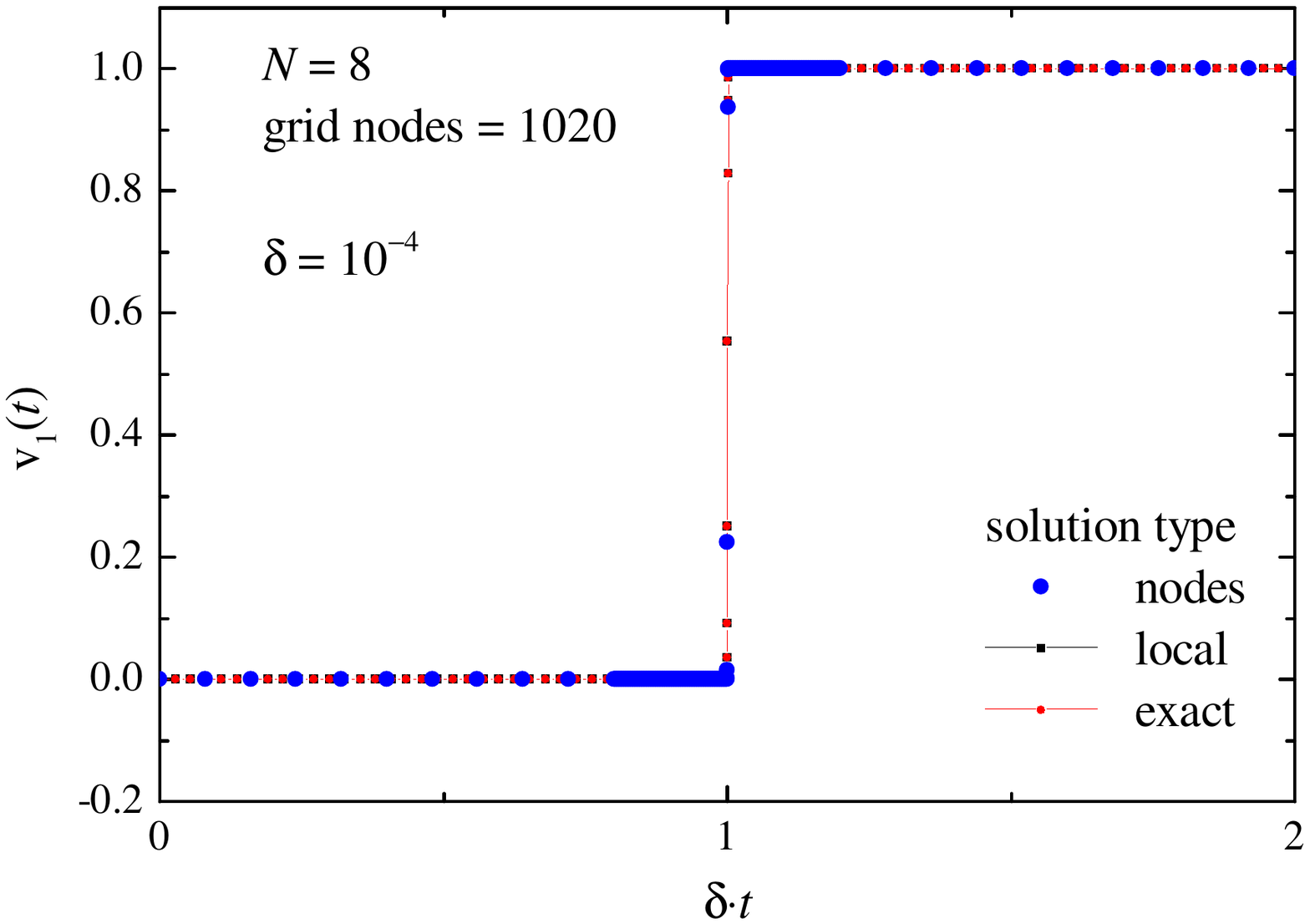}
\vspace{-8mm}\caption{\label{fig:shampine_flame_delta_10m4_sol_v_epss:a2}}
\end{subfigure}
\begin{subfigure}{0.320\textwidth}
\includegraphics[width=\textwidth]{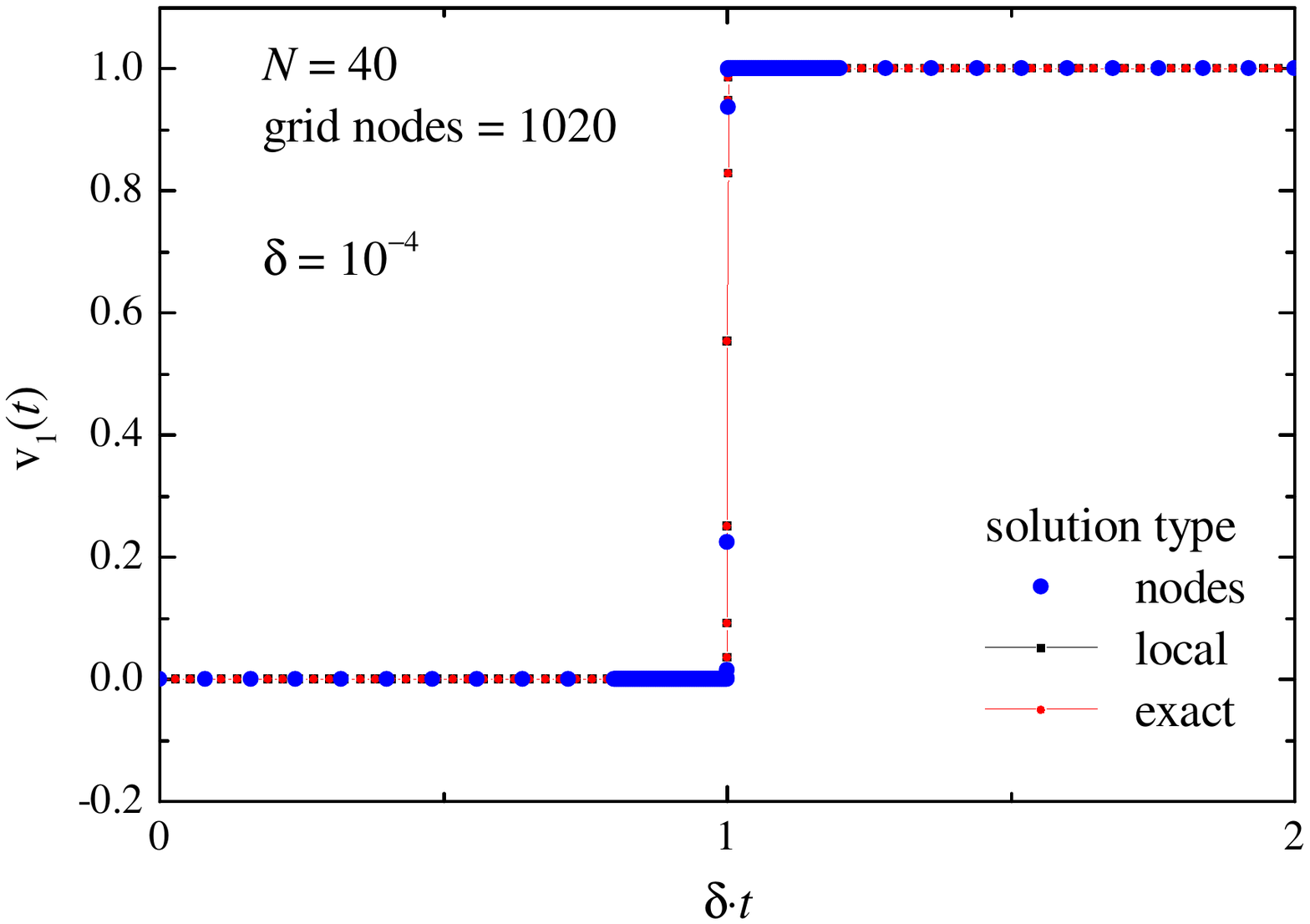}
\vspace{-8mm}\caption{\label{fig:shampine_flame_delta_10m4_sol_v_epss:a3}}
\end{subfigure}\\
\begin{subfigure}{0.320\textwidth}
\includegraphics[width=\textwidth]{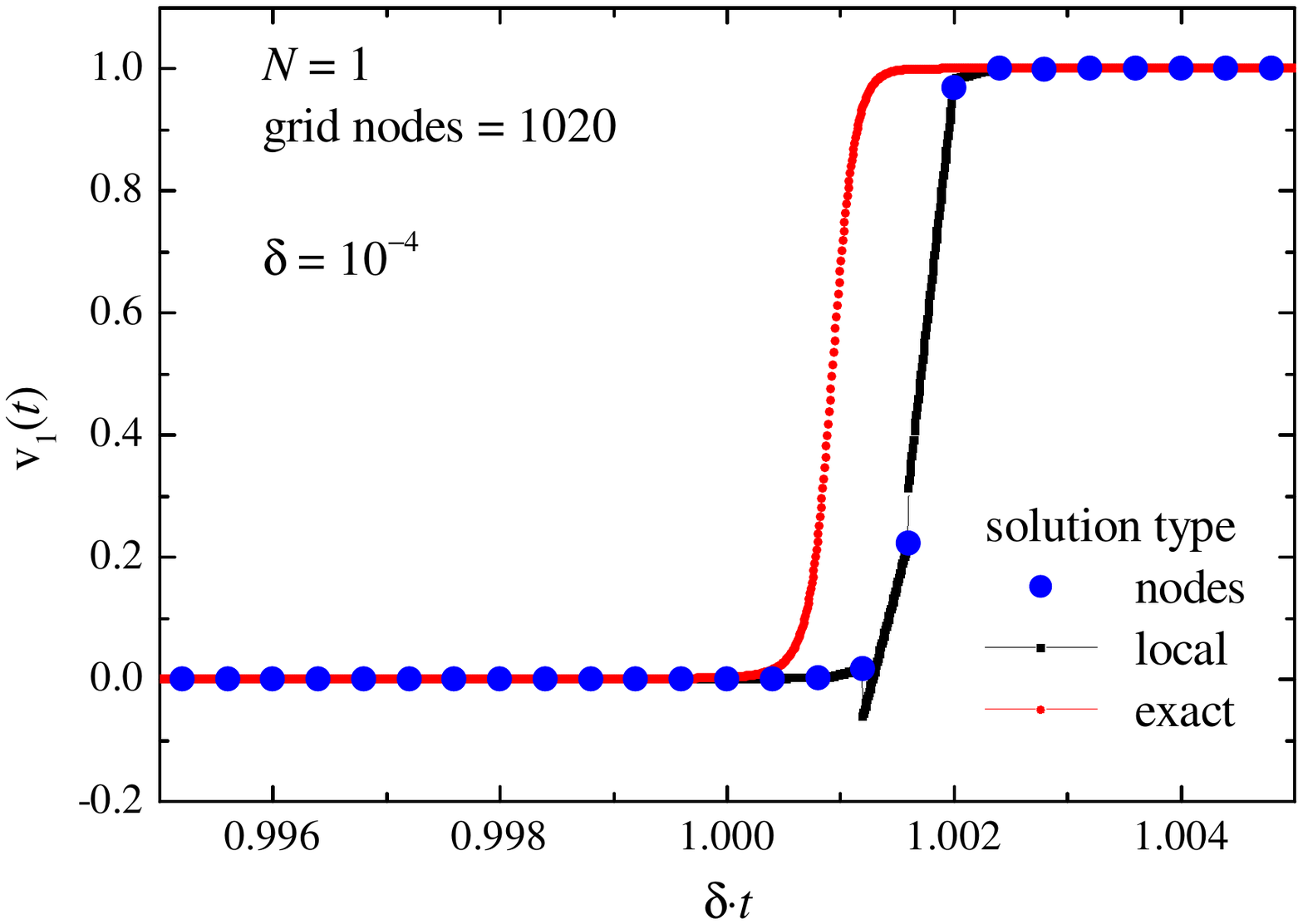}
\vspace{-8mm}\caption{\label{fig:shampine_flame_delta_10m4_sol_v_epss:b1}}
\end{subfigure}
\begin{subfigure}{0.320\textwidth}
\includegraphics[width=\textwidth]{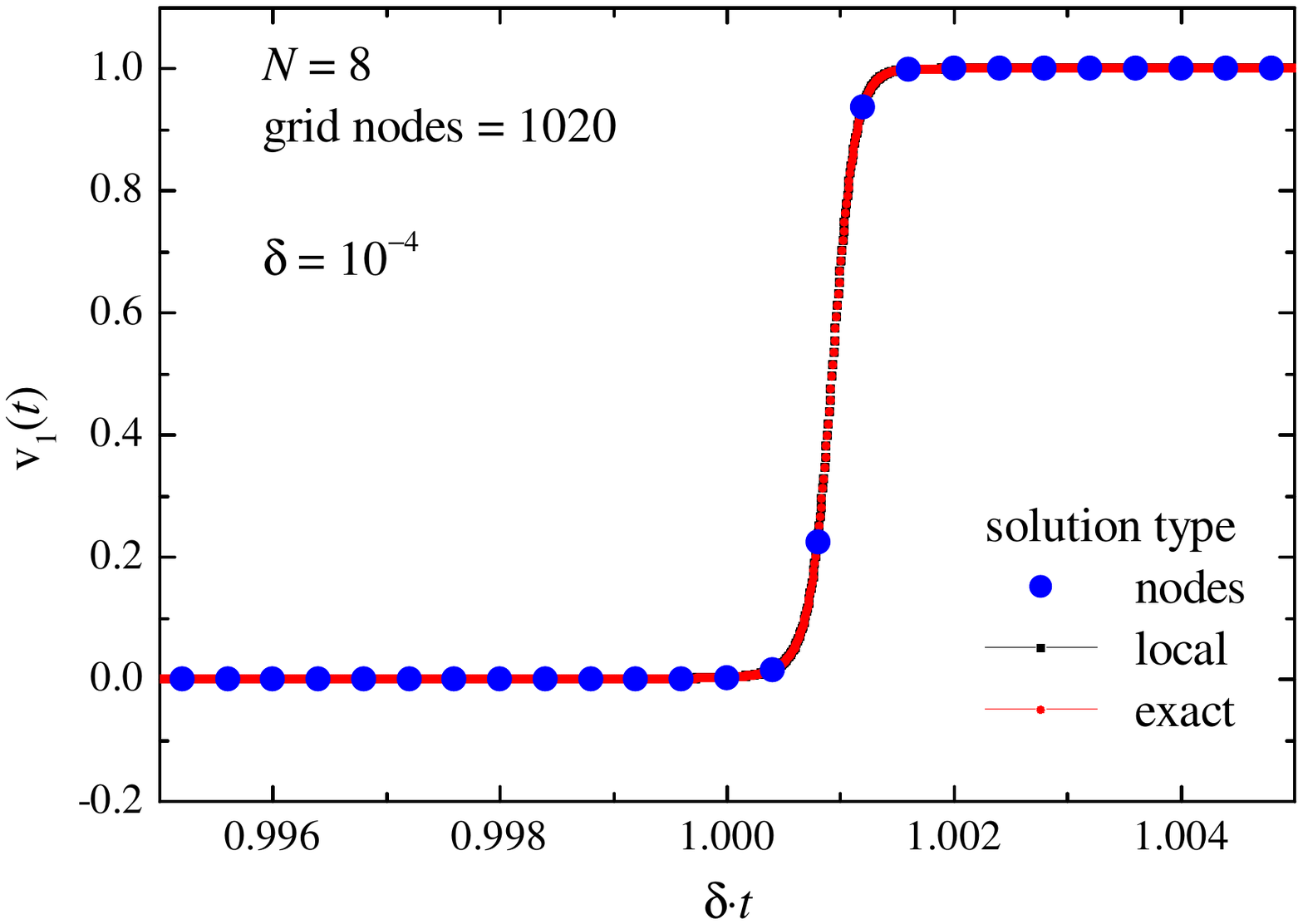}
\vspace{-8mm}\caption{\label{fig:shampine_flame_delta_10m4_sol_v_epss:b2}}
\end{subfigure}
\begin{subfigure}{0.320\textwidth}
\includegraphics[width=\textwidth]{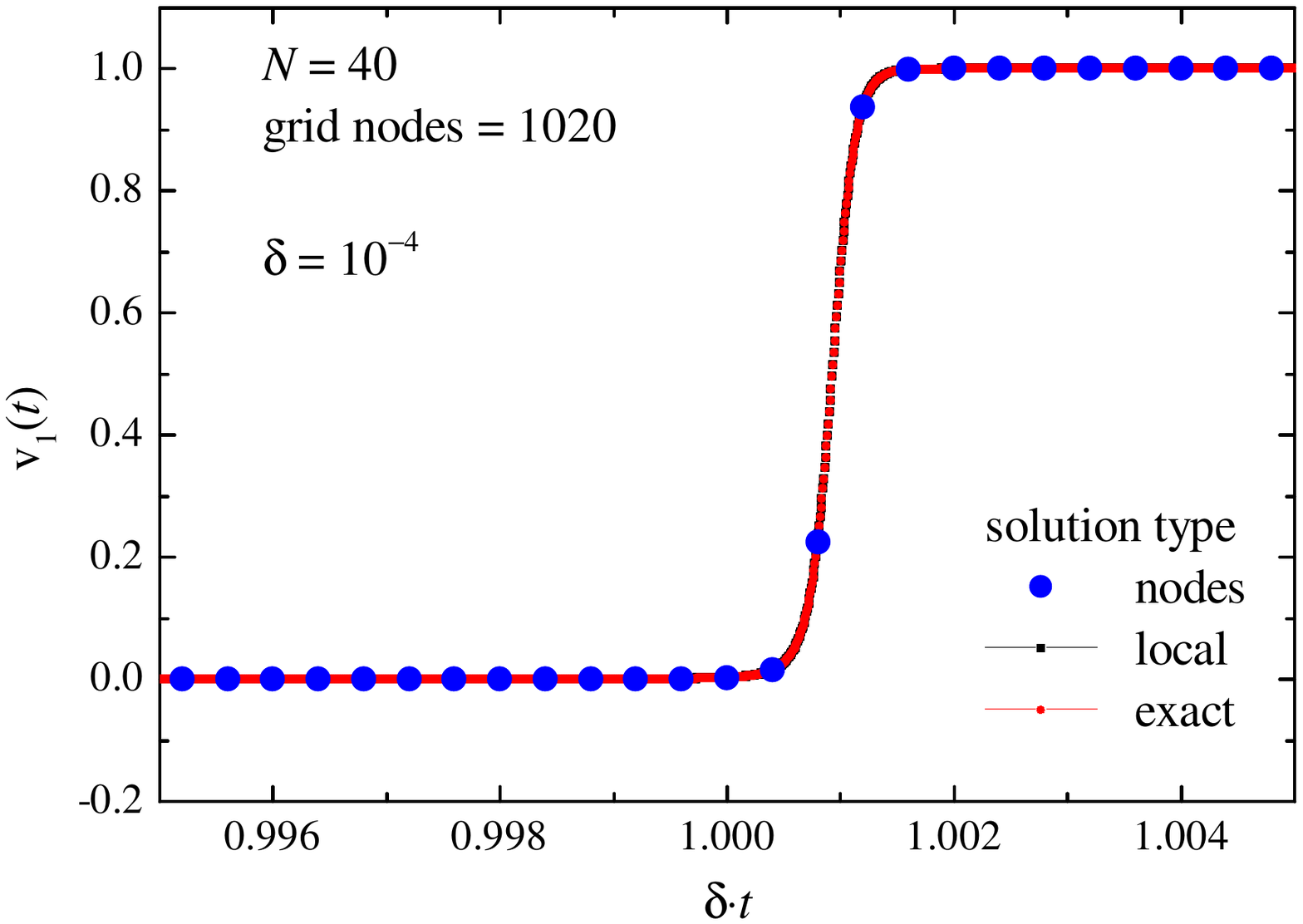}
\vspace{-8mm}\caption{\label{fig:shampine_flame_delta_10m4_sol_v_epss:b3}}
\end{subfigure}\\
\begin{subfigure}{0.320\textwidth}
\includegraphics[width=\textwidth]{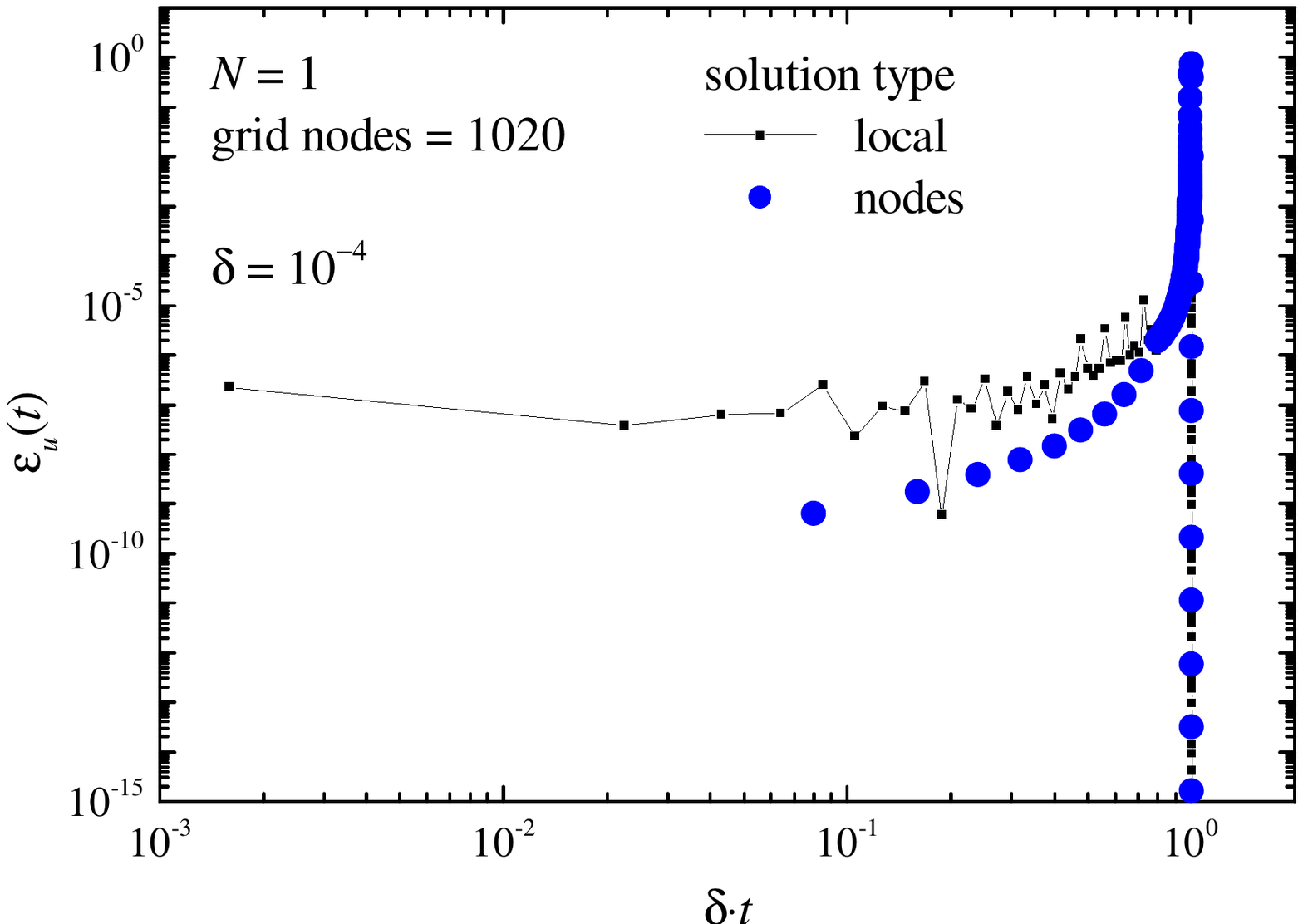}
\vspace{-8mm}\caption{\label{fig:shampine_flame_delta_10m4_sol_v_epss:c1}}
\end{subfigure}
\begin{subfigure}{0.320\textwidth}
\includegraphics[width=\textwidth]{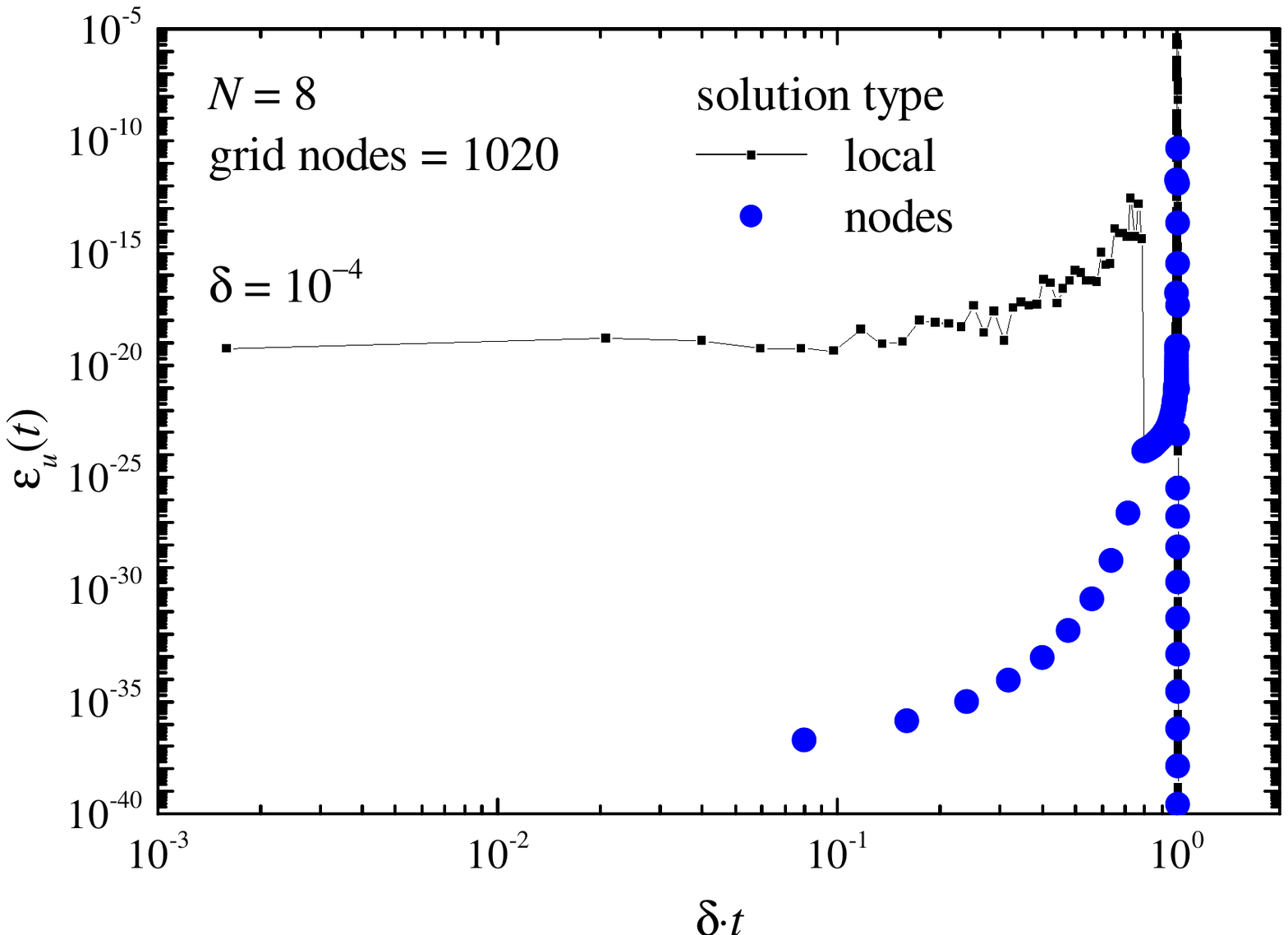}
\vspace{-8mm}\caption{\label{fig:shampine_flame_delta_10m4_sol_v_epss:c2}}
\end{subfigure}
\begin{subfigure}{0.320\textwidth}
\includegraphics[width=\textwidth]{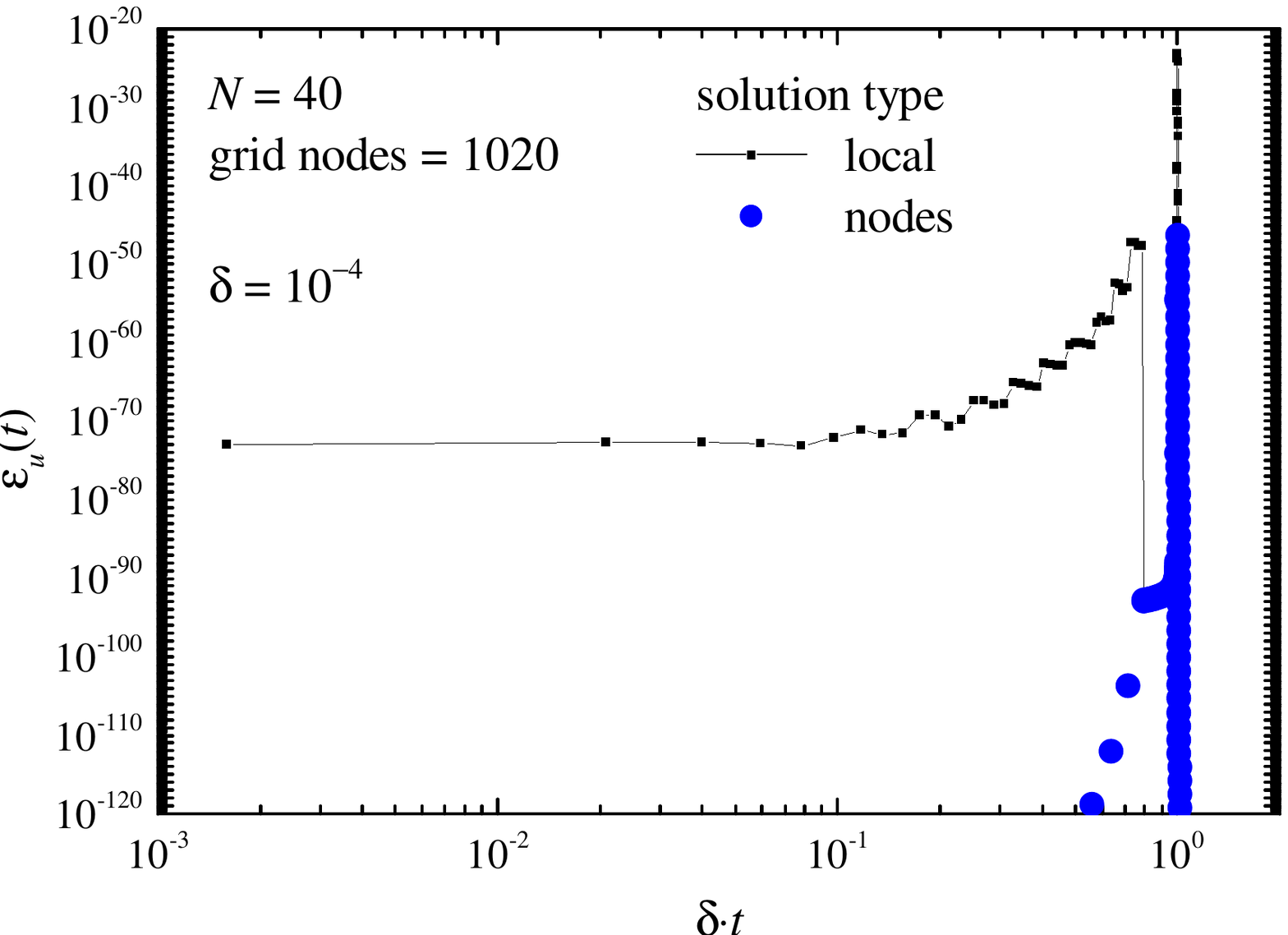}
\vspace{-8mm}\caption{\label{fig:shampine_flame_delta_10m4_sol_v_epss:c3}}
\end{subfigure}\\
\begin{subfigure}{0.320\textwidth}
\includegraphics[width=\textwidth]{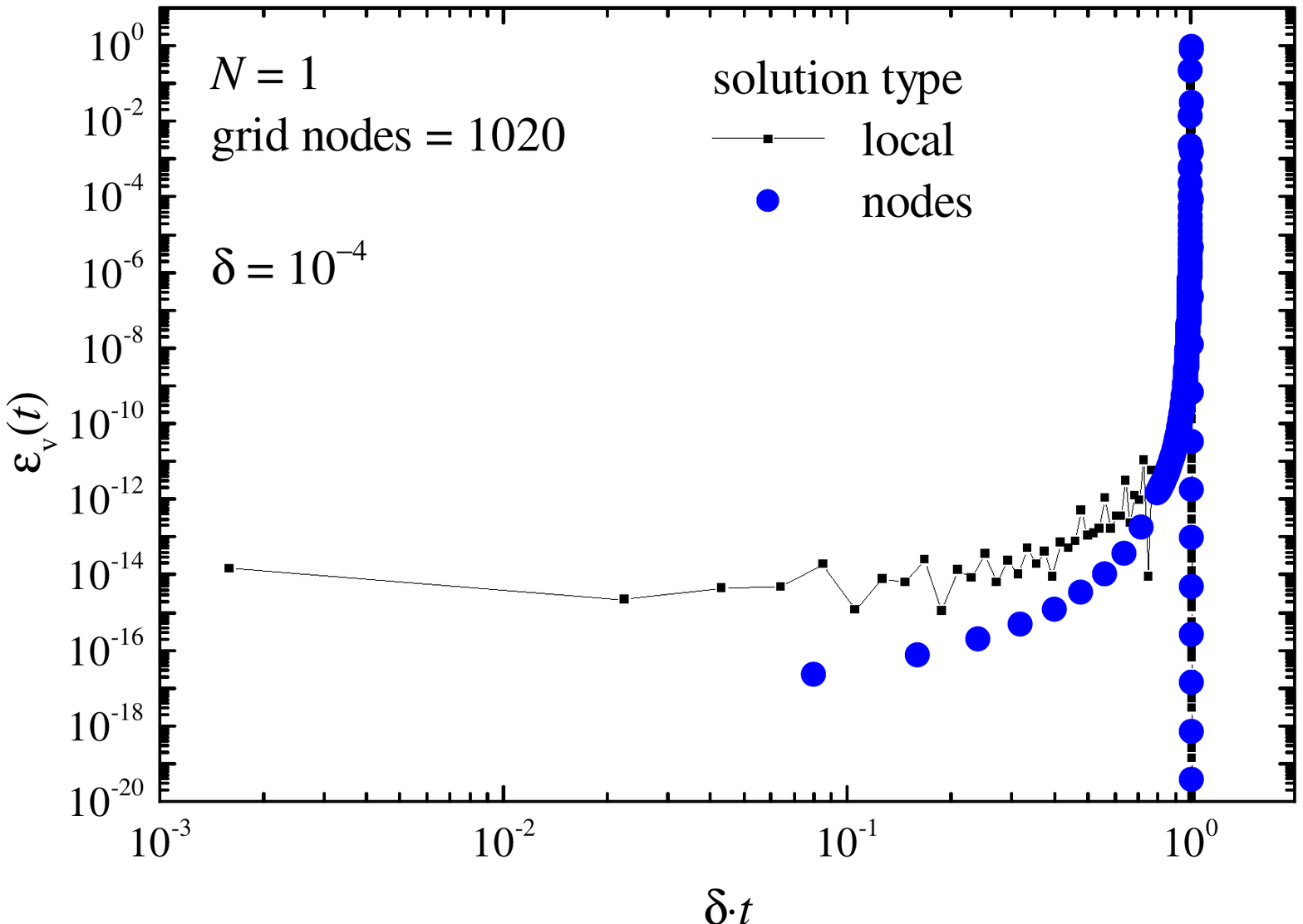}
\vspace{-8mm}\caption{\label{fig:shampine_flame_delta_10m4_sol_v_epss:d1}}
\end{subfigure}
\begin{subfigure}{0.320\textwidth}
\includegraphics[width=\textwidth]{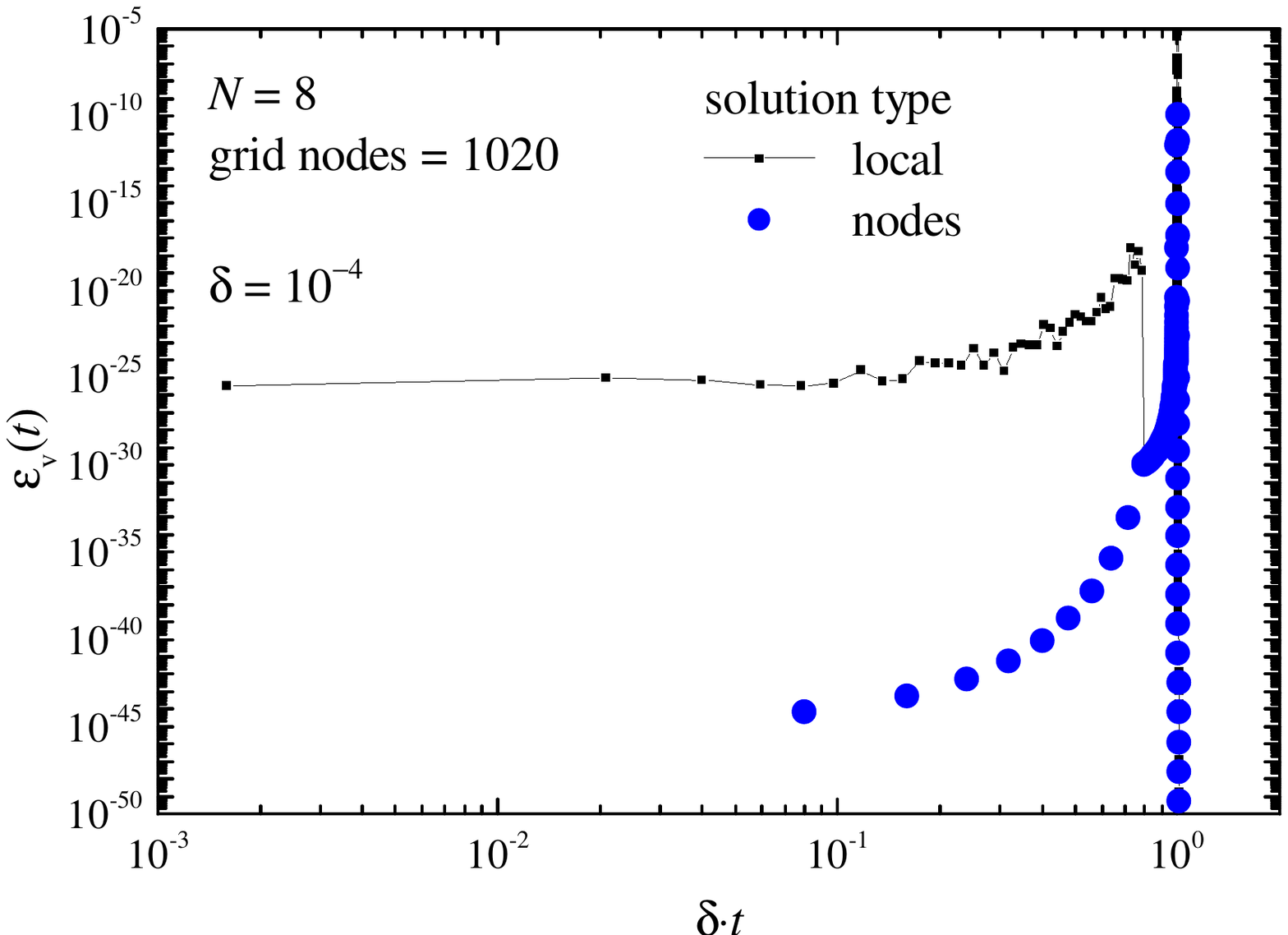}
\vspace{-8mm}\caption{\label{fig:shampine_flame_delta_10m4_sol_v_epss:d2}}
\end{subfigure}
\begin{subfigure}{0.320\textwidth}
\includegraphics[width=\textwidth]{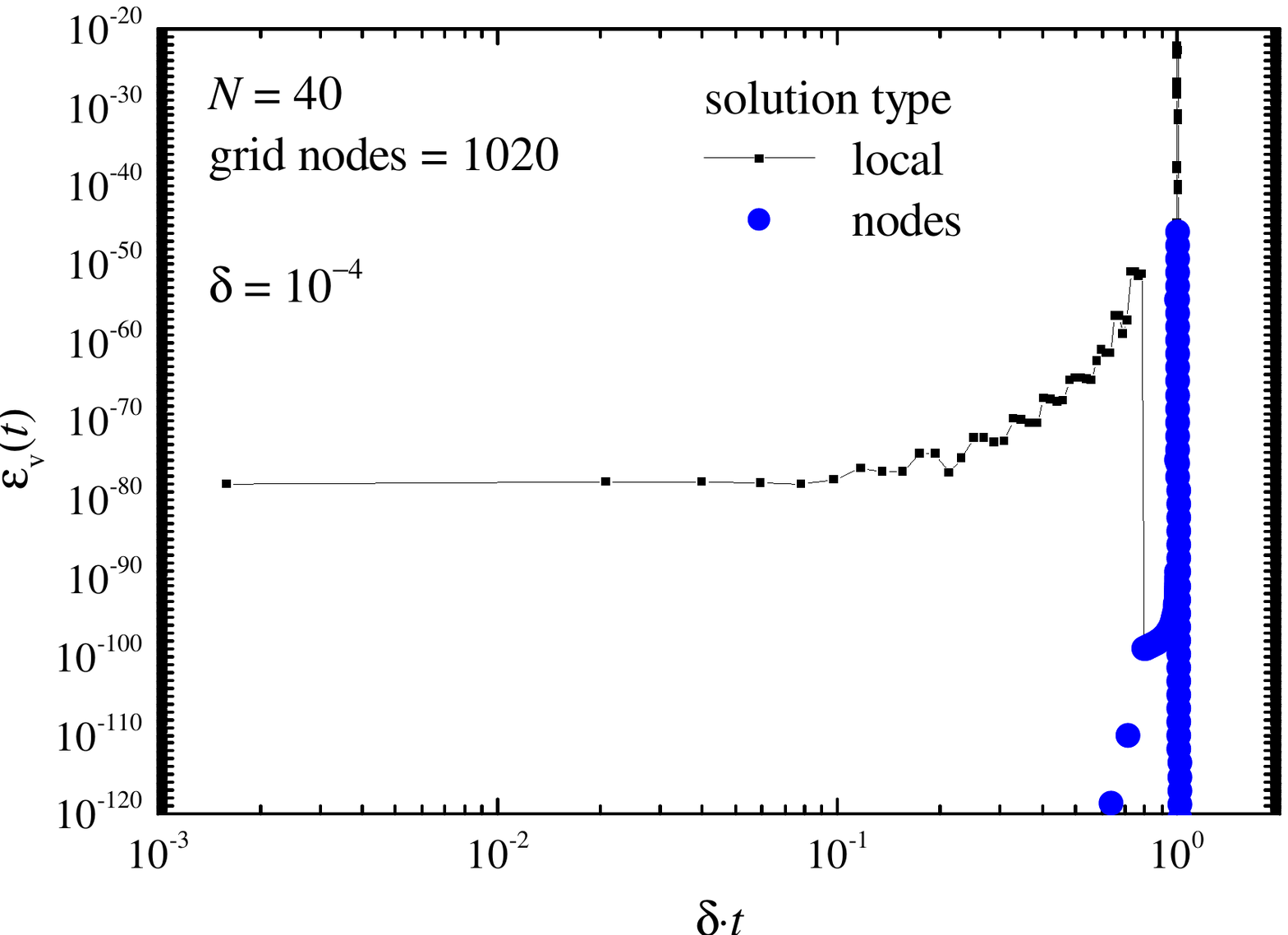}
\vspace{-8mm}\caption{\label{fig:shampine_flame_delta_10m4_sol_v_epss:d3}}
\end{subfigure}\\
\begin{subfigure}{0.320\textwidth}
\includegraphics[width=\textwidth]{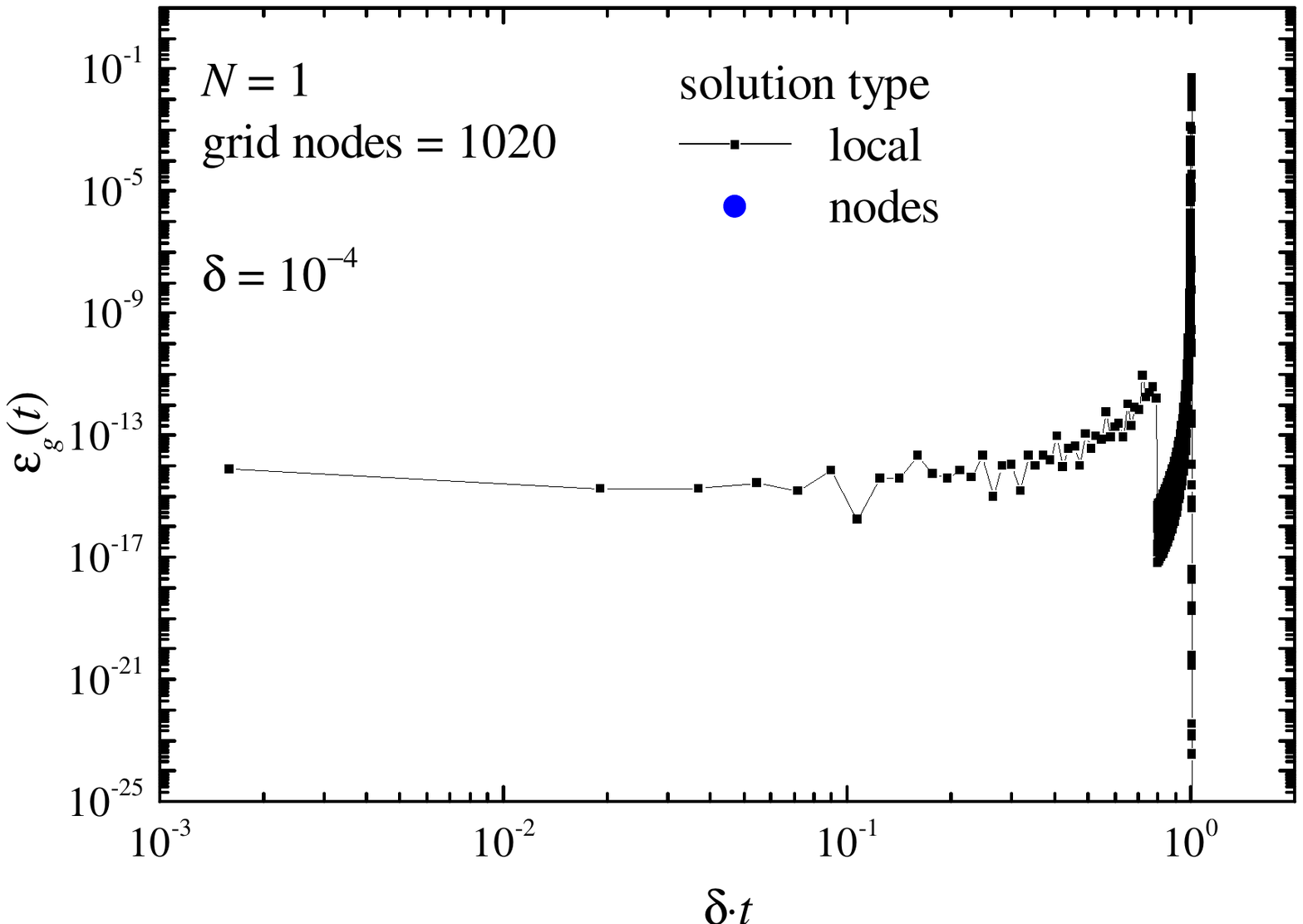}
\vspace{-8mm}\caption{\label{fig:shampine_flame_delta_10m4_sol_v_epss:e1}}
\end{subfigure}
\begin{subfigure}{0.320\textwidth}
\includegraphics[width=\textwidth]{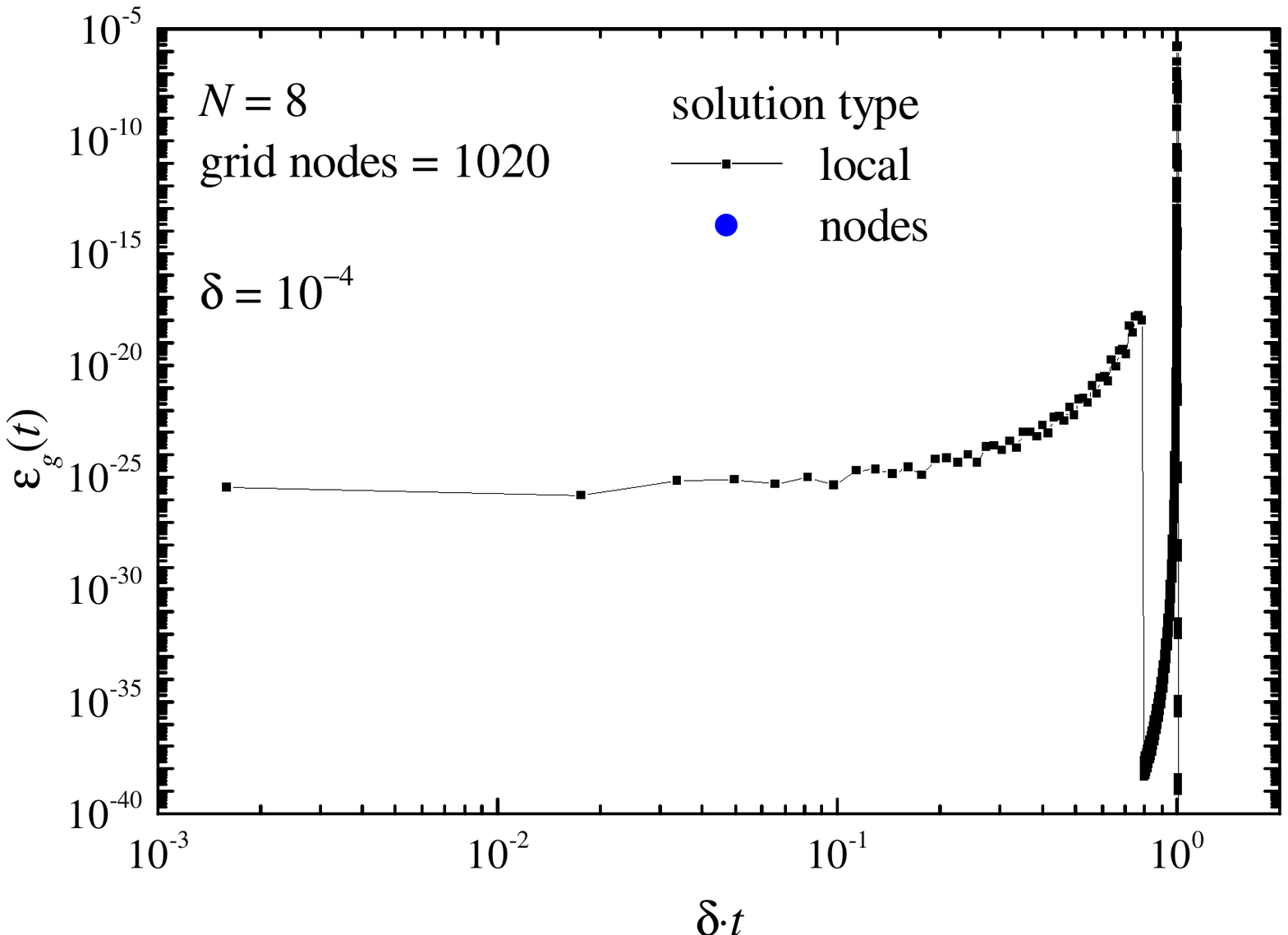}
\vspace{-8mm}\caption{\label{fig:shampine_flame_delta_10m4_sol_v_epss:e2}}
\end{subfigure}
\begin{subfigure}{0.320\textwidth}
\includegraphics[width=\textwidth]{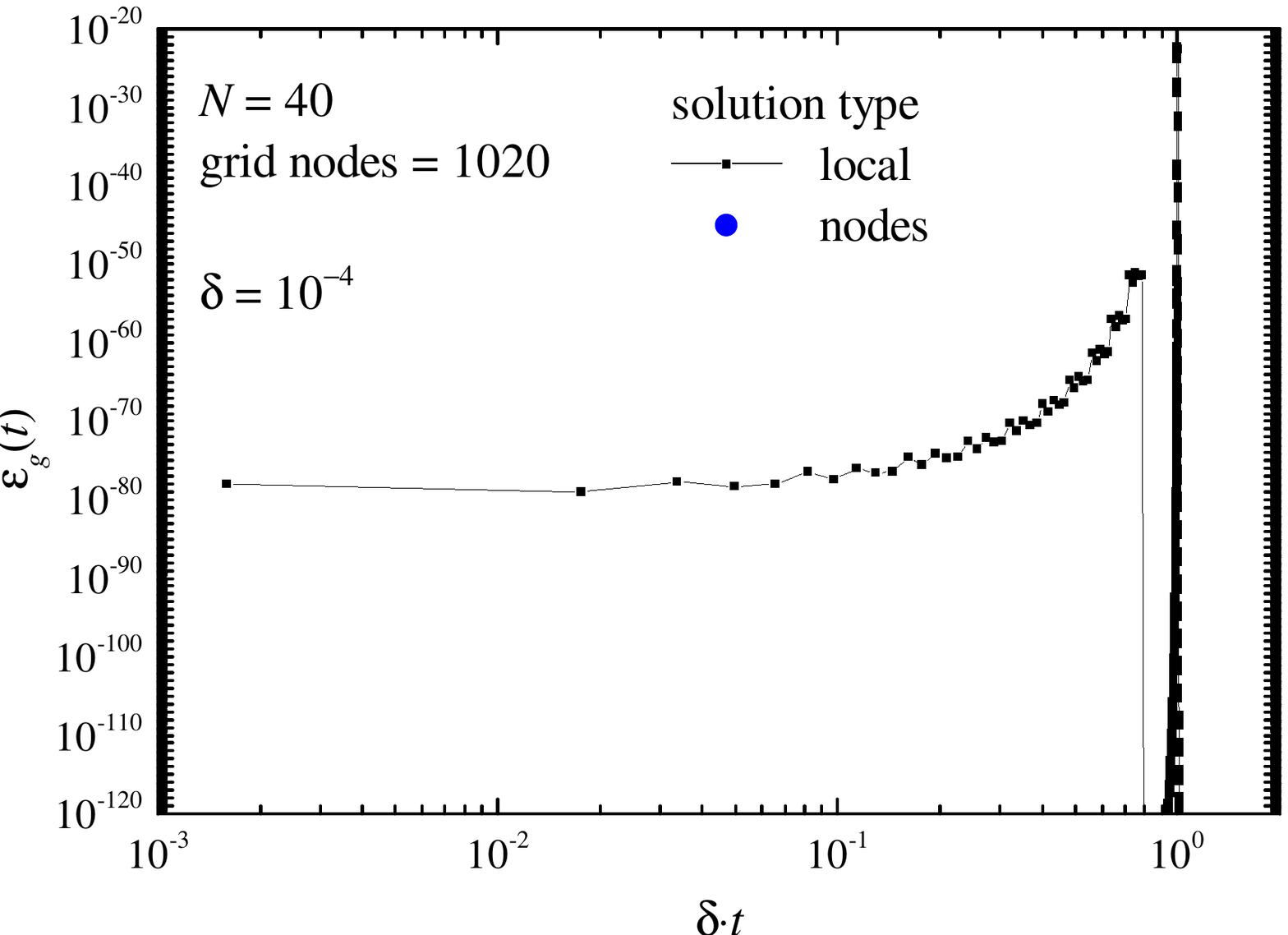}
\vspace{-8mm}\caption{\label{fig:shampine_flame_delta_10m4_sol_v_epss:e3}}
\end{subfigure}\\
\caption{%
Numerical solution of the stiff DAE system (\ref{eq:shampine_flame}) of index 1 with $\delta = 10^{-4}$. Comparison of the solution at nodes $\mathbf{v}_{n}$, the local solution $\mathbf{v}_{L}(t)$ and the exact solution $\mathbf{v}^{\rm ex}(t)$ for component $v_{1}$ (\subref{fig:shampine_flame_delta_10m4_sol_v_epss:a1}, \subref{fig:shampine_flame_delta_10m4_sol_v_epss:a2}, \subref{fig:shampine_flame_delta_10m4_sol_v_epss:a3}, \subref{fig:shampine_flame_delta_10m4_sol_v_epss:b1}, \subref{fig:shampine_flame_delta_10m4_sol_v_epss:b2}, \subref{fig:shampine_flame_delta_10m4_sol_v_epss:b3}), the errors $\varepsilon_{u}(t)$ (\subref{fig:shampine_flame_delta_10m4_sol_v_epss:c1}, \subref{fig:shampine_flame_delta_10m4_sol_v_epss:c2}, \subref{fig:shampine_flame_delta_10m4_sol_v_epss:c3}), $\varepsilon_{v}(t)$ (\subref{fig:shampine_flame_delta_10m4_sol_v_epss:d1}, \subref{fig:shampine_flame_delta_10m4_sol_v_epss:d2}, \subref{fig:shampine_flame_delta_10m4_sol_v_epss:d3}), $\varepsilon_{g}(t)$ (\subref{fig:shampine_flame_delta_10m4_sol_qug:e1}, \subref{fig:shampine_flame_delta_10m4_sol_qug:e2}, \subref{fig:shampine_flame_delta_10m4_sol_qug:e3}), obtained using polynomials with degrees $N = 1$ (\subref{fig:shampine_flame_delta_10m4_sol_qug:a1}, \subref{fig:shampine_flame_delta_10m4_sol_qug:b1}, \subref{fig:shampine_flame_delta_10m4_sol_qug:c1}, \subref{fig:shampine_flame_delta_10m4_sol_qug:d1}, \subref{fig:shampine_flame_delta_10m4_sol_qug:e1}), $N = 8$ (\subref{fig:shampine_flame_delta_10m4_sol_qug:a2}, \subref{fig:shampine_flame_delta_10m4_sol_qug:b2}, \subref{fig:shampine_flame_delta_10m4_sol_qug:c2}, \subref{fig:shampine_flame_delta_10m4_sol_qug:d2}, \subref{fig:shampine_flame_delta_10m4_sol_qug:e2}) and $N = 40$ (\subref{fig:shampine_flame_delta_10m4_sol_qug:a3}, \subref{fig:shampine_flame_delta_10m4_sol_qug:b3}, \subref{fig:shampine_flame_delta_10m4_sol_qug:c3}, \subref{fig:shampine_flame_delta_10m4_sol_qug:d3}, \subref{fig:shampine_flame_delta_10m4_sol_qug:e3}).
}
\label{fig:shampine_flame_delta_10m4_sol_v_epss}
\end{figure}

\begin{figure}[h!]
\captionsetup[subfigure]{%
	position=bottom,
	font+=smaller,
	textfont=normalfont,
	singlelinecheck=off,
	justification=raggedright
}
\centering
\begin{subfigure}{0.320\textwidth}
\includegraphics[width=\textwidth]{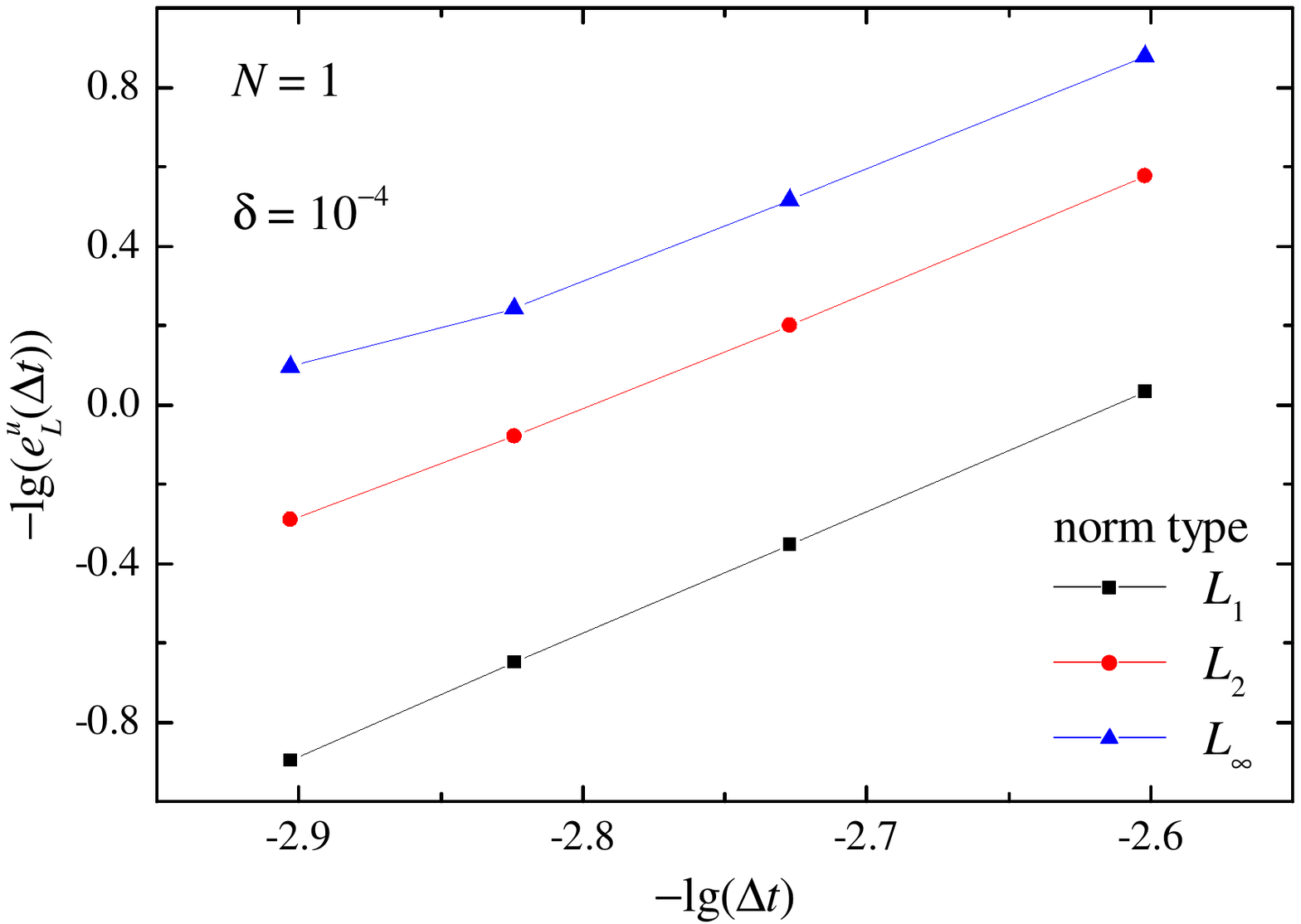}
\vspace{-8mm}\caption{\label{fig:shampine_flame_delta_10m4_errors:a1}}
\end{subfigure}
\begin{subfigure}{0.320\textwidth}
\includegraphics[width=\textwidth]{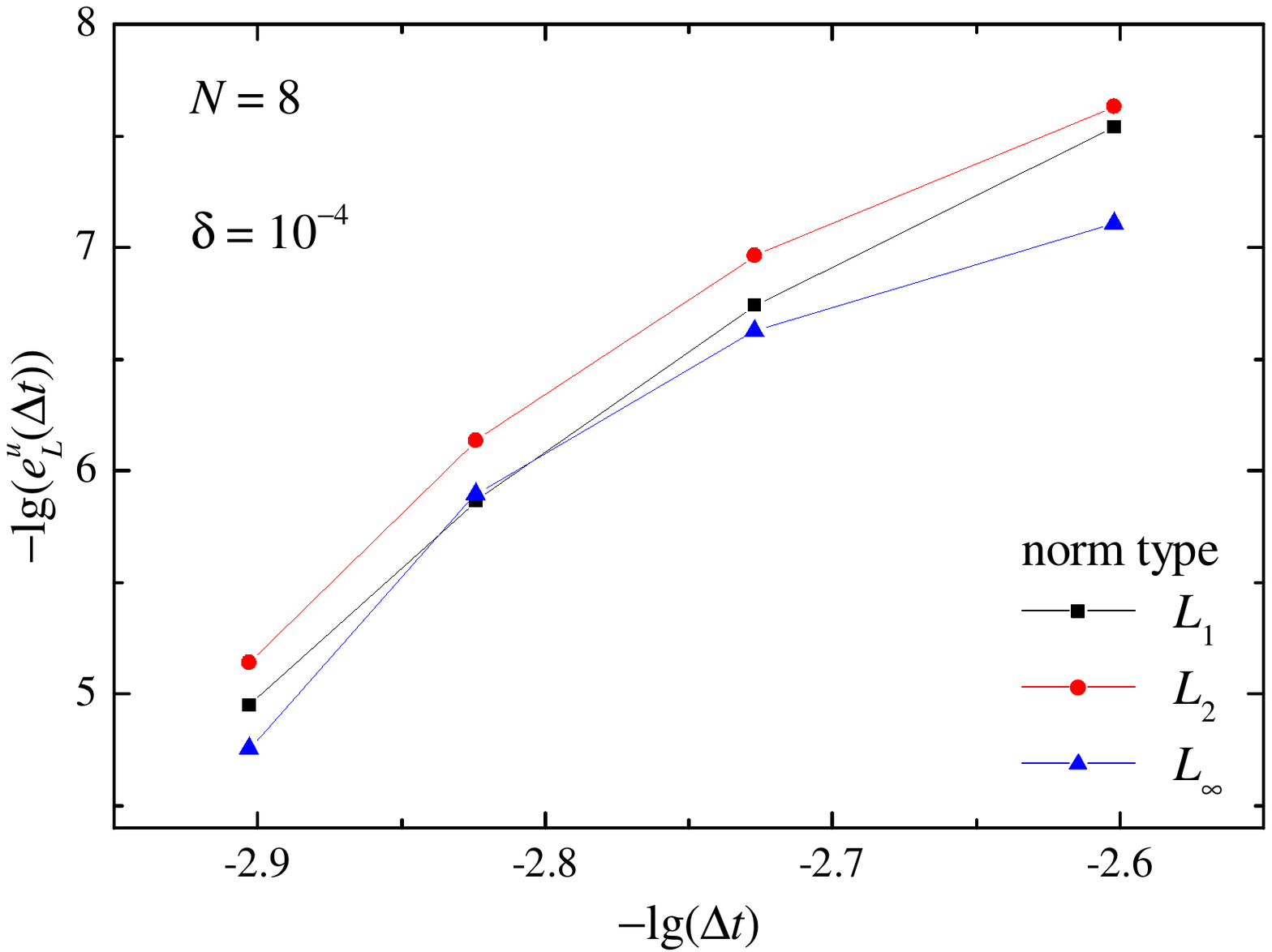}
\vspace{-8mm}\caption{\label{fig:shampine_flame_delta_10m4_errors:a2}}
\end{subfigure}
\begin{subfigure}{0.320\textwidth}
\includegraphics[width=\textwidth]{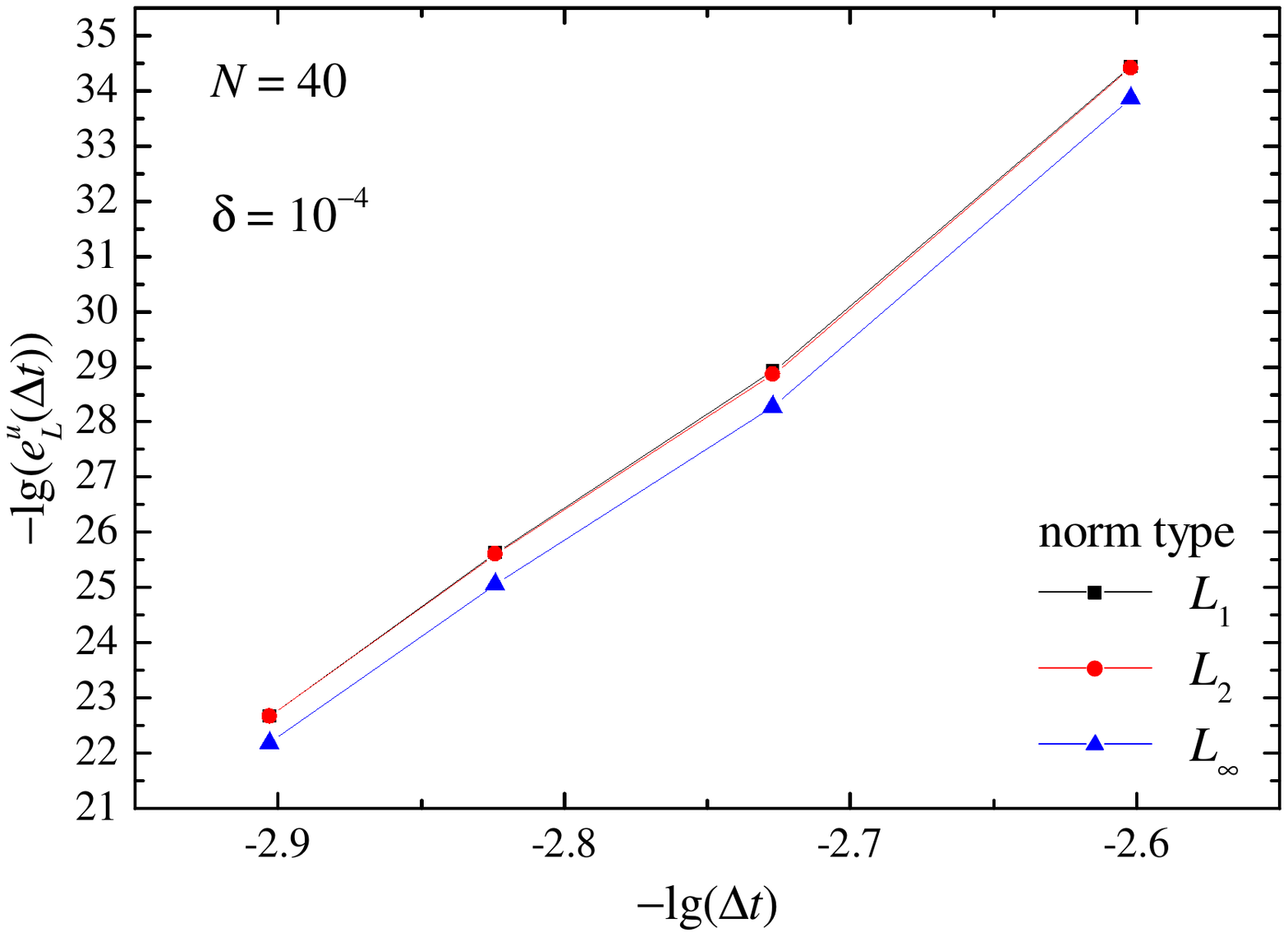}
\vspace{-8mm}\caption{\label{fig:shampine_flame_delta_10m4_errors:a3}}
\end{subfigure}\\
\begin{subfigure}{0.320\textwidth}
\includegraphics[width=\textwidth]{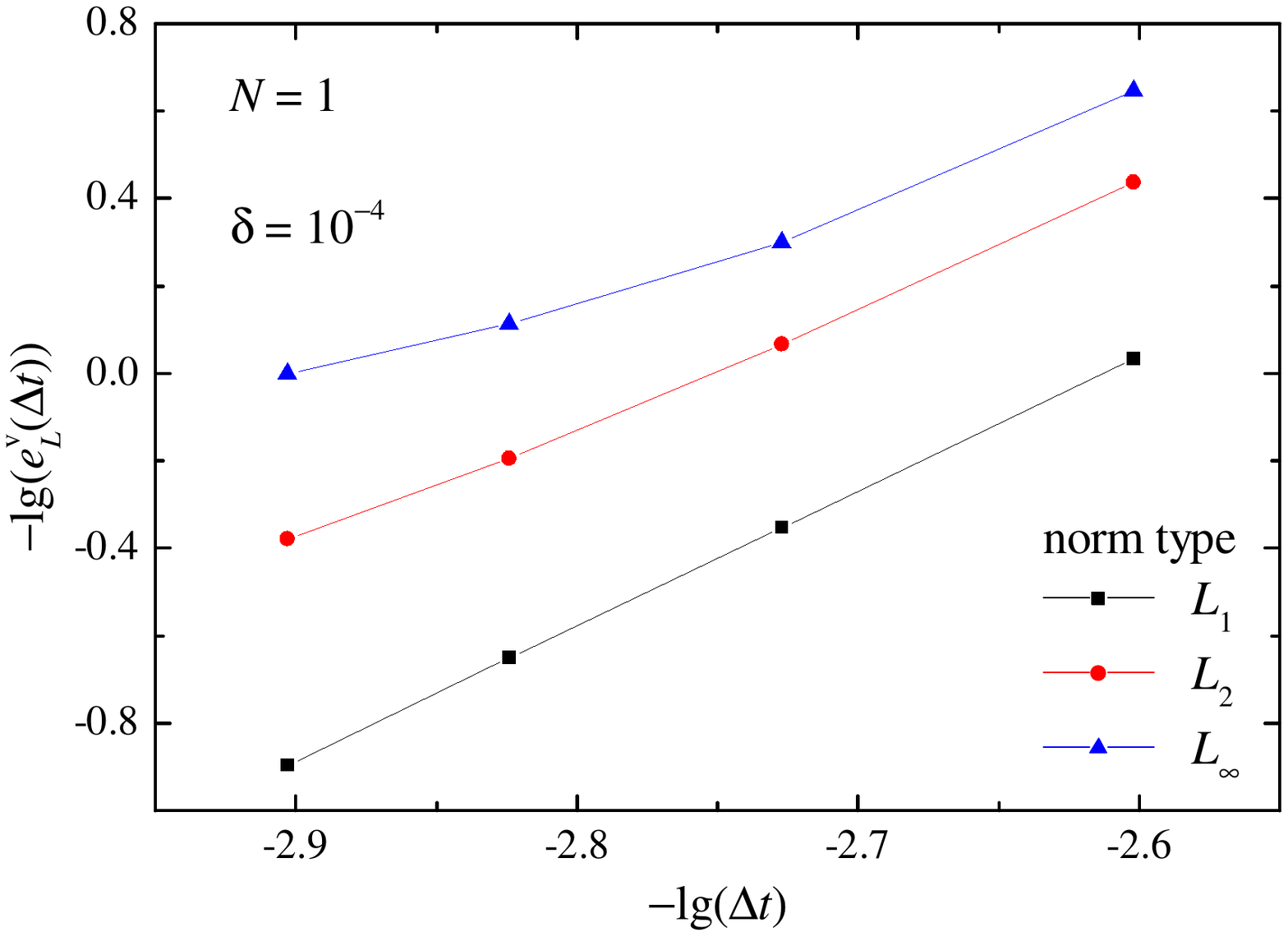}
\vspace{-8mm}\caption{\label{fig:shampine_flame_delta_10m4_errors:b1}}
\end{subfigure}
\begin{subfigure}{0.320\textwidth}
\includegraphics[width=\textwidth]{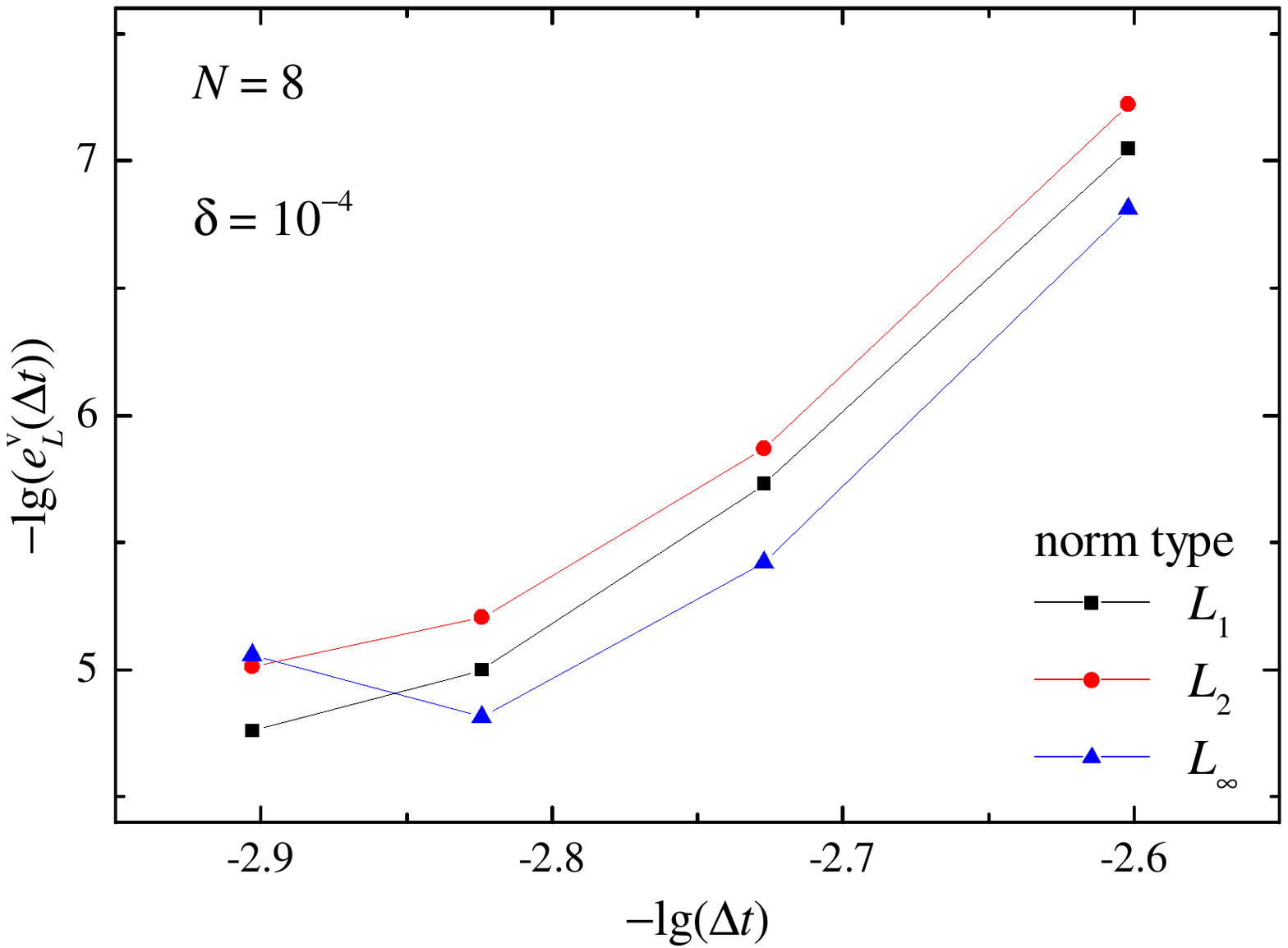}
\vspace{-8mm}\caption{\label{fig:shampine_flame_delta_10m4_errors:b2}}
\end{subfigure}
\begin{subfigure}{0.320\textwidth}
\includegraphics[width=\textwidth]{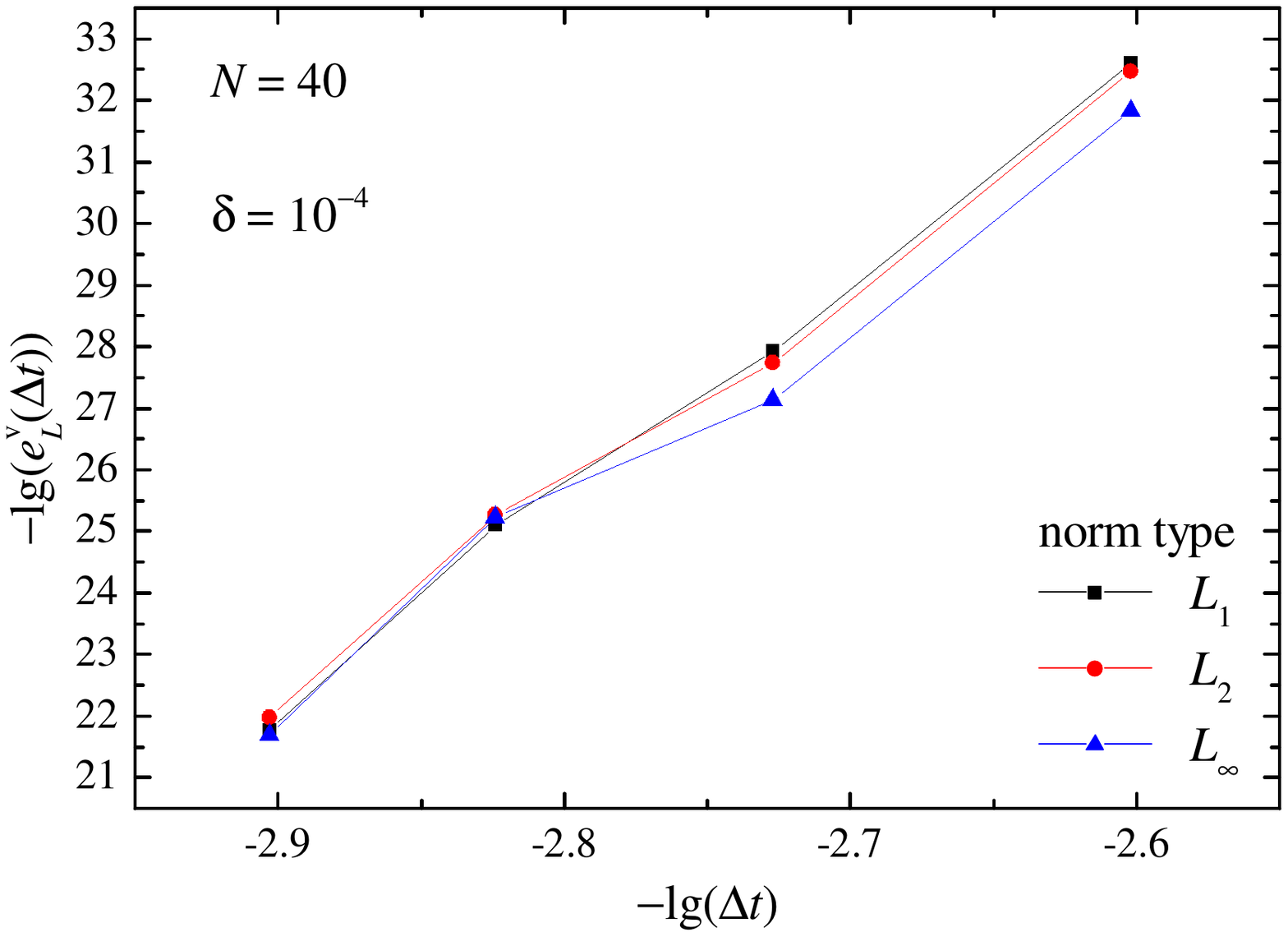}
\vspace{-8mm}\caption{\label{fig:shampine_flame_delta_10m4_errors:b3}}
\end{subfigure}\\
\begin{subfigure}{0.320\textwidth}
\includegraphics[width=\textwidth]{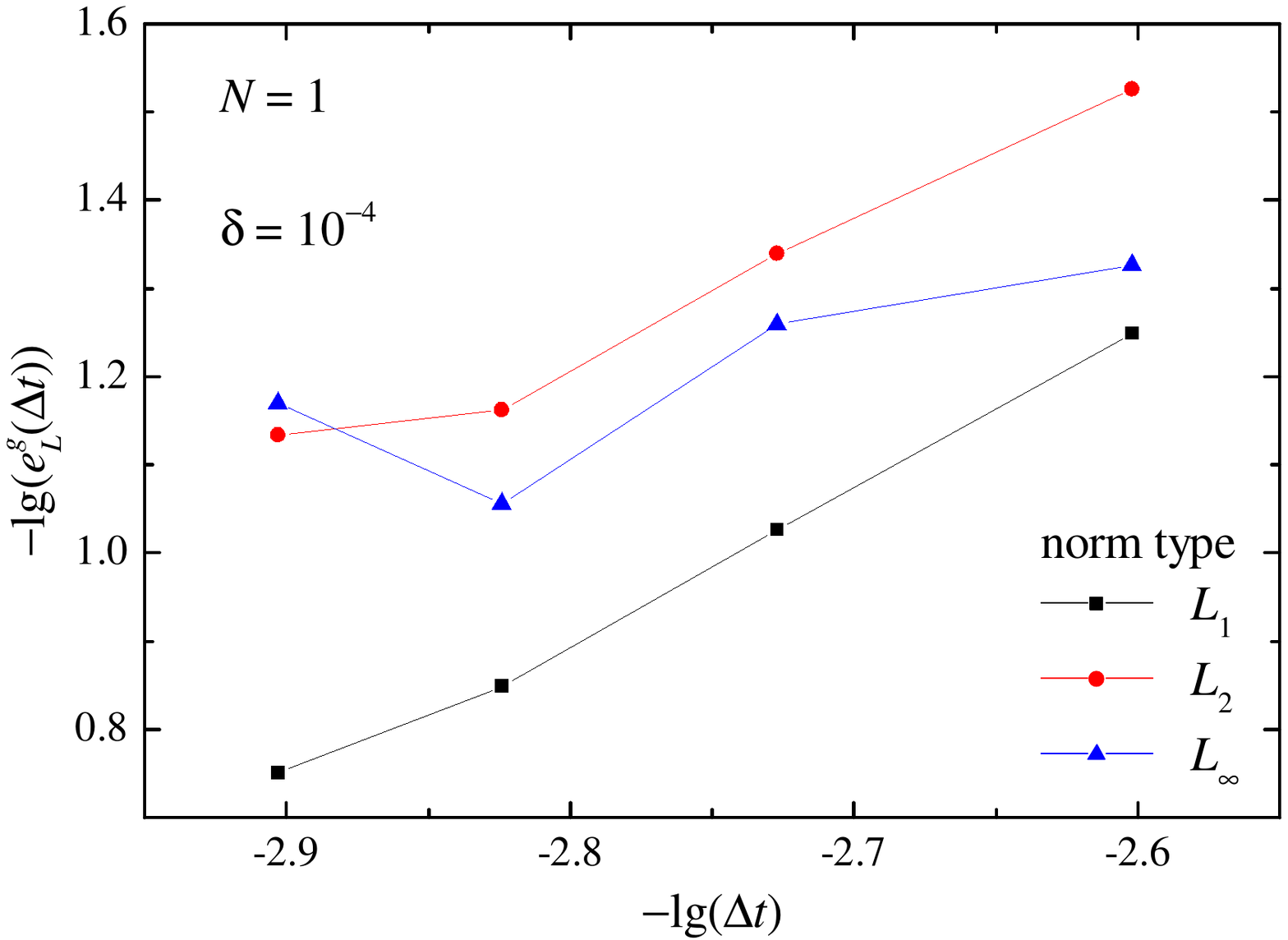}
\vspace{-8mm}\caption{\label{fig:shampine_flame_delta_10m4_errors:c1}}
\end{subfigure}
\begin{subfigure}{0.320\textwidth}
\includegraphics[width=\textwidth]{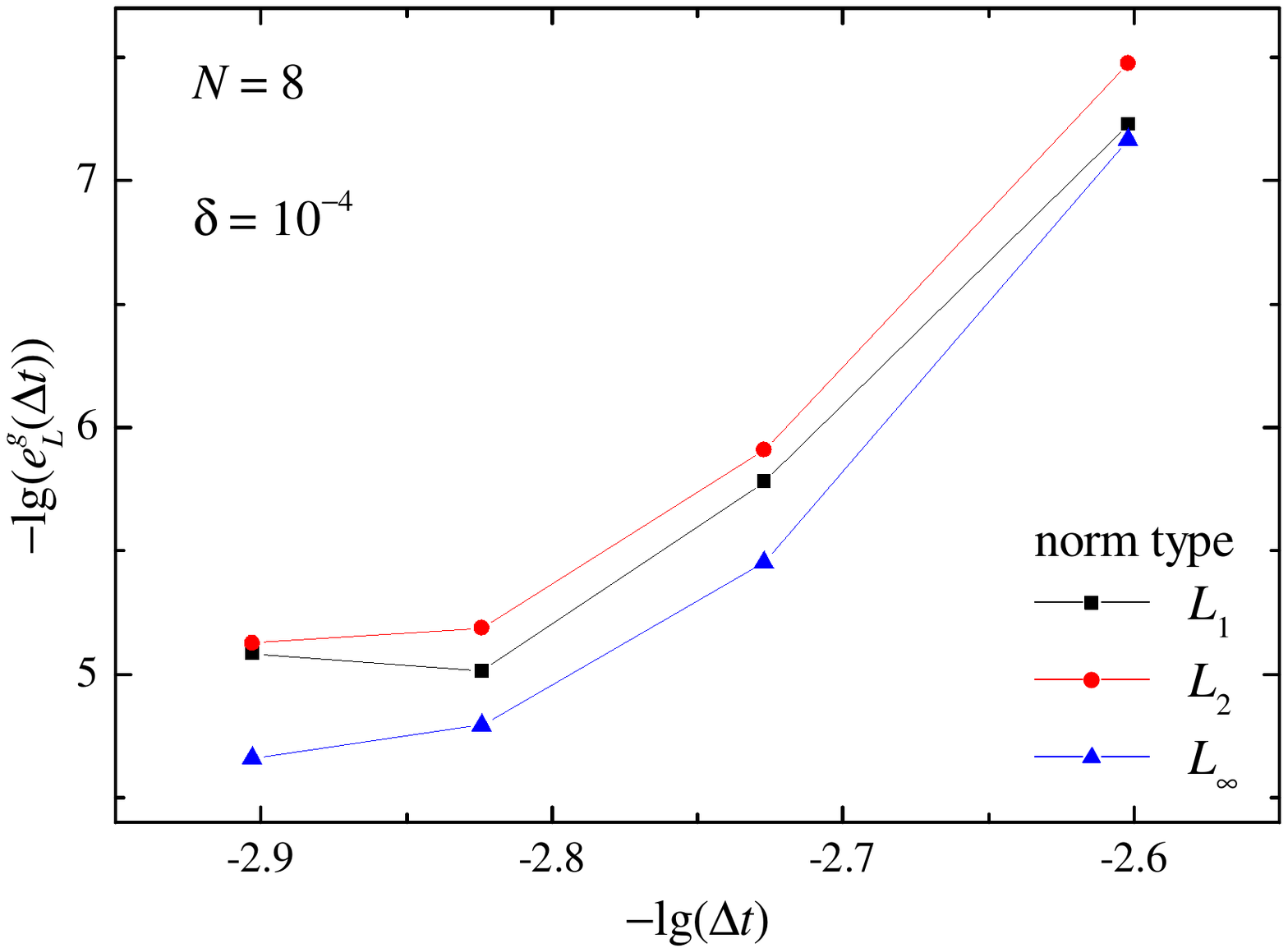}
\vspace{-8mm}\caption{\label{fig:shampine_flame_delta_10m4_errors:c2}}
\end{subfigure}
\begin{subfigure}{0.320\textwidth}
\includegraphics[width=\textwidth]{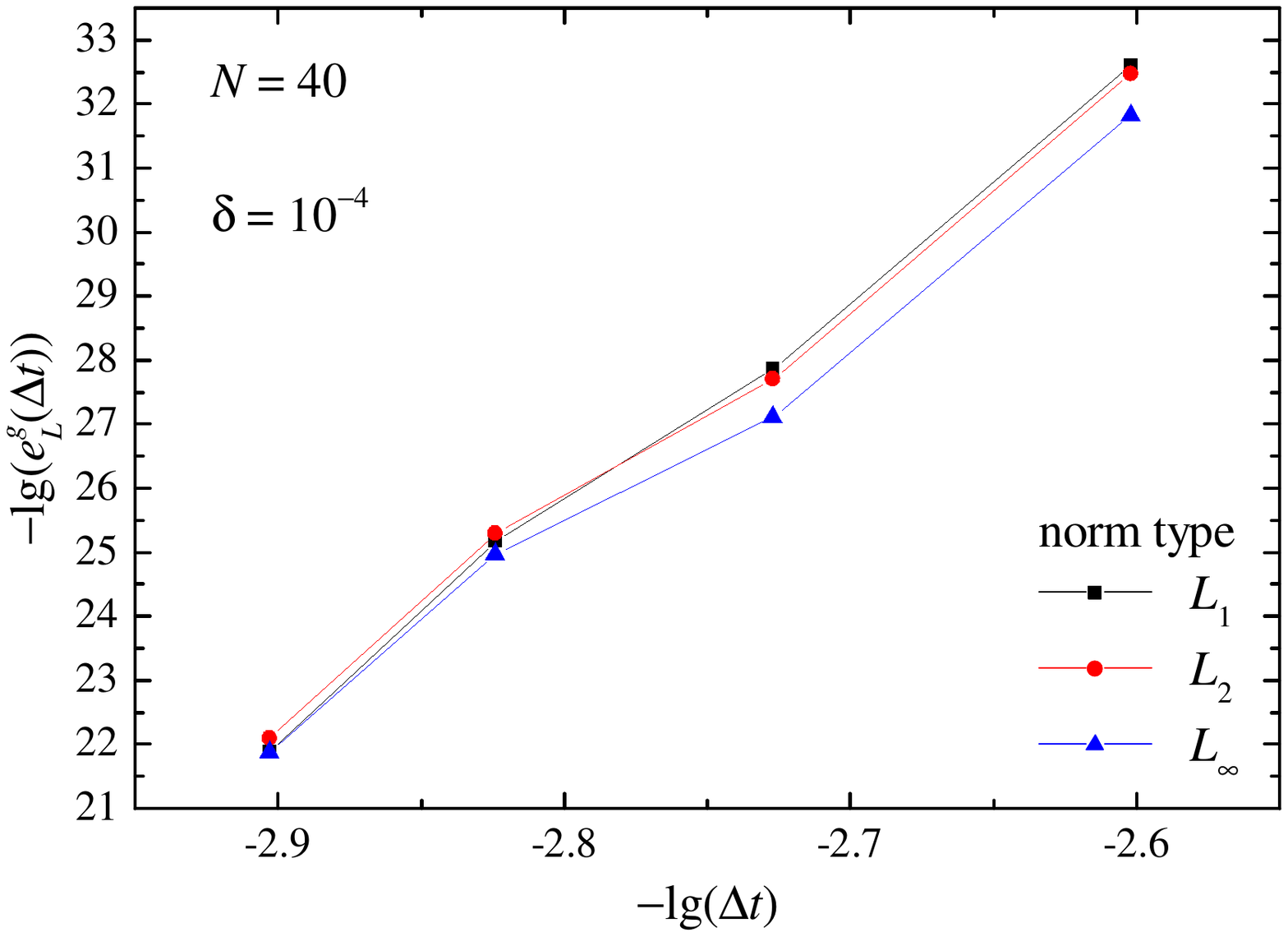}
\vspace{-8mm}\caption{\label{fig:shampine_flame_delta_10m4_errors:c3}}
\end{subfigure}\\
\begin{subfigure}{0.320\textwidth}
\includegraphics[width=\textwidth]{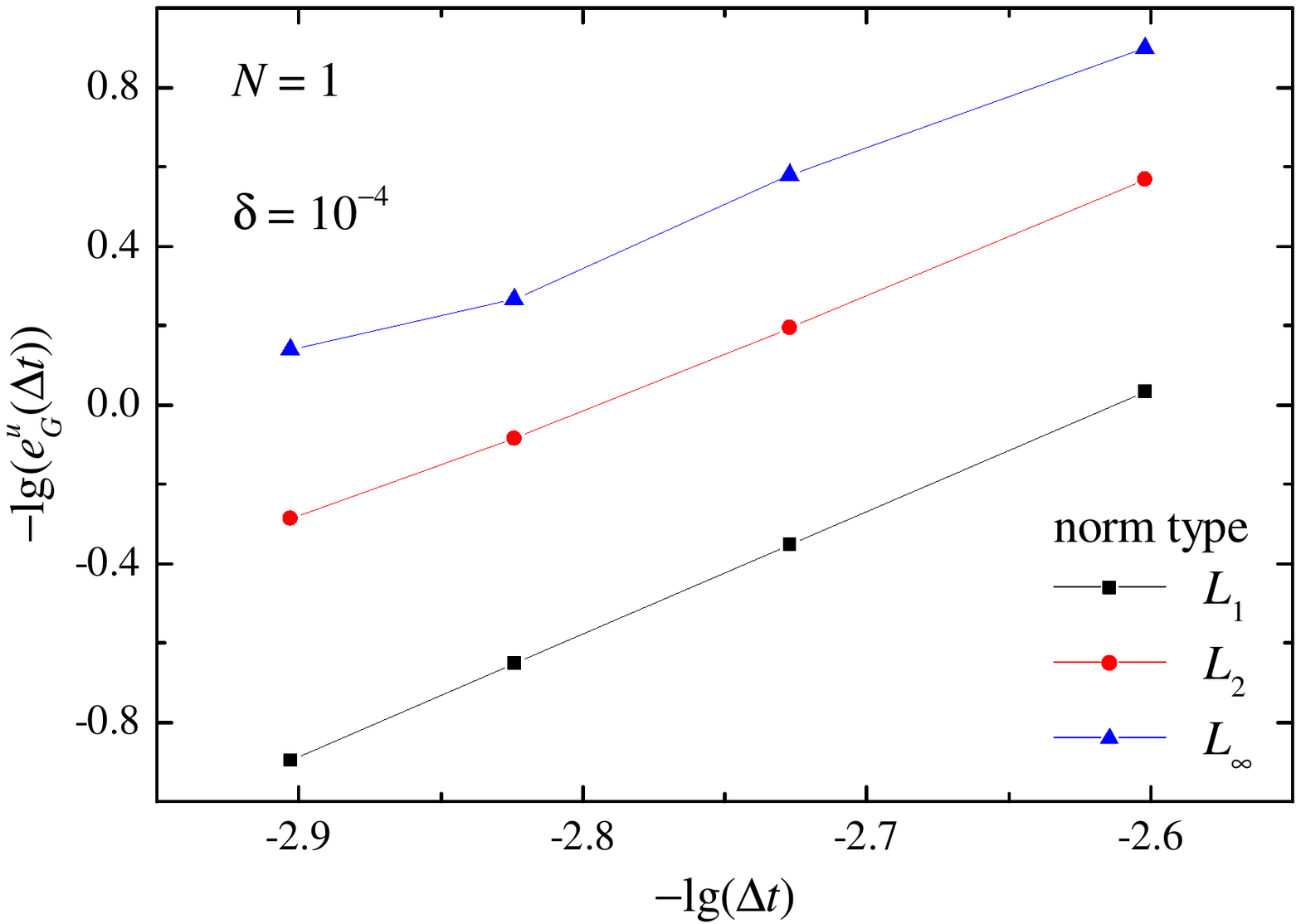}
\vspace{-8mm}\caption{\label{fig:shampine_flame_delta_10m4_errors:d1}}
\end{subfigure}
\begin{subfigure}{0.320\textwidth}
\includegraphics[width=\textwidth]{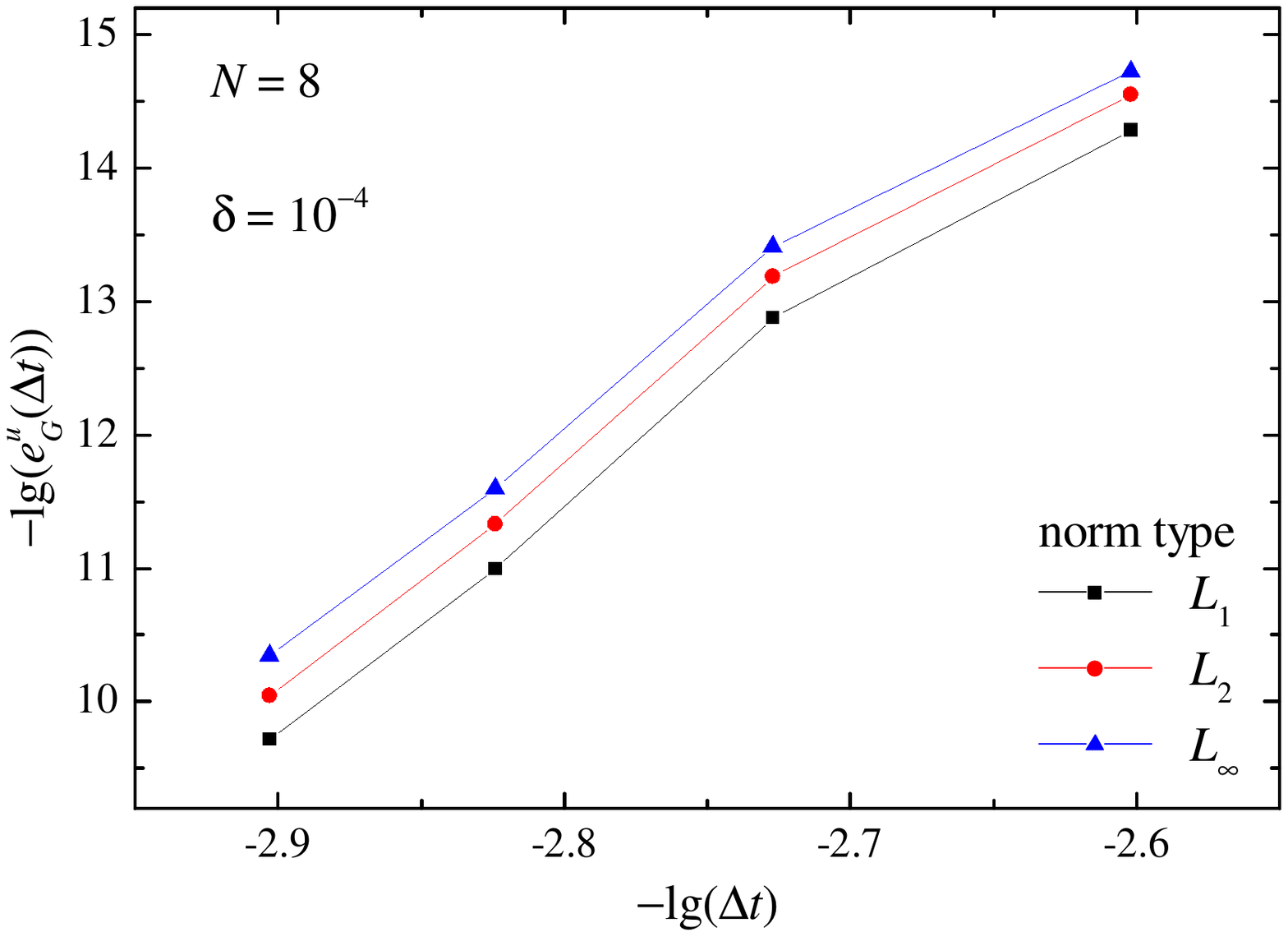}
\vspace{-8mm}\caption{\label{fig:shampine_flame_delta_10m4_errors:d2}}
\end{subfigure}
\begin{subfigure}{0.320\textwidth}
\includegraphics[width=\textwidth]{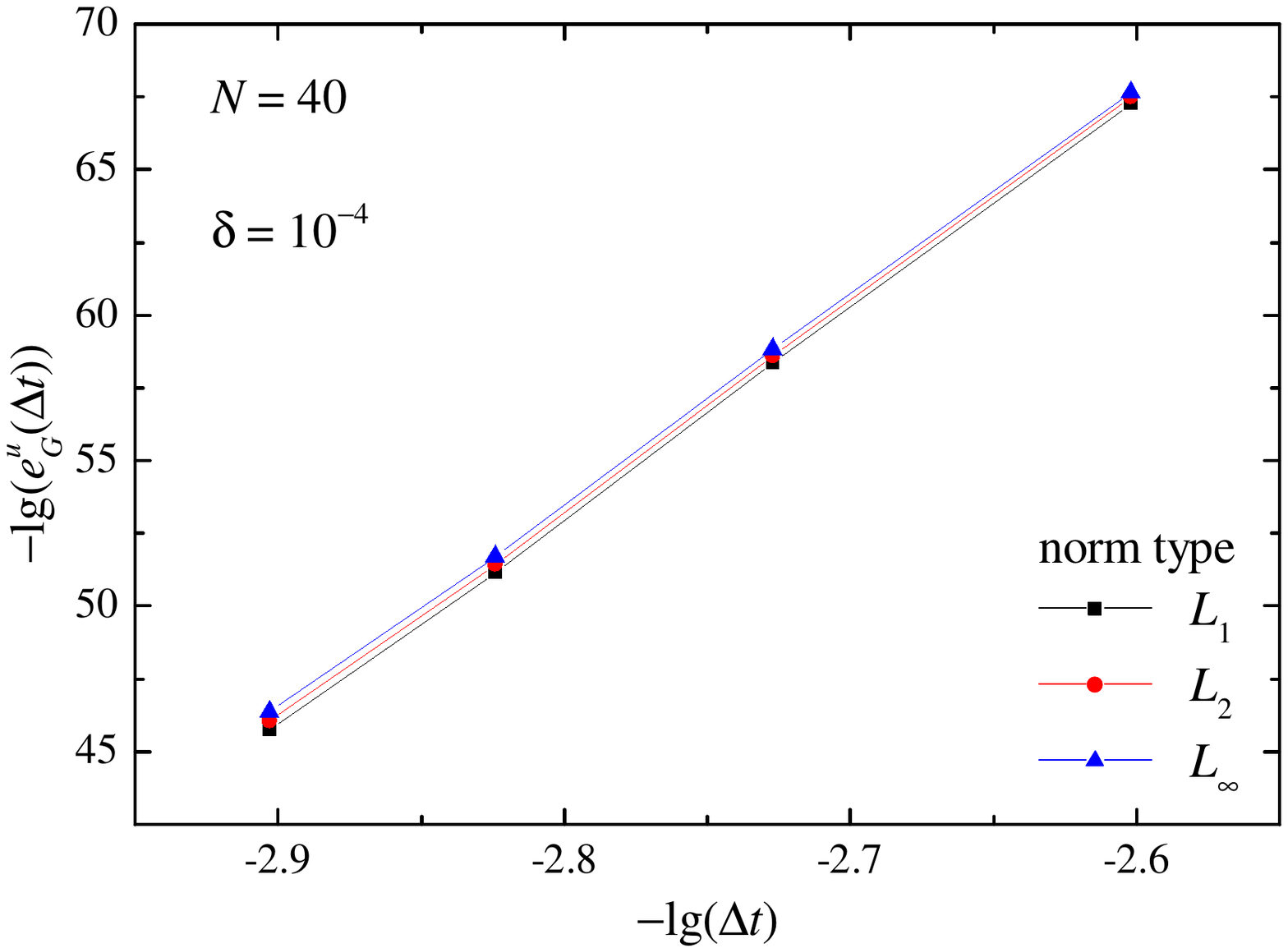}
\vspace{-8mm}\caption{\label{fig:shampine_flame_delta_10m4_errors:d3}}
\end{subfigure}\\
\begin{subfigure}{0.320\textwidth}
\includegraphics[width=\textwidth]{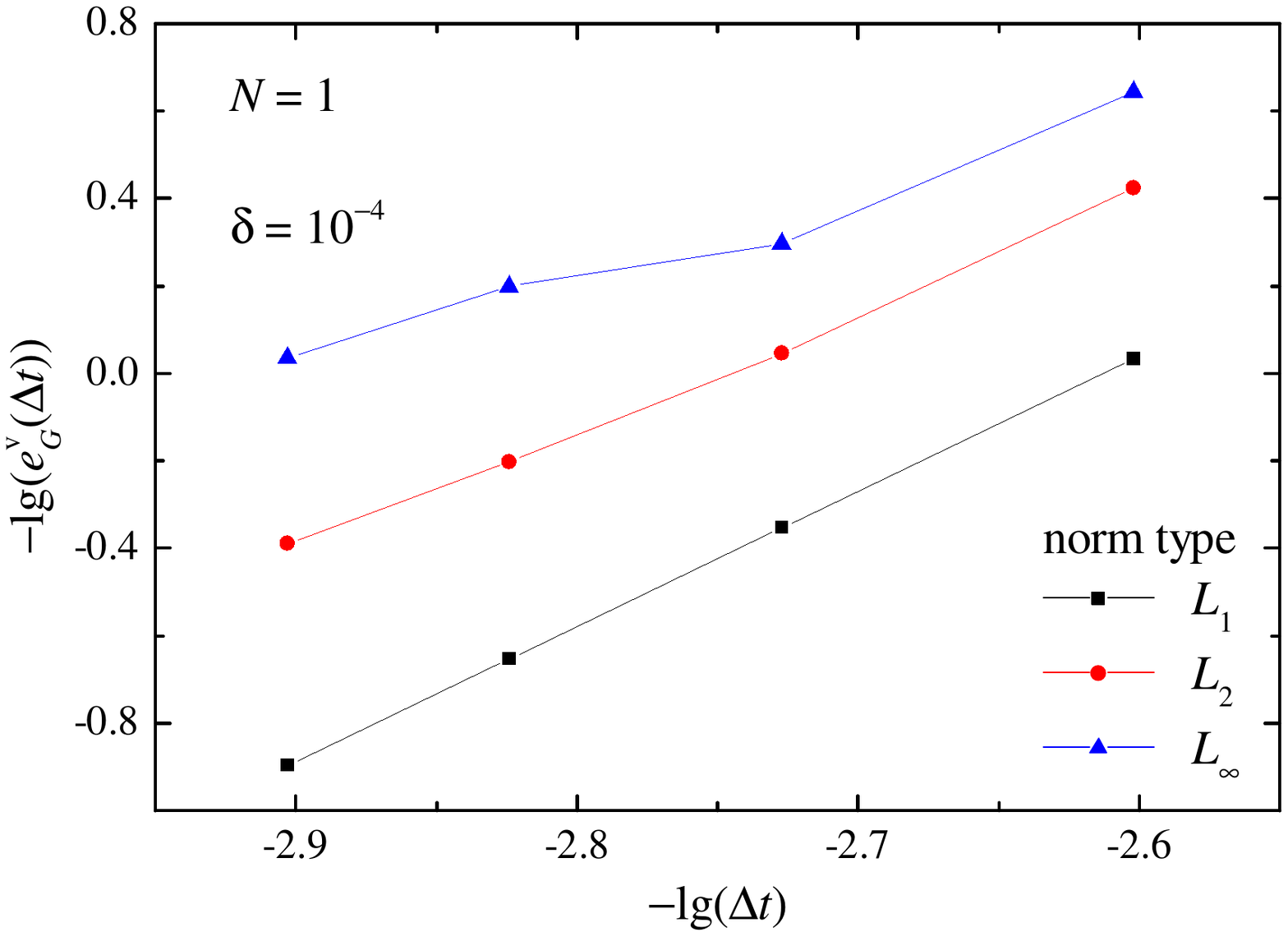}
\vspace{-8mm}\caption{\label{fig:shampine_flame_delta_10m4_errors:e1}}
\end{subfigure}
\begin{subfigure}{0.320\textwidth}
\includegraphics[width=\textwidth]{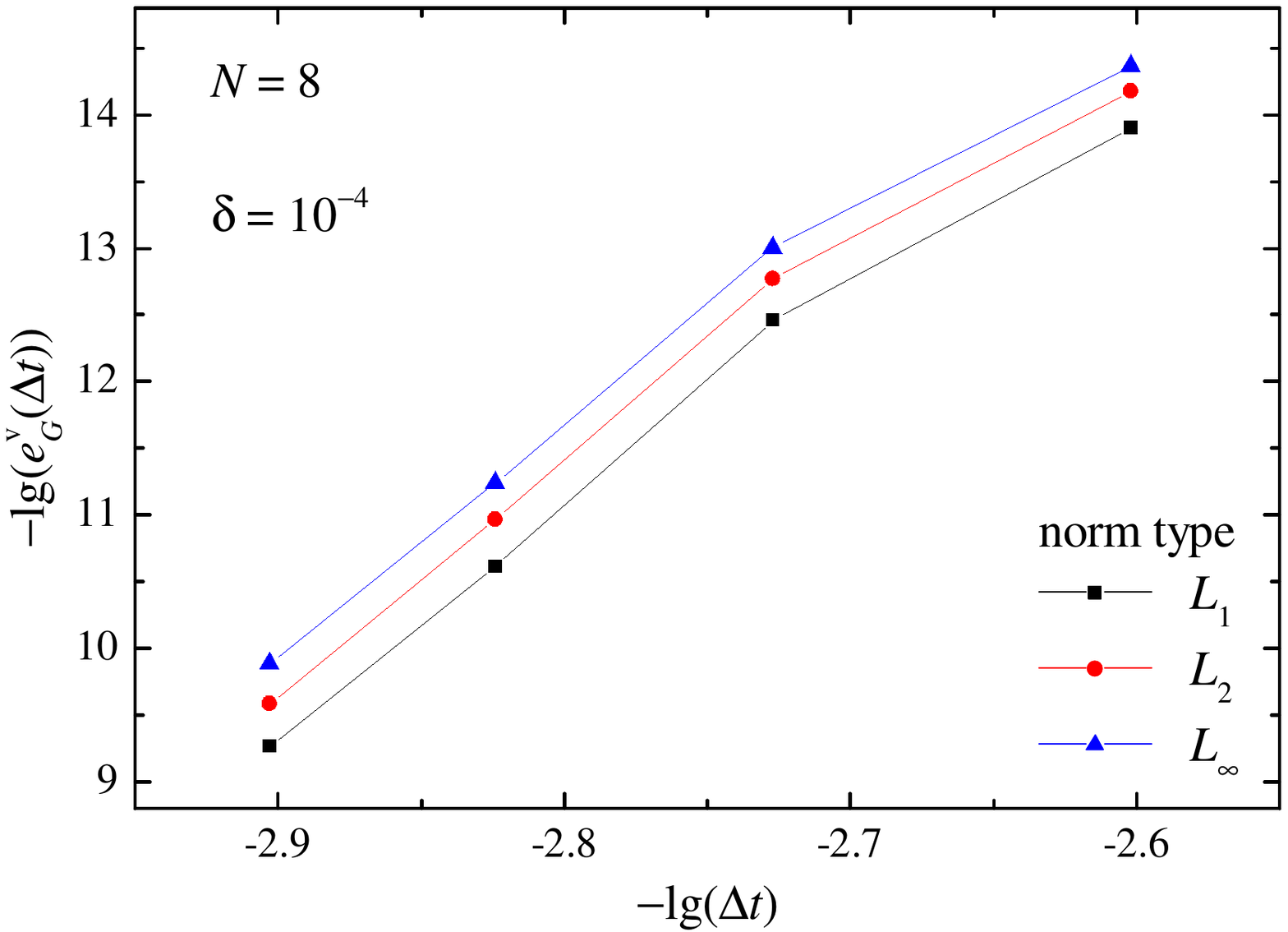}
\vspace{-8mm}\caption{\label{fig:shampine_flame_delta_10m4_errors:e2}}
\end{subfigure}
\begin{subfigure}{0.320\textwidth}
\includegraphics[width=\textwidth]{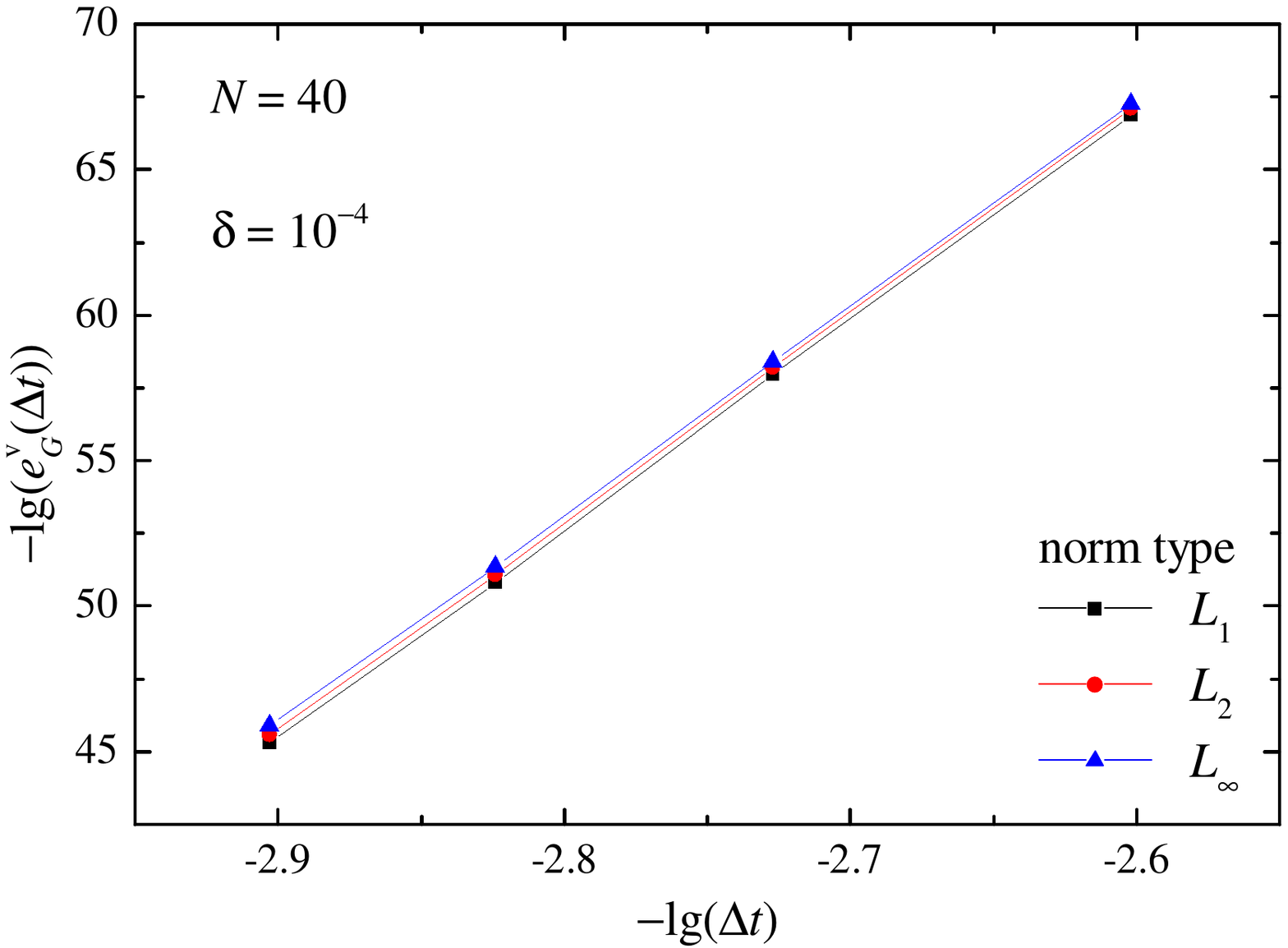}
\vspace{-8mm}\caption{\label{fig:shampine_flame_delta_10m4_errors:e3}}
\end{subfigure}\\
\caption{%
Log-log plot of the dependence of the global errors for the local solution $e_{L}^{u}$ (\subref{fig:shampine_flame_delta_10m4_errors:a1}, \subref{fig:shampine_flame_delta_10m4_errors:a2}, \subref{fig:shampine_flame_delta_10m4_errors:a3}), $e_{L}^{v}$ (\subref{fig:shampine_flame_delta_10m4_errors:b1}, \subref{fig:shampine_flame_delta_10m4_errors:b2}, \subref{fig:shampine_flame_delta_10m4_errors:b3}), $e_{L}^{g}$ (\subref{fig:shampine_flame_delta_10m4_errors:c1}, \subref{fig:shampine_flame_delta_10m4_errors:c2}, \subref{fig:shampine_flame_delta_10m4_errors:c3}) and the solution at nodes $e_{G}^{u}$ (\subref{fig:shampine_flame_delta_10m4_errors:d1}, \subref{fig:shampine_flame_delta_10m4_errors:d2}, \subref{fig:shampine_flame_delta_10m4_errors:d3}), $e_{G}^{v}$ (\subref{fig:shampine_flame_delta_10m4_errors:e1}, \subref{fig:shampine_flame_delta_10m4_errors:e2}, \subref{fig:shampine_flame_delta_10m4_errors:e3}) on the discretization step $\mathrm{\Delta}t$, obtained in the norms $L_{1}$, $L_{2}$ and $L_{\infty}$, by numerical solution of the stiff DAE system (\ref{eq:shampine_flame}) of index 1 with $\delta = 10^{-4}$, obtained using polynomials with degrees $N = 1$ (\subref{fig:shampine_flame_delta_10m4_errors:a1}, \subref{fig:shampine_flame_delta_10m4_errors:b1}, \subref{fig:shampine_flame_delta_10m4_errors:c1}, \subref{fig:shampine_flame_delta_10m4_errors:d1}, \subref{fig:shampine_flame_delta_10m4_errors:e1}), $N = 8$ (\subref{fig:shampine_flame_delta_10m4_errors:a2}, \subref{fig:shampine_flame_delta_10m4_errors:b2}, \subref{fig:shampine_flame_delta_10m4_errors:c2}, \subref{fig:shampine_flame_delta_10m4_errors:d2}, \subref{fig:shampine_flame_delta_10m4_errors:e2}) and $N = 40$ (\subref{fig:shampine_flame_delta_10m4_errors:a3}, \subref{fig:shampine_flame_delta_10m4_errors:b3}, \subref{fig:shampine_flame_delta_10m4_errors:c3}, \subref{fig:shampine_flame_delta_10m4_errors:d3}, \subref{fig:shampine_flame_delta_10m4_errors:e3}).
}
\label{fig:shampine_flame_delta_10m4_errors}
\end{figure}

The obtained results of the numerical solution of DAE system (\ref{eq:shampine_flame}) of index 1 with $\delta = 10^{-4}$ are presented in Figs.~\ref{fig:shampine_flame_delta_10m4_sol_qug}, \ref{fig:shampine_flame_delta_10m4_sol_v_epss}, \ref{fig:shampine_flame_delta_10m4_errors} and in Table~\ref{tab:conv_ords_shampine_flame_delta_10m4}. Fig.~\ref{fig:shampine_flame_delta_10m4_sol_qug} shows a comparison of the solution at nodes $\mathbf{u}_{n}$, the local solution $\mathbf{u}_{L}(t)$ and the exact solution $\mathbf{u}^{\rm ex}(t)$ for component $u_{1}$, quantitative satisfiability of the conditions $g_{1} = 0$, obtained using polynomials with degrees $N = 1$, $8$ and $40$. Fig.~\ref{fig:shampine_flame_delta_10m4_sol_v_epss} shows a comparison of the solution at nodes $\mathbf{v}_{n}$, the local solution $\mathbf{v}_{L}(t)$ and the exact solution $\mathbf{v}^{\rm ex}(t)$ for component $v_{1}$, the errors $\varepsilon_{u}(t)$, $\varepsilon_{v}(t)$, $\varepsilon_{g}(t)$, obtained using polynomials with degrees $N = 1$, $8$ and $40$. In the vicinity of point $t = 1.0/\delta$, a sharp change in the solution $(\mathbf{u},\, \mathbf{v})$ occurs, which leads to high stiffness of the DAE system (\ref{eq:shampine_flame}) with $\delta = 10^{-4}$, therefore the numerical solution at nodes $(\mathbf{u}_{n},\, \mathbf{v}_{n})$, the local solution $(\mathbf{u}_{L}(t),\, \mathbf{v}_{L}(t))$ and its comparison with the exact analytical solution $(\mathbf{u}^{\rm ex}(t),\, \mathbf{v}^{\rm ex}(t))$ (\ref{eq:shampine_flame:exact_solution}) in a small vicinity of point $t = 1.0/\delta$, specifically in the range $t \in [0.995/\delta,\, 1.005/\delta]$, was presented separately in Fig.~\ref{fig:shampine_flame_delta_10m4_sol_qug} (\subref{fig:shampine_flame_delta_10m4_sol_qug:d1}, \subref{fig:shampine_flame_delta_10m4_sol_qug:d2}, \subref{fig:shampine_flame_delta_10m4_sol_qug:d3}, \subref{fig:shampine_flame_delta_10m4_sol_qug:b1}, \subref{fig:shampine_flame_delta_10m4_sol_qug:b2}, \subref{fig:shampine_flame_delta_10m4_sol_qug:b3}) for differential variables $\mathbf{u}$ and in Fig.~\ref{fig:shampine_flame_delta_10m4_sol_v_epss} (\subref{fig:shampine_flame_delta_10m4_sol_v_epss:b1}, \subref{fig:shampine_flame_delta_10m4_sol_v_epss:b2}, \subref{fig:shampine_flame_delta_10m4_sol_v_epss:b3}) for algebraic variables $\mathbf{v}$. The comparison of the numerical solution and the exact analytical solution in this zoomed range $t \in [0.995/\delta,\, 1.005/\delta]$ allowed one to identify artifacts of the numerical solution obtained by the ADER-DG method with a local DG predictor more accurately and correctly, in the case of this stiff DAE system (\ref{eq:shampine_flame}) of index 1 with $\delta = 10^{-4}$. Fig.~\ref{fig:pend_ind2_errors} shows the log-log plot of the dependence of the global errors for the local solution $e_{L}^{u}$, $e_{L}^{v}$, $e_{L}^{g}$ and the solution at nodes $e_{G}^{u}$, $e_{G}^{v}$ on the discretization step $\mathrm{\Delta}t$, obtained in the norms $L_{1}$, $L_{2}$ and $L_{\infty}$, by numerical solution of the stiff DAE system (\ref{eq:shampine_flame}) of index 1 with $\delta = 10^{-4}$, obtained using polynomials with degrees $N = 1$, $8$ and $40$, which was used to calculate the empirical convergence orders $p$.

The dynamic dependencies of the numerical solution presented in Fig.~\ref{fig:shampine_flame_delta_10m4_sol_qug} and Fig.~\ref{fig:shampine_flame_delta_10m4_sol_v_epss} demonstrate a very high accuracy of the ADER-DG numerical method with a local DG predictor for solving the stiff problem (\ref{eq:shampine_flame}) with $\delta = 10^{-4}$. A point-by-point comparison of the dynamic dependencies of the solution at the nodes $\mathbf{u}_{n}$ and the local solution $\mathbf{u}_{L}(t)$ with the exact analytical solution $\mathbf{u}^{\rm ex}$ for differential variables $\mathbf{u}$ shows that in the case of using polynomials of degree $N = 1$, an expected artifact of the numerical solution is observed --- the range of a sharp change in the solution occurs in the region of larger values of the argument $t$ than in the exact analytical solution, and the shift occurs approximately by $\delta\cdot t \approx 0.002$ to the right. This effect has already been observed in the work~\cite{ader_dg_ivp_ode} when using the ADER-DG numerical method with a local DG predictor when solving the initial value problem for an ODE system that is completely equivalent to a DAE system (\ref{eq:shampine_flame}) of index 1 with $\delta = 10^{-4}$. Note that a similar effect also occurs when using the Runge-Kutta numerical methods of the Radau family, as well as the \texttt{BDF} and \texttt{LSODA} methods available in the \texttt{scipy} module of the \texttt{python} programming language (which was also noted in the work~\cite{ader_dg_ivp_ode}), while in the case of the \texttt{BDF} and \texttt{LSODA} methods this effect is much stronger. The point-wise comparison of the dynamic dependencies of the solution at the nodes $\mathbf{u}_{n}$ and the local solution $\mathbf{u}_{L}(t)$ with the exact analytical solution $\mathbf{u}^{\rm ex}$, in the case of using polynomials of degrees $N = 8$ and $40$, shows that this artifact of the numerical solution does not occur --- the numerical solution very accurately coincides point-wise with the exact analytical solution. It should be noted that the domain of a sharp change in the solution in the vicinity of the point $t \approx 1/\delta$ accounts for only $4$--$5$ grid nodes, while the comparison not only of the solution at the nodes $\mathbf{u}_{n}$, but also of the local solution $\mathbf{u}_{L}(t)$ available between the grid nodes $t_{n} \leqslant t \leqslant t_{n+1}$ in the discretization domains $\Omega_{n}$ shows high accuracy. The point-wise comparison of the dynamic dependencies of the solution at the nodes $\mathbf{v}_{n}$ and the local solution $\mathbf{v}_{L}(t)$ with the exact analytical solution $\mathbf{v}^{\rm ex}$ for the algebraic variables $\mathbf{v}$ shows similar behavior both in the region of smooth change of the solution and in the neighborhood of the point $t = 1/\delta$, which is shown in Fig.~\ref{fig:shampine_flame_delta_10m4_sol_v_epss} (\subref{fig:shampine_flame_delta_10m4_sol_v_epss:a1}, \subref{fig:shampine_flame_delta_10m4_sol_v_epss:a2}, \subref{fig:shampine_flame_delta_10m4_sol_v_epss:a3}, \subref{fig:shampine_flame_delta_10m4_sol_v_epss:b1}, \subref{fig:shampine_flame_delta_10m4_sol_v_epss:b2}, \subref{fig:shampine_flame_delta_10m4_sol_v_epss:b3}). 

Fig.~\ref{fig:shampine_flame_delta_10m4_sol_qug} (\subref{fig:shampine_flame_delta_10m4_sol_qug:e1}, \subref{fig:shampine_flame_delta_10m4_sol_qug:e2}, \subref{fig:shampine_flame_delta_10m4_sol_qug:e3}) shows the dynamic dependence of the error $|g_{1}|$ of satisfaction of the algebraic equation $g_{1} = 0$ (\ref{eq:shampine_flame}). The algebraic equation $g_{1} = 0$ is satisfied exactly (within the accuracy of representation of real numbers with a floating point) for the solution at the nodes $(\mathbf{u}_{n},\, \mathbf{v}_{n})$, therefore the error $|g_{1}|$ for the solution at the nodes is not shown in Fig.~\ref{fig:shampine_flame_delta_10m4_sol_qug}. The error $|g_{1}|$ of satisfaction of the algebraic equation $g_{1} = 0$ for the local solution $(\mathbf{u}_{L}(t),\, \mathbf{v}_{L}(t))$ in the region to the right of the region $t \approx 1/\delta$ of sharp change of the solution is also zero (within the accuracy of representation of real numbers with a floating point), but to the left of the region  $t \approx 1/\delta$ of sharp change of the solution this error is non-zero and is clearly visible in Fig.~\ref{fig:shampine_flame_delta_10m4_sol_qug}. In the case of using a polynomial of degree $N = 1$, the characteristic values of the error are $|g_{1}| \sim 10^{-18}$--$10^{-1}$, in the case of $N = 8$ --- $|g_{1}| \sim 10^{-39}$--$10^{-8}$, in the case of $N = 40$ --- $|g_{1}| \sim 10^{-150}$--$10^{-20}$. From the presented results it is evident that the error $|g_{1}|$ of satisfaction of the algebraic equation $g_{1} = 0$ for the local solution $(\mathbf{u}_{L}(t),\, \mathbf{v}_{L}(t))$ undergoes very significant changes, which is connected precisely with the high stiffness of the problem (\ref{eq:shampine_flame}) with $\delta = 10^{-4}$ --- the greatest error $|g_{1}|$ occurs in the vicinity of point $t \approx 1/\delta$, while to the left of this vicinity there is a minimum of the error $|g_{1}|$, in which the minimum value of the error $|g_{1}|$ for the entire domain of definition $\Omega$ is achieved.

\begin{table*}[h!]
\centering
\caption{%
Convergence orders $p_{L_{1}}$, $p_{L_{2}}$, $p_{L_{\infty}}$, calculated in norms $L_{1}$, $L_{2}$, $L_{\infty}$ of the ADER-DG method for the stiff DAE system (\ref{eq:shampine_flame}) of index 1 with $\delta = 10^{-4}$; $N$ is the degree of the basis polynomials $\varphi_{p}$. Orders $p^{n, u}$ are calculated for solution at the nodes $\mathbf{u}_{n}$; orders $p^{n, v}$ --- for solution at the nodes $\mathbf{v}_{n}$; orders $p^{l, u}$ --- for local solution $\mathbf{u}_{L}$; orders $p^{l, v}$ --- for local solution $\mathbf{v}_{L}$. The theoretical values of convergence order $p_{\rm th.}^{n} = 2N+1$ and $p_{\rm th.}^{l} = N+1$ are applicable for the ADER-DG method for ODE problems, and are presented for comparison.
}
\label{tab:conv_ords_shampine_flame_delta_10m4}
\begingroup
\setlength{\tabcolsep}{4pt}
\begin{adjustbox}{angle=90}
\begin{tabular}{@{}|l|lll|lll|c|lll|lll|c|@{}}
\toprule
$N$\hspace{-2mm} & $p_{L_{1}}^{n, u}$ & $p_{L_{2}}^{n, u}$ & $p_{L_{\infty}}^{n, u}$ & $p_{L_{1}}^{n, v}$ & $p_{L_{2}}^{n, v}$ & $p_{L_{\infty}}^{n, v}$ & $p_{\rm th.}^{n}$ & $p_{L_{1}}^{l, u}$ & $p_{L_{2}}^{l, u}$ & $p_{L_{\infty}}^{l, u}$ & $p_{L_{1}}^{l, v}$ & $p_{L_{2}}^{l, v}$ & $p_{L_{\infty}}^{l, v}$ & $p_{\rm th.}^{l}$ \\
\midrule
$1$	&	$3.08$	&	$2.85$	&	$2.61$	&	$3.09$	&	$2.70$	&	$1.95$	&	$3$	&	$3.08$	&	$2.88$	&	$2.64$	&	$3.09$	&	$2.73$	&	$2.16$	&	$2$\\
$2$	&	$4.49$	&	$4.26$	&	$3.70$	&	$4.73$	&	$4.56$	&	$4.07$	&	$5$	&	$2.47$	&	$2.24$	&	$1.26$	&	$2.62$	&	$2.35$	&	$1.73$	&	$3$\\
$3$	&	$8.55$	&	$8.59$	&	$8.21$	&	$8.50$	&	$8.53$	&	$8.13$	&	$7$	&	$3.89$	&	$3.66$	&	$3.27$	&	$4.66$	&	$4.58$	&	$4.12$	&	$4$\\
$4$	&	$6.68$	&	$6.10$	&	$5.56$	&	$6.82$	&	$6.29$	&	$5.79$	&	$9$	&	$4.68$	&	$4.51$	&	$3.95$	&	$3.29$	&	$3.43$	&	$3.51$	&	$5$\\
$5$	&	$12.25$	&	$12.27$	&	$11.97$	&	$12.36$	&	$12.30$	&	$11.89$	&	$11$	&	$6.02$	&	$5.74$	&	$5.44$	&	$6.44$	&	$6.22$	&	$5.57$	&	$6$\\
$6$	&	$12.21$	&	$11.77$	&	$11.22$	&	$12.49$	&	$12.04$	&	$11.45$	&	$13$	&	$6.59$	&	$6.49$	&	$6.07$	&	$5.32$	&	$5.21$	&	$5.24$	&	$7$\\
$7$	&	$14.53$	&	$14.50$	&	$14.23$	&	$14.60$	&	$14.50$	&	$14.15$	&	$15$	&	$7.35$	&	$7.28$	&	$6.92$	&	$7.99$	&	$7.72$	&	$7.05$	&	$8$\\
$8$	&	$15.44$	&	$15.23$	&	$14.79$	&	$15.58$	&	$15.43$	&	$15.03$	&	$17$	&	$8.52$	&	$8.14$	&	$7.56$	&	$7.76$	&	$7.49$	&	$6.19$	&	$9$\\
$9$	&	$17.49$	&	$17.24$	&	$16.76$	&	$17.66$	&	$17.46$	&	$16.99$	&	$19$	&	$9.60$	&	$9.34$	&	$8.06$	&	$8.94$	&	$8.76$	&	$8.05$	&	$10$\\
$10$	&	$19.18$	&	$18.91$	&	$18.44$	&	$19.36$	&	$19.13$	&	$18.67$	&	$21$	&	$9.87$	&	$9.69$	&	$9.17$	&	$10.49$	&	$10.26$	&	$9.76$	&	$11$\\
$11$	&	$19.92$	&	$19.68$	&	$19.21$	&	$20.09$	&	$19.89$	&	$19.44$	&	$23$	&	$11.93$	&	$11.98$	&	$12.09$	&	$10.75$	&	$10.70$	&	$10.57$	&	$12$\\
$12$	&	$20.24$	&	$19.98$	&	$19.51$	&	$20.42$	&	$20.20$	&	$19.75$	&	$25$	&	$12.09$	&	$11.68$	&	$11.36$	&	$12.80$	&	$12.82$	&	$12.75$	&	$13$\\
$13$	&	$22.23$	&	$21.97$	&	$21.50$	&	$22.39$	&	$22.18$	&	$21.73$	&	$27$	&	$13.87$	&	$14.20$	&	$14.26$	&	$11.98$	&	$12.16$	&	$12.17$	&	$14$\\
$14$	&	$24.66$	&	$24.41$	&	$23.94$	&	$24.83$	&	$24.63$	&	$24.18$	&	$29$	&	$13.51$	&	$13.09$	&	$12.76$	&	$15.55$	&	$15.44$	&	$15.95$	&	$15$\\
$15$	&	$25.82$	&	$25.56$	&	$25.09$	&	$25.98$	&	$25.78$	&	$25.33$	&	$31$	&	$16.06$	&	$15.97$	&	$15.81$	&	$13.52$	&	$13.16$	&	$12.87$	&	$16$\\
$16$	&	$27.69$	&	$27.45$	&	$26.98$	&	$27.86$	&	$27.66$	&	$27.21$	&	$33$	&	$14.82$	&	$14.24$	&	$12.83$	&	$16.81$	&	$16.85$	&	$16.50$	&	$17$\\
$17$	&	$30.93$	&	$30.69$	&	$30.23$	&	$31.10$	&	$30.91$	&	$30.46$	&	$35$	&	$17.39$	&	$17.32$	&	$16.97$	&	$16.03$	&	$15.24$	&	$14.14$	&	$18$\\
$18$	&	$32.22$	&	$31.98$	&	$31.51$	&	$32.39$	&	$32.20$	&	$31.75$	&	$37$	&	$17.09$	&	$16.56$	&	$16.09$	&	$18.47$	&	$18.26$	&	$17.75$	&	$19$\\
$19$	&	$34.42$	&	$34.17$	&	$33.70$	&	$34.59$	&	$34.39$	&	$33.94$	&	$39$	&	$19.01$	&	$18.78$	&	$18.36$	&	$18.07$	&	$17.46$	&	$16.90$	&	$20$\\
$20$	&	$35.89$	&	$35.64$	&	$35.17$	&	$36.06$	&	$35.86$	&	$35.41$	&	$41$	&	$19.07$	&	$18.65$	&	$18.19$	&	$19.81$	&	$19.35$	&	$18.84$	&	$21$\\
\midrule
$25$	&	$45.83$	&	$45.59$	&	$45.12$	&	$46.00$	&	$45.81$	&	$45.36$	&	$51$	&	$23.12$	&	$22.39$	&	$21.64$	&	$23.59$	&	$23.60$	&	$23.33$	&	$26$\\
$30$	&	$50.72$	&	$50.48$	&	$50.01$	&	$50.89$	&	$50.69$	&	$50.24$	&	$61$	&	$28.13$	&	$27.25$	&	$26.25$	&	$27.77$	&	$27.30$	&	$26.77$	&	$31$\\
$35$	&	$62.31$	&	$62.06$	&	$61.59$	&	$62.47$	&	$62.28$	&	$61.83$	&	$71$	&	$32.13$	&	$31.51$	&	$30.98$	&	$33.38$	&	$33.54$	&	$34.04$	&	$36$\\
$40$	&	$71.74$	&	$71.50$	&	$71.03$	&	$71.91$	&	$71.71$	&	$71.26$	&	$81$	&	$38.76$	&	$38.71$	&	$38.50$	&	$35.26$	&	$33.94$	&	$32.28$	&	$41$\\
\bottomrule
\end{tabular}
\end{adjustbox}
\endgroup
\end{table*}

\begin{figure}[h!]
\captionsetup[subfigure]{%
	position=bottom,
	font+=smaller,
	textfont=normalfont,
	singlelinecheck=off,
	justification=raggedright
}
\centering
\begin{subfigure}{0.320\textwidth}
\includegraphics[width=\textwidth]{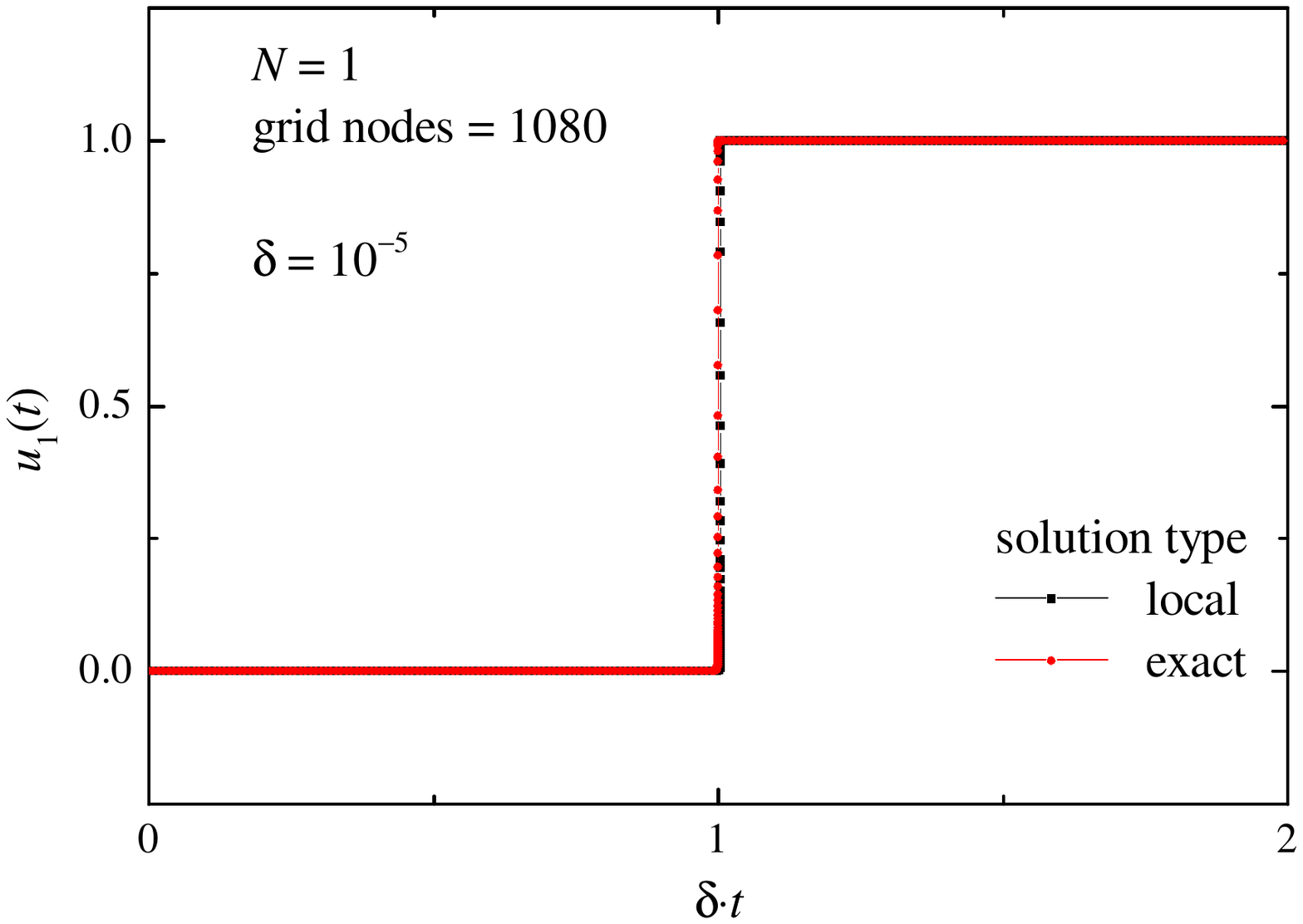}
\vspace{-8mm}\caption{\label{fig:shampine_flame_delta_10m5_sol_qug:a1}}
\end{subfigure}
\begin{subfigure}{0.320\textwidth}
\includegraphics[width=\textwidth]{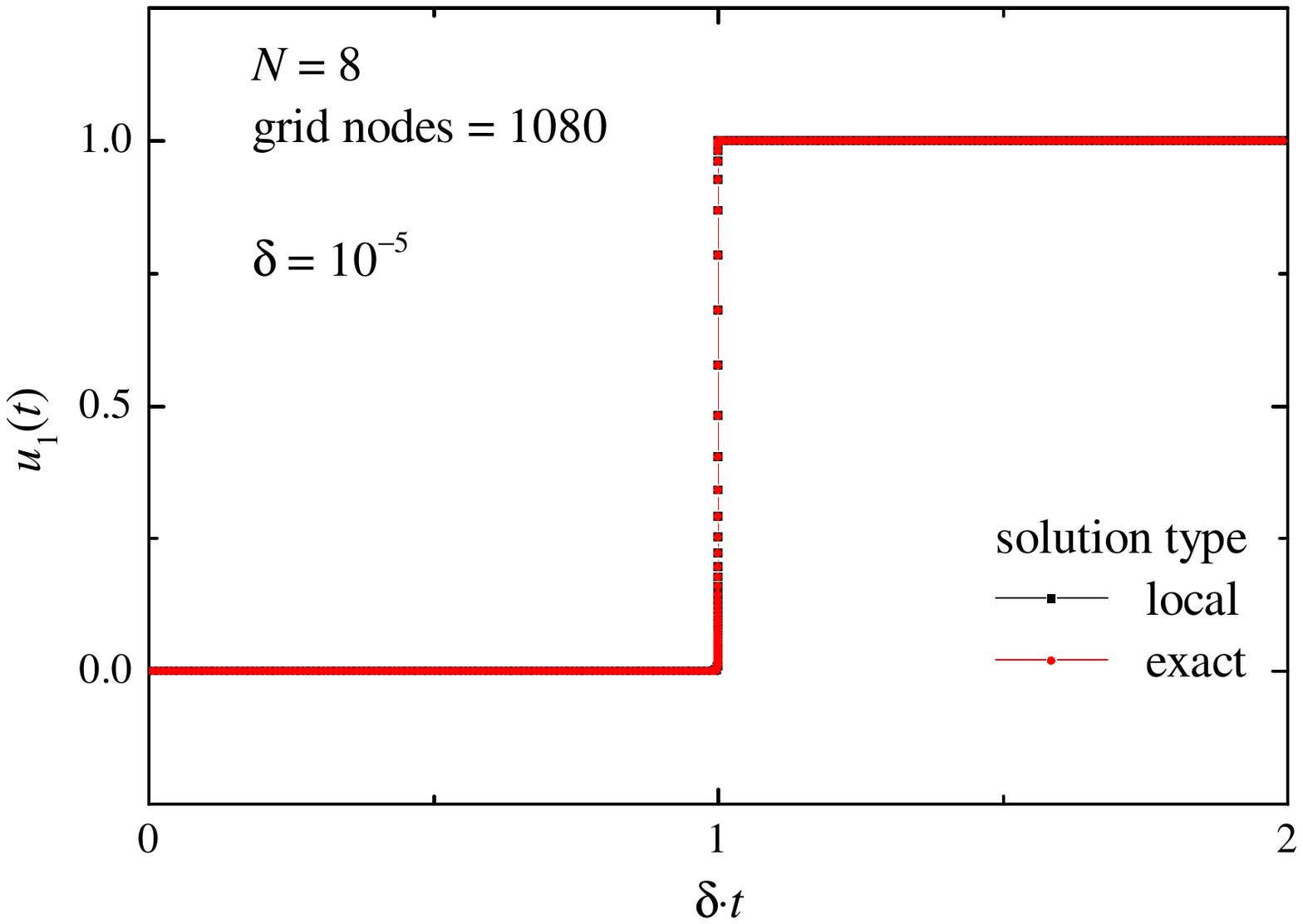}
\vspace{-8mm}\caption{\label{fig:shampine_flame_delta_10m5_sol_qug:a2}}
\end{subfigure}
\begin{subfigure}{0.320\textwidth}
\includegraphics[width=\textwidth]{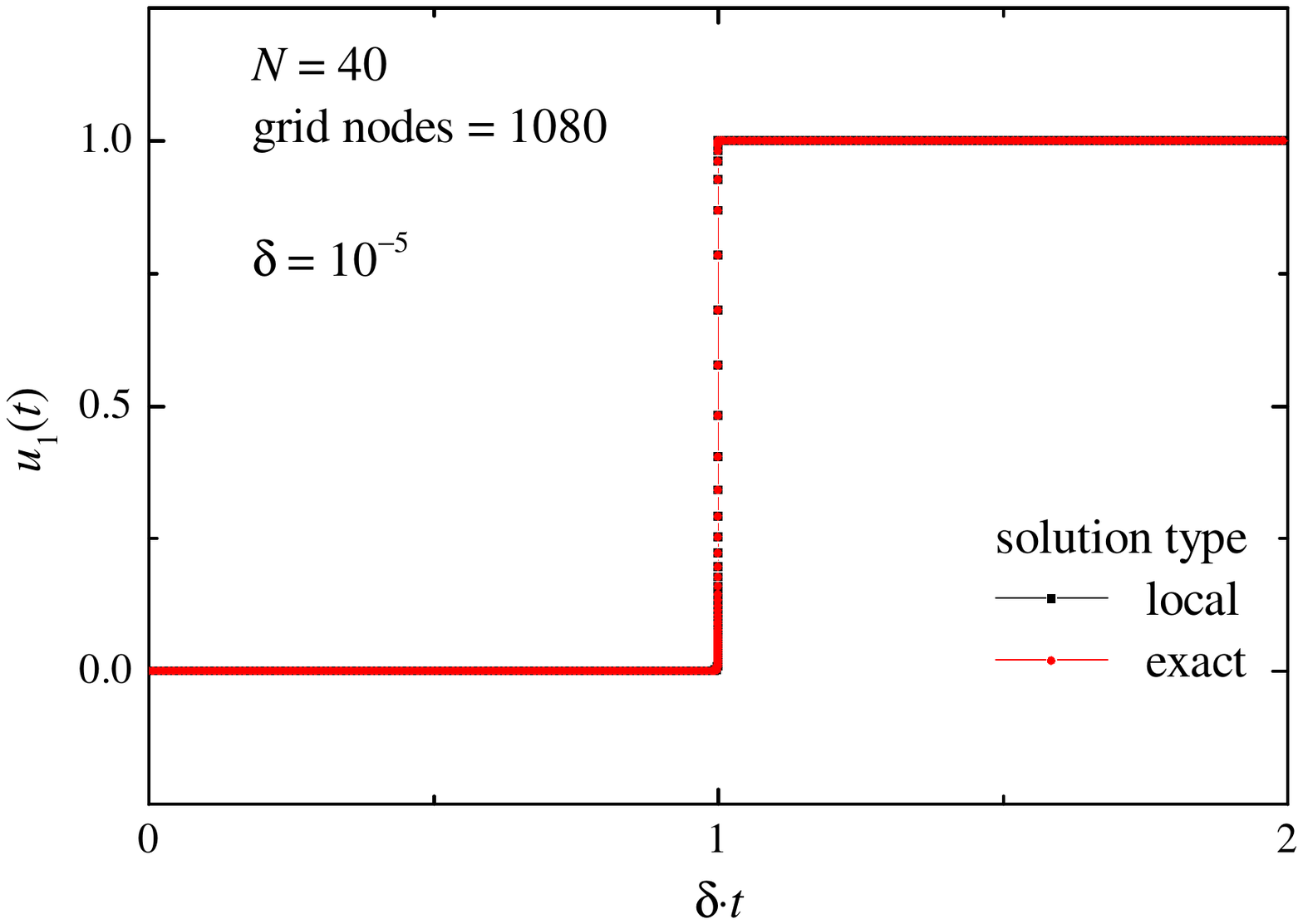}
\vspace{-8mm}\caption{\label{fig:shampine_flame_delta_10m5_sol_qug:a3}}
\end{subfigure}\\
\begin{subfigure}{0.320\textwidth}
\includegraphics[width=\textwidth]{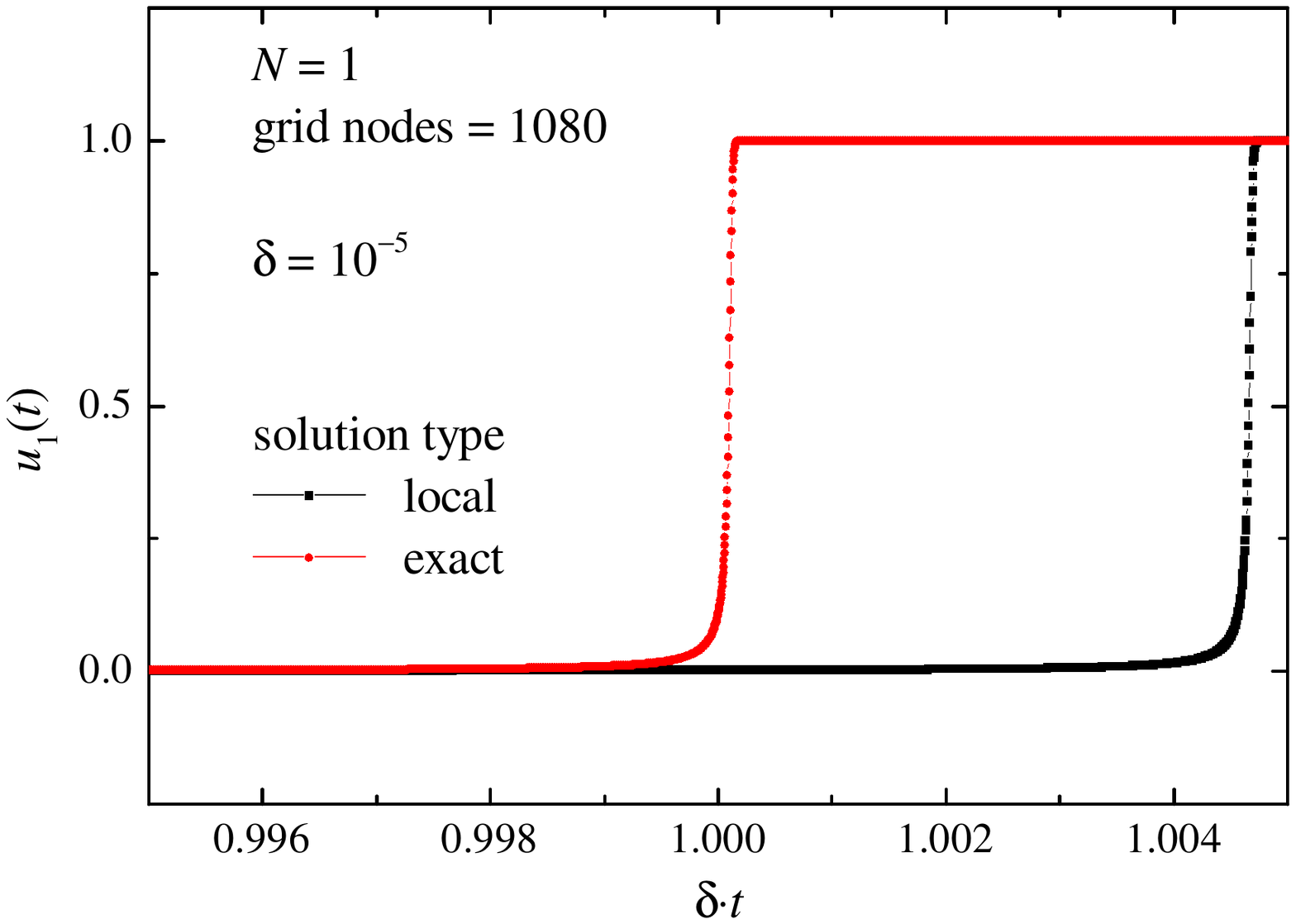}
\vspace{-8mm}\caption{\label{fig:shampine_flame_delta_10m5_sol_qug:b1}}
\end{subfigure}
\begin{subfigure}{0.320\textwidth}
\includegraphics[width=\textwidth]{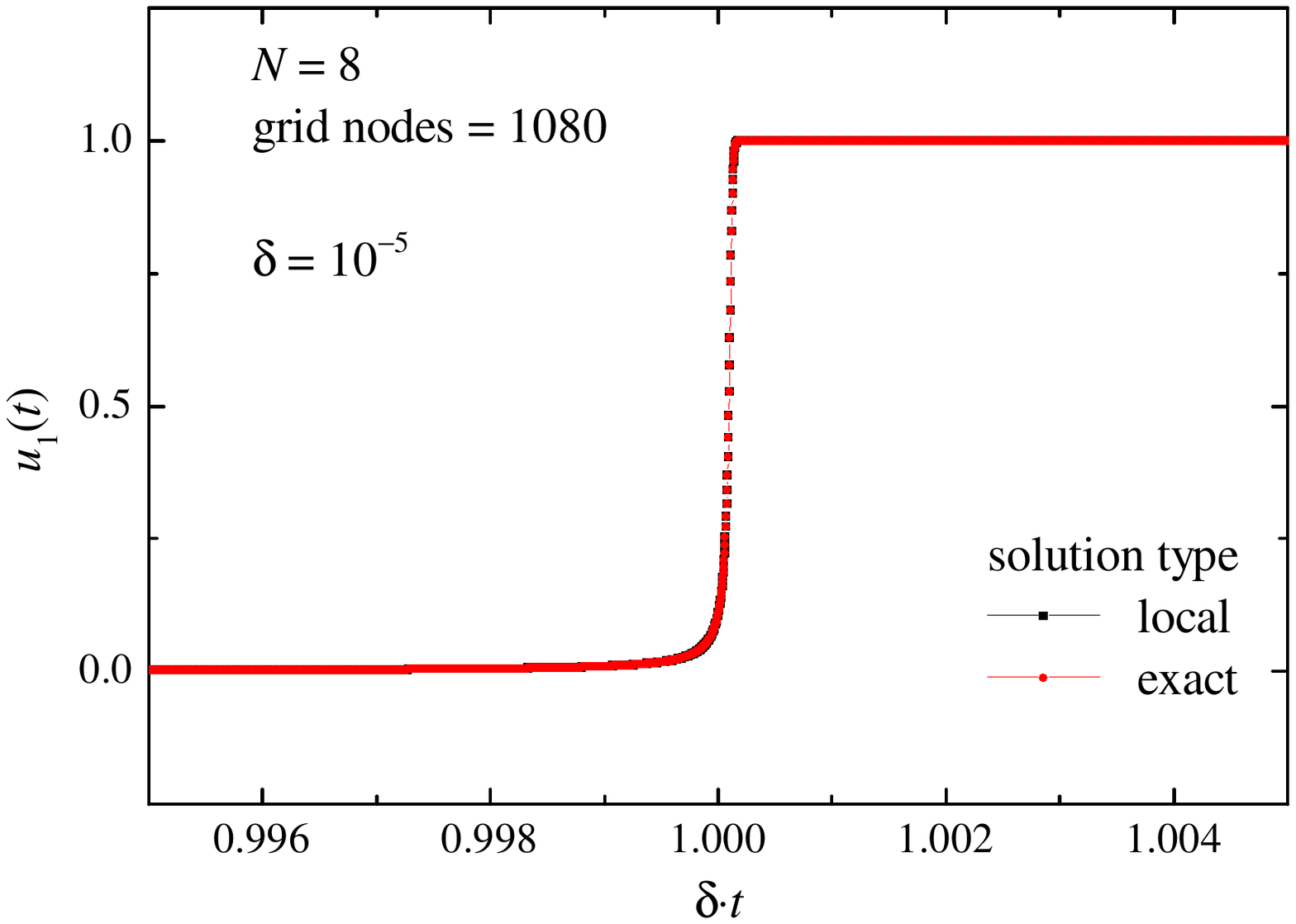}
\vspace{-8mm}\caption{\label{fig:shampine_flame_delta_10m5_sol_qug:b2}}
\end{subfigure}
\begin{subfigure}{0.320\textwidth}
\includegraphics[width=\textwidth]{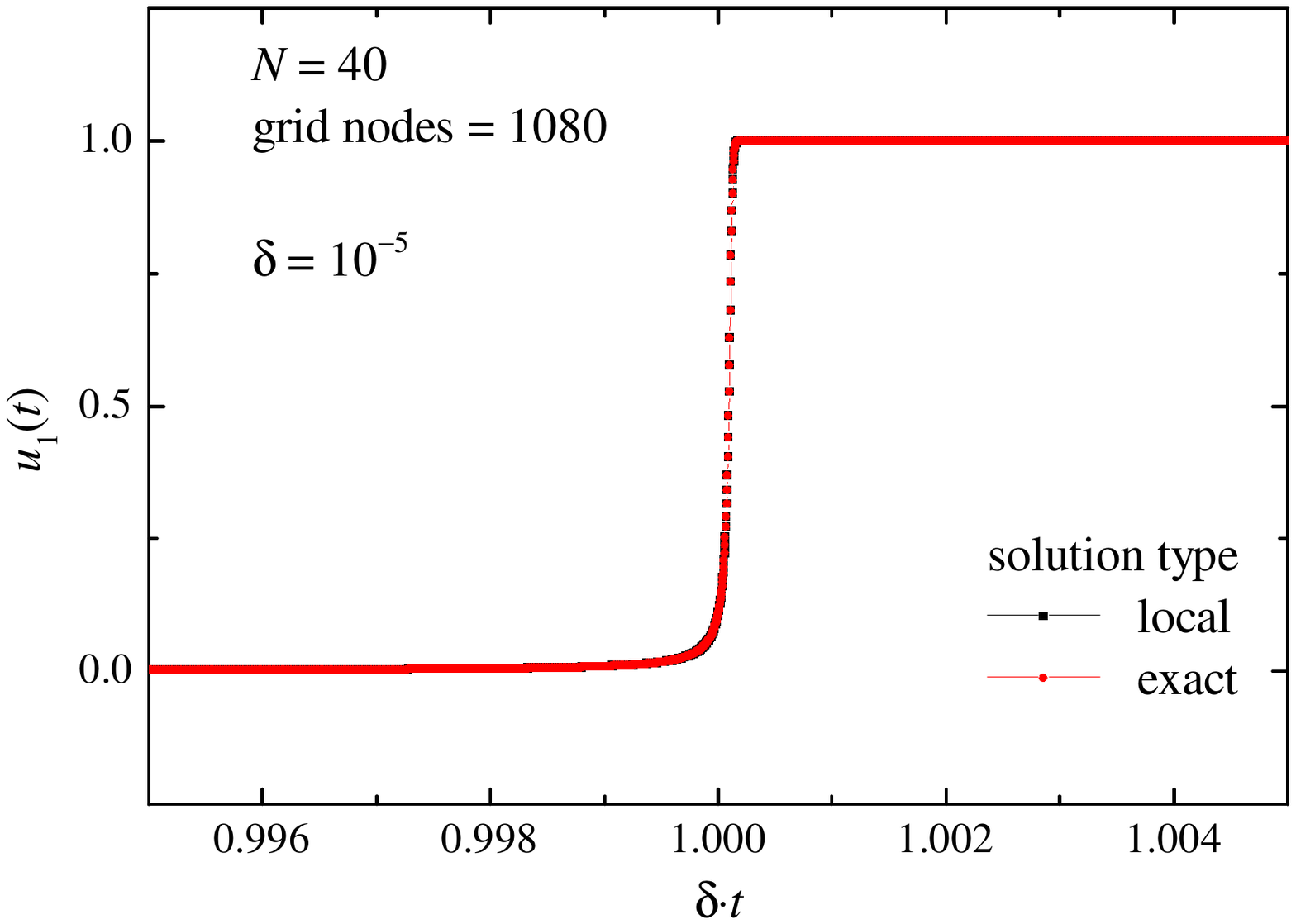}
\vspace{-8mm}\caption{\label{fig:shampine_flame_delta_10m5_sol_qug:b3}}
\end{subfigure}\\
\begin{subfigure}{0.320\textwidth}
\includegraphics[width=\textwidth]{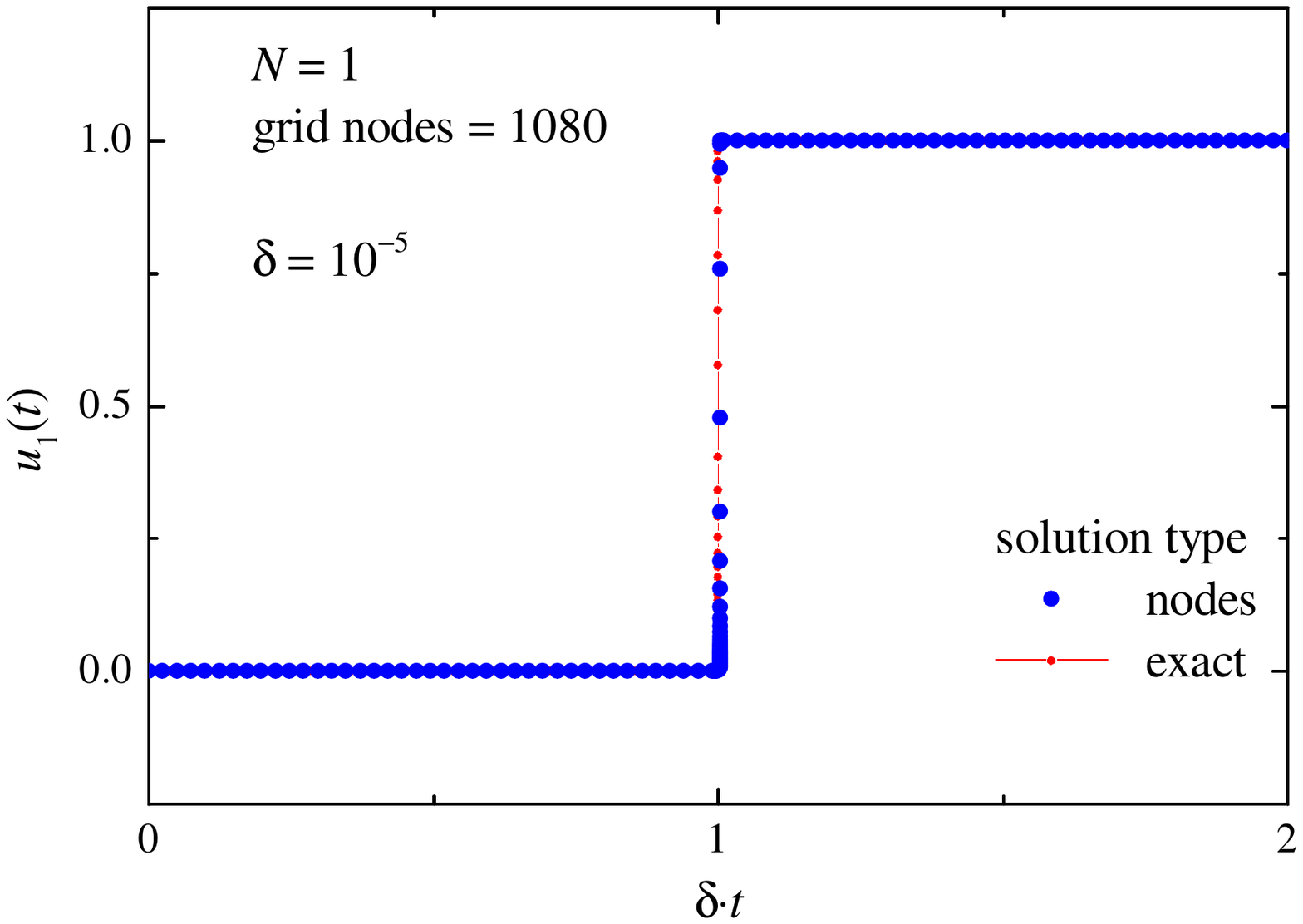}
\vspace{-8mm}\caption{\label{fig:shampine_flame_delta_10m5_sol_qug:c1}}
\end{subfigure}
\begin{subfigure}{0.320\textwidth}
\includegraphics[width=\textwidth]{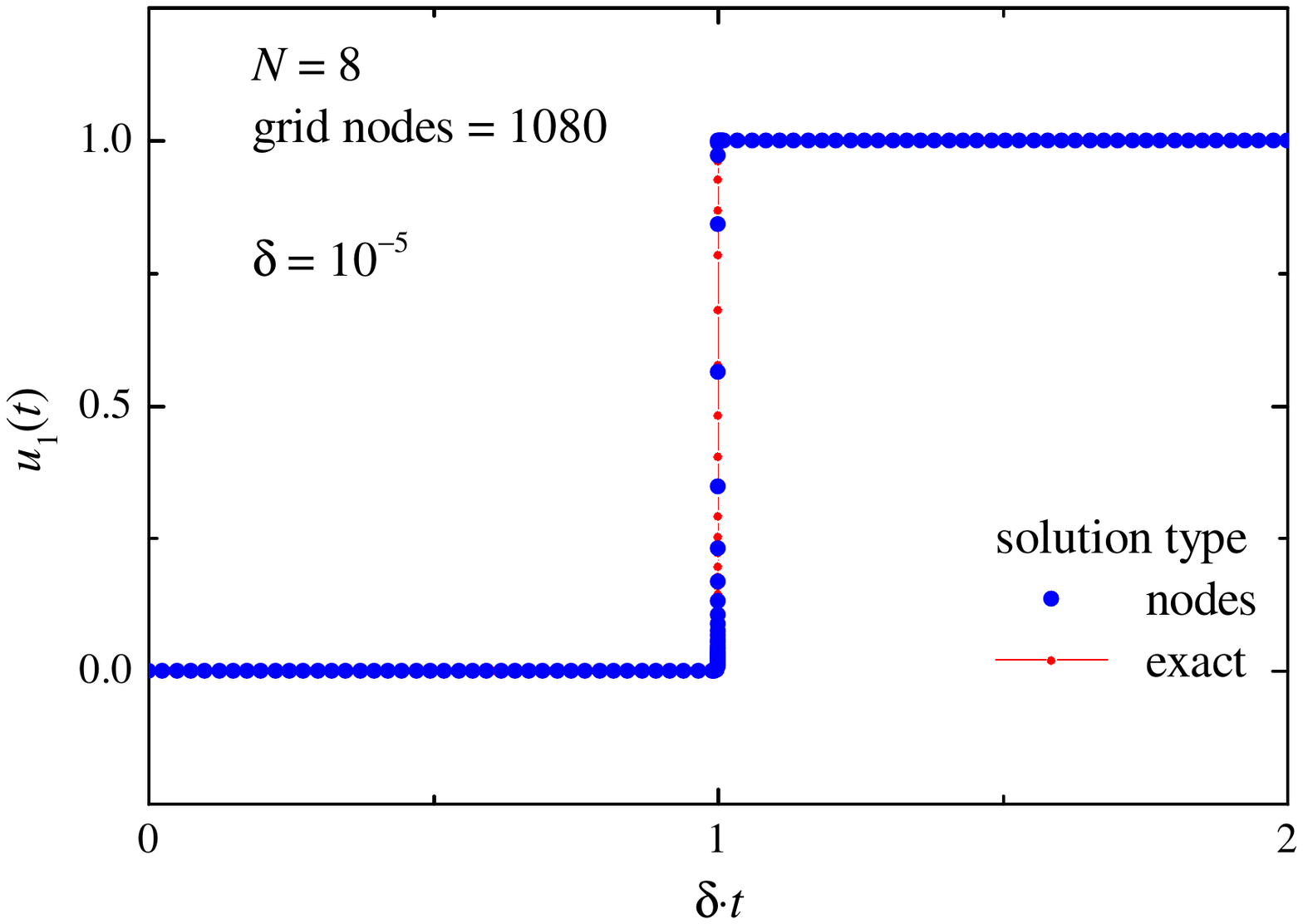}
\vspace{-8mm}\caption{\label{fig:shampine_flame_delta_10m5_sol_qug:c2}}
\end{subfigure}
\begin{subfigure}{0.320\textwidth}
\includegraphics[width=\textwidth]{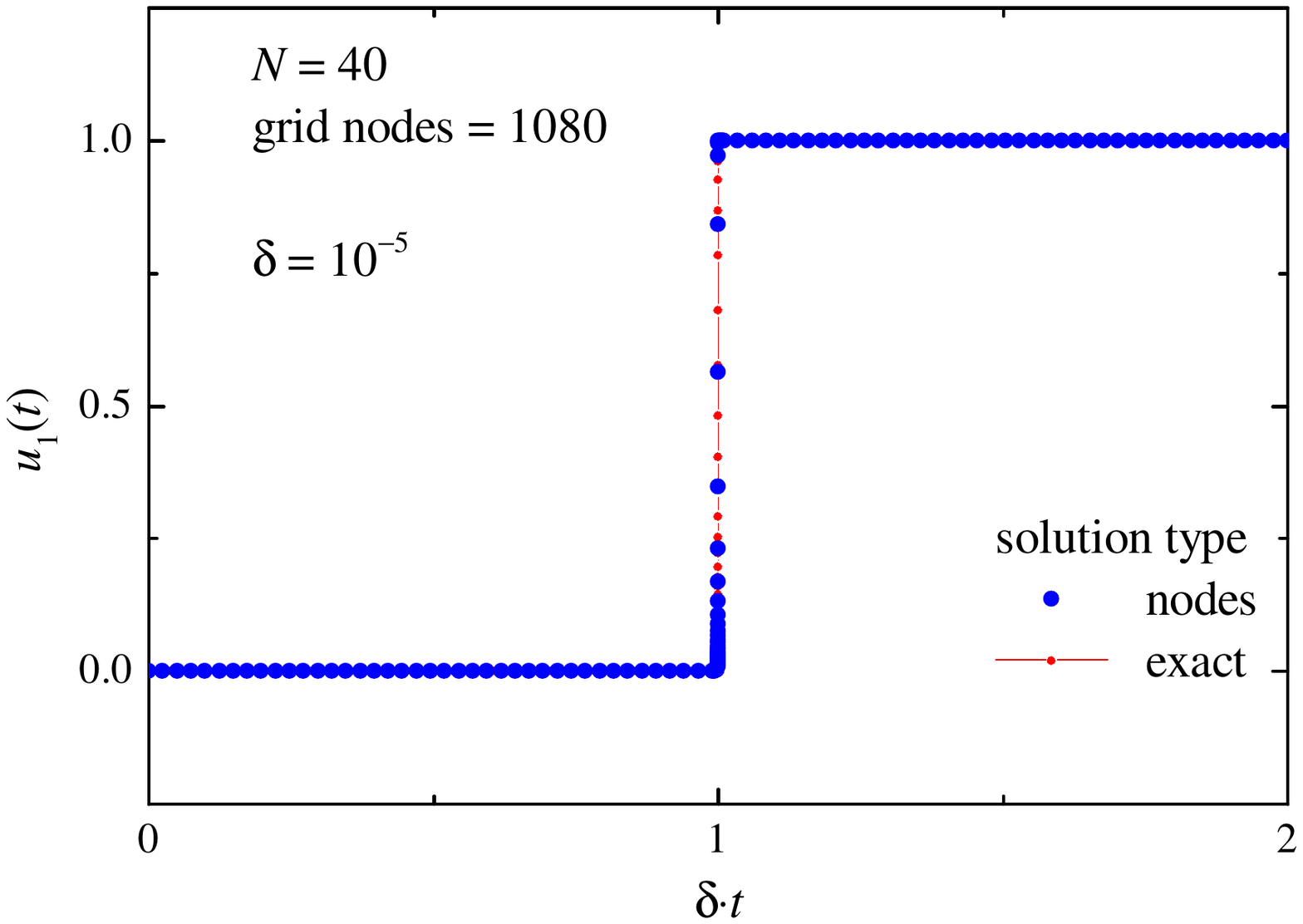}
\vspace{-8mm}\caption{\label{fig:shampine_flame_delta_10m5_sol_qug:c3}}
\end{subfigure}\\
\begin{subfigure}{0.320\textwidth}
\includegraphics[width=\textwidth]{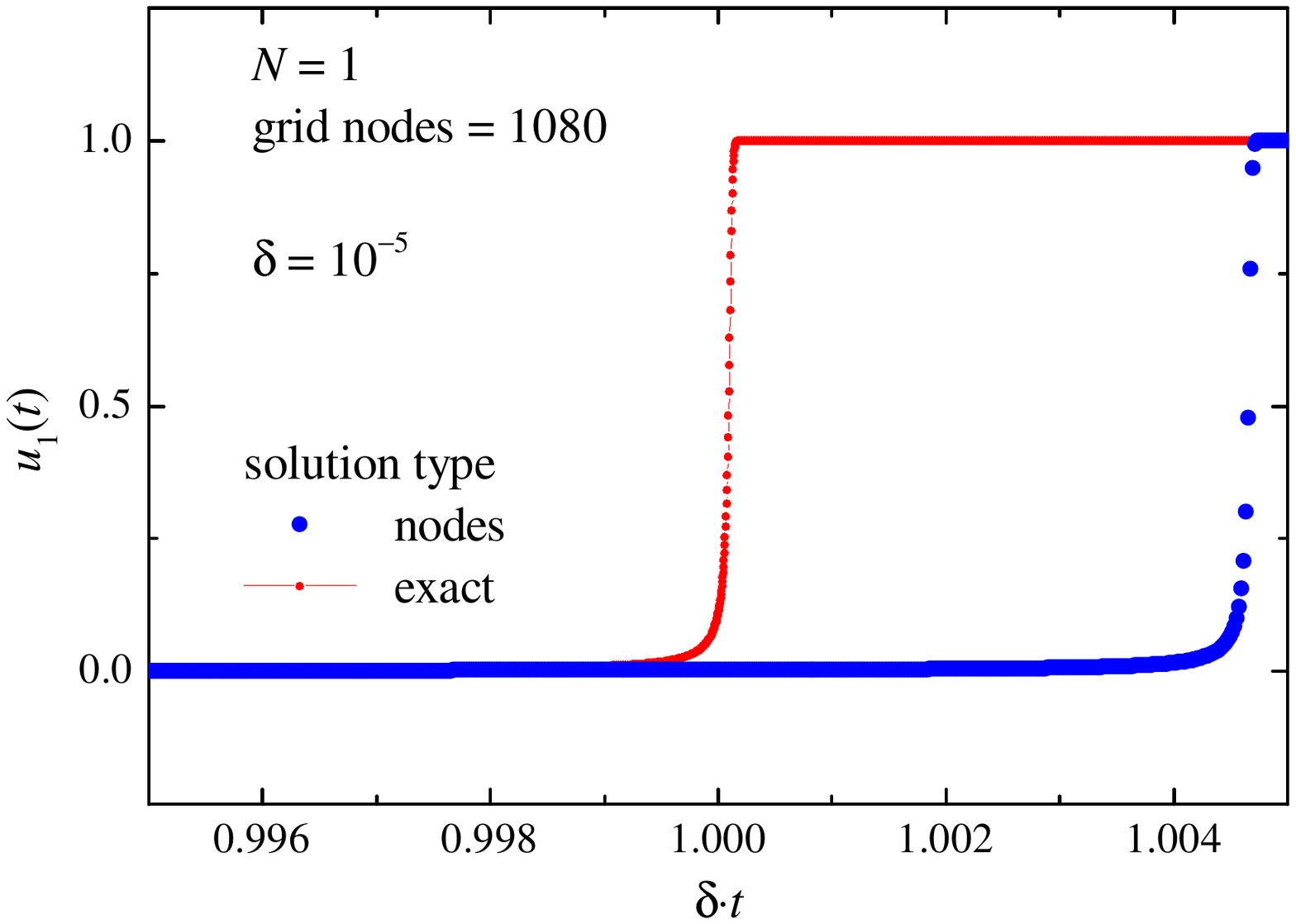}
\vspace{-8mm}\caption{\label{fig:shampine_flame_delta_10m5_sol_qug:d1}}
\end{subfigure}
\begin{subfigure}{0.320\textwidth}
\includegraphics[width=\textwidth]{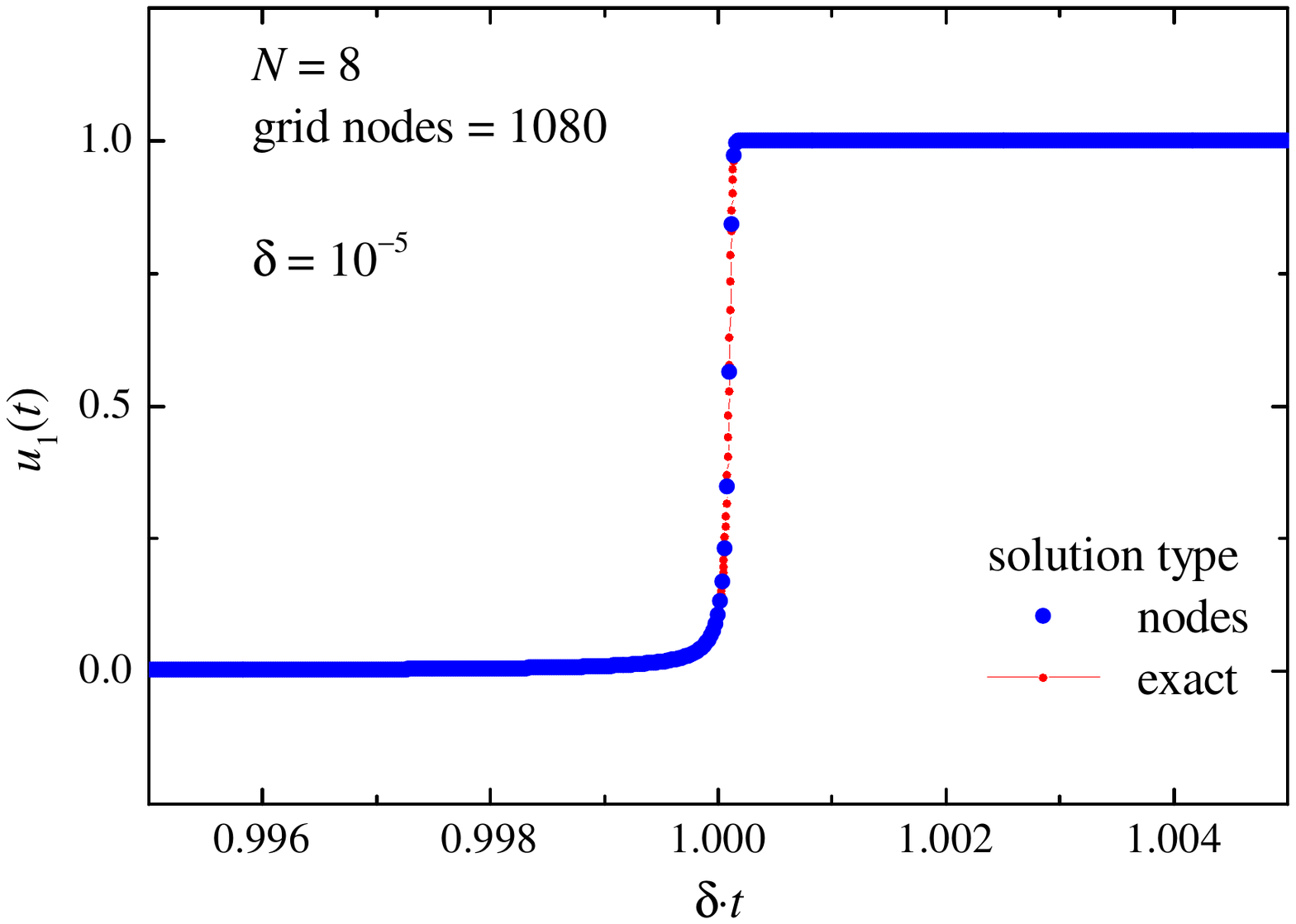}
\vspace{-8mm}\caption{\label{fig:shampine_flame_delta_10m5_sol_qug:d2}}
\end{subfigure}
\begin{subfigure}{0.320\textwidth}
\includegraphics[width=\textwidth]{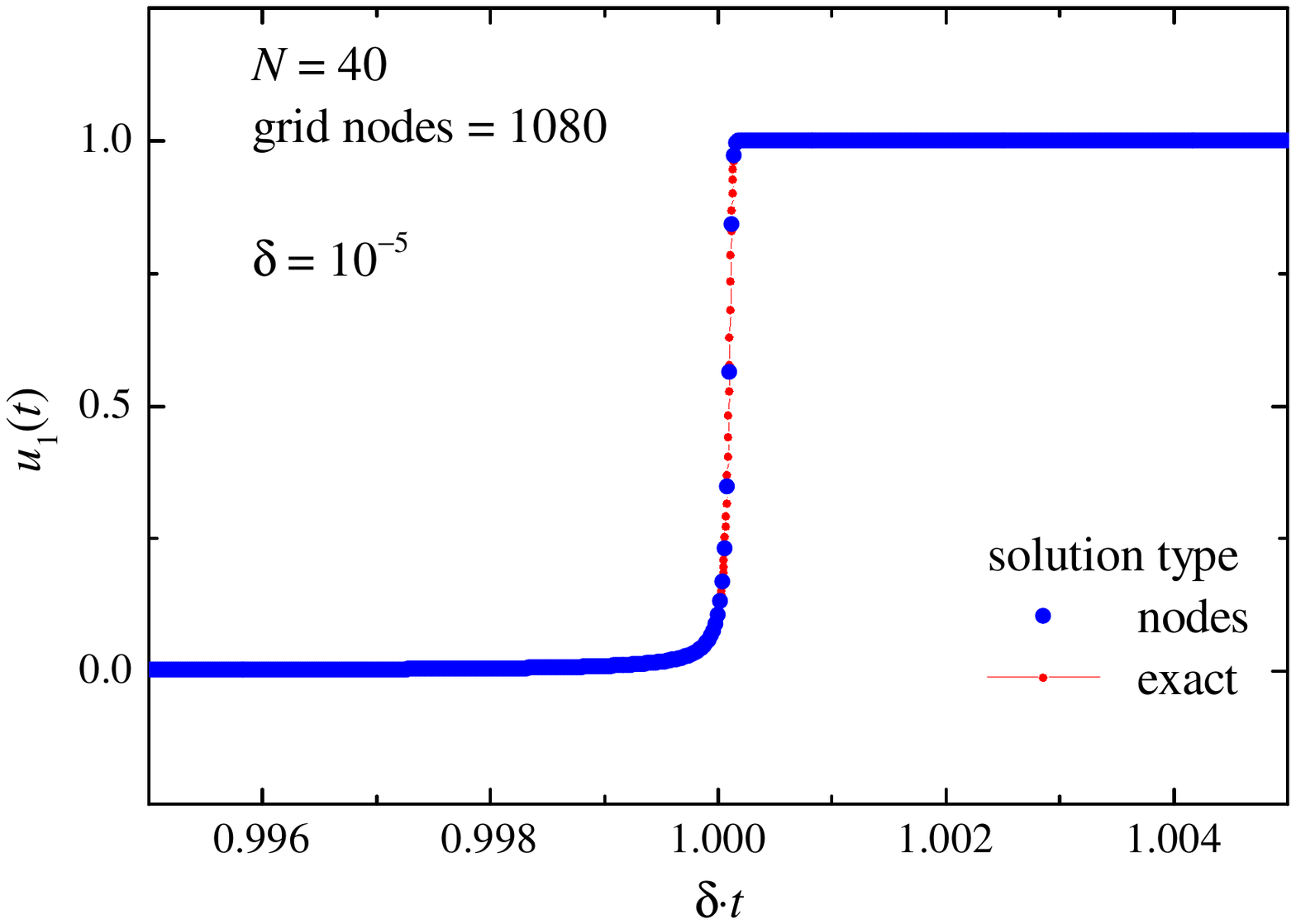}
\vspace{-8mm}\caption{\label{fig:shampine_flame_delta_10m5_sol_qug:d3}}
\end{subfigure}\\
\begin{subfigure}{0.320\textwidth}
\includegraphics[width=\textwidth]{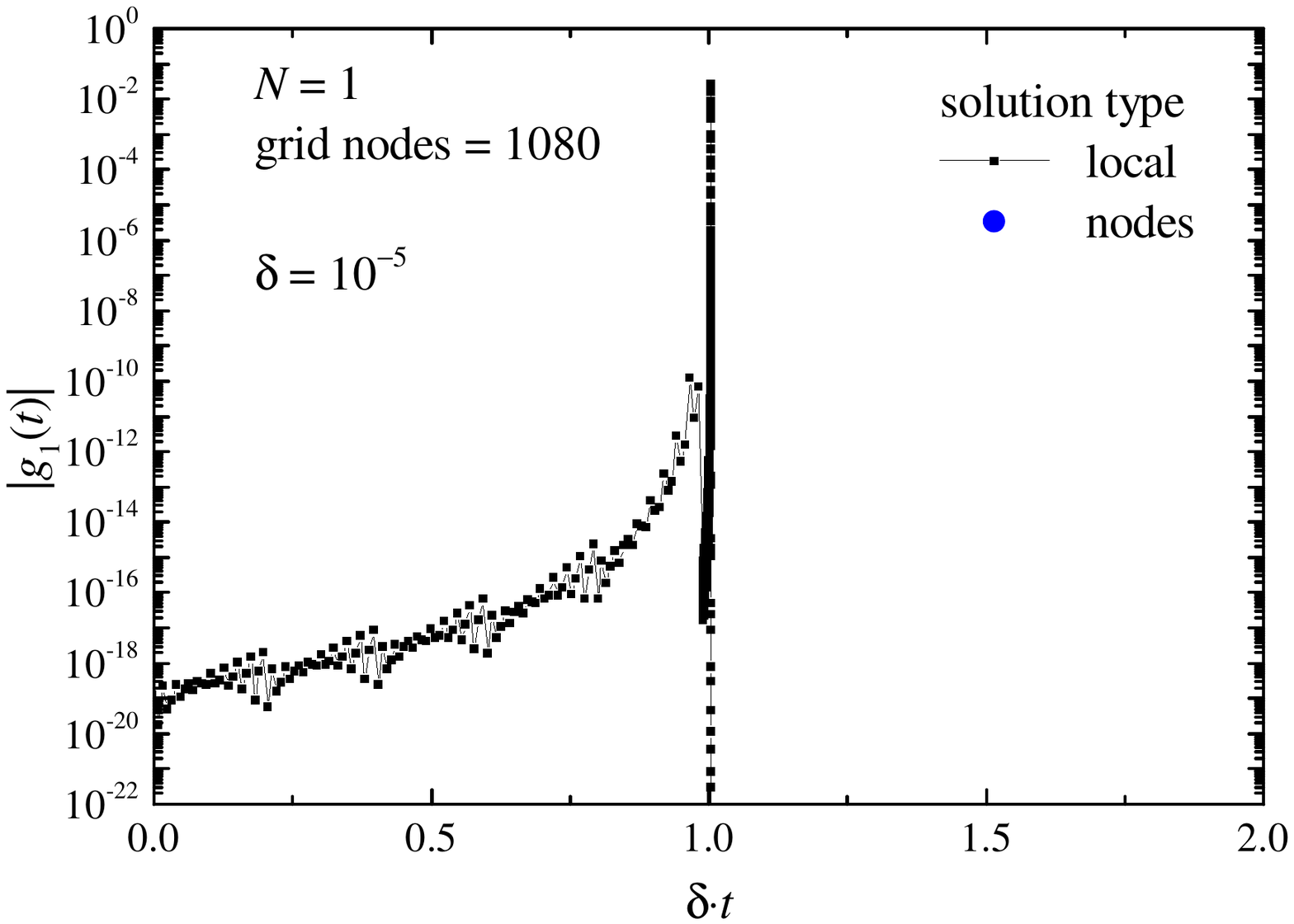}
\vspace{-8mm}\caption{\label{fig:shampine_flame_delta_10m5_sol_qug:e1}}
\end{subfigure}
\begin{subfigure}{0.320\textwidth}
\includegraphics[width=\textwidth]{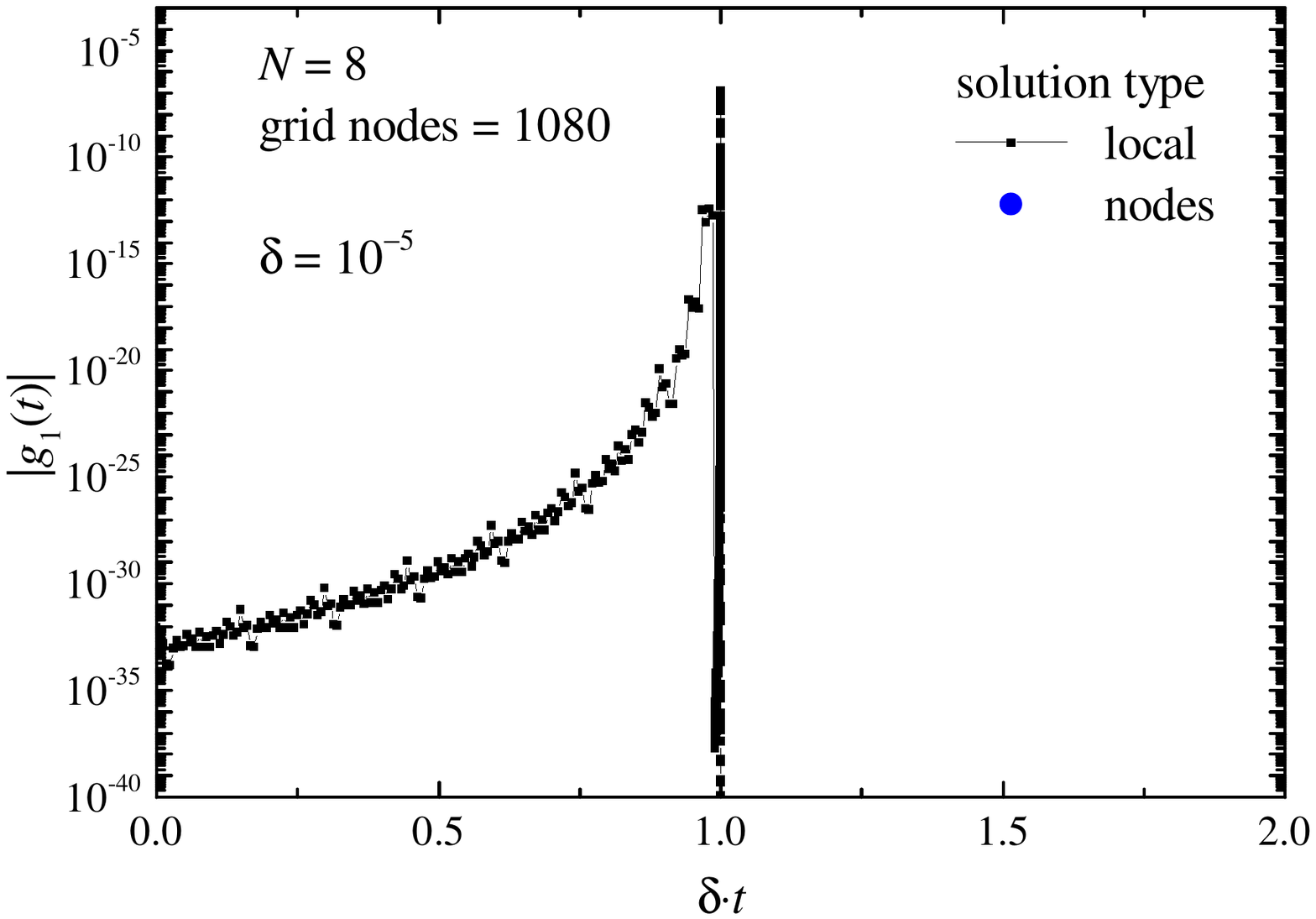}
\vspace{-8mm}\caption{\label{fig:shampine_flame_delta_10m5_sol_qug:e2}}
\end{subfigure}
\begin{subfigure}{0.320\textwidth}
\includegraphics[width=\textwidth]{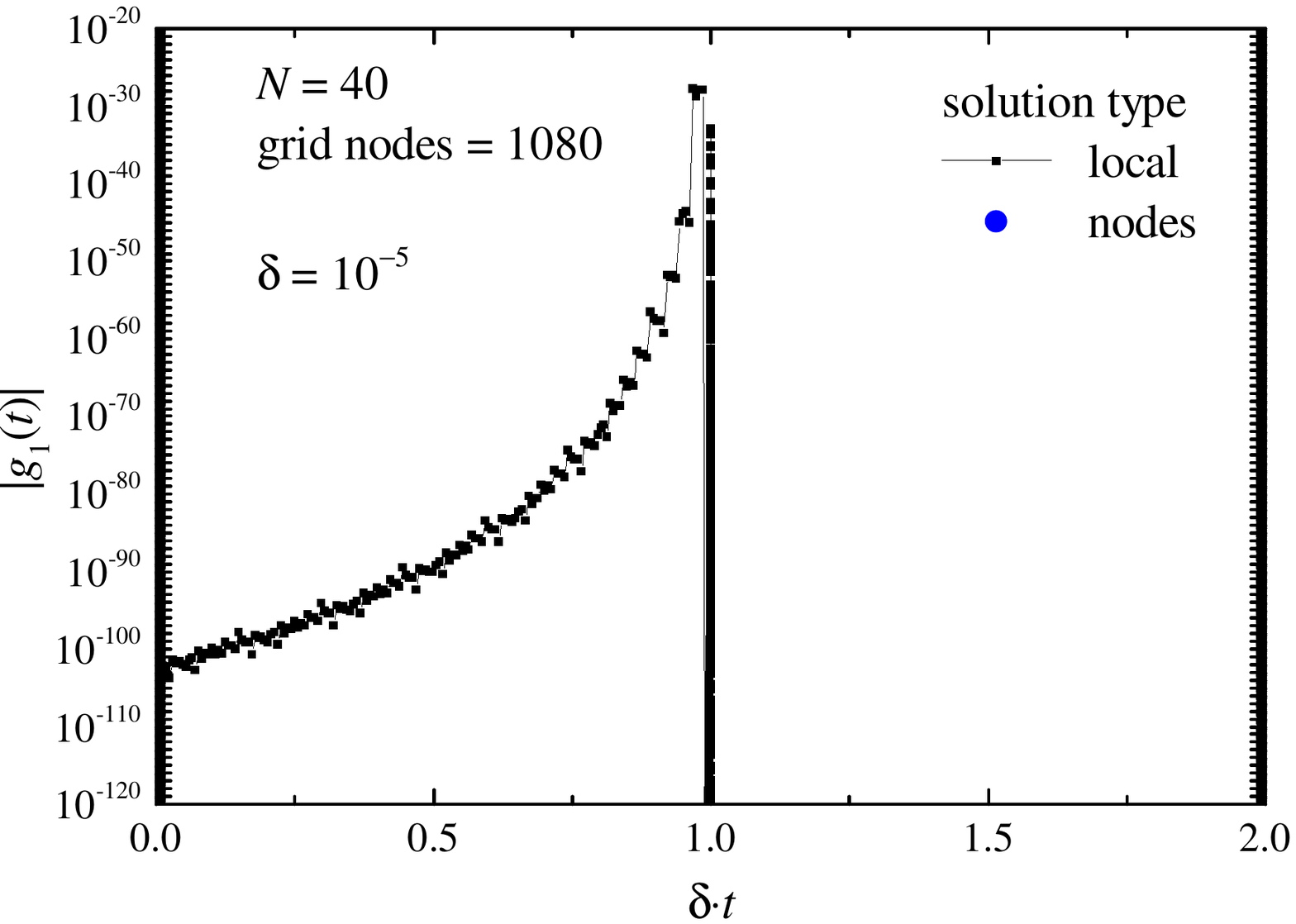}
\vspace{-8mm}\caption{\label{fig:shampine_flame_delta_10m5_sol_qug:e3}}
\end{subfigure}\\
\caption{%
Numerical solution of the stiff DAE system (\ref{eq:shampine_flame}) of index 1 with $\delta = 10^{-5}$. Comparison of the solution at nodes $\mathbf{u}_{n}$ (\subref{fig:shampine_flame_delta_10m5_sol_qug:c1}, \subref{fig:shampine_flame_delta_10m5_sol_qug:c2}, \subref{fig:shampine_flame_delta_10m5_sol_qug:c3}, \subref{fig:shampine_flame_delta_10m5_sol_qug:d1}, \subref{fig:shampine_flame_delta_10m5_sol_qug:d2}, \subref{fig:shampine_flame_delta_10m5_sol_qug:d3}), the local solution $\mathbf{u}_{L}(t)$ (\subref{fig:shampine_flame_delta_10m5_sol_qug:a1}, \subref{fig:shampine_flame_delta_10m5_sol_qug:a2}, \subref{fig:shampine_flame_delta_10m5_sol_qug:a3}, \subref{fig:shampine_flame_delta_10m5_sol_qug:b1}, \subref{fig:shampine_flame_delta_10m5_sol_qug:b2}, \subref{fig:shampine_flame_delta_10m5_sol_qug:b3}) and the exact solution $\mathbf{u}^{\rm ex}(t)$ for component $u_{1}$ (\subref{fig:shampine_flame_delta_10m5_sol_qug:a1}, \subref{fig:shampine_flame_delta_10m5_sol_qug:a2}, \subref{fig:shampine_flame_delta_10m5_sol_qug:a3}, \subref{fig:shampine_flame_delta_10m5_sol_qug:b1}, \subref{fig:shampine_flame_delta_10m5_sol_qug:b2}, \subref{fig:shampine_flame_delta_10m5_sol_qug:b3}, \subref{fig:shampine_flame_delta_10m5_sol_qug:c1}, \subref{fig:shampine_flame_delta_10m5_sol_qug:c2}, \subref{fig:shampine_flame_delta_10m5_sol_qug:c3}, \subref{fig:shampine_flame_delta_10m5_sol_qug:d1}, \subref{fig:shampine_flame_delta_10m5_sol_qug:d2}, \subref{fig:shampine_flame_delta_10m5_sol_qug:d3}), quantitative satisfiability of the conditions $g_{1} = 0$ (\subref{fig:shampine_flame_delta_10m5_sol_qug:e1}, \subref{fig:shampine_flame_delta_10m5_sol_qug:e2}, \subref{fig:shampine_flame_delta_10m5_sol_qug:e3}), obtained using polynomials with degrees $N = 1$ (\subref{fig:shampine_flame_delta_10m5_sol_qug:a1}, \subref{fig:shampine_flame_delta_10m5_sol_qug:b1}, \subref{fig:shampine_flame_delta_10m5_sol_qug:c1}, \subref{fig:shampine_flame_delta_10m5_sol_qug:d1}, \subref{fig:shampine_flame_delta_10m5_sol_qug:e1}), $N = 8$ (\subref{fig:shampine_flame_delta_10m5_sol_qug:a2}, \subref{fig:shampine_flame_delta_10m5_sol_qug:b2}, \subref{fig:shampine_flame_delta_10m5_sol_qug:c2}, \subref{fig:shampine_flame_delta_10m5_sol_qug:d2}, \subref{fig:shampine_flame_delta_10m5_sol_qug:e2}) and $N = 40$ (\subref{fig:shampine_flame_delta_10m5_sol_qug:a3}, \subref{fig:shampine_flame_delta_10m5_sol_qug:b3}, \subref{fig:shampine_flame_delta_10m5_sol_qug:c3}, \subref{fig:shampine_flame_delta_10m5_sol_qug:d3}, \subref{fig:shampine_flame_delta_10m5_sol_qug:e3}).
}
\label{fig:shampine_flame_delta_10m5_sol_qug}
\end{figure}

\begin{figure}[h!]
\captionsetup[subfigure]{%
	position=bottom,
	font+=smaller,
	textfont=normalfont,
	singlelinecheck=off,
	justification=raggedright
}
\centering
\begin{subfigure}{0.320\textwidth}
\includegraphics[width=\textwidth]{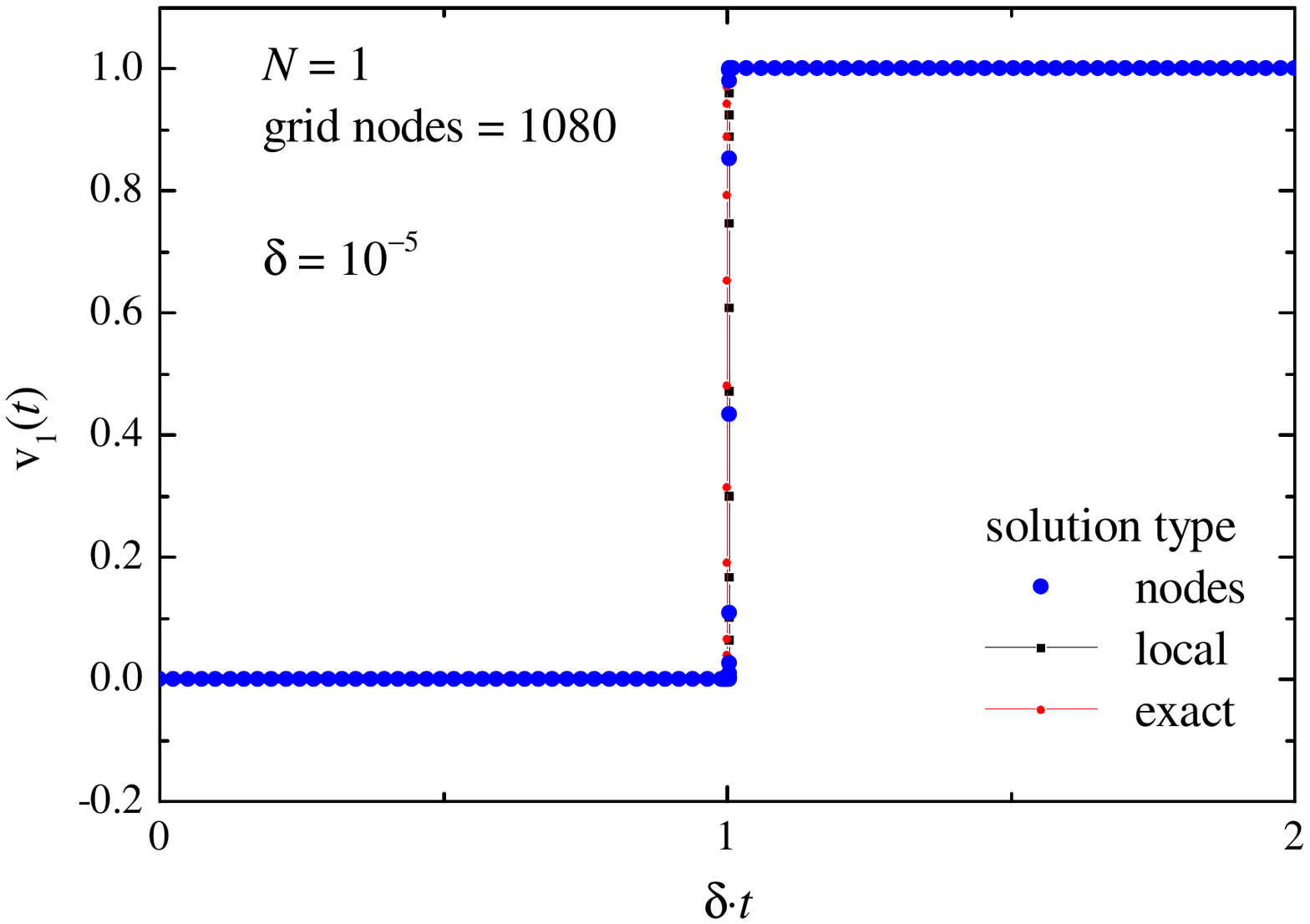}
\vspace{-8mm}\caption{\label{fig:shampine_flame_delta_10m5_sol_v_epss:a1}}
\end{subfigure}
\begin{subfigure}{0.320\textwidth}
\includegraphics[width=\textwidth]{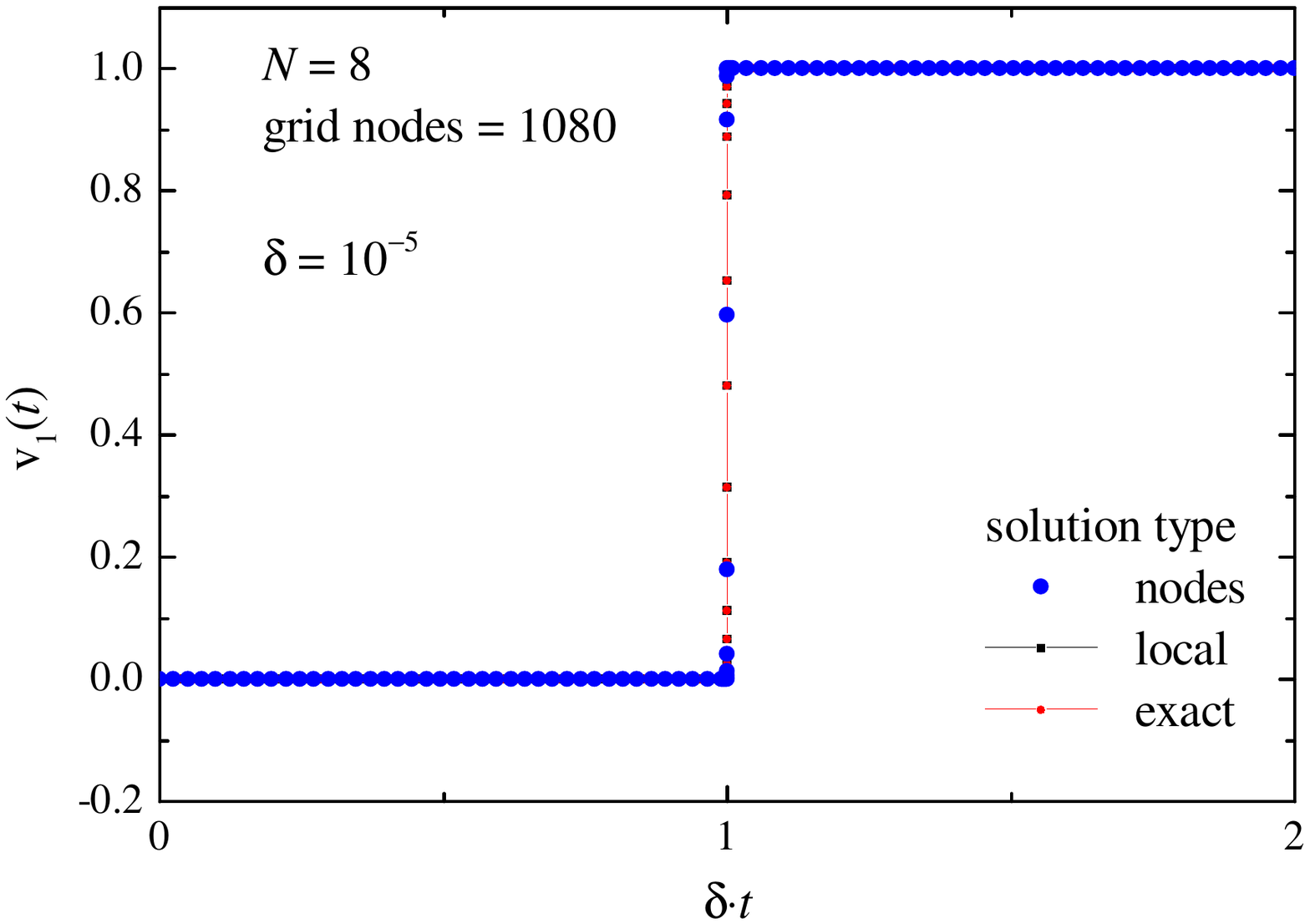}
\vspace{-8mm}\caption{\label{fig:shampine_flame_delta_10m5_sol_v_epss:a2}}
\end{subfigure}
\begin{subfigure}{0.320\textwidth}
\includegraphics[width=\textwidth]{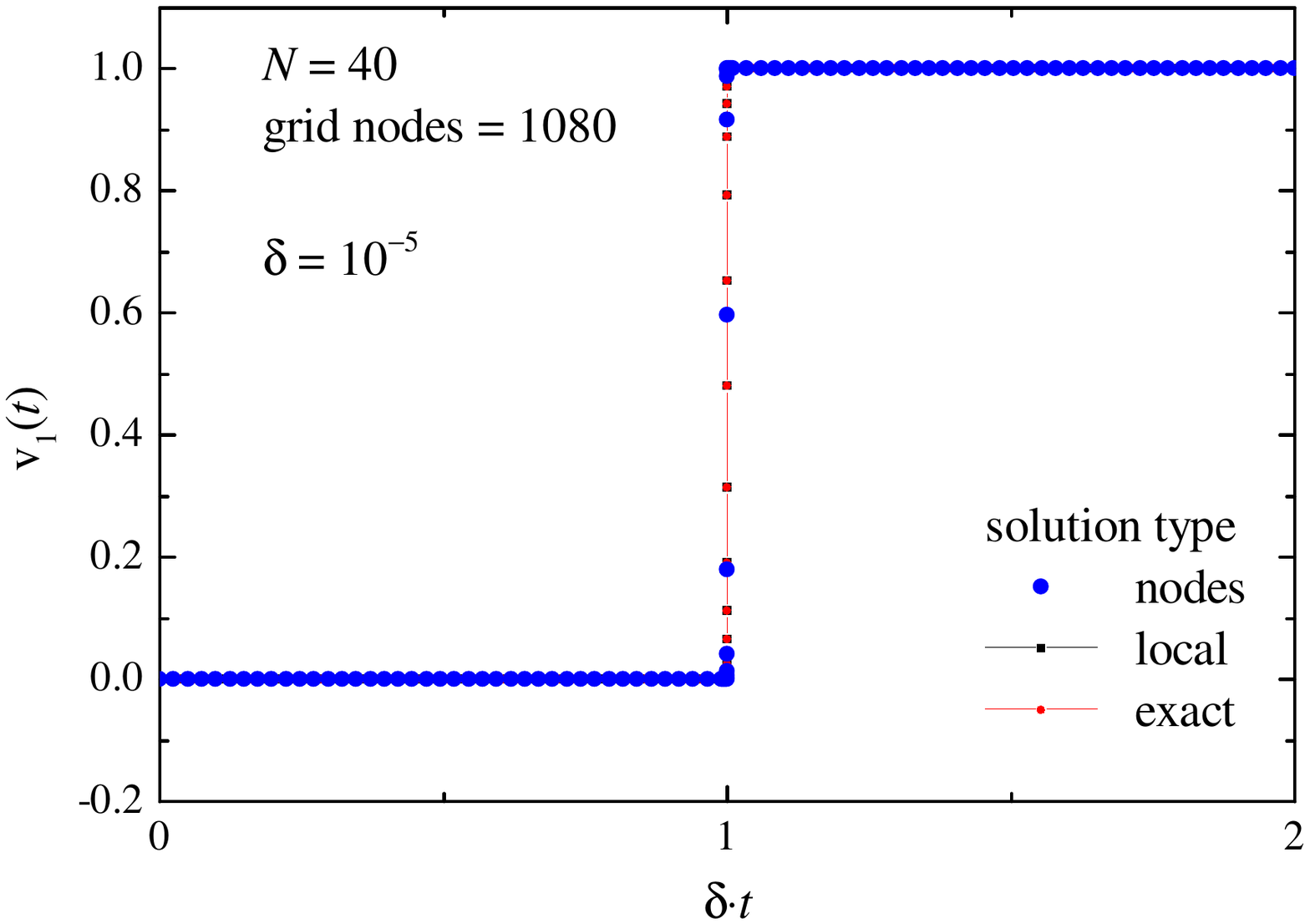}
\vspace{-8mm}\caption{\label{fig:shampine_flame_delta_10m5_sol_v_epss:a3}}
\end{subfigure}\\
\begin{subfigure}{0.320\textwidth}
\includegraphics[width=\textwidth]{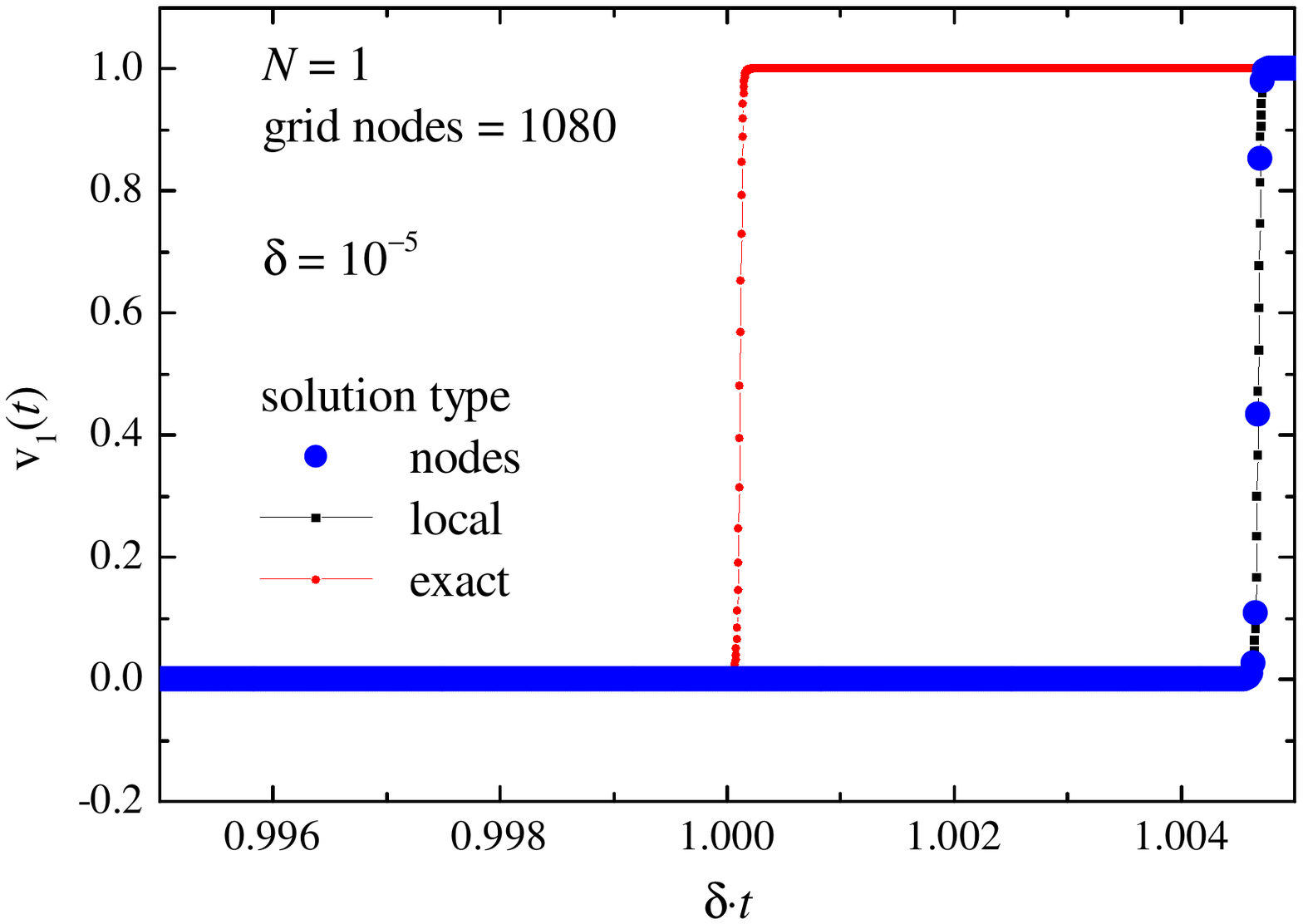}
\vspace{-8mm}\caption{\label{fig:shampine_flame_delta_10m5_sol_v_epss:b1}}
\end{subfigure}
\begin{subfigure}{0.320\textwidth}
\includegraphics[width=\textwidth]{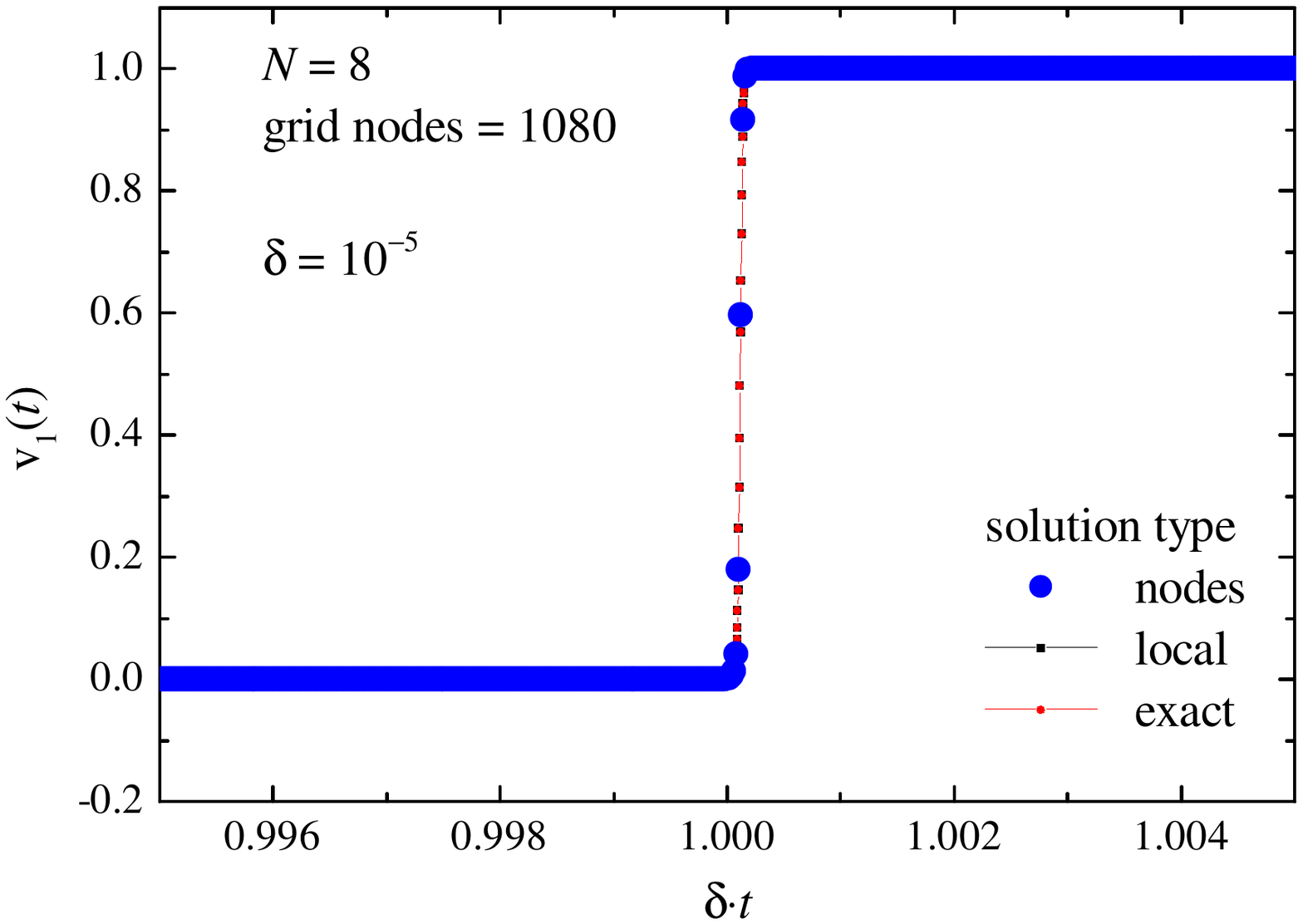}
\vspace{-8mm}\caption{\label{fig:shampine_flame_delta_10m5_sol_v_epss:b2}}
\end{subfigure}
\begin{subfigure}{0.320\textwidth}
\includegraphics[width=\textwidth]{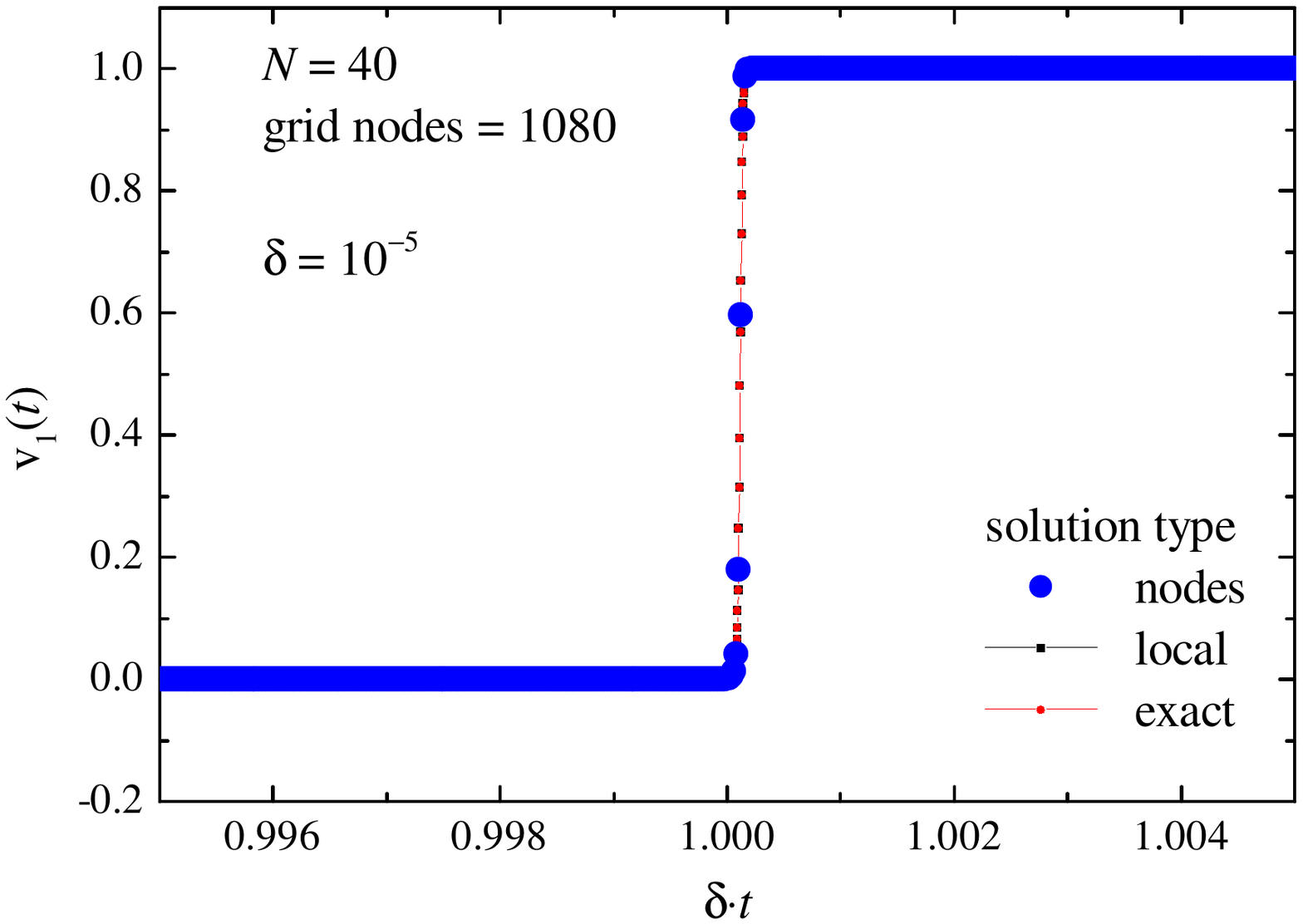}
\vspace{-8mm}\caption{\label{fig:shampine_flame_delta_10m5_sol_v_epss:b3}}
\end{subfigure}\\
\begin{subfigure}{0.320\textwidth}
\includegraphics[width=\textwidth]{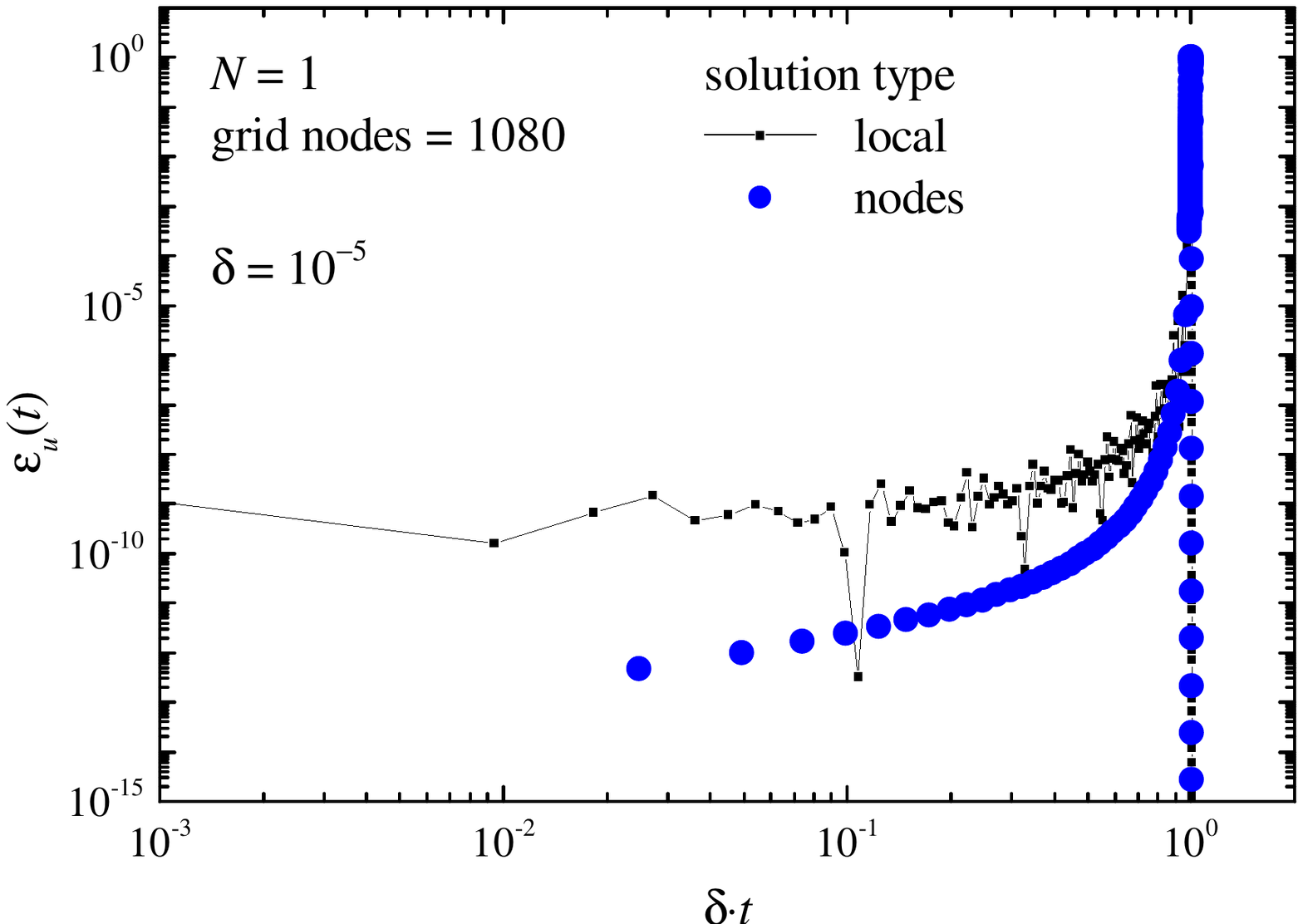}
\vspace{-8mm}\caption{\label{fig:shampine_flame_delta_10m5_sol_v_epss:c1}}
\end{subfigure}
\begin{subfigure}{0.320\textwidth}
\includegraphics[width=\textwidth]{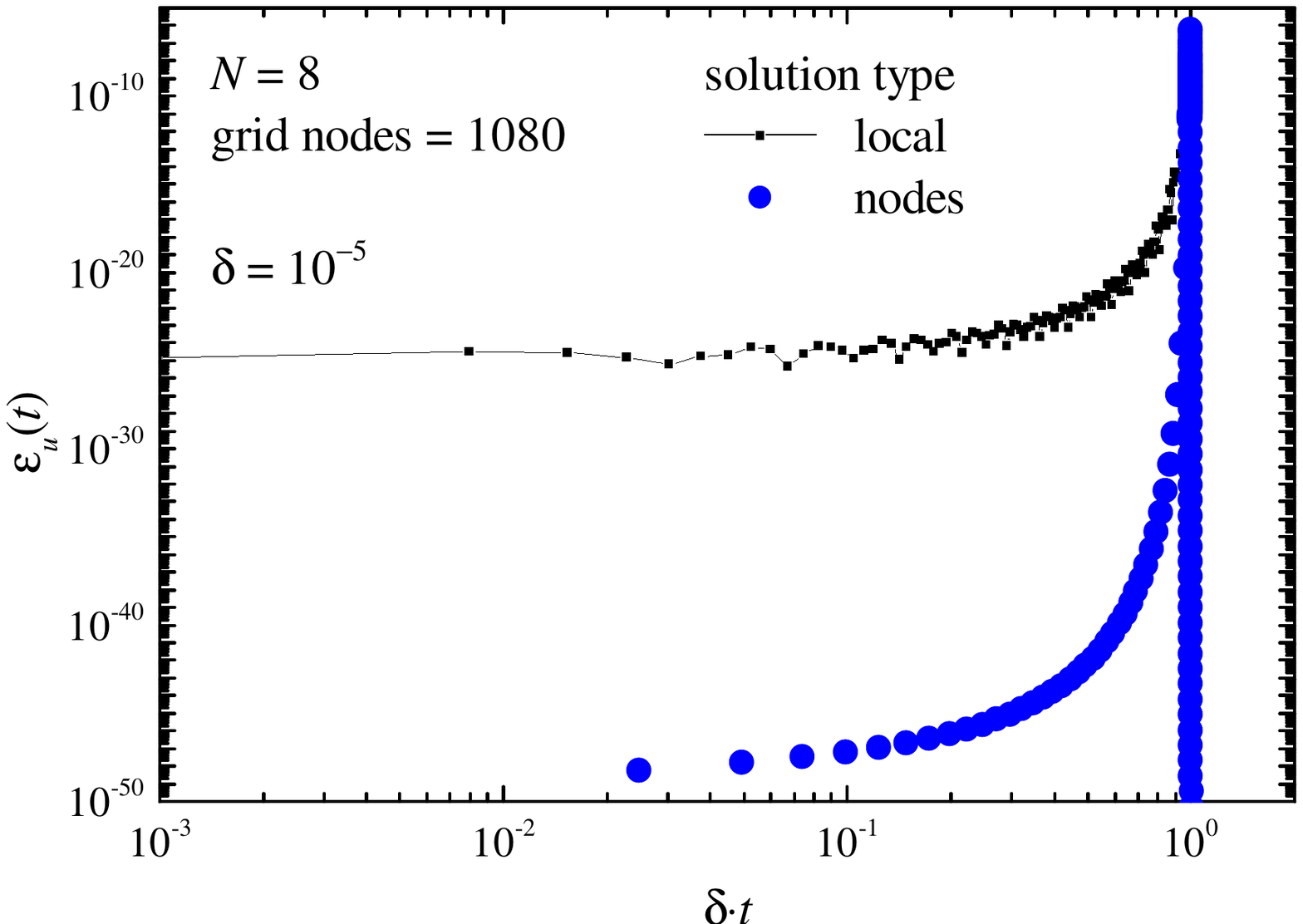}
\vspace{-8mm}\caption{\label{fig:shampine_flame_delta_10m5_sol_v_epss:c2}}
\end{subfigure}
\begin{subfigure}{0.320\textwidth}
\includegraphics[width=\textwidth]{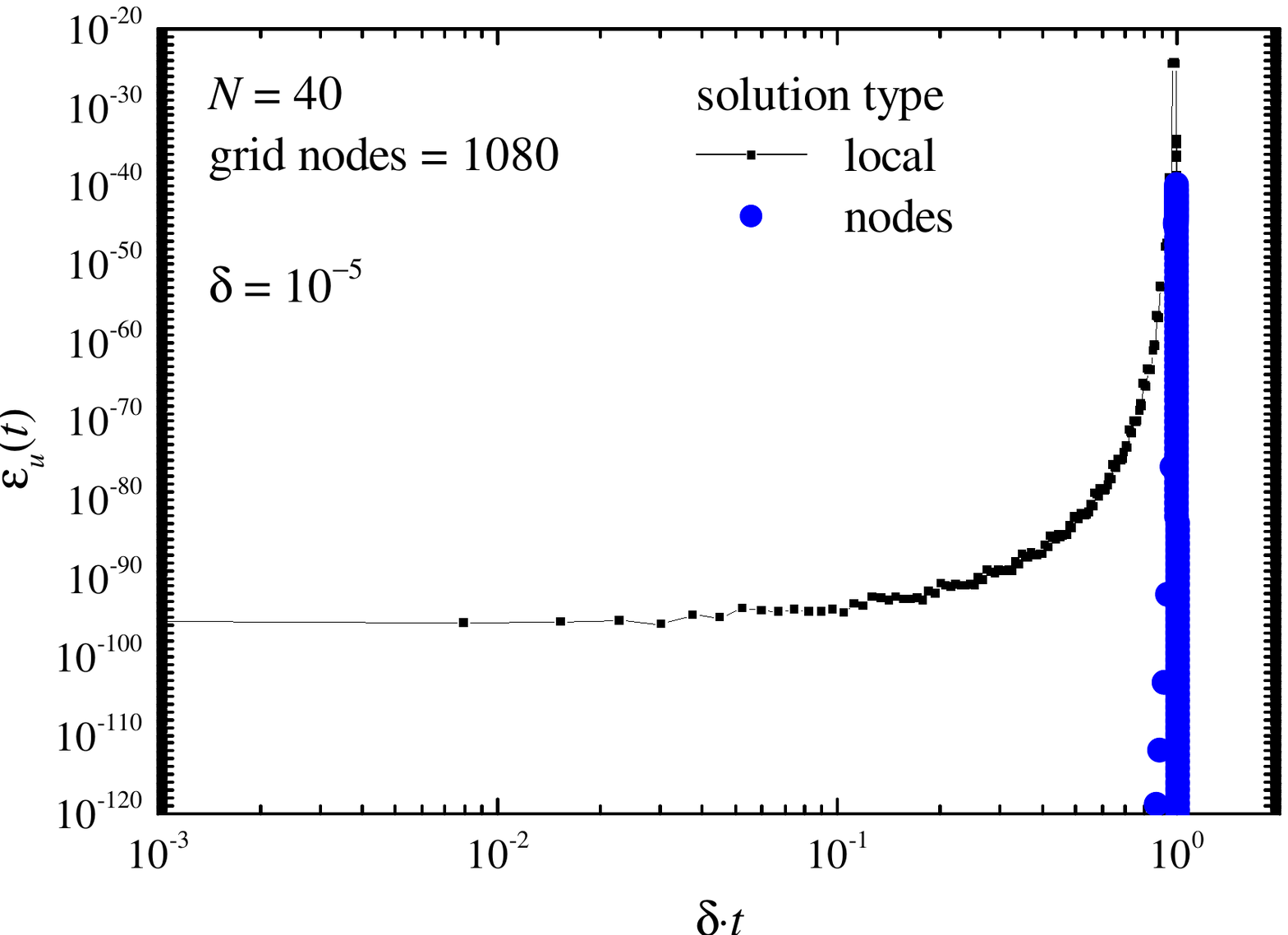}
\vspace{-8mm}\caption{\label{fig:shampine_flame_delta_10m5_sol_v_epss:c3}}
\end{subfigure}\\
\begin{subfigure}{0.320\textwidth}
\includegraphics[width=\textwidth]{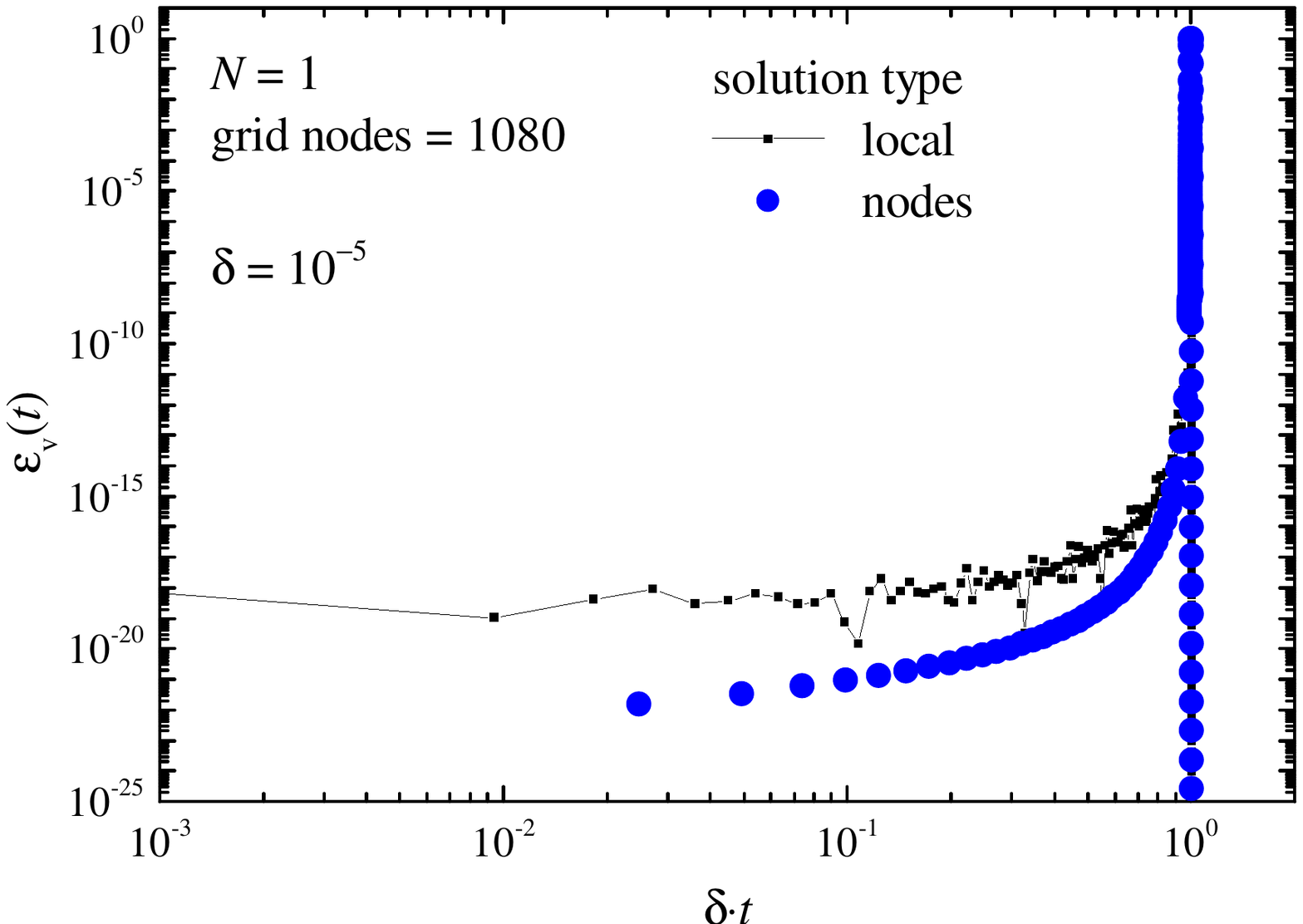}
\vspace{-8mm}\caption{\label{fig:shampine_flame_delta_10m5_sol_v_epss:d1}}
\end{subfigure}
\begin{subfigure}{0.320\textwidth}
\includegraphics[width=\textwidth]{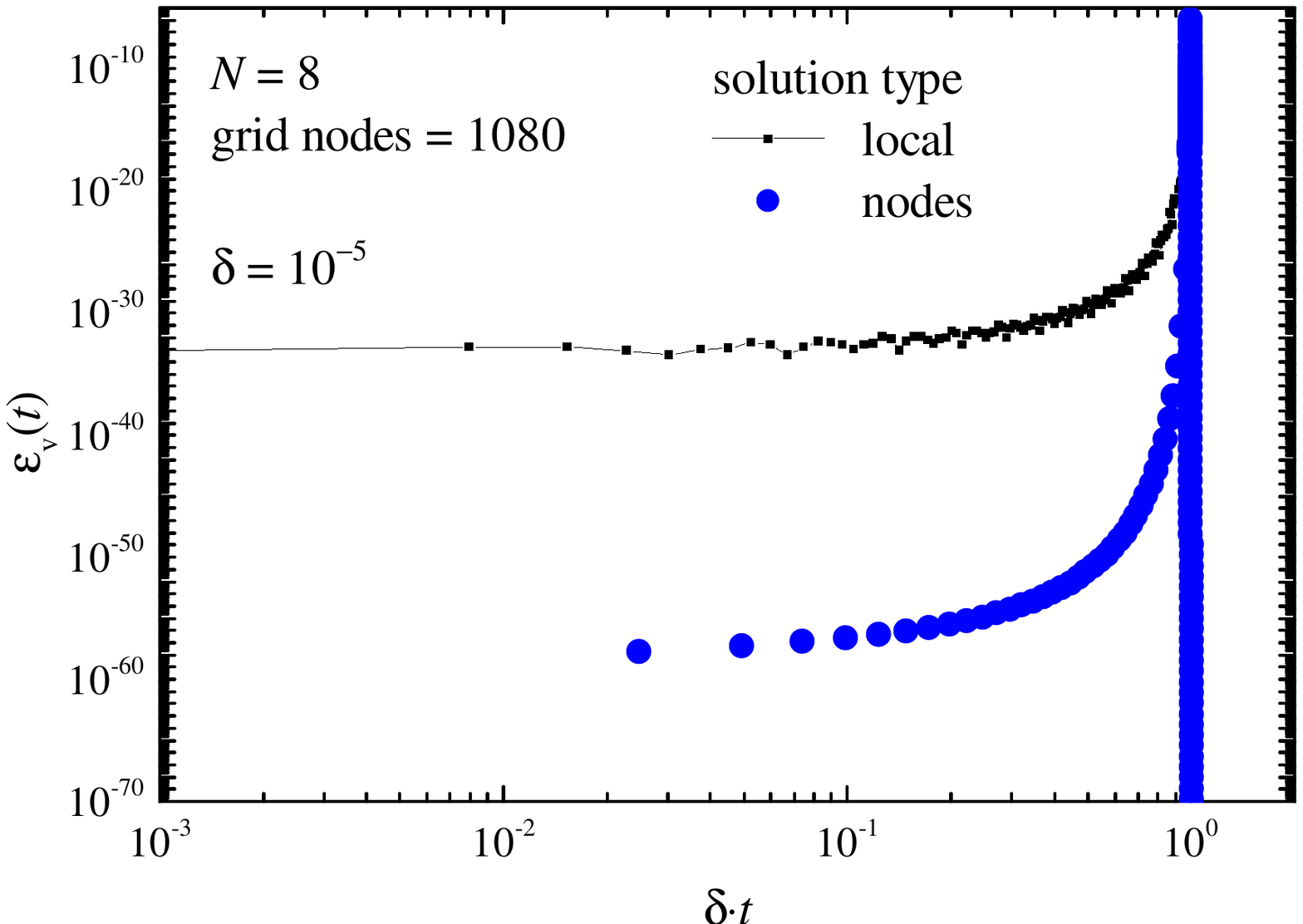}
\vspace{-8mm}\caption{\label{fig:shampine_flame_delta_10m5_sol_v_epss:d2}}
\end{subfigure}
\begin{subfigure}{0.320\textwidth}
\includegraphics[width=\textwidth]{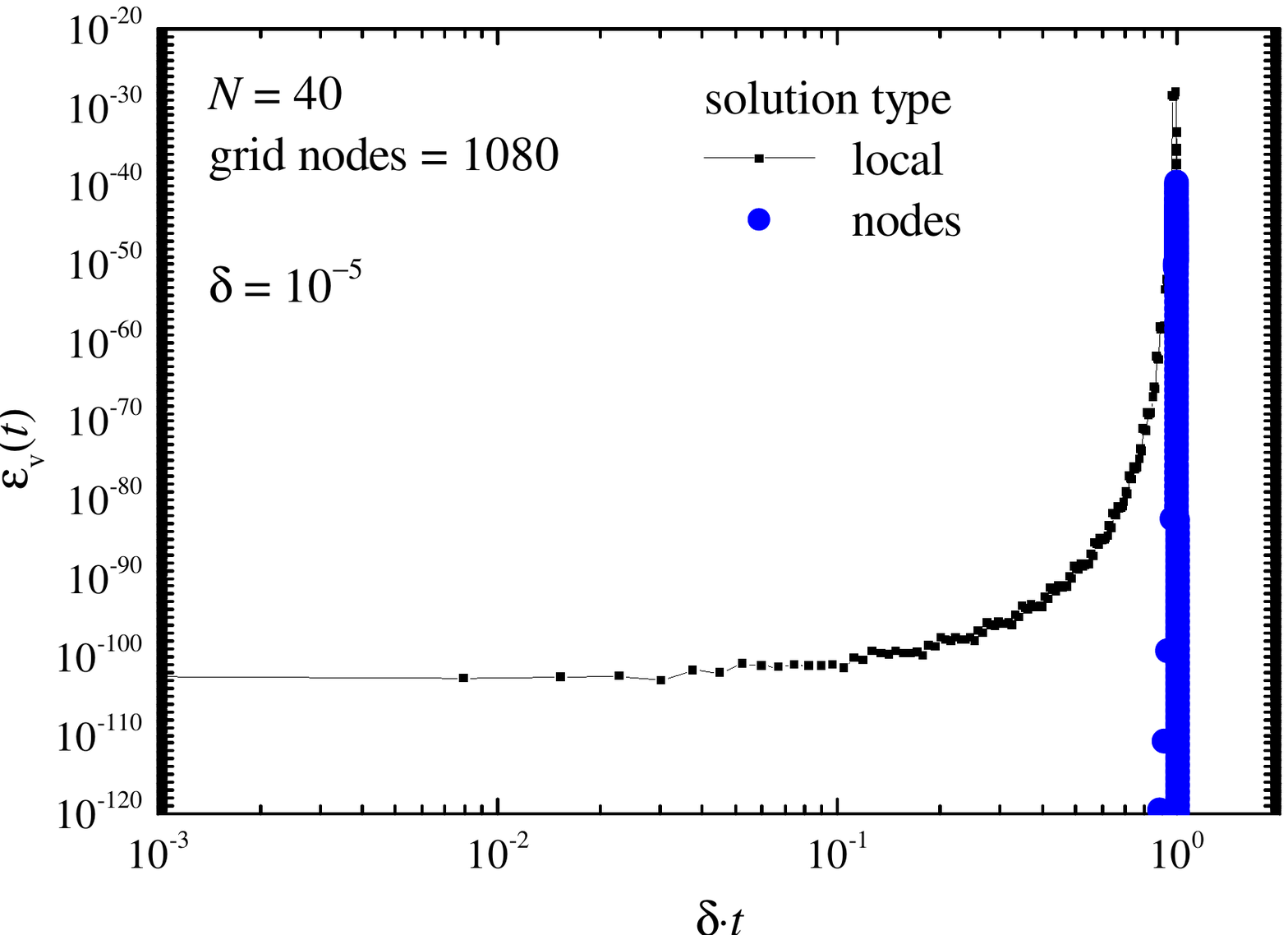}
\vspace{-8mm}\caption{\label{fig:shampine_flame_delta_10m5_sol_v_epss:d3}}
\end{subfigure}\\
\begin{subfigure}{0.320\textwidth}
\includegraphics[width=\textwidth]{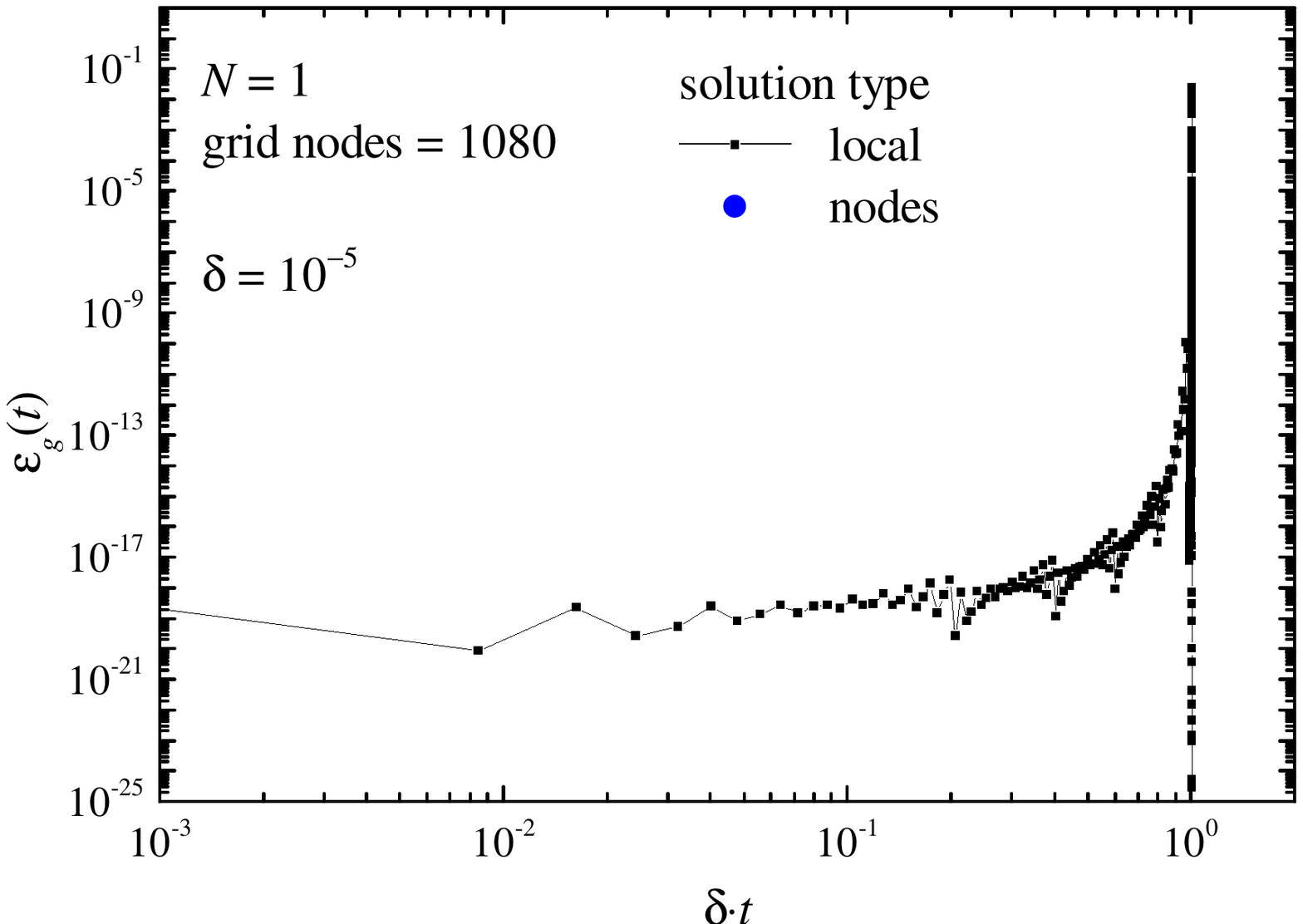}
\vspace{-8mm}\caption{\label{fig:shampine_flame_delta_10m5_sol_v_epss:e1}}
\end{subfigure}
\begin{subfigure}{0.320\textwidth}
\includegraphics[width=\textwidth]{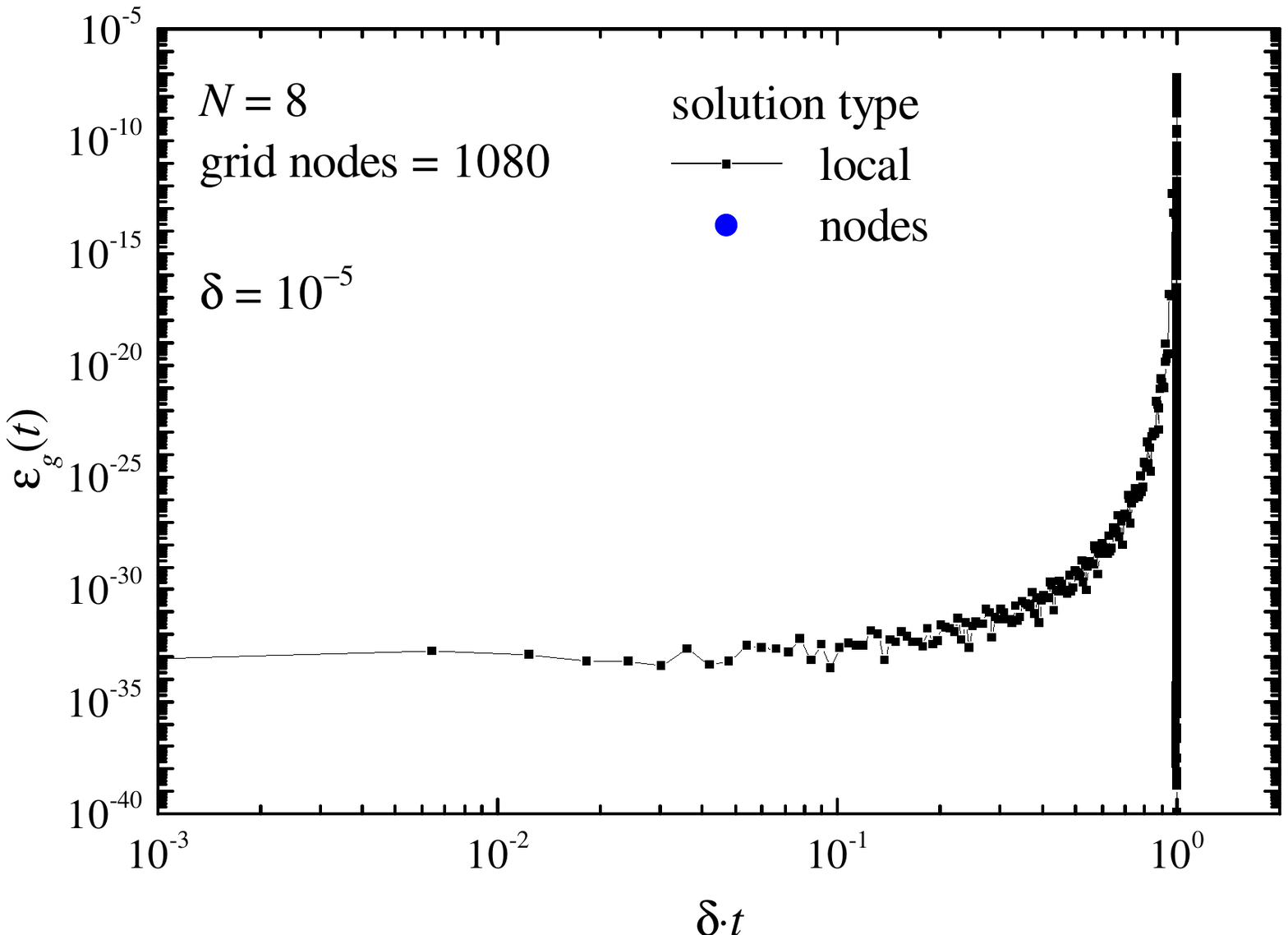}
\vspace{-8mm}\caption{\label{fig:shampine_flame_delta_10m5_sol_v_epss:e2}}
\end{subfigure}
\begin{subfigure}{0.320\textwidth}
\includegraphics[width=\textwidth]{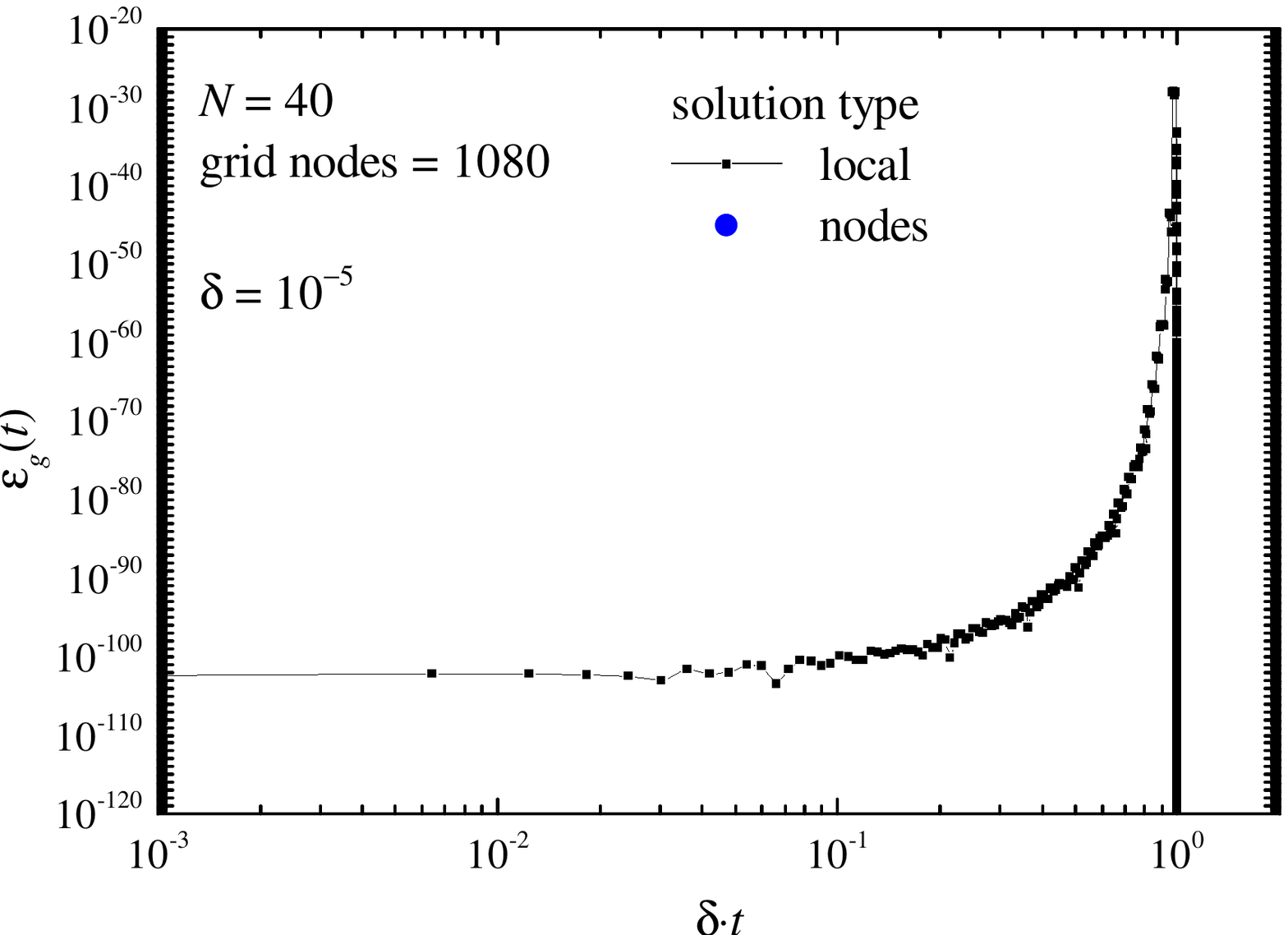}
\vspace{-8mm}\caption{\label{fig:shampine_flame_delta_10m5_sol_v_epss:e3}}
\end{subfigure}\\
\caption{%
Numerical solution of the stiff DAE system (\ref{eq:shampine_flame}) of index 1 with $\delta = 10^{-5}$. Comparison of the solution at nodes $\mathbf{v}_{n}$, the local solution $\mathbf{v}_{L}(t)$ and the exact solution $\mathbf{v}^{\rm ex}(t)$ for component $v_{1}$ (\subref{fig:shampine_flame_delta_10m5_sol_v_epss:a1}, \subref{fig:shampine_flame_delta_10m5_sol_v_epss:a2}, \subref{fig:shampine_flame_delta_10m5_sol_v_epss:a3}, \subref{fig:shampine_flame_delta_10m5_sol_v_epss:b1}, \subref{fig:shampine_flame_delta_10m5_sol_v_epss:b2}, \subref{fig:shampine_flame_delta_10m5_sol_v_epss:b3}), the errors $\varepsilon_{u}(t)$ (\subref{fig:shampine_flame_delta_10m5_sol_v_epss:c1}, \subref{fig:shampine_flame_delta_10m5_sol_v_epss:c2}, \subref{fig:shampine_flame_delta_10m5_sol_v_epss:c3}), $\varepsilon_{v}(t)$ (\subref{fig:shampine_flame_delta_10m5_sol_v_epss:d1}, \subref{fig:shampine_flame_delta_10m5_sol_v_epss:d2}, \subref{fig:shampine_flame_delta_10m5_sol_v_epss:d3}), $\varepsilon_{g}(t)$ (\subref{fig:shampine_flame_delta_10m5_sol_qug:e1}, \subref{fig:shampine_flame_delta_10m5_sol_qug:e2}, \subref{fig:shampine_flame_delta_10m5_sol_qug:e3}), obtained using polynomials with degrees $N = 1$ (\subref{fig:shampine_flame_delta_10m5_sol_qug:a1}, \subref{fig:shampine_flame_delta_10m5_sol_qug:b1}, \subref{fig:shampine_flame_delta_10m5_sol_qug:c1}, \subref{fig:shampine_flame_delta_10m5_sol_qug:d1}, \subref{fig:shampine_flame_delta_10m5_sol_qug:e1}), $N = 8$ (\subref{fig:shampine_flame_delta_10m5_sol_qug:a2}, \subref{fig:shampine_flame_delta_10m5_sol_qug:b2}, \subref{fig:shampine_flame_delta_10m5_sol_qug:c2}, \subref{fig:shampine_flame_delta_10m5_sol_qug:d2}, \subref{fig:shampine_flame_delta_10m5_sol_qug:e2}) and $N = 40$ (\subref{fig:shampine_flame_delta_10m5_sol_qug:a3}, \subref{fig:shampine_flame_delta_10m5_sol_qug:b3}, \subref{fig:shampine_flame_delta_10m5_sol_qug:c3}, \subref{fig:shampine_flame_delta_10m5_sol_qug:d3}, \subref{fig:shampine_flame_delta_10m5_sol_qug:e3}).
}
\label{fig:shampine_flame_delta_10m5_sol_v_epss}
\end{figure}

\begin{figure}[h!]
\captionsetup[subfigure]{%
	position=bottom,
	font+=smaller,
	textfont=normalfont,
	singlelinecheck=off,
	justification=raggedright
}
\centering
\begin{subfigure}{0.320\textwidth}
\includegraphics[width=\textwidth]{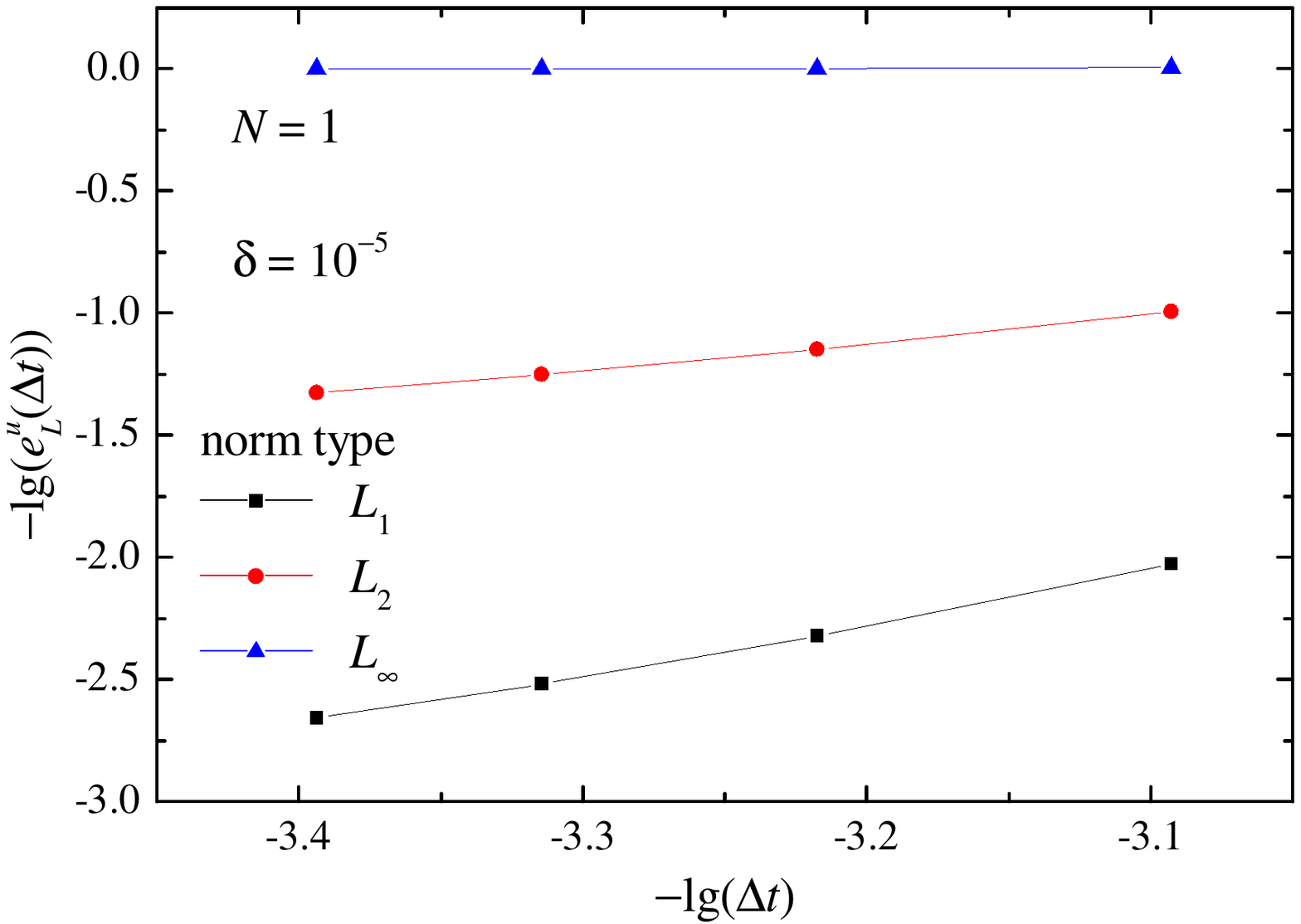}
\vspace{-8mm}\caption{\label{fig:shampine_flame_delta_10m5_errors:a1}}
\end{subfigure}
\begin{subfigure}{0.320\textwidth}
\includegraphics[width=\textwidth]{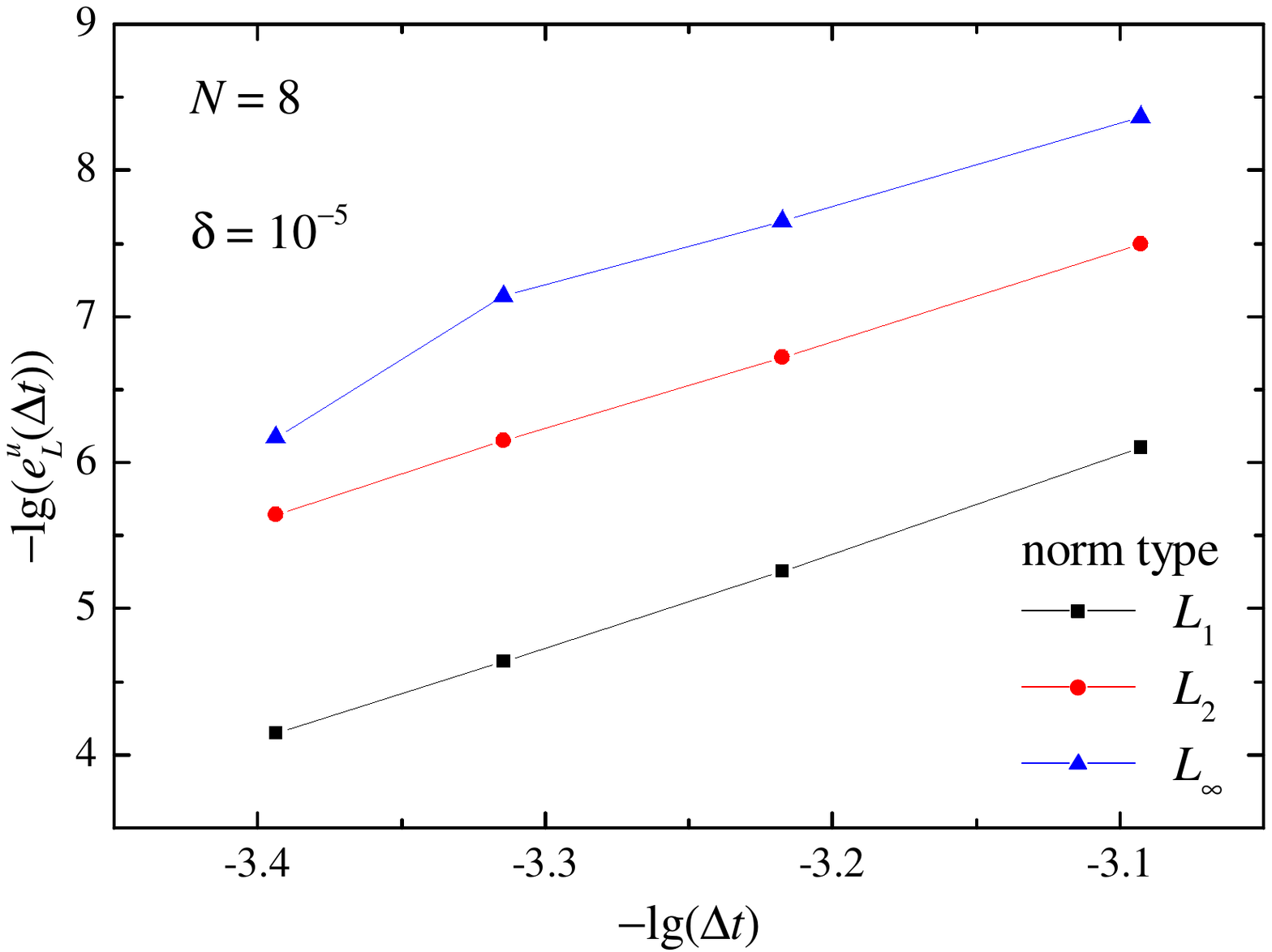}
\vspace{-8mm}\caption{\label{fig:shampine_flame_delta_10m5_errors:a2}}
\end{subfigure}
\begin{subfigure}{0.320\textwidth}
\includegraphics[width=\textwidth]{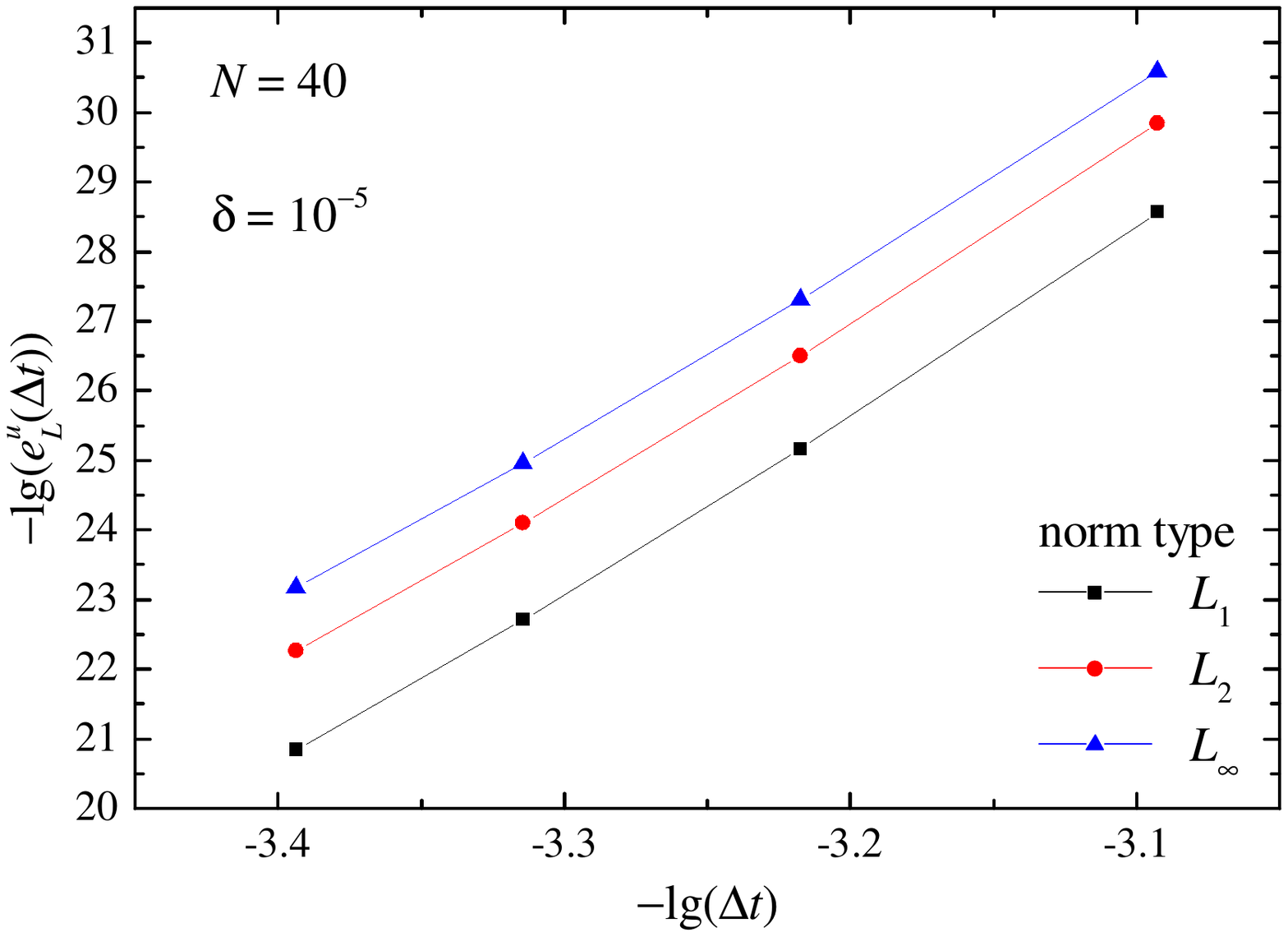}
\vspace{-8mm}\caption{\label{fig:shampine_flame_delta_10m5_errors:a3}}
\end{subfigure}\\
\begin{subfigure}{0.320\textwidth}
\includegraphics[width=\textwidth]{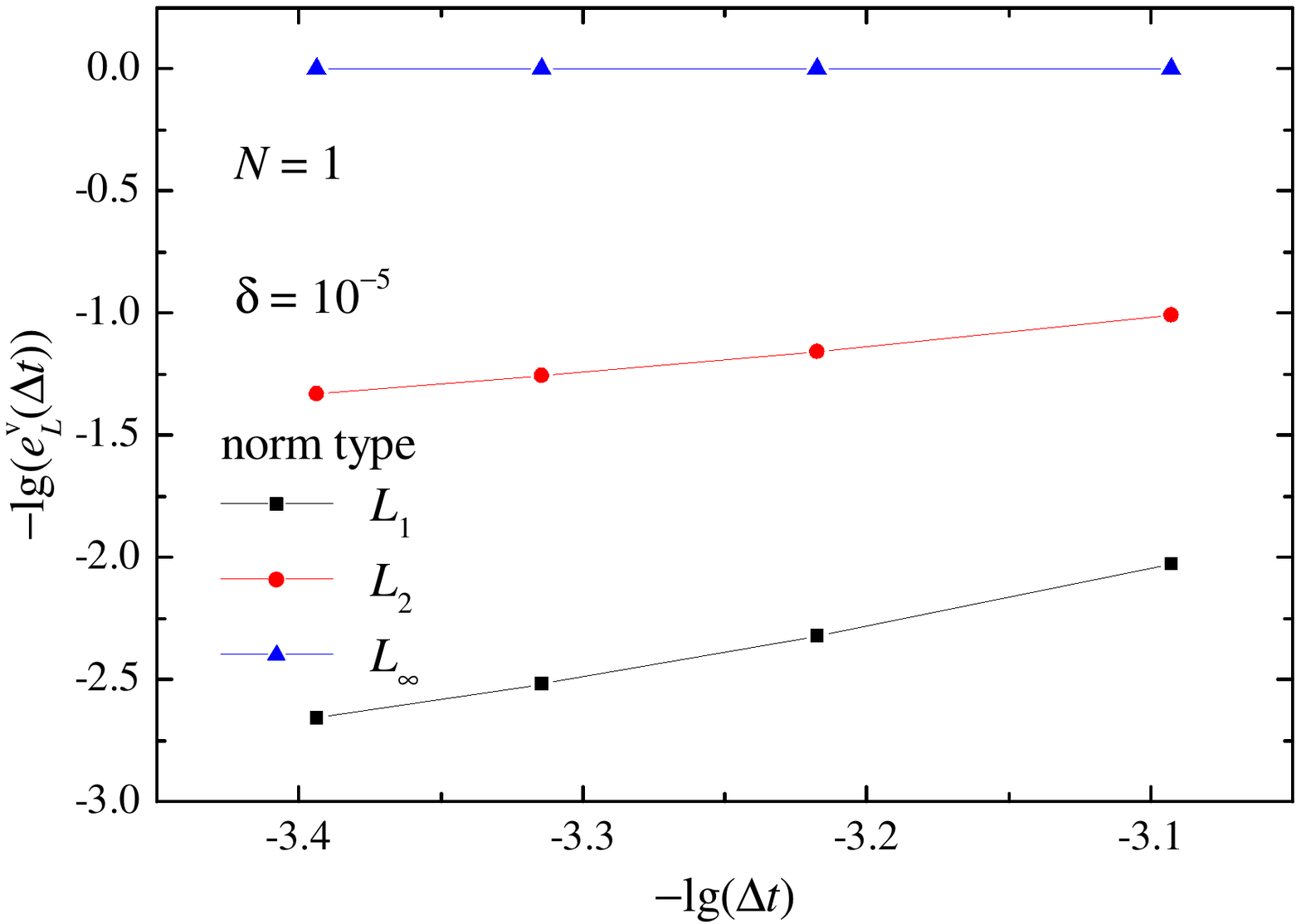}
\vspace{-8mm}\caption{\label{fig:shampine_flame_delta_10m5_errors:b1}}
\end{subfigure}
\begin{subfigure}{0.320\textwidth}
\includegraphics[width=\textwidth]{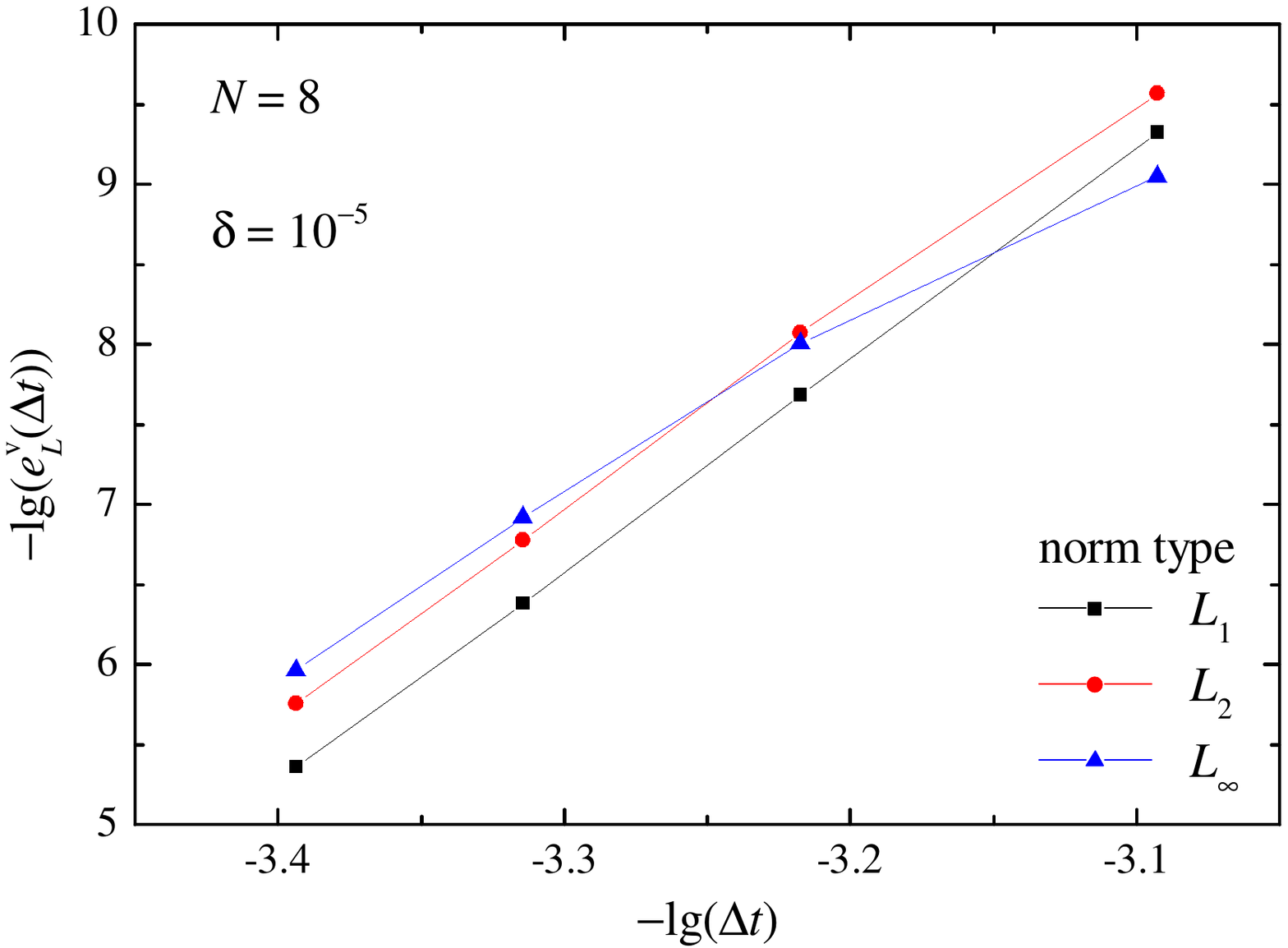}
\vspace{-8mm}\caption{\label{fig:shampine_flame_delta_10m5_errors:b2}}
\end{subfigure}
\begin{subfigure}{0.320\textwidth}
\includegraphics[width=\textwidth]{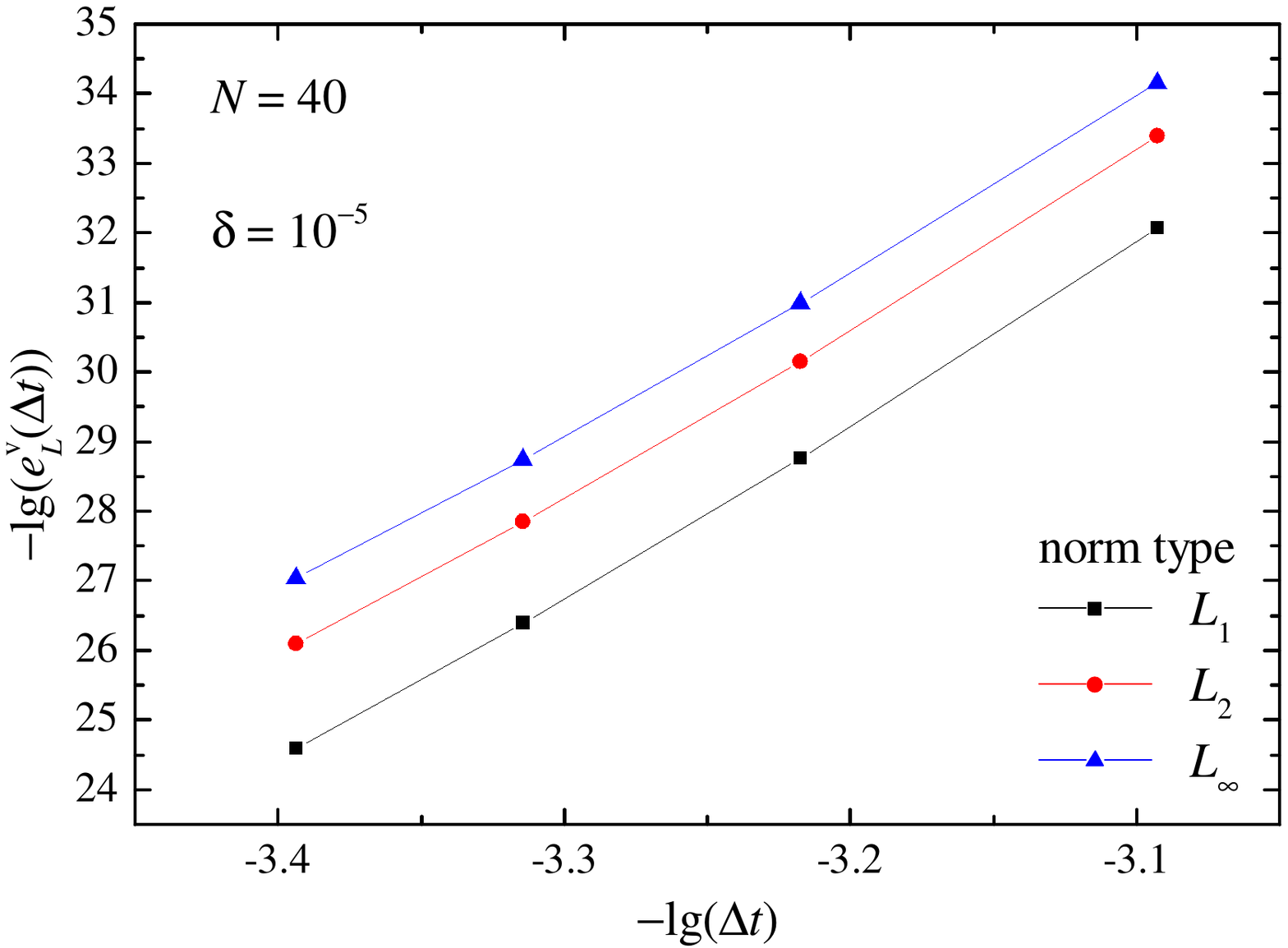}
\vspace{-8mm}\caption{\label{fig:shampine_flame_delta_10m5_errors:b3}}
\end{subfigure}\\
\begin{subfigure}{0.320\textwidth}
\includegraphics[width=\textwidth]{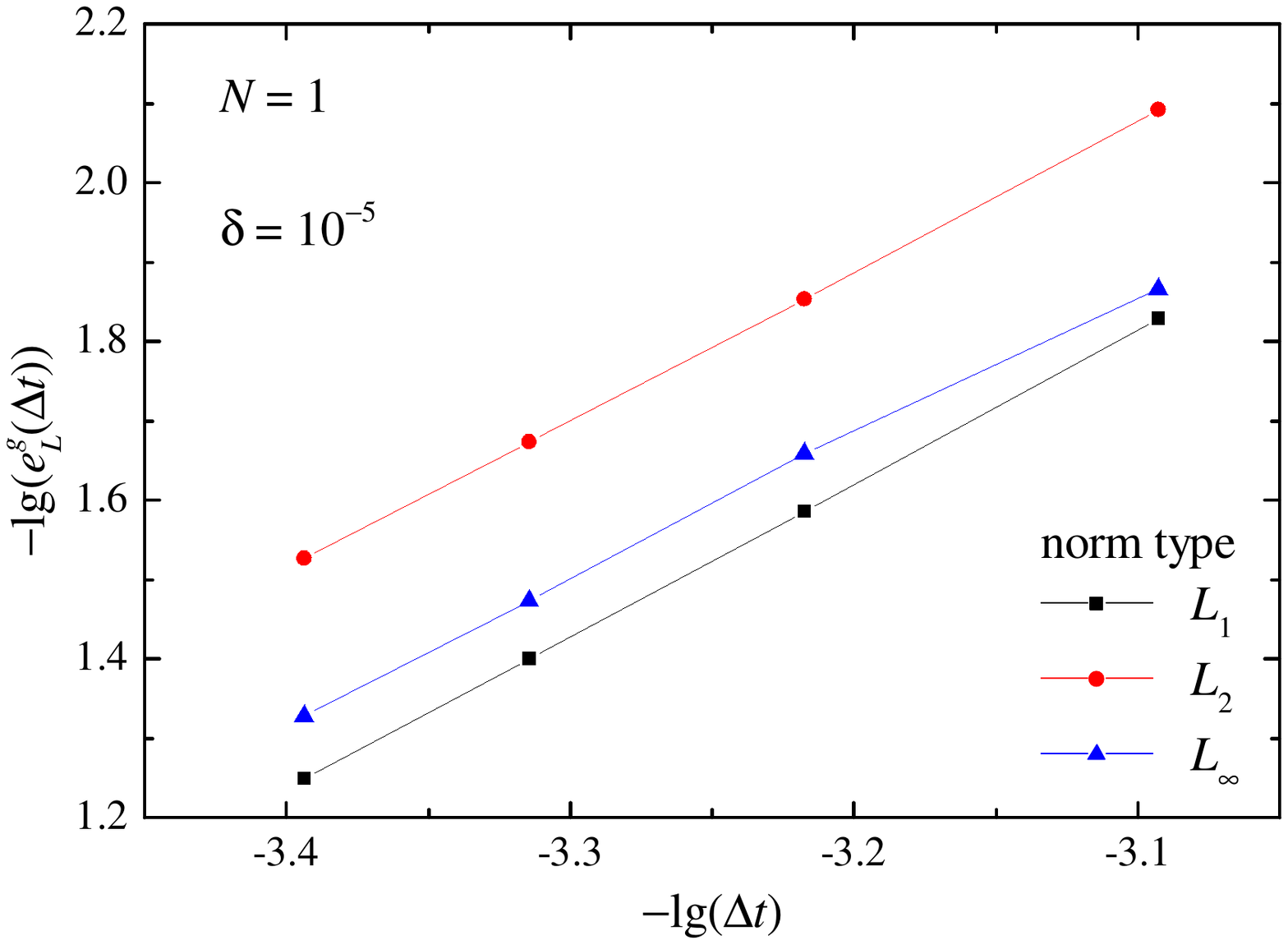}
\vspace{-8mm}\caption{\label{fig:shampine_flame_delta_10m5_errors:c1}}
\end{subfigure}
\begin{subfigure}{0.320\textwidth}
\includegraphics[width=\textwidth]{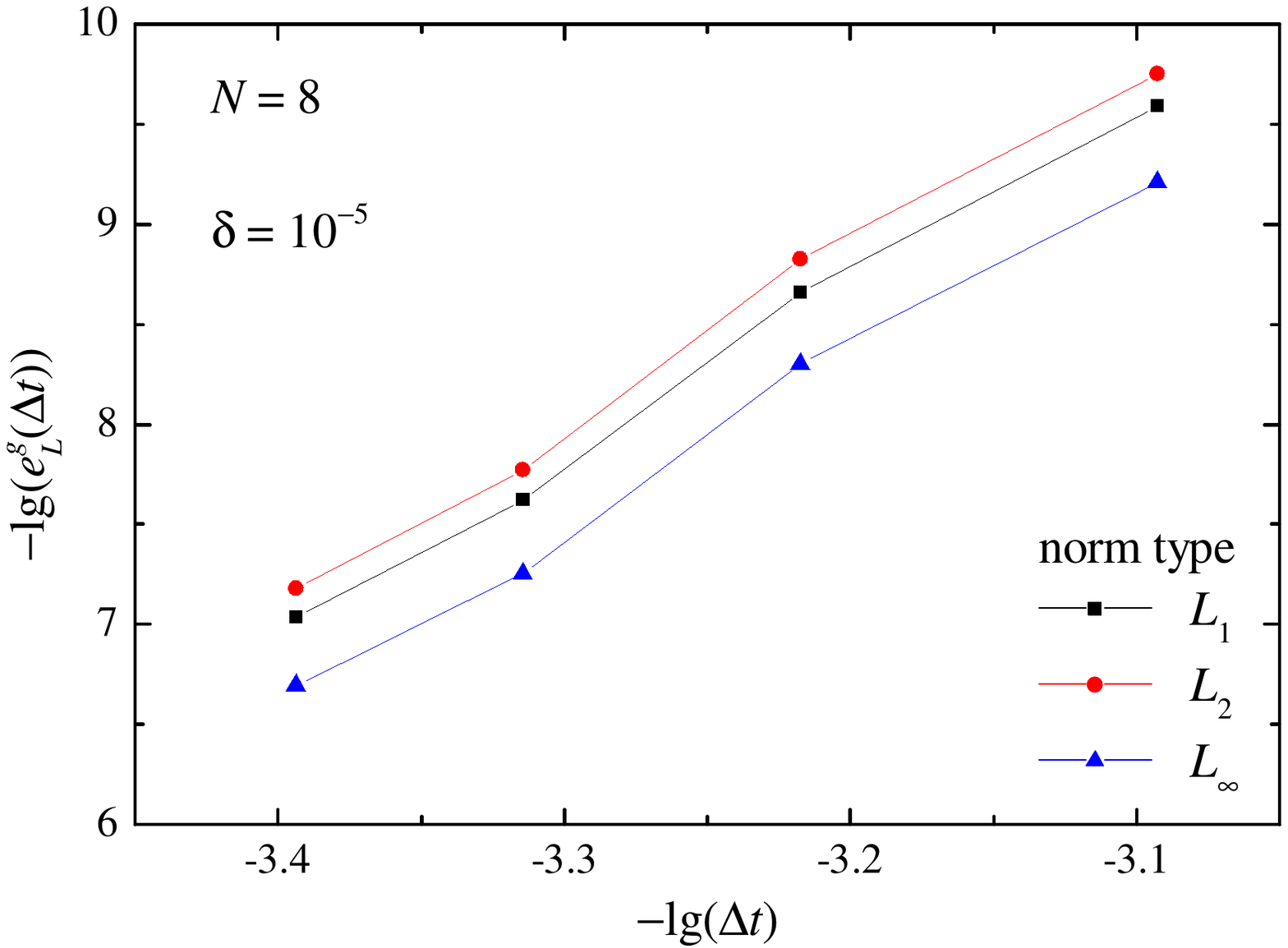}
\vspace{-8mm}\caption{\label{fig:shampine_flame_delta_10m5_errors:c2}}
\end{subfigure}
\begin{subfigure}{0.320\textwidth}
\includegraphics[width=\textwidth]{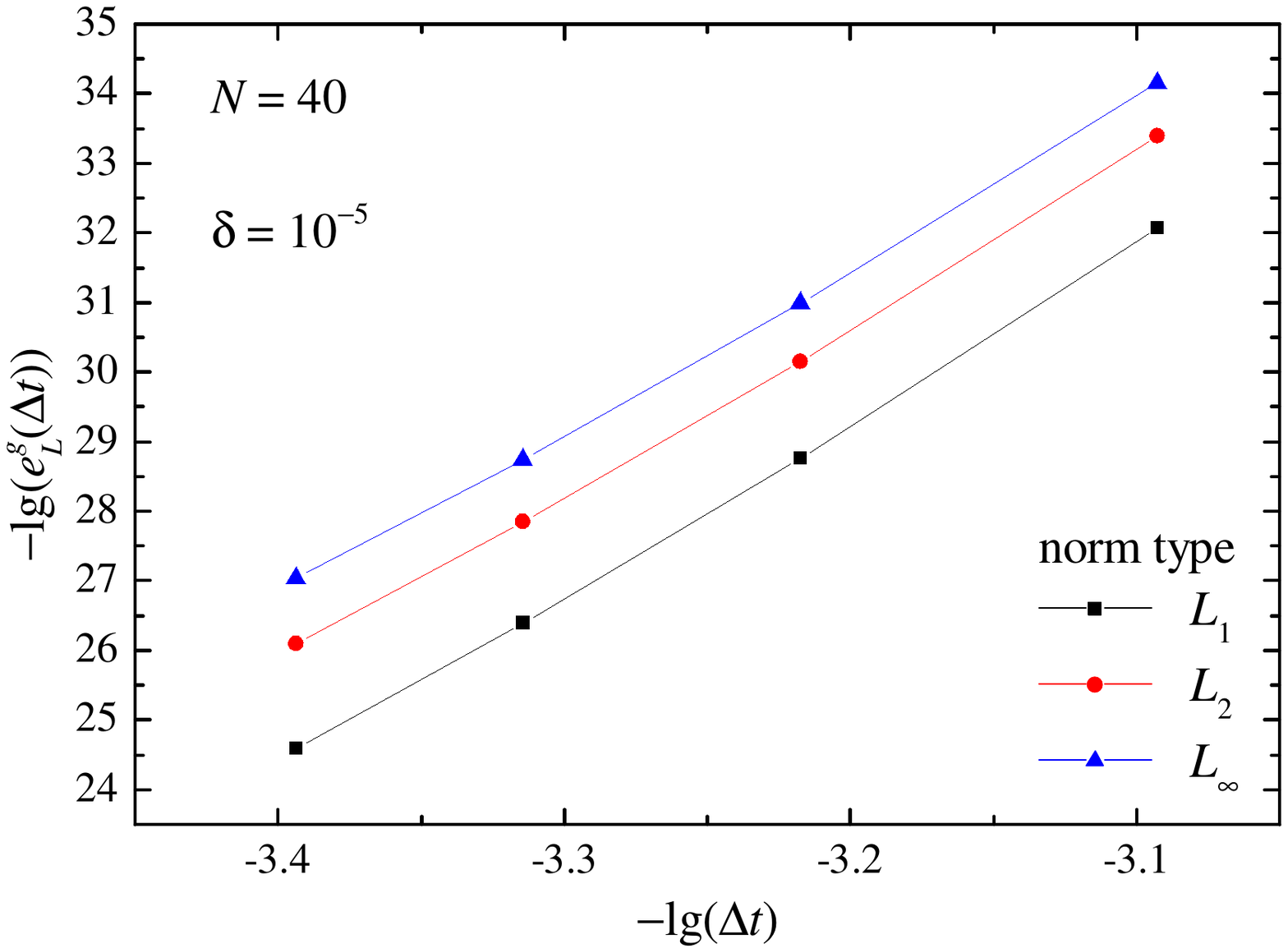}
\vspace{-8mm}\caption{\label{fig:shampine_flame_delta_10m5_errors:c3}}
\end{subfigure}\\
\begin{subfigure}{0.320\textwidth}
\includegraphics[width=\textwidth]{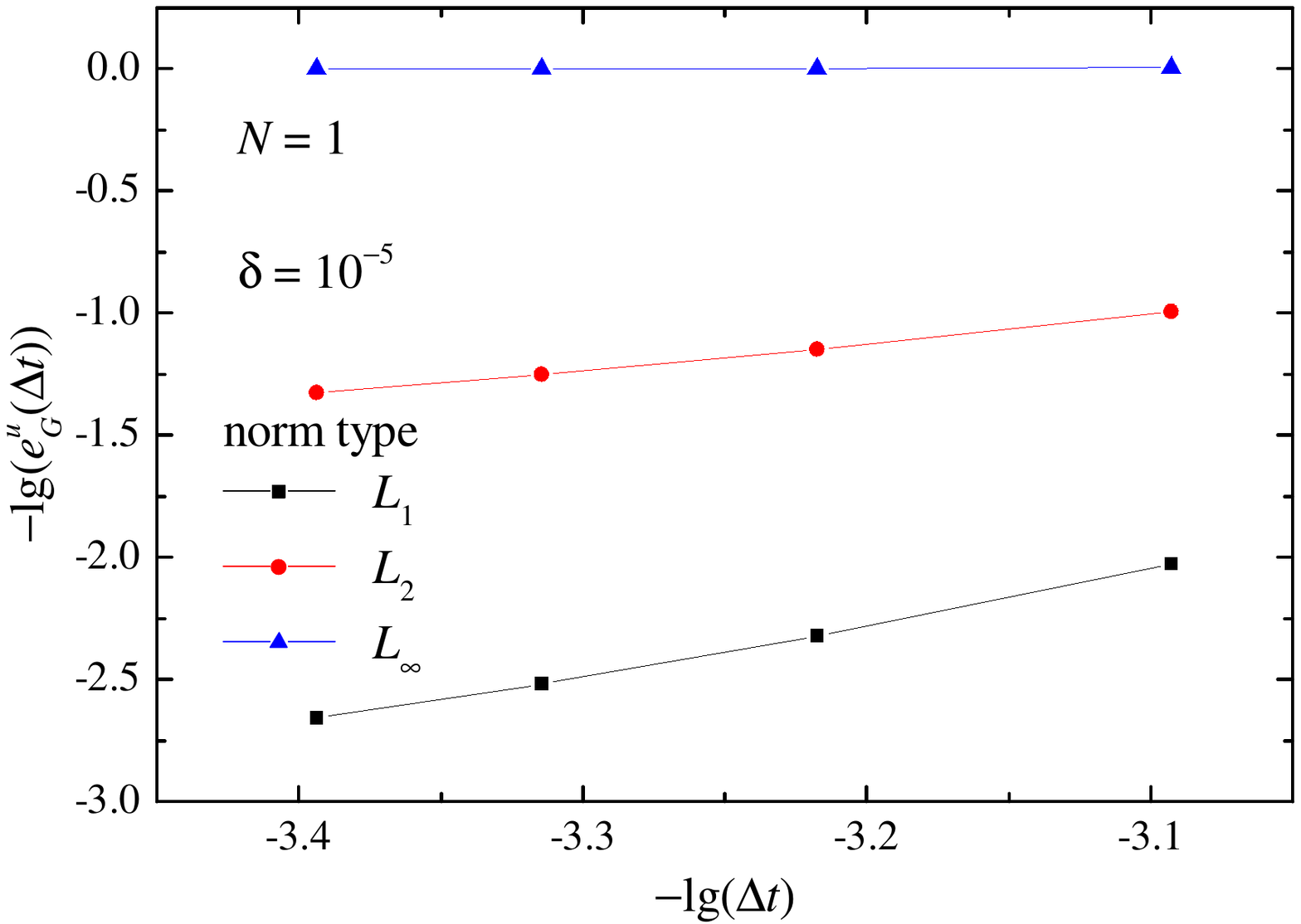}
\vspace{-8mm}\caption{\label{fig:shampine_flame_delta_10m5_errors:d1}}
\end{subfigure}
\begin{subfigure}{0.320\textwidth}
\includegraphics[width=\textwidth]{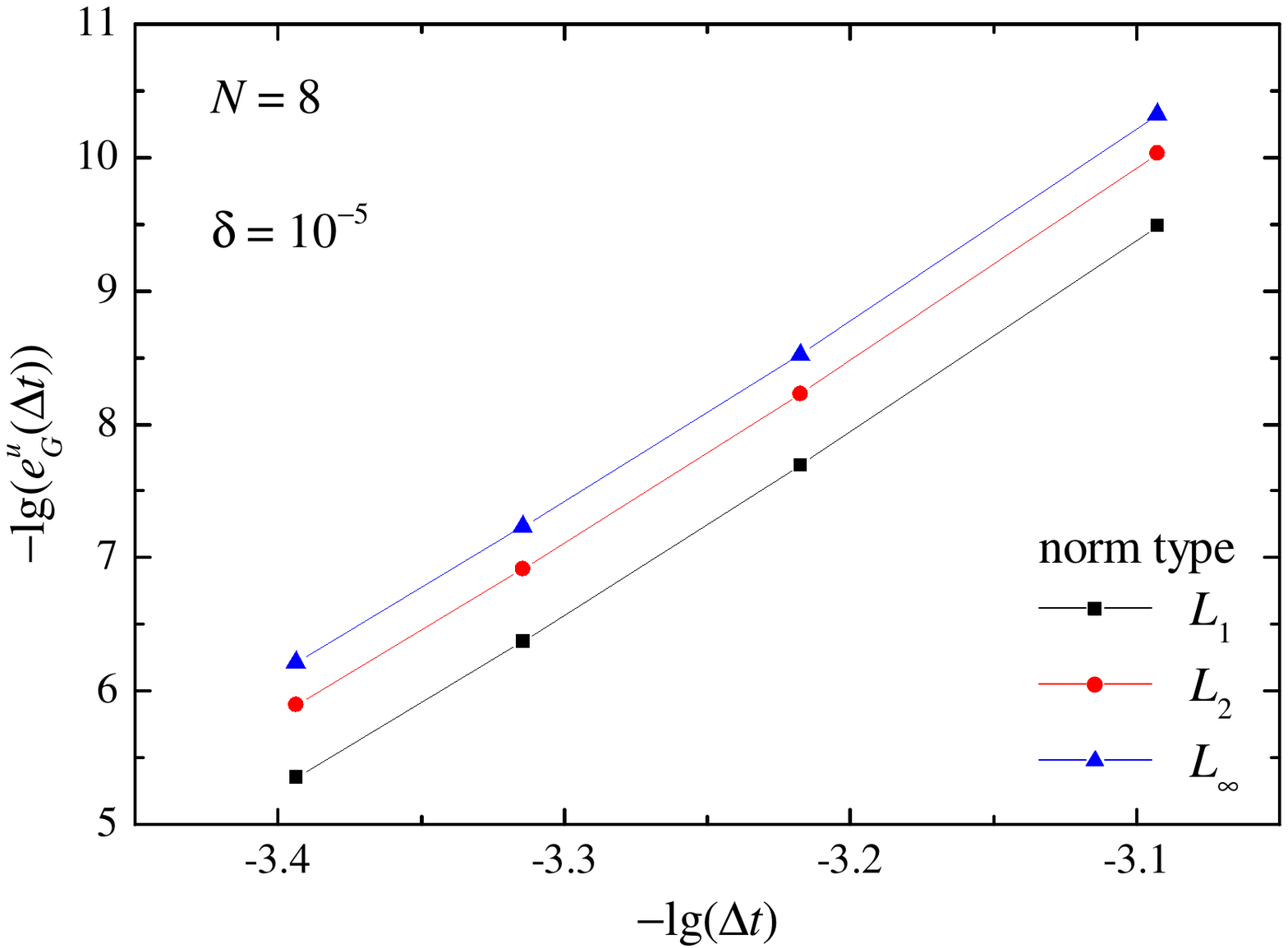}
\vspace{-8mm}\caption{\label{fig:shampine_flame_delta_10m5_errors:d2}}
\end{subfigure}
\begin{subfigure}{0.320\textwidth}
\includegraphics[width=\textwidth]{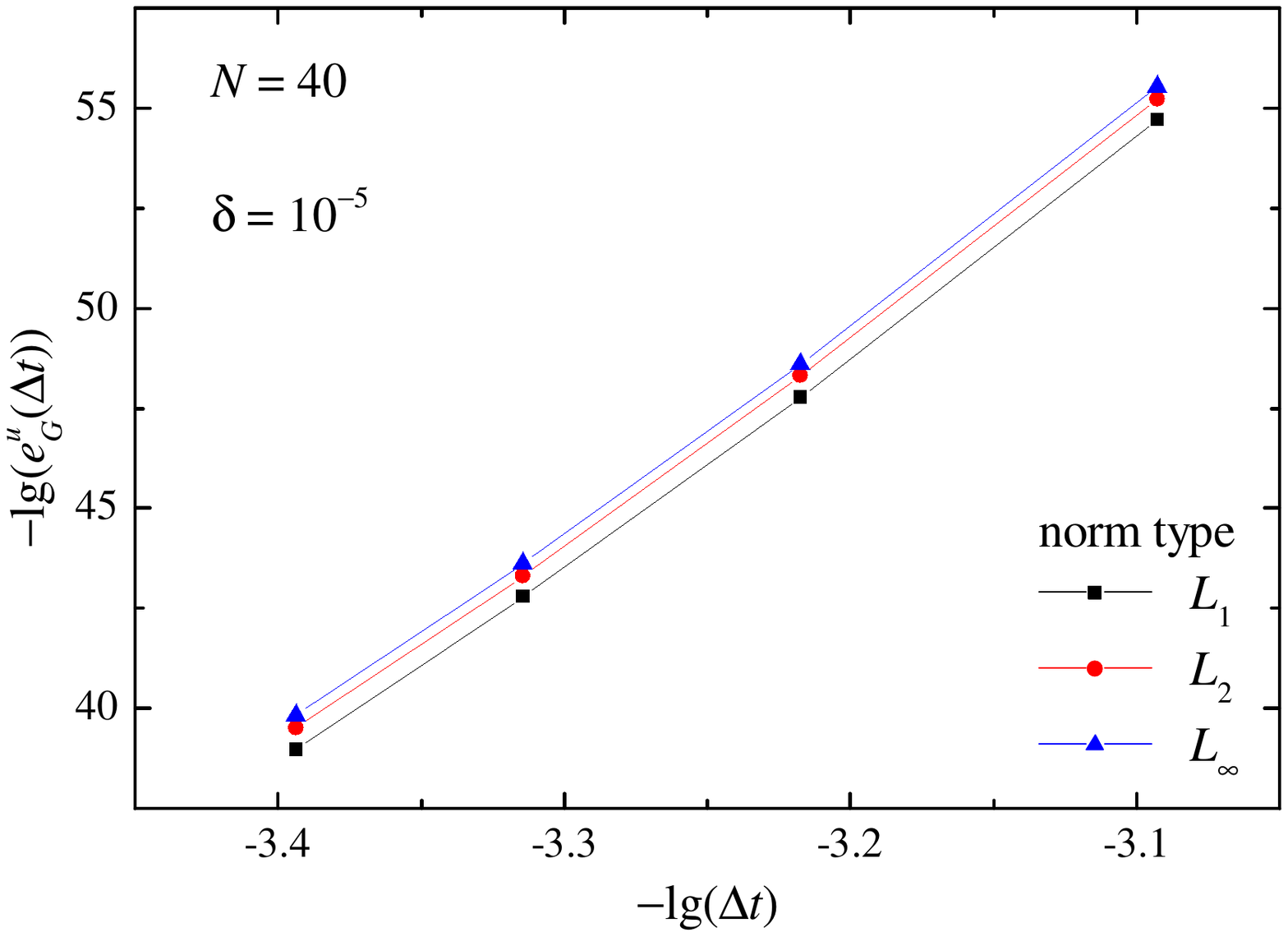}
\vspace{-8mm}\caption{\label{fig:shampine_flame_delta_10m5_errors:d3}}
\end{subfigure}\\
\begin{subfigure}{0.320\textwidth}
\includegraphics[width=\textwidth]{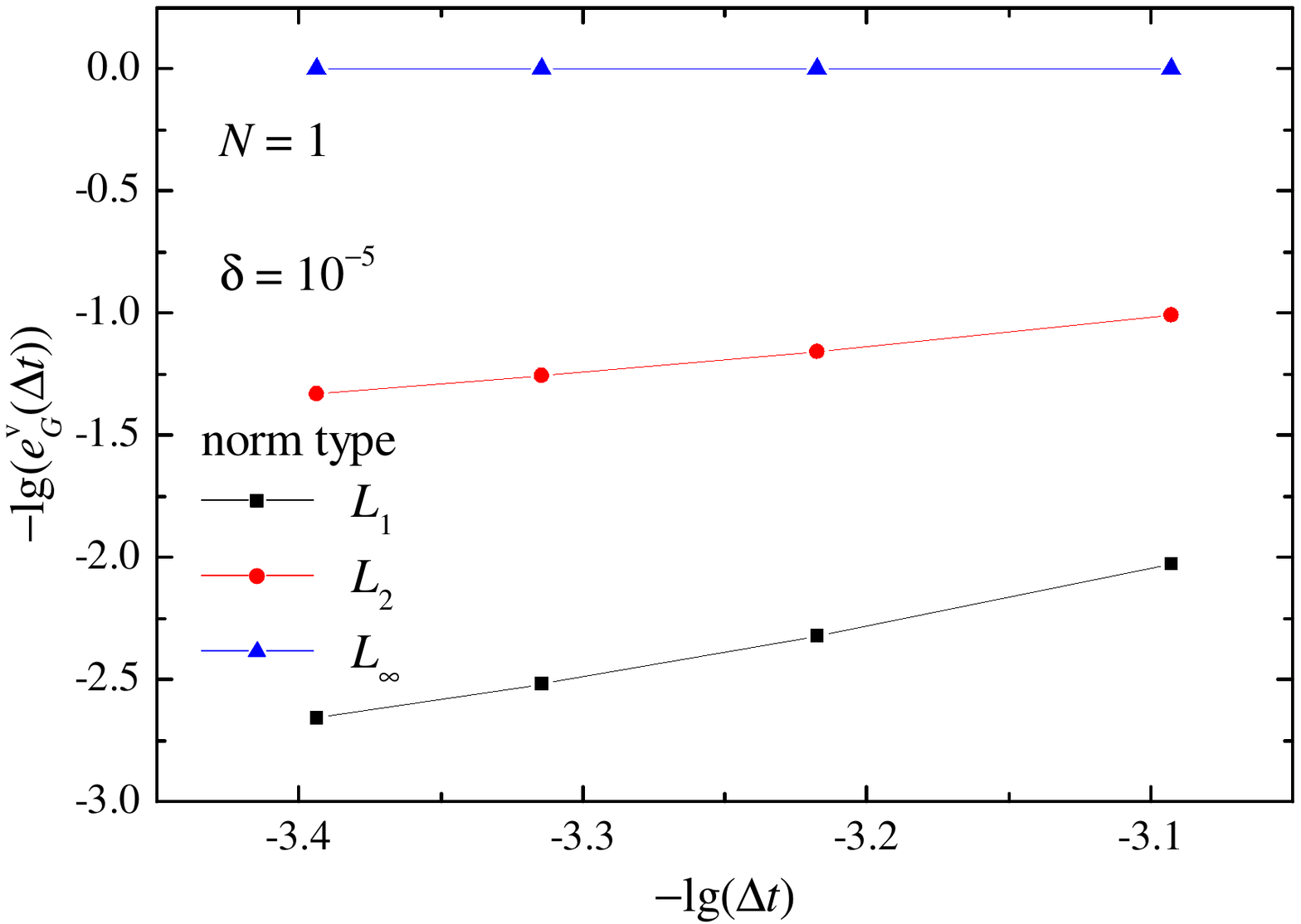}
\vspace{-8mm}\caption{\label{fig:shampine_flame_delta_10m5_errors:e1}}
\end{subfigure}
\begin{subfigure}{0.320\textwidth}
\includegraphics[width=\textwidth]{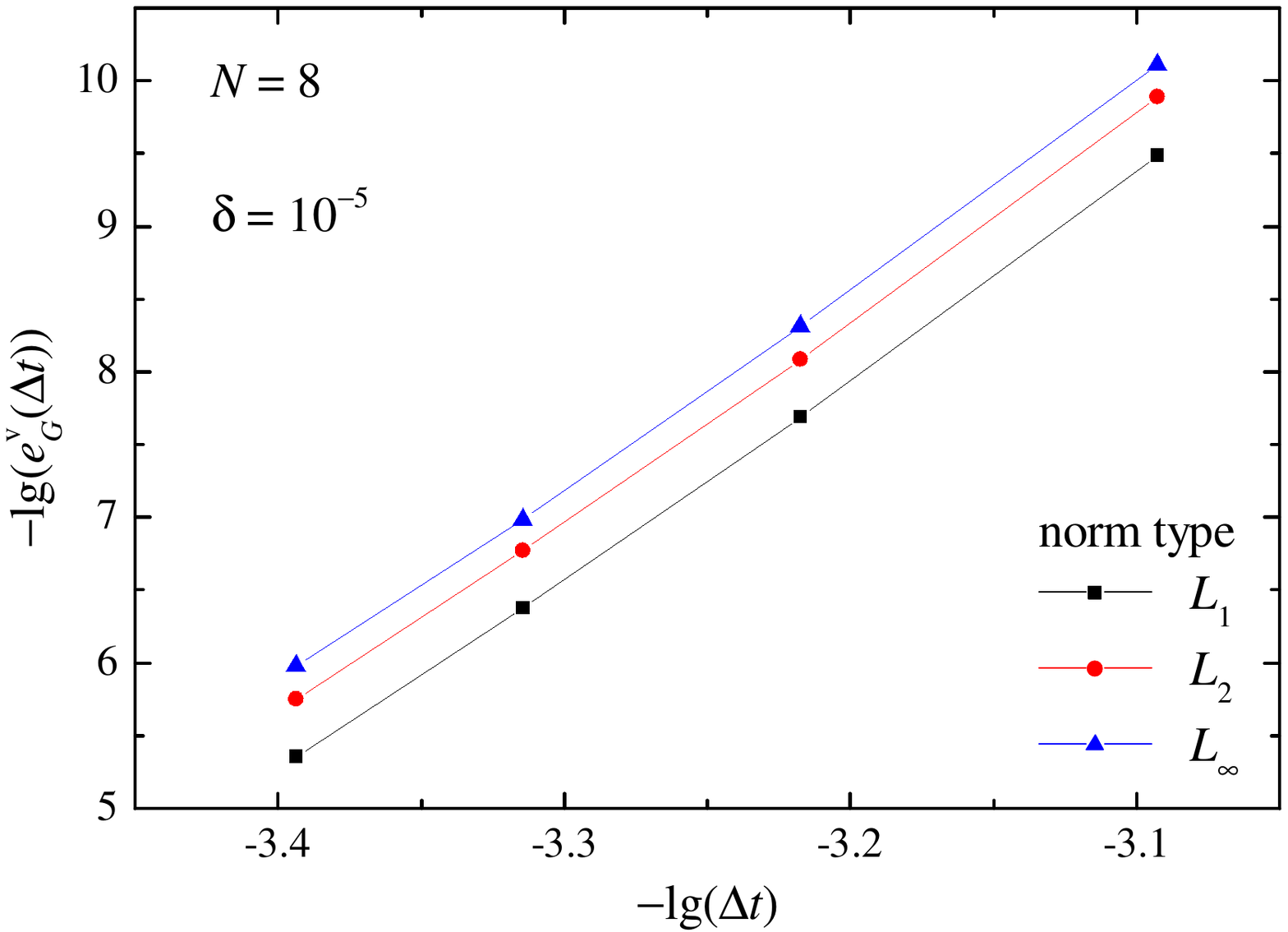}
\vspace{-8mm}\caption{\label{fig:shampine_flame_delta_10m5_errors:e2}}
\end{subfigure}
\begin{subfigure}{0.320\textwidth}
\includegraphics[width=\textwidth]{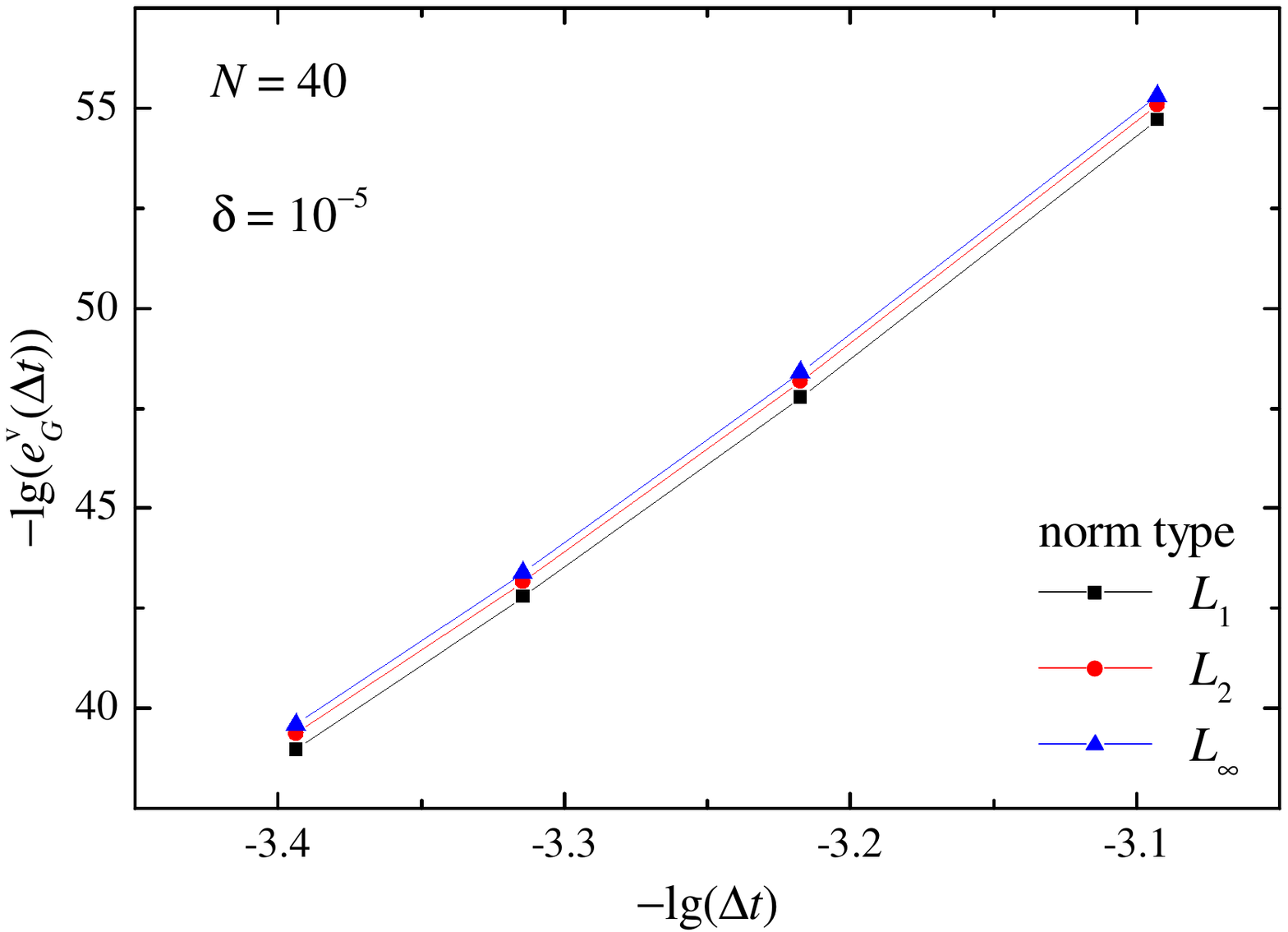}
\vspace{-8mm}\caption{\label{fig:shampine_flame_delta_10m5_errors:e3}}
\end{subfigure}\\
\caption{%
Log-log plot of the dependence of the global errors for the local solution $e_{L}^{u}$ (\subref{fig:shampine_flame_delta_10m5_errors:a1}, \subref{fig:shampine_flame_delta_10m5_errors:a2}, \subref{fig:shampine_flame_delta_10m5_errors:a3}), $e_{L}^{v}$ (\subref{fig:shampine_flame_delta_10m5_errors:b1}, \subref{fig:shampine_flame_delta_10m5_errors:b2}, \subref{fig:shampine_flame_delta_10m5_errors:b3}), $e_{L}^{g}$ (\subref{fig:shampine_flame_delta_10m5_errors:c1}, \subref{fig:shampine_flame_delta_10m5_errors:c2}, \subref{fig:shampine_flame_delta_10m5_errors:c3}) and the solution at nodes $e_{G}^{u}$ (\subref{fig:shampine_flame_delta_10m5_errors:d1}, \subref{fig:shampine_flame_delta_10m5_errors:d2}, \subref{fig:shampine_flame_delta_10m5_errors:d3}), $e_{G}^{v}$ (\subref{fig:shampine_flame_delta_10m5_errors:e1}, \subref{fig:shampine_flame_delta_10m5_errors:e2}, \subref{fig:shampine_flame_delta_10m5_errors:e3}) on the discretization step $\mathrm{\Delta}t$, obtained in the norms $L_{1}$, $L_{2}$ and $L_{\infty}$, by numerical solution of the stiff DAE system (\ref{eq:shampine_flame}) of index 1 with $\delta = 10^{-5}$, obtained using polynomials with degrees $N = 1$ (\subref{fig:shampine_flame_delta_10m5_errors:a1}, \subref{fig:shampine_flame_delta_10m5_errors:b1}, \subref{fig:shampine_flame_delta_10m5_errors:c1}, \subref{fig:shampine_flame_delta_10m5_errors:d1}, \subref{fig:shampine_flame_delta_10m5_errors:e1}), $N = 8$ (\subref{fig:shampine_flame_delta_10m5_errors:a2}, \subref{fig:shampine_flame_delta_10m5_errors:b2}, \subref{fig:shampine_flame_delta_10m5_errors:c2}, \subref{fig:shampine_flame_delta_10m5_errors:d2}, \subref{fig:shampine_flame_delta_10m5_errors:e2}) and $N = 40$ (\subref{fig:shampine_flame_delta_10m5_errors:a3}, \subref{fig:shampine_flame_delta_10m5_errors:b3}, \subref{fig:shampine_flame_delta_10m5_errors:c3}, \subref{fig:shampine_flame_delta_10m5_errors:d3}, \subref{fig:shampine_flame_delta_10m5_errors:e3}).
}
\label{fig:shampine_flame_delta_10m5_errors}
\end{figure}

\begin{table*}[h!]
\centering
\caption{%
Convergence orders $p_{L_{1}}$, $p_{L_{2}}$, $p_{L_{\infty}}$, calculated in norms $L_{1}$, $L_{2}$, $L_{\infty}$ of the ADER-DG method for the stiff DAE system (\ref{eq:shampine_flame}) of index 1 with $\delta = 10^{-5}$; $N$ is the degree of the basis polynomials $\varphi_{p}$. Orders $p^{n, u}$ are calculated for solution at the nodes $\mathbf{u}_{n}$; orders $p^{n, v}$ --- for solution at the nodes $\mathbf{v}_{n}$; orders $p^{l, u}$ --- for local solution $\mathbf{u}_{L}$; orders $p^{l, v}$ --- for local solution $\mathbf{v}_{L}$. The theoretical values of convergence order $p_{\rm th.}^{n} = 2N+1$ and $p_{\rm th.}^{l} = N+1$ are applicable for the ADER-DG method for ODE problems, and are presented for comparison.
}
\label{tab:conv_ords_shampine_flame_delta_10m5}
\begingroup
\setlength{\tabcolsep}{4pt}
\begin{adjustbox}{angle=90}
\begin{tabular}{@{}|l|lll|lll|c|lll|lll|c|@{}}
\toprule
$N$\hspace{-2mm} & $p_{L_{1}}^{n, u}$ & $p_{L_{2}}^{n, u}$ & $p_{L_{\infty}}^{n, u}$ & $p_{L_{1}}^{n, v}$ & $p_{L_{2}}^{n, v}$ & $p_{L_{\infty}}^{n, v}$ & $p_{\rm th.}^{n}$ & $p_{L_{1}}^{l, u}$ & $p_{L_{2}}^{l, u}$ & $p_{L_{\infty}}^{l, u}$ & $p_{L_{1}}^{l, v}$ & $p_{L_{2}}^{l, v}$ & $p_{L_{\infty}}^{l, v}$ & $p_{\rm th.}^{l}$ \\
\midrule
$1$	&	$2.10$	&	$1.10$	&	$0.01$	&	$2.10$	&	$1.06$	&	$0.00$	&	$3$	&	$2.10$	&	$1.10$	&	$0.01$	&	$2.10$	&	$1.06$	&	$0.00$	&	$2$\\
$2$	&	$7.02$	&	$6.21$	&	$4.91$	&	$7.02$	&	$5.86$	&	$4.26$	&	$5$	&	$6.94$	&	$6.21$	&	$4.86$	&	$7.02$	&	$5.86$	&	$4.17$	&	$3$\\
$3$	&	$7.95$	&	$7.94$	&	$7.91$	&	$7.95$	&	$7.92$	&	$7.75$	&	$7$	&	$7.44$	&	$7.94$	&	$7.91$	&	$7.95$	&	$7.93$	&	$7.88$	&	$4$\\
$4$	&	$9.06$	&	$9.06$	&	$8.98$	&	$9.06$	&	$9.07$	&	$9.07$	&	$9$	&	$7.26$	&	$8.98$	&	$8.92$	&	$9.06$	&	$9.06$	&	$8.81$	&	$5$\\
$5$	&	$10.21$	&	$10.21$	&	$10.12$	&	$10.21$	&	$10.22$	&	$10.23$	&	$11$	&	$6.52$	&	$9.06$	&	$9.87$	&	$10.21$	&	$10.20$	&	$9.57$	&	$6$\\
$6$	&	$11.38$	&	$11.38$	&	$11.28$	&	$11.38$	&	$11.39$	&	$11.39$	&	$13$	&	$6.03$	&	$7.83$	&	$9.93$	&	$11.36$	&	$11.30$	&	$10.54$	&	$7$\\
$7$	&	$12.56$	&	$12.55$	&	$12.46$	&	$12.56$	&	$12.56$	&	$12.57$	&	$15$	&	$6.10$	&	$6.55$	&	$8.56$	&	$12.52$	&	$12.44$	&	$10.98$	&	$8$\\
$8$	&	$13.74$	&	$13.74$	&	$13.65$	&	$13.74$	&	$13.75$	&	$13.76$	&	$17$	&	$6.49$	&	$6.12$	&	$6.96$	&	$13.21$	&	$12.70$	&	$10.25$	&	$9$\\
$9$	&	$14.94$	&	$14.93$	&	$14.84$	&	$14.94$	&	$14.94$	&	$14.95$	&	$19$	&	$7.04$	&	$6.51$	&	$5.99$	&	$14.18$	&	$13.86$	&	$11.45$	&	$10$\\
$10$	&	$16.13$	&	$16.13$	&	$16.04$	&	$16.13$	&	$16.14$	&	$16.14$	&	$21$	&	$7.61$	&	$7.09$	&	$6.56$	&	$13.54$	&	$13.18$	&	$11.03$	&	$11$\\
$11$	&	$17.33$	&	$17.32$	&	$17.23$	&	$17.33$	&	$17.33$	&	$17.34$	&	$23$	&	$8.21$	&	$7.69$	&	$7.16$	&	$13.08$	&	$13.62$	&	$11.57$	&	$12$\\
$12$	&	$18.53$	&	$18.52$	&	$18.43$	&	$18.53$	&	$18.53$	&	$18.54$	&	$25$	&	$8.82$	&	$8.29$	&	$7.76$	&	$12.15$	&	$13.77$	&	$12.37$	&	$13$\\
$13$	&	$19.73$	&	$19.72$	&	$19.63$	&	$19.73$	&	$19.74$	&	$19.74$	&	$27$	&	$9.43$	&	$8.89$	&	$8.36$	&	$10.95$	&	$13.44$	&	$11.92$	&	$14$\\
$14$	&	$20.93$	&	$20.93$	&	$20.84$	&	$20.93$	&	$20.94$	&	$20.95$	&	$29$	&	$10.02$	&	$9.49$	&	$8.97$	&	$10.61$	&	$14.15$	&	$14.77$	&	$15$\\
$15$	&	$22.14$	&	$22.13$	&	$22.04$	&	$22.14$	&	$22.14$	&	$22.15$	&	$31$	&	$10.62$	&	$10.09$	&	$9.57$	&	$10.21$	&	$12.14$	&	$13.10$	&	$16$\\
$16$	&	$23.34$	&	$23.34$	&	$23.25$	&	$23.34$	&	$23.35$	&	$23.36$	&	$33$	&	$11.23$	&	$10.69$	&	$10.17$	&	$10.72$	&	$12.71$	&	$16.00$	&	$17$\\
$17$	&	$24.55$	&	$24.54$	&	$24.45$	&	$24.55$	&	$24.55$	&	$24.56$	&	$35$	&	$11.83$	&	$11.30$	&	$10.77$	&	$11.08$	&	$11.60$	&	$16.17$	&	$18$\\
$18$	&	$25.75$	&	$25.75$	&	$25.66$	&	$25.75$	&	$25.76$	&	$25.77$	&	$37$	&	$12.42$	&	$11.90$	&	$11.38$	&	$11.66$	&	$11.95$	&	$16.53$	&	$19$\\
$19$	&	$26.96$	&	$26.95$	&	$26.86$	&	$26.96$	&	$26.97$	&	$26.97$	&	$39$	&	$13.04$	&	$12.50$	&	$11.98$	&	$12.22$	&	$11.85$	&	$15.88$	&	$20$\\
$20$	&	$28.17$	&	$28.16$	&	$28.07$	&	$28.17$	&	$28.17$	&	$28.18$	&	$41$	&	$13.65$	&	$13.11$	&	$12.59$	&	$12.82$	&	$12.28$	&	$14.97$	&	$21$\\
\midrule
$25$	&	$34.21$	&	$34.20$	&	$34.11$	&	$34.21$	&	$34.21$	&	$34.22$	&	$51$	&	$16.65$	&	$16.12$	&	$15.61$	&	$15.80$	&	$15.22$	&	$14.64$	&	$26$\\
$30$	&	$40.25$	&	$40.25$	&	$40.16$	&	$40.25$	&	$40.26$	&	$40.27$	&	$61$	&	$19.68$	&	$19.15$	&	$18.63$	&	$18.85$	&	$18.25$	&	$17.67$	&	$31$\\
$35$	&	$46.30$	&	$46.30$	&	$46.21$	&	$46.30$	&	$46.31$	&	$46.32$	&	$71$	&	$22.74$	&	$22.17$	&	$21.66$	&	$21.86$	&	$21.26$	&	$20.69$	&	$36$\\
$40$	&	$52.35$	&	$52.35$	&	$52.26$	&	$52.35$	&	$52.36$	&	$52.37$	&	$81$	&	$25.75$	&	$25.20$	&	$24.68$	&	$24.87$	&	$24.27$	&	$23.72$	&	$41$\\
\bottomrule
\end{tabular}
\end{adjustbox}
\endgroup
\end{table*}

\begin{figure}[h!]
\captionsetup[subfigure]{%
	position=bottom,
	font+=smaller,
	textfont=normalfont,
	singlelinecheck=off,
	justification=raggedright
}
\centering
\begin{subfigure}{0.320\textwidth}
\includegraphics[width=\textwidth]{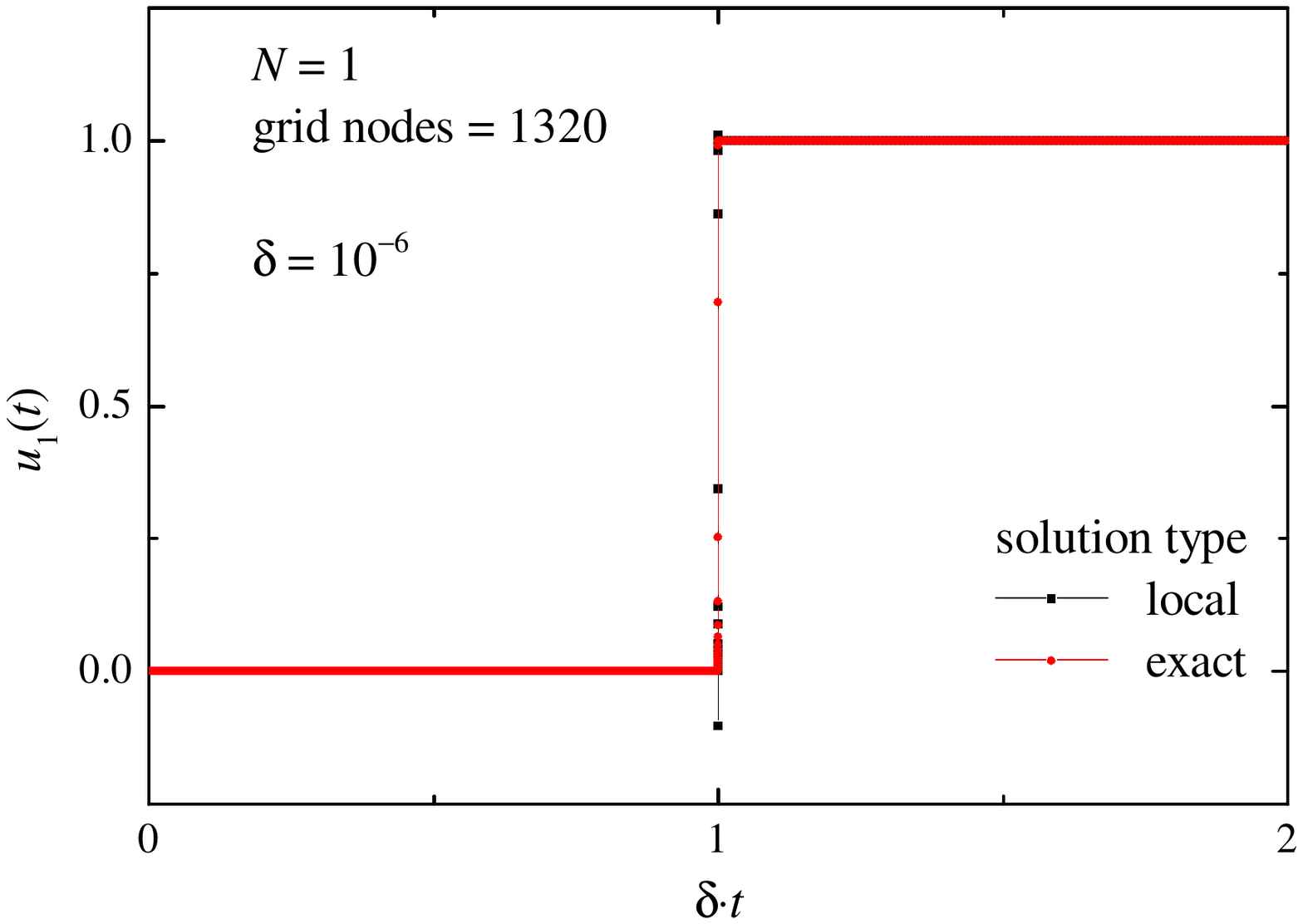}
\vspace{-8mm}\caption{\label{fig:shampine_flame_delta_10m6_sol_qug:a1}}
\end{subfigure}
\begin{subfigure}{0.320\textwidth}
\includegraphics[width=\textwidth]{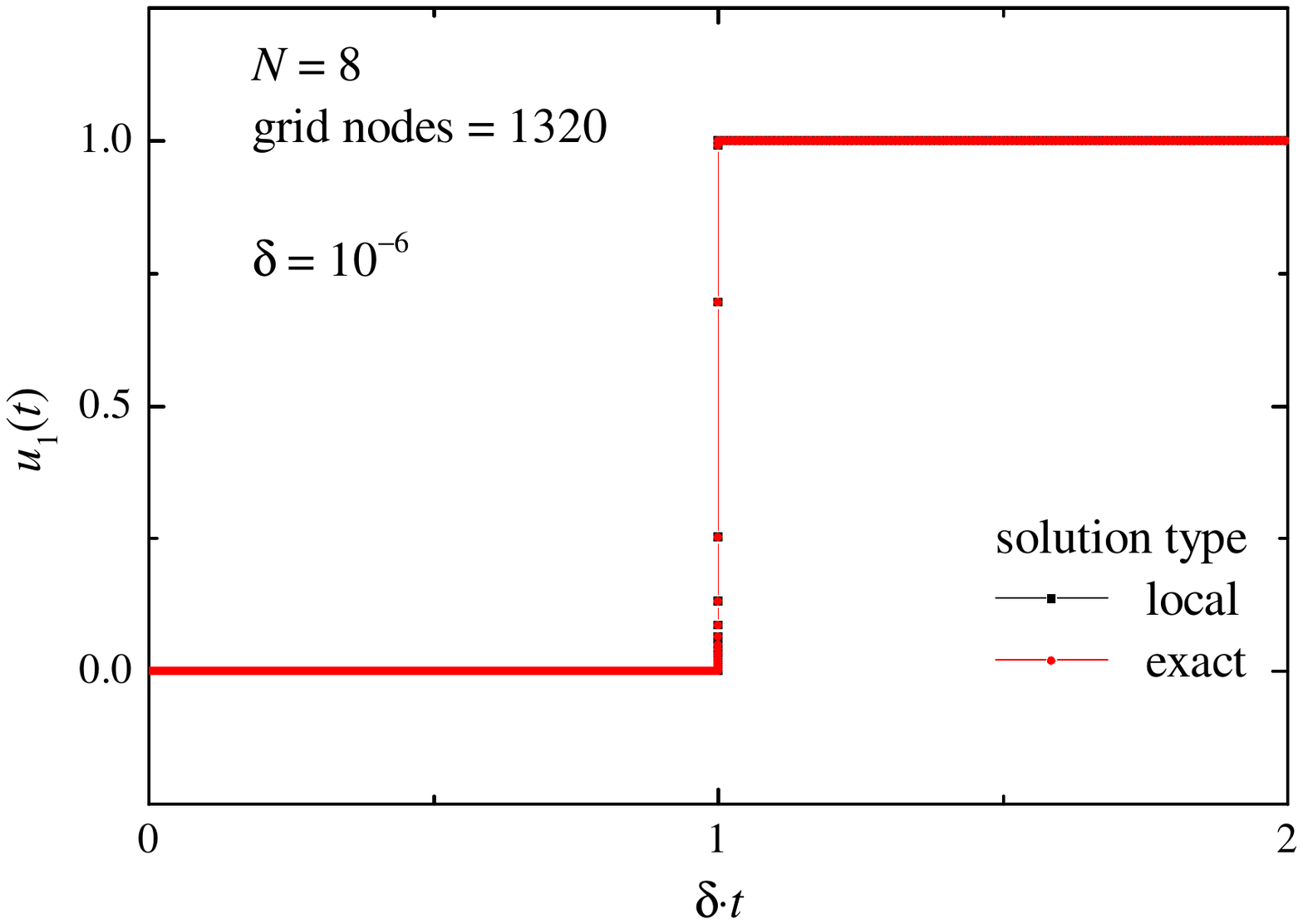}
\vspace{-8mm}\caption{\label{fig:shampine_flame_delta_10m6_sol_qug:a2}}
\end{subfigure}
\begin{subfigure}{0.320\textwidth}
\includegraphics[width=\textwidth]{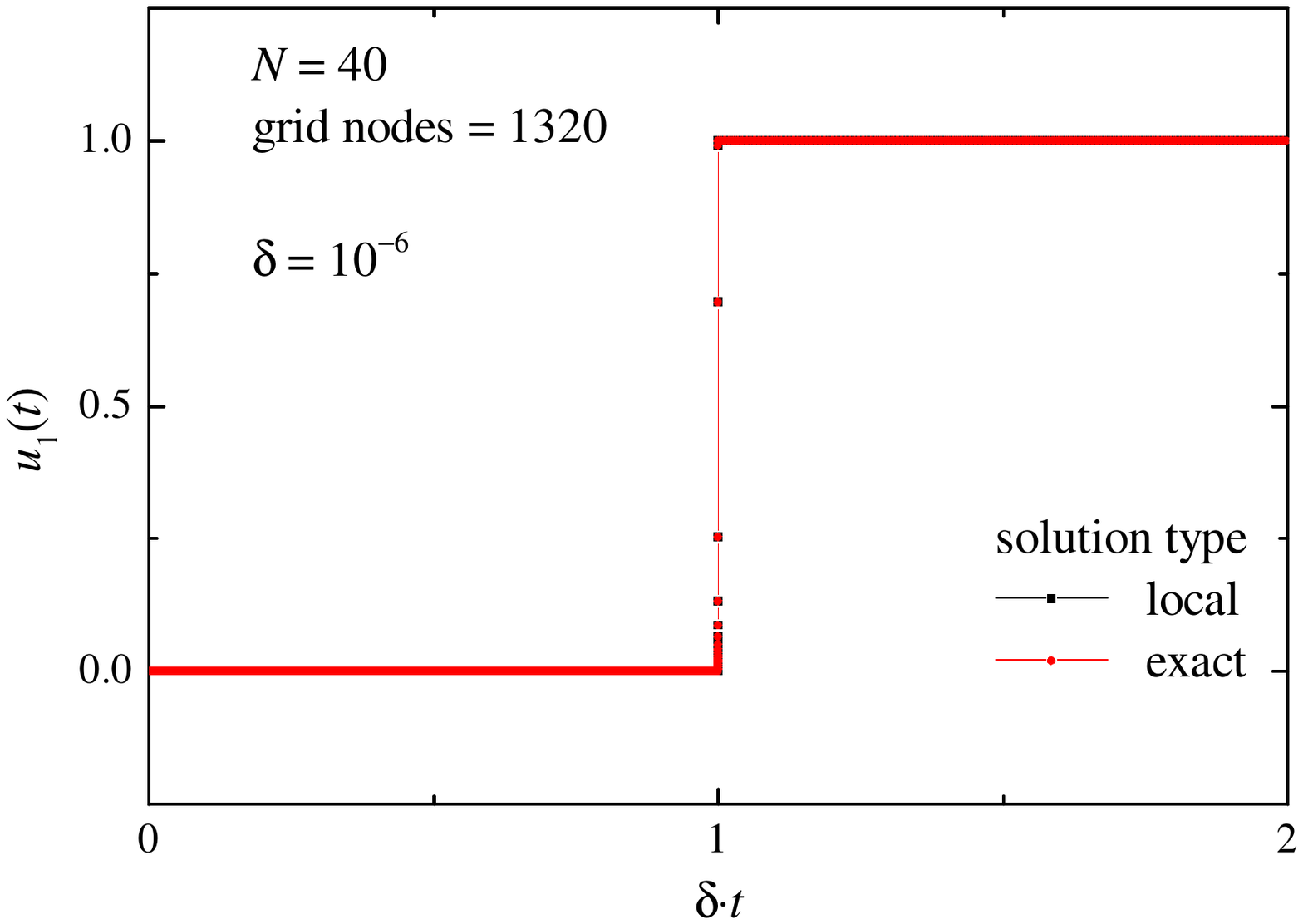}
\vspace{-8mm}\caption{\label{fig:shampine_flame_delta_10m6_sol_qug:a3}}
\end{subfigure}\\
\begin{subfigure}{0.320\textwidth}
\includegraphics[width=\textwidth]{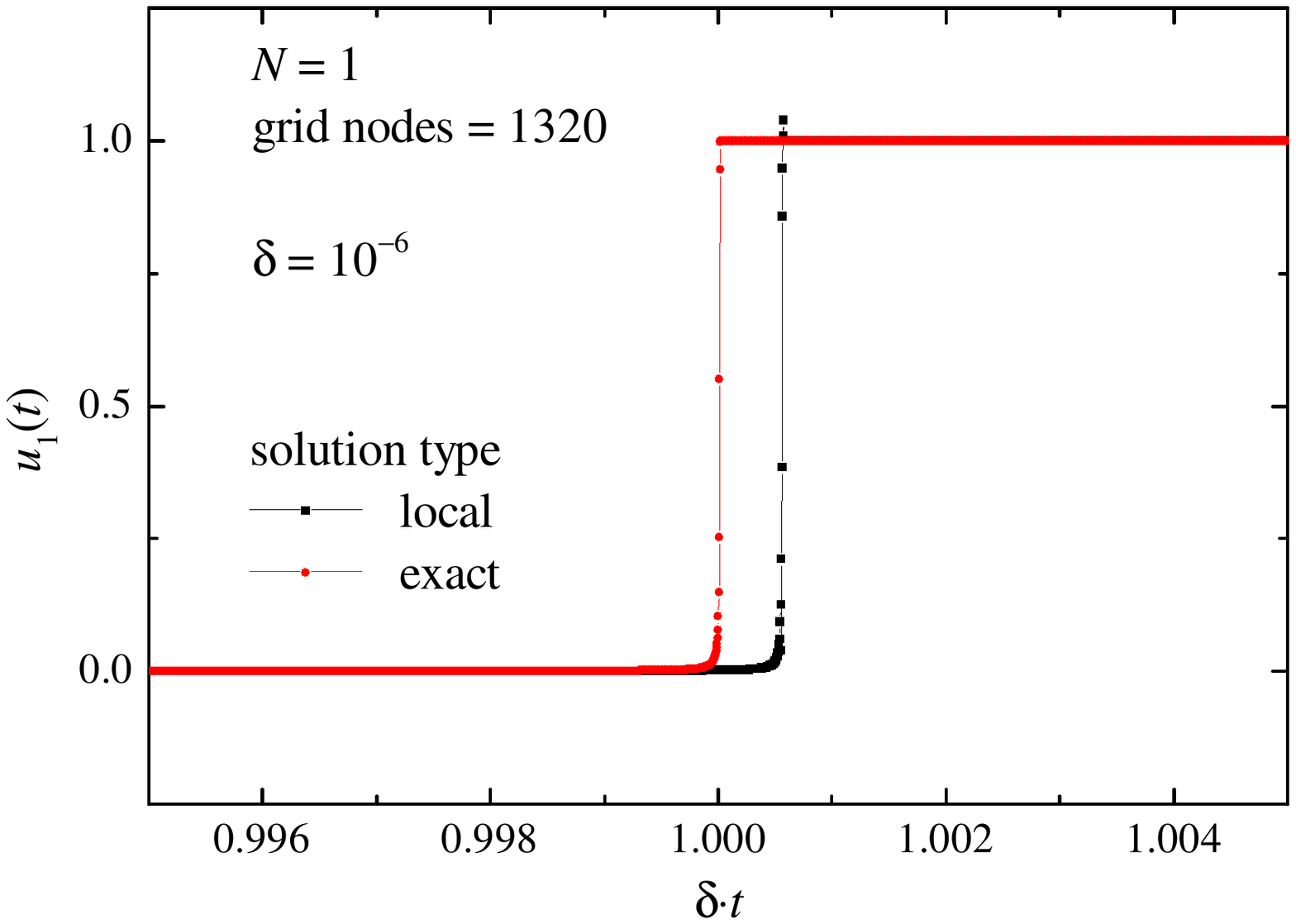}
\vspace{-8mm}\caption{\label{fig:shampine_flame_delta_10m6_sol_qug:b1}}
\end{subfigure}
\begin{subfigure}{0.320\textwidth}
\includegraphics[width=\textwidth]{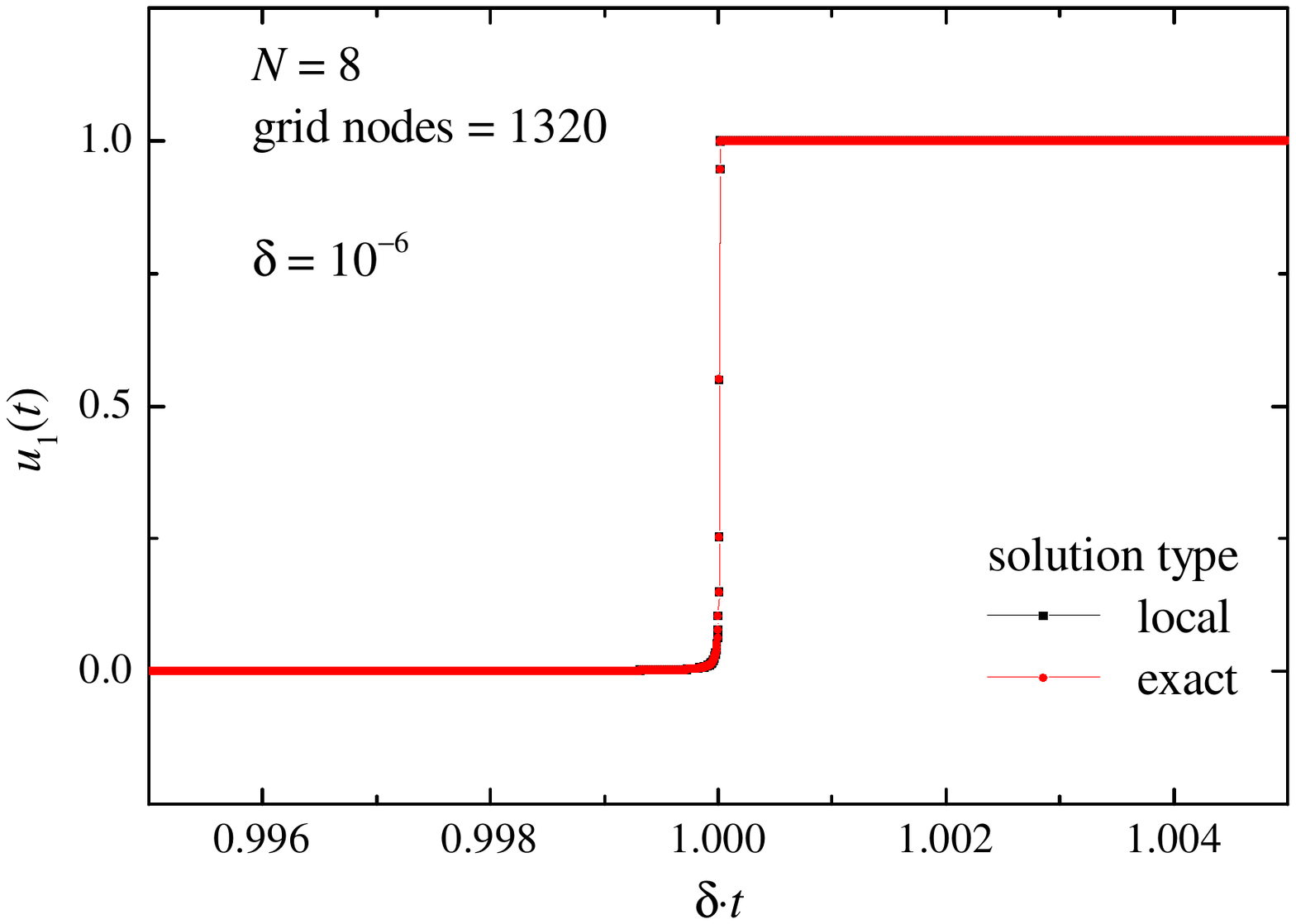}
\vspace{-8mm}\caption{\label{fig:shampine_flame_delta_10m6_sol_qug:b2}}
\end{subfigure}
\begin{subfigure}{0.320\textwidth}
\includegraphics[width=\textwidth]{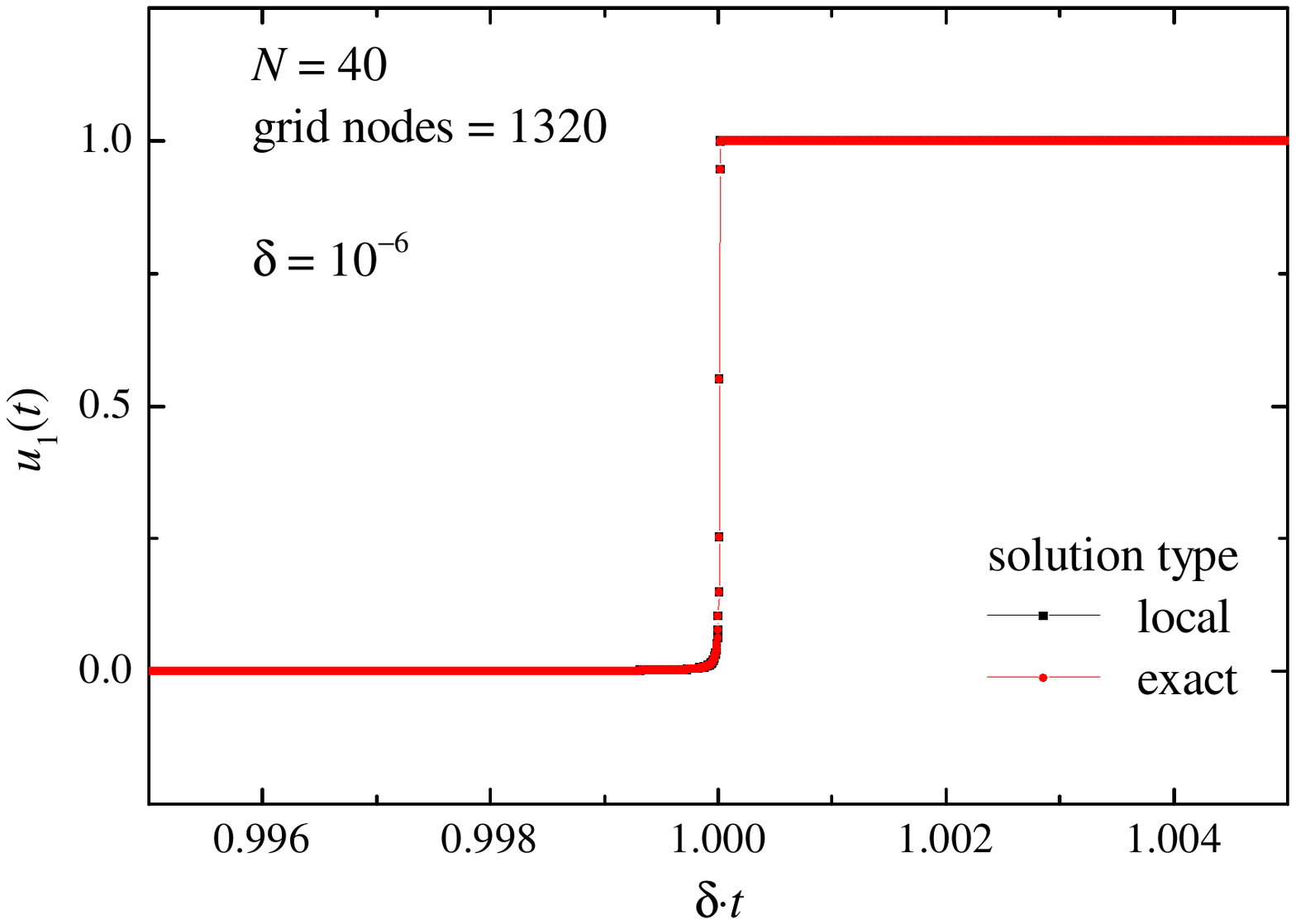}
\vspace{-8mm}\caption{\label{fig:shampine_flame_delta_10m6_sol_qug:b3}}
\end{subfigure}\\
\begin{subfigure}{0.320\textwidth}
\includegraphics[width=\textwidth]{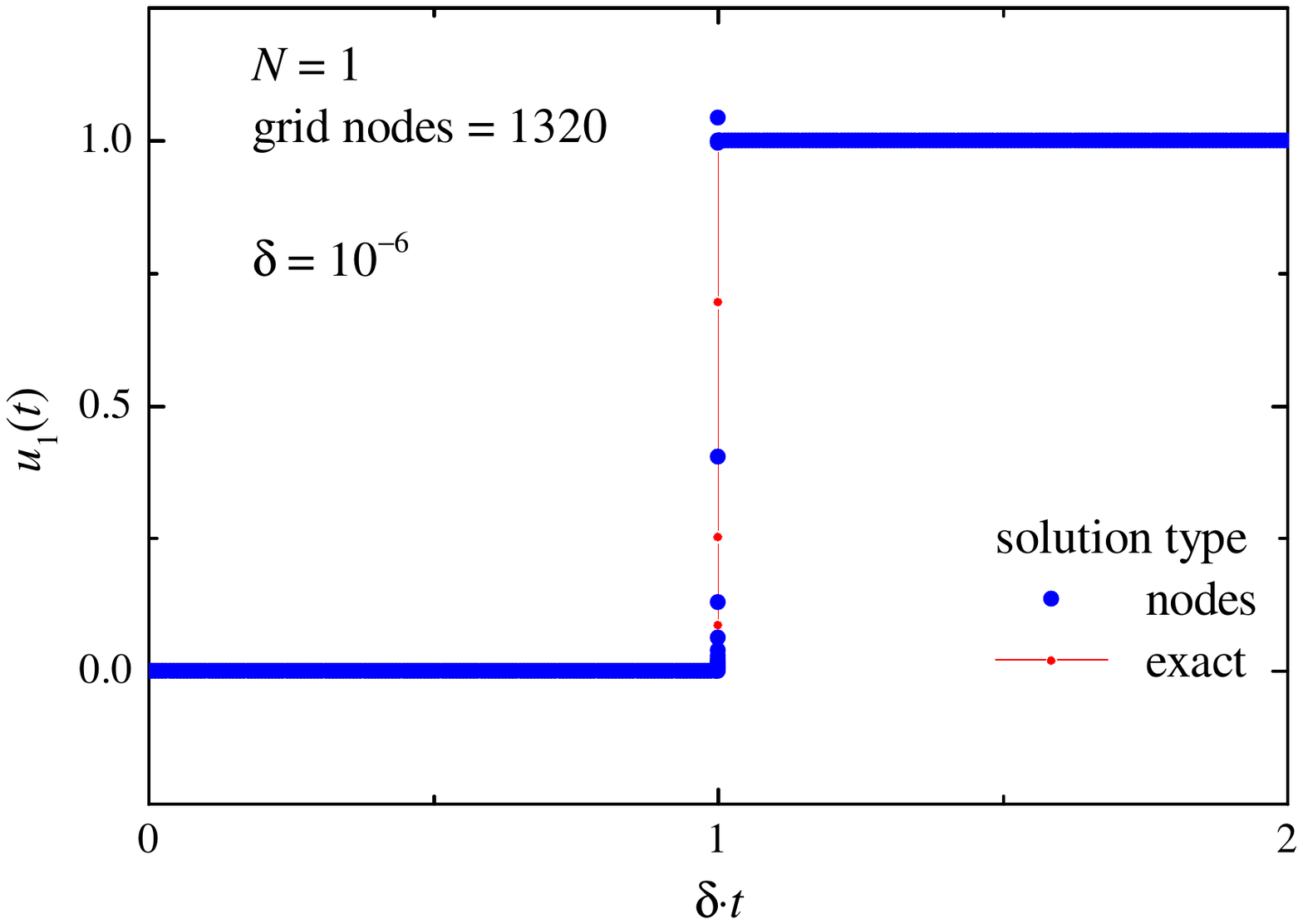}
\vspace{-8mm}\caption{\label{fig:shampine_flame_delta_10m6_sol_qug:c1}}
\end{subfigure}
\begin{subfigure}{0.320\textwidth}
\includegraphics[width=\textwidth]{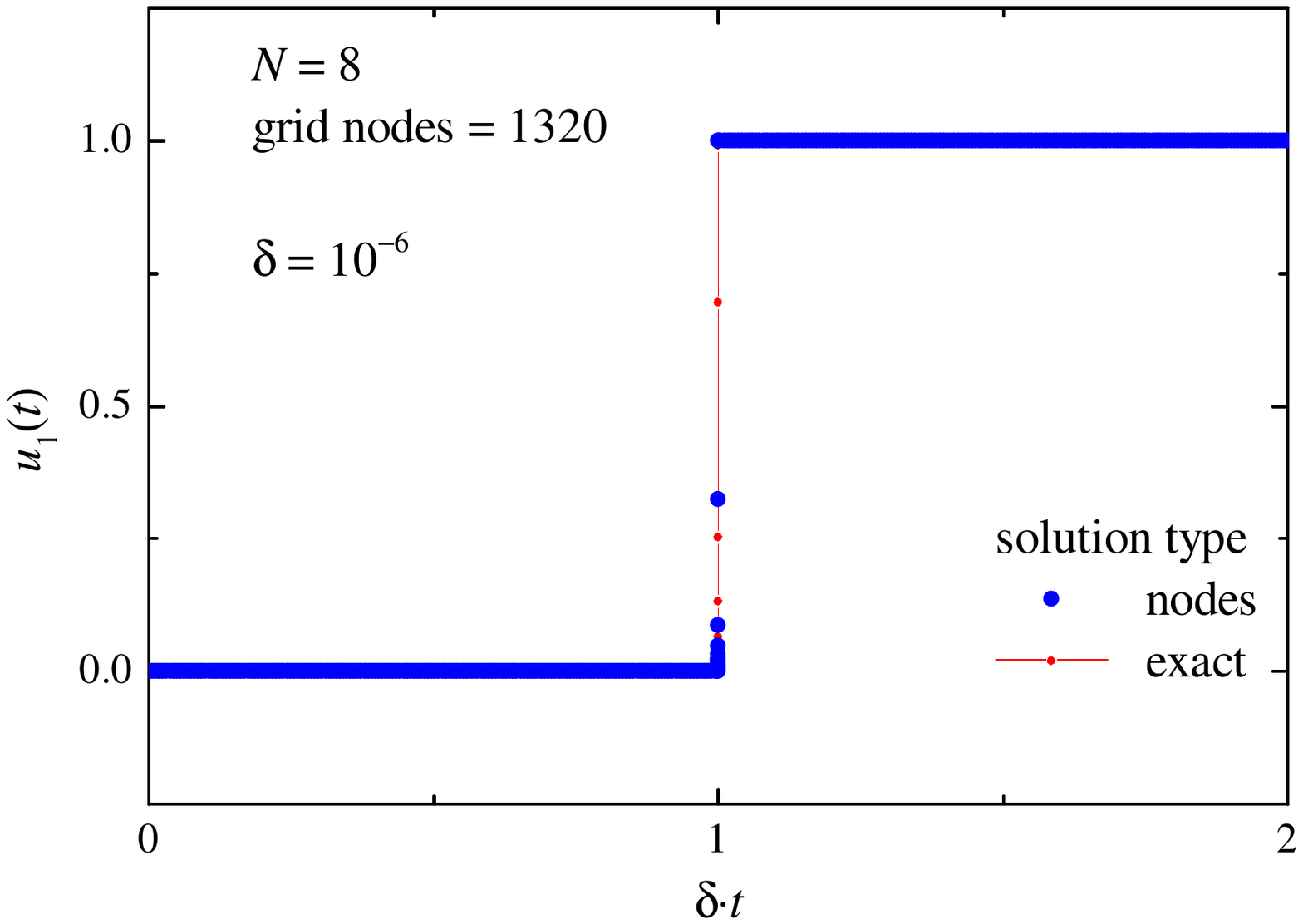}
\vspace{-8mm}\caption{\label{fig:shampine_flame_delta_10m6_sol_qug:c2}}
\end{subfigure}
\begin{subfigure}{0.320\textwidth}
\includegraphics[width=\textwidth]{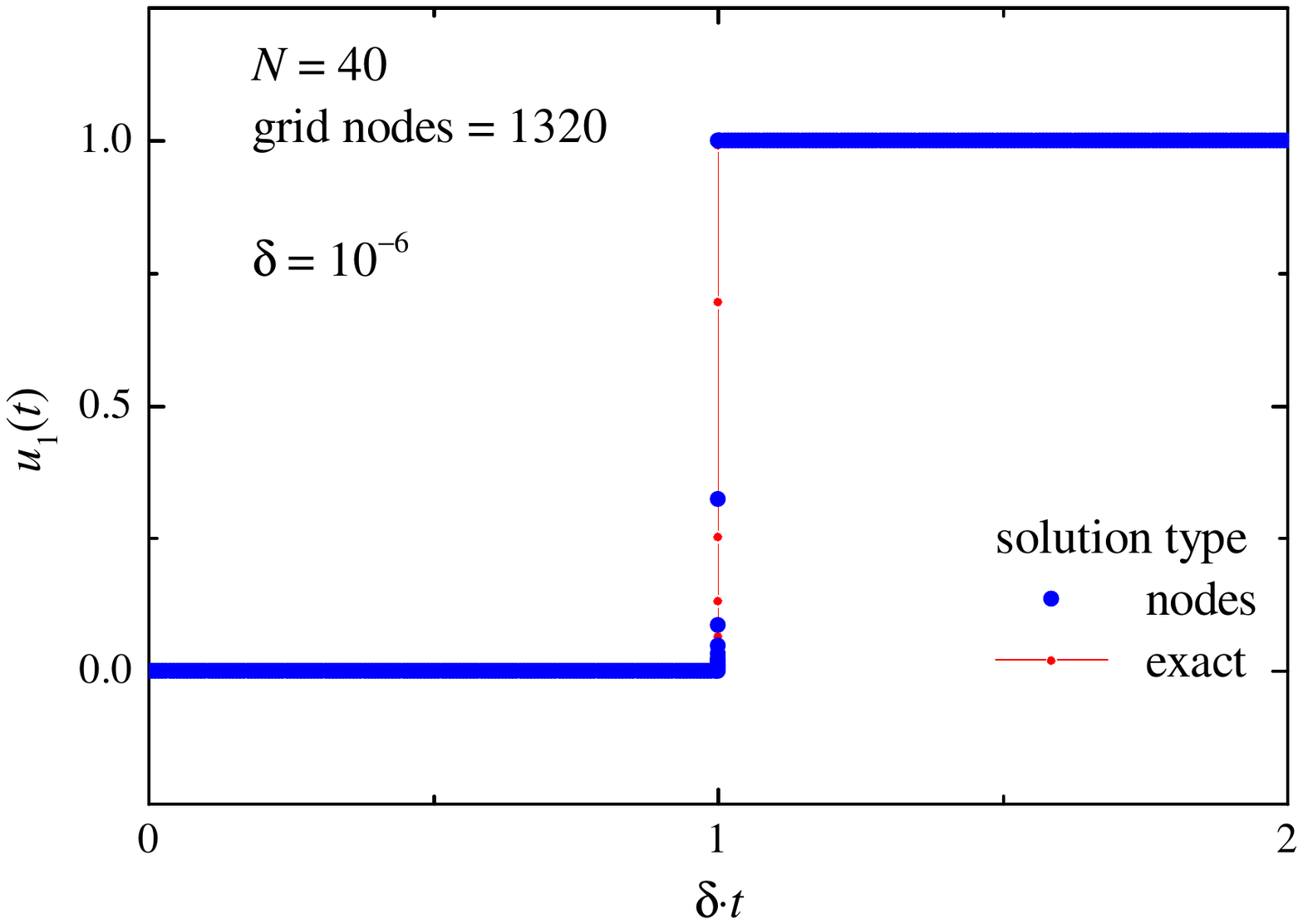}
\vspace{-8mm}\caption{\label{fig:shampine_flame_delta_10m6_sol_qug:c3}}
\end{subfigure}\\
\begin{subfigure}{0.320\textwidth}
\includegraphics[width=\textwidth]{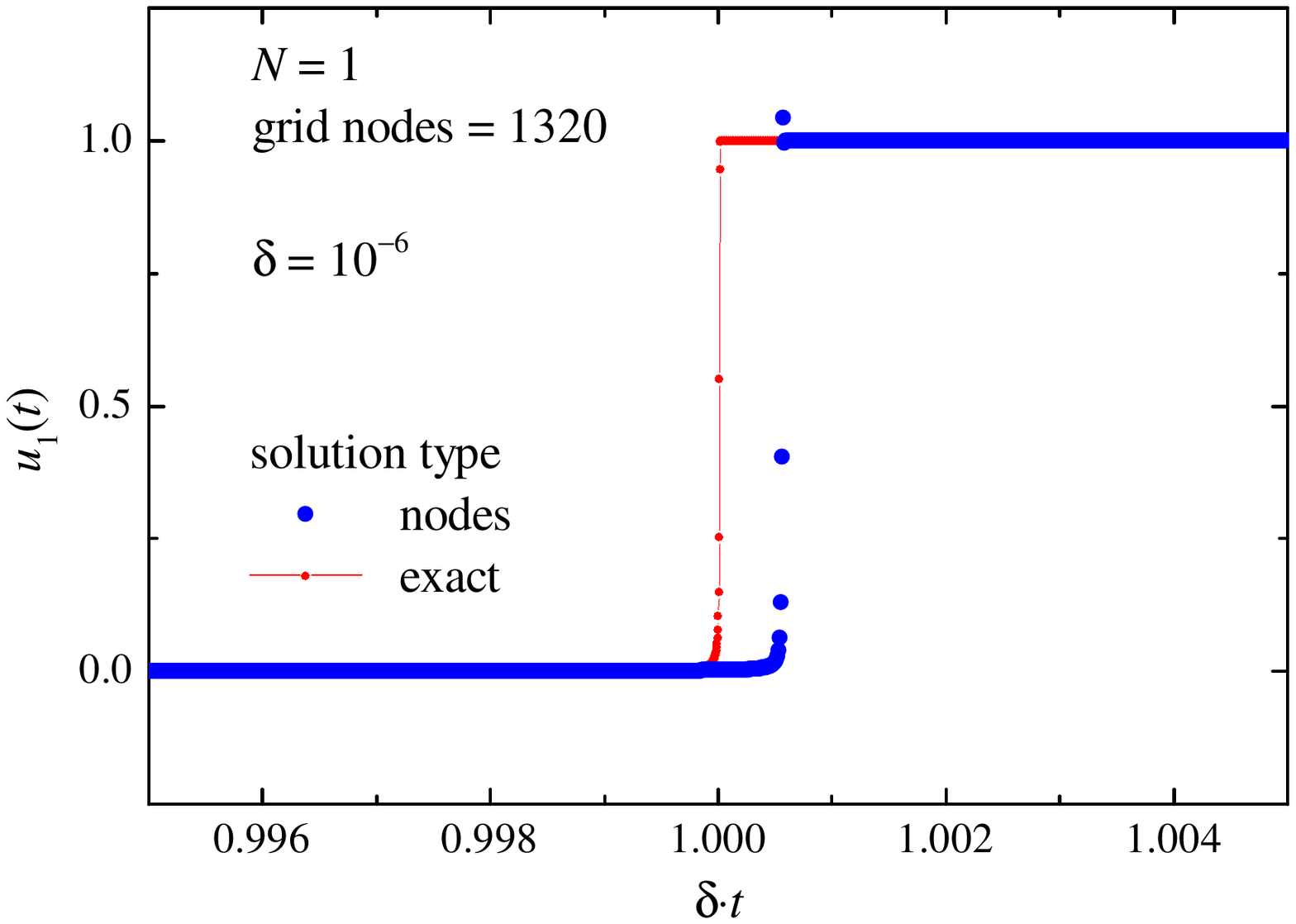}
\vspace{-8mm}\caption{\label{fig:shampine_flame_delta_10m6_sol_qug:d1}}
\end{subfigure}
\begin{subfigure}{0.320\textwidth}
\includegraphics[width=\textwidth]{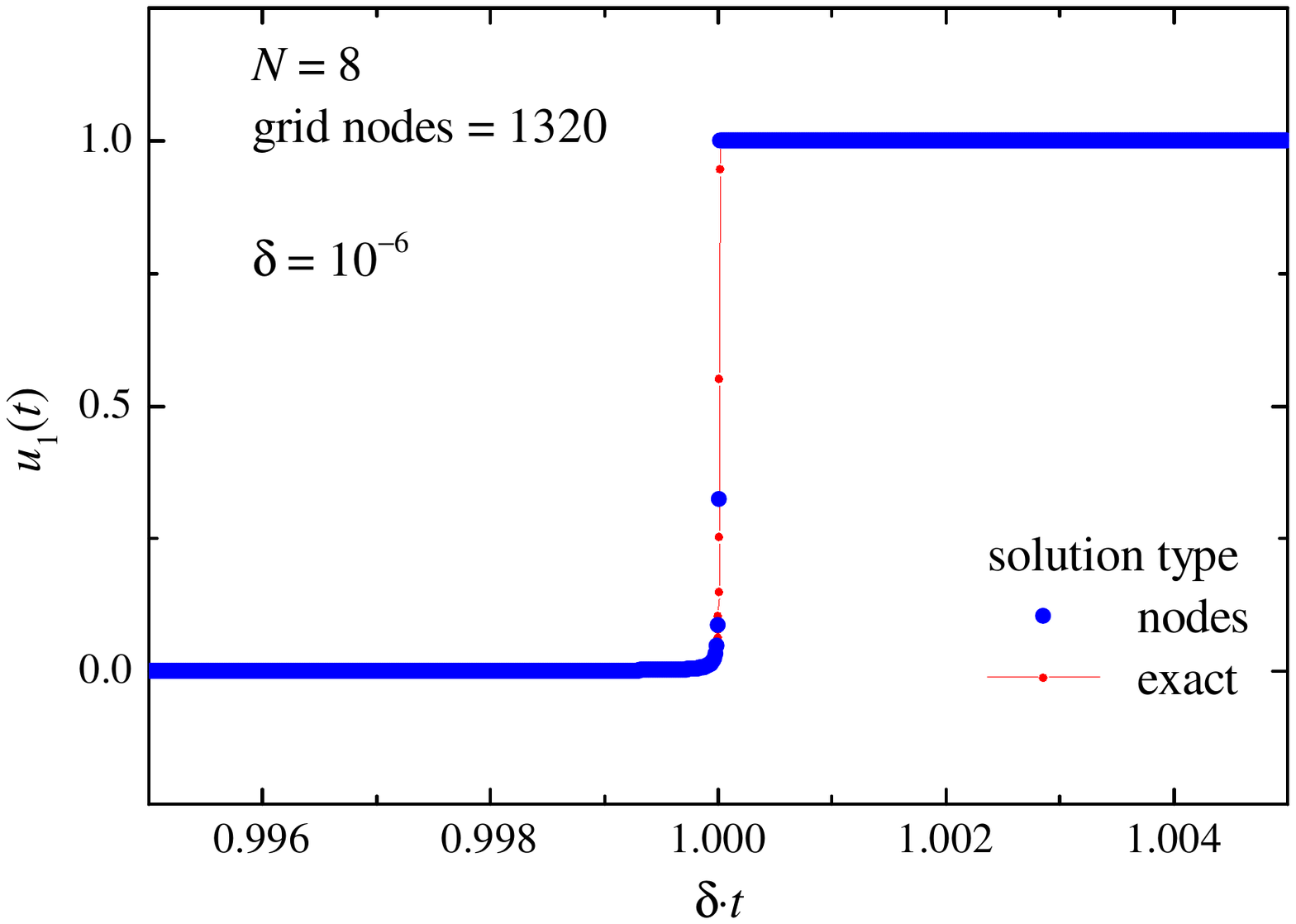}
\vspace{-8mm}\caption{\label{fig:shampine_flame_delta_10m6_sol_qug:d2}}
\end{subfigure}
\begin{subfigure}{0.320\textwidth}
\includegraphics[width=\textwidth]{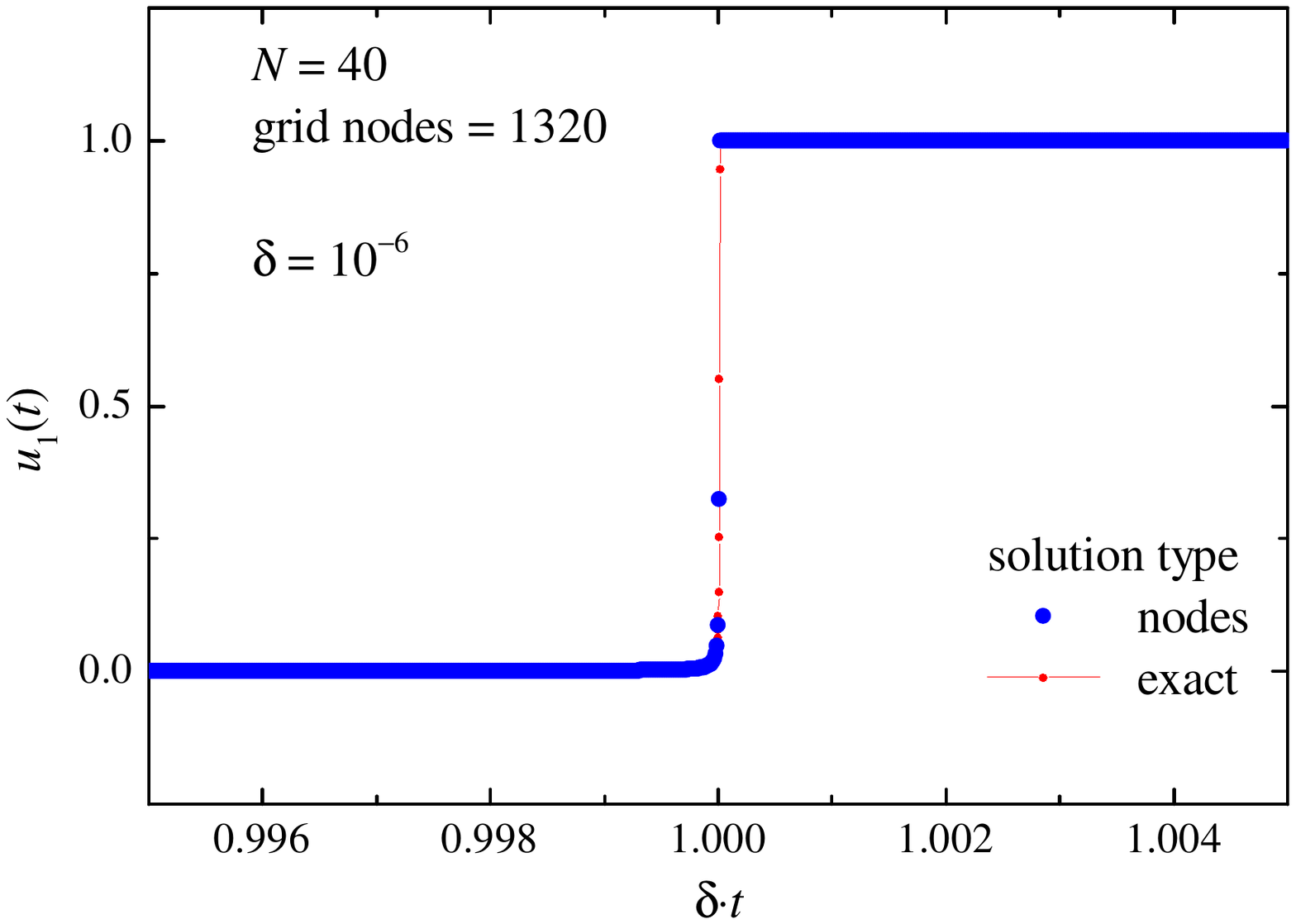}
\vspace{-8mm}\caption{\label{fig:shampine_flame_delta_10m6_sol_qug:d3}}
\end{subfigure}\\
\begin{subfigure}{0.320\textwidth}
\includegraphics[width=\textwidth]{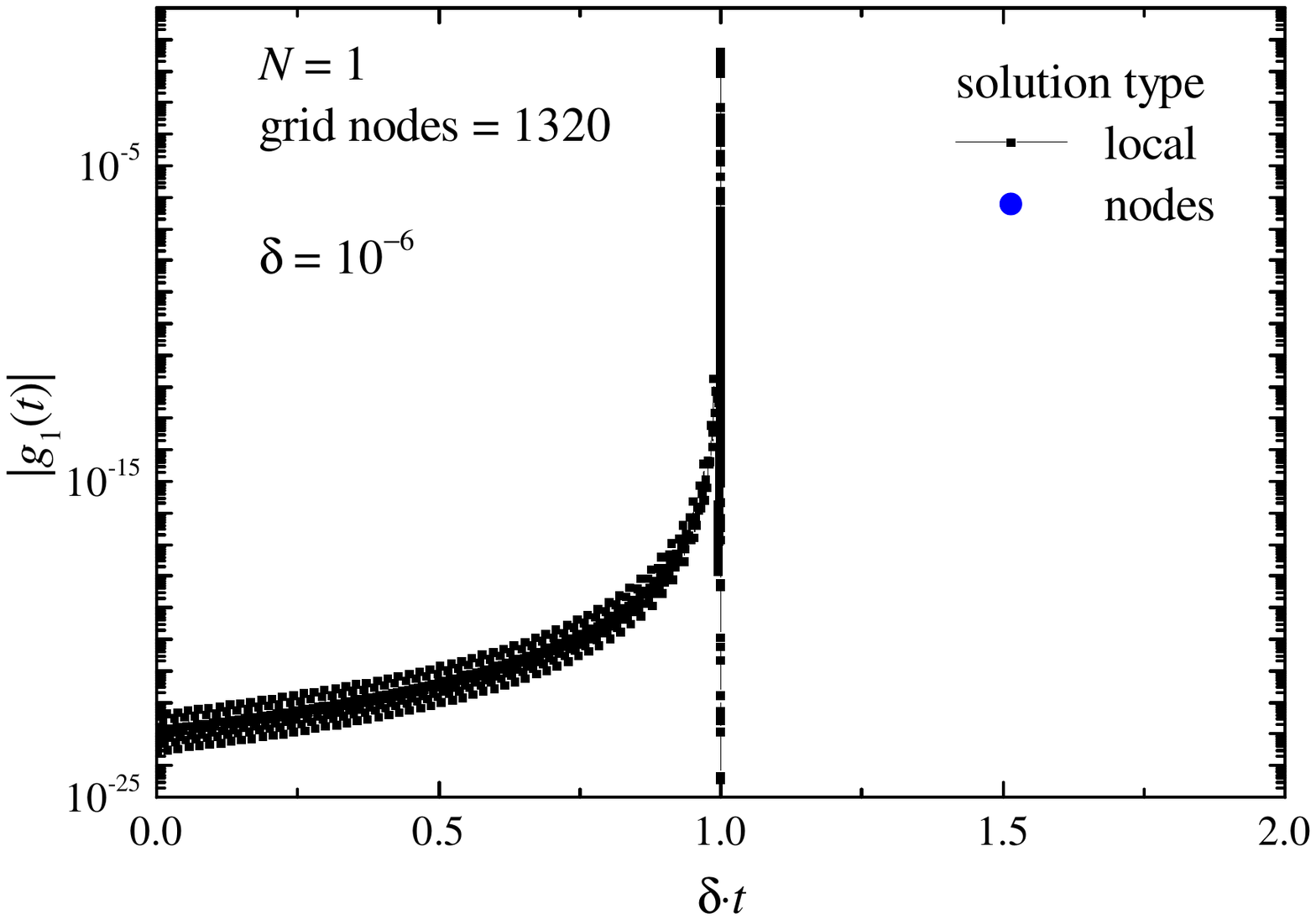}
\vspace{-8mm}\caption{\label{fig:shampine_flame_delta_10m6_sol_qug:e1}}
\end{subfigure}
\begin{subfigure}{0.320\textwidth}
\includegraphics[width=\textwidth]{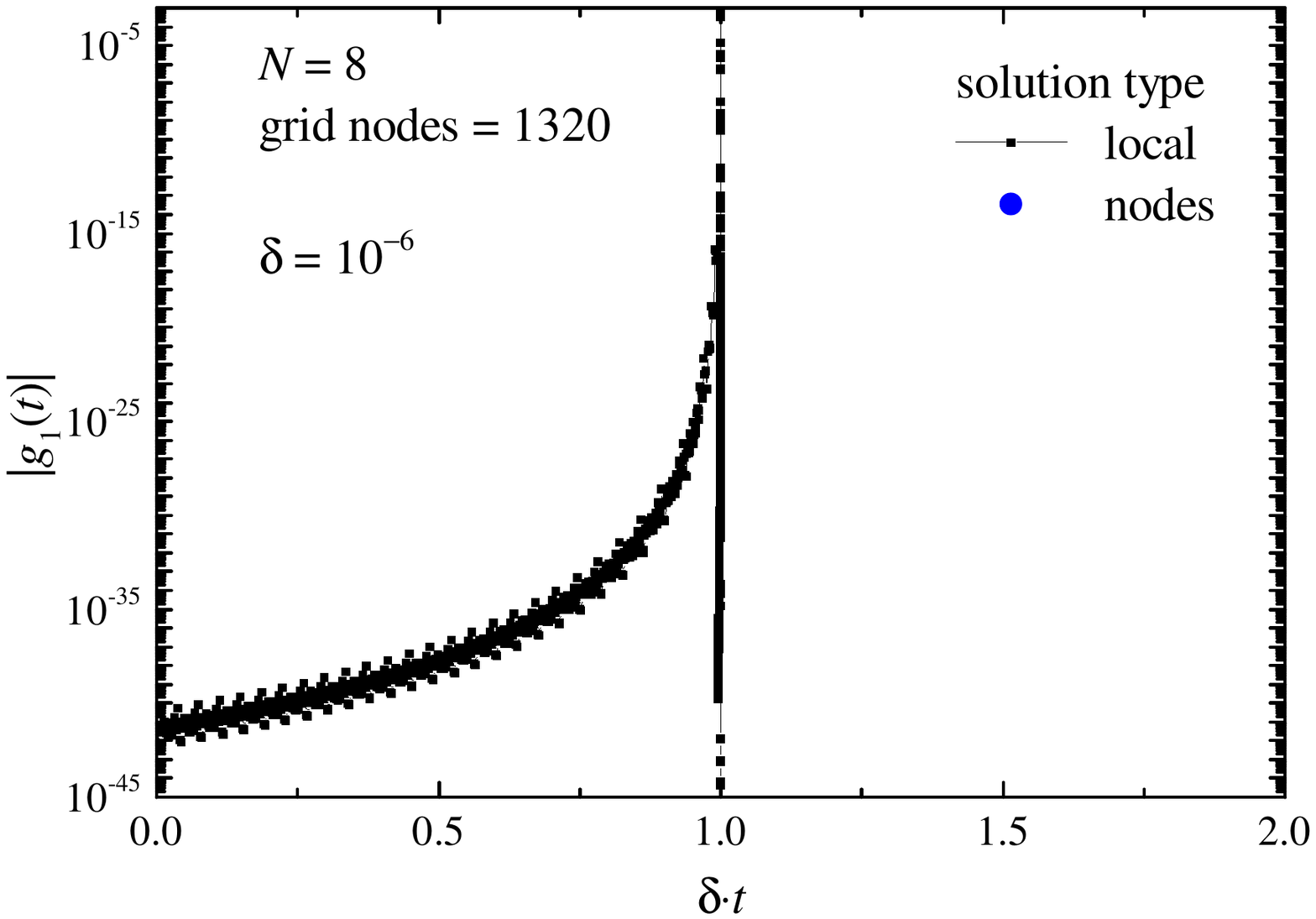}
\vspace{-8mm}\caption{\label{fig:shampine_flame_delta_10m6_sol_qug:e2}}
\end{subfigure}
\begin{subfigure}{0.320\textwidth}
\includegraphics[width=\textwidth]{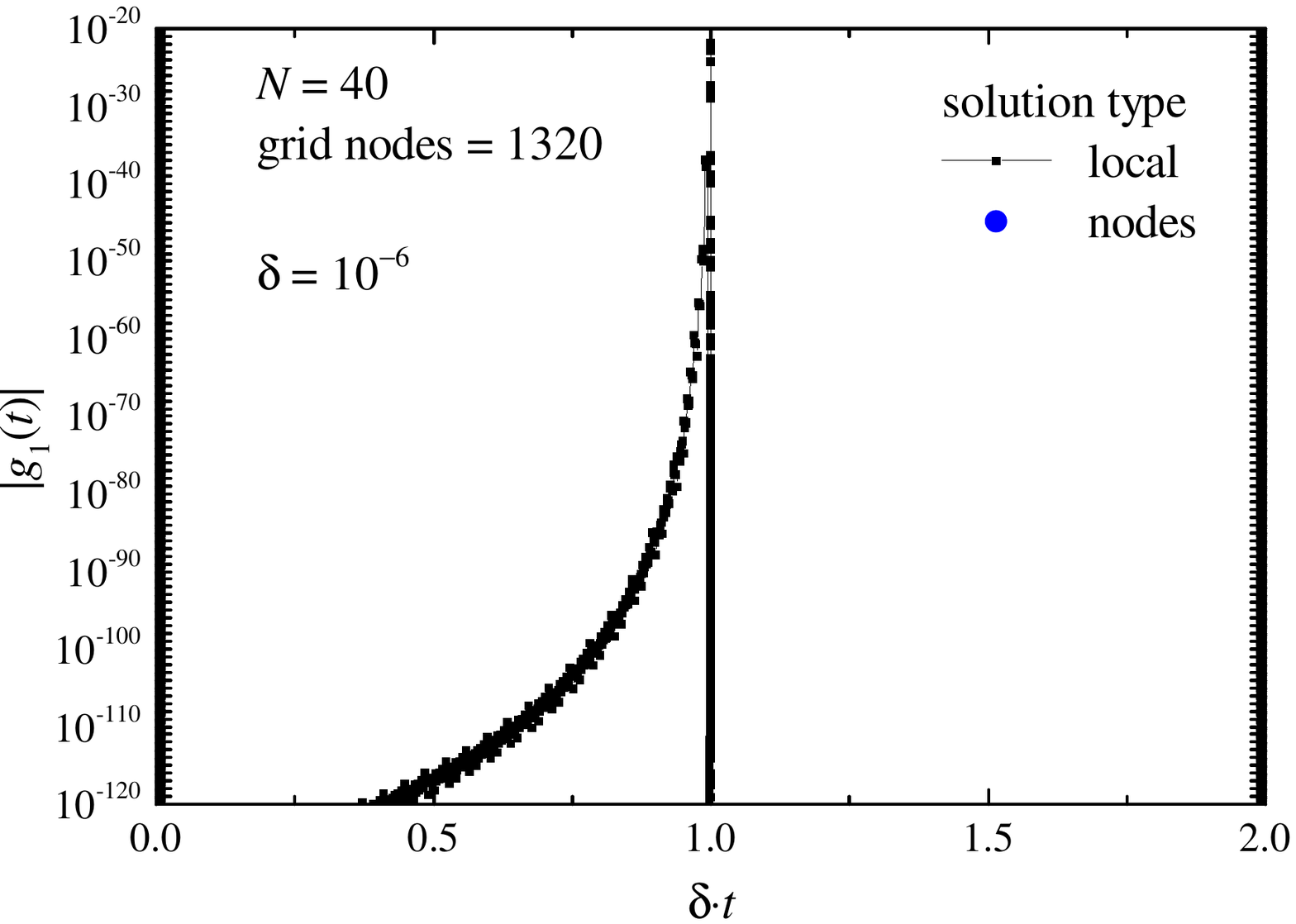}
\vspace{-8mm}\caption{\label{fig:shampine_flame_delta_10m6_sol_qug:e3}}
\end{subfigure}\\
\caption{%
Numerical solution of the stiff DAE system (\ref{eq:shampine_flame}) of index 1 with $\delta = 10^{-6}$. Comparison of the solution at nodes $\mathbf{u}_{n}$ (\subref{fig:shampine_flame_delta_10m6_sol_qug:c1}, \subref{fig:shampine_flame_delta_10m6_sol_qug:c2}, \subref{fig:shampine_flame_delta_10m6_sol_qug:c3}, \subref{fig:shampine_flame_delta_10m6_sol_qug:d1}, \subref{fig:shampine_flame_delta_10m6_sol_qug:d2}, \subref{fig:shampine_flame_delta_10m6_sol_qug:d3}), the local solution $\mathbf{u}_{L}(t)$ (\subref{fig:shampine_flame_delta_10m6_sol_qug:a1}, \subref{fig:shampine_flame_delta_10m6_sol_qug:a2}, \subref{fig:shampine_flame_delta_10m6_sol_qug:a3}, \subref{fig:shampine_flame_delta_10m6_sol_qug:b1}, \subref{fig:shampine_flame_delta_10m6_sol_qug:b2}, \subref{fig:shampine_flame_delta_10m6_sol_qug:b3}) and the exact solution $\mathbf{u}^{\rm ex}(t)$ for component $u_{1}$ (\subref{fig:shampine_flame_delta_10m6_sol_qug:a1}, \subref{fig:shampine_flame_delta_10m6_sol_qug:a2}, \subref{fig:shampine_flame_delta_10m6_sol_qug:a3}, \subref{fig:shampine_flame_delta_10m6_sol_qug:b1}, \subref{fig:shampine_flame_delta_10m6_sol_qug:b2}, \subref{fig:shampine_flame_delta_10m6_sol_qug:b3}, \subref{fig:shampine_flame_delta_10m6_sol_qug:c1}, \subref{fig:shampine_flame_delta_10m6_sol_qug:c2}, \subref{fig:shampine_flame_delta_10m6_sol_qug:c3}, \subref{fig:shampine_flame_delta_10m6_sol_qug:d1}, \subref{fig:shampine_flame_delta_10m6_sol_qug:d2}, \subref{fig:shampine_flame_delta_10m6_sol_qug:d3}), quantitative satisfiability of the conditions $g_{1} = 0$ (\subref{fig:shampine_flame_delta_10m6_sol_qug:e1}, \subref{fig:shampine_flame_delta_10m6_sol_qug:e2}, \subref{fig:shampine_flame_delta_10m6_sol_qug:e3}), obtained using polynomials with degrees $N = 1$ (\subref{fig:shampine_flame_delta_10m6_sol_qug:a1}, \subref{fig:shampine_flame_delta_10m6_sol_qug:b1}, \subref{fig:shampine_flame_delta_10m6_sol_qug:c1}, \subref{fig:shampine_flame_delta_10m6_sol_qug:d1}, \subref{fig:shampine_flame_delta_10m6_sol_qug:e1}), $N = 8$ (\subref{fig:shampine_flame_delta_10m6_sol_qug:a2}, \subref{fig:shampine_flame_delta_10m6_sol_qug:b2}, \subref{fig:shampine_flame_delta_10m6_sol_qug:c2}, \subref{fig:shampine_flame_delta_10m6_sol_qug:d2}, \subref{fig:shampine_flame_delta_10m6_sol_qug:e2}) and $N = 40$ (\subref{fig:shampine_flame_delta_10m6_sol_qug:a3}, \subref{fig:shampine_flame_delta_10m6_sol_qug:b3}, \subref{fig:shampine_flame_delta_10m6_sol_qug:c3}, \subref{fig:shampine_flame_delta_10m6_sol_qug:d3}, \subref{fig:shampine_flame_delta_10m6_sol_qug:e3}).
}
\label{fig:shampine_flame_delta_10m6_sol_qug}
\end{figure}

\begin{figure}[h!]
\captionsetup[subfigure]{%
	position=bottom,
	font+=smaller,
	textfont=normalfont,
	singlelinecheck=off,
	justification=raggedright
}
\centering
\begin{subfigure}{0.320\textwidth}
\includegraphics[width=\textwidth]{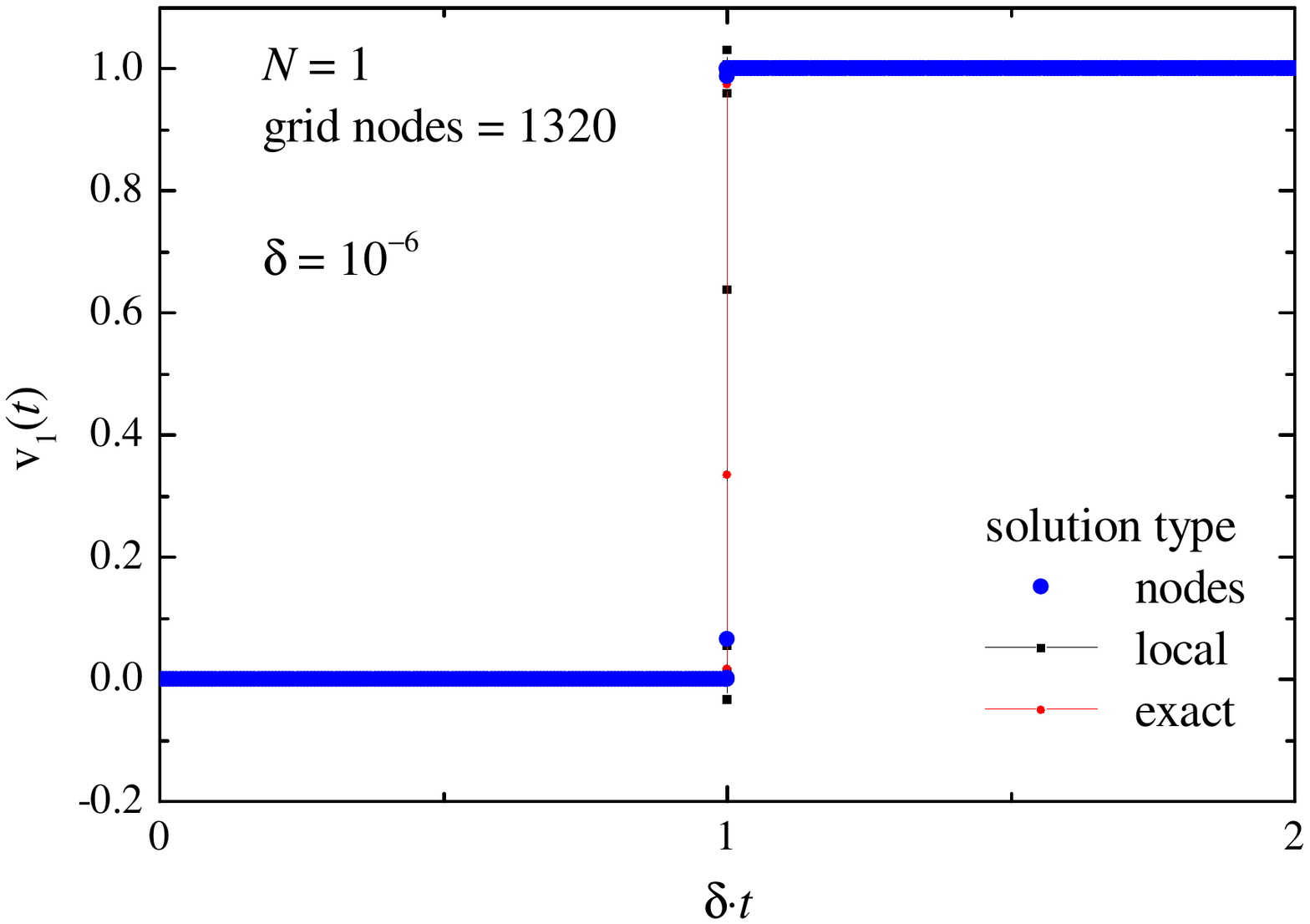}
\vspace{-8mm}\caption{\label{fig:shampine_flame_delta_10m6_sol_v_epss:a1}}
\end{subfigure}
\begin{subfigure}{0.320\textwidth}
\includegraphics[width=\textwidth]{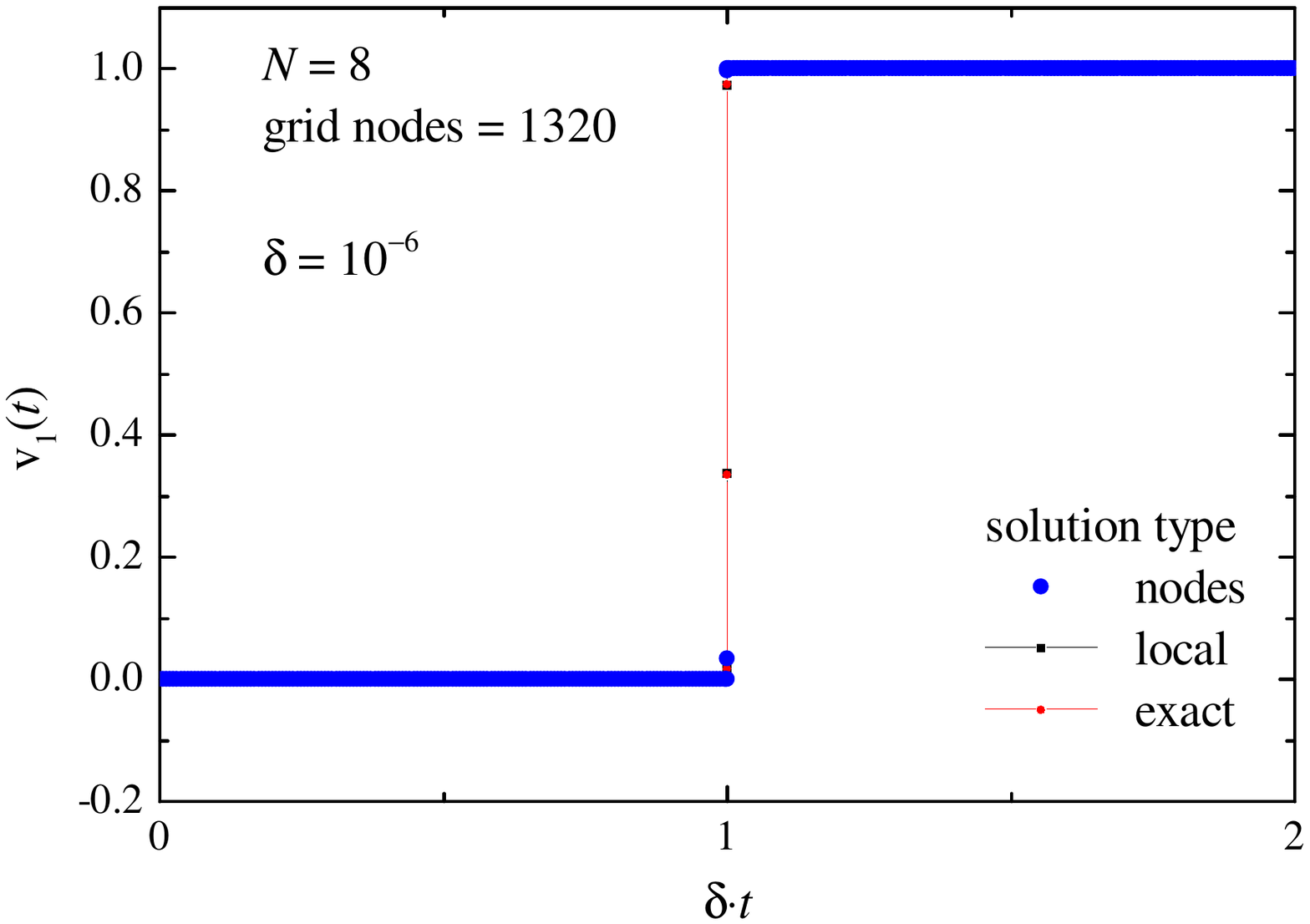}
\vspace{-8mm}\caption{\label{fig:shampine_flame_delta_10m6_sol_v_epss:a2}}
\end{subfigure}
\begin{subfigure}{0.320\textwidth}
\includegraphics[width=\textwidth]{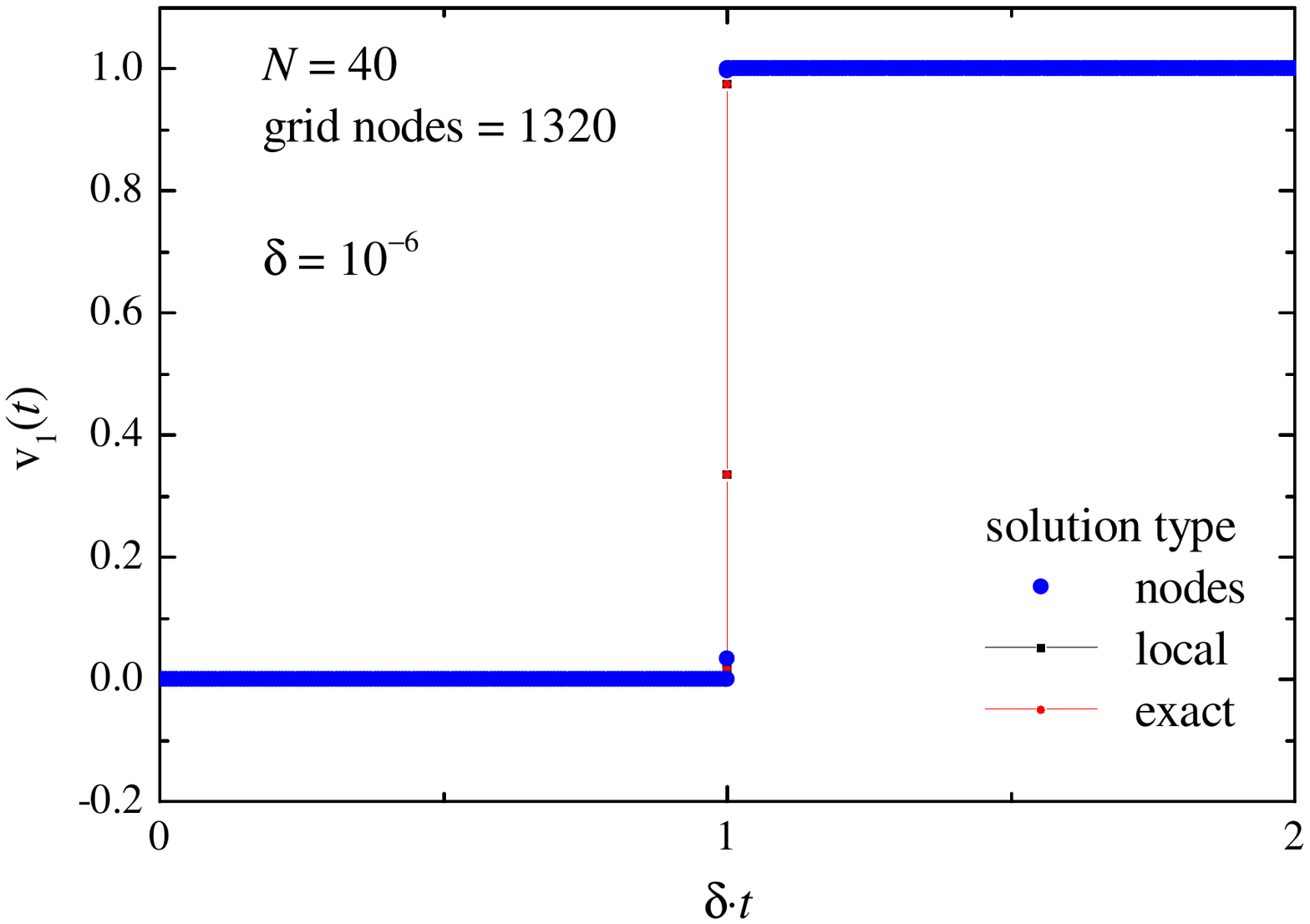}
\vspace{-8mm}\caption{\label{fig:shampine_flame_delta_10m6_sol_v_epss:a3}}
\end{subfigure}\\
\begin{subfigure}{0.320\textwidth}
\includegraphics[width=\textwidth]{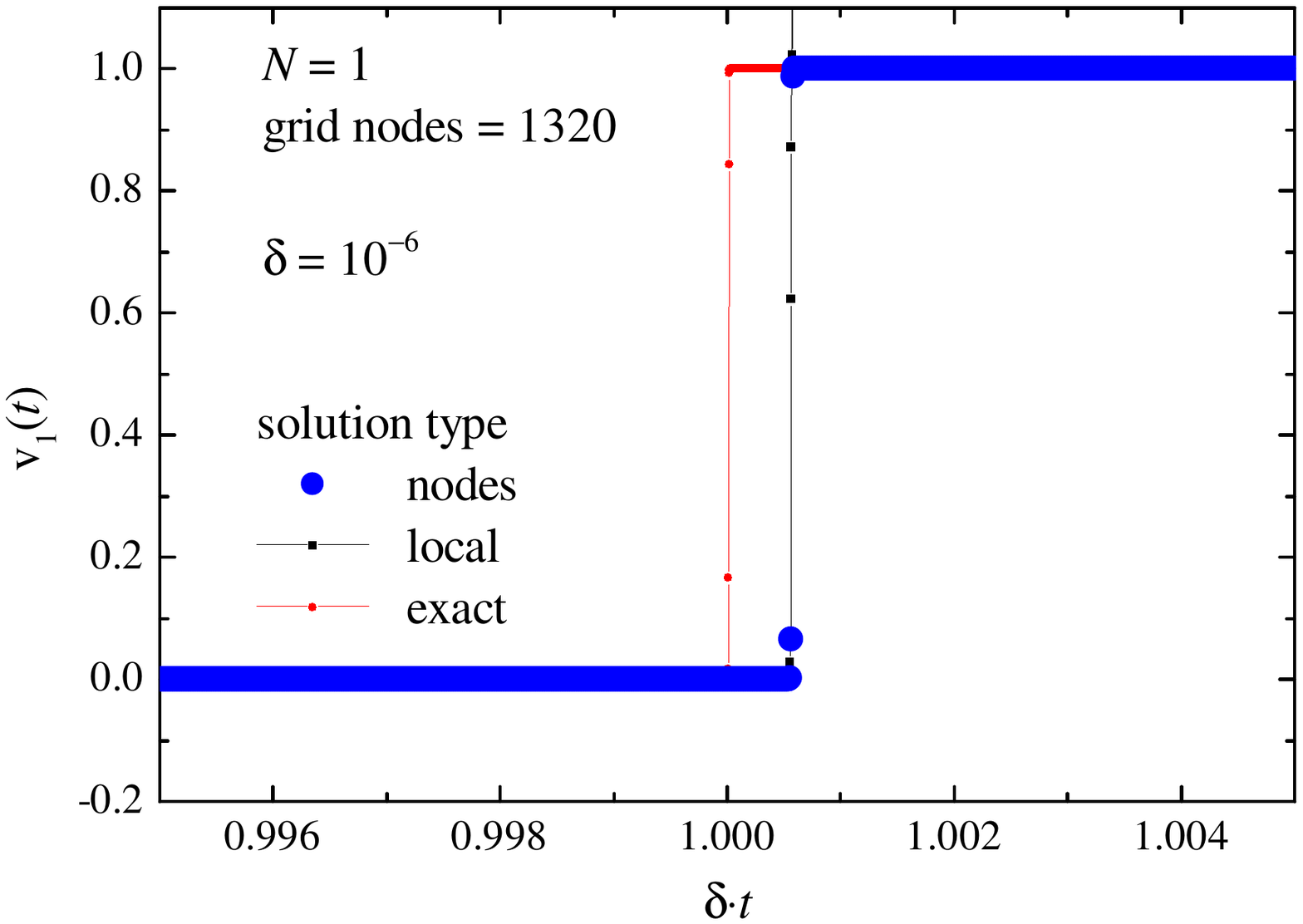}
\vspace{-8mm}\caption{\label{fig:shampine_flame_delta_10m6_sol_v_epss:b1}}
\end{subfigure}
\begin{subfigure}{0.320\textwidth}
\includegraphics[width=\textwidth]{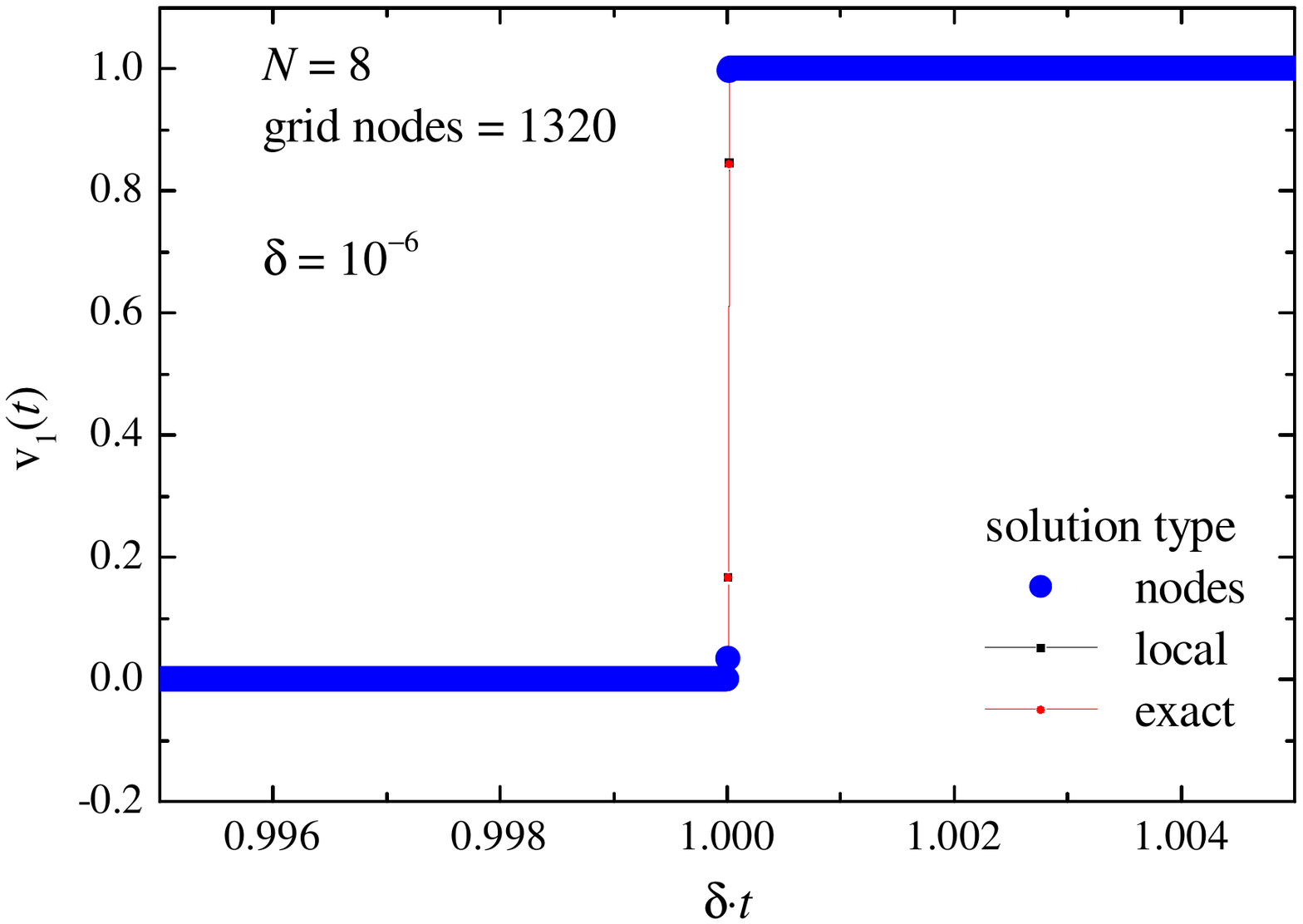}
\vspace{-8mm}\caption{\label{fig:shampine_flame_delta_10m6_sol_v_epss:b2}}
\end{subfigure}
\begin{subfigure}{0.320\textwidth}
\includegraphics[width=\textwidth]{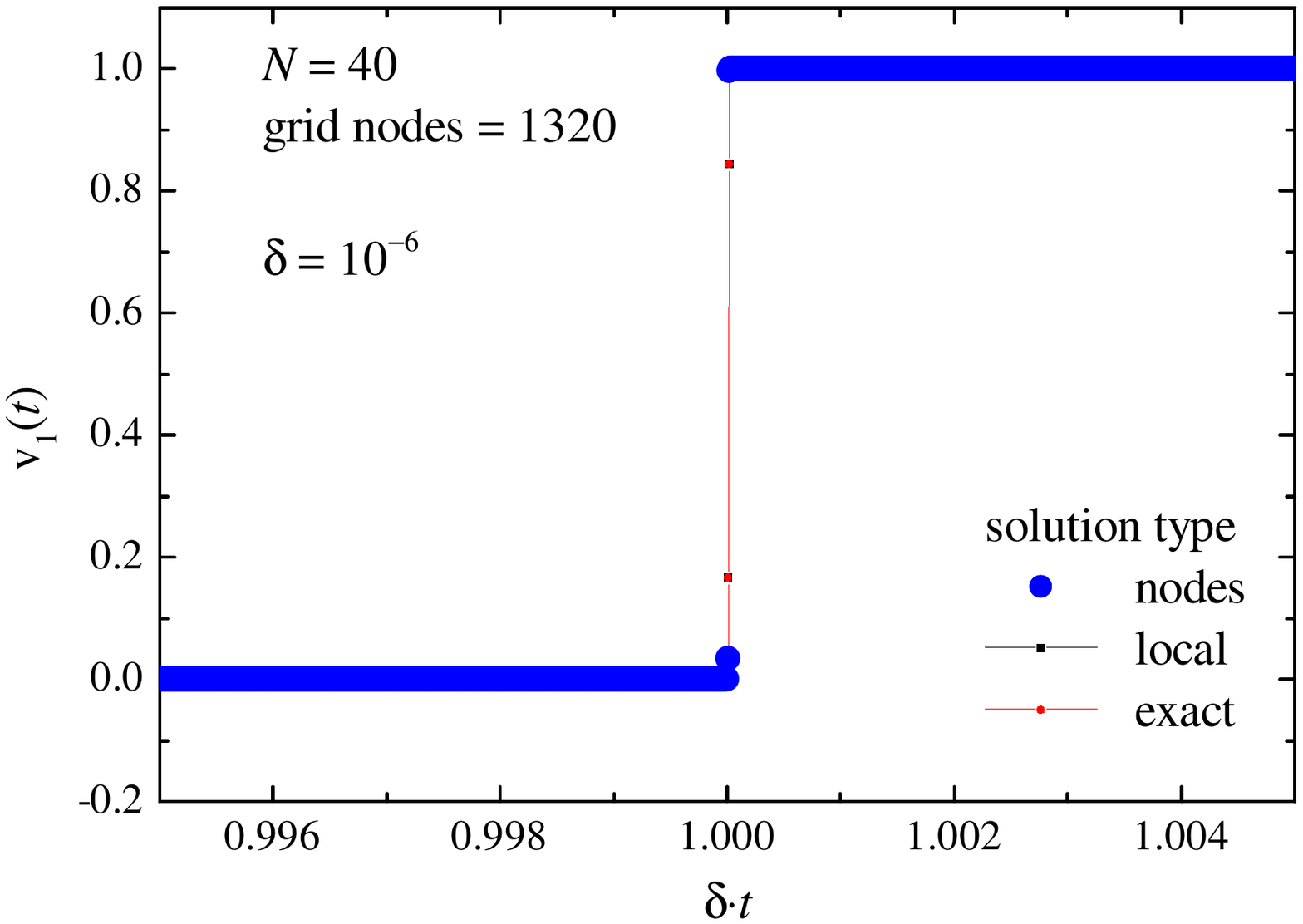}
\vspace{-8mm}\caption{\label{fig:shampine_flame_delta_10m6_sol_v_epss:b3}}
\end{subfigure}\\
\begin{subfigure}{0.320\textwidth}
\includegraphics[width=\textwidth]{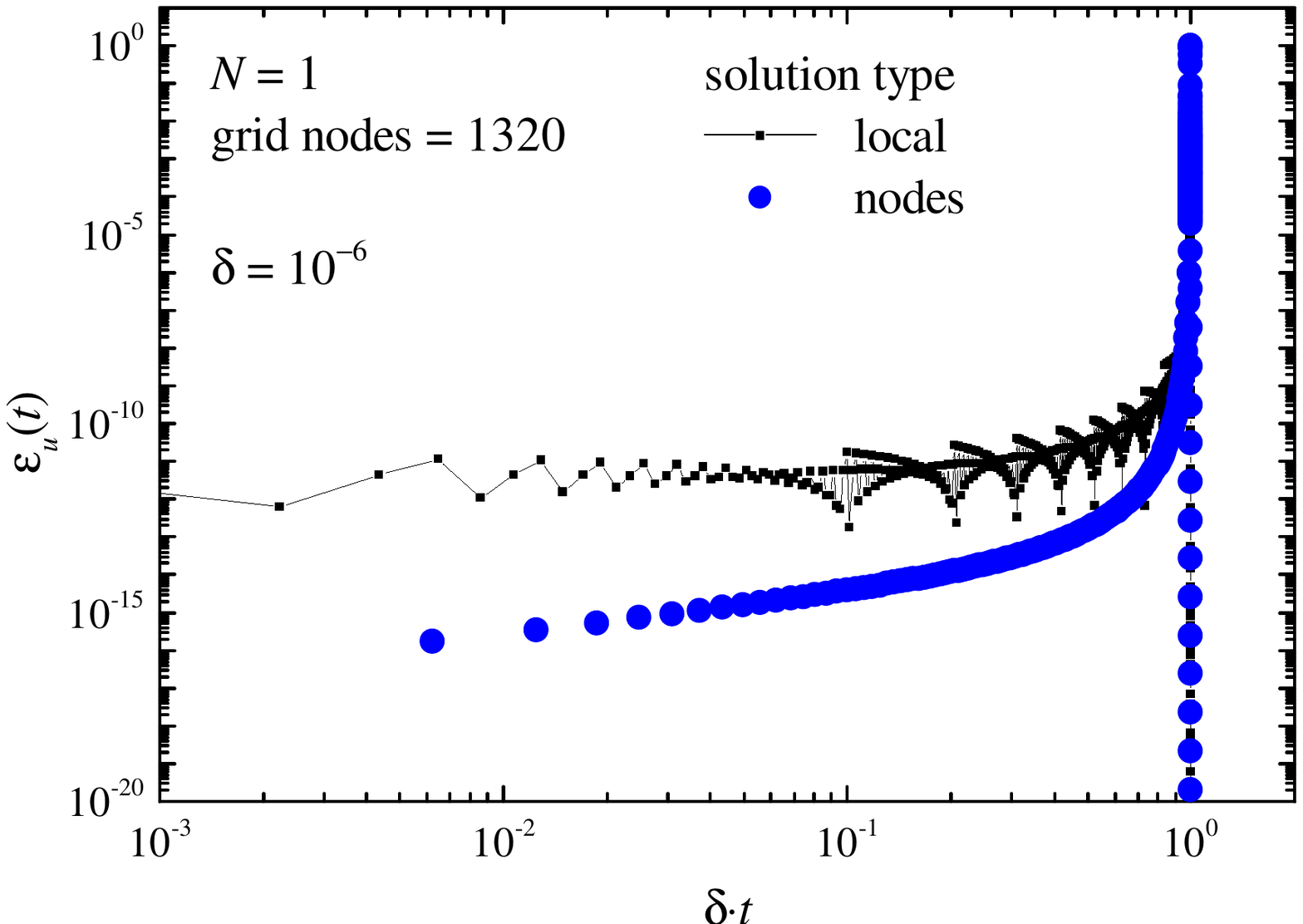}
\vspace{-8mm}\caption{\label{fig:shampine_flame_delta_10m6_sol_v_epss:c1}}
\end{subfigure}
\begin{subfigure}{0.320\textwidth}
\includegraphics[width=\textwidth]{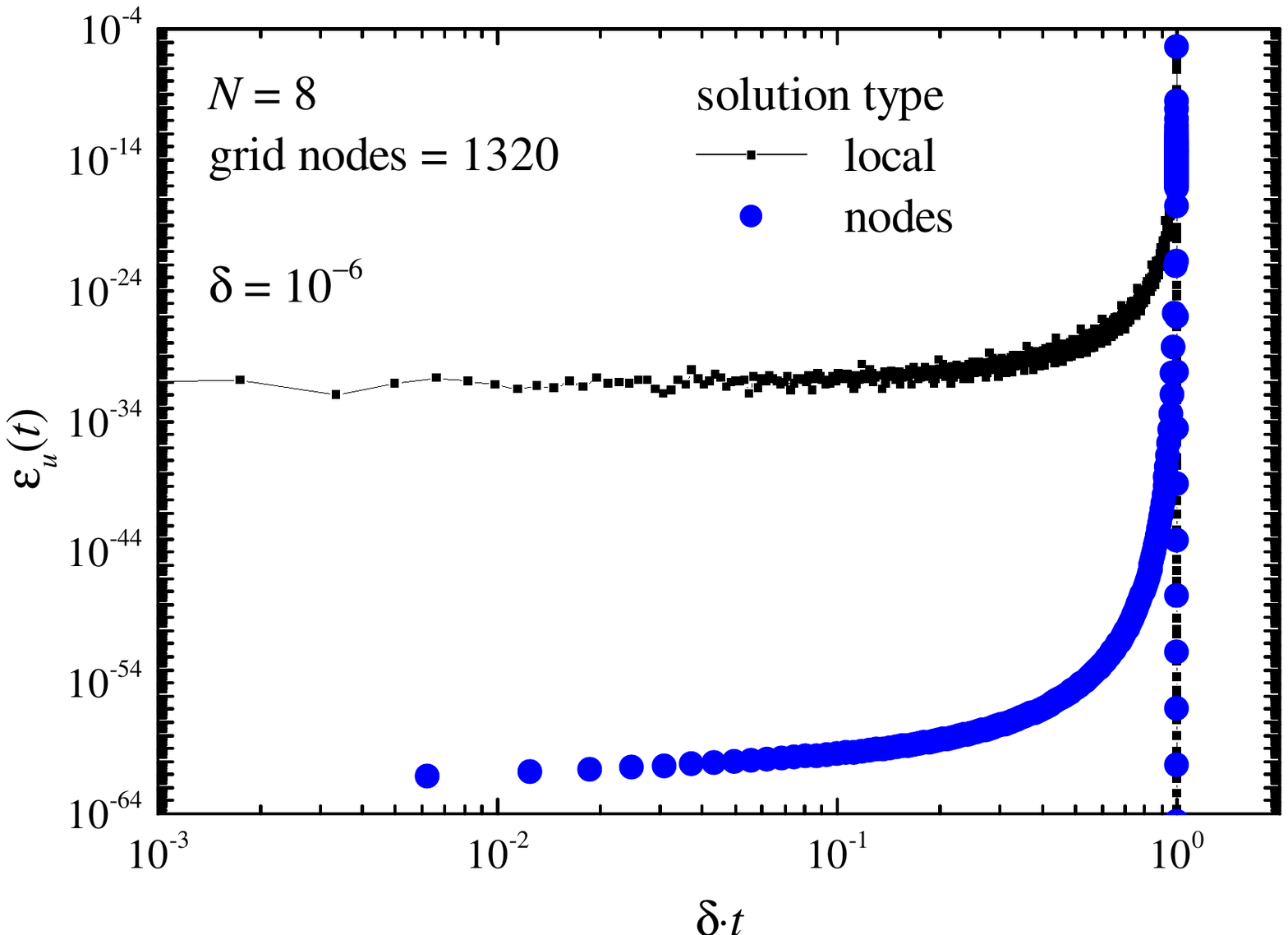}
\vspace{-8mm}\caption{\label{fig:shampine_flame_delta_10m6_sol_v_epss:c2}}
\end{subfigure}
\begin{subfigure}{0.320\textwidth}
\includegraphics[width=\textwidth]{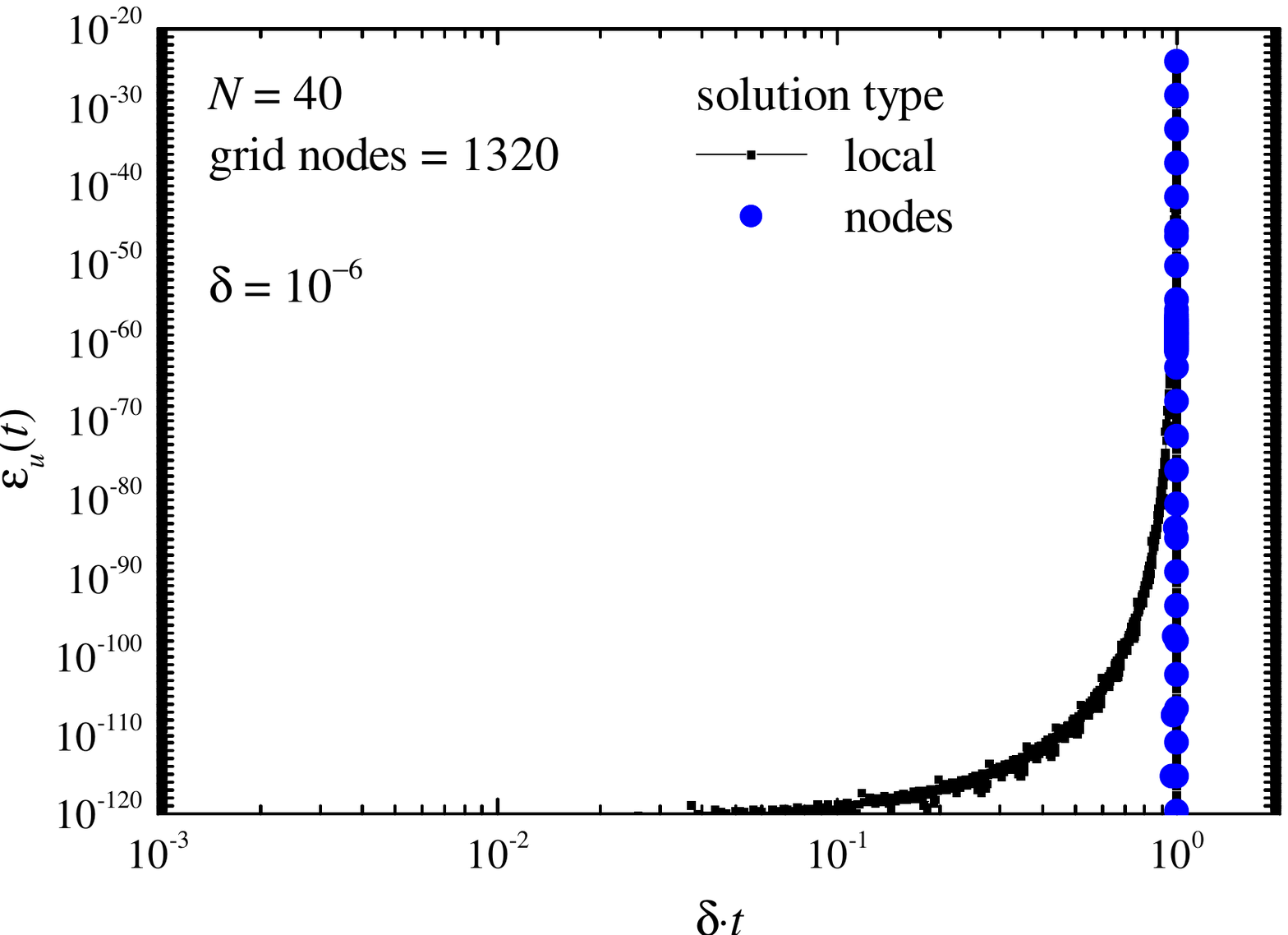}
\vspace{-8mm}\caption{\label{fig:shampine_flame_delta_10m6_sol_v_epss:c3}}
\end{subfigure}\\
\begin{subfigure}{0.320\textwidth}
\includegraphics[width=\textwidth]{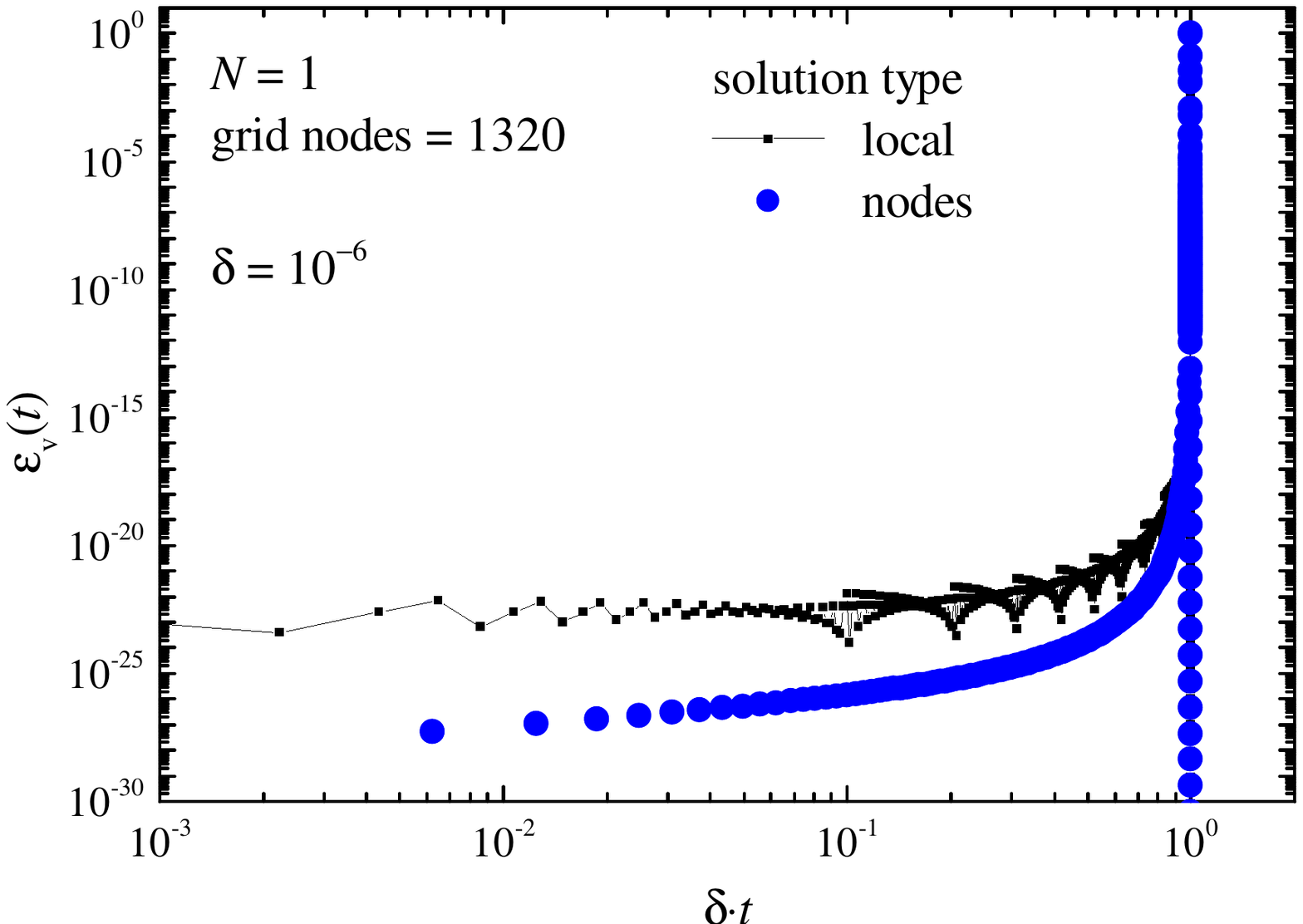}
\vspace{-8mm}\caption{\label{fig:shampine_flame_delta_10m6_sol_v_epss:d1}}
\end{subfigure}
\begin{subfigure}{0.320\textwidth}
\includegraphics[width=\textwidth]{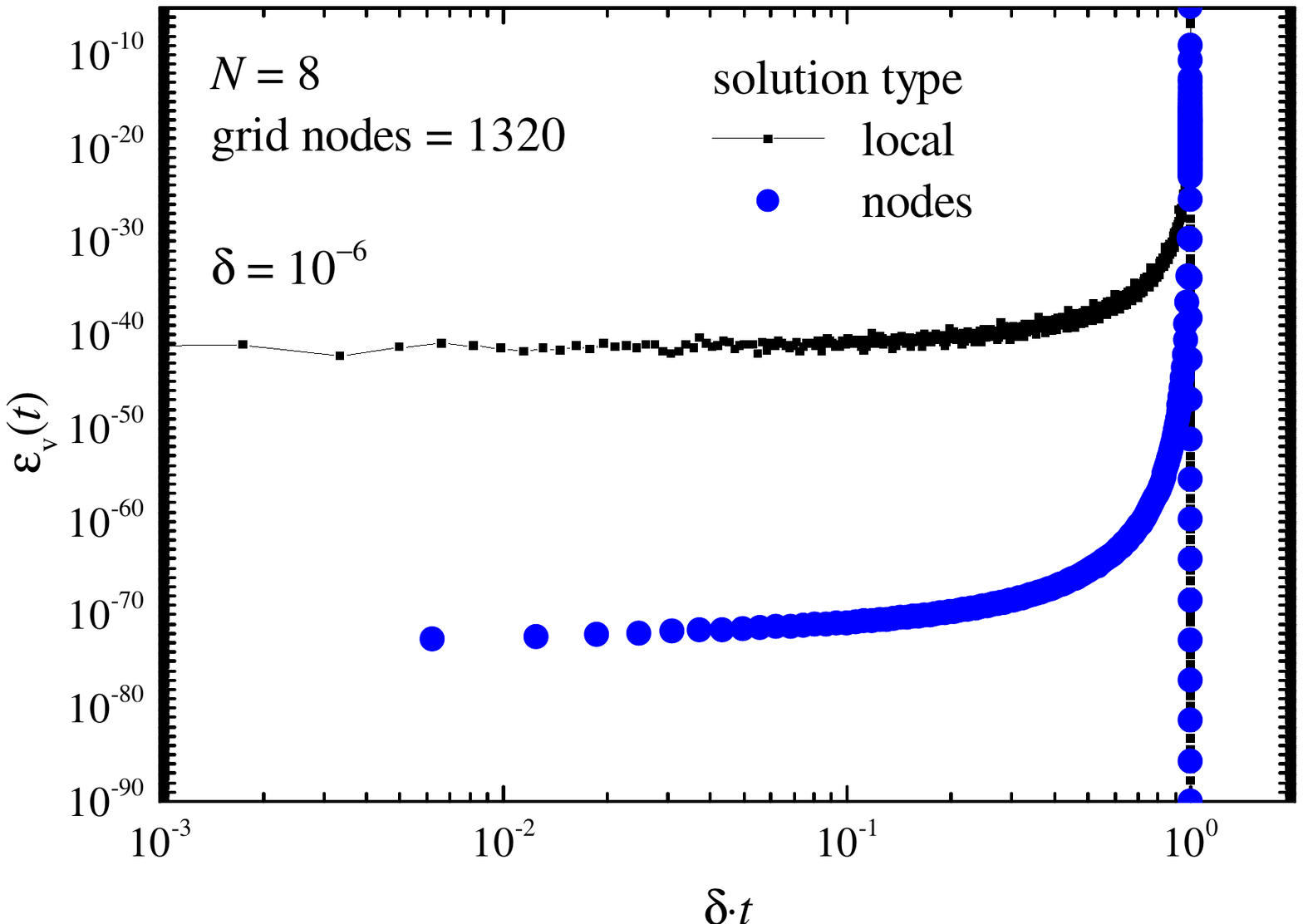}
\vspace{-8mm}\caption{\label{fig:shampine_flame_delta_10m6_sol_v_epss:d2}}
\end{subfigure}
\begin{subfigure}{0.320\textwidth}
\includegraphics[width=\textwidth]{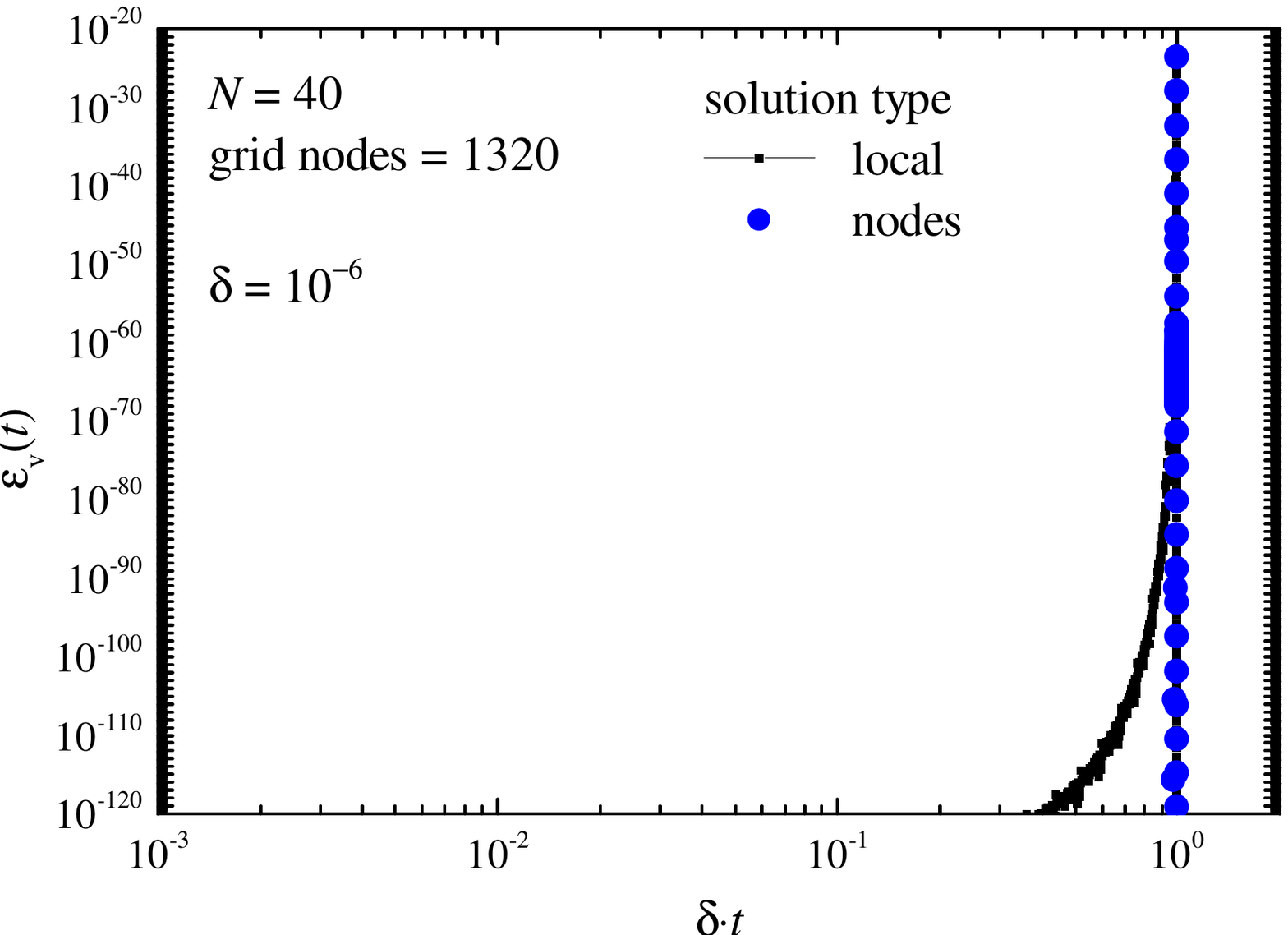}
\vspace{-8mm}\caption{\label{fig:shampine_flame_delta_10m6_sol_v_epss:d3}}
\end{subfigure}\\
\begin{subfigure}{0.320\textwidth}
\includegraphics[width=\textwidth]{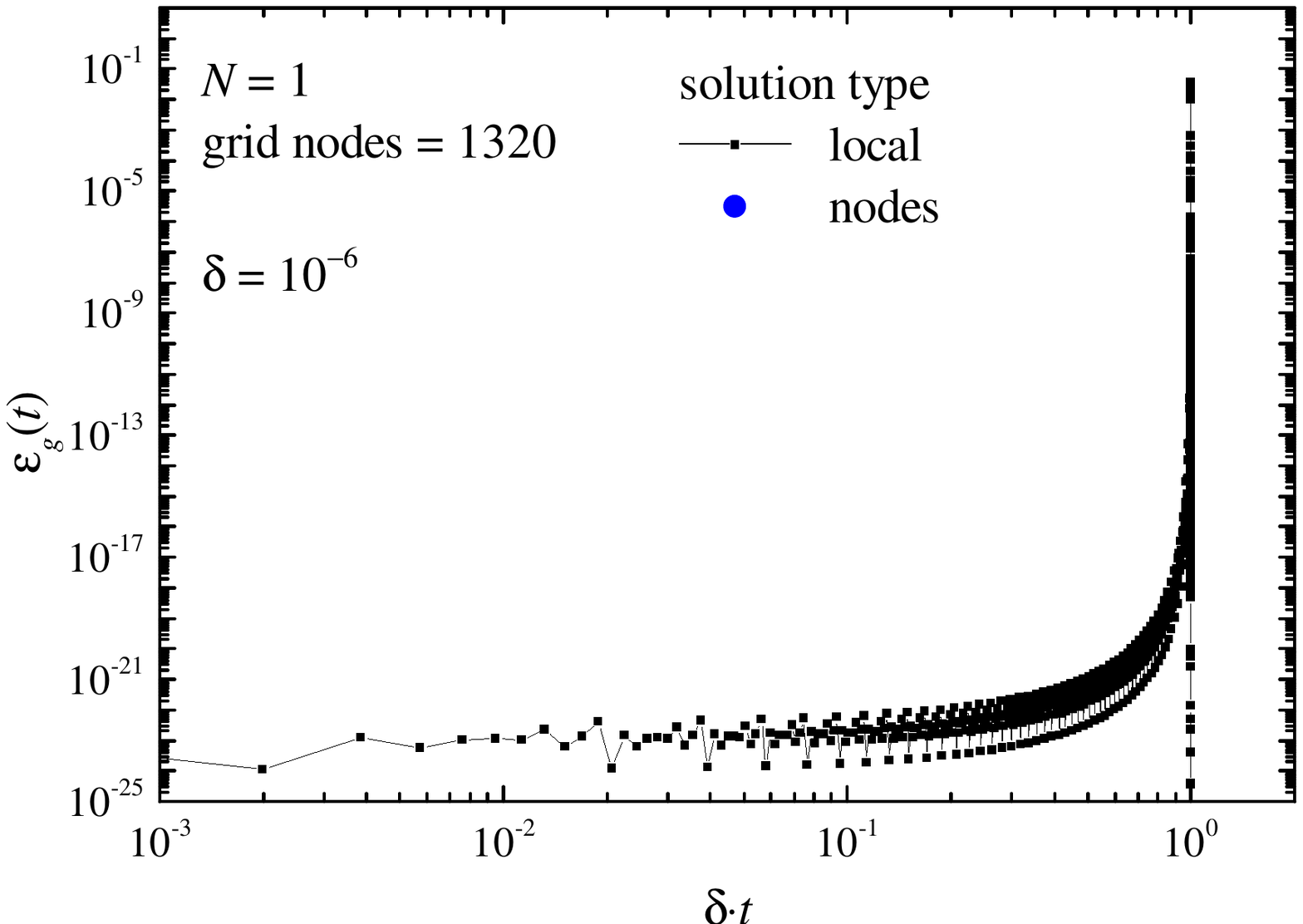}
\vspace{-8mm}\caption{\label{fig:shampine_flame_delta_10m6_sol_v_epss:e1}}
\end{subfigure}
\begin{subfigure}{0.320\textwidth}
\includegraphics[width=\textwidth]{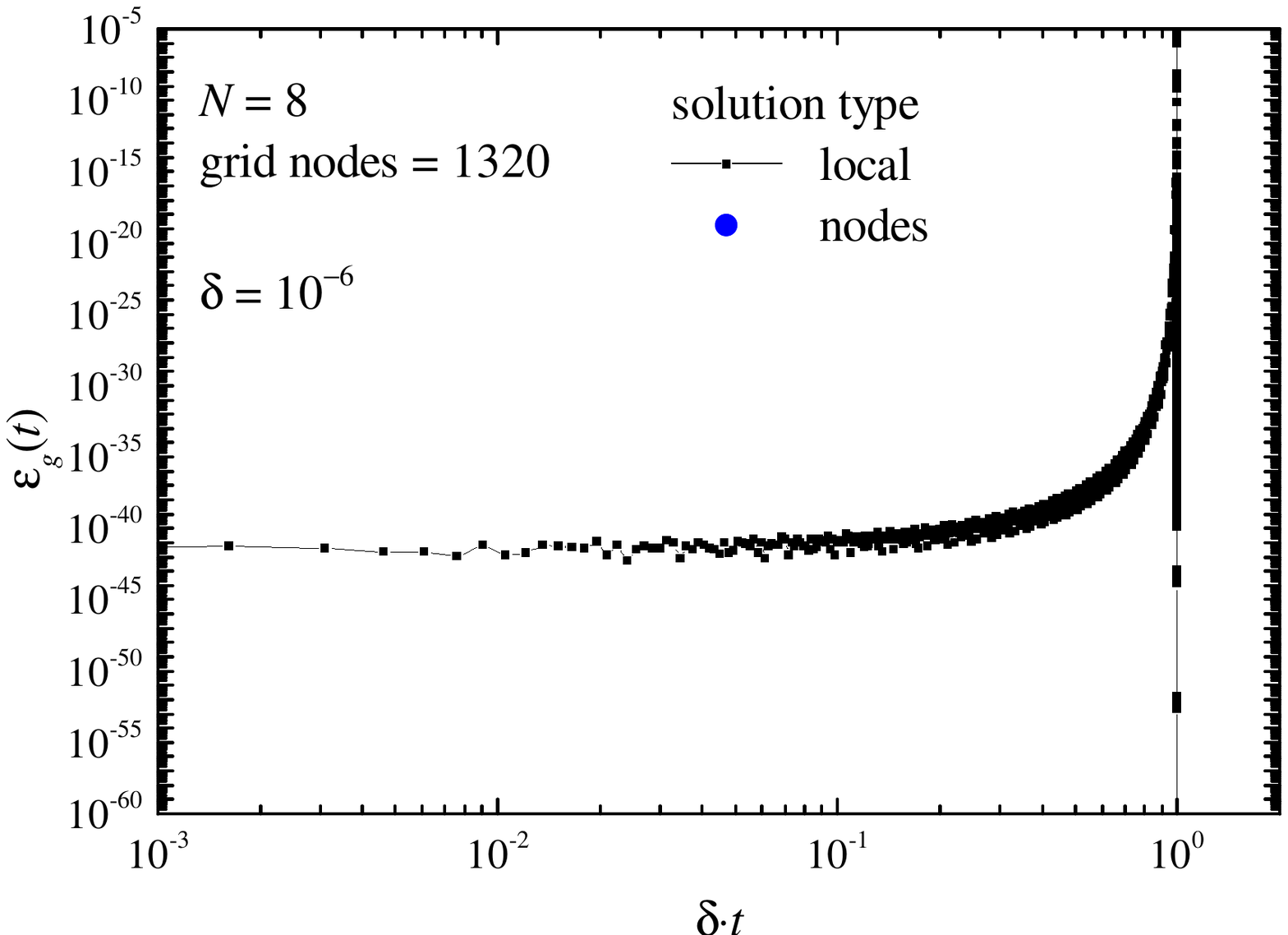}
\vspace{-8mm}\caption{\label{fig:shampine_flame_delta_10m6_sol_v_epss:e2}}
\end{subfigure}
\begin{subfigure}{0.320\textwidth}
\includegraphics[width=\textwidth]{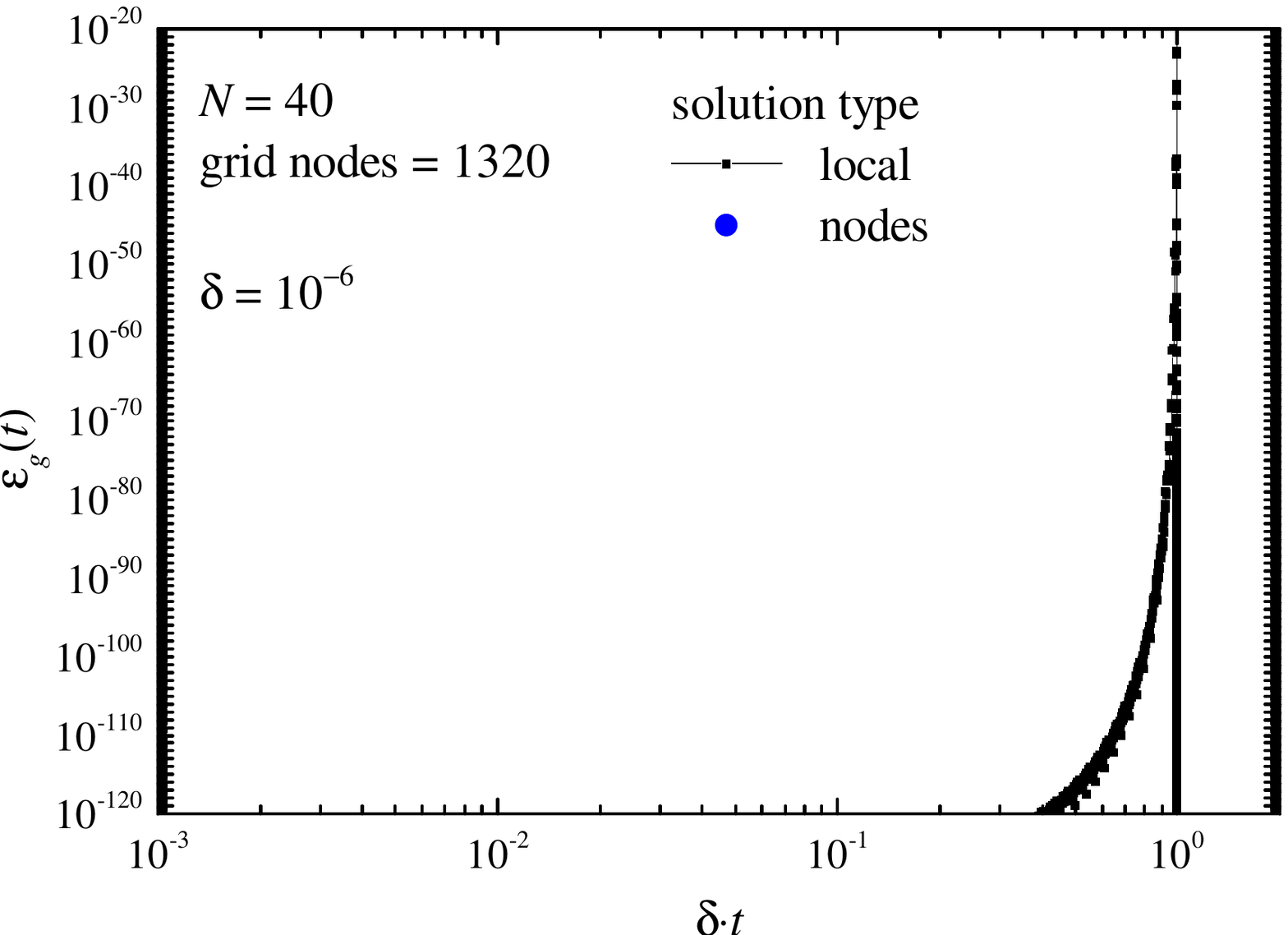}
\vspace{-8mm}\caption{\label{fig:shampine_flame_delta_10m6_sol_v_epss:e3}}
\end{subfigure}\\
\caption{%
Numerical solution of the stiff DAE system (\ref{eq:shampine_flame}) of index 1 with $\delta = 10^{-6}$. Comparison of the solution at nodes $\mathbf{v}_{n}$, the local solution $\mathbf{v}_{L}(t)$ and the exact solution $\mathbf{v}^{\rm ex}(t)$ for component $v_{1}$ (\subref{fig:shampine_flame_delta_10m6_sol_v_epss:a1}, \subref{fig:shampine_flame_delta_10m6_sol_v_epss:a2}, \subref{fig:shampine_flame_delta_10m6_sol_v_epss:a3}, \subref{fig:shampine_flame_delta_10m6_sol_v_epss:b1}, \subref{fig:shampine_flame_delta_10m6_sol_v_epss:b2}, \subref{fig:shampine_flame_delta_10m6_sol_v_epss:b3}), the errors $\varepsilon_{u}(t)$ (\subref{fig:shampine_flame_delta_10m6_sol_v_epss:c1}, \subref{fig:shampine_flame_delta_10m6_sol_v_epss:c2}, \subref{fig:shampine_flame_delta_10m6_sol_v_epss:c3}), $\varepsilon_{v}(t)$ (\subref{fig:shampine_flame_delta_10m6_sol_v_epss:d1}, \subref{fig:shampine_flame_delta_10m6_sol_v_epss:d2}, \subref{fig:shampine_flame_delta_10m6_sol_v_epss:d3}), $\varepsilon_{g}(t)$ (\subref{fig:shampine_flame_delta_10m6_sol_qug:e1}, \subref{fig:shampine_flame_delta_10m6_sol_qug:e2}, \subref{fig:shampine_flame_delta_10m6_sol_qug:e3}), obtained using polynomials with degrees $N = 1$ (\subref{fig:shampine_flame_delta_10m6_sol_qug:a1}, \subref{fig:shampine_flame_delta_10m6_sol_qug:b1}, \subref{fig:shampine_flame_delta_10m6_sol_qug:c1}, \subref{fig:shampine_flame_delta_10m6_sol_qug:d1}, \subref{fig:shampine_flame_delta_10m6_sol_qug:e1}), $N = 8$ (\subref{fig:shampine_flame_delta_10m6_sol_qug:a2}, \subref{fig:shampine_flame_delta_10m6_sol_qug:b2}, \subref{fig:shampine_flame_delta_10m6_sol_qug:c2}, \subref{fig:shampine_flame_delta_10m6_sol_qug:d2}, \subref{fig:shampine_flame_delta_10m6_sol_qug:e2}) and $N = 40$ (\subref{fig:shampine_flame_delta_10m6_sol_qug:a3}, \subref{fig:shampine_flame_delta_10m6_sol_qug:b3}, \subref{fig:shampine_flame_delta_10m6_sol_qug:c3}, \subref{fig:shampine_flame_delta_10m6_sol_qug:d3}, \subref{fig:shampine_flame_delta_10m6_sol_qug:e3}).
}
\label{fig:shampine_flame_delta_10m6_sol_v_epss}
\end{figure}

\begin{figure}[h!]
\captionsetup[subfigure]{%
	position=bottom,
	font+=smaller,
	textfont=normalfont,
	singlelinecheck=off,
	justification=raggedright
}
\centering
\begin{subfigure}{0.320\textwidth}
\includegraphics[width=\textwidth]{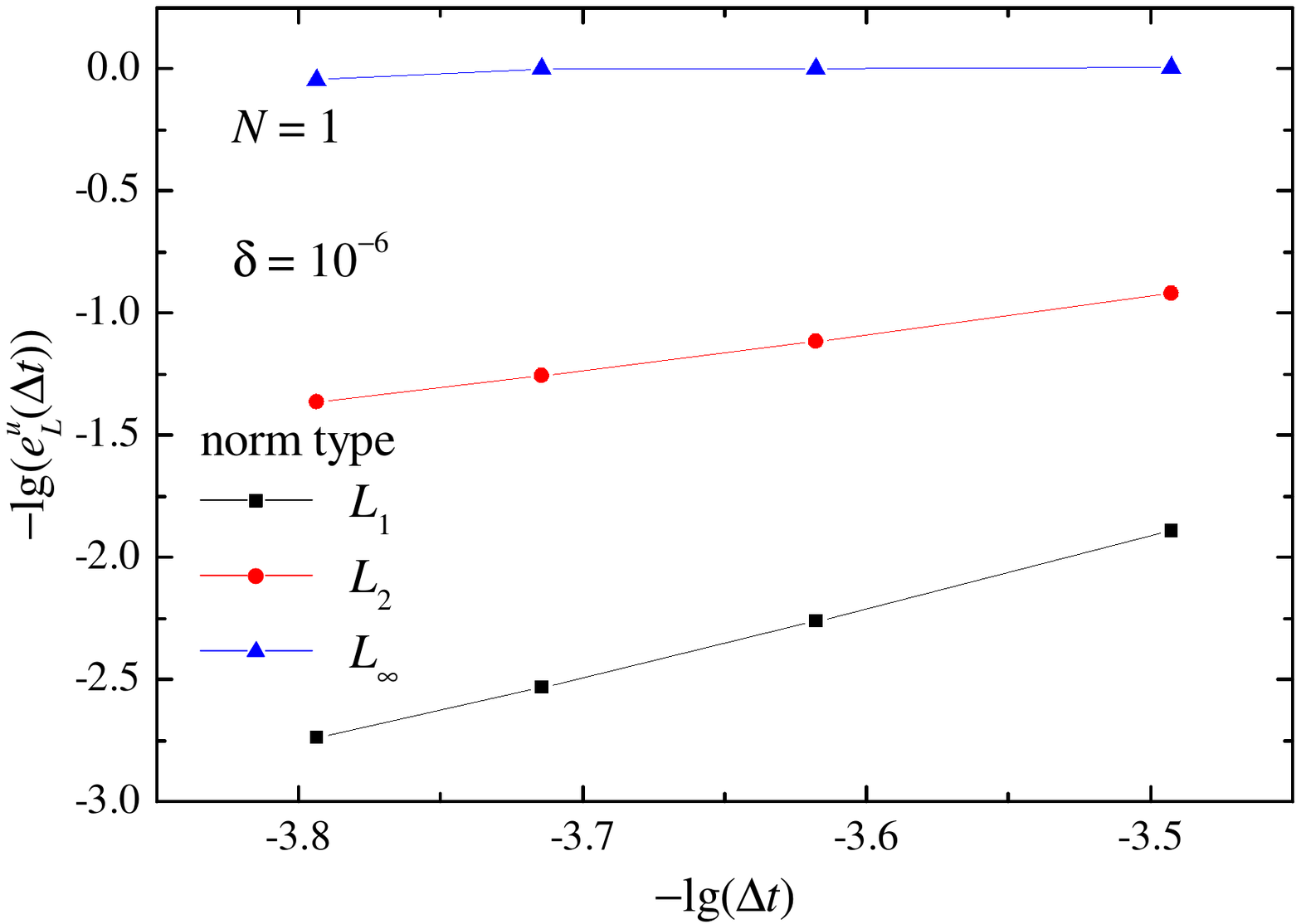}
\vspace{-8mm}\caption{\label{fig:shampine_flame_delta_10m6_errors:a1}}
\end{subfigure}
\begin{subfigure}{0.320\textwidth}
\includegraphics[width=\textwidth]{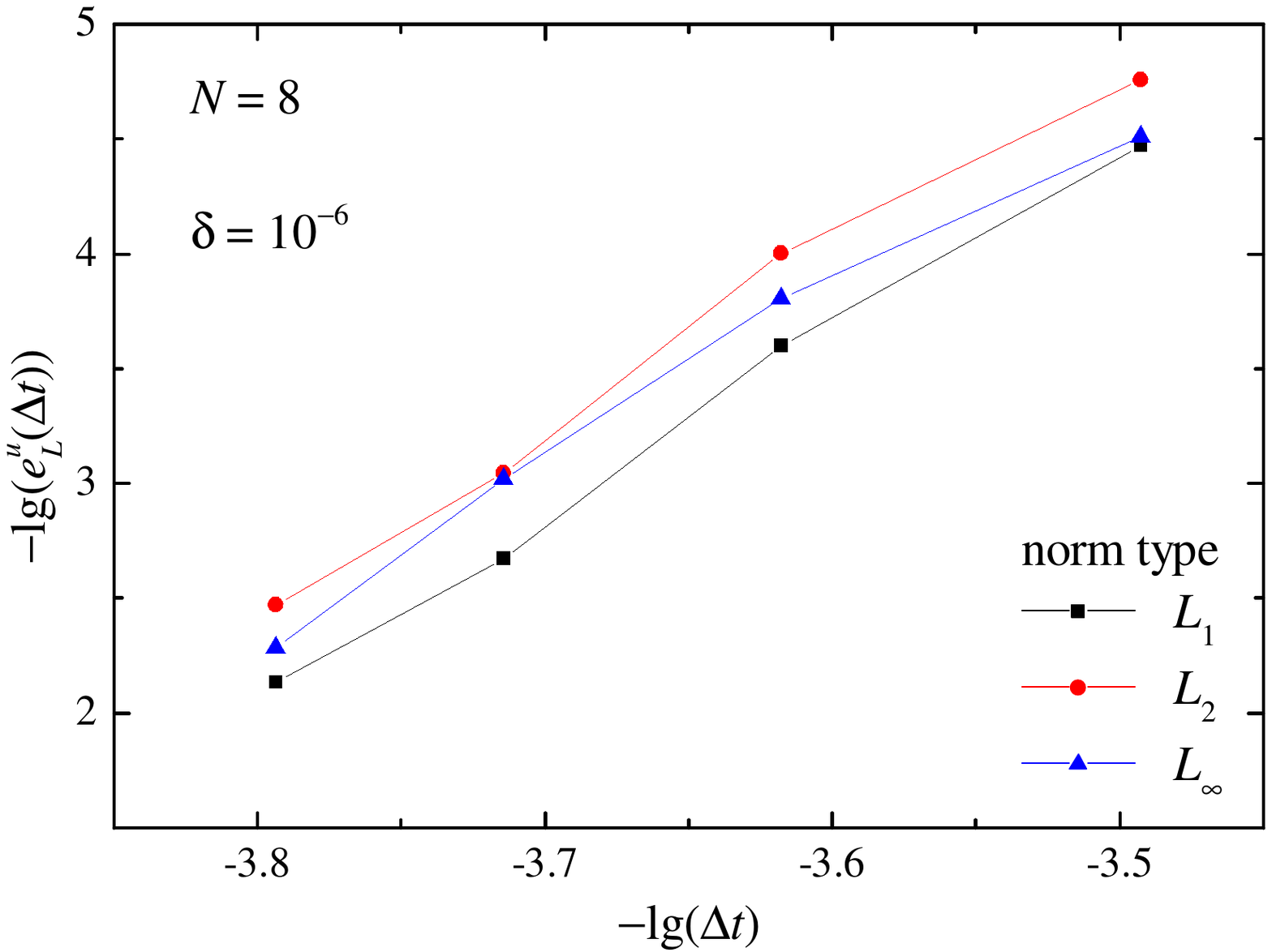}
\vspace{-8mm}\caption{\label{fig:shampine_flame_delta_10m6_errors:a2}}
\end{subfigure}
\begin{subfigure}{0.320\textwidth}
\includegraphics[width=\textwidth]{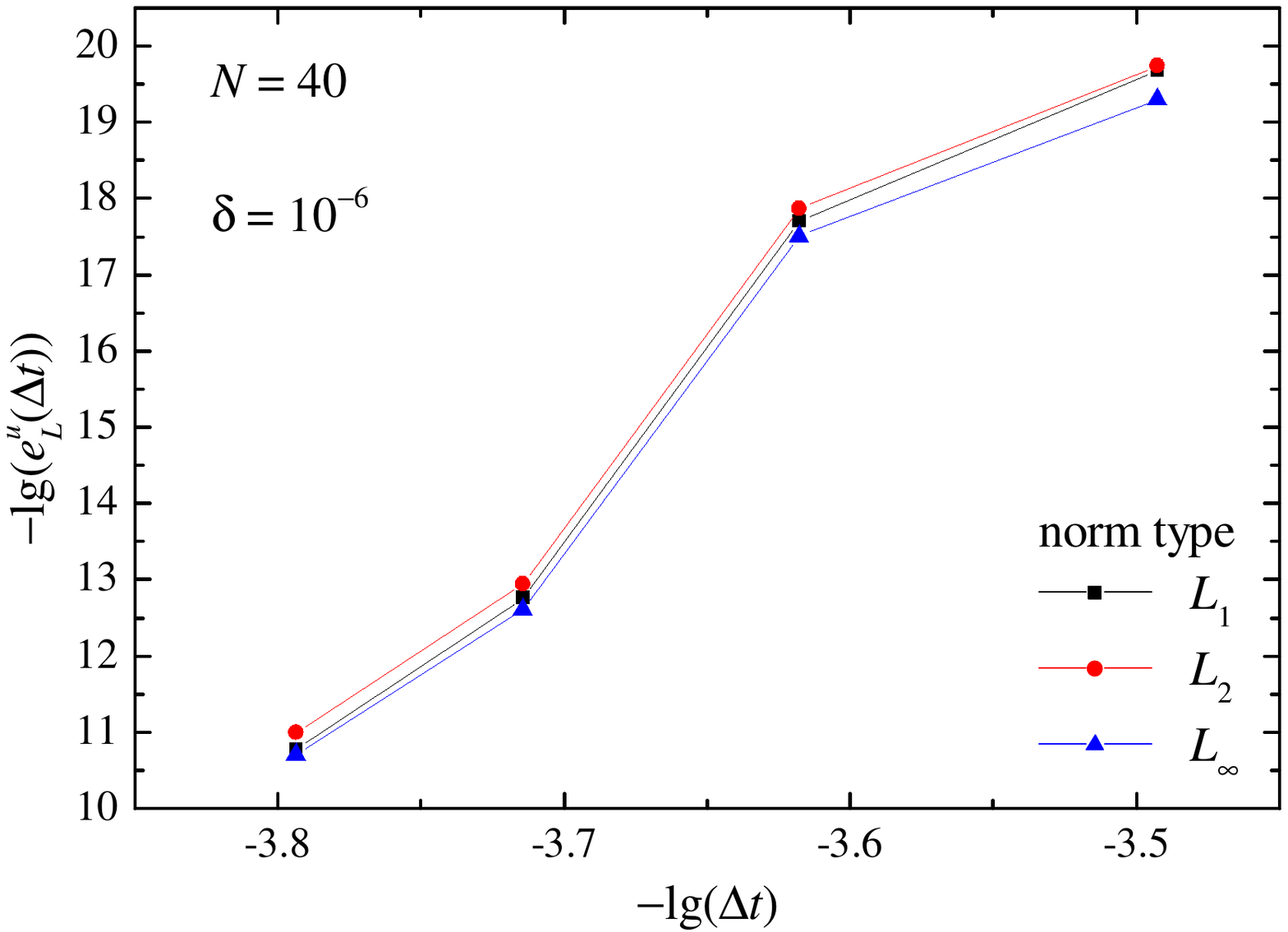}
\vspace{-8mm}\caption{\label{fig:shampine_flame_delta_10m6_errors:a3}}
\end{subfigure}\\
\begin{subfigure}{0.320\textwidth}
\includegraphics[width=\textwidth]{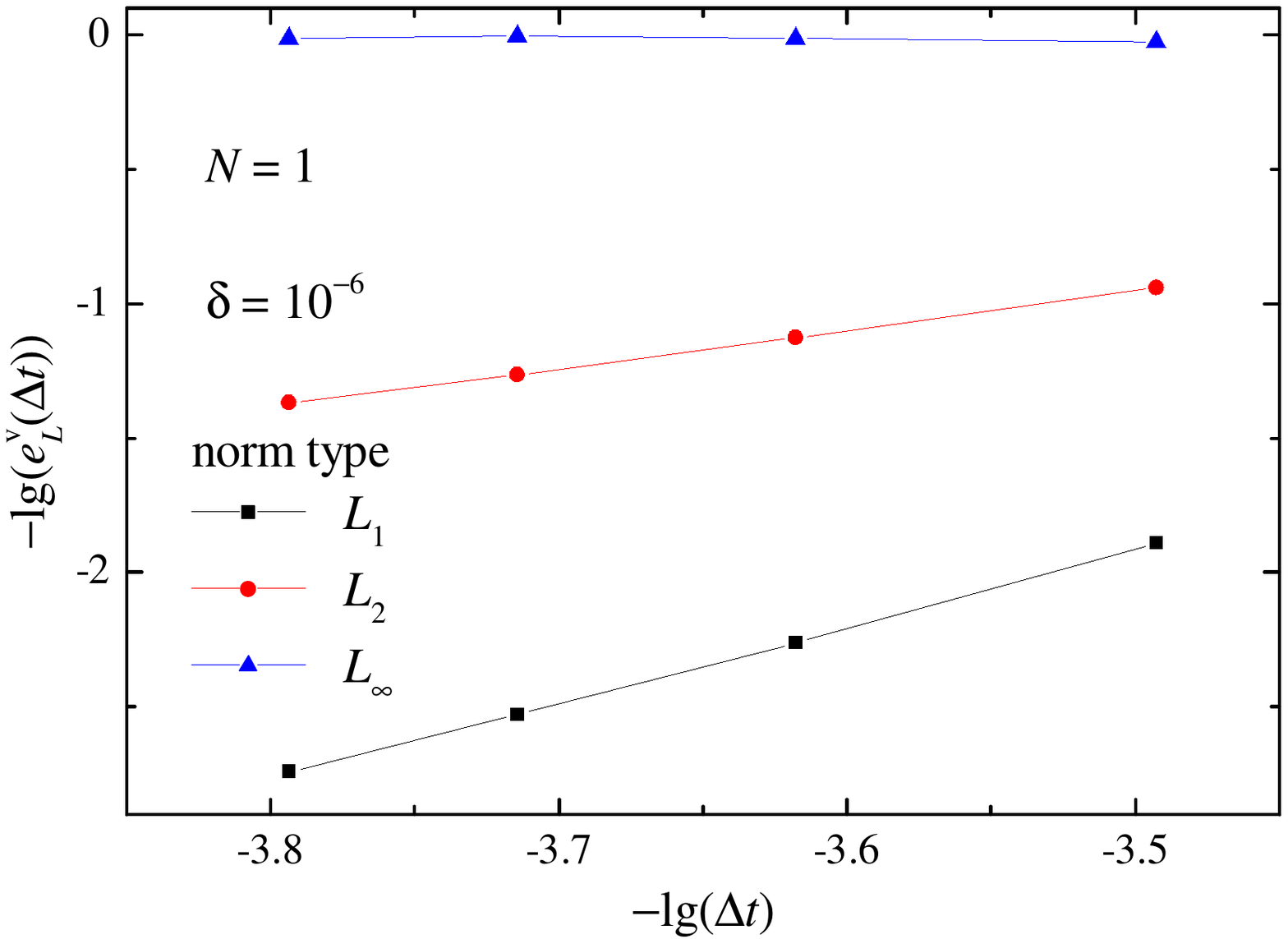}
\vspace{-8mm}\caption{\label{fig:shampine_flame_delta_10m6_errors:b1}}
\end{subfigure}
\begin{subfigure}{0.320\textwidth}
\includegraphics[width=\textwidth]{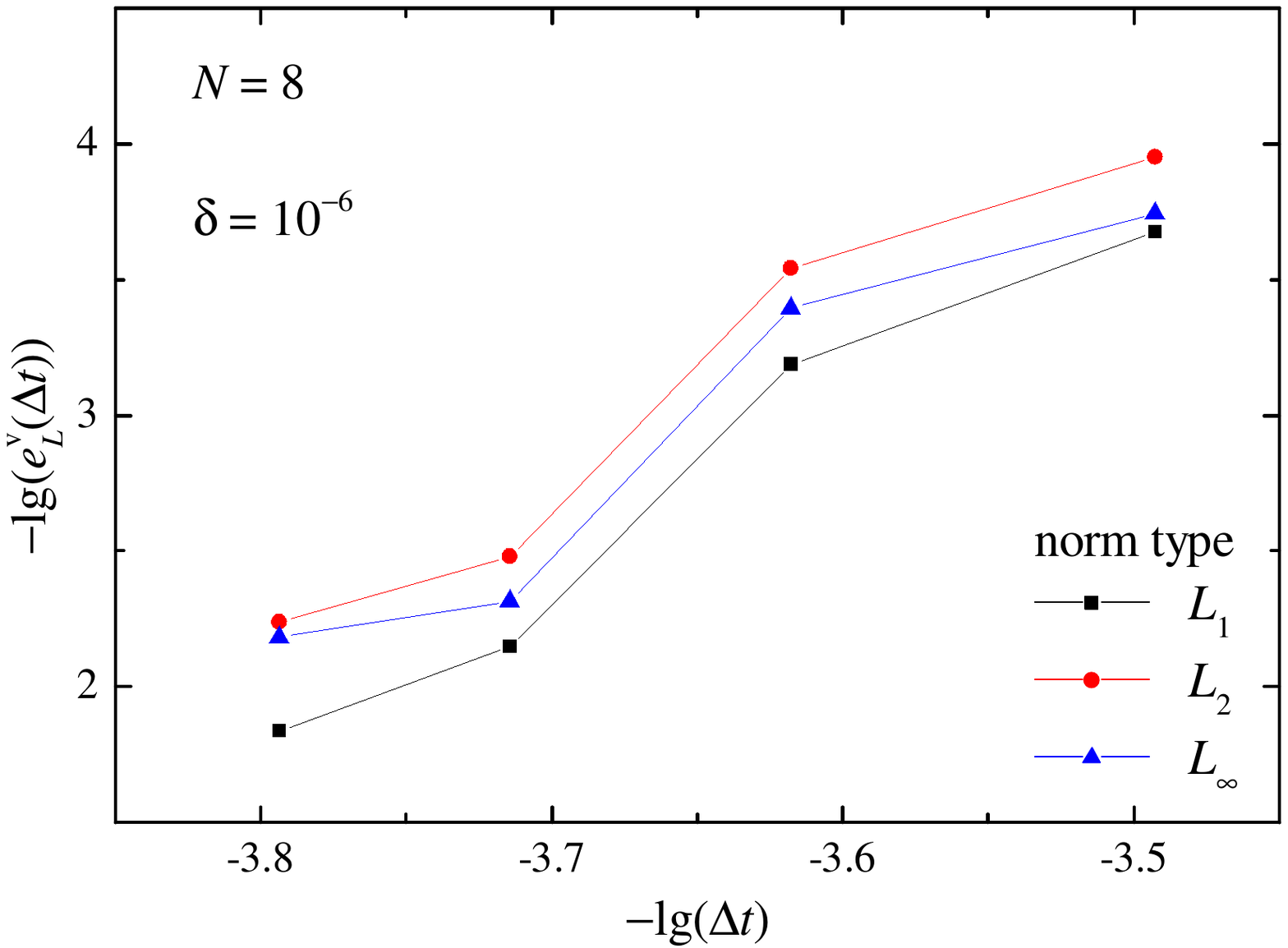}
\vspace{-8mm}\caption{\label{fig:shampine_flame_delta_10m6_errors:b2}}
\end{subfigure}
\begin{subfigure}{0.320\textwidth}
\includegraphics[width=\textwidth]{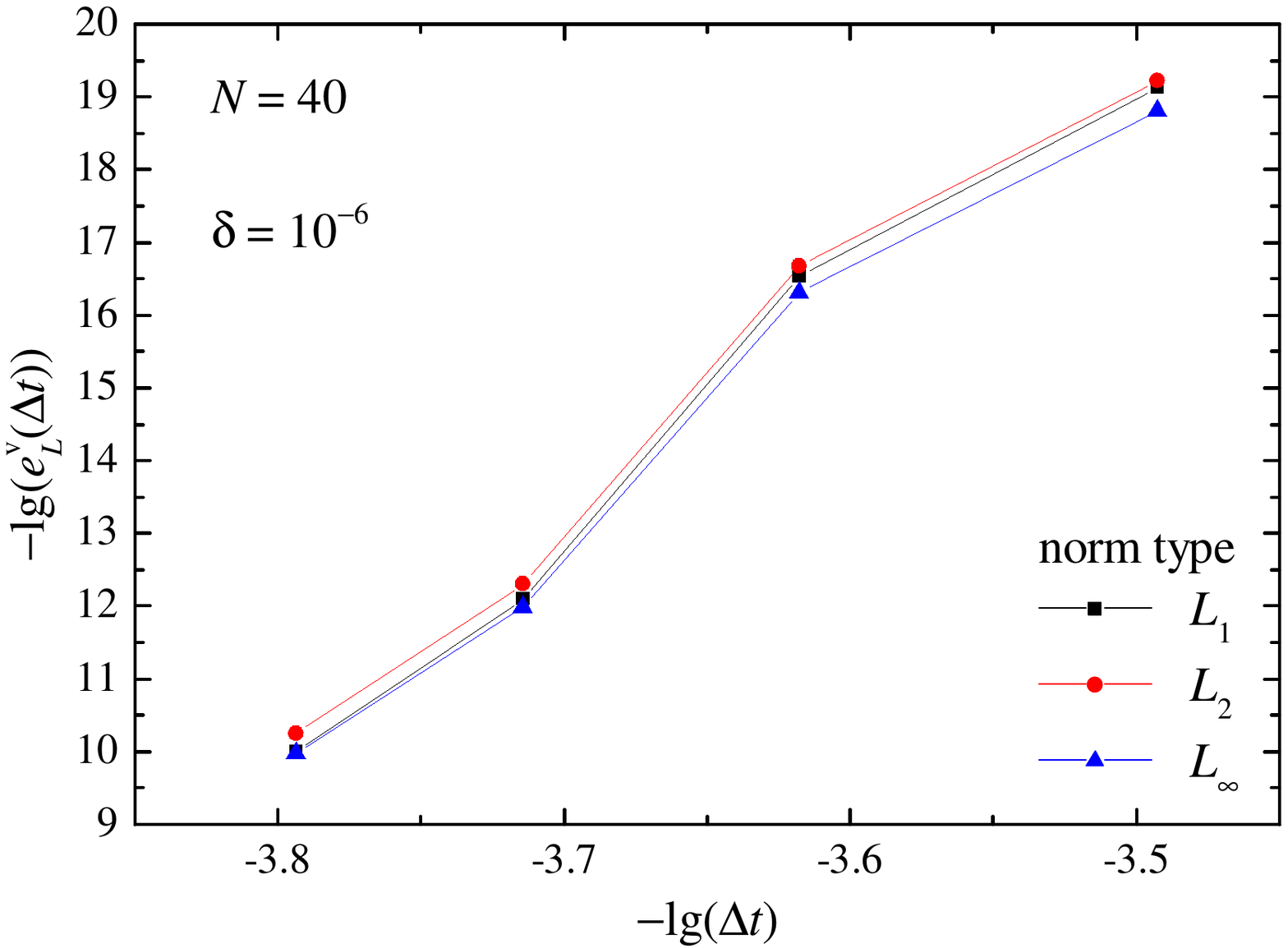}
\vspace{-8mm}\caption{\label{fig:shampine_flame_delta_10m6_errors:b3}}
\end{subfigure}\\
\begin{subfigure}{0.320\textwidth}
\includegraphics[width=\textwidth]{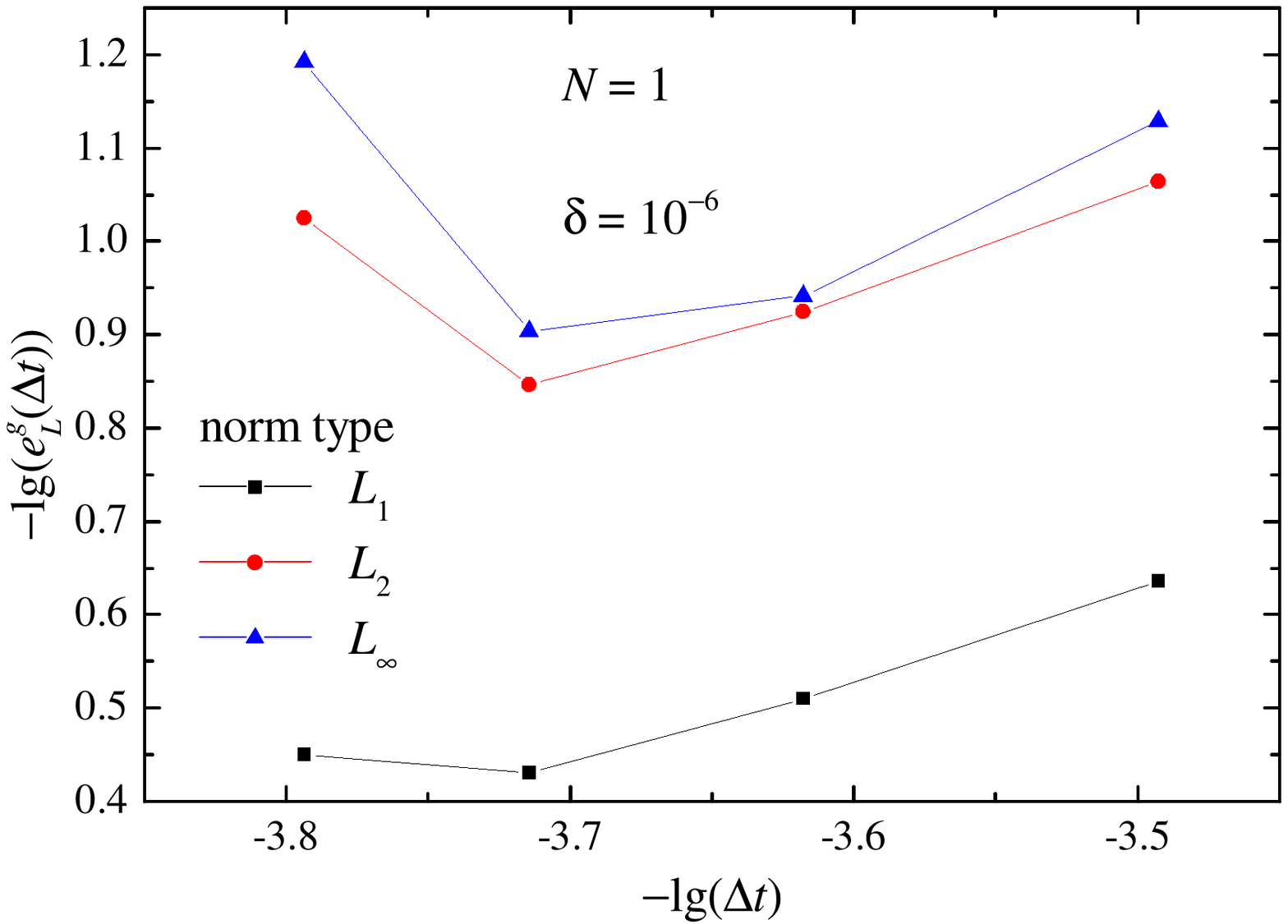}
\vspace{-8mm}\caption{\label{fig:shampine_flame_delta_10m6_errors:c1}}
\end{subfigure}
\begin{subfigure}{0.320\textwidth}
\includegraphics[width=\textwidth]{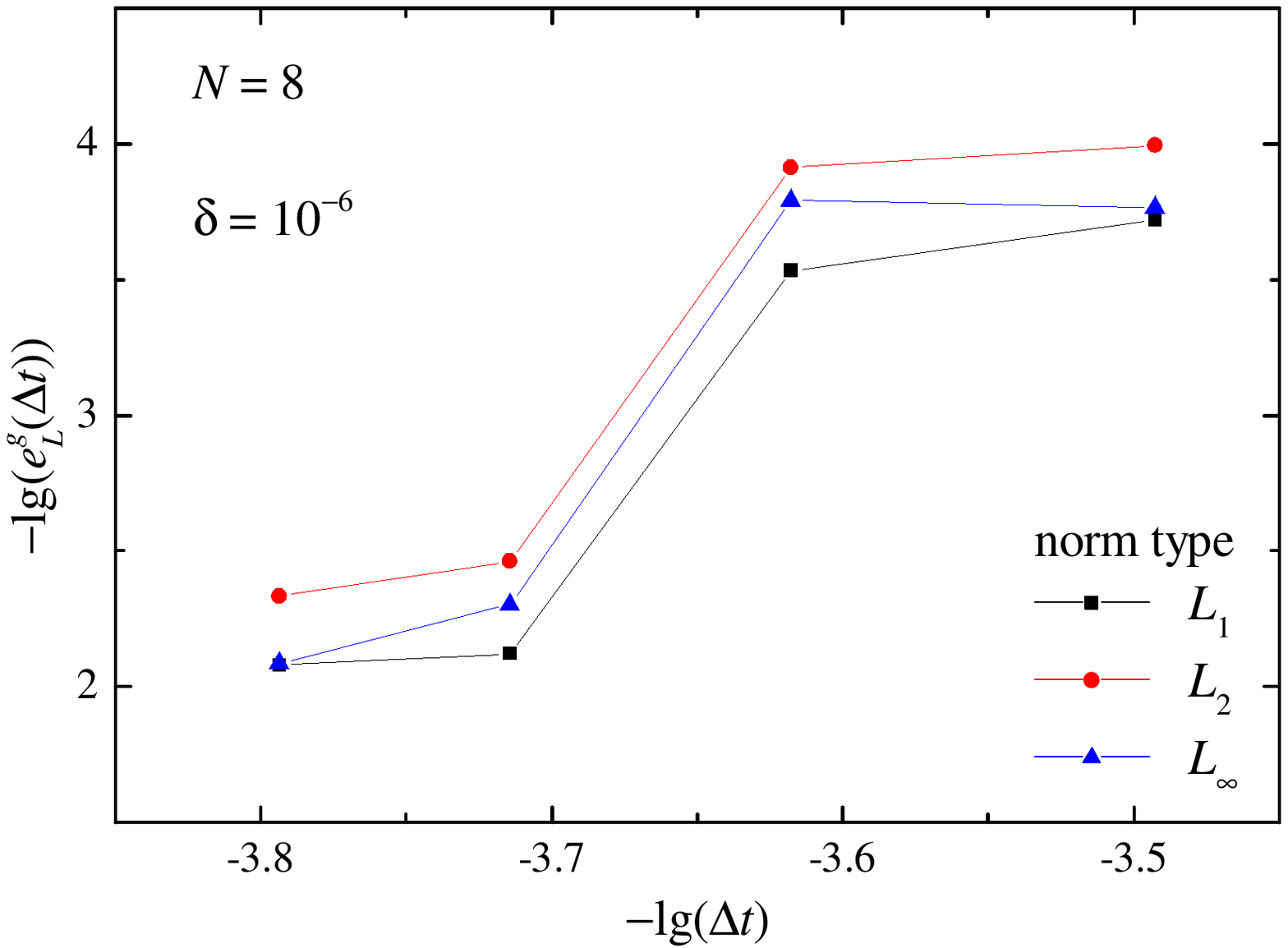}
\vspace{-8mm}\caption{\label{fig:shampine_flame_delta_10m6_errors:c2}}
\end{subfigure}
\begin{subfigure}{0.320\textwidth}
\includegraphics[width=\textwidth]{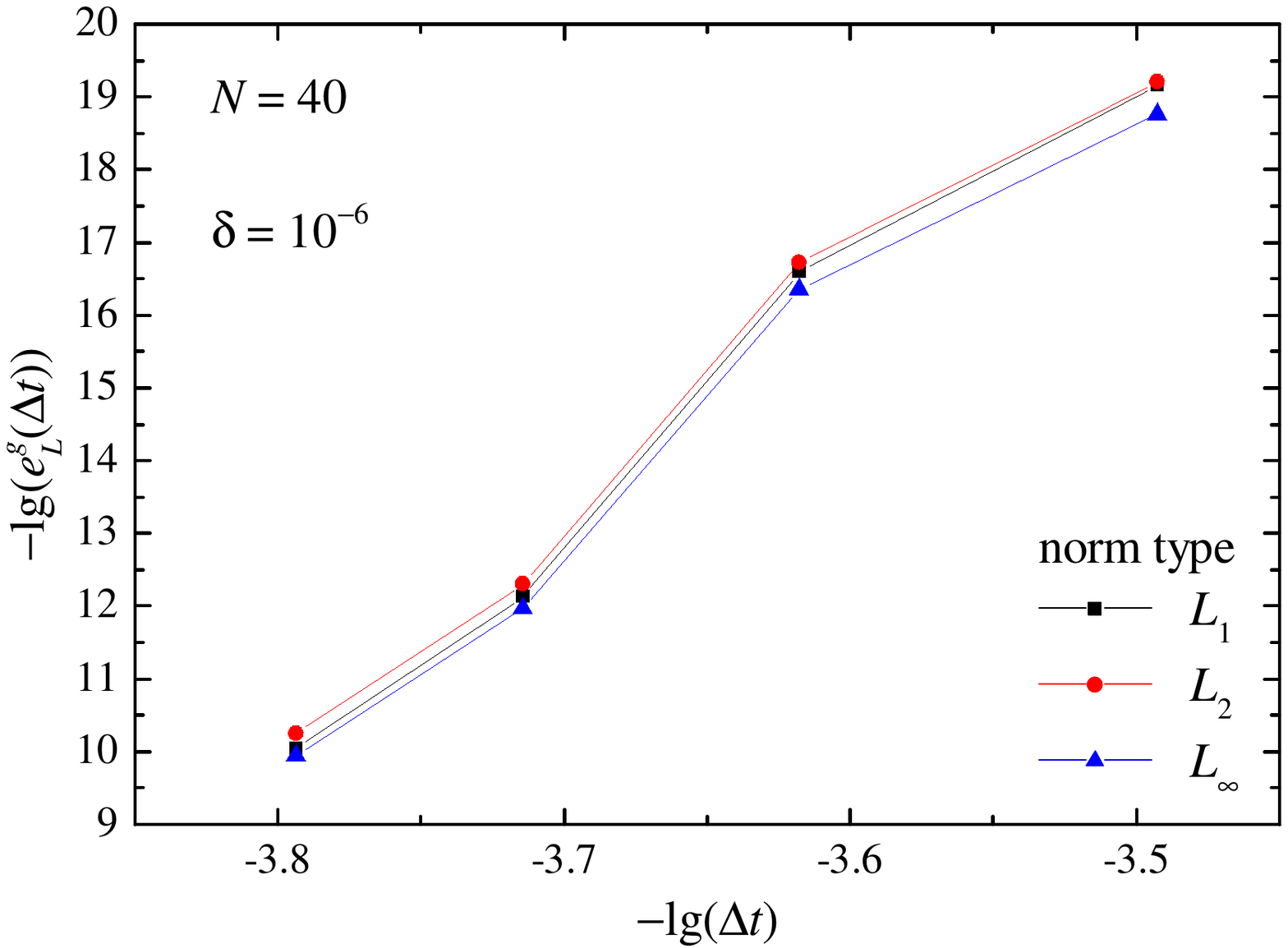}
\vspace{-8mm}\caption{\label{fig:shampine_flame_delta_10m6_errors:c3}}
\end{subfigure}\\
\begin{subfigure}{0.320\textwidth}
\includegraphics[width=\textwidth]{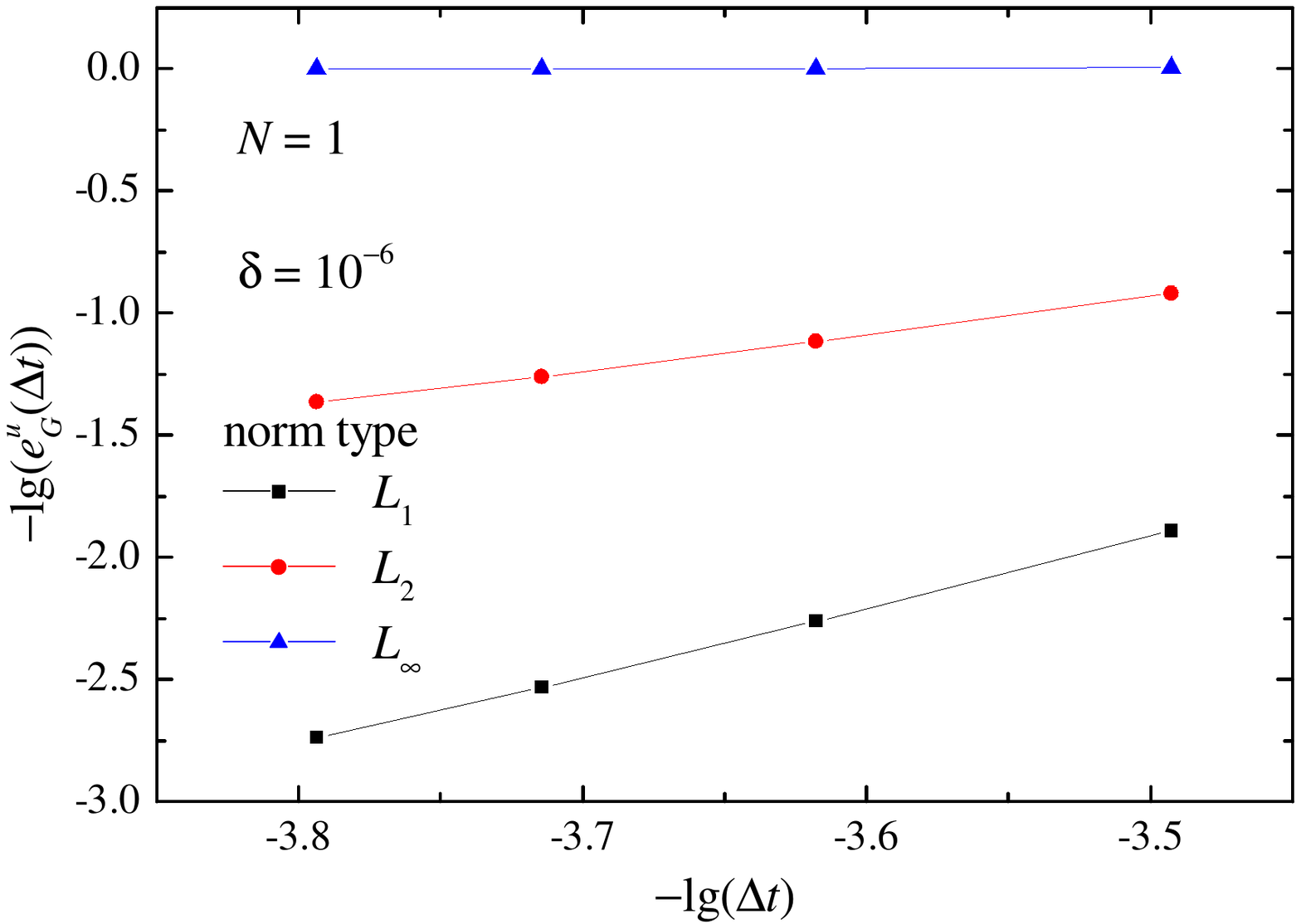}
\vspace{-8mm}\caption{\label{fig:shampine_flame_delta_10m6_errors:d1}}
\end{subfigure}
\begin{subfigure}{0.320\textwidth}
\includegraphics[width=\textwidth]{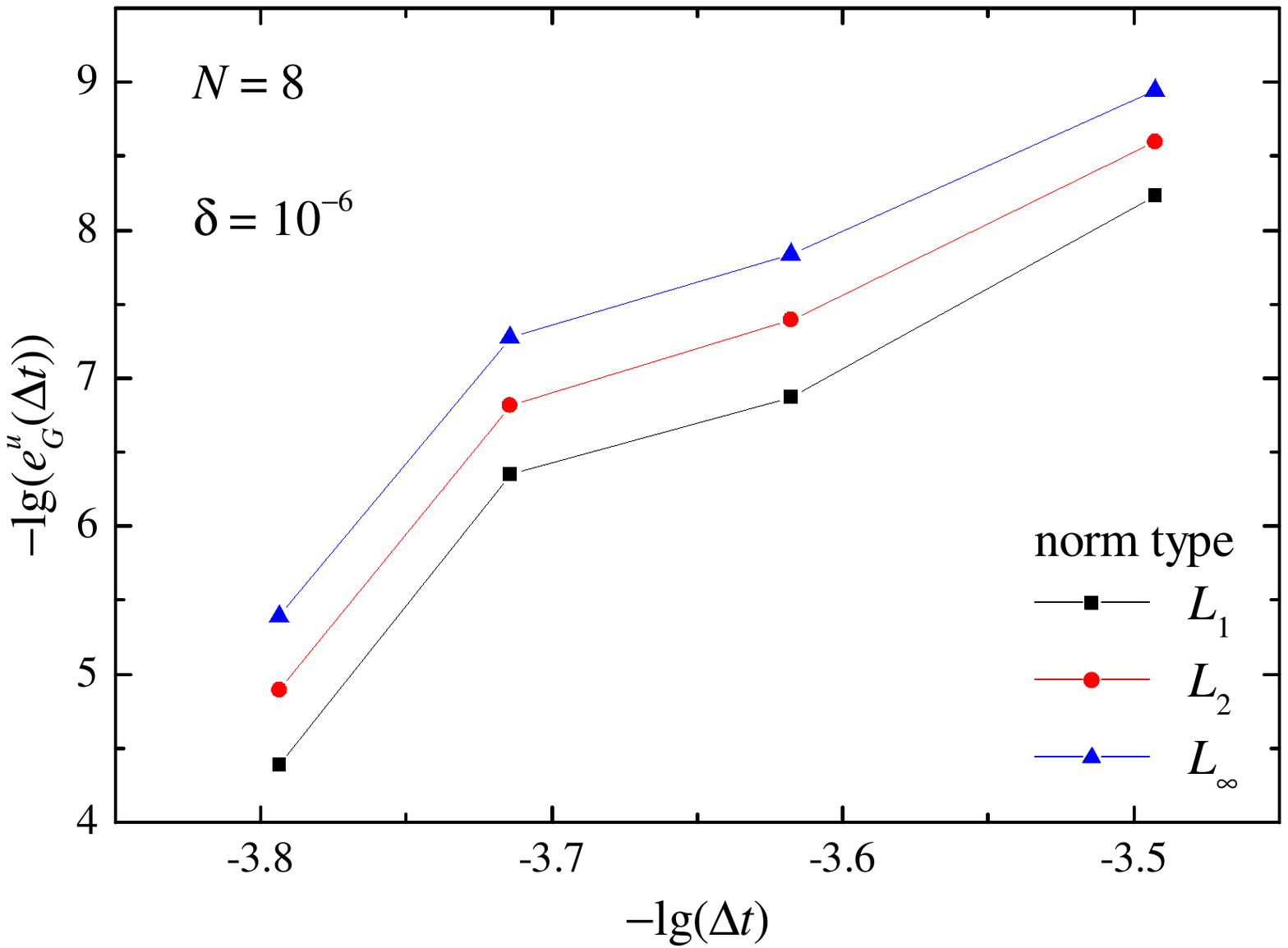}
\vspace{-8mm}\caption{\label{fig:shampine_flame_delta_10m6_errors:d2}}
\end{subfigure}
\begin{subfigure}{0.320\textwidth}
\includegraphics[width=\textwidth]{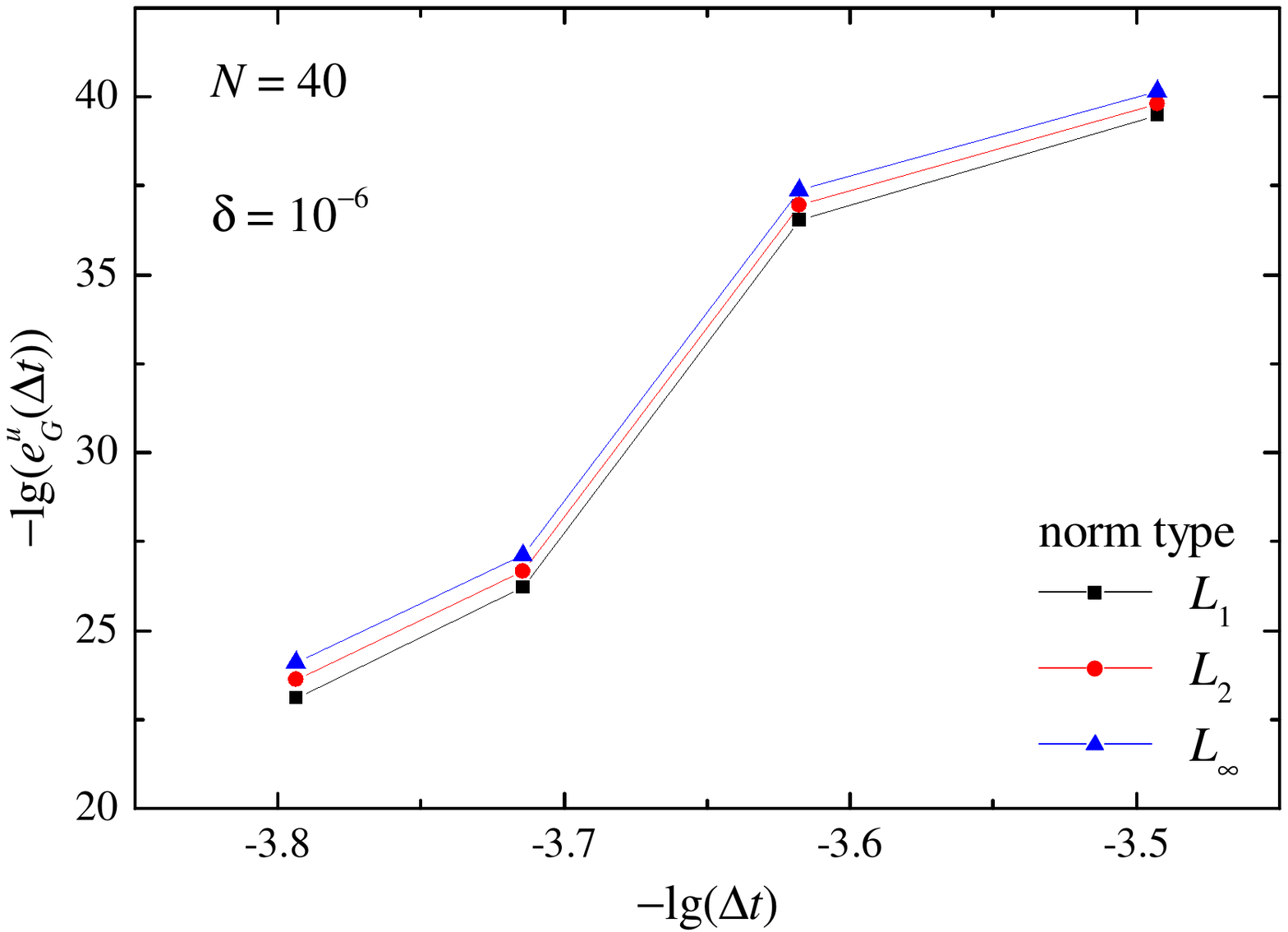}
\vspace{-8mm}\caption{\label{fig:shampine_flame_delta_10m6_errors:d3}}
\end{subfigure}\\
\begin{subfigure}{0.320\textwidth}
\includegraphics[width=\textwidth]{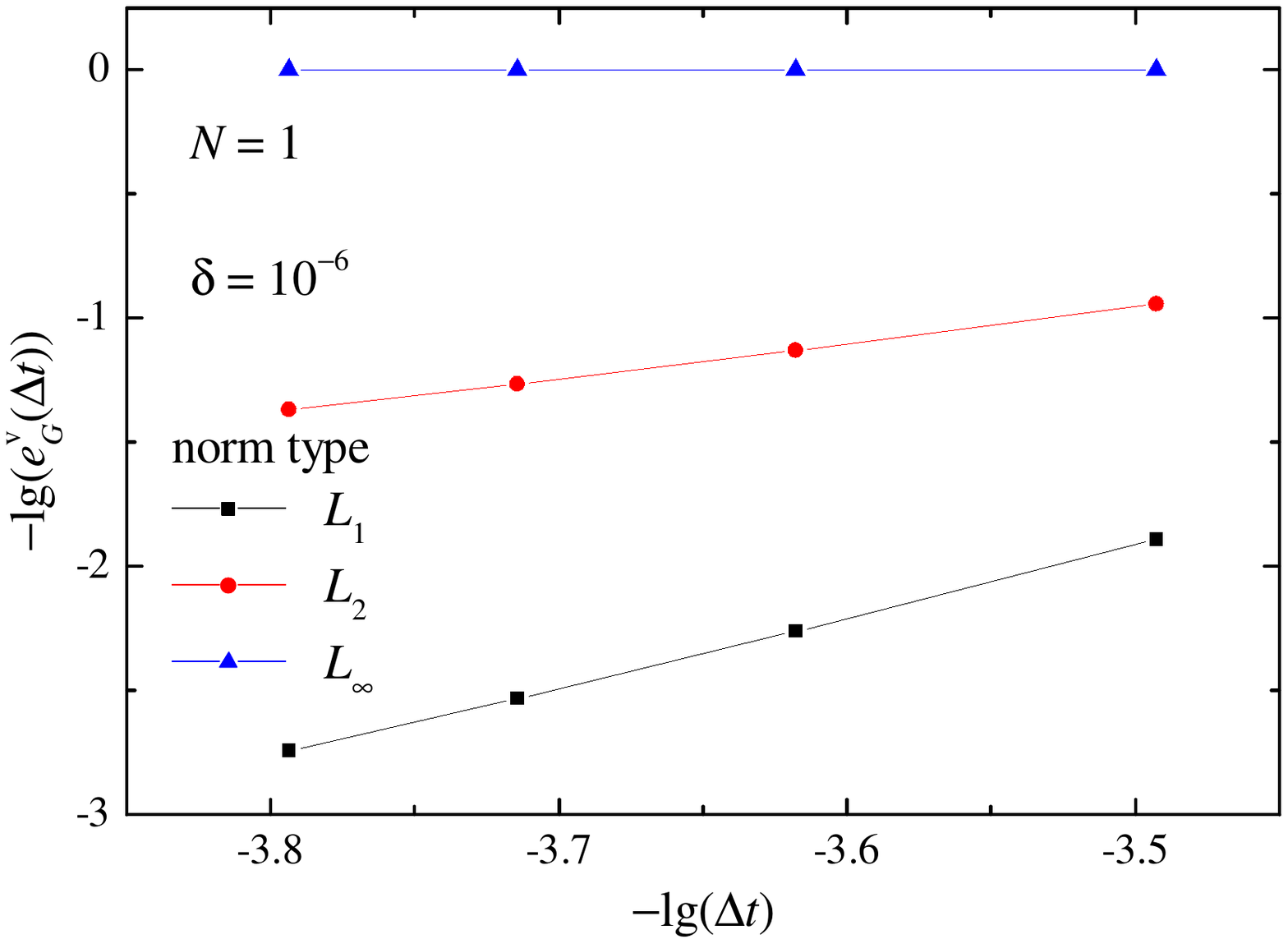}
\vspace{-8mm}\caption{\label{fig:shampine_flame_delta_10m6_errors:e1}}
\end{subfigure}
\begin{subfigure}{0.320\textwidth}
\includegraphics[width=\textwidth]{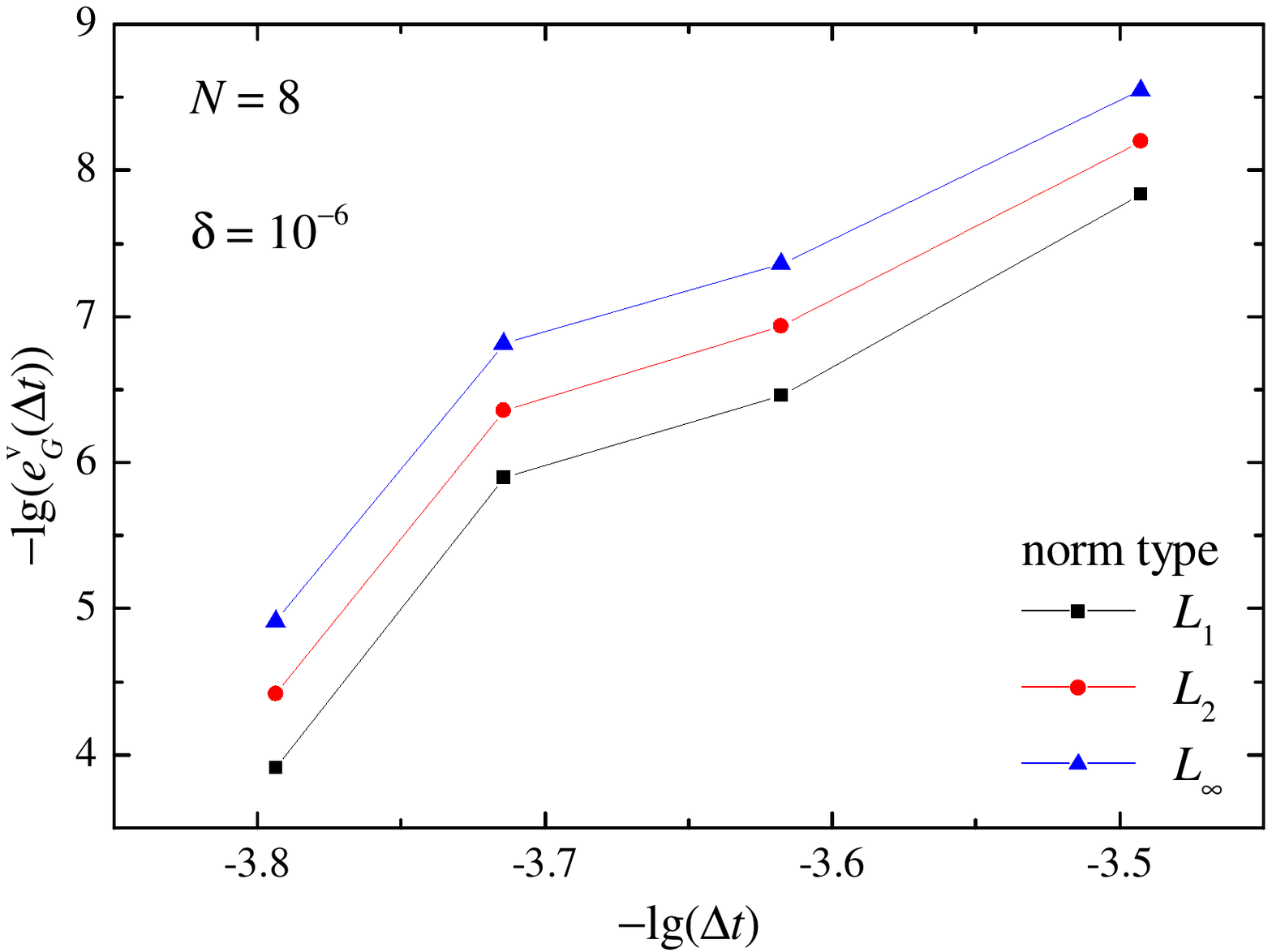}
\vspace{-8mm}\caption{\label{fig:shampine_flame_delta_10m6_errors:e2}}
\end{subfigure}
\begin{subfigure}{0.320\textwidth}
\includegraphics[width=\textwidth]{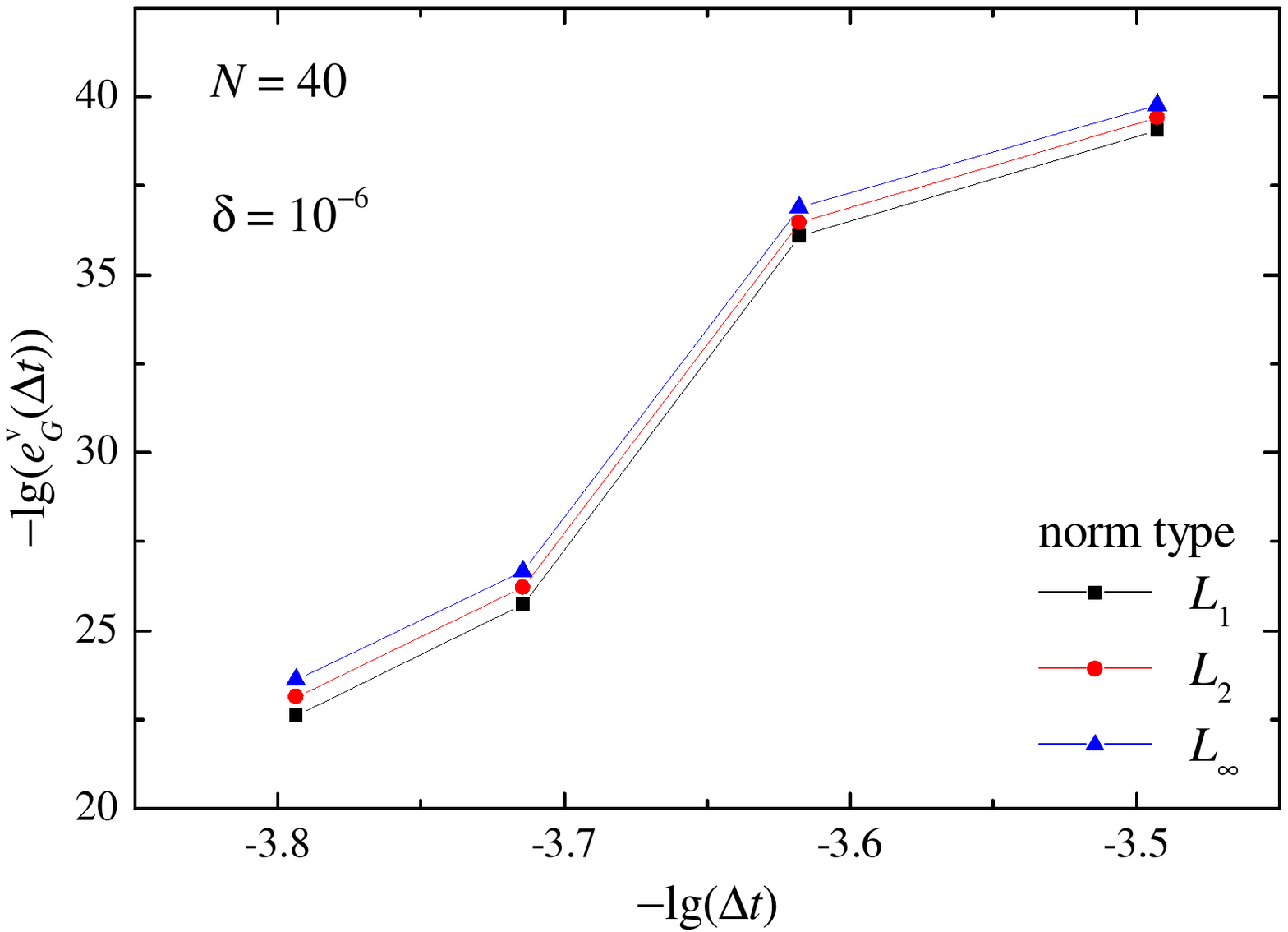}
\vspace{-8mm}\caption{\label{fig:shampine_flame_delta_10m6_errors:e3}}
\end{subfigure}\\
\caption{%
Log-log plot of the dependence of the global errors for the local solution $e_{L}^{u}$ (\subref{fig:shampine_flame_delta_10m6_errors:a1}, \subref{fig:shampine_flame_delta_10m6_errors:a2}, \subref{fig:shampine_flame_delta_10m6_errors:a3}), $e_{L}^{v}$ (\subref{fig:shampine_flame_delta_10m6_errors:b1}, \subref{fig:shampine_flame_delta_10m6_errors:b2}, \subref{fig:shampine_flame_delta_10m6_errors:b3}), $e_{L}^{g}$ (\subref{fig:shampine_flame_delta_10m6_errors:c1}, \subref{fig:shampine_flame_delta_10m6_errors:c2}, \subref{fig:shampine_flame_delta_10m6_errors:c3}) and the solution at nodes $e_{G}^{u}$ (\subref{fig:shampine_flame_delta_10m6_errors:d1}, \subref{fig:shampine_flame_delta_10m6_errors:d2}, \subref{fig:shampine_flame_delta_10m6_errors:d3}), $e_{G}^{v}$ (\subref{fig:shampine_flame_delta_10m6_errors:e1}, \subref{fig:shampine_flame_delta_10m6_errors:e2}, \subref{fig:shampine_flame_delta_10m6_errors:e3}) on the discretization step $\mathrm{\Delta}t$, obtained in the norms $L_{1}$, $L_{2}$ and $L_{\infty}$, by numerical solution of the stiff DAE system (\ref{eq:shampine_flame}) of index 1 with $\delta = 10^{-6}$, obtained using polynomials with degrees $N = 1$ (\subref{fig:shampine_flame_delta_10m6_errors:a1}, \subref{fig:shampine_flame_delta_10m6_errors:b1}, \subref{fig:shampine_flame_delta_10m6_errors:c1}, \subref{fig:shampine_flame_delta_10m6_errors:d1}, \subref{fig:shampine_flame_delta_10m6_errors:e1}), $N = 8$ (\subref{fig:shampine_flame_delta_10m6_errors:a2}, \subref{fig:shampine_flame_delta_10m6_errors:b2}, \subref{fig:shampine_flame_delta_10m6_errors:c2}, \subref{fig:shampine_flame_delta_10m6_errors:d2}, \subref{fig:shampine_flame_delta_10m6_errors:e2}) and $N = 40$ (\subref{fig:shampine_flame_delta_10m6_errors:a3}, \subref{fig:shampine_flame_delta_10m6_errors:b3}, \subref{fig:shampine_flame_delta_10m6_errors:c3}, \subref{fig:shampine_flame_delta_10m6_errors:d3}, \subref{fig:shampine_flame_delta_10m6_errors:e3}).
}
\label{fig:shampine_flame_delta_10m6_errors}
\end{figure}

\begin{table*}[h!]
\centering
\caption{%
Convergence orders $p_{L_{1}}$, $p_{L_{2}}$, $p_{L_{\infty}}$, calculated in norms $L_{1}$, $L_{2}$, $L_{\infty}$ of the ADER-DG method for the stiff DAE system (\ref{eq:shampine_flame}) of index 1 with $\delta = 10^{-6}$; $N$ is the degree of the basis polynomials $\varphi_{p}$. Orders $p^{n, u}$ are calculated for solution at the nodes $\mathbf{u}_{n}$; orders $p^{n, v}$ --- for solution at the nodes $\mathbf{v}_{n}$; orders $p^{l, u}$ --- for local solution $\mathbf{u}_{L}$; orders $p^{l, v}$ --- for local solution $\mathbf{v}_{L}$. The theoretical values of convergence order $p_{\rm th.}^{n} = 2N+1$ and $p_{\rm th.}^{l} = N+1$ are applicable for the ADER-DG method for ODE problems, and are presented for comparison.
}
\label{tab:conv_ords_shampine_flame_delta_10m6}
\begingroup
\setlength{\tabcolsep}{4pt}
\begin{adjustbox}{angle=90}
\begin{tabular}{@{}|l|lll|lll|c|lll|lll|c|@{}}
\toprule
$N$\hspace{-2mm} & $p_{L_{1}}^{n, u}$ & $p_{L_{2}}^{n, u}$ & $p_{L_{\infty}}^{n, u}$ & $p_{L_{1}}^{n, v}$ & $p_{L_{2}}^{n, v}$ & $p_{L_{\infty}}^{n, v}$ & $p_{\rm th.}^{n}$ & $p_{L_{1}}^{l, u}$ & $p_{L_{2}}^{l, u}$ & $p_{L_{\infty}}^{l, u}$ & $p_{L_{1}}^{l, v}$ & $p_{L_{2}}^{l, v}$ & $p_{L_{\infty}}^{l, v}$ & $p_{\rm th.}^{l}$ \\
\midrule
$1$	&	$2.82$	&	$1.48$	&	$0.02$	&	$2.82$	&	$1.43$	&	$0.00$	&	$3$	&	$2.83$	&	$1.48$	&	$0.14$	&	$2.83$	&	$1.43$	&	$-0.04$	&	$2$\\
$2$	&	$6.77$	&	$6.66$	&	$6.26$	&	$5.29$	&	$4.71$	&	$4.14$	&	$5$	&	$3.88$	&	$3.49$	&	$2.31$	&	$3.68$	&	$3.71$	&	$3.15$	&	$3$\\
$3$	&	$5.51$	&	$4.94$	&	$4.43$	&	$5.68$	&	$5.35$	&	$4.88$	&	$7$	&	$3.82$	&	$3.84$	&	$3.71$	&	$3.19$	&	$2.86$	&	$1.85$	&	$4$\\
$4$	&	$6.38$	&	$5.90$	&	$5.40$	&	$6.50$	&	$6.14$	&	$5.65$	&	$9$	&	$4.10$	&	$3.81$	&	$3.09$	&	$4.25$	&	$3.82$	&	$3.08$	&	$5$\\
$5$	&	$7.84$	&	$7.51$	&	$7.06$	&	$8.13$	&	$7.79$	&	$7.31$	&	$11$	&	$4.76$	&	$4.28$	&	$3.63$	&	$5.08$	&	$4.83$	&	$4.30$	&	$6$\\
$6$	&	$12.42$	&	$12.25$	&	$11.85$	&	$12.52$	&	$12.29$	&	$11.83$	&	$13$	&	$6.42$	&	$6.42$	&	$6.62$	&	$4.82$	&	$4.67$	&	$4.40$	&	$7$\\
$7$	&	$12.33$	&	$12.12$	&	$11.67$	&	$12.59$	&	$12.38$	&	$11.92$	&	$15$	&	$5.98$	&	$5.52$	&	$4.48$	&	$7.58$	&	$7.21$	&	$6.14$	&	$8$\\
$8$	&	$11.80$	&	$11.42$	&	$10.94$	&	$12.08$	&	$11.68$	&	$11.19$	&	$17$	&	$7.92$	&	$7.77$	&	$7.38$	&	$6.51$	&	$6.16$	&	$5.73$	&	$9$\\
$9$	&	$13.07$	&	$12.59$	&	$12.09$	&	$13.32$	&	$12.84$	&	$12.34$	&	$19$	&	$7.77$	&	$7.26$	&	$5.84$	&	$8.15$	&	$7.67$	&	$6.74$	&	$10$\\
$10$	&	$14.45$	&	$14.07$	&	$13.64$	&	$14.79$	&	$14.42$	&	$13.99$	&	$21$	&	$8.85$	&	$8.42$	&	$7.71$	&	$8.79$	&	$8.44$	&	$8.02$	&	$11$\\
$11$	&	$14.07$	&	$13.69$	&	$13.24$	&	$14.38$	&	$13.97$	&	$13.49$	&	$23$	&	$10.13$	&	$9.98$	&	$9.83$	&	$8.64$	&	$8.58$	&	$8.30$	&	$12$\\
$12$	&	$16.10$	&	$15.72$	&	$15.33$	&	$16.43$	&	$16.04$	&	$15.57$	&	$25$	&	$9.91$	&	$9.34$	&	$8.73$	&	$11.57$	&	$10.62$	&	$8.62$	&	$13$\\
$13$	&	$16.27$	&	$15.80$	&	$15.30$	&	$16.72$	&	$16.26$	&	$15.76$	&	$27$	&	$11.99$	&	$11.69$	&	$11.36$	&	$10.40$	&	$10.36$	&	$10.25$	&	$14$\\
$14$	&	$22.45$	&	$22.08$	&	$21.67$	&	$22.78$	&	$22.39$	&	$21.92$	&	$29$	&	$11.68$	&	$11.10$	&	$10.29$	&	$11.58$	&	$10.73$	&	$9.32$	&	$15$\\
$15$	&	$20.29$	&	$19.84$	&	$19.34$	&	$20.71$	&	$20.29$	&	$19.81$	&	$31$	&	$12.44$	&	$11.90$	&	$11.01$	&	$13.01$	&	$13.00$	&	$14.21$	&	$16$\\
$16$	&	$22.76$	&	$22.33$	&	$21.84$	&	$23.16$	&	$22.77$	&	$22.30$	&	$33$	&	$13.90$	&	$13.64$	&	$12.95$	&	$12.86$	&	$12.28$	&	$11.69$	&	$17$\\
$17$	&	$23.01$	&	$22.63$	&	$22.20$	&	$23.36$	&	$22.98$	&	$22.54$	&	$35$	&	$13.85$	&	$13.27$	&	$12.56$	&	$14.52$	&	$13.44$	&	$12.17$	&	$18$\\
$18$	&	$20.07$	&	$19.68$	&	$19.22$	&	$20.37$	&	$19.95$	&	$19.46$	&	$37$	&	$15.34$	&	$14.73$	&	$13.68$	&	$14.41$	&	$14.57$	&	$14.65$	&	$19$\\
$19$	&	$31.08$	&	$30.66$	&	$30.17$	&	$31.47$	&	$31.09$	&	$30.63$	&	$39$	&	$15.66$	&	$15.19$	&	$14.72$	&	$15.26$	&	$14.04$	&	$11.66$	&	$20$\\
$20$	&	$26.53$	&	$26.08$	&	$25.59$	&	$26.79$	&	$26.33$	&	$25.83$	&	$41$	&	$16.34$	&	$15.73$	&	$15.12$	&	$16.94$	&	$17.03$	&	$18.01$	&	$21$\\
\midrule
$25$	&	$33.91$	&	$33.45$	&	$32.96$	&	$34.34$	&	$33.91$	&	$33.42$	&	$51$	&	$19.71$	&	$18.74$	&	$17.23$	&	$21.42$	&	$21.42$	&	$22.27$	&	$26$\\
$30$	&	$43.53$	&	$43.05$	&	$42.55$	&	$43.78$	&	$43.29$	&	$42.79$	&	$61$	&	$23.55$	&	$22.69$	&	$21.62$	&	$24.73$	&	$24.37$	&	$24.31$	&	$31$\\
$35$	&	$52.14$	&	$51.69$	&	$51.19$	&	$52.40$	&	$51.94$	&	$51.44$	&	$71$	&	$27.54$	&	$26.84$	&	$26.26$	&	$28.28$	&	$27.41$	&	$25.86$	&	$36$\\
$40$	&	$58.64$	&	$58.15$	&	$57.65$	&	$58.88$	&	$58.40$	&	$57.90$	&	$81$	&	$31.29$	&	$30.76$	&	$30.25$	&	$31.58$	&	$31.05$	&	$30.57$	&	$41$\\
\bottomrule
\end{tabular}
\end{adjustbox}
\endgroup
\end{table*}

Fig.~\ref{fig:shampine_flame_delta_10m4_sol_v_epss} (\subref{fig:shampine_flame_delta_10m4_sol_v_epss:c1}, \subref{fig:shampine_flame_delta_10m4_sol_v_epss:c2}, \subref{fig:shampine_flame_delta_10m4_sol_v_epss:c3}, \subref{fig:shampine_flame_delta_10m4_sol_v_epss:d1}, \subref{fig:shampine_flame_delta_10m4_sol_v_epss:d2}, \subref{fig:shampine_flame_delta_10m4_sol_v_epss:d3}) shows the dynamic dependencies of the errors $\varepsilon_{u}(t)$, $\varepsilon_{v}(t)$ of the numerical solution at the nodes $(\mathbf{u}_{n},\, \mathbf{v}_{n})$ and the local solution $(\mathbf{u}_{L}(t),\, \mathbf{v}_{L}(t))$ for differential variables $\mathbf{u}$ and algebraic variables $\mathbf{v}$, as well as the dynamic dependencies of the errors $\varepsilon_{g}(t)$ (Fig.~\ref{fig:shampine_flame_delta_10m4_sol_v_epss} (\subref{fig:shampine_flame_delta_10m4_sol_qug:e1}, \subref{fig:shampine_flame_delta_10m4_sol_qug:e2}, \subref{fig:shampine_flame_delta_10m4_sol_qug:e3})) of satisfaction of algebraic equations $\mathbf{g} = \mathbf{0}$ (in this problem the decrease in the index of the DAE system was not carried out, therefore numerically this error $\varepsilon_{g}(t)$ is equal to the error $|g_{1}(t)|$ of satisfaction of the algebraic equation $g_{1} = 0$). The presented results clearly demonstrate that the error $\varepsilon_{u}(t)$ of the solution at the nodes $\mathbf{u}_{n}$ is less than the error $\varepsilon_{u}(t)$ of the local solution $\mathbf{u}_{L}(t)$ for all values of the argument $t$. This is observed both for differential variables $\mathbf{u}$ (which is shown in Fig.~\ref{fig:shampine_flame_delta_10m4_sol_v_epss} (\subref{fig:shampine_flame_delta_10m4_sol_v_epss:c1}, \subref{fig:shampine_flame_delta_10m4_sol_v_epss:c2}, \subref{fig:shampine_flame_delta_10m4_sol_v_epss:c3})) and algebraic variables $\mathbf{v}$ (which is shown in Fig.~\ref{fig:shampine_flame_delta_10m4_sol_v_epss} (\subref{fig:shampine_flame_delta_10m4_sol_v_epss:d1}, \subref{fig:shampine_flame_delta_10m4_sol_v_epss:d2}, \subref{fig:shampine_flame_delta_10m4_sol_v_epss:d3})). In the case of a local solution $(\mathbf{u}_{L}(t),\, \mathbf{v}_{L}(t))$, a minimum of errors $\varepsilon_{u}(t)$, $\varepsilon_{v}(t)$ is observed to the left of point $t \approx 1/\delta$, which is similar to the case of an error $|g_{1}|$ in satisfying an algebraic equation $g_{1} = 0$ (in Fig.~\ref{fig:shampine_flame_delta_10m4_sol_qug} (\subref{fig:shampine_flame_delta_10m4_sol_qug:e1}, \subref{fig:shampine_flame_delta_10m4_sol_qug:e2}, \subref{fig:shampine_flame_delta_10m4_sol_qug:e3})). However, in the case of a solution at nodes $(\mathbf{u}_{n},\, \mathbf{v}_{n})$, this feature of the numerical solution errors $\varepsilon_{u}(t)$, $\varepsilon_{v}(t)$ does not occur --- the errors to the left of point $t \approx 1/\delta$ smoothly increases from a zero value at the initial point $t_{0} = 0$, reaching a maximum value in the vicinity of point $t \approx 1/\delta$, and then quickly decreases to almost zero, within the accuracy of representation of real numbers with a floating point. This effect of the dynamic dependence of the solution errors $\varepsilon_{u}(t)$, $\varepsilon_{v}(t)$ at nodes $(\mathbf{u}_{n},\, \mathbf{v}_{n})$ occurs for both differential variables $\mathbf{u}$ and algebraic variables $\mathbf{v}$. Fig.~\ref{fig:shampine_flame_delta_10m4_errors} shows the log-log plots of the dependence of global errors $e$ of the numerical solution at nodes $(\mathbf{u}_{n},\, \mathbf{v}_{n})$, the local solution $(\mathbf{u}_{L}(t),\, \mathbf{v}_{L}(t))$ and of the satisfaction of algebraic equations $\mathbf{g} = \mathbf{0}$ on the discretization step ${\Delta t}$; in this case, the grid was chosen to be non-uniform, therefore, as a quantitative value of the discretization step ${\Delta t}$, the ``average'' value was chosen --- the length of the domain of definition $\Omega$, divided by the total number of grid nodes $L$: ${\Delta t} = 2/(\delta \cdot L)$, which, up to a constant factor, which turned into an additive constant after taking logarithm, coincides with the discretization step in any of the subdomains $(\Omega_{1},\, \Omega_{2},\, \Omega_{3})$, and this choice does not affect the value of the slope of the approximating line. The global errors $e$ were calculated separately for the norms $L_{1}$, $L_{2}$, $L_{\infty}$ (\ref{eq:norms_def}). The presented dependencies demonstrate the power law $e(\Delta t) \sim (\Delta t)^{p}$, from which the empirical convergence orders $p$ were calculated. However, in comparison with the previously studied test examples, in this case the approximations are not so good.

The calculated empirical convergence orders $p^{n, u}$, $p^{n, v}$ for the numerical solution at the nodes $(\mathbf{u}_{n},\, \mathbf{v}_{n})$ and $p^{l, u}$, $p^{l, v}$ for the local solution $(\mathbf{u}_{L}(t),\, \mathbf{v}_{L}(t))$ are presented in Table~\ref{tab:conv_ords_shampine_flame_delta_10m4}. The convergence orders were calculated separately for the norms $L_{1}$, $L_{2}$, $L_{\infty}$ (\ref{eq:norms_def}). The theoretical values of convergence order $p_{\rm th.}^{n} = 2N+1$ and $p_{\rm th.}^{l} = N+1$ (\ref{eq:expect_orders}) are applicable for the ADER-DG method for ODE problems, and are presented in Table~\ref{tab:conv_ords_shampine_flame_delta_10m4} for comparison. The obtained results show that the empirical convergence orders $p$ have values smaller than the expected convergence orders (\ref{eq:expect_orders}) for the ADER-DG method with a local DG predictor for the initial value problem for the ODE system. However, the obtained result has a \textit{fundamental feature} compared to the expected values. In the work~\cite{ader_dg_ivp_ode}, the initial value problem for the ODE system equivalent to the DAE system (\ref{eq:shampine_flame}) was solved using the ADER-DG method with a local DG predictor, and the empirical convergence orders $p$ for the case $\delta = 10^{-4}$ did not exceed $8.15$, $8.02$, $7.65$ in the norms of $L_{1}$, $L_{2}$, $L_{\infty}$, respectively, for the solution at the nodes (starting with the degree of polynomials $N = 3$, the empirical convergence orders $p$ for the solution at the nodes stopped showing an increase with the growth of the degree of polynomials $N$) and $6.50$, $6.45$, $5.93$ in the norms of $L_{1}$, $L_{2}$, $L_{\infty}$, respectively, for the local solution at the nodes (in the case of a local solution, the growth of the empirical convergence orders $p$ stopped with the degree of polynomials $N = 6$). In the present work, the empirical convergence orders $p$ are also smaller than the expected values (\ref{eq:expect_orders}) for non-stiff problems, however, the empirical convergence orders $p$ demonstrate a relatively monotonic (except for several local minima) growth with increasing polynomial degree $N$. In this case, the empirical convergence orders $p$ for the local solution $(\mathbf{u}_{L}(t),\, \mathbf{v}_{L}(t))$ are $1$--$3$ smaller than the expected values (\ref{eq:expect_orders}) in the region of low polynomial degrees $N$ and $3$--$5$ smaller than the expected values (\ref{eq:expect_orders}) in the region of high polynomial degrees $N$. The empirical convergence orders $p$ for the solution at the nodes $(\mathbf{u}_{n},\, \mathbf{v}_{n})$ are comparable or slightly higher than the expected values for polynomial degrees $N = 1$--$3$, $1$--$3$ smaller than the expected values in the region of low polynomial degrees $N$ and $5$--$10$ smaller than the expected values (\ref{eq:expect_orders}) in the region of high polynomial degrees $N$. This result is of \textit{fundamental importance} for the choice of basis functions when using the ADER-DG method with a local DG predictor to solve stiff problems --- Lagrange interpolation polynomials with nodal points at the roots of shifted Legendre polynomials as basis functions for representing a local discrete time solution (\ref{eq:qr_def_exp}) were used in the work~\cite{ader_dg_ivp_ode}, and set of Lagrange interpolation polynomials with nodal points at the roots of the right Radau polynomials were used to represent a local discrete time solution (\ref{eq:qr_def_exp}) in the present work. It turns out that the use of Lagrange interpolation polynomials with nodal points at the roots of the right Radau polynomials is much better suited for solving stiff problems.

The obtained results of the numerical solution of DAE system (\ref{eq:shampine_flame}) of index 1 with $\delta = 10^{-5}$ are presented in Figs.~\ref{fig:shampine_flame_delta_10m5_sol_qug}, \ref{fig:shampine_flame_delta_10m5_sol_v_epss}, \ref{fig:shampine_flame_delta_10m5_errors} and in Table~\ref{tab:conv_ords_shampine_flame_delta_10m5}. The arrangement of the obtained results in the Figures is similar to the previous presented case $\delta = 10^{-4}$: Fig.~\ref{fig:shampine_flame_delta_10m5_sol_qug} shows a comparison of the solution at nodes $\mathbf{u}_{n}$, the local solution $\mathbf{u}_{L}(t)$ and the exact solution $\mathbf{u}^{\rm ex}(t)$ for component $u_{1}$, quantitative satisfiability of the conditions $g_{1} = 0$; Fig.~\ref{fig:shampine_flame_delta_10m5_sol_v_epss} shows a comparison of the solution at nodes $\mathbf{v}_{n}$, the local solution $\mathbf{v}_{L}(t)$ and the exact solution $\mathbf{v}^{\rm ex}(t)$ for component $v_{1}$, the errors $\varepsilon_{u}(t)$, $\varepsilon_{v}(t)$, $\varepsilon_{g}(t)$. In the vicinity of point $t = 1.0/\delta$, a sharp change in the solution $(\mathbf{u},\, \mathbf{v})$ occurs, which leads to high stiffness of the DAE system (\ref{eq:shampine_flame}) with $\delta = 10^{-4}$, therefore the numerical solution at nodes $(\mathbf{u}_{n},\, \mathbf{v}_{n})$, the local solution $(\mathbf{u}_{L}(t),\, \mathbf{v}_{L}(t))$ and its comparison with the exact analytical solution $(\mathbf{u}^{\rm ex}(t),\, \mathbf{v}^{\rm ex}(t))$ (\ref{eq:shampine_flame:exact_solution}) in a small vicinity of point $t = 1.0/\delta$, specifically in the range $t \in [0.995/\delta,\, 1.005/\delta]$, was presented separately in Fig.~\ref{fig:shampine_flame_delta_10m5_sol_qug} (\subref{fig:shampine_flame_delta_10m5_sol_qug:d1}, \subref{fig:shampine_flame_delta_10m5_sol_qug:d2}, \subref{fig:shampine_flame_delta_10m5_sol_qug:d3}, \subref{fig:shampine_flame_delta_10m5_sol_qug:b1}, \subref{fig:shampine_flame_delta_10m5_sol_qug:b2}, \subref{fig:shampine_flame_delta_10m5_sol_qug:b3}) for differential variables $\mathbf{u}$ and in Fig.~\ref{fig:shampine_flame_delta_10m5_sol_v_epss} (\subref{fig:shampine_flame_delta_10m5_sol_v_epss:b1}, \subref{fig:shampine_flame_delta_10m5_sol_v_epss:b2}, \subref{fig:shampine_flame_delta_10m5_sol_v_epss:b3}) for algebraic variables $\mathbf{v}$. The comparison of the numerical solution and the exact analytical solution in this zoomed range $t \in [0.995/\delta,\, 1.005/\delta]$ allowed one to identify artifacts of the numerical solution obtained by the ADER-DG method with a local DG predictor more accurately and correctly, in the case of this stiff DAE system (\ref{eq:shampine_flame}) of index 1 with $\delta = 10^{-5}$. Fig.~\ref{fig:pend_ind2_errors} shows the log-log plot of the dependence of the global errors for the local solution $e_{L}^{u}$, $e_{L}^{v}$, $e_{L}^{g}$ and the solution at nodes $e_{G}^{u}$, $e_{G}^{v}$ on the discretization step $\mathrm{\Delta}t$, obtained in the norms $L_{1}$, $L_{2}$ and $L_{\infty}$, by numerical solution of the stiff DAE system (\ref{eq:shampine_flame}) of index 1 with $\delta = 10^{-5}$, which was used to calculate the empirical convergence orders $p$.

The dynamic dependencies of the numerical solution presented in Fig.~\ref{fig:shampine_flame_delta_10m5_sol_qug} and Fig.~\ref{fig:shampine_flame_delta_10m5_sol_v_epss} demonstrate a very high accuracy of the ADER-DG numerical method with a local DG predictor for solving the stiff problem (\ref{eq:shampine_flame}) with $\delta = 10^{-4}$. A point-by-point comparison of the dynamic dependencies of the solution at the nodes $(\mathbf{u}_{n},\, \mathbf{v}_{n})$ and the local solution $(\mathbf{u}_{L}(t),\, \mathbf{v}_{L}(t))$ with the exact analytical solution $(\mathbf{u}^{\rm ex},\, \mathbf{v}^{\rm ex})$ for differential variables $\mathbf{u}$ and algebraic variables $\mathbf{v}$ shows that in the case of using polynomials of degree $N = 1$, an expected artifact of the numerical solution is observed --- the range of a sharp change in the solution occurs in the region of larger values of the argument $t$ than in the exact analytical solution, and the shift occurs approximately by $\delta\cdot t \approx 0.0040$--$0.0045$ to the right. This behavior of the numerical solution is similar to case $\delta = 10^{-4}$ considered above. This effect, as in case $\delta = 10^{-4}$, has already been observed in the work~\cite{ader_dg_ivp_ode} when using the ADER-DG numerical method with a local DG predictor when solving the initial value problem for an ODE system that is completely equivalent to a DAE system (\ref{eq:shampine_flame}) of index 1 with $\delta = 10^{-5}$. Note that a similar effect also occurs when using the Runge-Kutta numerical methods of the Radau family, as well as the \texttt{BDF} and \texttt{LSODA} methods available in the \texttt{scipy} module of the \texttt{python} programming language (which was also noted in the work~\cite{ader_dg_ivp_ode}), while in the case of the \texttt{BDF} and \texttt{LSODA} methods this effect is much stronger. The point-wise comparison of the dynamic dependencies of the solution at the nodes $(\mathbf{u}_{n},\, \mathbf{v}_{n})$ and the local solution $(\mathbf{u}_{L}(t),\, \mathbf{v}_{L}(t))$ with the exact analytical solution $(\mathbf{u}^{\rm ex},\, \mathbf{v}^{\rm ex})$, in the case of using polynomials of degrees $N = 8$ and $40$, shows that this artifact of the numerical solution does not occur --- the numerical solution very accurately coincides point-wise with the exact analytical solution.

Despite the significantly stronger stiffness of the problem with $\delta = 10^{-5}$, compared to $\delta = 10^{-4}$, the dynamic dependencies of the error $|g_{1}|$ of satisfaction of the algebraic equation $g_{1} = 0$ for the local solution $(\mathbf{u}_{L}(t),\, \mathbf{v}_{L}(t))$ are at approximately the same level of values. The algebraic equation $g_{1} = 0$ is satisfied exactly (within the accuracy of representation of real numbers with a floating point) for the solution at the nodes $(\mathbf{u}_{n},\, \mathbf{v}_{n})$, therefore the error $|g_{1}|$ for the solution at the nodes is not shown in Fig.~\ref{fig:shampine_flame_delta_10m5_sol_qug}. However, an important difference is the situation with the minimum of the error, which in case $\delta = 10^{-4}$ was noticeably to the left of point $t = 1.0/\delta$ --- in case $\delta = 10^{-5}$, the area with the minimum of the error $|g_{1}|$ becomes very small, however, the minimum of the error $|g_{1}|$ is achieved exactly in it. Similar behavior is demonstrated by the dynamic dependencies of the errors $\varepsilon_{u}(t)$, $\varepsilon_{v}(t)$, $\varepsilon_{g}(t)$, obtained using polynomials with degrees $N = 1$, $8$ and $40$, taking into account the narrowing of the area of the minimum of the error to the left of point $t = 1.0/\delta$: in the case of a local solution $(\mathbf{u}_{L}(t),\, \mathbf{v}_{L}(t))$, a minimum of errors $\varepsilon_{u}(t)$, $\varepsilon_{v}(t)$ is observed to the left of point $t \approx 1/\delta$, which is similar to the case of an error $|g_{1}|$ in satisfying an algebraic equation $g_{1} = 0$ (in Fig.~\ref{fig:shampine_flame_delta_10m5_sol_qug} (\subref{fig:shampine_flame_delta_10m5_sol_qug:e1}, \subref{fig:shampine_flame_delta_10m5_sol_qug:e2}, \subref{fig:shampine_flame_delta_10m5_sol_qug:e3})); in the case of a solution at nodes $(\mathbf{u}_{n},\, \mathbf{v}_{n})$, this feature of the numerical solution errors $\varepsilon_{u}(t)$, $\varepsilon_{v}(t)$ does not occur --- the errors to the left of point $t \approx 1/\delta$ smoothly increases from a zero value at the initial point $t_{0} = 0$, reaching a maximum value in the vicinity of point $t \approx 1/\delta$, and then quickly decreases to almost zero, within the accuracy of representation of real numbers with a floating point. Fig.~\ref{fig:shampine_flame_delta_10m5_errors} shows the log-log plots of the dependence of global errors $e$ of the numerical solution at nodes $(\mathbf{u}_{n},\, \mathbf{v}_{n})$, the local solution $(\mathbf{u}_{L}(t),\, \mathbf{v}_{L}(t))$ and of the satisfaction of algebraic equations $\mathbf{g} = \mathbf{0}$ on the discretization step ${\Delta t}$. The global errors $e$ were calculated separately for the norms $L_{1}$, $L_{2}$, $L_{\infty}$ (\ref{eq:norms_def}). The presented dependencies demonstrate the power law $e(\Delta t) \sim (\Delta t)^{p}$, from which the empirical convergence orders $p$ were calculated. However, compared to case $\delta = 10^{-5}$, in this case $\delta = 10^{-5}$ the dependencies are much better approximated by linear functions.

The calculated empirical convergence orders $p^{n, u}$, $p^{n, v}$ for the numerical solution at the nodes $(\mathbf{u}_{n},\, \mathbf{v}_{n})$ and $p^{l, u}$, $p^{l, v}$ for the local solution $(\mathbf{u}_{L}(t),\, \mathbf{v}_{L}(t))$ are presented in Table~\ref{tab:conv_ords_shampine_flame_delta_10m5}. The convergence orders were calculated separately for the norms $L_{1}$, $L_{2}$, $L_{\infty}$ (\ref{eq:norms_def}). The theoretical values of convergence order $p_{\rm th.}^{n} = 2N+1$ and $p_{\rm th.}^{l} = N+1$ (\ref{eq:expect_orders}) are applicable for the ADER-DG method for ODE problems, and are presented in Table~\ref{tab:conv_ords_shampine_flame_delta_10m5} for comparison. In this case, qualitatively similar to case $\delta = 10^{-5}$ behavior of the dependence of the empirical convergence of orders $p$ on the degree of polynomials $N$ is observed. The obtained results show that the empirical convergence orders $p$ have values smaller than the expected convergence orders (\ref{eq:expect_orders}) for the ADER-DG method with a local DG predictor for the initial value problem for the ODE system. In the work~\cite{ader_dg_ivp_ode}, the initial value problem for the ODE system equivalent to the DAE system (\ref{eq:shampine_flame}) was solved using the ADER-DG method with a local DG predictor, and the empirical convergence orders $p$ for the case $\delta = 10^{-5}$ did not exceed $4.21$, $4.08$, $3.78$ in the norms of $L_{1}$, $L_{2}$, $L_{\infty}$, respectively, for the solution at the nodes (starting with the degree of polynomials $N = 5$, the empirical convergence orders $p$ for the solution at the nodes stopped showing an increase with the growth of the degree of polynomials $N$) and $4.22$, $4.20$, $3.40$ in the norms of $L_{1}$, $L_{2}$, $L_{\infty}$, respectively, for the local solution at the nodes (in the case of a local solution, the growth of the empirical convergence orders $p$ stopped with the degree of polynomials $N = 7$) --- therefore, with an increase in the stiffness of the problem, a significant decrease in the empirical convergence orders was observed. In the present work, the empirical convergence orders $p$ are also smaller than the expected values (\ref{eq:expect_orders}) for non-stiff problems, however, the empirical convergence orders $p$ demonstrate a relatively monotonic (unlike case $\delta = 10^{-4}$ --- strictly monotonic) growth with increasing polynomial degree $N$. In this case, the empirical convergence orders $p$ for the local solution $(\mathbf{u}_{L}(t),\, \mathbf{v}_{L}(t))$ are comparable or higher (by $1$--$4$) than the expected values for polynomial degrees $N = 5$--$11$ (for different orders and norms), $1$--$8$ smaller than the expected values in the region of low polynomial degrees $N$ and $10$--$20$ smaller than the expected values (\ref{eq:expect_orders}) in the region of high polynomial degrees $N$. The empirical convergence orders $p$ for the solution at the nodes $(\mathbf{u}_{n}(t),\, \mathbf{v}_{n}(t))$ are comparable or slightly higher than the expected values for polynomial degrees $N = 1$--$3$, $1$--$5$ smaller than the expected values in the region of low polynomial degrees $N$ and $10$--$28$ smaller than the expected values (\ref{eq:expect_orders}) in the region of high polynomial degrees $N$. This result is also related to the choice of basis functions when using the ADER-DG method to solve stiff problems --- Lagrange interpolation polynomials with nodal points at the roots of shifted Legendre polynomials as basis functions for representing a local discrete time solution (\ref{eq:qr_def_exp}) were used in the work~\cite{ader_dg_ivp_ode}, and set of Lagrange interpolation polynomials with nodal points at the roots of the right Radau polynomials were used to represent a local discrete time solution (\ref{eq:qr_def_exp}) in the present work.

The obtained results of the numerical solution of DAE system (\ref{eq:shampine_flame}) of index 1 with $\delta = 10^{-6}$ are presented in Figs.~\ref{fig:shampine_flame_delta_10m6_sol_qug}, \ref{fig:shampine_flame_delta_10m6_sol_v_epss}, \ref{fig:shampine_flame_delta_10m6_errors} and in Table~\ref{tab:conv_ords_shampine_flame_delta_10m6}. The arrangement of the obtained results in this Figures is similar to the previous presented cases $\delta = 10^{-4}$ and $\delta = 10^{-5}$. This case $\delta = 10^{-6}$ was not considered in the work~\cite{ader_dg_ivp_ode} on the numerical solution of the initial value problem for an ODE system equivalent to a DAE system (\ref{eq:shampine_flame}), and in principle is not presented in existing works, as far as the author knows. A point-by-point comparison of the dynamic dependencies of the solution at the nodes $(\mathbf{u}_{n},\, \mathbf{v}_{n})$ and the local solution $(\mathbf{u}_{L}(t),\, \mathbf{v}_{L}(t))$ with the exact analytical solution $(\mathbf{u}^{\rm ex},\, \mathbf{v}^{\rm ex})$ for differential variables $\mathbf{u}$ and algebraic variables $\mathbf{v}$ shows that in the case of using polynomials of degree $N = 1$, an expected artifact of the numerical solution is observed --- the range of a sharp change in the solution occurs in the region of larger values of the argument $t$ than in the exact analytical solution, and the shift occurs approximately by $\delta\cdot t \approx 0.00060$--$0.00075$ to the right. This behavior of the numerical solution is similar to cases $\delta = 10^{-4}$ and $\delta = 10^{-5}$ considered above. It is interesting to note that in units of $\delta\cdot t$ this shift is smaller than it was not only in case $\delta = 10^{-5}$, but also in case $\delta = 10^{-4}$, which are characterized by a significantly weaker stiffness, which can be seen from the results presented above in Figs.~\ref{fig:shampine_flame_delta_10m4_sol_qug}, \ref{fig:shampine_flame_delta_10m4_sol_v_epss}, \ref{fig:shampine_flame_delta_10m5_sol_qug}, \ref{fig:shampine_flame_delta_10m5_sol_v_epss}. Note that a similar effect also occurs when using the Runge-Kutta numerical methods of the Radau family, as well as the \texttt{BDF} and \texttt{LSODA} methods available in the \texttt{scipy} module of the \texttt{python} programming language, while in the case of the \texttt{BDF} and \texttt{LSODA} methods this effect is much stronger. The point-wise comparison of the dynamic dependencies of the solution at the nodes $(\mathbf{u}_{n},\, \mathbf{v}_{n})$ and the local solution $(\mathbf{u}_{L}(t),\, \mathbf{v}_{L}(t))$ with the exact analytical solution $(\mathbf{u}^{\rm ex},\, \mathbf{v}^{\rm ex})$, in the case of using polynomials of degrees $N = 8$ and $40$, shows that this artifact of the numerical solution does not occur --- the numerical solution very accurately coincides point-wise with the exact analytical solution. An important difference from the previous cases $\delta = 10^{-4}$ and $\delta = 10^{-5}$ is that in case $\delta = 10^{-6}$ the local errors $\varepsilon_{u}(t)$, $\varepsilon_{v}(t)$, $\varepsilon_{g}(t)$ minimum to the left of point $t \approx 1/\delta$ is not explicitly resolved on the selected grid sizes. Fig.~\ref{fig:shampine_flame_delta_10m5_errors} shows the log-log plots of the dependence of global errors $e$ of the numerical solution at nodes $(\mathbf{u}_{n},\, \mathbf{v}_{n})$, the local solution $(\mathbf{u}_{L}(t),\, \mathbf{v}_{L}(t))$ and of the satisfaction of algebraic equations $\mathbf{g} = \mathbf{0}$ on the discretization step ${\Delta t}$. The global errors $e$ were calculated separately for the norms $L_{1}$, $L_{2}$, $L_{\infty}$ (\ref{eq:norms_def}). The presented dependencies demonstrate the power law $e(\Delta t) \sim (\Delta t)^{p}$, from which the empirical convergence orders $p$ were calculated. However, compared to cases $\delta = 10^{-4}$ and $\delta = 10^{-5}$, in this case $\delta = 10^{-5}$ the approximations are not so good.

The calculated empirical convergence orders $p^{n, u}$, $p^{n, v}$ for the numerical solution at the nodes and $p^{l, u}$, $p^{l, v}$ for the local solution are presented in Table~\ref{tab:conv_ords_shampine_flame_delta_10m6}. The convergence orders were calculated separately for the norms $L_{1}$, $L_{2}$, $L_{\infty}$ (\ref{eq:norms_def}). The theoretical values of convergence order $p_{\rm th.}^{n} = 2N+1$ and $p_{\rm th.}^{l} = N+1$ (\ref{eq:expect_orders}) are applicable for the ADER-DG method for ODE problems, and are presented in Table~\ref{tab:conv_ords_shampine_flame_delta_10m5} for comparison. The empirical convergence orders $p$ obtained in this case are less than the expected values (\ref{eq:expect_orders}), however, in the region of large degrees of polynomials $N$, the empirical convergence orders $p$ are even higher than in the case $\delta = 10^{-5}$. This, of course, can be associated simply with a lower quality of approximation of dependencies $e(\Delta t)$ by $e(\Delta t) \sim (\Delta t)^{p}$.

Therefore, it can be concluded that the ADER-DG method with a local DG predictor allows obtaining a highly accurate numerical solution of the DAE system characterized by strong stiffness, however, a significant decrease in the order of convergence $p$ of the numerical solution is observed. An important result was obtained about the choice of basis functions for representing a local discrete time solution (\ref{eq:qr_def_exp}) --- Lagrange interpolation polynomials with nodal points at the roots of shifted Legendre polynomials as basis functions for representing a local discrete time solution (\ref{eq:qr_def_exp}) were used in the work~\cite{ader_dg_ivp_ode}, and set of Lagrange interpolation polynomials with nodal points at the roots of the right Radau polynomials were used to represent a local discrete time solution (\ref{eq:qr_def_exp}) in the present work. It turns out that the use of Lagrange interpolation polynomials with nodal points at the roots of the right Radau polynomials is much better suited for solving stiff problems. This result is of \textit{fundamental importance} for the choice of basis functions when using the ADER-DG method with a local DG predictor to solve stiff problems.

\section*{Conclusion}
\addcontentsline{toc}{section}{Conclusion}
\label{sec:conclusion}

In conclusion, it should be noted that in this work, a numerical ADER-DG method with a local DG predictor for solving a DAE system was presented. The developed method is a development of the numerical ADER-DG method presented in the work~\cite{ader_dg_ivp_ode} for solving the initial value problem for an ODE system, which was based on the formulation of the numerical ADER-DG methods using a local DG predictor presented in the works~\cite{ader_dg_ideal_flows, ader_dg_diss_flows, ader_dg_ale, ader_dg_grmhd, ader_dg_gr_prd, ader_dg_gr_z4_2024, ader_dg_simple_mod, ader_dg_PNPM, fron_phys, ader_dg_axioms, exahype, PNPM_DG_2009}. An important feature of the numerical ADER-DG method with a local DG predictor developed in this paper is the condition of $L_{2}$-orthogonality of the local solution error $\boldsymbol{\sigma}_{\rm v}$ to the basis functions $\varphi_{k}$, which was chosen in form (\ref{eq:lst_dg_error_cond}), and not in the form of a point-wise correspondence of the local solution $(\mathbf{q},\, \mathbf{r})$ to the subsystem of algebraic equations of the DAE system. The presented numerical ADER-DG method with a local DG predictor has an arbitrarily high convergence order.

Polynomials were chosen as basis functions $\varphi_{k}$, which differs from the polynomials chosen in the work~\cite{ader_dg_ivp_ode} and in the works~\cite{ader_dg_ideal_flows, ader_dg_diss_flows, ader_dg_ale, ader_dg_grmhd, ader_dg_gr_prd, ader_dg_gr_z4_2024, ader_dg_simple_mod, ader_dg_PNPM, fron_phys, ader_dg_axioms, exahype, PNPM_DG_2009}. The use of Lagrange polynomials with nodal points at the roots of Legendre polynomials $L_{N}(\tau)$ as basis functions $\varphi_{p}$ showed that such a basis is not very suitable for use in the ADER-DG method with a local DG predictor for solving the DAE system (\ref{eq:dae_chosen_form}), despite the fact that this basis was chosen in works~\cite{ader_stiff_1, PNPM_DG_2010, ader_dg_ivp_ode} for the numerical solution of the initial value problem for the ODE system and in works~\cite{ader_dg_ideal_flows, ader_dg_diss_flows, ader_dg_ale, ader_dg_grmhd, ader_dg_gr_prd, ader_dg_gr_z4_2024, ader_dg_simple_mod, ader_dg_PNPM, fron_phys, ader_dg_axioms, exahype, PNPM_DG_2009} for solving PDE problems. In this case, the local solution demonstrated very strong discontinuities at the grid nodes for DAE systems with a high index (2 and more). The choice of basis functions in the form of Lagrange polynomials with nodal points at the roots of right Radau polynomials $R_{N}(\tau)$ made it possible to eliminate this drawback of the ADER-DG method with a local DG predictor for solving DAE systems (\ref{eq:dae_chosen_form}). It was shown that the use of polynomials allows obtaining results similar to the results of the work~\cite{ader_dg_ivp_ode} in the case of using the DAE solver as the ODE solver. It was shown that the use of polynomials does not change the stability properties of the ADER-DG numerical method compared to the case in the work~\cite{ader_dg_ivp_ode} --- the numerical method remains $A$-stable and $L$-stable (more precisely, $L_{1}$-stable).

The presented ADER-DG method with a local DG predictor allows obtaining a highly accurate numerical solution even on very coarse grids, with a step larger than the main characteristic scale of solution change. In this case, the local discrete time solution $(\mathbf{q},\, \mathbf{r})$ (\ref{eq:qr_def_exp}) can be used as a numerical solution of the DAE system between grid nodes $t \in [t_{n},\, t_{n+1}]$, thereby providing subgrid resolution even in the case of very coarse grids. The local error $\varepsilon$ of the obtained local solution even on a very coarse grid reached the round-off error of representing real numbers by double-precision floating-point numbers. This property of the ADER-DG method presented in this paper allows using a local solution between grid nodes, specified by a high-precision representation in the form of an expansion in basis polynomials, to calculate the solution functionals within more complex problems, without using any methods of reconstructing the solution in the space between grid nodes, such as interpolation.

Using the developed numerical method ADER-DG with a local DG predictor, classical test examples were solved, including a simple DAE system of index 1 in Subsection~\ref{sec:2:ct:ex1}, Hessenberg DAE systems of index 1 and index 2 in Subsections~\ref{sec:2:ct:ex2} and~\ref{sec:2:ct:ex3}, the mathematical pendulum problem in Subsection~\ref{sec:2:ct:ex4}, and the double pendulum problem in Subsection~\ref{sec:2:ct:ex5}. Empirical convergence orders $p$ were calculated for the solution at the nodes $(\mathbf{u}_{n},\, \mathbf{v}_{n})$ and the local solution defined between the nodes $(\mathbf{u}_{L},\, \mathbf{v}_{L})$, using the basic functional norms $L_{1}$, $L_{2}$, $L_{\infty}$ (\ref{eq:norms_def}). The obtained values of the empirical convergence orders $p$ were compared with the expected theoretical values of convergence order $p_{\rm th.}^{n} = 2N+1$ and $p_{\rm th.}^{l} = N+1$ (\ref{eq:expect_orders}) are applicable for the ADER-DG method for ODE problems. It was found that in the case of numerically solving DAE systems of index 1, the empirical convergence orders $p$ approximately correspond to the expected convergence orders that occur in the case of solving ODE systems. With increasing index of the DAE system (when considering problems for DAE systems of indices 2 and 3), a decrease in the empirical convergence orders $p$ is observed, up to the same convergence orders for the solution at the nodes $p^{n}$ and the local solution $p^{l}$. Using the method of decreasing the index of the DAE system led to an increase in the empirical convergence orders. In general, these results, when the empirical convergence orders decrease with increasing index of the DAE system, were quite expected~\cite{Hairer_book_2}.

An unexpected result was obtained in the numerical solution of the DAE system of index 1 in Subsection~\ref{sec:2:st}, arising from the corresponding problem of the initial value for the stiffness system of ODEs. In this case, cases of extremely high stiffness for this problem were studied. The obtained empirical convergence orders significantly exceeded similar values obtained in the work~\cite{ader_dg_ivp_ode} when solving this problem, and in that work it was a solution of the ODE system, and not a DAE system (in this case, the DAE system can be represented as a DAE system of index 0). In particular, in the case of the first of the studied cases of the stiff DAE system, the empirical convergence orders were very close to the expected theoretical values, despite the stiffness of the problem and the index 1 of the DAE system. This result is of fundamental importance for the choice of basis functions when using the ADER-DG method with a local DG predictor to solve stiff problems --- Lagrange interpolation polynomials with nodal points at the roots of shifted Legendre polynomials as basis functions for representing a local discrete time solution (\ref{eq:qr_def_exp}) were used in the work~\cite{ader_dg_ivp_ode}, and set of Lagrange interpolation polynomials with nodal points at the roots of the right Radau polynomials were used to represent a local discrete time solution (\ref{eq:qr_def_exp}) in the present work. It turns out that the use of Lagrange interpolation polynomials with nodal points at the roots of the right Radau polynomials is much better suited for solving stiff problems.

Estimates of the computational costs showed that the computational costs of the ADER-DG numerical method with a local DG predictor are approximately comparable to the computational costs of implicit Runge-Kutta methods used to solve DAE systems. However, a case was noted in which the computational costs of obtaining a numerical solution using the ADER-DG method with a local DG predictor can exceed the computational costs of the implicit Runge-Kutta method. Methods (\ref{eq:lst_dg_newton_slae_v_opt}) and (\ref{eq:lst_dg_newton_slae_u_opt}) were proposed to reduce the computational costs of the ADER-DG numerical method with a local DG predictor based on the use of the properties of the matrix of the system of linear algebraic equations (\ref{eq:lst_dg_newton_slae}) that arises in the context of applying the classical Newton iterative method to solve the system of nonlinear algebraic equations of the DG predictor (\ref{eq:lst_dg_snae}).

\bmhead{Acknowledgments}
The reported study was supported by the Russian Science Foundation grant No. 21-71-00118 \url{https://rscf.ru/en/project/21-71-00118/}.

The author would like to thank A.P.~Popova for help in correcting the English text.

\subsection*{Data Availability Statement}
\addcontentsline{toc}{subsection}{Data Availability Statement}
Data will be made available on reasonable request.

%The datasets generated during and/or analysed during the current study are available from the corresponding author on reasonable request.

\section*{Declarations}
\addcontentsline{toc}{section}{Declarations}

\subsection*{Conflict of interest}
\addcontentsline{toc}{subsection}{Conflict of interest}
The author declares that he has no conflict of interest.

\addcontentsline{toc}{section}{References}


\end{document}